\documentclass[graybox,envcountchap,sectrefs,leqno]{svmono}
\usepackage{amssymb,amsbsy,amsmath,amsfonts,amssymb,amscd,times,
graphics,color,xy,fancyhdr,fancybox}

\usepackage{mathptmx}
\usepackage{helvet}
\usepackage{courier}
\usepackage{type1cm}
\usepackage{graphicx}     
\usepackage{multicol}        
\usepackage[bottom]{footmisc}

\usepackage[T1]{fontenc}
\newcommand{\C}{\mathbb{C}}\newcommand{\K}{\mathbb{K}}
\newcommand{\N}{\mathbb{N}}
\newcommand{\R}{\mathbb{R}}

\definecolor{black}{cmyk}{1.,1.,1.,1.}
\definecolor{blue}{cmyk}{1.,1.,0.,0.85}
\definecolor{red}{cmyk}{0.,1.,1.,0.85}
\definecolor{green}{cmyk}{1.,0.,1.,0.85}

\newcommand{\blue}{\textcolor{blue}}
\newcommand{\green}{\textcolor{green}}
\newcommand{\red}{\textcolor{red}}

\newcommand{\deutsch}[1]{{\normalsize\sc [#1]}}

\newcommand{\deutschplain}[1]{{\normalsize\sc #1}}

\newcommand{\emphasis}[1]{{\em #1}}

\newcommand{\Fill}{\fbox{\scriptsize\bf Fill~??}\,}

\newcommand{\linestop}{\nopagebreak\medskip\centerline{\bf 
-----------------}\medskip}

\newcommand{\leftbracket}{[}

\newcommand{\name}[1]{{\normalsize\sc #1}}

\newcommand{\pde}[1]{{\normalsize\sc #1}}

\newcommand{\plainstatement}[1]{\smallskip{\em #1}\medskip}

\newcommand{\rightbracket}{]}

\newcommand{\shortplainstatement}[1]{\smallskip{\em #1}}

\newcommand{\smallercharacters}[1]{\smallskip{\normalsize #1}\medskip}

\newcommand{\terminology}[1]{{\sl #1}}

\newcommand{\voir}{{\em see}}

\newcommand{\sectionengellie}[2]{%
\vspace{-0.75cm}
\def\thesection{}
\section{\!\!\!\!\!\!\!\!{\normalsize\sf #1}}
\sectionmark{\centerline{\scriptsize\sf #1}}
\vspace{-0.37cm}}

\makeatletter
\renewcommand{\@fnsymbol}[1]
{\ensuremath{\ifcase#1\or $\dag$ \or $\dag\dag$ \or $\dag\dag\dag$ \or
\else\@ctrerr\fi}}
\makeatother

\def\boiteepaisseavecuntitre#1{%
  \def\thickhrulefill{\leavevmode \leaders \hrule height 1pt \hfill \kern \z@}%
  \def\bkvz@before@breakbox{\ifhmode\par\fi\vskip\breakboxskip\relax}%
  \fboxrule=1pt
  \def\bkvz@set@linewidth{\advance\linewidth -2\fboxrule
                          \advance\linewidth -2\fboxsep}%
  \def\bkvz@left{\vrule \@width\fboxrule\hskip\fboxsep}%
  \def\bkvz@right{\hskip\fboxsep\vrule \@width\fboxrule}%
  \def\bkvz@top{\hbox to \hsize{%
      \vrule\@width\fboxrule\@height 1.2pt 
      \thickhrulefill{#1}\thickhrulefill
      \vrule\@width\fboxrule\@height 1.2pt}}%
  \def\bkvz@bottom{\hrule\@height\fboxrule}%
  \breakbox}

\newcommand{\encadre}{\smallskip\begin{boiteepaisseavecuntitre}{}
\normalsize\baselineskip=0.53cm\fboxrule=0.77pt}
\newcommand{\stopencadre}{\end{boiteepaisseavecuntitre}}

\begin{document}

\large

\author{Joël Merker}
\title{\LARGE Theory of Transformation Groups
\\
{\em by} Sophus Lie and Friedrich Engel
\\
(Vol.~I, 1888)}
\subtitle{\Large Modern Presentation and English Translation}
\maketitle

\frontmatter

\begin{center}

{\LARGE\bf
THEORIE}

\bigskip\medskip

{\Large\sc der}

\bigskip\medskip

{\LARGE\bf
TRANSFORMATIONSGRUPPEN}

\bigskip\medskip

{\bf ----------------}

\bigskip

{\large
ERSTER ABSCHNITT}

\bigskip

{\bf ----------------}

\bigskip\medskip

UNTER MITWIRKUNG VON Prof. Dr. FRIEDRICH ENGEL

\bigskip\medskip

{\Large\sc bearbeitet von}

\bigskip\medskip

{\Large\bf SOPHUS LIE,}

\medskip

{\sc weil. professor der geometrie and der universität leipzig}

{\sc und professor i transformasjonsgruppenes teori an der}

{\sc königlichen frederiks universität zu oslo}

\bigskip\bigskip

{\large\sc unveränderter neudruck}

{\large\sc mit unterstützung der königlichen}

{\large\sc frederiks universität zu oslo}

\bigskip\medskip

\vfill

{\large\bf 1930}

\medskip

{\large\bf VERLAG UND DRUCK VON B.G. TEUBNER IN LEIPZIG UND BERLIN}

\medskip

\end{center}

\newpage

\centerline{\sf Widmung.}

\medskip

\centerline{\bf ----------------}

\bigskip

{\Large
\blue{\sf
Die neuen Theorien, die in diesem Werke dargestellt sind, habe ich in
den Jahren von 1869 bis 1884 entwickelt, wo mir die Liberalität
meines}

\medskip

\centerline{\bf Geburtslandes Norwegen}

\medskip\noindent
\blue{\sf
gestattete, ungestört meine volle Arbeitskraft der Wissen\-schaft zu
widmen, die durch {\sc Abels} Werke in {\sl Norwegens}
wissenschaftlicher Literatur den ersten Platz erhalten hat.}

\blue{\sf Bei der Durchführung meiner Ideen im Einzelnen und bei ihrer
systematischen Darstellung genoss ich seit 1884 in der Grössten
Ausdehnung die unermüdliche Unterstützung des {\sc Professors}}

\medskip

\centerline{\bf Friedrich Engel}

\medskip\noindent
\blue{\sf
meines ausgezeichneten Kollegen and der Universität
{\sl Leipzig}.}

\blue{\sf Das in dieser Weise entstandene Werk widme ich
{\sl Frankreichs}}

\medskip

\centerline{\bf \'Ecole Normale Supérieure}

\medskip\noindent
\blue{\sf
deren unsterblicher Schüler {\sc Galois} zuerst die Bedeutung des
Begriffs {\sl discontinuirliche Gruppe} erkannte. Den hervorragenden
Lehrern dieses Instituts, besonders den Herren G. {\sc Darboux},
E. {\sc Picard} und J. {\sc Tannery} verdanke ich es, dass die
tüchtigsten jungen Mathematiker Frankreichs wetteifernd mit einer
Reihe junger {\sl deutscher} Mathematiker meine Untersuchungen über
{\sl continuirliche Gruppen}, über {\sl Geometrie} und über {\sl
Differentialgleichungen} studiren und mit glänzendem Erfolge
verwerthen.}

}

\bigskip
\hfill
{\bf Sophus Lie.}


\preface

This modernized English translation grew out of my old simultaneous
interest in the mathematics themselves and in the metaphysical
thoughts governing their continued development. I owe to the books of
Robert Hermann, of Peter Olver, of Thomas Hawkins, and of Olle
Stormark to have been introduced in Lie's original vast field.

Up to the end of the 18\textsuperscript{th} Century, the universal
language of Science was Latin, until its centre of gravity shifted to
German during the 19\textsuperscript{th} Century while
nowadays\,---\,and needless to say\,---\,English is widespread. Being
intuitively convinced that Lie's original works contain much more than
what is modernized up to now, I started three years ago to learn
German \emphasis{from scratch} 
just in order to read Lie, and with two main goals in mind:

\smallskip
$\square$\,\,
complete and modernize the Lie-Amaldi classification of
finite-dimensional local Lie group holomorphic actions on spaces of
complex dimensions 1, 2 and 3 for various applications in
complex and Cauchy-Riemann geometry;

\smallskip
$\square$\,\,
better understand the roots of \'Elie Cartan's achievements.

\smallskip
Then it gradually appeared to me that \emphasis{the mathematical
thought of Lie is universal and transhistorical}, hence deserves per
se to be translated. The present adapted English translation follows a
first monograph written in French and specially devoted to the
treatment by Engel and Lie of the so-called Riemann-Helmholtz problem
in the Volume~III of the
\emphasis{Theorie der Transformationsgruppen}, to appear in 2010
apud Hermann, Paris, cf. also the e-print: {\normalsize\sf
arxiv.org/abs/0910.0801/}.

A few observations are in order concerning the chosen format.  For
several reasons, it was essentially impossible to translate directly
the first few chapters in which Lie's intention was to set up the
beginnings of the theory in the highest possible generality,
especially in order to come up with the elimination of the axiom of
inverse, an aspect never dealt with in modern treatises. As a result,
I decided in the first four chapters to reorganize the material and to
reprove the concerned statements, without nevertheless removing
anything from the mathematical contents embraced. But starting with
Chap.~\ref{kapitel-5}, the exposition of Engel and Lie is so smooth,
so rigorous, so understandable, so systematic, so astonishingly well
organized\,---\,in one word: \emphasis{so beautiful for
thought}\,---\,that a pure translation is essential.

Lastly, the author is grateful to Gautam Bharali, to Philip Boalch, to
Egmont Porten for a few fine suggestions concerning the language, but is of
course the sole responsible for the lack of idiomatic English.

\vspace{\baselineskip}
\begin{flushright}\noindent
Paris, \'Ecole Normale Supérieure,\hfill {\it Joël Merker}\\
16 March 2010\hfill {\it $\:$}\\
\end{flushright}


\tableofcontents

\mainmatter


\chapter{Three Principles of Thought 
\\
Governing the Theory of Lie}
\label{three-principles-thought} 
\chaptermark{Three Principles of Thought Governing the Theory of Lie}

\setcounter{footnote}{0}

\abstract*{}

Let $x = (x_1, \dots, x_n )$ be coordinates on an $n$-dimensional real
or complex euclidean space $\C^n$ or $\R^n$, considered as a source
domain. The archetypal objects of Lie's Theory of Continuous
Transformation Groups are \terminology{point transformation equations}:
\[
x_i'
=
f_i(x_1,\dots\,x_n;\,a_1,\dots,a_r)
\ \ \ \ \ \ \ \ \ \ 
{\scriptstyle{(i\,=\,1\,\cdots\,n)}},
\]
parameterized by a finite number $r$ of real or complex parameters
$(a_1, \dots, a_r )$, namely each map $x' = f ( x; \, a) =: f_a ( x)$
is assumed to constitute a diffeomorphism from some 
domain\footnote{\,
By \terminology{domain}, it will always been meant a 
\emphasis{connected}, nonempty open set.
} 
in the source space into some domain in a target space of
the same dimension equipped with coordinates $(x_1', \dots, x_n')$.
Thus, the functional determinant:
\[
\aligned
\det{\rm Jac}(f)
=
\left\vert
\begin{array}{cccc}
\frac{\partial f_1}{\partial x_1}
&\cdots\cdots&
\frac{\partial f_1}{\partial x_n}
\\
\vdots
&\ddots&
\vdots
\\
\frac{\partial f_n}{\partial x_1}
&\cdots\cdots&
\frac{\partial f_n}{\partial x_n}
\end{array}
\right\vert
=
\sum_{\sigma\in{\sf Perm}_n}\,
{\rm sign}(\sigma)\,
\frac{\partial f_1}{\partial x_{\sigma(1)}}\,
\frac{\partial f_2}{\partial x_{\sigma(2)}}
\cdots
\frac{\partial f_n}{\partial x_{\sigma(n)}}
\endaligned
\]
vanishes at no point of the source domain. 

\label{vertauscht}
\begin{svgraybox}
\centerline{{\sf \S\,\,\,15. (\cite{ enlie1888-1}, pp.~25--26)}}

[The concepts of transformation $x' = f ( x)$ and of transformation
equations $x' = f ( x; \, a)$ are of purely analytic nature.]
However, these concepts can receive a graphic interpretation
\deutsch{anschaulich Auffasung}, 
when the concept of an $n$-times extended space
\deutsch{Raum} is introduced.

If we interpret $x_1, \dots, x_n$ as the coordinates for points
\deutsch{Punktcoordinaten} of such a space, then a transformation
$x_i' = f_i (x_1, \dots, x_n)$ appears as a point transformation
\deutsch{Punkttransformation}; consequently, this transformation can be
interpreted as an \emphasis{operation} which consists in that, every
point $x_i$ is transferred at the same time in the new position
$x_i'$. One expresses this as follows: 
the transformation in question is an operation
through which the points of the space $x_1, \dots, x_n$ are
permuted with each other.
\end{svgraybox}

Before introducing the
continuous group axioms in
Chap.~\ref{fundamental-differential}
below, the very first question to
be settled is: how many different transformations $x_i ' = f_i ( x;\,
a)$ do correspond to the $\infty^r$ different systems of values $(
a_1, \dots, a_r)$? Some parameters might indeed be superfluous, hence
they should be removed from the beginning, as will be achieved in
Chap.~\ref{essential-parameters}. 
For this purpose, it is crucial to formulate explicitly and
once for all \emphasis{three principles of thought concerning the
admission of hypotheses that do hold throughout the theory of
continuous groups developed by Lie}.

\subsubsection*{General Assumption of Analyticity}

Curves, surfaces, manifolds, groups, subgroups, coefficients of
infinitesimal transformations, etc., all mathematical objects of
the theory will be assumed to be \terminology{analytic}, i.e. their
representing functions will be assumed to be locally expandable in
convergent, univalent power series defined in a certain domain of an
appropriate $\R^m$.

\subsubsection*{Principle of Free Generic Relocalization}
\label{free-relocalization}

Consider a local mathematical object which is represented by functions
that are analytic in some domain $U_1$, and suppose that a
certain ``generic'' nice behaviour holds on $U_1 \backslash {\sf D}_1$
outside a certain proper closed analytic subset ${\sf D}_1 \subset
U_1$; for instance: the invertibility of a square matrix composed of
analytic functions holds outside the zero-locus of its determinant. Then
relocalize the considerations in some subdomain $U_2 \subset U_1
\backslash {\sf D}_1$.

\begin{center}
\input delocalization.pstex_t
\end{center}

\noindent
Afterwards, in $U_2$, further reasonings may demand to avoid another
proper closed analytic subset ${\sf D}_2$, hence to relocalize the
considerations into some subdomain $U_3 \subset U_2 \backslash {\sf
D}_2$, and so on. Most proofs of the \emphasis{Theorie der
Transformationsgruppen}, and especially the classification theorems,
do allow a great number of times such relocalizations, often without
any mention, such an \emphasis{act of thought} being considered as
implicitly clear, and free relocalization being justified by the
necessity of \emphasis{studying at first generic objects}.

\subsubsection*{Giving no Name to Domains or Neighbourhoods}

Without providing systematic notation, Lie and Engel commonly wrote 
\emphasis{the}
neighbourhood \deutsch{\emphasis{der} Umgebung} (of a point), 
similarly as
one speaks of \emphasis{the} neighbourhood of a house, or of 
\emphasis{the}
surroundings of a town, whereas contemporary topology conceptualizes
\emphasis{a} (say, sufficiently small) given neighbourhood amongst an
infinity. Contrary to what the formalistic, twentieth-century
mythology tells sometimes, Lie and Engel did emphasize the local
nature of the concept of transformation group in terms of narrowing 
down 
neighbourhoods; we shall illustrate this especially when presenting 
Lie's attempt to economize the axiom of inverse. Certainly,
it is true that most of Lie's results are stated without specifying
domains of existence, but in fact also, it is moreover quite plausible
that Lie soon realized that giving no name to neighbourhoods, and
avoiding superfluous denotation is efficient and expeditious in order
to perform far-reaching classification theorems.

\smallskip

Therefore, adopting the economical style of thought in Engel-Lie's
treatise, our ``modernization-translation'' of the theory will,
without nevertheless providing frequent reminders, presuppose
that:

\begin{itemize}

\smallskip\item[$\bullet$]
mathematical objects are analytic; 

\smallskip\item[$\bullet$]
relocalization is freely allowed;

\smallskip\item[$\bullet$]
open sets are often small, usually unnamed, and always
\emphasis{connected}.

\end{itemize}

\begin{svgraybox}
\noindent{\Large\bf Introduction} 
(\cite{ enlie1888-1}, pp.~1--8)

\linestop

If the variables $x_1', \dots, x_n'$ are determined as functions
of $x_1, \dots, x_n$ by $n$ equations solvable with respect
to $x_1, \dots, x_n$:
\[
x_i'
=
f_i(x_1,\dots,x_n)
\ \ \ \ \ \ \ \ \ \ \ \ \ {\scriptstyle{(i\,=\,1\,\cdots\,n)}},
\]
then one says that these equations represent a transformation
\deutsch{Transformation} between the variables $x$ and $x'$. In the
sequel, we will have to deal with such transformations; unless the
contrary is expressly mentioned, we will restrict ourselves to the
case where the $f_i$ are \emphasis{analytic \deutsch{analytisch}
functions} of their arguments. But because a not negligible portion
of our results is independent of this assumption, we will occasionally
indicate how various developments take shape by taking into
considerations functions of this sort.

\smallercharacters{

When the functions $f_i ( x_1, \dots, x_n)$ are analytic and are
defined inside a common region \deutsch{Bereich}, then according to
the known studies of \name{Cauchy}, \name{Weierstrass}, \name{Briot}
and \name{Bouquet}, one can always delimit, in the manifold of all
real and complex systems of values $x_1, \dots, x_n$, a region $(x)$
such that all functions $f_i$ are univalent in the complete extension
\deutsch{Ausdehnung} of this region, and so that as well, in the
neighbourhood \deutsch{Umgebung} of every system of values $x_1^0,
\dots, x_n^0$ belonging to the region $(x)$, the functions behave
regularly \deutsch{regulär verhalten}, that is to say, they can be
expanded in ordinary power series with respect to $x_1 - x_1^0$,
\dots, $x_n - x_n^0$ with only entire positive powers.

For the solvability of the equations $x_i' = f_i (x)$, a unique
condition is necessary and sufficient, namely the condition that the
functional determinant:
\[
\sum\,\pm\,
\frac{\partial f_1}{\partial x_1}\,\cdots\,
\frac{\partial f_n}{\partial x_n}
\]
should not vanish identically. If this condition is satisfied, then
the region $(x)$ defined above can specially be defined so that the
functional determinant does not take the value zero for any system of
values in the $(x)$. Under this assumption, if one lets the $x$ take
gradually all systems of values in the region $(x)$, then the
equations $x_i' = f_i (x)$ determine, in the domain \deutsch{Gebiete}
of the $x'$, a region of such a nature that $x_1, \dots, x_n$, in the
neighbourhood of every system of values ${x_1'}^0, \dots, {x_n'}^0$ in
this new region behave regularly as functions of $x_1', \dots, x_n'$,
and hence can be expanded as ordinary power series of $x_1' - {x_1'}^0$,
\dots, $x_n' - {x_n'}^0$. It is well known that from this,
it does not follow that the $x_i$ are univalent functions 
of $x_1', \dots, x_n'$ in the complete extension of
the new region; but when necessary, it is possible to narrow down 
the region $(x)$ defined above so that two different
systems of values $x_1, \dots, x_n$ of the region $(x)$ always
produce two, also different systems of values $x_1' = f_1 (x)$, 
\dots, $x_n' = f_n(x)$.

Thus, the equations $x_i' = f_i (x)$ establish a
univalent invertible relationship \deutsch{Beziehung} 
between regions in the domain of the $x$ and regions
in the domain of the $x'$; to every system of values
in one region, they associate one and only 
one system of values in the other region, and conversely.

}

If the equations $x_i' = f_i(x)$ are solved with respect to 
the $x$, then in turn, the resulting equations:
\[
x_k
=
F_k(x_1',\dots,x_n')
\ \ \ \ \ \ \ \ \ \ \ \ \ {\scriptstyle{(k\,=\,1\,\cdots\,n)}}
\]
again represent a transformation. The relationship between 
this transformation and the initial one is evidently 
a reciprocal relationship; accordingly, one says: 
the two transformations are \emphasis{inverse}
one to another. 
From this definition, it visibly follows:

\shortplainstatement{If one executes at first the transformation:
\[
x_i'
=
f_i(x_1,\dots,x_n)
\ \ \ \ \ \ \ \ \ \ \ \ \ {\scriptstyle{(i\,=\,1\,\cdots\,n)}}
\]
and afterwards the transformation inverse to it:
\[
x_i''
=
F_i(x_1',\dots,x_n')
\ \ \ \ \ \ \ \ \ \ \ \ \ {\scriptstyle{(i\,=\,1\,\cdots\,n)}},
\]
then one obtains the identity transformation:
\[
x_i''
=
x_i
\ \ \ \ \ \ \ \ \ \ \ \ \ {\scriptstyle{(i\,=\,1\,\cdots\,n)}}.
\]}

Here lies the real definition of the concept
\deutsch{Begriff} of two transformations inverse inverse
to each other.

In general, if one executes two arbitrary transformations:
\[
x_i'
=
f_i(x_1,\dots,x_n),\ \ \ \ \ \
x_i''
=
g_i(x_1',\dots,x_n')
\ \ \ \ \ \ \ \ \ \ \ \ \ {\scriptstyle{(i\,=\,1\,\cdots\,n)}}
\]
one after the other, then one obtains a new transformation, 
namely the following one:
\[
x_i''
=
g_i\big(f_1(x),\dots,f_n(x)\big)
\ \ \ \ \ \ \ \ \ \ \ \ \ {\scriptstyle{(i\,=\,1\,\cdots\,n)}}.
\]
In general, this new transformation naturally changes
when one changes the order \deutsch{Reihenfolge} of the
two transformations; however, it can also happen
that the order of the two transformations is indifferent. 
This case occurs when one has identically:
\[
g_i\big(f_1(x),\dots,f_n(x)\big)
\equiv
f_i\big(g_1(x),\dots,g_n(x)\big)
\ \ \ \ \ \ \ \ \ \ \ \ \ {\scriptstyle{(i\,=\,1\,\cdots\,n)}}\,;
\]
we then say, as in the process of the Theory of
Substitutions
\deutsch{Vorgang der Substitutionentheorie}: 
\emphasis{the two transformations:
\[
x_i'
=
f_i(x_1,\dots,x_n)
\ \ \ \ \ \ \ \ \ \ \ \ \ {\scriptstyle{(i\,=\,1\,\cdots\,n)}}
\]
and:
\[
x_i'
=
g_i(x_1,\dots,x_n)
\ \ \ \ \ \ \ \ \ \ \ \ \ {\scriptstyle{(i\,=\,1\,\cdots\,n)}}
\]
are interchangeable \deutsch{vertauschbar} one with the
other}.\,---

\renewcommand{\thefootnote}{\fnsymbol{footnote}}
\plainstatement{A finite or infinite family 
\deutsch{Schaar}
of transformations between
the $x$ and the $x'$ is called a \terminology{group of
transformations} or a \terminology{transformation group} when any two
transformations of the family executed one after the other give a
transformation which again belongs to the family}.\footnotemark[1]
\renewcommand{\thefootnote}{\arabic{footnote}}
\end{svgraybox}

\renewcommand{\thefootnote}{\fnsymbol{footnote}}
\footnotetext[1]{\,
Sophus \name{Lie}, Gesellschaft der Wissenschaften zu 
Christiania 1871, p.~243. \name{Klein}, Vergleichende Betrachtungen
über neuere geometrische Forschungen, Erlangen 1872. \name{Lie}, 
Göttinger Nachrichten 1873, 3. Decemb.
} 
\renewcommand{\thefootnote}{\arabic{footnote}}

\vspace{-1.2cm}

\begin{svgraybox}
A transformation group is called \terminology{discontinuous} when it
consists of a discrete number of transformations, and this number can
be finite or infinite. Two transformations of such a group are
finitely different from each other. The discontinuous groups belong
to the domain of the \emphasis{Theory of Substitutions}, so in the
sequel, they will remain out of consideration.

The discontinuous groups stand in opposition to the
\emphasis{continuous} transformation groups, which always contain
infinitely many transformations. A transformation group is called
\terminology{continuous} when it is possible, for every transformation
belonging to the group, to indicate certain other transformations
which are only infinitely little different from the transformation in
question, and when by contrast, it is not possible to reduce the
complete totality \deutsch{Inbegriff} of transformations contained in
the group to a single discrete family.

Now, amongst the continuous transformation groups, we consider again
two separate categories \deutsch{Kategorien} which, in the
nomenclature \deutsch{Benennung}, are distinguished as
\terminology{finite continuous} groups and as \terminology{infinite
continuous} groups. To begin with, we can only give provisional
definitions of the two categories, and these definitions will be
apprehended precisely later.

A \emphasis{finite continuous transformation group} will
be represented by \emphasis{one} system of $n$ equations:
\[
x_i'
=
f_i(x_1,\dots,x_n,\,a_1,\dots,a_r)
\ \ \ \ \ \ \ \ \ \ \ \ \ {\scriptstyle{(i\,=\,1\,\cdots\,n)}},
\]
where the $f_i$ denote analytic functions of the variables
$x_1, \dots, x_n$ and of the arbitrary parameters $a_1, \dots, a_r$. 
Since we have to deal with a group, two transformations:
\[
\aligned
x_i'
&
=
f_i(x_1,\dots,x_n,\,a_1,\dots,a_r)
\\
x_i''
&
=
f_i(x_1',\dots,x_n',\,a_1,\dots,a_r),
\endaligned
\]
when executed one after the other, must produce a transformation
which belongs to the group, hence which has the form:
\[
x_i''
=
f_i\big(
f_1(x,a),\dots,f_n(x,a),\,b_1,\dots,b_r\big)
=
f_i(x_1,\dots,x_n,\,c_1,\dots,c_r).
\]
Here, the $c_k$ are naturally independent of the $x$ and
so, are functions of only the $a$ and the $b$.

{\sf Example.}
A known group of this sort is the following:
\[
x'
=
\frac{x+a_1}{a_2\,x+a_3},
\]
which contains the three parameters $a_1$, $a_2$, $a_3$. 
If one executes the two transformations:
\[
x'
=
\frac{x+a_1}{a_2\,x+a_3},\ \ \ \ \ \
x''
=
\frac{x'+b_1}{b_2\,x'+b_3}
\]
one after the other, then one receives:
\[
x''
=
\frac{x+c_1}{c_2\,x+c_3},
\]
where $c_1$, $c_2$, $c_3$ are defined as functions
of the $a$ and the $b$ by the relations:
\[
c_1
=
\frac{a_1+b_1\,a_3}{1+b_1\,a_2},\ \ \ \ \ \ \
c_2
=
\frac{b_2+a_2\,b_3}{1+b_1\,a_2},\ \ \ \ \ \ \
c_3
=
\frac{b_2\,a_1+b_3\,a_3}{1+b_1\,a_2}.
\]

The following group with the $n^2$ parameters $a_{ ik}$
is not less known:
\[
x_i'
=
\sum_{k=1}^n\,a_{ik}\,x_k
\ \ \ \ \ \ \ \ \ \ \ \ \ {\scriptstyle{(i\,=\,1\,\cdots\,n)}}.
\]
If one sets here:
\[
x_\nu''
=
\sum_{i=1}^n\,b_{\nu i}\,x_i'
\ \ \ \ \ \ \ \ \ \ \ \ \ {\scriptstyle{(i\,=\,1\,\cdots\,n)}},
\]
then it comes:
\[
x_\nu''
=
\sum_{i,\,\,k}^{1\cdots n}\,
b_{\nu i}\,a_{ik}\,x_k
=
\sum_{k=1}^n\,c_{\nu k}\,x_k,
\]
where the $c_{\nu k}$ are determined by the equations:
\[
c_{\nu k}
=
\sum_{i=1}^n\,b_{\nu i}\,a_{ik}
\ \ \ \ \ \ \ \ \ \ \ \ \ 
{\scriptstyle{(\nu,\,\,k\,=\,1\,\cdots\,n)}}.\,-
\]

In order to arrive at a usable definition of a finite continuous
group, we want at first to somehow reshape the definition of finite
continuous groups. On the occasion, we use a proposition from the
theory of differential equations about which, besides, we will come
back later in a more comprehensive way (cf. Chap.~\ref{kapitel-10}).

Let the equations:
\[
x_i'
=
f_i(x_1,\dots,x_n,\,a_1,\dots,a_r)
\ \ \ \ \ \ \ \ \ \ \ \ \ {\scriptstyle{(i\,=\,1\,\cdots\,n)}}
\]
represent an arbitrary continuous group. According to the proposition
in question, it is then possible to define the functions $f_i$ through
a system of differential equations, insofar as they depend upon the
$x$. To this end, one only has to differentiate the equations $x_i' =
f_i ( x, a)$ with respect to $x_1, \dots, x_n$ sufficiently often and
then to set up all equations that may be obtained by elimination of
$a_1, \dots, a_r$. If one is gone sufficiently far by
differentiation, then by elimination of the $a$, one obtains a system
of differential equations for $x_1', \dots, x_n'$, whose most general
system of solutions is represented by the initial equations $x_i' =
f_i ( x, a)$ with the $r$ arbitrary parameters. Now, since by
assumption the equations $x_i' = f_i ( x, a)$ define a group, it
follows that the concerned system of differential equations possesses
the following remarkable property: if $x_i' = f_i ( x_1, \dots, x_n,
\, b_1, \dots, b_r)$ is a system of solutions of it, and if $x_i' =
f_i ( x, a)$ is a second system of solutions, then:
\[
x_i'
=
f_i
\big(
f_1(x,a),\dots,f_n(x,a),\,b_1,\dots,b_r
\big)
\ \ \ \ \ \ \ \ \ \ \ \ \ {\scriptstyle{(i\,=\,1\,\cdots\,n)}}
\]
is also a system of solutions.

From this, we see that the equations of an arbitrary finite continuous
transformation group can be defined by a system of differential
equations which possesses certain specific properties. Firstly, from
two systems of solutions of the concerned differential equations one
can always derive, in the way indicated above, a third system of
solutions: it is precisely in this that we have to deal with a group.
Secondly, the most general system of solutions of the concerned
differential equations depends only upon a finite number of arbitrary
constants: this circumstance expresses that our group is finite.

Now, we assume that there is a family of transformations $x_i' = f_i (
x_1, \dots, x_n)$ which is is defined by a system of differential
equations of the form:
\[
W_k
\bigg(
x_1',\dots,x_n',\,
\frac{\partial x_1'}{\partial x_1},\,\dots,\,
\frac{\partial^2x_1'}{\partial x_1^2},\,\dots
\bigg)
=
0
\ \ \ \ \ \ \ \ \ \ \ \ \ {\scriptstyle{(k\,=\,1,\,\,2\,\cdots)}}.
\]
Moreover, we assume that this system of differential equations
possesses the first one of the two mentioned properties, but not the
second one; therefore, with $x_i' = f_i (x_1, \dots, x_n)$ and $x_i' =
g_i ( x_1, \dots, x_n)$, then always, $x_i' = g_i \big( f_1 ( x),
\dots, f_n (x) \big)$ is also a system of solutions of these
differential equations, and the most general system of solutions of
them does not only depend upon a finite number of arbitrary constants,
but also upon higher sorts of elements, as for example, upon arbitrary
functions. Then the totality of all transformations which satisfy
the concerned differential equations evidently forms again a group, 
and in general, a continuous group, though no more a finite one, 
but one which we call \terminology{infinite continuous}. 

Straightaway, we give a few simple examples of infinite
continuous transformation groups. 

When the differential equations which define the concerned
infinite group reduce to the identity equation 
$0 = 0$, then the transformations of the group read:
\[
x_i'
=
\Pi_i(x_1,\dots,x_n)
\ \ \ \ \ \ \ \ \ \ \ \ \ {\scriptstyle{(i\,=\,1\,\cdots\,n)}}
\]
where the $\Pi_i$ denote arbitrary analytic functions
of their arguments.

Also the equations: 
\[
\frac{\partial x_i'}{\partial x_k}
=
0
\ \ \ \ \ \ \ \ \ \ \ \ \ 
{\scriptstyle{(i\,\neq\,k\,;\,\,\,
i,\,\,k\,=\,1\,\cdots\,n)}}
\]
define an infinite group, namely the following one:
\[
x_i'
=
\Pi_i(x_i)
\ \ \ \ \ \ \ \ \ \ \ \ \ {\scriptstyle{(i\,=\,1\,\cdots\,n)}}, 
\]
where again the $\Pi_i$ are absolutely arbitrary.

\smallercharacters{

Besides, it is to be observed that the concept of an infinite
continuous group can yet be understood more generally than what takes
place here. Actually, one could call infinite continuous any
continuous group which is not finite. However, this definition does
not coincide with the one given above.

For instance, the equations:
\[
x'
=
F(x),\ \ \ \ \ \ \
y'
=
F(y)
\]
in which $F$, for the two cases, denotes the same function of its
arguments, represent a group. This group is continuous, since all its
transformations are represented by a single system of
equations; in addition, it is obviously not finite. 
Consequently, it would be an infinite continuous group 
if we would interpret this concept in the more general sense
indicated above. But on the other hand, it is
not possible to define the family of the transformations:
\[
x'
=
F(x),\ \ \ \ \ \ \
y'
=
F(y)
\]
by differential equations that are free of arbitrary 
elements. Consequently, the definition stated first for
an infinite continuous group does not fit to this case. 
Nevertheless, we find suitable to consider only the
infinite continuous groups which can be defined by differential 
equations and hence, we always set as fundamental our first, 
tight definition.

}

We do not want to omit emphasizing that the concept of
``transformation group'' is still not at all exhausted by the difference
between discontinuous and continuous groups. Rather, there are
transformation groups which are subordinate to none of these two
classes but which have something in common
\label{mention-discontinuous} with each one of the two
classes. In the sequel, we must at least occasionally also treat this
sort of groups. Provisionally, two examples will suffice.

The totality of all coordinate transformations of a plane by which one
transfers an ordinary right-angled system of coordinates to another
one forms a group which is neither continuous, nor discontinuous.
Indeed, the group in question contains two separate categories of
transformations between which a continuous transition is not possible:
firstly, the transformations by which the old and the new systems of
coordinates are congruent, and secondly, the transformations by which
these two systems are not congruent.

The first ones have the form:
\[
x'-a
=
x\,\cos\alpha-y\,\sin\alpha,\ \ \ \ \
y'-b
=
x\,\sin\alpha+y\,\cos\alpha,
\]
while the analytic expression of the second ones reads:
\[
x'-a
=
x\,\cos\alpha+y\,\sin\alpha,\ \ \ \ \
y'-b=x\,\sin\alpha-y\,\cos\alpha.
\]
Each one of these systems of equations represents a continuous family
of transformations, hence the group is not discontinuous; but it is
also not continuous, because both systems of equations taken together
provide all transformations of the group; thus, the transformations of
the group decompose in two discrete families.
If one imagines the plane $x, y$ in ordinary space and if
one adds $z$ as a third right-angled coordinate, then 
one can imagine the totality of coordinate transformations
of the plane $z = 0$ is a totality of certain movements 
\deutsch{Bewegungen} of the space, namely the movements
during which the plane $z = 0$ keeps its position. 
Correspondingly, these movements separate in two classes,
namely in the class which only shifts the plane
in itself and in the class which 
turns the plane.

As a second example of such a group, one can yet
mention the totality of all projective and dualistic
transformations of the plane.

\linestop 

According to these general remarks about the concept of transformation
group, we actually turn ourselves to the consideration of the finite
continuous transformation groups which constitute the object of the
studies following. These studies are divided in three volumes
\deutsch{Abschnitte}. The \emphasis{first} volume treats finite
continuous groups in general. The \emphasis{second} volume treats the
finite continuous groups whose transformations are so-called
\terminology{contact transformations}
\deutsch{Beührungstransformationen}. Lastly, in the \emphasis{third}
volume, certain general problems of group theory will be carried out
in great details for a small number of variables.
\end{svgraybox}

\linestop


\chapter{Local Transformation Equations 
\\
and Essential Parameters}
\label{essential-parameters}
\chaptermark{Local Transformation Equations and Essential Parameters}

\setcounter{footnote}{0}

\abstract{
Let $\K = \R$ or $\C$, throughout.
As said in Chap.~\ref{three-principles-thought}, transformation
equations $x_i' = f_i ( x; \, a_1, \dots, a_r)$, $i =1, \dots, n$,
which are local, analytic diffeomorphisms of $\K^n$ parametrized by a
finite number $r$ of real or complex numbers $a_1, \dots, a_r$,
constitute the archetypal objects of Lie's theory. The preliminary
question is to decide whether the $f_i$ really depend upon
\emphasis{all} parameters, and also, to get rid of superfluous
parameters, if there are any.
\newline\indent Locally in a neighborhood of a fixed $x_0$, one
expands $f_i ( x; \, a) = \sum_{ \alpha \in \N^n}\, \mathcal{
U}_\alpha^i (a) \, ( x - x_0)^\alpha$ in power series and one looks at
the \terminology{infinite coefficient mapping}\, ${\sf U}_\infty: \, a
\longmapsto \big( \mathcal{ U}_\alpha^i ( a) \big)_{ \alpha \in
\N^n}^{ 1\leqslant i \leqslant n}$ from $\K^r$ to $\K^\infty$,
expected to tell faithfully the dependence with respect to $a$ in
question. If $\rho_\infty$ denotes the maximal, generic and locally
constant rank of this map, with of course $0 \leqslant \rho_\infty
\leqslant r$, then the answer says that locally in a neighborhood of a
generic $a_0$, there exist both a local change of parameters $a
\mapsto \big( {\sf u}_1 ( a), \dots, {\sf u}_{\rho_\infty} ( a) \big)
=: {\sf u}$ decreasing the number of parameters from $r$ down to
$\rho_\infty$, and new transformation equations:
\[
x_i'
=
g_i
\big(x;\,{\sf u}_1,\dots,{\sf u}_{\rho_\infty}\big)
\ \ \ \ \ \ \ \ \ \
{\scriptstyle{(i\,=\,1\,\cdots\,n)}}
\]
depending \emphasis{only} upon $\rho_\infty$ parameters 
which give again the old ones:
\[
g_i\big(x;\,{\sf u}(a)\big)
\equiv 
f_i(x;\, a)
\ \ \ \ \ \ \ \ \ \
{\scriptstyle{(i\,=\,1\,\cdots\,n)}}.
\]
At the end of this brief chapter, before introducing precisely the
local Lie group axioms, we present an example due to Engel which shows
that the axiom of inverse cannot be deduced from the axiom of
composition, contrary to one of Lie's \emphasis{Idées fixes}.
}

\section{Generic Rank of the Infinite Coefficient Mapping}
\label{generic-rank}

Thus, we consider local transformation equations:
\[
x_i'
=
f_i(x_1,\dots,x_n;a_1,\dots,a_r)
\ \ \ \ \ \ \ \ \ \
{\scriptstyle{(i\,=\,1\,\cdots\,n)}}.
\]
We want to illustrate how the principle of free generic relocalization
exposed just above on p.~\pageref{free-relocalization} helps to get
rid of superfluous parameters $a_k$. We assume that the $f_i$ are
defined and analytic for $x$ belonging to a certain (unnamed,
connected) domain of $\K^n$ and for $a$ belonging as well to some
domain of $\K^r$.

Expanding the $f_i$ of $x_i ' = f_i ( x; a)$ in power series with
respect to $x - x_0$ in some neighborhood of a point $x_0$:
\[
f_i(x;a)
=
\sum_{\alpha\in\N^n}\,\mathcal{U}_\alpha^i(a)\,
(x-x_0)^\alpha,
\]
we get an infinite number of analytic functions $\mathcal{ U}_\alpha^i
= \mathcal{ U}_\alpha^i (a)$ of the parameters that are defined in
some uniform domain of $\K^r$. Intuitively, this infinite collection
of coefficient functions $\mathcal{ U}_\alpha^i (a)$ should show how
$f ( x; a )$ does depend on $a$.

To make this claim precise, we thus consider the map:
\[
{\sf U}_\infty:
\ \ \ \ \
\K^r\ni\ a
\longmapsto
\big(
\mathcal{U}_\alpha^i(a)
\big)_{\alpha\in\N^n}^{1\leqslant i\leqslant n}\
\in\K^\infty.
\]
For the convenience of applying standard differential 
calculus in finite dimension, we simultaneously consider also 
all its $\kappa$-th truncations:
\[
{\sf U}_\kappa:
\ \ \ \ \ 
\K^r\ni\ a
\longmapsto
\big(
\mathcal{U}_\alpha^i(a)
\big)_{\vert\alpha\vert\leqslant\kappa}^{1\leqslant i\leqslant n}\
\in\K^{n\frac{(n+\kappa)!}{n!\,\,\kappa!}},
\] 
where $\frac{ ( n + \kappa ) ! }{ n! \, \, \kappa !}$ is the number of
multiindices $\alpha \in \N^n$ whose length $\vert \alpha \vert :=
\alpha_1 + \cdots + \alpha_n$ satisfies the upper bound $\vert \alpha
\vert \leqslant \kappa$. We call ${\sf U}_\kappa$, ${\sf U}_\infty$
the \terminology{(in)finite coefficient mapping(s)} of $x_i' = f_i (
x; \, a)$.

The \terminology{Jacobian matrix}\, of ${\sf U}_\kappa$ is the $r
\times \big( n\, \frac{(n+ \kappa)!}{n!\, \, \kappa!} \big)$ matrix:
\[
\Big(
{\textstyle{\frac{\partial\mathcal{U}_\alpha^i}{
\partial a_j}}}
(a)
\Big)_{1\leqslant j\leqslant r}^{
\vert\alpha\vert\leqslant\kappa,\, 
1\leqslant i\leqslant n}, 
\]
its $r$ rows being indexed by the partial derivatives. The
\terminology{generic rank}\, of ${\sf U}_\kappa$ is the largest
integer $\rho_\kappa \leqslant r$ such that there is a $\rho_\kappa
\times \rho_\kappa$ minor of ${\rm Jac}\, {\sf U}_\kappa$ which does
not vanish identically, but all $(\rho_\kappa + 1) \times (
\rho_\kappa + 1)$ minors do vanish identically. The uniqueness
principle for analytic functions then insures that the common zero-set
of all $\rho_\kappa \times \rho_\kappa$ minors is a
\emphasis{proper}\, closed analytic subset ${\sf D}_\kappa$ (of the
unnamed domain where the $\mathcal{ U}_\alpha^i$ are defined), so it
is stratified by a finite number of submanifolds of codimension
$\geqslant 1$ (\cite{ ma1967, bcr1987, ch1989, gu1990}), and in
particular, it is of empty interior, hence intuitively ``thin''.

So the set of parameters $a$ at which there is a least one
$\rho_\kappa \times \rho_\kappa$ minor of ${\sf Jac}\, {\sf U}_\kappa$
which does not vanish is open and \emphasis{dense}. Consequently,
``for a generic point $a$'', the map ${\sf U}_\kappa$ is of rank
$\geqslant \rho_\kappa$ at every point $a'$ sufficiently close to $a$
(since the corresponding $\rho_\kappa \times \rho_\kappa$ minor does
not vanish in a neighborhood of $a$), and because all $( \rho_\kappa +
1) \times ( \rho_\kappa + 1)$ minors of ${\rm Jac}\, {\sf U}_\kappa$
were assumed to vanish identically, the map ${\sf U}_\kappa$ happens
to be in fact of \emphasis{constant}\, rank ${\sf U}_\kappa$ in a
(small) neighborhood of every such a generic $a$.

\plainstatement{Insuring constancy of a
rank is one important instance of why free relocalization is useful: a
majority of theorems of the differential calculus and of the classical
theory of (partial) differential equations do hold under specific
local constancy assumptions}.

As $\kappa$ increases, the number of columns of ${\rm Jac}\, U_\kappa$
increases, hence $\rho_{ \kappa_1} \leqslant \rho_{\kappa_2}$ for
$\kappa_1 \leqslant \kappa_2$. Since $\rho_\kappa \leqslant r$ is
anyway bounded, the generic rank of ${\sf U}_\kappa$ becomes constant
for all $\kappa \geqslant \kappa_0$ bigger than some large enough
$\kappa_0$. Thus, let $\rho_\infty \leqslant r$ denote this maximal
possible generic rank.

\begin{definition}
The parameters $(a_1, \dots, a_r)$ of given point transformation
equations $x_i ' = f_i ( x; a)$ are called \terminology{essential}\,
if, after expanding $f_i ( x; a) = \sum_{ \alpha \in \N^n}\,
\mathcal{ U}_\alpha^i ( a) \, (x - x_0)^\alpha$ in power series at
some $x_0$, the generic rank $\rho_\infty$ of the coefficient mapping
$a \longmapsto \big( \mathcal{ U}_\alpha^i ( a) \big)_{ \alpha \in
\N^n}^{ 1 \leqslant i \leqslant n}$ is maximal, equal to the number
$r$ of parameters: $\rho_\infty = r$. 
\end{definition}

Without entering technical
details, we make a remark.
It is a consequence of the principle of analytic continuation and of
some reasonings with power series that the \emphasis{same}\, maximal
rank $\rho_\infty$ is enjoyed by the coefficient mapping $a \mapsto
\big( {\mathcal{ U}'}_\alpha^i (a) \big)_{ \alpha\in \N^n}^{
1\leqslant i \leqslant n}$ for the expansion of $f_i ( x; a) = \sum_{
\alpha \in \N^n}\, {\mathcal{ U }'}_\alpha^i (a) \, ( x - x_0 '
)^\alpha$ at another, arbitrary point $x_0'$. Also, one can prove that
$\rho_\infty$ is independent of the choice of coordinates $x_i$ and of
parameters $a_k$. These two facts will not be needed, and the
interested reader is referred to~\cite{ me2005} for proofs of quite
similar statements holding true in the context of
\emphasis{Cauchy-Riemann geometry}.

\section{Quantitative Criterion
\\
for the Number of Superfluous Parameters}
\sectionmark{Quantitative Criterion
for the Number of Superfluous Parameters}
\label{quantitative-criterion}

It is rather not practical to compute the generic rank of the infinite
Jacobian matrix ${\rm Jac}\, {\sf U }_\infty$. To check essentiality
of parameters in concrete situations, a helpful criterion due to Lie
is {\bf (iii)} below.

\label{Theorem-essential}
\begin{theorem}
The following three conditions are equivalent:

\begin{itemize}

\smallskip\item[{\bf \!(i)}]
In the transformation equations
\[
x_i'
=
f_i(x_1,\dots,x_n;\,a_1,\dots,a_r)
=
\sum_{\alpha\in\N^n}\,\mathcal{U}_\alpha^i(a)\,
(x-x_0)^\alpha
\ \ \ \ \ \ \ \ \ \
{\scriptstyle{(i\,=\,1\,\cdots\,n)}},
\]
the parameters $a_1, \dots, a_r$ are \emphasis{not}
essential.

\smallskip\item[{\bf \!\!\!(ii)}]
(By definition) The generic rank $\rho_\infty$ of the infinite
Jacobian matrix:
\[
{\rm Jac}\,{\sf U}_\infty(a)
=
\Big(
\frac{\partial\mathcal{U}_\alpha^i}{\partial a_j}(a)
\Big)_{1\leqslant j\leqslant r}^{\alpha\in\N^n,\,1\leqslant i\leqslant n}
\]
is strictly less than $r$.

\smallskip\item[{\bf \!\!\!\!(iii)}]
\label{iii-ess-param}
Locally in a neighborhood of {\rm every} $(x_0, a_0)$, 
there exists a not identically zero analytic vector field
on the parameter space:
\[
\mathcal{T}
=
\sum_{k=1}^n\,\tau_k(a)\,\frac{\partial}{\partial a_k}
\]
which annihilates all the $f_i ( x; a)$:
\[
0
\equiv
\mathcal{T}\,f_i
=
\sum_{k=1}^n\,\tau_k\,\frac{\partial f_i}{\partial a_k}
=
\sum_{\alpha\in\N^n}\,\sum_{k=1}^r\,\tau_k(a)\,
\frac{\partial\mathcal{U}_\alpha^i}{\partial a_k}(a)\,
(x-x_0)^\alpha
\ \ \ \ \ \ \ \ \ \ \ \ \
{\scriptstyle{(i\,=\,1\,\cdots\,n)}}.
\]
\end{itemize}\smallskip

\noindent
More generally, if $\rho_\infty$ denotes the generic rank of the
infinite coefficient mapping:
\[
{\sf U}_\infty:\ \ \ \ \ \ \
a\ \ 
\longmapsto\ \
\big(
\mathcal{U}_\alpha^i(a)
\big)_{\alpha\in\N^n}^{1\leqslant i\leqslant n},
\]
then locally in a neighborhood of {\rm every} $(x_0, a_0)$, there exist
exactly $r - \rho_\infty$, and no more, analytic vector fields:
\[
\mathcal{T}_\mu
=
\sum_{k=1}^n\,\tau_{\mu k}(a)\,
\frac{\partial}{\partial a_k}
\ \ \ \ \ \ \ \ \ \
{\scriptstyle{(\mu\,=\,1\,\cdots\,r-\rho_\infty)}},
\]
with the property that the dimension of ${\rm Span}\big( \mathcal{
T}_1 \big\vert_a$, \dots, $\mathcal{ T}_{ r - \rho_\infty} \big\vert_a
\big)$ is equal to $r - \rho_\infty$ at every parameter $a$ at which
the rank of ${\sf U}_\infty$ is maximal equal to $\rho_\infty$, such
that the derivations $\mathcal{ T}_\mu$ all annihilate the $f_i ( x;
a)$:
\[
0
\equiv
\mathcal{T}_\mu\,f_i
=
\sum_{k=1}^r\,\tau_{\mu k}(a)\,
\frac{\partial f_i}{\partial a_k}(x;a)
\ \ \ \ \ \ \ \ \ \
{\scriptstyle{(i\,=\,1\,\cdots\,n;\,\,\,
\mu\,=\,1\,\cdots\,r-\rho_\infty)}}.
\]
\end{theorem}

\begin{proof}\smartqed
Just by the chosen definition, we have {\bf (i)} $\Longleftrightarrow$
{\bf (ii)}. Next, suppose that condition {\bf (iii)} holds, in which
the coefficients $\tau_k ( a)$ of the concerned nonzero derivation
$\mathcal{ T}$ are locally defined. Reminding that the Jacobian matrix
${\rm Jac} \, {\sf U}_\infty$ has $r$ rows and an infinite number of
columns, we then see that the $n$ annihilation equations $0 \equiv
\mathcal{ T}\, f_i$, when rewritten in matrix form as:
\[
0
\equiv
\big(\tau_1(a),\dots,\tau_r(a)\big)\,
\Big(
\frac{\partial\mathcal{U}_\alpha^i}{\partial a_j}(a)
\Big)_{1\leqslant j\leqslant r}^{\alpha\in\N^n,\,1\leqslant i\leqslant n}
\]
just say that the transpose of ${\rm Jac }\, {\sf U }_\infty (a)$ has
nonzero kernel at each $a$ where the vector $\mathcal{ T} \big \vert_a
= \big( \tau_1 ( a), \dots, \tau_r ( a) \big)$ is nonzero.
Consequently, ${\rm Jac }\, {\sf U }_\infty$ has rank strictly less
than $r$ locally in a neighborhood of every $a_0$, hence in the whole
$a$-domain. So {\bf (iii)} $\Rightarrow$ {\bf (ii)}.

Conversely, assume that the generic rank $\rho_\infty$ of ${\rm Jac}\,
{\sf U}_\infty$ is $< r$. Then there exist $\rho_\infty <r$ ``basic''
coefficient functions $\mathcal{ U}_{ \alpha ( 1)}^{ i ( 1)}, \dots,
\mathcal{ U}_{ \alpha ( \rho_\infty)}^{ i ( \rho_\infty)}$ (there can
be several choices) such that the generic rank of the extracted map $a
\mapsto \big( \mathcal{ U}_{\alpha ( l)}^{ i ( l)} \big)_{ 1\leqslant
l \leqslant \rho_\infty}$ equals $\rho_\infty$ already. We abbreviate:
\[
{\sf u}_l(a)
:=
\mathcal{U}_{\alpha(l)}^{i(l)}(a)
\ \ \ \ \ \ \ \ \ \
{\scriptstyle{(l\,=\,1\,\cdots\,\rho_\infty)}}.
\]
The goal is to find vectorial local analytic solutions $\big( \tau_1 ( a),
\dots, \tau_r ( a) \big)$ to the infinite number of linear equations:
\[
0
\equiv
\tau_1(a)\,
\frac{\partial\mathcal{U}_\alpha^i(a)}{\partial a_1}(a)
+\cdots+
\tau_r(a)\,
\frac{\partial\mathcal{U}_\alpha^i(a)}{\partial a_r}(a)
\ \ \ \ \ \ \ \ \ \
{\scriptstyle{(i\,=\,1\,\cdots\,n\,;\,\,\,\alpha\,\in\,\N^n)}}.
\]
To begin with, we look for solutions of the finite, extracted linear
system of $\rho_\infty$ equations with the $r$ unknowns $\tau_k ( a)$:
\[
\left\{
\aligned
0
&
\equiv
\tau_1(a)\,\frac{\partial{\sf u}_1}{\partial a_1}(a)
+\cdots+
\tau_{\rho_\infty}(a)\,
\frac{\partial{\sf u}_1}{\partial a_{\rho_\infty}}(a)
+\cdots+
\tau_r(a)\,\frac{\partial{\sf u}_1}{\partial a_r}(a)
\\
&
\cdots\cdots\cdots\cdots\cdots\cdots\cdots\cdots\cdots\cdots\cdots\cdots
\cdots\cdots\cdots\cdots\cdots\cdots\cdots
\\
0
&
\equiv
\tau_1(a)\,\frac{\partial{\sf u}_{\rho_\infty}}{\partial a_1}(a)
+\cdots+
\tau_{\rho_\infty}(a)\,
\frac{\partial{\sf u}_{\rho_\infty}}{\partial a_{\rho_\infty}}(a)
+\cdots+
\tau_r(a)\,\frac{\partial{\sf u}_{\rho_\infty}}{\partial a_r}(a).
\endaligned\right.
\]
After possibly renumbering the variables
$(a_1, \dots, a_r)$, we can assume that the left
$\rho_\infty \times \rho_\infty$ minor of this system:
\[
\Delta(a)
:=
\det\,
\Big(
\frac{\partial {\sf u}_l}{\partial a_m}(a)
\Big)_{1\leqslant m\leqslant\rho_\infty}^{1\leqslant l\leqslant\rho_\infty}
\]
does not vanish identically. However, it can vanish at some points,
and while endeavoring to solve the above linear system by an
application of the classical rule of Cramer, the necessary division by
the determinant $\Delta ( a)$ introduces poles that are undesirable,
for we want the $\tau_k ( a)$ to be analytic. So, for any $\mu$ with
$1 \leqslant \mu \leqslant r - \rho_\infty$, we look for a solution
(rewritten as a derivation) under the specific form:
\[
\mathcal{T}_\mu
:=
-\Delta(a)\,\frac{\partial}{\partial a_{\rho_\infty+\mu}}
+
\sum_{1\leqslant k\leqslant\rho_\infty}\,
\tau_{\mu k}(a)\,\frac{\partial}{\partial a_k}
\ \ \ \ \ \ \ \ \ \
{\scriptstyle{(\mu\,=\,1\,\cdots\,r-\rho_\infty)}},
\]
in which we introduce in advance a factor $\Delta( a)$ designed to
compensate the unavoidable
division by $\Delta ( a)$. Indeed, such a $\mathcal{
T}_\mu$ will annihilate the ${\sf u}_l$:
\[
0
\equiv 
\mathcal{T}_\mu{\sf u}_1
\equiv\cdots\equiv
\mathcal{T}_\mu{\sf u}_{\rho_\infty}
\]
if and only its coefficients are solutions of the linear system:
\[
\left\{
\aligned
\Delta(a)\,
\frac{\partial{\sf u}_1}{\partial a_{\rho_\infty+\mu}}(a)
&
\equiv
\tau_{\mu 1}(a)\,
\frac{\partial{\sf u}_1}{\partial a_1}(a)
+\cdots+
\tau_{\mu\rho_\infty}(a)\,
\frac{\partial{\sf u}_1}{\partial a_{\rho_\infty}}(a)
\\
\cdots\cdots\cdots\cdots\cdots
&
\cdots\cdots\cdots\cdots\cdots\cdots\cdots\cdots\cdots\cdots\cdots\cdots
\cdots
\\
\Delta(a)\,
\frac{\partial{\sf u}_{\rho_\infty}}{\partial a_{\rho_\infty+\mu}}(a)
&
\equiv
\tau_{\mu 1}(a)\,
\frac{\partial{\sf u}_{\rho_\infty}}{\partial a_1}(a)
+\cdots+
\tau_{\mu\rho_\infty}(a)\,
\frac{\partial{\sf u}_{\rho_\infty}}{\partial a_{\rho_\infty}}(a).
\endaligned\right.
\]
Cramer's rule then yields the unique solution:
\[
\tau_{\mu k}(a)
=
\frac{1}{\Delta(a)}\,
\left\vert
\begin{array}{ccccc}
\frac{\partial{\sf u}_1}{\partial a_1}(a)
& \cdots & 
\Delta(a)\,\frac{\partial{\sf u}_1}
{\partial a_{\rho_\infty+\mu}}(a) 
& \cdots &
\frac{\partial{\sf u}_1}{\partial a_{\rho_\infty}}(a)
\\
\cdots\cdots
& \cdots & 
\cdots\cdots\cdots
& \cdots &
\cdots\cdots
\\
\frac{\partial{\sf u}_{\rho_\infty}}{\partial a_1}(a)
& \cdots & 
\Delta(a)\,\frac{\partial{\sf u}_{\rho_\infty}}
{\partial a_{\rho_\infty+\mu}}(a)
& \cdots &
\frac{\partial{\sf u}_{\rho_\infty}}{\partial a_{\rho_\infty}}(a)
\end{array}
\right\vert
\ \ \ \ \ \ \ \ \ \
{\scriptstyle{(k\,=\,1\,\cdots\,r)}},
\]
where, as predicted, the overall factor $\Delta ( a)$ of the $k$-th column
compensates the division by the determinant $\Delta ( a)$ of this
system, so that the $\mathcal{ T}_\mu$ all have \emphasis{analytic}
coefficients indeed.

Clearly, $\mathcal{ T}_1, \dots, \mathcal{ T}_{ r - \rho_\infty}$ are
linearly independent at all (generic) points $a$ where $\Delta ( a)
\neq 0$. It remains to show that the $\mathcal{ T}_\mu$ also
annihilate all other coefficient functions $\mathcal{ U}_\alpha^i$.

Thus, let $(i, \alpha) \neq ( i ( 1), \alpha ( 1) ), \dots, ( i (
\rho_\infty) , \alpha ( \rho_\infty) )$. By the very definition of
$\rho_\infty$, the generic rank of any $r \times ( 1+ \rho_\infty)$
extracted subJacobian matrix:
\[
\left(
\begin{array}{cccc}
\frac{\partial\mathcal{U}_\alpha^i}{\partial a_1}(a)
&
\frac{\partial{\sf u}_1}{\partial a_1}(a)
&
\cdots
&
\frac{\partial{\sf u}_{\rho_\infty}}{\partial a_1}(a)
\\
\cdots & \cdots & \cdots & \cdots
\\
\frac{\partial\mathcal{U}_\alpha^i}{\partial a_r}(a)
&
\frac{\partial{\sf u}_1}{\partial a_r}(a)
&
\cdots
&
\frac{\partial{\sf u}_{\rho_\infty}}{\partial a_r}(a)
\end{array}
\right)
\]
must always be equal to $\rho_\infty$, the generic rank of its last
$\rho_\infty$ columns. Consequently, in a neighborhood of every point
$a$ at which its top right $\rho_\infty \times \rho_\infty$ minor
$\Delta ( a)$ does not vanish, the first column is a certain linear
combination:
\[
\frac{\partial\mathcal{U}_i^\alpha}{\partial a_k}(a)
=
\lambda_1(a)\,\frac{\partial{\sf u}_1}{\partial a_k}(a)
+\cdots+
\lambda_{\rho_\infty}(a)\,
\frac{\partial{\sf u}_{\rho_\infty}}{\partial a_k}(a)
\ \ \ \ \ \ \ \ \ \
{\scriptstyle{(i\,=\,1\,\cdots\,n\,;\,\,\,\alpha\,\in\,\N^n\,;
\,\,\,k\,=\,1\,\cdots\,r)}}
\]
of the last $\rho_\infty$ columns in question, where, again thanks to
an appropriate application of Cramer's rule, the coefficients
$\lambda_l ( a)$ are \emphasis{analytic} in the concerned generic
neighborhood, for their denominator $\Delta ( a)$ is $\neq 0$ there.
It then follows immediately by appropriate scalar multiplication and
summation that:
\[
\aligned
\sum_{k=1}^r\,\tau_{\mu k}(a)\,
\frac{\partial\mathcal{U}_\alpha^i}{\partial a_k}(a)
&
\equiv
\lambda_1(a)\,
\sum_{k=1}^r\,\tau_{\mu k}(a)\,
\frac{\partial{\sf u}_1}{\partial a_k}(a)
+\cdots+
\lambda_{\rho_\infty}(a)\,
\sum_{k=1}^r\,\tau_{\mu k}(a)\,
\frac{\partial{\sf u}_{\rho_\infty}}{\partial a_k}(a)
\\
&
\equiv
\lambda_1(a)\,\mathcal{T}_\mu\,{\sf u}_1
+\cdots+
\lambda_{\rho_\infty}(a)\,\mathcal{T}_\mu\,{\sf u}_{\rho_\infty}
\\
&
\equiv
0
\ \ \ \ \ \ \ \ \ \ \ \ \ \ \ \ \ \ \ \ \ \ \ \ \ \ \ \ \ \
\ \ \ \ \ \ \ \ \ \ \ \ \ \ \ \ \ \ \ \ \ \ \ \ \ \ \ \ \ \
{\scriptstyle{(i\,=\,1\,\cdots\,n\,;\,\,\,\alpha\,\in\,\N^n)}}.
\endaligned
\]
But since these analytic equations hold on the dense open set where
$\Delta ( a) \neq 0$, we deduce by continuity that the equations:
\[
0
\equiv
\mathcal{T}_1\,\mathcal{U}_\alpha^i
\equiv\cdots\equiv
\mathcal{T}_{r-\rho_\infty}\,\mathcal{U}_\alpha^i 
\ \ \ \ \ \ \ \ \ \ \ \ \
{\scriptstyle{(i\,=\,1\,\cdots\,n\,;\,\,\,\alpha\,\in\,\N^n)}}
\]
do hold \emphasis{everywhere}, as desired. In conclusion, we have
shown the implication {\bf (ii)} $\Rightarrow$ {\bf (iii)}, and
simultaneously, we have established the last part of the theorem.
\qed\end{proof}

\begin{corollary}
Locally in a neighborhood of every generic point $a_0$ at which the
infinite coefficient mapping $a \mapsto {\sf U}_\infty ( a)$ has
maximal, locally constant rank equal to its generic rank
$\rho_\infty$, there exist both a local change of parameters $a
\mapsto \big( {\sf u}_1 ( a), \dots, {\sf u}_{\rho_\infty} ( a) \big)
=: {\sf u}$ decreasing the number of parameters from $r$ down to
$\rho_\infty$, and new transformation equations:
\[
x_i'
=
g_i
\big(x;\,{\sf u}_1,\dots,{\sf u}_{\rho_\infty}\big)
\ \ \ \ \ \ \ \ \ \
{\scriptstyle{(i\,=\,1\,\cdots\,n)}}
\]
depending \emphasis{only} upon $\rho_\infty$ parameters 
which give again the old ones:
\[
g_i\big(x;\,{\sf u}(a)\big)
\equiv 
f_i(x;\, a)
\ \ \ \ \ \ \ \ \ \
{\scriptstyle{(i\,=\,1\,\cdots\,n)}}.
\]
\end{corollary}

\begin{proof}\smartqed
Choose $\rho_\infty$ coefficients $\mathcal{ U}_{ \alpha ( l)}^{ i (
l)} ( a) =: {\sf u}_l ( a)$, $1 \leqslant l \leqslant \rho_\infty$,
with $ \Delta ( a) : = {\rm det}\, \big( \frac{\partial {\sf u}_l (
a)}{\partial a_m} (a) \big)_{ 1\leqslant m \leqslant \rho_\infty}^{
1\leqslant l
\leqslant \rho_\infty} \not \equiv 0$ as in the proof of the
theorem. Locally in some small neighborhood of any $a^0$ with $\Delta
( a_0) \neq 0$, the infinite coefficient map ${\sf U}_\infty$ has
constant rank $\rho_\infty$, hence the constant rank theorem provides,
for every $( i, \alpha)$, a certain function $\mathcal{ V}_\alpha^i$
of $\rho_\infty$ variables such that:
\[
\mathcal{U}_\alpha^i(a)
\equiv
\mathcal{V}_\alpha^i
\big(
{\sf u}_1(a),\dots,{\sf u}_{\rho_\infty}(a)
\big).
\]
Thus, we can work out the power series expansion:
\[
\aligned
f_i(x;a)
&
=
{\textstyle{\sum}}_{\alpha\in\N^n}\,
\mathcal{U}_\alpha^i(a)\,(x-x_0)^\alpha
\\
&
=
{\textstyle{\sum}}_{\alpha\in\N^n}\,
\mathcal{V}_\alpha^i({\sf u}_1(a),\dots,{\sf u}_{\rho_\infty}(a))\,
(x-x_0)^\alpha
\\
&
=:
g_i(x,{\sf u}_1(a),\dots,{\sf u}_{\rho_\infty}(a)).
\endaligned
\]
which yields the natural candidate for $g_i ( x; \, {\sf u} )$.
Lastly, one may verify that any Cauchy estimate for the growth
decrease of $\mathcal{ U }_\alpha^i (a)$ as $\vert \alpha \vert \to
\infty$ insures a similar Cauchy estimate for the growth decrease of
$b \mapsto \mathcal{ V}_\alpha^i ( {\sf u})$, whence each $g_i$ is
analytic, and in fact, termwise substitution was legitimate.
\qed\end{proof}

\section{Axiom of Inverse and Engel's Counter-Example}
\label{axiom-of-inverse}

Every analytic diffeomorphism of an $n$-times extended space permutes
all the points in a certain differentiable, invertible way. Although
they act on a set of infinite cardinality, diffeomorphisms can thus be
thought to be sorts of analogs of the \emphasis{substitutions} on a
finite set. In fact, in the years 1873--80, Lie's Idée fixe was to
build, in the geometric realm of $n$-dimensional continua, a
counterpart of the Galois theory of substitutions of roots of
algebraic equations (\cite{ h2001}).

As above, let $x ' = f ( x; \, a_1, \dots, a_r) = : f_a ( x)$ be a
family of (local) analytic diffeomorphisms parameterized by a finite
number $r$ of parameters. For Lie, the basic, single group axiom
should just require that such a family be \emphasis{closed under
composition}, namely that one always has $f_a \big( f_b ( x) \big)
\equiv f_c ( x)$ for some $c$ depending on $a$ and on $b$. More
precisions on this definition will given in the next chapter, but at
present, we ask whether one can really economize the other two group
axioms: existence of an identity element and existence of inverses.

\begin{lemma}
If $H$ is any subset of some abstract group $G$ with ${\rm
Card}\, H < \infty$ which is closed under group multiplication:
\[
h_1 h_2\in H 
\ \ \ \ \ \ \ \ \ \ \ \ 
\text{\rm whenever} 
\ \ \ \ \ \ \ 
h_1,\,h_2\in H, 
\]
then $H$ contains the identity element $e$ of $G$ and every $h \in H$
has an inverse in $H$, so that $H$ itself is a true subgroup of $G$.
\end{lemma}

\begin{proof}\smartqed
Indeed, picking $h \in H$ arbitrary, the infinite sequence $h, h^2,
h^3, \dots, h^k, \dots$ of elements of the finite set $H$ must become
eventually periodic: $h^a = h^{ a + n}$ for some $a \geqslant 1$ and
for some $n \geqslant 1$, whence $e = h^n$, so $e \in H$ and $h^{ n -
1}$ is the inverse of
$h$.
\qed\end{proof}

For more than thirteen years, Lie was convinced that a purely similar
property should also hold with $G = {\sf Diff }_n$ being the (infinite
continuous pseudo)group of analytic diffeomorphisms and with $H
\subset {\sf Diff }_n$ being any continuous
family closed under composition. We quote a characteristic 
excerpt of~\cite{
lie1880}, pp.~444--445.

\begin{svgraybox}
As is known, one shows in the theory of substitutions that the
permutations of a group can be ordered into pairwise inverse couples
of permutations. Now, since the distinction between a permutation
group and a transformation group only lies in the fact that the former
contains a finite and the latter an infinite number of operations, it
is natural to presume that the transformations of a transformation
group can also be ordered into pairs of inverse transformations. In
previous works, I came to the conclusion that this should actually be
the case. But because in the course of my investigations in question,
certain \emphasis{implicit} hypotheses have been made about the nature
of the appearing functions, then I think that it is necessary to
\emphasis{expressly add the requirement that the transformations of
the group can be ordered into pairs of inverse transformations}. In
any case, I conjecture \label{conjecture-composition}
that this is a necessary consequence of my
original definition of the concept \deutsch{Begriff} of transformation
group. However, it has been impossible for me to prove this in
general.
\end{svgraybox}

\smallskip
As a proposal of counterexample that Engel devised in the first year he
worked with Lie (1884), consider the family of transformation equations:
\[
x'
=
\zeta\,x,
\]
where $x, \, x' \in \C$ and the parameter $\zeta \in \C$ is restricted
to $\vert \zeta \vert < 1$. Of course, this family is closed under any
composition, say: $x' = \zeta_1 \, x$ and $x'' = \zeta_2 \, x' =
\zeta_1 \zeta_2 x$, with indeed $\vert \zeta_2 \,\zeta_1 \vert < 1$
when $\vert \zeta_1 \vert, \, \vert \zeta_2 \vert < 1$, but neither
the identity element nor any inverse transformation does belong to the
family. However, the requirement $\vert \zeta \vert < 1$ is here too
artificial: the family extends in fact trivially as the complete group
$\big( x' = \zeta\, x \big)_{ \zeta \in \C}$ of dilations of the
line. Engel's idea was to appeal to a Riemann map $\omega$ having $\{
\vert \zeta \vert = 1 \}$ as a frontier of nonextendability. 
The map used by Engel is the following\footnote{\,
In the treatise~\cite{ enlie1888-2}, this example is
presented at the end of Chap.~\ref{kapitel-9}, 
\voir~below p.~\pageref{univalent-odd}.
}. 
Let ${\sf od}_k$ denote the number of odd divisors (including $1$) 
of any integer $k\geqslant 1$. The theory of
holomorphic functions in one complex variables yields the
following. 

\begin{lemma}
The infinite series:
\[
\omega(a)
:=
\sum_{\nu\geqslant 1}\,
\frac{a^\nu}{1-a^{2\nu}}
=
\sum_{\nu\geqslant 1}\,
\big(
a^\nu+a^{3\nu}+a^{5\nu}+a^{7\nu}+\cdots
\big)
=
\sum_{k\geqslant 1}\,
{\sf od}_k\,a^k
\]
converges absolutely in every open disc $\Delta_\rho = \{ z\in \C: \,
\, \vert z \vert < \rho \}$ of radius $\rho < 1$ and defines a
univalent holomorphic function $\Delta \to \C$ from the unit disc
$\Delta := \{ \vert z \vert < 1 \}$ to $\C$ which \emphasis{does not}
extend holomorphically across any point of the unit circle $\partial
\Delta := \{ \vert z \vert = 1\}$.
\end{lemma}

\label{engel-counterexample}
In fact, any other similar Riemann biholomorphic map $\zeta
\longmapsto \omega ( \zeta) =: \lambda$ from the unit disc $\Delta$
onto some simply connected domain $\Lambda := \omega ( \Delta)$ having
fractal boundary not being a Jordan curve, as e.g. the Von Koch
Snowflake Island, would do the job\footnote{\, A concise presentation
of Carathéodory's theory may be found in Chap.~17 of~\cite{
mi2006}.}. Denote then by $\lambda \longmapsto \chi ( \lambda) =:
\zeta$ the inverse of such a map
and consider the family of transformation equations:
\[
\big(
x'
=
\chi(\lambda)\, x
\big)_{\lambda\in\Lambda}.
\]
By construction, $\vert \chi ( \lambda) \vert < 1$ for every $\lambda
\in \Lambda$. Any composition of $x' = \chi ( \lambda_1) \, x$ and of
$x'' = \chi ( \lambda_2) \, x'$ is of the form $x'' = \chi ( \lambda)
\, x$, with the uniquely defined parameter $\lambda := \omega \big(
\chi ( \lambda_1)\, \chi ( \lambda_2) \big)$, hence the group
composition axiom is satisfied. However there is again no identity
element, and again, none transformation has an inverse. And
furthermore crucially (and lastly), there does not exist any
prolongation of the family to a larger domain $\widetilde{ \Lambda}
\supset \Lambda$ together with a holomorphic prolongation $\widetilde{
\chi}$ of $\chi$ to $\widetilde{ \Lambda}$ so that $\widetilde{ \chi}
\big( \widetilde{ \Lambda} \big)$ contains a neighborhood of $\{ 1 \}$
(in order to catch the identity) or \emphasis{a fortiori}\, a
neighborhood of $\overline{ \Delta}$ (in order to catch inverses of
transformations $x' = \chi( \lambda) \, x$ with $\lambda \in
\Lambda$ close to $\partial \Lambda$).

\subsubsection*{Observation}
In Vol.~I of the \emphasis{Theorie der Transformationsgruppen}, this
example appears only in Chap.~9, on pp.~163--165, and it is written
in small characters. In fact, Lie still believed that a deep analogy
with substitution groups should come out as a theorem. Hence
\emphasis{the structure of the first nine chapters insist on setting
aside}, whenever possible, \emphasis{the two axioms of existence of
identity element and of existence of inverses}. To do justice to this
great treatise, we shall translate in Chap.~\ref{kapitel-9} 
how Master Lie managed
to produce the Theorem~26 on p.~\pageref{Theorem-26}, 
which he considered to provide
the sought analogy with finite group theory, after taking Engel's
counterexample into account.

\linestop


\chapter{Fundamental Differential Equations
\\
for Finite Continuous Transformation Groups}
\label{fundamental-differential}
\chaptermark{Fundamental Differential Equations
for Finite Continuous Transformation Groups}

\abstract{
A \terminology{finite continuous local transformation group}\, in the
sense of Lie is a family of local analytic diffeomorphisms $x_i ' =
f_i ( x; \, a)$, $i = 1 \dots, n$, parametrized by a finite number $r$
of parameters $a_1, \dots, a_r$ that is closed under composition and
under taking inverses:
\[
f_i\big(f(x;\,a);\,b)\big)
=
f_i\big(x;\,\text{\bf m}(a,b)\big)
\ \ \ \ \ \ \ 
\text{\rm and}
\ \ \ \ \ \ \
x_i
=
f_i\big(x';\,\text{\bf i}(a)\big),
\]
for some \terminology{group multiplication map}\, $\text{\bf m}$ and
for some \terminology{group inverse map}\, $\text{\bf i}$, both local
and analytic. Also, it is assumed that there exists $e = (e_1, \dots,
e_r)$ yielding the \terminology{identity}\, transformation: $f_i ( x;
\, e) \equiv x_i$.
\newline\indent Crucially, these requirements imply the existence of
\emphasis{fundamental partial differential equations}:
\[
\boxed{
\frac{\partial f_i}{\partial a_k}(x;\,a)
=
-
\sum_{j=1}^r\,\psi_{kj}(a)\,\frac{\partial f_i}{\partial a_j}(x;\,e)}
\ \ \ \ \ \ \ \ \ \ 
{\scriptstyle{(i\,=\,1\,\cdots\,n,\,\,\,k\,=\,1\,\cdots\,r)}}
\]
which, technically speaking, are cornerstones of the basic theory.
What matters here is that the group axioms
guarantee that the $r\times
r$ matrix $( \psi_{ kj})$ depends only on $a$ and it is locally
invertible near the identity. Geometrically speaking, these equations
mean that the $r$ infinitesimal transformations:
\[
X_k^a\big\vert_x
=
\frac{\partial f_1}{\partial a_k}(x;\,a)\,\frac{\partial}{\partial x_1}
+\cdots+
\frac{\partial f_n}{\partial a_k}(x;\,a)\,\frac{\partial}{\partial x_n}
\ \ \ \ \ \ \ \ \ \ \ \ \
{\scriptstyle{(k\,=\,1\,\cdots\,r)}}
\]
corresponding to an infinitesimal increment of the $k$-th parameter
computed at $a$:
\[
f(x;a_1,\dots,a_k+\varepsilon,\dots,a_r) 
- 
f(x;\,a_1,\dots,a_k,\dots,a_r)
\approx
\varepsilon
X_k^a\big\vert_x
\] 
are \emphasis{linear combinations}, with certain coefficients
$- \psi_{ kj} (a)$ \emphasis{depending only on the parameters}, of the
same infinitesimal transformations computed at the
identity:
\[
X_k^e\big\vert_x
=
\frac{\partial f_1}{\partial a_k}(x;\,e)\,\frac{\partial}{\partial x_1}
+\cdots+
\frac{\partial f_n}{\partial a_k}(x;\,e)\,\frac{\partial}{\partial x_n}
\ \ \ \ \ \ \ \ \ \ \ \ \
{\scriptstyle{(k\,=\,1\,\cdots\,r)}}.
\]
Remarkably, the process of removing superfluous parameters introduced
in the previous chapter applies to local Lie groups without the
necessity of relocalizing around a generic $a_0$, so that everything
can be achieved around the identity $e$ itself, 
without losing it.
}

\section{Concept of Local Transformation Group}
\label{concept-local-group}

Let $\K = \R$ or $\C$, let $n \geqslant 1$ be an integer and let $x =
(x_1, \dots, x_n) \in \K^n$ denote variables of an $n$-times extended
space. We shall constantly employ the sup-norm:
\[
\vert x\vert
:=
\max_{1\leqslant i\leqslant n}\,\vert x_i\vert,
\]
where $\vert \cdot \vert$ denotes the absolute value on $\R$, or the
modulus on $\C$. For various ``radii'' $\rho > 0$, we shall consider
the precise open sets centered at the origin that are defined by:
\[
\Delta_\rho^n
:=
\big\{x\in\K^n:\,
\vert x\vert<\rho
\big\};
\]
in case $\K = \C$, these are of course standard open polydiscs,
while in case $\K = \R$, these are open cubes.

On the other hand, again\footnote{\,
When $x\in \R^n$ is real, while studying local analytic Lie group
actions $x_i' = f_i ( x; a)$ below, we will naturally require that $a
\in \R^r$ also be real (unless we complexify both $x$ and $a$). When
$x\in \C^n$, we can suppose $a$ to be either real or complex. } with
$\K = \R$ or $\C$, let $r \geqslant 1$ be another integer and
introduce parameters $a = ( a_1, \dots, a_r)$ in $\K^r$, again
equipped with the sup-norm:
\[
\vert
a
\vert
:=
\max_{1\leqslant k\leqslant r}\,
\vert a_k\vert.
\]
For various $\sigma > 0$ similarly, we introduce the precise open
sets:
\[
\square_\sigma^r
:=
\big\{
a\in\K^r:\,\vert a\vert<\sigma
\big\}.
\]

\subsection{Transformation Group Axioms}

Let 
\[
x_i'
=
f_i(x_1,\dots,x_n;a_1,\dots,a_r)
\ \ \ \ \ \ \ \ \ \
{\scriptstyle{(i\,=\,1\,\cdots\,n)}}
\]
be a local transformation equations, as presented 
in Chapter~\ref{three-principles-thought}.
To fix the local character, the $f_i (x; a)$ will be assumed to be
defined when $\vert x \vert < \rho_1$ and when $\vert a \vert <
\sigma_1$, for some $\rho_1 > 0$ and for some $\sigma_1 > 0$.
We shall assume that for the parameter $a := e$
equal to the origin $0 \in
\K^r$, the transformation
corresponds to the identity, so that:
\[
f_i(x_1,\dots,x_n:0,\dots,0)
\equiv
x_i
\ \ \ \ \ \ \ \ \ \
{\scriptstyle{(i\,=\,1\,\cdots\,n)}}.
\]
Consequently, 
for the composition of two successive such transformations $x' = f (
x; a)$ and $x'' = f ( x' ; b)$ to be well defined,
it suffices to shrink $\rho_1$
to $\rho_2$ with $0 < \rho_2 < \rho_1$ and $\sigma_1$ to $\sigma_2$
with $0 < \sigma_2 < \sigma_1$ in order to insure that:
\[
\vert
f(x;a)
\vert
<
\rho_1
\ \ \ \ \ \ \ \
\text{\rm for all}
\ \ \ \ \
\vert x\vert<\rho_2\ \ \
\text{\rm and all}\ \ \
\vert a\vert<\sigma_2.
\]
This is clearly possible thanks to $f(x;0) = x$. Now, we can present
the {\em local} transformation group axioms, somehow with a rigorous
control of existence domains.

\medskip

\begin{center}
\input shrink-polydiscs.pstex_t
\end{center}

\noindent{\sf Group composition axiom.} For every $x\in
\Delta_{\rho_2}^n$, and $a, b \in \square_{\sigma_2}^r$, an arbitrary
composition:
\def\theequation{1}\begin{equation}
\label{f-f-composition}
x''
=
f\big(f(x;a);b\big)
=
f(x;c)
=
f(x;\text{\bf m}(a,b))
\end{equation}
always identifies to an element $f ( x; c)$ of the same family, for a
unique parameter $c = \text{\bf m} ( a,b)$ given by an auxiliary
\terminology{group-multiplication}\, local analytic map:
\[
\text{\bf m}:\ \ \ \
\square_{\sigma_1}^n\times\square_{\sigma_1}^n
\longrightarrow
\K^r
\]
which satisfies $\text{\bf m} \big( \square_{\sigma_2}^r \times
\square_{ \sigma_2}^r \big) \subset \square_{ \sigma_1}^r$ and
$\text{\bf m} ( a, e) \equiv \text{\bf m} ( e, a) \equiv a$.

For $a, b, c \in \square_{ \sigma_3}^n$ with $0 < \sigma_3 < \sigma_2
< \sigma_1$ small enough so that three successive compositions are
well defined, the associativity of diffeomorphism composition yields:
\[
f\big(x;\text{\bf m}(\text{\bf m}(a,b),c)\big)
=
f\big(f(f(x;a);b);c\big)
=
f
\big(
x;
\text{\bf m}(a,\text{\bf m}(b,c))\big),
\]
whence, thanks to the supposed uniqueness of $c = \text{\bf m} ( a,
b)$, it comes the group associativity: $\text{\bf m} \big( \text{\bf
m} ( a,b ), c\big) = \text{\bf m} \big( a, \text{\bf m} ( b, c)\big)$
for such restricted values of $a, b, c$.

Contrary to what his opponents sometimes claimed, e.g. Study, Slocum
and others, Lie was conscious of the necessity of emphasizing the
local character of transformation groups that is often required in
applications. Amongst the first 50 pages of the \emphasis{Theorie der
Transformationsgruppen}, at least 15 pages (written in small
characters) are devoted to rigorously discuss when and why domains of
definition should be shrunk. We translate for instance a relevant
excerpt (\cite{ enlie1888-3}, pp.~15--16) in which the symbol $( x)$,
due to Weierstrass and introduced by Engel and Lie just before,
denotes a region of the coordinate space and $(a)$ a region of the
parameter space, with the $f_i$ analytic there.

\begin{svgraybox}\indent
Here one has to observe that we have met the fixing about the
behavior of the functions $f_i ( x; a)$ only inside the regions $(x)$
and $(a)$.
\label{excerpt-shrinking}

Consequently, we have the permission to substitute the expression
$x_\nu ' = f_\nu ( x, a)$ in the equations $x_i'' = f_i ( x', b)$ only
when the system of values $x_1', \dots, x_n'$ lies in the region
$(x)$. That is why we are compelled to add, to the fixings met up to
now about the regions $(x)$ and $(a)$, yet the following assumption:
it shall be possible to indicate, inside the regions $(x)$ and $(a)$,
respective subregions $(\!( x )\!)$ and $(\!( a )\!)$ 
\label{S-16}
of such a nature
that the $x_i'$ always remain in the region $(x)$ when the $x_i$ run
arbitrarily in $(\!( x )\!)$ and when the $a_k$ run arbitrarily in
$(\!( a )\!)$; we express this briefly as: the region $x' = f\big(
(\!( x )\!) \, (\!( a )\!) \big)$ shall entirely fall into
the region $(x)$.

According to these fixings, if we choose $x_1, \dots, x_n$ in the region
$(\!( x )\!)$ and $a_1, \dots, a_r$ in the region $(\!( a )\!)$, we
then can really execute the substitution $x_k' = f_k ( x, a)$ in the
expression $f_i ( x_1', \dots, x_n', \, b_1, \dots, b_n)$; that is to
say, when $x_1^0, \dots, x_n^0$ means an arbitrary system of values in
the region $(\!( x)\!)$, the expression:
\[
f_i\big(
f_1(x,a),\dots,f_n(x,a),\,b_1,\dots,b_n
\big)
\] 
can be expanded, in the neighborhood of the system of values $x_k^0$,
as an ordinary power series in $x_1 - x_1^0, \dots, x_n - x_n^0$; the
coefficients of this power series are functions of $a_1, \dots, a_r$,
$b_1, \dots, b_r$ and behave regularly, when the $a_k$ are arbitrary
in $(\!( a )\!)$ and the $b_k$ are arbitrary in $(a)$.
\end{svgraybox}

\smallskip\noindent{\sf Existence of an inverse-element map.}
There exists a local analytic map:
\[
\text{\bf i}:\ \ \ 
\square_{\sigma_1}^r\longrightarrow\K^r
\]
with $\text{\bf i} ( e) = e$ (namely $\text{\bf i} ( 0) = 0$)
such that for every $a \in \square_{ \sigma_2}^r$:
\[
\aligned
&
e
=
\text{\bf m}(a,\text{\bf i}(a)\big)
=
\text{\bf m}\big(\text{\bf i}(a),a\big)
\\
&
\text{\rm whence in addition:}
\ \ \ \ \ 
x
=
f\big(f(x;a);\text{\bf i}(a)\big)
=
f\big(f(x;\text{\bf i}(a));a),
\endaligned
\]
for every $x \in \Delta_{\rho_2}^n$.

\subsection{Some Conventions}

In the sequel, the diffeomorphism $x \mapsto f( x; a)$ will
occasionally be written $x \mapsto f_a ( x)$. Also, we shall
sometimes abbreviate $\text{\bf m} ( a,b)$ by $a \cdot b$ and
$\text{\bf i} ( a)$ by $a^{ - 1}$. Also, a finite number of times, it
will be necessary to shrink again $\rho_2$ and $\sigma_2$. This will
be done automatically, without emphasizing it.

\section{Changes of Coordinates and of Parameters}

\begin{svgraybox}
In the variables $x_1, \dots, x_n$, let the equations:
\[
x_i'
=
f_i(x_1,\dots,x_n,\,a_1,\dots,a_r)
\ \ \ \ \ \ \ \ \ \ \ \ \ {\scriptstyle{(i\,=\,1\,\cdots\,n)}}
\]
of an $r$-term group be presented. 
Then there are various means to derive from these equations
other equations which represent again an $r$-term group.

On the first hand, in place of the $a$, we can introduce
$r$ arbitrary independent functions of them:
\[
\overline{a}_k
=
\beta_k(a_1,\dots,a_r)
\ \ \ \ \ \ \ \ \ \ \ \ \ {\scriptstyle{(k\,=\,1\,\cdots\,r)}}
\] 
as new parameters. By resolution with respect to 
$a_1, \dots, a_r$, one can obtain:
\[
a_k
=
\gamma_k
(\overline{a}_1,\dots,\overline{a}_r)
\ \ \ \ \ \ \ \ \ \ \ \ \ {\scriptstyle{(k\,=\,1\,\cdots\,r)}}
\]
and by substitution of these values, one may set:
\[
f_i(x_1,\dots,x_n,\,a_1,\dots,a_r)
=
\mathfrak{f}_i(x_1,\dots,x_n,\,
\overline{a}_1,\dots,\overline{a}_r).
\]
Then if we yet set:
\[
\beta_k(b_1,\dots,b_r)
=
\overline{b}_k,
\ \ \ \ \ \ \ \ \ \ \
\beta_k(c_1,\dots,c_r)
=
\overline{c}_k
\ \ \ \ \ \ \ \ \ \ \ \ \ {\scriptstyle{(k\,=\,1\,\cdots\,r)}},
\]
the composition equations:
\[
f_i\big(f_1(x,a),\,\dots,\,f_n(x,a),\,\,
b_1,\dots,b_r\big)
=
f_i(x_1,\dots,x_n,\,c_1,\dots,c_r)
\]
take, without effort, the form:
\[
\mathfrak{f}_i
\big(
\mathfrak{f}_1(x,
\overline{a}),\,\dots,\,\mathfrak{f}_n(x,
\overline{a}),\,\,
\overline{b}_1,\dots,\overline{b}_r
\big)
=
\mathfrak{f}_i(x_1,\dots,x_n,\,
\overline{c}_1,\dots,\overline{c}_r),
\]
from which it results that the equations:
\[
x_i'
=
\mathfrak{f}_i(x_1,\dots,x_n,\,
\overline{a}_1,\dots,\overline{a}_r)
\ \ \ \ \ \ \ \ \ \ \ \ \ {\scriptstyle{(i\,=\,1\,\cdots\,n)}}
\]
with the $r$ essential parameters $\overline{ a}_1,
\dots, \overline{ a}_r$ represent in the same
way an $r$-term group.

Certainly, the equations of this new group are different from the ones
of the original group, but these equations obviously represent exactly
the same transformations as the original equations $x_i' = f_i ( x,
a)$.  Consequently, the new group is fundamentally identical to the
old one.

On the other hand, we can also introduce new independent
variables $y_1, \dots, y_n$ in place of the $x$:
\[
y_i
=
\omega_i(x_1,\dots,x_n)
\ \ \ \ \ \ \ \ \ \ \ \ \ {\scriptstyle{(i\,=\,1\,\cdots\,n)}},
\]
or, if resolved:
\[
x_i
=
w_i(y_1,\dots,y_n)
\ \ \ \ \ \ \ \ \ \ \ \ \ {\scriptstyle{(i\,=\,1\,\cdots\,n)}}.
\]
Afterwards, we have to set:
\[
x_i'
=
w_i(y_1',\dots,y_n')
=
w_i',
\ \ \ \ \ \ \ \ \ \
x_i''
=
w_i(y_1'',\dots,y_n'')
=
w_i'',
\]
and we hence obtain, in place of the transformation equations:
\[
x_i'
=
f_i(x_1,\dots,x_n,\,a_1,\dots,a_r),
\]
the following ones:
\[
w_i(y_1',\dots,y_n')
=
f_i(w_1,\dots,w_n,\,a_1,\dots,a_r),
\]
or, by resolution:
\[
y_i'
=
\omega_i
\big(
f_1(w,a),\,\dots,\,f_n(w,a)
\big)
=
\mathfrak{F}_i(y_1,\dots,y_n,\,a_1,\dots,a_r).
\]

It is easy to prove that the equations:
\[
y_i'
=
\mathfrak{F}_i
(y_1,\dots,y_n,\,a_1,\dots,a_r)
\ \ \ \ \ \ \ \ \ \ \ \ \ {\scriptstyle{(i\,=\,1\,\cdots\,n)}}
\]
with the $r$ essential parameters $a_1, \dots, a_r$
again represent an $r$-term group. In fact, the known equations:
\[
f_i\big(
f_1(x,a),\,\dots,\,f_n(x,a),\,\,
b_1,\dots,b_r\big)
=
f_i(x_1,\dots,x_n,\,c_1,\dots,c_r)
\]
are transferred, after the introduction of the new variables, to:
\[
f_i\big(f(w,a),\,b\big)
=
f_i(w_1,\dots,w_n,\,c_1,\dots,c_r),
\]
which can also be written:
\[
f_i(w_1',\dots,w_n',\,b_1,\dots,b_r)
=
f_i(w_1,\dots,w_n,\,c_1,\dots,c_r)
=
w_i''\,;
\]
but from this, it comes by resolution with respect to
$y_1'', \dots, y_n''$:
\[
y_\nu''
=
\omega_\nu
\big(f_1(w',b),\,\dots,\,f_n(w',b)\big)
=
\omega_\nu
\big(
f_1(w,c),\dots,f_n(w,c)
\big),
\]
or, what is the same:
\[
y_\nu''
=
\mathfrak{F}_\nu(y_1',\dots,y_n',\,b_1,\dots,b_r)
=
\mathfrak{F}_\nu(y_1,\dots,y_n,\,c_1,\dots,c_r),
\]
that is to say: there exist the equations:
\[
\mathfrak{F}_\nu
\big(
\mathfrak{F}_1(y,a),\,\dots,\,\mathfrak{F}_n(y,a),\,\,
b_1,\dots,b_r
\big)
=
\mathfrak{F}_\nu(y_1,\dots,y_n,\,c_1,\dots,c_r),
\]
whence it is indeed proved that the equations $y_i' = \mathfrak{ F}_i
( y,a)$ represent a group.

Lastly, we can naturally introduce at the same time new parameters and
new variables in a given group; it is clear that in this way, we 
likewise obtain a new group from the original group.

At present, we set up the following definition:

\smallskip{\sf Definition.}
\label{S-24}
Two $r$-term groups:
\[
\aligned
x_i'
&
=
f_i(x_1,\dots,x_n,\,a_1,\dots,a_r)
\ \ \ \ \ \ \ \ \ \ \ \ \ {\scriptstyle{(i\,=\,1\,\cdots\,n)}}
\\
y_i'
&
=
\mathfrak{f}_i(y_1,\dots,y_n,\,b_1,\dots,b_r)
\ \ \ \ \ \ \ \ \ \ \ \ \ {\scriptstyle{(i\,=\,1\,\cdots\,n)}}
\endaligned
\]
in the same number of variables are \terminology{similar}
\deutsch{ähnlich} to each other as soon as the one converts into
the other by the introduction of appropriate new variables and of
appropriate new parameters.

Obviously, there is an unbounded number of groups which are similar to
a given one; but all these unboundedly numerous groups are known
simultaneously with the given one. For this reason, as it shall also
happen in the sequel, we can consider that two mutually similar groups
are not essentially distinct from each other.

\smallercharacters{

Above, we spoke about the introduction of new parameters and of new
variables without dealing with the assumptions by which we can
ascertain that all group-theoretic properties essential for us are
preserved here.  Yet a few words about this point.

For it to be permitted to introduce, in the group $x_i' = f_i ( x_1,
\dots, x_n, \, a_1, \dots, a_r)$, the new parameters $\overline{ a}_k
= \beta_k ( a_1, \dots, a_r)$ in place of the $a$, the $\overline{
a}_k$ must be univalent functions of the $a$ in the complete region
$(a)$ defined earlier on, and they must behave regularly everwhere in
it; the functional determinant $\sum \, \pm \,
\partial \beta_1 / \partial a_1 \, \cdots\, \, 
\partial \beta_r / \partial a_r$ should vanish
nowhere in the region $(a)$, and lastly, to two distinct systems of
values $a_1, \dots, a_r$ of this region, there must always be
associated two distinct systems of values $\overline{ a}_1, \dots,
\overline{ a}_r$.  In other words: in the region of the $\overline{
a}_k$, one must be able to delimit a region $( \overline{ a})$ on
which the systems of values of the region $(a)$ are represented in a
univalent way by the equations $\overline{ a}_k =
\beta_k ( a_1, \dots, a_r)$. 

On the other hand, for the introduction of the new variables $y_i =
\omega_i ( x_1, \dots, x_n)$ to be allowed, the $y$ must be univalent
and regular functions of the $x$ for all sytems of values $x_1, \dots,
x_n$ which come into consideration after establishing the
group-theoretic properties of the equations $x_i' = f_i ( x_1, \dots,
x_n, \, a_1, \dots, a_r)$; inside this region, the functional
determinant $\sum\, \pm \, \partial \omega_1 /
\partial x_1 \, \cdots\,\, \partial \omega_n / \partial x_n$
should vanish nowhere, and lastly, to two distinct systems of values
$x_1, \dots, x_n$ of this region, there must always be
associated two distinct systems of values $y_1, \dots, y_n$.  The
concerned system of values of the $x$ must therefore be represented
univalently onto a certain region of systems of values of the $y$.

If one would introduce, in the group $x_i' = f_i ( x, a)$, new
parameters or new variables without the requirements just explained
being satisfied, then it would be thinkable in any case that important
properties of the group, for instance the group composition property
itself, would be lost; a group with the identity transformation could
convert into a group which does not contain the identity
transformation, and conversely.

But it certain circumstances, the matter is only to study
the family of transformations $x_i' = f_i ( x, a)$ 
in the neighbourhood of a single point $a_1, \dots, a_r$ or
$x_1, \dots, x_n$. This study will often be facilitated 
by introducing new variables or new parameters
which satisfy the requirements mentioned above 
in the \emphasis{neighbourhood of the concerned
points}. 

In such a case, one does not need at all to deal with 
the question whether the concerned requirements are satisfied
in the whole extension of the regions $(x)$ and $(a)$. 

}

Now, if we have a family of $\infty^r$ transformations: 
\[
x_i'
=
f_i(x_1,\dots,x_n,\,a_1,\dots,a_r)
\]
which forms an $r$-term group, then there corresponds to this family a
family of $\infty^r$ operations by which the points of the space $x_1,
\dots, x_n$ are permuted. Evidently, any two of these $\infty^r$
operations, when executed one after the other, always produce an
operation which again belongs to the family.

Thus, if we actually call a family of operations of this sort a
group-operation
\deutsch{Operationsgruppe}, or shortly, a group, 
then we can say: \emphasis{every given $r$-term transformation group
can be interpreted as the analytic representation of a certain group
of $\infty^r$ permutations of the points $x_1, \dots, x_n$}.

Conversely, if a group of $\infty^r$ permutations of the points $x_1,
\dots, x_n$ is given, and if it is possible to represent these
permutations by \emphasis{analytic} transformation equations, then the
corresponding $\infty^r$ transformations naturally form a
transformation group.

Now, if one imagines that a determined group-operation 
is given, and in addition, that an analytic representation 
of it is given\,---\,hence if one has a transformation group\,---,
then this
representation has in itself two obvious 
incidental characters
\deutsch{Zufälligkeit}. 

The first incidental character is the choice of the parameters
$a_1, \dots, a_r$. It stands to reason that this choice
has in itself absolutely no influence on the group-operation, 
when we introduce in place of the $a$ the new
parameters $\overline{ a}_k = \beta_k ( a_1, \dots, a_r)$. 
Only the analytic expression for the group-operation
will be a different one on the occasion; therefore, this
expression represents a transformation group as
before. 

The second incidental character in the analytic representation of our
group-operation is the choice of the coordinates in the space $x_1,
\dots, x_n$. Every permutation of the points $x_1, \dots, x_n$ is
fully independent of the choice of the system of coordinates to which
one refers the points $x_1, \dots, x_n$; only the analytic
representation of the permutation changes with the concerned system of
coordinates.  Naturally, the same holds for any group of permutations.
From this, it results that by introducing new variables, that is to
say, by a change of the system of coordinates, one obtains, from a
transformation group, again a transformation group, for the
transformation equations that one receives after the introduction of
the new variables represent exactly the same group-operation as the
one of the initial transformation group, so that they form in turn a
transformation group.

Thus, the analytic considerations of the previous paragraphs are now
explained in a conceptual way.  Above all, it is at present clear why
two similar transformation groups are to be considered as not
essentially distinct from each other; namely, for the reason that they
both represent analytically one and the same group-operation.
\end{svgraybox}

\section{Geometric Introduction of Infinitesimal Transformations}

Now, for reasons of understandability, we shall present in advance the
basic geometric way in which infinitesimal transformations can be
introduced, a way which is knowingly passed over silence in the great
treatise~\cite{ enlie1888-3}.

Letting $\varepsilon$ denote either an infinitesimal quantity in the
sense of Leibniz, or a small quantity subjected to Weierstrass'
rigorous epsilon-delta formalism, for fixed $k \in \{ 1, 2, \dots, r
\}$, we consider all the points:
\[
\label{infinitesimal-transformations}
\aligned
x_i'
&
=
f_i
\big(x;\,e_1,\dots,e_k+\varepsilon,\dots,e_n\big)
\\
&
=
x_i
+
\frac{\partial f_i}{\partial a_k}(x;e)\,
\varepsilon
+\cdots
\ \ \ \ \ \ \ \ \ \
{\scriptstyle{(i\,=\,1\,\cdots\,n)}}
\endaligned
\]
that are infinitesimally pushed from the starting points $x = f ( x;
e)$ by adding the tiny increment $\varepsilon$ to only the $k$-th
identity parameter $e_k$. One may reinterpret this common spatial move
by introducing the vector field (and a new notation for its
coefficients):
\[
\aligned
X_k^e
:=
\sum_{i=1}^n\,
\frac{\partial f_i}{\partial a_k}(x;e)\,
\frac{\partial}{\partial x_i}
&
:=
\sum_{i=1}^n\,
\xi_{ki}(x)\,
\frac{\partial}{\partial x_i},
\endaligned
\] 
which is either written as a derivation in modern style, or considered
as a column vector:
\[
{}^\tau \, \big( 
{\textstyle{ \frac{ \partial f_1 }{\partial a_k}}}, \cdots,
{\textstyle{ 
\frac{ \partial f_n }{\partial a_k}}} \big) \Big\vert_x = 
{}^\tau \big( \xi_{ k1}, \dots, \xi_{ kn} \big)
\Big\vert_x 
\] 
based at $x$, where ${}^\tau ( \cdot)$ denotes a transposition,
yielding column vectors. Then
$x ' = x + \varepsilon\, X_k^e + \cdots$, or equivalently:
\[
x_i'
=
x_i
+
\varepsilon\,\xi_{ki}
+\cdots
\ \ \ \ \ \ \ \ \ \ \ \ \ \
{\scriptstyle{(i\,=\,1\,\cdots\,n)}},
\]
where the left out terms ``$+ \cdots$'' are of course an ${\rm O} (
\varepsilon^2)$, so that from the geometrical viewpoint, 
$x '$ is infinitesimally pushed along the vector $X_k^e 
\big\vert_x$ up to a
length $\varepsilon$.

\begin{center}
\input push-points.pstex_t
\end{center}
\label{push-points}

More generally, still starting from the identity parameter $e$, 
when we add to $e$ an arbitrary infinitesimal increment:
\[
\big(
e_1+\varepsilon\,\lambda_1,\dots,
e_k+\varepsilon\,\lambda_k,\dots,
e_r+\varepsilon\,\lambda_r
\big),
\]
where ${}^\tau(\lambda_1, \dots, \lambda_r) \big\vert_e$ is a fixed,
constant vector based at $e$ in the parameter space, it follows by
linearity of the tangential map, or else just by the chain rule in
coordinates, that:
\[
\aligned
f_i(x;\,e+\varepsilon\,\lambda)
&
=
x_i
+
\sum_{k=1}^n\,\varepsilon\,\lambda_k\,
\frac{\partial f_i}{\partial a_k}(x;e)
+\cdots
\\
&
=
x_i
+
\varepsilon\,
\sum_{k=1}^n\,\lambda_k\,\xi_{ki}(x)
+\cdots,
\endaligned
\]
so that all points $x ' = x + \varepsilon \, X + \cdots$
are infinitesimally and simultaneously pushed along the vector
field:
\[
X
:=
\lambda_1\,X_1^e
+\cdots+
\lambda_r\,X_r^e
\]
which is the general linear combination of the $r$ previous basic
vector fields $X_k^e$, $k=1, \dots, r$.

Occasionally, Lie wrote that such a vector field $X$
\terminology{belongs to the group} $x ' = f ( x; a)$, to mean that $X$
comes itself with the infinitesimal move $x ' = x + \varepsilon \, X$
it is supposed to perform (dots should now be suppressed in
intuition), and hence accordingly, Lie systematically called such an
$X$\, an \terminology{infinitesimal transformation}, viewing indeed $x
' = x + \varepsilon \, X$ as just a case of $x ' = f ( x, a)$.
Another, fundamental and \emphasis{very}\, deep reason why Lie said
that $X$ belongs to the group $x ' = f ( x, a)$ is that he showed that
local transformation group actions are in one-to-one correspondence
with the purely linear vector spaces:
\[
\text{\rm Vect}_\K
\big(
X_1,X_2,\dots,X_r
\big),
\]
of infinitesimal transformations, which in fact also inherit a crucial
additional \emphasis{algebraic}\, structure directly from the group
multiplication law. Without anticipating too much, let us come to the
{\em purely analytic} 
way how Engel and Lie introduce the infinitesimal
transformations.

\section{Derivation of Fundamental Partial Differential Equations}
\label{derivation-fundamental}

So we have defined the concept of a purely local Lie transformation
group, insisting on the fact that composition and inversion are both
represented by some precise local analytic maps defined around the
identity. The royal road towards the famous theorems of Lie is to
differentiate these finite data, namely to \emphasis{infinitesimalize}.

We start with the group composition law~\thetag{ 
\ref{f-f-composition}} which we rewrite
as follows:
\[
\label{beginning-computation}
x''
=
f\big(f(x;a);b\big)
=
f(x;a\cdot b)
=:
f(x;c).
\]
Here, $c := a\cdot b$ depends on $a$ and $b$, but instead of $a$ and
$b$, following~\cite{ enlie1888-3}, we want to consider $a$ and $c$ to be
the independent parameters, namely we rewrite $b = a^{ - 1} \cdot c =:
b(a,c)$ so that the equations:
\[
f_i\big(f(x;a);\,b(a,c)\big)
\equiv
f_i(x;c)
\ \ \ \ \ \ \ \ \ \
{\scriptstyle{(i\,=\,1\,\cdots\,n)}}
\]
hold identically for all $x$, all $a$ and all $c$. Next, we
differentiate these identities with respect to $a_k$, denoting $f_i'
\equiv f_i( x'; b)$ and $x_j' \equiv f_j(x;a)$:
\[
\frac{\partial f_i'}{\partial x_1'}\,
\frac{\partial x_1'}{\partial a_k}
+\cdots+
\frac{\partial f_i'}{\partial x_n'}\,
\frac{\partial x_n'}{\partial a_k}
+
\frac{\partial f_i'}{\partial b_1}\,
\frac{\partial b_1}{\partial a_k}
+\cdots+
\frac{\partial f_i'}{\partial b_r}\,
\frac{\partial b_r}{\partial a_k}
\equiv
0
\ \ \ \ \ \ \ \ \ \
{\scriptstyle{(i\,=\,1\,\cdots\,n)}}.
\]
Here of course again, the argument of $f'$ is $\big( f( x, a); b (
a,c) \big)$, the argument of $x'$ is $( x; a)$ and the argument of $b$
is $(a, c)$. Thanks to $x'' ( x'; e) \equiv x'$, the matrix
$\frac{\partial f_i'}{\partial x_k'} \big( f(x;e); b(e, e) \big)$ is
the identity $I_{ n \times n}$. So using Cramer's rule\footnote{\,
---\,\,and possibly also, shrinking $\rho_2$ and $\sigma_2$ if
necessary, and also below, without further mention. }, for each fixed
$k$, we can solve the preceding $n$ linear equations with respect to
the $n$ unknowns $\frac{ \partial x_1' }{ \partial a_k}, \dots, \frac{
\partial x_n' }{ \partial a_k}$, getting expressions of the form:
\def\theequation{2}\begin{equation}
\label{partial-Xi}
\aligned
\frac{\partial x_\nu'}{\partial a_k}
(x;\,a)
&
=
\Xi_{1\nu}(x',b)\,\frac{\partial b_1}{\partial a_k}(a,c)
+\cdots+
\Xi_{r\nu}(x',b)\,\frac{\partial b_r}{\partial a_k}(a,c)
\\
&
\ \ \ \ \ \ \ \ \ \ \ \ \ \ \ \
{\scriptstyle{(\nu\,=\,1\,\cdots\,n\,;\,k\,=\,1\,\cdots\,r)}},
\endaligned
\end{equation}
with some analytic functions $\Xi_{ j\nu} ( x', b)$ that are
independent of $k$.

On the other hand, in order to substitute the $\frac{ \partial b_j}{
\partial a_k}$, we differentiate with respect to $a_k$ the identically
satisfied identities: \label{c-m-a-b}
\[
c_\mu
\equiv
\text{\bf m}_\mu\big(a,b(a,c)\big)
\ \ \ \ \ \ \ \ \ \
{\scriptstyle{(\mu\,=\,1\,\cdots\,r)}}.
\]
We therefore get:
\[
0
\equiv
\frac{\partial\text{\bf m}_\mu}{\partial a_k}
+
\sum_{\pi=1}^r\,
\frac{\partial\text{\bf m}_\mu}{\partial b_\pi}\,
\frac{\partial b_\pi}{\partial a_k}
\ \ \ \ \ \ \ \ \ \
{\scriptstyle{(\mu\,=\,1\,\cdots\,r)}}
\]
But since the matrix $\frac{ \partial \text{\bf m}_\mu}{ \partial
b_\pi}$ is the identity $I_{ r \times r}$ for $\big( a, b (a, c) \big)
\big\vert_{ (a, c) = (e,e)} = ( e, e)$, just because $\text{\bf m} (e,
b) \equiv b$, Cramer's rule again enables one to solve this system
with respect to the $r$ unknowns $\frac{ \partial b_\pi}{ \partial
a_k}$, getting expressions of the form:
\[
\frac{\partial b_\pi}{\partial a_k}
(a,c)
=
\Psi_{k\pi}(a,c),
\]
for certain functions $\Psi_{ k\pi}$, defined on a possibly smaller
parameter product space $\square_{ \sigma_2 }^r \times \square_{
\sigma_2 }^r$. Putting again $(a, c) = ( e, e)$ in the same system
before solving it, we get in fact (notice the
minus sign):
\[
\Big(
\frac{\partial b_\pi}{\partial a_k}(e,e)
\Big)_{1\leqslant k\leqslant r}^{1\leqslant\pi\leqslant r}
=
-I_{r\times r},
\]
so that $\Psi_{ k \pi} ( e, e) = - \delta_k^\pi$.

We can therefore now insert in~\thetag{ \ref{partial-Xi}} the gained
value $\Psi_{ k \pi}$ of $\frac{ \partial b_\pi}{ \partial a_k}$,
obtaining (and quoting~\cite{ enlie1888-3}, p.~29) the following
crucial partial differential equations:

\begin{svgraybox}
\def\theequation{2'}\begin{equation}
\label{partial-Xi-prime}
\frac{\partial x_\nu'}{\partial a_k}(x;a)
=
\sum_{\pi=1}^r\,\Psi_{k \pi}(a,b)\,\Xi_{\pi\nu}(x',b)
\ \ \ \ \ \ \ \ \ \
{\scriptstyle{(\nu\,=\,1\,\cdots\,n,\,\,\,k\,=\,1\,\cdots\,r)}}
\end{equation}
\emphasis{These equations are of utmost importance \deutsch{\"ausserst
wichtig}, as we will see later.}
\end{svgraybox}

Here, we have replaced $c$ by $c = c ( a, b) = a \cdot b$, whence $b
\big( a, c (a, b) \big) \equiv b$, and we have reconsidered $(a, b)$
as the independent variables.

\subsection{Restricting Considerations to a Single System of 
Parameters} 

Setting $b:=e$ in~\thetag{ \ref{partial-Xi-prime}} above, the partial
derivatives of the group transformation equations $x_i' = x_i' ( x;
a)$ with respect to the parameters $a_k$ at $(x; a)$:
\def\theequation{2''}\begin{equation}
\label{fundamental-differential-equations}
\boxed{
\boxed{\fboxrule=0.25pt
\frac{\partial x_i'}{\partial a_k}(x;a)
=
\sum_{j=1}^r\,\psi_{kj}(a)\,
\xi_{ji}\big(x'(x;a)\big)
}}
\ \ \ \ \ \ \ \ \ \ 
{\scriptstyle{(i\,=\,1\,\cdots\,n,\,\,\,k\,=\,1\,\cdots\,r)}},
\end{equation}
appear to be \emphasis{linear combinations}, with certain coefficients
$\psi_{ kj} (a) := \Psi_{ kj} ( a, e)$ \emphasis{depending only upon
$a$}, of the quantities $\xi_{ ji} ( x') := \Xi_{ji} ( x', b)
\big\vert_{ b = e}$. But we in fact already know these quantities.

\subsection{Comparing Different Frames 
\\
of Infinitesimal Transformations} 

Indeed,
setting $a = e$ above, thanks to $\psi_{ kj} ( e) = - \delta_k^j$, we
get immediately:
\[
\xi_{ki}(x)
=
\xi_{ki}\big(x'(x;e)\big)
=
-
\frac{\partial x_i'}{\partial a_k}(x;e),
\]
whence the $\xi_{ ki} ( x)$ just coincide with the coefficients of the
$r$ infinitesimal transformations already introduced on
p.~\pageref{infinitesimal-transformations} (with an overall opposite
sign) and written as derivations:
\def\theequation{3}\begin{equation}
\label{minus-xi}
\aligned
X_k^e\big\vert_x
&
=
\frac{\partial f_1}{\partial a_k}(x;e)\,
\frac{\partial}{\partial x_1}
+\cdots+
\frac{\partial f_n}{\partial a_k}(x;e)\,
\frac{\partial}{\partial x_n}
\\
&
=:
-\xi_{k1}(x)\,\frac{\partial}{\partial x_1}
-\cdots-
\xi_{kn}(x)\,\frac{\partial}{\partial x_n}
\ \ \ \ \ \ \ \ \ \ \ \ \ \ \ \ \ \ \ \
{\scriptstyle{(k\,=\,1\,\cdots\,r)}}.
\endaligned
\end{equation}
\emphasis{Now at last}, after having reproduced some rather blind
computations by which Lie expresses his brilliant synthetic
thoughts, the crucial geometric interpretation of the twice-boxed
partial differential equations can be unveiled.

\begin{center}
\input geom-fund-diff-eqs.pstex_t
\end{center}

\noindent
Instead of differentiating with respect to $a_k$ only at $(x; e)$, we
must in principle, for reasons of generality, do the same at any
$(x;a)$, which yields the vector fields:
\[
X_k^a\big\vert_{f_a(x)}
:=
\frac{\partial f_1}{\partial a_k}(x;a)\,
\frac{\partial}{\partial x_1}
+\cdots+
\frac{\partial f_n}{\partial a_k}(x;a)\,
\frac{\partial}{\partial x_n}
\ \ \ \ \ \ \ \ \ \ \ \ \ \ \ \ \ \ \ \
{\scriptstyle{(k\,=\,1\,\cdots\,r)}},
\]
and then the fundamental partial differential equations say that
such new $r$ infinitesimal transformations:
\[
\boxed{
X_k^a\big\vert_{f_a(x)}
=
-\psi_{k1}(a)\,X_1^e\big\vert_{f_a(x)}
-\cdots-
\psi_{kr}(a)\,X_r^e\big\vert_{f_a(x)}}
\]
are just linear combinations, \emphasis{with coefficients depending
only upon the group parameters}, of the $r$ infinitesimal
transformations $X_1^e, \dots, X_r^e$ computed at the special
parameter $e$, \emphasis{and considered at the $a$-pushed points $f_a
( x)$}.

Finally, since the matrix $\psi_{kj} (a)$ has an inverse that we will
denote, as in~\cite{ enlie1888-3}, by $\alpha_{jk} ( a)$ which is
analytic in a neighborhood of $e$, we can also write the fundamental
differential equations under the reciprocal form, useful in the
sequel:
\def\theequation{4}\begin{equation}
\label{xi-reciprocal}
\xi_{ji}\big(x'(x;a)\big)
=
\sum_{k=1}^r\,\alpha_{jk}(a)\,
\frac{\partial x_i'}{\partial a_k}(x;a).
\ \ \ \ \ \ \ \ \ \ \ \ \
{\scriptstyle{(i\,=\,1\,\cdots\,n\,;\,\,\,j\,=\,1\,\cdots\,r)}}.
\end{equation}

\section{Essentializing the Group Parameters}
\label{essentializing}

As we have seen in the previous chapter, the (needed) suppression of
illusory parameters in the transformation group equations $x_i' = f_i
( x; a) = \sum_{ \alpha \in\N^n }\, \mathcal{ U}_\alpha^i ( a) \,
x^\alpha$ might require to relocalize considerations to some
neighborhood of some generic $a^0$, \emphasis{which might possibly not
include the identity $e$}. Fortunately, this does not occur: the group
property ensures that the rank of the infinite coefficient map ${\sf
U}_\infty : a \longmapsto \big( \mathcal{ U}_\alpha^i ( a) \big)_{
\alpha \in \N^n}^{ 1\leqslant i \leqslant n}$ is \emphasis{constant}
around $e$.

\begin{proposition}
\label{proposition-essential-group}
For a finite continuous local Lie transformation group $x_i' = f_i (
x; a) = \sum_{ \alpha \in \N^n}\, \mathcal{ U}_\alpha^i ( a) \,
x^\alpha$ expanded in power series with respect to $x$, the following
four conditions are equivalent:

\begin{itemize}

\smallskip\item[\!{\bf (i)}]
The parameters $(a_1, \dots, a_r)$ are \emphasis{not} essential, namely the
generic rank of $a \longmapsto \big( \mathcal{ U}_\alpha^i ( a)
\big)_{ \alpha \in \N^n}^{ 1 \leqslant i \leqslant n}$ is strictly
less than $r$.

\smallskip\item[\!\!{\bf (ii)}]
The rank at all $a$ near the identity $e = 0$ (and not only the
generic rank) of $a \longmapsto \big( \mathcal{ U}_\alpha^i ( a)
\big)_{ \alpha \in \N^n}^{ 1 \leqslant i \leqslant n}$ is
strictly less than $r$.

\smallskip\item[\!\!\!\!{\bf (iii)}]
There exists a vector field $\mathcal{ T} = \sum_{ k=1}^r \, \tau_k (
a) \, \frac{ \partial }{\partial a_k}$ having analytic coefficients
$\tau_k ( a)$ which vanishes nowhere near $e$ and annihilates all
coefficient functions: $0 \equiv \mathcal{ T} \mathcal{ U}_\alpha^i$
for all $i$ and all $\alpha$.

\smallskip\item[\!\!\!{\bf (iv)}]
There exist constants $e_1, \dots, e_r \in \K$ not all zero such that,
if
\[
X_k
=
\frac{\partial f_1}{\partial a_k}(x;e)\,
\frac{\partial}{\partial x_1}
+\cdots+
\frac{\partial f_n}{\partial a_k}(x;e)\,
\frac{\partial}{\partial x_n}
\ \ \ \ \ \ \ \ \ \
{\scriptstyle{(k\,=\,1\,\cdots\,n)}}
\]
denote the $r$ infinitesimal transformations associated
to the group, then:
\[
0
\equiv
e_1\,X_1
+\cdots+
e_r\,X_r.
\] 
\end{itemize}
\end{proposition}

\begin{proof}\smartqed
Suppose that {\bf (iv)} holds, namely via the 
expressions~\thetag{ \ref{minus-xi}}:
\[
0
\equiv
-e_1\,\xi_{1i}(x)
-\cdots
-e_r\,\xi_{ri}(x)
\ \ \ \ \ \ \ \ \ \ \ \ \
{\scriptstyle{(i\,=\,1\,\cdots\,n)}}.
\]
By applying the linear combination $\sum_{ j = 1}^r\, e_j ( \cdot)$ to
the partial differential equations~\thetag{ 
\ref{xi-reciprocal}}, the left member
thus becomes zero, and we get $n$
equations:
\[
0
\equiv
\sum_{j=1}^r\,\sum_{k=1}^r\,
e_j\,
\theta_{jk}(a)\,
\frac{\partial x_i'}{\partial a_k}(x;a)
\ \ \ \ \ \ \ \ \ \ \ \ \
{\scriptstyle{(i\,=\,1\,\cdots\,n)}}
\]
which just say that the \emphasis{not identically zero} vector field:
\[
\mathcal{T}
:=
\sum_{k=1}^r\,
\bigg(
\sum_{j=1}^r\,e_j\,\theta_{jk}(a)
\bigg)\,\frac{\partial}{\partial a_k}
\]
satisfies $0 \equiv \mathcal{ T}\, f_1 \equiv \cdots \equiv \mathcal{
T}\, f_n$, or equivalently $0 \equiv 
\mathcal{ T} \mathcal{ U}_\alpha^i$ for all
$i$ and all $\alpha$. Since the matrix $\theta (a)$ equals $- I_{ r
\times r}$ at the identity $e$, the vectors $\mathcal{ T} \big
\vert_a$ are \emphasis{nonzero} for all $a$ in a neighborhood of $a$. In
conclusion, {\bf (iv)} implies {\bf (iii)}, and then straightforwardly
{\bf (iii)} $\Rightarrow$ {\bf (ii)} $\Rightarrow$ {\bf (i)}.

\smallskip

It remains to establish the reverse implications: {\bf (i)}
$\Rightarrow$ {\bf (ii)} $\Rightarrow$ {\bf (iii)} $\Rightarrow$ {\bf
(iv)}. As a key fact, the generic rank in {\bf (i)} happens to be
constant.

\begin{lemma}
There exist two $\infty \times \infty$ matrices $\big( {\sf
Coeff}_\alpha^\beta ( a)\big)_{ \alpha\in \N^n}^{ \beta \in \N^n}$ and
$\big( \overline{ \sf Coeff}{\:\!}_\alpha^\beta ( a)\big)_{ \alpha\in
\N^n}^{ \beta \in \N^n}$ of analytic functions of the parameters $a_1,
\dots, a_r$ which enjoy appropriate Cauchy convergence
estimates\footnote{\, ---\,\,that are automatically satisfied and
hence can be passed over in silence\,\,---}, the second matrix being
the inverse of the first, such that:
\[
\aligned
\mathcal{U}_\alpha^i
\big(\text{\bf m}(a,b)\big)
\equiv
\sum_{\beta\in\N^n}\,
{\sf Coeff}_\alpha^\beta(a)\,
\mathcal{U}_\beta^i(b)
\ \ \ \ \ \ \ \ \ \ \ \ \
{\scriptstyle{(i\,=\,1\,\cdots\,n;\,\,\,\alpha\,\in\,\N^n)}},
\endaligned
\]
and inversely:
\[
\mathcal{U}_\alpha^i(b)
\equiv
\sum_{\beta\in\N^n}\,
\overline{\sf Coeff}{\:\!}_\alpha^\beta(a)\,
\mathcal{U}_\beta^i
\big(\text{\bf m}(a,b)\big)
\ \ \ \ \ \ \ \ \ \ \ \ \
{\scriptstyle{(i\,=\,1\,\cdots\,n;\,\,\,\alpha\,\in\,\N^n)}}.
\]
As a consequence, the rank of the infinite coefficient mapping $a
\longmapsto \big( \mathcal{ U}_\alpha^i ( a) \big)_{ \alpha \in
\N^n}^{ 1\leqslant i \leqslant n}$ is constant in a neighborhood of
$e$.
\end{lemma}

\begin{proof}\smartqed
By definition, the composition of $x' = \sum_{ \alpha \in \N^n} \,
\mathcal{ U}_\alpha (a) \, x^\alpha$ and of $x_i'' = \sum_{ \beta \in
\N^n}\, \mathcal{ U}_\beta^i ( b) \, (x')^\beta$ yields $x_i'' = f_i
\big( x; \, \text{\bf m} ( a, b)\big)$ and we must therefore compute
the composed Taylor series expansion in question, which we place in
the right-hand side:
\[
\sum_{\alpha\in\N^n}\,\mathcal{U}_\alpha^i
\big(\text{\bf m}(a,b)\big)\,x^\alpha
\equiv
\sum_{\beta\in\N^n}\,
\mathcal{U}_\beta^i(b)
\bigg(
\sum_{\alpha\in\N^n}\,\mathcal{U}_\alpha(a)\,x^\alpha
\bigg)^\beta.
\]
To this aim, we use the general formula for the expansion of a
$k$-th power:
\[
\bigg(
\sum_{\alpha\in\N^n}\,{\sf c}_\alpha\,x^\alpha
\bigg)^k
=
\sum_{\alpha\in\N^n}\,x^\alpha
\Big\{
\sum_{\alpha_1+\cdots+\alpha_k\,=\,\alpha}\,
{\sf c}_{\alpha_1}\cdots{\sf c}_{\alpha_k}
\Big\}
\]

\noindent
of a \emphasis{scalar} power series. Thus splitting ${x'}^\beta =
(x_1')^{\beta_1} \cdots (x_n')^{\beta_n}$ in scalar powers, we may
start to compute the wanted expansion:
\[
\aligned
\text{\sf Composition}
&
=
\sum_{(\beta_1,\dots,\beta_n)}\,
\mathcal{U}_\beta^i(b)\,
\bigg(
\sum_{\alpha^1\in\N^n}\,
\mathcal{U}_{\alpha^1}^1(a)\,x^{\alpha^1}
\bigg)^{\beta_1}
\cdots
\bigg(
\sum_{\alpha^n\in\N^n}\,
\mathcal{U}_{\alpha^n}^n(a)\,x^{\alpha^n}
\bigg)^{\beta_n}
\\
&
=
\sum_{(\beta_1,\dots,\beta_n)}\,
\mathcal{U}_\beta^i(b)\,
\bigg[
\sum_{\alpha^1\in\N^n}\,x^{\alpha^1}\,
\Big\{
\sum_{\alpha_1^1+\cdots+\alpha_{\beta_1}^1=\,\alpha^1}\,
\mathcal{U}_{\alpha_1^1}^1(a)\cdots
\mathcal{U}_{\alpha_{\beta_1}^1}^1(a)
\Big\}
\bigg]\cdots
\\
&
\ \ \ \ \ \ \ \ \ \ \ \ \ \ \
\cdots\bigg[
\sum_{\alpha^n\in\N^n}\,x^{\alpha^n}\,
\Big\{
\sum_{\alpha_1^n+\cdots+\alpha_{\beta_n}^n=\,\alpha^n}\,
\mathcal{U}_{\alpha_1^n}^n(a)\cdots\mathcal{U}_{\alpha_{\beta_n}^n}^n(a)
\Big\}
\bigg].
\endaligned
\]
To finish off, if we apply the expansion of a product of $n$ 
power series:
\[
\Big(
{\textstyle{\sum}_{\alpha^1}}\,{\sf c}_{\alpha^1}^1\,x^{\alpha^1}
\Big)
\cdots
\Big(
{\textstyle{\sum}_{\alpha^n}}\,{\sf c}_{\alpha^n}^n\,x^{\alpha^n}
\Big)
=
{\textstyle{\sum}}\,x^\alpha\,
\Big\{
{\textstyle{\sum}}_{\alpha^1+\cdots+\alpha^n\,=\,\alpha}\,\,\,
{\sf c}_{\alpha^1}^1\cdots{\sf c}_{\alpha^n}^n
\Big\},
\]
we then obtain straightforwardly:
\[
\aligned
&
\sum_{\alpha\in\N^n}\,\mathcal{U}_\alpha^i
\big(\text{\bf m}(a,b)\big)\,x^\alpha
=
\sum_{\alpha\in\N^n}\,x^\alpha\,
\bigg\{
\sum_{(\beta_1,\dots,\beta_n)}\,
\mathcal{U}_\beta^i(b)\,
\sum_{\alpha^1+\cdots+\alpha^n\,=\,\alpha}\,
\\
&
\Big(
\sum_{\alpha_1^1+\cdots+\alpha_{\beta_1}^1\,=\,\alpha^1}\,
\mathcal{U}_{\alpha_1^1}^1(a)\cdots\mathcal{U}_{\alpha_{\beta_1}^1}^1(a)
\Big)\cdots
\Big(
\sum_{\alpha_1^n+\cdots+\alpha_{\beta_n}^n\,=\,\alpha^n}\,
\mathcal{U}_{\alpha_1^n}^n(a)\cdots\mathcal{U}_{\alpha_{\beta_n}^n}^n(a)
\Big)
\bigg\}.
\endaligned
\]
By identifying the coefficients of $x^\alpha$ in both sides and by
abbreviating just as ${\sf Coeff}_\alpha^\beta ( a)$ the
$\sum_{\alpha^1 + \cdots + \alpha^n \, = \, \alpha}$ of the products
of the $n$ long sums appearing in the second line, we have thus gained
the first family of equations. The second one is obtained quite
similarly by just expanding the identities:
\[
f_i(x;\,b)
\equiv
f_i\big(f(x;\,\text{\bf i}(a));\,\text{\bf m}(a,b)\big)
\ \ \ \ \ \ \ \ \ \ \ \ \
{\scriptstyle{(i\,=\,1\,\cdots\,n)}}.
\] 
The two gained families of equations must clearly be inverses of each
other.

To show that the rank at $e$ of $a \mapsto {\sf U }_\infty (a)$ is the
same as it is at any $a$ near $e$, we apply $\frac{ \partial }{
\partial b_k} \big \vert_{ b = e}$ to the first, and also to the
second, family of equations, which gives:
\[
\aligned
\sum_{l=1}^r\,\frac{\partial\mathcal{U}_\alpha^i}{\partial a_l}(a)\,
\frac{\partial\text{\bf m}_l}{\partial b_k}(a,e)
&
\equiv
\sum_{\beta\in\N^n}\,
\text{\sf Coeff}{\:\!}_\alpha^\beta(a)
\,
\frac{\partial\mathcal{U}_\beta^i}{\partial b_k}(e)
\\
\frac{\partial\mathcal{U}_\alpha^i}{\partial b_k}(e)
&
\equiv
\sum_{\beta\in\N^n}\,
\overline{\sf Coeff}{\:\!}_\alpha^\beta(a)\,
\bigg(
\sum_{l=1}^r\,
\frac{\partial\mathcal{U}_\beta^i(a)}{\partial a_l}\,
\frac{\partial\text{\bf m}_l}{\partial b_k}(a,e)
\bigg)
\\
&
{\scriptstyle{(k\,=\,1\,\cdots\,r;\,\,\,i\,=\,1\,\cdots\,n;\,\,\,
\alpha\,\in\,\N^n)}}.
\endaligned
\]
Of course, the appearing $r \times r$ matrix:
\[
{\sf M}(a)
:=
\Big(
\frac{\partial\text{\bf m}_l}{\partial b_k}(a,e)
\Big)_{1\leqslant k\leqslant r}^{1\leqslant l\leqslant r}
\]
has nonzero determinant at every $a$, because left translations $b \mapsto
\text{\bf m} ( a, b)$ of the group are diffeomorphisms. Thus, if for
each fixed $i$ we denote by ${\rm Jac} \, {\sf U}_\infty^i (a )$ the
$r \times \infty$ Jacobian matrix $\big( \frac{\partial \mathcal{
U}_\alpha^i }{\partial a_k} ( a) \big)_{ 1 \leqslant k \leqslant r}^{
\alpha \in \N^n}$, the previous two families of identities when
written in matrix form:
\[
\aligned
{\sf M}(a)\,{\rm Jac}\,{\sf U}_\infty^i(a)
&
\equiv
{\sf Coeff}(a)\,{\rm Jac}\,{\sf U}_\infty^i(a)
\\
{\rm Jac}\,{\sf U}_\infty^i(e)
&
\equiv
\overline{\sf Coeff}(a)\,{\sf M}(a)\,{\rm Jac}\,{\sf U}_\infty^i(a)
\endaligned
\]
show well that the rank of ${\rm Jac}_\infty^i ( a)$ must be equal to
the rank of ${\rm Jac}_\infty^i ( e)$. This completes the proof of
the auxiliary lemma.
\qed\end{proof}

Thus, {\bf (i)} of the proposition implies {\bf (ii)} and moreover,
the last part of the theorem on p.~\pageref{Theorem-essential} yields
\emphasis{annihilating}
analytic vector fields $\mathcal{ T}_1, \dots, \mathcal{ T}_{ r -
\rho_\infty}$ on the parameter space with $\dim {\rm Vect} \big(
\mathcal{ T}_1 \big\vert_a, \dots, \mathcal{ T}_{\rho_\infty}
\big\vert_a \big) = r - \rho_\infty$ constant for all $a$ near $e$.
So {\bf (i)} $\Rightarrow$ {\bf (ii)} $\Rightarrow$ {\bf (iii)}.
Finally, assuming {\bf (iii)}, namely:
\[
0
\equiv
\mathcal{T}\,f_i(x;a)
=
\sum_{k=1}^r\,\tau_k(a)\,
\frac{\partial x_i'}{\partial a_k}(x;a)
\ \ \ \ \ \ \ \ \ \ \ \ \
{\scriptstyle{(i\,=\,1\,\cdots\,n)}}
\]
with $\mathcal{ T} \big\vert_e \neq 0$, and replacing $\frac{ \partial
x_i'}{\partial a_k}$ by its value~\thetag{
\ref{fundamental-differential-equations}} given by the fundamental
differential equations, we get:
\[
0
\equiv
\sum_{j=1}^r\,
\sum_{k=1}^r\,\tau_k(a)\,\psi_{kj}(a)\,\xi_{ji}
\big(x'(x;a)\big)
\ \ \ \ \ \ \ \ \ \ \ \ \
{\scriptstyle{(i\,=\,1\,\cdots\,n)}}.
\]
Setting $a$ to be the identity element in these equations and
introducing the $r$ constants (remind here that $\psi_{ kj} ( e) = -
\delta_k^j$):
\[
c_j
:=
{\textstyle{\sum_{k=1}^r}}\,
\tau_k(e)\,\psi_{kj}(e)
=
-\tau_j(e)
\ \ \ \ \ \ \ \ \ \ \ \ \
{\scriptstyle{(j\,=\,1\,\cdots\,r)}}
\]
that are not all zero, since $\mathcal{ T}\big\vert_e \neq 0$, we get
equations:
\[
0
\equiv
c_1\,\xi_{1i}(x)
+\cdots+
c_r\,\xi_{ri}(x)
\ \ \ \ \ \ \ \ \ \ \ \ \
{\scriptstyle{(i\,=\,1\,\cdots\,n)}}
\]
which, according to~\thetag{ \ref{minus-xi}}, 
are the coordinatewise expression of
{\bf (iv)}. This completes the proof of the proposition.
\qed\end{proof}

\section{First Fundamental Theorem}

Thanks to this proposition and to the corollary on
p.~\pageref{Theorem-essential}, \emphasis{no relocalization is
necessary to get rid of superfluous parameters in finite continuous
transformation groups $x' = f( x; a_1, \dots, a_r)$}. Thus, without
loss of generality parameters can (and surely will) \emphasis{always}
be assumed to be essential. We can now translate Theorem~3 on
pp.~33--34 of~\cite{ enlie1888-3} which summarizes all the preceding
considerations, and we add some technical precisions if this theorem
is to be interpreted as a {\em local} statement (though Lie has
something different in mind, cf. the next section).  As above, we
assume implicitly for simplicity that the group $x' = f ( x; \, a)$ is
local and contains the identity transformation $x' = f ( x; \, e) =
x$, but the theorem is in fact valid with less restrictive
assumptions, \voir~Sect.~\ref{substituting-axiom} and
especially Proposition~\ref{economic-differential}.

\begin{svgraybox}
\vspace{-0.4cm}
\def\thetheorem{3}\begin{theorem}
\label{Theorem-3-S-33} 
If the $n$ equations:
\[
x_i'
=
f_i(x_1,\dots,x_n,\,a_1,\dots,a_r)
\ \ \ \ \ \ \ \ \ \
{\scriptstyle{(i\,=\,1\,\cdots\,n)}}
\]
represent a finite continuous group, whose parameters $a_1, \dots,
a_r$ are all essential, then $x_1', \dots, x_n'$, considered as
functions of $a_1, \dots, a_r, \, x_1, \dots, x_n$ satisfy certain
differential equations of the form:
\def\theequation{5}\begin{equation}
\frac{\partial x_i'}{\partial a_k}
=
\sum_{j=1}^r\,\psi_{kj}(a_1,\dots,a_r)\,
\xi_{ji}(x_1',\dots,x_n')
\ \ \ \ \ \ \ \ \
{\scriptstyle{(i\,=\,1\,\cdots\,n;\,\,\,k\,=\,1\,\cdots\,r)}},
\end{equation}
which can also be written as:
\def\theequation{6}\begin{equation}
\xi_{ji}(x_1',\dots,x_n')
=
\sum_{k=1}^r\,\alpha_{jk}(a_1,\dots,a_r)\,
\frac{\partial x_i'}{\partial a_k}
\ \ \ \ \ \ \ \ \
{\scriptstyle{(i\,=\,1\,\cdots\,n;\,\,\,j\,=\,1\,\cdots\,r)}}.
\end{equation}
Here, neither the determinant of the $\psi_{ kj} (a)$, nor the one of
the $\alpha_{ jk} (a)$ vanishes identically\footnotemark; in
addition, it is impossible to indicate $r$ quantities $e_1, \dots,
e_r$ independent of $x_1', \dots, x_n'$ [constants] and not all
vanishing such that the $n$ expressions:
\[
e_1\,\xi_{1i}(x')
+\cdots+
e_r\,\xi_{ri}(x')
\ \ \ \ \ \ \ \ \ \
{\scriptstyle{(i\,=\,1\,\cdots\,n)}}
\]
vanish simultaneously.
\end{theorem}
\end{svgraybox}

\footnotetext{\,
The functions $\xi_{ ji} (x)$ are, up to an overall minus sign, just
the coefficients of the $r$ infinitesimal transformations obtained by
differentiation with respect to the parameters at the identity:
$\xi_{ji}(x) = -\frac{ \partial f_i}{ \partial a_j}( x; \,e)$, $i = 1,
\dots, n$, $j = 1, \dots, n$. Furthermore, with the purely local
assumptions we made above, we in fact have $\psi_{ kj} ( e) = -
\delta_k^j$, \label{psi-k-j-e}
so the mentioned determinants do not vanish for all $a$
in a neighborhood of $e$.
} 
 
Furthermore, the latter property is clearly equivalent to
the nonexistence of constants $e_1, \dots, e_r$ not all zero
such that:
\[
0
\equiv
e_1\,X_1+\cdots+e_r\,X_r, 
\]
that is to say: the $r$ infinitesimal transformations $X_1, \dots,
X_r$ are {\em linearly} independent. This property was shown to derive
from essentiality of parameters, and Theorem~8
p.~\pageref{Theorem-8-S-65} below will establish a satisfactory
converse.

Thus in the above formulation, Engel and
Lie do implicitly introduce the $r$ infinitesimal transformations:
\[
X_k
:=
\sum_{i=1}^n\,\xi_{ki}(x_1,\dots,x_n)\,
\frac{\partial}{\partial x_i}
\ \ \ \ \ \ \ \ \ \ \ \ \
{\scriptstyle{(k\,=\,1\,\cdots\,r)}}
\]
of an $r$-term continuous group \emphasis{not} as partial derivatives
at the identity element with respect to the parameters (which would be
the intuitively clearest way):
\[
X_k^e
:=
\sum_{i=1}^n\,
\frac{\partial f_i}{\partial a_k}(x;\,e)\,
\frac{\partial}{\partial x_i}
\equiv
-\,
X_k,
\] 
but rather indirectly as having coefficients $\xi_{ ki} (x)$ stemming
from the fundamental differential equations~\thetag{ 5}; 
in fact, we already
saw in~\ref{minus-xi} that both definitions agree modulo
sign. Slightly Later in the treatise, the introduction is 
done explicitly.

\begin{svgraybox}
\label{Satz-S-67} 
\noindent{\bf Proposition.} (\cite{ enlie1888-3}, p.~67)
{\em Associated to every $r$-term group $x_i' = f_i (x_1, \dots,x_n,$
$a_1, \dots, a_r)$, there are $r$ infinitesimal transformations:
\[
X_k(f)
=
\sum_{i=1}^n\,\xi_{ki}(x_1,\dots,x_n)\,
\frac{\partial f}{\partial x_i}
\ \ \ \ \ \ \ \ \ \ \ \ \
{\scriptstyle{(k\,=\,1\,\cdots\,r)}}
\]
which stand in such a relationship that equations of the form:
\[
\frac{\partial x_i'}{\partial a_k}
=
\sum_{j=1}^r\,\psi_{kj}(a_1,\dots,a_r)\,
\xi_{ji}(x_1',\dots,x_n')
\]
hold true, that can be resolved with respect to the 
$\xi_{ ji}$:}
\[
\xi_{ji}(x_1',\dots,x_n')
=
\sum_{k=1}^r\,\alpha_{jk}(a_1,\dots,a_r)\,
\frac{\partial x_i'}{\partial a_k}.
\]
\vspace{-0.4cm}
\end{svgraybox}

The gist of Lie's theory is to show that the datum of an $r$-term
continuous (local) transformation group is equivalent the datum of
$r$ infinitesimal transformations $X_1, \dots, X_r$ associated this
way. The next chapters are devoted to expose this (well known)
one-to-one correspondence between local Lie groups and Lie algebras of
local analytic vector fields (Chap.~\ref{kapitel-9}), by following, by
summarizing and by adapting the original presentation, but without
succumbing to the temptation of overformalizing some alternative
coordinate-free reasonings with the help of some available
contemporary views, because this would certainly impoverish the depth
of Lie's original thought.

\section{
Fundamental Differential Equations
\\
for the Inverse Transformations
}
\label{fundamental-inverse}
\sectionmark{Fundamental Differential Equations
for the Inverse Transformations}

According to the fundamental Theorem~3 just above
p.~\pageref{Theorem-3-S-33}, a general $r$-term continuous
transformation group $x_i ' = f_i ( x; a_1, \dots, a_r)$ satisfies
partial differential equations: $\partial f_i / \partial a_k = \sum_{
j = 1}^r\, \psi_{ kj} ( a) \,
\xi_{ ji} ( f_1, \dots, f_n)$ that are used everywhere in the basic
Lie theory. For the study of the adjoint group in
Chap.~\ref{kapitel-16} below, we must also know how to write precisely
the fundamental differential equations that are satisfied by the group
of \emphasis{inverse} transformations:
\[
x_i
=
f_i\big(x';\,\text{\bf i}(a)\big)
\ \ \ \ \ \ \ \ \ \ \ \ \
{\scriptstyle{(i\,=\,1\,\cdots\,n)}},
\]
and this is easy. Following an already known path, we must indeed
begin by differentiating these equations with respect to the
parameters $a_k$:
\[
\frac{\partial x_i}{\partial a_k}
=
\sum_{l=1}^r\,\frac{\partial f_i}{\partial a_l}
\big(x';\,\text{\bf i}(a)\big)\,
\frac{\partial\text{\bf i}_l}{\partial a_k}(a)
\ \ \ \ \ \ \ \ \ \ \ \ \
{\scriptstyle{(i\,=\,1\,\cdots\,n\,;\,\,\,
k\,=\,1\,\cdots\,r)}}.
\]
Naturally, we replace here the
$\partial f_i / \partial a_l$ by
their values $\sum_{ j=1 }^r\, \psi_{ lj} \, \xi_{ ji }$ given by the
fundamental differential equations
of Theorem~3, and we obtain a double sum:
\[
\aligned
\frac{\partial x_i}{\partial a_k}
&
=
\sum_{l=1}^r\,\sum_{j=1}^r\,
\psi_{lj}\big(\text{\bf i}(a)\big)\,
\xi_{ji}\big(
f(x';\text{\bf i}(a)
\big)\,
\frac{\partial\text{\bf i}_l}{\partial a_k}(a)
\\
&
=:
\sum_{j=1}^r\,\vartheta_{kj}(a)\,\xi_{ji}(x)
\ \ \ \ \ \ \ \ \ \ \ \ \
{\scriptstyle{(i\,=\,1\,\cdots\,n;\,\,\,k\,=\,1\,\cdots\,r)}},
\endaligned
\]
which we contract to a single sum by simply introducing the 
following new $r \times r$ auxiliary matrix of parameter functions:
\[
\vartheta_{kj}(a)
:=
\sum_{l=1}^r\,\psi_{lj}\big(\text{\bf i}(a)\big)\,
\frac{\partial\text{\bf i}_l}{\partial a_k}(a)
\ \ \ \ \ \ \ \ \ \ \ \ \
{\scriptstyle{(k,\,j\,=\,1\,\cdots\,r)}},
\]
whose precise expression in terms of
${\bf i} (a)$
will not matter anymore. It now remains to
check that this matrix $\big( \vartheta_{ kj} (a) \big)_{ 1 \leqslant
k \leqslant r}^{ 1 \leqslant j \leqslant r}$ is invertible for all $a$
in a neighborhood of the identity element $e = (e_1, \dots, e_r)$.
In fact, we claim that:
\[
\vartheta_{kj}(e)
=
\delta_k^j,
\]
which will clearly assure the invertibility in question. Firstly, we
remember from Theorem~3 on p.~\pageref{Theorem-3-S-33} that $\psi_{
lj} ( e) = - \delta_l^j$. Thus secondly, it remains now only to check
that $\frac{ \partial \text{\bf i}_l}{\partial a_k} (e) = -
\delta_k^l$.

To check this, we differentiate with respect to $a_k$ the 
identities: $e_j \equiv \text{ \bf m }_j \big(a, \text{ \bf i}
(a)\big)$, $j = 1, \dots, r$, 
which hold by definition, and we get:
\[
0
\equiv
\frac{\partial \text{\bf m}_j}{\partial a_k}(e,e)
+
\sum_{l=1}^r\,
\frac{\partial \text{\bf m}_j}{\partial b_l}(e,e)
\frac{\partial\text{\bf i}_l}{\partial a_k}(e)
\ \ \ \ \ \ \ \ \ \ \ \ \
{\scriptstyle{(j\,=\,1\,\cdots\,r)}}.
\]
From another side, by differentiating the two
families of $r$ identities $a_j \equiv
\text{ \bf m}_j ( a, e)$ and $b_j \equiv \text{ \bf m}_j ( e, b)$ with
respect to $a_k$ and with respect to $b_l$, we immediately get
two expressions:
\[
\frac{\partial\text{\bf m}_j}{\partial a_k}(e,e) 
=
\delta_k^j
\ \ \ \ \ \ \ \ \ \ \ \ \
\text{\rm and}
\ \ \ \ \ \ \ \ \ \ \ \ \
\frac{\partial\text{\bf m}_j}{\partial b_l}(e,e)
=
\delta_l^j
\]
which, when inserted just above, yield the announced $\frac{
\partial
\text{ \bf i}_l }{ \partial a_k } (e) = - \delta_k^l$. Sometimes, we
will write $g ( x; \, a)$ instead of $f \big( x; \, \text{\bf i} ( a)
\big)$. As a result:

\begin{lemma}
The finite continuous transformation group $x_i ' = f_i ( x; a)$
and its inverse transformations $x_i = g_i ( x'; a) := f_i \big( x'; \,
\text{ \bf i} ( a) \big)$ both satisfy fundamental partial
differential equations of the form:
\[
\left\{
\aligned
\frac{\partial x_i'}{\partial a_k}(x;\,a)
&
=
\sum_{j=1}^r\,\psi_{kj}(a)\,\xi_{ji}
\big(x'(x;a)\big)
\ \ \ \ \ \ \ \ \ \ \ \ \
{\scriptstyle{(i\,=\,1\,\cdots\,n\,;\,\,\,
k\,=\,1\,\cdots\,r)}},
\\
\frac{\partial x_i}{\partial a_k}(x';\,a)
&
=
\sum_{j=1}^r\,\vartheta_{kj}(a)\,\xi_{ji}
\big(x(x';a)\big)
\ \ \ \ \ \ \ \ \ \ \ \ \
{\scriptstyle{(i\,=\,1\,\cdots\,n\,;\,\,\,
k\,=\,1\,\cdots\,r)}},
\endaligned\right.
\]
where $\psi$ and $\vartheta$ are some two $r \times r$ matrices of
analytic functions with $- \psi_{ kj} ( e) = \vartheta_{ kj} ( e) =
\delta_k^j$, and where the functions $\xi_{ ji}$ appearing in 
\emphasis{both} systems of equations:
\[
\xi_{ji}(x)
:=
-\frac{\partial f_i}{\partial x_j}(x;e)
\ \ \ \ \ \ \ \ \ \ \ \ \
{\scriptstyle{(i\,=\,1\,\cdots\,n;\,\,\,j\,=\,1\,\cdots\,r)}}
\] 
are, up to an overall minus sign, just the coefficients of the $r$
infinitesimal transformations 
\[
X_1^e 
=
\frac{\partial f}{\partial a_1}(x;\,e), 
\dots\dots, 
X_r^e
=
\frac{\partial f}{\partial a_r}(x;\,e) 
\]
obtained by differentiating the finite equations with respect to the
parameters at the identity element.
\end{lemma}

\begin{svgraybox}\vspace{-0.4cm}
\def\thetheorem{4}\begin{theorem}
If, in the equations:
\[
x_i'
=
f_i(x_1,\dots,x_n,\,a_1,\dots,a_r)
\ \ \ \ \ \ \ \ \ \ \ \ \ {\scriptstyle{(i\,=\,1\,\cdots\,n)}}
\]
of a group with the $r$ essential parameters $a_1, \dots, a_r$,
one considers the $x_i$ as functions of $a_1, \dots, a_r$
and of $x_1', \dots, x_n'$, then there exist differential 
equations of the form:
\[
\frac{\partial x_i}{\partial a_k}
=
\sum_{j=1}^r\,\vartheta_{kj}(a_1,\dots,a_r)\,
\xi_{ji}(x_1,\dots,x_n)
\ \ \ \ \ \ \ \ \ \ \ \ \ 
{\scriptstyle{(i\,=\,1\,\cdots\,n\,;\,\,\,
k\,=\,1\,\cdots\,r)}}.
\]
\end{theorem}\vspace{-0.4cm}
\end{svgraybox}

\section{Transfer of Individual Infinitesimal Transformations 
\\
by the Group}
\label{transfer-infinitesimal-by-group}
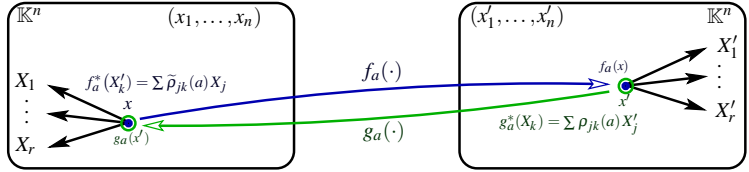
\sectionmark{Transfer of Infinitesimal Transformations by the Group}

With $x = g( x'; \, a)$
denoting the inverse of
$x' = f ( x; \, a)$, 
we now differentiate with respect to $a_k$ the identically satisfied
equations:
\[
x_i'
\equiv
f_i\big(g(x';\,a);\,a)
\ \ \ \ \ \ \ \ \ \ \ \ \
{\scriptstyle{(i\,=\,1\,\cdots\,n)}},
\]
which just say that an arbitrary transformation of the group followed by
its inverse gives again the identity transformation, 
and we immediately get:
\[
0
\equiv
\sum_{\nu=1}^n\,
\frac{\partial f_i}{\partial x_\nu}\,
\frac{\partial g_\nu}{\partial a_k}
+
\frac{\partial f_i}{\partial a_k}
\ \ \ \ \ \ \ \ \ \ \ \ \
{\scriptstyle{(i\,=\,1\,\cdots\,n;\,\,\,k\,=\,1\,\cdots\,r)}}.
\]
Thanks to the above two systems of partial differential equations, we
may then replace $\partial g_\nu / \partial a_k$ by its value from the
second equation of the lemma above, and also $\partial f_i / \partial
a_k$ by its value from the first equation in the same lemma:
\def\theequation{7}\begin{equation}
\aligned
0
&
\equiv
\sum_{\nu=1}^n\,
\Big\{
\sum_{j=1}^r\,\vartheta_{kj}(a)\,\xi_{j\nu}(g)
\Big\}\,
\frac{\partial f_i}{\partial x_\nu}
+
\sum_{j=1}^r\,\psi_{kj}(a)\,\xi_{ji}(f)
\\
&
\ \ \ \ \ \ \ \ \ \ \ \ \ \ \ \ \ \ \ \ \ \ \ \ \ \ \ \ \ \
{\scriptstyle{(i\,=\,1\,\cdots\,n;\,\,\,k\,=\,1\,\cdots\,r)}}.
\endaligned
\end{equation}
In order to bring these equations to a more symmetric form,
following~\cite{ enlie1888-3} pp.~44--45, we fix $k$ and we multiply,
for $i = 1$ to $n$, the $i$-th equation by $\frac{ \partial }{\partial
x_i'}$, we apply the summation $\sum_{ i = 1}^n$, we use the fact
that, through the diffeomorphism $x \mapsto f_a ( x) = x'$, the
coordinate vector fields transform as:
\[
\frac{\partial}{\partial x_\nu} 
= 
\sum_{i=1}^n\, 
\frac{\partial f_i}{\partial x_\nu}\,
\frac{ \partial }{\partial x_i'}, 
\]
which just means in contemporary
notation that:
\[
{\textstyle{
(f_a)_*
\big(
\frac{\partial}{\partial x_\nu}
\big)
=\sum_{i=1}^n
\frac{\partial f_i}{\partial x_\nu}\,
\frac{\partial}{\partial x_i'}}}
\ \ \ \ \ \ \ \ \ \ \ \ \
{\scriptstyle{(\nu\,=\,1\,\cdots\,n)}},
\]
and we obtain, thanks to this observation, completely symmetric
equations:
\[
0
\equiv
\sum_{j=1}^n\,
\vartheta_{kj}(a)\,
\sum_{\nu=1}^n\,\xi_{j\nu}(x)\,
\frac{\partial}{\partial x_\nu}
+
\sum_{j=1}^r\,\psi_{kj}(a)\,
\sum_{\nu=1}^r\,\xi_{j\nu}(x')\,
\frac{\partial}{\partial x_\nu'}
\ \ \ \ \ \ \ \ \ \ \ \ \
{\scriptstyle{(k\,=\,1\,\cdots\,r)}},
\]
in which the push-forwards $(f_a)_* \big( \partial / \partial
x_\nu \big)$ are now implicitly understood. It is easy to see that
exactly the {\rm same} equations, but with the opposite push-forwards $(
g_a)_* \big( \partial / \partial x_\nu' \big)$, can be obtained
by subjecting to similar calculations the reverse, identically
satisfied equations: $x_i \equiv g_i \big( f( x; \, a);\, a \big)$.
Consequently, we have obtained two families of equations:
\def\theequation{8}\begin{equation}
\label{two-fam-sym}
\left\{
\aligned
0
&
\equiv
\sum_{j=1}^n\,
\vartheta_{kj}(a)\,
\sum_{\nu=1}^n\,\xi_{j\nu}(x)\,
\frac{\partial}{\partial x_\nu}\,
\bigg\vert_{x\mapsto g_a(x')}
+
\sum_{j=1}^r\,\psi_{kj}(a)\,
\sum_{\nu=1}^r\,\xi_{j\nu}(x')\,
\frac{\partial}{\partial x_\nu'},
\\
0
&
\equiv
\sum_{\nu=1}^n\,
\vartheta_{kj}(a)\,
\sum_{\nu=1}^n\,\xi_{j\nu}(x)\,
\frac{\partial}{\partial x_\nu}
+
\sum_{j=1}^r\,\psi_{kj}(a)\,
\sum_{\nu=1}^r\,\xi_{j\nu}(x')\,
\frac{\partial}{\partial x_\nu'}\,
\bigg\vert_{x'\mapsto f_a(x)}
\
\\
&
\ \ \ \ \ \ \ \ \ \ \ \ \ \ \ \ \ \ \ \ \ \ \ \ \ \
\ \ \ \ \ \ \ \ \ \ \ \ \ \ \ \ \ \ \ \ \ \ \ \
{\scriptstyle{(k\,=\,1\,\cdots\,r)}},
\endaligned\right.
\end{equation}
in which we represent push-forwards of vector fields by the symbol of
variable replacement $x \mapsto g_a ( x')$ in the first line, and
similarly in the second line, by $x' \mapsto f_a ( x)$.

\begin{svgraybox}\vspace{-0.4cm}
\def\thetheorem{5}\begin{theorem}
\label{Theorem-5-S-45}
If the equations:
\[
x_i'
=
f_i(x_1,\dots,x_n,\,a_1,\dots,a_r)
\ \ \ \ \ \ \ \ \ \ \ \ \ {\scriptstyle{(i\,=\,1\,\cdots\,n)}}
\]
with the $r$ essential parameters $a_1, \dots, a_r$
represent an $r$-term transformation group, if, moreover, $F( x_1',
\dots, x_n')$ denotes an arbitrary function of $x_1', \dots, x_n'$ and
lastly, if the $\xi_{ j\nu} ( x)$, $\psi_{ kj} (a)$, $\vartheta_{kj}
(a)$ denote the same functions of their arguments as in the two
Theorems~3 and~4, then the relations:
\[
\aligned
\sum_{j=1}^r\,\vartheta_{kj}(a)\,
\sum_{\nu=1}^n\,\xi_{j\nu}(x)\,
&
\frac{\partial F}{\partial x_\nu}
+
\sum_{j=1}^r\,\psi_{kj}(a)\,
\sum_{\nu=1}^n\,\xi_{j\nu}(x')\,
\frac{\partial F}{\partial x_\nu'}
=
0 
\\
&\ \ \ \ \
{\scriptstyle{(k\,=\,1\,\cdots\,r)}}
\endaligned
\]
\end{theorem}
hold true after the substitution $x_1' = f_1 ( x, a)$, 
\dots, $x_n' = f_n ( x, a)$. 
\end{svgraybox}

\subsection{Synthetic, Geometric Counterpart of the Computations} 

To formulate the adequate interpretation
of the above considerations, we must introduce the two
systems of $r$ infinitesimal transformations ($1 \leqslant k \leqslant
r$):
\[
X_k
:=
\sum_{i=1}^n\,\xi_{ki}(x)\,
\frac{\partial}{\partial x_i}
\ \ \ \ \ \ \ \ \
\text{\rm and}
\ \ \ \ \ \ \ \ \ 
X_k'
:=
\sum_{i=1}^n\,\xi_{ki}(x')\,
\frac{\partial}{\partial x_i'},
\]
where the second ones are defined to be \emphasis{exactly the same vector
fields}\, as the first ones, though considered on the $x'$-space. This
target, auxiliary space $x'$ has in fact to be considered to be the
\emphasis{same} space as the $x$-space, because the considered
transformation group acts on a single individual space. So we can also
consider that $X_k'$ coincides with the value of $X_k$ at $x'$ and we
shall sometimes switch to another notation:
\[
\boxed{
X_k'
\equiv
X_k\big\vert_{x'}
}\,.
\]
Letting now $\alpha$ and $\widetilde{ \vartheta}$
be the inverse matrices of $\psi$ and of $\vartheta$, namely:
\[
\sum_{k=1}^r\,\alpha_{lk}(a)\,
\psi_{kj}(a)
=
\delta_l^j,
\ \ \ \ \ \ \ \ \ \ \ \ \ \ \ \ \
\sum_{k=1}^r\,\widetilde{\vartheta}_{lk}(a)\,
\vartheta_{kj}(a)
=
\delta_l^j,
\]
we can multiply the first (resp. the second) line 
of~\thetag{ \ref{two-fam-sym}}
by $\alpha_{ lk}(a)$ (resp. by $\widetilde{ \vartheta}_{
lk} (a)$) and then make summation over $k = 1, \dots, r$
in order to get resolved equations:
\[
\left\{
\aligned
0
&
\equiv
\sum_{k=1}^r\,\sum_{j=1}^r\,
\alpha_{lk}(a)\,\vartheta_{kj}(a)\,X_j
+
X_l'
\ \ \ \ \ \ \ \ \ \ \ \ \
{\scriptstyle{(k\,=\,1\,\cdots\,r)}},
\\
0
&
\equiv
X_l
+
\sum_{k=1}^r\,\sum_{j=1}^r\,
\widetilde{\vartheta}_{lk}(a)\,\psi_{kj}(a)\,X_j'
\ \ \ \ \ \ \ \ \ \ \ \ \
{\scriptstyle{(k\,=\,1\,\cdots\,r)}},
\endaligned\right.
\]
in which we have suppressed the push-forward symbols.
We can readily rewrite such equations under the contracted form:
\[
\aligned
X_k
=
\sum_{j=1}^r\,
\rho_{jk}(a)\,X_j'
\ \ \ \ \ \ \ \ \ \ \ \ 
&
\ \ \
\text{\rm and}
\ \ \ \ \ \ \ \ \ \ \ \ \ \ \
X_k'
=
\sum_{j=1}^r\,\widetilde{\rho}_{jk}(a)\,X_j
\\
&
{\scriptstyle{(k\,=\,1\,\cdots\,r)}},
\endaligned
\]
by introducing some two appropriate auxiliary $r \times r$ matrices
$\rho_{ jk} ( a) := - \sum_{ l=1}^r\, \widetilde{ \vartheta}_{ kl} (
a) \, \psi_{ lj} ( a)$ and $\widetilde{ \rho}_{ jk} ( a) := - \sum_{
l=1}^r\, \alpha_{ kl} ( a) \, \vartheta_{ lj} ( a)$ of
analytic functions (whose precise expression does not matter here)
which depend only upon $a$ and which, naturally, are inverses of each
other. A diagram illustrating what we have gained at that point is
welcome and intuitively helpful.

\begin{center}
\input transfer-infinitesimal.pstex_t
\end{center}

\begin{proposition}
\label{symmetric-X-X-prime}
If, in each one of the $r$ basic infinitesimal transformations of the
finite continuous transformation group $x' = f( x; \, a) = f_a ( x)$
having the inverse transformations $x = g_a ( x')$, namely if in the
vector fields:
\[
X_k
=
\sum_{i=1}^n\,
\xi_{ki}(x)\,\frac{\partial}{\partial x_i}
\ \ \ \ \ \ \ \ \ \ \ \ \
{\scriptstyle{(k\,=\,1\,\cdots\,r)}},
\ \ \ \ \ \ \ \ \ \ \ \ \
{\textstyle{
\xi_{ ki}(x)
:=
-\frac{\partial f_i}{\partial a_k}(x;\,e)}},
\]
one \emphasis{introduces} the new variables $x' = f_a ( x)$, that is
to say: \emphasis{replaces} $x$ by $g_a ( x')$ and $\frac{ \partial }{
\partial x_i}$ by $\sum_{ \nu = 1}^n \,\frac{ \partial f_\nu }{
\partial x_i }(x;\, a)\, \frac{ \partial}{ \partial x_\nu'}$, then one
necessarily obtains a linear combination of the same infinitesimal
transformations $X_l ' = \sum_{ i=1}^n\, \xi_{ ki} ( x') \, \frac{
\partial }{ \partial x_i'}$ at the point $x'$ with coefficients
depending only upon the parameters $a_1, 
\dots, a_r${\rm :}
\[
(f_a)_*
\big(
X_k\big\vert_x
\big)
=
(g_a)^*\big(
X_k\big\vert_{g_a(x')}
\big)
=
\sum_{l=1}^r\,\rho_{lk}(a_1,\dots,a_r)\,
X_l\big\vert_{x'}
\ \ \ \ \ \ \ \ \ \ \ \ \
{\scriptstyle{(k\,=\,1\,\cdots\,r)}}.
\]
Of course, through the inverse change of variable $x' \mapsto f_a (
x)$, the infinitesimal transformations $X_k'$ are subjected to similar
linear substitutions:
\[
(g_a)_*
\big(
X_k'\big\vert_{x'}
\big)
=
(f_a)^*
\big(
X_k'\big\vert_{f_a(x)}
\big)
=
\sum_{l=1}^r\,\widetilde{\rho}_{lk}(a)\,
X_l\big\vert_x
\ \ \ \ \ \ \ \ \ \ \ \ \
{\scriptstyle{(k\,=\,1\,\cdots\,r)}}.
\]
\end{proposition}

\subsection{Transfer of General Infinitesimal Transformations}

Afterwards, thanks to the linearity of the tangent
map, we deduce that the general infinitesimal
transformation of our group:
\[
X
:=
e_1\,X_1
+\cdots+
e_r\,X_r,
\]
coordinatized in the basis $\big( X_k\big)_{ 1\leqslant k \leqslant
r}$ by means of some $r$ arbitrary constants $e_1, \dots, e_r \in \K$,
then transforms as:
\[
\aligned
(g_a)^*
\Big(
e_1\,X_1+\cdots+e_r\,X_r
\big\vert_{g_a(x')}
\Big)
&
=
\sum_{k=1}^r\,e_k\,\sum_{l=1}^r\,\rho_{lk}(a)\,X_l\big\vert_{x'}
\\
&
=:
e_1'(e;\,a)\,X_1
\big\vert_{x'}
+\cdots+
e_r'(e;\,a)\,X_r
\big\vert_{x'}
\endaligned
\]
and hence we obtain that \emphasis{the change of variables $x' = f_a (
x)$ caused by a general transformation of the group acts linearly
on the space $\simeq \K^r$ of its infinitesimal transformations}:
\[
e_k'(e;\,a)
:=
\sum_{l=1}^r\,\rho_{kl}(a)\,e_l
\ \ \ \ \ \ \ \ \ \ \ \ \
{\scriptstyle{(k\,=\,1\,\cdots\,r)}},
\]
\emphasis{by just multiplying the coordinates $e_l$ by the matrix
$\rho_{ kl} ( a)$}. Inversely, we have:
\[
\label{e-e-prime-a}
e_k(e';\,a)
=
\sum_{l=1}^r\,\widetilde{\rho}_{kl}(a)\,e_l'
\ \ \ \ \ \ \ \ \ \ \ \ \
{\scriptstyle{(k\,=\,1\,\cdots\,r)}},
\]
where $\widetilde{ \rho} (a)$ is the inverse matrix of $\rho ( a)$.

\begin{svgraybox}
\def\theproposition{4}\begin{proposition}\vspace{-0.4cm}
\text{\rm (\cite{enlie1888-3}, p.~81)}
\label{Satz-4-S-81}
If the equations $x_i' = f_i(x_1,\dots, x_n, \, a_1, \dots, a_r)$
represent an $r$-term group and if this group contains the $r$
independent infinitesimal transformations:
\[
X_k(f)
=
\sum_{i=1}^n\,\xi_{ki}(x_1,\dots,x_n)\,
\frac{\partial f}{\partial x_i}
\ \ \ \ \ \ \ \ \ \ \ \ \ {\scriptstyle{(k\,=\,1\,\cdots\,r)}},
\]
then after the introduction of the new variables $x_i' = f_i(x, a)$,
the general infinitesimal transformation:
\[
e_1\,X_1(f)
+\cdots+
e_r\,X_rf
\]
keeps its form in so far as, for every system of values $e_1, \dots,
e_r$, there is a relation of the form:
\[
\sum_{k=1}^r\,e_k\,X_k(f)
=
\sum_{k=1}^r\,e_k'\,X_k'(f),
\]
where $e_1', \dots, e_r'$ are independent linear homogeneous functions
of $e_1, \dots, e_r$ with coefficients which are functions of $a_1,
\dots, a_r$.
\end{proposition}\vspace{-0.3cm}
\end{svgraybox}

\subsection{ Towards the Adjoint Action}

Afterwards, thanks to the linearity of the tangent
map, we deduce that the general transformation of our group:
\[
X
:=
e_1\,X_1
+\cdots+
e_r\,X_r,
\]
coordinatized in the basis $\big( X_k\big)_{ 1\leqslant k \leqslant
r}$ by means of some $r$ arbitrary constants $e_1, \dots, e_r \in \K$,
then transforms as:
\[
\aligned
(g_a)^*
\Big(
e_1\,X_1+\cdots+e_r\,X_r
\big\vert_{g_a(x')}
\Big)
&
=
\sum_{k=1}^r\,e_k\,\sum_{l=1}^r\,\rho_{lk}(a)\,X_l\big\vert_{x'}
\\
&
=:
e_1'(e;\,a)\,X_1
\big\vert_{x'}
+\cdots+
e_r'(e;\,a)\,X_r
\big\vert_{x'}
\endaligned
\]
and hence we obtain that \emphasis{the change of variables $x' = f_a (
x)$ performed by a general transformation of the group then acts linearly
on the space $\simeq \K^r$ of its infinitesimal transformations}:
\[
\boxed{
e_k'(e;\,a)
:=
\sum_{l=1}^r\,\rho_{kl}(a)\cdot e_l}
\ \ \ \ \ \ \ \ \ \ \ \ \
{\scriptstyle{(k\,=\,1\,\cdots\,r)}},
\]
\emphasis{by just multiplying the coordinates $e_l$ by the matrix $\rho_{ kl} (
a)$}.

\smallskip

Nowadays, the adjoint action is defined as an action of an abstract
Lie group on its Lie algebra. But in Chap.~\ref{kapitel-16}
below, Lie defines it in the more general context of a transformation
group, as follows. Employing the present way of expressing, 
one considers the general infinitesimal
transformation $X \big\vert_{ x'} = e_1 \, X_1 + \cdots + e_r X_r
\big\vert_{ x'}$ of the group as being based at the point $x'$, and
one computes the adjoint action ${\rm Ad}\, f_a \big( X \big \vert_{
x'} \big)$ of $f_a$ on $X \big \vert_{ x'}$ by differentiating at $t =
0$ the composition\footnote{\,
One terms groups $(t, x) \longmapsto
\exp (tX) ( x)$ are introduced in Chap.~\ref{one-term-groups}.
}: 
\[
f_a 
\circ 
\exp(tX)
\circ 
f_a^{-1}
\]
which represents the action of
the interior automorphism associated to $f_a$ on the one-parameter
subgroup $(t, w)
\longmapsto \exp (t X ) ( x)$ generated by $X$:
\[
\aligned
{\rm Ad}\,f_a\big(X\big\vert_{x'}\big)
:=
&\
\frac{\D}{\D\,t}
\Big(
f_a\circ\exp(tX)(\cdot)\circ f_a^{-1}
(x')
\Big)
\Big\vert_{t=0}
\\
=
&\
(f_a)_*
\frac{\D}{\D\,t}\,
\Big(
\exp(tX)
\big(f_a^{-1}(x')\big)
\Big)
\Big\vert_{t=0}
\\
=
&\
(f_a)_*\big(
X\big\vert_{f_a^{-1}(x')}\big)
\\
=
&\
(g_a)^*
\big(
X\big\vert_{g_a(x')}
\big)
\\
=
&\
(g_a)^*
\Big(
e_1\,X_1+\cdots+e_r\,X_r
\big\vert_{g_a(x')}
\Big)
\\
=
&\
e_1'(e;\,a)\,X_1\big\vert_{x'}
+\cdots+
e_r'(e;\,a)\,X_r\big\vert_{x'}.
\endaligned
\]
We thus recover \emphasis{exactly} the linear action $e_k' = e_k' ( e;
\, a_1, \dots, a_r)$ boxed above. 
A diagram is welcome: from the left, the point $x$ is sent by $f_a$ 
to $x' = f_a ( x)$\,---\,or inversely $x = g_a ( x')$\,---,
the vector $X\vert_x = \sum\, e_k\, X_k \big\vert_{ g_a ( x')}$ is sent
to $(g_a)^* \big( X \big\vert_{ g_a ( x')} \big)$ 
by means of the differential $(f_a)_* = (g_a)^*$ 
and Lie's theorem says that
this transferred vector is a linear combination of the existing vectors
$X_k \vert_{ x'}$ based at $x'$, with coefficients
depending on the $e_j$ and on the $a_j$.  

\begin{center}
\input Ad-fa.pstex_t
\end{center}

\section{Substituting the Axiom of Inverse
\\
for a Differential Equations Assumption}
\label{substituting-axiom}
\sectionmark{Substituting the Axiom of Inverse
for a Differential Equations Assumption}

Notwithstanding the counterexample discovered by Engel in 1884
(p.~\pageref{engel-counterexample}), Lie wanted in his systematic
treatise to avoid as much as possible employing both the axiom of
inverse and the axiom of identity. As a result, the first nine
fundamental chapters of~\cite{ enlie1888-3} regularly emphasize what
can be derived from only the axiom of composition.

In this section, instead of purely local hypotheses valuable in small
polydiscs (Sect.~\ref{concept-local-group}), we present the
semi-global topological hypotheses accurately made by Lie and Engel at
the very beginning of their book. After a while, they emphasize that
the fundamental differential equations $\frac{ \partial x_i'}{
\partial a_k} = \sum_{ j = 1}^r\, \psi_{ kj} ( a) \, \xi_{ ji} ( x')$,
which can be deduced from only the composition axiom, should take
place as being the main continuous group assumption. A technically
central theorem states that if some transformation equations $x '= f(
x; \, a)$ with essential parameters satisfy such kind of differential
equations $\frac{ \partial x_i'}{ \partial a_k} = \sum_{ j = 1}^r\,
\psi_{ kj} ( a) \, \xi_{ ji} ( x')$, for $x$ and $a$ running in
appropriate domains $\mathcal{ X} \subset \K^n$ and $\mathcal{ A}
\subset \K^r$, then every transformation $x' = f ( x; \, a)$ whose
parameter $a$ lies in a small neighbourhood of some fixed $a^0 \in
\mathcal{ A}$ can be obtained by firstly performing the initial
transformation $\overline{ x} = f ( x; \, a^0)$ and then secondly by
performing a certain transformation:
\[
x_i'
=
\exp\big(t\lambda_1X_1+\cdots+t\lambda_rX_r\big)(\overline{x}_i)
\ \ \ \ \ \ \ \ \ \ \ \ \
{\scriptstyle{(i\,=\,1\,\cdots\,n)}}
\]
of the one-term group 
(Chap.~\ref{one-term-groups} below)
generated by some suitable linear combination
$\lambda_1 X_1 + \cdots + \lambda_r X_r$ of the $n$ infinitesimal
transformations $X_k := \sum_{ i = 1}^n\, \xi_{ ki} ( x) \, \frac{
\partial }{ \partial x_i}$. This theorem will be of crucial use when
establishing the so-called \terminology{Second Fundamental Theorem}:
\emphasis{To any Lie algebra of local analytic vector fields is
associated a local Lie transformation group containing the identity
element} (Chap.~\ref{kapitel-9} below).

Also postponed to Chap.~\ref{kapitel-9} below, Lie's answer to Engel's
counterexample will show that, starting from transformation equations
$x' = f ( x; \, a)$ that are only 
assumed to be closed under composition, one can
always catch the identity element and all the inverses
of transformations near the identity
by appropriately changing coordinates in the parameter space
(Theorem~26 p.~\pageref{Theorem-26}).

\subsection{Specifying Domains of Existence}
\label{specifying-domains-existence} 

Thus, we consider $\K$-analytic transformation equations
$x_i' = f_i ( x; \, a)$ defined on a more general domain than a
product $\Delta_\rho^n \times \square_\sigma^r$ of two small polydiscs
centered at the origin. Here is how Lie and Engel specify their domains
of existence on p.~14 of~\cite{ enlie1888-3}, and these domains might
be global.

\begin{svgraybox}
\centerline{{\sf \S\,\,\,2.}}

In the transformation equations:
\def\theequation{1}\begin{equation}
x_i'
=
f_i(x_1,\dots,x_n;\,a_1,\dots,a_r)
\ \ \ \ \ \ \ \ \ \ \ \ \
{\scriptstyle{(i\,=\,1\,\cdots\,n)}},
\end{equation}
let now all the parameters $a_1, \dots, a_r$ be essential.

\smallskip

\smallercharacters{
Since the $f_i$ are analytic functions of their arguments, in the
domain \deutsch{Gebiete} of all systems of values $x_1, \dots, x_n$ and
in the domain of all systems of values $a_1, \dots, a_r$, we can
choose a region \deutsch{Bereich} $(x)$ and, respectively, a region
$(a)$ such that the following holds:

{\sf Firstly}. The $f_i ( x,\, a)$ are single-valued \deutsch{eindeutig}
functions of the $n + r$ variables $x_1, \dots, x_n$, $a_1, \dots,
a_r$ in the complete extension \deutsch{Ausdehnung} of the two regions
$(x)$ and $(a)$.

{\sf Secondly}. The $f_i (x, \, a)$ behave regularly in the
neighbourhood of every system of values $x_1^0, \dots, x_n^0$, $a_1^0,
\dots, a_r^0$, hence are expandable in ordinary power series with
respect to $x_1 - x_1^0$, \dots, $x_n - x_n^0$, as soon as $x_1^0,
\dots, x_n^0$ lies arbitrarily in the domain $(x)$, and
$a_1^0, \dots, a_r^0$ lies arbitrarily in the domain $(a)$.

{\sf Thirdly}. The functional determinant:
\[
\sum\,\pm\,\frac{\partial f_1}{\partial x_1}\,\cdots\,
\frac{\partial f_n}{\partial x_n}
\]
vanishes for no combination of systems of values of $x_i$ and of $a_k$
in the two domains $(x)$ and $(a)$, respectively.

{\sf Fourthly}. If one gives to the parameters $a_k$ in the equations
$x_i' = f_i( x, \, a)$ any of the values $a_k^0$ in domain $(a)$, then
the equations:
\[
x_i'
=
f_i(x_1,\dots,x_n;\,a_1^0,\dots,a_r^0)
\ \ \ \ \ \ \ \ \ \ \ \ \
{\scriptstyle{(i\,=\,1\,\cdots\,n)}}
\]
always produce two different systems of values $x_1', \dots, x_n'$ for
two different systems of values $x_1, \dots, x_n$ of the domain $(x)$.

We assume that the two regions $(x)$ and $(a)$ are chosen in such a
way that these four conditions are satisfied. If we give to the
variables $x_i$ in the equations $x_i' = f_i ( x, \, a)$ all
possible values in $(x)$ and to the parameters $a_k$ all possible
values in $(a)$, then in their domain, the $x_i'$ run throughout a
certain region, which we can denote symbolically by the equation $x' =
f \big( (x) (a) \big)$. This new domain has the following properties:

{\sf Firstly}. 
If $a_1^0, \dots, a_r^0$ is an arbitrary system of values of $(a)$ and
${x_1 '}^0, \dots, {x_n '}^0$ an arbitrary system of values of the
subregion $x' = f \big( (x) (a^0) \big)$, then in the neighbourhood of
the system of values ${ x_i'}^0$, $a_k^0$, the $x_1, \dots, x_n$ can
be expanded in ordinary power series with respect to $x_1' - {
x_1'}^0, \dots, x_n' - {x_n'}^0$, $a_1 - a_1^0, \dots, a_n - a_n^0$.

{\sf Secondly}. If one gives to the $a_k$ fixed values $a_k^0$ in
domain $(a)$, then in the equations:
\[
x_i'
=
f_i(x_1,\dots,x_n;\,a_1^0,\dots,a_r^0)
\ \ \ \ \ \ \ \ \ \ \ \ \
{\scriptstyle{(i\,=\,1\,\cdots\,n)}},
\]
the quantities $x_1, \dots, x_n$ will be single-valued functions of
$x_1', \dots, x_n'$, which behave regularly in the complete extension
of the region $x' = f \big( (x) (a^0) \big)$.

}
\vspace{-0.2cm}
\end{svgraybox}

\smallskip
Referring also to the
excerpt p.~\pageref{excerpt-shrinking}, 
we may rephrase these basic assumptions as follows. 
The $f_i (x; \, a)$ are defined for $(x, a)$ belonging to the
product:
\[
\mathcal{X}\times\mathcal{A}\,
\subset\,\K^n\times\K^r
\]
of two domains $\mathcal{ X} \subset \K^n$ and $\mathcal{ A} \subset
\K^r$. These functions are $\K$-analytic in both variables, hence
expandable in Taylor series at every point $(x^0, a^0)$. Furthermore,
for every fixed $a^0$, the map $x \mapsto f ( x; \, a^0)$ is assumed to
constitute a $\K$-analytic diffeomorphism of $\mathcal{ X} \times \{
a^0 \}$ onto its image. Of course, the inverse map is also locally
expandable in power series, by virtue of the $\K$-analytic inverse
function theorem.

To insure 
that the composition of two transformations exists, one requires
that there exist nonempty subdomains:
\[
\mathcal{X}^1\subset\mathcal{X}
\ \ \ \ \ \ \ \ \ 
\text{\rm and}
\ \ \ \ \ \ \ \ \ 
\mathcal{A}^1\subset\mathcal{A},
\]
with the property that for every fixed $a^1 \in \mathcal{ A}^1$:
\[
f\big(\mathcal{X}^1\times\{a^1\}\big)
\subset
\mathcal{X},
\]
so that for every such an $a^1 \in \mathcal{ A}^1$ and 
for every fixed $b \in \mathcal{ A}$, the composed map:
\[
x
\longmapsto
f\big(f(x;\,a^1);\,b\big)
\]
is well defined for all $x \in \mathcal{ X}^1$ and moreover, is
a $\K$-analytic diffeomorphism
onto its image. In fact, it is even $\K$-analytic with
respect to all the variables $(x, a^1, b)$ in $\mathcal{ X}^1 \times
\mathcal{ A}^1 \times \mathcal{ A}$. Lie's fundamental
and unique
\terminology{group composition axiom} may then
be expressed as follows. 

\smallskip\begin{itemize}

\item[{\bf \!\!\!(A1)}] \ \ \
There exists a $\K^m$-valued 
$\K$-analytic map $\varphi = \varphi ( a, b)$ defined in
$\mathcal{ A}^1 \times \mathcal{ A}^1$ 
with $\varphi \big( \mathcal{ A}^1
\times \mathcal{ A}^1\big) \subset \mathcal{ A}$ such that:
\[
f\big(f(x;\,a);\,b\big)
\equiv
f\big(x;\,\varphi(a,b)\big)
\ \ \ \ \ \ \ \ \ \ 
\text{\rm for all}\ \
x\in\mathcal{X}^1,\ 
a\in\mathcal{A}^1,\
b\in\mathcal{A}^1.
\]

\end{itemize} 

Here are two further specific unmentioned assumptions that Lie
presupposes, still with the goal of admitting neither the identity
element, nor the existence of inverses.

\smallskip\begin{itemize}

\item[{\bf \!\!\!(A2)}] \ \ \ \label{A-2}
There is a $\K^m$-valued $\K$-analytic map $(a, c) 
\longmapsto {\sf b} = {\sf b} ( a, c)$
defined for $a$ running in a certain (nonempty) subdomain $\mathcal{
A}^2 \subset \mathcal{ A}^1$ and for $c$ running in a certain
(nonempty) subdomain $\mathcal{ C}^2 \subset \mathcal{ A}^1$ with
${\sf b} \big( \mathcal{ A}^2 \times \mathcal{ C}^2 \big) \subset
\mathcal{ A}^1$ which solves $b$ in terms of $(a, c)$ in the equations
$c_k = \varphi_k ( a, b)$, namely which satisfies identically:
\[
c
\equiv
\varphi\big(a,\,{\sf b}(a,c)\big)
\ \ \ \ \ \ \ \ \ \ 
\text{\rm for all}\ \ 
a\in\mathcal{A}^2,\
c\in\mathcal{C}^2.
\]
Inversely, $\varphi ( a,b)$ solves $c$ in terms of $(a, b)$ in the
equations $b_k = {\sf b}_k 
( a, c)$, namely more precisely: there exists a
certain (nonempty) subdomain $\mathcal{ A}^3 \subset 
\mathcal{ A}^2$ ($\subset \mathcal{ A}^1$)
and a certain (nonempty) subdomain $\mathcal{ B}^3 \subset \mathcal{
A}^1$ with $\varphi \big( \mathcal{ A}^3 \times \mathcal{ B}^3 \big)
\subset \mathcal{ C}^2$ such that one has identically:
\[
b
\equiv
{\sf b}\big(a,\,\varphi(a,b)\big)
\ \ \ \ \ \ \ \ \ \ 
\text{\rm for all}\ \ 
a\in\mathcal{A}^3,\
b\in\mathcal{B}^3.
\]

\end{itemize}\smallskip

{\sf Example.}
In Engel's counterexample of the group $x' = \chi ( \lambda) \,x$ with
a Riemann uniformizing map $\omega : \Delta \to \Lambda$ as on
p.~\pageref{engel-counterexample} having inverse $\chi : \Lambda \to
\Delta$, these three requirements are satisfied, and in addition, we
claim that one may even take $\mathcal{ X} = \K$ and $\mathcal{ A}^1 =
\mathcal{ A} = \Lambda$, with no shrinking, for composition happens to
hold in fact without restriction in this case.  Indeed, starting from
the general composition:
\[
x'' 
= 
\chi(\lambda_2)\, 
x' 
= 
\chi(\lambda_2)\,\chi(\lambda_1)\,x,
\] 
that is to say, from $x'' = \chi (\lambda_2) \, \chi (\lambda_1) \,
x$, in order to represent it under the specific form $x'' = \chi
(\lambda_3) \, x$, it is necessary and sufficient to solve $\chi
(\lambda_3) = \chi ( \lambda_1) \chi(\lambda_2)$, hence we may take
for $\varphi$:
\[
\lambda_3
=
\omega\big(\chi(\lambda_1)\chi(\lambda_2)\big)
=:
\varphi(\lambda_1,\lambda_2),
\]
without shrinking the domains, for the two inequalities $\vert \chi
(\lambda_1) \vert < 1$ and $\vert \chi ( \lambda_2) \vert < 1$ readily
imply that $\vert \chi (\lambda_1)\chi ( \lambda_2) \vert < 1$ too so
that $\omega \big( \chi(\lambda_1)\chi(\lambda_2) \big)$ is defined.
On the other hand, for solving $\lambda_2$ in terms of $(\lambda_1,
\lambda_3)$ in the above equation, we are naturally led to define:
\[
{\sf b}(\lambda_1,\lambda_3)
:=
\omega\big(\chi(\lambda_3)\big/\chi(\lambda_1)\big),
\]
and then ${\sf b} = {\sf b} ( \lambda_1, \lambda_3)$ is defined under
the specific restriction that $\vert \chi (\lambda_3) \vert < \vert
\chi ( \lambda_1) \vert$.

\smallskip
However, the axiom {\bf (A2)} happens to be still
incomplete for later use, and one should add the following axiom
in order to be able to also solve $a$ in $c = \varphi ( a, b)$. 

\smallskip\begin{itemize}

\item[{\bf \!\!\!(A3)}] \ \
There is a $\K^m$-valued $\K$-analytic map $(b, c) \longmapsto {\sf a}
= {\sf a} ( b, c)$ defined in $\mathcal{ B}^4 \times \mathcal{ C}^4$
with $\mathcal{ B}^4 \subset \mathcal{ A}^1$ and $\mathcal{ C}^4
\subset \mathcal{ A}^1$, and with ${\sf a} \big( \mathcal{ B}^4 \times
\mathcal{ C}^4 \big) \subset \mathcal{ A}^1$, such that one has
identically:
\[
c
\equiv
\varphi\big({\sf a}(b,c),\,b\big)
\ \ \ \ \ \ \ \ \ \ 
\text{\rm for all}\ \ 
b\in\mathcal{B}^4,\
c\in\mathcal{C}^4. 
\]
Inversely, $\varphi ( a, b)$ solves $c$ in the equations $a_k = {\sf
a}_k ( b, c)$, namely more precisely: there exist $\mathcal{ B}^5
\subset \mathcal{ B}^4$ and $\mathcal{ A}^5 \subset \mathcal{ A}^1$
with $\varphi \big( \mathcal{ A}^5 \times \mathcal{ B}^5 \big) \subset
\mathcal{ C}^4$ such that one has identically:
\[
a
\equiv
{\sf a}
\big(b,\,\varphi(a,b)\big)
\ \ \ \ \ \ \ \ \ \ 
\text{\rm for all}\ \ 
a\in\mathcal{A}^5,\
b\in\mathcal{B}^5, 
\]
and furthermore in addition, with
${\sf a} \big( \mathcal{ B}^5 \times
\mathcal{ C}^4 \big) \subset \mathcal{ A}^2$ and
with ${\sf b} \big( \mathcal{ A}^5 \times 
\mathcal{ C}^4 \big) \subset \mathcal{ B}^4$
such that one also has identically: 
\[
\aligned
b
&
\equiv
{\sf b}\big({\sf a}(b,c),\,c\big)
\ \ \ \ \ \ \ \ \ \ 
\text{\rm for all}\ \ 
b\in\mathcal{B}^5,\
c\in\mathcal{C}^4
\\
a
&
\equiv
{\sf a}\big({\sf b}(a,c),\,a\big) 
\ \ \ \ \ \ \ \ \ \ 
\text{\rm for all}\ \ 
a\in\mathcal{A}^5,\
c\in\mathcal{C}^4.
\endaligned
\]
\end{itemize}

The introduction of the numerous (nonempty) domains
$\mathcal{ A}^2$, $\mathcal{ C}^2$, 
$\mathcal{ A}^3$, $\mathcal{ B}^3$, 
$\mathcal{ B}^4$, $\mathcal{ C}^4$, 
$\mathcal{ A}^5$, $\mathcal{ B}^5$ which 
appears slightly unnatural and seems to
depend upon the order in which 
the solving maps ${\sf a}$ and ${\sf b}$ are considered
can be avoided by requiring 
from the beginning only that there exist two 
subdomains $\mathcal{ A}^3 \subset \mathcal{ A}^2
\subset \mathcal{ A}^1$ such that one has uniformly: 
\[
\aligned
c
&
\equiv
\varphi\big(a,\,{\sf b}(a,c)\big) 
\ \ \ \ \ \ \ \ \ \ 
\text{\rm for all}\ \ 
a\in\mathcal{A}^3,\
c\in\mathcal{A}^3
\\
c
&
\equiv
\varphi\big({\sf a}(b,c),\,b\big) 
\ \ \ \ \ \ \ \ \ \ 
\text{\rm for all}\ \ 
b\in\mathcal{A}^3,\
c\in\mathcal{A}^3
\\
b
&
\equiv
{\sf b}\big({\sf a}(b,c),\,c\big) 
\ \ \ \ \ \ \ \ \ \ \
\text{\rm for all}\ \ 
b\in\mathcal{A}^3,\
c\in\mathcal{A}^3
\\
b
&
\equiv
{\sf b}\big(a,\,\varphi(a,b)\big)
\ \ \ \ \ \ \ \ \ \,
\text{\rm for all}\ \ 
a\in\mathcal{A}^3,\
b\in\mathcal{A}^3
\\
a
&
\equiv
{\sf a}\big({\sf b}(a,c),\,c\big)
\ \ \ \ \ \ \ \ \ \ \
\text{\rm for all}\ \ 
a\in\mathcal{A}^3,\
c\in\mathcal{A}^3
\\
a
&
\equiv
{\sf a}\big(b,\,\varphi(a,b)\big)
\ \ \ \ \ \ \ \ \ \ 
\text{\rm for all}\ \ 
a\in\mathcal{A}^3,\
b\in\mathcal{A}^3.
\endaligned
\]
We will adopt these axioms in the next subsection.  Importantly, we
would like to point out that, although ${\sf b} (a, c)$ seems to
represent the group product $a^{ -1} \cdot c = 
{\bf m} \big( {\bf i} (a), \, c \big)$, the assumption {\bf (A2)}
neither reintroduces inverses, nor the identity element, it just means
that one may solve $b$ by means of the implicit function theorem in
the parameter composition equations $c_k = \varphi_k ( a_1, \dots,
a_r, \, b_1, \dots, b_r)$.

\subsection{Group Composition Axiom
And Fundamental Differential Equations}

As was said earlier on, the fundamental Theorem~3,
p.~\pageref{Theorem-3-S-33} about differential equations satisfied by
a transformation group was in fact stated and proved in~\cite{
enlie1888-3} under semi-global assumptions essentially equivalent to
the ones we just formulated above with $\mathcal{ A}^3 \subset \mathcal{
A}^2 \subset \mathcal{ A}^1$, and now, we can restitute it really.

\begin{proposition}
\label{economic-differential}
Under these assumptions, there is an $r \times r$ matrix of functions
$( \psi_{ kj} ( a) )_{ 1 \leqslant k \leqslant r}^{ 1
\leqslant j \leqslant r}$ which is
$\K$-analytic and invertible in $\mathcal{ A}^3$, and there are
certain functions $\xi_{ j i} ( x)$, $\K$-analytic in $\mathcal{ X}$
such that the following differential equations:
\def\theequation{2}\begin{equation}
\frac{\partial x_i'}{\partial a_k}(x;\,a)
=
\sum_{j=1}^r\,\psi_{kj}(a)\,\xi_{j i}(x')
\ \ \ \ \ \ \ \ \ \ \ \ \
{\scriptstyle{(i\,=\,1\,\cdots\,n\,;\,\,\,k\,=\,1\,\cdots\,r)}}
\end{equation}
are identically satisfied for all $x \in \mathcal{ X}^1$ and all $a
\in \mathcal{ A}^3$ after replacing $x'$ by $f( x; \, a)$.  Here, the
functions $\xi_{ ji} (x')$ are defined by choosing arbitrarily some
fixed $b^0 \in \mathcal{ A}^3$ and by setting:
\[
\xi_{ji}(x_1',\dots,x_n')
=
\bigg[
\frac{\partial x_i'}{\partial b_j}
\bigg]_{b=b^0}
=
\bigg[
\sum_{k=1}^r\,
\frac{\partial x_i'}{\partial a_k}\,
\frac{\partial{\sf a}_k}{\partial b_j}
\bigg]_{b=b^0},
\]
and moreover, in the equations inverse to~\thetag{ 2}:
\[
\xi_{ji}(x')
=
\sum_{k=1}^r\,
\alpha_{jk}(a)\,
\frac{\partial x_i'}{\partial a_k}(x;\,a),
\]
the (inverse) coefficients $\alpha_{ jk} (a)$
are defined by:
\[
\alpha_{jk}(a)
=
\bigg[
\frac{\partial{\sf a}_k}{\partial b_j}
\bigg]_{b=b^0}.
\] 
\end{proposition}

\begin{proof}\smartqed
The computations which we have already conducted on
p.~\pageref{beginning-computation} for a purely local transformation
group can here be performed in a similar way. In brief, using $b =
{\sf b} ( a, c)$ from {\bf (A2)}, we consider the identities:
\[
f_i\big(f(x;\,a);\,{\sf b}(a,c)\big)
\equiv
f_i(x;\,c)
\ \ \ \ \ \ \ \ \ \ \ \ \
{\scriptstyle{(i\,=\,1\,\cdots\,n)}}
\]
and we differentiate them with respect to $a_k$; if we
shortly denote $f_i' \equiv f_i \big( x'; 
{\sf b} (a, c) \big)$ and $x_j' \equiv f_j(x;a)$, 
this gives:
\[
\sum_{\nu=1}^n\,
\frac{\partial f_i'}{\partial x_\nu'}\,
\frac{\partial x_\nu'}{\partial a_k}
+
\sum_{j=1}^r\,
\frac{\partial f_i'}{\partial b_j}\,
\frac{\partial {\sf b}_j}{\partial a_k}
\equiv
0
\ \ \ \ \ \ \ \ \ \ \ \ \
{\scriptstyle{(i\,=\,1\,\cdots\,n)}}.
\]
By the diffeomorphism assumptions, the matrix $\big( \frac{ \partial
f_i'}{ \partial x_\nu'} (x'; \, b) \big)_{ 1 \leqslant i \leqslant
n}^{ 1 \leqslant \nu \leqslant n}$ has a $\K$-analytic inverse for all
$(x', b) \in \mathcal{ X} \times \mathcal{ A}$, so an application of
Cramer's rule yields a resolution of the form:
\[
\aligned
\frac{\partial x_\nu'}{\partial a_k}
(x;\,a)
&
=
\Xi_{1\nu}\big(x',{\sf b}(a,c)\big)\,
\frac{\partial {\sf b}_1}{\partial a_k}(a,c)
+\cdots+
\Xi_{r\nu}\big(x',{\sf b}(a,c)\big)\,
\frac{\partial {\sf b}_r}{\partial a_k}(a,c)
\\
&
\ \ \ \ \ \ \ \ \ \ \ \ \ \ \ \ \ \ \ \ \ \ \ \
{\scriptstyle{(\nu\,=\,1\,\cdots\,n\,;\,\,\,
k\,=\,1\,\cdots\,r)}}.
\endaligned
\]
Here of course, $a,\, c \in \mathcal{ A}^3$.  Then we replace $c$ by
$\varphi ( a, b)$, with $a,\, b \in \mathcal{ A}^3$, we confer to $b$
any fixed value, say $b^0$ (remind that the identity $e$ is not
available), and we get the desired differential equations with $\xi_{
ji} ( x') := \Xi_{ ji} ( x', b^0)$ and with $\psi_{ k j} ( a) :=
\frac{\partial {\sf b}_j}{ \partial a_k} ( a, b^0)$.  Naturally, the
invertibility of the matrix $\psi_{ kj} (a)$ comes from {\bf (A2)}.

Next, we multiply each equation just obtained
(changing indices):  
\def\theequation{2'}\begin{equation}
\frac{\partial x_i'}{\partial a_k}(x;\,a)
\equiv
\sum_{\nu=1}^r\,
\frac{\partial{\sf b}_\nu}{\partial a_k}(a,c)\,
\Xi_{\nu i}
\big(x',\,{\sf b}(a,c)\big)
\end{equation}
by $\frac{ \partial {\sf a}_k}{ \partial b_l} ( b, c)$, where $l \in
\{ 1, \dots, r\}$ is fixed, we consider $(b, c)$, instead of $(a, c)$,
as the $2\, r$ independent variables while $a = {\sf a} ( b, c)$, and
we sum with respect to $k$ for $k = 1, \dots, r$:
\[
\sum_{k=1}^r\,
\frac{\partial x_i'}{\partial a_k}\,
\frac{\partial{\sf a}_k}{\partial b_l}
=
\sum_{\nu=1}^r\,\sum_{k=1}^r\,
\frac{\partial{\sf a}_k}{\partial b_l}\,
\frac{\partial{\sf b}_\nu}{\partial a_k}\,
\Xi_{\nu i}. 
\]
Now, the chain rule and the fact that $\partial {\sf a}_k / \partial
b_l$ is the inverse matrix of $\partial {\sf b}_\nu / 
\partial a_k$ enables
us to simplify both sides (interchanging members):
\def\theequation{3}\begin{equation}
\sum_{k=1}^r\,
\frac{\partial x_i'}{\partial a_k}\,
\frac{\partial{\sf a}_k}{\partial b_l}
=
\frac{\partial x_i'}{\partial b_l}
=
\Xi_{li}
\end{equation}
Specializing $b := b^0 \in \mathcal{ A}^3$, we get the announced
representation: 
\[
\xi_{ji}(x') 
=
\Xi_{ji}(x',\,b^0)
=
\partial x_i'/\partial b_j
\big\vert_{b=b_0},
\]
and by identification, we also obtain at the same time the
representation $\alpha_{ jk} ( a) = \partial {\sf a}_k /
\partial b_j \big\vert_{b = b_0}$.
\qed\end{proof}

We end up by observing that, similarly as we did on
p.~\pageref{c-m-a-b}, one could proceed to some further computations,
although it would not really be needed for the proposition. We
may indeed differentiate the equations:
\[
c_\mu
\equiv
\varphi_\mu
\big({\sf a}(b,c),\,b)
\ \ \ \ \ \ \ \ \ \ \ \ \ 
{\scriptstyle{(\mu\,=\,1\,\cdots\,r)}}
\]
with respect to $b_l$, and for this, 
we translate a short passage of~\cite{
enlie1888-3}, p.~20.
\begin{svgraybox}
Hence one has:
\[
\sum_{k=1}^r\,
\frac{\partial\varphi_\mu}{\partial a_k}\,
\frac{\partial {\sf a}_k}{\partial b_l}
+
\frac{\partial\varphi_\mu}{\partial b_l}
=
0
\ \ \ \ \ \ \ \ \ \ \ \ \ 
{\scriptstyle{(\mu,\,\,l\,=\,1\,\cdots\,r)}},
\]
whence it comes:
\[
\frac{\partial {\sf a}_k}{\partial b_l}
=
-\,
\frac{\sum\,\pm\,
\frac{\partial\varphi_1}{\partial a_1}
\,\cdots\,\,
\frac{\partial\varphi_{k-1}}{\partial a_{k-1}}\,
\frac{\partial\varphi_k}{\partial b_l}\,
\frac{\partial\varphi_{k+1}}{\partial a_{k+1}}
\,\cdots\,\,
\frac{\partial\varphi_r}{\partial a_r}}{
\sum\,\pm\,
\frac{\partial\varphi_1}{\partial a_1}
\,\cdots\,\,
\frac{\partial\varphi_r}{\partial a_r}}
=
A_{l k}(a_1,\dots,a_r,\,b_1,\dots,b_r).
\]
When we insert these values in~\thetag{ 3}, we obtain the equations:
\def\theequation{3'}\begin{equation}
\Xi_{li}(x',b)
=
\sum_{k=1}^r\,A_{l k}(a,\,b)\,
\frac{\partial x_i'}{\partial a_k}
\ \ \ \ \ \ \ \ \ \ \ \ \ 
{\scriptstyle{(i\,=\,1\,\cdots\,n\,;\,\,\,
l\,=\,1\,\cdots\,r)}}.
\end{equation}\vspace{-0.4cm}
\end{svgraybox}
\noindent
Of course by identification with~\thetag{ 3}, we must then have $A_{
lk} ( a, \, b) = \partial {\sf a}_k / \partial b_l$ here.  In
conclusion, we have presented the complete thought of Lie and Engel,
who did not necessarily consider the axioms of groups to be strictly
local.

\subsection{The Differential Equations Assumption
\\
And Its Consequences}

Now, we would like to emphasize that instead of axioms {\bf (A1-2-3)},
in his answer to Engel's counterexample and in several other places as
well, Lie says that he wants to set as a fundamental hypothesis the
existence of a system of differential equations as the one above, with
an invertible matrix $\psi_{ k j} ( a)$.

\begin{svgraybox}
\centerline{{\sf \S\,\,\,17. (\cite{ enlie1888-3}, pp.~67--68)}}

\smallskip
For the time being, we want to retain from assuming that the equations
$x_i' = f_i ( x_1, \dots, x_n, \, a_1, \dots, a_r)$ should represent
an $r$-term group. Rather, about the equations~\thetag{ 1}, we want
only to assume: \emphasis{firstly}, that they represent a family of
$\infty^r$ different transformations, hence that the $r$ parameters
$a_1, \dots, a_r$ are all essential, and \emphasis{secondly}, that they
satisfy differential equations of the specific form~\thetag{ 2}.
\end{svgraybox}

\smallskip

So renaming the domain $\mathcal{ A}^3$ considered above simply as
$\mathcal{ A}^1$, we will fundamentally assume in
Sect.~\ref{application-theorem-9} of Chap.~\ref{one-term-groups} that
differential equations of the specific form~\thetag{ 2} hold for all
$x \in \mathcal{ X}^1$ and all $a \in \mathcal{ A}^1$, 
\emphasis{forgetting
completely about composition}, and most importantly, 
\emphasis{without assuming
neither the existence of the identity element, nor assuming existence
of inverse transformations}.  At first,
as explained by Engel and Lie, one can easily deduce from
such new economical assumptions two basic nondegeneracy conditions.

\begin{lemma}
\label{lemma-two-nondegeneracies}
Consider transformation equations $x_i' = f_i ( x; \, a)$ defined for
$x \in \mathcal{ X}$ and $a \in \mathcal{ A}$ having essential
parameters $a_1, \dots, a_r$ which satisfy differential equations of
the form:
\def\theequation{2}\begin{equation}
\frac{\partial x_i'}{\partial a_k}(x;\,a)
=
\sum_{j=1}^r\,\psi_{kj}(a)\,\xi_{j i}(x')
\ \ \ \ \ \ \ \ \ \ \ \ \
{\scriptstyle{(i\,=\,1\,\cdots\,n\,;\,\,\,k\,=\,1\,\cdots\,r)}},
\end{equation}
for all $x \in \mathcal{ X}^1$ and all
$a \in \mathcal{ A}^1$. 
Then the determinant of the $\psi_{ kj} (a)$ does not vanish
identically and furthermore, the $r$ infinitesimal transformations:
\[
X_k'
=
\sum_{i=1}^n\,\xi_{ki}(x')\,
\frac{\partial}{\partial x_i'}
\ \ \ \ \ \ \ \ \ \ \ \ \
{\scriptstyle{(k\,=\,1\,\cdots\,r)}}
\]
are independent of each other.
\end{lemma}

\begin{proof}\smartqed
If the determinant of the $\psi_{ kj} (a)$ would vanish identically,
there would exist a locally defined (for $a$ running in the locus
where $\psi_{ kj} ( a)$ is of maximal, locally constant rank) nonzero
$\K$-analytic vector $\big(\tau_1 ( a), \dots, \tau_r ( a)\big) $ in
the kernel of $\psi_{ kj} ( a)$. Consequently, after
multiplying~\thetag{ 2} by $\tau_k ( a)$, we would derive the
equations: $\sum_{ k=1}^r\, \tau_k ( a) \, \frac{ \partial x_i'}{
\partial a_k} ( x; \, a) \equiv 0$ which would then contradict
essentiality of parameters, according to the theorem on
p.~\pageref{Theorem-essential}.

As a result, for any $a$ belonging to open set where $\det \psi_{ kj}
( a) \neq 0$, we can locally
invert the differential equations~\thetag{ 2} and
write them down under the form:
\[
\xi_{ji}(x')
=
\sum_{k=1}^r\,\alpha_{jk}(a)\,
\frac{\partial x_i'}{\partial a_k}(x;\,a)
\ \ \ \ \ \ \ \ \ \ \ \ \
{\scriptstyle{(i\,=\,1\,\cdots\,n;\,\,\,j\,=\,1\,\cdots\,r)}}.
\] 
If there would exist constants $e_1', \dots, e_r'$ not all zero with
$e_1' X_1' + \cdots + e_r' X_r ' = 0$, we would then deduce the
relation:
\[
0
\equiv 
\sum_{k=1}^r\,
\sum_{j=1}^r\,e_j'\,
\alpha_{jk}(a)\, 
\frac{\partial x_i'}{\partial a_k}
(x;\,a) 
\]
which would again contradict essentiality of parameters.
\qed\end{proof}

\subsection{Towards the Theorem~26}

At the end of the next chapter in Sect.~\ref{application-theorem-9},
we shall be in a position to pursue the restoration of further
refined propositions towards Lie's Theorem~26 (translated in
Chap.~\ref{kapitel-29}, p.~\pageref{Theorem-26}) which will answer
fully Engel's counterexample. In brief, and by anticipation, this
theorem states what follows.

Let $x_i' = f_i ( x; \, a_1, \dots, a_r)$ be a
family of transformation equations
which is only assumed to be closed under composition, a
\terminology{finite continuous transformation group}, in the sense of
Lie. According to Proposition~\ref{economic-differential} above, 
there exists a system of fundamental 
differential equations of the form: 
\[
\frac{\partial x_i'}{\partial a_k}
=
\sum_{j=1}^r\,\psi_{kj}(a)\cdot\xi_{ji}(x')
\ \ \ \ \ \ \ \ \ \ \ \ \
{\scriptstyle{(i\,=\,1\,\cdots\,n\,;\,k\,=\,1\,\cdots\,r)}}
\]
which is identically satisfied by the functions $f_i ( x, a)$,
where the $\psi_{ kj}$ are certain analytic functions
of the parameters
$(a_1, \dots, a_r)$.
If one introduces the $r$ infinitesimal transformations:
\[
\sum_{i=1}^n\,\xi_{ki}(x)\,\frac{\partial f}{\partial x_i}
=:
X_k(f)
\ \ \ \ \ \ \ \ \ \ \ \ \
{\scriptstyle{(k\,=\,1\,\cdots\,r)}},
\]
and if one forms the so-called 
\terminology{canonical} finite equations: 
\[
\aligned
x_i'
&
=
\exp\big(\lambda_1X_1+\cdots+\lambda_kX_k\big)(x)
\\
&
=:
g_i(x;\,\lambda_1,\dots,\lambda_r)
\ \ \ \ \ \ \ \ \ \ \ \ \
{\scriptstyle{(i\,=\,1\,\cdots\,n)}}
\endaligned
\]
of the $r$-term group which is generated by these $r$ infinitesimal
transformations\,---\,\voir~the next Chap.~\ref{one-term-groups} for
$\exp X(x)$\,---\,then this group contains the identity element,
namely $g( x; \, 0)$, and its transformations are ordered as inverses
by pairs, namely: $g(x; \, - \lambda) = g (x; \,
\lambda)^{ -1}$. Lastly, 
the Theorem~26 in questions states that in these finite equations
$x_i' = g_i ( x; \, \lambda)$, it is possible to introduce
\emphasis{new} local parameters $\overline{ a}_1,
\dots, \overline{ a}_r$
in place of $\lambda_1, \dots, \lambda_r$ so that the resulting
transformation equations:
\[
\aligned
x_i'
&
=
g_i\big(x;\,\lambda_1(\overline{a}),\dots,\lambda_r(\overline{a})\big)
\\
&
=:
\overline{f}_i
(x_1,\dots,x_n,\,\overline{a}_1,\dots,\overline{a}_r)
\ \ \ \ \ \ \ \ \ \ \ \ \
{\scriptstyle{(i\,=\,1\,\cdots\,n)}}
\endaligned
\]
represent a family of $\infty^r$ transformations which embraces,
\emphasis{possibly after analytic prolongation}, all the $\infty^r$
initial transformations:
\[
x_i'
=
f_i(x_1,\dots,x_n,\,a_1,\dots,a_r).
\]
In this way, Lie not only answers Engel's counterexample $x' = \chi (
\lambda) \, x$ by saying that one has to restitute the plain
transformations $x' = \zeta\, x$ by changing appropriately coordinates
in the parameter space, but also, Lie really establishes the
conjecture he suspected (quotation
p.~\pageref{conjecture-composition}), modulo the fact that the
conjecture was not true without a suitable change of coordinates in
the parameter space. To our knowledge, no modern treatise restitutes
this theorem of Lie, although a great deal of the first 170 pages of
the Theorie der Transformationsgruppen is devoted to economize the
axiom of inverse.  We end up this chapter by a quotation
(\cite{enlie1888-3}, pp.~81--82) motivating the introduction of
one-term groups $\exp (tX) (x)$ in the next chapter.

\subsection{Metaphysical Links With Substitution Theory}

We conclude this chapter with a brief quotation motivating
what will follow.

\renewcommand{\thefootnote}{\fnsymbol{footnote}}
\begin{svgraybox}\indent
The concepts and the propositions of the theory of the continuous
transformation groups often have their analogues in the
\emphasis{theory of substitutions}\footnotemark[1], that is to say, in
the theory of the discontinuous groups. In the course of our studies,
we will not emphasize this analogy every time, but we will more often
remember it by translating the terminology of the theory of
substitutions into the theory of transformation groups, and this shall
take place as far as possible.

Here, we want to point out that the \emphasis{one-term groups} in the
theory of transformation groups play the same rôle as the
\emphasis{groups generated by a single substitution} in the theory of
substitutions.

\medskip

In a way, we shall consider the one-term groups, or their 
infinitesimal transformations, as the elements of the $r$-term
group. In the studies about $r$-term groups, it is, 
almost in all circumstances, 
advantageous to direct at first the attention towards
the infinitesimal transformations of the 
concerned group and to choose them as
the object of study. 
\end{svgraybox}

\footnotetext[1]{\,
\name{C. Jordan}, Traité des substitutions, Paris 1870. 
} 

\renewcommand{\thefootnote}{\arabic{footnote}}

\linestop


\setcounter{footnote}{0}

\chapter{One-Term Groups 
\\
and Ordinary Differential Equations}
\label{one-term-groups}
\chaptermark{One-Term Groups and Ordinary Differential Equations}

\abstract{
The flow $x' = \exp ( tX) ( x)$ of a single, arbitrary vector field $X
= \sum_{ i = 1}^n\, \xi_i ( x) \, \frac{ \partial }{
\partial x_i}$ with analytic coefficients
$\xi_i ( x)$ always generates a one-term (local) continuous
transformation group:
\[
\exp\big(t_1X\big(\exp(t_2X)(x)\big)\big)
=
\exp\big((t_1+t_2)X\big)(x)
\ \ \ \ \ \ \
\text{\rm and}
\ \ \ \ \ \ \
\left[
\exp(tX)(\cdot)
\right]^{-1}
=
\exp(-tX)(\cdot).
\]
In a neighbourhood of any point at which $X$ does not vanish, an
appropriate local diffeomorphism $x \mapsto y$ may straighten $X$ to
just $\frac{ \partial }{ \partial y_1}$, hence its flow becomes
$y_1' = y_1 + t$, $y_2 ' = y_2, \dots, y_n' = y_n$.
\newline\indent
In fact, in the analytic category (only), computing a general flow
$\exp( tX) ( x)$ amounts to adding the differentiated terms
appearing in the formal expansion of \emphasis{Lie's 
exponential series}:
\[
\exp(tX)(x_i)
=
\sum_{k\geqslant 0}\,
\frac{(tX)^k}{k!}(x_i)
=
x_i
+
t\,X(x_i)
+\cdots+
\frac{t^k}{k!}\,
\underbrace{X\big(\cdots\big(
X\big(X}_{k\,\,\text{\rm times}}(x_i)\big)\big)\cdots\big)
+\cdots,
\]
that have been studied extensively by Gr\"obner in~\cite{ gr1960}. 
\newline\indent
The famous Lie bracket is introduced by looking at the way how a
vector field $X = \sum_{ i = 1}^n \,
\xi_i ( x) \frac{ \partial }{ \partial x_i}$ is
perturbed, to first order, while introducing the new coordinates $x' =
\exp ( tY) ( x) =: \varphi ( x)$ provided by the flow of another
vector field $Y$:
\[
\varphi_*(X)
=
X'
+
t\,\left[X',\,Y'\right]
+\cdots,
\]
with $X ' = \sum_{ i=1}^n \, \xi_i ( x') \,
\frac{ \partial }{\partial x_i'}$ and
$Y' = \sum_{ i = 1}^n \, \eta_i ( x') \,
\frac{ \partial }{ \partial x_i'}$ denoting the two
vector fields in the target space $x'$ having the 
\emphasis{same} coefficients as
$X$ and $Y$. Here, the analytical expression of the Lie bracket is:
\[
\left[X',\,Y'\right]
=
\sum_{i=1}^n\,
\bigg(
\sum_{l=1}^n\,
\xi_l(x')\,\frac{\partial\eta_i}{\partial x_l'}(x')
-
\eta_l(x')\,\frac{\partial\xi_i}{\partial x_l'}(x')
\bigg)\,
\frac{\partial}{\partial x_i'}.
\]
\newline\indent
An $r$-term group $x' = f ( x; \, a)$ satisfying his fundamental
differential equations $\frac{ \partial x_i'}{ \partial a_k} = \sum_{
j = 1}^r \,
\psi_{ kj} ( a) \, \xi_{ ji} ( x')$ can, alternatively,
be viewed as being generated by its infinitesimal transformations $X_k
= \sum_{ i = 1}^n \,
\xi_{ ki} ( x) \, \frac{ \partial }{\partial x_i}$ 
in the sense that the totality of the transformations $x' = f ( x; \,
a)$ is identical with the totality of all transformations:
\[
\aligned
x_i'
&
=
\exp\big(
\lambda_1\,X_1
+\cdots+
\lambda_r\,X_r\big)(x_i)
\\
&
=
x_i
+
\sum_{k=1}^r\,\lambda_k\,\xi_{ki}(x)
+
\sum_{k,\,j}^{1\dots r}\,
\frac{\lambda_k\,\lambda_j}{1\cdot 2}\,
X_k(\xi_{ji})
+
\cdots
\ \ \ \ \ \ \ \ \ \ \ \ \
{\scriptstyle{(i\,=\,1\,\cdots\,n)}}
\endaligned
\]
obtained as the time-one map of the one-term group $\exp \big( t
\sum\, \lambda_i X_i \big) ( x)$ generated by the general linear
combination of the infinitesimal transformations.
\newline\indent
A beautiful idea of analyzing the (diagonal) action ${x^{( \mu)}}' = f
\big( x^{ ( \mu)}; \, a\big)$ induced on $r$-tuples of points $\big(
x^{(1)}, \dots, x^{ ( r)} \big)$ in general position enables Lie to
show that for every collection of $r$ linearly independent vector
fields $X_k = \sum_{ i = 1}^n\, \xi_{ ki} ( x) \, \frac{ \partial }{
\partial x_i}$, the parameters $\lambda_1, \dots, \lambda_r$ in the
finite transformation equations $x' = \exp \big( \lambda_1 \, X_1 +
\cdots + \lambda_r \, X_r \big) ( x)$ are all essential.
}

\section{Mechanical and Mental Images}
\label{mechanical-mental}

\begin{svgraybox}
\centerline{{\sf \S\,\,\,15. (\cite{ enlie1888-4})}}
The concept of infinitesimal transformation and likewise the one of
one-term group gain a certain graphic nature when one makes use of
geometric and mechanical images.

In the infinitesimal transformation:
\[
x_i'
=
x_i+\xi_i(x_1,\dots,x_n)\,\delta t
\ \ \ \ \ \ \ \ \ \ \ \ \
{\scriptstyle{(i\,=\,1\,\cdots\,n)}},
\]
or in:
\[
X(f)
=
\sum_{i=1}^n\,\xi_i\,
\frac{\partial f}{\partial x_i},
\]
we interpret the variables $x$ as Cartesian coordinates of an
$n$-times extended space. Then the transformation can obviously be
interpreted in such a way that each point of coordinates $x_1, \dots,
x_n$ is transferred to an infinitesimally neighbouring point of
coordinates $x_1', \dots, x_n'$. {\em So the transformation attaches
to every point at which not all the $\xi_i$ vanish a certain direction
of progress \label{fortsch}
\deutsch{Fortschreitungsrichtung}, and
along this direction of progress, a certain infinitely small line}
\deutsch{Strecke}; the direction of progress is determined by the
proportion: $\xi_1 ( x) \colon \xi_2 ( x) \colon
\cdots \colon \xi_n ( x)$; the infinitely small
line has the length: $\sqrt{ \xi_1^2 + \cdots + \xi_n^2} \, \delta t$.
If at a point $x_1, \dots, x_n$ all the $\xi$ are zero, then no
direction of progress is attached to the point by the infinitesimal
transformation.

If we imagine that the whole space is filled with a compressible
fluid, then we can interpret the infinitesimal transformation $X ( f)$
simply as an infinitely small movement \deutsch{Bewegung} of this
fluid, and $\delta t$ as the infinitely small time interval
\deutsch{Zeitabschnitt} during which this movement proceeds. 
Then evidently, the quantities $\xi_1 ( x) \colon \cdots \colon \xi_n
( x)$ are the components of the velocity of the fluid particle
\deutsch{Flüssigkeitstheilchens}
which is located precisely in the point $x_1, \dots, x_n$.

The finite transformations of the one-term group $X ( f)$ come into
being \deutsch{entstehen} by repating infinitely many times
\deutsch{durch unendlichmalige Wiederholung} 
the infinitesimal transformation $x_i' = x_i + \xi_i \, \delta t$. In
order to arrive at a finite transformation, we must hence imagine that
the infinitely small movement represented by the infinitesimal
transformation is repeated during infinitely many time intervals
$\delta t$; in other words, we must follow the movement of the fluid
particle during a finite time interval. To this end, it is necessary
to integrate the differential equations of this movement, that is to
say the simultaneous system:
\[
\frac{\D\,x_1'}{\D\,\xi_1(x')}
=\cdots=
\frac{\D\,x_n'}{\xi_n(x')}
=
\D\,t.
\]
A fluid particle which, for the time $t = 0$, is located in the point
$x_1, \dots, x_n$ will, after a lapse of time $t$, reach the point
$x_1', \dots, x_n'$; so the integration must be executed in such a way
that for $t = 0$, one has: $x_i' = x_i$. In reality, we found earlier
the general form of a finite transformation of our group exactly in
this manner.

The movement of our fluid is one which is a so-called stationary
movement, because the velocity components $\D\, x_i' / \D\, 
t$ are free of
$t$. From this, it follows that in one and the same point, the
movement is the same at each time, and consequently that the whole
process of movement always takes the same course, whatever the time at
which one starts to consider it. The proof that the finite
transformations produced by the infinitesimal transformation $X ( f)$
constitute a group already lies fundamentally in this observation.

If we consider a determined fluid particle, say the one which, for the
time $t = 0$ is located in the point $x_1^0, \dots, x_n^0$, then we
see that it moves on a curve which passes through the point $x_1^0,
\dots, x_n^0$. So the entire space is decomposed only in curves of
a constitution such that each particle remains on the curve on which
it is located once. We want to call these curves the
\terminology{integral curves
\deutsch{Bahncurven} of the 
infinitesimal transformation $X(f)$}. Obviously, there are $\infty^{
n-1}$ such integral curves.

For a given point $x_1^0, \dots, x_n^0$, it is easy to display the
equations of the integral curve passing through it.
The \emphasis{integral curve is nothing but the locus of
all points in which the point $x_1^0, \dots, x_n^0$ is
transferred, when the $\infty^1$ transformations:
\[
x_i'
=
x_i
+
\frac{t}{1}\,
\xi_i(x)
+\cdots
\ \ \ \ \ \ \ \ \ \ \ \ \
{\scriptstyle{(i\,=\,1\,\cdots\,n)}}
\]
of the one-term group are executed on it}.
Consequently, the equations:
\[
x_i'
=
x_i^0
+
\frac{t}{1}\,\xi_i(x^0)
+
\frac{t^2}{1\cdot 2}\,
\big(
X(\xi_i)
\big)_{x=x_0}
+
\cdots
\ \ \ \ \ \ \ \ \ \ \ \ \
{\scriptstyle{(i\,=\,1\,\cdots\,n)}}
\]
represent the integral curve in question, when
one considers $t$ as an independent variable.
\end{svgraybox}

\section{Straightening of Flows and the Exponential Formula}
\label{straightening-flows}

Let $x_i' = f_i ( x; \, a)$ be a \emphasis{one-term} local
transformation group, with $a \in \K$ ($= \R$ or $\C$) a scalar and
with the identity $e$ corresponding to the origin $0 \in \K$ as
usual. Its fundamental differential equations
(p.~\pageref{fundamental-differential-equations}):
\[
\frac{\D\,x_i'}{\D\,a}
=
\psi(a)\,\xi_i(x_1',\dots,x_n')
\ \ \ \ \ \ \ \ \ \ \ \ \
{\scriptstyle{(i\,=\,1\,\cdots\,n)}}
\]
then consist of a complete first order
\pde{pde} system. But
by introducing the new parameter:
\[
t
=
t(a)
:=
\int_0^a\,
\psi(a_1)\,\D\,a_1,
\]
we immediately transfer these fundamental differential equations to
the \emphasis{time-independent} system of $n$ 
ordinary differential equations:
\def\theequation{1}\begin{equation}
\frac{\D\,x_i'}{\D \,t}
=
\xi_i(x_1',\dots,x_n')
\ \ \ \ \ \ \ \ \ \ \ \ \
{\scriptstyle{(i\,=\,1\,\cdots\,n)}},
\end{equation}
the integration of which amounts to computing the so-called
\emphasis{flow} of the vector field $X := \sum_{ i = 1}^n\, \xi_i ( x
) \, \frac{ \partial }{ \partial x_i }$. 

With the same letters $f_i$, we will write $f_i ( x; \, t)$ instead of
$f_i ( x; \, a(t))$. Of course, the (unique) solution of the
system~\thetag{ 1} with the initial condition $x_i' ( x; \, 0) = x_i$
is nothing but $x_i' = f_i ( x; \, t)$: the flow was in fact known
from the beginning. Furthermore, uniqueness of the flow and the fact
that the $\xi_i$ are independent of $t$ both imply that the group
composition property corresponds just to addition of time 
parameters (\cite{ arno1978, rao1981}):
\[
f_i\big(f(x;\,t_1);\,t_2)
\equiv
f_i(x;\,t_1+t_2)
\ \ \ \ \ \ \ \ \ \ \ \ \
{\scriptstyle{(i\,=\,1\,\cdots\,n)}}.
\]
It is classical that one may (locally) straighten $X$ to $\frac{
\partial }{\partial y_n }$.

\begin{theorem}
\label{theorem-straightening}
Every one-term continuous transformation group:
\[
x_i'
=
f_i(x_1,\dots,x_n;\,t)
\ \ \ \ \ \ \ \ \ \ \ \ \
{\scriptstyle{(i\,=\,1\,\cdots\,n)}}
\]
satisfying differential equations of the form:
\[
\frac{\D\,f_i}{\D \,t}
=
\xi_i(f_1,\dots,f_n)
\ \ \ \ \ \ \ \ \ \ \ \ \
{\scriptstyle{(i\,=\,1\,\cdots\,n)}}
\]
is locally equivalent, through a suitable change of variables $y_i =
y_i ( x)$, to a group of translations:
\[
y_1'
=
y_1+t,\dots\dots,\ \ \ 
y_2'
=
y_2,\ \ \ \ \ \ 
y_n'
=
y_n.
\]
\end{theorem}

\begin{center}
\input flow-straightening.pstex_t
\end{center}

\begin{proof}\smartqed
We may suppose from the beginning that the coordinates $x_1, \dots,
x_n$ had been chosen so that $\xi_1 ( 0) = 1$ and $\xi_2 ( 0 ) =
\cdots = \xi_n ( 0) = 0$. In an auxiliary space $y_1, \dots, y_n$
drawn on the left of the figure, we consider all points 
$\big( 0, \widehat{ y} \,\big)
:= (0, y_2, \dots, y_n)$ near the origin lying on the coordinate
hyperplane which is complementary to the $y_1$-axis, and we introduce
the diffeomorphism:
\[
y
\longmapsto
x
=
x(y)
:=
f\big(0,\widehat{y};\,y_1\big)
=:
\Phi(y),
\]
defined by following the flow up to time $y_1$ from
$\big( 0, \widehat{ y} \,\big)$; this is indeed a diffeomorphism
fixing the origin thanks to $\frac{ \partial \Phi_1}{ \partial y_1 } (
0) = \frac{ \partial f_1 }{ \partial t} ( 0) = \xi_1 ( 0) = 1$, to
$\frac{ \partial \Phi_k }{ \partial y_1} ( 0) = \frac{ \partial f_k}{
\partial t} ( 0) = \xi_k(0) = 0$ and to $\Phi_k (0, \widehat{ y })
\equiv \widehat{ y }_k$ for $k = 2, \dots, n$. Consequently, we get
that the (wavy) flow represented on the right figure side has been
straightened on the left side to be just a uniform translation 
directed by the $y_1$-axis, because by substituting:
\[
x'
=
f\big(
0,\widehat{y}';\,y_1'\big)
=
f\big(
0,\widehat{y};\,y_1+t\big)
=
f\big(f\big(
0,\widehat{y};\,y_1
\big);\,t\big)
=
f(x;\,t),
\]
we recover the uniquely defined flow $x' = f ( x; \, t)$ when assuming
that $\widehat{ y}' = \widehat{ y}$ and $y_1' = y_1 + t$.
\qed\end{proof}

\begin{svgraybox}\vspace{-0.5cm}
\def\thetheorem{6}\begin{theorem}
\label{Theorem-6-S-49}
If a one-term group:
\[
x_i'
=
f_i(x_1,\dots,x_n,\,a_1,\dots,a_r)
\ \ \ \ \ \ \ \ \ \ \ \ \ {\scriptstyle{(i\,=\,1\,\cdots\,n)}}
\]
contains the identity transformation, then its transformations are
interchangeable \deutsch{vertauschbar} one with another and they
can be ordered as inverses by pairs. Every group 
of this sort is equivalent to a group of
translations:
\[
y_1'=y_1,
\,\,\,\dots,\,\,\,
y_{n-1}'=y_{n-1},\,\,\,
y_n'=y_n+t.
\]
\end{theorem}\vspace{-0.4cm}
\end{svgraybox}

\subsection{The Exponential Analytic Flow Formula} 

In fact, because we
have universally assumed analyticity
of all data, the solution $x' = x' ( x; \,
t)$ to the \pde{pde} system~\thetag{ 1} can be sought by expanding the
unknown $x'$ in power series with respect to $t$:
\[
x_i'(x;\,t)
=
\sum_{k\geqslant 0}\,\Xi_{ik}(x)\,t^k
=
x_i
+
t\,\xi_i(x)
+\cdots
\ \ \ \ \ \ \ \ \ \ \ \ \
{\scriptstyle{(i\,=\,1\,\cdots\,n)}}.
\]
So, what are the coefficient functions $\Xi_{ ik} ( x)$?
Differentiating once more and twice more~\thetag{ 1}, we get for
instance:
\[
\aligned
\frac{\D^2x_i'}{\D\,t^2}
&
=
\sum_{k=1}^n\,\frac{\partial\xi_i}{\partial x_k}\,
\frac{\D\,x_k'}{\D\,t}
=
\sum_{k=1}^n\,\frac{\partial\xi_i}{\partial x_k}\,\xi_k
=
X(\xi_i)
\ \ \ \ \ \ \ \ \ \ \ \ \
{\scriptstyle{(i\,=\,1\,\cdots\,n)}}
\\
\frac{\D^3x_i'}{\D\,t^3}
&
=
\sum_{k=1}^n\,\frac{\partial X(\xi_i)}{\partial x_k}\,
\frac{\D\,x_k'}{\D\,t}
=
\sum_{k=1}^n\,\frac{\partial X(\xi_i)}{\partial x_k}\,
\xi_k
=
X\big(X(\xi_i)\big),
\endaligned
\]
etc., and hence generally by a straightforward induction:
\[
\frac{\D^kx_i'}{\D\,t^k}
=
\underbrace{X\big(\cdots\big(
X}_{k-1\,\,\text{\rm times}}(\xi_i)\big)\cdots\big)
=
\underbrace{X\big(\cdots\big(
X\big(X}_{k\,\,\text{\rm times}}(x_i')\big)\big)\cdots\big),
\]
for every nonnegative integer $k$, with the convention $X^0 x_i =
x_i$. Setting $t = 0$, we therefore get:
\[
k!\,\,\Xi_{ik}(x)
\equiv
\underbrace{X\big(\cdots\big(
X\big(X}_{k\,\,\text{\rm times}}(x_i)\big)\big)\cdots\big).
\]
Thus rather strikingly, computing a flow boils down in the
analytic category to summing up differentiated terms.

\begin{proposition}
The unique solution $x' ( x; \, t)$ to a local analytic system of
ordinary differential equations $\frac{ \D\,x_i'}{ \D\,t} = \xi_i ( x_1',
\dots, x_n')$ with initial condition $x_i' ( x; \, 0) = x_i$ is
provided by the power series expansion:
\def\theequation{2}\begin{equation}
x_i'(x;\,t)
=
x_i
+
t\,X(x_i)
+\cdots+
\frac{t^k}{k!}\,
\underbrace{X\big(\cdots\big(
X\big(X}_{k\,\,\text{\rm times}}(x_i)\big)\big)\cdots\big)
+\cdots
\ \ \ \ \ \ \ \ \ \ \ \ \
{\scriptstyle{(i\,=\,1\,\cdots\,n)}},
\end{equation}
which can also be written, quite adequately, by means of an
exponential denotation:
\def\theequation{2'}\begin{equation}
\label{flow-exponential-formula}
\aligned
x_i'
=
\exp\big(t\,X\big)(x_i)
=
\sum_{k\geqslant 0}\,\frac{(t\,X)^k}{k!}(x_i)
\ \ \ \ \ \ \ \ \ \ \ \ \
{\scriptstyle{(i\,=\,1\,\cdots\,n)}}.
\endaligned
\end{equation}
\end{proposition}

\subsection{Action on Functions} 

Letting $f = f ( x_1, \dots, x_n)$
be an arbitrary analytic function, we 
might compose $f$ with the above flow:
\[
f'
:=
f(x_1',\dots,x_n')
=
f\big(x_1'(x;\,t),\dots,x_n'(x;\,t)\big),
\]
and we should then expand the
result in power series with respect to $t$:
\[
f'
=
\big(f'\big)_{t=0}
+
{\textstyle{\frac{t}{1!}}}\,
\big(
{\textstyle{\frac{\D\,f'}{\D\,t}}}
\big)_{t=0}
+
{\textstyle{\frac{t^2}{2!}}}\,
\big(
{\textstyle{\frac{\D^2f'}{\D\,t^2}}}
\big)_{t=0}
+\cdots.
\]
Consequently, we need to 
compute the differential quotients $\frac{ \D\,f'}{
\D\,t}$, $\frac{ \D^2 f'}{ \D\,t^2}$, \dots, and
if we set $\xi_i' := \xi_i ( x_1', \dots, x_n')$ and 
$X' := \sum_{ i=1}^n \, \xi_i' \frac{ \partial }{ \partial x_i'}$, 
\[
\aligned
\frac{\D\,f'}{\D\,t}
&
=
\sum_{i=1}^n\,\xi_i'\,
\frac{\partial f'}{\partial x_i'}
=
X'(f'),
\\
\frac{\D^2f'}{\D\,t^2}
&
=
X'\bigg(
\sum_{i=1}^n\,\xi_i'\,
\frac{\partial f'}{\partial x_i'}
\bigg)
=
X'\big(
X'(f')\big),
\endaligned
\]
and so on. After setting $t = 0$, the $x_i'$ become $x_i$, the $f'$
becomes $f$, the $X' ( f')$ becomes $X (f)$, and so on, whence we
obtain the expansion\footnote{\,
Changing $t$ to $-t$ exchanges the rôles of
$x' = \exp (tX) (x)$ and of 
$x = \exp ( -t X) (x')$, hence we also have: 
\def\theequation{3a}\begin{equation}
\label{t-3a}
f(x_1,\dots,x_n)
=
f(x_1',\dots,x_n')
-
\frac{t}{1!}\,X(f)
+
\cdots
+
(-1)^k\,
\frac{t^k}{k!}\,
X\big(\cdots\big(X
(f)\big)\cdots\big)
+\cdots.
\end{equation}
}: 
\def\theequation{3}\begin{equation}
\label{S-52-eq-7}
f(x_1',\dots,x_n')
=
f(x_1,\dots,x_n)
+
\frac{t}{1!}\,X(f)
+
\cdots
+
\frac{t^k}{k!}\,
\underbrace{X\big(\cdots\big(X}_{k\,\,\text{\rm times}}
(f)\big)\cdots\big)
+\cdots.
\end{equation}

\begin{svgraybox}
\centerline{{\sf \S\,\,\,13. (\cite{ enlie1888-4})}}

\smallskip
Amongst the $\infty^1$ transformations of the one-term group~\thetag{
2}, those whose parameter $t$ has an infinitely small
value, say the value $\delta t$, play an outstanding role. We now want
to consider more precisely these ``\terminology{infinitely small}'' or
``\terminology{infinitesimal}'' transformations of the group.
 
If we take into account only the first power of $\delta t$, whereas we
ignore the second and all the higher ones, then we obtain
from~\thetag{ 2} the desired infinitesimal transformation under the
form:
\def\theequation{4}\begin{equation}
x_i'
=
x_i
+
\xi_i(x_1,\dots,x_n)\,\delta t
\ \ \ \ \ \ \ \ \ \ \ \ \
{\scriptstyle{(i\,=\,1\,\cdots\,n)}};
\end{equation}
on the other hand, if we use the equation~\thetag{ 3}, we get the
last $n$ equations condensed in the single one:
\[
f'
=
f
+X(f)\,\delta t,
\]
or written in greater length:
\[
f(x_1',\dots,x_n')
=
f(x_1,\dots,x_n)
+
\delta t\,
\sum_{i=1}^n\,\xi_i\,
\frac{\partial f}{\partial x_i}.
\]
It is convenient to introduce a specific naming 
for the difference $x_i' - x_i$, that is
to say, for the expression $\xi_i \, \delta t$.
Occasionally, we want to call $\xi_i \delta t$ the
``\terminology{increase}'' \deutsch{Zuwachs}, 
\label{S-53} or the
``\terminology{increment}''
\deutsch{Increment}, or also the ``\terminology{variation}'' 
\deutsch{Variation}
of $x_i$, and write for that: $\delta x_i$. Then we can also represent
the infinitesimal transformation under the form:
\[
\delta x_1
=
\xi_1\,\delta t,\dots,\,
\delta x_n
=
\xi_n\,\delta t.
\]

Correspondingly, we will call the difference $f' - f$, or the
expression $X ( f) \, \delta t$ the 
\terminology{increase}, or the \terminology{variation} of the
function $f ( x_1, \dots, x_n)$ and we shall write:
\[
f'
-
f
=
X(f)\,\delta t
=
\delta f.
\]

It stands to reason that the expression, completely alone:
\[
X(f)
=
\sum_{i=1}^n\,\xi_i\,
\frac{\partial f}{\partial x_i}
\]
already fully determines the infinitesimal transformation $\delta x_i
= \xi_i \, \delta t$, when one understands by $f ( x_1, \dots, x_n)$
some undetermined function of its arguments. Indeed, all the $n$
functions $\xi_1, \dots, \xi_n$ are individually given at the same
time with $X (f)$.

\smallskip

\emphasis{This is why we shall introduce the expression $X ( f) = \frac{
\delta f}{ \delta t}$ as being the symbol of the infinitesimal
transformation~\thetag{ 4}}, so we will really speak of the 
``\terminology{infinitesimal transformation $X ( f )$}''. 
However, we want to point
out just now that the symbol of the infinitesimal
transformation~\thetag{ 4} is basically determined only up to an
arbitrary remaining constant factor. In fact, when we multiply the
expression $X ( f)$ by any finite constant $c$, then the resulting
expression $c \, X (f)$ is also to be considered as the symbol of the
infinitesimal transformation~\thetag{ 4}. Indeed, according to the
concept of an infinitely small quantity, it makes no difference,
when we substitute in~\thetag{ 4} the infinitely small quantity
$\delta t$ by $c \, \delta t$.

The introduction of the symbol $X ( f)$ for the infinitesimal
transformation~\thetag{ 4} presents many advantages. \emphasis{Firstly},
it is very convenient that the $n$ equations $x_i' = x_i + \xi_i \,
\delta t$ of the transformation are replaced by the single expression
$X ( f)$. \emphasis{Secondly}, it is convenient that in the symbol $X
(f)$, we have to deal with only \emphasis{one} series of variables, not
with the two series: $x_1, \dots, x_n$ and $x_1', \dots, x_n'$.
Lastly and \emphasis{thirdly}, the symbol $X (f)$ establishes the 
connection
between infinitesimal transformations and linear partial differential
\label{cited-S-82}
equations; because in the latter theory, expressions such as $X (f)$
do indeed play an important role. We shall go into more details later
about this connection (cf. Chap.~6).

The preceding developments show that a one-term group with the 
identity
transformation always comprises a well determined infinitesimal
transformation $x_i' = x_i + \xi_i \, \delta t$, or briefly $X ( f)$.
But it is also clear that conversely, the one-term group in question
is perfectly determined, as soon as one knows its infinitesimal
transformation. Indeed, the infinitesimal transformation $x_i' = x_i
+ \xi_i \, \delta t$ is, so to speak, only another way of writing the
simultaneous system~\thetag{ 1} from which are derived the
equations~\thetag{ 2} of the one-term group.

Thus, since every one-term group with the identity transformation
is completely determined by its infinitesimal transformation, 
we shall, for the sake of convenience, introduce
the following way of expressing.

\plainstatement{Every transformation of the one-term group:
\[
x_i'
=
x_i
+
\frac{t}{1}\,\xi_i
+
\frac{t^2}{1\cdot 2}\,X(\xi_i)
+\cdots
\ \ \ \ \ \ \ \ \ \
{\scriptstyle{(i\,=\,1\,\cdots\,n)}}
\]
is obtained by repeating infinitely many times
the infinitesimal transformation}:
\[
x_i'
=
x_i+\xi_i\,\delta t
\ \ \ \ \
\text{\rm or}
\ \ \ \ \
X(f)
=
\xi_1\,\frac{\partial f}{\partial x_1}
+\cdots+
\xi_n\,\frac{\partial f}{\partial x_n}.
\]

\noindent
Or yet more briefly:

\plainstatement{The one-term group in question is 
generated by its infinitesimal
transformations. }

In contrast to the infinitesimal transformation $X ( f)$, we call the
equations:
\[
x_i'
=
x_i
+
\frac{t}{1}\,\xi_i
+
\frac{t^2}{1\cdot 2}\,X(\xi_i)
+\cdots
\]
the \terminology{finite} equations of the one-term group in question.

Now, we may enunciate briefly as follows the connection found
earlier between the simultaneous system~\thetag{ 1} and the one-term
group~\thetag{ 2'}:

\def\theproposition{1}\begin{proposition}
Every one-term group which contains the identity transformation
is generated by a well determined infinitesimal transformation.
\end{proposition}

\smallskip
And conversely:

\def\theproposition{2}\begin{proposition}
Every infinitesimal transformation generates a completely 
determined one-term group.
\end{proposition}
\vspace{-0.3cm}
\end{svgraybox}

\section{Exponential Change of Coordinates and Lie Bracket}

Let $X = \sum_{ i=1}^n\, \xi_i ( x) \, \frac{ \partial }{\partial
x_i}$ be an infinitesimal transformation which generates the one-term
group $x' = \exp ( t X) ( x)$. What happens if the variables $x_1,
\dots, x_n$ are subjected to an analytic diffeomorphism $y_\nu =
\varphi_\nu ( x)$ which transfers naturally $X$ to the vector field:
\[
Y
:=
\varphi_*(X)
=
\sum_{\nu=1}^n\,X(y_\nu)\,
\frac{\partial}{\partial y_\nu}
=:
\sum_{\nu=1}^n\,
\eta_\nu(y_1,\dots,y_n)\,
\frac{\partial}{\partial y_\nu}
\]
in the new variables $y_1, \dots, y_n$?

\begin{proposition}
\label{S-57}
The new one-term group $y' = \exp ( tY) ( y)$ associated to $Y =
\varphi_* ( X)$ can be recovered from the old one $x' = \exp ( t X) (
x)$ thanks to the formula:
\[
\exp(tY)(y)
=
\varphi\big(
\exp(tX)(x)\big)
\big\vert_{x=\varphi^{-1}(y)}.
\]
\end{proposition}

\begin{proof}\smartqed
Since we work in the analytic category, we are allowed to deal with
power series expansions. Through the introduction of the new variables
$y_\nu = \varphi_\nu ( x)$, an arbitrary function $f ( x_1, \dots,
x_n)$ is transformed to the function $F = F ( y)$ defined by the
identity:
\[
f(x)
\equiv
F\big(\varphi(x)\big).
\]
With $x' = \exp ( tX) ( x)$ by assumption, we may also 
define $y_\nu' := \varphi_\nu ( x')$ so that:
\[
F(y')\big\vert_{y'=\varphi(x')}
=
f(x').
\]
On the other hand, the Jacobian matrix of $\varphi$ induces a
transformation between vector fields; equivalently, this
transformation $X \mapsto \varphi_* ( X) =: Y$ can be defined by the
requirement that for any function $f$:
\[
Y(F)\big\vert_{y=\varphi(x)}
=
X(f).
\] 
By a straightforward induction, it follows for any integer $k\geqslant
1$ that we have:
\[
Y\big(Y(F)\big)
\big\vert_{y=\varphi(x)}
=
X\big(X(f)\big),
\dots\dots,\,\,
Y^k(F)
\big\vert_{y=\varphi(x)}
=
X^k(f).
\]
Consequently, in the expansion~\thetag{ 3} of $f ( x') = f \big( \exp
( tX) (x)\big)$ with respect to the powers of $t$, we may perform
replacements:
\[
F(y')\big\vert_{y'=\varphi(x')}
=
f(x')
=
\sum_{k\geqslant 0}\,
\frac{t^k}{k!}\,
X^k(f)
=
\sum_{k\geqslant 0}\,
\frac{t^k}{k!}\,
Y^k(F)\big\vert_{y=\varphi(x)}.
\]
Removing the two replacements 
$\vert_{y'=\varphi(x')}$ and $\vert_{y=\varphi(x)}$, 
we get an identity in terms of the variables $y$ and $y'$:
\[
F(y')
=
\sum_{k\geqslant 0}\,
\frac{t^k}{k!}\,Y^k(F),
\]
which, because the function $f$\,---\,and hence $F$ 
too\,---\,was arbitrary,
shows that $y'$ must coincide with $\exp ( t Y) ( y)$, namely:
\[
y_i'
=
\sum_{k\geqslant 0}\,
\frac{t^k}{k!}\,Y^k(y_i)
\ \ \ \ \ \ \ \ \ \ \ \ \
{\scriptstyle{(i\,=\,1\,\cdots\,n)}}.
\]
Consequently, after replacing here $y'$ by $\varphi ( x')$ and $y$ by
$\varphi ( x)$, we obtain:
\[
\varphi\big(\exp(tX)(x)\big)
=
\varphi(x')
=
y'
=
\sum_{k\geqslant 0}\,
\frac{t^k}{k!}\,
Y^k(y)\big\vert_{y=\varphi(x)}
=
\exp(tY)(y)\big\vert_{y=\varphi(x)},
\]
and this is the same relation between flows as the one stated in
the proposition, but viewed in the $x$-space.
\qed\end{proof}

\begin{svgraybox}\vspace{-0.4cm}
\def\theproposition{3}\begin{proposition}
\label{S-58}
If, after the introduction of the new independent variables:
\[
y_i
=
\varphi_i(x_1,\dots,x_n)
\ \ \ \ \ \ \ \ \ \ \ \ \ {\scriptstyle{(i\,=\,1\,\cdots\,n)}},
\]
the symbol:
\[
X(f)
=
\sum_{i=1}^n\,\xi_i\,
\frac{\partial f}{\partial x_i}
\]
of the infinitesimal transformation:
\[
x_i'
=
x_i
+
\xi_i\,\delta t
\ \ \ \ \ \ \ \ \ \ \ \ \ {\scriptstyle{(i\,=\,1\,\cdots\,n)}}
\]
receives the form:
\[
X(f)
=
\sum_{\nu=1}^n\,X(y_\nu)\,
\frac{\partial f}{\partial y_\nu}
=
\sum_{\nu=1}^n\,
\eta_\nu(y_1,\dots,y_n)
=
Y(f),
\]
then the finite transformations generated by $X(f)$:
\[
x_i'
=
x_i
+
\frac{t}{1}\,\xi_i
+
\frac{t^2}{1\cdot 2}\,X(\xi_i)
+\cdots
\ \ \ \ \ \ \ \ \ \ \ \ \ {\scriptstyle{(i\,=\,1\,\cdots\,n)}}
\]
are given, in the new variables, the shape:
\[
y_i'
=
y_i
+
\frac{t}{1}\,\eta_i
+
\frac{t^2}{1\cdot 2}\,Y(\eta_i)
+\cdots
\ \ \ \ \ \ \ \ \ \ \ \ \ {\scriptstyle{(i\,=\,1\,\cdots\,n)}},
\]
where the parameter $t$ possesses the same value
in both cases.
\end{proposition}\vspace{-0.3cm}
\end{svgraybox}

\subsection{Flows as Changes of Coordinates} 

Suppose we are given an arbitrary infinitesimal transformations:
\[
X
=
\sum_{i=1}^n\,\xi_i(x)\,
\frac{\partial}{\partial x_i}.
\]
For later use, we want to know how $X$ transforms through the change
of coordinates represented by the flow diffeomorphism:
\[
x_i'
=
\exp(tY)(x_i)
=
x_i
+
\frac{t}{1}\,\eta_i
+
\frac{t^2}{1\cdot 2}\,Y(\eta_i)
+
\cdots
\ \ \ \ \ \ \ \ \ \ \ \ \
{\scriptstyle{(i\,=\,1\,\cdots\,n)}}
\]
of {\em another} infinitesimal transformation:
\[
Y
=
\sum_{i=1}^n\,
\eta_i(x_1,\dots,x_n)\,
\frac{\partial}{\partial x_i}.
\]
At least, we would like to control how $X$ transforms modulo second
order terms in $t$. By definition, $X$ is transformed to the vector
field:
\[
\varphi_*(X)
=
\sum_{i=1}^n\,X(x_i')\,
\frac{\partial}{\partial x_i'},
\]
where the expressions $X ( x_i')$ should still be expressed in terms
of the target coordinates $x_1', \dots, x_n'$. Replacing $x_i'$ by its
above expansion and neglecting the second and the higher powers of
$t$, we thus get:
\[
X(x_i')
=
X(x_i)
+
t\,X(\eta_i)
+
\cdots
\ \ \ \ \ \ \ \ \ \ \ \ \
{\scriptstyle{(i\,=\,1\,\cdots\,n)}}.
\]
To express the first term of the right-hand side in terms of the
target coordinates, we begin by inverting $x' = \exp ( tY ) ( x)$,
getting $x = \exp ( -t Y') (x')$, with of course $Y' := \sum_{
i=1}^n\, \eta_i ( x') \, \frac{ \partial }{\partial y_i'}$ denoting
the same field as $Y$ but written in the $x'$-space, and we compute
$X(x_i)$ as follows:
\[
\aligned
X(x_i)
=
\xi_i(x)
&
=
\xi_i\big(
\exp(-tY')(x')
\big)
\\
&
=
\xi_i
\big(
x_1'-t\,\eta_1(x'),\dots\dots,x_n'-t\,\eta_n(x')
\big)
\\
&
=
\xi_i(x')
-
t\,\sum_{l=1}^n\,
\frac{\partial\xi_i}{\partial x_l'}(x')\,
\eta_l(x')
+\cdots
\\
&
=
\xi_i(x')
-
t\,Y'\big(\xi_i(x')\big)
+\cdots
\ \ \ \ \ \ \ \ \ \ \ \ \
{\scriptstyle{(i\,=\,1\,\cdots\,n)}}.
\endaligned
\]
If we now abbreviate $\xi_i' := \xi_i ( x')$ and $\eta_i' :=
\eta_i ( x')$, because ${\rm O} (t^2)$ is neglected,
no computation is needed to express the second term in terms of $x'$:
\[
t\,X(\eta_i)
=
t\,X'(\eta_i')
-
\cdots
\ \ \ \ \ \ \ \ \ \ \ \ \
{\scriptstyle{(i\,=\,1\,\cdots\,n)}},
\]
and consequently by adding we get:
\[
X(x_i')
=
\xi_i'
+
t\,\big(
X'(\eta_i')
-
Y'(\xi_i')
\big)
+\cdots
\ \ \ \ \ \ \ \ \ \ \ \ \
{\scriptstyle{(i\,=\,1\,\cdots\,n)}}.
\]
As the $n$ coefficients $X' ( \eta_i') - Y' ( \xi_i')$, either one
\emphasis{recognizes} the coefficients of the
\terminology{Lie bracket}:
\[
\big[X',\,Y'\big]
=
\bigg[
\sum_{i=1}^n\,\xi_i(x')\,\frac{\partial}{\partial x_i'},\,\,
\sum_{i=1}^n\,\eta_i(x')\,\frac{\partial}{\partial x_i'}
\bigg]
:=
\sum_{i=1}^n\,
\Big(
X'(\eta_i')-Y'(\xi_i')
\Big)
\frac{\partial}{\partial x_i'},
\]
between $X' := \sum_{ i=1}^n\, \xi_i ( x') \, \frac{ \partial
}{\partial x_i'}$ and $Y' := \sum_{ i = 1}^n \,
\eta_i ( x') \, \frac{ \partial }{\partial x_i'}$,
or one chooses to \emphasis{define} once for all the Lie bracket in
such a way. In fact, Engel and Lie will mainly introduce brackets in
the context of the Clebsch-Frobenius theorem, Chap.~\ref{kapitel-5}
below. At present, let us state the gained result in a self-contained
manner.

\begin{lemma}
\label{S-141}
{\rm (\cite{enlie1888-4}, p.~141)}
If, in the infinitesimal transformation $X = \sum_{ i=1}^n\, \xi_i
( x) \, \frac{ \partial }{\partial x_i}$ of
the space $x_1, \dots, x_n$, one introduces as new variables:
\[
x_i'
=
\exp(tY)(x_i)
=
x_i
+
t\,Y(x_i)
+
\cdots
\ \ \ \ \ \ \ \ \ \ \ \ \
{\scriptstyle{(i\,=\,1\,\cdots\,n)}},
\]
those
induced by the one-term group generated by another infinitesimal
transformation $Y = \sum_{ i = 1}^n \, \eta_i ( x) \, \frac{ \partial
}{\partial x_i}$, then setting $X' := \sum_{ i=1}^n\, \xi_i (
x') \, \frac{ \partial }{\partial x_i'}$ and $Y' := \sum_{ i = 1}^n \,
\eta_i ( x') \, \frac{ \partial }{\partial x_i'}$, one obtains 
a transformed vector field:
\[
\varphi_*(X)
=
X'
+
t\,
\big[
X',\,Y'
\big]
+\cdots,
\]
with a first order perturbation which
is the Lie bracket $\big[ X', \, Y'\big]$.
\end{lemma}

\section{Essentiality of Multiple Flow Parameters}

At present, we again consider $r$ arbitrary vector fields with
analytic coefficients defined on a certain, unnamed domain of $\K^n$
which contains the origin:
\[
X_k
=
\sum_{i=1}^n\,\xi_{ki}(x_1,\dots,x_n)\,
\frac{\partial}{\partial x_i}
\ \ \ \ \ \ \ \ \ \ \ \ \
{\scriptstyle{(k\,=\,1\,\cdots\,r)}}.
\]
Though the collection $X_1, \dots, X_k$ does not necessarily stem from
an $r$-term group, each individual $X_k$ nonetheless generates
the one-term continuous transformation group $x' = \exp ( t X_k) ( x)$
with corresponding infinitesimal transformations $x_i' = x_i +
\varepsilon\, \xi_{ ki} ( x)$; this is the reason why we shall hence
sometimes alternatively 
refer to such general vector fields $X_k$ as being
\terminology{infinitesimal transformations}.

\begin{definition}
The $r$ infinitesimal transformations $X_1, \dots, X_r$ will be
called \terminology{independent} (\terminology{of each other}) 
if they are linearly independent, namely if the $n$ equations:
\[
0
\equiv
e_1\,\xi_{1i}(x)
+\cdots+
e_r\,\xi_{ri}(x)
\ \ \ \ \ \ \ \ \ \ \ \ \
{\scriptstyle{(i\,=\,1\,\cdots\,n)}}
\]
in which $e_1, \dots, e_r$ are \emphasis{constants}, 
do imply $e_1 = \cdots = e_r = 0$.
\end{definition}

For instance, Theorem~3 on 
p.~\pageref{Theorem-3-S-33} states that an $r$-term continuous
local transformation
group $x_i' = f_i ( x; \, a_1, \dots, a_r)$ whose parameters $a_k$ are
all \emphasis{essential} always gives rise to the $r$ infinitesimal
transformations $X_k := - \frac{ \partial f}{ \partial a_k} (x; \,
e)$, $k = 1, \dots, r$, which are \emphasis{independent of each other}.

Introducing $r$ arbitrary auxiliary constants $\lambda_1, \dots,
\lambda_r$, one may consider the one-term group generated by the
general linear combination:
\[
C
:=
\lambda_1\,X_1
+\cdots+
\lambda_r\,X_r,
\] 
of $X_1, \dots, X_r$, namely the flow:
\[
\aligned
x_i'
=
\exp(tC)(x_i)
&
=
x_i
+
\frac{t}{1}\,C(x_i)
+
\frac{t^2}{1\cdot 2}\,
C\big(C(x_i)\big)
+
\cdots
\\
&
=
x_i
+
t\,\sum_{k=1}^r\,\lambda_k\,\xi_{ki}
+
t^2\,\sum_{k,\,j}^{1\dots r}\,
\frac{\lambda_k\,\lambda_j}{1\cdot 2}\,
X_k(\xi_{ji})
+\cdots
\\
&
=:
h_i(x;\,t,\lambda_1,\dots,\lambda_r)
\ \ \ \ \ \ \ \ \ \ \ \ \
{\scriptstyle{(i\,=\,1\,\cdots\,n)}}.
\endaligned
\]
If the $X_k \equiv - \frac{ \partial f}{ \partial a_k} ( x; \, e)$ do
stem from an $r$-term continuous group $x_i' = f_i( x; \, a_1, \dots,
a_r)$, a natural question 
\label{question-compare}
is then to compare the above integrated
finite equations $h_i ( x; \, t, \lambda_1, \dots, \lambda_r)$ to the
original transformation equations $f_i ( x; \, a_1, \dots,
a_r)$. Before studying this question together with Lie and Engel, we
focus our attention on a subquestion whose proof shows a beautiful,
synthetical, geometrical idea: that
of prolongating the action jointly to finite
sets of points. 

At first,
without assuming that the $X_k$ stem from an $r$-term continuous
group, it is to be asked (subquestion) whether the parameters
$\lambda_1, \dots, \lambda_r$ in the above integrated transformation
equations $x_i' = h_i ( x; \, \lambda_1, \dots, \lambda_r)$ are all
essential. In the formula just above, we notice that the $r+1$
parameters only appear in the form $t\, \lambda_1$, \dots, $t\,
\lambda_r$, hence because the $\lambda_k$ are arbitrary, there is no
restriction to set $t = 1$. We shall then simply write $h_i ( x; \,
\lambda_1, \dots, \lambda_r)$ instead of $h_i ( x; \, 1, \lambda_1,
\dots, \lambda_r)$.

\begin{svgraybox}
\vspace{-0.4cm}
\def\thetheorem{8}\begin{theorem}
\label{Theorem-8-S-65}
If the $r$ independent infinitesimal transformations:
\[
X_k(f)
=
\sum_{i=1}^n\,\xi_{ki}(x_1,\dots,x_n)\,
\frac{\partial f}{\partial x_i}
\ \ \ \ \ \ \ \ \ \ \ \ \
{\scriptstyle{(k\,=\,1\,\cdots\,r)}}
\]
are independent of each other, and if furthermore $\lambda_1, \dots,
\lambda_r$ are arbitrary parameters, then the totality of all one-term
groups $\lambda_1 \, X_1 ( f) + \cdots + \lambda_r \, X_r (f)$ forms a
family of transformations:
\def\theequation{5}\begin{equation}
x_i'
=
x_i
+
\sum_{k=1}^r\,\lambda_k\,\xi_{ki}
+
\sum_{k,\,j}^{1\dots r}\,
\frac{\lambda_k\,\lambda_k}{1\cdot 2}\,
X_k(\xi_{ji})
+\cdots
\ \ \ \ \ \ \ \ \ \ \ \ \
{\scriptstyle{(i\,=\,1\,\cdots\,n)}},
\end{equation}
in which the $r$ parameters $\lambda_1, \dots, \lambda_r$ are all 
\emphasis{essential}, hence a family of $\infty^r$ different 
transformations.
\end{theorem}
\vspace{-0.3cm}
\end{svgraybox}

\begin{proof}\smartqed
Here exceptionally, we observed a harmless technical incorrection in
Engel-Lie's proof (\cite{ enlie1888-4}, pp.~62--65) about the link
between the generic rank of $X_1 \big\vert_x, \dots, X_r \big\vert_x$
and a lower bound for the number of essential 
parameters\footnote{\,
On page~63, it is said that if the number $r$ of the independent
infinitesimal transformations $X_k$ is $\leqslant n$, then the $r
\times n$ matrix $\big( \xi_{ ki} ( x) \big)_{ 1\leqslant k \leqslant
r}^{ 1 \leqslant i \leqslant n}$ of their coefficients is of generic
rank equal to $r$, although this claim is contradicted with $n = r =
2$ by the two vector fields $x \, \frac{ \partial }{\partial x} + y \,
\frac{ \partial }{\partial y}$ and $xx \, \frac{ \partial }{\partial
x} + xy \, \frac{ \partial }{\partial y}$. Nonetheless, the ideas and
the arguments of the written proof (which does not really needs such a
fact) are perfectly correct.
}. 

However, Lie's main idea is clever and pertinent: it consists in the
introduction of exactly $r$ (the number of $\lambda_k$'s) copies of
the same space $x_1, \dots, x_n$ whose coordinates are labelled as
$x_1^{ ( \mu)}, \dots, x_n^{ ( \mu)}$ for $\mu = 1, \dots, r$ and to
consider the family of transformation equations induced by the 
\emphasis{same transformation equations}:
\[
{x_i^{(\mu)}}' 
=
\exp(C)\big(x_i^{(\mu)}\big)
= 
h_i\big(x^{(\mu)};\,
\lambda_1,\dots,\lambda_r\big)
\ \ \ \ \ \ \ \ \ \ \ \ \
{\scriptstyle{(i\,=\,1\,\cdots\,n;\,\,\,\mu\,=\,1\,\cdots\,r)}}
\]
on each copy of space, again with $t = 1$. Geometrically, one thus
views how the initial transformation equations $x_i' = h_i ( x; \,
\lambda_1, \dots, \lambda_r)$ act \emphasis{simultaneously} on $r$-tuples
of points. Written in greater length, these transformations read:
\def\theequation{5'}\begin{equation}
\aligned
{x_i^{(\mu)}}'
=
x_i^{(\mu)}
+
&
\sum_{k=1}^r\,\lambda_k\,\xi_{ki}^{(\mu)}
+
\sum_{k,\,j}^{1\dots r}\,
\frac{\lambda_k\,\lambda_j}{1\cdot 2}\,
X_k^{(\mu)}\big(\xi_{ji}^{(\mu)}\big)
+
\cdots
\\
&
\ \ \ \ \ \ \ \ \ \ \ \ \
{\scriptstyle{(i\,=\,1\,\cdots\,n;\,\,\,\mu\,=\,1\,\cdots\,r)}},
\endaligned
\end{equation}
where we have of course set: $\xi_{ ki}^{ ( \mu)} := \xi_{ ki} ( x^{ (
\mu)})$ and $X_k^{ ( \mu)} := \sum_{ i=1}^n\, \xi_{ ki} ( x^{ ( \mu)})
\, \frac{ \partial }{\partial x_i^{ (\mu)}}$.
Such an idea also reveals to be fruitful in 
other contexts.

According to the theorem stated 
on p.~\pageref{Theorem-essential}, in order to
check that the parameters $\lambda_1, \dots, \lambda_r$ are essential,
one only has to expand $x'$ in power series with respect to the powers of
$x$ at the origin:
\[
x_i'
=
\sum_{\alpha\in\N^n}\,
\mathcal{U}_\alpha^i(\lambda)\,
x^\alpha
\ \ \ \ \ \ \ \ \ \ \ \ \
{\scriptstyle{(i\,=\,1\,\dots\,n)}},
\]
and to show that the generic rank of the infinite coefficient mapping
$\lambda \longmapsto \big( \mathcal{ U}_\alpha^i ( \lambda) \big)_{
\alpha \in \N^n}^{ 1 \leqslant i \leqslant n}$ is maximal possible
equal to $r$. Correspondingly and immediately, we get the expansion of
the $r$-times copied transformation equations:
\def\theequation{5''}\begin{equation}
{x_i^{(\mu)}}'
=
\sum_{\alpha\in\N^n}\,
\mathcal{U}_\alpha^{i,(\mu)}(\lambda)\,(x^{(\mu)})^\alpha
\ \ \ \ \ \ \ \ \ \ \ \ \
{\scriptstyle{(i\,=\,1\,\dots\,n;\,\,\,\mu\,=\,1\,\dots\,r)}},
\end{equation}
with, for each $\mu = 1, \dots, r$, the \emphasis{same} 
coefficient functions:
\[
\mathcal{U}_\alpha^{i,(\mu)}(\lambda)
\equiv 
\mathcal{U}_\alpha^i
(\lambda)
\ \ \ \ \ \ \ \ \ \ \ \ \
{\scriptstyle{(i\,=\,1\,\dots\,n;\,\,\,\alpha\,\in\,\N^n;
\,\,\,\mu\,=\,1\,\dots\,r)}}. 
\]
So the generic rank of the corresponding infinite
coefficient matrix, which is just an $r$-times copy of the same
mapping $\lambda \longmapsto \big( \mathcal{ U}_\alpha^i ( \lambda)
\big)_{ \alpha \in \N^n}^{ 1 \leqslant i \leqslant n}$, does neither
increase nor decrease.

\plainstatement{Thus, 
the parameters $\lambda_1, \dots, \lambda_r$ for the
transformation equations $x' = h ( x; \, \lambda)$ are essential if
and only if they are essential for the diagonal transformation
equations ${ x^{ ( \mu )}} ' = h \big( x^{ ( \mu)}; \, \lambda\big)$,
$\mu = 1, \dots, r$, induced on the $r$-fold copy of the space $x_1,
\dots, x_n$. }

Therefore, we are left with the purpose of showing that the generic
rank of the $r$ times copy of the infinite coefficient matrix $\lambda
\longmapsto \big( \mathcal{ U}_\alpha^{ i, (\mu)} ( \lambda) \big)_{
\alpha \in \N^n}^{ 1 \leqslant i \leqslant n, \, \, 1 \leqslant \mu
\leqslant r}$ is equal to $r$. We shall in fact establish more, 
namely that the rank at $\lambda = 0$ of this map already equals $r$,
or equivalently, that the infinite constant matrix:
\[
\bigg(
\frac{\partial\mathcal{U}_\alpha^{i,(\mu)}}{\partial\lambda_k}(0)
\bigg)_{1\leqslant k\leqslant r}^{
1\leqslant i\leqslant n,\,\alpha\in\N^n,\,
1\leqslant\mu\leqslant r},
\]
whose $r$ lines are labelled with respect to partial derivatives, 
has rank equal to $r$.

To prepare this infinite matrix, if we differentiate the
expansions~\thetag{ 5'} which identify to~\thetag{ 5''} with respect
to $\lambda_k$ at $\lambda = 0$, and if we expand the coefficients of
our infinitesimal transformations:
\[
\xi_{ki}(x^{(\mu)})
=
\sum_{\alpha\in\N^n}\,
\xi_{ki\alpha}\,(x^{(\mu)})^\alpha
\ \ \ \ \ \ \ \ \ \ \ \ \
{\scriptstyle{(i\,=\,1\,\dots\,n;\,\,\,k\,=\,1\,\dots\,r;
\,\,\,\mu\,=\,1\,\dots\,r)}}
\]
with respect to the powers of $x_1, \dots, x_n$, we obtain a more
suitable expression of it:
\[
\aligned
\bigg(
\frac{\partial\mathcal{U}_\alpha^{i,(\mu)}}{\partial\lambda_k}(0)
\bigg)_{
1\leqslant k\leqslant r}^{1\leqslant i\leqslant n,\,\alpha\in\N^n,\,
1\leqslant\mu\leqslant r}
&
\equiv
\Big(
\big(
\xi_{ki\alpha}
\big)_{1\leqslant k\leqslant r}^{
1\leqslant i\leqslant n,\,\alpha\in\N^n}
\,\,\cdots\,\,
\big(
\xi_{ki\alpha}
\big)_{1\leqslant k\leqslant r}^{
1\leqslant i\leqslant n,\,\alpha\in\N^n}
\Big)
\\
&
=:
\Big(
{\sf J}^\infty\Xi(0)
\,\,\cdots\,\,
{\sf J}^\infty\Xi(0)
\Big).
\endaligned
\]
As argued up to now, it thus suffices to show that this matrix has
rank $r$. Also, we observe that this matrix identifies with the
\emphasis{infinite 
jet matrix} ${\sf J }^\infty \Xi(0)$ of Taylor coefficients of the
$r$-fold copy of the same $r \times n$ matrix of coefficients of the
vector fields $X_k$:
\[
\Xi(x)
:=
\left(
\begin{array}{ccc}
\xi_{11}(x) & \cdots & \xi_{1n}(x)
\\
\cdots & \cdots & \cdots
\\
\xi_{r1}(x) & \cdots & \xi_{rn}(x)
\end{array}
\right).
\]
This justifies the symbol ${\sf J}^\infty$ introduced just above. At
present, we can formulate an auxiliary lemma which will enable us to
conclude.

\begin{lemma}
Let $n\geqslant 1$, $q \geqslant 1$, $m \geqslant 1$ be
integers, let ${\sf x} \in \K^n$ and let:
\[
{\sf A}({\sf x})
=
\left(
\begin{array}{ccc}
a_{11}({\sf x}) & \cdots & a_{1m} ({\sf x})
\\
\cdots & \cdots & \cdots
\\
a_{q1}({\sf x}) & \cdots & a_{qm}({\sf x})
\end{array}
\right)
\]
be an arbitrary $q \times m$ matrix of analytic functions:
\[
a_{ij}({\sf x})
=
\sum_{\alpha\in\N^n}\,a_{ij\alpha}\,{\sf x}^\alpha
\ \ \ \ \ \ \ \ \ \ \ \ \
{\scriptstyle{(i\,=\,1\,\dots\,q;\,\,\,j\,=\,1\,\dots\,m)}}
\]
that are all defined in a fixed neighborhood of the origin in $\K^n$,
and introduce the $q \times \infty$ constant matrix of Taylor
coefficients:
\[
{\sf J}^\infty{\sf A}(0)
:=
\big(
a_{ij\alpha}
\big)_{1\leqslant i\leqslant q}^{
1\leqslant j\leqslant m,\,\alpha\in\N^n}
\]
whose $q$ lines are labelled by the index $i$. Then the following
inequality between (generic) ranks holds true:
\[
{\rm rk}\,{\sf J}^\infty{\sf A}(0)
\geqslant
{\rm genrk}\,{\sf A}({\sf x}).
\]
\end{lemma}

\begin{proof}\smartqed
Here, our infinite matrix ${\sf J}^\infty {\sf A} ( 0)$ will be
considered as acting by \emphasis{left} multiplication on \emphasis{horizontal
vectors} $u = (u_1, \dots, u_q)$, so that $u \, {\sf J}^\infty {\sf
A} ( 0)$ is an $\infty \times 1$ matrix, namely an infinite horizontal
vector. Similarly, ${\sf A} ( {\sf x})$ will act on horizontal vectors of
analytic functions $(u_1 ( {\sf x}), \dots, u_r ( {\sf x}))$.

Supposing that $u = (u_1, \dots, u_q) \in \K^q$ is any nonzero vector
in the kernel of ${\sf J}^\infty {\sf A} ( 0)$, namely: $0 = u \,
{\sf J}^\infty {\sf A} ( 0)$, or else in greater length:
\[
0
=
u_1\, a_{1j\alpha}
+\cdots+
u_q\, a_{qj\alpha}
\ \ \ \ \ \ \ \ \ \ \ \ \
{\scriptstyle{(j\,=\,1\,\dots\,m;\,\,\,\alpha\,\in\,\N^n)}},
\]
we then immediately deduce, after multiplying each such equation by
${\sf x}^\alpha$ and by summing up over all $\alpha \in \N^n$:
\[
0
\equiv
u_1\, a_{1j}({\sf x})
+\cdots+
u_q\, a_{qj}({\sf x})
\ \ \ \ \ \ \ \ \ \ \ \ \
{\scriptstyle{(j\,=\,1\,\dots\,m)}},
\]
so that the same constant vector $u = (u_1, \dots, u_q)$ also
satisfies $0 \equiv u \, A ({\sf x})$. It follows that the
dimension of the kernel of ${\sf J}^\infty {\sf A} ( 0)$ is smaller
than or equal to the dimension of the kernel of $A ( {\sf x})$ (at a
generic ${\sf x}$): this is just equivalent to the above inequality
between (generic) ranks.
\qed\end{proof}

Now, for each $q = 1, 2, \dots, r$, we want to apply the lemma with
the matrix ${\sf A} ({\sf x})$ being the $q$-fold copy of matrices
$\big( \Xi(x^{(1)}) \, \cdots \, \Xi (x^{ ( q)}) \big)$, or
equivalently in greater length:
\[\label{q-fold-extended-coefficient-matrix}
\Xi_q
\big(
\widetilde{\sf x}_q
\big)
:=
\left(
\begin{array}{cccccccccc}
\xi_{11}^{(1)} & \cdots & \xi_{1n}^{(1)} &
\xi_{11}^{(2)} & \cdots & \xi_{1n}^{(2)} &
\cdots\cdots &
\xi_{11}^{(q)} & \cdots & \xi_{1n}^{(q)}
\\
\cdots & \cdots & \cdots &
\cdots & \cdots & \cdots &
\cdots\cdots &
\cdots & \cdots & \cdots
\\ 
\xi_{r1}^{(1)} & \cdots & \xi_{rn}^{(1)} &
\xi_{r1}^{(2)} & \cdots & \xi_{rn}^{(2)} &
\cdots\cdots &
\xi_{r1}^{(q)} & \cdots & \xi_{rn}^{(q)}
\end{array}
\right),
\]
where we have abbreviated: 
\[
\widetilde{\sf x}_q 
:= 
\big(x^{(1)},\dots, 
x^{(q)}\big).
\]

\begin{lemma}
\label{crucial-assertion}
It is a consequence of the fact that $X_1, X_2, \dots, X_r$ are
linearly independent of each other that for every $q = 1, 2 ,\dots,
r$, one has:
\[
{\rm genrk}\,
\Big(
\Xi\big(x^{(1)}\big)\ \ \
\Xi\big(x^{(2)}\big)\ \ \
\cdots\ \ \
\Xi\big(x^{(q)}\big)
\Big)
\geqslant q.
\]
\end{lemma}

\begin{proof}\smartqed
Indeed, for $q = 1$, it is at first clear that ${\rm genrk}\, \big(
\Theta ( x^{ (1)} ) \big) \geqslant 1$, just because not all the
$\xi_{ ki} (x)$ vanish identically.

We next establish by induction that, as long as they remain $< r$,
generic ranks do increase of at least one unity at each step:
\[
{\rm genrk}
\Big(
\Xi_{q+1}
\big(\widetilde{\sf x}_{q+1}\big)
\Big)
\geqslant 
1
+
{\rm genrk}
\Big(
\Xi_q\big(\widetilde{\sf x}_q\big)
\Big),
\]
a fact which will immediately yield the lemma.

Indeed, if on the contrary, the generic ranks would stabilize, and
still be $< r$, then locally in a neighborhood of a generic, fixed
$\widetilde{
\sf x}_{ q+1}^0$, both matrices $\Xi_{ q+1}$ and $\Xi_q$ would have the
same, locally constant rank. Consequently, the solutions $\big(
\vartheta_1 ( \widetilde{ \sf x}_q) \, \, \cdots \, \, \vartheta_r (
\widetilde{ \sf x}_q ) \big)$ to the (kernel-like) system of linear
equations written in matrix form:
\[
0
\equiv
\big(
\vartheta_1(\widetilde{\sf x}_q)
\,\,\cdots\,\,
\vartheta_r(\widetilde{\sf x}_q)
\big)
\,
\Xi_q\big(\widetilde{\sf x}_q\big),
\]
which are analytic near $\widetilde{ \sf x}_q^0$ thanks to an application
of Cramer's rule and thanks to constancy of rank, would be
automatically also solutions of the extended system:
\[
0
\equiv
\big(
\vartheta_1(\widetilde{\sf x}_q)
\,\,\cdots\,\,
\vartheta_r(\widetilde{\sf x}_q)
\big)
\,
\big(
\Xi_q(\widetilde{\sf x}_q)\,\,\,\,
\Xi(x^{(q+1)})
\big),
\]
whence there would exist \emphasis{nonzero} solutions $(\vartheta_1, \dots,
\vartheta_r)$ to the linear dependence equations:
\[
0
=
\big(
\vartheta_1\,\,\cdots\,\,\vartheta_r
\big)\,
\Xi
\big(x^{(q+1)}\big)
\]
which are \emphasis{constant} with respect to the variable $x^{ ( q+1)}$,
since they only depend upon $\widetilde{ \sf x}_q$. This exactly
contradicts the assumption that $X_1^{ (q+1)}, \dots, X_r^{ (q+1)}$
are independent of each other.
\qed\end{proof}

Lastly, we may chain up a series of (in)equalities that are now
obvious consequences of the lemma and of the assertion:
\[
{\rm rank}\,
\Big(
{\sf J}^\infty\Xi(0)
\,\,\cdots\,\,
{\sf J}^\infty\Xi(0)
\Big)
=
{\rm rank}\,
{\sf J}^\infty\Xi_r(0)
\geqslant
{\rm genrk}\,
\Xi_r\big(\widetilde{\sf x}_r\big)
=
r,
\]
and since all ranks are anyway $\leqslant r$, we get the promised rank
estimation:
\[
r
=
{\rm rank}\,\Big(
{\sf J}^\infty\Xi(0)
\,\,\cdots\,\,
{\sf J}^\infty\Xi(0)
\Big),
\]
which finally completes the proof of the theorem.
\qed\end{proof}

In order to keep a memory track of the trick of extending the group
action to an $r$-fold product of the base space, we also translate a
summarizing proposition which is formulated on p.~66 of~\cite{
enlie1888-4}.

\begin{svgraybox}
\vspace{-0.4cm}
\def\theproposition{5}\begin{proposition}
\label{Satz-S-66}
If the $r$ infinitesimal transformations:
\[
X_k(f)
=
\sum_{i=1}^n\,\xi_{ki}(x_1,\dots,x_n)\,
\frac{\partial f}{\partial x_i}
\ \ \ \ \ \ \ \ \ \ \ \ \
{\scriptstyle{(k\,=\,1\,\dots\,r)}}
\]
are independent of each other, if furthermore:
\[
x_1^{(\mu)},\dots,x_n^{(\mu)}
\ \ \ \ \ \ \ \ \ \ \ \ \
{\scriptstyle{(\mu\,=\,1\,\dots\,r)}}
\]
are $r$ different systems of $n$ variables, and if lastly one sets for
abbreviation:
\[
X_k^{(\mu)}(f)
=
\sum_{i=1}^n\,
\xi_{ki}
\big(
x_1^{(\mu)},\dots,x_n^{(\mu)}
\big)\,
\frac{\partial f}{\partial x_i^{(\mu)}}
\ \ \ \ \ \ \ \ \ \ \ \ \
{\scriptstyle{(k,\,\mu\,=\,1\,\dots\,r)}},
\]
then the $r$ infinitesimal transformations:
\[
W_k(f)
=
\sum_{\mu=1}^r\,X_k^{(\mu)}(f)
\ \ \ \ \ \ \ \ \ \ \ \ \
{\scriptstyle{(k\,=\,1\,\dots\,r)}}
\]
in the $nr$ variables $x_i^{ ( \mu)}$ satisfy \emphasis{no} relation
of the form:
\[
\sum_{k=1}^n\,
\chi_k
\big(
x_1^{(1)},\dots,x_n^{(1)},\dots\dots,
x_1^{(r)},\dots,x_n^{(r)}
\big)\,
W_k(f)
\equiv
0.
\]
\end{proposition}
\vspace{-0.3cm}
\end{svgraybox}

\smallskip

Also, we remark for later use as in~\cite{ enlie1888-4}, 
p.~65, that during
the proof of the theorem~8 above, it did not really matter that the
equations~\thetag{ 5} represented the finite equations of a family of
one-term groups. In fact, we only considered the terms of first order
with respect to $\lambda_1, \dots, \lambda_r$ in the finite
equations~\thetag{ 5}, and the crucial 
Lemma~\ref{crucial-assertion} 
emphasized during the
proof was true under the only assumption that the infinitesimal
transformations $X_1, \dots, X_r$ were mutually
independent. Consequently, Theorem~8 can be somewhat
generalized as
follows.

\begin{svgraybox}
\vspace{-0.4cm}
\def\theproposition{4}\begin{proposition}
\label{Satz-4-S-65} 
If a family of transformations contains the $r$ arbitrary 
parameters $e_1, \dots, e_r$ and if its equations, 
when they are expanded with respect to powers
of $e_1, \dots, e_r$, possess the form:
\[
x_i'
=
x_i
+
\sum_{k=1}^r\,e_k\,\xi_{ki}(x_1,\dots,x_n)
+\cdots
\ \ \ \ \ \ \ \ \ \ \ \ \
{\scriptstyle{(i\,=\,1\,\cdots\,n)}},
\]
where the neglected terms in $e_1, \dots, e_r$ are of second and of
higher order, and lastly, if the functions $\xi_{ ki} ( x)$ have the
property that the $r$ infinitesimal transformations made up with them:
\[
X_k(f)
=
\sum_{i=1}^n\,\xi_{ki}(x_1,\dots,x_n)\,
\frac{\partial f}{\partial x_i}
\ \ \ \ \ \ \ \ \ \ \ \ \
{\scriptstyle{(k\,=\,1\,\cdots\,r)}}
\]
are independent of each other, then those transformation
equations represent $\infty^r$ different transformations,
or what is the same: the $r$ parameters $e_1, \dots, e_r$
are essential.
\end{proposition}
\vspace{-0.3cm}
\end{svgraybox}

\section{Generation of an $r$-Term Group by its One-Term Subgroups}
\label{generation-by-one-term}

After these preparations, we may now come back to our
question formulated
on p.~\pageref{question-compare}: 
how to compare the equations $x_i' = f_i ( x; \, a)$ of a
given finite continuous transformation group to the equations:
\[
\aligned
x_i'
&
=
\exp\big(
t\,\lambda_1\,X_1
+\cdots+
t\,\lambda_r\,X_r\big)(x_i)
\\
&
=:
h_i\big(x;\,t,\lambda_1,\dots,\lambda_r\big)
\ \ \ \ \ \ \ \ \ \ \ \ \ \ \ \ \ \ \ \ \ \ \
{\scriptstyle{(i\,=\,1\,\dots\,n)}}
\endaligned
\]
obtained by integrating the general linear combination of its $r$
infinitesimal transformations $X_k = - \frac{ \partial f_i}{ \partial
a_k} ( x; \, e)$? Sometimes, such equations will be called 
as in~\cite{ enlie1888-4} the
\terminology{canonical finite equations} of the group.

Abbreviating $\lambda_1 \, X_1 + \cdots + \lambda_r \, X_r$ as the
infinitesimal transformation $C := \sum_{ i = 1}^n\, \xi_i ( x) \,
\frac{ \partial }{\partial x_i}$, whose coefficients are given by:
\[
\xi_i(x)
:=
\sum_{j=1}^n\,\lambda_j\,\xi_{ji}(x)
\ \ \ \ \ \ \ \ \ \ \ \ \
{\scriptstyle{(i\,=\,1\,\dots\,n)}},
\]
then by definition of the flow $x_i' = \exp ( tC) ( x)$, the functions
$h_i$ satisfy the first order system of ordinary differential
equations $\frac{ \D\,h_i}{ \D\,t} = \xi_i ( h_1, \dots, h_n)$, or
equivalently:
\def\theequation{6}\begin{equation}
\frac{\D\,h_i}{\D\,t}
=
\sum_{j=1}^r\,\lambda_j\,\xi_{ji}(h_1,\dots,h_n)
\ \ \ \ \ \ \ \ \ \ \ \ \
{\scriptstyle{(i\,=\,1\,\dots\,n)}},
\end{equation}
with of course the initial condition $h ( x; 0, \lambda) = x$ at $t =
0$. On the other hand, according to Theorem~3 on
p.~\pageref{Theorem-3-S-33},
we remember that the $f_i$ satisfy the fundamental differential
equations:
\def\theequation{7}\begin{equation}
\xi_{ji}(f_1,\dots,f_n)
=
\sum_{k=1}^r\,\alpha_{jk}(a)\,
\frac{\partial f_i}{\partial a_k}
\ \ \ \ \ \ \ \ \ \ \ \ \ \
{\scriptstyle{(i\,=\,1\,\dots\,n;\,\,\,j\,=\,1\,\dots\,r)}}.
\end{equation}

\begin{proposition}
\label{S-69}
If the parameters $a_1, \dots, a_r$ are the unique solutions $a_k ( t,
\lambda)$ to the system of first order ordinary differential
equations:
\[
\frac{\D\,a_k}{\D\,t}
=
\sum_{j=1}^r\,\lambda_j\,
\alpha_{jk}(a)
\ \ \ \ \ \ \ \ \ \ \ \ \
{\scriptstyle{(k\,=\,1\,\dots\,r)}}
\]
with initial condition $a ( 0, \lambda) = e$ being the identity
element, then the following identities hold:
\[
\aligned
f_i\big(x;\,a(t,\lambda)\big)
&
\equiv
\exp
\big(
t\,\lambda_1\,X_1
+\cdots+
t\,\lambda_r\,X_r
\big)
(x_i)
=
h_i\big(x;\,t,\lambda_1,\dots,\lambda_r\big)
\\
&
\ \ \ \ \ \ \ \ \ \ \ \ \ \ \ \ \ \ \ \ \ \ \ \ \ \
{\scriptstyle{(i\,=\,1\,\dots\,n)}}
\endaligned
\]
and they show how the $h_i$ are recovered from the $f_i$.
\end{proposition}

\begin{proof}\smartqed
Indeed, multiplying the equation~\thetag{ 7} by $\lambda_j$ and
summing over $j$ for $j$ equals $1$ up to $r$, we get:
\[
\sum_{k=1}^r\,
\frac{\partial f_i}{\partial a_k}\,
\sum_{j=1}^r\,\lambda_j\,\alpha_{jk}(a)
=
\sum_{j=1}^r\,\lambda_j\,\xi_{ji}(f_1,\dots,f_n)
\ \ \ \ \ \ \ \ \ \ \ \ \
{\scriptstyle{(i\,=\,1\,\dots\,n)}}.
\]
Thanks to the assumption about the $a_k$, we can replace the second
sum of the left-hand side by $\frac{ \D\,a_k}{ \D\,t}$, which yields
identities:
\[
\sum_{k=1}^r\,
\frac{\partial f_i}{\partial a_k}\,
\frac{\D\,a_k}{\D\,t}
\equiv
\sum_{j=1}^r\,\lambda_j\,\xi_{ji}(f_1,\dots,f_n)
\ \ \ \ \ \ \ \ \ \ \ \ \
{\scriptstyle{(i\,=\,1\,\dots\,n)}}
\]
in the left hand side of which we recognize just a plain derivation
with respect to $t$:
\[
\frac{\D\,f_i}{\D\,t}
=
\frac{\D}{\D\,t}
\big[
f_i\big(x;\,a(t,\lambda)\big)
\big]
\equiv
\sum_{j=1}^r\,\lambda_j\,\xi_{ji}(f_1,\dots,f_n)
\ \ \ \ \ \ \ \ \ \ \ \ \
{\scriptstyle{(i\,=\,1\,\dots\,n)}}.
\]
But since $f \big( x; \, a ( 0, \lambda) \big) = f ( x; \, e) = x$ has
the same initial condition $x$ at $t = 0$ as the solution $h \big( x;
\, t, \lambda\big)$ to~\thetag{ 6}, the uniqueness of solutions to
systems of first order ordinary differential equations immediately
gives the asserted coincidence $f \big( x; \, a ( t, \lambda) \big)
\equiv h\big( x; \, t, \lambda\big)$.
\qed\end{proof}

\section{Applications to the Economy of Axioms}
\label{application-theorem-9}

We now come back to the end of Chap.~\ref{fundamental-differential}, 
Sect.~\ref{substituting-axiom}, where
the three standard group axioms (composition; identity element;
existence of inverses) were superseded by the hypothesis of existence
of differential equations. We remind that in the proof of
Lemma~\ref{lemma-two-nondegeneracies}, a relocalization was needed to
assure that $\det
\psi_{ kj} ( a) \neq 0$. For the sake of clarity and of
rigor, we will, in the hypotheses, explicitly mention the subdomain
$\mathcal{ A}^1 \subset \mathcal{ A}$ where the determinant of the
$\psi_{ kj} ( a)$ does not vanish.

The following (apparently technical) theorem which is a mild
modification of the Theorem~9 on p.~72 of~\cite{ enlie1888-4}, will be
used in an essential way by Lie to derive his famous three fundamental
theorems in Chap.~\ref{kapitel-9} below.

\begin{svgraybox}
\vspace{-0.3cm}
\def\thetheorem{9}\begin{theorem}
\label{Theorem-9-S-72}
If, in the transformation equations defined for $(x, a) \in \mathcal{
X} \times \mathcal{ A}${\rm :}
\def\theequation{1}\begin{equation}
x_i'
=
f_i(x_1,\dots,x_n;\,a_1,\dots,a_r)
\ \ \ \ \ \ \ \ \ \ \ \ \
{\scriptstyle{(i\,=\,1\,\cdots\,n)}},
\end{equation}
the $r$ parameters $a_1, \dots, a_r$ are all essential and if in
addition, certain differential equations of the form:
\def\theequation{2}\begin{equation}
\frac{\partial x_i'}{\partial a_k}
=
\sum_{j=1}^r\,\psi_{kj}(a_1,\dots,a_r)\,
\xi_{ji}(x_1',\dots,x_n')
\ \ \ \ \ \ \ \ \ \ \ \ \
{\scriptstyle{(i\,=\,1\,\cdots\,n\,;\,\,\,k\,=\,1\,\cdots\,r)}}
\end{equation}
are identically satisfied by $x_1' = f_1 (x; \, a), \dots, x_n' = f_n
( x; \, a)$, where the matrix $\psi_{ kj} (a)$ is holomorphic and
invertible in some nonempty subdomain $\mathcal{ A}^1 \subset
\mathcal{ A}$, and where the functions $\xi_{ ji} ( x')$ are
holomorphic in $\mathcal{ X}$, then by introducing the $r$
infinitesimal transformations:
\[
X_k
:=
\sum_{i=1}^n\,\xi_{ki}(x)\,
\frac{\partial}{\partial x_i},
\]
it holds true that every transformation $x_i' = f_i ( x; \, a)$ whose
parameters $a_1, \dots, a_r$ lie in a small neighborhood of some fixed
$a^0 \in \mathcal{ A}^1$ can be obtained by firstly performing the
transformation:
\[
\overline{x}_i
=
f_i(x_1,\dots,x_n;\,a_1^0,\dots,a_r^0)
\ \ \ \ \ \ \ \ \ \ \ \ \
{\scriptstyle{(i\,=\,1\,\cdots\,n)}},
\]
and then secondly, by performing a certain transformation:
\[
x_i'
=
\exp\big(t\lambda_1X_1+\cdots+t\lambda_rX_r\big)(\overline{x}_i)
\ \ \ \ \ \ \ \ \ \ \ \ \
{\scriptstyle{(i\,=\,1\,\cdots\,n)}}
\]
of the one-term group generated by some suitable linear combination of
the $X_k$, where $t$ and $\lambda_1, \dots, \lambda_r$ are small
complex numbers.
\end{theorem}
\vspace{-0.2cm}
\end{svgraybox}

Especially, this technical statement will be useful later to show that
whenever $r$ infinitesimal transformations $X_1, \dots, X_r$ form a
Lie algebra, the composition of two transformations of the form $x' =
\exp \big( t \lambda_1 X_1 + \cdots +
t \lambda_r X_r \big)$ is again of the same form, hence the totality
of these transformations
truly constitutes a group.

\begin{proof}\smartqed
The arguments are essentially the same as the ones developed at the
end of the previous section (p.~\pageref{S-69}) for a genuinely local
continuous transformation group, except that the identity parameter
$e$ (which does not necessarily exist here) should be replaced by
$a^0$.

For the sake of completeness, let us perform the proof. On 
the first hand, we fix $a^0 \in
\mathcal{ A}^1$ and we introduce the solutions $a_k = a_k ( t,
\lambda_1, \dots, \lambda_r)$ of the following system of ordinary
differential equations:
\[
\frac{\D\,a_k}{\D\,t}
=
\sum_{j=1}^r\,\lambda_j\,\alpha_{jk}(a)
\ \ \ \ \ \ \ \ \ \ \ \ \
{\scriptstyle{(k\,=\,1\,\cdots\,r)}},
\]
with initial condition $a_k ( 0, \lambda_1, \dots, \lambda_r) = a_k^0$,
where $\lambda_1, \dots, \lambda_r$ are small complex parameters and
where, as before, $\alpha_{ jk} ( a)$ denotes the inverse
matrix of $\psi_{ jk} ( a)$, which is holomorphic in the whole of
$\mathcal{ A}_1$.

On the second hand, we introduce the local flow:
\[
\exp\big(t\lambda_1X_1+\cdots+t\lambda_rX_r\big)(\overline{x})
=:
h\big(\overline{x};\,t,\lambda\big)
\]
of the general linear combination $\lambda_1 X_1 + \cdots + \lambda_r
X_r$ of the $r$ infinitesimal transformations $X_k = \sum_{ i = 1}^n
\, \xi_{ ki} ( x) \, \frac{ \partial}{ \partial x_i}$, where 
$\overline{ x}$ is assumed to run in $\mathcal{ A}^1$. Thus by its
very definition, this flow integrates the ordinary differential
equations:
\[
\frac{\D\,h_i}{\D\,t}
=
\sum_{j=1}^r\,\lambda_j\,\xi_{ji}(h_1,\dots,h_n)
\ \ \ \ \ \ \ \ \ \ \ \ \
{\scriptstyle{(i\,=\,1\,\cdots\,n)}}
\]
with the initial condition $h \big( \overline{ x}; \, 0, \lambda \big) 
= \overline{ x}$.

On the third hand, we solve first the $\xi_{ ji}$ in the fundamental
differential equations~\thetag{ 2} using the inverse matrix
$\alpha$:
\[
\xi_{ji}(f_1,\dots,f_n)
=
\sum_{k=1}^r\,\alpha_{jk}(a)\,
\frac{\partial f_i}{\partial a_k}
\ \ \ \ \ \ \ \ \ \ \ \ \
{\scriptstyle{(i\,=\,1\,\cdots\,n\,;\,\,\,j\,=\,1\,\cdots\,r)}}.
\]
Then we multiply by $\lambda_j$, we sum and we recognize
$\frac{ \D\,a_k}{ \D\,t}$, which we then substitute:
\[
\aligned
\sum_{j=1}^r\,\lambda_j\,\xi_{ji}(f_1,\dots,f_n)
&
=
\sum_{k=1}^r\,
\frac{\partial f_i}{\partial a_k}\,
\sum_{j=1}^r\,\lambda_j\,\alpha_{jk}(a)
\\
&
=
\sum_{k=1}^r\,
\frac{\partial f_i}{\partial a_k}\,
\frac{\D\,a_k}{\D\,t}
\\
&
=
\frac{\D}{\D\,t}
\big[
f_i\big(x;\,a(t,\lambda)\big)
\big]
\ \ \ \ \ \ \ \ \ \ \ \ \
{\scriptstyle{(i\,=\,1\,\cdots\,n)}}.
\endaligned
\]
So the $f_i \big( x; \, a ( t, \lambda)\big)$ satisfy the
\emphasis{same} differential equations as the $h_i \big( \overline{
x}; \, t, \lambda
\big)$, and in addition, if we set $\overline{ x}$ equal to $f ( x; \,
a^0)$, both collections of solutions will have the
\emphasis{same} initial
value for $t = 0$, namely $f(x; \, a^0)$. In conclusion, by observing
that the $f_i$ and the $h_i$ satisfy the same equations, the
uniqueness property enjoyed by first order ordinary differential
equations yields the identity:
\[
\boxed{
f\big(x;\,a(t,\lambda)\big)
\equiv
\exp\big(t_1\lambda_1 X_1
+\cdots+
t\lambda_rX_r\big)\big(f(x;\,a^0)\big)
}
\]
expressing that every transformation $x' = f( x; \, a)$ for $a$ in a
neighbourhood of $a^0$ appears to be the composition of the fixed
transformation $\overline{ x} = f ( x; \, a^0)$ followed by a certain
transformation of the one-term group $\exp \big( t \lambda_1 X_1 +
\cdots + t \lambda_r X_r \big) ( \overline{ x})$.
\qed\end{proof}

\begin{svgraybox}
\centerline{\sf \S\,\,\,18.}

\smallskip
We now apply the preceding general developments to the
peculiar case where the $\infty^r$ transformations 
$x_i' = f_i ( x_1, \dots, x_n,\, a_1, \dots, a_r)$ 
consitute an $r$-term group.

If the equations~\thetag{ 1} represent an $r$-term group, then 
according to Theorem~3 p.~\pageref{Theorem-3-S-33}, 
there always are differential equations of the
form~\thetag{ 2}; so we do not need to specially enunciate
this requirement.

Moreover, we observe that all infinitesimal transformations
of the form
$\sum_{ i=1}^n \Big\{ \sum_{ j=1}^r\, \lambda_j \xi_{ ji} ( x)
\Big\} \, \frac{ \partial f}{ \partial x_i}$ can be
linearly expressed by means of the following
$r$ infinitesimal transformations:
\[
X_k(f)
=
\sum_{i=1}^n\,\xi_{ki}(x_1,\dots,x_n)\,
\frac{\partial f}{\partial x_i}
\ \ \ \ \ \ \ \ \ \ \ \ \
{\scriptstyle{(k\,=\,1\,\cdots\,r)}},
\]
since all the infinitesimal transformations~\thetag{ 9}
are contained in the expression:
\[
\sum_{k=1}^r\,\lambda_k\,X_k(f).
\]
Here, as we have underscored already in the introduction of this
chapter, the infinitesimal transformations $X_1 (f), \dots, X_r ( f)$
are independent of each other.

Consequently, we can state the following theorem about
arbitrary $r$-term groups:

\def\thetheorem{10}\begin{theorem}
To every $r$-term group:
\[
x_i'
=
f_i(x_1,\dots,x_n,\,a_1,\dots,a_r)
\ \ \ \ \ \ \ \ \ \ \ \ \
{\scriptstyle{(i\,=\,1\,\cdots\,n)}}
\]
are associated $r$ independent infinitesimal transformations:
\[
X_k(f)
=
\sum_{i=1}^n\,\xi_{ki}(x_1,\dots,x_n)\,
\frac{\partial f}{\partial x_i}
\ \ \ \ \ \ \ \ \ \ \ \ \
{\scriptstyle{(i\,=\,1\,\cdots\,n)}},
\]
which stand in the following relationship to the finite
transformations of the group: if $x_i' = f_i ( x_1, \dots, x_n; \,
a_1^0,
\dots, a_r^0 )$ is any transformation of the group, then 
every transformation $x_i' = f_i ( x, a)$ whose parameter lies in a
certain neighbourhood of $a_1^0, \dots, a_r^0$ can be obtained by
firstly executing the transformation $\overline{ x}_i = f_i ( x_1,
\dots, x_n, \, a_1^0, \dots, a_r^0)$ and secondly a certain
transformation $x_i' = \omega_i (
\overline{ x}_1, \dots, \overline{ x}_n )$
of a one-term group, the infinitesimal transformation of which has the
form $\lambda_1 X_1 ( f) + \cdots +
\lambda_r X_r ( f)$, where $\lambda_1, \dots, \lambda_r$ denote
certain suitably chosen constants.
\end{theorem}

If we not only know that the equations $x_i' = f_i ( x, a)$ 
represent an $r$-term group, but also that this
group contains the identity transformation, and 
lastly also that the parameters $a_k^0$ represent
a system of values of the identity transformation
in the domain $(\!( a_k )\!)$, then we can still 
say more. Indeed, if we in particular choose
for the transformation $\overline{ x}_i = f_i ( x, a^0)$ the
identity transformation, we then realize immediately 
that the transformations of our group are nothing but 
the transformations of those
one-term groups that are generated by the
infinitesimal transformations:
\[
\lambda_1\,X_1(f)
+\cdots+
\lambda_r\,X_r(f).
\]
Thus, if for abbreviation we set:
\[
\sum_{k=1}^r\,\lambda_k\,X_k(f)
=
C_k(f),
\]
then the equations:
\[
x_i'
=
x_i
+
\frac{t}{1}\,C(x_i)
+
\frac{t^2}{1\cdot 2}\,
C\big(C(x_i)\big)
+\cdots
\ \ \ \ \ \ \ \ \ \ \ \ \
{\scriptstyle{(i\,=\,1\,\cdots\,n)}}
\]
represent the $\infty^r$ transformations of the group. 
The fact that the $r +
1$ parameters: $\lambda_1, \dots, \lambda_r, t$ appear is just
fictitious here, for they are indeed only found in the $r$
combinations $\lambda_1 t, \dots,
\lambda_r t$. We can therefore quietly set $t$ equal to $1$. 
In addition, if we remember the representation of a one-term group by
a single equation given in eq.~\thetag{ 3} 
on p.~\pageref{S-52-eq-7}, 
then we realize that the equations of our $r$-term group may be
condensed in the single equation:
\[
f(x_1',\dots,x_n')
=
f(x_1,\dots,x_n)
+
C(f)
+
\frac{1}{1\cdot 2}\,
C\big(C(f)\big)
+\cdots.
\]
That it is possible to order the transformations
of the group as inverses in pairs
thus requires hardly any mention. 

We can briefly state as follows
the above result about the $r$-term groups
with identity transformation.

\def\thetheorem{11}\begin{theorem}
\label{Theorem-11}
If an $r$-term group:
\[
x_i'
=
f_i(x_1,\dots,x_n,\,a_1,\dots,a_r)
\ \ \ \ \ \ \ \ \ \ \ \ \
{\scriptstyle{(i\,=\,1\,\cdots\,n)}}
\]
contains the identity transformation, then its $\infty^r$
transformations can be organized in $\infty^{ r - 1}$ families of
$\infty^1$ transformations in such a way that each family amongst these
$\infty^{ r - 1}$ families consists of all the transformations of a
certain one-term group with the identity transformation. In order to
find these one-term groups, one forms the known equations:
\[
\frac{\partial x_i'}{\partial a_k}
=
\sum_{j=1}^r\,\psi_{kj}(a_1,\dots,a_r)\,
\xi_{ji}(x_1',\dots,x_n')
\ \ \ \ \ \ \ \ \ \ \ \ \
{\scriptstyle{(i\,=\,1\,\cdots\,n\,;\,
k\,=\,1\,\cdots\,r)}},
\]
which are identically satisfied after substituting $x_i' = f_i ( x, a)$, 
one further sets:
\[
\sum_{i=1}^n\,\xi_{ki}(x_1,\dots,x_n)\,
\frac{\partial f}{\partial x_i}
=
X_k(f)
\ \ \ \ \ \ \ \ \ \ \ \ \
{\scriptstyle{(k\,=\,1\,\cdots\,r)}},
\]
hence the expression:
\[
\sum_{k=1}^r\,\lambda_k\,X_k(f)
\]
with the $r$ arbitrary parameters $\lambda_1, \dots, \lambda_r$
represents the infinitesimal transformations of these $\infty^{ r - 1}$
one-term groups, and their finite equations have the form:
\[
x_i'
=
x_i
+
\sum_{k=1}^r\,\lambda_k\,\xi_{ki}(x)
+
\sum_{k,\,j}^{1\dots r}\,
\frac{\lambda_k\,\lambda_j}{1\cdot 2}\,
X_k(\xi_{ji})
+\cdots
\ \ \ \ \ \ \ \ \ \ \ \ \
{\scriptstyle{(i\,=\,1\,\cdots\,n)}}.
\]
The totality of all these finite transformations is identical with
the totality of all transformations of the group $x_i' = f_i ( x, a)$.
Besides, the transformations of this group can be ordered as inverses
in pairs.
\end{theorem}

\medskip
\centerline{\sf \S\,\,\,19.}

\smallskip
In general, if an $r$-term group contains all transformations of some
one-term group and if, in the sense discussed earlier, this one-term
group is generated by the infinitesimal transformation $X ( f)$, then
we say that the $r$-term group \terminology{contains the infinitesimal
transformation $X ( f)$}. Now, we have just seen that every $r$-term
group with the identity transformation can be brought to the form:
\[
\aligned
f(x_1',\dots,x_n')
&
=
f(x_1,\dots,x_n)
+
\sum_{k=1}^r\,\lambda_k\,X_k(f)
\\
&\ \ \ \ \
+
\frac{1}{1\cdot 2}\,
\sum_{k,\,j}^{1\cdots r}\,
\lambda_k\,\lambda_j\,
X_k\big(X_j(f)\big)
+\cdots,
\endaligned
\]
where $\lambda_1, \dots, \lambda_r$ denote arbitrary constants, while
$X_1 ( f), \dots, X_r ( f)$ stand for mutually independent
infinitesimal transformations. So we can say that \emphasis{every such
$r$-term group contains $r$ independent infinitesimal
transformations}.

One is very close to presume that an $r$-term group 
cannot contain more than $r$ independent infinitesimal
transformations.

\label{S-76}
In order to clarify this point, we want to consider directly the
question of when the infinitesimal transformation:
\[
Y(f)
=
\sum_{i=1}^n\,
\eta_i(x_1,\dots,x_n)\,
\frac{\partial f}{\partial x_i}
\]
is contained in the $r$-term group with the
$r$ independent infinitesimal transformations
$X_1 ( f) , \dots, X_r ( f)$.

If $Y ( f)$ belongs to the $r$-term group in question,
then the same is also true of the transformations:
\[
x_i'
=
x_i
+
\frac{\tau}{1}\,\eta_i(x)
+
\frac{\tau^2}{1\cdot 2}\,Y(\eta_i)
+\cdots
\ \ \ \ \ \ \ \ \ \ \ \ \
{\scriptstyle{(i\,=\,1\,\cdots\,n)}}
\]
of the one-term group generated by $Y (f)$. Hence if we execute
first an arbitrary transformation of this one-term group and after an
arbitrary transformation:
\[
x_i''
=
x_i'
+
\sum_{k=1}^r\,\lambda_k\,\xi_{ki}'
+
\sum_{k,\,j}^{1\cdots r}\,
\frac{\lambda_k\,\lambda_j}{1\cdot 2}\,
X_j'\big(\xi_{ki}'\big)
+\cdots
\]
of the $r$-term group, we then must obtain a transformation which also
belongs to the $r$-term group. By a calculation, we find that this
new transformation has the form:
\[
x_i''
=
x_i
+
\tau\,\eta_i(x)
+
\sum_{k=1}^r\,\lambda_k\,\xi_{ki}(x)
+\cdots
\ \ \ \ \ \ \ \ \ \ \ \ \
{\scriptstyle{(i\,=\,1\,\cdots\,n)}},
\]
where all the left out terms are of second and of higher order with
respect to $\lambda_1, \dots, \lambda_r, \tau$. For arbitrary
$\lambda_1, \dots, \lambda_r, \tau$, this transformation must belong
to the $r$-term group. Now, if the $r + 1$ infinitesimal
transformations $X_1 ( f), \dots, X_r ( f), Y ( f)$ were independent
of each other, then according to the Proposition~4 
p.~\pageref{Satz-4-S-65}, the 
last written equations would represent $\infty^{ r+1}$
transformations; but this is impossible, for the $r$-term group
contains in general only $\infty^r$ transformations. Consequently,
$X_1 ( f), \dots, X_r ( f), Y( f)$ are not independent of each other,
but since $X_1 ( f), \dots, X_r ( f)$ are so, then $Y (f)$ must be
linearly deduced from $X_1 ( f), \dots, X_r ( f)$, hence have the
form:
\[
Y(f)
=
\sum_{k=1}^r\,l_k\,X_k(f),
\] 
where $l_1, \dots, l_r$ denote appropriate constants. That is to say,
the following holds:

\def\theproposition{2}\begin{proposition}
If an $r$-term group contains the identity transformation, then it
contains $r$ independent infinitesimal transformations $X_1 ( f),
\dots, X_r (f)$ and every infinitesimal transformation contained in it
possesses the form $\lambda_1 X_1 ( f) + \cdots + \lambda_r \, X_r (
f)$, where $\lambda_1, \dots, \lambda_r$ denote constants.
\end{proposition}

We even saw above that \emphasis{every} infinitesimal transformation
of the form $\lambda_1 X_1 ( f) + \cdots + \lambda_r X_r (f)$ belongs
to the group; hence in future, we shall call the expression $\lambda_1
X_1 ( f) + \cdots + \lambda_r X_r (f)$ with the $r$ arbitrary
constants $\lambda_1, \dots, \lambda_r$ \terminology{the general
infinitesimal transformation} of the $r$-term group in question.

From the preceding considerations, it also comes the following
certainly special, but nevertheless important:

\def\theproposition{2}\begin{proposition}
If an $r$-term group contains the $m \leqslant r$ mutually independent
infinitesimal transformations $X_1 ( f), \dots, X_r ( f)$, then it
also contains every infinitesimal transformation of the following
form: $\lambda_1 \, X_1 ( f) + \cdots + \lambda_m X_m (f)$, where
$\lambda_1, \dots, \lambda_m$ denote completely arbitrary constants.
\end{proposition}

Of course, the researches done so far give the means to determine the
infinitesimal transformations of an $r$-term group $x_i' = f_i ( x, \,
a)$ with the identity transformation. But it is possible to reach the
objective more rapidly.

Let the identity transformation of our group go with the parameters:
$a_1^0, \dots, a_r^0$, and let $a_1^0,
\dots, a_r^0$ lie in the domain $( \! ( a ) \! )$,
so that the determinant $\sum\, \pm \psi_{ 11} ( a^0) \cdots
\psi_{ rr} ( a^0)$ is hence certainly different
from zero. Now, we have:
\[
\sum_{j=1}^r\,\psi_{kj}(a_1,\dots,a_r)\,
\xi_{ji}\big(f_1(x,a),\dots,f_n(x,a)\big)
\equiv
\frac{\partial}{\partial a_k}\,
f_i(x,\,a),
\]
hence when one sets $a_k = a_k^0$:
\[
\sum_{j=1}^r\,\psi_{kj}(a^0)\,\xi_{ji}(x_1,\dots,x_n)
\equiv
\left[
\frac{\partial}{\partial a^k}\,f_i(x,\,a)
\right]_{a=a^0}.
\]
We multiply this equation by $\frac{ \partial f}{\partial x_i}$
and we sum for $i$ from $1$ to $n$, which then gives:
\[
\sum_{j=1}^r\,\psi_{kj}(a^0)\, X_j(f)
=
\sum_{i=1}^n\,
\left[
\frac{\partial}{\partial a_k}\,f_i(x,\,a)
\right]_{a=a^0}\,
\frac{\partial f}{\partial x_i}
\ \ \ \ \ \ \ \ \ \ \ \ \
{\scriptstyle{(k\,=\,1\,\cdots\,r)}}
\]
Now, since $X_1 ( f), \dots, X_r ( f)$ are independent infinitesimal
transformations and since in addition the determinant of the $\psi_{
kj} ( a^0)$ is different from zero, the right-hand sides of the latter
equations represent $r$ independent infinitesimal transformations of
our group.

The following method is even somewhat simpler.

One sets $a_k = a_k^0 + \delta t_k$, where it 
is understood that $\delta t_1, \dots,
\delta_r$ are infinitely small quantities.
Then it comes:
\[
\aligned
x_i'
&
=
f_i\big(x_1,\dots,x_n,\,a_1^0+\delta t_1,\dots,a_r^0+\delta t_r\big)
\\
&
=
x_i
+
\sum_{k=1}^r\,
\left[
\frac{\partial}{\partial a_k}\,
f_i(x,\,a)
\right]_{a=a^0}\,\delta t_k
+\cdots,
\endaligned
\]
where the left out terms are of second and of higher order with
respect to the $\delta t_k$. Here, it is now immediately apparent
that our group contains the $r$ infinitesimal transformations:
\[
\aligned
x_i'
=
x_i
+
&
\left[\frac{\partial}{\partial a_k}f_i(x,a)\right]_{a=a_0}\,
\delta t_k\ \ \ \ \ \ \ \ \ 
{\scriptstyle{(i\,=\,1\,\cdots\,n)}}
\\
&
\ \ \ \ \ \ \ \ \ \ \ \ \
{\scriptstyle{(k\,=\,1\,\cdots\,r)}}
\endaligned
\]
However, the question whether these $r$ infinitesimal transformations
are independent of each other requires in each individual case yet a
specific examination, if one does not know from the beginning that
the determinant of the $\psi$ does not vanish for $a_k = a_k^0$.

\smallskip{\sf Example.}
We consider the general projective group:
\[
x'
=
\frac{x+a_1}{a_2x+a_3}
\]
of the once-extended manifold.

The infinitesimal transformations of this group
are obtained very easily by means of the method shown right now.
Indeed, one has $a_1^0 = 0$, $a_2^0 = 0$, $a_3^0 = 1$, hence we 
have:
\[
x'
=
\frac{x+\delta t_1}{x\,\delta t_2+1+\delta t_3}
=
x+\delta t_1
-
x\, \delta t_3
-
x^2\,\delta t_2
+\cdots,
\]
that is to say, our group contains the three mutually 
independent infinitesimal transformations:
\[
X_1(f)
=
\frac{\D\,f}{\D\,x},
\ \ \ \ \ \ \ \ \ \ \ \
X_2(f)
=
x\,\frac{\D\,f}{\D\,x},
\ \ \ \ \ \ \ \ \ \ \ \ 
X_3(f)
=
x^2\,\frac{\D\,f}{\D\,x}.
\]
The general infinitesimal transformation of our group has the
form:
\[
\big(\lambda_1+\lambda_2x+\lambda_3x^2\big)\,
\frac{\D\,f}{\D\,x},
\]
hence we obtain its finite transformations by integrating the ordinary
differential equation:
\[
\frac{\D\,x'}{\lambda_1+\lambda_2x'+\lambda_3{x'}^2}
=
\D\,t,
\]
adding the initial condition: $x' = x$ for $t = 0$.

In order to carry out this integration, we bring the differential 
equation to the form:
\[
\frac{\D\,x'}{x'-\alpha}
-
\frac{\D\,x'}{x'-\beta}
=
\gamma\,\D\,t,
\]
by setting:
\[
\lambda_1
=
\frac{\alpha\beta\gamma}{\alpha-\beta},
\ \ \ \ \ \ \ \
\lambda_2
=
-\frac{\alpha+\beta}{\alpha-\beta}\,\gamma,
\ \ \ \ \ \ \ \
\lambda_3
=
\frac{\gamma}{\alpha-\beta},
\]
whence:
\[
\aligned
2\alpha
=
-\frac{\lambda_2}{\lambda_3}
&
+
\frac{\sqrt{\lambda_2^2-4\lambda_1\lambda_3}}{\lambda_3},
\ \ \ \ \ \ \ 
2\beta
=
-\frac{\lambda_2}{\lambda_3}
-
\frac{\sqrt{\lambda_2^2-4\lambda_1\lambda_3}}{\lambda_3},
\\
&
\ \ \ \ \ \ \ \ \ \ \ \
\gamma
=
\sqrt{\lambda_2^2-4\lambda_1\lambda_3}.
\endaligned
\]
By integration, we find:
\[
l(x'-\alpha)
-
l(x'-\beta)
=
\gamma\,t
+
l(x-\alpha)
-
l(x-\beta),
\]
or:
\[
\frac{x'-\alpha}{x'-\beta}
=
e^{\gamma t}\,
\frac{x-\alpha}{x-\beta},
\]
and now there is absolutely no difficulty to express $\alpha, \beta,
\gamma$ in terms of $\lambda_1, \lambda_2, \lambda_3$, in order
to receive in this way the $\infty^3$ transformations of our
three-term group arranged in $\infty^2$ one-term groups, exactly as is
enunciated in Theorem~11.

Besides, a simple known form of our group is obtained if one keeps the
two parameters $\alpha$ and $\beta$, while one introduces the new
parameter $\overline{ \gamma}$ instead of $e^{ \gamma t}$; then our
group appears under the form:
\[
\frac{x'-\alpha}{x'-\beta}
=
\overline{\gamma}\,\frac{x-\alpha}{x-\beta}.
\]
\end{svgraybox}

\linestop


\chapter{Complete Systems of Partial Differential Equations}
\label{kapitel-5}
\chaptermark{Complete Systems of Partial Differential Equations}

\setcounter{footnote}{0}

\begin{svgraybox}
\abstract{
Any infinitesimal transformation $X = \sum_{ i = 1}^n \, \xi_i ( x) \,
\frac{ \partial }{ \partial x_i}$ can be considered as the first order
analytic partial differential equation $X \omega = 0$ with the unknown
$\omega$. After a relocalization, a renumbering and a rescaling, one
may supppose $\xi_n (x) \equiv 1$. Then the general solution $\omega$
happens to be any (local, analytic) function $\Omega \big( \omega_1,
\dots, \omega_{ n-1} \big)$ of the $(n-1)$ functionally independent
solutions defined by the formula:
\[
\omega_k(x)
:=
\exp\big(-x_nX\big)(x_k)
\ \ \ \ \ \ \ \ \ \ \ \ \
{\scriptstyle{(k\,=\,1\,\cdots\,n-1)}}.
\]
What about first order \emphasis{systems} $X_1 \omega = \cdots = X_q
\omega = 0$ of such differential equations? Any solution $\omega$
trivially satisfies also $X_i \big( X_k ( \omega) \big) - X_k \big(
X_i ( \omega) \big) = 0$. But it appears that the subtraction in the
Jacobi commutator $X_i \big( X_k ( \cdot ) \big) - X_k \big( X_i (
\cdot ) \big)$ kills all the second-order differentiation terms, so
that one may freely add such supplementary first-order differential
equations to the original system, continuing again and again, until
the system, still denoted by $X_1 \omega = \cdots = X_q \omega = 0$,
becomes \emphasis{complete} in the sense of Clebsch, namely satisfies,
locally in a neighborhood of a generic point $x^0$:
\newline\indent
\smallskip{\bf (i)} 
for all indices $i, k = 1, \dots, q$, there are appropriate functions
$\chi_{ ik \mu} ( x)$ so that $X_i\big( X_k ( f) \big) - X_k \big( X_i
( f ) \big) = \chi_{ ik1} ( x) \, X_1 ( f) +
\cdots + \chi_{ ik q} ( x) \, X_q ( f)$;
\newline\indent
\smallskip{\bf (ii)}
the rank of the vector space generated by the $q$ vectors
$X_1\big\vert_x, \dots, X_q \big\vert_x$ is constant equal to $q$ for
all $x$ near the central point $x^0$.
\newline\noindent
Under these assumptions, it is shown in this chapter that there are $n
- q$ functionally independent solutions $x_1^{ ( q)},
\dots, x_{ n- q}^{ (q)}$ of the system that are
analytic near $x_0$ such that any other solution is a suitable
function of these $n-q$ fundamental solutions.
}

\subsubsection*{First Order Scalar Partial Differential Equation}

As a prologue, we ask what are the general solutions $\omega$ of a
first order partial differential equation $X \omega = 0$ naturally
associated to a local analytic vector field $X = \sum_{ i=1}^n\, \xi_i
( x) \, \frac{ \partial }{\partial x_i}$. Free relocalization 
being always allowed in the theory of Lie, we may assume, after
possibly renumbering the variables, that $\xi_n$ does not vanish in a
small neighborhood of some point at which we center the origin of the
coordinates. Dividing then by $\xi_n (x)$, it is equivalent to seek
functions $\omega$ that are annihilated by the differential operator:
\[
X
=
\sum_{i=1}^{n-1}\,
\frac{\xi_i(x)}{\xi_n(x)}\,
\frac{\partial}{\partial x_i}
+
\frac{\partial}{\partial x_n},
\]
still denoted by $X$ and which now satisfies $X ( x_n) \equiv 1$. We
recall that the corresponding system of ordinary differential equations
which defines curves that are everywhere tangent to 
$X$, namely the system:
\[
\frac{\D\,{\sf x}_1}{\D\,t}
=
\frac{\xi_1\big({\sf x}(t)\big)}{\xi_n\big({\sf x}(t)\big)},\ \
\dots\dots,\ \ \ \
\frac{\D\,{\sf x}_{n-1}}{\D\,t}
=
\frac{\xi_{n-1}\big({\sf x}(t)\big)}{\xi_n\big({\sf x}(t)\big)},\ \ \ \ 
\frac{\D\,{\sf x}_n(t)}{\D\,t}
=
1,
\]
with initial condition for $t = 0$ being an arbitrary 
point of the hyperplane $\{ x_n = 0 \}$:
\[
{\sf x}_1(0)
=
x_1,\dots\dots,
{\sf x}_{n-1}(0)
=
x_{n-1},\ \ \
{\sf x}_n(0)
=
0
\]
is \emphasis{solvable} and has a unique vectorial solution $( {\sf x}_1,
\dots, {\sf x}_{ n-1}, {\sf x}_n)$ which is analytic in a neighborhood
of the origin. In fact, ${\sf x}_n (t) = t$ by an obvious integration,
and the $(n-1)$ other ${\sf x}_k (t)$ are given by the marvelous
exponential formula already shown on
p.~\pageref{flow-exponential-formula}:
\[
{\sf x}_k(t)
=
\exp(tX)(x_k)
=
\sum_{l\geqslant 0}\,
\frac{t^l}{l!}\,
X^l(x_k)
\ \ \ \ \ \ \ \ \ \ \ \ \
{\scriptstyle{(k\,=\,1\,\cdots\,n-1)}}.
\]
We then set $t := -x_n$ in this formula (minus sign will be crucial)
and we define the $(n-1)$ functions that are relevant to us:
\[
\aligned
\omega_k(x_1,\dots,x_n)
:=
{\sf x}_k(-x_n)
&
=
\exp\big(-x_nX\big)(x_k)
\\
&
=
\sum_{l\geqslant 0}\,
(-1)^l\,
\frac{(x_n)^l}{l!}\,
X^l(x_k).
\endaligned
\]

\begin{proposition}
The $(n-1)$ so defined functions $\omega_1, \dots, \omega_{ n-1}$ are
functionally independent solutions of the partial differential
equation $X \omega = 0$ with the rank of their Jacobian matrix $\big(
\frac{ \partial
\omega_k}{\partial x_i} \big)_{ 1 \leqslant i \leqslant n}^{ 1
\leqslant k \leqslant n-1}$ being equal to $n-1$ at the origin.
Furthermore, for every other solution $\omega$ of $X \omega = 0$,
there exists a local analytic function $\Omega = \Omega
\big( \omega_1, \dots, \omega_{ n-1} \big)$
defined in a neighborhood of the origin in $\K^{ n-1}$ such that:
\[
\omega(x)
\equiv
\Omega\big(\omega_1(x),\dots,\omega_{n-1}(x)\big).
\]
\end{proposition}

\begin{proof}\smartqed
Indeed, when applying $X$ to the above series defining the $\omega_k$,
we see that all terms do cancel out, just thanks to an application
Leibniz' formula developed in the form:
\[
X\big[
(x_n)^l\,X^l(x_k)
\big]
=
l\,(x_n)^{l-1}\,
X^l(x_k)
+
(x_n)^l\,X^{l+1}(x_k).
\]
Next, the assertion that the map $x \longmapsto
\big( \omega_1 (x), \dots, \omega_{ n-1} (x) \big)$ has rank 
$n-1$ is clear, for $\omega_k ( x_1, \dots, x_{ n-1}, 0) \equiv x_k$
by construction. Finally, after straightening $X$ to $X' :=
\frac{ \partial }{ \partial x_n'}$ in some new coordinates 
$(x_1', \dots, x_n')$ thanks to the theorem on
p.~\pageref{theorem-straightening}, the general solution $\omega ' (
x')$ to $X' \omega' = 0$ happens trivially to just be any function
$\Omega' (x_1', \dots, x_{ n-1}')$ of $x_1' \equiv \omega_1',
\dots, x_{ n-1}' \equiv \omega_{n-1}'$. 
\vspace{-0.4cm}
\qed\end{proof}

\end{svgraybox}

\bigskip

\centerline{\large\sf C\,h\,a\,p\,t\,e\,r\,\,\,5.}

\medskip
\begin{center}
{\large\bf
The Complete Systems.}
\end{center}

\medskip
We assume that the theory of the integration of an individual first
order linear partial differential equation:
\[
X(f)
=
\sum_{i=1}^n\,\xi_i(x_1,\dots,x_n)\,
\frac{\partial f}{\partial x_i}
=
0,
\]
or of the equivalent simultaneous system of ordinary differential
equations:
\[
\frac{\D\,x_1}{\xi_1}
=
\cdots
=
\frac{\D\,x_n}{\xi_n},
\]
is known; nonetheless, as an introduction, we compile without
demonstration a few related propositions. Based on these propositions,
we shall very briefly derive the theory of the integration of
simultaneous linear partial differential equations of the first
order. In the next chapter, we shall place in a new light \deutsch{in
ein neues Licht setzen} this theory due in the main whole to
\name{Jacobi} and to \name{Clebsch}, by explaining more closely the
connection between the concepts \deutsch{Begriffen} of ``linear
partial differential equation'' and of ``infinitesimal
transformation'', a connection that we have already mentioned earlier
(Chap.~\ref{one-term-groups}, p.~\pageref{cited-S-82}).

\sectionengellie{\S\,\,\,21.}

One can suppose that $\xi_1, \dots, \xi_n$ behave regularly in the
neighborhood of a determinate system of values $x_1^0, \dots, x_n^0$,
and as well that $\xi_n ( x_1^0, \dots, x_n^0)$ is different from
zero. Under these assumptions, one can determine $x_1, \dots, x_{ n -
1}$ as analytic functions of $x_n$ in such a way that by sustitution
of these functions, the simultaneous system:
\[
\frac{\D\,x_1}{\D\,x_n}
=
\frac{\xi_1}{\xi_n},\dots,
\frac{\D\,x_{n-1}}{\D\,x_n}
=
\frac{\xi_{n-1}}{\xi_n}
\]
is identically satisfied, and that in addition $x_1, \dots, x_{ n-1}$
for $x_n = x_n^0$ take certain prescribed \emphasis{initial values}
\deutsch{Anfangswerthe} $x_1', \dots, x_{ n-1}'$. These initial values
have to be interpreted as the integration constants.

The equations which, in the concerned way, represent $x_1, \dots, x_{
n - 1}$ as functions of $x_n$ are called the \terminology{complete
integral equations of the simultaneous system}; they can receive the
form:
\[
\aligned
x_k
=
x_k'
+
(x_n^0-x_n)\,
&
\mathfrak{P}_k
\big(
x_1'-x_1^0,\dots,x_{n-1}'-x_{n-1}^0,x_n^0-x_n
\big)
\\
&
\ \ \ \
{\scriptstyle{(k\,=\,1\,\cdots\,n\,-\,1)}},
\endaligned
\]
where the $\mathfrak{ P}_k$ denote ordinary power series in their
arguments. By inversing the relation
between $x_1, \dots, x_{ n-1}, x_n$ and
$x_1', \dots, x_{ n-1}', x_n^0$, one again obtains the integral
equations, resolved with respect to only the initial values $x_1',
\dots, x_{ n - 1}'$:
\[
\aligned
x_k'
=
x_k
+
(x_n-x_n^0)\,
&
\mathfrak{P}_k
\big(x_1-x_1^0,\dots,x_n-x_n^0\big)
=
\omega_k(x_1,\dots,x_n)
\\
&
\ \ \ \
{\scriptstyle{(k\,=\,1\,\cdots\,n\,-\,1)}}.
\endaligned
\]
Here, the functions $\omega_k$ are the so-called \terminology{integral
functions} of the simultaneous system, since the differentials of these
functions:
\[
\D\,\omega_k
=
\sum_{i=1}^n\,
\frac{\partial\omega_k}{\partial x_i}\,
\D\,x_i
\ \ \ \ \ \ \ \ \ \ \ \ \ \
{\scriptstyle{(k\,=\,1\,\cdots\,n\,-\,1)}}
\]
all vanish identically by virtue of the simultaneous system, and every
function of this sort is called an integral function of the
simultaneous system. But every such integral function is at the same
time a solution of the linear partial differential equation $X ( f ) =
0$, whence $\omega_1, \dots, \omega_{ n - 1}$ are solutions of $X ( f)
= 0$, and in fact, they are obviously independent. In a certain
neighborhood of $x_1^0, \dots, x_n^0$ these solutions behave
regularly; in addition, they reduce for $x_n = x_n^0$ to $x_1, \dots,
x_{ n-1}$ respectively; that is why they are called the
\terminology{general solutions of the equation $X ( f) = 0$ relative
to $x_n = x_n^0$}.

If one knows altogether $n-1$ independent solutions:
\[
\psi_1(x_1,\dots,x_n),\dots\dots,
\psi_{n-1}(x_1,\dots,x_n)
\]
of the equation $X ( f) = 0$, then the most general solution of it has
the form $\Omega ( \psi_1, \dots, \psi_{ n-1})$, where $\Omega$
denotes an arbitrary analytic function of its arguments.

\sectionengellie{\S\,\,\,22.}

If a function $\psi( x_1, \dots, x_n)$ satisfies the two equations:
\[
X_1(f)
=
0,
\ \ \ \ \ \ \ \ \
X_2(f)
=
0,
\]
then it naturally also satisfies the two differential equations of
second order:
\[
X_1\big(X_2(f)\big)
=
0,
\ \ \ \ \ \ \ \ \ 
X_2\big(X_1(f)\big)
=
0,
\]
and in consequence of that, also the equation:
\[
X_1\big(X_2(f)\big)
-
X_2\big(X_1(f)\big)
=
0,
\]
which is obtained by subtraction from the last two written ones.

If now one lets:
\[
X_k(f)
=
\sum_{i=1}^n\,
\xi_{ki}(x_1,\dots,x_n)\,
\frac{\partial f}{\partial x_i}
\ \ \ \ \ \ \ \ \ \ \ \ \
{\scriptstyle{(k\,=\,1,\,2)}},
\]
then it comes:
\[
X_1\big(X_2(f)\big)
-
X_2\big(X_1(f)\big)
=
\sum_{i=1}^n\,
\big\{
X_1(\xi_{2i})
-
X_2(\xi_{1i})
\big\}\,
\frac{\partial f}{\partial x_i},
\]
because all terms which contain second order differential quotients
are cancelled. Thus, the following holds:

\def\theproposition{1}\begin{proposition}
\label{Satz-1-S-84}
If a function $\psi ( x_1, \dots, x_n)$ satisfies the two
differential equations of first order:
\[
X_k(f)
=
\sum_{i=1}^n\,\xi_{ki}(x_1,\dots,x_n)\,
\frac{\partial f}{\partial x_i}
=
0
\ \ \ \ \ \ \ \ \ \ \ \ \
{\scriptstyle{(k\,=\,1,\,2)}},
\]
then it also satisfies the equation:
\[
X_1\big(X_2(f)\big)
-
X_2\big(X_1(f)\big)
=
\sum_{i=1}^n\,
\big\{
X_1(\xi_{2i})
-
X_2(\xi_{1i})
\big\}\,
\frac{\partial f}{\partial x_i}
=
0,
\]
which, likewise, is of first order.
\end{proposition}

It is of great importance to know how the expression $X_1 \big( X_2 (
f) \big) - X_2 \big( X_1 ( f) \big)$ behaves when, in place of $x_1,
\dots, x_n$, new independent variables $y_1, \dots, y_n$ are
introduced.

\smallskip

We agree that by introduction of the $y$, it arises:
\[
\aligned
X_k(f)
=
\sum_{i=1}^n\,X_k(y_i)\,\frac{\partial f}{\partial y_i}
=
&
\sum_{i=1}^n\,\eta_{ki}(y_1,\dots,y_n)\,\frac{\partial f}{\partial y_i}
=
Y_k(f)
\\
&
\ 
{\scriptstyle{(k\,=\,1,\,2)}}.
\endaligned
\]
Since $f$ denotes here a completely arbitrary function of $x_1, \dots,
x_n$, we can substitute $X_1 ( f)$ or $X_2 ( f)$ in place of $f$, so
we have:
\[
\aligned
X_1\big(X_2(f)\big)
&
=
Y_1\big(X_2(f)\big)
=
Y_1\big(Y_2(f)\big)
\\
X_2\big(X_1(f)\big)
&
=
Y_2\big(X_1(f)\big)
=
Y_2\big(Y_1(f)\big),
\endaligned
\]
and consequently:
\[
X_1\big(X_2(f)\big)
-
X_2\big(X_1(f)\big)
=
Y_1\big(Y_2(f)\big)
-
Y_2\big(Y_1(f)\big).
\]
Thus we have the

\def\theproposition{2}\begin{proposition}
\label{Satz-2-S-84}
If, by the introduction of a new independent variable, the expressions
$X_1 ( f)$ and $X_2 (f)$ are transferred to $Y_1 (f)$ and respectively
to $Y_2 (f)$, then the expression $X_1 \big( X_2 (f) \big) - X_2 \big(
X_1 (f) \big)$ is transferred to $Y_1 \big( Y_2 (f) \big) - Y_2 \big(
Y_1 ( f) \big)$.
\end{proposition}

This property of the expression $X_1 \big( X_2 ( f) \big) - X_2 \big(
X_1 (f) \big)$ will be frequently used in the course of our study. The
same proposition can be stated more briefly: \emphasis{the expression $X_1
\big( X_2 (f) \big) - X_2 \big( X_1 (f) \big)$ behaves invariantly
through the introduction of a new variable}.

\smallskip

We now consider the $q$ equations:
\def\theequation{1}\begin{equation}
X_1(f)
=
0,\dots\dots,
X_q(f)
=
0,
\end{equation}
and we ask about its possible joint solutions.

It is thinkable that between the expressions $X_k (f)$, there are
relations of the form:
\def\theequation{2}\begin{equation}
\sum_{k=1}^q\,\chi_k(x_1,\dots,x_n)\,
X_k(f)
\equiv
0.
\end{equation}
If this would be the case, then certain amongst our equations would be
a consequence of the remaining ones, and they could easily be left out
while taking for granted the solution of the stated problem.
Therefore it is completely legitimate to make the assumption that
there are no relations of the form~\thetag{ 2}, hence that the
equations~\thetag{ 1} are solvable with respect to the $q$ of the
differential quotients $\frac{ \partial f }{ \partial x_i}$. It is to
be understood in this sense, when we refer to the 
\emphasis{equations}~\thetag{ 1} as 
\emphasis{independent} of each other
\footnote{\,
This is a typical place where a relocalization is in general required
in order to insure that the vectors $X_1 \big\vert_x, \dots, X_q
\big\vert_x$ are locally linearly independent, and the authors, as
usual, understand it mentally.
}. 

According to what has been said above about the two equations $X_1 (
f) = 0$ and $X_2 ( f) = 0$, it is clear that the possible joint
solutions of our $q$ equations do also satisfy all equations of the
form:
\[
X_i\big(X_k(f)\big)
-
X_k\big(X_i(f)\big)
=
0.
\]
And now, two cases can occur.

Firstly, the equations obtained this way can be a consequence
of the former, when for every 
$i$ and $k \leqslant q$, a relation of the following
form:
\[
\aligned
&
X_i\big(X_k(f)\big)
-
X_k\big(X_i(f)\big)
=
\\
\ \ \ \ \ \ \ \
&
=
\chi_{ik1}(x_1,\dots,x_n)\,X_1(f)
+\cdots+
\chi_{ikq}(x_1,\dots,x_n)\,X_q(f).
\endaligned
\]
holds. With \name{Clebsch}, we then say that \emphasis{the $q$
independent equations $X_1 (f) = 0$, \dots, $X_q ( f) = 0$ form a}
\terminology{$q$-term complete system} \deutsch{$q$-gliedrig
vollständig System}.

However, in general the second possible case will occur;
amongst the new formed equations:\label{S-85}
\[
X_i\big(X_k(f)\big)
-
X_k\big(X_i(f)\big)
=
0,
\]
will be found a certain number which are independent of each other,
and from the $q$ presented ones. We add them, say:
\[
X_{q+1}(f)=0,
\dots\dots,
X_{q+s}(f)
=
0
\]
to the $q$ initial ones and we now treat the obtained $q + s$ 
equations exactly as the $q$ ones given earlier.
So we continue; but since we cannot come to more than
$n$ equations $X_i ( f) = 0$ that are independent of each other,
we must finally reach a complete system which consists of 
$n$ or less independent equations.
Therefore we have the proposition:

\def\theproposition{3}\begin{proposition}
The determination of the joint solutions of $q$ given linear partial
differential equations of first order $X_1 ( f) = 0, \dots, X_q ( f) =
0$ can always be reduced, by differentiation and elimination, to the
integration of a complete system.
\end{proposition}

We now assume that $X_1 ( f) = 0$, \dots, $X_q ( f) = 0$ form a
$q$-term complete system. Obviously, these equations may be
replaced by $q$ other ones:
\[
Y_k(f)
=
\sum_{j=1}^q\,\psi_{kj}(x_1,\dots,x_n)\,X_j(f)
=
0
\ \ \ \ \ \ \ \ \ \ \ \ \
{\scriptstyle{(k\,=\,1\,\cdots\,q)}}.
\]
In the process, it is only required that the determinant of the
$\psi_{ kj}$ does not vanish identically, whence the $X_j ( f)$ can
also be expressed linearly in terms of the $Y_k ( f)$.
Visibly, relations of the form:
\def\theequation{3}\begin{equation}
Y_i\big(Y_k(f)\big)
-
Y_k\big(Y_i(f)\big)
=
\sum_{j=1}^q\,\omega_{ikj}(x_1,\dots,x_n)\,
Y_j(f)
\end{equation}
then hold true; so the equations $Y_k (f) = 0$ too do form a $q$-term
complete system and hence are totally equivalent to the equations $X_k
( f) = 0$.

As it has been pointed out for the first time by \name{Clebsch}, some
$\psi_{ kj}$ are always available for which all the $\omega_{ ikj}$
vanish. In order to attain this in the simplest way, we select with
A.~\name{Mayer} the $\psi_{ kj}$ so that the $Y_k (f) = 0$ appear to be
resolved with respect to $q$ of the differential quotients, for
instance with respect to $\frac{ \partial f}{ \partial x_n}$, \dots,
$\frac{ \partial f}{ \partial x_{ n-q+1}}$:
\def\theequation{4}\begin{equation}
Y_k(f)
=
\frac{\partial f}{\partial x_{n-q+k}}
+
\sum_{i=1}^{n-q}\,
\eta_{ki}\,\frac{\partial f}{\partial x_i}
\ \ \ \ \ \ \ \ \ \ \ \ \
{\scriptstyle{(k\,=\,1\,\cdots\,q)}}.
\end{equation}
Then the expressions $Y_i \big( Y_k (f) \big) - Y_k \big( Y_i (f)
\big)$ will all be free of $\frac{ \partial f}{ \partial x_{ n - q +
1}}$, \dots, $\frac{ \partial f}{ \partial x_n}$, and consequently
they can have the form $\sum_j \, \omega_{ ikj}\, Y_j (f)$ only if all
the $\omega_{ ikj}$ vanish. Thus:

\def\theproposition{4}\begin{proposition}
If one solves a $q$-term complete system:
\[
X_1(f)=0,\dots\dots,
X_q(f)=0
\]
with respect to $q$ of the differential quotients, 
then the resulting equations:
\def\theequation{4}\begin{equation}
Y_k(f)
=
\frac{\partial f}{\partial x_{n-q+k}}
+
\sum_{i=1}^{n-q}\,
\eta_{ki}\,
\frac{\partial f}{\partial x_i}
=
0
\ \ \ \ \ \ \ \ \ \ \ \ \
{\scriptstyle{(k\,=\,1\,\cdots\,q)}}
\end{equation}
stand pairwise in the relationships:
\def\theequation{5}\begin{equation}
Y_i\big(Y_k(f)\big)
-
Y_k\big(Y_i(f)\big)
=
0
\ \ \ \ \ \ \ \ \ \ \ \ \
{\scriptstyle{(i,\,k\,=\,1\,\cdots\,q)}}.
\end{equation}
\end{proposition}

\sectionengellie{\S\,\,\,23.}

We now imagine that a given $q$-term complete system
is brought to the above-mentioned form:
\def\theequation{4}\begin{equation}
Y_k(f)
=
\frac{\partial f}{\partial x_{n-q+k}}
+
\sum_{i=1}^{n-q}\,\eta_{ki}\,
\frac{\partial f}{\partial x_i}
=
0
\ \ \ \ \ \ \ \ \ \ \ \ \
{\scriptstyle{(k\,=\,1\,\cdots\,q)}}.
\end{equation}
It will be shown that this system possesses $n-q$ independent
solutions, for the determination of which it
suffices to integrate $q$ individual linear partial differential
equations one after the other.

At first, we integrate the differential equation:
\[
Y_q(f)
=
\frac{\partial f}{\partial x_n}
+
\sum_{i=1}^{n-q}\,\eta_{qi}\,\frac{\partial f}{\partial x_i}
=
0;
\]
amongst its $n-1$ independent solutions $x_1', \dots, x_n'$, the
following $q-1$, namely:
\[
x_{n-q+1}'
=
x_{n-q+1},\dots\dots,
x_{n-1}'
=
x_{n-1}
\]
are known at once. If the $n-1$ expressions $x_1', \dots, x_{ n-1}'$
together with the quantity $x_n$ which is independent of them, are
introduced as new variables, then our complete system receives the
form:
\[
\aligned
Y_q(f)
=
\frac{\partial f}{\partial x_n}
=
0,
&
\ \ \ \ \ \ \ \
Y_k(f)
=
\frac{\partial f}{\partial x_{n-q+1}'}
+
\sum_{i=1}^{n-q}\,
\eta_{ki}'\,\frac{\partial f}{\partial x_i'}
=
0
\\
&
\ \ \ \ \ \ \ \ \ \ \ \ \
{\scriptstyle{(k\,=\,1\,\cdots\,q\,-\,1)}}.
\endaligned
\]
Now, since the expressions $Y_i \big( Y_k (f) \big) - Y_k \big( Y_i
(f) \big)$ behave invariantly \label{S-87} through the introduction
of new variables (cf. Proposition~2 of this chapter), they must now
again also vanish, from which it follows that all $\eta_{ ki}'$ must
be functions of only $x_1', \dots, x_{ n-1}'$, and be free of
$x_n$. The initial problem of integration is therefore reduced to
finding out the joint solutions of the $q-1$ equations:
\[
\aligned
Y_k'(f)
=
\frac{\partial f}{\partial x_{n-q+k}'}
+
&
\sum_{i=1}^{n-q}\,
\eta_{ki}'(x_1',\dots,x_{n-1}')\,
\frac{\partial f}{\partial x_i'}
=
0
\\
&
{\scriptstyle{(k\,=\,1\,\cdots\,q\,-\,1)}},
\endaligned
\]
and in fact these equations, which depend upon only $n-1$ independent
variables, namely $x_1', \dots, x_{ n-1}'$, do stand again pairwise in
the relationships: $Y_i' \big( Y_k' (f) \big) - Y_k' \big( Y_i' (f)
\big) = 0$.

We formulate this result in the following way.

\def\theproposition{5}\begin{proposition}
The joint solutions of the equations of a $q$-term complete system in
$n$ variables can also be defined as the joint solutions of the
equations of a $(q-1)$-term complete system in $n-1$ variables. In
order to be able to set up this new complete system, one only has to
integrate a single linear partial differential equation of first order
in $(n - q + 1)$ variables. 
\end{proposition}

The new complete system again appears in resolved form; 
we can hence at once continue our above-mentioned process 
and by applying it $(q-1)$ times, we obtain the following:

\def\theproposition{6}\begin{proposition}
The joint solutions of the equations of a $q$-term complete system in
$n$ variables can also be defined as the solutions of a single linear
partial differential equation of first order in $n-q+1$ variables. In
order to be able to set up this equation, it suffices to 
integrate $q-1$ individual equations of this sort one
after the other.
\end{proposition}

From this it follows easily:

\def\theproposition{7}\begin{proposition}
A $q$-term complete system in $n$ independent variables
always possesses $n - q$ independent solutions.
\end{proposition}

But conversely, the following also holds:

\def\theproposition{8}\begin{proposition}
\label{Satz-8-S-88}
If $q$ independent linear partial differential equations of 
first order in $n$ independent variables:
\[
X_1(f)=0,\dots\dots,
X_q(f)=0
\]
have exactly $n-q$ independent solutions in common, then they form a
$q$-term complete system.
\end{proposition}

For the proof, we remark that according to what precedes, the
equations:
\[
X_1(f)=0,\dots\dots,
X_q(f)=0
\]
determine a complete system with $q$ or more terms; now, if this
complete system would contain more than $q$ independent equations,
then it would possess not $n-q$, but only a smaller number of
independent solutions; under the assumptions of the proposition,
it is therefore $q$-term, that is to say, it is
constituted by the equations $X_1 (f) = 0$, \dots, 
$X_q (f) = 0$ themselves.

If $n-q$ independent functions $\omega_1, \dots, \omega_{ n-q}$
of $x_1, \dots, x_n$ are presented, then one always can
set up $q$ independent linear partial differential equations
which are identically satisfied by all $\omega$.
Indeed, if we take for granted that the determinant:
\[
\sum\,
\pm\,
\frac{\partial\omega_1}{\partial x_1}\,\dots\,
\frac{\partial\omega_{n-q}}{\partial x_{n-q}}
\]
does not vanish identically, which we can achieve
without restriction, 
then the $q$ equations:
\[
\left\vert
\begin{array}{cccc}
\frac{\partial f}{\partial x_1} & \cdots &
\frac{\partial f}{\partial x_{n-q}} &
\frac{\partial f}{\partial x_{n-q+k}}
\\
\frac{\partial\omega_1}{\partial x_1} & \cdots &
\frac{\partial\omega_1}{\partial x_{n-q}} &
\frac{\partial\omega_1}{\partial x_{n-q+k}}
\\
\cdot & \cdots & \cdot & \cdot
\\
\frac{\partial\omega_{n-q}}{\partial x_1} & \cdots &
\frac{\partial\omega_{n-q}}{\partial x_{n-q}} &
\frac{\partial\omega_{n-q}}{\partial x_{n-q+k}}
\end{array}
\right\vert
=
0
\ \ \ \ \ \ \ \ \ \ \ \ \
{\scriptstyle{(k\,=\,1\,\cdots\,q)}}
\]
are differential equations of that kind; because after the
substitution $f = \omega_j$, they go to an identity, and because they
are independent of each other, as shows their form. According to the
last proposition, these equations form a $q$-term complete system.
Thus:

\def\theproposition{9}\begin{proposition}
\label{Satz-9-S-89}
If $n-q$ independent functions of $n$ variables $x_1, \dots, x_n$ are
presented, then there always exists a determinate $q$-term complete
system in $x_1, \dots, x_n$ of which these functions constitute a
system of solutions.
\end{proposition}

\sectionengellie{\S\,\,\,24.}

However, for our goal, it does not suffice to have proved the
existence of the solutions to a complete system, it is on the contrary
rather necessary to deal more closely with the analytic properties of
these solutions.

On this purpose, we shall assume that in the presented $q$-term
complete system:
\[
Y_k(f)
=
\frac{\partial f}{\partial x_{n-q+k}}
+
\sum_{i=1}^{n-q}\,
\eta_{ki}\,\frac{\partial f}{\partial x_i}
=
0
\ \ \ \ \ \ \ \ \ \ \ \ \
{\scriptstyle{(k\,=\,1\,\cdots\,q)}},
\]
the analytic functions $\eta_{ ki}$ behave regularly in the
neighborhood of:
\[
x_1=\cdots=x_n=0.
\]
As quantities $x_1', \dots, x_{ n-1}'$, we now choose amongst the
solutions of the equation:
\[
\frac{\partial f}{\partial x_n}
+
\sum_{i=1}^{n-q}\,\eta_{qi}\,
\frac{\partial f}{\partial x_i}
=
0
\]
the formerly defined general solutions of this equation, 
relative to $x_n = 0$.

We know that these general solutions $x_1', \dots, x_{ n-1}'$ are
ordinary power series with respect to $x_1, \dots, x_n$ in a certain
neighborhood of $x_1 = \cdots = x_n = 0$, and that for $x_n = 0$, they
reduce to $x_1, \dots, x_{ n-1}$ respectively. The solutions already
mentioned above:
\[
x_{n-q+1}'
=
x_{n-q+1},\dots\dots,
x_{n-1}'
=
x_{n-1}
\]
are therefore general solutions.

Conversely, according to the observations in \S\,\,21, $x_1, \dots,
x_{ n-1}$ are also analytic functions of $x_1', \dots, x_{ n-1}', x_n$
and they behave regularly in the neighborhood of the system of values
$x_1 ' = \cdots = x_{ n-1} ' = x_n = 0$. Now, since in the new
complete system $Y_k' (f) = 0$, the coefficients $\eta_{ ki}'$ behave
regularly as functions of $x_1, \dots, x_n$, they will also be
ordinary power series with respect to $x_1', \dots, x_{ n-1}', x_n$ in
a certain neighborhood of $x_1' = 0, \dots, x_{ n-1}' = 0$, $x_n = 0$,
and in fact they will be, as we know, free of $x_n$.

Next, we determine the general solutions of the equation:
\[
\frac{\partial f}{\partial x_{n-1}'}
+
\sum_{i=1}^{n-q}\,\eta_{q-1,i}'\,
\frac{\partial f}{\partial x_i'}
=
0
\]
relative to $x_{ n-1} ' = 0$. These solutions, that we may call
$x_1'', \dots, x_{ n-2} ''$, behave regularly as functions of $x_1',
\dots, x_{ n-1}'$ in the neighborhood of $x_1' = \cdots = x_{ n-1} ' =
0$ and are hence also regular in the neighborhood of $x_1 = \cdots =
x_n = 0$, as functions of $x_1, \dots, x_n$. After the substitution
$x_{ n-1}' = 0$, the functions $x_1'', \dots, x_{ n-2}''$ reduce to
$x_1', \dots, x_{ n-2}'$, whence they reduce to $x_1, \dots, x_{ n-2}$
after the substitution $x_{ n-1} = x_n = 0$. The coefficients of the
next $(q-2)$-term complete system are naturally ordinary power series
in $x_1'', \dots, x_{ n-2}''$.

After iterating $q$ times these considerations, we obtain at the end
$n - q$ independent solutions: $x_1^{ (q)}, \dots, x_{ n-q}^{ (q)}$ of
our complete system. These are ordinary power series with respect to
the $x_i^{ (q-1)}$ in a certain neighborhood of $x_1^{ (q-1)} =
\cdots = x_{n-q + 1}^{ (q)} = 0$, and just in the same
way also, series with respect to the $x_i$ in a certain neighborhood
of $x_1 = \cdots = x_n = 0$. For $x_{ n-q + 1}^{ (q-1)} = 0$, the
$x_1^{ (q)}, \dots, x_{ n-q}^{ (q)}$ reduce respectively to $x_1^{
(q-1)}, \dots, x_{ n-q}^{ (q-1)}$ and hence to $x_1, \dots, x_{ n-q}$
for $x_{ n-q + 1} = \cdots = x_n = 0$. Therefore one has:
\[
x_i^{(q)}
=
x_i
+
\mathfrak{P}_i(x_1,\dots,x_n)
\ \ \ \ \ \ \ \ \ \ \ \ \
{\scriptstyle{(i\,=\,1\,\cdots\,n\,-\,q)}},
\]
where the $\mathfrak{P}_i$ all vanish for $x_{ n-q +1} =\cdots = x_n =
0$.

We call the solutions $x_1^{ (q)}, \dots, x_{ n-q}^{ (q)}$ of our
complete system its \terminology{general solutions relative to $x_{
n-q+1} = 0$, \dots, $x_n = 0$}.

We may state the gained result in a somewhat more general form by
introducing a general system of values: $x_1^0, \dots, x_n^0$ in place
of the special one:
\[
x_1=x_2=\cdots=x_n=0.
\]
Then we can say:

\def\thetheorem{12}\begin{theorem}
\label{Theorem-12-S-91}
Every $q$-term complete system:
\[
\frac{\partial f}{\partial x_{n-q+k}}
+
\sum_{i=1}^{n-q}\,
\eta_{ki}(x_1,\dots,x_n)\,
\frac{\partial f}{\partial x_i}
=
0
\ \ \ \ \ \ \ \ \ \ \ \ \
{\scriptstyle{(k\,=\,1\,\cdots\,q)}}
\]
whose coefficients $\eta_{ ki}$ behave regularly in the neighborhood
of $x_1 = x_1^0$, \dots, $x_n = x_n^0$ possesses $n - q$ independent
solutions $x_1^{ (q)}, \dots, x_{ n-q}^{ (q)}$ which behave regularly
in a certain neighborhood of $x_1 = x_1^0$, \dots, $x_n = x_n^0$ and
which in addition reduce to $x_1, \dots, x_{ n-q}$ respectively after
the substitution $x_{ n-q + 1} = x_{ n-q + 1}^0$, 
\dots, $x_n = x_n^0$.
\end{theorem}

The main theorem of the theory of complete systems has not been stated
in this precise version, neither by \name{Jacobi}, nor by
\name{Clebsch}. Nevertheless, this theorem is implicitly contained in
the known studies due to \name{Cauchy}, \name{Weierstrass},
\name{Briot} and \name{Bouquet}, \name{Kowalevsky} and \name{Darboux}
and which treat the question of existence of solutions to given
differential equations.

\sectionengellie{\S\,\,\,25.}

The theory of a single linear partial differential equation:
\[
X(f)
=
\sum_{i=1}^n\,\xi_i(x_1,\dots,x_n)\,
\frac{\partial f}{\partial x_i}
=
0
\]
stands, as already said, in the closest connection 
with the theory of the simultaneous system:
\[
\frac{\D\,x_1}{\xi_1}
=\cdots=
\frac{\D\,x_n}{\xi_n}.
\]
Moreover, something completely analogous takes place also for systems
of linear partial differential
\renewcommand{\thefootnote}{\fnsymbol{footnote}}
equations\footnote[1]{\,
This connection has been explained in detail for the
first time by \name{Boole}. Cf. also \name{A.~Mayer} about
unrestricted integrable differential equations, Math. Ann. Vol.~V.
}.

\label{S-91-sq}
Consider 
\renewcommand{\thefootnote}{\arabic{footnote}}
$q$ independent linear partial differential equations of
first order which, though, need not constitute a complete system.
For the sake of simplicity, we imagine that the equations are solved
with respect to $q$ of the differential quotients:
\def\theequation{4'}\begin{equation}
Y_k(f)
=
\frac{\partial f}{\partial x_{n-q+k}}
+
\sum_{i=1}^{n-q}\,\eta_{ki}(x_1,\dots,x_n)\,
\frac{\partial f}{\partial x_i}
=
0
\ \ \ \ \ \ \ \ \ \ \ \ \
{\scriptstyle{(k\,=\,1\,\cdots\,q)}}.
\end{equation}
If $\omega ( x_1, \dots, x_n)$ is a general solution of
these equations, then it holds:
\[
\frac{\partial\omega}{\partial x_{n-q+k}}
\equiv
-
\sum_{i=1}^{n-q}\,
\eta_{ki}\,\frac{\partial\omega}{\partial x_i}
\ \ \ \ \ \ \ \ \ \ \ \ \
{\scriptstyle{(k\,=\,1\,\cdots\,q)}},
\]
whence:
\[
\D\,\omega
=
\sum_{i=1}^{n-q}\,
\frac{\partial\omega}{\partial x_i}
\bigg\{
dx_i
-
\sum_{k=1}^q\,\eta_{ki}\,dx_{n-q+k}
\bigg\}.
\]
Consequently, the differential $d\omega$ vanishes identically
by virtue of the $n- q$ total differential equations:
\def\theequation{6}\begin{equation}
\D\,x_i
-
\sum_{k=1}^q\,
\eta_{ki}\,dx_{n-q+k}
=
0
\ \ \ \ \ \ \ \ \ \ \ \ \
{\scriptstyle{(i\,=\,1\,\cdots\,n\,-\,q)}}.
\end{equation}
But one calls every function of this nature an 
\terminology{integral function} of this
system of total differential equations.
Thus we can say:

\plainstatement{Every joint solution of the $q$ linear
partial differential equations~\thetag{ 4'} is an integral
function of the system of $n-q$ total differential 
equations~\thetag{ 6}. 

But conversely also, every integral function of the
system~\thetag{ 6} is a joint solution of the equations~\thetag{ 4'}. } 

Indeed, if $w (x_1, \dots, x_n)$ is an integral function 
of~\thetag{ 6}, then the expression:
\[
\D\,w
=
\sum_{i=1}^n\,\frac{\partial w}{\partial x_i}\,\D\,x_i
\]
vanishes identically by virtue of~\thetag{ 6}, that it is to say, 
one has:
\[
\sum_{k=1}^q\,
\bigg\{
\frac{\partial w}{\partial x_{n-q+k}}
+
\sum_{i=1}^{n-q}\,\eta_{ki}\,\frac{\partial w}{\partial x_i}
\bigg\}\,
\D\,x_{n-q+k}
\equiv
0,
\]
and from this it is evident that $w$ satisfies the equations~\thetag{
4} identically.

But now, the integration of the system~\thetag{ 6} is synonymous to
the finding \deutsch{Auffindung} of all its integral functions;
for, what does it mean to integrate the system~\thetag{ 6}?
Nothing else, as is known, but to determine all possible functions
$\rho_1, \dots, \rho_{ n-q}$ of $x_1, \dots, x_n$
which make the expression:
\[
\sum_{i=1}^{n-q}\,
\rho_i(x_1,\dots,x_n)
\bigg\{
\D\,x_i
-
\sum_{k=1}^q\,\eta_{ki}\,\D\,x_{n-q+k}
\bigg\}
\]
to become a complete differential, hence to be the 
differential of an integral function.

In other words:

\plainstatement{The integration of the system~\thetag{ 4'} of $q$
linear partial differential equations is accomplished if one
integrates the system of the $n-q$ total differential
equations~\thetag{ 6}, and conversely. }

This connection between the two systems~\thetag{ 4'} and~\thetag{ 6}
naturally presupposes that they (respectively) possess integral
functions; however, there still exists a certain connection also when
there are no joint solutions to~\thetag{ 4'}, and hence too, no
integral functions for~\thetag{ 6}. That will be considered on
another occasion (Chap.~6, p.~\pageref{S-105}).
\label{S-93}

Since the equations~\thetag{ 4'} have at most $n - q$ independent
solution in common, the system~\thetag{ 6} has at most $n - q$
integral functions. If it possess precisely $n - q$ such independent
solutions, the system~\thetag{ 6} is called
\terminology{unrestrictedly integrable}
\deutsch{unbeschränkt integrabel}; so 
this case occurs only when the $q$ equations~\thetag{ 4'} form a
$q$-term complete system.

We suppose that the system~\thetag{ 6} is unrestrictedly integrable,
or what is the same, that all the expressions $Y_k \big( Y_j ( f) \big) -
Y_j \big( Y_k ( f) \big)$ vanish identically.
Furthermore, we imagine that the $n - q$ general solutions:
\[
\omega_1(x_1,\dots,x_n),\,\dots,\,\,
\omega_{n-q}(x_1,\dots,x_n)
\]
of the complete system~\thetag{ 4'}
relative to:
\[
x_{n-q+1}=x_{n-q+1}^0,\dots,
x_n=x_n^0
\]
are determined. Then the equations:
\[
\omega_1(x_1,\dots,x_n)
=
a_1,\,\dots,\,\,
\omega_{n-q}(x_1,\dots,x_n)
=
a_{n-q}
\]
with the $n-q$ arbitrary constants $a_1, \dots, a_{ n-q}$ are called
the \terminology{complete integral equations} of the system~\thetag{
6}. These integral equations are obviously solvable with respect to
$x_1, \dots, x_{ n-q}$, for $\omega_1,
\dots, \omega_{ n-q}$ reduce to $x_1, \dots, 
x_{ n-q}$ (respectively) for:
\[
x_{n-q+1}
=
x_{n-q+1}^0,\,\dots,\,\,
x_n
=
x_n^0.
\]
Hence we obtain:
\[
x_i
=
\psi_i\big(x_{n-q+1},\dots,x_n,\,a_1,\dots,a_{n-q}\big)
\ \ \ \ \ \ \ \ \ \ \ \ \
{\scriptstyle{(i\,=\,1\,\cdots\,n\,-\,q)}}.
\]

One can easily see that the equations~\thetag{ 6}
become identities after the substitution $x_1 = \psi_1, \dots, \,
x_{ n-q} = \psi_{ n-q}$. Indeed, we at first have:
\[
x_i
-
\psi_i\big(x_{n-q+1},\dots,x_n,\,\omega_1,\dots,\omega_{n-q}\big)
\equiv
0
\ \ \ \ \ \ \ \ \ \ \ \ \
{\scriptstyle{(i\,=\,1\,\cdots\,n\,-\,q)}};
\]
so if we introduce $x_i - \psi_i$ in place of $f$
in $Y_k (f)$, we naturally again obtain an identically 
vanishing expression, and since all $Y_k ( \omega_1), \dots, 
Y_k ( \omega_{ n-q})$ are identically zero, it comes:
\[
\aligned
\eta_{ki}(x_1,\dots,x_n)
&
-
\frac{\partial}{\partial x_{n-q+k}}
\psi_i\big(
x_{n-q+1},\dots,x_n,\,\omega_1,\dots,\omega_{n-q}\big)
\equiv
0
\\
&
\ \ \ \ \ \
\ \ \ \ \ \ \ \ \ \ \ \ \
{\scriptstyle{(k\,=\,1\,\cdots\,q;\,\,i\,=\,1\,\cdots\,n\,-\,q)}}.
\endaligned
\]
If we here make the substitution $x_1 = \psi_1, \dots, \, 
x_{ n-q} = \psi_{ n-q}$, it comes:
\[
\aligned
\eta_{ki}\big(\psi_1,\dots,\psi_{n-q},\,x_{n-q+1},\dots,x_n\big)
-
\frac{\partial}{\partial x_{n-q+k}}\,
\psi_i
&
\big(x_{n-q+1},\dots,x_n,\,
\\
&\ \ \
a_1,\dots,a_{n-q}\big)
\equiv
0.
\endaligned
\]
We multiply this by $dx_{ n-q+k}$, we sum with respect to $k$ from $1$
to $q$ and we then realize that the expression:
\[
\D\,x_i
-
\sum_{k=1}^q\,\eta_{ki}(x_1,\dots,x_n)\,\D\,x_{n-q+k}
\]
does effectively vanish identically after the substitution $x_1 =
\psi_1, \dots, \, x_{ n-q} = \psi_{ n-q }$. If we yet add that, from
the equations $\omega_i = a_i$, we can always determine the $a_i$ so
that the variables $x_1, \dots, x_{ n-q }$ take prescribed initial
values $\overline{ x }_1, \dots, \overline{ x }_{ n-q}$ for $x_{ n-q +
1} = x_{ n-q + 1}^0, \dots, \, x_n = x_n^0$, then we can say:

\plainstatement{If the system of total differential 
equations~\thetag{ 6} is unrestrictedly integrable, then it is always
possible to determine analytic functions $x_1, \dots, x_{ n-q}$ of
$x_{ n-q+1}, \dots, x_n$ in such a way that the system~\thetag{ 6} is
identically satisfied and that $x_1, \dots, x_{ n-q}$ take prescribed
initial values for $x_{ n-q+1} = x_{ n-q+1}^0, \dots, \, x_n =
x_n^0$. }

\sectionengellie{\S\,\,\,26.}

\label{S-94}
For reasons of convenience, we introduce a few abbreviations useful
in the future.

Firstly, the parentheses around the $f$ in $X ( f)$ shall from now on
be frequently left out.

Further, since expressions of the form $X \big( Y ( f) \big) - Y \big(
X ( f) \big)$ shall occur always more frequently
\deutsch{immer haüfiger}, we want to write:
\[
X\big(Y(f)\big)
-
Y\big(X(f)\big)
=
XYf
-
YXf
=
\big\leftbracket X,\,Y\big\rightbracket;
\]
also, we shall not rarely employ the following language
\deutsch{Redeweise bedienen}: the expression, or the 
infinitesimal transformation 
$\big\leftbracket X, \, Y \big\rightbracket$ arises as the
``\terminology{composition}'' \deutsch{Zusammensetzung}, or the
``\terminology{combination}'' \deutsch{Combination}, of
$Xf$ and $Yf$.

Yet the following observation may also find its position at this
place: 

Between any three expressions $Xf, Yf, Zf$ there always
exists the following identity:
\def\theequation{7}\begin{equation}
\label{jacobi-identity}
\big\leftbracket\big\leftbracket X,\,
Y\big\rightbracket,\,Z\big\rightbracket
+
\big\leftbracket\big\leftbracket Y,\,
Z\big\rightbracket,\,X\big\rightbracket
+
\big\leftbracket\big\leftbracket Z,\,
X\big\rightbracket,\,Y\big\rightbracket
\equiv
0.
\end{equation}
This identity is a special case of the so-called \terminology{Jacobi
identity}, which we shall get to know later. Here, we want to content
ourselves with verifying the correctness of the special
identity~\thetag{ 7}; later at the concerned place
(cf. Volume~2), we shall enter into
the sense of the Jacobi identity.

One has obviously:
\[
\big\leftbracket\big\leftbracket X,\,Y\big\rightbracket,\,Z\big\rightbracket
=
XYZf
-
YXZf
-
ZXYf
+
ZYXf;
\]
if one permutes here circulary $Xf, Yf, Zf$ and then sums the three
obtained relations, one gets rid of all the terms in the 
right-hand side, and one receives the identity pointed out above.

This special Jacobi identity turns out to be extremely important 
in all the researches on transformation groups.

\smallercharacters{The above simple verification of the
\emphasis{special} Jacobi identity could be given for the first time
by \name{Engel}. }

\linestop


\chapter{New Interpretation of the Solutions of a Complete System}
\label{kapitel-6}
\chaptermark{New Interpretation of the Solutions of a Complete System}

\setcounter{footnote}{0}

\abstract*{??}

If, in the developments of the preceding chapter, we interpret the
expression $X ( f)$ as the symbol of an infinitesimal transformation,
or what is the same, as the symbol of a one-term group, then all what
has been said 
receives a new sense. If on the other hand, 
we interpret the variables $x_1,
\dots, x_n$ as point coordinates \deutsch{Punktcoordinaten} 
in a space of $n$ dimensions, the results obtained at that time also
receive a certain graphic nature
\deutsch{Anschaulichkeit}.

The goal of the present chapter is to present these two aspects
in details and then to put them in association; 
but for that, the introduction of various new concepts
turns out to be 
\renewcommand{\thefootnote}{\fnsymbol{footnote}}
necessary\footnote[1]{\, 
The formations of concepts \deutsch{Begriffsbildungen}
presented in this chapter have been developed by Lie in the Memoirs of
the Scientific Society of Christiania 1872, 1873, 1874 and 19 february
1875. Cf. also Math. Ann. Vols. VIII, IX and XI. 
}. 
\renewcommand{\thefootnote}{\arabic{footnote}}

\sectionengellie{\S\,\,\,27.}

Let $x_i' = f_i ( x_1, \dots, x_n)$ be a transformation
in the variables $x_1, \dots, x_n$ and let
$\Phi ( x_1, \dots, x_n)$ be an arbitrary 
function; now, if by chance this function is constituted
in such a way that the relation:
\[
\Phi\big(f_1(x),\dots,f_n(x)\big)
=
\Phi(x_1,\dots,x_n)
\]
holds identically, then we say:
\terminology{the function $\Phi ( x_1, \dots, x_n)$ admits
\deutsch{gestattet}
the transformation $x_i' = f_i ( x_1, \dots, x_n)$}, or: it
\terminology{allows} \deutsch{zulässt} this transformation; we also
express ourselves as follows:
\terminology{the function $\Phi ( x_1, \dots, x_n)$ remains
invariant through the mentioned transformation, it behaves as an
invariant with respect to this transformation}.

If a function $\Phi ( x_1, \dots, x_n)$ admits all the $\infty^r$
transformations $x_i' = f_i ( x_1, \dots, x_n, \, a_1, \dots, a_r)$ of
an $r$-term group, we say that it 
\label{S-95}
\terminology{remains invariant by
this group}, and that \terminology{it admits this group}; at the same
time, we call $\Phi$ an \terminology{absolute invariant}, or briefly
an \terminology{invariant} of the group.

Here, we restrict ourselves to one-term groups.
So consider a on-term group:
\[
X(f)
=
\sum_{i=1}^n\,\xi_i(x_1,\dots,x_n)\,
\frac{\partial f}{\partial x_i},
\]
the finite transformations of which have the form:
\[
x_i'
=
x_i
+
\frac{t}{1}\,\xi_i
+
\frac{t^2}{1\cdot 2}\,X(\xi_i)
+\cdots
\ \ \ \ \ \ \ \ \ \ \ \ \
{\scriptstyle{(i\,=\,1\,\cdots\,n)}}.
\]
We ask for all invariants of this group.

If $\Phi ( x_1, \dots, x_n)$ is such an invariant, 
then for every $t$, the equation:
\[
\Phi\big(x_1+t\xi_1+\cdots,\,\,\dots\dots,\,\,
x_n+t\xi_n+\cdots\big)
=
\Phi(x_1,\dots,x_n)
\]
must hold identically. If we expand here the left-hand side in a
series of powers of $t$ according to the general formula~\thetag{ 7}
in Chap.~\ref{one-term-groups}, p.~\pageref{S-52-eq-7}, 
then we obtain the condition:
\[
\Phi(x_1,\dots,x_n)
+
\frac{t}{1}\,X(\Phi)
+
\frac{t^2}{1\cdot 2}\,
X\big(X(\Phi)\big)
+\cdots
\equiv
\Phi(x_1,\dots,x_n)
\]
for every $t$. From this, it immediately follows that the expression
$X ( \Phi)$ must vanish identically, if $\Phi$ admits our one-term
group; as a result, we thus have a \emphasis{necessary} criterion for
the invariance of the function $\Phi$ by the one-term group $X ( f)$.

As we have seen earlier (cf. p.~\pageref{S-53}), the expression $X (
\Phi)$ determines the increase \deutsch{Zuwachs} that the function
$\Phi$ undergoes by the infinitesimal transformation $X ( f)$. Now,
since this increase $\delta \Phi = X ( \Phi) \delta t$ vanishes
together with $X ( \Phi)$, it is natural to introduce the following
language: \terminology{when the expression $X ( \Phi)$ vanishes
identically, we say that the function $\Phi$ admits the infinitesimal
transformation $X (f )$}.

Thus, our result above can also be enunciated as follows:

\plainstatement{For a function $\Phi$ of $x_1, \dots, x_n$ to
admit all transformations of the one-term group $X ( f)$, it is a
necessary condition that it admits the infinitesimal transformation $X
( f)$ of the concerned group}.
 
But it is easy to see that this necessary condition is at the same
time sufficient. Indeed, together with $X ( \Phi)$, all the
expressions $X \big( X ( \Phi) \big)$, $X \big( X \big( X ( \Phi )
\big) \big)$, etc. also vanish identically, and consequently,
the equation:
\[
\Phi(x_1',\dots,x_n')
=
\Phi(x_1,\dots,x_n)
+
\frac{t}{1}\,X(\Phi)
+\cdots
\]
reduces to $\Phi ( x') = \Phi ( x)$ for every value of $t$, hence with
that, it is proved that the function $\Phi ( x)$ admits all
transformations of the one-term group $X ( f)$. Now, since the
functions $\Phi ( x_1, \dots, x_n)$ for which the expression $X (
\Phi)$ vanishes identically are nothing but the solutions of the
linear partial differential equation $X ( f) = 0$, we therefore can
state the following theorem.

\def\thetheorem{13}\begin{theorem}
\label{Theorem-13-S-97}
The solutions of the linear partial differential equation:
\[
X(f)
=
\sum_{i=1}^n\,\xi_i(x_1,\dots,x_n)\,\frac{\partial f}{\partial x_i}
=
0
\]
are invariant by the one-term group $X ( f)$, and in fact, they 
are the only invariants by $X (f)$.
\end{theorem}

Of course, it should not be forgotten that the invariants of the
one-term group $X ( f)$ are also at the same time the invariants of
any one-term group of the form:
\[
\rho(x_1,\dots,x_n)\,\sum_{i=1}^n\,\xi_i(x_1,\dots,x_n)\,
\frac{\partial f}{\partial x_i}
=
\rho\,X(f),
\]
whichever $\rho$ can be as a function of its arguments.

This simply follows from the fact that one can multiply by any
arbitrary function $\rho$ of $x_1, \dots, x_n$ the equation $X ( f) =
0$ which defines the invariants in question. We can also express this
fact as follows: if a function of $x_1, \dots, x_n$ admits the
infinitesimal transformation $X ( f)$, then at the same time, it
admits every infinitesimal transformation $\rho ( x_1, \dots, x_n)
\, X ( f)$.

One sees that the two concepts of one-term group $X ( f)$ and of
infinitesimal transformation $X ( f)$ are more special
\deutsch{specieller} than the concept of linear partial differential
equation $X ( f) = 0$.

On the basis of our developments mentioned above, we can now express
in the following way our former observation that the joint solutions
of two equations $X_i ( f) = 0$ and $X_k ( f) = 0$ simultaneously
satisfy the third equation $X_i \big( X_k ( f) \big) - X_k \big( X_i (
f) \big) = 0$.

\def\theproposition{1}\begin{proposition}
\label{Satz-1-S-97}
If a function of $x_1, \dots, x_n$ admits the two infinitesimal
transformations $X_i ( f)$ and $X_k ( f)$ in these variables, then it
also admits the infinitesimal transformation $X_i \big( X_k ( f) \big)
- X_k \big( X_i ( f) \big)$.
\end{proposition}

Expressed differently:

\def\theproposition{2}\begin{proposition}
If a function of $x_1, \dots, x_n$ admits the two one-term groups
$X_i ( f)$ and $X_k ( f)$, then it also admits the one-term group $X_i
\big( X_k ( f) \big) - X_k \big( X_i ( f) \big)$.
\end{proposition}

If $\psi_1, \dots, \psi_{ n-q}$ is a system of independent solutions
of the $q$-term complete system $X_1 ( f) = 0, \dots, X_q ( f) = 0$,
then $\Omega ( \psi_1, \dots, \psi_{ n-q} )$ is the general form of a
solution of this complete system, whence $\Omega ( \psi_1, \dots,
\psi_{ n-q} )$ generally admits every infinitesimal transformation of
the shape:
\def\theequation{1}\begin{equation}
\chi_1(x_1,\dots,x_n)\, X_1(f)
+\cdots+
\chi_q(x_1,\dots,x_n)\, X_q(f),
\end{equation}
whatever functions $\chi_1, \dots, \chi_q$ one can choose. It is even
clear that aside from the ones written just now, there are no
infinitesimal transformations by which all functions of the form
$\Omega ( \psi_1, \dots, \psi_{ n-q} )$ remain \label{S-98-bis}
invariant; for we know
that the $q$-term complete system $X_1 ( f) = 0,
\dots, \, X_q ( f) = 0$ is characterized by its solutions
$\psi_1, \dots, \psi_{ n-q}$.

Naturally, the functions $\Omega$ also admit all finite
transformations of the one-term group~\thetag{ 1}. 
\label{S-98}
Besides, one can
indicate all finite transformations through which all functions
$\Omega ( \psi_1, \dots, \psi_{ n-q} )$ remain simultaneously
invariant. The form of these transformations obviously is:
\[
\aligned
&
\psi_k(x_1',\dots,x_n')
=
\psi_k(x_1,\dots,x_n)
\ \ \ \ \ \ \ \ \ \ \ \ \
{\scriptstyle{(k\,=\,1\,\cdots\,n\,-\,q)}}
\\
&
\mathfrak{U}_j(x_1',\dots,x_n',\,x_1,\dots,x_n)
=
0
\ \ \ \ \ \ \ \ \ \ \ \ \
{\scriptstyle{(j\,=\,1\,\cdots\,q)}},
\endaligned
\]
where the $\mathfrak{ U}_j$ are subjected to the condition that really
a transformation arises. Moreover, the $x_k'$ and the $x_k$ as well
must stay inside a certain region.

We want to say that \terminology{the system of equations:
\[
\pi_1(x_1,\dots,x_n)=0,
\,\dots,\,\,
\pi_m(x_1,\dots,x_n)=0
\]
admits the transformation $x_i' = f_i ( x_1, \dots, x_n)$ 
when the system of equations:
\[
\pi_1(x_1',\dots,x_n')=0,
\,\dots,\,\,
\pi_m(x_1',\dots,x_n')=0
\]
is equivalent to $\pi_1 ( x) = 0, \, \dots, \, \pi_m ( x) = 0$ after
the substitution $x_i' = f_i ( x_1, \dots, x_n)$, 
\label{S-98-ter}
hence when every
system of values $x_1, \dots, x_n$ which satisfies the $m$ equations
$\pi_\mu ( x) = 0$ also satisfies the $m$ equations: }
\[
\pi_1\big(f_1(x),\dots,f_n(x)\big)
=
0,\,\dots,\,\,
\pi_m\big(f_1(x),\dots,f_m(x)\big)
=
0.
\]

With the introduction of this definition, it is not even necessary to
assume that the $m$ equations $\pi_1 = 0, \dots, \pi_m = 0$ are
independent of each other, though this assumption will always be done
in the sequel, unless the contrary is expressly admitted.

From what precedes, it now comes immediately the

\def\theproposition{3}\begin{proposition}
If $W_1, W_2, \dots, W_m$ $(m \leqslant n - q)$ are arbitrary
solutions of the $q$-term complete system $X_1 (f ) = 0, \dots, X_q (
f) = 0$ in the $n$ independent variables $x_1, \dots, x_n$ and if
furthermore $a_1, \dots, a_m$ are arbitrarily chosen constants, then
the system of equation:
\[
W_1=a_1,\,\dots,\,\,W_m=a_m
\]
admits any one-term group of the form:
\[
\sum_{k=1}^q\,\chi_k(x_1,\dots,x_n)\,X_k(f),
\]
where it is understood that $\chi_1, \dots, \chi_q$ are arbitrary
functions of their arguments. 
\end{proposition}

\sectionengellie{\S\,\,\,28.}

In the previous paragraph, we have shown in a new light the
theory of the
integration of linear partial differential equations, in such a
way that we brought to connection the infinitesimal transformations
and the one-term groups. At present, we want to take another
route, we want to attempt to make accessible 
the clear, illustrated conception of this
theory of integration (and of what is linked with it), 
by means of manifold considerations
\deutsch{Mannigfaltigkeitsbetrachtungen}.

If we interpret $x_1, \dots, x_n$ as coordinates in an $n$-times
extended space $R_n$, then the simultaneous system:
\[
\frac{\D\,x_1}{\xi_1}
=\cdots=
\frac{\D\,x_n}{\xi_n}
\]
receives a certain illustrative sense; namely
it attaches to each point $x_1, \dots, x_n$ of the $R_n$ a certain 
direction \deutsch{Richtung}.

The integral equations of the simultaneous system determine $n-1$ of
the variables $x_1, \dots, x_n$, hence for instance $x_1, \dots, x_{
n-1}$, as functions of the $n$-th: $x_n$ and of the initial values
$\overline{ x}_1, \dots,
\overline{ x}_n$; 
consequently, after a definite choice of initial values, these
integral equations represent a determinate once-extended manifold
which we call an \terminology{integral curve} of the simultaneous
system. \emphasis{Every such integral curve comes into contact
\deutsch{berührt}, in each one of its points, with the direction
attached to the point}.

There is in total $\infty^{ n-1}$ different integral curves
of the simultaneous system, and in fact, through every point of the
$R_n$, there passes in general one integral curve.

If $\psi_1, \dots, \psi_{ n-1}$ are independent integral functions
of the simultaneous system, then all the integral curves are also
represented by means of the $n - 1$ equations:
\[
\psi_k(x_1,\dots,x_n)
=
C_k
\ \ \ \ \ \ \ \ \ \ \ \ \
{\scriptstyle{(k\,=\,1\,\cdots\,n\,-\,1)}},
\]
with the $n-1$ arbitrary constants $C_1, \dots, C_{ n-1}$. If one
sets an arbitrary integral function $\Omega ( \psi_1,
\dots, \psi_{ n-1} )$ to be equal to an arbitrary 
constant:
\[
\Omega(\psi_1,\dots,\psi_{n-1})
=
A,
\]
then one gets the equation of $\infty^1$ $(n-1)$-times extended
manifolds, which are entirely constituted of integral curves and in
fact, every such manifold is constituted of $\infty^{ n-2}$ different
integral manifolds. Lastly, if one sets in general $m \leqslant n-1$
independent integral functions $\Omega_1, \dots,
\Omega_m$ to be equal to arbitrary constants:
\[
\Omega_\mu(\psi_1,\dots,\psi_{n-1})
=
A_\mu
\ \ \ \ \ \ \ \ \ \ \ \ \
{\scriptstyle{(\mu\,=\,1\,\cdots\,m)}},
\]
then one obtains the analytic expression of a
family of $\infty^m$ $(n-m)$-times extended manifolds, 
each one of which consists of $\infty^{ n-m-1 }$ integral curves.

The integral functions of our simultaneous system are at the
same time the solutions of the linear partial differential 
equations:
\[
\sum_{i=1}^n\,\xi_i(x_1,\dots,x_n)\,
\frac{\partial f}{\partial x_i}
=
X(f)
=
0.
\]

Occasionally, we also call the integral curves of the simultaneous
system the \terminology{characteristics of the linear partial
differential equation $X ( f) = 0$}, if we take up again a
terminology introduced by \name{Monge} for $n = 3$. Using this way
of expressing, we can also say: every solution of the linear partial
differential equation $X ( f) = 0$ represents, when it is set equal to
a constant, a family of $\infty^1$ manifolds which consists of
$\infty^{ n-2}$ characteristics of $X ( f) = 0$.

Now, we imagine that we are given two linear partial differential
equations, say $X_1 ( f) = 0$ and $X_2 ( f) = 0$.

It is possible that the two equations have characteristics in
common. This case happens when there is an identity of the form:
\[
\chi_1(x_1,\dots,x_n)\,X_1(f)
+
\chi_2(x_1,\dots,x_n)\,X_2(f)
\equiv
0,
\]
without $\chi_1$ and $\chi_2$ both vanishing. Then obviously every
solution of $X_1 (f) = 0$ satisfies also $X_2 ( f) = 0$ and
conversely.

If the two equations $X_1 ( f) = 0$ and $X_2 ( f) = 0$ have distinct
characteristics, then they do not have all their solutions in common;
then the question is whether they in general possess solutions in
common, or, as we can now express: whether the characteristics of $X_1
( f) = 0$ can be gathered as manifolds which consist of
characteristics of $X_2 ( f) = 0$.

This question can be answered directly when one knows the 
characteristics of the two equations $X_1 ( f)= 0$ and
$X_2 ( f)= 0$; however, we do not want to 
halt here. 
In the sequel, we shall restrict ourselves to 
expressing in the language of the
theory of manifolds
\deutsch{Mannigfaltigkeitslehre} 
the former results which have been deduced by means of the analytic
method
\deutsch{durch analytische Methoden}.

Let the $q$ mutually independent equations:
\[
X_1(f)=0,\,
\dots,\,\,
X_q(f)=0
\]
form a $q$-term complete system; 
let $\psi_1, \dots, \psi_{ n-q}$ be independent
solution of it. Then the equations:
\[
\psi_1(x_1,\dots,x_n)=C_1,
\,\dots,\,\,
\psi_{n-q}(x_1,\dots,x_n)=C_{n-q}
\]
with the $n-q$ arbitrary constants $C_k$ represent a
family of $\infty^{ n-q}$ $q$-times extended manifolds, 
each one of which consists of $\infty^{ q-1}$ characteristics
of each individual equation amongst the $q$ equations
$X_1 (f) = 0, \dots, X_q ( f) = 0$. 
We call these $\infty^{ n-q}$ manifolds the \label{S-101}
\terminology{characteristic manifolds of the complete
system}.

If one sets any $n - q - m$ independent functions of
$\psi_1, \dots, \psi_{ n-q}$ equal to 
arbitrary constants:
\[
\Omega_1(\psi_1,\dots,\psi_{n-q})
=
A_1,\,\dots,\,\,
\Omega_{n-q-m}(\psi_1,\dots,\psi_{n-q})
=
A_{n-q-m},
\]
then one gets the analytic expression of a family of
$\infty^{ n-q-m}$ $(q+m)$-times extended manifolds, 
amongst which each individual one consists of $\infty^m$
characteristic manifolds. 

The equations of the $\infty^{ n-q}$ characteristic manifolds show
that every point of the $R_n$ belongs to one and to only one
characteristic manifold. Consequently, we can say that
\terminology{the whole $R_n$ is decomposed
\deutsch{zerlegt} in $\infty^{ n-q}$ $q$-times extended
manifolds, hence that our complete system defines a decomposition
\label{S-101-bis}
\deutsch{Zerlegung} of the space}.

Conversely, every decomposition of the $R_n$ in $\infty^{ n-q}$
$q$-times extended manifolds:
\[
\varphi_1(x_1,\dots,x_n)
=
A_1,\,\dots,\,\,
\varphi_{n-q}(x_1,\dots,x_n)
=
A_{n-q}
\]
can be defined by means of a $q$-term complete system; for $\varphi_1,
\dots, \varphi_{ n-q}$ necessarily are independent functions, whence
according to Proposition~9, p.~\pageref{Satz-9-S-89}, there is
a $q$-term complete system the most general solutions of which is an
arbitrary function of $\varphi_1, \dots, \varphi_{ n-q}$; this
complete system then defines the decomposition in question.

An individual linear partial differential equation $X ( f) = 0$
attaches to every point of the $R_n$ a certain direction. If one
has several such equations, for instance the following ones, which can
be chosen in a completely arbitrary way:
\[
X_k(f)
=
\sum_{i=1}^n\,\xi_{ki}(x_1,\dots,x_n)\,
\frac{\partial f}{\partial x_i}
=
0
\ \ \ \ \ \ \ \ \ \ \ \ \
{\scriptstyle{(k\,=\,1\,\cdots\,q)}},
\]
then each one of these equations associates to every point of the
space a direction of progress
\deutsch{Fortschreitungsrichtung}. For instance, 
the $q$ directions associated to the point $x_1^0, \dots, x_n^0$ are
determined by:
\[
\aligned
&
\delta x_1^0\,:\,\delta x_2^0\,:\,\cdots\,:\,\delta x_n^0
=
\xi_{k1}(x^0)\,:\,\xi_{k2}(x^0)\,:\,\cdots\,:\,\xi_{kn}(x^0)
\\
&
\ \ \ \ \ \ \ \ \ \ \ \ \ \ \ \ \ 
\ \ \ \ \ \ \ \ \ \ \ \ \ \ \ \ \ \
{\scriptstyle{(k\,=\,1\,\cdots\,q)}}.
\endaligned
\]
We call these $q$ directions in the chosen point
\terminology{independent of each other} 
\label{S-102}
if none of it can be linearly
deduced from the remaining ones, that is to say: if it is not possible
to indicate $q$ numbers $\lambda_1, \dots, \lambda_q$ which do not all
vanish although the $n$ equations:
\[
\lambda_1\xi_{1i}^0
+\cdots+
\lambda_q\xi_{qi}^0
\ \ \ \ \ \ \ \ \ \ \ \ \
{\scriptstyle{(i\,=\,1\,\cdots\,n)}}
\]
are satisfied.

From this, it follows that the $q$ equations $X_1 ( f) = 0, \dots, X_q
( f) = 0$ associate to \emphasis{every point} in general position $q$
independent directions when they are themselves mutually independent,
thus when the equation:
\[
\sum_{k=1}^q\,\chi_k(x_1,\dots,x_n)\,X_k(f)
=
0
\]
can be identically satisfied only for $\chi_1 = 0, \,
\dots, \, \chi_q = 0$.

If one wants to visualize geometrically what should be understood by
``independent directions'', one should best start in ordinary,
thrice-extended space $R_3$; for it is then really obvious. In a
point of the $R_3$ one calls two directions independent of each other
when they are generally distinct; three directions are independent
when they do not fall in a same plane passing through the point; there
is in general no more than three mutually independent directions in a
point of the $R_3$.

Accordingly, $q$ directions in a point of the $R_n$ are independent of
each other if and only if, when collected together, they are not
contained in any smooth manifold through this point which has less
than $q$ dimensions.

Every possible common solution of the $q$ equations:
\[
X_1(f)=0,\,\dots,\,\,
X_q(f)=0
\]
also satisfies all equations of the form:
\def\theequation{2}\begin{equation}
\sum_{k=1}^q\,\chi_k(x_1,\dots,x_n)\,X_k(f)
=
0.
\end{equation}
The totality of all these equations associates a whole family of
directions to any point of the $R_n$. If we assume from the beginning
that the equations $X_1 ( f) = 0, \, \dots, \, X_q ( f) = 0$ are
mutually independent, we do not restrict the generality of the
investigation; thus through the equations~\thetag{ 2}, we are given a
family of $\infty^{ q-1}$ different directions which are attached to
each point $x_1, \dots, x_n$. One easily realizes that these $\infty^{
q-1}$ directions in every point form a smooth bundle, and hence
determine a smooth $q$-times extended manifold passing through this
point, namely the smallest smooth manifold through the point which
contains the $q$ independent directions of the $q$ equations $X_1 ( f)
= 0, \dots, X_q ( f) = 0$.

Every possible joint solution of the equations:
\[
X_1(f)=0,\,\dots,\,\,X_q(f)=0
\]
satisfies also all equations of the form:
\[
X_k\big(X_j(f)\big)
-
X_j\big(X_k(f)\big)
=
0.
\]
These equations also do attach to any point $x_1, \dots, x_n$ certain
directions, but in general, the directions in question shall only
exceptionally belong to the above-mentioned bundle of $\infty^{ q-1}$
directions at the point $x_1, \dots, x_n$. It is only in one case
that at each point of the space, the directions attached to all
equations:
\[
X_k\big(X_j(f)\big)
-
X_j\big(X_k(f)\big)
=
0
\]
belong to the bundle in question, namely only 
if for every $k$ and $j$ there exists a
relation of the form:
\[
X_k\big(X_j(f)\big)
-
X_j\big(X_k(f)\big)
=
\sum_{s=1}^q\,\omega_{kjs}(x_1,\dots,x_n)\,X_s(f),
\]
that is to say, when the equations $X_1 ( f) = 0, \, 
\dots, \, X_q ( f) = 0$ do coincidentally form a 
$q$-term complete system. ---

We can also define by the equations:
\[
\aligned
\delta x_1\,:\,\cdots\,:\,\delta x_n
&
=
\sum_{k=1}^q\,\chi_k(x)\,\xi_{k1}(x)\,\colon\,\cdots\,\colon\,
\sum_{k=1}^q\,\chi_k(x)\,\xi_{kn}(x)
\\
&
\ \ \ \ \ \ \ \ \ \ \ \ \ \ \ \ 
{\scriptstyle{(k\,=\,1\,\cdots\,q)}}
\endaligned
\]
the bundle of directions which is determined at every point $x_1,
\dots, x_n$ by the equations $X_1 f = 0, \, \dots, \, X_q f = 0$. By
eliminating from this the arbitrary functions $\chi_1, \dots, \chi_q$,
or, what is the same, by setting equal to zero all the determinants in
$q+1$ columns of the matrix:
\[
\left\vert
\begin{array}{cccc}
\D\,x_1 & \D\,x_2 & \cdot\,\,\cdot & \D\,x_n
\\
\xi_{11} & \xi_{12} & \cdot\,\,\cdot & \xi_{1n}
\\
\cdot & \cdot & \cdot\,\,\cdot & \cdot
\\
\xi_{q1} & \xi_{q2} & \cdot\,\,\cdot & \xi_{qn}
\end{array}
\right\vert
\]
we obtain a system of $n-q$ independent total differential equations.
This system attaches to every point $x_1, \dots, x_n$ exactly the same
smooth bundle of $\infty^{ q-1}$ directions as did the
equations~\thetag{ 2}. It follows conversely that the totality of the
linear partial differential equations~\thetag{ 2} is also completely
determined by the system of total differential equations just
introduced.

The formulas become particularly convenient when one replaces the $q$
equations $X_1 ( f) = 0, \,
\dots, \, X_q ( f) = 0$ by $q$ other equations which are
resolved with respect to $q$ of the differential quotients $\frac{
\partial f}{ \partial x_i}$, for instance by the following $q$
equations:
\[
Y_k(f)
=
\frac{\partial f}{\partial x_{n-q+k}}
+
\sum_{i=1}^{n-q}\,
\eta_{ki}(x_1,\dots,x_n)\,
\frac{\partial f}{\partial x_i}
=
0
\ \ \ \ \ \ \ \ \ \ \ \ \ \ \ \ 
{\scriptstyle{(k\,=\,1\,\cdots\,q)}}.
\]

The totality of all the equations~\thetag{ 2} is equivalent
to the totality of all the equations:
\[
\sum_{k=1}^q\,\chi_k(x_1,\dots,x_n)\,
Y_k(f)
=
0;
\]
hence the directions which are attached to the point
$x_1, \dots, x_n$ are also represented by the equations:
\[
\aligned
&
\D\,x_1\,:\,\cdots\,:\,\D\,x_{n-q}\,:\,\D\,x_{n-q+1}\,:\,\cdots\,:\,
\D\,x_n=
\\
&
=
\sum_{k=1}^q\,\chi_k\,\eta_{k1}\,\colon\,\cdots\,\colon\,
\sum_{k=1}^q\,\chi_k\,\eta_{k,\,n-q}\,\colon\,\chi_1\,
\colon\,\cdots\,\colon\,\chi_q.
\endaligned
\]
If we therefore eliminate $\chi_1, \dots, \chi_q$, we obtain the
following system of total differential equations:
\def\theequation{3}\begin{equation}
\D\,x_i
-
\sum_{k=1}^q\,\eta_{ki}(x)\,\D\,x_{n-q+k}
=
0
\ \ \ \ \ \ \ \ \ \ \ \ \
{\scriptstyle{(i\,=\,1\,\cdots\,n\,-\,q)}}.
\end{equation}

We have already seen in the preceding chapter, 
from p.~\pageref{S-91-sq}
up to p.~\pageref{S-93}, that there is a connection between the system
of the linear partial differential equations $Y_1 (f) = 0, \, \dots,
\, Y_q ( f) = 0$ and the above system of total differential equations.
But at that time, we limited ourselves to the special case where the
$q$ equations:
\[
Y_1(f)=0,\,\dots,\,\,Y_q(f)=0
\]
possess solutions in common, and we showed that the determination of
these joint solutions amounts to the integration of the above total
differential equations.

However, in the developments carried out just now, the integrability
of the concerned system of differential equations is out of the
question. With that, the connection between the system of the linear
partial differential equations $Y_1 (f) = 0, \, \dots, \, Y_q ( f) =
0$ and the system of the total differential equations~\thetag{ 3} is
completely independent of the integrability of these
two systems; this connection is just based on the fact that the two
systems attach a single and the same smooth bundle of $\infty^{ q-1}$
directions. \label{S-105}

To conclude, we still make the following remark which seems most
certainly obvious, but nevertheless has to be done: if the $q$
equations $X_1 f = 0, \, \dots, \, X_q f = 0$ constitute a $q$-term
complete system, then at each point $x_1, \dots, x_n$, the
characteristic manifold of the complete system comes into contact with
the $\infty^{ q-1}$ directions that all equations of the form~\thetag{
2} attach to this point.

\sectionengellie{\S\,\,\,29.}

Lastly, it is advisable to combine the manifold-type considerations
\deutsch{Mannigfaltigkeitsbetrachtungen} of the
previous paragraph with the developments of
\S\,\,27. But before, we still pursue a bit the
manifold considerations.

Every transformation $x_i' = f_i (x_1, \dots, x_n)$ can be interpreted
as an operation which exchanges \deutsch{vertauscht} the points of the
$R_n$, as it transfers each point $x_1, \dots, x_n$ in the new
position $x_1' = f_1 ( x), \, \dots, \, x_n' = f_n ( x)$
(Chap.~\ref{three-principles-thought}, p.~\pageref{vertauscht}).

A system of $m$ independent equations: 
\[
\Omega_1(x_1,\dots,x_n)=0,
\,\dots,\,\,
\Omega_m(x_1,\dots,x_n)=0
\]
represents a $(n-m)$-times extended manifold of the $R_n$. We say
that this manifold \terminology{admits the transformation $x_i' = f_i
( x_1, \dots, x_n)$} if the system of equations $\Omega_1 = 0, \,
\dots, \, \Omega_m = 0$ admits this transformation. According to
\S\,\,27 it is the case when every system of values $x_1, \dots, x_n$
which satisfies the equations $\Omega_1 ( x) = 0, \, \dots, \,
\Omega_m ( x) = 0$ satisfies at the same time the equations:
\[
\Omega_1\big(f_1(x),\dots,f_n(x)\big)=0,
\,\dots,\,\,
\Omega_m\big(f_1(x),\dots,f_n(x)\big)=0.
\]
Hence we can also express ourselves as follows:

\plainstatement{The manifold:
\[
\Omega_1(x_1,\dots,x_n)=0,\,\dots,\,\,
\Omega_m(x_1,\dots,x_n)=0
\]
admits the transformation:
\[
x_1'
=
f_1(x_1,\dots,x_n),\,\dots,\,\,
x_n'
=
f_n(x_1,\dots,x_n)
\]
if every point $x_1, \dots, x_n$ of the manifold 
is transferred by this transformation to a
point $x_1', \dots, x_n'$ which likewise belongs
to the manifold. }

In Chapter~3, p.~\pageref{fortsch} sq.,
we have seen that through every point $x_1, \dots, x_n$ in general
position of $R_n$ there passes an integral curve
\deutsch{Bahncurve} of the infinitesimal transformation:
\[
X(f)
=
\sum_{i=1}^n\,\xi_i(x_1,\dots,x_n)\,
\frac{\partial f}{\partial x_i}.
\]
There, we defined this integral curve as the totality of all the
positions that the point $x_1, \dots, x_n$ can take by means of all
the $\infty^1$ transformations of the one-term group $X ( f)$;
besides, the point $x_1, \dots, x_n$ on the integral curve in question
could be chosen completely arbitrarily. From this, it results that, 
through every transformation of the one-term group $X ( f)$, 
every point $x_1, \dots, x_n$ stays on the integral curve passing
through it, whence 
\emphasis{every integral curve of the infinitesimal
transformation $X ( f)$ remains invariant by the $\infty^1$
transformations of the one-term group $X ( f)$. The same naturally
holds true for every manifold which consists of integral curves}.

But the integral curves of the infinitesimal transformation $X ( f)$
are nothing else than the integral curves of the simultaneous system:
\[
\frac{\D\,x_1}{\xi_1}
=\cdots=
\frac{\D\,x_n}{\xi_n}
\]
hence from this it again follows that the above-mentioned
characteristics of the linear partial differential equation
$X ( f) = 0$ coincide with the integral curves of the infinitesimal
transformation $X ( f)$. 

Earlier on (Chap.~\ref{one-term-groups}, p.~\pageref{fortsch}), 
we have emphasized that an
infinitesimal transformation $Xf$ attaches to every point $x_1, \dots,
x_n$ in general position a certain direction of progress, 
namely the one with which comes
into contact the integral curve passing through. Obviously this
direction of progress coincides with the direction that the equation
$X f$ associates to the point in question.

Several, say $q$, infinitesimal transformations:
\[
X_kf
=
\sum_{i=1}^n\,\xi_{ki}(x_1,\dots,x_n)\,
\frac{\partial f}{\partial x_i}
\ \ \ \ \ \ \ \ \ \ \ \ \
{\scriptstyle{(k\,=\,1\,\cdots\,q)}}
\]
determine $q$ different directions of progress at each 
point $x_1, \dots, x_n$ in general position; 
in accordance with what precedes, we call these
directions of progress \terminology{independent of
each other} if the equations $X_1 f = 0, \, \dots, 
\, X_q f = 0$ are mutually independent.

From this, it follows immediately that the $q$ infinitesimal
transformations $X_1 f, \dots, X_qf$ attach 
precisely $h \leqslant q$ independent directions of
progress when all the $(h+1)\times (h+1)$ determinants
of the matrix:
\[
\left\vert
\begin{array}{ccc}
\xi_{11} & \cdot\,\,\cdot & \xi_{1n}
\\
\cdot & \cdot\,\,\cdot & \cdot
\\
\xi_{q1} & \cdot\,\,\cdot & \xi_{qn}
\end{array}
\right\vert
\]
vanish identically, without all its $h \times h$ determinants 
doing so. ---

What we have in addition to say here can be summarized as a statement.

\def\theproposition{4}\begin{proposition}
If $q$ infinitesimal transformations $X_1 f, \dots, X_q f$ of the
$n$-times extended space $x_1, \dots, x_n$ are constituted in such a
way that the $q$ equations $X_1 f = 0, \, \dots, \, X_q f = 0$ are
independent of each other, then $X_1f, \dots, X_q f$ attach to every
point $x_1, \dots, x_n$ in general position $q$ independent directions
of progress; moreover, if the equations $X_1 f = 0, \, \dots, \, X_q f
= 0$ form a $q$-term complete system, then $X_1 f, \dots, X_q f$
determine a decomposition of the space in $\infty^{ n-q}$ $q$-times
extended manifolds, the characteristic manifolds of the complete
system. Each one of these manifolds comes into contact in each of its
points with the directions that $X_1 f, \dots, X_qf$ associate to the
point; each such manifold can be engendered by the $\infty^{ q-1}$
integral curves of an arbitrary infinitesimal transformation of the
form:
\[
\chi_1(x_1,\dots,x_n)\,X_1f
+\cdots+
\chi_q(x_1,\dots,x_n)\,X_qf,
\]
where it is understood that $\chi_1, \dots, \chi_q$ are arbitrary
functions of their arguments; lastly, each one of the discussed
manifolds admits all the transformations of any one-term group:
\[
\chi_1\,X_1f
+\cdots+
\chi_q\,X_qf.
\]
\end{proposition}

\linestop


\chapter{Determination of All Systems of Equations
\\
Which Admit Given Infinitesimal Transformations
}
\label{kapitel-7}
\chaptermark{Determination of All Systems of Equations
Which Admit Given Infinitesimal Transformations}

\setcounter{footnote}{0}

\abstract*{??}

At first, we shall define what should mean the phrase that
the system of equations:
\[
\Omega_1(x_1,\dots,x_n)=0,
\,\dots,\,\,
\Omega_{n-m}(x_1,\dots,x_n)=0
\]
admits the infinitesimal transformation $X ( f)$. Afterwards, we
shall settle the extremely important problem of determining all
systems of equations which admit given infinitesimal 
\renewcommand{\thefootnote}{\fnsymbol{footnote}}
transformations\footnote[1]{\,
Cf. Lie, 
Scientific Society of Christiania 1872--74, as also
Math. Ann. Vol.~XI, Vol.~XXIV, pp.~542--544.
}. 
\renewcommand{\thefootnote}{\arabic{footnote}}

But beforehand, we still want to observe the following:

We naturally consider only such system of equations:
\label{S-107-sq}
\[
\Omega_1=0,\,\dots,\,\,
\Omega_{n-m}=0
\]
that are really satisfied by certain systems of values $x_1, \dots,
x_n$; at the same time, we \emphasis{always} restrict ourselves to
systems of values $x_1, \dots, x_n$ in the neighbourhood of which the
functions $\Omega_1, \dots, \Omega_{ n-m}$ behave regularly. In
addition, we want \emphasis{once for all} agree on the following:
unless the contrary is expressly allowed, every system of equations
$\Omega_1 = 0, \, \dots, \, \Omega_{ n-m} = 0$ which we consider
should be constituted in such a way that not all $(n-m)\times (n-m)$
determinants of the matrix:
\def\theequation{1}\begin{equation}
\left\vert
\begin{array}{ccc}
\frac{\partial\Omega_1}{\partial x_1} & \cdot\,\,\cdot &
\frac{\partial\Omega_1}{\partial x_n}
\\
\cdot & \cdot\,\,\cdot & \cdot
\\
\frac{\partial\Omega_{n-m}}{\partial x_1} & \cdot\,\,\cdot &
\frac{\partial\Omega_{n-m}}{\partial x_n}
\end{array}
\right\vert
\end{equation}
vanish by means of $\Omega_1 = 0, \, \dots, \,
\Omega_{ n-m} = 0$. 
It is permitted to make this assumption, since a system of equations
which does not possess the demanded property can always be brought to
a form in which it satisfies the stated requirement.

\sectionengellie{\S\,\,\,30.}

In the variables $x_1, \dots, x_n$, let an infinitesimal 
transformation:
\[
Xf
=
\sum_{i=1}^n\,\xi_i(x_1,\dots,x_n)\,
\frac{\partial f}{\partial x_i}
\]
be given.
In the studies which relate to such an infinitesimal 
transformation, 
\emphasis{we shall always restrict ourselves to 
systems of values for which the $\xi_i$ behave regularly}.

If a system of equations $\Omega_1 = 0, \, \dots, \, \Omega_{
n-m} = 0$ admits all finite transformations:
\[
x_i'
=
x_i
+
e\,Xx_i
+\cdots
\ \ \ \ \ \ \ \ \ \ \ \ \
{\scriptstyle{(i\,=\,1\,\cdots\,n)}}
\]
of the one-term group $Xf$,
then the system of equations:
\[
\aligned
\Omega_k\big(x_1+e\,Xx_1
&
+\cdots,\,\,\dots,\,\,
x_n+e\,Xx_n+\cdots\big)
=
0
\\
&
\ \ \ \ \ \ \ \ \ \
{\scriptstyle{(k\,=\,1\,\cdots\,n\,-\,m)}},
\endaligned
\]
or, what is the same, the system:
\[
\Omega_k+e\,X\Omega_k+\cdots
=
0
\ \ \ \ \ \ \ \ \ \
{\scriptstyle{(k\,=\,1\,\cdots\,n\,-\,m)}}
\]
must be equivalent to the system of equations:
\[
\Omega_1=0,
\,\dots,\,\,
\Omega_{n-m}=0
\]
for all values of $e$. To this end it is obviously
\emphasis{necessary} that all $X \Omega_k$ vanish for the system of
values of $\Omega_1 = 0, \, \dots, \, \Omega_{ n-m} = 0$, hence that
the increment $X \Omega_k \delta t$ that $\Omega_k$ undergoes by the
infinitesimal transformation $x_i' = x_i + \xi_i \delta t$, vanishes
by means of $\Omega_1 = 0, \, \dots, \, \Omega_{ n-m} = 0$.

These considerations conduct us to set up the following definition:
\label{S-109-sq}

\plainstatement{A system of equations:
\label{vanish-by-means}
\[
\Omega_1(x_1,\dots,x_n)=0,
\,\dots,\,\,
\Omega_{n-m}=0
\]
admits the infinitesimal transformation $Xf$ as soon as all the $n -
m$ expressions $X \Omega_k$ vanish by means\footnote{\,
An example appearing p.~\pageref{x-2-0} below illustrates this
condition. 
} 
of the system of
equations. }

Then on the basis of this definition, the following obviously holds:

\def\theproposition{1}\begin{proposition}
If a system of equations admits all transformations of
the one-term group $Xf$, then in any case, it must
admit the infinitesimal transformation $Xf$.
\end{proposition}

In addition, we immediately realize that
\emphasis{a system of equations which admits the
infinitesimal transformation $Xf$ allows at the same time every
infinitesimal transformation of the form $\chi (x_1, \dots, x_n) \,
Xf$, provided of course that the function $\chi$ behaves regularly for
the system of values which comes into consideration in the concerned
system of equations. }

Without difficulty, one can see that the above definition is
independent of the choice of coordinates, so that every system of
equations $\Omega_1 = 0, \, \dots, \,
\Omega_{ n-m} = 0$ which admits the infinitesimal
transformation $Xf$ in the sense indicated above, must also admit it,
when new independent variables $y_1, \dots, y_n$ are introduced in
place of the $x$. The fact that this really holds true follows
immediately from the behaviour of the symbol $Xf$ after the
introduction of new variables. Here as always, it is assumed that the
$y$ are ordinary power series in the $x$, and that the $x$ are
ordinary power series in the $y$, for all the systems of values $x_1,
\dots, x_n$ and $y_1, \dots, y_n$ coming into consideration.

It yet remains to show that the definition set up above is also
independent of the form of the system of equations $\Omega_1 = 0, \,
\dots, \, \Omega_{ n-m} = 0$. Only when we will have proved this fact
shall the legitimacy of the definition be really established.

Now, in order to be able to perform this proof, we provide
at first a few general developments which are actually
already important and which will later find several applications.

Suppose that the system of equations $\Omega_1 = 0, \, \dots, \, 
\Omega_{ n-m} = 0$ admits the infinitesimal transformation
$Xf$. Since not all $m \times m$ determinants of the matrix~\thetag{ 1}
vanish by means of $\Omega_1 = 0, \, \dots, \, \Omega_{
n-m} = 0$, we can assume that the determinant:
\[
\sum\,\pm\,
\frac{\partial\Omega_1}{\partial x_1}
\,\dots\,
\frac{\partial\Omega_{n-m}}{\partial x_{n-m}}
\]
belongs to the nonvanishing ones. 
Then it is possible to resolve the equations 
$\Omega_k = 0$ with respect to $x_1, \dots, x_m$, and
this naturally delivers a system of equations:
\[
x_1
=
\varphi_1(x_{n-m+1},\dots,x_n),\,\dots,\,\,
x_{n-m}
=
\varphi_{n-m}(x_{n-m+1},\dots,x_n)
\]
which is analytically equivalent to the system $\Omega_1 = 0, \,
\dots, \, \Omega_{ n - m} = 0$.
Therefore, if by the sign $[ \, \, ]$ we denote the substitution $x_1
= \varphi_1, \, \dots, \, x_{ n-m} = \varphi_{ n-m}$, we have:
\[
\big[\Omega_1\big]
\equiv
0,
\,\dots,\,\,
\big[\Omega_{n-m}\big]
\equiv
0;
\]
and moreover, the fact that the system of equations $\Omega_1 = 0, \,
\dots, \, \Omega_{ n-m} = 0$ admits the infinitesimal transformation
$Xf$ is expressed by the identities:
\[
\big[X\Omega_1\big]
\equiv
0,
\,\dots,\,\,
\big[X\Omega_{n-m}\big]
\equiv 
0.
\]

Now, let $\Phi ( x_1, \dots, x_n)$ be an arbitrary function which 
behaves regularly for the system of values $x_1, \dots, x_n$
coming into consideration.
Then one has:
\[
\big[X\Phi\big]
=
\sum_{i=1}^n\,\big[Xx_i\big]\,
\bigg[
\frac{\partial\Phi}{\partial x_i}
\bigg],
\]
and on the other hand:
\[
\big[X[\Phi]\big]
=
\sum_{k=1}^{n-m}\,
\big[X\varphi_k\big]\,
\bigg[
\frac{\partial\Phi}{\partial x_k}
\bigg]
+
\sum_{\mu=1}^m\,
\big[Xx_{n-m+\mu}\big]\,
\bigg[
\frac{\partial\Phi}{\partial x_{n-m+\mu}}
\bigg],
\]
whence:
\def\theequation{2}\begin{equation}
\big[X\Phi\big]
=
\big[X[\Phi]\big]
+
\sum_{k=1}^{n-m}\,
\big[X(x_k-\varphi_k)\big]\,
\bigg[
\frac{\partial\Phi}{\partial x_k}
\bigg].
\end{equation}
If in place of $\Phi$ we insert one after the other the functions
$\Omega_1, \dots, \Omega_{ n-m}$, and if we take into account that
$\big[ \Omega_j \big]$ and also $X \big[ \Omega_j \big]$ plus $\big[ X
[ \Omega_j ] \big]$ vanish identically, then we find:
\[
\aligned
\big[X\Omega_j\big]
&
=
\sum_{k=1}^{n-m}\,
\big[X(x_k-\varphi_k)\big]\,
\bigg[
\frac{\partial\Omega_j}{\partial x_k}
\bigg]
\\
&\ \ \ \ \ \ \ \ \ \ \ \ \
{\scriptstyle{(j\,=\,1\,\cdots\,n\,-\,m)}}.
\endaligned
\]
Now, because $\big[ X \Omega_j \big]$ vanishes identically
in any case, while the determinant:
\[
\sum\,\pm\,
\bigg[
\frac{\partial\Omega_1}{\partial x_1}
\bigg]
\,\dots\,
\bigg[
\frac{\partial\Omega_{n-m}}{\partial x_{n-m}}
\bigg]
\]
does not vanish identically, it follows:
\[
\big[X(x_k-\varphi_k)\big]
\equiv
0 \ \ \ \ \ \ \ \ \ \ \ \ \
{\scriptstyle{(k\,=\,1\,\cdots\,n\,-\,m)}},
\]
whence the equation~\thetag{ 2} takes the form:
\label{relation-3-S-110}
\def\theequation{3}\begin{equation}
\big[X\Phi\big]
\equiv
\big[X[\Phi]\big].
\end{equation}
This formula which is valid for any function $\Phi ( x_1,
\dots, x_n)$ will later be very useful. 
Here, we need it only in the special case where $\Phi$ vanishes by
means of $\Omega_1 = 0, \, \dots, \, \Omega_{ n-m} = 0$; then $\big[
\Phi \big]$ is identically zero and likewise $\big[ X [
\Phi] \big]$; our formula hence shows that also $\big[ X \Phi \big]$
vanishes identically. In words, we can express this result as
follows:

\def\theproposition{2}\begin{proposition}
\label{Satz-2-S-111}
If a system of equations:
\[
\Omega_1(x_1,\dots,x_n)=0,
\,\dots,\,\,
\Omega_{n-m}(x_1,\dots,x_n)=0
\]
admits the infinitesimal transformation:
\[
Xf
=
\sum_{i=1}^n\,\xi_i(x_1,\dots,x_n)\,
\frac{\partial f}{\partial x_i}
\]
and if $V (x_1, \dots, x_n)$ is a function which vanishes by means of
this system of equation, then the function $X V$ also vanishes by
means of $\Omega_1 = 0, \, \dots, \, \Omega_{ n-m} = 0$.
\end{proposition}

Now, if $V_1 = 0, \, \dots, \, V_{ n-m} = 0$ is an arbitrary system
analytically equivalent to:
\[
\Omega_1=0,
\,\dots,\,\,
\Omega_{n-m}=0,
\]
then according to the proposition just stated, all the $n-m$
expressions $X V_k$ vanish by means of $\Omega_1 = 0, \, \dots, \,
\Omega_{ n-m} = 0$ and hence they also vanish by means of $V_1 = 0, 
\, \dots, \, V_{ n-m} = 0$; in other words:
the system of equations $V_1 = 0, \, \dots, \, V_{ n-m} = 0$ too
admits the infinitesimal transformation $Xf$.

Finally, as a result, it is established that our above definition for
the invariance of a system of equations by an infinitesimal
transformation is also independent of the form of this system of
equations. Therefore, the introduction of this definition is
completely natural
\deutsch{naturgemäss}.

We know that a system of equations can admit all transformations of
the one-term group $Xf$ only when it admits the infinitesimal
transformation $Xf$. But this necessary condition is at the same time
sufficient; indeed, it can be established that every system of
equations which admits the infinitesimal transformation $Xf$ generally
allows all the transformations of the one-term group $Xf$.

In fact, let the system of equation $\Omega_1 = 0, \, \dots, \,
\Omega_{ n-m} = 0$ admit the infinitesimal transformation $Xf$;
moreover, let $x_1 = \varphi_1, \, \dots, \, x_{ n-m} = \varphi_{
n-m}$ be a resolved form for the system of equations $\Omega_1 = 0, \,
\dots, \, \Omega_{ n-m} = 0$; 
lastly, let the substitution $x_\mu = \varphi_\mu$ be again denoted by
the sign $\big[ \, \, \big]$.

Under these assumptions, we have at first $\big[ \Omega_k \big] \equiv 
0$, then $\big[ X\Omega_k \big] \equiv 0$ and from the
Proposition~2 stated just now we obtain furthermore:
\[
\big[XX\Omega_k\big]
\equiv
0,
\ \ \ \ \ \ \ \ \ \
\big[XXX\Omega_k\big]
\equiv
0,\,\cdots.
\]
Consequently, the infinite series:
\[
\Omega_k
+
\frac{e}{1}\,X\Omega_k
+
\frac{e^2}{1\cdot 2}\,XX\Omega_k
+\cdots
\]
vanishes identically after the substitution $x_\mu = \varphi_\mu$,
whichever value the parameter $e$ can have. So for any $e$, the system of
equations:
\[
\Omega_k
+
e\,X\Omega_k
+\cdots
=
0
\ \ \ \ \ \ \ \ \ \ \ \ \
{\scriptstyle{(k\,=\,1\,\cdots\,n\,-\,m)}}
\]
will be satisfied by the systems of values of the system
of equations $\Omega_1 = 0, \, \dots, \,
\Omega_{ n-m} = 0$, and according to what
has been said earlier, this does mean nothing but: the system of
equations:
\[
\Omega_1=0,
\,\dots,\,\,
\Omega_{n-m}=0
\]
admits all transformations:
\[
x_i'
=
x_i
+
e\,Xx_i+\cdots
\ \ \ \ \ \ \ \ \ \ \ \ \
{\scriptstyle{(i\,=\,1\,\cdots\,n)}}
\]
of the one-term group $Xf$. With that, the assertion made 
above is proved; as a result, we have the

\def\thetheorem{14}\begin{theorem}
\label{Theorem-14-S-112}
The system of equations:
\[
\Omega_1(x_1,\dots,x_n)=0,
\,\dots,\,\,
\Omega_{n-m}=0
\]
admits all transformations of the one-term group $Xf$ if and only if
it admits the infinitesimal transformation $Xf$, that is to say, when
all the $n -m$ expressions $X \Omega_k$ vanish by means of $\Omega_1 =
0, \, \dots, \, \Omega_{ n-m} =0$.
\end{theorem}

This theorem is proved under the assumption, which we always make
unless something else is expressly notified, under the assumption
namely that not all $(n-m) \times (n - m)$ determinants of the
matrix~\thetag{ 1} vanish by means of $\Omega_1 = 0, \, \dots, \,
\Omega_{ n-m} = 0$. In addition, as already said, both the $\Omega_k$
and the $\xi_i$ must behave regularly for the systems of values $x_1,
\dots, x_n$ coming into consideration.

It can be seen that Theorem~14 does not hold true anymore when this
assumption about the determinant of the matrix~\thetag{ 1} is not
fulfilled. Indeed, we consider for instance the system of equations:
\label{x-2-0}
\[
\Omega_1=x_1^2=0,
\,\dots,\,\,
\Omega_{n-m}=x_{n-m}^2=0,
\]
by means of which all the $(n-m) \times ( n-m)$ determinants of the
matrix~\thetag{ 1} vanish. We find for these equations: $X \Omega_k =
2\, x_k \,X x_k$, hence all the $X \Omega_k$ vanish by means of
$\Omega_1 = 0, \, \dots, \,
\Omega_{ n-m} = 0$, whatever form $Xf$ can have. 
Consequently, if the Theorem~14 would be also true here, then the
system of equations:
\[
x_1^2
=
0,\,\dots,\,\,
x_{n-m}^2=0
\]
would admit any arbitrary one-term group $Xf$, which
obviously is not the case. 

From this we conclude the following: when the system of equations
$\Omega_1 = 0, \, \dots, \, \Omega_{ n-m} = 0$ brings to zero all the
$(n - m) \times ( n - m)$ determinants of the matrix~\thetag{ 1}, the
vanishing of all the $X \Omega_k$ by means of $\Omega_1 = 0, \, \dots,
\, \Omega_{ n-m} =0$ is of course necessary in order that this system
of equations admits the one-term group $Xf$, but however, it is not
sufficient.

Nevertheless, some general researches frequently
conduct one to systems of equations for which one has no
means to decide whether the repeatedly mentioned requirement is met.
Then how should one recognize that the system of equations
in question admits, or does not admit, a given one-term 
group?

In such circumstances, there is a criterion which 
is frequently of great help and which we now want to develop.

Let:
\[
\Delta_1(x_1,\dots,x_n)=0,
\,\dots,\,\,
\Delta_s(x_1,\dots,x_n)=0
\]
be a system of equations. We assume that the functions $\Delta_1,
\dots, \Delta_s$ behave regularly inside a certain region $B$, in the
neighbourhood of those systems of values $x_1, \dots, x_n$ which
satisfy the system of equations. However, we assume nothing about the
behaviour of the functional determinants of the $\Delta$'s; we do not
demand anymore that our $s$ equations are independent of each other,
so the number $s$ can even be larger than $n$ in certain
circumstances.

Moreover, let:
\[
Xf
=
\sum_{i=1}^n\,\xi_i(x_1,\dots,x_n)\,
\frac{\partial f}{\partial x_i}
\]
be an infinitesimal transformation and suppose that
amongst the system of values $x_1, \dots, x_n$ for which
$\xi_1, \dots, \xi_n$ behave regularly, there
exist some which satisfy the equations 
$\Delta_1 = 0, \, \dots, \, \Delta_s = 0$ and
which in addition belong to the domain $B$. 

Now, if under these assumptions the $s$ expressions
$X \Delta_\sigma$ can be represented as:
\[
X\Delta_\sigma
\equiv
\sum_{\tau=1}^s\,
\rho_{\sigma\tau}(x_1,\dots,x_n)\,
\Delta_\tau(x_1,\dots,x_n)\ \ \ \ \ \ \ \ \ \ \ \ \
{\scriptstyle{(\sigma\,=\,1\,\cdots\,s)}},
\]
and if at the same time all the $\rho_{ \sigma \tau}$ behave
regularly for the concerned system of values of:
\[
\Delta_1=0,
\,\dots,\,\,
\Delta_s=0,
\]
then our system of equations admits every transformation:
\[
x_i'
=
x_i
+
\frac{e}{1}\,Xx_i
+\cdots\ \ \ \ \ \ \ \ \ \ \ \ \
{\scriptstyle{(i\,=\,1\,\cdots\,n)}}
\]
of the one-term group $Xf$.

The proof of that is very simple. We have:
\[
\Delta_\sigma(x_1',\dots,x_n')
=
\Delta_\sigma(x_1,\dots,x_n)
+
\frac{e}{1}\,X\Delta_\sigma
+\cdots;
\]
but it comes:
\[
XX\Delta_\sigma
\equiv
\sum_{\tau=1}^s\,
\bigg\{
X\rho_{\sigma\tau}
+
\sum_{\pi=1}^s\,\rho_{\sigma\pi}\,\rho_{\pi\tau}
\bigg\}\,\Delta_\tau,
\]
where in the right-hand side the coefficients of the $\Delta$ again
behave regularly for the system of values of $\Delta_1 = 0, \, \dots,
\, \Delta_s = 0$.
In the same way, the $XXX \Delta_\sigma$ express linearly in terms of
$\Delta_1, \dots, \Delta_s$, and so on. In brief, we find:
\[
\Delta_\sigma(x_1',\dots,x_n')
=
\sum_{\tau=1}^s\,\psi_{\sigma\tau}
(x_1,\dots,x_n,\,e)\,
\Delta_\tau(x)\ \ \ \ \ \ \ \ \ \ \ \ \
{\scriptstyle{(\sigma\,=\,1\,\cdots\,s)}},
\]
where the $\psi_{ \sigma\tau}$ are ordinary power series in $e$ and
behave regularly for the system of values of $\Delta_1 = 0, \, \dots,
\, \Delta_s = 0$. From this, it results that every system of values
$x_1, \dots, x_n$ which satisfies the equations $\Delta_1 ( x) = 0, \,
\dots, \,
\Delta_s ( x) = 0$ also satisfies the equations:
\[
\Delta_\sigma
\Big(
x_1+\frac{e}{1}\,Xx_1+\cdots,
\,\,\dots,\,\,
x_n+\frac{e}{1}\,Xx_n+\cdots
\Big)\ \ \ \ \ \ \ \ \ \ \ \ \
{\scriptstyle{(\sigma\,=\,1\,\cdots\,s)}}
\]
so that the system of equations $\Delta_1 = 0, \, \dots, 
\, \Delta_s = 0$ really admits the one-term group $Xf$. 

As a result, we have the

\def\theproposition{3}\begin{proposition}
\label{X-Delta-sigma}
If, in the variables $x_1, \dots, x_n$, a system of equations:
\[
\Delta_1(x_1,\dots,x_n)=0,
\,\dots,\,\,
\Delta_s(x_1,\dots,x_n)=0
\]
is given, about which it is not assumed that its equations are
mutually independent, and even less that the $s \times s$ determinants
of the matrix:
\[
\left\vert
\begin{array}{ccc}
\frac{\partial\Delta_1}{\partial x_1} & \cdot\,\,\cdot &
\frac{\partial\Delta_1}{\partial xi_1}
\\
\cdot & \cdot\,\,\cdot & \cdot
\\
\frac{\partial\Delta_s}{\partial x_1} & \cdot\,\,\cdot &
\frac{\partial\Delta_s}{\partial x_n}
\end{array}
\right\vert
\]
vanish or do not vanish by means of $\Delta_1 = 0, \,
\dots, \, \Delta_s = 0$, then this system of equations
surely admits all the transformations of the one-term group $Xf$ when
the $s$ expressions $X \Delta_\sigma$ can be represented under the
form:
\[
X\Delta_\sigma
\equiv
\sum_{\tau=1}^s\,\rho_{\sigma\tau}(x_1,\dots,x_n)\,\Delta_\tau
\ \ \ \ \ \ \ \ \ \ \ \ \
{\scriptstyle{(\sigma\,=\,1\,\cdots\,s)}},
\]
and when at the same time the $\rho_{ \sigma \tau}$ behave regularly
for those systems of values $x_1, \dots, x_n$ which satisfy the system
of equations $\Delta_1 = 0, \, \dots, \,
\Delta_s = 0$.
\end{proposition}

\sectionengellie{\S\,\,\,31.}

In the preceding paragraph, we have shown that the determination of
all systems of equations which admit the one-term group $Xf$ amounts
to determining all systems of equations which admit the infinitesimal
transformation $Xf$. Hence the question arises to ask for 
\deutsch{Es entsteht daher die Frage nach} all systems
of equations $\Omega_1 = 0, \, \dots, \,
\Omega_{ n-m} = 0$ $(m \leqslant n)$ 
which admit the infinitesimal transformation\footnote{\,
(in the sense of the definition p.~\pageref{vanish-by-means})
}: 
\[
Xf
=
\sum_{i=1}^n\,\xi_i(x_1,\dots,x_n)\,
\frac{\partial f}{\partial x_i}.
\]

This question shall find its answer in the present paragraph.

Two cases must be distinguished, namely either
not all functions $\xi_1, \dots, \xi_n$ vanish
by means of:
\[
\Omega_1=0,
\,\dots,\,\,
\Omega_{n-m}=0
\]
or the equations $\xi_1 = 0, \, \dots, \, \xi_n = 0$ are
a consequence of:
\[
\Omega_1=0,
\,\dots,\,\,
\Omega_{n-m}=0.
\]

To begin with, we treat the first case.

Suppose to fix ideas that $\xi_n$ does not vanish by
means of $\Omega_1 = 0, \, \dots, \Omega_{ n-m} = 0$.
Then the concerned system of equations also admits the
infinitesimal transformation:
\[
Yf
=
\frac{1}{\xi_n}\,Xf
=
\frac{\xi_1}{\xi_n}\,
\frac{\partial f}{\partial x_1}
+\cdots+
\frac{\xi_{n-1}}{\xi_n}\,
\frac{\partial f}{\partial x_{n-1}}
+
\frac{\partial f}{\partial x_n}.
\]
If now $x_1^0, \dots, x_n^0$ is a system of values which satisfies the
equations $\Omega_k = 0$ and for which $\xi_n$ does not vanish, then
we may think that the general solutions of $Xf = 0$ relative to $x_n =
x_n^0$, or, what is the same, of $Yf$, are determined; these general
solutions, that we may call $y_1, \dots, y_{ n-1}$, behave regularly
in the neighbourhood of $x_1^0, \dots, x_n^0$ and are independent of
$x_n$. Hence if we introduce the new independent variables $y_1,
\dots, y_{ n-1}, \, y_n = x_n$ in place of the $x$, this will be an
allowed transformation. Doing so, $Yf$ receives the form $\frac{
\partial f}{\partial y_n}$, and the system of equations $\Omega_1 = 0,
\, \dots, \, \Omega_{ n-m} = 0$ is transferred to a new one:
\[
\overline{\Omega}_1(y_1,\dots,y_n)=0,
\,\dots,\,\,
\overline{\Omega}_{n-m}(y_1,\dots,y_n)=0
\]
which admits the infinitesimal transformation $\frac{ \partial f}{
\partial y_n}$.

From this, it follows that the system of equations $\overline{
\Omega}_k = 0$ is not solvable with respect to $y_n$. Indeed, if it
would be solvable with respect to $y_n$, and hence would yield $y_n -
\varphi ( y_1, \dots, y_{ n-1} ) = 0$, then the expression:
\[
Y(y_n-\psi)
=
\frac{\partial}{\partial y_n}
(y_n-\psi)
=
1
\]
would vanish by means of $\overline{ \Omega}_1 = 0, \,
\dots, \, \overline{ \Omega}_{ n-m} = 0$, which is
nonsensical. Consequently, $y_n$ can at most appear purely formally
in the equations $\overline{ \Omega}_k = 0$, that is to say, these
equations can in all circumstances be brought to a form such as they
represent relations between $y_1, \dots, y_{ n-1}$ alone\footnote{\,
For instance, the system of two equations $y_1 y_n = 0$ and
$y_1 = 0$ invariant by $\frac{ \partial }{ \partial y_n}$
amounts to just $y_1 = 0$. 
}. 
Here, the
form of these relations is subjected to no further restriction.

If we now return to the initial variables, we immediately realize that
the system of equations $\Omega_k = 0$ can be expressed by means of
relations between the $n-1$ independent solutions $y_1, \dots, y_{
n-1}$ of the equation $Xf = 0$. This outcome is obviously independent
of the assumption that $\xi_n$ itself should not vanish by means of
$\Omega_k= 0$; we therefore see that every system of equations which
admits the infinitesimal transformation $Xf$ and which does not
annihilate $\xi_1, \dots, \xi_n$ is represented by relations between
the solutions of $Xf = 0$. On the other hand, we know that completely
arbitrary relations between the solutions of $Xf = 0$ do represent a
system of equations which admits not only the infinitesimal
transformation, but also all transformations of the one-term group
$Xf$ (Chap.~\ref{kapitel-6}, 
Proposition~3, p.~\pageref{S-98}). Consequently, this
confirms the previously established result that our system of
equations $\Omega_k = 0$ admits the one-term group $Xf$.

We now come to the second of the above two distinguished cases;
naturally, this case can occur only when there are in general systems
of values $x_1, \dots, x_n$ for which all the $n$ functions $\xi_i$
vanish.

If $x_1^0, \dots, x_n^0$ is an arbitrary system of values for which
all $\xi_i$ vanish, then the transformation of our one-term group:
\[
x_i'
=
x_i
+
\frac{e}{1}\,\xi_i
+
\frac{e^2}{1\cdot 2}\,X\xi_i
+\cdots\ \ \ \ \ \ \ \ \ \ \ \ \
{\scriptstyle{(i\,=\,1\,\cdots\,n)}}
\]
reduces after the substitution $x_i = x_i^0$ to:
\[
x_1'
=
x_1^0,
\,\dots,\,\,
x_n=x_n^0,
\]
and we express this as: the system of values $x_1^0, \dots, x_n^0$
remains invariant by all transformations of the one-term group $Xf$.
For this reason, it comes that every system of equations of the form:
\[
\xi_1=0,
\,\dots,\,\,\xi_n=0,
\ \ \ \ \ \
\psi_1(x_1,\dots,x_n)=0,
\ \ \
\psi_2(x_1,\dots,x_n)=0,\,\dots,
\]
admits the one-term group $Xf$, whatever systems
of values $x_1, \dots, x_n$ are involved in it.
But such a form embraces\footnote{\,
In Lie's thought, \label{explain-embrace} a first system of equations
\terminology{embraces} (verb: \deutschplain{umfassen}) a second system
of equations when the first zero-set is larger than the second one so
that the first system \emphasis{implies} the second one, at least
locally and generically, and perhaps after some allowed algebraic
manipulations. Nothing really more precise about this notion will
come up later, and certainly nothing approaching either the
Nullstellensatz or some of the concepts of the so-called theory of
complex spaces.
} 
\deutsch{umfasst} all systems of equations which bring $\xi_1, \dots,
\xi_n$ to zero; so as a result, the second one of the two previously
distinguished cases is settled.

We summarize the gained result in the

\def\thetheorem{15}\begin{theorem}
\label{Theorem-15-S-117}
There are two sorts of systems of equations which admit the
infinitesimal transformation:
\[
Xf
=
\sum_{i=1}^n\,\xi_i(x_1,\dots,x_n)\,
\frac{\partial f}{\partial x_i},
\]
and hence in general admit all transformations of the one-term group
$Xf$. The systems of equations of the first sort are represented by
completely arbitrary relations between the solutions of the linear
partial differential equation $Xf = 0$. The systems of equations of
the second sort have the form:
\[
\xi_1=0,
\,\dots,\,\,\xi_n=0,
\ \ \ \ \ \
\psi_1(x_1,\dots,x_n)=0,
\ \ \
\psi_2(x_1,\dots,x_n)=0,\,\dots,
\]
in which the $\psi$ are absolutely arbitrary, 
provided of course that there are systems of
values $x_1, \dots, x_n$ which satisfy the
equations in question.
\end{theorem} 

We yet make a brief remark on this.

Let $C_1, \dots, C_{ n-m}$ be arbitrary constants, and let $\Omega_1,
\dots, \Omega_{ n-m}$ be functions of $x$, which however are
free of the $C$; lastly, suppose that each system of equations of the
form:
\[
\Omega_1(x_1,\dots,x_n)=C_1
\,\dots,\,\,
\Omega_{n-m}(x_1,\dots,x_n)=C_{n-m}
\]
admits the infinitesimal transformation $Xf$, which is free of $C$. 
Under these assumptions, the $n -m$ expression $X \Omega_k$
must vanish by means of:
\[
\Omega_1=C_1,
\,\dots,\,\,
\Omega_{n-m}=C_{n-m},
\]
and certainly, for all values of the $C$. But since the $X\Omega_k$
are all free of the $C$, this is only possible when the $X \Omega_k$
vanish identically, that is to say, when the $\Omega_k$ are solutions
of the equation $Xf = 0$. Thus, the following holds.

\def\theproposition{4}\begin{proposition}
If the equations:
\[
\Omega_1(x_1,\dots,x_n)=C_1,
\,\dots,\,\,
\Omega_{n-m}(x_1,\dots,x_n)=C_{n-m}
\]
with the arbitrary constants $C_1, \dots, C_{ n-m}$ represent a system
of equations which admits the infinitesimal transformation $Xf$, then
$\Omega_1, \dots, \Omega_{ n-m}$ are solutions of the differential
equation $Xf$, which is free of the $C$.
\end{proposition}

\sectionengellie{\S\,\,\,32.}

\plainstatement{At present, consider $q$ arbitrary infinitesimal
transformations:
\[
\label{S-118-sq}
X_kf
=
\sum_{i=1}^n\,\xi_{ki}(x_1,\dots,x_n)\,
\frac{\partial f}{\partial x_i}
\ \ \ \ \ \ \ \ \ \ \ \ \
{\scriptstyle{(k\,=\,1\,\cdots\,q)}}
\]
and ask for the systems of equations $\Omega_1 = 0, \, \dots, \,
\Omega_2 = 0, \, \dots$ which admit all these infinitesimal
transformations}. As always, we restrict ourselves here to systems of
values $x_1, \dots, x_n$ for which all the $\xi_{ ki}$ behave
regularly.

At first, it is clear that every sought system of equations
also admits all infinitesimal transformations
of the form:
\[
\sum_{k=1}^q\,\chi_k(x_1,\dots,x_n)\,X_kf,
\]
provided that $\chi_1, \dots, \chi_q$ behave regularly for the systems
of values $x_1, \dots, x_n$ which satisfy the system of equations.
From the previous paragraph, we moreover see that each such system of
equations also admits all finite transformations of the one-term groups
that arise from the discussed infinitesimal 
transformations\footnote{\,
In fact, Proposition~3 p.~\pageref{X-Delta-sigma} gives a more
precise statement. 
}. 

Furthermore, we remember Chap.~\ref{kapitel-6}, Proposition~1,
p.~\pageref{Satz-1-S-97}. At that time, we saw that every
function of $x_1, \dots, x_n$ which admits the two infinitesimal
transformations $X_1 f$ and $X_2 f$ also allows the transformation
$X_1 X_2 f - X_2 X_1 f$. Exactly the same property also holds true for
every system of equations which allows the two infinitesimal
transformations $X_1 f$ and $X_2 f$.

In fact, assume that the system of equations:
\[
\Omega_k(x_1,\dots,x_n)=0\ \ \ \ \ \ \ \ \ \ \ \ \
{\scriptstyle{(k\,=\,1\,\cdots\,n\,-\,m)}}
\]
admits the two infinitesimal transformations $X_1 f$ and $X_2 f$, so
that all the expressions $X_1 \Omega_j$ and $X_2 \Omega_j$ vanish by
means of the system $\Omega_1 = 0, \, \dots, \, \Omega_{ n-m} = 0$.
Then according to Proposition~2 p.~\pageref{Satz-2-S-111}, the
same holds true for all expressions $X_1 X_2 \Omega_j$ and $X_2 X_1
\Omega_j$, so that each $X_1 X_2 \Omega_j - X_2 X_1 \Omega_j$ vanishes
by means of the system of equations $\Omega_1 = 0, \, \dots, \,
\Omega_{ n-m} = 0$. As a result, we have the

\def\theproposition{5}\begin{proposition}
\label{Satz-5-S-118}
If a system of equations:
\[
\Omega_1(x_1,\dots,x_n)=0,
\,\dots,\,\,
\Omega_{n-m}(x_1,\dots,x_n)=0
\]
admits the two infinitesimal transformations $X_1 f$ and $X_2 f$, then
it also admits the infinitesimal transformation $X_1 X_2 f - X_2 X_1
f$.
\end{proposition}

We now apply this proposition similarly as we did at a previous time
in the Proposition~1 on p.~\pageref{Satz-1-S-84}, where the
question was to determine the joint solutions of $q$ given equations
$X_1 f = 0, \, \dots, \, X_q f = 0$. At that time, we reduced the
stated problem to the determination of the solutions of a complete
system. Now, we proceed as follows.

We form all infinitesimal transformations:
\[
X_kX_jf
-
X_jX_kf
=
\big[X_k,\,X_j\rightbracket
\]
and we ask whether the linear partial differential equations $\leftbracket
X_k, \, X_j \rightbracket = 0$ are a consequence of $X_1 f = 0, \, \dots, \,
X_q f = 0$. If this is not the case, then we add all transformations
$\leftbracket X_k, \, X_j \rightbracket$ to the infinitesimal 
\renewcommand{\thefootnote}{\fnsymbol{footnote}}
transformations $X_1 f,
\dots, X_q f$\footnote[1]{\,
In practice, one will in general not add all infinitesimal
transformations $\leftbracket X_k, \, X_j \rightbracket$, but only a certain number
amongst them.
}, 
which is permitted, 
\renewcommand{\thefootnote}{\arabic{footnote}}
since every system of equations which admits $X_1f, \dots, X_q f$ also
allows $\leftbracket X_k, \, X_j \rightbracket$. At present, we treat the
infinitesimal transformations taken together:
\[
X_1f,\,\dots,\,X_qf,\,\,
\leftbracket X_k,\,X_j\rightbracket\ \ \ \ \ \ \ \ \ \ \ \ \
{\scriptstyle{(k,\,j\,=\,1\,\cdots\,q)}}
\]
exactly as we did at first with $X_1f, \dots, X_q f$, that is to say,
we form all infinitesimal transformations:
\[
\leftbracket\leftbracket X_k,\,X_j\rightbracket,\,X_l\rightbracket,
\ \ \ \ \ \ \ \ \ \ \
\leftbracket \leftbracket X_k,\,X_j\rightbracket,\,\leftbracket X_h,\,X_l\rightbracket\rightbracket,
\]
and we ask whether the equations obtained by setting these expressions
equal to zero are a consequence of $X_k f = 0$, $\leftbracket X_k, \,
X_j\rightbracket = 0$ ($k, \, j = 1, \dots, q$). If this is not the case, we
add all the found infinitesimal transformations to $X_k f$, $\leftbracket 
X_k, \, X_j \rightbracket$ ($k, \, j = 1, \dots, q$).

We continue in this way, and so at the end we must obtain a series of
infinitesimal transformations:
\[
X_1f,\,\dots,\,X_qf,\,\,
X_{q+1}f,\,\dots,\,X_{q'}f
\ \ \ \ \ \ \ \ \ \ \ \ \
{\scriptstyle{(q'\,\geqslant\,q)}}
\]
which is constituted in such a way that every equation:
\[
\leftbracket X_k,\,X_j\rightbracket
=
0
\ \ \ \ \ \ \ \ \ \ \ \ \
{\scriptstyle{(k,\,j\,=\,1\,\cdots\,q')}}
\]
is a consequence of $X_1 f = 0, \, \dots, \, X_{ q'} f = 0$. The
equations:
\[
X_1f=0,
\,\dots,\,\,
X_{q'}f=0
\]
then define a complete system with $q'$ or less terms. Hence we have
the

\def\thetheorem{16}\begin{theorem}
The problem of determining all systems of equations which admit $q$
given infinitesimal transformations $X_1 f, \dots, X_q f$ can always
be led back to the determination of all systems of equations which,
aside from $X_1 f, \dots, X_q f$, admit yet certain further
infinitesimal transformations:
\[
X_{q+1}f,\,\dots,\,X_{q'}f
\ \ \ \ \ \ \ \ \ \ \ \ \
{\scriptstyle{(q'\,\geqslant\,q)}},
\] 
where now the equations:
\[
X_1f=0,\,\dots,\,\,
X_qf=0,\ \
X_{q+1}f=0,\,\dots,\,\,
X_{q'}f=0
\]
define a complete system which has as many terms as
there are independent equations in it.
\end{theorem}

Thus we can from now on limit ourselves to the following
more special problem:

\plainstatement{Consider $q$ infinitesimal transformations
in the variables $x_1, \dots, x_n$:
\[
X_kf
=
\sum_{i=1}^n\,\xi_{ki}(x_1,\dots,x_n)\,
\frac{\partial f}{\partial x_i}\ \ \ \ \ \ \ \ \ \ \ \ \
{\scriptstyle{(k\,=\,1\,\cdots\,q)}}
\]
with the property that amongst the $q$ equations:
\[
X_1f=0,\,\dots,\,\,X_qf=0,
\]
there are exactly $p \leqslant q$ equations which are mutually
independent, and with the property that any $p$ independent ones amongst
these $q$ equations form a $p$-term complete system which belongs to
$X_1 f = 0, \dots, X_q f = 0$. To seek all systems of equations in
$x_1, \dots, x_n$ which admit the infinitesimal transformations $X_1f,
\dots, X_q f$. }

It is clear that we can suppose without loss of generality that $X_1f,
\dots, X_q f$ are independent of each other as \emphasis{infinitesimal
transformations}. Furthermore, we notice that under the assumptions
of the problem, the $(p+1) \times (p+1)$ determinants of the matrix:
\def\theequation{4}\begin{equation}
\left\vert
\begin{array}{ccc}
\xi_{11} & \cdot\,\,\cdot & \xi_{1n}
\\
\cdot & \cdot\,\,\cdot & \cdot
\\
\xi_{q1} & \cdot\,\,\cdot & \xi_{qn}
\end{array}
\right\vert
\end{equation}
vanish \emphasis{all identically}, whereas not all $p \times p$
determinants do.

The first step towards the solution of our problem is to distribute
the systems of equations which admit the $q$ infinitesimal
transformations $X_1f, \dots, X_q f$ in two separate
classes \deutsch{in zwei getrennte Classen}; as a principle of
classification, we here take the behaviour of the $p \times p$
determinants of~\thetag{ 4}.

\plainstatement{In the first class, we reckon all systems of equations
by means of which not all $p \times p$ determinants of the
matrix~\thetag{ 4} vanish. \label{S-120}

In the second class, we reckon all systems of equations 
by means of which the $p \times p$ determinants in 
question all vanish. }

We now examine the two classes one after the other.

\sectionengellie{\S\,\,\,33.}

Amongst the $q$ equations $X_1 f = 0, \, 
\dots, \, X_q f = 0$, we choose any $p$ equations
which are independent of each other; to fix ideas, 
let $X_1 f = 0, \, \dots, \, X_p f = 0$ be such
equations, so that not all $p \times p$ determinants of the
matrix:
\def\theequation{5}\begin{equation}
\left\vert
\begin{array}{ccc}
\xi_{11} & \cdot\,\,\cdot & \xi_{1n}
\\
\cdot & \cdot\,\,\cdot & \cdot
\\
\xi_{p1} & \cdot\,\,\cdot & \xi_{pn}
\end{array}
\right\vert
\end{equation}
vanish identically. We now seek, amongst the systems of equations
$\Omega_1 = 0, \, \dots, \, \Omega_{ n-m} = 0$ 
of the first species, those which do not cancel
the independence of the equations $X_1 f = 0, \, \dots, \, X_p f = 0$,
hence those which do not annihilate all $p \times p$ determinants of
the matrix~\thetag{ 5}. Here, we want at first to make the special
assumption that a definite $p \times p$ determinant of the
matrix~\thetag{ 5}, say the following one:
\[
D
=
\left\vert
\begin{array}{ccc}
\xi_{1,n-p+1} & \cdot\,\,\cdot & \xi_{1n}
\\
\cdot & \cdot\,\,\cdot & \cdot
\\
\xi_{p,n-p+1} & \cdot\,\,\cdot & \xi_{pn}
\end{array}
\right\vert
\]
neither vanishes identically, nor vanishes by means of
$\Omega_1 = 0, \, \dots, \, \Omega_{ n-m} = 0$.

Under the assumptions made, there are identities of the form:
\[
X_{p+j}f
\equiv
\sum_{\pi=1}^p\,\chi_{j\pi}(x_1,\dots,x_n)\,
X_\pi f
\ \ \ \ \ \ \ \ \ \ \ \ \
{\scriptstyle{(j\,=\,1\,\cdots\,q\,-\,p)}}.
\]
For the determination of the functions $\chi_{ j\pi}$, we
have here the equations: \label{S-121}
\[
\sum_{\pi=1}^p\,
\xi_{\pi\nu}\,\chi_{j\pi}
=
\xi_{p+j,\,\nu}
\ \ \ \ \ \ \ \ \ \ \ \ \
{\scriptstyle{(\nu\,=\,1\,\cdots\,n;\,j\,=\,1\,\cdots\,q\,-\,p)}};
\]
now since the determinant $D$ does not vanish by means of $\Omega_1 =
0, \, \dots, \, \Omega_{ n-m} = 0$, we realize that the $\chi_{
j\pi}$ behave regularly for the system of values of $\Omega_1 = 0, \,
\dots, \, \Omega_{ n-m} = 0$, so that 
we can, under the assumptions made, leave out the infinitesimal
transformations $X_{ p+1} f, \dots, X_q f$; because, if the system of
equations $\Omega_1 = 0, \, \dots, \,
\Omega_{ n - m} = 0$ admits the transformations 
$X_1 f, \dots, X_p f$, then it automatically admits also $X_{ p+1} f ,
\dots, X_q f$.

We replace the infinitesimal transformations $X_1 f, \dots, 
X_p f$ by $p$ other infinitesimal transformations 
of the specific form:
\[
Y_\pi f
=
\frac{\partial f}{\partial x_{n-p+\pi}}
+
\sum_{i=1}^{n-p}\,
\eta_{\pi i}(x_1,\dots,x_n)\,
\frac{\partial f}{\partial x_i}
\ \ \ \ \ \ \ \ \ \ \ \ \
{\scriptstyle{(\pi\,=\,1\,\cdots\,p)}}.
\]
We are allowed to do that, because for the determination of
$Y_1 f, \dots, Y_p f$ we obtain the equations:
\[
\sum_{\pi=1}^p\,\xi_{j,\,n-p+\pi}\,Y_\pi f
=
X_jf
\ \ \ \ \ \ \ \ \ \ \ \ \
{\scriptstyle{(j\,=\,1\,\cdots\,p)}},
\]
which are solvable, and which provide for the $Y_\pi$ expressions
of the form:
\[
Y_\pi f
=
\sum_{j=1}^p\,\rho_{\pi j}(x_1,\dots,x_n)\,X_jf
\ \ \ \ \ \ \ \ \ \ \ \ \
{\scriptstyle{(\pi\,=\,1\,\cdots\,p)}},
\]
where the coefficients $\rho_{ \pi j}$ behave regularly
for the systems of values of:
\[
\Omega_1=0,
\,\dots,\,\,
\Omega_{n-m}=0.
\]
Hence if the system of equations $\Omega_k = 0$ admits the
infinitesimal transformations $X_1 f, \dots, X_p f$, it also admits
$Y_1 f, \dots, Y_p f$, and conversely.

Let $x_1^0, \dots, x_n^0$ be any system of values which satisfies the
equations $\Omega_1 = 0, \, \dots, \,
\Omega_{ n-m} = 0$ and for which the determinant 
$D$ is nonzero. Then the $\eta_{ \pi i}$ behave regularly in the
neighbourhood of $x_1^0, \dots, x_n^0$. Now, the equations $Y_1 f =
0, \, \dots, \, Y_p f = 0$ constitute a $p$-term complete system just
as the equations $X_1 f = 0, \, \dots, \, X_p f = 0$, hence according
to Theorem~12 p.~\pageref{Theorem-12-S-91}, they possess $n-p$ general
solutions $y_1, \dots, y_{ n-p}$ which behave regularly in the
neighbourhood of $x_1^0, \dots, x_n^0$ and which reduce to $x_1,
\dots, x_{ n-p}$ (respectively) for $x_{ n-p+1} = x_{ n-p + 1}^0, \,
\dots, \, x_n = x_n^0$.

Thus, if we still set $y_{ n-p+ 1} = x_{ n-p+1}, \, \dots, \, 
y_n = x_n$, we can introduce $y_1, \dots, y_n$ as new
variables in place of the $x$. 
At the same time, the infinitesimal transformations
$Y_\pi f$ receive the form:
\[
Y_1f
=
\frac{\partial f}{\partial y_{n-p+1}},\,\dots,\,\,
Y_pf
=
\frac{\partial f}{\partial y_n};
\]
but the system of equations $\Omega_k = 0$ is transferred to:
\[
\overline{\Omega}_1(y_1,\dots,y_n)=0,
\,\dots,\,\,
\overline{\Omega}_{n-m}(y_1,\dots,y_n)=0,
\]
and now these new equations must admit the infinitesimal
transformations $\frac{ \partial f}{ \partial y_{ n-p+1}}, \dots,
\frac{ \partial f}{ \partial y_n}$. From this, it follows that the
equations $\overline{ \Omega}_k = 0$ are solvable with respect to none
of the variables $y_{ n- p + 1}, \dots, y_n$, that they contain these
variables at most fictitiously and that they can be reshaped so as to
represent only relations between $y_1, \dots, y_{ n-p}$ alone.

If we now return to the original variables $x_1, \dots, x_n$, we then
see that the equations $\Omega_1 = 0, \,
\dots, \, \Omega_{ n-m} = 0$ are nothing but relations
between the solutions of the complete system $Y_1 f = 0, \, \dots, \,
Y_p f = 0$, or, what is the same, of the complete system $X_1 f = 0,
\, \dots, \, X_p f = 0$. Moreover, this result is independent of the
assumption that precisely the determinant $D$ should not vanish by
means of $\Omega_1 = 0, \, \dots, \, \Omega_{ n-m} = 0$; thus
it holds always true when not all $p \times p$ determinants of the
matrix~\thetag{ 5} vanish by means of 
$\Omega_1 = 0, \, \dots, \, \Omega_{ n-m} = 0$.

Consequently, we can state the following theorem:

\def\thetheorem{17}\begin{theorem}
\label{Theorem-17-S-123}
If $q$ infinitesimal transformations $X_1f, \dots, X_q f$ in the
variables $x_1, \dots, x_n$ provide exactly $p \leqslant q$
independent equations when equated to zero, say $X_1 f = 0, \, \dots,
\, X_p f = 0$, and if the latter equations form a $p$-term complete
system, then every system of equations:
\[
\Omega_1(x_1,\dots,x_n)=0,
\,\dots,\,\,
\Omega_{n-m}(x_1,\dots,x_n)=0
\]
which admits the $q$ infinitesimal transformations $X_1f, \dots, 
x_q f$ without cancelling the independence 
of the equations:
\[
X_1f=0,\,\dots,\,\,X_pf=0,
\]
is represented by relations between the solutions
of the $p$-term complete system $X_1 f = 0, \, \dots, 
\, X_p f = 0$.
\end{theorem}

Thanks to this theorem, the determination of all systems of equations
which belong to our first class is accomplished. In place of $X_1f =
0, \, \dots, \, X_p f = 0$, we only have to insert in the theorem one
after the other all the systems of $p$ independent equations amongst the
$q$ equations $X_1 f = 0, \, \dots, \, X_q f = 0$.

\sectionengellie{\S\,\,\,34.}

We now come to the second class of systems of equations which 
\label{S-123}
admit the infinitesimal transformations $X_1f, \dots, X_q f$, 
namely to the systems of equations which bring to zero all $p \times p$
determinants of the matrix~\thetag{ 4}.

\plainstatement{Here, a series of subcases must immediately 
be distinguished. Namely it is possible that 
aside from the $p \times p$ determinants of the matrix~\thetag{ 4},
the system of equations:
\[
\Omega_1=0,
\,\dots,\,\,\Omega_{n-m}=0
\]
also brings to zero yet all $(p-1) \times ( p-1)$ determinants, 
all $(p-2) \times (p-2)$ determinants, and so on. }

We hence see that to every sought system of equations is associated a
determinate number $h < p$ with the property that the concerned system
of equations brings to zero all $p \times p$, all $(p-1) \times (
p-1)$, \dots, all $(h+1) \dots ( h+1)$ determinants of the
matrix~\thetag{ 4}, \label{S-123-bis}
but not all $h \times h$ determinants. Thus, we
must go in details through all the various possible values $1, 2,
\dots, p-1$ of $h$ and for each one of these values, we must set up
the corresponding systems of equations that are admitted by $X_1 f,
\dots, X_q f$.

Let $h$ be any of the numbers $1, 2, \dots, p-1$. A system of
equations which brings to zero all $(h+1) \times (h+1)$ determinants
of the matrix~\thetag{ 4}, but not all the $h \times h$ ones, contains
in any case all equations which are obtained by equating to zero all
$( h+1) \times ( h+1)$ determinants $\Delta_1, \dots, \Delta_s$
of~\thetag{ 4}. Now, if the equations $\Delta_1 = 0, \, \dots, \,
\Delta_s = 0$ would absolutely not be satisfied by the systems of
values $x_1, \dots, x_n$ for which the $\xi_{ ki} (x)$ behave
regularly, or else, if the equations $\Delta_1 = 0, \, \dots, \,
\Delta_s = 0$ would also bring to zero all $h \times h$ determinants
of~\thetag{ 4}, then this would be a sign that there is absolutely no
system of the demanded nature to which is associated the chosen number
$h$. So we assume that none of these two cases occurs.

At first, it is to be examined whether the system of equations
$\Delta_1 = 0, \, \dots,\, \Delta_s = 0$ is reducible. When this is
the case, one should consider every irreducible\footnote{\,
This means smooth after stratification and not decomposable further in
an invariant way.
} 
\deutsch{irreducibel}
system of equations
which comes from $\Delta_1 = 0, \, \dots, \, \Delta_s = 0$, and which
does not bring to zero all $h \times h$ determinants of~\thetag{ 4}.

Let $W_1 = 0, \, \dots, \, W_l = 0$ be one of the found irreducible
systems of equations, and suppose that it is already brought to a form
such that not all $l \times l$ determinants of the matrix:
\[
\left\vert
\begin{array}{ccc}
\frac{\partial W_1}{\partial x_1} & \cdot\,\,\cdot &
\frac{\partial W_1}{\partial x_n}
\\ 
\cdot & \cdot\,\,\cdot & \cdot
\\
\frac{\partial W_l}{\partial x_1} & \cdot\,\,\cdot &
\frac{\partial W_l}{\partial x_n}
\end{array}
\right\vert
\]
vanish by means of $W_1 = 0, \, \dots, \, W_l = 0$. Thus we must
determine all systems of equations which are admitted by $X_1f, \dots,
X_qf$, which contain the equations $W_1 = 0, \, \dots, \, W_l = 0$ and
which at the same time do not bring to zero all $h \times h$
determinants of ~\thetag{ 4}. When we execute this for each
individual irreducible system obtained from $\Delta_1 = 0, \, \dots,
\, \Delta_s = 0$, we find all systems of 
equations which admit $X_1 f, \dots, X_q f$ and to which is associated
the number $h$.

In general, \label{S-124-sq}
the system of equations $W_1 = 0, \, \dots, W_l = 0$ will
in fact not admit the infinitesimal transformations $X_1f, \dots, X_q
f$. To a system of equations which contains $W_1 = 0, \, \dots, \,
W_l = 0$ and which in addition admits $X_1 f , \dots, X_q f$ there
must belong in any case yet the equations:
\[
X_kW_\lambda=0,
\ \ \ \ \ \ \
X_jX_kW_\lambda=0
\ \ \ \ \ \ \ \ \ \ \ \ \
{\scriptstyle{(k,\,j\,=\,1\,\cdots\,q;\,\,\lambda\,=\,1\,\cdots\,l)}},
\]
and so on. We form these equations and we examine whether they become
contradictory with each other, or with $W_1 = 0, \, \dots, \, W_l =
0$, and whether they possibly bring to zero all $h \times h$
determinants of~\thetag{ 4}. If one of these two cases occurs, then
there is no system of equations of the demanded constitution; if none
occurs, then the independent equations amongst the equations:
\[
W_\lambda=0,
\ \ \ \ \
X_kW_\lambda=0,
\ \ \ \ \
X_jX_kW_\lambda=0,\dots
\ \ \ 
{\scriptstyle{(\lambda\,=\,1\,\cdots\,l;\,\,j,\,k\,=\,1\,\cdots\,q)}}
\]
represent a system of equations which is admitted by $X_1 f, \dots,
X_q f$. Into this system of equations which can naturally be
reducible, we stick all systems of equations which admit $X_1f, \dots,
X_q f$, which contain the equations $W_1 = 0, \, \dots, \, W_l = 0$,
but which embrace no smaller system of equations of the same nature.

The corresponding systems of equations are obviously the
\emphasis{smallest} systems of equations which are admitted by $X_1f,
\dots, X_qf$ and which bring to zero all $(h+1) \times ( h+1)$
determinants of~\thetag{ 4}, though not all the $h \times h$ ones. 

Let now:
\def\theequation{6}\begin{equation}
W_1=0,
\,\dots,\,\,
W_{n-m}=0
\ \ \ \ \ \ \ \ \ \ \ \ \
{\scriptstyle{(n\,-\,m\,\geqslant\,l)}}
\end{equation}
be one of the found irreducible systems of equations; then the
question is to add, in the most general way, new equations to this
system so that one obtains a system of equations which is admitted by
$X_1f, \dots, X_q f$ but which does not bring to zero all $h \times h$
determinants of~\thetag{ 4}.

Since not all $h \times h$ determinants of~\thetag{ 4} vanish
by means of:
\def\theequation{6}\begin{equation}
W_1=0,
\,\dots,\,\,
W_{n-m}=0,
\end{equation}
we can assume that for instance the determinant:
\[
\Delta
=
\left\vert
\begin{array}{ccc}
\xi_{1,\,n-h+1} & \cdot\,\,\cdot & \xi_{1n}
\\
\cdot & \cdot\,\,\cdot & \cdot
\\
\xi_{h,\,n-h+1} & \cdot\,\,\cdot & \xi_{hn}
\end{array}
\right\vert
\]
belongs to the nonvanishing ones and we can set ourselves the problem
of determining all systems of equations which are admitted by $X_1f,
\dots, X_q f$ and which at the same time do not make $\Delta$ equal to
zero. When we carry out this problem for every individual $h \times h$
determinant of~\thetag{ 4} which does not vanish already by means
of~\thetag{ 6}, we then evidently obtain all the sought systems of
equations which embrace $W_1 = 0, \, \dots, \, W_{ n-m} = 0$.

The independence of the equations $X_1 f = 0, \, \dots, \, X_h f = 0$
is, under the assumptions made, not cancelled by $W_1 = 0, \, \dots,
\, W_{ n-m} = 0$, and in fact, for the systems of values of~\thetag{
6}, there are certain relations of the form:
\[
X_{h+j}f
=
\sum_{k=1}^h\,\psi_{jk}(x_1,\dots,x_n)\,X_kf
\ \ \ \ \ \ \ \ \ \ \ \ \
{\scriptstyle{(j\,=\,1\,\cdots\,q\,-\,h)}}
\]
where the $\psi_{ jk}$ are to be determined out from the
equations:
\[
\xi_{h+j,\,\nu}
=
\sum_{k=1}^h\,\psi_{jk}\,\xi_{k\nu}
\ \ \ \ \ \ \ \ \ \ \ \ \
{\scriptstyle{(j\,=\,1\,\cdots\,q\,-\,h;\,\,\nu\,=\,1\,\cdots\,n)}}.
\]
Since all $(h+1) \times ( h+1)$ determinants of~\thetag{ 4} vanish by
means of~\thetag{ 6} while $\Delta$ does not, the functions $\psi_{
jk}$ are perfectly determined and they behave regularly for the
systems of values of~\thetag{ 6}. Without loss of generality, we are
hence allowed to leave out the infinitesimal transformations $X_{ h+1}
f, \dots, X_q f$; because every system of equations which contains the
equations~\thetag{ 6}, which at the same time does not bring to zero
the determinant $\Delta$, and lastly, which is admitted by $X_1f, \dots,
X_h f$, is also automatically admitted by $X_{ h+1} f, \dots, X_q f$.

\label{S-126}
One can easily see that no relation between $x_{ n-h+1}, 
\dots, x_n$ alone can be derived from the equations~\thetag{ 6}. 
Indeed, if one obtained such a relation, say:
\[
x_n
-
\omega(x_{n-h+1},\dots,x_{n-1})
=
0,
\]
then the $h$ expressions:
\[
X_k(x_n-\omega)
=
\xi_{kn}
-
\sum_{j=1}^{h-1}\,
\frac{\partial\omega}{\partial x_{n-h+j}}\,
\xi_{k,\,n-h+j}
\ \ \ \ \ \ \ \ \ \ \ \ \
{\scriptstyle{(k\,=\,1\,\cdots\,h)}}
\]
would vanish by means of~\thetag{ 6}. From this, we draw
the conclusion that the number $n - m$ is in any case 
not larger than $n - h$ and that the equations~\thetag{ 6}
can be resolved with respect to $n-m$ of the variables
$x_1, \dots, x_{ n-h}$, say with respect to $x_1, \dots, 
x_{ n-m}$:
\def\theequation{7}\begin{equation}
x_k
=
\varphi_k(x_{n-m+1},\dots,x_n)
\ \ \ \ \ \ \ \ \ \ \ \ \
{\scriptstyle{(k\,=\,1\,\cdots\,n\,-\,m)}} 
\ \ \ \ \ \ \ \
{\scriptstyle{(n\,-\,m\,\leqslant\,n\,-\,h)}}.
\end{equation}
We can therefore replace the system~\thetag{ 6} by these equations.

Every system of equations which contains the equations~\thetag{ 6}
or equivalently the equations~\thetag{ 7} can obviously be brought to 
a form such that, aside from the equations~\thetag{ 7}, it yet
contains a certain number of relations:
\[
V_j(x_{n-m+1},\dots,x_n)
=
0
\ \ \ \ \ \ \ \ \ \ \ \ \
{\scriptstyle{(j\,=\,1,\,2,\,\cdots)}} 
\]
between $x_{ n-m+1}, \dots, x_n$ alone. Now, if the system of
equations in question is supposed to
admit the infinitesimal transformations
$X_1 f, \dots, X_h f$, then all the expressions $X_k V_j$ must vanish
by means of~\thetag{ 7} and of $V_1 = 0, \, \dots$; or, if we again
denote by the sign $\big[ \,\, \big]$ the substitution $x_1 =
\varphi_1, \, \dots, \, x_{ n-m} = \varphi_{ n-m}$: 
the expressions:
\[
\big[X_kV_j\big]
=
\sum_{\mu=1}^m\,
\big[\xi_{k,\,n-m+\mu}\big]\,
\frac{\partial V_j}{\partial x_{n-m+\mu}}
\ \ \ \ \ \ \ \ \ \ \ \ \
{\scriptstyle{(k\,=\,1\,\cdots\,h;\,\,j\,=\,1,\,2,\,\cdots)}}
\]
vanish by means of $V_1 = 0, \, \dots$. But we can also express this
as follows: the system of equations $V_1 = 0, \, \dots$ in the
variables $x_{ n-m+1}, \dots, x_n$ must admit the $h$ reduced
infinitesimal transformations:
\[
\overline{X}_kf
=
\sum_{\mu=1}^m\,
\big[\xi_{k,\,n-m+\mu}\big]\,
\frac{\partial f}{\partial x_{n-m+\mu}}
\ \ \ \ \ \ \ \ \ \ \ \ \
{\scriptstyle{(k\,=\,1\,\cdots\,h)}}
\]
in these variables. In addition, the system of equations 
$V_1 = 0, \dots$ should naturally not bring to zero
the determinant $\big[ \Delta \big]$. 

Conversely, every system of equations $V_1 = 0, \, \dots$ which
possesses the property discussed just now provides, together
with~\thetag{ 7}, a system of equations which does not make $\Delta$
be zero and which in addition admits $X_1 f, \dots, X_h f$ and also
$X_{ h+1} f, \dots, X_q f$ as well.

As a result, our initial problem is reduced to the simpler problem of
determining all systems of equations in $m < n$ variables which admit
the $h \leqslant m$ infinitesimal transformations $\overline{ X}_1 f,
\dots, \overline{ X}_h f$ and for which the determinant:
\[
\sum\,\pm\,
\big[\xi_{1,\,n-h+1}\big]\cdots\big[\xi_{h,\,n-h+h}\big]
\]
is not made zero.

We summarize what was done up to now in the

\def\thetheorem{18}\begin{theorem}
If $q$ infinitesimal transformations:
\[
X_kf
=
\sum_{i=1}^n\,\xi_{ki}(x_1,\dots,x_n)\,
\frac{\partial f}{\partial x_i}
\ \ \ \ \ \ \ \ \ \ \ \ \
{\scriptstyle{(k\,=\,1\,\cdots\,q)}}
\]
are constituted in such a way that
all $(p+1) \times ( p+1)$ determinants, but not all
$p \times p$ determinants of the matrix:
\[
\left\vert
\begin{array}{ccc}
\xi_{11} & \cdot\,\,\cdot & \xi_{1n}
\\
\cdot & \cdot\,\,\cdot & \cdot
\\
\xi_{q1} & \cdot\,\,\cdot & \xi_{qn}
\end{array}
\right\vert
\]
\label{S-127}
vanish identically, and that any $p$ independent equations amongst the
equations $X_1 f = 0, \, \dots, \, X_q f = 0$ form a $p$-term complete
system, then one finds in the following manner all systems of
equations in $x_1, \dots, x_n$ which admit $X_1f, \dots, X_q f$ and
which at the same time bring to zero all $( h+1) \times ( h+1)$
determinants of the above matrix, but not all $h \times h$ ones: by
forming determinants
\deutsch{durch Determinantenbildung},
one seeks at first the smallest system of equations for which all
$(h+1) \times ( h+1)$ determinants of the matrix vanish, but not all
$h \times h$ ones. If there exists a system of this sort, and if $W_1
= 0, \, \dots, \, W_l = 0$ is one such, then one
forms the equations $X_k W_i = 0$, $X_j X_k W_i = 0$, \dots, and one
determines in this way the possibly existing smallest system of
equations which embraces $W_1 = 0, \, \dots, \, W_l = 0$, which is
admitted by $X_1 f, \dots, X_q f$ and which does not make equal to
zero all $h \times h$ determinants of the matrix; if $W_1 = 0, \,
\dots, \, W_{ n-m} = 0$ $(n-m \geqslant l)$ is such a system of
equations, which does not bring to zero for instance the determinant:
\[
\Delta
=
\sum\,\pm\,
\xi_{1,\,n-h+1}\cdots\xi_{h,\,n-h+h},
\]
then $h$ is $\leqslant m$ and the equations $W_1 = 0, \, \dots, \, W_{
n-m} = 0$ can be resolved with respect to $n - m$ of the variables
$x_1, \dots, x_{ n-h}$, say as follows:
\[
x_k
=
\varphi_k(x_{n-m+1},\dots,x_n)
\ \ \ \ \ \ \ \ \ \ \ \ \
{\scriptstyle{(k\,=\,1\,\cdots\,n\,-\,m)}}.
\]
Lastly, one determines all systems of equations in the $m$
variables $x_{ n- m + 1}, \dots, x_n$ which admit the
$h$ reduced infinitesimal transformations:
\[
\aligned
\overline{X}_kf
&
=
\sum_{\mu=1}^m\,\xi_{k,\,n-m+\mu}(\varphi_1,\dots,\varphi_{n-m},\,\,
x_{n-m+1},\dots,x_n)\,
\frac{\partial f}{\partial x_{n-m+\mu}}
\\
&
=
\sum_{\mu=1}^m\,
\big[\xi_{k,\,n-m+\mu}\big]\,
\frac{\partial f}{\partial x_{n-m+\mu}}
\ \ \ \ \ \ \ \ \ \ \ \ \
{\scriptstyle{(k\,=\,1\,\cdots\,h)}}
\endaligned
\]
and which do not bring to zero the determinant:
\[
\big[\Delta\big]
=
\sum\pm\,\big[\xi_{1,\,n-h+1}\big]
\cdots
\big[\xi_{h,\,n_h+h}\big].
\]
Each one of these systems of equations represents, after adding the
equations $x_1 = \varphi_1, \, \dots, \, x_{ n-m} = \varphi_{ n-m}$, a
system of equations of the demanded constitution. By carrying out the
indicated developments in all possible cases, one obtains all systems
of equations of the demanded constitution.
\end{theorem}

It frequently happens that one already knows a system of equations
$U_1 = 0, \, U_2 = 0, \, \dots, \, U_{ n-s} = 0$ which admits the
infinitesimal transformations $X_1f, \dots, X_q f$ and which brings to
zero all $(h+1) \times ( h+1)$ determinants of the matrix~\thetag{ 4},
and which nevertheless does not bring to zero the $h \times h$
determinant $\Delta$. Then one can ask for all systems of equations
which comprise the equations $U_1 = 0, \, \dots, \, U_{ n-s} = 0$ and
which likewise do not bring $\Delta$ to zero.

The determination of all these systems of equations can be carried
out exactly as in the special case above, where the
smallest system of equations $W_1 = 0, \, \dots, \, W_{ n-m} = 0$
of the relevant nature was known, 
and where one was looking for all systems of
equations which were admitted by $X_1 f, \dots, X_q f$, which
comprised the system $W_1 = 0, \, \dots, \, W_{ n-m} = 0$
and which did not bring $\Delta$ to zero.

Indeed, in exactly the same way as above, one shows at first
that $h \leqslant s$ and that the equations $U_1 = 0, \, 
\dots, \, U_{ n-s} = 0$ can be resolved with respect to 
$n - s$ of the variables $x_1, \dots, x_{ n-h}$, say 
as follows:
\[
x_k
=
\psi_k(x_{n-s+1},\dots,x_n)
\ \ \ \ \ \ \ \ \ \ \ \ \
{\scriptstyle{(k\,=\,1\,\cdots\,n\,-\,s)}}.
\]
Now, in order to find the sought systems of equations, 
one forms the $h$ reduced infinitesimal transformations:
\[
\aligned
\widehat{X}_kf
&
=
\sum_{\sigma=1}^s\,
\xi_{k,\,n-s+\sigma}
\big(\psi_1,\dots,\psi_{n-s},\,x_{n-s+1},\dots,x_n)\,
\frac{\partial f}{\partial x_{n-s+\sigma}}
\\
&
=
\sum_{\sigma=1}^s\,
\widehat{\xi}_{k,\,n-s+\sigma}\,
\frac{\partial f}{\partial x_{n-s+\sigma}}
\ \ \ \ \ \ \ \ \ \ \ \ \
{\scriptstyle{(k\,=\,1\,\cdots\,h)}}
\endaligned
\]
and next, one determines all systems of equations in the $s$ variables
$x_{ n-s+1}, \dots, x_n$ which admit $\widehat{ X}_1 f, \dots,
\widehat{ X}_h f$ and which do not bring to zero the determinant:
\[
\widehat{\Delta}
=
\sum\,\pm\,\widehat{\xi}_{1,\,n-h+1}\cdots
\widehat{\xi}_{h,\,n-h+h}.
\]
By adding these equations one after the other to the equations $x_1 -
\psi_1 = 0, \, \dots, \, x_{ n-s} - \psi_{ n-s} = 0$, one obtains all
systems of equations of the demanded constitution.

It is not necessary to explain more precisely what has been just said.

\sectionengellie{\S\,\,\,35.}

The problem which has been set up at the beginning of
\S\,\,32 is now basically settled. 
Indeed, thanks to the latter theorem, this problem is reduced to the
determination of all systems of equations in the $m$ variables $x_{
n-m+1}, \dots, x_n$ which admit the $h$ infinitesimal transformations
$\overline{ X}_1 f, \dots, \overline{ X}_h f$. But this is a problem
of the same type as the original one, which is simplified only
inasmuch as the number $m$ of the variables is smaller than $n$.

Now, about the reduced problem, the same considerations as those about
the initial problem can be made use of. That is to say: if the
equations $\overline{ X}_1 f = 0, \, \dots, \,
\overline{ X}_h f = 0$ do not actually form an $h$-term
complete system, then one has to set up for $k,\, j = 1, \dots, h$ the
infinitesimal transformations $\leftbracket \overline{ X}_k, \, \overline{
X}_j \rightbracket$ and to ask whether the independent equations amongst the
equations $\overline{ X}_k f = 0$, $\leftbracket \overline{ X}_k, \,
\overline{ X}_j \rightbracket = 0$ for a complete system, and so on, in brief,
one proceeds just as for the original problem. The only difference in
comparison to the former situation is that from the beginning one
looks for systems of equations which do not cancel the independence of
the equations:
\[
\overline{X}_1f=0,
\,\dots,\,\,
\overline{X}_hf=0,
\]
and especially, which do not bring the determinant $\big[ \Delta
\big]$ to zero.

Exactly as above for the original problem, the reduced problem can be
dealt with in parts, and be partially reduced to a problem in yet less
variables, and so forth. Briefly, one sees that the complete
resolution of the original problem can be attained after a finite
number of steps.

We therefore need not to address this issue further,
but we want to consider somehow into more details a special
particularly important case. 

Let as above the system of equations:
\def\theequation{7'}\begin{equation}
x_1-\varphi_1(x_{n-m+1},\dots,x_n)=0,
\,\dots,\,\,
x_{n-m}-\varphi_{n-m}(x_1,\dots,x_n)=0
\end{equation}
be constituted in such a way that it admits the infinitesimal
transformations $X_1f, \dots, X_q f$ and assume that it brings to zero
all $(h+1) \times ( h+1)$ determinants of the matrix~\thetag{ 4},
while it does not bring to zero the determinant $\Delta$. In
comparison, by the system of equations~\thetag{ 7'}, we want to
understand \emphasis{not just the smallest one}, but a
\emphasis{completely arbitrary one} of the demanded constitution.

Amongst the equations $X_1 f = 0, \, \dots, \, X_q f = 0$, assume that
any $p$ independent amongst them, say for instance $X_1 f = 0, \, \dots,
\, X_p f = 0$, form a $p$-complete system; thus in any case, there are
relations of the form:
\[
\leftbracket X_k,\,X_j\rightbracket
=
\sum_{\sigma=1}^q\,
\omega_{kj\sigma}(x_1,\dots,x_n)\,X_\sigma f
\ \ \ \ \ \ \ \ \ \ \ \ \
{\scriptstyle{(k,\,j\,=\,1\,\cdots\,q)}}.
\]
For the systems of values of~\thetag{ 7'} 
there remain only $h$ equations:
\[
X_1f=0,
\,\dots,\,\,
X_hf=0
\]
that are independent of each other, whereas $X_{ h+1}f, \dots, 
X_q f$ can be represented under the form:
\[
X_{h+j}f
=
\sum_{\tau=1}^h\,
\psi_{j\tau}(x_1,\dots,x_n)\,X_\tau f
\ \ \ \ \ \ \ \ \ \ \ \ \
{\scriptstyle{(j\,=\,1\,\cdots\,q\,-\,h)}},
\]
where the $\psi_{ jk}$ behave regularly for the
concerned system of values. 

Now, \label{S-130}
we want to make the specific assumption that also all the
coefficients $\omega_{ kjs}$ behave regularly 
for the
system of values of~\thetag{ 7}; then evidently for the
concerned system of values, equations of the form:
\[
\leftbracket
X_k,\,X_j
\rightbracket
=
\sum_{\sigma=1}^h\,
\bigg\{
\omega_{kj\sigma}
+
\sum_{\tau=1}^{q-h}\,
\omega_{kjh+\tau}\,\psi_{\tau\sigma}
\bigg\}\,X_\sigma f
\]
hold true. In this special case one can indicate in one stroke
all systems of equations which embrace the equations~\thetag{ 7'}, 
which admit $X_1f, \dots, X_q f$ and which do not bring
$\Delta$ to zero. This should now be shown.

For all systems of values $x_1, \dots, x_n$ 
which satisfy the equations~\thetag{ 7'}, there are relations of the 
form:
\[
\leftbracket X_k,\,X_j\rightbracket
=
\sum_{\sigma=1}^h\,w_{kj\sigma}(x_1,\dots,x_n)\,X_\sigma f
\ \ \ \ \ \ \ \ \ \ \ \ \
{\scriptstyle{(k,\,j\,=\,1\,\cdots\,h)}}.
\]
Here, the $w_{ kj \sigma}$ behave regularly and likewise
the functions $[ w_{ kjs} ]$, if by the sign
$[ \, \, ]$ we denote as before the substitution:
\[
x_1=\varphi_1,
\,\dots,\,\,x_{n-m}=\varphi_{n-m}.
\]

We decompose the above relations in the following ones:
\[
\aligned
X_k\xi_{j\nu}
&
-
X_j\xi_{k\nu}
=
\sum_{s=1}^h\,w_{kjs}(x_1,\dots,x_n)\,\xi_{s\nu}
\\
&
\ \ \ \ \ \ \ \ \ \ \ \ \
{\scriptstyle{(k,\,j\,=\,1\,\cdots\,h;\,\,\nu\,=\,1\,\cdots\,n)}},
\endaligned
\]
and they naturally hold identically after the substitution $[ \, \,
]$, so that we have:
\def\theequation{8}\begin{equation}
\big[X_k\xi_{j\nu}\big]
-
\big[X_j\xi_{k\nu}\big]
\equiv
\sum_{s=1}^h\,[w_{kjs}]\,[\xi_{s\nu}].
\end{equation}
But now the system of equations $x_k - \varphi_k = 0$ admits the
infinitesimal transformations $X_1 f, \dots, X_h f$, so the
relation~\thetag{ 3} derived on p.~\pageref{relation-3-S-110}:
\[
\big[X_k\Phi\big]
\equiv
\big[X_k[\Phi]\big]
\]
holds true, in which $\Phi ( x_1, \dots, x_n)$ is a completely
arbitrary function of its arguments. We can write this relation
somewhat differently, if we remember the infinitesimal
transformations: \label{S-131}
\[
\overline{X}_kf
=
\sum_{\mu=1}^m\,
[\xi_{k,\,n-m+\mu}]\,
\frac{\partial f}{\partial x_{n-m+\mu}}
\ \ \ \ \ \ \ \ \ \ \ \ \
{\scriptstyle{(k\,=\,1\,\cdots\,h)}};
\]
indeed, one evidently has:
\[
\big[X_k[\Phi]\big]
\equiv
\overline{X}_k[\Phi],
\]
whence:
\[
\big[X_k\Phi\big]
\equiv
\overline{X}_k[\Phi].
\]
From this, it follows that the identities~\thetag{ 8}
can be replaced by the following ones:
\[
\overline{X}_k[\xi_{j\nu}]
-
\overline{X}_j[\xi_{k\nu}]
\equiv
\sum_{s=1}^h\,[w_{kjs}]\,[\xi_{s\nu}].
\]
In other words, there are identities of the form:
\[
\leftbracket \overline{X}_k,\,\overline{X}_j\rightbracket
=
\overline{X}_k\overline{X}_jf
-
\overline{X}_j\overline{X}_kf
\equiv
\sum_{s=1}^h\,[w_{kjs}]\,\overline{X}_sf,
\]
that is to say, the equations $\overline{ X}_1 f = 0, \,
\dots, \, \overline{ X}_h f = 0$ form an $h$-term complete 
system in the $m$ independent variables $x_{ n-m + 1}, \dots, x_n$.

Now, we remember the observations that we have linked up with
Theorem~18. They showed that our problem stated above comes down to
determining all systems of equations in $x_{ n-m+1},
\dots, x_n$ which admit $\overline{ X}_1 f, \dots, 
\overline{ X}_h f$ and do not cancel the independence
of the equations $\overline{ X}_1 f = 0, \,
\dots, \, \overline{ X}_h f = 0$. But since in our case 
the equations $\overline{ X}_1 f = 0, \,
\dots, \, \overline{ X}_h f = 0$ form an $h$-term complete
system, we can at once apply the Theorem~17. Thanks to it, we see
that the sought systems of equations in $x_{ n-m+1 }, \dots, x_n$ are
represented by relations between the solutions of the complete system
$\overline{ X}_1 f = 0, \, \dots, \, \overline{ X}_h f = 0$. 

Consequently, if we add arbitrary relations between the solutions of
this complete system to the equations $x_k = \varphi_k$, we obtain the
general form of a system of equations which admits the infinitesimal
transformations $X_1 f, \dots, X_q f$, which comprises the equations
$x_k = \varphi_k$ and which does not bring $\Delta$ to zero.

It appears superfluous to formulate as a proposition the result
obtained here in its full generality. By contrast, it is useful for
the sequel to expressly state the following theorem, which corresponds
to the special assumption $q = p = h$. 

\def\thetheorem{19}\begin{theorem}
\label{Theorem-19-S-132}
If a system of $n - m$ independent equations in the
variables $x_1, \dots, x_n$ admits the $h$ infinitesimal 
transformations:
\[
X_kf
=
\sum_{i=1}^n\,\xi_{ki}(x_1,\dots,x_n)\,
\frac{\partial f}{\partial x_i}
\ \ \ \ \ \ \ \ \ \ \ \ \
{\scriptstyle{(k\,=\,1\,\cdots\,h)}}
\]
and if at the same time the determinant:
\[
\Delta
=
\sum\,\pm\,
\xi_{1,\,n-h+1}\cdots\xi_{h,\,n-h+h}
\]
neither vanishes identically, nor vanishes by means
of the system of equations, then $h$ is $\leqslant m$ and the
system of equations can be resolved with respect to 
$n - m$ of the variables $x_1, \dots, x_{ n-h}$, 
say as follows:
\[
x_1
=
\varphi_1(x_{n-m+1},\dots,x_n),
\,\dots,\,\,
x_{n-m}
=
\varphi_{n-m}(x_{n-m+1},\dots,x_n).
\]
Now if, for the systems of values
$x_1, \dots, x_n$ of these equations, all the expressions $\leftbracket 
X_k, \, X_j \rightbracket$ can be represented in the form:
\[
\leftbracket X_k,\,X_j\rightbracket
=
\sum_{s=1}^h\,w_{kjs}(x_1,\dots,x_n)\,X_sf
\ \ \ \ \ \ \ \ \ \ \ \ \
{\scriptstyle{(k,\,j\,=\,1\,\cdots\,h)}},
\]
where the $w_{ kjs}$ behave regularly for the 
concerned systems of values, then one finds as follows
all systems of equations which comprise the equations:
\[
x_1-\varphi_1=0,
\,\dots,\,\,
x_{n-m}-\varphi_{n-m}=0,
\]
which admit the infinitesimal transformations $X_1f, \dots, 
X_hf$ and which do not bring to zero the determinant $\Delta$:
one sets up the $h$ reduced infinitesimal 
transformations:
\[
\aligned
\overline{X}_kf
=
\sum_{\mu=1}^m\,
\xi_{k,\,n-m+\mu}
&
\big(\varphi_1,\dots,\varphi_{n-m},\,x_{n-m+1},\dots,x_n\big)\,
\frac{\partial f}{\partial x_{n-m+\mu}}
\\
&
\ \ \ \ \ \ \ \ \ \ \ \ \
{\scriptstyle{(k\,=\,1\,\cdots\,h)}}; 
\endaligned
\]
then the $h$ mutually independent equations $\overline{ X}_1 f = 0, \,
\dots, \, \overline{ X}_h f = 0$ form an $h$-term
complete system in the independent variables $x_{ n - m + 1}, \dots,
x_n$; if:
\[
u_1(x_{n-m+1},\dots,x_n),\dots,\,
u_{m-h}(x_{n-m+1},\dots,x_n)
\] 
are independent solutions of this complete system, then:
\[
\aligned
x_1-\varphi_1=0,\,\dots,\,\,
&
x_{n-m}-\varphi_{n-m}=0,\ \ \ \ \
\Phi_i(u_1,\dots,u_{m-h})=0
\\
&
\ \ \ \ \ \ \ \ \ \ \ 
{\scriptstyle{(i\,=\,1,\,2,\,\cdots)}}
\endaligned
\]
is the general form of the sought system of equations; 
by the $\Phi_i$ here are to be understood arbitrary function 
of their arguments.
\end{theorem}
\label{S-133}

\sectionengellie{\S\,\,\,36.}

The analytical developments of the present chapter receive a certain
clarity and especially a better transparency when one applies the
ideas and the forming of concepts of the theory of manifolds
\deutsch{die Vorstellungen und Begriffsbildungen der
Mannigfaltigkeitslehre}. We now want to do that. The sequel stands in
comparison to the \S\S\,\,30 to~35 in exactly the same relationship as
Chap.~6 stands in comparison to Chap.~5.

Every system of equations in the variables $x_1, \dots, x_n$
represents a manifold of the $n$-times extended manifold $R_n$. If
the system of equations admits the one-term group $Xf$, then according
to the previously introduced way of expressing, the corresponding
manifold also admits $Xf$; so if a point $x_1, \dots, x_n$ belongs to
the manifold, then all points into which $x_1, \dots, x_n$ is
transferred by all the transformations of the one-term group $Xf$ also
lie in the same manifold.

If now $x_1^0, \dots, x_n^0$ is any point of the space, two cases can
occur: either $\xi_1, \dots, \xi_n$ do not all vanish for $x_i =
x_i^0$, or the quantities $\xi_1 ( x^0),
\dots, \xi_n ( x^0)$ are all equal to zero. 
In the first case, the point $x_1^0, \dots, x_n^0$ takes infinitely
many positions by the $\infty^1$ transformations of the one-term group
$Xf$, and as we know, the totality of these positions is invariant by
the one-term group $Xf$ and it constitutes an integral curve of the
infinitesimal transformation $Xf$. In the second case, $x_1^0, \dots,
x_n^0$ keeps its position \label{S-134}
through all transformations of the one-term
group $Xf$; the integral curve passing by $x_1^0, \dots, x_n^0$
shrinks to the point itself.

So if a manifold admits the one-term group $Xf$ and if it consists in
general of points for which $\xi_1, \dots, \xi_n$ do not all vanish,
then it is constituted of integral curves of the infinitesimal
transformation $Xf$. By contrast, if the manifold in question
consists only of points for which $\xi_1, \dots, \xi_n$ all vanish,
then each point of the manifold keeps its position through all
transformations of the one-term group $Xf$. Evidently, every smaller
manifold contained in this manifold then admits also the one-term
group $Xf$.

With these words, the conceptual content
\deutsch{der begriffliche Inhalt} of the Theorem~15
is clearly stated, and in fact fundamentally, the preceding
considerations can be virtually regarded as a providing a new
demonstration of the Theorem~15.


We still have to explain what does it means conceptually for a system
of equations $\Omega_1 = 0, \, \dots, \,
\Omega_{ n - m} =0$ to admit the infinitesimal 
transformation $Xf$.

To this end, we remember that the infinitesimal transformation $Xf$
attaches to each point $x_1, \dots, x_n$ at which not all $\xi$ vanish
the completely determined direction of progress:
\[
\delta x_1\,:\,\cdots\,:\,\delta x_n
=
\xi_1\,:\,\cdots\,:\,\xi_n,
\]
while it attaches no direction of progress to a point for which
$\xi_1 = \cdots = \xi_n = 0$. We can therefore say:

\plainstatement{The system of equations 
$\Omega_1 = 0, \, \dots, \, \Omega_{ n- m} = 0$ admits the
infinitesimal transformation $Xf$ when the latter attaches to each
point of the manifold $\Omega_1 = 0, \, \dots, \, \Omega_{ n - m} = 0$
either absolutely no direction of progress, or a direction of progress
which satisfies the $n - m$ equations:
\[
\frac{\partial\Omega_k}{\partial x_1}\,\delta x_1
+\cdots+
\frac{\partial\Omega_k}{\partial x_n}\,\delta x_n
\ \ \ \ \ \ \ \ \ \ \ \ \
{\scriptstyle{(k\,=\,1\,\cdots\,n\,-\,m)}},
\]
hence which comes into contact with the manifold.}

This definition is visibly independent both of the choice of the
variables and of the form of the system of equations:
\[
\Omega_1=0,
\,\dots,\,\,\Omega_{n-m}=0.
\]

The Theorem~14 can now receive the following visual version:

\plainstatement{If the infinitesimal transformation $Xf$ 
attaches to each point of a manifold either absolutely no direction of
progress, or a direction of progress which comes into contact with the
manifold, then the manifold admits all transformations of the one-term
group $Xf$. }

In conclusion, if we yet introduce the language: ``the manifold
$\Omega_1 = 0, \, \dots, \, \Omega_{ n-m} = 0$ admits the infinitesimal 
transformation $Xf$, we can express as follows the Theorem~14.

\plainstatement{A manifold admits all transformations of the
one-term group $Xf$ if and only if it admits the infinitesimal
transformation $Xf$. }

In the \S\S\,\,32 and~34, we gave a classification
of all systems of equations which admit the $q$ infinitesimal
transformations $X_1f, \dots, X_q f$. There, we
took as a starting point the behaviour of the 
determinants of the matrix~\thetag{ 4}:
\[
\left\vert
\begin{array}{ccc}
\xi_{11} & \cdot\,\,\cdot & \xi_{1n}
\\
\cdot & \cdot\,\,\cdot & \cdot
\\
\xi_{q1} & \cdot\,\,\cdot & \xi_{qn}
\end{array}
\right\vert.
\]
We assumed that all $(p+1) \times (p+1)$ determinants of this matrix
vanished identically, whereas the $p \times p$ determinants did not.
We then reckoned as belonging to the same class all systems of
equations by means of which all $(h+1) \times ( h+1)$ determinants of
the matrix vanish, but not all $h \times h$ determinants; in the
process, the integer $h$ could have the $p+1$ different values: $p,
p-1, \dots, 2, 1, 0$.

Already in Chap.~6 we observed that, under the assumptions made just
now, the $q$ infinitesimal transformations:
\[
X_1f,\dots,X_qf
\]
attach exactly $p$ independent directions of progress to any point
$x_1, \dots, x_n$ in general position. 
It stands to reason in the
corresponding way that $X_1f, \dots, X_q f$ attach exactly $h$
independent directions of progress to a point $x_1, \dots, x_n$, when
the $h \times h$ determinants of the above matrix do not all vanish at
the point in question, 
\label{S-135}
while by contrast all $(h+1) \times ( h+1)$
determinants take the zero value. Now, since every system of
equations which admits the infinitesimal transformations $X_1f, \dots,
X_q f$ represents a manifold having the same property, our previous
classification of the systems of equations immediately provides a
classification of the manifolds. Indeed, amongst the manifolds which
admit $X_1f, \dots, X_q f$ we always reckon as belonging to one class
those systems, to the points of which are associated the same number
$h \leqslant p$ of independent directions of progress by the
infinitesimal transformations $X_1 f, \dots, X_q f$.

If a manifold admits the infinitesimal transformations $X_1f, \dots,
X_q f$, then at each one of its points, it comes into contact with the
directions of progress that are attached to the point by $X_1f, \dots,
X_qf$. Now, if $X_1 f, \dots, X_q f$ determine precisely $h$
independent directions at each point of the manifold, then the
manifold must obviously have at least $h$ dimensions. Hence we can
state the proposition:

\def\theproposition{6}\begin{proposition}
If $q$ infinitesimal transformations $X_1f, \dots, X_q f$ attach
precisely $h$ independent directions of progress to a special point
$x_1^0, \dots, x_n^0$, then there is in any case no manifold of
smaller dimension than $h$ which contains the point $x_1^0, \dots,
x_n^0$ and which admits the $q$ infinitesimal transformations $X_1f,
\dots, X_q f$.
\end{proposition}

Basically, this proposition is only another formulation of a former
result. In \S\,\,34 indeed, we considered the systems of equations
which admit $X_1f, \dots, X_qf$ and which at the same time leave only
$h$ mutually independent equations amongst the $q$ equations $X_1 f =
0, \, \dots, \, X_q f = 0$. On the occasion, we saw that such a
system of equations consists of at most $n-h$ independent equations,
so that it represents a manifold having at least $h$ 
\label{S-136}
dimensions.

\linestop


\chapter{Complete Systems Which Admit 
\\
All Transformations
\\
of a One-term Group
}
\label{kapitel-8}
\chaptermark{Complete Systems Which Admit All Transformations
of a One-term Group}

\setcounter{footnote}{0}

\abstract*{??}

If, in a $q$-term complete system:
\[
X_kf
=
\sum_{i=1}^r\,\xi_{ki}(x_1,\dots,x_n)\,
\frac{\partial f}{\partial x_i}
\ \ \ \ \ \ \ \ \ \ \ \ \
{\scriptstyle{(k\,=\,1\,\cdots\,q)}}
\]
one introduces new independent variables $x_1 ' = F_1 ( x_1, \dots,
x_n), \dots, x_n' = F_n$, then as was already observed earlier (cf.
Chap.~\ref{kapitel-5}, p.~\pageref{S-87}), one again obtains a
$q$-term complete system in $x_1', \dots, x_n'$. Naturally, this new
complete system has in general a form different from the initial one;
nonetheless, it can also happen that the two complete systems do
essentially not differ in their form, when relationships of the shape:
\[
X_kf
=
\sum_{j=1}^q\,\psi_{kj}(x_1',\dots,x_n')\,
\sum_{i=1}^n\,\xi_{ji}(x_1',\dots,x_n')\,
\frac{\partial f}{\partial x_i'}
\ \ \ \ \ \ \ \ \ \ \ \ \
{\scriptstyle{(k\,=\,1\,\cdots\,q)}}
\]
hold, in which of course the determinant of the $\psi_{ kj}$ does not
vanish identically. In this case, we say: \terminology{the complete
system:
\[
X_1f=0,\,\,\dots,\,\,X_qf=0
\]
admits the transformation $x_i' = F_i ( x_1, \dots, x_n)$}, or:
\terminology{it remains invariant through this
transformation}.

By making use of the abbreviated notations:
\[
\aligned
\xi_{ki}(x_1',\dots,x_n')
&
=
\xi_{ki}',
\\
\sum_{i=1}^n\,\xi_{ki}'\,
\frac{\partial f}{\partial x_i'}
&
=
X_k'f,
\endaligned
\]
we can set up the following definition:

\renewcommand{\thefootnote}{\fnsymbol{footnote}}
\plainstatement{The $q$-term complete system:
\[
X_kf
=
\sum_{i=1}^n\,\xi_{ki}(x_1,\dots,x_n)\,
\frac{\partial f}{\partial x_i}
=
0
\ \ \ \ \ \ \ \ \ \ \ \ \
{\scriptstyle{(k\,=\,1\,\cdots\,q)}}
\]
admits the infinitesimal transformation $x_i' = F_i ( x_1, \dots,
x_n)$ if and only if, for every $k$, there is a relation of the
form}:\footnote[1]{\,
Lie, Scientific Society of Christiania, February 1875.
} 
\def\theequation{1}\begin{equation}
X_kf
=
\sum_{j=1}^q\,\psi_{kj}(x_1',\dots,x_n')\, X_j'f.
\end{equation}
\renewcommand{\thefootnote}{\arabic{footnote}}

The sense of this important definition will best be 
elucidated \deutsch{verdeutlichen} by means of a simple example.

The two equations:
\[
X_1f
=
\frac{\partial f}{\partial x_1}
=
0,
\ \ \ \ \ \ \ \ \ \ \ \
X_2f
=
\frac{\partial f}{\partial x_2}
=
0
\]
in the three variables $x_1, x_2, x_3$ 
form a two-term complete system. By introducing the new variables:
\[
x_1'
=
x_1+x_2,
\ \ \ \ \ \ \
x_2'=x_1-x_2,
\ \ \ \ \ \ \
x_3'=x_3
\]
in place of the $x$, we obtain the new complete system:
\[
\frac{\partial f}{\partial x_1}
=
\frac{\partial f}{\partial x_1'}
+
\frac{\partial f}{\partial x_2'}
=
0,
\ \ \ \ \ \ \
\frac{\partial f}{\partial x_2}
=
\frac{\partial f}{\partial x_1'}
-
\frac{\partial f}{\partial x_2'}
=
0
\]
which is equivalent to the system $\frac{ \partial f}{ \partial x_1'}
= 0$, $\frac{\partial f}{ \partial x_2'} = 0$. 
So there are relations of the form:
\[
X_1f
=
X_1'f+X_2'f,
\ \ \ \ \ \ \ \ \ \ \ \ \ \
X_2f
=
X_1'f-X_2'f,
\]
whence the complete system $\frac{ \partial f}{ \partial x_1} = 0$,
$\frac{ \partial f}{ \partial x_2} = 0$ admits the transformation:
$x_1' = x_1 + x_2$, $x_2' = x_1 - x_2$.

\sectionengellie{\S\,\,\,37.}

Let the $q$-term complete system $X_1 f = 0, \, \dots, \, X_q f = 0$
admit the transformation $x_i ' = F_i ( x_1,
\dots, x_n)$, so that for every $f$ there are 
relations of the form~\thetag{ 1}. Now, if $\varphi ( x_1, \dots,
x_n)$ is a solution of the complete system, the right-hand side
of~\thetag{ 1} vanishes identically after the substitution $f =
\varphi ( x_1', \dots, x_n')$, whence the left-hand side also vanishes
identically after the substitution $f = \varphi \big( F_1 ( x), \dots,
F_n ( x) \big)$; consequently, as $\varphi ( x)$ itself, $\varphi
\big( F_1 ( x), \dots, F_n ( x) \big)$ constitutes at the same time a
solution, or, as one can express this fact: the transformation $x_i '
= F_i ( x)$ transfers every solution of the complete system $X_1f = 0,
\, \dots, \, X_q f = 0$ to a solution of the same complete system.

But conversely the following also holds true: if every solution of the
$q$-term complete system $X_1f = 0, \, \dots, \, X_q f = 0$ is
transferred, by means of a transformation $x_i' = F_i ( x_1, \dots,
x_n)$, to a solution, then the complete system admits the
transformation in question. Indeed, by the introduction of the
variables $x_i$ in place of the variables $x_i'$, the equations $X_1'
f = 0, \, \dots, \, X_q' f = 0$ convert into a $q$-term complete
system which has all its solutions in common with the $q$-term system
$X_1f = 0, \, \dots, \, X_q f = 0$; but from this it follows that
relations of the form~\thetag{ 1} hold, hence that the complete system
$X_1f = 0, \, \dots, \, X_q f = 0$ really admits the transformation
$x_i ' = F_i ( x)$.

From all that, it results that 
\label{S-138}
\emphasis{a $q$-term complete system
admits the transformation $x_i' = F_i ( x_1, \dots, x_n)$ if and only
if this transformation transfers every solution of the complete system
into a solution. Naturally, for this to hold, it is only necessary
that the transformation transfers to solutions any $n - q$ independent
solutions of the system.}

Next, the present line of thought will completely correspond to the
one followed in Chaps.~6 and~7 if we now ask the question: how can one
realize that a $q$-term complete system $X_1 f = 0, \, \dots, \, X_q f
= 0$ admits all transformations of the one-term group $Yf$? Indeed,
this question actually belongs to the general researches about
differential equations which admit one-term groups, and that is why 
this question
will also be took up again in a subsequent chapter of this Division,
in the Chapter on differential invariants, and on the basis of the
general theory developed there, it will be settled. But before, we
need criteria by means of which we can recognize whether a given
complete system admits or not all transformations of a given one-term
group. That is why we want to derive such criteria already now, with
somewhat simpler expedients.

Let us denote any $n-q$ independent solutions of the $q$-term
complete system $X_1 f = 0$, \dots, $X_q f = 0$ by $\varphi_1, 
\dots, \varphi_{ n-q}$. 
If now the complete system admits all transformations: 
\[
x_i'
=
x_i
+
\frac{t}{1}\,Yx_i
+
\frac{t^2}{1\cdot 2}\,YYx_i
+\cdots
\ \ \ \ \ \ \ \ \ \ \ \ \
{\scriptstyle{(i\,=\,1\,\cdots\,n)}}
\]
of the one-term group $Yf$, then the $n-q$
independent functions: 
\[
\aligned
\varphi_k\big(x+t\,Yx+\cdots\big)
&
=
\varphi_k(x)
+
\frac{t}{1}\,Y\varphi_k
+
\frac{t^2}{1\cdot 2}\,YY\varphi_k
+\cdots
\\
&
\ \ \ \ \ \ 
{\scriptstyle{(k\,=\,1\,\cdots\,n\,-\,q)}}
\endaligned
\]
must also be solutions of the system, and for every value of $t$.
From this, we deduce that the $n-q$ expressions $Y
\varphi_k$ are in any case 
solutions of the system, hence that relations of the
form:
\def\theequation{2}\begin{equation}
Y\varphi_k
=
\omega_k(\varphi_1,\dots,\varphi_{n-q})
\ \ \ \ \ \ \ \ \ \ \ \ 
{\scriptstyle{(k\,=\,1\,\cdots\,n\,-\,q)}}
\end{equation}
must hold. This condition is necessary; but at the
same time it is also sufficient, since if it
is satisfied, then all $YY \varphi_k$, 
$YYY \varphi_k$, \dots, will be functions of 
$\varphi_1, \dots, \varphi_{ n-q}$ only, whence
the expressions $\varphi_k ( x + t \, Yx + \cdots )$ 
will be solutions of the complete
system, and from this it follows that this system 
effectively admits all transformations of the one-term
group $Yf$. 

Consequently, the following holds.

\def\theproposition{1}\begin{proposition}
\label{Satz-1-S-139}
A $q$-term complete system $X_1f = 0$, \dots, $X_q f = 0$
with the $n-q$ independent solutions $\varphi_1, \dots, 
\varphi_{ n-q}$ admits all the transformations of the
one-term group:
\[
Yf
=
\sum_{i=1}^n\,\eta_i(x_1,\dots,x_n)\,
\frac{\partial f}{\partial x_i},
\]
if and only if $n-q$ relations of the form: 
\[
Y\varphi_k
=
\omega_k(\varphi_1,\dots,\varphi_{n-q})
\ \ \ \ \ \ \ \ \ \ \ \ 
{\scriptstyle{(k\,=\,1\,\cdots\,n\,-\,q)}}
\]
hold.
\end{proposition}

Saying this, the gained criterion naturally 
is practically applicable only when the complete
system is already integrated. 
But it is very easy to deduce from this criterion
another one which does not presupposes that one
knows the solutions of the complete system. 

If in the identity:
\[
X_k\big(Y(f)\big)
-
Y\big(X_k(f)\big)
\equiv
\sum_{i=1}^n\,
\big(
X_k\eta_i-Y\xi_{ki}
\big)\,
\frac{\partial f}{\partial x_i}
\]
one sets in place of $f$ any solution $\varphi$ of the
complete system, then one obtains:
\[
X_k\big(Y(f)\big)
\equiv
\sum_{i=1}^n\,
\big(
X_k\eta_i-Y\xi_{ki}
\big)
\frac{\partial\varphi}{\partial x_i}.
\]

Now, if the complete system allows all transformations of the
one-term group $Yf$, then by the above, 
$Y( \varphi)$ also a solution of the system, hence
the left-hand side of the last equation vanishes identically; 
naturally, 
the right-hand side does the same, so every solution of the 
complete system also satisfies the $q$ equations:
\[
\sum_{i=1}^n\,
\big(X_k\eta_i-Y\xi_{ki}\big)\,
\frac{\partial f}{\partial x_i}
=
0\,;
\]
but from this it follows that $q$ identities of the
form:
\[
\sum_{i=1}^n\,
\big(X_k\eta_i-Y\xi_{ki}\big)\,
\frac{\partial f}{\partial x_i}
\equiv
\sum_{j=1}^q\,\chi_{kj}(x_1,\dots,x_n)\,X_jf
\ \ \ \ \ \ \ \ \ \ \ \ \
{\scriptstyle{(k\,=\,1\,\cdots\,q)}}
\]
hold. 

On the other hand, we assume that identities of 
this form exist and we again understand by $\varphi$, 
any solution of the complete system. 
Then it follows immediately that the
$q$ expressions $X_k \big( Y ( f) \big)
- Y \big( X_k ( f) \big)$ vanish identically 
after the substitution $f = \varphi$; 
but from that, it comes: 
$X_k \big( Y ( \varphi) \big) \equiv 0$, 
that is to say $Y \varphi$ is a solution
of the system, so there exist $n-q$ relations
of the form:
\def\theequation{2}\begin{equation}
Y\varphi_k
=
\omega_k(\varphi_1,\dots,\varphi_{n-q})
\ \ \ \ \ \ \ \ \ \ \ \ \
{\scriptstyle{(k\,=\,1\,\cdots\,n\,-\,q)}},
\end{equation}
and they show that the complete system admits
all transformations of the one-term group $Y f$. 

As a result, we have the

\renewcommand{\thefootnote}{\fnsymbol{footnote}}
\def\thetheorem{20}\begin{theorem}
\footnote[1]{\,
Lie, Gesellschaft d. W. zu Christiania, 1874.
} 
A $q$-term complete system:
\label{Theorem-20-S-140}
\[
X_kf
=
\sum_{i=1}^n\,\xi_{ki}(x_1,\dots,x_n)\,
\frac{\partial f}{\partial x_i}
=
0
\ \ \ \ \ \ \ \ \ \ \ \ \
{\scriptstyle{(k\,=\,1\,\cdots\,q)}}
\]
in the variables $x_1, \dots, x_n$ admits all transformations
of the one-term group:
\[
Yf
=
\sum_{i=1}^n\,\eta_i(x_1,\dots,x_n)\,
\frac{\partial f}{\partial x_i}
\]
if and only of all $X_k \big( Y ( f) \big)
- Y \big( X_k ( f) \big)$ can be expressed linearly 
in terms of $X_1f, \dots, X_q f$:
\def\theequation{3}\begin{equation}
\aligned
\leftbracket X_k,\,Y\rightbracket
=
X_k\big(Y(f)\big)
&
-
Y\big(X_k(f)\big)
=
\sum_{j=1}^q\,\chi_{kj}(x_1,\dots,x_n)\,X_jf
\\
&
\ \ \ \ \ \ \ \ \ 
{\scriptstyle{(k\,=\,1\,\cdots\,q)}}.
\endaligned
\end{equation}
\end{theorem}
\renewcommand{\thefootnote}{\arabic{footnote}}

The existence of relations of the form~\thetag{ 3}
is therefore necessary and sufficient in order that
the system $X_1f = 0$, \dots, $X_q f = 0$
admits all transformations of the one-term group
$Yf$. For various fundamental reasons, 
it appears to be desirable to still prove the {\em necessity}
of the relation~\thetag{ 3} also by means of a direct method. 

The transformations of the one-term group $Yf$ have the form:
\[
x_i'
=
x_i
+
\frac{t}{1}\,\eta_i
+
\frac{t^2}{1\cdot 2}\,Y\eta_i
+\cdots
\ \ \ \ \ \ \ \ \ \ \ \ \
{\scriptstyle{(i\,=\,1\,\cdots\,n)}}.
\]
If we now introduce the new variables $x_1', \dots, x_n'$ in the
expressions $X_kf$ by means of this formula in place
of $x_1, \dots, x_n$, we receive:
\[
X_kf
=
\sum_{i=1}^n\,X_kx_i'\,\,
\frac{\partial f}{\partial x_i'}. 
\]
Here, the $X_k x_i'$ have still to be expressed in terms of $x_1',
\dots, x_n'$. By leaving out the second and the higher powers of
$t$, it comes immediately:
\def\theequation{4}\begin{equation}
X_kx_i'
=
X_kx_i+t\,X_k\eta_i+\cdots.
\end{equation}
Furthermore, one has (cf. Chap.~\ref{one-term-groups}, Eq.~\thetag{3a},
p.~\pageref{t-3a}):
\[
\aligned
&
X_kx_i
=
\xi_{ki}(x)
=
\xi_{ki}'
-
\frac{t}{1}\,Y'\xi_{ki}'
+\cdots
\\
&
t\,X_k\eta_i
=
t\,X_k'\eta_i'
-\cdots.
\endaligned
\]
Hence it comes:
\[
X_kx_i'
=
\xi_{ki}'+t\big(X_k'\eta_i'-Y'\xi_{ki}'\big)
+\cdots,
\]
and lastly:
\def\theequation{5}\begin{equation}
X_kf
=
X_k'f
+
\frac{t}{1}\,
\big(X_k'Y'f-Y'X_k'f\big)
+\cdots,
\end{equation}
where the left out terms are of second or of higher order in $t$. 

Now, if the complete system $X_1 f = 0$, \dots, $X_q f = 0$
admits all transformations of the one-term group $Yf$, then
the system of the $q$ equations:
\[
X_kf
+
\frac{t}{1}\,\leftbracket X_k,\,Y\rightbracket
+\cdots
=
0
\ \ \ \ \ \ \ \ \ \ \ \ \
{\scriptstyle{(k\,=\,1\,\cdots\,q)}}
\]
must be equivalent to $X_1 f = 0$, \dots, $X_q f = 0$ and
this, for all values of $t$. 
Consequently, the coefficients of $t$ must
in any case be expressible by means of $X_1f, \dots, X_qf$, that
is to say, relations of the form~\thetag{ 3}:
\[
\aligned
X_k\big(Y(f)\big)
-
Y\big(X_k(f)\big)
&
=
\leftbracket X_k,\,Y\rightbracket
=
\sum_{j=1}^q\,\chi_{kj}(x_1,\dots,x_n)\,X_jf
\\
&
\ \ \ \ \ 
{\scriptstyle{(k\,=\,1\,\cdots\,q)}}
\endaligned
\]
must exist. 

As a result, we have shown directly that the existence
of these relations is necessary; but we do not want to 
show once more its sufficiency, and we want to 
content ourselves with what has been said. 

We want to speak of complete systems which admit an infinitesimal
transformation $Yf$ in the same way as we spoke of the systems of
equations:
\[
\Omega_1(x_1,\dots,x_n)=0,\,\,\dots,\,\,
\Omega_{n-m}(x_1,\dots,x_n)=0
\]
which do the same. We shall say: \terminology{the $q$-term complete
system $X_1f = 0$, \dots, $X_q f = 0$ admits the infinitesimal
transformation $Yf$ when relations of the form~\thetag{ 3} hold}.

With the use of this language, we can also express the Theorem~20 as
follows:

\plainstatement{The $q$-term complete system $X_1f = 0$, \dots, $X_q f
= 0$ admits all transformations of the one-term group $Yf$ if and only
if it admits the infinitesimal transformation $Yf$}.

When a complete system admits all transformations of a one-term group,
we say briefly \terminology{that it admits this one-term group}.

The conditions of the Theorem~20 are in particular satisfied
when $Yf$ has the form:
\[
Yf
=
\sum_{j=1}^q\,\rho_j(x_1,\dots,x_n)\,X_jf,
\] 
where it is understood that the $\rho_j$ are arbitrary functions of
the $x$. Indeed, there always exist relations of the form~\thetag{ 3}
and in addition, it holds that $Y \varphi \equiv 0$. This is in
accordance with the developments of Chap.~\ref{kapitel-6},
p.~\pageref{S-98-bis} which showed that every solution $\Omega (
\varphi_1, \dots, \varphi_{ n-q})$ of the complete system $X_1f = 0,
\dots, X_q f = 0$ remains invariant by the transformations of all
one-term groups of the form $\sum \rho_j X_j f$.

On the other hand, if an infinitesimal transformation $Yf$ which does
not have the form $\sum \rho_j ( x) \, X_jf$ satisfies the conditions
of Theorem~20, then the finite transformations of the one-term group
$Yf$
do not at all leave invariant each individual solution of the complete
system, but instead, they leave invariant the
totality of all these solutions. 

Amongst the infinitesimal transformations that are admitted by a given
complete system $X_1 f = 0$, \dots, $X_q f = 0$, there are those of
the form $\sum\, \rho_j ( x) \, X_j f$ which are given simultaneously
with the complete system, and which, for this reason, have to
considered as trivial. By contrast, one cannot
in general indicate
the remaining infinitesimal transformations that the
system admits before one has integrated the system. 

As we have seen above, 
when the complete system $X_1 f = 0$, \dots, 
$X_q f = 0$ admits the infinitesimal transformation
$Yf$, every solution $\Omega ( \varphi_1, \dots, 
\varphi_{ n-q})$ of the complete system
is transferred to a solution by every transformation:
\[
x_i'
=
x_i
+
\frac{t}{1}\,Yx_i
+\cdots
\ \ \ \ \ \ \ \ \ \ \ \ \
{\scriptstyle{(i\,=\,1\,\cdots\,n)}}
\]
of the one-term group $Yf$. Hence if we interpret the $x$ and the
$x'$ as coordinates for the points of an $n$-times extended space and
the transformation just written as an operation by which the point
$x_1, \dots, x_n$ takes the new position $x_1', \dots, x_n'$, and if
we remember in addition that, with the constants $a_1, \dots, a_{
n-q}$, the equations $\varphi_1 = a_1$, \dots, $\varphi_{ n-q} = a_{
n-q}$ represent a characteristic manifold of the complete system $X_1
f = 0$, \dots, $X_q f = 0$ (cf. Chap.~\ref{kapitel-6},
p.~\pageref{S-101}), then we realize immediately that the
transformations of the one-term group $Yf$ send every characteristic
manifold of the complete system to a characteristic manifold.
\emphasis{The \label{S-143}
characteristic manifolds of our complete system are
therefore permuted with each other by the transformations of the
one-term group $Yf$, hence they form, as we want to express this, a
family \deutsch{Schaar} which is invariant by the one-term group
$Yf$.} Since the mentioned characteristic manifolds determine,
according to Chap.~\ref{kapitel-6}, p.~\pageref{S-101-bis}, a
decomposition of the space, we can also say that \emphasis{this
decomposition remains invariant by the one-term group $Yf$}.

\sectionengellie{\S\,\,\,38.}

Still a few simple statements about complete systems
which admit one-term groups:

\def\theproposition{2}\begin{proposition}
If a $q$-term complete system admits every transformation of 
the two one-term groups $Yf$ and $Zf$, then
it also admits every transformation of the
one-term group 
$Y \big( Z ( f) \big) 
- Z \big ( Y ( f) \big)
= \leftbracket Y, \, Z \rightbracket$.
\end{proposition}

Let $X_k f = 0$ be the equations of the complete system. 
One then forms the Jacobi identity:
\[
\big\leftbracket\leftbracket 
Y,\,Z\rightbracket,\,X_k\big\rightbracket
+
\big\leftbracket\leftbracket 
Z,\,X_k\rightbracket,\,Y\big\rightbracket
+
\big\leftbracket\leftbracket 
X_k,\,Y\rightbracket,\,Z\big\rightbracket
=
0,
\]
and one takes into consideration that $\leftbracket Y, \, X_k
\rightbracket$ and $\leftbracket Z, \, X_k \rightbracket$ can,
according to the assumption, be expressed linearly in terms of $X_1f,
\dots, X_q f$, hence one realizes that the same property is enjoyed by
$\big\leftbracket \leftbracket Y, \, Z \rightbracket, \, X_k
\rightbracket$. But as a result, our proposition is proved.

\def\theproposition{3}\begin{proposition}
If the equations $A_1 f = 0$, \dots, $A_q f = 0$ and likewise the
equations $B_1 f = 0$, \dots, $B_s f = 0$ form a complete system, then
a complete system is also formed by the totality of all possible
equations $Cf = 0$ which the two complete systems share, in the sense
that relations of the form:
\[
Cf
=
\sum_{j=1}^q\,\alpha_j(x_1,\dots,x_n)\,A_jf
=
\sum_{k=1}^s\,\beta_k(x_1,\dots,x_n)\,B_kf
\]
hold. If the two complete systems $A_k f = 0$ and $B_k f = 0$ admit a
certain one-term group $Xf$, then the complete system of equations $Cf
= 0$ also admits this group.
\end{proposition}

{\sf Proof.}
Amongst the equations $Cf = 0$, one can select
exactly, say $m$ but not more, equations:
\[
C_1f=0,\,\,\dots,\,\,C_mf=0
\]
which are independent of each other. Next, there exist for every $\mu
= 1, \dots, m$ relations of the form:
\[
C_\mu f
=
\sum_{j=1}^q\,\alpha_{\mu j}(x)\,A_jf
=
\sum_{k=1}^s\,\beta_{\mu k}(x)\,B_kf\,;
\]
consequently, every $C_\mu \big( C_\nu (f) \big) - C_\nu \big( C_\mu
(f) \big) = \leftbracket C_\mu, \, C_\nu \rightbracket$ can be
expressed both in terms of the $Af$ and in terms of the $Bf$, that is
to say, every $\leftbracket C_\mu, \, C_\nu \rightbracket$ expresses
linearly in terms of $C_1 f, \dots, C_m f$. As a result, the first
part of our proposition is proved. Further, every $\leftbracket X, \,
C_\mu \rightbracket$ expresses both in terms of the $Af$ and in terms
of the $Bf$, hence there there exist relations of the form:
\[
\leftbracket X,\,C_\mu\rightbracket
=
\sum_{\nu=1}^m\,\gamma_{\mu\nu}(x_1,\dots,x_n)\,C_\nu f
\ \ \ \ \ \ \ \ \ \ \ \ \
{\scriptstyle{(\mu\,=\,1\,\cdots\,m)}}.
\]
This is the second part of the proposition. 

If two complete systems $A_1f = 0$, \dots, $A_q f = 0$ and $B_1 f =
0$, \dots, $B_s f = 0$ are given, then all possible solutions that are
common to the to systems can be defined by means of a complete system
which, under the guidance of Chap.~\ref{kapitel-5}, p.~\pageref{S-85},
one can derive from the equations:
\[
A_1f=0,\,\,\dots,\,\,A_qf=0,\ \ \ \ \
B_1f=0,\,\,\dots,\,\,B_sf=0.
\]
What is more, the following holds.

\def\theproposition{3}\begin{proposition}
If the two complete systems $A_1 f = 0$, \dots, $A_q f = 0$ and $B_1 f
= 0$, \dots, $B_s f = 0$ admit the one-term group $Xf$, then the
complete system which defines the common solution of all the equations
$A_k f = 0$ and $B_k f = 0$ also admits the one-term group $Xf$.
\end{proposition}

{\sf Proof.}
The identity:
\[
\big\leftbracket\leftbracket 
A_j,\,B_k\rightbracket,\,X\big\rightbracket
+
\big\leftbracket\leftbracket 
B_k,\,X\rightbracket,\,A_j\big\rightbracket
+
\big\leftbracket\leftbracket 
X,\,A_j\rightbracket,\,B_k\big\rightbracket
=
0
\]
shows that all $\big\leftbracket \leftbracket A_j,\, B_k
\rightbracket,\, X\big \rightbracket$ can be expressed linearly in
terms of the $Af$, of the $Bf$ and of the $\leftbracket A, \, B
\rightbracket$. Hence, when the equations $\leftbracket A_j, \, B_k
\rightbracket = 0$, together with the $Af = 0$ and the $Bf = 0$,
already produce a complete system, the assertion of our proposition is
proved. Otherwise, one treats the system of the equations $Af = 0$,
$Bf = 0$, $\leftbracket A, \, B \rightbracket = 0$ exactly in the same
way as the system $A f = 0$, $Bf = 0$ just considered, that is to say:
one forms the Jacobi identity with $Xf$ for any two expressions
amongst the $Af$, $Bf$, $\leftbracket A, \, B \rightbracket$, and so
on.

\medskip

The propositions~3 and~4 can be given a simple conceptual sense when
$x_1, \dots, x_n$ are interpreted as coordinates for the points of a
space $R_n$.

At first, we remember that the $\infty^{ n-q}$ $q$-times extended
characteristic manifolds $M_q$ of the complete system $A_1 f = 0$,
\dots, $A_q f = 0$ form a family that remains invariant
by the one-term group $Xf$. Next, we remark that also the $\infty^{
n-s}$ $s$-times extended characteristic manifolds ${\sf M}_s$ of the
complete system $B_1 f = 0$, \dots, $B_s f = 0$ form such an
\emphasis{invariant family}.

Let us agree that the number $s$ is at least equal to $q$. Then every
${\sf M}_s$ of general position will be decomposed in a family of
$\infty^{s - q + h}$ $(q-h)$-extended manifolds, those in which it is
\label{S-145}
cut by the $M_q$, where one understands that $h$ is a determined
number amongst $0, 1, \dots, q$. Therefore in this way, the whole
$R_n$ is decomposed in a family of $\infty^{ n-q + h}$
$(q-h)$-extended manifolds. Naturally, the totality of these manifolds
remain invariant by the one-term group $Xf$, for it is the cutting of
the totality of all $M_q$ with the totality of all ${\sf M}_s$, and
these two totalities are, as already said, invariant by the group
$Xf$.

The considered $(q -h)$-times extended manifolds are nothing but the
characteristic manifolds of the complete system $Cf = 0$ which appears
in Proposition~3, and under the assumptions made, this complete system
is $(q-h)$-term.

On the other hand, one can ask for the smallest manifolds which
consist both of $M_q$ and of \label{S-146}
${\sf M}_s$.  When there are manifolds of this kind, the totality of
them naturally remains invariant by the one-term group $Xf$; they are
the characteristic manifolds of the complete system which is defined
in Proposition~4.

\linestop


\chapter{Characteristic Relationships
\\
Between the Infinitesimal Transformations of a Group}
\label{kapitel-9}
\chaptermark{Characteristic Relationships
Between the Infinitesimal Transformations of a Group}

\setcounter{footnote}{0}

\abstract*{??}

In Chap.~\ref{fundamental-differential} and in Chap.~4 (Proposition~1,
p.~\pageref{Satz-S-67}) it has been shown that to every
$r$-term group, there belong $r$ independent infinitesimal
transformations which stand in a characteristic relationship to the
group in question. Now, we want at first to derive certain important
relations which exist between these infinitesimal
transformations. Afterwards, we shall prove the equally important
proposition that $r$ independent infinitesimal transformations which
satisfy the concerned relations do always determine an $r$-term group
with the identity transformation.

\sectionengellie{\S\,\,\,39.}

Instead of considering an $r$-term group, we at the moment want to
take the somewhat more general point of view of considering a family
of $\infty^r$ different transformations:
\[
x_i'
=
f_i(x_1,\dots,x_n,\,a_1,\dots,a_r)
\ \ \ \ \ \ \ \ \ \ \ \ \
{\scriptstyle{(i\,=\,1\,\cdots\,n)}}
\]
which satisfy differential equations
of the form:
\[
\label{S-146-sq}
\aligned
\frac{\partial x_i'}{\partial a_k}
=
\sum_{j=1}^r\,
&
\psi_{kj}(a_1,\dots,a_r)\,\xi_{ji}(x_1',\dots,x_n')
\\
&
{\scriptstyle{(i\,=\,1\,\cdots\,n;\,\,\,k\,=\,1\,\cdots\,r)}}.
\endaligned
\]
We then know (cf. Chap.~\ref{fundamental-differential}, 
p.~\pageref{lemma-two-nondegeneracies})
that the $r$ infinitesimal transformations:
\[
X_k'(f)
=
\sum_{i=1}^n\,\xi_{ki}(x')\,
\frac{\partial f}{\partial x_i'}
\ \ \ \ \ \ \ \ \ \ \ \ \
{\scriptstyle{(k\,=\,1\,\cdots\,r)}}
\]
are independent of each other, and that the determinant of the $\psi_{
kj} (a)$ do not vanish identically; consequently, as we have already
done earlier, we can also write down 
the above differential equations as:
\def\theequation{1}\begin{equation}
\aligned
&
\xi_{ji}(x_1',\dots,x_n')
=
\sum_{k=1}^r\,\alpha_{jk}(a_1,\dots,a_r)\,
\frac{\partial x_i'}{\partial a_k}
\\
&
\ \ \ \ \ \ \ \ \ \ \ \ \ \ \ \ \ \ \
{\scriptstyle{(i\,=\,1\,\cdots\,n;\,\,\,j\,=\,1\,\cdots\,r)}}.
\endaligned
\end{equation}
Here naturally, the determinant of the $\alpha_{ jk}
(a)$ does
not vanish identically.

If on the other hand, we imagine that the equations:
\[
x_i'
=
f_i(x_1,\dots,x_n,\,a_1,\dots,a_r)
\ \ \ \ \ \ \ \ \ \ \ \ \
{\scriptstyle{(i\,=\,1\,\cdots\,n)}}
\]
are solved with respect to $x_1, \dots, x_n$:
\[
x_i
=
F_i(x_1',\dots,x_n',\,a_1,\dots,a_r)
\ \ \ \ \ \ \ \ \ \ \ \ \
{\scriptstyle{(i\,=\,1\,\cdots\,n)}},
\]
then we can very easily derive certain differential equations
that are satisfied by $F_1, \dots, F_n$. 
We simply differentiate the identity:
\[
F_i\big(f_1(x,\,a),\dots,f_n(x,\,a),\,
a_1,\dots,a_r\big)
\equiv
x_i
\]
with respect to $a_k$; then we have:
\[
\sum_{\nu=1}^n\,
\frac{\partial F_i(x',a)}{\partial x_\nu'}\,
\frac{\partial f_\nu(x,a)}{\partial a_k}
+
\frac{\partial F_i(x',a)}{\partial a_k}
\equiv 
0,
\]
provided that one sets $x_\nu ' = f_\nu (x, \, a)$ 
everywhere. We multiply this identity by $\alpha_{ jk} (a)$
and we sum it for $k$ equals $1$ to $r$; 
then on account of:
\[
\sum_{k=1}^r\,\alpha_{jk}(a)\,
\frac{\partial f_\nu(x,a)}{\partial a_k}
\equiv
\xi_{j\nu}(f_1,\dots,f_n),
\]
we obtain the following equations:
\[
\aligned
\sum_{\nu=1}^n\,\xi_{j\nu}(x_1',\dots,x_n')\,
&
\frac{\partial F_i}{\partial x_\nu'}
+
\sum_{k=1}^r\,\alpha_{jk}(a_1,\dots,a_r)\,
\frac{\partial F_i}{\partial a_k}
=
0
\\
&
{\scriptstyle{(i\,=\,1\,\cdots\,n;\,\,\,j\,=\,1\,\cdots\,r)}}.
\endaligned
\]

According to their derivation, these equations hold identically
when one makes the substitution $x_\nu ' = f_\nu (x, \, a)$ 
in them; 
but since they do not at all contain $x_1, \dots, x_n$, they
must actually hold identically, that is to say: 
$F_1, \dots, F_n$ are all solutions of the following
linear partial differential equations:
\def\theequation{2}\begin{equation}
\aligned
\Omega_j(F)
=
\sum_{\nu=1}^n\,\xi_{j\nu}(x')\,
&
\frac{\partial F}{\partial x_\nu'}
+
\sum_{\mu=1}^r\,\alpha_{j\mu}(a)\,
\frac{\partial F}{\partial a_\mu}
=
0
\\
&
{\scriptstyle{(j\,=\,1\,\cdots\,r)}}.
\endaligned
\end{equation}
These $r$ equations contain $n + r$ variables, namely $x_1', \dots,
x_n'$ and $a_1, \dots, a_r$; in addition, they are independent of each
other, for the determinant of the $\alpha_{ j\mu} (a)$ does
not vanish identically, and hence a resolution with respect to the $r$
differential quotients $\frac{ \partial F}{ \partial a_1}, \dots,
\frac{ \partial F}{ \partial a_r}$ is possible. But on the other
hand, the equations~\thetag{ 2} have $n$ independent solutions in
common, namely just the functions $F_1 (x', \, a), \dots, F_n( x', \,
a)$ whose functional determinant with respect to the $x'$:
\[
\sum\,\pm\,
\frac{\partial F_1}{\partial x_1'}\,\dots\,
\frac{\partial F_n}{\partial x_n'}
=
\frac{1}{
\sum\,\pm\,
\frac{\partial f_1}{\partial x_1}\,\dots\,
\frac{\partial f_n}{\partial x_n}
}
\]
does not vanish identically, because the equations $x_i' = f_i( x, \,
a)$ represent transformations by assumption. Therefore, the hypotheses
of the Proposition~8 in Chap.~\ref{kapitel-5}, 
p.~\pageref{Satz-8-S-88} are met by the equations~\thetag{
2}, that is to say, these equations constitute an $r$-term complete
system.

If we set:
\[
\sum_{k=1}^r\,\alpha_{jk}(a)\,
\frac{\partial F}{\partial a_k}
=
A_j(F)
\]
and furthermore:
\[
\sum_{\nu=1}^n\,\xi_{j\nu}(x')\,
\frac{\partial F}{\partial x_\nu'}
=
X_j'(F),
\]
in accordance with a designation employed earlier, 
then the equations~\thetag{ 2} receive the form:
\[
\Omega_j(F)
=
X_j'(F)
+
A_j(F)
=
0
\ \ \ \ \ \ \ \ \ \ \ \ \
{\scriptstyle{(j\,=\,1\,\cdots\,r)}}.
\]

As we know, the fact that they constitute a complete system is found
in their expressions: certain equations of the form:
\[
\aligned
\Omega_k\big(\Omega_j(F)\big)
-
\Omega_j\big(\Omega_k(F)\big)
&
=
\sum_{s=1}^r\,
\vartheta_{kjs}(x_1',\dots,x_n',a_1,\dots,a_r)\,
\Omega_s(F)
\\
&
{\scriptstyle{(k,\,j\,=\,1\,\cdots\,r)}}
\endaligned
\]
must hold identically\footnote{\,
These are identities between vector fields.
}, 
whichever $F$ can be as a function of $x_1', \dots, x_n', \, a_1,
\dots, a_r$. But since these identities can also be written as:
\[
\aligned
X_k'\big(X_j'(F)\big)
&
-
X_j'\big(X_k'(F)\big)
+
A_k\big(A_j(F)\big)
-
A_j\big(A_k(F)\big)
=
\\
&
=
\sum_{s=1}^r\,\vartheta_{kjs}\,X_s'(F)
+
\sum_{s=1}^r\,\vartheta_{kjs}\,A_s(F),
\endaligned
\]
we can immediately split them in two identities:
\def\theequation{3}\begin{equation}
\left\{
\aligned
X_k'\big(X_j'(F)\big)
-
X_j'\big(X_k'(F)\big)
&
=
\sum_{s=1}^r\,\vartheta_{kjs}\,X_s'(F)
\\
A_k\big(A_j(F)\big)
-
A_j\big(A_k(F)\big)
&
=
\sum_{s=1}^r\,\vartheta_{kjs}\,A_s(F),
\endaligned\right.
\end{equation}
and here, the second series of equations can yet again
be decomposed in:
\[
A_k\big(\alpha_{j\mu}\big)
-
A_j\big(\alpha_{k\mu}\big)
=
\sum_{s=1}^r\,\vartheta_{kjs}\,\alpha_{s\mu}
\ \ \ \ \ \ \ \ \ \ \ \ \
{\scriptstyle{(k,\,j,\,\mu\,=\,1\,\cdots\,r)}}.
\]

Now, because the determinant of the $\alpha_{ s\mu}$ does
not vanish identically, then the $\vartheta_{ kjs}$ in these last
conditions are completely determined, and it comes out that the
$\vartheta_{ kjs}$ can only depend upon $a_1, \dots, a_r$, whereas
they are in any case free of $x_1', \dots, x_n'$. But it can be
established that the $\vartheta_{ kjs}$ are also free of $a_1, \dots,
a_r$. For, if in the first series of the identities~\thetag{ 3}, we
consider $F$ as an arbitrary function of only $x_1', \dots, x_n'$,
then we obtain by differentiating with respect to $a_\mu$ the
following identically satisfied equations:
\[
0
\equiv
\sum_{s=1}^r\,\frac{\partial\vartheta_{kjs}}{\partial a_\mu}\,
X_s'(F)
\ \ \ \ \ \ \ \ \ \ \ \ \
{\scriptstyle{(k,\,j,\,\mu\,=\,1\,\cdots\,r)}}.
\]
But since $X_1' (F), \dots, X_r' ( F)$ are independent infinitesimal
transformations, and since in addition the $\frac{ \partial
\vartheta_{ kjs}}{\partial a_\mu}$ do not depend upon $x_1', \dots,
x_n'$, then all the $\frac{ \partial \vartheta_{ kjs }}{ \partial
a_\mu}$ vanish identically; that is to say, the $\vartheta_{ kjs}$ are
also free of $a_1, \dots, a_r$, they are numerical constants.

Thus, we have the

\def\thetheorem{21}\begin{theorem}
\label{Theorem-21-S-149}
If a family of $\infty^r$ transformations:
\[
x_i'
=
f_i(x_1,\dots,x_n,\,a_1,\dots,a_r)
\ \ \ \ \ \ \ \ \ \ \ \ \
{\scriptstyle{(i\,=\,1\,\cdots\,n)}}
\]
satisfies certain differential equations of the specific
form:
\[
\frac{\partial x_i'}{\partial a_k}
=
\sum_{j=1}^r\,\psi_{kj}(a_1,\dots,a_r)\,
\xi_{ji}(x_1',\dots,x_n')
\ \ \ \ \ \ \ \ \ \ \ \ \
{\scriptstyle{(i\,=\,1\,\cdots\,n;\,\,\,k\,=\,1\,\cdots\,r)}}
\]
and if one writes, what is always possible, these
differential equations in the form:
\[
\xi_{ji}(x_1',\dots,x_n')
=
\sum_{k=1}^r\,\alpha_{jk}(a_1,\dots,a_r)\,
\frac{\partial x_i'}{\partial a_k}
\ \ \ \ \ \ \ \ \ \ \ \ \
{\scriptstyle{(i\,=\,1\,\cdots\,n;\,\,\,j\,=\,1\,\cdots\,r)}},
\]
then there exist between the $2r$ independent infinitesimal
transformations:
\[
\aligned
X_j'(F)
&
=
\sum_{i=1}^n\,\xi_{ji}(x_1',\dots,x_n')\,
\frac{\partial F}{\partial x_i}
\ \ \ \ \ \ \ \ \ \ \ \ \
{\scriptstyle{(j\,=\,1\,\cdots\,r)}}
\\
A_j(F)
&
=
\sum_{\mu=1}^r\,\alpha_{j\mu}(a_1,\dots,a_r)\,
\frac{\partial F}{\partial a_\mu}
\ \ \ \ \ \ \ \ \ \ \ \ \
{\scriptstyle{(j\,=\,1\,\cdots\,r)}}
\endaligned
\]
relationships of the form:
\def\theequation{4}\begin{equation}
\left\{
\aligned
X_k'\big(X_j'(F)\big)
-
X_j'\big(X_k'(F)\big)
&
=
\sum_{s=1}^r\,c_{kjs}\,X_s'(F)
\ \ \ \ \ \ \ \ \ \ \ \ \
{\scriptstyle{(k,\,j\,=\,1\,\cdots\,r)}},
\\
A_k\big(A_j(F)\big)
-
A_j\big(A_k(F)\big)
&
=
\sum_{s=1}^r\,c_{kjs}\,A_s(F)
\ \ \ \ \ \ \ \ \ \ \ \ \
{\scriptstyle{(k,\,j\,=\,1\,\cdots\,r)}},
\endaligned\right.
\end{equation}
where the $c_{ kjs}$ denote numerical constants. In consequence of
that, the $r$ equations:
\[
X_j'(F)
+
A_j(F)
=
0
\ \ \ \ \ \ \ \ \ \ \ \ \
{\scriptstyle{(k\,=\,1\,\cdots\,r)}},
\]
which are
solvable with respect to $\frac{ \partial F}{ \partial a_1}$, \dots,
$\frac{ \partial F}{ \partial a_r}$, constitute an $r$-term complete
system in the $n+r$ variables $x_1', \dots, x_n'$, $a_1, \dots, a_r$;
if one solves the $n$ equations $x_i' = f_i ( x, \, a)$ with respect
to $x_1, \dots, x_n$:
\[
x_i
=
F_i(x_1',\dots,x_n',\,a_1,\dots,a_r)
\ \ \ \ \ \ \ \ \ \ \ \ \
{\scriptstyle{(i\,=\,1\,\cdots\,n)}},
\]
then $F_1 (x', \, a)$, \dots, $F_n( x', \, a)$ are independent
solutions of this complete system.
\end{theorem}

This theorem can now be immediately applied to all $r$-term groups,
whether or not they contain the identity transformation.

When applied to the case of an $r$-term group with the identity
transformation, the theorem gives us certain relationships which exist
between the infinitesimal transformations of this group. We therefore
obtain the important

\renewcommand{\thefootnote}{\fnsymbol{footnote}}
\def\thetheorem{22}\begin{theorem}
If an $r$-term group in the variables $x_1, \dots, x_n$ contains
the $r$ independent infinitesimal transformations:
\[
X_k(f)
=
\sum_{i=1}^n\,\xi_{ki}(x_1,\dots,x_n)\,
\frac{\partial f}{\partial x_i}
\ \ \ \ \ \ \ \ \ \ \ \ \
{\scriptstyle{(k\,=\,1\,\cdots\,r)}},
\]
then there exist between these infinitesimal transformations
pairwise relationships of the form:
\[
X_k\big(X_j(f)\big)
-
X_j\big(X_k(f)\big)
=
\sum_{s=1}^r\,c_{kjs}\,X_s(f),
\]
where the $c_{ kjs}$ denote numerical constants\footnote[1]{\,
Lie, Math. Ann. Vol. 8, p.~303; Göttinger Nachr. 1874.
}. 
\end{theorem}
\renewcommand{\thefootnote}{\arabic{footnote}}

From this, it follows in particular the following
important

\def\theproposition{1}\begin{proposition}
If a finite continuous group contains the two
infinitesimal transformations:
\[
X(f)
=
\sum_{i=1}^n\,\xi_i(x_1,\dots,x_n)\,
\frac{\partial f}{\partial x_i},
\ \ \ \ \ \ 
Y(f)
=
\sum_{i=1}^n\,\eta_i(x_1,\dots,x_n)\,
\frac{\partial f}{\partial x_i},
\]
then it also contains the infinitesimal transformation:
\[
X\big(Y(f)\big)
-
Y\big(X(f)\big).
\]
\end{proposition}

\sectionengellie{\S\,\,\,40.}

Conversely, we imagine that $r$ independent
infinitesimal transformations in $x_1', \dots, x_n'$:
\[
X_j'(F)
=
\sum_{i=1}^n\,\xi_{ji}(x_1',\dots,x_n')\,
\frac{\partial F}{\partial x_i'}
\ \ \ \ \ \ \ \ \ \ \ \ \
{\scriptstyle{(j\,=\,1\,\cdots\,r)}}
\]
are presented, which stand pairwise in relationships of the form:
\[
X_k'\big(X_j'(F)\big)
-
X_j'\big(X_k'(F)\big)
=
\sum_{s=1}^r\,c_{kjs}\,X_s'(F),
\]
where the $c_{ kjs}$ are numerical constants.
In addition, we imagine that $r$ infinitesimal
transformations in $a_1, \dots, a_r$:
\[
A_j(F)
=
\sum_{\mu=1}^r\,\alpha_{j\mu}(a_1,\dots,a_r)\,
\frac{\partial F}{\partial a_\mu}
\ \ \ \ \ \ \ \ \ \ \ \ \
{\scriptstyle{(j\,=\,1\,\cdots\,r)}}
\]
are given, which satisfy analogous
relations in pairs of the form:
\[
A_k\big(A_j(F)\big)
-
A_j\big(A_k(F)\big)
=
\sum_{s=1}^r\,c_{kjs}\,A_s(F),
\]
with the same $c_{ kjs}$ and whose determinant $\sum\, \pm 
\alpha_{ 11} ( a) \cdots \alpha_{ rr} ( a)$
does not vanish identically. We will show that under these
assumptions, the infinitesimal transformations $X_1 ' ( F), \dots, X_r
' ( F)$ generate a completely determined $r$-term group with the
identity transformation.

To this aim, we form the equations:
\[
\Omega_j(F)
=
X_j'(F)
+
A_j(F)
=
0
\ \ \ \ \ \ \ \ \ \ \ \ \
{\scriptstyle{(j\,=\,1\,\cdots\,r)}},
\]
which, according to the assumptions made, constitute an $r$-term
complete system; indeed, there exist relations of the form:
\[
\Omega_k\big(\Omega_j(F)\big)
-
\Omega_j\big(\Omega_k(F)\big)
=
\sum_{s=1}^r\,c_{kjs}\,\Omega_s(F)
\] 
and in addition, the equations $\Omega_1 ( F) = 0$, \dots, $\Omega_r (
F) = 0$ are solvable with respect to $\frac{ \partial F}{\partial
a_1}$, \dots, $\frac{ \partial F}{ \partial a_r}$.

Now, let $a_1^0, \dots, a_r^0$ be a system of values of the $a$, in a
neighbourhood of which the $\alpha_{ jk} (a)$ behave
regularly and for which the determinant $\sum \, \pm \alpha_{ 
11} ( a^0) \cdots \alpha_{ rr} ( a^0)$ is different
from zero. Then according to Theorem~12, Chap.~\ref{kapitel-5},
p.~\pageref{Theorem-12-S-91}, the complete system $\Omega_j ( F) = 0$
possesses $n$ solutions $F_1 ( x',\, a), \dots, F_n ( x', \, a)$ which
reduce to $x_1', \dots, x_n'$ respectively for $a_k = a_k^0$; they are
the so-called general solutions of the complete system relative to
$a_k = a_k^0$. We imagine that these general solutions are given, we
form the $n$ equations:
\[
x_i
=
F_i(x_1',\dots,x_n',\,a_1,\dots,a_r)
\ \ \ \ \ \ \ \ \ \ \ \ \
{\scriptstyle{(i\,=\,1\,\cdots\,n)}},
\]
and we resolve them with respect to 
$x_1', \dots, x_n'$, which is always possible, 
for $F_1, \dots, F_n$ are obviously independent of
each other, as far as $x_1', \dots, x_n'$ are
concerned.
The equations obtained in this way:
\[
x_i'
=
f_i(x_1,\dots,x_n,\,a_1,\dots,a_r)
\ \ \ \ \ \ \ \ \ \ \ \ \
{\scriptstyle{(i\,=\,1\,\cdots\,n)}}
\]
represent, as we shall now show, an $r$-term group
and in fact naturally, a group with the identity
transformation; because for $a_k = a_k^0$, one
gets $x_i' = x_i$.

At first, we have identically:
\def\theequation{5}\begin{equation}
\label{S-152}
\aligned
&
\sum_{\nu=1}^n\,\xi_{j\nu}(x')\,\frac{\partial F_i}{\partial x_\nu'}
+
\sum_{\mu=1}^r\,\alpha_{j\mu}(a)\,
\frac{\partial F_i}{\partial a_\mu}
=
0
\\
&
\ \ \ \ \ \ \ \ \ \ \ \ \
{\scriptstyle{(i\,=\,1\,\cdots\,n;\,\,\,j\,=\,1\,\cdots\,r)}}.
\endaligned
\end{equation}
On the other hand, by differentiating $x_i = F_i (x', \, a)$
with respect to $a_\mu$, one obtains the equation:
\[
0
=
\sum_{\nu=1}^n\,\frac{\partial F_i}{\partial x_\nu'}\,
\frac{\partial x_\nu'}{\partial a_\mu}
+
\frac{\partial F_i}{\partial a_\mu}
\ \ \ \ \ \ \ \ \ \ \ \ \
{\scriptstyle{(i\,=\,1\,\cdots\,n;\,\,\,\mu\,=\,1\,\cdots\,r)}},
\]
which holds identically after the substitution $x_\nu ' = f_\nu ( x,
\, a)$. We multiply this equation by $\alpha_{ j\mu} ( a)$
and we sum for $\mu$ equals 1 to $r$, hence we obtain an equation which,
after using~\thetag{ 5}, goes to:
\[
\aligned
&
\sum_{\nu=1}^n\,
\frac{\partial F_i}{\partial x_\nu'}\,
\Big(
\sum_{\mu=1}^r\,\alpha_{j\mu}(a)\,
\frac{\partial x_\nu'}{\partial a_\mu}
-
\xi_{j\nu}(x')
\Big)
=
0
\\
&
\ \ \ \ \ \ \ \ \ \ \ \ \ \ \ \ \ \ \ \ \ \ \
{\scriptstyle{(i\,=\,1\,\cdots\,n;\,\,\,\mu\,=\,1\,\cdots\,r)}}.
\endaligned
\]
But since the determinant $\sum \, \pm \frac{ \partial F_1}{\partial
x_1'} \, \dots \, \frac{ \partial F_n}{\partial x_n'}$ does not vanish
identically, we therefore obtain:
\[
\sum_{\mu=1}^r\,\alpha_{j\mu}(a)\,
\frac{\partial x_\nu'}{\partial a_\mu}
=
\xi_{j\nu}(x'),
\]
a system that we can again resolve with respect to the $\frac{
\partial x_\nu'}{ \partial a_\mu}$, for the determinant of the
$\alpha_{ j\mu} (a)$ does not vanish, indeed. Thus, we obtain
finally that equations of the form:
\def\theequation{6}\begin{equation}
\aligned
&
\frac{\partial x_\nu'}{\partial a_\mu}
=
\sum_{j=1}^r\,\psi_{\mu j}(a_1,\dots,a_r)\,
\xi_{j\nu}(x_1',\dots,x_n')
\\
&
\ \ \ \ \ \ \ \ \ \ \ \ \ \ \ \ \ \ \ \ \ \ \
{\scriptstyle{(\nu\,=\,1\,\cdots\,n;\,\,\,\mu\,=\,1\,\cdots\,r)}}
\endaligned
\end{equation}
do hold true, which naturally reduce to identities after the
substitution $x_\nu ' = f_\nu ( x, \, a)$.

At this point, the demonstration that the equations $x_i' = f_i (x, \,
a)$ represent an $r$-term group is not at all difficult.

Indeed, it is at first easy to see that the equations $x_i' = f_i( x,
\, a)$ represent $\infty^r$ distinct transformations, hence that the
parameters $a_1, \dots, a_r$ are all essential. Otherwise
indeed, all functions $f_1 ( x, \, a)$, \dots, $f_n( x, \, a)$ should
satisfy a linear partial differential equation of the form 
(cf. p.~\pageref{Theorem-essential}):
\[
\sum_{k=1}^r\,\chi_k(a_1,\dots,a_r)\,
\frac{\partial f}{\partial a_k}
=
0,
\]
where the $\chi_k$ would be free of $x_1, \dots, x_n$.
On account of~\thetag{ 6}, we would then have:
\[
\sum_{k,\,j}^{1\dots r}\,
\chi_k(a)\,\psi_{kj}(a)\,\xi_{j\nu}(f_1,\dots,f_n)
\equiv 
0
\ \ \ \ \ \ \ \ \ \ \ \ \
{\scriptstyle{(\nu\,=\,1\,\cdots\,n)}},
\]
whence, since $X_1' ( F), \dots, X_r' ( F)$ are
independent infinitesimal transformations:
\[
\sum_{k=1}^r\,\chi_k(a)\,\psi_{kj}(a)
=
0
\ \ \ \ \ \ \ \ \ \ \ \ \
{\scriptstyle{(j\,=\,1\,\cdots\,r)}}\,;
\]
but from this, it follows immediately: $\chi_1 ( a) = 0$, \dots,
$\chi_r ( a) = 0$, because the determinant of the $\psi_{ kj} ( a)$
does not vanish identically.

Thus, the equations $x_i' = f_i (x, \, a)$ effectively represent a
family of $\infty^r$ different transformations. But now this family
satisfies certain differential equations of the specific form~\thetag{
6}; we hence can immediately apply 
the Theorem~9 of Chap.~\ref{one-term-groups} on
p.~\pageref{Theorem-9-S-72}. 
According to it, the following holds
true: if $\overline{ a}_1, \dots, \overline{ a}_r$ is a system of
values of the $a$ for which the the determinant $\sum \, \pm \, \psi_{
11} ( \overline{ a} ) \cdots \psi_{ rr} ( \overline{ a})$ does not
vanish, and the $\psi_{ kj} (a)$ behave regularly, then every
transformation $x_i' = f_i( x, \, a)$ whose parameters $a_1, \dots,
a_r$ lie in a certain neighbourhood of $\overline{ a}_1, \dots,
\overline{ a}_r$, can be obtained by firstly executing the
transformation:
\[
\overline{x}_i
=
f_i(x_1,\dots,x_n,\,\overline{a}_1,\dots,\overline{a}_r)
\ \ \ \ \ \ \ \ \ \ \ \ \
{\scriptstyle{(i\,=\,1\,\cdots\,n)}}
\]
and then a transformation: 
\[
x_i'
=
\overline{x}_i
+
\sum_{k=1}^r\,\lambda_k\,\xi_{ki}(\overline{x})
+\cdots
\ \ \ \ \ \ \ \ \ \ \ \ \
{\scriptstyle{(i\,=\,1\,\cdots\,n)}}
\]
of a one-term group $\lambda_1 \, X_1 (f) + \cdots + \lambda_r \, X_r
( f)$, where it is understood that $\lambda_1, \dots, \lambda_r$ are
appropriate constants. If we set in particular $\overline{ a}_k =
a_k^0$, we then get $\overline{ x}_i = x_i$, hence we see at first
that the family of the $\infty^r$ transformations $x_i' = f_i (x, \,
a)$ coincides, in a certain neighbourhood of $a_1^0, \dots, a_r^0$,
with the family of the transformations: 
\def\theequation{7}\begin{equation}
\aligned
x_i'
=
x_i
&
+
\sum_{k=1}^r\,\lambda_k\,\xi_{ki}(x)
+\cdots
\\
&
{\scriptstyle{(i\,=\,1\,\cdots\,n)}}.
\endaligned
\end{equation}

If, on the other hand, we choose $\overline{ a}_1, \dots, \overline{
a}_r$ arbitrary in a certain neighbourhood of $a_1^0, \dots, a_r^0$,
then the transformation $\overline{ x}_i = f_i ( x, \,
\overline{ a})$ always belongs to the family~\thetag{ 7}. But if we first
execute the transformation $\overline{ x}_i = f_i(x, \, \overline{
a})$ and then an appropriate transformation:
\[
x_i'
=
\overline{x}_i
+
\sum_{k=1}^r\,\lambda_k\,\xi_{ki}(\overline{x})
+\cdots
\]
of the family~\thetag{ 7}, then by what has been said above, we obtain
a transformation $x_i' = f_i ( x, \, a)$, where $a_1, \dots, a_r$ can
take all values in a certain neighbourhood of $\overline{ a}_1, \dots,
\overline{ a}_r$. In particular, if we choose $a_1, \dots, a_r$ in the
neighbourhood of $a_1^0, \dots, a_r^0$ mentioned earlier on, which is
always possible, then again the transformation $x_i' = f_i ( x, \, a)$
also belongs to the family~\thetag{ 7}; Consequently\footnote{\,
The present goal is mainly to establish that if $r$ infinitesimal
transformations $X_1, \dots, X_r$ stand pairwise in the relationships
$\big[ X_k, \, X_j \big] = \sum_{ j=1}^r \, c_{ kjs}\, X_s$, where the
$c_{ kjs}$ are constants, then the totality of transformations $x' =
\exp \big( \lambda_1 X_1 + \cdots + \lambda_r X_r \big) ( x)$
constitutes an $r$-term continuous local Lie group; Theorem~24 below
will conclude such a fundamental statement. Especially, the
exponential family $x' = \exp \big( \lambda_1 X_1 + \cdots + \lambda_r
X_r \big) ( x)$ will be shown to be (after appropriate shrinkings)
closed under composition, a property that we may abbreviate informally
by $\exp \circ \exp \equiv \exp$.

However the Theorem~9 on p.~\pageref{Theorem-9-S-72} only said that $f
\circ \exp \equiv f$, or in greater details:
\def\theequation{$*$}\begin{equation}
\left(
\aligned
&
\overline{x}=f\big(x;\,\overline{a}\big)
\\
&\ \ 
\overline{a}\,\,\text{\rm near}\,\,a^0
\endaligned
\right)
\circ
\left(
\aligned
&
x'=\exp\big(\lambda\,X\big)(\overline{x})
\\
&\ \ \ \ \ \ \
\lambda\,\,\text{\rm near}\,\,0
\endaligned
\right)
\equiv
\left(
\aligned
&
x'=f(x;\,a)
\\
&\ \ 
a\,\,\text{\rm near}\,\,\overline{a}
\endaligned
\right).
\end{equation}
But if we now apply this statement with $\overline{ a} = a^0$ being
the system of values introduced while solving the complete system
$\Omega_1 ( F) = \cdots = \Omega_r ( f) = 0$, then by construction, $a^0$
yields the identity transformation $\overline{ x} = f ( x; \, a^0 ) =
x$ and we hence get in particular:
\[
\left(
\aligned
&
x'=\exp\big(\lambda\,X\big)(\overline{x})
\\
&\ \ \ \ \ \ \
\lambda\,\,\text{\rm near}\,\,0
\endaligned
\right)
\equiv
\left(
\aligned
&
x'=f(x;\,a)
\\
&\ \ 
a\,\,\text{\rm near}\,\,a^0
\endaligned
\right).
\]
We can therefore replace by exponentials the two occurrences of the
family $f (x; \, a)$ in~\thetag{ $*$} to see that the
exponential family is indeed closed under composition, as
was shortly claimed in the text.
}, 
we see that two transformations of the family~\thetag{ 7}, when
executed one after the other, do once again yield a transformation of
this family. As a result, this family, and naturally also the family
$x_i' = f_i (x, \, a)$ which identifies with it, forms an $r$-term
group, a group which contains the identity transformation and whose
transformations can be ordered as inverses in pairs. We can state the
gained result as follows:

\def\thetheorem{23}\begin{theorem}
\label{Theorem-23-S-154}
If $r$ independent infinitesimal transformations:
\[
X_k'(f)
=
\sum_{i=1}^n\,\xi_{ki}(x_1',\dots,x_n')\,
\frac{\partial f}{\partial x_i'}
\ \ \ \ \ \ \ \ \ \ \ \ \
{\scriptstyle{(k\,=\,1\,\cdots\,r)}}
\]
in the variables $x_1', \dots, x_n'$ satisfy conditions in pairs of
the form:
\[
X_k'\big(X_j'(f)\big)
-
X_j'\big(X_k'(f)\big)
=
\sum_{s=1}^r\,c_{kjs}\,X_s'(f),
\]
if furthermore $r$ independent infinitesimal transformations:
\[
A_k(f)
=
\sum_{\mu=1}^r\,
\alpha_{k\mu}(a_1,\dots,a_r)\,
\frac{\partial f}{\partial a_\mu}
\ \ \ \ \ \ \ \ \ \ \ \ \
{\scriptstyle{(k\,=\,1\,\cdots\,r)}}
\]
in the variables $a_1, \dots, a_r$ satisfy the analogous 
conditions:
\[
A_k\big(A_j(f)\big)
-
A_j\big(A_k(f)\big)
=
\sum_{s=1}^r\,c_{kjs}\,A_s(f)
\]
with the same $c_{ kjs}$, and if, in addition, the determinant $\sum \,
\pm \, \alpha_{ 11} (a) \cdots \alpha_{ rr} (
a)$ does not vanish identically, then one obtains in the following way
the equations of an $r$-term group: one forms the $r$-term complete
system:
\[
X_k'(f)+A_k(f)
=
0
\ \ \ \ \ \ \ \ \ \ \ \ \
{\scriptstyle{(k\,=\,1\,\cdots\,r)}}
\]
and one determines its general solutions relative to a suitable system
of values $a_k = a_k^0$. If $x_i = F_i ( x_1', \dots, x_n', \, a_1,
\dots, a_r)$ are these general solutions, then the equations $x_i' =
f_i ( x_1, \dots, x_n, \, a_1, \dots, a_r)$ which arise from
resolution represent an $r$-term continuous transformation group.
This group contains the identity transformation and for each one of
its transformations, yet the inverse transformation;
it is generated by the $\infty^{ r - 1}$ infinitesimal
transformations:
\[
\lambda_1\,X_1'(f)
+\cdots+
\lambda_r\,X_r'(f),
\]
where $\lambda_1, \dots, \lambda_r$ denote arbitrary constants. By
introducing new parameters in place of the $a_k$, the equations of the
group can therefore be brought to the form:
\[
x_i'
=
x_i
+
\sum_{k=1}^r\,\lambda_k\,\xi_{ki}(x)
+
\sum_{k,\,j}^{1\dots r}\,
\frac{\lambda_k\,\lambda_j}{1\,2}\,
X_j(\xi_{ki})
+\cdots
\ \ \ \ \ \ \ \ \ \ \ \ \
{\scriptstyle{(i\,=\,1\,\cdots\,n)}}.
\]
\end{theorem}

Obviously, even the equations $x_i = F_i ( x', \, a)$ 
appearing in this theorem do represent a group, and
in fact, just the group $x_i' = f_i ( x, \, a)$.

\sectionengellie{\S\,\,\,41.}

The hypotheses which are made in the important Theorem~23 can be
simplified in an essential way.

The theorem expresses that the $2r$ infinitesimal transformations $X_k
(f)$ and $A_k (f)$ determine a certain $r$-term group in the 
$x$-space; but
at the same time, there is a representation of this group which is
absolutely independent of the $A_k (f)$; indeed, according to the
cited theorem, the group in question identifies with the family of the
$\infty^{ r-1}$ one-term groups $\lambda_1 \, X_1 (f) + \cdots +
\lambda_r \, X_r (f)$, and this family is already completely
determined by the $X_k (f)$ alone. This circumstance conducts us to
the conjecture that the family of the $\infty^{ r-1}$ one-term groups
$\lambda_1 \, X_1 (f) + \cdots + \lambda_r \, X_r( f)$ always forms an
$r$-term group if and only if the independent infinitesimal
transformations $X_1 (f), \dots, X_r (f)$ stand pairwise in relations
of the form:
\def\theequation{8}\begin{equation}
X_k\big(X_j(f)\big)
-
X_j\big(X_k(f)\big)
=
\left[
X_k,\,X_j
\right]
=
\sum_{s=1}^r\,c_{kjs}\,X_s(f).
\end{equation}
According to the Theorem~22, this condition is necessary for the
$\infty^{ r-1}$ one-term groups $\sum \, \lambda_k \, X_k (f)$ to form
an $r$-term group. Thus our conjecture amounts to the fact that this
necessary condition is also sufficient.

This presumption would be brought to certainty if we would succeed in
producing, for every system of the discussed nature, $r$ independent
infinitesimal transformations in $(a_1, \dots, a_r)$:
\[
\label{S-156}
A_k(f)
=
\sum_{\mu=1}^r\,\alpha_{k\mu}(a_1,\dots,a_r)\,
\frac{\partial f}{\partial a_\mu}
\ \ \ \ \ \ \ \ \ \ \ \ \
{\scriptstyle{(k\,=\,1\,\cdots\,r)}}
\]
which satisfy the corresponding relations:
\[
A_k\big(A_j(f)\big)
-
A_j\big(A_k(f)\big)
=
\sum_{s=1}^r\,c_{kjs}\,A_s(f),
\]
while, however, the determinant $\sum \, \pm\, \alpha_{ 11}
\cdots \alpha_{ rr }$ does not vanish identically, or
expressed differently, while no relation of the form:
\[
\sum_{k=1}^r\,\chi_k(a_1,\dots,a_r)\,A_k(f)
=
0
\]
holds identically.

With the help of the Proposition~5 in Chap.~\ref{one-term-groups},
p.~\pageref{Satz-S-66}, we can now in fact always manage to
produce such a system of infinitesimal transformations $A_k ( f)$.
Similarly as at that previous time, we set:
\[
X_k^{(\mu)}(f)
=
\sum_{i=1}^n\,\xi_{ki}\big(x_1^{(\mu)},\dots,x_n^{(\mu)}\big)\,
\frac{\partial f}{\partial x_i^{(\mu)}},
\]
and we make up the $r$ infinitesimal transformations:
\[
W_k(f)
=
\sum_{\mu=1}^r\,X_k^{(\mu)}(f).
\]
According to the stated proposition, these infinitesimal
transformations have the property that no relation of the form:
\[
\sum_{k=1}^r\,\psi_k\big(
x_1',\dots,x_n',x_1'',\dots,x_n'',\cdots\cdots,
x_1^{(r)},\dots,x_n^{(r)}
\big)\,W_k(f)
=
0
\]
holds. Now, since in addition we have:
\[
W_k\big(W_j(f)\big)
-
W_j\big(W_k(f)\big)
=
\sum_{s=1}^r\,c_{kjs}\,W_s(f),
\]
the $r$ equations, independent of each other: \label{S-157}
\[
W_1(f)=0,\cdots\cdots,W_r(f)=0
\]
form an $r$-term complete system in the $rn$ variables $x_1', \dots,
x_n', \cdots, x_1^{(r)}, \dots, x_n^{ (r)}$. This complete system
possesses $r ( n-1)$ independent solutions, which can be called $u_1,
u_2, \dots, u_{ rn-r}$. Hence, if we select $r$ functions $y_1, \dots,
y_r$ of the $rn$ quantities $x_i^{ ( \mu)}$ that are independent of
each other and independent of $u_1, \dots, u_{ rn - r}$, we can
introduce the $y$ and the $u$ as new independent variables in place of
the $x_i^{ ( \mu)}$. By this, we obtain:
\[
W_k(f)
=
\sum_{\pi=1}^r\,
W_k(y_\pi)\,
\frac{\partial f}{\partial y_\pi}
+
\sum_{\tau=1}^{rn-r}\,W_k(u_\tau)\,\frac{\partial f}{\partial
u_\tau},
\]
or, since all $W_k ( u_\tau)$ vanish identically:
\[
W_k(f)
=
\sum_{\pi=1}^r\,\omega_{k\pi}
(y_1,\dots,y_r,u_1,\dots,u_{rn-r})\,
\frac{\partial f}{\partial y_\pi},
\]
where $W_1 ( f), \dots, W_r ( f)$ are linked by no relation 
of the form:
\[
\sum_{k=1}^r\,\varphi_k(y_1,\dots,y_r,\,
u_1,\dots,u_{rn-r})\,
W_k(f)
=
0.
\]
This property of the $W_k( f)$ remains naturally also true when we
confer to the $u_\tau$ appropriate fixed values $u_\tau^0$. If we
then set $\omega_{ k\pi} ( y, u^0) = \omega_{ k\pi}^0 ( y)$, 
the $r$ independent infinitesimal transformations 
in the independent variables $y_1, \dots, y_r$: 
\[
V_k(f)
=
\sum_{\pi=1}^r\,\omega_{k\pi}^0(y_1,\dots,y_r)\,
\frac{\partial f}{\partial y_\pi}
\]
stand pairwise in the relationships:
\[
V_k\big(V_j(f)\big)
-
V_j\big(V_k(f)\big)
=
\sum_{s=1}^r\,c_{kjs}\,V_s(f)
\]
and in addition, are linked by no relation of the form:
\[
\sum_{k=1}^r\,\varphi_k(y_1,\dots,y_r)\,V_k(f)
=
0.
\]
Consequently, the $V_k ( f)$ are infinitesimal transformations of the
required constitution. So we can immediately apply the Theorem~23
p.~\pageref{Theorem-23-S-154} to the $2r$ infinitesimal
transformations $X_1 ( f), \dots, X_r ( f)$, $V_1 ( f), \dots, V_r (
f)$ and as a result, we have proved that the $\infty^{ r - 1}$
one-term groups $\sum \, \lambda_k \, X_k ( f)$ constitute an $r$-term
group. Therefore the following holds true:

\renewcommand{\thefootnote}{\fnsymbol{footnote}}
\def\thetheorem{24}\begin{theorem}
\label{Theorem-24-S-158}
If $r$ independent infinitesimal transformations:
\[
X_k(f)
=
\sum_{i=1}^n\,\xi_{ki}(x_1,\dots,x_n)\,
\frac{\partial f}{\partial x_i}
\ \ \ \ \ \ \ \ \ \ \ \ \
{\scriptstyle{(k\,=\,1\,\cdots\,r)}}
\]
stand pairwise in the relationships:
\def\theequation{8}\begin{equation}
X_k\big(X_j(f)\big)
-
X_j\big(X_k(f)\big)
=
\left[
X_k,\,X_j
\right]
=
\sum_{s=1}^r\,c_{kjs}\,X_s(f),
\end{equation}
where the $c_{ kjs}$ are constants, then the totality
of the $\infty^{ r - 1}$ one-term groups:
\[
\lambda_1\,X_1(f)
+\cdots+
\lambda_r\,X_r(f)
\]
forms an $r$-term continuous group, which contains
the identity transformation and whose transformations
are ordered as inverses in pairs\footnote[1]{\,
Lie, Math. Annalen Vol.~8, p.~303, 1874; 
G\"ottinger Nachrichten, 1874, p.~533 and 540; 
Archiv for Math. og Naturv. Christiania 1878.
}. 
\end{theorem}

\renewcommand{\thefootnote}{\arabic{footnote}}
When we have $r$ independent infinitesimal transformations:
\[
X_1(f),\dots,X_r(f)
\]
which match the hypotheses of the above theorem, we shall henceforth
say that $X_1 ( f), \dots, X_r ( f)$ generate an $r$-term group, and
\terminology{we shall also virtually speak of the $r$-term group $X_1
( f), \dots, X_r ( f)$.}

\sectionengellie{\S\,\,\,42.}

Let $X_1 ( f), \dots, X_r ( f)$ be an $r$-term group in the variables
$x_1, \dots, x_n$ and $Y_1,
\dots Y_m ( f)$ be an $m$-term group in 
the same variables. The relations between the $X_k ( f)$ and the
$Y_\mu ( f)$ respectively may have the form:
\[
\aligned
X_k\big(X_j(f)\big)
-
X_j\big(X_k(f)\big)
&
=
\sum_{s=1}^r\,c_{kjs}\,X_s(f)
=
\left[X_k,\,X_j\right]
\\
Y_\mu\big(Y_\nu(f)\big)
-
Y_\nu\big(Y_\mu(f)\big)
&
=
\sum_{s=1}^m\,c_{\mu\nu s}'\,Y_s(f)
=
\left[Y_\mu,\,Y_\nu\right]
\\
&
\!\!\!\!\!\!\!\!\!\!\!\!\!\!\!\!
{\scriptstyle{(k,\,j\,=\,1\,\cdots\,r;\,\,\,\mu,\,\nu\,=\,1\,\cdots\,m)}}.
\endaligned
\]
Now, it can happen that these two groups have certain infinitesimal
transformations in common. We suppose that they have exactly $l$
independent such transformations in common, say:
\[
Z_\lambda(f)
=
\sum_{k=1}^r\,g_{\lambda k}\,X_k(f)
=
\sum_{\mu=1}^m\,h_{\lambda\mu}\,Y_\mu(f)
\ \ \ \ \ \ \ \ \ \ \ \ \
{\scriptstyle{(\lambda\,=\,1\,\cdots\,l)}},
\]
where the $g_{ \lambda k}$ and the $h_{ \lambda \mu}$ denote
constants. Then every other infinitesimal transformation contained in
the two groups can be linearly deduced from $Z_1 ( f), \dots, Z_l (
f)$. But if we form the expressions:
\[
Z_\lambda\big(Z_\nu(f)\big)
-
Z_\nu\big(Z_\lambda(f)\big)
=
\left[Z_\lambda,\,Z_\nu\right],
\]
we realize that they can be deduced linearly from $X_1 ( f), \dots,
X_r (f)$ and also from $Y_1 ( f), \dots, Y_m ( f)$ as well, hence that
they are common to the two groups. Consequently, relations of the
form:
\[
Z_\lambda\big(Z_\nu(f)\big)
-
Z_\nu\big(Z_\lambda(f)\big)
=
\left[Z_\lambda,\,Z_\nu\right]
=
\sum_{s=1}^l\,d_{\lambda\nu s}\,Z_s(f)
\]
hold true, that is to say $Z_1 ( f), \dots, Z_l ( f)$ generate an
$l$-term group.

As a result, we have the

\def\theproposition{2}\begin{proposition}
\label{Satz-2-S-159}
If the two continuous groups: $X_1 ( f), \dots, X_r ( f)$ and $Y_1 (
f), \dots, Y_m ( f)$ in the same variables have exactly $l$ and not
more independent infinitesimal transformations in common, then these
transformations generate, as far as they are concerned, an $l$-term
continuous group.
\end{proposition}

\sectionengellie{\S\,\,\,43.}

If the equations $x_i' = f_i ( x_1, \dots, x_n, \, a_1, \dots, a_r)$
represent a family of $\infty^r$ transformations and if in addition
they satisfy differential equations of the specific form:
\[
\frac{\partial x_i'}{\partial a_k}
=
\sum_{j=1}^r\,\psi_{kj}(a_1,\dots,a_r)\,
\xi_{ji}(x_1',\dots,x_n')
\ \ \ \ \ \ \ \ \ \ \ \ \
{\scriptstyle{(i\,=\,1\,\cdots\,n;\,\,\,k\,=\,1\,\cdots\,r)}},
\]
then as we know, the $r$ infinitesimal transformations:
\[
X_k(f)
=
\sum_{i=1}^n\,\xi_{ki}(x_1,\dots,x_n)\,
\frac{\partial f}{\partial x_i}
\ \ \ \ \ \ \ \ \ \ \ \ \
{\scriptstyle{(k\,=\,1\,\cdots\,r)}}
\]
are independent of each other, and in addition, 
according to Theorem~21 on p.~\pageref{Theorem-21-S-149}, 
they are linked together through relations of the form:
\[
X_k\big(X_j(f)\big)
-
X_j\big(X_k(f)\big)
=
\left[X_k,\,X_j\right]
=
\sum_{s=1}^r\,c_{kjs}\,X_s(f).
\]
Hence
the family of the $\infty^{ r-1}$ one-term groups:
\[
\lambda_1\,X_1(f)
+\cdots+
\lambda_r\,X_r(f)
\]
forms an $r$-term group with the identity transformation. 
Consequently, we can state the Theorem~9 on p.~\pageref{Theorem-9-S-72}
as follows:

\def\thetheorem{25}\begin{theorem}
\label{Theorem-25-S-160}
If a family of $\infty^r$ transformations:
\[
x_i'
=
f_i(x_1,\dots,x_n,\,a_1,\dots,a_r)
\ \ \ \ \ \ \ \ \ \ \ \ \
{\scriptstyle{(i\,=\,1\,\cdots\,n)}}
\]
satisfies certain differential equations of the form:
\[
\frac{\partial x_i'}{\partial a_k}
=
\sum_{j=1}^r\,\psi_{kj}(a_1,\dots,a_r)\,
\xi_{ji}(x_1',\dots,x_n')
\ \ \ \ \ \ \ \ \ \ \ \ \
{\scriptstyle{(i\,=\,1\,\cdots\,n\,;\,\,\,k\,=\,1\,\cdots\,r)}},
\]
and if $a_1^0, \dots, a_r^0$ is a system of values of the $a$ for
which the $\psi_{ kj} ( a)$ behave regularly and for which in
addition, the determinant:
\[
\sum\,\pm\,\psi_{11}(a)\cdots\psi_{rr}(a)
\]
is different from zero, then every transformation $x_i' = f_i ( x, \,
a)$ whose parameters $a_1, \dots, a_r$ lie in a certain neighbourhood
of $a_1^0, \dots, a_r^0$ can be thought to be produced by performing
firstly the transformation $\overline{ x}_i = f_i ( x, \, a^0)$ and
secondly, a completely determined transformation:
\[
x_i'
=
\overline{x}_i
+
\sum_{k=1}^r\,\lambda_k\,\xi_{ki}(\overline{x})
+\cdots
\ \ \ \ \ \ \ \ \ \ \ \ \
{\scriptstyle{(i\,=\,1\,\cdots\,n)}}
\]
of the $r$-term group which, under the
assumptions made, is generated by the $r$ independent infinitesimal
transformations:
\[
X_k(f)
=
\sum_{i=1}^n\,\xi_{ki}(x_1,\dots,x_n)\,
\frac{\partial f}{\partial x_i}
\ \ \ \ \ \ \ \ \ \ \ \ \
{\scriptstyle{(k\,=\,1\,\cdots\,r)}}.
\]
\end{theorem}

\smallercharacters{The above theorem is then 
of special interest when the equations $x_i'
= f_i ( x, a)$ represent an $r$-term group which does not contain the
identity transformation at least in the domain $( \! ( a ) \! )$. In
this case, we will yet derive a few important conclusions.

Thus, let $x_i' = f_i ( x_1, \dots, x_n, \, a_1, \dots, a_r)$
be an $r$-term group without identity transformation, 
and assume that the two transformations:
\[
\aligned
x_i'
&
=
f_i(x_1,\dots,x_n,\,a_1,\dots,a_r)
\\
x_i''
&
=
f_i(x_1',\dots,x_n',\,b_1,\dots,b_r)
\endaligned
\]
executed one after the other produce the transformation:
\[
x_i''
=
f_i(x_1,\dots,x_n,\,c_1,\dots,c_r)
=
f_i\big(x_1,\dots,x_n,\,\varphi_1(a,b),\dots,\varphi_r(a,b)\big).
\]
Here, if we employ the previous notation, $x_1, \dots, x_n$
lie arbitrarily in the domain $( \! ( x ) \! )$, $a_1, \dots, 
a_r$ and $b_1, \dots, b_r$ in the domain $( \! ( a ) \! )$, while
the positions of the $x_i'$, $x_i''$ and of the $c_k$ are determined
by the indicated equations. In addition, there are
still differential equations of the specific form:
\[
\aligned
\frac{\partial x_i'}{\partial a_k}
&
=
\sum_{j=1}^r\,\psi_{kj}(a_1,\dots,a_r)\,
\xi_{ji}(x_1',\dots,x_n')
\\
&
\ \ \ \ \ \ \ \ \ \ \ \ \
{\scriptstyle{(i\,=\,1\,\cdots\,n\,;\,k\,=\,1\,\cdots\,r)}}.
\endaligned
\]

In what follows, $a_1^0, \dots, a_r^0$ and likewise $b_1^0,
\dots, b_r^0$ should now denote a determined point \deutsch{Stelle} of
the domain $( \! ( a ) \! )$ and $\varphi_k ( a^0, b^0)$ should be
equal to $c_k^0$. By contrast, by $\overline{ a}_1,
\dots, \overline{ a}_r$ we want to understand
an arbitrary point in the domain $( a)$, so that the equations:
\[
\overline{x}_i
=
f_i(x_1,\dots,x_n,\,\overline{a}_1,\dots,\overline{a}_r)
\]
should represent any transformation of the given group.

Every transformation of the form $\overline{ x}_i = f_i ( x,
\overline{ a})$ can be obtained by performing at first the
transformation $x_i' = f_i ( x, a^0)$ and afterwards a certain second
transformation. In order to find the latter transformation, we solve
the equations $x_i' = f_i (x, a^0)$ with respect to $x_1, \dots, x_n$:
\[
x_i
=
F_i(x_1',\dots,x_n',\,a_1^0,\dots,a_r^0)
\]
and we introduce these values of the $x_i$ in $\overline{ x}_i = f_i (
x, \overline{ a})$. In this way, we obtain for the sought
transformation an expression of the form:
\def\theequation{9}\begin{equation}
\overline{x}_i
=
\Phi_i
\big(x_1',\dots,x_n',\,\overline{a}_1,\dots,\overline{a}_r)
\ \ \ \ \ \ \ \ \ \ \ \ \
{\scriptstyle{(i\,=\,1\,\cdots\,n)}}\,;
\end{equation}
we here do not write down the $a_k^0$, because we want to consider
them as numerical constants.

The transformation~\thetag{ 9} is well defined for all systems of
values $\overline{ a}_k$ in the domain $(a)$ and its expression can be
analytically continued to the whole domain of such systems of values;
this follows from the hypotheses that we have made previously about
the nature of the functions $f_i$ and $F_i$.

We now claim that for certain values of the parameters $\overline{
a}_k$, the transformations of the family $\overline{ x}_i = \Phi_i (
x', \, \overline{ a})$ belong to the initially given group $x_i' = f_i
( x, a)$, whereas by contrast, for certain other values of the
$\overline{ a}_k$, they belong to the group $X_1 f, \dots, X_r f$ with
identity transformation.

We establish as follows the first part of the claim stated just now. We
know that the two transformations:
\[
x_i'
=
f_i(x_1,\dots,x_n,\,a_1^0,\dots,a_r^0),
\ \ \ \ \ \ \ \ \ \ \
\overline{x}_i
=
f_i(x_1',\dots,x_n',\,b_1,\dots,b_r)
\]
executed one after the other produce the transformation $\overline{
x}_i = f_i ( x, c)$, where $c_k = \varphi_k ( a^0, b)$; here, we may
set for $b_1, \dots, b_r$ any system of values of the domain $( \! ( a
) \! )$, while the system of values $c_1, \dots, c_r$ then lies
in the domain $(a)$, in a certain neighbourhood of $c_1^0, \dots,
c_r^0$. But according to what has been said earlier, the
transformation $\overline{ x}_i = f_i ( x, c)$ is also obtained when
the two transformations:
\[
x_i'
=
f_i(x_1,\dots,x_n,\,a_1^0,\dots,a_r^0),
\ \ \ \ \ \ \ \ \ \ 
\overline{x}_i
=
\Phi_i(x_1',\dots,x_n',\,\overline{a}_1,\dots,\overline{a}_r)
\]
are executed one after the other and when one chooses $\overline{ a}_k =
c_k$. Consequently, the transformation $\overline{ x}_i =
\Phi_i ( x', \overline{ a})$ is identical, after
the substitution $\overline{ a}_k =
\varphi_k ( a^0, b)$, to the transformation $\overline{ x}_i = 
f_i ( x', b)$, that is to say: \emphasis{all transformations $\overline{
x}_i = \Phi_i ( x', \overline{ a})$ whose parameters $\overline{ a}_k$
lie in a certain neighbourhood of $c_1^0, \dots, c_r^0$ defined
through the equation $\overline{ a}_k =
\varphi_k ( a^0, b)$, belong to the presented group $x_i' = 
f_i( x, a)$}.

In order to establish the second part of our claim stated above, we
remember the Theorem~25.
If $\overline{ a}_1, \dots, \overline{ a}_r$ lie in a certain 
neighbourhood of $a_1^0, \dots, a_r^0$, then according to this
theorem, the transformation $\overline{ x}_i = 
f_i ( x, \overline{ a})$ can be obtained 
by performing at first the transformation:
\[
x_i'
=
f_i(x_1,\dots,x_n,\,a_1^0,\dots,a_r^0)
\]
and afterwards a completely determined transformation:
\def\theequation{10}\begin{equation}
\overline{x}_i
=
x_i'
+
\sum_{k=1}^r\,\lambda_k\,\xi_{ki}(x')
+\cdots
\end{equation}
of the $r$-term group that is generated by the
$r$ independent infinitesimal transformations:
\[
X_k(f)
=
\sum_{i=1}^n\,\xi_{ki}(x_1,\dots,x_n)\,
\frac{\partial f}{\partial x_i}
\ \ \ \ \ \ \ \ \ \ \ \ \
{\scriptstyle{(k\,=\,1\,\cdots\,r)}}.
\]
From earlier considerations (cf. Chap.~\ref{one-term-groups}, proof of
Theorem~9, p.~\pageref{Theorem-9-S-72}), we know in addition that one
finds the transformation~\thetag{ 10} in question when one chooses in
an appropriate way $\overline{ a}_1, \dots, \overline{ a}_r$ as
independent functions of $\lambda_1, \dots, \lambda_r$ and when one
then determines by resolution $\lambda_1, \dots, \lambda_r$ as
functions of $\overline{ a}_1, \dots, \overline{ a}_r$. But on the
other hand, we also obtain the transformation $\overline{ x}_i = f_i (
x, \overline{ a} )$ when we at first execute the transformation $x_i'
= f_i ( x, a^0)$ and afterwards the transformation $\overline{ x}_i =
\Phi_i ( x', \overline{ a})$. Consequently, the transformation
$\overline{ x}_i = \Phi_i ( x', \overline{ a})$ belongs to the group
generated by $X_1f, \dots, X_r f$ as soon as the system of values
$\overline{ a}_1, \dots, \overline{ a}_r$ lies in a certain
neighbourhood of $a_1^0, \dots, a_r^0$. Expressed differently: the
equations $\overline{ x}_i = \Phi_i ( x', \overline{ a})$ are
transferred to the equations~\thetag{ 10} when $\overline{ a}_1,
\dots, \overline{ a}_r$ are replaced by the functions of $\lambda_1,
\dots, \lambda_r$ discussed above.

With these words, our claim stated above is completely proved.

Thus the transformation equations $\overline{ x}_i =
\Phi_i ( x', \overline{ a })$
possess the following important property: if in place of the
$\overline{ a}_k$, the new parameters $b_1, \dots, b_r$ are introduced
by means of the equations $\overline{ a}_k = \varphi_k ( a^0, b)$,
then for a certain domain of the variables, the equations $\overline{
x}_i = \Phi_i ( x', \overline{ a})$ take the form $\overline{ x}_i =
f_i ( x', b)$; on the other hand, if in place of the $\overline{
a}_k$, the new parameters $\lambda_1, \dots, \lambda_r$ are
introduced, then for a certain domain, the equations $\overline{ x}_i
= \Phi_i ( x', \overline{ a})$ convert into:
\[
x_i
=
x_i'
+
\sum_{k=1}^r\,\lambda_k\,\xi_{ki}(x')
+\cdots
\ \ \ \ \ \ \ \ \ \ \ \ \
{\scriptstyle{(i\,=\,1\,\cdots\,n)}}
\]

Here lies an important feature of the initially given group $x_i' =
f_i ( x, a)$. Namely, when we introduce in the equations $x_i' = f_i
( x, a)$ the new parameters $\overline{ a}_1, \dots, \overline{ a}_r$
in place of the $a_k$ by means of $\overline{ a}_k =
\varphi_k ( a^0, a)$, 
then we obtain a
system of transformation equations $x_i' = \Phi_i ( x,
\overline{ a})$ which, by performing in addition 
its analytic continuation, represents a family of transformations to
which belong all transformations of some $r$-term group with identity
transformation.

We can also express this as follows.

\def\thetheorem{26}\begin{theorem}
\label{Theorem-26}
Every $r$-term group $x_i' = f_i ( x_1, \dots, x_n, \, a_1, \dots,
a_r)$ which is not generated by $r$ independent infinitesimal
transformations can be derived from an $r$-term group with $r$
independent infinitesimal transformations in the following way: one
sets up at first the differential equations:
\[
\frac{\partial x_i'}{\partial a_k}
=
\sum_{j=1}^r\,\psi_{kj}(a)\,\xi_{ji}(x')
\ \ \ \ \ \ \ \ \ \ \ \ \
{\scriptstyle{(i\,=\,1\,\cdots\,n\,;\,k\,=\,1\,\cdots\,r)}},
\]
which are satisfied by the equations $x_i' = 
f_i ( x, a)$, then one sets:
\[
\sum_{i=1}^n\,\xi_{ki}(x)\,\frac{\partial f}{\partial x_i}
=
X_k(f)
\ \ \ \ \ \ \ \ \ \ \ \ \
{\scriptstyle{(k\,=\,1\,\cdots\,r)}}
\]
and one forms the finite equations:
\[
x_i'
=
x_i
+
\sum_{k=1}^r\,\lambda_k\,\xi_{ki}(x)
+\cdots
\ \ \ \ \ \ \ \ \ \ \ \ \
{\scriptstyle{(i\,=\,1\,\cdots\,n)}}
\]
of the $r$-term group with identity transformation which is generated
by the $r$ independent infinitesimal transformations $X_1 f, \dots,
X_r f$. Then it is possible, in these finite equations, to
introduce new parameters $\overline{ a}_1, \dots, \overline{ a}_r$ in
place of $\lambda_1, \dots, \lambda_r$ in such a way that the
resulting transformation equations:
\[
x_i'
=
\Phi_i(x_1,\dots,x_n,\,\overline{a}_1,\dots,\overline{a}_r)
\ \ \ \ \ \ \ \ \ \ \ \ \
{\scriptstyle{(i\,=\,1\,\cdots\,n)}}
\]
represent a family of $\infty^r$ transformations which embrace, after
analytic continuation, all the $\infty^r$ transformations:
\[
x_i'
=
f_i(x_1,\dots,x_n,\,a_1,\dots,a_r)
\]
of the group.
\end{theorem}

\sectionengellie{\S\,\,\,44.}

Without difficulty, it can be shown that there really exist groups
which do not contain the identity transformation and also whose
transformations are not ordered as inverses by pairs.

The equation:
\[
x'
=
ax
\]
with the arbitrary parameter $a$ represents a one-term group. If one
executes two transformations:
\[
x'
=
ax,\ \ \ \ \ \ \ \
x''
=
bx'
\]
of this group one after the other, then one gets the transformation:
\[
x''=abx,
\]
which belongs to the group as well. From this, it results that the
family of all transformations of the form $x' = ax$, in which the
absolute value of $a$ is smaller than $1$, constitutes a group
too\footnote{\,
In the sought economy of axioms, what matters is only closure under
composition.
}. 
Obviously, neither this family contains the identity transformation,
nor its transformations order as inverses by pairs.

Hence, if there would exist an analytic expression which would
represent only those transformations of the form $x' = ax$ in which
the absolute value of $a$ is smaller than $1$, then with it, we would
have a finite continuous group without the identity transformation and
without inverse transformations.

Now, an analytic expression of the required constitution can
effectively be indicated.

It is known that the function\footnote{\,
A presentation of this passage has been anticipated in
Sect.~\ref{axiom-of-inverse}.
}: 
\[
\label{univalent-odd}
\sum_{\nu=1}^\infty\,
\frac{a^\nu}{1-a^{2\nu}}
\]
can be expanded, in the neighbourhood of $a = 0$, as an ordinary power
series with respect to $a$ which converges as long as the absolute
value $\vert a\vert$ of $a$ is smaller than $1$. The power series in
question has the form:
\[
\sum_{\mu=1}^\infty\,
k_\mu\,a^\mu
=
\omega(a),
\]
where the $k_\mu$ denote entire numbers depending on the index
$\mu$. Hence if we interpret the complex values of $a$ as points in a
plane, then $\omega ( a)$, as an analytic function of $a$, is defined
in the interior of the circle of radius $1$ which can be described
\deutsch{beschreiben} around the point $a = 0$.

Furthermore, it is known that the function $\omega ( a)$ has no more
sense for the $a$ whose absolute value equals $1$, so that the circle
in question around the point $a = 0$ constitutes the natural frontier
for $\omega ( a)$, across which this function cannot by analytically
continued.

We not set $\omega ( a) = \lambda$, and moreover, let $\vert a^0 \vert
< 1$ and $\omega (a^0) = \lambda^0$; then we can solve the equation
$\omega ( a) = \lambda$ with respect to $a$, that is to say, we can
represent $a$ as an ordinary power series in $\lambda - \lambda^0$ in
such a way that it gives: $a = a^0$ for $\lambda = \lambda^0$ and that
the equation $\omega ( a) = \lambda$ is identically satisfied after
substitution of this expression for $a$.

Let $a = \chi ( \lambda)$; then $\chi ( \lambda)$ is an analytic
function which takes only values whose absolute value are smaller than
$1$; this holds true not only for the found function element which is
represented in the neighbourhood of $\lambda = \lambda^0$ by an
ordinary power series in $\lambda - \lambda^0$, but also for every
analytic continuation of this function element.

Hence if we set: 
\[
x'
=
\chi(\lambda)\,x,
\]
we obtain the desired analytic expression for all transformations $x'
= ax$ in which $\vert a \vert$ is smaller than $1$. If now we have:
\[
x'=\chi(\lambda_1)\,x,
\ \ \ \ \ \ \ \ \
x''=\chi(\lambda_2)\,x',
\]
then it comes:
\[
x''=\chi(\lambda_1)\,\chi(\lambda_2)\,x\,;
\]
but this equation can always be brought to the form: 
\[
x''=\chi(\lambda_3)\,x\,;
\]
indeed, we have $\vert \chi ( \lambda_1) \, \chi ( \lambda_2) \vert
<1$; so if we set $\chi ( \lambda_1) \, \chi ( \lambda_2) = \alpha$,
we simply receive: $\lambda_3 = \omega ( \alpha)$.

With this, it is shown that the equation $x' = \chi ( \lambda) \, x$
with the arbitrary parameter $\lambda$ represents a group. This group
is continuous and finite, but lastly, neither it contains the
identity transformation, nor its transformations can be ordered as
inverses by pairs. Our purpose: the proof that there are groups of
this kind, is attained with it. Besides, it is easy to see that in a
similar manner, one can form arbitrarily many groups having this
constitution.

\smallskip
{\sf Remarks.}
In his first researches about finite continuous transformation groups,
\name{Lie} has attempted to show that {\em every} $r$-term group
contains the identity transformation plus infinitesimal
transformations, and is generated by the latter (cf. notably the two
articles in Archiv for Math. og Naturvid., Bd. 1, Christiania 1876).
However he soon realized that in his proof were made certain implicit
assumptions about the constitution of the occurring functions; as a
consequence, he restricted himself expressly to groups whose
transformations can be ordered as inverses by pairs and he showed that
in any case, the mentioned statement was correct for such groups
(Math. Ann. Vol.~16, p.~441 sq.).

Later, in the year 1884, \name{Engel} succeeded to make up a finite
continuous group which does not contain the identity transformation
and whose transformations do not order as inverses by pairs; this was
the group set up in the preceding paragraph.

Finally, \name{Lie} found that the equations of an {\em arbitrary}
finite continuous group with $r$ parameters can in any case be
derived, after introduction of new parameters and analytic
continuation, from the equations of an $r$-term group which contains
the identity transformation and $r$ independent infinitesimal
transformations, while its finite transformations can be ordered as
inverses by pairs (Theorem~26).

\sectionengellie{\S\,\,\,45.}

In Chap.~\ref{one-term-groups}, p.~\pageref{Theorem-11}, we found that
every $r$-term group which contains $r$ independent infinitesimal
transformations has the property that its finite transformations can
be ordered as inverses by pairs. On the other hand, it was mentioned
at the end of the previous paragraph that this statement can be
reversed, hence that every $r$-term group whose transformations order
as inverses by pairs, contains the identity transformation and is
produced by $r$ infinitesimal transformations. We will indicate how
the correctness of this assertion can be seen.

Let the equations:
\[
x_i'
=
f_i(x_1,\dots,x_n,\,a_1,\dots,a_n)
\ \ \ \ \ \ \ \ \ \ \ \ \
{\scriptstyle{(i\,=\,1\,\cdots\,n)}}
\]
with $r$ essential parameters $a_1, \dots, a_r$ represent an $r$-term
group with inverse transformations by pairs. By resolution with
respect to $x_1, \dots, x_n$, one can obtain:
\[
x_i
=
F_i(x_1',\dots,x_n',\,a_1,\dots,a_r)
\ \ \ \ \ \ \ \ \ \ \ \ \
{\scriptstyle{(i\,=\,1\,\cdots\,n)}}
\]

When the system of values $\varepsilon_1, \dots, \varepsilon_r$ lies
in a certain neighbourhood of $\varepsilon_1 = 0$, \dots,
$\varepsilon_r = 0$, then we can execute the two transformations:
\[
\aligned
x_i
&
=
F_i(x_1',\dots,x_n',\,a_1,\dots,a_r)
\\
x_i''
&
=
f_i(x_1,\dots,x_n,\,a_1+\varepsilon_1,\,\,\dots,\,\,
a_r+\varepsilon_r)
\endaligned
\]
one after the other and we thus obtain a transformation:
\[
x_i''
=
f_i\big(
F_1(x',a),\dots,F_n(x',a),\,\,
a_1+\varepsilon_1,\,\,\dots,\,\,a_r+\varepsilon_r\big)
\]
which likewise belongs to our group and which can be expanded in power
series with respect to $\varepsilon_1, \dots, \varepsilon_r$:
\def\theequation{11}\begin{equation}
x_i''
=
x_i'
+
\sum_{k=1}^r\,
\varepsilon_k\,
\bigg[
\frac{\partial f_i(x,a)}{\partial a_k}
\bigg]_{x=F(x',a)}
+\cdots.
\end{equation}

If we set here all $\varepsilon_k$ equal to zero, then we get the
identity transformation, which hence appears to our group. If on the
other hand we choose all the $\varepsilon_k$ infinitely small, we
obtain transformations of our group which are infinitely little
different from the identity transformation.

For brevity, we set:
\def\theequation{12}\begin{equation}
\bigg[
\frac{\partial f_i(x,a)}{\partial a_k}
\bigg]_{x=F(x',a)}
=
\eta_{ki}(x',a),
\end{equation}
so that the transformation~\thetag{ 11} hence receives the form:
\[
x_i''
=
x_i'
+
\sum_{k=1}^r\,\varepsilon_k\,\eta_{ki}(x',a)
+\cdots
\ \ \ \ \ \ \ \ \ \ \ \ \
{\scriptstyle{(i\,=\,1\,\cdots\,n)}}.
\]

Then, in the variables $x_1', \dots, x_n'$, we form the $r$
infinitesimal transformations:
\[
Y_k'f
=
\sum_{i=1}^n\,\eta_{ki}(x',a)\,
\frac{\partial f}{\partial x_k'}
\ \ \ \ \ \ \ \ \ \ \ \ \
{\scriptstyle{(k\,=\,1\,\cdots\,r)}},
\]
which are certainly independent from each other for undetermined
values of the $a_k$. Indeed, on the contrary case, the would be $r$
not all vanishing quantities $\rho_1, \dots, \rho_r$ which would not
depend upon $x_1', \dots, x_n'$ such that the equation:
\[
\sum_{k=1}^r\,\rho_k\,Y_k'f
=
0
\]
would be identically satisfied; 
then the $n$ relations: 
\[
\sum_{k=1}^r\,\rho_k\,\eta_{ki}(x',a)
=
0
\ \ \ \ \ \ \ \ \ \ \ \ \
{\scriptstyle{(i\,=\,1\,\cdots\,n)}}
\]
would follow and in turn, these relations would, after the
substitution $x_i' = f_i ( x, a)$, be transferred to:
\[
\sum_{k=1}^r\,\rho_k\,
\frac{\partial f_i(x,a)}{\partial a_k}
=
0
\ \ \ \ \ \ \ \ \ \ \ \ \
{\scriptstyle{(i\,=\,1\,\cdots\,n)}}\,;
\]
but such relations cannot exist, since according to the assumption,
the parameters are essential in the equations $x_i' = f_i ( x, a)$
(cf. Chap.~\ref{essential-parameters}, p.~\pageref{iii-ess-param}).

From this, we then conclude that $Y_1' f, \dots, Y_r ' f$ also remain
independent from each other, when one inserts for $a_1, \dots, a_r$
some determined system of values in general position. If
$\overline{ a}_1, \dots, \overline{ a}_r$ is
such a system of values, we want to write:
\[
\eta_{ki}(x',\overline{a})
=
\xi_{ki}(x')\,;
\]
then the $r$ infinitesimal transformations:
\[
X_k'f
=
\sum_{i=1}^n\,\xi_{ki}(x')\,
\frac{\partial f}{\partial x_i'}
\ \ \ \ \ \ \ \ \ \ \ \ \
{\scriptstyle{(k\,=\,1\,\cdots\,r)}}
\]
are also independent of each other. 
It remains to show that our group is generated
by the $r$ infinitesimal transformations
$X_k' f$. 

We execute two transformations of our group one after the
other, namely firstly the transformation:
\[
x_i''
=
x_i'
+
\sum_{k=1}^r\,\varepsilon_k\,\eta_{ki}(x',a)
+\cdots,
\]
and secondly the transformation:
\[
\aligned
x_i'''
&
=
x_i''+\sum_{k=1}^r\,\vartheta_k\,\eta_{ki}(x'',\overline{a})
+\cdots
\\
&
=
x_i''+\sum_{k=1}^r\,\vartheta_k\,\xi_{ki}(x'')
+\cdots.
\endaligned
\]
If, as up to now, we only take into consideration
the first-order terms, we then obtain in the 
indicated way the transformation:
\[
x_i'''
=
x_i'
+
\sum_{k=1}^r\,\varepsilon_k\,\eta_{ki}(x',a)
+
\sum_{k=1}^r\,\vartheta_k\,\xi_{ki}(x')
+\cdots,
\]
which belongs naturally to our group, and this, inside a certain
region for all values of the parameters $a$, $\varepsilon$,
$\vartheta$.

If there would be, amongst the infinitesimal transformations $Y_1'f,
\dots, Y_r' f$, also a certain number which would be independent of
$X_1' f, \dots, X_r' f$, then the last written equations would,
according to Chap.~\ref{one-term-groups}, Proposition~4,
p.~\pageref{Satz-4-S-65}, represent at least $\infty^{ r+1}$
different transformations, whereas our group nevertheless contain only
$\infty^r$ different transformations. Consequently, each
transformation $Y_k' f$ must be linearly expressible in terms of $X_1'
f, \dots, X_r' f$, whichever values the $a$ can have. By
considerations similar to those of Chap.~\ref{one-term-groups}
one now realizes that $r$ identities 
of the form:
\[
Y_k'f
\equiv
\sum_{j=1}^r\,
\psi_{kj}(a_1,\dots,a_r)\,X_j'f
\ \ \ \ \ \ \ \ \ \ \ \ \
{\scriptstyle{(k\,=\,1\,\cdots\,r)}}
\]
hold true, where the $\psi_{ kj}$ behave regularly in a certain
neighbourhood of $a_k = \overline{ a}_k$; in addition, the determinant
of the $\psi_{ kj}$ does not vanish identically, since otherwise $Y_1'
f, \dots, Y_r'f$ would not anymore be independent infinitesimal
transformations.

At present, it is clear that the $\eta_{ ki} (x', a)$ can be expressed
as follows in terms of the $\xi_{ ji} ( x')$:
\[
\eta_{ki}(x',a)
\equiv
\sum_{j=1}^r\,\psi_{kj}(a)\,\xi_{ji}(x').
\]

Finally, if we remember the equations~\thetag{ 12} which define the
functions $\eta_{ ki} (x', a)$, we realize that the differential
equations:
\def\theequation{13}\begin{equation}
\aligned
\frac{\partial x_i'}{\partial a_k}
&
=
\sum_{j=1}^r\,\psi_{kj}(a)\,\xi_{ji}(x')
\\
&
\ \,
{\scriptstyle{(i\,=\,1\,\cdots\,n\,;\,\,\,
k\,=\,1\,\cdots\,r)}}
\endaligned
\end{equation}
are identically satisfied after the substitution $x_i' = f_i( x, a)$.

As a result, it has been directly shown that every group with inverse
transformations by pairs satisfies certain differential equations of
the characteristic form~\thetag{ 13}; insofar, we have reached the
starting point for the developments of
Chap.~\ref{fundamental-differential}, Sect.~\ref{substituting-axiom}.

\emphasis{Now, if it would be possible to prove that, for the
parameter values $a_1^0, \dots, a_r^0$ of the identity transformation,
the determinant of the $\psi_{ kj} ( a)$ has a value distinct from
zero, then it would follow from the stated developments that the group
$x_i' = f_i ( x, a)$ is generated by the $r$ infinitesimal
transformations $X_1f, \dots, X_r f$}. But now, it is not in the
nature of things that one can prove that the determinant $\sum \pm
\psi_{ 11} (a^0) \cdots \psi_{ rr} ( a^0)$ is distinct from zero. One
can avoid this trouble as follows\footnote{\,
For (local) continuous finite transformation groups 
containing the identity transformation, this property 
has already been seen in Sect.~\ref{derivation-fundamental}.
}. 

One knows that the equations:
\def\theequation{14}\begin{equation}
\aligned
x_i'
&
=
x_i
+
\sum_{k=1}^r\,\varepsilon_k\,\xi_{ki}(x)
+\cdots
\\
&
\ \ \ \ \ \ \ \ \ \ \ \ \
{\scriptstyle{(i\,=\,1\,\cdots\,n)}}
\endaligned
\end{equation}
represent transformations of our group $x_i' = f_i( x, a)$ as soon as
the $\varepsilon_k$ lie in a certain neighbourhood of $\varepsilon_1 =
0$, \dots, $\varepsilon_r = 0$; besides, one can show that one obtains
in this way \emphasis{all} transformations $x_i' = f_i ( x, a)$ the
parameters of which lie in a certain neighbourhood of $a_1^0, \dots,
a_r^0$. Now, if one puts the equations~\thetag{ 14} at the
foundation, one easily realizes that the $x'$, interpreted as
functions of the $\varepsilon$ and of the $x$, satisfy differential
equations of the form:
\[
\aligned
\frac{\partial x_i'}{\partial\varepsilon_k}
&
=
\sum_{j=1}^r\,
\chi_{kj}(\varepsilon_1,\dots,\varepsilon_r)\,
\xi_{ji}(x_1',\dots,x_n')
\\
&
\ \ \ \ \ \ \ \ \ \ \ \ \
{\scriptstyle{(i\,=\,1\,\cdots\,n\,;\,\,
k\,=\,1\,\cdots\,r)}},
\endaligned
\]
where now the determinant of the $\chi_{kj} ( \varepsilon)$ for
$\varepsilon_1 = 0$, \dots, $\varepsilon_r = 0$ does not vanish. In
this way, one comes finally to the following result:

\plainstatement{\label{S-169}
Every $r$-term group with transformations inverse by pairs contains
the identity transformation, and in addition $r$ independent
infinitesimal transformations by which it is generated.}

}

\sectionengellie{\S\,\,\,46.}

\label{S-169-sq}
Consider $r$ independent infinitesimal transformations:
\[
X_kf
=
\sum_{\nu=1}^n\,\xi_{k\nu}(x_1,\dots,x_n)\,
\frac{\partial f}{\partial x_\nu}
\ \ \ \ \ \ \ \ \ \ \ \ \
{\scriptstyle{(k\,=\,1\,\cdots\,r)}}
\]
which satisfy relations in pairs of the form:
\def\theequation{8}\begin{equation}
\aligned
X_i\big(X_k(f)\big)
&
-
X_k\big(X_i(f)\big)
=
\big[X_i,\,X_k\big]
=
\sum_{s=1}^r\,c_{iks}\,X_s(f)
\\
&\ \ \ \ \ \ \ \ \ \ \ \ \ \ \ \ \ \ \ 
{\scriptstyle{(i,\,k\,=\,1\,\cdots\,r)}},
\endaligned
\end{equation}
with certain constants $c_{ iks}$, so that
according to Theorem~24 p.~\pageref{Theorem-24-S-158}, 
the totality of all one-term groups of the
form:
\[
\lambda_1X_1(f)+\cdots+\lambda_rX_r(f)
\]
constitutes an $r$-term group. We will show that, as far as they are
concerned, the constants $c_{ iks}$ in the above relations are then
tied up together with certain equations.

To begin with, we have $\big[ X_i, \, X_k \big] = - \big[ X_k, \, X_i
\big]$, from which it comes immediately: $c_{ iks} = - c_{ kis}$.
When the number $r$ is greater than $2$, we find still other
relations. Indeed, in this case, there is the Jacobi identity
(Chap.~\ref{kapitel-5}, 
\S\,\,26, p.~\ref{jacobi-identity}) which holds between
any three $X_i f$, $X_k f$, $X_j f$ amongst the $r$ infinitesimal
transformations $X_1 f, \dots, X_r f$:
\[
\aligned
\big[\big[X_i,\,X_k\big],\,X_j\big]
&
+
\big[\big[X_k,\,X_j\big],\,X_i\big]
+
\big[\big[X_j,\,X_i\big],\,X_k\big]
=
0
\\
&
\ \ \ \ \ \ \ \ \ \ \ \ \
{\scriptstyle{(i,\,k,\,j\,=\,1\,\cdots\,n)}}.
\endaligned
\]

By making use of the above relation~\thetag{ 8}, we obtain firstly
from this:
\[
\sum_{s=1}^r\,
\Big\{
c_{iks}\,\big[X_s,\,X_j\big]
+
c_{kjs}\,\big[X_s,\,X_i\big]
+
c_{jis}\,\big[X_s,\,X_k\big]
\Big\}
=
0,
\]
and then by renewed applications of this relation:
\[
\sum_{s,\,\tau}^{1\cdots r}\,
\Big\{
c_{iks}\,c_{sj\tau}
+
c_{kjs}\,c_{si\tau}
+
c_{jis}\,c_{sk\tau}
\Big\}\,
X_\tau f
=
0.
\]
But since the infinitesimal transformations $X_\tau f$ are independent
of each other, this equation decomposes in the following
$r$ equations:
\def\theequation{15}\begin{equation}
\aligned
\sum_{s=1}^r\,
&
\big\{
c_{iks}\,c_{sj\tau}
+
c_{kjs}\,c_{si\tau}
+
c_{jis}\,c_{sk\tau}
\big\}
=
0
\\
&\ \ \ \ \ \ \ \ \ \ \ \ \ \ \ \ \ \
{\scriptstyle{(\tau\,=\,1\,\cdots\,r)}}.
\endaligned
\end{equation}

Thus the following holds.
\renewcommand{\thefootnote}{\fnsymbol{footnote}}

\def\thetheorem{27}\begin{theorem}\footnote[1]{\,
Lie, Archiv for Math. og Naturv.
Bd. 1, p.~192, Christiania 1876.
} 
If $r$ independent infinitesimal transformations $X_1f, \dots, X_r f$
are constituted in such a way that they satisfy relations in pairs of
the form:
\label{Theorem-27-S-170}
\def\theequation{8}\begin{equation}
\aligned
X_i\big(X_k(f)\big)
&
-
X_k\big(X_i(f)\big)
=
\big[X_i,\,X_k\big]
=
\sum_{s=1}^r\,c_{iks}\,X_s(f)
\\
&\ \ \ \ \ \ \ \ \ \ \ \ \ \ \ \ \ \ \ 
{\scriptstyle{(i,\,k\,=\,1\,\cdots\,r)}},
\endaligned
\end{equation}
with certain constants $c_{ iks}$, so that the totality of all
$\infty^{ r - 1}$ one-term groups of the form:
\[
\lambda_1\,X_1(f)
+\cdots+
\lambda_r\,X_r(f)
\]
forms an $r$-term group, then between the constants $c_{ iks}$, 
there exist the following relations:
\def\theequation{16}\begin{equation}
\label{S-170}
\left\{
\aligned
&\ \ \ \ \ \ \ \ \ \ \ \ \
c_{ik\tau}+c_{ki\tau}
=
0
\\
&
\sum_{s=1}^r\,
\big\{
c_{iks}\,c_{sj\tau}+c_{kjs}\,c_{si\tau}+c_{jis}\,c_{sk\tau}
\big\}
=
0
\\
&
\ \ \ \ \ \ \ \ \ \ \ \ \
{\scriptstyle{(i,\,k,\,j,\,\tau\,=\,1\,\cdots\,r)}}.
\endaligned\right.
\end{equation}
\end{theorem}

\renewcommand{\thefootnote}{\arabic{footnote}}
The equations~\thetag{ 16} are completely independent of the number of
the variables $x_1, \dots, x_n$ in the infinitesimal transformations
$X_1 ( f), \dots, X_r ( f)$. Hence from the the preceding statement,
we can still conclude what follows:

Even if the number $n$ of the independent variables $x_1, \dots, x_n$
can be chosen in an completely arbitrary way, one nevertheless cannot
associate to every system of constants $c_{ iks}$ ${\scriptstyle{ ( i,
\, k, \, s \, = \, 1 \, \cdots \, r )}}$ $r$ independent infinitesimal
transformations:
\[
X_k(f)
=
\sum_{\nu=1}^n\,\xi_{k\nu}(x_1,\dots,x_n)\,
\frac{\partial f}{\partial x_\nu}
\ \ \ \ \ \ \ \ \ \ \ \ \
{\scriptstyle{(k\,=\,1\,\cdots\,r)}}
\]
which pairwise satisfy the relations:
\[
\big[X_i,\,X_k\big]
=
\sum_{s=1}^r\,c_{iks}\,X_s(f)
\ \ \ \ \ \ \ \ \ \ \ \ \
{\scriptstyle{(i,\,k\,=\,1\,\cdots\,r)}}.
\]
Rather, for the existence of such infinitesimal transformations, the
existence of the equations~\thetag{ 16} is necessary; but it also is
sufficient, as we will see later.

\linestop

\plainstatement{Unless 
the contrary is specially notified, in all the subsequent
studies, we shall restrict ourselves to the $r$-term groups which
contain $r$ independent infinitesimal transformations:
\[
X_kf
=
\sum_{i=1}^n\,\xi_{ki}(x_1,\dots,x_n)\,
\frac{\partial f}{\partial x_i},
\]
and so, which are generated by the same transformations. Here, we
shall always consider only systems of values $x_1, \dots, x_n$ for
which all the $\xi_{ ki}$ behave regularly.

Further, we stress once more that in the future, we shall often call
shortly an $r$-term group with the independent infinitesimal
transformations $X_1f, \dots, X_r f$ by ``\terminology{the group
$X_1f, \dots, X_r f$}''. Amongst the various forms the finite equations
of the $r$-term group $X_1f, \dots, X_r f$ can be given, we shall call
the following:
\[
x_i'
=
x_i
+
\sum_{k=1}^r\,e_k\,\xi_{ki}(x)
+\cdots
\ \ \ \ \ \ \ \ \ \ \ \ \
{\scriptstyle{(i\,=\,1\,\cdots\,n)}}
\]
by ``a \terminology{canonical form}'' 
of the group.}

\linestop


\chapter{
Systems of Partial Differential Equations
\\
the General Solution of Which
\\
Depends Only Upon a Finite Number
\\
of Arbitrary Constants
}
\label{kapitel-10}
\chaptermark{On Certain Higher-Order
Systems of Partial Differential Equations}

\setcounter{footnote}{0}

\abstract*{??}

In the variables $x_1, \dots, x_n$, $z_1, \dots, z_m$ 
let us imagine that a system of partial differential 
equations of arbitrary order is given. Aside from 
$x_1, \dots, x_n$, $z_1, \dots, z_m$, such a system
can only contain differential quotients of 
$z_1, \dots, z_m$ with respect to $x_1, \dots, x_n$, 
so that we have hence to consider $x_1, \dots, x_n$
as independent of each other, while $z_1, \dots, z_m$ 
are to be determined as functions of $x_1, \dots, x_n$
in such a way that the system be identically satisfied. 

Our system of differential equations shall not at all
be completely arbitrary, but it will possess certain 
special properties. We want to assume that, in the form
in which it is presented, it satisfies the following 
conditions.
\label{S-171}

\smallskip{\sf Firstly.}
If $s$ is the order of the highest differential quotient occurring in
the system, then by {\em resolution} of the equations of the system,
all $s$-th order differential quotients of $z_1, \dots, z_m$ with
respect to $x_1, \dots, x_n$ are supposed to 
be expressible in terms of the
differential quotients of the first order up to the $(s-1)$-th, and in
terms of $z_1, \dots, z_m$, $x_1, \dots, x_n$. By contrast, it shall
not be possible in this way to express all differential quotients of
order $(s-1)$ in terms of those of lower order and in terms of $z_1,
\dots, z_m$, $x_1, \dots, x_n$. 

\smallskip{\sf Secondly.}
By differentiating once the given system with respect to the
individual variables $x_1, \dots, x_n$ and by combining the obtained
equations, it shall result only relations between $x_1, \dots, x_n$,
$z_1, \dots, z_m$ and the differential quotients of order $(s-1)$
which already follow from the given system.

\smallskip
We make these special assumptions about the form of the given system
for reasons of convenience. Naturally
\deutsch{Selbstverständlich},
all subsequent considerations can be applied on the whole to every
system of partial differential equations which can, trough
differentiations and elimination, be given the form just described.

According to the known theory of differential equations, it results
without difficulty that every system of differential equations which
possesses the properties discussed just now is integrable and that the
most general functions $z_1, \dots, z_m$ of $x_1, \dots, x_n$ which
satisfy the system depend only on a finite number of constants.

We want now to give somewhat different reasons to this proposition, by
reducing the discussed problem of integration to the problem of
finding systems of equations which admit a certain number of given
infinitesimal transformations; on the basis of the developments of
Chap.~\ref{kapitel-7}, we can indeed solve the latter problem
straightaway.

On the other hand, we will also show that the mentioned proposition
can be inverted: we will show that every system of equations:
\[
z_\mu
=
Z_\mu(x_1,\dots,x_n,\,a_1,\dots,a_r)
\ \ \ \ \ \ \ \ \ \ \ \ \
{\scriptstyle{(\mu\,=\,1\,\cdots\,m)}}
\]
which contains a finite number $r$ of arbitrary parameters $a_k$
represents the most general system of solutions for a certain system
of partial differential equations.

For us, this second proposition is the more important;
we shall use it again in the next chapter.

Therefore, it appears to be desirable to derive the two propositions
in an independent way, without marching into the theory of total
differential equations.

\sectionengellie{\S\,\,\,47.}

Let the number of differential quotients of order $k$ of $z_1, \dots,
z_m$ with respect to $x_1, \dots, x_n$ be denoted $\varepsilon_k$. For
the differential quotients of order $k$ themselves, we introduce the
notation $p_i^{ (k)}$, where $i$ has to run through the values $1, 2,
\dots, \varepsilon_k$; but we reserve us the right to write simply
$z_1, \dots, z_m$ instead of $p_1^{ (0)}$, \dots, $p_{
\varepsilon_0}^{(0)}$. Lastly, we also set:
\[
\frac{\partial p_i^{(k)}}{\partial x_j}
=
p_{ij}^{(k)},
\]
so that $p_{ ij}^{ (k)}$ hence denotes one of the
$\varepsilon_{ k+1}$ differential quotients
of order $(k+1)$. 

According to these settlements, we can write as 
follows:
\def\theequation{1}\begin{equation}
\left\{\aligned
&
W_1(x,\,z,\,p^{(1)},\,\dots,\,p^{(s-1)})
=
0,
\,\,\,\,\dots,\,\,\,\,
W_q(x,\,z,\,p^{(1)},\,\dots,\,p^{(s-1)})
=
0
\\
&
p_{ij}^{(s-1)}
=
P_{ij}(x,\,z,\,p^{(1)},\,\dots,\,p^{(s-1)})
\ \ \ \ \ \ \ \ \ \ \ \ \
{\scriptstyle{(i\,=\,1\,\cdots\,\varepsilon_{s\,-\,1}\,;\,\,\,
j\,=\,1\,\cdots\,n)}}
\endaligned\right.
\end{equation}
the system of differential equations to be studied.

We assume here that the equations $W_1 = 0$, \dots, $W_q = 0$ are
independent of each other. The $n \varepsilon_{ s-1}$ equations $p_{
ij}^{ (s-1)} = P_{ ij}$ are independent of the $W = 0$, but they are
not of each other, because indeed the $n \varepsilon_{ s-1}$
expressions $p_{ ij}^{ (s-1)}$ do not represent only distinct
differential quotients of order $s$. However for what follows the
above way of writing is more convenient than if we had written the
system of equations~\thetag{ 1} under the form:
\[
W_1=0,\,\,\,\dots,\,\,\,
W_q=0,\ \ \ \ \
p_i^{(s)}
=
P_i(x,\,z,\,p^{(1)},\,\dots,\,p^{(s-1)})
\ \ \ \ \ \ \ \ \ \ \ \ \
{\scriptstyle{(i\,=\,1\,\cdots\,\varepsilon_s)}}.
\]

Now, by hypothesis, our system of equations~\thetag{ 1} has the
property that by differentiating it once with respect to the $x$, no
new relation between the $x$, $z$, $p^{ (1)}$, \dots, $p^{ (s-1)}$ is
produced. All relations between the $x$, $z$, $p^{ (1)}$, \dots, $p^{
(s-1)}$ which come out by differentiating once~\thetag{ 1} must hence
be a consequence of $W_1 = 0$, \dots, $W_q = 0$.

Obviously, we find the relations in question by
differentiating~\thetag{ 1} with respect to the $n$ variables $x_1,
\dots, x_n$ and afterwards, by taking away all differential quotients
of order $(s+1)$, and by substituting all differential quotients of
order $s$ by means of~\thetag{ 1}.

By differentiation of $W_k = 0$ and then by 
elimination of the differential quotients of
order $s$, we receive the equations:
\[
\aligned
\frac{\partial W_k}{\partial x_\nu}
&
+
\sum_{i=1}^m\,p_{i\nu}^{(0)}\,
\frac{\partial W_k}{\partial z_i}
+
\sum_{i=1}^{\varepsilon_1}\,p_{i\nu}^{(1)}\,
\frac{\partial W_k}{\partial p_i^{(1)}}
+\cdots+
\\
&
+
\sum_{i=1}^{\varepsilon_{s-2}}\,
p_{i\nu}^{(s-2)}\,
\frac{\partial W_k}{\partial p_i^{(s-2)}}
+
\sum_{i=1}^{\varepsilon_{s-1}}\,
P_{i\nu}\,
\frac{\partial W_k}{\partial p_i^{(s-1)}}
=
0
\\
&
\ \ \ \ \ \ \ \ \ \ \ \ \
{\scriptstyle{(k\,=\,1\,\cdots\,q\,;\,\,\,
\nu\,=\,1\,\cdots\,n)}}.
\endaligned
\]
By the above, these equations are a consequence of $W_1 = 0$, \dots,
$W_q = 0$. In other words: \emphasis{the system of equations $W_1 =
0$, \dots, $W_q = 0$ in the variables $x$, $z$, $p^{ (1)}$, \dots,
$p^{ (s-1)}$ admits the $n$ infinitesimal transformations:
\def\theequation{2}\begin{equation}
\left\{
\aligned
\Omega_\nu f
=
\frac{\partial f}{\partial x_\nu}
&
+
\sum_{i=1}^m\,p_{i\nu}^{(0)}\,
\frac{\partial f}{\partial z_i}
+\cdots+
\\
&
+
\sum_{i=1}^{\varepsilon_{s-2}}\,
p_{i\nu}^{(s-2)}\,
\frac{\partial f}{\partial p_i^{(s-2)}}
+
\sum_{i=1}^{\varepsilon_{s-1}}\,
P_{i\nu}\,
\frac{\partial f}{\partial p_i^{(s-1)}}
\\
\ \ \ \ \ \ \ \ \ \ \ \ \
{\scriptstyle{(\nu\,=\,1\,\cdots\,n)}}.
\endaligned\right.
\end{equation}
in these variables (cf. Chap.~\ref{kapitel-7}, 
p.~\pageref{vanish-by-means})}.

If on the other hand, one differentiates the equations $p_{ ij}^{
(s-1)} = P_{ ij}$ with respect to $x_\nu$ and then eliminates all
differential quotients of order $s$, then one obtains:
\[
\frac{\partial}{\partial x_\nu}\,
p_{ij}^{(s-1)}
=
\frac{\partial^2}{\partial x_\nu\partial x_j}\,
p_i^{(s-1)}
=
\Omega_\nu(P_{ij}). 
\]
One has still to take away all the differential quotients of order
$(s+1)$ from these equations. One easily realizes that only the
following equations come out:
\def\theequation{3}\begin{equation}
\aligned
&
\Omega_\nu(P_{ij})-\Omega_j(P_{i\nu})
=
0
\\
& \ \
{\scriptstyle{(\nu,\,j\,=\,1\,\cdots\,n\,;\,\,\,
i\,=\,1\,\cdots\,\varepsilon_{s\,-\,1})}},
\endaligned
\end{equation}
which likewise must hence be a consequence of
$W_1 = 0$, \dots, $W_q = 0$. 

With this, the properties of the system~\thetag{ 1} demanded
in the introduction are formulated analytically. 

Now, we imagine that an arbitrary system of
solutions:
\[
z_1
=
\varphi_1(x_1,\dots,x_n),\,\,\,\dots,\,\,\,
z_m
=
\varphi_m(x_1,\dots,x_n)
\]
of the differential equations~\thetag{ 1} is presented. 
By differentiating this system we obtain that
the $p^{ (1)}$, \dots, $p^{ (s-1)}$, 
$p^{ (s)}$ are represented as functions of
$x_1, \dots, x_n$:
\[
\aligned
p_{i_1}^{(1)}
=
\varphi_{i_1}^{(1)}(x_1,\dots,x_n),
&
\,\,\,\dots,\,\,\,
p_{i_{s-1}}^{(s-1)}
=
\varphi_{i_{s-1}}^{(s-1)}(x_1,\dots,x_n),\ \ \
p_{i_{s-1},\nu}^{(s-1)}
=
\frac{\partial\varphi_{i_{s-1}}^{(s-1)}}{\partial x_\nu}
\\
&
{\scriptstyle{(i_k\,=\,1,\,2\,\cdots\,\varepsilon_k\,;\,\,\,
\nu\,=\,1\,\cdots\,n)}},
\endaligned
\]
and when we insert these expressions, and
the expressions for $z_1, \dots, z_m$ as well, 
inside the equations~\thetag{ 1}, we
naturally receive nothing but identities. 
From this, it results that
the equations $W_1 = 0$, \dots, $W_q = 0$ in their
turn convert into identities
after the substitution:
\def\theequation{4}\begin{equation}
\aligned
z_\mu
=
\varphi_\mu(x_1,\dots,x_n),\ \ \
&
p_{i_1}^{(1)}
=
\varphi_{i_1}^{(1)}(x_1,\dots,x_n),\,\,\,\dots,\,\,\,
p_{i_{s-1}}^{(s-1)}
=
\varphi_{i_{s-1}}^{(s-1)}(x_1,\dots,x_n)
\\
&\ \ \ \ 
{\scriptstyle{(\mu\,=\,1\,\cdots\,m\,;\,\,\,
i_k\,=\,1\,\cdots\,\varepsilon_k)}}\,;
\endaligned
\end{equation}
clearly, we can also express this as follows: the system of
equations~\thetag{ 4} embraces \deutsch{umfasst} the equations $W_1 =
0$, \dots, $W_q = 0$.

Furthermore, we claim that the system of equations~\thetag{ 4} admits
the infinitesimal transformations $\Omega_1 f, \dots, \Omega_nf$
discussed above.

Indeed at first, all the expressions:
\[
\aligned
\Omega_\nu(z_\mu-\varphi_\mu)
&
=
p_{\mu\nu}^{(0)}
-
\frac{\partial\varphi_\mu}{\partial x_\nu},
\\
\Omega_\nu(p_{i_1}^{(1)}-\varphi_{i_1}^{(1)})
&
=
p_{i,\nu}^{(1)}
-
\frac{\partial\varphi_{i_1}^{(1)}}{\partial x_\nu},
\\
\cdots\cdots\cdots\cdots\cdots\cdots
&
\cdots\cdots\cdots\cdots\cdots\cdots
\\
\Omega_\nu\big(
p_{i_{s-2}}^{(s-2)}-\varphi_{i_{s-2}}^{(s-2)}
\big)
&
=
p_{i_{s-2}\nu}^{(s-2)}
-
\frac{\partial\varphi_{i_{s-2}}^{(s-2)}}{
\partial x_\nu}
\endaligned
\]
vanish by means of~\thetag{ 4}, since the equations~\thetag{ 4} come
from $z_\mu - \varphi_\mu = 0$ by differentiation with respect to
$x_1, \dots, x_n$.
But also the expressions: 
\[
\Omega_\nu\big(
p_{i_{s-1}}^{(s-1)}-\varphi_{i_{s-1}}^{(s-1)}\big)
=
P_{i_{s-1}\nu}
-
\frac{\partial\varphi_{i_{s-1}}^{(s-1)}}{
\partial x_\nu}
\]
vanish by means of~\thetag{ 4}, since the equations:
\[
p_{i_{s-1}\nu}^{(s-1)}
=
P_{i_{s-1}\nu}
\]
are, as already said above, identically 
satisfied after the substitution:
\[
z_\mu
=
\varphi_\mu,\ \ \
p_{i_1}^{(1)}=\varphi_{i_1}^{(1)},
\,\,\,\dots,\,\,\,
p_{i_{s-1}}^{(s-1)}
=
\varphi_{i_{s-1}}^{(s-1)},\ \ \
p_{i_{s-1}\nu}^{(s-1)}
=
\frac{\partial\varphi_{i_{s-1}}^{(s-1)}}{
\partial x_\nu}.
\]

As a result, our claim stated above is proved.

Conversely, imagine now that we are given a system of equations of the
form~\thetag{ 4} about which moreover we do know neither whether it
admits the infinitesimal transformations $\Omega_1 f, \dots,
\Omega_nf$, nor whether it embraces the equations $W_1 = 0$, \dots,
$W_q = 0$.

Clearly, the equations~\thetag{ 4} in question are all created by
differentiating the equations $z_\mu = \varphi_\mu$ with respect to
$x_1, \dots, x_n$. In addition, since under the made assumptions, the
equations:
\[
\Omega_\nu
\big(
p_{i_{s-1}}^{(s-1)}
-
\varphi_{i_{s-1}}^{(s-1)}
\big)
=
P_{i_{s-1}\nu}
-
\frac{\partial\varphi_{i_{s-1}}^{(s-1)}}{
\partial x_\nu}
=
0
\]
are a consequence of~\thetag{ 4} and likewise, of the equations $W_1 =
0$, \dots, $W_q = 0$, then the system of equations~\thetag{ 1} will be
identically satisfied when one executes the substitution~\thetag{ 4} in
it, and then sets also $p_{ i\nu}^{ (s-1)} = \frac{ \partial
\varphi_i^{ (s-1)}}{\partial x_\nu}$. Consequently, the equations
$z_\mu = \varphi_\mu$ represent a system of solutions for the
differential equations~\thetag{ 1}.

From this, we see: every solution $z_\mu = \varphi_\mu$ of the
differential equations~\thetag{ 1} provides a completely determined
system of equations of the form~\thetag{ 4} which admits the
infinitesimal transformations $\Omega_1 f, \dots, \Omega_n f$ and
which embraces the equations $W_1 = 0$, \dots, $W_q = 0$; conversely,
every system of equations of the form~\thetag{ 4} which possesses the
properties just indicated provides a completely determined system of
solutions for the differential equations~\thetag{ 1}.
\emphasis{Consequently, the problem of determining all systems of
solutions of the differential equations~\thetag{ 1} is equivalent to
the problem of determining all systems of equations of the
form~\thetag{ 4} which admit the infinitesimal transformations
$\Omega_1 f, \dots, \Omega_nf$ and which in addition embrace the
equations $W_1 = 0$, \dots, $W_q = 0$. If one knows the most general
solution for one of these two problems, then at the same time, the
most general solution of the other system is given.}

But we can solve this new problem on the basis of the
developments of Chap.~\ref{kapitel-7}, p.~\pageref{S-118-sq} sq.

According to Chap.~\ref{kapitel-7}, Proposition~5,
p.~\ref{Satz-5-S-118}, every sought system of equations admits,
simultaneously with $\Omega_1 f, \dots, \Omega_nf$, also all
infinitesimal transformations of the form $\Omega_\nu \big( \Omega_j (
f) \big) - \Omega_j \big( \Omega_\nu ( f) \big)$. By computation, one
verifies that:
\[
\Omega_\nu\big(\Omega_j(f)\big)
-
\Omega_j\big(\Omega_\nu(f)\big)
=
\sum_{i=1}^{\varepsilon_s-1}\,
\big\{
\Omega_\nu(P_{ij})-\Omega_j(P_{i\nu})
\big\}\,
\frac{\partial f}{\partial p_i^{(s-1)}}
\ \ \ \ \ \ \ \ \ \ \ \ \
{\scriptstyle{(\nu,\,\,j\,=\,1\,\cdots\,n)}},
\]
since the expressions:
\[
\Omega_\nu(p_{ij}^{(k)})
-
\Omega_j(p_{i\nu}^{(k)})
\]
all vanish identically, as long as $k$ is smaller than $s - 1$. But
now, we have seen above that the equations~\thetag{ 3}:
\[
\Omega_\nu(P_{ij})
-
\Omega_j(P_{i\nu})
=
0
\]
are a consequence of $W_1 = 0$, \dots, $W_q = 0$. Consequently, for
the systems of values $x$, $z$, $p^{(1)}$, \dots, $p^{ (s-1)}$ of the
system of equations $W_1 = 0$, \dots, $W_q = 0$, there exist relations
of the form:
\[
\Omega_\nu\big(\Omega_j(f)\big)
-
\Omega_j\big(\Omega_\nu(f)\big)
=
\sum_{\tau=1}^n\,
\omega_{\nu j\tau}
\big(
x,\,z,\,p^{(1)},\dots,p^{(s-1)}
\big)\cdot
\Omega_\tau f
\ \ \ \ \ \ \ \ \ \ \ \ \
{\scriptstyle{(\nu,\,\,j\,=\,1\,\cdots\,n)}},
\]
where the functions $\omega_{ \nu j \tau}$ are all
equal to zero and behave regularly for the
systems of values of $W_1 = 0$, \dots, $W_q = 0$. 

In addition, it is still to be underlined that amongst
the $n \times n$ determinants of the matrix:
\[
\left\vert
\begin{array}{ccccccccccc}
1 & 0 & \cdots & 0 & p_{11}^{(0)} & \cdots
& p_{\varepsilon_{s-2},1}^{(s-2)} & P_{11} 
& \cdots & P_{\varepsilon_s-1,1}
\\
\cdot & \cdot & \cdots & \cdot & \cdots & \cdot & \cdot
& \cdot & \cdots & \cdot
\\
0 & 0 & \cdots & 1 & p_{1n}^{(0)} & \cdots
& p_{\varepsilon_{s-2},n}^{(s-2)} & P_{1n} 
& \cdots & P_{\varepsilon_{s-1},n}
\end{array}
\right\vert,
\]
one has the value $1$, and thus cannot be brought to zero for any of
the sought systems of equations.

From this, we see that we have in front of us the special case whose
handling was given by Theorem~19 in Chap.~\ref{kapitel-7},
p.~\pageref{Theorem-19-S-132}. With the help of this theorem, we can
actually set up all systems of equations which admit $\Omega_1 f,
\dots, \Omega_n f$ and which embrace the equations $W_1 = 0$, \dots,
$W_q = 0$. There is no difficulty to identify such systems of
equations which can be given the form~\thetag{ 4}.

From one system of equations which admits the infinitesimal
transformations $\Omega_1 f, \dots, \Omega_nf$, one can never derive a
relation between the variables $x_1, \dots, x_n$. This follows from
the mentioned theorem and also can easily be seen directly. As a
result, it is possible to solve the equations $W_1 = 0$, \dots, $W_q =
0$ with respect to $q$, amongst $\varepsilon_0 + \varepsilon_1 +
\cdots + \varepsilon_{ s-1}$, of the quantities $z$, $p^{ (1)}$,
\dots, $p^{ (s-1)}$. When we do that, we receive $q$ of the $\sum
\varepsilon_k$ variables $z$, $p^{ (1)}$, \dots, $p^{ ( s-1)}$
expressed by means of the $\sum \varepsilon_k - q$ remaining ones and
by means of $x_1, \dots, x_n$.

Thanks to these preparations, we can form the reduced
infinitesimal transformations about which it is question
in the mentioned theorem. We obtain them when
we leave out all differential quotients 
of $f$ with respect to the $q$ considered variables
amongst $z$, $p^{ (1)}$, \dots, $p^{ (s-1)}$
and afterwards, by substituting, in the remaining
terms, each one of the $q$ variables with its expression
by means of the $\sum \varepsilon_k - q$ left ones and 
the $x$. 

The so obtained $n$ reduced infinitesimal transformations, 
which we can denote by $\overline{ \Omega}_1 f, \dots, 
\overline{ \Omega}_n f$, contain $n - q + \sum \, \varepsilon_k$
independent variables and stand in addition 
pairwise in the relationships:
\[
\overline{\Omega}_\nu\big(
\overline{\Omega}_j(f)\big)
-
\overline{\Omega}_j\big(
\overline{\Omega}_\nu(f)\big)
\equiv
0
\ \ \ \ \ \ \ \ \ \ \ \ \
{\scriptstyle{(\nu,\,\,j\,=\,1\,\cdots\,n)}},
\]
according to the mentioned theorem. 

Consequently, the $n$ mutually independent equations
$\overline{ \Omega}_1 f = 0$, \dots, $\overline{ 
\Omega}_n f = 0$ form an $n$-term complete
system with the $\sum \, \varepsilon_k - q$
independent solutions: 
$u_1$, $u_2$, \dots, $u_{\, 
\sum \varepsilon_k - q}$. 
If these solutions are determined, then one
can indicate all systems of equations which 
admit the infinitesimal transformations
$\Omega_1 f, \dots, \Omega_n f$ and
which embrace the equations $W_1 = 0, \dots, 
W_q = 0$. The general form
of a system of this kind is: 
$W_1 = 0$, \dots, $W_q = 0$ together
with arbitrary relations between the $\sum \, \varepsilon_k - q$
solutions $u$. 

However now, the matter is not about all systems of equations
of this kind, but only about those which can 
be brought to the form~\thetag{ 4}. 
Every system of equations of this constitution
contains exactly $\sum \, \varepsilon_k$ equations, 
hence can be given the form:
\def\theequation{5}\begin{equation}
W_1=0,\,\,\,\dots,\,\,\,\,W_q=0,\ \ \ \ \
u_1=a_1,\,\,\,\dots,\,\,\,
u_{\,\sum\varepsilon_k-q}
=
a_{\,\sum\varepsilon_k-q},
\end{equation}
where the $a$ denotes constants. But every system of the form just
indicated can be solved with respect to the $\sum \, \varepsilon_k$
variables $z$, $p^{ (1)}$, \dots, $p^{ (s-1)}$, since it admits the
infinitesimal transformations $\Omega_1 f$, \dots, $\Omega_nf$. If we
perform the resolution in question, we then obtain a system of
equations:
\def\theequation{5'}\begin{equation}
\aligned
z_\mu
&
=
Z_\mu(x_1,\dots,x_n,\,a_1,a_2,\dots),\ \ \ 
p_{i_1}^{(1)}
=
\Pi_{i_1}^{(1)}
(x_1,\dots,x_n,\,a_1,a_2,\dots),
\\
&
\ \ \ \ \ \ \ \ \ \ \
\cdots\,
p_{i_{s-1}}^{(s-1)}
=
\Pi_{i_{s-1}}(x_1,\dots,x_n,\,a_1,a_2,\dots)
\\
&
\ \ \ \ \ \ \ \ \ \ \ \ \ \ \ \ \ \ \ \ \ \ \ \ \ \ \ \ \ \ \
{\scriptstyle{(\mu\,=\,1\,\cdots\,m\,;\,\,\,
i_k\,=\,1\,\cdots\,\varepsilon_k)}}
\endaligned
\end{equation} 
which satisfies all the stated requirements:
it admits the infinitesimal transformations
$\Omega_1 f, \dots, \Omega_n f$, it
embraces the equations
$W_1 = 0$, \dots, $W_q = 0$, and
it possesses the form~\thetag{ 4}. 
At present, when we consider the $\sum \, \varepsilon_k
-q$ constants $a$ as arbitrary constants, 
we therefore obviously have
the most general system of equations of the
demanded constitution.

From this, it follows that according to what
has been said, the equations:
\def\theequation{6}\begin{equation}
z_\mu
=
Z_\mu(x_1,\dots,x_n,\,a_1,a_2,\dots) 
\ \ \ \ \ \ \ \ \ \
{\scriptstyle{(\mu\,=\,1\,\cdots\,m)}}
\end{equation}
represent a system of solutions for the differential
equations~\thetag{ 1} and in fact, the most general system of
solutions. We now claim that in this system of solutions, the $\sum
\, \varepsilon_k - q$ arbitrary constants $a$ are all essential.

To prove our claim, we remember that the equations~\thetag{ 5'} can be
obtained from the equations $z_\mu = Z_\mu$ by differentiation with
respect to $x_1, \dots, x_n$, provided that the $p^{ (1)}$, $p^{
(2)}$, \dots\, are again interpreted as differential quotients of the
$z$ with respect to the $x$. Now, if the $\sum \, \varepsilon_k - q$
parameters $a$ in the equations~\thetag{ 6} were not essential, then
the number of parameters could be lowered by introducing appropriate
functions of them. But with this, according to what precedes, the
number of parameters in the equations~\thetag{ 5'} would at the same
time be lowered, and this is impossible, since the equations~\thetag{
5'} can be brought to the form~\thetag{ 5}, from which it follows
immediately that the parameters $a$ in~\thetag{ 5'} are all essential.
This is a contradiction, so the assumption made a short while ago is
false and the parameters $a$ in the equations~\thetag{ 6} are all
essential.

We can hence state the following proposition.

\def\theproposition{1}\begin{proposition}
\label{Satz-1-S-179}
If a system of partial differential equations of
order $s$ of the form:
\[
F_\sigma
\bigg(
x_1,\dots,x_n,\,\,z_1,\dots,z_m,\,\,
\frac{\partial z_1}{\partial x_1},\,\dots,\,
\frac{\partial^2z_1}{\partial x_1^2},\,\dots,\,
\frac{\partial^s z_m}{\partial x_n^s}
\bigg)
=
0
\ \ \ \ \ \ \ \ \ \ \ \ \
{\scriptstyle{(\sigma\,=\,1,\,2,\,\cdots)}}
\]
possesses the property that all differential quotients of order $s$ of
the $z$ with respect to the $x$ can be expressed by means of the
differential quotients of lower order, 
by means of $z_1, \dots, z_m$ and by means of
$x_1, \dots, x_n$, while the corresponding property does not hold in
any case for the differential quotients of order $(s-1)$, and if in
addition, by differentiating it once with respect to the $x$, the system
produces only relations between $x_1, \dots, x_n$, $z_1, \dots, z_n$
and the differential quotients of the first order up to the $(s-1)$-th
which follow from the system itself, then the most general system of
solutions:
\[
z_\mu
=
\varphi_\mu(x_1,\dots,x_n)
\ \ \ \ \ \ \ \ \ \ \ \ \
{\scriptstyle{(\mu\,=\,1\,\cdots\,m)}}
\]
of the concerned system of differential equations contains only a
finite number of arbitrary constants. The number of these arbitrary
constants is equal to $\varepsilon_0 + \varepsilon_1 + \cdots +
\varepsilon_{ s-1} - q$, where $\varepsilon_k$ denotes the number of
differential quotients of order $k$ of $z_1, \dots, z_m$ with respect
to $x_1, \dots, x_n$ and where $q$ is the number of independent
relations which the system in question yields between $x_1, \dots,
x_n$, $z_1, \dots, z_n$ and the differential quotients up to order
$(s-1)$. One finds the most general system of solutions $z_\mu =
\varphi_\mu$ itself by integrating an $n$-term complete
system in $n - q + \varepsilon_0 + \cdots + \varepsilon_{ s-1}$
independent variables.
\end{proposition}

For the proof of the above proposition, we imagined the equations $W_1
= 0$, \dots, $W_q = 0$ solved with respect to $q$ of the quantities
$z$, $p^{ (1)}$, \dots, $p^{ (s-1)}$. In principle, it is completely
indifferent with respect to which amongst these quantities the
equations are solved; but now, since the equations $W_1 = 0$, \dots,
$W_q = 0$ are differential equations and since the $p^{ (1)}$, $p^{
(2)}$, \dots\, denote differential quotients, it is advisable to
undertake the concerned resolution in a specific way, which we want to
now explain.

At first, we eliminate all differential quotients $p^{ (1)}$, $p^{
(2)}$, \dots, $p^{ (s-1)}$ from the equations $W_1 = 0$, \dots, $W_q =
0$; then we obtain, say, $\nu_0$ independent equations between the $x$
and $z$ alone, and so we can represent $\nu_0$ of the $z$ as functions
of the $\varepsilon_0 - \nu_0 = m - \nu_0$ remaining ones, and of the
$x$.

We insert the expressions for these $\nu_0$ quantities $z$ in the
equations $W_1 = 0$, \dots, $W_q = 0$, which now reduce to $q - \nu_0$
mutually independent equations. Afterwards, from these $q - \nu_0$
equations, we remove all the differential quotients $p^{ (2)}$, \dots,
$p^{ (s-1)}$ and we obtain, say $\nu_1$ mutually independent equations
by means of which we can express $\nu_1$ of the quantities $p^{ (1)}$
in terms of the $\varepsilon_1 - \nu_1$ remaining ones, in terms of
the $\varepsilon_0 - \nu_0$ of the $z$ and in terms of $x_1, \dots,
x_n$.

If we continue in the described way, then at the end, the system of
equations $W_1 = 0$, \dots, $W_q = 0$ will be resolved, 
\label{S-180} and to be
precise, it will be resolved with respect to $\nu_k$ amongst the
$\varepsilon_k$ differential quotients $p^{ (k)}$, where it is
understood that $k$ is any of the numbers $0$, $1$, $2$, \dots, $s-1$.
Here, the concerned $\nu_k$ amongst the $p^{ (k)}$
are each time expressed in terms of the 
$\varepsilon_k - \nu_k$ remaining ones, in 
terms of certain of the differential quotients:
$p^{ (k-1)}$, \dots, $p^{ (1)}$, $p^{ (0)}$, and
in terms of the $x$. Naturally, 
the sum $\nu_0 + \nu_1 + \cdots + 
\nu_{ s-1}$ has the value $q$. 
Lastly, it can be shown that $\nu_k$ is always smaller than
$\varepsilon_k$. Indeed at first, 
$\nu_{ s-1}$ is certainly smaller than $\varepsilon_{ s-1}$, 
because we have assumed that in our system
of differential equations, not
all differential quotients of order $(s-1)$
can be expressed in terms of those of lower order, and
in terms of the $x$. 
But if any other numbers $\nu_k$ would be equal
to $\varepsilon_k$, we would have:
\[
p_i^{(k)}
=
F_i(x,\,z,\,p^{(1)},\,\dots,\,p^{(k-1)})
\ \ \ \ \ \ \ \ \ \ \ \ \
{\scriptstyle{(i\,=\,1,\,2\,\cdots\,\varepsilon_k)}}\,;
\]
now, since the system of equations
$W_1 = 0$, \dots, $W_q = 0$ admits the infinitesimal
transformations $\Omega_1 f, \dots, \Omega_n f$, then
all systems of equations:
\[
p_{i\nu}^{(k)}
=
\Omega_\nu(F_i)
\ \ \ \ \ \ \ \ \ \ \ \ \
{\scriptstyle{(i\,=\,1\,\cdots\,\varepsilon_k\,;\,\,\,
\nu\,=\,1\,\cdots\,n)}}
\]
would be consequences of $W_1 = 0$, \dots, $W_q = 0$, so one would
also have $\nu_{ k+1} = \varepsilon_{ k+1}$ and in the same way $\nu_{
k+2} = \varepsilon_{ k+2}$, \dots, $\nu_{ s-1} = \varepsilon_{ s-1}$,
and this is impossible according to what has been said.

Once we have resolved the equations $W_1 = 0$, \dots, $W_q = 0$, then
with this at the same time, the system of differential
equations~\thetag{ 1} is resolved in a completely determined way,
namely with respect to $\nu_k$ of the $\varepsilon_k$ differential
quotients of order $k$, where it is understood that $k$ is any of the
numbers $0$, $1$, $2$, \dots, $s$. At each time, the $\nu_k$
differential quotients of order $k$ are expressed in terms of the
$\varepsilon_k - \nu_k$ remaining ones of order $k$, in terms of those
of lower order, and in terms of the $x$.

If one knows all the numbers $\varepsilon_k$ and $\nu_k$, then
one can immediately indicate the number
of arbitrary constants in the most general 
system of solutions of the differential 
equations~\thetag{ 1}. Indeed, 
according to the above proposition, 
this number is equal to:
\[
\label{S-181}
\varepsilon_0+\varepsilon_1
+\cdots+
\varepsilon_{s-1}-q
=
(\varepsilon_0-\nu_0)+(\varepsilon_1-\nu_1)
+\cdots+
(\varepsilon_{s-1}-\nu_{s-1}).
\]

\sectionengellie{\S\,\,\,48.}

Conversely, consider a system of equations of the form:
\def\theequation{7}\begin{equation}
\aligned
z_\mu
&
=
Z_\mu(x_1,\dots,x_n,\,a_1,\dots,a_r)
\\
&
\ \ \ \ \ \ \ \ \ \ \ \ \
{\scriptstyle{(\mu\,=\,1\,\cdots\,m)}},
\endaligned
\end{equation}
in which $x_1, \dots, x_n$ are interpreted as
independent variables, and $a_1, \dots, a_r$ as
arbitrary parameters. 
The $r$ parameters $a_1, \dots, a_r$ whose number
is finite can be assumed to be all essential.

We will prove that there is a system
of partial differential equations which is free of
$a_1, \dots, a_r$ whose most general solution is
represented by the equations~\thetag{ 7}. 

By differentiating the equations~\thetag{ 7} with 
respect to $x_1, \dots, x_n$, we obtain the ones
after the others equations of the form:
\def\theequation{$7_1$}\begin{equation}
p_i^{(1)}
=
Z_i^{(1)}(x_1,\dots,x_n,\,a_1,\dots,a_r)
\ \ \ \ \ \ \ \ \ \ \ \ \
{\scriptstyle{(i\,=\,1\,\cdots\,\varepsilon_1)}}
\end{equation}
\def\theequation{$7_2$}\begin{equation}
p_i^{(2)}
=
Z_i^{(2)}(x_1,\dots,x_n,\,a_1,\dots,a_r)
\ \ \ \ \ \ \ \ \ \ \ \ \
{\scriptstyle{(i\,=\,1\,\cdots\,\varepsilon_2)}},
\end{equation}
and so on. 

Now, by means of the equations~\thetag{ 7}, a certain number, say
$\nu_0$, of the $a$ can be expressed in terms of the $r - \mu_0$
remaining ones, and in terms of the $x$ and the $z$. By taking
together the equations~\thetag{ 7} and~\thetag{ $7_1$}, we can express
$\mu_0 + \mu_1$ of the $a$ in terms of the $r - \mu_0 - \mu_1$
remaining ones, in terms of the $p^{ (1)}$, the $z$ and the $x$. In
general, if we take the equations~\thetag{ 7}, \thetag{ $7_1$}, \dots,
\thetag{ $7_k$} together, we can then express $\mu_0 + \mu_1 + \cdots
+ \mu_k$ of the quantities $a$ in terms of the $r - \mu_0 - \cdots -
\mu_k$ remaining ones and in terms of the $p^{ (k)}$, $p^{ (k-1)}$,
\dots, $p^{ (1)}$, $z$, $x$. Here, $\mu_0$, $\mu_1$, \dots\, are
entirely determined positive entire numbers.

Naturally, the sum $\mu_0 + \mu_1 + \cdots + \mu_k$ is at most equal
to $r$; furthermore, the number $\mu_0$ is in any case different from
zero, because we can assume that $r$ is bigger than
zero. Consequently, there must exist a finite entire number $s
\geqslant 1$ having the property that $\mu_0$, $\mu_1$, \dots, $\mu_{
s-1}$ are all different from zero, but the number $\mu_s$ vanishes.
Then if the said quantities amongst $a_1, \dots, a_r$ are determined
from the equations~\thetag{ 7}, \thetag{ $7_1$}, \dots, \thetag{ $7_{
s-1}$} and if their values are inserted in the equations~\thetag{
$7_s$}, then one obtains only equations which are free of $a_1, \dots,
a_r$; because on the contrary case, one could determine more than
$\mu_0 + \cdots + \mu_{ s-1}$ of the quantities $a$ from the
equations~\thetag{ 7}, \thetag{ $7_1$}, \dots, \thetag{ $7_s$}, so one
would have $\mu_s > 0$, which is not the case according to the above.

As a result, by eliminating $a_1, \dots, a_r$ from the 
equations~\thetag{ 7}, \thetag{ $7_1$}, \dots, \thetag{ $7_s$},
we obtain the following equations:

\smallskip{\sf Firstly}, $\varepsilon_s$ equations which express all
the $\varepsilon_s$ differential quotients $p^{ (s)}$ in terms of the
$p^{ (s-1)}$, \dots, $p^{ (1)}$, $z$, $x$, and:

\smallskip{\sf Secondly:} yet $\varepsilon_0 - \mu_0 + \varepsilon_1 -
\mu_1 + \cdots + \varepsilon_{ s-1} - \mu_{ s-1}$ mutually independent
equations between the $x$, $z$, $p^{ (1)}$, \dots, $p^{ (s-1)}$.

\smallskip
It is easy to see that the so obtained system of
differential equations of order $s$ possesses all 
properties which were ascribed
in Proposition~1, p.~\pageref{Satz-1-S-179} 
to the differential equations $F_\sigma \big(
x, \, z, \, \frac{ \partial z}{ \partial x}, \dots \big)$. 

Indeed, all $\varepsilon_s$ differential quotients of order
$s$ are functions of the differential quotients of
lower order, and of the $x$; 
however, the corresponding property does not
hold true for the differential quotients of
order $(s-1)$, because the number $\mu_{ s-1}$
discussed above is indeed different from zero. 
Lastly, by differentiation with respect to the $x$
and by combination of the obtained equations,
it results only relations between the $x$, $z$, 
$p^{ (1)}$, \dots, $p^{ (s-1)}$
which follow from the above-mentioned relations. Indeed, 
the equations~\thetag{ 7}, 
\thetag{ $7_1$}, \dots, \thetag{ $7_{ s-1}$}
and the ones which follow from them are
the only finite relations through which
the quantities $a_1, \dots, a_r$, $x$, $z$, 
$p^{ (1)}$, \dots, $p^{ (s-1)}$ are
linked; hence when we eliminate the $a$, we obtain 
the only finite relations which exist between 
the $x$, $z$, $p^{ (1)}$, \dots, 
$p^{ (s-1)}$, namely the $\varepsilon_0 - \mu_0 + \cdots + 
\varepsilon_{ s-1} - \mu_{ s-1}$ relations mentioned
above. 

The Proposition~1 p.~\pageref{Satz-1-S-179}
can therefore easily be applied to our
system of differential equations of
order $s$. The numbers $\nu_k$ defined at that time are equal
to $\varepsilon_k - \mu_k$, and the number $q$ has hence
in the present case the value:
\[
\varepsilon_0-\mu_0
+\cdots+
\varepsilon_{s-1}-\mu_{s-1},
\]
and therefore the most general system of solutions of
our differential equations contains precisely 
$\mu_0 + \mu_1 + \cdots + \mu_{ s-1}$
arbitrary constants. Now, since on the other hand
the equations~\thetag{ 7} also represent a
system of solutions for our differential equations, and
a system with the $r$ essential parameters
$a_1, \dots, a_r$ as arbitrary constants, 
it follows that the number $\mu_0 + \cdots + 
\mu_{ s-1}$, which is at most equal to $r$, 
must precisely be equal to $r$. In other
words: the equations~\thetag{ 7} represent
the most general system of solutions of our differential equations.

As a result, we have the

\def\theproposition{2}\begin{proposition}
If $z_1, \dots, z_m$ are given functions of the variables
$x_1, \dots, x_n$ and of a finite number of parameters
$a_1, \dots, a_r$:
\[
z_\mu
=
Z_\mu(x_1,\dots,x_n,\,a_1,\dots,a_r)
\ \ \ \ \ \ \ \ \ \ \ \ \
{\scriptstyle{(\mu\,=\,1\,\cdots\,m)}},
\]
then there always exists an integrable system of partial differential
equations which determines the $z$ as functions of the $x$, and whose
most general solutions are represented by the equations $z_\mu = Z_\mu
( x, \, a)$.
\end{proposition}

If we relate this proposition with the Proposition~1
p.~\pageref{Satz-1-S-179}, we immediately obtain the

\def\theproposition{3}\begin{proposition}
\label{Satz-3-S-183}
If an integrable system of partial differential equations:
\[
F_\sigma
\bigg(
x_1,\dots,x_n,\,z_1,\dots,z_m,\,\,
\frac{\partial z_1}{\partial x_1},\,\dots,\,
\frac{\partial^2z_1}{\partial x_1^2},\,\dots\,
\bigg)
=
0\ \ \ \ \ \ \ \ \ \ \ \ \
{\scriptstyle{(\sigma\,=\,1,\,2\,\cdots)}}
\]
is constituted in such a way that its most general solutions depend
only on a finite number of arbitrary constants, then by means of
differentiation and of elimination, it can always be brought to a form
which possesses the following two properties: firstly, all
differential quotients of a certain order, say $s$, can be expressed
in terms of those of lower order, and in terms of $z_1, \dots, z_m$,
$x_1, \dots, x_n$, whereas the corresponding property does not, in any
case, hold true for all differential quotients of order
$(s-1)$. Secondly, by differentiating once with respect to the $x$, it
only comes relations between the $x$, $z$ and the differential
quotients of orders $1$ up to $(s-1)$ which follow from the already
extant equations. 
\end{proposition}

Besides, the developments of the present $\S$ provide a simple method
to answer the question of how many parameters $a_1, \dots, a_r$ are
essential amongst the ones of a given system of equations:
\[
z_\mu
=
Z_\mu(x_1,\dots,x_n,\,a_1,\dots,a_r)
\ \ \ \ \ \ \ \ \ \ \ \ \
{\scriptstyle{(\mu\,=\,1\,\cdots\,m)}}.
\]

Indeed, in order to be able to answer this question, we only need to
compute the entire numbers $\mu_0$, $\mu_1$, \dots, $\mu_{ s-1}$
defined above; the sum $\mu_0 + \mu_1 + \cdots + \mu_{ s-1}$ then
identifies the number of essential parameters amongst $a_1, \dots,
a_r$, because the equations $z_\mu = Z_\mu ( x, a)$ represent the most
general solutions of a system of differential equations the most
general solutions of which, according to the above developments,
contain precisely $\mu_0 + \cdots + \mu_{ s-1}$ essential parameters.

\linestop


\chapter{The Defining Equations
\\
for the Infinitesimal Transformations of a Group
}
\label{kapitel-11}
\chaptermark{The Defining Equations of a Group}

\setcounter{footnote}{0}

\abstract*{??}

In Chap.~\ref{kapitel-9}, we have reduced the finding of all $r$-term
groups to the determination of all systems of $r$ independent
infinitesimal transformations $X_1 f, \dots, X_r f$ which satisfy
relations of the form:
\[
X_iX_kf
-
X_kX_if
=
\leftbracket
X_i,\,X_k
\rightbracket
=
\sum_{s=1}^r\,c_{iks}\,X_sf, 
\]
with certain constants $c_{ iks}$. Only later will we find means of
treating this reduced problem; temporarily, we must restrict ourselves
to admitting systems $X_1 f, \dots, X_r f$ of the concerned nature and
to study their properties.

\renewcommand{\thefootnote}{\fnsymbol{footnote}}
In the present chapter we begin with an application of the
developments of the preceding chapter; from this chapter, 
we conclude that the general infinitesimal
transformation:
\[
e_1\,X_1f
+\cdots+
e_r\,X_rf
\]
of a given $r$-term group $X_1f, \dots, X_r f$ can be defined by means
of certain linear partial differential equations, which we call the
\terminology{defining equations} of the group. From that, further
conclusions will then be drawn\footnote[1]{\,
\name{Lie}, Archiv for Mathematik og Naturvidenskab Vol.~3, 
Christiania 1878 and Vol.~8, 1883; Gesellschaft der
Wissenschaften zu Christiania 1883, No.~12. 
}. 
\renewcommand{\thefootnote}{\arabic{footnote}}

\sectionengellie{\S\,\,\,49.}

Consider $r$ infinitesimal transformations:
\[
X_kf
=
\sum_{i=1}^n\,\xi_{ki}(x_1,\dots,x_n)\,
\frac{\partial f}{\partial x_i}
\ \ \ \ \ \ \ \ \ \ \ \ \
{\scriptstyle{(k\,=\,1\,\cdots\,r)}}
\]
which generate an $r$-term group. Then the general infinitesimal
transformation of this group has the form:
\[
\sum_{i=1}^n\,\xi_i\,\frac{\partial f}{\partial x_i}
=
\sum_{i=1}^n\,\sum_{k=1}^r\,
e_k\,\xi_{ki}\,
\frac{\partial f}{\partial x_i},
\]
where the $e_k$ denote arbitrary constants. For the
$\xi_i$, it comes from this the expressions:
\def\theequation{1}\begin{equation}
\xi_i
=
\sum_{k=1}^r\,e_k\,\xi_{ki}(x)
\ \ \ \ \ \ \ \ \ \ \ \ \
{\scriptstyle{(i\,=\,1\,\cdots\,n)}},
\end{equation}
in which the $\xi_{ ki} (x)$ are given functions of
$x_1, \dots, x_n$. 

Now, according to Proposition~2 of the preceding chapter, the
expressions $\xi_1, \dots, \xi_n$ just written are the most general
system of solutions of a certain system of partial differential
equations which is free of the arbitrary constants $e_1, \dots, e_r$.
In order to set up this system, we proceed as in the introduction to
\S\,\,48; we differentiate the equations~\thetag{ 1} with respect to
every variable $x_1, \dots, x_n$, next we differentiate in the same
way the obtained equations with respect to $x_1, \dots, x_n$, and so
on. Then when we have computed all differential quotients of the $\xi$
up to a certain order (considered in more details in \S\,\,48), we
take away from the found equations the arbitrary constants $e_1,
\dots, e_r$ and we receive in this way the desired system of
differential equations, the most general solutions of which are just
the expressions $\xi_1, \dots, \xi_n$.

We can always arrange that all equations of the discussed system in
the $\xi_i$ and in their differential quotients are linear homogeneous;
indeed, from any two systems of solutions:
\[
\xi_{k\nu},\ \ \
\xi_{j\nu}
\ \ \ \ \ \ \ \ \ \ \ \ \
{\scriptstyle{(\nu\,=\,1\,\cdots\,n)}},
\]
one can always derive another system of solutions $e_k \xi_{ k\nu} +
e_j \xi_{ j\nu}$ of the concerned differential equations with two
arbitrary constants $e_k$ and $e_j$.

The system of the differential equations which define $\xi_1, \dots,
\xi_n$ therefore has the form:
\[
\sum_{\nu=1}^n\,A_{\mu\nu}(x_1,\dots,x_n)\,\xi_\nu
+
\sum_{\nu,\,\,\pi}^{1\cdots n}\,
B_{\mu\nu\pi}(x_1,\dots,x_n)\,
\frac{\partial\xi_\nu}{\partial x_\pi}
+\cdots
=
0,
\]
where the $A$, $B$, \dots, are free of the arbitrary constants $e_1,
\dots, e_r$.

We briefly \deutsch{kurtzweg} call these differential equations the
\terminology{defining equations of the group}, \emphasis{since they
completely define the totality of all infinitesimal transformations of
this group and therefore, they define the group itself}.

What has been said can be illustrated precisely at this place by means
of a couple of examples.

In the general infinitesimal transformation:
\[
\xi_1\,\frac{\partial f}{\partial x_1}
+
\xi_2\,\frac{\partial f}{\partial x_2}
\]
of the six-term linear group:
\[
\aligned
x_1'
&
=
a_1\,x_1+a_2\,x_2+a_3
\\
x_2'
&
=
a_4\,x_1+a_5\,x_2+a_6,
\endaligned
\]
$\xi_1$ and $\xi_2$ have the form: 
\[
\xi_1
=
e_1+e_2\,x_1+e_3\,x_2,
\ \ \ \ \ \ \ \ \
\xi_2
=
e_4+e_5\,x_1+e_6\,x_2.
\]

From this, it comes that the defining equations of the group are the
following:
\[
\aligned
\frac{\partial^2\xi_1}{\partial x_1^2}
=
0,\ \ \ \ \
\frac{\partial^2\xi_1}{\partial x_1\partial x_2}
=
0,\ \ \ \ \
\frac{\partial^2\xi_1}{\partial x_2^2}
=
0,
\\
\frac{\partial^2\xi_2}{\partial x_1^2}
=
0,\ \ \ \ \
\frac{\partial^2\xi_2}{\partial x_1\partial x_2}
=
0,\ \ \ \ \
\frac{\partial^2\xi_2}{\partial x_2^2}
=
0.
\endaligned
\]

The known eight-term group:
\[
x_1'
=
\frac{a_1\,x_1+a_2\,x_2+a_3}{a_7\,x_1+a_8\,x_2+1},
\ \ \ \ \ \ \ \ \ \
x_2'
=
\frac{a_4\,x_1+a_5\,x_2+a_6}{a_7\,x_1+a_8\,x_2+1}
\]
of all projective transformations of the plane serves as a second
example. Its general infinitesimal transformation:
\[
\aligned
\xi_1
&
=
e_1+e_2\,x_1+e_3\,x_2+e_7\,x_1^2+e_8\,x_1x_2
\\
\xi_2
&
=
e_4+e_5\,x_1+e_6\,x_1+e_7\,x_1x_2+e_8\,x_2^2
\endaligned
\]
will be defined by means of relations between the differential
quotients of second order of the $\xi$, namely by means of:
\[
\frac{\partial^2\xi_1}{\partial x_1^2}
-
2\,\frac{\partial^2\xi_2}{\partial x_1\partial x_2}
=
0,\ \ \ \ \
\frac{\partial^2\xi_2}{\partial x_2^2}
-
2\,\frac{\partial^2\xi_1}{\partial x_1\partial x_2}
=
0,\ \ \ \ \
\frac{\partial^2\xi_1}{\partial x_2^2}
=
0,\ \ \ \ \
\frac{\partial^2\xi_2}{\partial x_1^2}
=
0.
\]
By renewed differentiation, one finds 
that all third-order differential quotients of
$\xi_1$ and of $\xi_2$ vanish. 

\sectionengellie{\S\,\,\,50.}

Conversely, when does a system of linear homogeneous
differential equations:
\[
\sum_{\nu=1}^n\,A_{\mu\nu}(x)\,\xi_\nu
+
\sum_{\nu,\,\,\pi}^{1\cdots n}\,
B_{\mu\nu\pi}(x)\,
\frac{\partial\xi_\nu}{\partial x_\pi}
+\cdots
=
0
\ \ \ \ \ \ \ \ \ \ \ \ \
{\scriptstyle{(\mu\,=\,1,\,2\,\cdots)}}
\]
define the general infinitesimal transformation of a finite continuous
group?

Naturally, the first condition is that the most general solutions
$\xi_1, \dots, \xi_n$ of the system depend only on a finite number of
arbitrary constants. Let this condition be fulfilled. Then according
to Proposition~3, p.~\pageref{Satz-3-S-183} of the preceding
chapter, it is always possible, by differentiation and elimination, to
bring the system to a certain specific form in which all differential
quotients of the highest, say the $s$-th, order can be expressed in
terms of those of lower order, and in terms of $x_1, \dots, x_n$,
whereas the corresponding property does not, in any case, hold true
for all differential quotients of order $(s-1)$;
in addition, differentiating yet once with respect to the $x$ produces
no new relation between $x_1, \dots, x_n$, $\xi_1, \dots, \xi_n$ and
the differential quotients of first order up to the $(s-1)$-th. We
assume that the system has received this form and we imagine that it
has been solved in the way indicated on p.~\pageref{S-180}; then for
$k = 0, 1, \dots, s$, we obtain at each time that a certain number of
differential quotients, say $\nu_k$, amongst the $\varepsilon_k$ of
the $k$-th order differential quotients of the $\xi$, are represented
as linear homogeneous functions of the remaining differential
quotients of $k$-th order, and of certain differential quotients of
$(k-1)$-th, \dots, first, zeroth order, with coefficients which depend
only upon $x_1, \dots, x_n$, where $\nu_s = \varepsilon_s$, but else,
$\nu_k$ is always smaller than $\varepsilon_k$.

As was shown in the preceding chapter, p.~\pageref{S-181}, under the
assumptions made, the most general system of solutions $\xi_1, \dots,
\xi_n$ of our differential equations comprises precisely:
\[
(\varepsilon_0-\nu_0)+(\varepsilon_1-\nu_1)
+\cdots+
(\varepsilon_{s-1}-\nu_{s-1})
=
r
\]
arbitrary constants; this most general system of solutions can be
deduced from $r$ particular systems of solutions $\xi_{ 1i}, \dots,
\xi_{ ri}$ with the help of $r$ constants of integration as follows:
\[
\xi_i
=
\sum_{k=1}^r\,e_k\,\xi_{ki}\ \ \ \ \ \ \ \ \ \ \ \ \
{\scriptstyle{(i\,=\,1\,\cdots\,n)}},
\]
and here, the particular systems of solutions in question must only be
such that the $r$ expressions:
\[
\sum_{i=1}^n\,\xi_{1i}\,\frac{\partial f}{\partial x_i},
\,\,\,\,\dots,\,\,\,
\sum_{i=1}^n\,\xi_{ri}\,\frac{\partial f}{\partial x_i}
\]
represent as many independent infinitesimal transformations.

Now according to Theorem~24 p.~\pageref{Theorem-24-S-158}, for $\xi_1
\, \frac{\partial f}{ \partial x_1} + \cdots + \xi_n \, \frac{
\partial f}{ \partial x_n}$ to be the general infinitesimal
transformation of an $r$-term group, a certain condition is necessary
and sufficient, namely: when $\xi_{ k1}, \dots, \xi_{ kn}$ and $\xi_{
j1}, \dots, \xi_{ jn}$ are two particular systems of solutions, then
the expression:
\[
\sum_{\nu=1}^n\,
\bigg(
\xi_{k\nu}\,
\frac{\partial\xi_{ji}}{\partial x_\nu}
-
\xi_{j\nu}\,
\frac{\partial\xi_{ki}}{\partial x_\nu}
\bigg)
\ \ \ \ \ \ \ \ \ \ \ \ \
{\scriptstyle{(i\,=\,1\,\cdots\,n)}}
\]
must always represent a system of solutions. 
As a result, we have the

\def\thetheorem{28}\begin{theorem}\label{Theorem-28}
If $\xi_1, \dots, \xi_n$, as functions of $x_1, \dots, x_n$, 
are determined by certain linear and homogeneous
partial differential equations:
\[
\sum_{\nu=1}^n\,
A_{\mu\nu}(x)\,\xi_\nu
+
\sum_{\nu,\,\,\pi}^{1\cdots n}\,
B_{\mu\nu\pi}(x)\,
\frac{\partial\xi_\nu}{\partial x_\pi}
+\cdots
=
0
\ \ \ \ \ \ \ \ \ \ \ \ \
{\scriptstyle{(\mu\,=\,1,\,2\,\cdots)}},
\]
then the expression $\xi_1 \, \frac{ \partial f}{ \partial x_1} +
\cdots + \xi_n \, \frac{ \partial f}{\partial x_n}$ represents the
general infinitesimal transformation of a finite continuous group if
and only if: firstly the most general system of solutions of these
differential equations depends only on a finite number of arbitrary
constants, and secondly from any two particular systems of solutions
$\xi_{ k1}, \dots, \xi_{ kn}$ and $\xi_{ j1}, \dots, \xi_{ jn}$, by
the formation of the expression:
\[
\sum_{\nu=1}^n\,
\bigg(
\xi_{k\nu}\,
\frac{\partial\xi_{ji}}{\partial x_\nu}
-
\xi_{j\nu}\,
\frac{\partial\xi_{ki}}{\partial x_\nu}
\bigg)
\ \ \ \ \ \ \ \ \ \ \ \ \
{\scriptstyle{(i\,=\,1\,\cdots\,n)}},
\]
one always obtains a new system of solutions.
\end{theorem}

If $\xi_1 \, \frac{ \partial f}{ \partial x_1} + \cdots + \xi_n \,
\frac{ \partial f}{\partial x_n}$ really is the general infinitesimal
transformation of a finite group, then naturally, the above
differential equations are the defining equations of this group.

If the defining equations of a finite group are
given, the numbers $\nu_0$, $\nu_1$, \dots, 
$\nu_{ s-1}$ discussed earlier on can be determined
by differentiation and by elimination; 
the number:
\[
r
=
(\varepsilon_0-\nu_0)
+
(\varepsilon_1-\nu_1)
+\cdots+
(\varepsilon_{s-1}-\nu_{s-1})
\]
indicates how many parameters the group contains.

In the first one of the two former examples, one has:
\[
s=2,\ \ \ \ \
\nu_0=0,\ \ \ \ \
\nu_1=0,\ \ \ \ \
\varepsilon_0=2,\ \ \ \ \
\varepsilon_1=4,
\]
and therefore $r = 6$. 
In the second example, one has:
\[
s=3,\ \ \ \ \
\nu_0=0,\ \ \ \ \
\nu_1=0,\ \ \ \ \
\nu_2=4,\ \ \ \ \
\varepsilon_0=2,\ \ \ \ \
\varepsilon_1=4,\ \ \ \ \
\varepsilon_2=6,
\]
whence $r = 8$. 

\sectionengellie{\S\,\,\,51.}

\label{S-188-sq}
Now, let again: 
\[
X_kf
=
\sum_{i=1}^n\,\xi_{ki}(x_1,\dots,x_n)\,
\frac{\partial f}{\partial x_i}
\ \ \ \ \ \ \ \ \ \ \ \ \
{\scriptstyle{(k\,=\,1\,\cdots\,r)}}
\]
be independent infinitesimal transformations of an $r$-term group. We
imagine that the defining equations of this group are set up and are
brought to the form discussed above, hence resolved with respect to
$\nu_k$ of the $\varepsilon_k$ differential quotients of order $k$ of
the $\xi$ ($k=0, 1, \dots, s$); here, as said earlier on, $\nu_k$ is
always $< \varepsilon_k$, except in the case $k = s$, in which $\nu_s
= \varepsilon_s$.

The coefficients in the resolved defining equations are visibly
\emphasis{rational} functions of the $\xi$ together with their
differential quotients of first order, up to th $s$-th. Now, since as
a matter of principle (cf. p.~\pageref{S-171}), we restrict ourselves
to systems of values $x_1, \dots, x_n$ for which all $\xi_{ ki}$
behave regularly, then the meant coefficient will in general behave
regularly too for the systems of values $x_1, \dots, x_n$ coming into
consideration, but obviously only in general: there can well exist
points $x_1, \dots, x_n$ in which all $\xi_{ ki}$ indeed behave
regularly, but in which not all coefficients of the solved defining
equations do so. In what follows, one must everywhere pay heed 
to this distinction between the different points $x_1, \dots, x_n$. 

Let $x_1^0, \dots, x_n^0$ be \label{S-189-sq}
\emphasis{a point for which all
coefficients in the solved defining equations behave regularly}; then
$\xi_1, \dots, \xi_n$ are, in a certain neighbourhood of $x_1^0,
\dots, x_n^0$, ordinary power series with respect to the $x_k -
x_k^0$:
\[
\xi_i
=
g_i^0
+
\sum_{\nu=1}^n\,g_{i\nu}'\,(x_\nu-x_\nu^0)
+
\sum_{\nu,\,\,\pi}^{1\cdots n}\,
g_{i\nu\pi}''\,(x_\nu-x_\nu^0)(x_\pi-x_\pi^0)
+\cdots,
\]
where, always, the terms of the same order are thought to be combined
together. Here, the coefficients $g^0$, $g'$, \dots\, must be
determined in such a way that the given differential equations are
identically satisfied after insertion of the power series expansions
for $\xi_1, \dots, \xi_n$.

If we now remember that our group is $r$-term, then we immediately
realize that on the whole, $r$ of the coefficients $g^0$, $g'$,
\dots\, remain undetermined, hence that certain amongst the initial
values which the $\xi$ and its differential quotients take for $x_1 =
x_1^0$, \dots, $x_n = x_n^0$, can be chosen arbitrarily. Now, since
our differential equations show that all differential quotients of
order $s$ and of higher order can be expressed in terms of the
differential quotients of orders zero up to $(s-1)$, and in terms of
$x_1, \dots, x_n$, it follows that exactly $r$ amongst the initial
values $g^0$, $g'$, $g''$, \dots, $g^{ (s-1)}$ can be chose
arbitrarily, or, what amounts to the same: the $\varepsilon_0 +
\varepsilon_1 + \cdots + \varepsilon_{ s-1}$ mentioned initial values
must be linked together by certain $\sum \, \varepsilon_k - r$
independent relations. At present, we want to set up these relations.

Or system of differential equations produces all relations which
exist, at an {\em arbitrary} $x$, between the $\xi$ and its
differential quotients of first order up to the $(s-1)$-th. Hence,
when we make the substitution $x_i = x_i^0$ in our differential
equations, we obtain certain relationships by which the initial values
$g^0$, $g'$, \dots, $g^{ (s-1)}$ are linked together. In this manner,
it results $\nu_0$ linear homogeneous relations between the $g^0$
alone, furthermore $\nu_1$ relations between the $g'$ and certain
$g^0$, and in general $\nu_k$ relations between the $g^{ (k)}$ and
certain $g^{ (k-1)}$, \dots, $g'$, $g^0$. The relations in question
are resolved, and to be precise, they are resolved with respect to
$\nu_0$ of the $g^0$, with respect to $\nu_1$ of the $g'$, and so on,
and in total, there are $\sum\, \nu_k = \sum\, \varepsilon_k - r$
independent relations. Now, since according to what has been said,
there do not exist more than $\sum \, \varepsilon_k - r$ independent
relations between $g^0$, $g'$, \dots, $g^{ (s-1)}$, we therefore have
found all relations by means of which these $\sum \, \varepsilon_k$
initial values are linked together; at the same time, we obtain $\nu_0
+ \cdots + \nu_{ s-1}$ of the quantities $g^0, \dots, g^{ s-1}$
represented as \emphasis{linear homogeneous} functions of the $\sum (
\varepsilon_k - \nu_k) = r$ remaining ones, which stay entirely
arbitrary.

\medskip

We use the preceding observations in order to draw conclusions from
them about the constitution of the particular system of solutions to
our differential equations.

At first, it can be shown that there does not exist a particular
system of solutions $\xi_1, \dots, \xi_n$ whose power series expansion
with respect to the $x_i - x_i^0$ begins with terms of order $s$ or
yet higher order. Indeed, in the power series expansion of such a
system of solutions, the coefficients $g^0$, $g'$, \dots, $g^{ (s-1)}$
would all be equal to zero, so all $g^{ (s)}$, $g^{ (s+1)}$, \dots\,
would also vanish and the system of solutions would therefore reduce
to: $\xi_1 = 0$, \dots, $\xi_n = 0$. This system of solutions
certainly satisfies the given differential equations, but
alone it does not deliver us any infinitesimal
transformation of our group, and is therefore useless.

But for that, there is a certain number of particular systems of
solutions $\xi_1, \dots, \xi_n$ whose power series expansions begin
with terms of order lower than the $s$-th, let us say with terms of
$k$-th order. The coefficients $g^0$, $g'$, \dots, $g^{ (k-1)}$ can
then all be chosen equal to zero; so the existing relations between
these are all satisfied, and it yet remains only $\nu_k$ relations
between the $g^{ (k)}$, $\nu_{ k+1}$ relations between the $g^{
(k+1)}$ and certain $g^{ (k)}$, \dots, and lastly, $\nu_{ s-1}$
relations between the $g^{ (s-1)}$ and certain $g^{ (s-2)}$, \dots,
$g^{ (k)}$. As a result, there are in sum still $(\varepsilon_k -
\nu_k) + \cdots + (\varepsilon_{ s-1} - \nu_{ s-1})$ of the constants
$g^{ (k)}$, \dots, $g^{ (s-1)}$ which can be chosen arbitrarily, and
when one disposes of these constants in such a way that not all the
$\varepsilon_k$ quantities $g^{ (k)}$ vanish, then one always obtains
a particular system of solutions $\xi_1, \dots, \xi_n$ whose power
series expansions with respect to the $x_i - x_i^0$ contain terms of
the $k$-th order, but no terms of lower order.

If $\xi_1, \dots, \xi_n$ is a particular system of solutions
of our differential equations, then the infinitesimal
transformation:
\[
\xi_1\,\frac{\partial f}{\partial x_1}
+\cdots+
\xi_n\,\frac{\partial f}{\partial x_n}
\]
belongs to our group. From what has been said above, it results that
this group always contains infinitesimal transformations whose power
series expansion with respect to the $x_i - x_i^0$ begins with terms
of order $k$, only as soon as $k$ is one of the numbers $0$, $1$, $2$,
\dots, $s-1$. By contrast, there are no infinitesimal transformations
in the group whose power series expansions begin with terms of order
$s$ or of higher order. Naturally, all of this is proved only under
the assumption that the coefficients of the resolved differential
equations behave regularly in the point $x_1^0, \dots, x_n^0$.


In order to meet, at least up some extent, the requirements for
conciseness of the expression, we want henceforth to say: \emphasis{an
infinitesimal transformation is of the $k$-th order in the $x_i -
x_i^0$ when its power series expansion with respect to the $x_i -
x_i^0$ begins with terms of order $k$}. Then we can enunciate the
preceding result as follows:

If $x_1^0, \dots, x_n^0$ is a point for which the coefficients in the
resolved defining equations of the group $X_1 f, \dots, X_r f$ behave
regularly, then the group contains certain infinitesimal
transformations of zeroth order in the $x_i - x_i^0$, certain of the
first order, and in general, certain of the $k$-th order, where $k$
means one arbitrary number amongst $0$, $1$, \dots, $s-1$; by
contrast, the group contains no infinitesimal transformation of order
$s$ or of higher order in the $x_i - x_i^0$.

It is clear that two infinitesimal transformations of different orders
in the $x_i - x_i^0$ are always independent of each other. Actually,
for the examination whether several given infinitesimal
transformations are independent of each other or not, the
consideration of the terms of lowest order in their power series
expansions already settles the question many times; indeed, if the
terms of lowest order, taken for themselves, determine independent
infinitesimal transformations, then the given infinitesimal
transformations are also independent of each other.

The general expression of an infinitesimal transformation which
belongs to our group and which is of order $k$ with respect to the
$x_i - x_i^0$ contains, as we know, $(\varepsilon_k - \nu_k) + \cdots
+ (\varepsilon_{ s-1} - \nu_{ s-1}) = \rho_k$ arbitrary and essential
constants, namely the $\varepsilon_k - \nu_k$ which can be chosen
arbitrarily amongst the $g^{ (k)}$, the $\varepsilon_{ k+1} - \nu_{
k+1}$ arbitrary amongst the $g^{ (k+1)}$, and so on; however here, the
arbitrariness of these $\varepsilon_k - \nu_k$ quantities $g^{ (k)}$
is restricted inasmuch as not all $g^{ (k)}$ are allowed to vanish
simultaneously. From this, it follows that $\rho_k$ independent
infinitesimal transformations of our group can be exhibited which are
of order $k$ in the $x_i - x_i^0$; but it is easy to see that from
these $\rho_k$ infinitesimal transformations, one can derive in total
$\rho_{ k+1}$ independent ones which are of the $(k+1)$-th order, or
yet of higher order. The general expression of an infinitesimal
transformation which is linearly deduced from these $\rho_k$ ones
indeed contains exactly the same arbitrary constants as the general
expression of an infinitesimal transformation of order $k$ in the $x_i
- x_i^0$, with the only difference that in the first expression, all
the $\varepsilon_k - \nu_k$ available $g^{ (k)}$ can be set equal to
zero, which gives always an infinitesimal transformation of order
$(k+1)$ or of higher order. Consequently, amongst these $\rho_k$
infinitesimal transformations of order $k$, there are only $\rho_{
k+1} - \rho_k = \varepsilon_k - \nu_k$ which are independent of each
other and out of which no infinitesimal transformation of $(k+1)$-th
order or of higher order in the $x_i - x_i^0$ can be linearly deduced.

We recapitulate the present result in the

\def\thetheorem{29}\begin{theorem}
\label{Theorem-29-S-192} 
To every $r$-term
group $X_1 f, \dots, X_r f$ in $n$ variables $x_1, \dots, x_n$ is
associated a completely determined entire number $s \geqslant 1$ of
such a nature that, in the neighbourhood of a point $x_i^0$ for which
the coefficients of the resolved defining equations behave regularly,
the group contains certain infinitesimal transformations of zeroth, of
first, \dots, of $(s-1)$-th order in the $x_i - x_i^0$, but none of
$s$-th or of higher order. In particular, one can always select $r$
independent infinitesimal transformations of the group such that, for
each one of the $s$ values $0$, $1$, \dots, $s-1$ of the number $k$,
exactly $\varepsilon_k - \nu_k$ \label{S-192}
mutually independent infinitesimal
transformations of order $k$ in the $x_i - x_i^0$ are extant out of
which no infinitesimal transformation of order $(k+1)$ or of higher
order can be linearly deduced. At the same time, the number $\nu_k$
can be determined from the defining equations for the general
infinitesimal transformation $\xi_1 \, \frac{ \partial f}{ \partial
x_1} + \cdots + \xi_n \, \frac{ \partial f}{ \partial x_n}$ of the
group, and from $\varepsilon_k$, which denotes the number of all
differential quotients of order $k$ of the $\xi_1, \dots, \xi_n$ with
respect to $x_1, \dots, x_n$ and is always larger than $\nu_k$.
\end{theorem}

{\sf Example.} Earlier on, we have already mentioned the equations:
\[
\frac{\partial^2\xi_1}{\partial x_1^2}
=
\frac{\partial^2\xi_1}{\partial x_1\partial x_2}
=
\frac{\partial^2\xi_1}{\partial x_2^2}
=
\frac{\partial^2\xi_2}{\partial x_1^2}
=
\frac{\partial^2\xi_2}{\partial x_1\partial x_2}
=
\frac{\partial^2\xi_2}{\partial x_2^2}
=
0
\]
as the defining equations of the six-term linear group:
\[
x_1'
=
a_1\,x_1+a_2\,x_2+a_3,
\ \ \ \ \ \ \ \ \
x_2'
=
a_4\,x_1+a_5\,x_2+a_6.
\]

These defining equations are already presented under the resolved
form; all the appearing coefficients are equal to zero, hence they
behave regularly. Amongst the infinitesimal transformations of the
group, we can select the following six mutually independent ones:
\[
\frac{\partial f}{\partial x_1},\ \ \
\frac{\partial f}{\partial x_2},\ \ \
(x_1-x_1^0)\,\frac{\partial f}{\partial x_1},\ \ \
(x_2-x_2^0)\,\frac{\partial f}{\partial x_1},\ \ \
(x_1-x_1^0)\,\frac{\partial f}{\partial x_2},\ \ \
(x_2-x_2^0)\,\frac{\partial f}{\partial x_2}\,;
\]
the first two of them are of zeroth order, and the last four are of
first order in the $x_i - x_i^0$. 

\medskip

For the calculations with infinitesimal transformations, 
the expressions of the form:
\[
X\big(Y(f)\big)
-
Y\big(X(f)\big)
=
\leftbracket
X,\,Y
\rightbracket
\]
play an important rôle. So if $Xf$ and $Yf$ are expanded in the
neighbourhood of the point $x_i^0$ with respect to the power of the
$x_i - x_i^0$, the related question is how does the transformation
$\leftbracket X, \, Y \rightbracket$ behave in this point.

Let the power series expansion of $Xf$ begin with terms of order
$\mu$, let that of $Yf$ begin with terms of order $\nu$, that is to
say, let:
\[
Xf
=
\sum_{k=1}^n\,
(\xi_k^{(\mu)}+\cdots)\,
\frac{\partial f}{\partial x_k},
\ \ \ \ \ \ \ \
Yf
=
\sum_{j=1}^n\,
(\eta_j^{(\nu)}+\cdots)\,
\frac{\partial f}{\partial x_j},
\]
where the $\xi^{ (\mu)}$ and the $\eta^{ (\nu)}$ denote homogeneous
functions of order $\mu$ and of order $\nu$, respectively, in the $x_i
- x_i^0$, while the terms of higher order in the $x_i - x_i^0$ are
left out. Under these assumptions, the power series expansion for
$\leftbracket X, \, Y \rightbracket$ is, if one only considers terms
of the lowest order, the following:
\[
\leftbracket
X,\,Y
\rightbracket
=
\sum_{j=1}^n\,
\bigg\{
\sum_{k=1}^n\,
\bigg(
\xi_k^{(\mu)}\,
\frac{\partial\eta_j^{(\nu)}}{\partial x_k}
-
\eta_k^{(\nu)}\,
\frac{\partial\xi_j^{(\mu)}}{\partial x_k}
\bigg)
+\cdots
\bigg\}
\frac{\partial f}{\partial x_j}.
\]

So the terms of lowest order in $\leftbracket X, \, Y \rightbracket$
are of order $\mu + \nu - 1$ and they stem solely and only from the
terms of orders $\mu$ and $\nu$ in $Xf$ and in $Yf$, respectively. 

\def\thetheorem{30}\begin{theorem}\label{Theorem-30}
If $Xf$ and $Yf$ are two infinitesimal transformations whose power
series expansions with respect to the powers of $x_1 - x_1^0$, \dots,
$x_n - x_n^0$ begin with terms of orders $\mu$ and $\nu$,
respectively, then the power series expansion of the infinitesimal
transformation $XY f - YX f = \leftbracket X, \, Y \rightbracket$
begins with terms of order $(\mu + \nu - 1)$ which are entirely
determined by the terms of orders $\mu$ and $\nu$ in $Xf$ and in $Yf$,
respectively. If these terms of order $(\mu + \nu - 1)$ vanish, then
about the power series expansion of $\leftbracket X, \, Y
\rightbracket$, it can only be said that it starts with terms of order
$(\mu + \nu)$, or of higher order.
\end{theorem}

If the two numbers $\mu$ and $\nu$ are larger than one, then the
number $\mu + \nu - 1$ is larger than both of them. This remark is
often of great utility for the calculations with infinitesimal
transformations of various orders.

For the derivation of the Theorem~30, it is not at all assumed that
the two infinitesimal transformations $Xf$ and $Yf$ belong to a group;
the only assumption is that both $Xf$ and $Yf$ can be expanded in
powers of $x_k - x_k^0$.

\sectionengellie{\S\,\,\,52.}

Let the defining equations of an $r$-term group be given in the form
discussed earlier on, hence resolved with respect to $\nu_k$ of the
$\varepsilon_k$ differential quotients of order $k$ of $\xi_1, \dots,
\xi_n$. Moreover, let $x_i^0$ be a point in which the coefficients
of the resolved defining equations behave regularly.

Under these assumptions, we can expand the infinitesimal
transformations of our group in ordinary power series
of the $x_i - x_i^0$. 
We even know that our group contains a completely determined
number of independent
infinitesimal transformations, namely $\varepsilon_0 
- \nu_0$, which are of zeroth order in the $x_i - x_i^0$ and
out of which no infinitesimal transformation of
first order or of higher order can be linearly deduced; 
moreover, the group contains a completely determined
number of independent infinitesimal transformations, namely 
$\varepsilon_1 - \nu_1$, of first order in the $x_i - x_i^0$ 
out of which none of second order or of 
higher order can be linearly deduced, and so on.

Shortly, our group associates to every point of the
indicated nature a series of $s$ entire numbers
$\varepsilon_0 - \nu_0$, 
$\varepsilon_1 - \nu_1$, \dots, 
$\varepsilon_{ s-1} - \nu_{ s-1}$
and these entire numbers are the same for
all points of this kind.

Now, there can also be points $\overline{ x}_i$ in
special position, hence points in the neighbourhood of
which the coefficients of the resolved defining equations
do not behave regularly anymore, while by contrast,
all infinitesimal transformations of the group
can be expanded in ordinary power series in the 
$x_i - \overline{ x}_i$. 
If $\overline{ x}_1, \dots, 
\overline{ x}_n$ is a determined point of this sort, then
naturally, there is in our group a completely
determined number of infinitesimal transformations
of zeroth order in the $x_i - \overline{ x}_i$
out of which no infinitesimal transformation
of higher order can be linearly deduced, and so on. 

Consequently, our group also associates to every point
in special position a determined, obviously 
finite series of entire numbers; 
frequently, to two different points in special position 
there will be associated two also different series of
entire numbers. 

An example will best make clear the matter.

The defining equations of the two-term group
$\frac{ \partial f}{ \partial x_1}$, 
$x_2^2 \, \frac{ \partial f}{ \partial x_1}$ 
read in the resolved form:
\[
\aligned
\xi_2=0,\ \ \ \ \
\frac{\partial\xi_1}{\partial x_1}
&
=
-\,\frac{\partial\xi_2}{\partial x_1}
=
\frac{\partial\xi_2}{\partial x_2}
,=
0
\\
\frac{\partial^2\xi_1}{\partial x_1^2}
&
=
\frac{\partial^2\xi_1}{\partial x_1\partial x_2}
=
0,\ \ \ \ \
\frac{\partial^2\xi_1}{\partial x_2^2}
=
\frac{1}{x_2}\,\frac{\partial\xi_1}{\partial x_2},
\\
\frac{\partial^2\xi_2}{\partial x_1^2}
&
=
\frac{\partial^2\xi_2}{\partial x_1\partial x_2}
=
\frac{\partial^2\xi_2}{\partial x_2^2}
=
0.
\endaligned
\]

The coefficients appearing here behave regularly for all points $x_1,
x_2$ located in the finite, except only for the points of the line
$x_2 = 0$.

At first, let us consider a point $x_1^0$, $x_2^0$ with nonvanishing
$x_2^0$. We have $s = 2$, and moreover $\varepsilon_0 = 2$,
$\varepsilon_1 = 4$, $\nu_0 = 1$, $\nu_1 = 3$, hence to the point
$x_1^0, x_2^0$ are associated the two numbers $1, 1$. All
infinitesimal transformations of the group can be linearly deduced
from the two:
\[
\frac{\partial f}{\partial x_1},\ \ \ \ \ \ \
\Big(
x_2-x_2^0
+
\frac{1}{2x_0}\,(x_2-x_2^0)^2
\Big)\,
\frac{\partial f}{\partial x_1}
\]
amongst which the first is of zeroth order in the $x_i - x_i^0$, and
the second of first order.

Next, let us consider a point $\overline{ x}_1, \overline{ x}_2 = 0$.

To such a point, the group associates the three numbers $1, 0, 1$,
since amongst its infinitesimal transformations there are none of
first order in the $x_i - \overline{ x}_i$, but one of second order,
hence of $s$-th order, namely: $x_2^2 \, \frac{ \partial f}{ \partial
x_1}$. ---

If $x_1^0, \dots, x_n^0$ is a point for which the coefficients of the
resolved defining equations behave regularly, then according to
Theorem~29, the group certainly contains infinitesimal transformations
of zeroth, of first, \dots, of $(s-1)$-th orders in the $x_i - x_i^0$,
but none of $s$-th or of higher order. Now, our example
discussed just now shows that for a point
$\overline{ x}_i$ in which not all the coefficients in 
question behave regularly, no general 
statement of this kind holds anymore: the group
can very well contain infinitesimal transformations
of $s$-th order in the $x_i - \overline{ x}_i$, 
and perhaps also some of higher order; 
on the other hand, it can occur that for one number
$k < s$, the group actually contains no infinitesimal 
transformation of $k$-th order in the $x_i - 
\overline{ x}_i$. 

\medskip

\label{S-196-sq}
If $x_1^0, \dots, x_n^0$ denotes an arbitrary point in which all the
$\xi$ behave regularly, then as already said, the infinitesimal
transformations of our group can be classified according to their
orders in the $x_i - x_i^0$. It is of great importance that this
classification stays obtained when in place of the $x$, new variables
$y_1, \dots, y_n$ are introduced. Of course, the concerned change of
variables must, in the neighbourhood of the place $x_1^0, \dots,
x_n^0$, possess the following properties: $y_1, \dots, y_n$ must
firstly be ordinary power series in the $x_i - x_i^0$:
\def\theequation{2}\begin{equation}
y_k
=
y_k^0
+
\sum_{i=1}^n\,a_{ki}\,(x_i-x_i^0)
+\cdots
\ \ \ \ \ \ \ \ \ \ \ \ \
{\scriptstyle{(k\,=\,1\,\cdots\,n)}}\,;
\end{equation}
and secondly, $x_1, \dots, x_n$ must also be representable as ordinary
power series in the $y_k - y_k^0$ and in fact, so that every $x_i$ for
$y_1 = y_1^0$, \dots, $y_n = y_n^0$ must take the value $x_i^0$. If
the first one of these two requirements is satisfied, it is known that
the second one is then always satisfied, when the determinant $\sum \,
\pm \, a_{ 11} \cdots a_{ nn}$ is different from zero.

Now, in order to prove that the discussed classification stays
obtained after the transition to the variables $y_1, \dots, y_n$, we
need only to show that every infinitesimal transformation of $\mu$-th
order in the $x_i - x_i^0$ converts, by the introduction of the new
variables $y_k - y_k^0$, into an infinitesimal transformation of
$\mu$-th order in the $y_k - y_k^0$. But this is not difficult.

The general form of an infinitesimal transformation of the $\mu$-th
order in the $x_i - x_i^0$ is:
\[
Xf
=
\sum_{j=1}^n\,
(\xi_j^{(\mu)}+\cdots)\,
\frac{\partial f}{\partial x_j}\,;
\]
here, $\xi_1^{ (\mu)}, \dots, \xi_n^{ (\mu)}$ denote entire rational
functions\footnote{\,
\deutschplain{ganze rationale Functionen}, that is to say,
polynomials. }
which are homogeneous of order $\mu$ and do not all vanish; the terms
of orders $(\mu+1)$ and higher are left out.

By the introduction of $y_1, \dots, y_n$ it comes:
\[
Xf
=
\sum_{k=1}^n\,Xy_k\,
\frac{\partial f}{\partial y_k}\,;
\]
here at first, the $Xy_k$ are ordinary power series in the
$x_i - x_i^0$:
\[
Xy_k
=
\sum_{j=1}^n\,a_{kj}\,\xi_j^{(\mu)}
+\cdots
\]
and they begin with terms of order $\mu$.
These terms of order $\mu$ do not all vanish, 
since otherwise one would have:
\[
\sum_{j=1}^n\,a_{kj}\,
\xi_j^{(\mu)}
=
0
\ \ \ \ \ \ \ \ \ \ \ \ \
{\scriptstyle{(k\,=\,1\,\cdots\,n)}},
\]
which is impossible, because the determinant
$\sum \, \pm \, a_{ 11} \cdots a_{ nn}$
is different from zero, and because $\xi_1^{ (\mu)}, 
\dots, \xi_n^{ (\mu)}$ do not all vanish.
Now, if in $Xy_1, \dots, Xy_n$ we express
the $x_i$ in terms of the $y_i$, we
obtains $n$ ordinary power series in the
$y_i - y_i^0$. These power series
likewise begin with terms of order $\mu$
which do not all vanish. Indeed, one
obtains the terms of order $\mu$ in question by 
substituting, in the $n$ expressions:
\[
\sum_{j=1}^n\,
a_{kj}\,\xi_j^{(\mu)}
\ \ \ \ \ \ \ \ \ \ \ \ \
{\scriptstyle{(k\,=\,1\,\cdots\,n)}},
\]
the $x$ for the $y$ by means of the equations:
\[
y_k
=
y_k^0
+
\sum_{i=1}^n\,a_{ki}\,(x_i-x_i^0)
\ \ \ \ \ \ \ \ \ \ \ \ \
{\scriptstyle{(k\,=\,1\,\cdots\,n)}}\,;
\]
but since the $n$ shown expressions do not all vanish, 
then they also do not all vanish after
introduction of the $y$.

Consequently, the infinitesimal transformation
$Xf$ is transferred, by the introduction of the $y$, 
to an infinitesimal transformation which is
of order $\mu$ in the $y_i - y_i^0$. But this
was to be shown.

As a result, we have the

\def\theproposition{1}\begin{proposition}
\label{S-197}
If, in an infinitesimal transformation
$Xf$ which is of order $\mu$ in
$x_1 - x_1^0$, \dots, $x_n - x_n^0$, one introduces
new variables:
\[
\aligned
y_k
=
y_k^0
+
\sum_{i=1}^n\,a_{ki}(x_i-x_i^0)
&
+
\sum_{i,\,\,j}^{1\cdots n}\,
a_{kij}\,(x_i-x_i^0)\,(x_j-x_j^0)
+\cdots
\\
&
{\scriptstyle{(k\,=\,1\,\cdots\,n)}},
\endaligned
\]
where the determinant $\sum \, \pm \, a_{ 11} \cdots a_{nn}$ is
different from zero, then $Xf$ converts into an infinitesimal
transformation of order $\mu$ in $y_1 - y_1^0$, \dots, $y_n - y_n^0$.
\end{proposition} 

From this, it immediately follows the somewhat
more specific

\def\theproposition{2}\begin{proposition}
If, in the neighbourhood of the point $x_1^0, \dots, x_n^0$, an
$r$-term group contains exactly $\tau_\mu$ independent infinitesimal
transformations of order $\mu$ in the $x_i - x_i^0$ out of which none
of higher order can be linearly deduced, and if new variables:
\[
\aligned
y_k
=
y_k^0
+
\sum_{i=1}^n\,a_{ki}(x_i-x_i^0)
&
+\cdots
\ \ \ \ \ \ \ \ \ \
{\scriptstyle{(k\,=\,1\,\cdots\,n)}},
\endaligned
\]
are introduced in this group, where the determinant $\sum \, \pm \,
a_{ 11} \cdots a_{ nn}$ is different from zero, then in turn in the
neighbourhood of the point $y_k^0$, the new group which one finds in
this way contains exactly $\tau_\mu$ independent infinitesimal
transformations of order $\mu$ in the $y_k - y_k^0$ out of which none
of higher order can be linearly deduced.
\end{proposition} 

One therefore sees: the series of entire numbers which the initial
group associates to the point $x_1^0, \dots, x_n^0$ is identical to
the series of entire numbers which the new groups associates to the
point $y_1^0, \dots, y_n^0$.

\sectionengellie{\S\,\,\,53.}

If one knows the defining equations of an $r$-term group and if one
has resolved them in the way discussed earlier on, then as we have
seen, one can immediately identify the numbers $\varepsilon_k - \nu_k$
defined above. For every point $x_1^0, \dots, x_n^0$ in which the
coefficients of the resolved defining equations behave regularly, one
therefore knows the number of all independent infinitesimal
transformations of the group which are of order $k$ in the $x_i -
x_i^0$ and which possess the property that out of them, no
infinitesimal transformation of order $(k+1)$ or of higher order can
be linearly deduced.

Naturally, one can compute the numbers in question also for the points
$x_1^0, \dots, x_n^0$ in which the coefficients of the resolved
defining equations do not behave regularly. For that, the knowledge
of the defining equations already suffices, however it is incomparably
more convenient when $r$ arbitrary independent infinitesimal
transformations are already given, which is what we will assume in the
sequel. Then one proceeds as follows.

At first, one determines how many independent
infinitesimal transformations of order $k$ or higher 
in the $x_i - x_i^0$ the group contains.
To this aim, one expands the general infinitesimal
transformation:
\[
e_1\,X_1f
+\cdots+
e_r\,X_rf
\]
with respect to the powers of the $x_i - x_i^0$ and
then in the $n$ expressions: 
\[
e_1\,\xi_{1i}
+\cdots+
e_r\,\xi_{ri}
\ \ \ \ \ \ \ \ \ \ \ \ \
{\scriptstyle{(i\,=\,1\,\cdots\,n)}},
\]
one sets equal to zero all coefficients of zeroth, of first, \dots, of
$(k-1)$-th order. In this way, one obtains a certain number of linear
homogeneous equations between $e_1, \dots, e_r$; one then easily
determines how many independent infinitesimal transformations are
extant amongst these equations by calculating certain determinants; if
$r - \omega_k$ is the number of independent equations, then it follows
that the group contains exactly $\omega_k$ independent infinitesimal
transformations of order $k$ or higher in the $x_i - x_i^0$. So
obviously, $\omega_k - \omega_{ k+1}$ is the number of independent
infinitesimal transformations of order $k$ out of which no
infinitesimal transformation of higher order can be linearly deduced.

It is hardly not necessary to make the observation that the
operations just indicated remain applicable also to
every point $x_1^0, \dots, x_n^0$ for which the coefficients
of the resolved defining equations behave regularly. 

Somewhat more precisely, we want to occupy ourselves with the
infinitesimal transformations $\sum \, e_j \, X_j f$ of the group $X_1
f, \dots, X_r f$ whose power series expansion in the $x_i - x_i^0$
contain only terms of the first and higher orders, but none of the
zeroth. At first, we shall examine how many independent infinitesimal
transformations of this nature there are and we shall show how one can
set up them in a simple manner. Here, by $x_1^0, \dots, x_n^0$, we
understand a completely arbitrary, though determined, point.

Evidently, such infinitesimal transformations are characterized by the
fact that themselves, and as well the one-term groups generated by
them, do leave at rest the point $x_i = x_i^0$
(cf. Chap.~\ref{kapitel-7}, p.~\pageref{S-134}), or, what amounts to
the same, by the fact that they are the only ones amongst the
infinitesimal transformations $\sum \, e_j \, X_jf$ which do not
attach any direction to the point $x_i = x_i^0$.

Analytically, the most general transformation $\sum\, e_j \, X_j f$ of
the concerned constitution will be determined by the equations:
\[
e_1\,\xi_{1i}(x_1^0,\dots,x_n^0)
+\cdots+
e_r\,\xi_{ri}(x_1^0,\dots,x_n^0)
=
0
\ \ \ \ \ \ \ \ \ \ \ \ \
{\scriptstyle{(i\,=\,1\,\cdots\,n)}}.
\]
Now, if in the matrix:
\def\theequation{3}\begin{equation}
\label{S-199}
\left\vert
\begin{array}{cccc}
\xi_{11}(x^0) & \,\cdot\, & \,\cdot\, & \xi_{1n}(x^0)
\\
\cdot & \,\cdot\, & \,\cdot\, & \cdot
\\
\xi_{r1}(x^0) & \,\cdot\, & \,\cdot\, & \xi_{rn}(x^0)
\end{array}
\right\vert,
\end{equation}
all $(r+1) \times (r+1)$ determinants vanish, but not all $h \times h$
ones, then $h$ of the quantities $e_1, \dots, e_r$ can be represented
as linear homogeneous functions of the $r - h$ left ones, which remain
completely arbitrary. As a result, we obtain the following simple but
important result:

\def\theproposition{3}\begin{proposition}
\label{Satz-3-S-200}
If all $(h+1) \times (h+1)$ determinants of the
matrix:
\[
\left\vert
\begin{array}{cccc}
\xi_{11}(x) & \,\cdot\, & \,\cdot\, & \xi_{1n}(x)
\\
\cdot & \,\cdot\, & \,\cdot\, & \cdot
\\
\xi_{r1}(x) & \,\cdot\, & \,\cdot\, & \xi_{rn}(x)
\end{array}
\right\vert,
\]
vanish for $x_1 = x_1^0$, \dots, $x_n = x_n^0$, but not all $h \times
h$ ones vanish, then the $r$-term group:
\[
X_kf
=
\sum_{i=1}^n\,\xi_{ki}(x_1,\dots,x_n)\,
\frac{\partial f}{\partial x_i}
\ \ \ \ \ \ \ \ \ \ \ \ \
{\scriptstyle{(k\,=\,1\,\cdots\,r)}}
\]
contains exactly $r - h$ independent infinitesimal 
transformations which, when expanded in power series
in $x_1 - x_1^0$, \dots, $x_n - x_n^0$, 
contain no term of zeroth order\,\,---\,\,which, in 
other words, leave at rest the point $x_1^0, \dots, x_n^0$.
\end{proposition}

At the same time, it yet comes from what has been said
the following

\def\theproposition{4}\begin{proposition}
\label{Satz-4-S-200}
When all $(h+1) \times (h+1)$ determinants of
the matrix:
\[
\left\vert
\begin{array}{cccc}
\xi_{11}(x) & \,\cdot\, & \,\cdot\, & \xi_{1n}(x)
\\
\cdot & \,\cdot\, & \,\cdot\, & \cdot
\\
\xi_{r1}(x) & \,\cdot\, & \,\cdot\, & \xi_{rn}(x)
\end{array}
\right\vert,
\]
are set to zero, then the resulting equations determine the locus of
all points $x_1, \dots, x_n$ which admit\footnote{\,
---\,\,in the sense the manifold constituted of such a point is left
invariant, i.e. at rest\,\,--- 
} 
at least $r - h$ independent infinitesimal transformations of the
group $X_1 f, \dots, X_r f$; amongst the found points, those which
do not bring to zero all $h \times h$ determinants of the matrix
admit exactly $r - h$ independent infinitesimal transformations of
the group.
\end{proposition}

Earlier on (p.~\pageref{S-135}), we have underlined that
the infinitesimal transformations $X_1f, \dots, X_r f$ associate
to a determined point $x_1, \dots, x_n$
precisely $h$ independent directions when all $(h+1) \times (h+1)$
determinants of the matrix~\thetag{ 3} vanish, 
while by contrast not all $h \times h$ determinants do. 
From this, we see that the found result can also 
be expressed as follows.

\def\theproposition{5}\begin{proposition}
If an $r$-term group $X_1f, \dots, X_r f$ of the space
$x_1, \dots, x_n$ contains exactly $r - h$ independent
infinitesimal transformations which leave at rest
a determined point $x_1^0, \dots, x_n^0$, then
the infinitesimal transformations of the group associate
to this point exactly $h$ independent directions.
\end{proposition}

At present, we continue one step further to really set up all
infinitesimal transformations $\sum \, e_j \, X_j f$ which leave at
rest a determined point $x_1^0, \dots, x_n^0$.

We want to suppose that in the matrix~\thetag{ 3}, all $(h+1) \times
(h+1)$ determinants vanish, but not all $h \times h$ ones, and
specifically, that in the smaller matrix:
\[
\left\vert
\begin{array}{cccc}
\xi_{11}(x^0) & \,\cdot\, & \,\cdot\, & \xi_{1n}(x^0)
\\
\cdot & \,\cdot\, & \,\cdot\, & \cdot
\\
\xi_{h1}(x^0) & \,\cdot\, & \,\cdot\, & \xi_{hn}(x^0)
\end{array}
\right\vert,
\]
not all $h \times h$ determinants are equal to zero.

Under these assumptions, there are obviously no infinitesimal
transformations of the form $e_1 \, X_1 f + \cdots + e_h \, X_h f$
which leave at rest the point $x_1^0, \dots, x_n^0$; by contrast, the
$r-h$ infinitesimal transformations:
\[
X_{h+k}f+\lambda_{k1}\,X_1f
+\cdots+
\lambda_{kh}\,X_hf
\ \ \ \ \ \ \ \ \ \ \ \ \
{\scriptstyle{(k\,=\,1\,\cdots\,r\,-\,h)}}
\]
do leave it at rest, as soon as one chooses the constants $\lambda$ in
an appropriate way, which visibly is possible only in one single way.
Saying this, $r - h$ independent infinitesimal transformations are
found whose power series expansions contain no term of order zero;
naturally, out of these $r - h$ transformation, 
every other of the
same constitution can be linearly deduced.
From this, it follows that amongst the infinitesimal
transformations of our group which are
of zeroth order in the $x_i - x_i^0$, there are
only $h$ independent ones out of which
no transformation of first order or of higher order
can be linearly deduced; of course, $X_1 f, \dots, X_r f$
are transformations of zeroth order of this nature; 
hence they attach to the point $x_1^0, \dots, x_n^0$ exactly
$h$ independent directions. 

With these words, we have the

\def\theproposition{6}\begin{proposition}
If the $r$ infinitesimal transformations:
\[
X_kf
=
\sum_{i=1}^n\,\xi_{ki}(x_1,\dots,x_n)\,
\frac{\partial f}{\partial x_i}
\ \ \ \ \ \ \ \ \ \ \ \ \
{\scriptstyle{(k\,=\,1\,\cdots\,r)}}
\]
of an $r$-term group of the space $x_1, \dots, x_n$ are
constituted in such a way that for $x_1 = x_1^0$, 
\dots, $x_n = x_n^0$, all $(h+1) \times (h+1)$
determinants, but not all $h \times h$ determinants, 
of the matrix: 
\[
\left\vert
\begin{array}{cccc}
\xi_{11}(x) & \,\cdot\, & \,\cdot\, & \xi_{1n}(x)
\\
\cdot & \,\cdot\, & \,\cdot\, & \cdot
\\
\xi_{r1}(x) & \,\cdot\, & \,\cdot\, & \xi_{rn}(x)
\end{array}
\right\vert
s\]
vanish, and especially, if not all $h \times h$
determinants of the matrix:
\[
\left\vert
\begin{array}{cccc}
\xi_{11}(x) & \,\cdot\, & \,\cdot\, & \xi_{1n}(x)
\\
\cdot & \,\cdot\, & \,\cdot\, & \cdot
\\
\xi_{h1}(x) & \,\cdot\, & \,\cdot\, & \xi_{hn}(x)
\end{array}
\right\vert
\]
are zero for $x_i = x_i^0$, then firstly: all infinitesimal
transformations:
\[
e_1\,X_1f
+\cdots+
e_h\,X_hf
\]
are of zeroth order in the $x_i - x_i^0$ and they attach to the point
$x_1^0, \dots, x_n^0$ exactly $h$ independent directions, and
secondly: one can always choose $h ( r - h)$ constants $\lambda_{
kj}$, but only in one single way, so that in the $r - h$ independent
infinitesimal transformations:
\[
X_{h+k}f+\lambda_{k1}\,X_1f
+\cdots+
\lambda_{kh}\,X_hf
\ \ \ \ \ \ \ \ \ \ \ \ \
{\scriptstyle{(k\,=\,1\,\cdots\,r\,-\,h)}}
\]
all terms of zeroth order in the $x_i - x_i^0$ are missing; then out
of these $r - h$ infinitesimal transformations, one can linearly
deduce all infinitesimal transformations of the group $X_1 f, \dots,
X_r f$ which are of the first order in the $x_i - x_i^0$, or of higher
order.
\end{proposition}

For the sequel, it is useful to state this proposition
in a somewhat more specific way. 

We want to assume that all $(q+1) \times (q+1)$ determinants
of the matrix:
\def\theequation{4}\begin{equation}
\left\vert
\begin{array}{cccc}
\xi_{11}(x) & \,\cdot\, & \,\cdot\, & \xi_{1n}(x)
\\
\cdot & \,\cdot\, & \,\cdot\, & \cdot
\\
\xi_{r1}(x) & \,\cdot\, & \,\cdot\, & \xi_{rn}(x)
\end{array}
\right\vert
\end{equation}
vanish {\em identically}, but that this is not the case
for all $q \times q$ determinants and specially, that
not all $q \times q$ determinants of the matrix:
\def\theequation{5}\begin{equation}
\left\vert
\begin{array}{cccc}
\xi_{11}(x) & \,\cdot\, & \,\cdot\, & \xi_{1n}(x)
\\
\cdot & \,\cdot\, & \,\cdot\, & \cdot
\\
\xi_{q1}(x) & \,\cdot\, & \,\cdot\, & \xi_{qn}(x)
\end{array}
\right\vert
\end{equation}
are identically zero.

Under these assumptions, it is impossible to exhibit $q$ not all
vanishing functions $\chi_1 (x), \dots, \chi_q (x)$ which make
identically equal to zero the expression $\chi_1 (x) \, X_1 f + \cdots
+ \chi_q (x) \, X_q f$. By contrast, one can determine $q ( r-q)$
functions $\varphi_{ jk} (x)$ so that the $r - q$ equations:
\[
\aligned
X_{q+j}f
=
\varphi_{j1}(x_1,\dots,x_n)\,
&
X_1f
+\cdots+
\varphi_{jq}(x_1,\dots,x_n)\,X_qf
\\
&
{\scriptstyle{(j\,=\,1\,\cdots\,r\,-\,q)}}
\endaligned
\]
are identically satisfied; indeed, every $\varphi_{ jk}$ will be equal
to a quotient whose numerator is a certain $q \times q$ determinant of
the matrix~\thetag{ 4} and whose denominator is a not identically
vanishing $q \times q$ determinant of the matrix~\thetag{ 5} (cf. the
analogous developments in Chap.~\ref{kapitel-7}, p.~\pageref{S-121}).

Now, let $x_1^0, \dots, x_n^0$ be a point in general position, or if
said more precisely, a point for which not all $q \times q$
determinants of~\thetag{ 5} vanish. 
\label{S-203}
Then the expressions $\varphi_{
jk} (x^0)$ are determined, finite constants, and at the same time, the
$r - q$ infinitesimal transformations:
\[
X_{q+j}f
-
\varphi_{j1}(x^0)\,X_1f
-\cdots-
\varphi_{jq}(x^0)\,X_qf
\ \ \ \ \ \ \ \ \ \ \ \ \
{\scriptstyle{(j\,=\,1\,\cdots\,r\,-\,q)}}
\]
belong to our group. These infinitesimal transformations are clearly
independent of each other and in addition, they possess the property
that their power series expansions with respect to the $x_i - x_i^0$
lack of all zeroth order terms. Hence according to Proposition~3,
p.~\pageref{Satz-3-S-200}, every infinitesimal transformation
of our group whose power series expansion with respect to the $x_i -
x_i^0$ \label{S-203-bis}
only contains terms of first order or of higher order must be
linearly expressible by means of the $r - q$ infinitesimal
transformations just found.

Thus, the following holds true.

\def\theproposition{7}\begin{proposition}
\label{Satz-7-S-203}
If the first $q$ of the infinitesimal transformations $X_1f, \dots,
X_r f$ of an $r$-term group are not linked by linear relations of the
form:
\[
\chi_1(x_1,\dots,x_n)\,X_1f
+\cdots+
\chi_q(x_1,\dots,x_n)\,X_qf
=
0,
\]
while $X_{ q+1} f, \dots, X_r f$ can be linearly expressed in terms of
$X_1 f, \dots, X_q f$:
\[
X_{q+j}f
\equiv
\sum_{k=1}^q\,
\varphi_{jk}(x_1,\dots,x_n)\,X_kf
\ \ \ \ \ \ \ \ \ \ \ \ \
{\scriptstyle{(j\,=\,1\,\cdots\,r\,-\,q)}},
\]
then in the neighbourhood of every point $x_1^0, \dots, x_n^0$ in
general position, the group contains exactly $q$ independent
infinitesimal transformations, for example $X_1f, \dots, X_q f$, which
are of zeroth order and out of which no infinitesimal transformation
of first order or of higher order in the $x_i - x_i^0$ can be linearly
deduced. By contrast, in the neighbourhood of $x_i^0$, the group
contains exactly $r - q$ independent infinitesimal transformations,
for instance:
\[
\label{S-203-ter}
X_{q+j}f
-
\sum_{k=1}^q\,\varphi_{jk}(x_1^0,\dots,x_n^0)\,X_kf
\ \ \ \ \ \ \ \ \ \ \ \ \
{\scriptstyle{(j\,=\,1\,\cdots\,r\,-\,q)}},
\]
which contain no terms of zeroth order in the $x_i - x_i^0$, hence
which leave at rest the point $x_1^0, \dots, x_n^0$.
\end{proposition}

We have pointed out above more precisely what is to be understood, in
this proposition, for a point in general position.

\linestop


\chapter{Determination of All Subgroups of an $r$-term Group}
\label{kapitel-12}
\chaptermark{Determination of All Subgroups of a Group}

\setcounter{footnote}{0}

\abstract*{??}

If all the transformations of a $\rho$-term group are contained in a
group with more than $\rho$ parameters, say with $r$ parameters, then
the $\rho$-term group is called a \terminology{subgroup} of the
$r$-term group.

The developments of Chap.~\ref{one-term-groups} 
already gave us an example
of subgroups of an $r$-term group; indeed, these developments showed
that every $r$-term group contains $\infty^{ r-1}$ one-term
subgroups. In the present chapter, we will at first develop a few
specific methods which enable one to find all subgroups of a given
group. Then we submit us to the question of how one should proceed in
order to determine all subgroups of a given group. About it, we
obtain the important result that the determination of all continuous
subgroups of an $r$-term group can always be achieved by resolution of
\emphasis{algebraic} equations.

\sectionengellie{\S\,\,\,54.}

In the preceding chapter, we imagined the infinitesimal
transformations of a given $r$-term group expanded with respect to the
powers of the $x_i - x_i^0$, where it is understood that $x_i^0$ is a
system of values for which all these transformations behave regularly.

For the infinitesimal transformations of the group, there resulted in
this way a classification which will now conduct us towards the
existence of certain subgroups. However, the considerations of this
paragraph find an application only to groups which in any case for
certain points $x_i^0$, contain not only infinitesimal transformations
of zeroth order, but also some of higher order in the $x_i - x_i^0$.
In the neighbourhood of the point $x_1^0, \dots, x_n^0$, let an
$r$-term group of the space $x_1, \dots, x_n$ contain exactly
$\omega_k$ independent infinitesimal transformations:
\[
Y_1f,\dots,Y_{\omega_k}f,
\]
whose power series expansions with respect to the $x_i - x_i^0$ start
with terms of order $k$ or of higher order.

We want to suppose that $k$ is $\geqslant 1$. Then if we combine one
with the other two infinitesimal transformations $Y_i f$ and $Y_j f$,
we obtain (Theorem~30, p.~\pageref{Theorem-30}) an infinitesimal
transformation $\leftbracket Y_i, \, Y_j \rightbracket$ of order
$(2k-1)$ or higher, hence at least of order $k$. Consequently,
$\leftbracket Y_i, \, Y_j \rightbracket$ must be linearly expressible
in terms of the $Yf$:
\[
\leftbracket
Y_i,\,Y_j
\rightbracket
=
\sum_{\nu=1}^{\omega_k}\,
d_{ij\nu}\,Y_\nu f,
\]
or, what is the same: the $Y_1f, \dots, Y_{ \omega_k} f$ generate an
$\omega_k$-term subgroup of the given group. Hence the following holds
true.

\def\theproposition{1}\begin{proposition}
\label{Satz-1-S-205}
If an $r$-term group of the space $x_1, \dots, x_n$ contains, in the
neighbourhood of $x_1^0, \dots, x_n^0$, exactly $\omega_k$ independent
infinitesimal transformations of order $k$ or of higher order and if
$k$ is $\geqslant 1$ here, then these $\omega_k$ infinitesimal
transformations generate an $\omega_k$-term subgroup of the group in
question.
\end{proposition}

If the point $x_1^0, \dots, x_n^0$ is procured so that, for it the
coefficients of the resolved defining equations of the group $X_1f,
\dots, X_r f$ behave regularly, then $\omega_k$ has the value:
\[
(\varepsilon_k-\nu_k)
+\cdots+
(\varepsilon_{s-1}-\nu_{s-1})
\]
(cf. Chap.~\ref{kapitel-11}, p.~\pageref{S-192}).

The case $k = 1$ is particularly important, hence we want to
yet dwell on it.

If:
\[
X_jf
=
\sum_{i=1}^n\,\xi_{ji}(x_1,\dots,x_n)\,
\frac{\partial f}{\partial x_i}
\ \ \ \ \ \ \ \ \ \ \ \ \
{\scriptstyle{(j\,=\,1\,\cdots\,r)}}
\]
are independent infinitesimal transformations of the $r$-term group,
then according to Chap.~\ref{kapitel-11}, p.~\pageref{S-199} sq., one
finds the number $\omega_1$ by examining the determinants of the
matrix:
\[
\left\vert
\begin{array}{cccc}
\xi_{11}(x^0) & \,\cdot\, & \,\cdot\, & \xi_{1n}(x^0)
\\
\cdot & \,\cdot\, & \,\cdot\, & \cdot
\\
\xi_{r1}(x^0) & \,\cdot\, & \,\cdot\, & \xi_{rn}(x^0)
\end{array}.
\right\vert
\]

Moreover, we remember (cf. p.~\pageref{S-199} sq.)
that all infinitesimal transformations of the group
which contain only terms of first order or of higher
order in the $x_i - x_i^0$ are characterized by the
fact that they leave at rest the point $x_1^0, \dots, 
x_n^0$. Hence if $k=1$, we can also enunciate the above
proposition as follows.

\def\theproposition{2}\begin{proposition}
\label{Satz-2-S-205}
If, in a group of the space $x_1, \dots, x_n$, there
are precisely $\omega_1$ independent infinitesimal transformations 
which leave invariant a determined point
$x_1^0, \dots, x_n^0$, then these transformations
generate an $\omega_1$-term subgroup of the
concerned group.
\end{proposition}

It is clear that in the variables $x_1, \dots, x_n$, there are no more
than $n$ infinitesimal transformations which are of zeroth order in
the $x_i - x_i^0$ and out of which no infinitesimal transformation of
first order or of higher order can be linearly deduced. From this, we
conclude that every $r$-term group in $n < r$ variables contains at
least $r - n$ independent infinitesimal transformations which are of
the first order or of higher order in the $x_i - x_i^0$. We therefore
have the

\def\theproposition{3}\begin{proposition}
Every $r$-term group in $n < r$ variables
contains subgroups with at least $r - n$
parameters.
\end{proposition}

From the Proposition~7 of the preceding chapter
(p.~\pageref{Satz-7-S-203}), we finally obtain for the points
$x_1^0, \dots, x_n^0$ in general position yet the

\def\theproposition{4}\begin{proposition}
If the infinitesimal transformations $X_1f, \dots, X_qf, \dots, X_r f$
of an $r$-term group in the space $x_1, \dots, x_n$ are constituted in
such a manner that $X_1f, \dots, X_qf$ are linked by no linear
relation of the form:
\[
\chi_1(x_1,\dots,x_n)\,X_1f
+\cdots+
\chi_q(x_1,\dots,x_n)\,X_qf
\equiv
0,
\]
while by contrast $X_{ q+1} f, \dots, X_r f$ can be expressed linearly
in terms of $X_1f, \dots, X_q f$:
\[
\aligned
X_{q+j}f
\equiv
\varphi_{j1}(x_1,\dots,x_n)\,
&
X_1f
+\cdots+
\varphi_{jq}(x_1,\dots,x_n)\,X_qf
\\
&
{\scriptstyle{(j\,=\,1\,\cdots\,r\,-\,q)}},
\endaligned
\]
and if in addition $x_1^0, \dots, x_n^0$ is a point in general
position, then the $r - q$ infinitesimal transformations:
\[
X_{q+j}f
-
\sum_{\mu=1}^q\,
\varphi_{j\mu}(x_1^0,\dots,x_n^0)\,X_\mu f
\ \ \ \ \ \ \ \ \ \ \ \ \
{\scriptstyle{(j\,=\,1\,\cdots\,r\,-\,q)}}
\]
are all of the first order, or of higher order, in the $x_i - x_i^0$
and they generate an $(r - q)$-term subgroup whose transformations are
characterized by the fact that they leave invariant 
\label{S-206}
the point $x_1^0,
\dots, x_n^0$.
\end{proposition}

By a point in general position,
as on p.~\pageref{S-203}, we understand here, a point 
which does not bring to zero all
$q \times q$ determinants of the matrix:
\[
\left\vert
\begin{array}{cccc}
\xi_{11}(x) & \,\cdot\, & \,\cdot\, & \xi_{1n}(x)
\\
\cdot & \,\cdot\, & \,\cdot\, & \cdot
\\
\xi_{q1}(x) & \,\cdot\, & \,\cdot\, & \xi_{qn}(x)
\end{array}
\right\vert.
\]

\sectionengellie{\S\,\,\,55.}

Here equally, we yet want to draw attention on a somewhat more general
method which often conducts, in a very simple way, to the
determination of certain subgroups of a given group. This method is
founded on the following

\def\thetheorem{31}\begin{theorem}
\label{Theorem-31-S-207}
If the $r$-term group $X_1f, \dots, X_r f$
contains some infinitesimal transformations
for which a given system of equations:
\[
\Omega_i(x_1,\dots,x_n)
=
0
\ \ \ \ \ \ \ \ \ \ \ \ \
{\scriptstyle{(i\,=\,1,\,2\,\cdots)}}
\]
remains invariant, and if every infinitesimal 
transformation of this nature can be linearly
deduced from the $m$ infinitesimal transformations:
\[
Y_kf
=
\sum_{\nu=1}^r\,h_{k\nu}\,X_\nu f
\ \ \ \ \ \ \ \ \ \ \ \ \
{\scriptstyle{(k\,=\,1\,\cdots\,m)}},
\]
then $Y_1 f, \dots, Y_m f$ generate an $m$-term subgroup of the group
$X_1f, \dots, X_r f$. 
\end{theorem}

The correctness of this theorem follows almost immediately from the
Proposition~5 in Chap.~\ref{kapitel-7},
p.~\pageref{Satz-5-S-118}. Indeed, according to this
proposition the system $\Omega_i = 0$ admits all infinitesimal
transformations of the form $\leftbracket Y_k, \, Y_j
\rightbracket$. Since the $\leftbracket Y_k, \, Y_j \rightbracket$
also belong to the group $X_1f, \dots, X_r f$, then under the
assumptions made, none of these infinitesimal transformations can be
independent of $Y_1f, \dots, Y_mf$, hence on the contrary, there must
exist relations of the form:
\[
\leftbracket
Y_k,\,Y_j
\rightbracket
=
\sum_{\mu=1}^m\,h_{kj\mu}\,Y_\mu f,
\]
that is to say, the $Y_k f$ really generate a group.

Clearly, we can also state the above theorem as follows.

\def\theproposition{5}\begin{proposition}
\label{Satz-5-S-207}
If, amongst the infinitesimal transformations of an $r$-term group in
the variables $x_1, \dots, x_n$, there are precisely $m$ independent
transformations by which a certain manifold of the space $x_1, \dots,
x_n$ remains invariant, then these $m$ infinitesimal transformations
generate an $m$-term subgroup of the $r$-term group.
\end{proposition}

The simplest case is the case where the invariant system 
of equations represents an invariant point, 
so that it has the form:
\[
x_1=x_1^0,
\,\,\,\dots,\,\,\,
x_n=x_n^0.
\]

The subgroup which corresponds to this system of equations
is of course generated by all infinitesimal transformations
whose power series expansions with respect to 
the $x_i - x_i^0$ start with terms of
first order, or of higher order. 
Thus, we arrive here at one of the subgroups
that we have already found in the preceding $\S$.

As a second example, we consider a subgroup of the eight-term general
projective group of the plane. The equation of a nondegenerate conic
section admits exactly three independent infinitesimal projective
transformations of the plane; hence these three infinitesimal
transformations generate a three-term subgroup of the general
projective group.

Lastly, yet another example, got from the ten-term group of all
conformal point transformations of the $R_3$. In this group, there
are exactly six independent infinitesimal transformations which leave
invariant an arbitrarily chosen sphere. These transformations
generate a six-term subgroup of the ten-term group.

\smallercharacters{

The Theorem~31 is only a special case of the
following more general

\def\thetheorem{32}\begin{theorem}
\label{Theorem-32-S-208}
If, in the variables $x_1, \dots, x_n$, an arbitrary group is given,
finite or infinite, continuous or not continuous, then the totality of
all transformations contained in it which leave invariant an arbitrary
system of equations in $x_1, \dots, x_n$, also forms a group.
\end{theorem}

The proof of this theorem is very simple. Any two infinitesimal
transformations of the group which, when executed one after the other, 
leave invariant the system of
equations give a transformation
which belongs again to the group and which at the same time leaves
invariant the system of equations. As a result, the proof that the
totality of transformations defined in the theorem effectively forms a
group is produced.

Instead of a system of equations, one can naturally consider also a
system of differential equations, and then the theorem would still
remain also valid.

Besides, if the given group is continuous, then the subgroup
which is defined through the invariant system of equations
can very well be discontinuous.

}

\sectionengellie{\S\,\,\,56.}

After we have so far got to know a few methods
in order to discover individual subgroups of a given group, 
we now turn to the
more general problem of determining 
\emphasis{all} continuous subgroups of a
given $r$-term group $X_1f, \dots, X_r f$.

Some $m$ arbitrary independent infinitesimal
transformations:
\[
Y_\mu f
=
\sum_{\rho=1}^r\,h_{\mu\rho}\,X_\rho f
\ \ \ \ \ \ \ \ \ \ \ \ \
{\scriptstyle{(\mu\,=\,1\,\cdots\,m)}}
\] 
of our group generate an $m$-term subgroup if and only
if all:
\[
\leftbracket
Y_\mu,\,Y_\nu
\rightbracket
=
\sum_{\rho,\,\,\sigma}^{1\cdots r}\,
h_{\mu\rho}\,h_{\nu\sigma}\,
\leftbracket
X_\rho,\,X_\sigma
\rightbracket
\]
express by means only of $Y_1f, \dots, Y_mf$.
If we insert here the values:
\[
\leftbracket
X_\rho,\,X_\sigma
\rightbracket
=
\sum_{\tau=1}^r\,
c_{\rho\sigma\tau}\,X_\tau f,
\] 
then it comes:
\[
\leftbracket
Y_\mu,\,Y_\nu
\rightbracket
=
\sum_{\rho\sigma\tau}^{1\cdots r}\,
h_{\mu\rho}\,h_{\nu\sigma}\,c_{\rho\sigma\tau}\,
X_\tau f,
\]
and it is demanded that these equations take the form:
\[
\leftbracket
Y_\mu,\,Y_\nu
\rightbracket
=
\sum_{\pi=1}^m\,l_{\mu\nu\pi}\,Y_\pi f
=
\sum_{\tau=1}^r\,\sum_{\pi=1}^m\,
l_{\mu\nu\pi}\,h_{\pi\tau}\,X_\tau f.
\]
For this to hold, it is necessary and sufficient that the
equations:
\def\theequation{1}\begin{equation}
\aligned
\sum_{\rho,\,\,\sigma}^{1\cdots r}\,
&
h_{\mu\rho}\,h_{\nu\sigma}\,c_{\rho\sigma\tau}
=
\sum_{\pi=1}^m\,l_{\mu\nu\pi}\,h_{\pi\tau}
\\
&
\ \ \ \ \ \ 
{\scriptstyle{(\mu,\,\nu\,=\,1\,\cdots\,m\,;\,\,\,
\tau\,=\,1\,\cdots\,r)}}
\endaligned
\end{equation}
can be satisfied; for this to be possible, all the $(m + 1) \times
(m+1)$ determinants of the matrix:
\def\theequation{2}\begin{equation}
\left\vert
\begin{array}{cccc}
\sum_{\rho,\,\,\sigma}^{1\cdots r}\,
h_{\mu\rho}h_{\nu\sigma}c_{\rho\sigma 1} \ \
& \ h_{11} \ & \,\cdot\, & \ h_{m1}
\\
\cdot
& \ \cdot \ & \,\cdot\, & \ \cdot
\\
\sum_{\rho,\,\,\sigma}^{1\cdots r}\,
h_{\mu\rho}h_{\nu\sigma}c_{\rho\sigma r} \ \
& \ h_{1r} \ & \,\cdot\, & \ h_{mr}
\end{array}
\right\vert
\end{equation}
must vanish. 

As a result, we have a series of algebraic equations
for the determination of the $mr$ unknowns $h_{ \pi \rho}$. 
But because $Y_1 f, \dots, Y_m f$
should be independent infinitesimal transformations,
one should from the beginning exclude
every system of values $h_{ \pi \rho}$ which
brings to zero all $m \times m$ determinants
of the matrix:
\[
\left\vert
\begin{array}{cccc}
h_{11} & \,\cdot\, & \,\cdot\, & h_{m1}
\\
\cdot & \,\cdot\, & \,\cdot\, & \cdot
\\
h_{1r} & \,\cdot\, & \,\cdot\, & h_{mr}
\end{array}
\right\vert.
\]

If the $h_{ \pi \rho}$ are determined so that all these conditions are
satisfied, then for every pair of numbers $\mu$, $\nu$, the
equations~\thetag{ 1} reduce to exactly $m$ equations that determine
completely the unknown constants $l_{ \mu\nu 1}, \dots, l_{\mu\nu
m}$. Therefore, every system of solutions $h_{ \pi \rho}$ provides an
$m$-term subgroup, and it is clear that in this way, one finds all
$m$-term subgroups.

With this, we have a general method for the determination of all
subgroups of a given $r$-term group; however in general, this method
is practically applicable only when the number $r$ is not too large;
nonetheless it shows that the problem of determining all these
subgroups necessitates only algebraic operations, which actually is
already a very important result.

\renewcommand{\thefootnote}{\fnsymbol{footnote}}
\def\thetheorem{33}\begin{theorem}
\label{Theorem-33-S-210}
The determination of all continuous subgroups of a given $r$-term
group $X_1f, \dots, X_rf$ necessitates only algebraic operations; the
concerned operations are completely determined in terms of the
constants $c_{ iks}$ in the relationships\footnote[1]{\,
\name{Lie}, Archiv for Mathematik of Naturv., Vol.~1, 
Christiania 1876.
}: 
\[
\leftbracket
X_i,\,X_k
\rightbracket
=
\sum_{s=1}^r\,c_{iks}\,X_s f
\ \ \ \ \ \ \ \ \ \ \ \ \
{\scriptstyle{(i,\,k\,=\,1\,\cdots\,r)}}.
\]
\end{theorem}
\renewcommand{\thefootnote}{\arabic{footnote}}

In specific cases, the determination of all subgroups of a given group
will often be facilitated by the fact that one knows from the
beginning certain subgroups and actually also certain properties of
the concerned group; naturally, there is also a simplification when
one has already settled the corresponding problem for one subgroup of
the given group. In addition, we shall see later that the matter is
actually not of really setting up {\em all} subgroups, but
rather, that it suffices to identify certain of these subgroups
(cf. the studies on types of subgroups, Chap.~23).

\sectionengellie{\S\,\,\,57.}

Let the $m$ independent infinitesimal transformations:
\[
Y_\mu f
=
\sum_{k=1}^r\,h_{\mu k}\,X_kf
\ \ \ \ \ \ \ \ \ \ \ \ \
{\scriptstyle{(\mu\,=\,1\,\cdots\,m)}}
\]
generate an $m$-term subgroup of the $r$-term group $X_1f, \dots,
X_rf$. Then the general infinitesimal transformation of this subgroup
is:
\[
\sum_{\mu=1}^m\,\alpha_\mu\,Y_\mu f
=
\sum_{\mu=1}^m\,\sum_{k=1}^r\,\alpha_\mu\,h_{\mu k}\,
X_kf,
\]
where the $\alpha_\mu$ denote arbitrary parameters.

As a consequence of this, all infinitesimal transformations $e_1 \,
X_1 f + \cdots + e_r\, X_r f$ of the $r$-term group which belong to
the $m$-term subgroup $Y_1 f, \dots, Y_m f$ are defined by the
equations:
\[
e_k
=
\sum_{\mu=1}^m\,\alpha_\mu\,h_{\mu k}
\ \ \ \ \ \ \ \ \ \ \ \ \
{\scriptstyle{(k\,=\,1\,\cdots\,r)}}.
\]
Hence if we imagine here that the $m$ arbitrary parameters
$\alpha_\mu$ are eliminated, we obtain exactly $r - m$ independent
linear homogeneous equations between $e_1, \dots, e_r$, so that we can
say:

\def\theproposition{6}\begin{proposition}
\label{Satz-6-S-211}
If $e_1 \, X_1 f + \cdots + e_r\, X_r f$ is the general infinitesimal
transformation of an $r$-term group, then the infinitesimal
transformations of an arbitrary $m$-term subgroup of this group can be
defined by means of $r - m$ independent linear homogeneous relations
between $e_1, \dots, e_r$.
\end{proposition}

The infinitesimal transformations which are common to two distinct
subgroups of an $r$-term group $X_1f, \dots, X_r f$ generate in turn a
subgroup; indeed, according to Chap.~\ref{kapitel-9}, Proposition~2,
p.~\pageref{Satz-2-S-159}, the infinitesimal transformations
common to the two groups generate a group, which, naturally, is
contained in $X_1f, \dots, X_r f$ as a subgroup.

Now, if we assume that one of the two groups is $m$-term, and the
other $\mu$-term, then their common infinitesimal transformations will
be defined by means of $r - m + r - \mu$ linear homogeneous equations
between the $e$, some equations which need not, however, be mutually
independent.

From this, we conclude that amongst the common infinitesimal
transformations, there are at least $r - (2r - m - \mu) = m + \mu - r$
which are independent. Therefore, we have the statement:

\def\theproposition{7}\begin{proposition}
\label{Satz-7-S-211}
If an $r$-term group contains two subgroups with $m$ and $\mu$
parameters, respectively, then these two subgroups have at least $m +
\mu - r$ independent infinitesimal transformations in common. The
infinitesimal transformations of the $r$-term group which are actually
common to the two subgroups generate in turn a subgroup.
\end{proposition}

This proposition can evidently be generalized:

Actually, if amongst the infinitesimal transformations $e_1 \, X_1 f +
\cdots + e_r \, X_r f$ of an $r$-term group, two
\emphasis{families} are sorted, the one by means of $r - m$ linear
homogeneous equations, the other by means of $r - \mu$ such equations,
then there are at least $m + \mu - r$ independent infinitesimal
transformations which are common to the two families.

\linestop


\chapter{Transitivity, Invariants, Primitivity
}
\label{kapitel-13}
\chaptermark{Transitivity, Invariants, Primitivity}

\setcounter{footnote}{0}

\abstract*{??}

\renewcommand{\thefootnote}{\fnsymbol{footnote}}
The concepts of transitivity and of primitivity which play a so broad
rôle in the theory of substitutions, shall be extended
\deutsch{ausgedehnt} here to finite continuous transformation
groups. In passing, let us mention that these concepts can
actually be extended to all groups, namely to finite groups
and to infinite groups as well, to continuous
groups and to not continuous groups as well\footnote[1]{\,
After that \name{Lie} had integrated in 1869 a few differential
equations with known continuous groups, in 1871 and in 1872, he stated
in conjunction with \name{Klein} the problem of translating the
concepts of the theory of substitutions as far as possible into the
theory of the continuous transformation groups. \name{Lie} gave the
settlement of this problem in details; on the basis of the
presentation and of the concepts exhibited through here, as early as
in 1874, he developed the fundamentals of a general theory of
integration of the complete systems which admit known infinitesimal
transformations (Verh. d. G. d. W. zu Christiania, 1874). He reduced
this problem to the case where the known infinitesimal transformations
generate a finite continuous group which is imprimitive, since its
transformations permute the characteristic manifolds of the complete
system. Amongst other things, he determined all cases where the
integration of the complete system can be performed by means of
quadratures.
}. 
\renewcommand{\thefootnote}{\arabic{footnote}}

\sectionengellie{\S\,\,\,58.}

A finite continuous group in the variables $x_1, \dots, x_n$ is called
\terminology{transitive} when, in the space $(x_1, \dots, x_n)$, there
is an $n$-times extended domain inside which each point can be
transferred to any other point by means of at least one transformation
of the group. One calls \terminology{intransitive} every group which
is not transitive.

According to this definition, an $r$-term group:
\[
x_i'
=
f_i(x_1,\dots,x_n,\,a_1,\dots,a_r)
\ \ \ \ \ \ \ \ \ \ \ \ \
{\scriptstyle{(i\,=\,1\,\cdots\,n)}}
\] 
is transitive when in general\footnote{\,
As usual (cf. Chap.~\ref{three-principles-thought}), one reasons
\emphasis{generically}, hence the concept of transitivity is
essentially considered for (sub)domains and for generic values of the
$x$ and of the $x'$.
}, 
to every system of values $x_1, \dots, x_n$, $x_1', \dots, x_n'$ at
least one system of values $a_1, \dots, a_r$ can be determined so that
the equations $x_i' = f_i ( x, \, a)$ are satisfied by the concerned
values of $x$, $a$, $x'$. In other words: \emphasis{The equations
$x_i' = f_i ( x, \, a)$ of a \emphasis{transitive group} can be
resolved with respect to $n$ of the $r$ parameters $a_1, \dots,
a_r$. If by contrast such a resolution if impossible, and if rather,
from from the equations $x_i' = f_i ( x, \, a)$ of the group, one can
derive equations which are free of the parameters $a$ and which
contain only the variables $x_1, \dots, x_n$, $x_1', \dots, x_n'$,
then the group is not transitive, it is \emphasis{intransitive}}.

From this, we see that every transitive group of the space $x_1,
\dots, x_n$ contains at least $n$ essential parameters.

If a \emphasis{transitive} group in $n$ variables has exactly $n$
essential parameters, then in general, it contains one, but only one,
transformation which transfers an arbitrary point of the space to
another arbitrary point; hence in particular, aside from the identity
transformation, it contains no transformation which leaves invariant a
point in general position. We call 
\label{S-212}
\terminology{simply transitive}
\deutsch{einfach transitiv} every group of this nature.

The general criterion for transitivity and, respectively, for
intransitivity of a group given above is practically applicable only
when one knows the finite equations of the group. But should not there
be a criterion the application of which would require only the
knowledge of the infinitesimal transformations of the group? We will
show that such a criterion can indeed be exhibited. At the same time,
we will find means in order to recognize how many and which relations
there are between the $x$ and the $x'$ only for an intransitive group
with known infinitesimal transformations.

Let us be given $r$ independent infinitesimal
transformations:
\[
X_kf
=
\sum_{i=1}^n\,\xi_{ki}(x_1,\dots,x_n)\,
\frac{\partial f}{\partial x_i}
\ \ \ \ \ \ \ \ \ \ \ \ \
{\scriptstyle{(k\,=\,1\,\cdots\,r)}}
\]
which generate an $r$-term group. This group, the finite equations of
which we can imagine written down in the form:
\[
\Phi_i
=
-\,x_i'
+
x_i
+
\sum_{k=1}^r\,e_k\,\xi_{ki}(x)
+\cdots
\ \ \ \ \ \ \ \ \ \ \ \ \
{\scriptstyle{(i\,=\,1\,\cdots\,n)}},
\]
is, according to what precedes, transitive if and only if the $n$
equations $\Phi_1 = 0$, \dots, $\Phi_n = 0$ are solvable with respect
to $n$ of the $r$ parameters $e_1, \dots, e_r$. Consequently, for the
transitivity of our group, it is necessary and sufficient that not all
$n \times n$ determinants of the matrix:
\def\theequation{1}\begin{equation}
\left\vert
\begin{array}{cccc}
\frac{\partial\Phi_1}{\partial e_1} 
& \,\cdot\, & \,\cdot\, &
\frac{\partial\Phi_n}{\partial e_1}
\\
\cdot 
& \,\cdot\, & \,\cdot\, &
\cdot
\\
\frac{\partial\Phi_1}{\partial e_r} 
& \,\cdot\, & \,\cdot\, &
\frac{\partial\Phi_n}{\partial e_r}
\end{array}
\right\vert
\end{equation}
vanish identically. From this, it follows that the group is certainly
transitive when the corresponding determinants do not all vanish for
$e_1 = 0, \dots, e_r = 0$, hence when not all $n \times n$
determinants of the matrix:
\def\theequation{2}\begin{equation}
\left\vert
\begin{array}{ccccc}
\xi_{11} & \cdot & \cdot & \cdot & \xi_{1n}
\\
\cdot & \cdot & \cdot & \cdot & \cdot
\\
\xi_{r1} & \cdot & \cdot & \cdot & \xi_{rn}
\end{array}
\right\vert
\end{equation}
are identically zero.

As a result, we have found a sufficient condition for the transitivity
of our group. We now want to study what happens when this condition
is not satisfied.

So, let all $n \times n$ determinants of the above-written
matrix~\thetag{ 2} be identically zero; in order to embrace all
possibilities in one case, we in addition want to assume that all
determinants of sizes $(n-1)$, $(n-2)$, \dots, $(q+1)$ vanish
identically, whereas not all $q \times q$ determinants do this.

Under these circumstances, it is certain that not all $q \times q$
determinants of the matrix~\thetag{ 1} vanish, hence we can conclude
that from the $n$ equations $\Phi_1 = 0$, \dots, $\Phi_n = 0$, at most
$n - q$ equations between the $x$ and $x'$ only that are free of $e_1,
\dots, e_r$ and are mutually independent can be deduced.

But now, under the assumptions made, the $r$ equations $X_1 f = 0$,
\dots, $X_r f = 0$ reduce to $q < n$ independent ones, say to: $X_1 f
= 0$, \dots, $X_q f = 0$, while $X_{ q+1} f, \dots, X_r f$ can be
expressed as follows:
\[
X_{q+\nu}f
\equiv
\varphi_{\nu 1}(x)\,X_1f
+\cdots+
\varphi_{\nu q}(x)\,X_qf
\ \ \ \ \ \ \ \ \ \ \ \ \
{\scriptstyle{(\nu\,=\,1\,\cdots\,r\,-\,q)}}.
\] 
Consequently, for all values of $i$ and $k$, one will have:
\[
\leftbracket
X_i,\,X_k
\rightbracket
=
\sum_{j=1}^q\,
\bigg(
c_{ikj}
+
\sum_{\nu=1}^{r-q}\,c_{ik,\,q+\nu}\,\varphi_{\nu j}
\bigg)\,
X_jf,
\]
that is to say: the $q$ equations $X_1 f = 0$, \dots, $X_q f = 0$ form
a $q$-term complete system with $n - q$ independent solutions, which
can be denoted by $\Omega_1 ( x), \dots, \Omega_{ n-q} ( x)$. These
solutions admit every infinitesimal transformation of the form $e_1 \,
X_1 f + \cdots + e_r \, X_r f$, hence every infinitesimal
transformation, and in consequence of that, also every finite
transformation of our $r$-term group $X_1f, \dots, X_r f$
(cf. Chap.~\ref{kapitel-6}, p.~\pageref{S-98}). Analytically, this
expresses by saying that between the variables $x$ and $x'$, which
appear in the transformation equations of our group, the following $n
- q$ equations free of the $e$ are extant:
\[
\Omega_1(x_1',\dots,x_n')=\Omega_1(x_1,\dots,x_n),
\,\,\,\,\dots\dots,\,\,\,\,
\Omega_{n-q}(x_1',\dots,x_n')=\Omega_{n-q}(x_1,\dots,x_n).
\]

Above, we said that between the $x$ and the $x'$ alone, there could
exist at most $n - q$ independent relations, and therefore we have
found all the relations in question.

In particular, we realize that the group $X_1f, \dots, X_r f$ is
intransitive as soon as all $n \times n$ determinants of the
matrix~\thetag{ 2} vanish identically. With this, it is shown that
the sufficient condition found a short while ago for the transitivity
of the group $X_1 f, \dots, X_r f$ is not only sufficient, and that it
is also necessary.

We formulate the gained results as propositions. At the head, we
state the

\def\thetheorem{34}\begin{theorem}
The $r$-term group $X_1f, \dots, X_r f$ of the space $x_1, \dots, x_n$
is transitive when amongst the $r$ equations $X_1f = 0, \dots, X_r f =
0$, exactly $n$ mutually independent ones are found, and it is
intransitive in the opposite case.
\end{theorem}

Then a proposition can follow: \label{S-215}

\def\theproposition{1}\begin{proposition}
From the finite equations:
\[
x_i'
=
f_i(x_1,\dots,x_n;\,a_1,\dots,a_r)
\ \ \ \ \ \ \ \ \ \ \ \ \
{\scriptstyle{(i\,=\,1\,\cdots\,n)}}
\]
of an $r$-term group with the infinitesimal transformations $X_1f,
\dots, X_r f$, one can eliminate the $r$ parameters $a_1, \dots, a_r$
only when the group is intransitive; in this case, one obtains between
the $x$ and the $x'$ a certain number of relations that can be brought
to the form:
\[
\Omega_k(x_1',\dots,x_n')
=
\Omega_k(x_1,\dots,x_n)
\ \ \ \ \ \ \ \ \ \ \ \ \
{\scriptstyle{(k\,=\,1,\,2\,\cdots)}}\,;
\]
here, $\Omega_1 ( x)$, $\Omega_2 ( x)$, \dots\, are a system of
independent solutions of the complete system which is determined by
the $r$ equations $X_1f = 0$, \dots, $X_r f = 0$.
\end{proposition}

In Chap.~\ref{kapitel-6}, p.~\pageref{Theorem-13-S-97}, we saw that
the solutions of the linear partial differential equations $Xf = 0$
are the only invariants of the one-term group $Xf$; accordingly, the
common solutions of the equations $X_1f = 0$, \dots, $X_r f = 0$ are
the only invariants of the group $X_1f, \dots, X_r f$. We can hence
state the proposition:

\def\theproposition{2}\begin{proposition}
If the $r$-term group $X_1f, \dots, X_r f$ of the space $x_1, \dots,
x_n$ is transitive, then it has no invariant; if it is intransitive,
then the common solutions of the equations $X_1f = 0$, \dots, $X_r f =
0$ are its only invariants.
\end{proposition}

In order to get straight \deutsch{klarstellen} the conceptual sense
\deutsch{begrifflichen Sinn} of the gained analytic results, we now
want to interpret $x_1, \dots, x_n$ as point coordinates of a space of
$n$ dimensions.

Let the complete system mentioned in the theorem be $q$-term and let
the number $q$ be smaller than $n$, so that the group $X_1f, \dots,
X_r f$ is intransitive. Let the functions $\Omega_1 ( x_1, \dots,
x_n)$, \dots, $\Omega_{ n-q} ( x_1, \dots, x_n)$ be independent
solutions of the complete system in question, and let $C_1, \dots, C_{
n-q}$ denote arbitrary constants. \emphasis{Then the equations: 
\label{S-216}
\[
\Omega_1=C_1,
\,\,\,\dots,\,\,\,
\Omega_{n-q}=C_{n-q}
\]
decompose the whole space in $\infty^{ n-q}$ different $q$-times
extended subsidiary domains \deutsch{Theilgebiete} which all remain
invariant by the group $X_1f, \dots, X_r f$}. Every point of the
space belongs to a completely determined subsidiary domain and can be
transferred only to points of the same subsidiary domain by
transformations of the group. Still, the points of a subsidiary
domain are transformed transitively, that is to say, every point in
general position in the concerned subsidiary domain can be transferred
to any other such point by means of at least one transformation of the
group.

If $x_1, \dots, x_n$ is a point in general position, then we also call
the functions $\Omega_1 ( x), \dots, \Omega_{ n-q} ( x)$ the
invariants of the point $x_1, \dots, x_n$ with respect to the group
$X_1f, \dots, X_r f$. The number of these invariants indicates the
degree of intransitivity of our group, since the larger the number of
invariants is, the smaller is the dimension number of the subsidiary
domains, inside which the point $x_1, \dots, x_n$ stays through the
transformations of the group.

Relatively to a transitive group, a point in general position clearly
has no invariant.

\medskip

The theorem stated above enables to determine whether a given group
$X_1f, \dots, X_r f$ is transitive or not. At present, one can give
various other versions of the criterion contained there for the
transitivity or the intransitivity of a group.

At first, with a light change of the way of expressing, we can say:

An $r$-term group $X_1f, \dots, X_rf$ in $x_1, \dots, x_n$ is
transitive if and only if amongst its
infinitesimal transformations, there are
exactly $n$\,\,---\,\,say $X_1f, \dots, X_n f$\,\,---\,\,which
are linked by no relation of the form:
\[
\chi_1(x_1,\dots,x_n)\,X_1f
+\cdots+
\chi_n(x_1,\dots,x_n)\,X_nf
=
0,
\]
whereas $X_{ n+1} f, \dots, X_r f$ express
as follows in terms of $X_1f, \dots, X_nf$:
\[
\aligned
X_{n+j}f
=
\varphi_{j1}(x)\,
&
X_1f
+\cdots+
\varphi_{jn}(x)\,X_nf
\\
&
{\scriptstyle{(j\,=\,1\,\cdots\,r\,-\,n)}}.
\endaligned
\]

If there are no infinitesimal transformations of this constitution,
then the group is intransitive.

Now, if we remember that every infinitesimal transformation $X_kf$
attaches a direction to each point of the space $x_1, \dots, x_n$,
and if we yet add what has been said in Chap.~\ref{kapitel-6},
p.~\pageref{S-102} about the independence of such directions which
pass through the same point, then we can state also as follows the
criterion for the transitivity of a group:

\def\theproposition{4}\begin{proposition}
\label{Satz-4-S-217}
A group $X_1f, \dots, X_r f$ in the variables $x_1, \dots, x_n$ is
transitive when it contains $n$ infinitesimal transformations which
attach to every point in general position $n$ independent
directions; if the group contains no infinitesimal transformation of
this constitution, then it is intransitive.
\end{proposition}

On the other hand, let us remember the discussions in
Chap.~\ref{kapitel-11}, p.~\pageref{S-203-bis}, where we imagined the
infinitesimal transformations of the group expanded with respect to
the powers of the $x_i - x_i^0$ in the neighbourhood of a point
$x_i^0$ in general position. Since a transitive group $X_1f, \dots,
X_r f$ in the variables $x_1, \dots, x_n$ contains $n$ infinitesimal
transformations, say: $X_1f, \dots, X_nf$ which are linked by no
relation $\chi_1 ( x) \, X_1f + \cdots + \chi_n ( x) \, X_nf = 0$, we
obtain the following proposition:

\def\theproposition{5}\begin{proposition}
An $r$-term group $X_1f, \dots, X_r f$ in the $n$ variables $x_1,
\dots, x_n$ is transitive if, in the neighbourhood of a point $x_i^0$
in general position, it contains exactly $n$ independent infinitesimal
transformations of zeroth order in the $x_i - x_i^0$ out of which no
infinitesimal transformation of first order or of higher order can be
linearly deduced. If the number of such infinitesimal transformations
of zeroth order is smaller than $n$, then the group is intransitive.
\end{proposition}

From this, one sees \emphasis{that one needs only know the defining
equations of the group $X_1f, \dots, X_r f$ in order to settle its
transitivity or its intransitivity}.

Finally, we can also state as follows the first part
of the Proposition~5:

The group $X_1f, \dots, X_r f$ is transitive when in the neighbourhood
of a point $x_i^0$ in general position, it contains exactly $r - n$
independent infinitesimal transformations whose power series
expansions with respect to the $x_i - x_i^0$ start with terms of of
first or of higher order, hence when the group contains exactly $r -
n$, and not more, independent infinitesimal transformations which
leave invariant a point in general position. ---

\medskip

If one only knows the infinitesimal transformations $X_1f, \dots, X_r
f$ of an intransitive group, then as we have seen, one finds the
invariants of the group by integrating the complete system $X_1f = 0,
\dots, X_r f = 0$. Now, it is of great importance that this
integration is not necessary when the finite equations of the group
are known, so that in this case, the invariants of the group are
rather found by plain elimination.

In order to prove this, we imagine that the finite equations of an
intransitive group are given:
\[
x_i'
=
f_i(x_1,\dots,x_n\,a_1,\dots,a_r)
\ \ \ \ \ \ \ \ \ \ \ \ \
{\scriptstyle{(i\,=\,1\,\cdots\,n)}}, 
\]
and then we eliminate the parameters $a_1, \dots, a_r$ from
them. According to what precedes, it must be possible to bring the $n
- q$ independent equations obtained in this way:
\def\theequation{3}\begin{equation}
W_\mu(x_1,\dots,x_n,\,x_1',\dots,x_n')
=
0
\ \ \ \ \ \ \ \ \ \ \ \ \
{\scriptstyle{(\mu\,=\,1\,\cdots\,n\,-\,q)}},
\end{equation}
to the form:
\def\theequation{4}\begin{equation}
\Omega_\mu(x_1',\dots,x_n')
=
\Omega_\mu(x_1,\dots,x_n)
\ \ \ \ \ \ \ \ \ \ \ \ \
{\scriptstyle{(\mu\,=\,1\,\cdots\,n\,-\,q)}},
\end{equation}
where the $\Omega_\mu ( x)$ are the sought invariants. Hence, when we
solve the equations~\thetag{ 3} with respect to $n - q$ of the
variables $x_1', \dots, x_n'$:
\[
x_\mu'
=
\Pi_\mu(x_1,\dots,x_n,\,x_{n-q+1}',\dots,x_n')
\ \ \ \ \ \ \ \ \ \ \ \ \
{\scriptstyle{(\mu\,=\,1\,\cdots\,n\,-\,q)}},
\]
which is always possible, then we obtain $n - q$ functions $\Pi_1,
\dots, \Pi_{ n-q}$ in which the variables $x_1, \dots, x_n$ occur only
the combinations $\Omega_1 ( x), \dots, \Omega_{ n-q} (x)$.
Consequently, the $n - q$ expressions:
\[
\Pi_\mu(x_1,\dots,x_n,\,\alpha_{n-q+1},\dots,\alpha_n)
\ \ \ \ \ \ \ \ \ \ \ \ \
{\scriptstyle{(\mu\,=\,1\,\cdots\,n\,-\,q)}}
\]
in which $\alpha_{ n-q+1}, \dots, \alpha_n$ denote
constants, represent invariants of our group, and in fact
clearly, $n - q$ independent invariants. ---

The following therefore holds true.

\def\thetheorem{35}\begin{theorem}
If one knows the finite transformations:
\[
x_i'
=
f_i(x_1,\dots,x_n,\,a_1,\dots,a_r)
\ \ \ \ \ \ \ \ \ \ \ \ \
{\scriptstyle{(i\,=\,1\,\cdots\,n)}}
\]
of an intransitive group, then one can find the invariants of this
group by means of elimination.
\end{theorem}

\sectionengellie{\S\,\,\,59.}

When we studied how a point in general position behaves relatively to
the transformations of an $r$-term group, we were conducted with
necessity to the concepts of transitivity and of intransitivity
\deutsch{wurden wir mit Nothwendigkeit auf die Begriffe Transitivität
und Intransitivität geführt}. We obtained in this way a division
\deutsch{Eintheilung} of all $r$-term groups of a space of $n$
dimensions in two different classes, exactly the same way as in the
theory of substitutions; but at the same time, we yet obtained a
division of the intransitive groups too, namely according to the
number of the invariants that a point in general position possesses
with respect to the concerned group.

Correspondingly to the process of the theory of substitutions, we now
can also go further and study the behaviour of two or more points in
general position relatively to an $r$-term group. This gives us a new
classification of the groups of the $R_n$.

Let:
\[
y_i
=
f_i(x_1,\dots,x_n,\,a_1,\dots,a_r)
\ \ \ \ \ \ \ \ \ \ \ \ \
{\scriptstyle{(i\,=\,1\,\cdots\,n)}}
\]
be an $r$-term group and let $X_1f, \dots, X_r f$ be $r$ independent
infinitesimal transformations of it.

At first, we want to consider two points $x_1', \dots, x_n'$ and
$x_1'', \dots, x_n''$ and to seek their invariants relatively to our
group, that is to say: we seek all functions of $x_1', \dots, x_n'$,
$x_1'', \dots, x_n''$ which remain invariant by the transformations of
our group.

To this end, we write the infinitesimal transformations $X_kf$ once in
the $x'$ under the form $X_k'f$ and once in the $x''$ under the form
$X_k'' f$; then simply, the sought invariants are the invariants of
the $r$-term group:
\def\theequation{5}\begin{equation}
X_k'f+X_k''f
\ \ \ \ \ \ \ \ \ \ \ \ \
{\scriptstyle{(k\,=\,1\,\cdots\,r)}}
\end{equation}
in the variables $x_1', \dots, x_n'$, $x_1'', \dots, x_n''$.

If $J_1 ( x), \dots, J_{ \rho_1} (x)$ are the invariants of the group
$X_1f, \dots, X_rf$, then without effort, the $2 \rho_1$ functions:
\[
J_1(x'),\dots,J_{\rho_1}(x'),\ \
J_1(x''),\dots,J_{\rho_1}(x'')
\]
are invariants, and in fact, independent invariants of the
group~\thetag{ 5}; but in addition, there can yet be a certain number,
say $\rho_2$, of invariants:
\[
J_1'(x_1',\dots,x_n',\,x_1'',\dots,x_n''),
\,\,\,\dots,\,\,\,
J_{\rho_2}'(x_1',\dots,x_n',\,x_1'',\dots,x_n'')
\]
which are mutually independent and are independent of the $2\rho_1$
above invariants. So in this case, two points in general position
have $2 \rho_1 + \rho_2$ independent invariants relatively to the
group $X_1f, \dots, X_r f$, amongst which however, only $\rho_2$ have
to be considered as essential, because each one of the two points
already has $\rho_1$ invariants for itself. Under these assumptions,
from the equations:
\[
\aligned
y_i'
=
f_i(x_1',\dots,x_n',\,a_1,\dots,
&
a_r),
\ \ \ \ \ \ \
y_i''
=
f_i(x_1'',\dots,x_n'',\,a_1,\dots,a_r)
\\
&
{\scriptstyle{(i\,=\,1\,\cdots\,n)}}
\endaligned
\]
there result the following relations free of the $a$:
\[
\aligned
J_k(y')
&
=
J_k(x'),
\ \ \ \ \ \
J_k(y'')
=
J_k(x'')
\ \ \ \ \ \ \ \ \ \ \ \ \ {\scriptstyle{(k\,=\,1\,\cdots\,\rho_1)}}
\\
&
J_j'(y',\,y'')
=
J_j'(x',\,x'')
\ \ \ \ \ \ \ \ \ \ \ \ \ {\scriptstyle{(j\,=\,1\,\cdots\,\rho_2)}}.
\endaligned
\]

Therefore, if we imagine that the quantities $x_1', \dots, x_n'$ and
$y_1', \dots, y_n'$ are chosen fixed, then the totality of all
positions $y_1'', \dots, y_n''$ which the point $x_1'', \dots, x_n''$
can take are determined by the equations:
\[
\aligned
J_k(y'')
=
&
J_k(x''),
\ \ \ \ \ \
J_j'(y',y'')
=
J_j'(x',x'')
\\
&\
{\scriptstyle{(k\,=\,1\,\cdots\,\rho_1\,;\,\,\,
j\,=\,1\,\cdots\,\rho_2)}}\,;
\endaligned
\]
so there are $\infty^{ n - \rho_1 - \rho_2}$ distinct positions of
this sort.

In a similar way, one can determine the invariants that three, four
and more points have relatively to the group. \emphasis{So one finds
a series of entire number $\rho_1$, $\rho_2$, $\rho_3$, \dots\, which
are characteristic of the group and which are independent of the
choice of variables}. If one computes these numbers one after the
other, one always comes to a number $\rho_m$ which vanishes, while at
the same time all numbers $\rho_{ m+1}$, $\rho_{ m+2}$, \dots\, are
equal to zero.

We do not want to address further the issue about these behaviours,
but it must be observed that analogous considerations can be
endeavoured for every family of $\infty^r$ transformations, shall this
family constitute a group or not.

\sectionengellie{\S\,\,\,60.}

Above, we have seen that an intransitive group $X_1f, \dots, X_r f$
decomposes the entire space $(x_1, \dots, x_n)$ in a continuous family
of $q$-times extended manifolds of points:
\[
\Omega_1(x_1,\dots,x_n)
=
C_1,
\,\,\,\dots,\,\,\,
\Omega_{n-q}(x_1,\dots,x_n)
=
C_{n-q}
\]
which remain all invariant by the transformations of the group. Here,
the $\Omega$ denote independent solutions of the $q$-term complete
system which is determined by the equations $X_1f = 0, \dots, X_r f =
0$.

Each point of the space belongs to one and to only one of the
$\infty^{ n-q}$ manifolds $\Omega_1 = a_1$, \dots, $\Omega_{ n-q} =
a_{ n-q}$, so in the sense provided by Chap.~6, \pageref{S-101-bis},
we are dealing with a decomposition of the space. This decomposition
remains invariant by all transformations of the group $X_1f, \dots,
X_r f$; and at the same time, each one of the individual subsidiary
domains in which the space is decomposed stays invariant.

\label{S-220-sq}
It can also happen for transitive groups that there exists a
decomposition of the space in $\infty^{ n-q}$ $q$-times extended
manifolds $\Omega_1 = a_1$, \dots, $\Omega_{ n-q} = a_{ n-q}$ which
remain invariant by the group. But naturally, each individual
manifold amongst the $\infty^{ n-q}$ manifolds need not remain
invariant, since otherwise the group would be intransitive; these
$\infty^{ n-q}$ manifolds must rather be permuted by the group, while
the totality of them remains invariant.

Now, a group of the $R_n$ is called \terminology{imprimitive}
\deutsch{imprimitiv} when it determines at least one invariant
decomposition of the space by $\infty^{ n-q}$ $q$-times extended
manifolds; a group for which there appears absolutely no invariant
decomposition is called \terminology{primitive} \deutsch{primitiv}.
That the values $q = 0$ and $q = n$ are excluded requires hardly any
mention here.

The intransitivity is obviously a special case of imprimitivity: every
intransitive group is at the same time also imprimitive. On the other
hand, every primitive group is necessarily transitive.

Now, in order to obtain an analytic definition for the imprimitivity
of an $r$-term group $X_1f, \dots, X_r f$, we need only to remember
that every decomposition of the space in $\infty^{ n-q}$ $q$-times
extended manifolds $y_1 = {\rm const.}$, \dots, $y_{ n-q} = {\rm
const.}$ is analytically defined by the $q$-term complete system $Y_1
f = 0$, \dots, $Y_q f = 0$, the solutions of which are the $y_k$. The
fact that the concerned decomposition remains invariant by the group
$X_1f, \dots, X_r f$ amounts to the fact that the corresponding
$q$-term complete system admits all transformations of the group.

The $q$-term complete system $Y_k f = 0$ admits our group as soon it
admits the general one-term group $\lambda_1 \, X_1 f + \cdots +
\lambda_r \, X_r f$. According to Theorem~20, Chap.~\ref{kapitel-8},
p.~\pageref{Theorem-20-S-140}, this is the case when, between the
$X_if$ and the $Y_kf$, there exist relationships of the following
form:
\[\label{S-221-bis}
\leftbracket
X_i,\,Y_k
\rightbracket
=
\sum_{\nu=1}^q\,\psi_{ik\nu}(x)\,Y_\nu f.
\]

Consequently, it is a necessary and sufficient condition for the
imprimitivity \label{S-221}
of the group $X_1f, \dots, X_r f$ that there exists a
$q$-term complete system:
\[
Y_1f=0,
\,\,\,\dots,\,\,\,
Y_qf=0
\ \ \ \ \ \ \ \ \ \ \ \ \ {\scriptstyle{(q\,<\,n)}}
\]
which stands in such relationships with the $X_kf$.

Besides, the group $X_1f, \dots, X_rf$ can also be imprimitive in
several manners, that is to say, there can exist many, and even
infinitely many systems that the group admits.

Later, we will develop a method for the setting up of all complete
systems which remain invariant by a given group. Hence in particular,
we will also be able to determine whether the group in question is
primitive or not. Naturally, the latter question requires a special
examination only for transitive groups.

Now, yet a brief remark.

Let the equations:\label{S-222}
\[
y_1
=
C_1,
\,\,\,\dots,\,\,\,
y_{n-q}
=
C_{n-q}
\]
represent a decomposition of the space $x_1, \dots, x_n$ which is
invariant by the group $X_1f, \dots, X_r f$. Then if we introduce
$y_1, \dots, y_{ n-q}$ together with $q$ other appropriate functions
$z_1, \dots, z_q$ of $x_1, \dots, x_n$ as new independent variables,
the infinitesimal transformations $X_k f$ receive the specific form:
\[
\aligned
X_kf
=
\sum_{\mu=1}^{n-q}\,
\omega_{k\mu}(y_1,\dots,y_{n-q})\,
&
\frac{\partial f}{\partial y_\mu}
+
\sum_{j=1}^q\,\zeta_{kj}(y_1,\dots,y_{n-q},\,z_1,\dots,z_q)\,
\frac{\partial f}{\partial z_j}
\\
&
\ \ {\scriptstyle{(k\,=\,1\,\cdots\,r)}}.
\endaligned
\]

\label{S-222-bis}
Now, according to Chap.~\ref{kapitel-12},
p.~\pageref{Satz-5-S-207}, all infinitesimal transformations
$e_1 \, X_1 f + \cdots + e_r\, X_r f$ which leave invariant a
determined manifold $y_1 = y_1^0$, \dots, $y_{ n-q} = y_{ n-q}^0$
generate a subgroup. If we want to find the infinitesimal
transformations in question, then we only have to seek all systems of
values $e_1, \dots, e_r$ which satisfy the $n - q$ equations:
\[
\sum_{k=1}^r\,e_k\,\omega_{k\mu}(y_1^0,\dots,y_{n-q}^0)
=
0
\ \ \ \ \ \ \ \ \ \ \ \ \ 
{\scriptstyle{(\mu\,=\,1\,\cdots\,n\,-\,q)}}.
\]
If $r > n-q$, then there always are systems of values $e_1, \dots,
e_r$ of this constitution, and consequently in this case, the group
$X_1f, \dots, X_r f$ certainly contains subgroups with a least $r - n
+ q$ parameters.

Besides, it is clear that the $r$ \terminology{reduced}
\deutsch{verkürzten} infinitesimal transformations:
\[
\overline{X}_kf
=
\sum_{\mu=1}^{n-q}\,
\omega_{k\mu}(y_1,\dots,y_{n-q})\,
\frac{\partial f}{\partial y_\mu}
\ \ \ \ \ \ \ \ \ \ \ \ \ {\scriptstyle{(k\,=\,1\,\cdots\,r)}}
\]
in the $n-q$ variables $y_1, \dots, y_{ n-q}$ generate a group;
however, this group possibly contains not $r$, but only a smaller
number of essential parameters. Clearly, the calculations indicated
just now amount to the determination of all infinitesimal
transformations $e_1 \, \overline{ X}_1 f + \cdots + e_r\, \overline{
X}_r f$ which leave invariant the point $y_1^0, \dots, y_{ n-q}^0$ of
the $(n-q)$-times extended space $y_1, \dots, y_{ n-q}$.

\linestop


\chapter{Determination of All Systems of Equations
\\
Which Admit a Given $r$-term Group
}
\label{kapitel-14}
\chaptermark{Determination of Invariant Systems of Equations}

\setcounter{footnote}{0}

\abstract*{??}

If a system of equations remains invariant by all transformations of
an $r$-term group $X_1f, \dots, X_rf$, we say that it admits the group
in question. Every system of equations of this constitution admits
all transformations of the general one-term group $\sum \, e_k \,
X_kf$, and therefore specially, all $\infty^{ r-1}$
\emphasis{infinitesimal} transformations $\sum\, e_k \, X_kf$ of the
$r$-term group.

Now on the other hand, we have shown earlier on that every system of
equations which admits the $r$ \emphasis{infinitesimal}
transformations $X_1f, \dots, X_rf$ and therefore also, all $\infty^{
r-1}$ infinitesimal transformations $\sum \, e_k \, X_kf$, allows at
the same time all \emphasis{finite} transformations of the one-term
group $\sum\, e_k \, X_kf$, that is to say, all transformations of the
group $X_1f, \dots, X_rf$ (cf. Theorem~14,
p.~\pageref{Theorem-14-S-112}). Hence if all systems of equations
which admit the $r$-term group $X_1f, \dots, X_rf$ 
are to be
determined, this shall be a problem which is completely settled by the
developments of the Chap.~\ref{kapitel-7}. 
\label{S-223}
Indeed, the problem of
setting up all systems of equations which admit $r$ given
infinitesimal transformations is solved in complete generality there.

However, the circumstance where the $X_1f, \dots, X_rf$ which are
considered here generate an $r$-term group means a really major
simplification in comparison to the general case. Hence it appears to
be completely legitimate that we settle independently the special case
where the $X_k f$ generate a group.

The treatment of the addressed problem turns out to be not
inessentially different, whether or not one also knows the finite
equations of the concerned group. In the first case, no integration is
required. But in the second case, one does not make it in general
without integration; however, some operations which were necessary for
the general problem of the chapter~\ref{kapitel-7} drop.

\renewcommand{\thefootnote}{\fnsymbol{footnote}}
We shall treat these two cases one after the other, above all because
of the applications, but also in order to afford a deeper insight into
the matter; in fact, the concerned developments complement one another
mutually\footnote[1]{\,
\name{Lie}, Math. Ann. Vol. XI, pp.~510--512, Vol.~XVI, p.~476.
Archiv for Math. og Nat., Christiania 1878, 1882, 1883. Math. Ann.
Vol. XXIV.
}. 
\renewcommand{\thefootnote}{\arabic{footnote}}

Lastly, let us yet mention that from now on, we shall frequently
translate the common symbolism of the theory of substitutions into the
theory of the transformation groups. So for example, we denote by
$S$, $T$, \dots\, individual transformations, and by $S^{ -1}$, $T^{
-1}$, \dots\, the corresponding inverse transformations. By $S\,T$,
we understand the transformation which is obtained when the
transformation $S$ is executed first, and then the transformation $T$.
From this, it follows that expressions of
the form $S S^{ -1}$, $T T^{ -1}$
mean the identity transformation.

\sectionengellie{\S\,\,\,61.}

We consider an arbitrary point $P$ of the space. The totality of all
positions that this point takes by the $\infty^r$ transformations of
the group forms a certain manifold $M$; we shall show that this
manifold admits the group, or in other words, that every point of $M$
is transferred to a point again of $M$ by every transformation of the
group.

Indeed, let $P'$ be any point of $M$, and let $P'$ come from $P$ by
the transformation $S$ of our group, what we want to express by means
of the symbolic equation:
\[
(P)\,S
=
(P').
\]
Next, if $T$ is a completely arbitrary transformation of the group,
then by the execution of $T$, $P'$ is transferred to:
\[
(P')\,T
=
(P)\,S\,T\,;
\]
but since the transformation $S \, T$ belongs to the group as well,
$(P) \, S \, T$ is also a point of $M$, and our claim is therefore
proved.

Obviously, every manifold invariant by the group which contains the
point $P$ must at the same time contain the manifold $M$. That is why
we can also say: $M$ is the \emphasis{smallest} manifold invariant by
the group to which the point $P$ belongs.

But still, there is something more. It can be shown that with the help
of transformations of the group, every point of $M$ can be transferred
to any other point of this manifold. Indeed, if $P'$ and $P''$ are
any two points of $M$ and if they are obtained from $P$ by means of
the transformations $S$ and $U$, respectively, one has the
relations:
\[
(P)\,S
=
(P'),
\ \ \ \ \ \ \ \ \
(P)\,U
=
(P'')\,;
\]
from the first one, it follows:
\[
(P')\,S^{-1}
=
(P)\,S\,S^{-1}
=
(P)\,;
\]
hence with the help of the second one, it comes:
\[
(P')\,S^{-1}\,U
=
(P''),
\]
that is to say, by the transformation $S^{ -1} U$ which likewise
belongs to the group, the point $P'$ is transferred to the point
$P''$. As a result, the assertion stated above is proved.

From this, we realize that the manifold $M$ can also be defined as the
totality of all positions which any of its other points, not just $P$,
take by the $\infty^r$ transformations of the group.

Consequently, the following holds true. 

\def\thetheorem{36}\begin{theorem}
\label{Theorem-36-S-225}
If, on a point $P$ of the space $(x_1, \dots, x_n)$, one executes all
$\infty^r$ transformations of an $r$-term group of this space, then
the totality of all positions that the point takes in this manner
forms a manifold invariant by the group; this manifold contains no
smaller subsidiary domain invariant by the group, and it is itself
contained in all invariant manifolds in which the point $P$ lies.
\end{theorem}
 
If one assumes that the finite equations $x_i' = f_i ( x_1, \dots,
x_n, \, a_1, \dots, a_r)$ of the $r$-term group are known, then
without difficulty, \label{S-225-ter}
one can indicate for every point $x_1^0, \dots,
x_n^0$ the smallest invariant manifold to which it belongs. Indeed by
the above, the manifold in question consists of the totality of all
positions $x_1, \dots, x_n$ that the point $x_1^0, \dots, x_n^0$ take
by the transformations of the group; but evidently, the totality of
these positions is represented by the $n$ equations:
\[
x_i
=
f_i(x_1^0,\dots,x_n^0,\,a_1,\dots,a_r)
\ \ \ \ \ \ \ \ \ \ \ \ \ {\scriptstyle{(i\,=\,1\,\cdots\,n)}}
\]
in which the parameters $a$ are to be interpreted as independent
variables. If one eliminates the $a$, one obtains the sought manifold
represented by equations between the $x$ alone.

\smallercharacters{

Here, it is to be recalled that in the equations $x_i = f_i ( x^0,
a)$, the $a$ are not completely arbitrary; indeed, all systems of
values $a_1, \dots, a_r$ for which the determinant:
\[
\label{S-225-bis}
\sum\,\pm\,
\bigg[
\frac{\partial f_1(x,a)}{\partial x_1}
\bigg]_{x=x^0}
\,\,\,\,\cdots\,\,\,\,\,\,
\bigg[
\frac{\partial f_n(x,a)}{\partial x_n}
\bigg]_{x=x^0}
\]
vanish are excluded from the beginning, because we always use only
transformations which are solvable. From this, it follows that in
certain circumstances, one obtains, by elimination of the $a$, a
manifold which contains, aside from the points to which $x_1^0, \dots,
x_n^0$ is transferred by the \emphasis{solvable} transformations of
the group, yet other points; then as one easily sees, the latter
points form in turn an invariant manifold.

}

Furthermore, it is to be remarked that the discussed elimination can
take different shapes for different systems of values $x_k^0$; indeed,
the elimination of the $a$ need not conduct always to the same number
of relations between the $x$, which again means that the smallest
invariant manifolds in question need not have all the same dimension
number.

If, amongst all smallest invariant manifolds of the same dimension
number, one takes infinitely many such invariant manifolds according
to an arbitrary analytic rule, then their totality also forms an
invariant manifold. In this way, all invariant manifolds can
obviously be obtained.

So we have the

\def\thetheorem{37}\begin{theorem}
\label{Theorem-37-S-226}
If one knows the finite equations of an $r$-term group 
$X_1f, \dots, X_rf$ in the variables $x_1, \dots, x_n$, 
then without integration, one can find all invariant
systems of equations invariant by the group, or, what is the
same, all manifolds invariant by it. 
\end{theorem}

\sectionengellie{\S\,\,\,62.}

In Chap.~\ref{kapitel-12}, p.~\pageref{S-206}, we have seen that an
$r$-term group $G_r$ in $x_1, \dots, x_n$ associates to every point of
the space a completely determined subgroup, namely the subgroup which
consists of all transformations of the $G_r$ which leave the point
invariant.

Let the point $P$ be invariant by an $m$-term subgroup of the $G_r$,
but by no subgroup with more terms; let $S$ be the general symbol of a
transformation of this $m$-term subgroup, so that one hence has: $(P)
\, S = (P)$. Moreover, let $T$ be a transformation which transfers
the point $P$ to the new position $P'$:
\[
(P)\,T
=
P'.
\]

Now, if ${\sf T}$ is an arbitrary transformation of the
$G_r$ which transfers $P$ to $P'$ as well, we have:
\[
(P)\,{\sf T}
=
(P')
=
(P)\,T,
\]
hence it comes:
\[
(P)\,{\sf T}\,T^{-1}
=
(P).
\]

From this, it is evident that ${\sf T} T^{ -1}$
belongs to the transformations $S$, hence that: 
\[
{\sf T}
=
S\,T
\]
is the general form of a transformation of the same constitution as
${\sf T}$. Now, since there are precisely $\infty^m$ transformations
$S$, we see that:

\plainstatement{The $G_r$ contains exactly $\infty^m$ transformations
which transfer $P$ to $P'$.}

On the other hand, if we ask for all transformations $S'$ of the $G_r$
which leave invariant the point $P'$, then we have to fulfill the
condition $(P') \, S' = (P')$. From it, we see that:
\[
(P)\,T\,S'
=
(P)\,T
\ \ \ \ \ \ \
\text{\rm and}
\ \ \ \ \ \ \
(P)\,T\,S'\,T^{-1}
=
(P),
\]
and consequently $T\, S' \, T^{ -1}$ is a transformation $S$,
that is to say $S'$ has the form:
\[
S'
=
T^{-1}\,S\,T.
\]

One sees with easiness here that $S$ can be a completely arbitrary
transformation of the subgroup associated to the point $P$; as a
result, our group contains exactly $\infty^m$ different
transformations $S'$, and in turn now, they obviously form an $m$-term
subgroup of the $G_r$.

The results of this paragraph obtained up to now 
can be summarized as follows. 

\def\theproposition{1}\begin{proposition}
\label{Satz-1-S-227}
If an $r$-term group $G_r$ of the $R_n$ contains exactly $\infty^m$
and not more transformations $S$ which leave the point $P$ invariant,
and if in addition it contains at least one transformation $T$ which
transfers the point $P$ to the point $P'$, then it contains on the
whole $\infty^m$ different transformations which transfer $P$ to $P'$;
the general form of these transformations is: $S\, T$. In addition,
the $G_r$ contains exactly $\infty^m$ transformations which leave
invariant the point $P'$; their general form is: $T^{-1} \, S \, T$.
\end{proposition}

From the second part of this proposition, it follows that the points
which admit exactly $\infty^m$ transformations of our group are
permuted by the transformations of the group, while their totality
remains invariant.

Hence the following hold true.

\def\thetheorem{38}\begin{theorem}
\label{Theorem-38-S-227}
The totality of all points which admit the same number, say
$\infty^m$, and not more transformations of an $r$-term group, remains
invariant by all transformations of the group.
\end{theorem}

We have proved this theorem by applying considerations which are 
borrowed 
from the theory of substitutions. But at the same time, we want to
show how one can conduct the proof in case one abstains from such
considerations, or from a more exact language.

A point $x_1^0, \dots, x_n^0$ which allows
$\infty^m$ transformations of the $r$-term group:
\[
X_kf
=
\sum_{i=1}^n\,\xi_{ki}(x_1,\dots,x_n)\,
\frac{\partial f}{\partial x_i}
\ \ \ \ \ \ \ \ \ \ \ \ \ {\scriptstyle{(k\,=\,1\,\cdots\,r)}}
\]
admits precisely $m$ independent infinitesimal transformations of this
group. Hence the group contains, in the neighbourhood of $x_1^0,
\dots, x_n^0$, exactly $m$ independent infinitesimal transformations
whose power series expansions with respect to the $x_i - x_i^0$ begin
with terms of first order or of higher order. Now, if we imagine that
new variables:
\[
\aligned
\overline{x}_i
=
\overline{x}_i^0
+
&
\sum_{k=1}^n\,\alpha_{ki}\,(x_k-x_k^0)
+\cdots
\ \ \ \ \ \ \ \ \ \ \ \ \ {\scriptstyle{(i\,=\,1\,\cdots\,n)}}
\\
&
{\textstyle{\sum}}\,\pm\,\alpha_{11}\,\cdots\,\alpha_{nn}\neq 0,
\endaligned
\]
are introduced in the group, then according to Chap.~\ref{kapitel-11},
p.~\pageref{S-197}, we obtain a new group in the $\overline{ x}_i$
which, in the neighbourhood of $\overline{ x}_i^0$, contains in the
same way exactly $m$ independent infinitesimal transformations of
first order or of higher order. In particular, if we imagine that the
transition from the $x_i$ to the $\overline{ x}_i$ is a transformation
of the group $X_1f, \dots, X_r f$, then the group in the $\overline{
x}_i$ is simply identical to the group:
\[
\overline{X}_kf
=
\sum_{i=1}^n\,\xi_{ki}(\overline{x}_1,\dots,\overline{x}_n)\,
\frac{\partial f}{\partial\overline{x}_i}
\ \ \ \ \ \ \ \ \ \ \ \ \ {\scriptstyle{(k\,=\,1\,\cdots\,r)}}
\]
(cf. Chap.~\ref{fundamental-differential}, 
p.~\pageref{symmetric-X-X-prime}). 
In other words, if, by a transformation of our group, the point
$x_i^0$ is transferred to the point $\overline{ x}_i^0$, then this
point also admits precisely $m$ independent infinitesimal
transformations of the group. But with this, the Theorem~38 is
visibly proved.

From the Theorem~38 it comes immediately that the following
proposition also holds true:

\def\theproposition{2}\begin{proposition}
The totality of all points $x_1, \dots, x_n$ which admit $m$ or more
independent infinitesimal transformations of the $r$-term group $X_1f,
\dots, X_rf$ remains invariant by this group.
\end{proposition}

If we associate this proposition with the developments in
Chap.~\ref{kapitel-11}, Proposition~4,
p.~\pageref{Satz-4-S-200}, we obtain a new important result.
At that time, we indeed saw that the points $x_1, \dots, x_n$ which
admit $m$ or more independent infinitesimal transformations $e_1 \,
X_1f + \cdots + e_r \, X_rf$ of the $r$-term group $X_1f, \dots, X_rf$
are characterized by the fact that all $(r-m+1) \times (r - m+1)$
determinants of a certain matrix are brought to vanishing. At
present, we recognize that the system of equations which is obtained
by equating to zero these $(r - m + 1) \times (r - m + 1)$
determinants admits all transformations of the group $X_1f, \dots, X_r
f$. As a result, we have the

\def\thetheorem{39}\begin{theorem}
\label{Theorem-39-S-228}
If $r$ independent infinitesimal transformations:
\[
X_kf
=
\sum_{i=1}^n\,\xi_{ki}(x_1,\dots,x_n)\,
\frac{\partial f}{\partial x_i}
\ \ \ \ \ \ \ \ \ \ \ \ \ {\scriptstyle{(k\,=\,1\,\cdots\,r)}}
\]
generate an $r$-term group, then by equating to zero all $(r - m + 1)
\times (r - m + 1)$ determinants of the matrix:
\def\theequation{1}\begin{equation}
\left\vert
\begin{array}{cccc}
\xi_{11}(x) & \,\cdot\, & \,\cdot\, & \xi_{1n}
\\
\cdot & \,\cdot\, & \,\cdot\, & \cdot
\\
\xi_{r1}(x) & \,\cdot\, & \,\cdot\, & \xi_{rn}
\end{array}
\right\vert,
\end{equation}
one always obtains a system of equations which admits all
transformations of the group $X_1f, \dots, X_rf$; this holds true for
every number $m \leqslant r$, provided only that there actually exist
systems of values $x_1, \dots, x_n$ which bring to zero all the $(r -
m + 1) \times (r - m + 1)$ determinants in question.
\end{theorem}

Later in the course of this chapter (\S\S\,\,66 and~67,
p.~\pageref{S-66} and~\pageref{S-67} resp.), we will give yet two
different purely analytic proofs \label{S-229}
of the above important theorem.
Temporarily, we observe only the following:

The Theorem~39 shows that there is an essential difference 
\label{S-229-bis}
between the
problem of the Chap.~\ref{kapitel-7} p.~\pageref{S-123} up
to~\pageref{S-133} and the one of the present chapter.

If $X_1f, \dots, X_rf$ generate an $r$-term group, then by equating to
zero all $(r - m + 1) \times (r - m + 1)$ determinants of the
matrix~\thetag{ 1}, one always obtains an invariant system of
equations, only as soon as all these determinants really can vanish at
the same time. But this is not anymore true when it is only assumed
about the $X_kf$ that, when set to zero, they constitute a complete
system consisting of of $r$, or less equations. In this case, it is
certainly possible that there are invariant systems of equations which
embrace the equations obtained by equating to zero the determinants in
question, but it is not at all always the case that one obtains an
invariant system of equations by equating to zero these determinants,
just like that. To get this, further operations are rather necessary
in general, as it is explained in Chap.~\ref{kapitel-7},
p.~\pageref{S-124-sq} sq.

In the last paragraph of this chapter, p.~\pageref{S-243-sq}, we will
study this point in more details.

\sectionengellie{\S\,\,\,63.}

In the preceding paragraph, we have seen that, \emphasis{from the
infinitesimal transformations} of an $r$-term group, one can derive
\emphasis{without integration} certain systems of equations which
remain invariant by the concerned group. Now, we will show how one
finds \emphasis{all} systems of equations which admit an $r$-term
group \emphasis{with given infinitesimal transformations}:
\[
X_kf
=
\sum_{i=1}^n\,\xi_{ki}(x_1,\dots,x_n)\,
\frac{\partial f}{\partial x_i}
\ \ \ \ \ \ \ \ \ \ \ \ \ {\scriptstyle{(k\,=\,1\,\cdots\,r)}}.
\]

We want to suppose that amongst the $r$ equations $X_1 f = 0$, \dots,
$X_r f = 0$, exactly $q$ mutually independent are extant, hence that
in the matrix~\thetag{ 1}, all $(q+1) \times (q+1)$ determinants
vanish identically, but not all $q \times q$ ones.

Exactly as in Chap.~\ref{kapitel-7}, p.~\pageref{S-120}
and~\pageref{S-123-bis} \emphasis{we can then distribute in $q+1$
different classes the systems of equations which admit the $r$
infinitesimal transformations $X_1f, \dots, X_r f$. In one and the
same class, we reckon here the systems of equations by virtue of which
all $(p+1) \times (p+1)$ determinants of the matrix~\thetag{ 1}
vanish, but not all $p \times p$ determinants, where it is understood
that $p$ is one of the $q + 1$ numbers $q$, $q-1$, \dots, $1$, $0$}.
If we prefer to apply the way of expressing of the theory of
manifolds, we must say: to one and the same class belong the invariant
manifolds whose points admit the same number, say exactly $r - p$, of
infinitesimal transformations $e_1 \, X_1 f + \cdots + e_r \, X_r f$.
To every point of such a manifold, the infinitesimal transformations
$X_1f, \dots, X_r f$ attach exactly $p$ independent directions, which
in turn are in contact with the manifold
(cf. Chap.~\pageref{kapitel-7}, p.~\pageref{S-135}).

The usefulness of this classification is that it makes it possible to
consider every individual class for itself and to determine the
systems of equations which belong to it, or, respectively, the
manifolds.

If the number $p$ equals $q$, then the determination of all invariant
systems of equations which belong to the concerned class is achieved
by the Theorem~17 in Chap.~\ref{kapitel-7},
p.~\pageref{Theorem-17-S-123}. Every such system of equations can be
represented by relations between the common solutions
\label{S-230}
of the equations $X_1f = 0$, \dots, $X_r f = 0$. Since these $r$
equations determine a $q$-term complete system, they will of course
possess common solutions only when $q$ is smaller than $n$.

At present, we can disregard the case $p = q$. Hence we assume from
now on that $p$ is one of the numbers $0$, $1$, \dots, $q-1$ and we
state the problem of determining all manifolds invariant by the group
$X_1f, \dots, X_r f$ which belong to the class defined by the number
$p$.

The first step for solving this problem is the determination of the
locus of all points for which all $(p+1) \times (p+1)$ determinants of
the matrix~\thetag{ 1} vanish, whereas not all $p \times p$
determinants do. Indeed, the corresponding locus clearly contains all
manifolds invariant by the group which belong to our class; besides,
according to Theorem~38 p.~\pageref{Theorem-38-S-227}, this locus
itself constitutes an invariant manifold.

In order to find the sought locus, we at first study the totality of
all points for which all $(p+1) \times (p+1)$ determinants of the
matrix~\thetag{ 1} vanish, that is to say, we calculate all the $(p+1)
\times (p+1)$ determinants of the matrix in question\,\,---\,\,they
can be denoted $\Delta_1$, $\Delta_2$, \dots,
$\Delta_\rho$\,\,---\,\,and we set them equal to zero:
\[
\Delta_1=0,
\,\,\,\dots,\,\,\,
\Delta_\rho=0.
\]
The so obtained equation then represent a manifold which, according to
Theorem~39, p.~\pageref{Theorem-39-S-228}, is invariant by the group
and which contains the sought locus.

If there is in fact no system of values $x_1, \dots, x_n$ which brings
to zero all the $\Delta$, or if the $(p+1) \times (p+1)$ determinants
can vanish only in such a way that all $p \times p$ determinants also
vanish at the same time, then it is clear that actually no manifold
invariant by the group $X_1f, \dots, X_r f$ belongs to the class which
is defined by $p$. Consequently, we see that not exactly each one of
our $q+1$ classes need to be represented by manifolds which belong to
it.

We assume that for the $p$ chosen by us, none of the two exceptional
cases discussed just now occurs, so that there really are systems of
values $x_1, \dots, x_n$ for which all $(p+1) \times (p+1)$
determinants of the matrix~\thetag{ 1}, but not all $p \times p$ ones
vanish. 


As we have already observed, the manifold $\Delta_1 = 0$, \dots,
$\Delta_\rho = 0$ remains invariant by the group $X_1f, \dots, X_r
f$. Now, if this manifold is reducible, it therefore consists of a
discrete number of finitely many different manifolds, so it decomposes
without difficulty in as many different individual invariant
manifolds. Indeed, the group $X_1f, \dots, X_rf$ is generated by
infinitesimal transformations; hence when it leaves invariant the
totality of finitely many manifolds, then each individual manifold
must stay at rest\footnote{\,
In fact, the argument is that each stratum is kept invariant because
the group acts close to the identity.
}. 

Let $M_1$, $M_2$, \dots\, be the individual irreducible, and so
invariant by the group, manifolds in which the manifold $\Delta_1 =
0$, \dots, $\Delta_\rho = 0$ decomposes. Then possibly amongst these
manifolds, there are some, for the points of which all $p \times p$
determinants of the matrix~\thetag{ 1} also vanish. When we exclude
all manifolds of this special constitution, we still keep certain
manifolds ${\sf M}_1$, ${\sf M}_2$, \dots, the totality of which
clearly forms the locus of all points for which all $(p+1) \times
(p+1)$ determinants of the matrix~\thetag{ 1} vanish, but not all $p
\times p$ ones.

With this, we have found the sought locus; at the same time, we see
that this locus can consist of a discrete number of individual
invariant manifolds ${\sf M}_1$, ${\sf M}_2$, \dots\, which,
naturally, belong all to the class defined by $p$.

Clearly, each manifold invariant by our group which belongs to the
class defined by $p$ is contained in one of the manifolds ${\sf M}_1$,
${\sf M}_2$, \dots\, So in order to find all such manifolds, we need
only to examine each individual manifold ${\sf M}_1$, ${\sf M}_2$,
\dots\, and to determine the invariant manifolds contained in them
which belong to the said class. According to a remark made earlier on
(Chap.~\ref{kapitel-7}, Proposition~6, p.~\pageref{S-136}), each
invariant manifold belonging to the class $p$ is at least $p$-times
extended.

\sectionengellie{\S\,\,\,64.}

The problem, to which we have been conducted at the end of the
preceding paragraph is a special case of the following general
problem:

\plainstatement{\label{S-232}
Suppose the equations of an irreducible manifold $M$ which remain
invariant by the transformations of the $r$-term group $X_1f, \dots,
X_rf$ are given. In general, the infinitesimal transformations $X_1f,
\dots, X_r f$ attach exactly $p$ independent directions to the points
of the manifold with which they are naturally in contact, and the
manifold is at least $p$-times extended. To seek all invariant
subsidiary domains contained in $M$, to the points of which the
transformations $X_1f, \dots, X_r f$ attach exactly $p$ independent
directions.}

We now want to solve this problem.

We imagine the equations of $M$ presented under the resolved form:
\[
x_{s+i}
=
\varphi_{s+i}(x_1,\dots,x_s)
\ \ \ \ \ \ \ \ \ \ \ \ \ {\scriptstyle{(i\,=\,1\,\cdots\,n\,-\,s)}},
\]
but here it should not be forgotten that by the choice of a
\emphasis{determined} resolution, we exclude all systems of values
$x_1, \dots, x_n$ for which precisely this resolution is not possible.
It is thinkable that we exclude in this manner certain invariant
subsidiary domains of $M$ which are caught by another resolution.

The case $p = 0$ requires no special treatment, because obviously, the
manifold $M$ then consists only of invariant points.

In order to be able to solve the problem for the remaining values of
$p$, we must begin by mentioning a few remarks that stand in tight
connection to the analytic developments in Chap.~\ref{kapitel-7},
P.~\pageref{S-126} and~\pageref{S-127}, and which already possess in
principle a great importance.

Since the manifold $M$ remains invariant by the transformations of our
group, its points are permuted by the transformations of the
group. Hence, if we disregard all points lying outside of $M$, then our
group $X_1f, \dots, X_r f$ determines a certain group of
transformations of the points of $M$. However, this new group need
not contain $r$ essential parameters, since it can happen that a
subgroup of the group $X_1f, \dots, X_rf$ leaves all points of $M$
individually fixed.
 
We want to summarize at first what has been said:

\def\thetheorem{40}\begin{theorem}
\label{Theorem-40-S-233}
The points of a manifold that remains invariant by an $r$-term group
of the space $(x_1, \dots, x_n)$ are in their turn transformed by a
continuous group with $r$ or less parameters.
\end{theorem}

Since we have assumed that the invariant manifold $M$ is irreducible,
we can consider it as being a space for itself. The analytic
expression of the transformation group by which the points of this
space are transformed must therefore be obtained by covering the
points of $M$ by means of a related coordinate system and by
establishing how these coordinates are transformed by the group $X_1
f, \dots, X_r f$.

Under the assumptions made, the group which transforms the points of
$M$ can be immediately indicated. Indeed, we need only to interpret
$x_1, \dots, x_s$ as coordinates of the points of $M$, and in the
finite equations $x_i' = f_i ( x, a)$ of the group $X_1f, \dots, X_r
f$, to replace the $x_{ s+1}, \dots, x_n$ by $\varphi_{ s+1}, \dots,
\varphi_n$ and to leave out $x_{ s+1}', \dots, x_n'$;
then we obtain the equations of the concerned group, 
they are the following:
\[
x_i'
=
f_i(x_1,\dots,x_s,\,
\varphi_{s+1},\dots,\varphi_n;\,
a_1,\dots,a_r)
\ \ \ \ \ \ \ \ \ \ \ \ \ {\scriptstyle{(i\,=\,1\,\cdots\,n)}}.
\]

One could convince oneself directly that one really has to face with a
group in the variables $x_1, \dots, x_n$. For that, one would only
need to execute two transformations one after another\footnote{\,
Remind from Chaps.~\ref{one-term-groups} and~\ref{kapitel-9} that only
closure under composition counts for Lie. 
} 
under the form just written, and then to take account of two facts:
firstly, that the transformations $x_i' = f_i ( x_1, \dots, x_n, \,
a_1, \dots, a_r)$ form a group, and secondly that the system of
equations $x_{ s+i} = \varphi_{ s+i}$ admits this group.

\label{S-233-sq}
If on the other hand one would want to know the
\emphasis{infinitesimal} transformations of the group in $x_1, \dots,
x_s$, then one would only need to leave out the terms with $\partial f
/ \partial x_{ s+1}$, \dots, $\partial f / \partial x_n$ in the $X_k
f$ and to replace $x_{ s+1}, \dots, x_n$ by the $\varphi$ in the
remaining terms. One finds the $r$ infinitesimal transformations:
\[
\label{S-234}
\aligned
\overline{X}_kf
=
\sum_{\nu=1}^s\,
\xi_{k\nu}
&
(x_1,\dots,x_s,\,\varphi_{s+1},\dots,\varphi_n)\,
\frac{\partial f}{\partial x_\nu}
\\
&\ \ \ \
{\scriptstyle{(k\,=\,1\,\cdots\,r)}},
\endaligned
\]
which, however, need not be independent of each other.

We will verify directly that the reduced infinitesimal transformations
$\overline{ X}_kf$ generate a group. The concerned computation is
mostly similar to the one executed in Chap.~\ref{kapitel-7},
p.~\pageref{S-131}.

At that time, we indicated the execution of the substitution $x_{ s+i
} = \varphi_{ s+i}$ by means of the symbol $[ \, \, \, ]$. So we have
at first:
\[
\overline{X}_kf
=
\sum_{\nu=1}^s\,
[\xi_{k\nu}]\,
\frac{\partial f}{\partial x_\nu}.
\]

Furthermore as before, we see through 
the Eq.~\thetag{ 3} of the
Chap.~\ref{kapitel-7}, 
p.~\pageref{relation-3-S-110}, that:
\[
\overline{X}_k[\Omega]
\equiv
[X_k\Omega],
\]
where it is understood that $\Omega$ is a completely arbitrary
function of $x_1, \dots, x_n$. From this, it therefore comes:
\[
\leftbracket
\overline{X}_k,\,\overline{X}_j
\rightbracket
=
\sum_{\nu=1}^s\,
\big\{
[X_k\xi_{j\nu}]-[X_j\xi_{k\nu}]
\big\}
\frac{\partial f}{\partial x_\nu}\,;
\]
and since relations of the form:
\[
\leftbracket
X_k,\,X_j
\rightbracket
=
\sum_{\pi=1}^r\,
c_{kj\pi}\,X_\pi f,
\]
or, what is the same, of the form: 
\[
X_k\xi_{j\nu}-X_j\xi_{k\nu}
=
\sum_{\pi=1}^r\,c_{kj\pi}\,\xi_{\pi\nu}
\]
hold true, then we obtain simply:
\[
\leftbracket
\overline{X}_k,\,\overline{X}_j
\rightbracket
=
\sum_{\pi=1}^r\,
c_{kj\pi}\,\overline{X}_\pi f.
\]

As a result, it is proved in a purely analytic way that $\overline{
X}_1 f, \dots, \overline{ X}_r f$ really generate a group.

The infinitesimal transformations $X_1f, \dots, X_r f$ of the space
$x_1, \dots, x_n$ attach exactly $p$ independent directions $\D\, x_1
\colon \cdots \colon \D\, x_n$ to every point in general position on the
manifold $M$, and these directions, as is known, are in contact with
the manifold. It can be foreseen that the infinitesimal
transformations $\overline{ X}_1 f, \dots, \overline{ X}_r f$ of the
space $x_1, \dots, x_s$, or what is the same, of the manifold $M$,
also attach to every point $x_1, \dots, x_s$ in general position
exactly $p$ independent directions $\D\, 
x_1 \colon \cdots \colon \D\, x_s$.
We shall verify that this is really so.

Under the assumptions made, after the substitution $x_{ s+i } =
\varphi_{ s+i}$, all $(p+1) \times (p+1)$ determinants of the
matrix~\thetag{ 1} vanish, but not all $p \times p$ determinants, and
therefore, amongst the $r$ equations:
\def\theequation{2}\begin{equation}
[\xi_{k1}]\,\frac{\partial f}{\partial x_1}
+\cdots+
[\xi_{kn}]\,\frac{\partial f}{\partial x_n}
\ \ \ \ \ \ \ \ \ \ \ \ \ {\scriptstyle{(k\,=\,1\,\cdots\,r)}},
\end{equation}
exactly $p$ independent ones are extant. From this, it follows that
amongst the $r$ equations $\overline{ X}_k f = 0$, there are at most
$p$ independent ones; our problem is to prove that there are exactly
$p$. This is not difficult.

Since the system of equations $x_{ s+i} - \varphi_{ s+i} = 0$
admits the infinitesimal transformations $X_kf$, 
we have identically:
\[
[X_k(x_{s+i}-\varphi_{s+i})]
\equiv
0,
\]
or if written in greater length:
\[
[\xi_{k,s+i}]
\equiv
\sum_{\nu=1}^s\,
[\xi_{k\nu}]\,
\frac{\partial\varphi_{s+i}}{\partial x_\nu}.
\]

Hence, if by $\chi_1, \dots, \chi_r$ we denote
arbitrary functions of $x_1, \dots, x_s$, we then have:
\[
\sum_{k=1}^r\,
\chi_k\,[\xi_{k,s+i}]
\equiv
\sum_{k=1}^r\,\chi_k\,
\overline{X}_k\varphi_{s+i}.
\]

Now, if there are $r$ functions $\psi_1, \dots, \psi_r$ of $x_1,
\dots, x_n$ not all vanishing such that the equation:
\[
\sum_{k=1}^r\,\psi_k(x_1,\dots,x_s)\,
\overline{X}_kf
\equiv
0
\]
is identically satisfied, then we have:
\[
\sum_{k=1}^r\,\psi_k\,[\xi_{k,s+i}]
\equiv
\sum_{k=1}^r\,\psi_k\,\overline{X}_k\,\varphi_{s+i}
\equiv
0,
\]
and consequently also:
\[
\sum_{k=1}^r\,\psi_k(x_1,\dots,x_s)\,
\sum_{\nu=1}^n\,[\xi_{k\nu}]\,
\frac{\partial f}{\partial x_\nu}
\equiv
0.
\]

As a result, it is proved that amongst the equations:
\[
\overline{X}_1f=0,
\,\,\,\dots,\,\,\,
\overline{X}_rf=0,
\]
there are exactly as many independent equations as there are amongst
the equations~\thetag{ 2}, that is to say, exactly $p$ independent
ones.

At present, we are at last in a position to settle the problem posed
in the beginning of the paragraph, on p.~\pageref{S-232}.

\label{S-236}
The thing is the determination of certain subsidiary domains invariant
by the group $X_1f, \dots, X_r f$ in the invariant manifold $M$,
namely the subsidiary domains to the points of which the infinitesimal
transformations $X_1f, \dots, X_r f$ attach exactly $p$ independent
directions. According to what precedes, these subsidiary domains can
be defined as certain manifolds of the space $x_1, \dots, x_s$; as
such, they are characterized by the fact that they admit the group
$\overline{ X}_1f, \dots, \overline{ X}_r f$ and that, to their
points, are attached exactly $p$ independent directions by the
infinitesimal transformations $\overline{ X}_1f, \dots, \overline{
X}_r f$. Consequently, our problem amounts to the following:

In a space $M$ of $s$ dimensions, let the group $\overline{ X}_1f,
\dots, \overline{ X}_r f$ be given, whose infinitesimal
transformations attach, to the points of this space in general
position, exactly $p \leqslant s$ independent directions. To seek all
invariant manifolds contained in $M$ having the same constitution.

But we already have solved this problem above (p.~\pageref{S-230});
only at that time we had the group $X_1f, \dots, X_r f$ in place of
the group $\overline{ X}_1 f, \dots, \overline{ X}_r f$, the number
$n$ in place of the number $s$, the number $q$ in place of the number
$p$. Thus, the wanted manifolds are represented by means of relations
between the solutions of the $p$-term complete system that the
equations $\overline{ X}_1f = 0$, \dots, $\overline{ X}_r f = 0$
determine. If one adds these relations to the equations of $M$, then
one obtains the equations of the invariant subsidiary domains of $M$
in terms of the initial variables $x_1, \dots, x_n$.

Naturally, there are invariant subsidiary domains in $M$ of the
demanded sort only when $s$ is larger than $p$, and there are none,
when the numbers $s$ and $p$ are equal one to another.

With this, we therefore have at first the following important result:

\def\thetheorem{41}\begin{theorem}
If an $s$-times extended manifold of the space $x_1, \dots, x_n$
admits the $r$-term group $X_1f, \dots, X_rf$ and if, to the points of
this manifold, the infinitesimal transformations attach exactly $p$
independent directions which then surely fall into the
manifold\footnote{\,
Act of ``intrinsiqueness'': directions attached to $M$ inside
the ambient space happen to in fact be intrinsically attached
to $M$. 
}, 
then $s$ is $\geqslant p$; in case $s > p$, the manifold decomposes in
$\infty^{ s-p}$ $p$-times extended subsidiary domains, each of which
admits the group $X_1f, \dots, X_r f$.
\end{theorem}

At the same time, the problem to which we were conducted 
at the end of the
preceding paragraph (p.~\pageref{S-232}) is also completely solved,
and with this, the determination of all systems of equations that the
group $X_1f, \dots, X_r f$ admits, is achieved. In view of
applications, we put together yet once more the guidelines
\deutsch{Massregeln} required for that.

\def\thetheorem{42}\begin{theorem}
\label{Theorem-42-S-237}
If the $r$ independent infinitesimal transformations:
\[
X_kf
=
\sum_{i=1}^n\,\xi_{ki}(x_1,\dots,x_n)\,
\frac{\partial f}{\partial x_i}
\ \ \ \ \ \ \ \ \ \ \ \ \ {\scriptstyle{(k\,=\,1\,\cdots\,r)}}
\]
generate an $r$-term group and if
at the same time, all $(q+1) \times (q+1)$
determinants of the matrix:
\[
\left\vert
\begin{array}{ccccc}
\xi_{11} & \,\cdot\, & \,\cdot\, & \,\cdot\, & \xi_{1n}
\\
\cdot & \,\cdot\, & \,\cdot\, & \,\cdot\, & \cdot
\\
\xi_{r1} & \,\cdot\, & \,\cdot\, & \,\cdot\, & \xi_{rn}
\end{array}
\right\vert
\]
vanish identically, whereas not all $q \times q$ determinants do, then
one finds as follows all systems of equations, or, what is the same,
all manifolds that the group admits.

One distributes the systems of equations or the manifolds in question
in $q+1$ different classes by reckoning to be always in the same class
the systems of equations by virtue of which\footnote{\,
This just means systems of equations {\em including} the equations of
$(p+1) \times (p+1)$ minors. 
} 
all $(p+1) \times (p+1)$ determinants of the above matrix, but not all
$p \times p$ ones, vanish, where it is understood that $p$ is one of
the numbers $q$, $q-1$, \dots, $1$, $0$.

Then in order to find all invariant systems of equations which belongs
to a determined class, one forms all $(p+1) \times (p+1)$ determinants
$\Delta_1$, $\Delta_2$, \dots, $\Delta_\rho$ of the matrix and one
sets them equal to zero. If there is no system of values $x_1, \dots,
x_n$ which brings to zero all the $\rho$ determinants $\Delta_i$ at
the same time, then actually, the class which is defined by the number
$p$ contains absolutely no invariant manifold; and the same evidently
holds also for the classes with the numbers $p-1$, $p-2$, \dots, $1$,
$0$. On the other hand, if all systems of values $x_1, \dots, x_n$
which make $\Delta_1, \dots, \Delta_\rho$ equal to zero would at the
same time bring to zero all $p \times p$ determinants of the matrix,
then also in this case, the class with the number $p$ would absolutely
not be present as manifolds. If none of these two cases occurs, then
the system of equations:
\[
\Delta_1=0,
\,\,\,\dots,\,\,\,
\Delta_\rho=0
\]
represents the manifold ${\sf M}$ invariant by the group inside which
all invariant manifolds with the class number $p$ are contained. If
${\sf M}$ decomposes in a discrete number of manifolds ${\sf M}_1$,
${\sf M}_2$, \dots, then these manifolds remain individually
invariant, but in each one of them, infinitely many invariant
subsidiary domains can yet be contained which belong to the same class
as ${\sf M}$. In order to find these subsidiary domains, one sets the
equations of, say ${\sf M}_1$, in resolved form:
\[
x_{s+i}
=
\varphi_{s+i}(x_1,\dots,x_s)
\ \ \ \ \ \ \ \ \ \ \ \ \ {\scriptstyle{(i\,=\,1\,\cdots\,n\,-\,s)}},
\]
where the entire number $s$ is at least equal to $p$. Lastly, one
forms the $r$ equations:
\[
\overline{X}_kf
=
\sum_{\nu=1}^s\,
\xi_{k\nu}(x_1,\dots,x_s,\,
\varphi_{s+1},\dots,\varphi_n)\,
\frac{\partial f}{\partial x_\nu}
=
0,
\]
and one computes any $s - p$ independent solutions:
\[
\omega_1(x_1,\dots,x_s),
\,\,\,\dots,\,\,\,
\omega_{s-p}(x_1,\dots,x_s)
\]
of the $p$-term complete system determined by these equations. The
general analytic expression for the sought invariant subsidiary
domains of ${\sf M}_1$ is then:
\[
\aligned
x_{s+i}
&
-\varphi_{s+i}(x_1,\dots,x_n)
=
0,\ \ \ \ \ \ \ \
\psi_j(\omega_1,\dots,\omega_{s-p})=0
\\
&
\ \ \ \ \ \ \ \ \ \ \ \ \ \ \ \ \ \ \ 
{\scriptstyle{(i\,=\,1\,\cdots\,n\,-\,s\,;\,\,\,
j\,=\,1\,\cdots\,m)}},
\endaligned
\]
where the $m \leqslant s - p$ relations $\psi_j = 0$ are completely
arbitrary. ---

Naturally, also ${\sf M}_2$, \dots\, must be treated in the same
manner as ${\sf M}_1$. In addition, for $p$, one has to insert one
after the other all the $q+1$ numbers $q$, $q-1$, \dots, $1$, $0$.
\end{theorem}\label{S-238}

\sectionengellie{\S\,\,\,65.}

In order to apply the preceding researches to an example, we consider
the three-term group:
\[
\aligned
X_1f
&
=
\frac{\partial f}{\partial y}
+
x\,\frac{\partial f}{\partial z},
\ \ \ \ \ \ \ \ \
X_2f
=
y\,\frac{\partial f}{\partial y}
+
z\,\frac{\partial f}{\partial z},
\\
X_3f
&
=
(-z+xy)\,\frac{\partial f}{\partial x}
+
y^2\,\frac{\partial f}{\partial y}
+
yz\,\frac{\partial f}{\partial z}
\endaligned
\]
of the ordinary space. The group is transitive, because the determinant:
\[
\Delta
=
\left\vert
\begin{array}{ccc}
0 & 1 & x
\\
0 & y & z
\\
-z\!+\!xy\, & y^2 & yz
\end{array}
\right\vert
=
-\,(z-xy)^2
\]
does not vanish identically. From this, we conclude that the surface
of second degree: $z - xy = 0$ remains invariant by the group, and
else that no further surface does.

For the points of the surface $z - xy = 0$, but also only for these
points, all the $2 \times 2$ subdeterminants of $\Delta$ even vanish,
while its $1 \times 1$ subdeterminants cannot vanish simultaneously.
From this, it follows that the invariant surface decomposes in
$\infty^1$ invariant curves, but that apart from these, no other
invariant curves exist, while there are actually no invariant points.

I order to find the $\infty^1$ curves on the surface
$z - xy = 0$, we choose $x$ and $y$ as coordinates for the
points of the surface and we form, according
to the instructions given above, the
reduced infinitesimal transformations in $x, y$:
\[
\overline{X}_1f
=
\frac{\partial f}{\partial y},
\ \ \ \ \ \ \ \
\overline{X}_2f
=
y\,\frac{\partial f}{\partial y},
\ \ \ \ \ \ \ \
\overline{X}_3f
=
y^2\,\frac{\partial f}{\partial y}.
\]

The three equations $\overline{X}_k f = 0$ reduce to a single
one whose solutions is $x$. Therefore, the sought curves
are represented by the equations: 
\[
z-xy=0,
\ \ \ \ \ \ \ \ \ \
x={\rm const.},
\]
that is to say, all individuals of a family of generatrices 
on the surface of second degree remain invariant. 

\sectionengellie{\S\,\,\,66.}

Here, we give \label{S-66}
one of the two new proofs promised on p.~\pageref{S-229} for the
important Theorem~39.

As before, we denote by $\Delta_1 ( x), \dots, \Delta_\rho ( x)$
all the $(p+1) \times (p+1)$ determinants of the matrix:
\def\theequation{3}\begin{equation}
\left\vert
\begin{array}{cccc}
\xi_{11}(x) & \,\cdot\, & \,\cdot\, & \xi_{1n}(x)
\\
\cdot & \,\cdot\, & \,\cdot\, & \cdot
\\
\xi_{r1}(x) & \,\cdot\, & \,\cdot\, & \xi_{rn}(x)
\end{array}
\right\vert.
\end{equation}

In addition, we assume that there are systems of values $x_1, \dots,
x_n$ which bring to zero all the $\rho$ determinants $\Delta$. Then
it is to be proved that the system of equations:
\[
\Delta_1(x)=0,
\,\,\,\dots,\,\,\,
\Delta_\rho(x)=0
\]
admits all transformations:
\[
x_i'
=
f_i(x_1,\dots,x_n,\,a_1,\dots,a_r)
\ \ \ \ \ \ \ \ \ \ \ \ \ {\scriptstyle{(i\,=\,1\,\cdots\,n)}}
\]
of the $r$-term group $X_1f, \dots, X_r f$. 

According to Chap.~\ref{kapitel-6}, p.~\pageref{S-98-ter}, 
this only amounts to prove that the system of equations:
\[
\Delta_1(x')=0,
\,\,\,\dots,\,\,\,
\Delta_\rho(x')=0
\]
is equivalent, after the substitution $x_i' = f_i ( x, \, a)$,
to the system of equations:
\[
\Delta_1(x)=0,
\,\,\,\dots,\,\,\,
\Delta_\rho(x)=0\,;
\]
here, it is completely indifferent whether or not the equations
$\Delta_1 = 0$, \dots, $\Delta_\rho = 0$ are mutually independent.

In order to prove that our system of equations really possesses the
property in question, we proceed as follows:

In Chap.~\ref{fundamental-differential}, 
p.~\pageref{e-e-prime-a}\footnote{\,
The matrix denoted $\widetilde{ \rho} (a)$ is
denoted here $\omega ( a)$
}, 
we have seen that by
virtue of the equations $x_i' = f_i ( x, a)$, a relation of the form:
\def\theequation{4}\begin{equation}
\sum_{k=1}^r\,e_k'\,X_k'f
=
\sum_{k=1}^r\,e_k\,X_kf
\end{equation}
holds, in which the $e_k'$ are related to the $e_k$
by the $r$ equations:
\[
e_j
=
-\,\sum_{\pi,\,\,k}^{1\cdots r}\,
\vartheta_{kj}(a)\,
\alpha_{\pi k}(a)\,
e_\pi'
=
\sum_{\pi=1}^r\,\omega_{j\pi}(a)\,e_\pi'.
\]
If we insert the just written expression of the
$e_k$ in~\thetag{ 4} and if we compare the
coefficients of the two sides, we obtain $r$
relations:
\[
X_k'f
=
\sum_{j=1}^r\,\omega_{jk}(a)\,X_jf
\ \ \ \ \ \ \ \ \ \ \ \ \ {\scriptstyle{(k\,=\,1\,\cdots\,r)}}
\]
which clearly reduce to identities
as soon as one expresses the $x'$ in terms of
the $x$ by means of the equations $x_i' = f_i ( x, a)$. 

By inserting the function $x_i'$ in place of $f$ 
in the equations just found, we obtain the equations:
\[
\aligned
X_k'x_i'
=
\xi_{ki}(x')
&
=
\sum_{j=1}^r\,\omega_{jk}(a)\,X_jx_i'
\\
&
=
\sum_{j=1}^r\,\omega_{jk}(a)\,
\sum_{\nu=1}^n\,\xi_{j\nu}(x)\,
\frac{\partial f_i(x,a)}{\partial x_\nu}
\endaligned
\]
which express directly the $\xi_{ki}(x')$ as functions
of the $x$ and $a$. 
Thanks to this, we are in a position to 
study the behaviour of the equations $\Delta ( x') = 0$
after the substitution $x_i' = f_i ( x, a)$. 

The determinants $\Delta_1 ( x'), \dots, \Delta_\rho (x)$
are made up from the matrix:
\[
\left\vert
\begin{array}{cccc}
\xi_{11}(x') & \,\cdot\, & \,\cdot\, & \xi_{1n}(x')
\\
\cdot & \,\cdot\, & \,\cdot\, & \cdot
\\
\xi_{r1}(x') & \,\cdot\, & \,\cdot\, & \xi_{rn}(x')
\end{array}
\right\vert
\]
in the same way as the determinants $\Delta_1 ( x), \dots, 
\Delta_\rho ( x)$ are made up from the matrix~\thetag{ 3}. 
Now, if we imagine that the values found a while ago: 
\[
\xi_{ki}(x')
=
\sum_{j=1}^r\,\omega_{jk}\,X_jx_i'
\]
are inserted in the matrix just written and then
that the determinants $\Delta ( x')$ are computed, we
realize at first what follows: 
the determinants $\Delta_\sigma ( x')$ have the
form: 
\[
\Delta_\sigma(x')
=
\sum_{\tau=1}^\rho\,\chi_{\sigma\tau}(a)\,D_\tau
\ \ \ \ \ \ \ \ \ \ \ \ \ {\scriptstyle{(\sigma\,=\,1\,\cdots\,\rho)}},
\]
where the $\chi_{\sigma\tau}$ are certain determinants
formed with the $\omega_{ jk}(a)$, while
$D_1, \dots, D_\rho$ denote all the
$(p+1) \times (p+1)$ determinants of the matrix:
\[
\left\vert
\begin{array}{cccc}
X_1x_1' & \,\cdot\, & \,\cdot\, & X_1x_n'
\\
\cdot & \,\cdot\, & \,\cdot\, & \cdot
\\
X_rx_1' & \,\cdot\, & \,\cdot\, & X_rx_n'
\end{array}
\right\vert.
\]

Lastly, if the replace each $X_k x_i'$ by its value:
\[
X_kx_i'
=
\sum_{\nu=1}^n\,\xi_{k\nu}(x)\,
\frac{\partial f_i(x,a)}{\partial x_\nu},
\]
we obtain for the determinants $D_\tau$ expressions of the form:
\[
D_\tau
=
\sum_{\mu=1}^\rho\,\psi_{\tau\mu}(x,a)\,\Delta_\mu(x),
\]
where the $\psi_{ \tau \mu}$ are certain 
determinants formed with the 
$\frac{ \partial f_i ( x, a)}{ \partial x_\nu}$. 

\smallskip

With this, it is proved that, after
the substitution $x_i' = f_i ( x,a)$, the $\Delta_\sigma (x')$
take the form:
\[
\Delta_\sigma(x')
=
\sum_{\tau,\,\,\mu}^{1\cdots\rho}\,
\chi_{\sigma\tau}(a)\,\psi_{\tau\mu}(x,a)\,\Delta_\mu(x)
\ \ \ \ \ \ \ \ \ \ \ \ \ {\scriptstyle{(\sigma\,=\,1\,\cdots\,\rho)}}.
\]

Now, since the functions $\chi_{ \sigma \tau} (a)$, $\psi_{ \tau \mu}
( x, a)$ behave regularly for all systems of values $x, a$ coming into
consideration, it is clear that the system of equations $\Delta_\sigma
( x') = 0$ is equivalent, after the substitution $x_i' = f_i ( x,a)$,
to the system of equations $\Delta_\sigma ( x) = 0$, hence that the
latter system of equations admits all transformations $x_i' = f_i ( x,
a)$. But this is what was to be proved.

\sectionengellie{\S\,\,\,67.}

The Theorem~39 \label{S-67}
is so important that it appears not to be superfluous
to produce yet a third proof of it.

According to the Proposition~3 of the Chap.~\ref{kapitel-7}
(p.~\pageref{X-Delta-sigma}), the system of equations $\Delta_1 = 0$,
\dots, $\Delta_\rho = 0$ certainly admits all transformations of the
$r$-term group $X_1f, \dots, X_r f$ when there exist relations of the
form:
\[
X_k\Delta_\sigma
=
\sum_{\tau=1}^\rho\,\omega_{\sigma\tau}(x_1,\dots,x_n)\,
\Delta_\tau
\ \ \ \ \ \ \ \ \ \ \ \ \ {\scriptstyle{(k\,=\,1\,\cdots\,r\,;\,\,\,
\sigma\,=\,1\,\cdots\,\rho)}}
\]
and when in addition the functions $\omega_{ \sigma \tau}$
behave regularly for the systems of values
which satisfy the equations $\Delta_1 = 0$, \dots, 
$\Delta_\rho = 0$. 
Now, it is in our case not more difficult to 
prove that the system of equations
$\Delta_\sigma = 0$
satisfies this property. 
But in order not to be too much extensive, we
want to execute this proof only in a special 
case. In such a way, one will see well
how to treat the most general case. 

We will firstly assume that our group is
simply transitive. So it contains $n$ independent infinitesimal
transformations and in addition, the determinant:
\[
\Delta
=
\left\vert
\begin{array}{ccccc}
\xi_{11}(x) & \,\cdot\, & \,\cdot\, & \,\cdot\, & \xi_{1n}(x)
\\
\cdot & \,\cdot\, & \,\cdot\, & \,\cdot\, & \cdot
\\
\xi_{n1}(x) & \,\cdot\, & \,\cdot\, & \,\cdot\, & \xi_{nn}(x)
\end{array}
\right\vert
\]
does not vanish identically.

Furthermore, we will restrict ourselves to establish
that there are $n$ relations of the form:
\[
X_i\Delta
=
\omega_i(x_1,\dots,x_n)\,\Delta
\ \ \ \ \ \ \ \ \ \ \ \ \ {\scriptstyle{(i\,=\,1\,\cdots\,n)}},
\]
so that the equation $\Delta = 0$ admits
all transformations of the group. 
By contrast, we will not consider the
invariant systems of equations which are
obtained by equating to zero all 
subdeterminants of the determinant $\Delta$. 

If we express the $(n-1) \times (n-1)$ subdeterminants
of $\Delta$ as partial differential quotients
of $\Delta$ with respect to the $\xi_{ \mu \nu}$, 
it comes for $X_i \Delta$ the expression:
\[
X_i\Delta
=
\sum_{\mu,\,\,\nu}^{1\cdots n}\,
X_i\,\xi_{\mu\nu}\,
\frac{\partial\Delta}{\partial\xi_{\mu\nu}}.
\]

Now, $X_1f, \dots, X_n f$ generate an $n$-term group, so
there are relations of the form:
\[
\leftbracket
X_i,\,X_\mu
\rightbracket
=
\sum_{s=1}^n\,c_{i\mu s}\,X_sf,
\]
or else, if written with more details:
\[
X_i\,\xi_{\mu\nu}
-
X_\mu\,\xi_{i\nu}
=
\sum_{s=1}^n\,c_{i\mu s}\,\xi_{s\nu}.
\]
Consequently, for $X_i\, \xi_{ \mu\nu}$, it results the following
expression:
\[
X_i\,\xi_{\mu\nu}
=
\sum_{s=1}^n\,
\bigg(
\xi_{\mu s}\,
\frac{\partial\xi_{i\nu}}{\partial x_s}
+
c_{i\mu s}\,\xi_{s\nu}
\bigg).
\]
If, in this expression, we insert the equation above for 
$X_i \Delta$, then it comes:
\[
X_i\,\Delta
=
\sum_{\mu,\,\,\nu,\,\,s}^{1\cdots n}\,
\bigg(
\xi_{\mu s}\,
\frac{\partial\xi_{i\nu}}{\partial x_s}
+
c_{i\mu s}\,\xi_{s\nu}
\bigg)
\frac{\partial\Delta}{\partial\xi_{\mu\nu}}.
\]

Here according to a known proposition about determinants, the
coefficients of $\partial \xi_{ i\nu} / \partial x_s$ and of $c_{ i\mu
s}$ can be expressed in terms of $\Delta$. Namely, one has:
\[
\aligned
\sum_{\mu=1}^n\,\xi_{\mu s}\,
\frac{\partial\Delta}{\partial\xi_{\mu\nu}}
&
=
\varepsilon_{s\nu}\,\Delta,
\\
\sum_{\nu=1}^n\,\xi_{s\nu}\,
\frac{\partial\Delta}{\partial\xi_{\mu\nu}}
&
=
\varepsilon_{s\mu}\Delta,
\endaligned
\]
where the quantities $\varepsilon_{\pi\rho}$ vanish as soon as $\pi$
and $\rho$ are different from each other, while $\varepsilon_{ \pi
\pi}$ always has the value $1$.
Using these formulas, we obtain:
\[
X_i\Delta
=
\Delta
\bigg\{
\sum_{\nu,\,\,s}^{1\cdots n}\,
\varepsilon_{s\nu}\,
\frac{\partial\xi_{i\nu}}{\partial x_s}
+
\sum_{\mu,\,\,s}^{1\cdots n}\,
\varepsilon_{s\mu}\,c_{i\mu s}
\bigg\},
\]
and from this, it follows lastly that:
\def\theequation{5}\begin{equation}
X_i\Delta
=
\Delta\,\sum_{s=1}^n\,
\bigg(
\frac{\partial\xi_{is}}{\partial x_s}
+
c_{iss}
\bigg)
\ \ \ \ \ \ \ \ \ \ \ \ \ {\scriptstyle{(i\,=\,1\,\cdots\,n)}}.
\end{equation}

Since, as always, only systems of values $x_1, \dots, x_n$ for which
all $\xi_{ ki} (x)$ behave regularly are took into consideration, then
clearly, the factor of $\Delta$ in the right-hand side behaves
regularly for the considered systems of values $x_1, \dots, x_n$.
Hence, if the equation $\Delta = 0$ can be satisfied by such systems
of values $x_1, \dots, x_n$, then according to the Proposition~3,
p.~\pageref{X-Delta-sigma}, it admits all transformations of the group
$X_1f, \dots, X_nf$.

\sectionengellie{\S\,\,\,68.}

As was already underlined on p.~\pageref{S-223}, the developments of
the present chapter \label{S-243-sq} have great similarities with
those of the Chap.~\ref{kapitel-7}, p.~\pageref{S-120} sq. Therefore,
it is important to be conscious of the differences between the two
theories.

We have already mentioned the first difference on
page~\pageref{S-229-bis}. It consists in what follows:

When the $r$ independent infinitesimal transformations $X_1f, \dots,
X_r f$ generate an $r$-term group, then each one of systems of
equations $\Delta_1 = 0$, \dots, $\Delta_\rho = 0$ obtained by forming
determinants and mentioned more than enough admits all infinitesimal
transformations $X_1f, \dots, X_r f$. By contrast, when the $r$
infinitesimal transformations are only subjected to the restriction
that the independent equations amongst the equations $X_1f = 0$,
\dots, $X_r f = 0$ form a complete system, then in general, none of
the systems of equations $\Delta_1 = 0$, \dots, $\Delta_\rho = 0$
needs to admit the infinitesimal transformations $X_1f, \dots, X_rf$.

\smallercharacters{ 

Amongst certain conditions, one can in fact also
be sure from the beginning, 
in the second one of the two cases just mentioned, 
that a system of equations $\Delta_1 = 0$, \dots,
$\Delta_\rho = 0$ obtained by forming determinants admits the
infinitesimal transformations $X_1f, \dots, X_r f$.

Let the $r$ independent infinitesimal transformations $X_1f , \dots, 
X_r f$ be constituted in such a way that the
independent equations amongst the $r$ equations
$X_1f = 0$, \dots, $X_r f = 0$
form a complete system, so that
there are relations of the form:
\def\theequation{6}\begin{equation}
\aligned
\leftbracket
X_i,\,X_k
\rightbracket
=
\gamma_{ik1}(x_1,\dots,x_n)\,
&
X_1f
+\cdots+
\gamma_{ikr}(x_1,\dots,x_n)\,X_rf
\\
&
{\scriptstyle{(i,\,\,k\,=\,1\,\cdots\,r)}}
\endaligned
\end{equation}
Let the $(p+1) \times (p+1)$ determinants
of the matrix:
\def\theequation{7}\begin{equation}
\left\vert
\begin{array}{cccc}
\xi_{11}(x) & \,\cdot\, & \,\cdot\, & \xi_{1n}(x)
\\
\cdot & \,\cdot\, & \,\cdot\, & \cdot
\\
\xi_{r1}(x) & \,\cdot\, & \,\cdot\, & \xi_{rn}(x)
\end{array}
\right\vert
\end{equation}
be denoted by $\Delta_1, \dots, \Delta_\rho$. If there are systems of
values $x_1, \dots, x_n$ for which all determinants $\Delta_1, \dots,
\Delta_\rho$ vanish and if \emphasis{all functions $\gamma_{ ikj}
(x_1, \dots, x_n)$ behave regularly for these systems of values}, then
it can be shown that the system of equations $\Delta_1 = 0$, \dots,
$\Delta_\rho = 0$ admits the infinitesimal transformations $X_1f,
\dots, X_r f$.

In what follows, this proposition plays no rôle; it will therefore
suffice that we prove it only in a special simple case; the proof for
the general proposition can be executed in an entirely similar way.

We want to assume that $r = n$ and that the $n$ equations $X_1f = 0$,
\dots, $X_n f = 0$ are mutually independent; in addition, let $p =
n-1$. The mentioned matrix then reduces to the not identically
vanishing determinant:
\[
\Delta
=
\sum\,\pm\,\xi_{11}\cdots\xi_{nn},
\]
and it contains only a single $(p+1) \times (p+1)$ determinant, namely
itself. We will show that the equation $\Delta = 0$ then certainly
admits the infinitesimal transformations $X_1f, \dots, X_r f$ when the
functions $\gamma_{ ikj}$ in the equations:
\[
\leftbracket
X_i,\,X_k
\rightbracket
=
\sum_{j=1}^n\,\gamma_{ikj}\,X_jf
\ \ \ \ \ \ \ \ \ \ \ \ \ {\scriptstyle{(i,\,\,k\,=\,1\,\cdots\,n)}}
\]
behave regularly for the systems of values $x_1, \dots, x_n$ which
bring $\Delta$ to zero.

According to Chap.~\ref{kapitel-7}, Proposition~3,
p.~\pageref{X-Delta-sigma}, we need only to show that each $X_k
\Delta$ can be represented under the form $\omega_k (x_1, \dots, x_n)
\, \Delta$ and that the $\omega_k$ behave regularly for the systems of
values of $\Delta = 0$. This proof succeeds in the same way as in the
preceding paragraph. We simply compute the expressions $X_k \Delta$
and we find in the same way as before:
\[
\aligned
X_k\,\Delta
=
\Delta\,\sum_{\nu=1}^n\,
&
\bigg\{
\frac{\partial\xi_{k\nu}}{\partial x_\nu}
+
\gamma_{k\nu\nu}(x_1,\dots,x_n)
\bigg\}
\\
&\ \
{\scriptstyle{(k\,=\,1\,\cdots\,n)}}.
\endaligned
\]

The computation necessary for that is exactly the previous one,
although the constants $c_{ ikj}$ are replaced by the functions
$\gamma_{ ikj} (x)$; but still, that there occurs no difference has
its reason in the fact that in the preceding paragraphs, it was made
no use of the constancy property of the $c_{ ikj}$.

Clearly, the factors of $\Delta$ in the right-hand side of the above
equations behave regularly for the systems of values of $\Delta = 0$,
hence we see that the equation $\Delta = 0$ really admits the
infinitesimal transformations $X_1f, \dots, X_nf$.

}

A second important difference between the case of an $r$-term
\emphasis{group} $X_1 f, \dots, X_rf$ and the more general case of the
Chap.~\ref{kapitel-7} comes out as soon as one already knows a
manifold $M$ which admits all infinitesimal transformations $X_1f,
\dots, X_r f$.

We want to assume that $X_1f, \dots, X_r f$ attach exactly $p$
independent directions to the points of $M$ and especially, that
$X_1f, \dots, X_pf$ determine $p$ independent directions. Under these
assumptions, for the points of $M$, there are relations of the form:
\[
\aligned
X_{p+k}f
=
\varphi_{k1}(x_1,\dots,x_n)\,
&
X_1f
+\cdots+
\varphi_{kp}(x_1,\dots,x_n)\,X_pf
\\
&
{\scriptstyle{(k\,=\,1\,\cdots\,r\,-\,p)}},
\endaligned
\]
where the $\varphi_{ kj}$ behave regularly; on the other hand, there
is no relation of the form:
\[
\chi_1(x_1,\dots,x_n)\,X_1f
+\cdots+
\chi_p(x_1,\dots,x_n)\,X_pf
=
0.
\]

Now, if $X_1f, \dots, X_rf$ generate an $r$-term group, 
then for the points of $M$, all $\leftbracket X_i,\,X_k \rightbracket$
can be represented under the form:
\[
\leftbracket
X_i,\,X_k
\rightbracket
=
\sum_{j=1}^p\,\psi_{ikj}(x_1,\dots,x_n)\,X_jf
\ \ \ \ \ \ \ \ \ \ \ \ \ 
{\scriptstyle{(i,\,\,k\,=\,1\,\cdots\,r)}},
\]
and here, \emphasis{the $\psi_{ ikj}$ behave regularly}. By contrast,
if $X_1f, \dots, X_rf$ only possess the property that the independent
equations amongst the equations $X_1f = 0$, \dots, $X_r f = 0$ form a
complete system, then such a representation of the $\leftbracket
X_i,\,X_k \rightbracket$ for the points of $M$ is not possible in all
cases; but it is always possible when the functions $\gamma_{ ikj}$ in
the equations~\thetag{ 6} behave regularly for the systems of values
$x_1, \dots, x_n$ on $M$. We have already succeeded to make
use of this condition in Chap.~\ref{kapitel-7}, p.~\pageref{S-130} sq.

\linestop


\chapter{Invariant Families of Infinitesimal Transformations}
\label{kapitel-15}
\chaptermark{Invariant Families of Infinitesimal Transformations}

\setcounter{footnote}{0}

\abstract*{??}

\begin{svgraybox}
One studies in this chapter the general linear combination:
\[
e_1\,X_1
+\cdots+
e_q\,X_q
\]
of $q \geqslant 1$ given arbitrary local infinitesimal transformations:
\[
X_k
=
\sum_{i=1}^n\,\xi_{ki}(x)\,
\frac{\partial}{\partial x_i}
\ \ \ \ \ \ \ \ \ \ \ \ \
{\scriptstyle{(k\,=\,1\,\cdots\,q)}}
\]
having analytic coefficients $\xi_{ ki} ( x)$ and which are assumed to
be independent of each other. When one introduces new variables $x_i'
= \varphi_i ( x_1, \dots, x_n)$ in place of the $x_l$, every
transformation $X_k$ of this general combination receives another
form, but it may sometimes happen under certain circumstances that the
complete family in its wholeness remains unchanged, namely that there
are functions $e_k' = e_k' ( e_1, \dots, e_q)$ such that:
\[
\varphi_*
\big(
e_1\,X_1
+\cdots+
e_q\,X_q
\big)
=
e_1'(e)\,X_1'
+\cdots+
e_q'(e)\,X_q',
\]
where, as in previous circumstances, the $X_k ' = \sum_{ i = 1}^n\,
\xi_{ ki} (x') \frac{ \partial }{ \partial x_i'}$ denote the same
vector fields, but viewed in the target space $x_1', \dots,
x_n'$.

\begin{definition}
The family $e_1 \, X_1 + \cdots + e_q \, X_q$ of infinitesimal
transformations is said to \terminology{remain invariant after the
introduction of the new variables $x' = \varphi ( x)$} if there are
functions $e_k' = e_k' ( e_1, \dots, e_q)$ depending on $\varphi$ such
that:
\def\theequation{1}\begin{equation}
\varphi_*
\big(
e_1\,X_1
+\cdots+
e_q\,X_q
\big)
=
e_1'(e)\,X_1'
+\cdots+
e_q'(e)\,X_q';
\end{equation}
alternatively, one says that the family \terminology{admits} 
the transformation
which is represented by the concerned change of variables.
\end{definition}

\begin{proposition}
Then the functions $e_k' ( e)$ in question necessarily are linear:
\[
e_k'
=
\sum_{j=1}^q\,
\rho_{kj}\, e_j
\ \ \ \ \ \ \ \ \ \ \ \ \
{\scriptstyle{(k\,=\,1\,\cdots\,q)}},
\]
with the constant matrix $\big( \rho_{ kj} \big)_{1 \leqslant k
\leqslant q}^{ 1 \leqslant j \leqslant q}$ being invertible: $e_k =
\sum_{ j = 1}^q \, \widetilde{ \rho}_{ kj} \, e_j'$.
\end{proposition}

\begin{proof}\smartqed
Indeed, through the change of coordinates $x' = \varphi (x)$, if we
write that the vector fields $X_k$ are transferred to:
\[
\varphi_*(X_k)
=
\sum_{i=1}^n\,X_k(x_i')\,
\frac{\partial}{\partial x_i'}
=:
\sum_{i=1}^n\,
\eta_{ki}(x_1',\dots,x_n')\,
\frac{\partial}{\partial x_i'}
\ \ \ \ \ \ \ \ \ \ \ \ \
{\scriptstyle{(k\,=\,1\,\cdots\,q)}},
\]
with their coefficients $\eta_{ ki} = \eta_{ ki} (x')$ being expressed
in terms of the target coordinates, and if we substitute the resulting
expression into~\thetag{ 1}, we get the following linear relations:
\def\theequation{1'}\begin{equation}
\sum_{k=1}^q\,e_k'\,\xi_{ki}(x')
=
\sum_{k=1}^q\,e_k\,\eta_{ki}(x')
\ \ \ \ \ \ \ \ \ \ \ \ \
{\scriptstyle{(i\,=\,1\,\cdots\,n)}}.
\end{equation}
The idea is to substitute here for $x'$ exactly the same number $q$ of
different systems of fixed values:
\[
x_1^{(1)},\dots,x_n^{(1)},
x_1^{(2)},\dots,x_n^{(2)},
\dots\dots,
x_1^{(q)},\dots,x_n^{(q)}
\]
that are mutually in general position and considered will be
considered as constant. In fact, according to the proposition on
p.~\pageref{Satz-S-66}, or equivalently, according to the
assertion formulated just below the long matrix located on
p.~\pageref{q-fold-extended-coefficient-matrix}, the linear
independence of $X_1, \dots, X_q$ insures that for most such $q$
points, the long $q \times qn$ matrix in question:
\[
\left(
\begin{array}{cccccccccc}
\xi_{11}^{(1)} & \cdots & \xi_{1n}^{(1)} &
\xi_{11}^{(2)} & \cdots & \xi_{1n}^{(2)} &
\cdots\cdots &
\xi_{11}^{(q)} & \cdots & \xi_{1n}^{(q)}
\\
\cdots & \cdots & \cdots &
\cdots & \cdots & \cdots &
\cdots\cdots &
\cdots & \cdots & \cdots
\\ 
\xi_{q1}^{(1)} & \cdots & \xi_{qn}^{(1)} &
\xi_{q1}^{(2)} & \cdots & \xi_{qn}^{(2)} &
\cdots\cdots &
\xi_{q1}^{(q)} & \cdots & \xi_{qn}^{(q)}
\end{array}
\right)
\]
has rank equal to $q$, where we have set $\xi_{ ki}^{ ( \nu )} := \xi_{
xi} \big( x^{ ( \nu )} \big)$. Consequently, while considering the values
of $\xi_{ ki} \big( x^{ ( \nu)} \big)$ and of $\eta_{ ki} \big( x^{ (
\nu)} \big)$ as \emphasis{constant}, the linear system above is solvable
with respect to the unknowns $e_k'$ and we obtain:
\[
e_k'
=
\sum_{j=1}^q\,\rho_{kj}\, e_j
\ \ \ \ \ \ \ \ \ \ \ \ \
{\scriptstyle{(k\,=\,1\,\cdots\,q)}},
\]
for some constants $\rho_{ kj}$. In addition, we claim that the
determinant of the matrix $\big( \rho_{ kj}\big)_{ 1 \leqslant k
\leqslant q}^{ 1 \leqslant j \leqslant q}$ is in fact nonzero. Indeed,
the linear independence of $X_1, \dots, X_q$ being obviously
equivalent to the linear independence of $\varphi_* (X_1), \dots,
\varphi_* (X_q)$, the other corresponding long matrix:
\[
\left(
\begin{array}{cccccccccc}
\eta_{11}^{(1)} & \cdots & \eta_{1n}^{(1)} &
\eta_{11}^{(2)} & \cdots & \eta_{1n}^{(2)} &
\cdots\cdots &
\eta_{11}^{(q)} & \cdots & \eta_{1n}^{(q)}
\\
\cdots & \cdots & \cdots &
\cdots & \cdots & \cdots &
\cdots\cdots &
\cdots & \cdots & \cdots
\\ 
\eta_{q1}^{(1)} & \cdots & \eta_{qn}^{(1)} &
\eta_{q1}^{(2)} & \cdots & \eta_{qn}^{(2)} &
\cdots\cdots &
\eta_{q1}^{(q)} & \cdots & \eta_{qn}^{(q)}
\end{array}
\right)
\]
then also has rank equal to $q$ and we therefore can also solve
symmetrically:
\[
e_k
=
\sum_{j=1}^q\,\widetilde{\rho}_{kj}\, e_j'
\ \ \ \ \ \ \ \ \ \ \ \ \
{\scriptstyle{(k\,=\,1\,\cdots\,q)}},
\]
with coefficients $\widetilde{ \rho}_{ kj}$ which necessarily coincide
with the elements of the inverse matrix.
\qed\end{proof}

\section{Families Invariant by One-Term Subgroups}
For an important application to the study of the adjoint group in
Chap.~L, we now want to study families $e_1\, X_1 + \cdots + e_q \,
X_q$ that are invariant when the transition from $x$ to a new variable
$x'$ is performed by an arbitrary transformation of some one-term
group $x' = \exp (tY) (x)$, where $Y$ is any (local, analytic) vector
field.

\end{svgraybox}

As usual, we understand by: 
\[
X_kf
=
\sum_{i=1}^n\,\xi_{ki}\,\frac{\partial f}{\partial x_i}
\ \ \ \ \ \ \ \ \ \ \ \ \ {\scriptstyle{(k\,=\,1\,\cdots\,q)}}
\]
\emphasis{independent} infinitesimal transformations; temporarily,
this shall be the only assumption which we make about the $X_k f$.

We consider the family of $\infty^{ q-1}$ infinitesimal
transformations which is represented by the expression:
\[
e_1\,X_1f
+\cdots+
e_q\,X_qf
\]
with the $q$ arbitrary parameters $e_1, \dots, e_q$. When we
introduce, in this expression, new independent variables $x'$ in place
of the $x$, then each infinitesimal transformations in our family
takes another form; evidently, we then obtain in general a completely
new family of $\infty^{ q-1}$ infinitesimal transformations. However,
in certain circumstances, it can happen that the new family does not
essentially differ in its form from the original family, when for
arbitrary values of the $e$, there is a relation of the shape:
\def\theequation{1}\begin{equation}
\sum_{k=1}^q\,e_k\,\sum_{i=1}^n\,
\xi_{ki}(x_1,\dots,x_n)\,
\frac{\partial f}{\partial x_i}
=
\sum_{k=1}^q\,e_k'\,
\sum_{i=1}^n\,\xi_{ki}(x_1',\dots,x_n')\,
\frac{\partial f}{\partial x_i'},
\end{equation}
in which the $e_k'$ do not depend upon the $x$, 
but only\footnote{\,
In fact, as will be see in a while, the dependence can then only be
linear.
} 
upon $e_1, \dots, e_q$. 

If there is such a relation, which we can also write shortly as:
\def\theequation{2}\begin{equation}
\sum_{k=1}^q\,e_k\,X_kf
=
\sum_{k=1}^q\,e_k'\,X_k'f,
\end{equation}
then we say: \terminology{the family of the infinitesimal
transformations $\sum\, e_k \, X_k f$ remains invariant after the
introduction of the new variables $x'$}, or: \terminology{it admits
the transformation which is represented by the concerned change of
variables}.

\sectionengellie{\S\,\,\,69.}

Let the family of the $\infty^{ q-1}$ infinitesimal transformations
$\sum \, e_k \, X_k f$ remain invariant through the transition to the
variables $x'$, so that there is a relation of the form:
\def\theequation{2}\begin{equation}
\sum_{k=1}^q\,e_k\,X_kf
=
\sum_{k=1}^q\,e_k'\,X_k'f,
\end{equation}
in which the $e'$ are certain functions of only the $e$. To begin
with, we study this relationship of dependence between the $e$ and the
$e'$; in this way, we reach the starting point for the more precise
study of such families of infinitesimal transformations.

The expressions $X_kf$ can be written as:
\[
X_kf
=
\sum_{i=1}^n\,X_k\,x_i'\,\,
\frac{\partial f}{\partial x_i'},
\]
or, when one expresses the $X_k\, x_i'$ in terms
of the $x'$, as:
\[
X_kf
=
\sum_{i=1}^n\,\eta_{ki}(x_1',\dots,x_n')\,
\frac{\partial f}{\partial x_i'}.
\]

If we insert these values in the equation~\thetag{ 2}, we can equate
the coefficients of the $\partial f / \partial x_i'$ in the two sides,
and so we obtain the following linear relations between the $e$ and
the $e'$:
\def\theequation{2'}\begin{equation}
\sum_{k=1}^q\,e_k'\,\xi_{ki}(x')
=
\sum_{k=1}^q\,e_k\,\eta_{ki}(x')
\ \ \ \ \ \ \ \ \ \ \ \ \ {\scriptstyle{(i\,=\,1\,\cdots\,n)}}
\end{equation}

According to our assumption, it is possible to enter for the $e'$
functions of the $e$ alone so that the equations~\thetag{ 2'} are
satisfied for all values of the $x'$.  It can be shown that the
concerned functions of the $e$ are completely determined.

Since the equations~\thetag{ 2'} 
are supposed to hold true for all values of the
$x'$, then they must also be satisfied yet when we replace $x_1',
\dots, x_n'$ by any other system of variables.
We want to do this, and to write down the equations~\thetag{ 2'}
in exactly $q$ different systems of variables
$x_1', \dots, x_n'$, 
$x_1'', \dots, x_n''$, $\cdots$, 
$x_1^{ (q)}, \dots, x_n^{ (q)}$:
\[
\aligned
\sum_{k=1}^q\,e_k'\,\xi_{ki}(x^{(\nu)})
&
=
\sum_{k=1}^q\,e_k\,\eta_{ki}(x^{(\nu)})
\ \ \ \ \ \ \ \ \ \ \ \ \ {\scriptstyle{(i\,=\,1\,\cdots\,n)}}
\\
&
{\scriptstyle{(\nu\,=\,1\,\cdots\,q)}}.
\endaligned
\]

The so obtained equations are solvable with respect to $e_1', \dots,
e_q'$, because under the assumptions made, according to the
developments of the Chap.~\ref{one-term-groups},
p.~\pageref{crucial-assertion}, not all $q \times q$ determinants of
the matrix:
\[
\left\vert
\begin{array}{cccccccccccc}
\xi_{11}' & \,\cdot\, & \,\cdot\, & \xi_{1n}'& \xi_{11}''
& \,\cdot\, & \,\cdot\, & \xi_{1n}'' 
& \,\cdot\, & \,\cdot\, & \,\cdot\, & \xi_{1n}^{(q)}
\\
\cdot & \,\cdot\, & \,\cdot\, & \cdot & \cdot
& \,\cdot\, & \,\cdot\, & \cdot
& \,\cdot\, & \,\cdot\, & \,\cdot\, & \cdot
\\
\xi_{q1}' & \,\cdot\, & \,\cdot\, & \xi_{qn}'& \xi_{q1}''
& \,\cdot\, & \,\cdot\, & \xi_{qn}'' 
& \,\cdot\, & \,\cdot\, & \,\cdot\, & \xi_{qn}^{(q)}
\end{array}
\right\vert
\]
vanish. 

In addition, since the said equations are certainly compatible with
each other\footnote{\,
---\,and since, furthermore, the Lemma on
p.~\pageref{crucial-assertion} insures that, with a suitable choice of
generic fixed points $x_1', \dots, x_n'$, $x_1'', \dots, x_n''$,
$\cdots$, $x_1^{ (q)}, \dots, x_n^{ (q)}$, the rank of the considered
matrix of $\xi$'s is maximal equal to $q$\,---
}, 
we obtain the $e'$ represented as linear homogeneous functions of the
$e$:
\[
e_k'
=
\sum_{j=1}^q\,\rho_{kj}\,e_j
\ \ \ \ \ \ \ \ \ \ \ \ \ {\scriptstyle{(k\,=\,1\,\cdots\,q)}}.
\]

Here naturally, the $\rho_{kj}$ are independent of the $x'$, $x''$,
\dots, $x^{ (q)}$ and hence are absolute constants; the determinant of
the $\rho_{ kj}$ is different from zero, because visibly the $e_k$
can, in exactly the same way, be represented as linear homogeneous
functions of the $e'$.

Although, under the assumptions made, the \emphasis{family} of the
infinitesimal transformations $\sum \, e_k \, X_kf$ remains invariant
after the introduction of the $x'$, in general, its individual
transformations are permuted.  However, there always exists at least
one infinitesimal transformation $\sum \, e_k^0 \, X_kf$ which remains
itself invariant, since the condition which the coefficients $e_k^0$
of such an infinitesimal transformation:
\[
\sum_{k=1}^q\,e_k^0\,X_kf
=
\omega\,
\sum_{k=1}^q\,e_k^0\,X_k'f
\]
must satisfy can be replaced by the $q$ equations:
\[
\omega\,e_k^0
=
\sum_{j=1}^q\,\rho_{kj}\,e_j^0
\ \ \ \ \ \ \ \ \ \ \ \ \ {\scriptstyle{(k\,=\,1\,\cdots\,q)}}
\]
and these last equations can always be satisfied
without all the $e_k^0$ being zero. 

For closer illustration of what has been said, 
an \emphasis{example} is more suitable. 

In the family of the $\infty^3$ transformations:
\[
e_1\,\frac{\partial f}{\partial x_1}
+
e_2\,\frac{\partial f}{\partial x_2}
+
e_3
\bigg(
x_1^2\,\frac{\partial f}{\partial x_1}
+
x_1x_2\,\frac{\partial f}{\partial x_2}
\bigg)
+
e_4
\bigg(
x_1x_2\,\frac{\partial f}{\partial x_1}
+
x_2^2\,\frac{\partial f}{\partial x_2}
\bigg),
\]
we introduce new variables by setting:
\[
x_1'
=
a_1\,x_1+a_2x_2,
\ \ \ \ \ \ \ \ \ 
x_2'
=
a_3\,x_1+a_4\,x_2.
\]
At the same time, the family receives the new form:
\[
e_1'\,\frac{\partial f}{\partial x_1'}
+
e_2'\,\frac{\partial f}{\partial x_2'}
+
e_3'
\bigg(
{x_1'}^2\,\frac{\partial f}{\partial x_1'}
+
x_1'x_2'\,\frac{\partial f}{\partial x_2'}
\bigg)
+
e_4'
\bigg(
x_1'x_2'\,\frac{\partial f}{\partial x_1'}
+
{x_2'}^2\,\frac{\partial f}{\partial x_2'}
\bigg),
\]
where $e_1', \dots, e_4'$ express as follows:
\[
e_1'
=
a_1\,e_1+a_2\,e_2,\ \ \ \ \
e_2'
=
a_3\,e_1+a_4\,e_2,\ \ \ \ \
e_3'
=
\frac{a_4\,e_3-a_3\,e_4}{a_1a_4-a_2a_3},\ \ \ \ \
e_4'
=
\frac{a_1\,e_4-a_2\,e_3}{a_1a_4-a_2a_3}.
\]

Consequently, the family remains invariant in the above sense.  If one
would want to know which individual infinitesimal transformations of
the family remain invariant, one would only have to determine $\omega$
from the equation:
\[
\left\vert
\begin{array}{cc}
a_1-\omega & a_2
\\
a_3 & a_4-\omega
\end{array}
\right\vert
\cdot
\left\vert
\begin{array}{cc}
a_1\omega-1 & a_2
\\
a_3 & a_4\omega-1
\end{array}
\right\vert
=0,
\]
and to choose $e_1, \dots, e_4$ so that $e_k' = \omega \, e_k$;
the concerned systems of values of the $e_k$ 
provide the invariant infinitesimal transformations.

At present, by coming back to the general case, we want yet to
specialize in a certain direction the assumptions made above.  Namely,
we want to suppose that the transition from the $x$ to the $x'$ is a
completely arbitrary transformation of a determined group.
Correspondingly, we state from now on the following question:

\plainstatement{Under which conditions does the family $\sum\, e_k \,
X_k f$ remain invariant through every transformation $x_i' = f_i (
x_1, \dots, x_n, t)$ of the one-term group $Yf$, 
that is to say, under which conditions does a relation:
\[
\sum_{k=1}^q\,e_k\, X_kf
=
\sum_{k=1}^q\,e_k'\, X_k'f,
\]
hold for all systems of values $e_1, \dots e_q, t$, in which
the $e_k'$, aside from upon the $e_j$, 
yet only depend upon $t$?}

When, in order to introduce new variables in $X_k f$, we apply the
general transformation:
\[
x_i'
=
x_i
+
t\, Y\,x_i
+\cdots
\ \ \ \ \ \ \ \ \ \ \ \ \
{\scriptstyle{(i\,=\,1\,\cdots\,n)}}
\]
of the one-term group $Yf$, we obtain according to 
Chap.~\ref{one-term-groups}, p.~\pageref{S-141}:
\[
X_kf
=
X_k'f
+
t\big(
X_k'Y'f
-
Y'X_k'f
\big)
+\cdots;
\]
hence also inversely:
\def\theequation{3}\begin{equation}
X_k'f
=
X_kf
+
t\leftbracket Y,\,X_k\rightbracket
+\cdots,
\end{equation}
which is more convenient for what follows.

Now, if every infinitesimal transformation $X_k f + t \leftbracket Y,
\, X_k \rightbracket + \cdots$ should belong to the family $e_1 X_1 f +
\cdots + e_q X_q f$, and in fact so for every value of $t$, then
obviously every infinitesimal transformation $\leftbracket Y, \, X_k
\rightbracket$ would also be contained in this family. As a result,
certain necessary conditions for the invariance of our family would be
found, some conditions which amount to the fact that $q$ relations of
the form:
\def\theequation{4}\begin{equation}
\label{Satz-4-S-259}
\leftbracket Y,\,X_k\rightbracket
=
\sum_{j=1}^q\,g_{kj}\, X_jf
\ \ \ \ \ \ \ \ \ \ \ \ \
{\scriptstyle{(k\,=\,1\,\cdots\,q)}}
\end{equation}
should hold, in which the $g_{ kj}$ denote absolute constants.

If the family of the infinitesimal transformations:
\[
e_1\, X_1f
+\cdots+
e_q\, X_qf
\]
is constituted so that for every $k$, a relation of the form~\thetag{
4} holds true, then we want to say that \emphasis{the family admits the
infinitesimal transformation $Yf$}. By this settlement 
of terminology, we
can state as follows the result just obtained:

\plainstatement{If 
the family of the infinitesimal transformations:
\[
e_1\, X_1f
+\cdots+
e_q\, X_qf
\]
admits all transformations of the one-term group
$Yf$, then it also admits the infinitesimal transformation
$Yf$. }

But the converse too holds true, as we will now show.

We want to suppose that the family of the transformations $\sum \, e_k
\, X_kf$ admits the infinitesimal transformation $Yf$, hence that
relations of the form~\thetag{ 4} hold true. If now the family
$\sum\, e_k \, X_kf$ is supposed to 
simultaneously admit all finite
transformations of the one-term group $Yf$, then it must be possible
to determine $e_1 ' , \dots, e_q'$ as functions of $e_1, \dots, e_q$
in such a way that the equation:
\[
\sum_{k=1}^q\,e_k'\, X_k'f
=
\sum_{k=1}^q\,e_k\, X_kf
\]
is identically satisfied, as soon as one introduces the variable $x$
in place of $x'$ in the $X_k' f$. Consequently, if $X_k'f$ takes the
form:
\[
X_k'f
=
\sum_{i=1}^n\,\zeta_{ki}(x_1,\dots,x_n,t)\,
\frac{\partial}{\partial x_i}
\]
after the introduction of the $x$, then 
one must be able to determine the $e_k'$ 
so that the expression:
\[
\sum_{k=1}^q\,e_k'\, X_k'f
=
\sum_{k=1}^q\,\sum_{i=1}^n\,
e_k'\,
\zeta_{ki}(x_1,\dots,x_n,t)\,
\frac{\partial f}{\partial x_i}
\]
is free of $t$, hence so that the differential quotient:
\[
\frac{\partial}{\partial t}\,
\sum_{k=1}^q\,e_k'\, X_k'f
=
\sum_{i=1}^n\,\frac{\partial f}{\partial x_i}\,
\frac{\partial}{\partial t}\,
\sum_{k=1}^q\,e_k'\,
\zeta_{ki}(x_1,\dots,x_n,t)
\]
vanishes\footnote{\,
Indeed, differentiation with respect to $t$ of $e_1 \, X_1 + \cdots +
e_r \, X_r$ yields: $0 \equiv \frac{ \partial }{\partial t}\, \sum_{ k
= 1}^q\, e_k\, X_k$.
}; 
but at the same time, the $e$ must still also satisfy the initial
condition: $e_k' = e_k$ for $t = 0$.

In order to be able to show that under the assumptions made there
really are functions $e'$ of the required constitution, we must at
first calculate the differential quotient:
\[
\frac{\partial}{\partial t}\,X_k'f
=
\sum_{i=1}^n\,
\frac{\partial\zeta_{ki}(x_1,\dots,x_n,t)}{\partial t}\,
\frac{\partial f}{\partial x_i};
\]
for this, we shall take an indirect route.

Above, we saw that $X_k'f$ can be expressed in the following
way in terms of $x_1, \dots, x_n$ and $t$:
\[
X_k'f
=
X_kf
+
t\leftbracket Y,\,X_k\rightbracket
+\cdots,
\]
when the independent variables $x'$ entering the $X_k'$ are determined
by the equations $x_i' = f_i ( x_1, \dots, x_n, t)$ of the one-term
group $Yf$. So the desired differential quotient can be obtained by
differentiation of the infinite power series in $t$ lying in the
right-hand side, or differently enunciated: it is the coefficient of
$\tau^1$ in the expansion of the expression:
\[
X_kf
+
(t+\tau)\leftbracket Y,\,X_k\rightbracket
+\cdots
=
\sum_{i=1}^n\,\xi_{ki}(x_1'',\dots,x_n'')\,
\frac{\partial f}{\partial x_i''}
=
X_k''f
\]
with respect to powers of $\tau$. Here, 
the $x''$ mean the quantities:
\[
x_i''
=
f_i(x_1,\dots,x_n,\,t+\tau).
\]

However, the expansion coefficient \deutsch{Entwickelungscoefficient}
discussed just above appears at first as an infinite series of powers
of $t$; nevertheless, there is no difficulty to find a finite closed
expression for it.

As we know, the transition from the variables $x$ to the variables
$x_i' = f_i (x_1, \dots, x_n, t)$ occurs through a transformation of
the one-term \emphasis{group} $Yf$, and to be precise, through a
transformation with the parameter $t$. One comes from the $x$ to the
$x_i'' = f_i ( x_1, \dots, x_n, \, t + \tau)$ through a transformation
of \emphasis{the same} group, namely through the transformation with
the parameter $t + \tau$. But this transformation can be substituted
for the succession of two transformations, of which the first one
possesses the parameter $t$, and the second one the parameter $\tau$;
consequently, the transition from the $x'$ to the $x''$ is likewise
got through a transformation of the one-term group $Yf$, namely
through the transformation whose parameter is $\tau$:
\[
x_i''
=
f_i(x_1',\dots,x_n',\tau).
\]

From this, we conclude that the series expansion
of $X_k''f$ with respect to powers of $\tau$ reads:
\[
X_k''f
=
X_k'f
+
\tau\leftbracket Y',\,X_k'\rightbracket
+\cdots.
\]
As a result, we have then found a finite closed expression for the
expansion coefficient mentioned a short while ago; the sought
differential quotient $\frac{ \partial (X_k' f)}{\partial t}$ is
hence:
\def\theequation{5}\begin{equation}
\frac{\partial}{\partial t}\,
X_k'f
=
\leftbracket Y',\,X_k'\rightbracket
=
Y'\,X_k'\,f
-
X_k'\,Y'\,f.
\end{equation}

Naturally, this formula holds generally, whatever also one can choose
as the two infinitesimal transformations $X_kf$ and $Yf$. However, in
our specific case, $X_1 f, \dots, X_qf, Yf$ are not absolutely
arbitrary, but they are linked together through the relations~\thetag{
4}. So under the assumptions made above, we receive:
\def\theequation{6}\begin{equation}
\frac{\partial(X_k'f)}{\partial t}
=
\sum_{\nu=1}^q\,g_{k\nu}\,
X_\nu'f
\ \ \ \ \ \ \ \ \ \ \ \ \
{\scriptstyle{(k\,=\,1\,\cdots\,q)}}.
\end{equation}

Now, if we form the differential quotient of $\sum \, e_k' \, X_k'f$
with respect to $t$, we obtain:
\[
\aligned
\frac{\partial}{\partial t}\,
\sum_{k=1}^q\,
e_k'\, X_k'f
&
=
\sum_{k=1}^q\,
\frac{\D\,e_k'}{\D\,t}\,X_k'f
+
\sum_{k=1}^q\,e_k'\,
\sum_{\nu=1}^q\,g_{k\nu}\,
X_\nu'f
\\
&
=
\sum_{k=1}^q\,
\Big\{
\frac{\D\,e_k'}{\D\,t}
+
\sum_{\nu=1}^q\,g_{\nu k}\,e_\nu'
\Big\}\,
X_k'f.
\endaligned
\]
Obviously, this expression vanishes only when the $e_k'$ satisfy the
differential equations:
\def\theequation{7}\begin{equation}
\frac{\D\,e_k'}{\D\,t}
+
\sum_{\nu=1}^q\,g_{\nu k}\,e_\nu'
=
0
\ \ \ \ \ \ \ \ \ \ \ \ \
{\scriptstyle{(k\,=\,1\,\cdots\,q)}}.
\end{equation}
But from this the $e_k'$ can be determined as functions of $t$ in such
a way that for $t = 0$, each $e_k'$ converts into the corresponding
$e_k$; in addition, the $e'$ are linear homogeneous functions of the
$e$.

If one puts the values of the $e'$ in question in the expression $\sum
\, e_k' \, X_k' f$ and returns afterwards from the $x'$ to the initial
variables $x_1, \dots, x_n$, then $\sum \, e_k' \, X_k'$ will be
independent of $t$, that is to say, it will be equal to $\sum \, e_k
\, X_kf$. Consequently, the family of the infinitesimal
transformations $\sum \, e_k \, X_kf$ effectively remains invariant by
the change of variables in question.

As a result, we can state the following theorem:

\renewcommand{\thefootnote}{\fnsymbol{footnote}}
\def\thetheorem{43}\begin{theorem}
\label{Theorem-43-S-252}
A family of $\infty^{ q - 1}$ infinitesimal transformations $e_1 \,
X_1 f + \cdots + e_q \, X_q f$ remains invariant, through the
introduction of new variables $x'$ which are defined by the equations
of a one-term group:
\[
x_i'
=
x_i
+
t\, Yx_i
+\cdots
\ \ \ \ \ \ \ \ \ \ \ \ \
{\scriptstyle{(i\,=\,1\,\cdots\,n)}},
\]
if and only if, between $Yf$ and the $X_kf$ there are $q$ relations of
the form:
\def\theequation{4}\begin{equation}
\leftbracket
Y,\,X_k
\rightbracket
=
\sum_{\nu=1}^q\,g_{k\nu}\, X_\nu f
\ \ \ \ \ \ \ \ \ \ \ \ \
{\scriptstyle{(k\,=\,1\,\cdots\,q)}},
\end{equation}
in which the $g_{ k\nu}$ denote constants. If these conditions are
satisfied, then by the concerned change of variables, $\sum\, e_k\,
X_kf$ receives the form $\sum\, e_k' \, X_k'f$, where $e_1', \dots,
e_q'$ determine themselves through the differential equations:
\[
\frac{\D\,e_k'}{\D\,t}
+
\sum_{\nu=1}^q\,g_{\nu k}\,e_\nu'
=
0
\ \ \ \ \ \ \ \ \ \ \ \ \
{\scriptstyle{(k\,=\,1\,\cdots\,q)}},
\]
by taking account of the initial conditions: $e_k' = e_k$ for $t =
0$\footnote[1]{\,
\name{Lie}, Archiv for Mathematik og Naturvidenskab Vol.~3, 
Christiania 1878.
}. 
\end{theorem}
\renewcommand{\thefootnote}{\arabic{footnote}}

If one performs the integration of which the preceding theorem speaks, 
hence determines $e_1', \dots, e_r'$ from the differential
equations:
\[
\frac{\D\,e_k'}{\D\,t}
=
-
\sum_{\nu=1}^q\,g_{\nu k}\,e_\nu'
\ \ \ \ \ \ \ \ \ \ \ \ \
{\scriptstyle{(k\,=\,1\,\cdots\,q)}}
\]
taking as a basis the initial conditions: $e_k' = e_k$ for $t = 0$,
then one obtains equations of the form:
\[
e_k'
=
\sum_{j=1}^q\,d_{kj}(t)\, e_j
\ \ \ \ \ \ \ \ \ \ \ \ \
{\scriptstyle{(k\,=\,1\,\cdots\,q)}}.
\]
It is clear that these equations represent the finite transformations
of a certain one-term group, namely the one which is generated by the
infinitesimal transformation:
\[
\sum_{k=1}^q\,
\Big\{
\sum_{\nu=1}^q\,g_{\nu k}\,e_\nu
\Big\}\,
\frac{\partial f}{\partial e_k}
\]
(cf. Chap.~\ref{one-term-groups}). 

\medskip

From the theorem just stated, we want to derive an important
proposition which is certainly closely suggested.

If the family of $\infty^{ q-1}$ infinitesimal transformations $e_1 \,
X_1 f + \cdots + e_q \, X_q$ admits the two infinitesimal
transformations $y_1f$ and $Y_2 f$, then it also admits at the same
time each transformation $c_1 \, Y_1f + c_2 \, Y_2f$ which is linearly
deduced from $Y_1f$ and $Y_2f$; this follows immediately from the fact
that the infinitesimal transformation:
\[
\big\leftbracket
c_1\,Y_1f+c_2\,Y_2f,\,\,X_kf
\big\rightbracket
=
c_1\,\leftbracket
Y_1,\,X_k
\rightbracket
+
c_2\,\leftbracket
Y_2,\,X_k
\rightbracket
\]
can, in our case, be linearly expressed in terms
of the $X_if$. But
our family $e_1 \, X_1f + \cdots + e_q \, X_qf$ also
admits the infinitesimal transformation
$\leftbracket Y_1,\,Y_2 \rightbracket$. 
Indeed, one forms the Jacobi identity:
\[
\big\leftbracket
\leftbracket
Y_1,\,Y_2
\rightbracket,\,X_k
\big\rightbracket
+
\big\leftbracket
\leftbracket
Y_2,\,X_k
\rightbracket,\,Y_1
\big\rightbracket
+
\big\leftbracket
\leftbracket
X_k,\,Y_1
\rightbracket,\,Y_2
\big\rightbracket
=
0,
\]
and one takes into account that $\leftbracket Y_1,\,X_k \rightbracket$
and $\leftbracket Y_2,\,X_k \rightbracket$ can be linearly expressed
in terms of the $X_if$, so one realizes that this is also the case for
$\big\leftbracket \leftbracket Y_1, \, Y_2 \rightbracket, \, X_k
\big\rightbracket$.

By combining these two observations, one obtains the announced

\def\theproposition{1}\begin{proposition}
If the most general infinitesimal transformation
which leaves invariant a family of $\infty^{ q-1}$
infinitesimal transformations:
\[
e_1\,X_1f
+\cdots+
e_q\,X_qf
\]
can be linearly deduced from a bounded number of infinitesimal
transformations, say from $Y_1 f, \dots, Y_m f$, then the $Y_kf$
generate an $m$-term group.
\end{proposition}

This proposition can yet be generalized; 
indeed, it is evident that that the
totality of \emphasis{all finite} transformations
which leave invariant the family $\sum \, e_k \, X_k f$
always forms a group. 

\sectionengellie{\S\,\,\,70.}

Let
the conditions of the latter theorem be satisfied, 
namely let the family of the $\infty^{ q-1}$ infinitesimal
transformations $\sum \, e_k \, X_kf$
be invariant by all transformations of the one-term group $Y_f$. 

Now, according to Chap.~\ref{one-term-groups}, p.~\pageref{S-57}, 
the following holds true: 
if, after the
introduction of new
variables, the infinitesimal transformation $Xf$ is
transferred to $Zf$, then at the
same time, the transformations
of the one-term group $Xf$
are transferred to the transformations
of the one-term group $Zf$. So we
deduce that under the assumptions 
of the Theorem~43, not only the family of
the $\infty^{ q - 1}$ infinitesimal transformations
$e_1 \, X_1f + \cdots + e_r \, X_r f$
remains invariant, but also
the family of the $\infty^{ q-1}$
one-term groups generated by these
infinitesimal transformations, and
naturally also, the totality
of the $\infty^q$ finite transformations which
belong to these one-term groups. 

But we yet want to go a step further: 
we want to study how the analytic expression
of the individual finite transformations of the
one-term groups $e_1 X_1 f + \cdots + e_q \, X_q f$
behave, when the new variables:
\[
x_i'
=
x_i
+
t\,Yx_i
+\cdots
\]
are introduced in place of the $x$.

The answer to this question is given
by the Proposition~3 in Chap.~\ref{one-term-groups}, 
p.~\pageref{S-58}. Indeed, from this proposition, 
it easily results that, after
the introduction of the new variables
$x'$, every finite transformation:
\def\theequation{8}\begin{equation}
\aligned
\overline{x}_i
=
x_i
+
\sum_k^{1\cdots q}\,
e_k\,
&
X_k\,x_i
+
\sum_{k,\,\,j}^{1\cdots q}\,
\frac{e_k\,e_j}{1\cdot 2}\,
X_k\,X_j\,x_i
+\cdots
\\
&\ \ \
{\scriptstyle{(i\,=\,1\,\cdots\,n)}}
\endaligned
\end{equation}
receives the shape:
\def\theequation{9}\begin{equation}
\overline{x}_i'
=
x_i'
+
\sum_k^{1\cdots q}\,e_k'\,X_k'\,x_i'
+
\sum_{k,\,\,j}^{1\cdots q}\,
\frac{e_k'\,e_j'}{1\cdot 2}\,
X_k'\,X_j'\,x_i'
+\cdots,
\end{equation}
where the connection between the $e_k$ and the $e_k'$
is prescribed through the relation: 
\def\theequation{10}\begin{equation}
\sum_{k=1}^q\,e_k\,X_kf
+
\sum_{k=1}^q\,e_k'\,X_k'f.
\end{equation}
Consequently, we see directly that the concerned family of $\infty^q$
finite transformations in the new variables $x'$ possesses exactly the
same form as in the initial variables $x$.  But in addition, we remark
that a finite transformation which has the parameters $e_1, \dots,
e_q$ in the $x$ possesses, after the introduction of the new variables
$x'$, the parameters $e_1', \dots, e_q'$.

Now, as said just now, the connection between the $e$ and the $e'$
through the identity~\thetag{ 10} is completely prescribed; hence this
identity is absolutely sufficient when the question is to determine
the new form which an arbitrary finite transformation~\thetag{ 8}
takes after the transition to the $x'$.

In order to be as distinct as possible, we bring this circumstance to
expression when we \emphasis{interpret \label{S-255}
$\sum\, e_k\, X_kf$ virtually
as the symbol of the finite transformation:
\[
x_i'
=
x_i
+
\sum_{k=1}^q\,e_k\,X_k\,x_i
+\cdots
\ \ \ \ \ \ \ \ \ \ \ \ \ {\scriptstyle{(i\,=\,1\,\cdots\,n)}},
\]
where the absolute\footnote{\,
---\,namely the values themselves, but not the `absolute values'
$\vert e_k \vert$ in the modern sense\,---
} 
values of the $e_k$ then come into consideration, not only their
ratio}.  Then we can simply say:

\plainstatement{After the introduction of the new variables $x'$, the
finite transformation $\sum\, e_k\, X_kf$ is transferred to the finite
transformation $\sum\, e_k'\, X_k'f$.}

In the next studies of this chapter, the symbol
$\sum\, e_k\, X_kf$ will be employed 
sometimes as the symbol of a finite transformation, 
sometimes as the symbol of an infinitesimal 
transformation.
Hence in each individual case, we shall underline 
which one of the two interpretations of the symbol is meant.

\sectionengellie{\S\,\,\,71.}

Let the family of the $\infty^q$ transformations $e_1\, X_1f + \cdots
+ e_q \, X_q f$ remain invariant by all transformations of the
one-term group $Yf$. It can happen that $Yf$ itself is an 
infinitesimal transformation of the family $\sum \, e_k \, X_kf$; 
indeed, the case where $Yf$ is one arbitrary of the
$\infty^{ q-1}$ infinitesimal transformations $\sum\, e_k\, X_kf$
is of special interest. 
This will occur if, but also, only if
between the $X_kf$, there are relations of the form:
\[
\leftbracket
X_i,\,X_k
\rightbracket
=
\sum_{s=1}^q\,g_{iks}\,X_sf.
\]
Hence from the theorem of the preceding paragraph, we obtain
the following more special theorem in which
we permit ourselves to write $r$ instead $q$
and $c_{ iks}$ instead of $g_{ iks}$.

\def\thetheorem{44}\begin{theorem}
\label{Theorem-44-S-256}
For a family of $\infty^r$ finite transformations:
\[
e_1\,X_1f
+\cdots+
e_r\,X_rf
\ \ \ \ \ \ \
\text{or:}
\ \ \ \ \ \ \
\overline{x}_i
=
\mathfrak{P}_i
(x_1,\dots,x_n,\,e_1,\dots,e_r)
\]
to remain invariant by every transformation which
belongs to it\,---\,so that, after the introduction of the
new variables:
\[
x_i'
=
\mathfrak{P}_i(x_1,\dots,x_n,\,h_1,\dots,h_r),
\ \ \ \ \ \ \ 
\overline{x}_i'
=
\mathfrak{P}_i(\overline{x}_1,\dots,\overline{x}_n,\,h_1,\dots,h_r)
\]
in place of the $x$ and the $\overline{ x}$, it takes
the form:
\[
\overline{x}_i'
=
\mathfrak{P}_i(x_1',\dots,x_n',\,l_1,\dots,l_r)
\]
where the $l$ only depend upon $e_1, \dots, e_r$ and
$h_1, \dots, h_r$\,---\,it is necessary 
and sufficient that the $Xf$
stand pairwise in the relationships:
\[
\leftbracket
X_i,\,X_k
\rightbracket
=
\sum_{s=1}^r\,c_{iks}\,X_sf,
\]
where the $c_{ iks}$ are absolute constants.
\end{theorem}

This theorem states an important property that the family of the
$\infty^r$ finite transformations $\sum \, e_k \, X_kf$ possesses as
soon as relations of the form $\leftbracket X_i,\,X_k \rightbracket$
exist.  It is noteworthy that, for the proof of this theorem, we have
used the studies of the preceding chapter only for the smallest part;
moreover, we have used no more than a few developments of the
Chaps.~\ref{essential-parameters}, \ref{fundamental-differential},
\ref{one-term-groups} and~\ref{kapitel-8}. Namely, one should observe
that we have made no use of the Theorem~24, Chap.~\ref{kapitel-9},
p.~\pageref{Theorem-24-S-158}.

If one assumes that the latter theorem is known, then one can shorten
the proof of the Theorem~44 as follows: one shows at first, as above,
that the relations $\leftbracket X_i,\,X_k \rightbracket = \sum\,
c_{iks}\, X_sf$ are necessary; then from Theorem~24,
p.~\pageref{Theorem-24-S-158} it comes that the $\infty^{ r-1}$
infinitesimal transformations $\sum\, e_k\, X_kf$ generate an $r$-term
group.  If $x_i' = \mathfrak{ P}_i (x_1, \dots, x_n, \, h_1, \dots,
h_r)$ are the finite equations of this group, then according to
Theorem~5, p.~\pageref{Theorem-5-S-45}, there is an identity of the
form:
\[
\sum_{k=1}^r\,e_k\,X_kf
=
\sum_{k=1}^r\,e_k'\,X_k'f\,;
\]
with this, the proof of the Theorem~44 is produced.

One does not even need to refer to Theorem~5
p.~\pageref{Theorem-5-S-45}, but one can conclude in the following
way:

The equations $x_i' = \mathfrak{ P}_i ( x, h)$ of our group, when
resolved, give a transformation of the shape:
\[
x_i
=
\mathfrak{P}_i
\big(x_1',\dots,x_n',\,\chi_1(h),\dots,\chi_r(h)\big),
\]
that is to say, the transformation which is inverse to the
transformation with the parameters $h_1, \dots, h_r$.  Now, if one
imagines that these values of the $x_i$ are inserted in the equations
$\overline{ x}_i = \mathfrak{ P}_i ( x, e)$ and if one takes into
consideration that one has to deal with a group, then one realizes
that there exist certain equations of the form:
\[
\overline{x}_i
=
\mathfrak{P}_i
\big(x_1',\dots,x_n',\,\psi_1(h,e),\dots,\psi_r(h,e)\big).
\]

Lastly, if one inserts these expressions for the $\overline{ x}_i$ in
the equations $\overline{ x}_i' = \mathfrak{ P}_i (\overline{ x}, h)$,
then one obtains:
\[
\overline{x}_i'
=
\mathfrak{P}_i(x_1',\dots,x_n',\,l_1,\dots,l_r).
\]
This is the new form that the transformations $\overline{ x}_i =
\mathfrak{ P}_i ( x, e)$ take after the introduction of the new
variables $x'$.  Here evidently, the $l$ are functions of only the $e$
and the $h$, exactly as it is claimed by the Theorem~44,
p.~\pageref{Theorem-44-S-256}.

As a result, the connection which exists between the Theorem~24,
p.~\pageref{Theorem-24-S-158} and the Theorem~44 of the
present chapter is clarified.

\sectionengellie{\S\,\,\,72.}

In order to be able to state the gained results more briefly, or, if
one wants, more clearly, we will, as earlier on, translate the
symbolism of the theory of substitutions to the theory of
transformation groups.

We shall denote all finite transformations $e_1\, X_1f + \cdots + e_r
X_r f$ by the common symbol $T$, and the individual transformations by
marking an appended index, so that for instance
the symbol $T_{ (a)}$ denotes the finite transformation:
\[
a_1\,X_1f
+\cdots+
a_r\,X_rf.
\]

Using this way of expressing, we can 
at first state the Theorem~24, p.~\pageref{Theorem-24-S-158}
as follows:

\def\theproposition{2}\begin{proposition}
If $r$ independent infinitesimal transformations
$X_1f, \dots, X_rf$ stand pairwise in
relationships of the form:
\[
\leftbracket
X_i,\,X_k
\rightbracket
=
\sum_{s=1}^r\,c_{iks}\,X_sf,
\]
then the family of all finite transformations $\sum\, e_k\, X_kf$, or
$T_{ (e)}$, contains, simultaneously with the two transformations $T_{
(a)}$ and $T_{ (b)}$, also the transformation $T_{ (a)} \, T_{ (b)}$;
hence there is a symbolic equation of the form:
\[
T_{(a)}\,T_{(b)}
=
T_{(c)},
\]
in which the parameters $c$ are functions of the
$a$ and of the $b$. 
\end{proposition}

Correspondingly, from the Theorem~44 of the present chapter, we obtain
the following proposition, which, however, does not exhaust the
complete content of the theorem:

\def\theproposition{3}\begin{proposition}
If $r$ independent infinitesimal transformations $X_1f, \dots, X_rf$
stand pairwise in the relationships $\leftbracket X_i,\,X_k
\rightbracket$, then the family of the $\infty^r$ finite
transformations $\sum\, e_k\, X_kf$, or $T_{ (e)}$, contains,
simultaneously with the transformations $T_{ (a)}$ and $T_{ (b)}$,
also the transformation $T_{ (a)}^{-1} \, T_{ (b)}\, T_{ (a)}$, whence
there exists an equation of the form:
\[
T_{(a)}^{-1}\,T_{(b)}\,T_{(a)}
=
T_{(c')}
\]
in which the parameters $c'$ are functions of the $a$ and of the $b$.
\end{proposition}

Obviously, the existence of the symbolic relation: $T_{ (a)}^{ -1}\,
T_{(b)}\, T_{(a)} = T_{(c')}$ is a consequence of the former relation:
$T_{ (a)}\, T_{ (b)} = T_{ (c)}$, hence the latter proposition is also
a consequence of the preceding one, as we already have seen in the
previous paragraphs. 

Finally, by combining the two Theorems~44 and~24, 
p.~\pageref{Theorem-24-S-158}, we yet obtain the following
the following curious result. 

\def\thetheorem{45}\begin{theorem}
If a family of $\infty^r$ finite transformations: $a_1\, X_1f + \cdots
+ a_r\, X_rf$, or shortly $T_{ (a)}$, possesses the property that the
transformation $T_{ (a)}^{-1} \, T_{ (b)} \, T_{ (a)}$ always belongs
to the family, whatever values the parameters $a_1, \dots, a_r$, $b_1,
\dots, b_r$ can have, then the family of $\infty^r$ transformations in
question forms an $r$-term group, that is to say: $T_{ (a)} \, T_{
(b)}$ is always a transformation which also belongs to the family.
\end{theorem}

\sectionengellie{\S\,\,\,73.}

\label{S-258}
If the transformation $T_{ (a)}^{ -1} \, T_{ (b)} \, T_{ (a)}$
coincides with the transformation $T_{ (b)}$, a fact that we express
by means of the symbolic equation:
\[
T_{(a)}^{-1}\,T_{(b)}\,T_{(a)}
=
T_{(b)},
\]
then we say:
\terminology{the transformation $T_{(b)}$ remains
invariant by the transformation $T_{(a)}$}.

But in this case, we also have:
\[
T_{(b)}^{-1}\,T_{(a)}\,T_{(b)}
=
T_{(a)},
\]
\label{S-259}
that is to say, the transformation $T_{ (a)}$ remains invariant by the
transformation $T_{ (b)}$; on the other hand:
\[
T_{(a)}\,T_{(b)}
=
T_{(b)}\,T_{(a)}
\]
is an equation which expresses that the two transformations $T_{ (a)}$
and $T_{ (b)}$ are interchangeable one with another.

We already remarked in Theorem~6, p.~\pageref{Theorem-6-S-49} that the
transformations of an arbitrary one-term group are interchangeable by
pairs. At present, we can also settle the more general question of
when the transformations of two \emphasis{different} one-term groups
$Xf$ and $Yf$ are interchangeable one with another.

In the general finite transformation $e\, Xf$ of the one-term
group $Xf$, we introduce the new variables $x_i'$
which are defined by the finite equations $x_i' = 
x_i + t\, Y x_i + \cdots$ of the one-term group $Yf$. 
The transformations of our two one-term groups will
then be interchangeable
if and only if every transformation of
the form $e\, Xf$
remains invariant after the introduction of the $x'$, 
whence $e\, Xf$ is equal to $e\, X'f$. 

According to the Theorem~43, p.~\pageref{Theorem-43-S-252}, 
for the existence of an equation of the form:
\[
e\,Xf
=
e'\,Xf,
\]
it is necessary and sufficient that the infinitesimal
transformations $Xf$ and $Yf$ satisfy a relation:
\[
\leftbracket
Y,\,X
\rightbracket
=
g\,Xf\,;
\]
at the same time, the $e'$ determines itself
through the differential equation:
\[
\frac{\D\,e'}{\D\,t}
+
g\,e'
=
0.
\]

But in our case, $e'$ is supposed, for every value of $t$, to be equal
to $e$, hence $e'$ depends absolutely not on $t$, that is to say, the
differential quotient $\D\, e' / \D\, t$ vanishes, and with it, the
quantity $g$ too. At the same time, this condition is 
evidently necessary and sufficient. Consequently, 
the following holds true:

\def\theproposition{4}\begin{proposition}
The finite transformations of two one-term groups
$Xf$ and $Yf$ are interchangeable one with another
if and only if the expression 
$\leftbracket X,\,Y \rightbracket$ vanishes identically.
\end{proposition}

It stands to reason to call \terminology{interchangeable two
infinitesimal transformations $Xf$ and $Yf$ which stand in the
relationship $\leftbracket X,\,Y \rightbracket \equiv 0$}.

If we introduce this way of expressing, we can say that the finite
transformations of two one-term groups are interchangeable by pairs if
and only if the infinitesimal transformations of the two groups are
so.

Moreover, from the latter proposition, it yet comes the

\renewcommand{\thefootnote}{\fnsymbol{footnote}}
\def\thetheorem{46}\begin{theorem}
The finite transformations of an $r$-term group $X_1f, \dots, X_rf$
are pairwise interchangeable if and only if all expressions
$\leftbracket X_i,\,X_k \rightbracket$ vanish identically, or stated
differently, if and only if the infinitesimal transformations $X_1f,
\dots, X_r f$ are interchangeable by pairs\footnote[1]{\,
\name{Lie}, Gesellschaft der Wissenschaften zu Christiania 1872;
Archiv for Mathematik og Naturvidenskab Vol.~8, p.~180, 1882;
Math. Annalen Vol.~24, p.~557, 1884.
}. 
\end{theorem}
\renewcommand{\thefootnote}{\arabic{footnote}}

\sectionengellie{\S\,\,\,74.}

From the general developments of the \S\S\,\,69 and~70, we will now
yet draw a few further consequences that are of importance.

We again assume that the $r$ independent infinitesimal transformations
$X_1f, \dots, X_r f$ stand pairwise in the relationships:
\[
\leftbracket
X_i,\,X_k
\rightbracket
=
\sum_{s=1}^r\,c_{iks}\,X_sf,
\]
whence according to the Theorem~24 of the Chap.~\ref{kapitel-9},
p.~\pageref{Theorem-24-S-158}, the $\infty^r$ finite transformations
$\sum\, e_k\, X_kf$ form an $r$-term group.  Now, if the family of the
transformations $\sum\, e_k\, X_kf$ remains invariant by all
transformations of a one-term group $Yf$, then according to
Theorem~43, p.~\pageref{Theorem-43-S-252}, there exist $r$ equations
of the form: 
\def\theequation{4}\begin{equation}
\leftbracket
Y,\,X_k
\rightbracket
=
\sum_{s=1}^r\,g_{ks}\,X_sf.
\end{equation}

These equations show that the $r+1$ infinitesimal transformations
$X_1f, \dots, X_rf$, $Yf$ generate an $(r+1)$-term group to which
$X_1f, \dots, X_rf$ belongs as an $r$-term subgroup.

In addition, it is clear that the relations~\thetag{ 4} do not
essentially change their form when one inserts in place of $Yf$ a
completely arbitrary infinitesimal transformation of the $(r+1)$-term
group.

Hence, if $x_i' = \psi_i ( x_1, \dots, x_n, \, a_1, \dots, a_{ r+1})$
are the finite equations of the $(r+1)$-term group, then the family of
the finite transformations $\sum\, e_k\, X_kf$ remains invariant when
one introduces in the same way the new variables $x'$ in place of the
$x$, that is to say, only the parameters vary in the analytic
expression of the $r$-term group.

A similar property would hold true if, instead of the single
transformation $Yf$, one would have several, say $m$, such
transformations all of which would satisfy relations of the
form~\thetag{ 4}; in addition, we yet want to add the assumption that
these $m$ infinitesimal transformations $Y_1 f, \dots, Y_m f$,
together with $X_1f, \dots, X_rf$, generate an $(r+m)$-term group.
Then if we introduce new variables in the group $X_1f, \dots, X_rf$ by
means of an arbitrary transformation of the $(r+m)$-term group, the
family of finite transformations of our $r$-term group remains
invariant, although the parameters are changed in their analytic
representation.

We want to express this relationship between the two groups by saying
shortly: \label{S-261}
\terminology{if the $r$-term group $X_1f, \dots, X_rf$
remains invariant by all transformations of the $(r+m)$-term group, it
is an invariant subgroup of it}.

If we translate the way of expressing used commonly in the theory of
substitutions, we can also interpret as follows the definition of the
invariants subgroups:

If $T$ is the symbol of an arbitrary transformation of the
$(r+m)$-term group $G$, and if $S$ is an arbitrary transformation of a
subgroup of $G$, then this subgroup is invariant in $G$ when the
transformation $T^{ -1} \, S \, T$ always belongs also to the subgroup
in question.

In what has been said above, the analytic conditions for the
invariance of a subgroup are completely exhibited; so we need only to
summarize them once again:

\def\thetheorem{47}\begin{theorem}
If the $r$-term group $X_1f, \dots, X_rf$ is contained in an
$(r+m)$-term group $X_1f, \dots, X_rf, \, Y_1f, \dots, Y_m f$, then it
is an invariant subgroup of it when every $\leftbracket Y_i,\,X_k
\rightbracket$ expresses in terms of $X_1f, \dots, X_rf$ linearly with
constant coefficients.
\end{theorem}

On can give the Theorem~47 in the following more 
general\,---\,though only in a formal sense\,---\,version:

\def\theproposition{5}\begin{proposition}
\label{Satz-5-S-261}
If the totality of all infinitesimal transformations
$e_1\, X_1f + \cdots + e_m\, X_mf$ forms an invariant
family in the $r$-term group $X_1f, \dots, X_m f, \dots, 
X_rf$, then $X_1f, \dots, X_mf$ generate an $m$-term invariant
subgroup of the $r$-term group.
\end{proposition}

Indeed, since all $\leftbracket X_i,\,X_k \rightbracket$, in which $i$
is $\leqslant m$, can be linearly deduced from $X_1f, \dots, X_mf$,
then in particular the same holds true of all $\leftbracket X_i,\,X_k
\rightbracket$ in which both $i$ and $k$ are $\leqslant m$.  As a
result, $X_1f, \dots, X_m f$ generate an $m$-term subgroup to which
the Theorem~47 can immediately be applied.

\medskip

At first, a few examples of invariant subgroups. 

\renewcommand{\thefootnote}{\fnsymbol{footnote}}
\def\theproposition{6}\begin{proposition}\footnote[1]{\,
In the Archiv for Mathematik og Naturvidenskab Vol.~8, p.~390,
Christiania 1883, \name{Lie} observed that the $\leftbracket X_i,\,X_k
\rightbracket$ form an invariant subgroup.  \name{Killing} has
realized this independently in the year 1886. }
If \label{Satz-6-S-261}
the $r$ independent infinitesimal transformations $X_1f, \dots,
X_rf$ generate an $r$-term group, then the totality of all
infinitesimal transformations $\leftbracket X_i,\,X_k \rightbracket$
also generate a group; if the latter group contains $r$ parameters,
then it is identical to the group $X_1f, \dots, X_rf$; if it contains
less than $r$ parameters, then it is an invariant subgroup of the
group $X_1f, \dots, X_rf$; if one adds to the $\leftbracket X_i,\,X_k
\rightbracket$ arbitrarily many mutually independent infinitesimal
transformations $e_1\, X_1f + \cdots + e_r\,X_rf$ that are also
independent of the $\leftbracket X_i,\,X_k \rightbracket$, then one
always obtains again an invariant subgroup of the group $X_1f, \dots,
X_rf$.
\end{proposition}
\renewcommand{\thefootnote}{\arabic{footnote}}

It is clear that the $\leftbracket X_i,\,X_k \rightbracket$ can at
most generate an $r$-term group, since they all belong to the group
$X_1f, \dots, X_rf$; the fact that they effectively generate a group
comes immediately from the relations:
\[
\leftbracket
X_i,\,X_k
\rightbracket
=
\sum_{s=1}^r\,c_{iks}\,X_sf,
\]
for one indeed has:
\[
\big\leftbracket
\leftbracket
X_i,\,X_k
\rightbracket,\,\,
\leftbracket
X_j,\,X_l
\rightbracket
\big\rightbracket
=
\sum_{s,\,\sigma}^{1\cdots r}\,
c_{iks}\,c_{jl\sigma}\,
\leftbracket
X_s,\,X_\sigma
\rightbracket.
\]

The claim that the group generated by the $\leftbracket X_i,\,X_k
\rightbracket$ in the group $X_1f, \dots, X_rf$ is invariant becomes
evident from the equations:
\[
\big\leftbracket
X_j,\,
\leftbracket
X_i,\,X_k
\rightbracket
\big\rightbracket
=
\sum_{s=1}^r\,c_{iks}\,
\leftbracket
X_j,\,X_s
\rightbracket.
\]
The last part of the proposition does not require any closer
explanation.

The following generalization of the proposition just proved is
noteworthy:

\def\theproposition{7}\begin{proposition}
\label{Satz-7-S-262}
If $m$ infinitesimal transformations $Z_1f, \dots, Z_mf$ of the
$r$-term group $X_1f, \dots, X_rf$ generate an $m$-term subgroup of
this group, and if the subgroup in question is invariant in the
$r$-term group, then the subgroup which is generated by all the
infinitesimal transformations $\leftbracket Z_\mu,\,Z_\nu
\rightbracket$ is also invariant in the $r$-term group.
\end{proposition}

The proof of that is very simple.  Under the assumptions of the
proposition, there are relations of the shape:
\[
\leftbracket
X_k,\,Z_\mu
\rightbracket
=
\sum_{\lambda=1}^m\,
h_{k\mu\lambda}\,Z_\lambda f
\ \ \ \ \ \ \ \ \ \ \ \ \ 
{\scriptstyle{(k\,=\,1\,\cdots\,r\,;
\,\,\,\mu\,=\,1\,\cdots\,m)}},
\]
where the $h_{k\mu\lambda}$ are constants. 
Next, if we form the Jacobi identity
(cf. Chap.~\ref{kapitel-5}, p.~\pageref{jacobi-identity}):
\[
\big\leftbracket
X_k,\,
\leftbracket
Z_\mu,\,Z_\nu
\rightbracket
\big\rightbracket
+
\big\leftbracket
Z_\mu,\,
\leftbracket
Z_\nu,\,X_k
\rightbracket
\big\rightbracket
+
\big\leftbracket
Z_\nu,\,
\leftbracket
X_k,\,Z_\mu
\rightbracket
\big\rightbracket
=
0, 
\] 
and if we insert in it the expressions written 
above for $\leftbracket Z_\mu,\,X_k \rightbracket$ and
$\leftbracket Z_\nu,\,X_k \rightbracket$, we obtain
the equations:
\[
\big\leftbracket
X_k,\,
\leftbracket
Z_\mu,\,Z_\nu
\rightbracket
\big\rightbracket
=
\sum_{\lambda=1}^m\,
\big\{
h_{k\nu\lambda}\,
\leftbracket
Z_\mu,\,Z_\lambda
\rightbracket
-
h_{k\mu\lambda}\,
\leftbracket
Z_\nu,\,Z_\lambda
\rightbracket
\big\},
\]
from which it comes that the subgroup of the group $X_1f, \dots, X_rf$
generated by the $\leftbracket Z_\mu,\,Z_\nu \rightbracket$ is
invariant in the former group.  But this is what was to be
proved.\,---

Let the $r$-term group $X_1f, \dots, X_rf$, or shortly $G_r$, contain
an $(r-1)$-term invariant subgroup and let $Y_1f, \dots, Y_rf$ be $r$
independent infinitesimal transformations of the $G_r$ selected in a
such a way that $Y_1f, \dots, Y_{ r-1}f$ is this invariant subgroup.
Then there exist relations of the form\footnote{\,
Indeed, the invariance yields such relations for all $i = 1, \dots,
r-1$ and all $k = 1, \dots, r$, and also by skew-symmetry for all $i =
1, \dots, r$ and all $k = 1, \dots, r-1$; it yet remains only
$\leftbracket Y_n,\,Y_n \rightbracket$ which, anyway, is zero. 
}: 
\[
\leftbracket
Y_i,\,Y_k
\rightbracket
=
c_{ik1}\,Y_1f
+\cdots+
c_{ik,r-1}Y_{r-1}f
\ \ \ \ \ \ \ \ \ \ \ \ \ {\scriptstyle{(i,\,\,k\,=\,1\,\cdots\,r)}},
\]
whence all $\leftbracket Y_i,\,Y_k \rightbracket$ and also all
$\leftbracket X_i,\,X_k \rightbracket$ can be linearly deduced from
only $Y_1f, \dots, Y_{ r-1}f$.  From this, we conclude that every
$(r-1)$-term invariant subgroup of the $G_r$ contains all the
infinitesimal transformations $\leftbracket X_i,\,X_k \rightbracket$
and so, we realize that the following proposition holds true:

\def\theproposition{8}\begin{proposition}
\label{Satz-8-S-263}
In the $r$-term group $X_1f, \dots, X_rf$, there is an $(r-1)$-term
invariant subgroup if and only if the infinitesimal transformations
$\leftbracket X_i,\,X_k \rightbracket$ generate a group with less than
$r$ parameters; if there are, amongst the $\leftbracket X_i,\,X_k
\rightbracket$ exactly $r_1 < r$ mutually independent infinitesimal
transformations, then one obtains all $(r-1)$-term invariant subgroups
of the group $X_1f, \dots, X_rf$ by adding to the $\leftbracket
X_i,\,X_k \rightbracket$ in the most general way $r - r_1 - 1$
infinitesimal transformations $e_1\, X_1f + \cdots + e_r\, X_rf$ which
are mutually independent and are independent of the $\leftbracket
X_i,\,X_k \rightbracket$.
\end{proposition}

\smallskip

The Proposition~1 in Chap.~\ref{kapitel-12}, p.~\pageref{Satz-1-S-205}
provides us with another example of invariant subgroup.

Namely, if, in the neighbourhood of a point $x_1^0, \dots, x_n^0$, a
group contains infinitesimal transformations of first or of higher
order in the $x_i - x_i^0$, then each time, all infinitesimal
transformations of order $k$ ($k >0$) and higher generate a subgroup.
Now, by bracketing \deutsch{Klammeroperation}, two infinitesimal
transformations of respective orders $k$ and $k+\nu$ produce a
transformation $\leftbracket X,\,Y \rightbracket$ of the order $2k +
\nu -1$.  Because $k$ is $>0$, one has $2k + \nu -1 \geqslant k +
\nu$, hence $\leftbracket X,\,Y \rightbracket$ must be linearly
expressible in terms of the infinitesimal transformations of orders
$k+\nu$ and higher.  In other words: all the infinitesimal
transformations of orders $k+\nu$ and higher generate a group which is
invariant in the group generated by the infinitesimal transformations
of orders $k$ and higher.  And as said, all of that holds true under
the only assumption that the number $k$ is larger than zero.

In particular, if the numbers $k$ and $k+\nu$ can be chosen in such a
way that the group contains no infinitesimal transformation of orders
$(2k+\nu-1)$ and higher in the neighbourhood of $x_1^0, \dots, x_n^0$,
then the expression $\leftbracket X,\,Y \rightbracket$ must vanish
identically.  Hence the following holds true:

\def\theproposition{9}\begin{proposition}
\label{Satz-9-S-264}
If, in the neighbourhood of a point $x_1^0, \dots, x_n^0$, a group
contains no infinitesimal transformation of $(s+1)$-th order, or of
higher order, and if, by contrast, it contains transformations of
$k$-th order, where $k$ satisfies the condition $2k-1 > s$, then all
infinitesimal transformations of the group which are of the $k$-th
order and of higher order generate a group with pairwise 
interchangeable transformations.
\end{proposition}

Here, the point $x_1^0, \dots, x_n^0$ needs absolutely not be such
that the coefficients of the resolved defining equations of the group
(cf. Chap.~\ref{kapitel-11}, p.~\pageref{S-189-sq} sq.)  behave
regularly.\,---

Consequently \deutsch{consequenterweise} one must say that each finite
continuous group contains two invariant subgroups, namely firstly
itself and secondly the identity transformation.  One realizes this by
setting $m$, in the Theorem~47, firstly equal to $r$, and secondly
equal to zero; in the two cases the condition for the invariance of
the subgroup $X_1f, \dots, X_rf$ is satisfied by itself, only as soon
as $X_1f, \dots, X_rf$, $Y_1f, \dots, Y_mf$ is an $(r+m)$-term group.

The groups which contain absolutely no invariant subgroup,
disregarding the two which are always present, are of special
importance.  That is why these groups are also supposed to have a
special name, and they should be called \label{S-264}
\terminology{simple}
\deutsch{einfach}.  In contrast to this, a group is called
\terminology{compound} \deutsch{zusammengesetzt} when, aside from the
two invariant subgroups indicated above, it yet contains other
invariant subgroups.

In conclusion, here are still two propositions about invariant
subgroups:

\def\theproposition{10}\begin{proposition}
\label{Satz-10-S-264}
The transformations which are common to two invariant subgroups of a
group $G$ form in the same way a subgroup which is invariant in $G$.
\end{proposition}

The transformations in question certainly form a subgroup of $G$
(cf. Chap.~\ref{kapitel-9}, Proposition~2, p.~\pageref{Satz-2-S-159});
this subgroup must be invariant in $G$, since by all transformations
of $G$, it is transferred to a group which belongs to the two
invariant subgroups, that is to say, to itself.

\def\theproposition{11}\begin{proposition}
\label{Satz-11-S-264}
If two invariant subgroups $Y_1f, \dots, Y_mf$ and $Z_1f, \dots, Z_pf$
of a group $G$ have no infinitesimal transformations in common, then
all expressions $\leftbracket Y_i,\,Z_k \rightbracket$ vanish
identically, that is to say, every transformation of one subgroup is
interchangeable with every other transformation of the other subgroup.
\end{proposition}

Indeed, under the assumptions made, every expression $\leftbracket
Y_i,\,Z_k \rightbracket$ must be expressible with constant
coefficients both in terms of $Y_1f, \dots, Y_mf$ and in terms of
$Z_1f, \dots, Z_p f$; but since the two subgroups have no
infinitesimal transformations in common, then the thing is nothing
else but that all expressions $\leftbracket Y_i,\,Z_k \rightbracket$
vanish identically.  The rest follows from the Proposition~4,
p.~\pageref{Satz-4-S-259}.

\sectionengellie{\S\,\,\,75.}

\renewcommand{\thefootnote}{\fnsymbol{footnote}}
There are $r$-term groups for which one can select $r$ independent
infinitesimal transformations: $Y_1f, \dots, Y_rf$ so that for every
$i < r$, the $i$ independent infinitesimal transformations $Y_1f,
\dots, Y_if$ generate an $i$-term group which is invariant in the
$(i+1)$-term group $Y_1f, \dots, Y_{ i+1} f$.  Then between $Y_1f,
\dots, Y_r f$, there are relations of the form:
\def\theequation{11}\begin{equation}\label{S-265}
\aligned
\leftbracket
Y_i,\,Y_{i+k}
\rightbracket
&
=
\overline{c}_{i,i+k,1}\,Y_1f
+\cdots+
\overline{c}_{i,i+k,i+k-1}\,Y_{i+k-1}f
\\
&
\ \ \ \ \ \ \
{\scriptstyle{(i\,=\,1\,\cdots\,r\,-\,1\,;\,\,\,
k\,=\,1\,\cdots\,r\,-\,i)}}.
\endaligned
\end{equation}
In the integration theory of those systems of differential equations
which admit finite groups, it comes out that the groups of the just
defined specific constitution 
occupy a certain outstanding position in
comparison to all other groups\footnote[1]{\,
\name{Lie}, Ges. der Wiss. zu Christiania, 1874, p.~273. Math.
Ann. Vol.~XI, p.~517 and 518. Archiv for Math. og Nat.  Vol.~3, 1878,
p.~105 sq., Vol.~8, 1883.
}. 
\renewcommand{\thefootnote}{\arabic{footnote}}

Later, in the chapter about linear homogeneous groups, we will occupy
ourselves more accurately with this special category
of groups; at
present, we want only to show in which way one can determine whether a
given $r$-term group $X_1f, \dots, X_rf$ belongs, or does not belong,
to the category in question.

In order that it be possible to choose, amongst the infinitesimal
transformations $e_1\, X_1f + \cdots + e_r\, X_r f$, $r$ mutually
independent ones: $Y_1f, \dots, Y_rf$ which stand in relationships of
the form~\thetag{ 11}, there must above all exist an $(r-1)$-term
invariant subgroup in the $G_r$: $X_1f, \dots, X_rf$.  Thanks to
Proposition~8, p.~\pageref{Satz-8-S-263}, we are in a position to
determine whether this is the case: according to this proposition, the
group $X_1f, \dots, X_rf$ contains an $(r-1)$-term invariant subgroup
only when the group generated by all $\leftbracket X_i,\,X_k
\rightbracket$ contains less than $r$ parameters, say $r_1$; if this
condition is satisfied, then one obtains all $(r-1)$-term invariant
subgroups of the $G_r$ by adding to the $\leftbracket X_i,\,X_k
\rightbracket$ in the most general way $r - r_1 -1$ infinitesimal
transformations $e_1\, X_1f + \cdots + e_r\, X_rf$ that are mutually
independent and are independent of the $\leftbracket X_i,\,X_k
\rightbracket$.

However, not every $r$-term group which contains an $(r-1)$-term
invariant subgroup does belong to the specific category defined above;
for it to belong to this category, it must contain an $(r-1)$-term
invariant subgroup, which in turn must contain an $(r-2)$-term invariant
subgroup, and again the latter subgroup must contain an $(r-3)$-term
invariant subgroup, and so on.

From this, it follows how we have to proceed with the group $X_1f,
\dots, X_rf$: amongst all $(r-1)$-term invariant subgroups of the
$G_r$, we must select those which contain at least an $(r-2)$-term
invariant subgroup and we must determine all their $(r-2)$-term
invariant subgroups; according to what precedes, this presents no
difficulty.  Afterwards, amongst the found $(r-2)$-term subgroups, we
must select those which contain $(r-3)$-term invariant subgroups, and
so forth.

It is clear that in this manner, we arrive at the answer
to the question which we have asked about the group $X_1f, \dots, 
X_rf$. Either we realize that this group does not belong
to the discussed specific category, or we find $r$ independent
infinitesimal transformations $Y_1f, \dots, Y_rf$ of our group
which are linked by relations of the form~\thetag{ 11}.

\medskip

In certain circumstances, the computations just indicated will be made
more difficult by the fact that the subgroups to be studied contain
arbitrary parameters which specialize themselves in the course of the
computations. So we will indicate yet another process which conducts
as well to answering our question, but which avoids all computations
with arbitrary parameters.

We want to suppose that amongst the $\leftbracket X_i,\,X_k
\rightbracket$, one finds exactly $r_1 \leqslant r$ 
\label{S-266} which are
independent and that all $\leftbracket X_i,\,X_k \rightbracket$ can be
linearly deduced from $X_1'f, \dots, X_{ r_1}'f$, and furthermore
correspondingly, that all $\leftbracket X_i',\,X_k' \rightbracket$ can
be linearly deduced from the $r_2 \leqslant r_1$ independent
transformations $X_1''f, \dots, X_{ r_2}''f$, all $\leftbracket
X_i'',\,X_k'' \rightbracket$ from the $r_3 \leqslant r_2$ independent
$X_1'''f, \dots, X_{ r_3}'''f$, and so on.  Then according to
Proposition~6, p.~\pageref{Satz-6-S-261}, $X_1'f, \dots, X_{ r_1}'f$
generate an $r_1$-term invariant subgroup $G_{ r_1}$ of the $G_r$:
$X_1f, \dots, X_rf$, and furthermore $X_1''f, \dots, X_{ r_2}''f$
generate an $r_2$-term invariant subgroup $G_{ r_2}$ of the $G_{
r_1}$, and so forth; briefly, we obtain a series of subgroups $G_{
r_1}$, $G_{ r_2}$, $G_{r_3}$, \dots, of the $G_r$ in which each
subgroup is contained in all the preceding ones and is invariant in
the immediately preceding one, in all cases.  But now, according to
the Proposition~7, p.~\pageref{Satz-7-S-262}, the $G_{ r_2}$ is at
first invariant not only in the $G_{ r_1}$, but also in the $G_r$, and
furthermore, according to the same proposition, the $G_{ r_3}$ is not
only invariant in the $G_{ r_2}$, but also in the $G_{ r_1}$ and even
in the $G_r$ itself, and so forth. One sees that in the series of the
groups: $G_r$, $G_{ r_1}$, $G_{ r_2}$, \dots, each individual group is
contained in all the groups preceding, and is invariant in all the
groups preceding.

In the series of the entire numbers $r$, $r_1$, $r_2$, \dots, there is
none which is larger than the preceding one, and on the
other hand, none which is smaller than zero. 
Consequently, there must exist a positive number $q$
of such a nature that $r_{ q+1}$ is equal to 
$r_q$, while for $j < q$, it always holds true that: 
$r_{ j+1} < r_j$. Evidently, one then has:
\[
r_q
=
r_{q+1}
=
r_{q+2}
=
\cdots,
\]
so actually:
\[
r_{q+k}
=
r_q
\ \ \ \ \ \ \ \ \ \ \ \ \ {\scriptstyle{(k\,=\,1,\,2\,\cdots)}}.
\]

Now, there are two cases to be distinguished, according to whether the
number $r_q$ has the value zero, or is larger than 
zero\footnote{\,
In this case, $r_q \geqslant 2$ in fact, for a one-term
group is always solvable.
}. 

In the case $r_q = 0$, it is always possible, as we will show, to
select $r$ independent infinitesimal transformations $Y_1f, \dots,
Y_rf$ of the $G_r$ which stand in relationships of the form~\thetag{
11}.

Indeed, we choose as $Y_rf$, $Y_{ r-1} f$, \dots, $Y_{ r_1+1}f$ any $r
- r_1$ infinitesimal transformations $e_1\, X_1f + \cdots + e_r\,
X_rf$ of the $G_r$ that are mutually independent and are independent
of $X_1'f, \dots, X_{ r_1}'f$; as $Y_{ r_1}f$, $Y_{ r_1 -1}f$, \dots,
$Y_{ r_2+1}f$, we choose any $r_1 - r_2$ infinitesimal transformations
$e_1' X_1'f + \cdots + e_{ r_1}' X_{ r_1}'f$ of the $G_{ r_1}$ that
are mutually independent and are independent of $X_1''f, \dots, X_{
r_2}''f$, etc.; lastly, as $Y_{ r_{q-1}}f$, $Y_{ r_{q-1}-1}f$,
\dots, $Y_1f$, we choose any $r_{ q-1}$ independent infinitesimal
transformations $e_1^{ (q-1)} X_1^{ (q-1)}f + \cdots + e_{ r_{
q-1}}^{(q-1)} X_{ r_{ q-1}}^{ (q-1)} f$ of the $G_{ r_{ q-1}}$.
In this way, we obviously obtain $r$ mutually independent
infinitesimal transformations $Y_1f, \dots, Y_rf$ of our $G_r$; about
them, we claim that for every $i < r$, the $i$ transformations $Y_1f,
\dots, Y_if$ generate an $i$-term group which is invariant in the
$(i+1)$-term group: $Y_1 f, \dots, Y_{ i+1}f$. If we succeed to prove
this claim, then we have proved also that $Y_1f, \dots, Y_rf$ stand in
the relationships~\thetag{ 11}.

Let $j$ be any of the numbers $1, 2, \dots, q$.  Then it is clear that
the $r_j$ mutually independent infinitesimal transformations $Y_1f,
\dots, Y_{ r_j}f$ generate an $r_j$-term group, namely the group $G_{
r_j}$ defined above; certainly, it must be remarked that in the case
$j = q$, the $G_{ r_j}$ reduces to the identity transformation.

As we already have remarked earlier on, according to Proposition~6,
p.~\pageref{Satz-6-S-261},
the $G_{ r_j}$ is invariant in the group
$G_{ r_{ j-1}}$; 
\label{S-267}
but from this proposition, one can conclude even
more, namely one can conclude that we always obtain an invariant
subgroup of the $G_{ r_{ j-1}}$ when, to the infinitesimal
transformations $Y_1f, \dots, Y_{ r_j}f$ of the $G_{ r_j}$, we add any
infinitesimal transformations of the group $G_{ r_{ j-1}}$ that are
mutually independent and are independent of $Y_1f, \dots, Y_{ r_j}f$.
Now, since $Y_{ r_j+1}f, \dots, Y_{ r_{ j-1}}f$ belong to the $G_{ r_{
j-1}}$, and furthermore, since they are mutually independent and
independent of $Y_1f, \dots, Y_{ r_j}f$, then it comes that each one
of the following systems of infinitesimal transformations:
\[
\aligned
Y_1f,\dots,Y_{r_j}f,\
&
Y_{r_j+1}f,\dots,Y_{r_{j-1}-2}f,Y_{r_{j-1}-1}f
\\
Y_1f,\dots,Y_{r_j}f,\
&
Y_{r_j+1}f,\dots,Y_{r_{j-1}-2}f
\\
\cdots\cdots\cdots\cdots\cdots
&
\cdots\cdots\cdots\cdots\cdots\cdot
\\
Y_1f,\dots,Y_{r_j}f,\
&
Y_{r_j+1}f,\,Y_{r_j+2}f
\\
Y_1f,\dots,Y_{r_j}f,\
&
Y_{r_j+1}f
\endaligned
\]
generates an invariant subgroup of the $G_{ r_{ j-1}}$.

In this way, between the two groups $G_{ r_{j-1}}$ and $G_{ r_j}$,
there are certain groups\,---\,let us call them $\Gamma_{ r_{j-1}
-1}$, $\Gamma_{r_{ j-1} -2}$, \dots, $\Gamma_{ r_j+1}$\,---\,which
interpolate them and which possess the following properties: each one
of them has a number of terms exactly one less than that of the group
just preceding in the series, each one of them is contained in the
$G_{ r_{ j-1}}$ and in all the groups preceding in the series, and
each one of them is invariant in the $G_{ r_{ j-1}}$, hence also
invariant in all the groups preceding in the series, and in
particular, invariant in the group immediately preceding.  About that,
it yet comes that the $G_{ r_j}$ is contained in the $\Gamma_{ r_j+1}$
as an invariant subgroup.

What has been said holds for all values: $1, 2, \dots, q$ of the
number $j$, and consequently, we have effectively proved that $Y_1f,
\dots, Y_rf$ are independent infinitesimal transformations of the
group $X_1f, \dots, X_rf$ such that, for any $i < r$, $Y_1f, \dots,
Y_if$ always generate an $i$-term group which is invariant in the
$(i+1)$-term group $Y_1f, \dots, Y_{ i+1}f$. But this is what we
wanted to prove.

The case where the entire number $r_q$ defined above vanishes is now
settled, so it still remains the case where $r_q$ is larger than zero.
We shall show that it is impossible in this case
to select $r$ infinitesimal transformations $Y_1f, \dots, Y_rf$
in the group $X_1f, \dots, X_rf$ which stand in relationships
of the form~\thetag{ 11}. 

If $r_q >0$, then the $r_q$-term group $G_{ r_q}$ which is generated
by the $r_q$ independent infinitesimal transformations $X_1^{ (q)}f,
\dots, X_{ r_q}^{(q)}f$ certainly contains no $(r_q-1)$-term invariant
subgroup. Indeed, since $r_{ q+1} = r_q$, then amongst the
infinitesimal transformations $\leftbracket X_i^{ (q)},\,X_k^{(q)}
\rightbracket$, one finds exactly $r_q$ that are mutually independent,
whence by taking account of the Proposition~8,
p.~\pageref{Satz-8-S-263}, the property of the group $G_{ r_q}$ just
stated follows.

At present, we assume\footnote{\,
(reasoning by contradiction)
} 
that in the group $X_1f, \dots, X_rf$, there are $r$ independent
infinitesimal transformations $Y_1f, \dots, Y_rf$ which are linked by
relations of the form~\thetag{ 11}, and we denote by $\mathfrak{ G}_i$
the $i$-term group which, under this assumption, is generated by
$Y_1f, \dots, Y_if$.

According to p.~\pageref{S-267}, the group $G_{ r_q}$ is invariant in
the group $X_1f, \dots, X_rf$, but as we saw just now, it contains no
$(r_q - 1)$-term invariant subgroup. Now, since each $\mathfrak{ G}_i$
contains an $(i-1)$-term invariant subgroup, namely the group
$\mathfrak{ G}_{ i-1}$, it then follows immediately that the $G_{
r_q}$ cannot coincide with the group $\mathfrak{ G}_{ r_q}$, and at
the same time, it also comes that there exists an integer number $m$
which is at least equal to $r_q$ and is smaller than $r$, such that
the $G_{ r_q}$ is contained in none of the groups $\mathfrak{ G}_{
r_q}$, $\mathfrak{ G}_{ r_{q+1}}$, \dots, $\mathfrak{ G}_m$, while by
contrast, it is contained in all groups: $\mathfrak{ G}_{ m+1}$,
$\mathfrak{ G}_{ m+2}$, \dots, $\mathfrak{ G}_r$.

As a result, we have an $m$-term group $\mathfrak{ G}_m$
and an $r_q$-term group $G_{ r_q}$ which are both contained
in the $(m+1)$-term group $\mathfrak{ G}_{ m+1}$ as subgroups,
and to be precise, which are both evidently contained in it as 
invariant\footnote{\,
Since, as we saw, $G_{ r_q}$ is invariant in $G_r = \mathfrak{ G}_r$,
it is then trivially invariant in $\mathfrak{ G}_{ m+1} \subset
\mathfrak{ G}_r$.
} 
subgroups. According to Chap.~\ref{kapitel-12}, Proposition~7,
p.~\pageref{Satz-7-S-211}, the transformations common to the
$\mathfrak{ G}_m$ and to the $G_{ r_q}$ form a group $\Gamma$ which
has at least $r_q - 1$ parameters, and which, according to
Chap.~\ref{kapitel-15}, Proposition~10, p.~\pageref{Satz-10-S-264}, is
invariant in the $\mathfrak{ G}_{ m+1}$.  Now, since under the
assumptions made, the $G_{ r_q}$ is not contained in the $\mathfrak{
G}_m$, it comes that $\Gamma$ is exactly $(r_q-1)$-term and is at the
same time invariant in the $G_{ r_q}$.

This is a contradiction, since according to what precedes, the $G_{
r_q}$ contains absolutely no invariant $(r_q - 1)$-term subgroup.
Consequently, the assumption which we took as a starting point is
false, namely the assumption that in the group $X_1f, \dots, X_rf$,
one can indicate $r$ independent infinitesimal transformations $Y_1f,
\dots, Y_rf$ which stand mutually in relationships of the
form~\thetag{ 11}.  From this, we see that in the case $r_q >0$, there
are no infinitesimal transformations $Y_1f, \dots, Y_rf$ of this
constitution.

In the preceding developments, a simple process has been provided by
means of which one can realize whether a given $r$-term group $X_1f,
\dots, X_rf$ belongs, or does not belong, to the specific category
defined on page~\pageref{S-265}.

\renewcommand{\thefootnote}{\fnsymbol{footnote}}
\def\theproposition{12}\begin{proposition}
\footnote[1]{\,
In this proposition, $r_1$ naturally has the same meaning
as on p.~\pageref{S-266}, and likewise $r_2$, $r_3$, etc.
} 
If \label{Satz-12-S-270}
$X_1f, \dots, X_r f$ are independent infinitesimal transformations
of an $r$-term group, if $X_1'f, \dots, X_{ r_1}'f$ ($r_1 \leqslant
r$) are independent infinitesimal transformations from which all
$\leftbracket X_i,\,X_k \rightbracket$ can be linearly deduced, if
furthermore $X_1''f, \dots, X_{ r_2}''f$ ($r_2 \leqslant r_1$) are
independent infinitesimal transformations from which all $\leftbracket
X_i',\,X_k' \rightbracket$ can be linearly deduced, and if one defines
in a corresponding way $r_3 \leqslant r_2$ mutually independent
infinitesimal transformations $X_1'''f, \dots, X_{r_3}'''f$, and so
on, then the $X'f$ generate an $r_1$-term group, the $X''f$ an
$r_2$-term group, the $X'''f$ an $r_3$-term group, and so on, and to
be precise, each one of these groups is invariant in all the preceding
groups, and also in the group $X_1f, \dots, X_rf$.\,--- In the series
of the numbers $r$, $r_1$, $r_2$, \dots\, there is a number, say
$r_q$, which is equal to all the numbers following: $r_{ q+1}$, $r_{
q+2}$, \dots, while by contrast the numbers $r$, $r_1$, \dots, $r_q$
are all distinct one from another.  Now, if $r_q = 0$, then it is
always possible to indicate, in the group $X_1f, \dots, X_rf$, $r$
mutually independent infinitesimal transformations $Y_1f, \dots, Y_rf$
such that for every $i < r$, the transformations $Y_1f, \dots, Y_if$
generate an $i$-term group which is invariant in the $(i+1)$-term
group: $Y_1f, \dots, Y_{ i+1}f$, so that there exist relations of the
specific form:
\[
\aligned
\leftbracket
Y_i,\,Y_{i+k}
\rightbracket
&
=
c_{i,i+k,1}\,Y_1f
+\cdots+
c_{i,i+k,i+k-1}\,Y_{i+k-1}f
\\
&
\ \ \ \ \ \ \
{\scriptstyle{(i\,=\,1\,\cdots\,r\,-\,1\,;\,\,\,
k\,=\,1\,\cdots\,r\,-\,i)}}.
\endaligned
\]
However, if $r_q > 0$, it is not possible to determine
$r$ independent infinitesimal transformations $Y_1f, \dots, Y_rf$
in the group $X_1f, \dots, X_rf$ having the constitution defined 
above.
\end{proposition}
\renewcommand{\thefootnote}{\arabic{footnote}}

\linestop


\chapter{The Adjoint Group}
\label{kapitel-16}
\chaptermark{The Adjoint Group}

\setcounter{footnote}{0}

\abstract*{??}

Let $x_i' = f_i ( x_1, \dots, x_n, \, a_1, \dots, a_r)$ be an $r$-term
group with the $r$ infinitesimal transformations:
\[
X_kf
=
\sum_{i=1}^n\,\xi_{ki}(x)\,\frac{\partial f}{\partial x_i}
\ \ \ \ \ \ \ \ \ \ \ \ \
{\scriptstyle{(k\,=\,1\,\cdots\,r)}}.
\]
If one introduces the $x_i'$ as new variables in the expression
$\sum\, e_k \, X_kf$, then as it has been already shown in
Chap.~\ref{fundamental-differential}, 
Proposition~4, p.~\pageref{Satz-4-S-81}, 
one gets for all values of the $e_k$ an equation of the form:
\[
\sum_{k=1}^r\,e_k\, X_kf
=
\sum_{k=1}^r\,e_k'\, X_k'f.
\]
Here, the $e_k'$ are certain linear, homogeneous functions of the
$e_k$ with coefficients that depend upon $a_1, \dots, a_r$:
\def\theequation{1}\begin{equation}
e_k'
=
\sum_{j=1}^r\,\rho_{kj}(a_1,\dots,a_r)\, e_j.
\end{equation}

If one again introduces in $\sum\, e_k' \, X_k' f$ the new variables
$x_i'' = f_i ( x, \, b)$, then one receives:
\[
\sum_{k=1}^r\,e_k'\, X_k'f
=
\sum_{k=1}^r\,e_k''\, X_k''f,
\]
where:
\def\theequation{1'}\begin{equation}
e_k''
=
\sum_{j=1}^r\,\rho_{kj}(b_1,\dots,b_r)\, e_j'.
\end{equation}
But now, because the equations $x_i' = f_i ( x, \, a)$ represent a
group, the $x''$ are consequently linked with the $x$ through
relations of the form $x_i'' = f_i ( x, \, c)$ in which the $c$ depend
only upon $a$ and $b$:
\[
c_k
=
\varphi_k(a_1,\dots,a_r,\,b_1,\dots,b_r).
\]
Hence if one passes directly from the $x$ to the $x''$,
one finds:
\[
\sum_{k=1}^r\,e_k\, X_kf
=
\sum_{k=1}^r\,e_k''\, X_k''f,
\]
and to be precise, one has:
\def\theequation{1''}\begin{equation}
e_k''
=
\sum_{j=1}^r\,\rho_{kj}(c_1,\dots,c_r)\, e_j
=
\sum_{j=1}^r\,\rho_{kj}\big(\varphi_1(a,b),\dots,\varphi_r(a,b)\big)
\, e_j.
\end{equation}
From this, it can be deduced that the totality of all transformations
$e_k' = \sum \, \rho_{ kj} ( a) \, e_j$ forms a group. Indeed, by
combination of the equations~\thetag{ 1} and~\thetag{ 1'} it comes:
\[
e_k''
=
\sum_{j,\,\,\nu}^{1\cdots r}\,
\rho_{kj}(b_1,\dots,b_r)\,
\rho_{j\nu}(a_1,\dots,a_r)\, e_\nu,
\]
what must naturally coincide with the equations~\thetag{ 1''} and in
fact, for all values of the $e$, of the $a$ and of the $b$. 
Consequently, there are the $r^2$ identities:
\[
\rho_{k\nu}\big(
\varphi_1(a,b),\dots,\varphi_r(a,b)\big)
\equiv
\sum_{j=1}^r\,\rho_{j\nu}(a_1,\dots,a_r)\,
\rho_{kj}(b_1,\dots,b_r),
\]
from which it results that the family of the transformations $e_k' =
\sum\, \rho_{ kj} ( a) \, e_j$ effectively forms a group.

\renewcommand{\thefootnote}{\fnsymbol{footnote}}
To every $r$-term group $x_i' = f_i ( x, \, a)$ therefore belongs
a fully determined linear homogeneous group:
\def\theequation{1}\begin{equation}
e_k'
=
\sum_{j=1}^r\,\rho_{kj}(a_1,\dots,a_r)\, e_j
\ \ \ \ \ \ \ \ \ \ \ \ \
{\scriptstyle{(k\,=\,1\,\cdots\,r)}},
\end{equation}
which we want to call the 
\terminology{adjoint group}\footnote[1]{\,
\name{Lie}, Archiv for Math., Vol.~1, Christiania 1876.
} 
\deutsch{adjungirte Gruppe} 
of the group $x_i' = f_i ( x, \, a)$.
\renewcommand{\thefootnote}{\arabic{footnote}}

We consider \emphasis{for example} the two-term group $x' = ax + b$ with
the two independent infinitesimal transformations: 
$\frac{ \D\, f}{ \D\,x}$,
$x\, \frac{ \D\,f}{ \D\,x}$. We find:
\[
e_1\,\frac{\D\,f}{\D\,x}
+
e_2\,x\,\frac{\D\,f}{\D\,x}
=
e_1\,a\,\frac{\D\,f}{\D\,x'}
+
e_2\,(x'-b)\,\frac{\D\,f}{\D\,x'}
=
e_1'\,\frac{\D\,f}{\D\,x'}
+
e_2'\,x'\,\frac{\D\,f}{\D\,x'},
\]
whence we obtain for the adjoint group of the group $x' = ax + b$ the
following equations:
\[
e_1'
=
a\,e_1-b\,e_2,
\ \ \ \ \ \ \ \ 
e_2'
=
e_2,
\]
which visibly really represent a group.

The adjoint group of the group $x_i' = f_i ( x, \, a)$ contains, under
the form in which it has been found above, precisely $r$ arbitrary
parameters: $a_1, \dots, a_r$. But for every individual group $x_i' =
f_i ( x, \, a)$, a special research is required to investigate whether
the parameters $a_1, \dots, a_r$ are all essential in the adjoint
group. Actually, we shall shortly see that there are $r$-term groups
whose adjoint group does not contain $r$ essential parameters.

Besides, in all circumstances, one special
transformation appears in the adjoint
group of the group $x_i' = f_i ( x, \, a)$, namely the identity
transformation; for if one sets for $a_1, \dots, a_r$ in the
equations~\thetag{ 1} the system of values which produces the identity
transformation $x_i' = x_i$ in the group $x_i' = f_i ( x, \, a)$, then
one obtains the transformation: $e_1' = e_1, \dots, e_r' = e_r$, which
hence is always present in the adjoint group. However, as we shall
see, it can happen that the adjoint group consists only of the
identity transformation: $e_1' = e_1, \dots, e_r ' = e_r$.

\sectionengellie{\S\,\,\,76.}

In order to make accessible the study of the adjoint group, we must
above all determine its infinitesimal transformations. We easily
reach this end by an application of the Theorem~43, 
Chap.~\ref{kapitel-15}, p.~\pageref{Theorem-43-S-252}; yet
we must in the process replace the equations $x_i' = f_i ( x, \, a)$
of our group by the equivalent \emphasis{canonical}\, equations:
\def\theequation{2}\begin{equation}
x_i'
=
x_i
+
\frac{t}{1}\,\sum_{k=1}^r\,
\lambda_k\, X_k\,x_i
+\cdots
\ \ \ \ \ \ \ \ \ \ \ \ \
{\scriptstyle{(i\,=\,1\,\cdots\,n)}},
\end{equation}
which represent the $\infty^{ r - 1}$ one-term subgroups of the group
$x_i' = f_i( x, \, a)$. According to Chap.~\ref{one-term-groups},
p.~\pageref{S-69}, the $a_k$ are defined here as functions of $t$ and
$\lambda_1, \dots, \lambda_r$ by the simultaneous system:
\def\theequation{3}\begin{equation}
\frac{\D\,a_k}{\D\,t}
=
\sum_{j=1}^r\,\lambda_j\,
\alpha_{jk}(a_1,\dots,a_r)
\ \ \ \ \ \ \ \ \ \ \ \ \
{\scriptstyle{(k\,=\,1\,\cdots\,r)}}.
\end{equation}

By means of the equations~\thetag{ 2}, we therefore have to
introduce the new variables $x_i'$ in $\sum \, e_k\, X_kf$ and
we must then obtain a relation of the
form:
\[
\sum_{k=1}^r\,e_k\, X_kf
=
\sum_{k=1}^r\,e_k'\, X_k'f.
\]
The infinitesimal transformation denoted by $Yf$ in 
Theorem~43 on p.~\pageref{Theorem-43-S-252}
now writes: $\lambda_1 \, X_1 f + \cdots + \lambda_r \, X_r f$; we
therefore receive in our case:
\[
\aligned
Y\big(X_k(f)\big)
-
X_k\big(Y(f)\big)
&
=
\sum_{\nu=1}^r\,\lambda_\nu\,
\leftbracket
X_\nu,\,X_k\rightbracket
\\
&
=
\sum_{s=1}^r\,
\Big\{
\sum_{\nu=1}^r\,\lambda_\nu\,c_{\nu ks}
\Big\}\,X_sf.
\endaligned
\]
Consequently, we obtain the following differential equations for
$e_1', \dots, e_r'$:
\def\theequation{4}\begin{equation}
\label{S-273}
\frac{\D\,e_s'}{\D\,t}
+
\sum_{\nu=1}^r\,\lambda_\nu\,
\sum_{k=1}^r\,c_{\nu ks}\,e_k'
=
0
\ \ \ \ \ \ \ \ \ \ \ \ \
{\scriptstyle{(s\,=\,1\,\cdots\,r)}}.
\end{equation}

We consider the integration of these differential equations as an
executable operation, \label{S-273-bis}
for it is known that it requires only the
resolution of an algebraic equation of $r$-th degree. So if we
perform the integration on the basis of the initial condition: $e_k' =
e_k$ for $t = 0$, we obtain $r$ equations of the form:
\def\theequation{5}\begin{equation}
e_k'
=
\sum_{j=1}^r\,\psi_{kj}(\lambda_1t,\dots,\lambda_rt)\, e_j
\ \ \ \ \ \ \ \ \ \ \ \ \
{\scriptstyle{(k\,=\,1\,\cdots\,r)}},
\end{equation}
which are equivalent to the equations~\thetag{ 1}, as soon as the
$a_k$ are expressed as functions of $\lambda_1 t, \dots, \lambda_r t$
in the latter.

From this, 
it follows that the equations~\thetag{ 5} represent the
adjoint group too. But now we have derived the equations~\thetag{ 5}
in exactly the same way as if we would have wanted to determine all
finite transformations which are generated by the infinitesimal
transformations:
\[
\sum_{\nu=1}^r\,\lambda_\nu\,
\sum_{k,\,\,s}^{1\cdots r}\,c_{k\nu s}\,e_k\,
\frac{\partial f}{\partial e_s}
=
\sum_{\nu=1}^r\,\lambda_\nu\, E_\nu f
\]
(cf. p.~51 above). \emphasis{Consequently, we conclude that the
adjoint group~\thetag{ 1} consists of the totality of all one-term
groups of the form $\lambda_1 \, E_1 f + \cdots + \lambda_r \, E_r
f$.}

If amongst the family of all infinitesimal transformations $\lambda_1
\, E_1 f + \cdots + \lambda_r \, E_r f$ there are exactly $\rho$
transformations and no more which are independent, say $E_1 f, \dots,
E_\rho f$, then all the finite transformations of the one-term groups
$\lambda_1 \, E_1 f + \cdots + \lambda_r \, E_r f$ are already
contained in the totality of all finite transformations of the
$\infty^{ \rho - 1}$ groups $\lambda_1 \, E_1 f + \cdots +
\lambda_\rho \, E_\rho f$. The totality of these $\infty^\rho$ finite
transformations forms the adjoint group: $e_k' = \sum \, \rho_{ kj} (
a) \, e_j$, which therefore contains only $\rho$ essential
parameters (Chap.~\ref{one-term-groups}, Theorem~8, 
p.~\pageref{Theorem-8-S-65}).

According to what precedes, it is to be supposed that $E_1 f, \dots,
E_\rho f$ are linked together by relations of the form:
\[
\leftbracket
E_\mu,\,E_\nu
\rightbracket
=
\sum_{s=1}^\rho\,
g_{\mu\nu s}\, E_sf;
\]
we can also confirm this by a computation. By a direct calculation,
it comes:
\[
E_\mu\big(E_\nu(f)\big)
-
E_\nu\big(E_\mu(f)\big)
=
\sum_{\sigma,\,\,k,\,\,\pi}^{1\cdots r}\,
\big(
c_{\pi\mu k}\,c_{k\nu\sigma}
-
c_{\pi\nu k}\,c_{k\mu\sigma}
\big)\,
e_\pi\,
\frac{\partial f}{\partial e_\sigma}.
\]
But between the $c_{ ik s}$, there exist the relations:
\[
\sum_{k=1}^r\,
\big(
c_{\pi\mu k}\,c_{k\nu\sigma}
+
c_{\mu\nu k}\,c_{k\pi\sigma}
+
c_{\nu\pi k}\,c_{k\mu\sigma}
\big)
=
0,
\]
which we have deduced from the Jacobi identity some time ago
(cf. Chap.~\ref{kapitel-9}, Theorem~27,
p.~\pageref{Theorem-27-S-170}). If we yet use for this that $c_{ \nu
\pi k} = - c_{ \pi \nu k}$ and $c_{ k \pi \sigma} = - c_{ \pi k
\sigma}$, we can bring the right hand-side of our equation for
$\leftbracket E_\mu, \, E_\nu \rightbracket$ to the form:
\[
\sum_{k=1}^r\,c_{\mu\nu k}\,
\sum_{\sigma,\,\,\pi}^{1\cdots r}\,
c_{\pi k\sigma}\,e_\pi\,
\frac{\partial f}{\partial e_\sigma},
\]
whence it comes:
\[
\label{S-274}
\leftbracket
E_\mu,\,E_\nu
\rightbracket
=
\sum_{k=1}^r\,c_{\mu\nu k}\, E_kf.
\]
Lastly, under the assumptions made above, the right hand side can be
expressed by means of $E_1 f, \dots, E_\rho f$ alone, so that
relations of the form:
\[
\leftbracket
E_\mu,\,E_\nu
\rightbracket
=
\sum_{s=1}^\rho\,
g_{\mu\nu s}\, E_s f
\]
really hold, in which the $g_{ \mu \nu s}$ denote constants.

Before we continue, we want yet to recapitulate in cohesion
\deutsch{im Zusammenhange wiederholen} the results of the chapter
obtained until now.

\def\thetheorem{48}\begin{theorem}
\label{Theorem-48-S-275}
If, in the general infinitesimal transformation $e_1 \, X_1 f + \cdots
+ e_r \, X_r f$ of the $r$-term group $x_i' = f_i ( x, \, a)$, one
introduces the new variables $x'$ in place of the $x$, then one
obtains an expression of the form:
\[
e_1'\, X_1'f
+\cdots+
e_r'\, X_r'f;
\]
in the process, the $e'$ are linked to the $e$ by
equations of the shape:
\[
e_k'
=
\sum_{j=1}^r\,\rho_{kj}(a_1,\dots,a_r)\, e_j
\ \ \ \ \ \ \ \ \ \ \ \ \
{\scriptstyle{(k\,=\,1\,\cdots\,r)}},
\]
which represent a group in the variables $e$, the so-called adjoint
group of the group $x_i' = f_i ( x, \, a)$. This adjoint group
contains the identity transformation and is generated by certain
infinitesimal transformations; if, between $X_1 f, \dots, X_r f$,
there exist the Relations:
\[
\leftbracket X_i,\,X_k\rightbracket
=
\sum_{s=1}^r\,c_{iks}\, X_sf
\ \ \ \ \ \ \ \ \ \ \ \ \
{\scriptstyle{(i,\,k\,=\,1\,\cdots\,r)}},
\]
and if one sets:
\[
\label{S-275}
E_\mu f
=
\sum_{k,\,\,j}^{1\cdots r}\,
c_{j\mu k}\,e_j\,\frac{\partial f}{\partial e_k}
\ \ \ \ \ \ \ \ \ \ \ \ \
{\scriptstyle{(\mu\,=\,1\,\cdots\,r)}},
\]
then $\lambda_1 \, E_1 f + \cdots + \lambda_r \, E_r f$ is the general
infinitesimal transformation of the adjoint group and between $E_1 f,
\dots, E_r f$, there are at the same time the relations:
\[
\leftbracket
E_i,\,E_k
\rightbracket
=
\sum_{s=1}^r\,c_{iks}\, E_sf
\ \ \ \ \ \ \ \ \ \ \ \ \
{\scriptstyle{(i,\,k\,=\,1\,\cdots\,r)}}.
\]
\end{theorem}

If two $r$-term groups $X_1 f, \dots, X_r f$ and $Y_1 f, \dots, Y_r f$
are constituted in such a way that one simultaneously has:
\[
\leftbracket X_i,\,X_k\rightbracket
=
\sum_{s=1}^r\,c_{iks}\, X_s,
\ \ \ \ \ \ \ \ \ \ \ \ \ \ \ \ \ \ \
\leftbracket Y_i,\,Y_k\rightbracket
=
\sum_{s=1}^r\,c_{iks}\, Y_sf,
\]
with the same constants $c_{ iks}$ in the two cases, then both groups
obviously have the same adjoint group. Later, we will see that in
certain circumstances, also certain groups which do not possess an
equal number of terms can nonetheless have the same adjoint group.

\sectionengellie{\S\,\,\,77.}

Now, by what can one recognize how many independent infinitesimal
transformation there are amongst $E_1 f, \dots, E_r f$?

If $E_1 f, \dots, E_r f$ are not all independent of each other, 
then there is at least one infinitesimal transformation
$\sum \, g_\mu \, E_\mu f$ which vanishes identically.
From the identity:
\[
\sum_{\mu=1}^r\,g_\mu\,
\sum_{k,\,\,j}^{1\cdots r}\,
c_{j\mu k}\,e_j\,
\frac{\partial f}{\partial e_k}
\equiv
0,
\]
it comes:
\[
\sum_{\mu=1}^r\,g_\mu\,c_{j\mu k}
=
0
\]
for all values of $j$ and $k$, and consequently the expression:
\[
\bigg\leftbracket
X_j,\,\sum_{\mu=1}^r\,g_\mu\, X_\mu f
\bigg\rightbracket
=
\sum_{k=1}^r\,\Big\{
\sum_{\mu=1}^r\,g_\mu\,c_{j\mu k}
\Big\}\,X_kf
\]
vanishes, that is to say: the infinitesimal transformation $\sum\,
g_\mu \, X_\mu f$ is interchangeable with all the $r$ infinitesimal
transformations $X_jf$. Conversely, if the group $X_1f, \dots, X_r f$
contains an infinitesimal transformation $\sum\, g_\mu \, X_\mu f$
which is interchangeable with all the $X_k f$, then if follows in the
same way that the infinitesimal transformation $\sum \, g_\mu \, E_\mu
f$ vanishes identically.

In order to express this relationship in a manner
which is as brief as possible, we introduce the following 
terminology:

\plainstatement{
An infinitesimal transformation $\sum \, g_\mu \, X_\mu f$ of the
$r$-term group $X_1f, \dots, X_r f$ is called an 
\label{S-276}
\terminology{excellent}
infinitesimal transformation of this group if it is interchangeable with
all the $X_k f$. }

Incidentally, the excellent infinitesimal transformations of the group
$X_1f, \dots, X_r f$ are also characterized by the fact that they keep
their form through the introduction of the new variables $x_i' = f_i (
x, \, a)$, whichever values the parameters $a_1, \dots, a_r$ can
have. Indeed, if the infinitesimal transformation $\sum\, g_\mu \,
X_\mu \, f$ is excellent, then according to Chap.~15, p.~259, there is
a relation of the form:
\[
{\textstyle{
\sum\,g_\mu\, X_\mu f
=
\sum\,g_\mu\, X_\mu'f}}.
\]
In addition, the cited developments show that each finite
transformation of the one-term group $\sum \, g_\mu \, X_\mu f$ is
interchangeable with every finite transformation of the group $X_1 f,
\dots, X_r f$.

According to what has been said above, to every excellent
infinitesimal transformation of the group $X_1f, \dots, X_rf$, it
corresponds a linear relation between $E_1f, \dots, E_rf$. If $\sum\,
g_\mu\, X_\mu f$ is an excellent infinitesimal transformation, then
there exists between the $Ef$ simply the relation: $\sum\, g_\mu\,
E_\mu f = 0$. Consequently, between $E_1f, \dots, E_rf$ there are
exactly as many independent relations of this sort as there are
independent excellent infinitesimal transformations in the group
$X_1f, \dots, X_rf$. If there are exactly $m$ and not more such
independent transformations, then amongst the infinitesimal
transformations $E_1f, \dots, E_rf$ there are exactly $r-m$ and not
more which are independent, and likewise, the adjoint group contains
the same number of essential parameters.

We therefore have the

\renewcommand{\thefootnote}{\fnsymbol{footnote}}
\def\thetheorem{49}\begin{theorem}
\label{Theorem-49-S-277}
The adjoint group $e_k' = \sum\, \rho_{ kj}(a)\, e_j$ of an $r$-term
group $X_1f, \dots, X_rf$ contains $r$ essential parameters if and
only if none of the $\infty^{ r-1}$ infinitesimal transformations
$\sum\, g_\mu\, X_\mu f$ is excellent inside the group $X_1f, \dots,
X_rf$; by contrast, the adjoint group has less than $r$ essential
parameters, namely $r - m$, when the group $X_1f, \dots, X_rf$
contains exactly $m$ and not more independent excellent infinitesimal
transformations.\footnote{\,
\name{Lie}, Math. Ann. Vol. XXV, p.~94.
}
\end{theorem}
\renewcommand{\thefootnote}{\arabic{footnote}}

Let us take for example the group:
\[
\frac{\partial f}{\partial x_1},
\,\,\,\dots,\,\,\,
\frac{\partial f}{\partial x_r}
\ \ \ \ \ \ \ \ \ \ \ \ \ {\scriptstyle{(r\,\,\leqslant\,n)}}.
\]
Its infinitesimal transformations are all excellent, since all $c_{
iks}$ are zero. So all $r$ expressions $E_1f, \dots, E_rf$ vanish
identically, and the adjoint group reduces to the identity
transformation.

On the other hand, let us take the four-term group:
\[
x\,\frac{\partial f}{\partial x},\ \ \ \ \
y\,\frac{\partial f}{\partial x},\ \ \ \ \
x\,\frac{\partial f}{\partial y},\ \ \ \ \
y\,\frac{\partial f}{\partial y}.
\]
This group contains a single excellent infinitesimal transformation,
namely:
\[
x\,\frac{\partial f}{\partial x}
+
y\,\frac{\partial f}{\partial y}\,;
\]
its adjoint group is hence only three-term.

\medskip

If the group $X_1f, \dots, X_rf$ contains the excellent infinitesimal
transformation $\sum\, g_\mu\, X_\mu f$, then $\sum\, g_\mu\, X_\mu f$
is a linear partial differential equation which remains invariant by
the group; in consequence of that, the group is imprimitive
(cf. p.~\pageref{S-221}). From this, we conclude that the following
proposition holds:

\def\theproposition{1}\begin{proposition}
\label{Satz-1-S-277}
If the group $X_1f, \dots, X_rf$ in the variables $x_1, \dots, x_n$
is primitive, then it contains no excellent infinitesimal
transformation.
\end{proposition}

If we combine this proposition with the latter theorem, we yet obtain 
the

\def\theproposition{2}\begin{proposition}
The adjoint group of an $r$-term primitive group contains
$r$ essential parameters.
\end{proposition}

In a later place, we shall give a somewhat more general version of the
last two propositions (cf. the Chapter~\ref{kapitel-24} about systatic
groups).

According to what precedes, the infinitesimal transformations of an
$r$-term group $X_1f, \dots, X_rf$ and those of the adjoint group
$E_1f, \dots, E_rf$ can be mutually ordered in such a way that, to
every infinitesimal transformation $\sum\, e_k\, X_k$ which is not
excellent there corresponds a nonvanishing infinitesimal
transformation $\sum\, e_k\, E_kf$, while to every excellent
infinitesimal transformation $\sum\, e_k\, X_kf$ is associated the
\emphasis{identity} transformation in the adjoint group. In addition,
if one takes account of the fact that the two systems of equations:
\[
\leftbracket
X_i,\,X_k
\rightbracket
=
\sum_{s=1}^r\,c_{iks}\,X_sf,
\ \ \ \ \ \ \ \ \ \ \
\leftbracket
E_i,\,E_k
\rightbracket
=
\sum_{s=1}^r\,c_{iks}\,E_sf
\]
have exactly the same form, then one could be led to the presumption
that the adjoint group cannot contain excellent infinitesimal
transformations. Nevertheless, this presumption would be false; that
is shown by the group $\partial f / \partial x_2$, $x_1\, \partial f /
\partial x_2$, $\partial f / \partial x_1$ whose adjoint group
consists of two interchangeable infinitesimal transformations and
therefore contains two independent excellent infinitesimal
transformations.

\sectionengellie{\S\,\,\,78.}

The starting point of our study was the remark that the expression
$\sum\, e_k\, X_kf$ takes the similar form $\sum\, e_k'\, X_k'f$ after
the introduction of the new variables $x_i' = f_i ( x, a)$. But now,
according to Chap.~\ref{kapitel-15}, p.~\pageref{S-255}, the
expression $\sum\, e_k\, X_kf$ can be interpreted as the symbol of the
general finite transformation of the group $X_1f, \dots, X_rf$; here,
the quantities $e_1, \dots, e_r$ are to be considered as the
parameters of the finite transformations of the group $X_1f, \dots,
X_rf$. Consequently, we can also say: after the transition to the
variables $x'$, the finite transformations of the group $X_1f, \dots,
X_rf$ are permuted with each other, while their totality remains
invariant. Thanks to the developments of the preceding paragraph, we
can add that the concerned permutation \deutsch{Vertauschung} is
achieved by a transformation of the adjoint group.

But when it is spoken of a ``\emphasis{permutation 
of the finite transformations $\sum\, e_k\, X_kf$}'', 
the interpretation fundamentally lies in the fact that one
imagines these transformations as \emphasis{individuals}
\deutsch{Individuen}; 
we now want to pursue this interpretation 
somewhat in details. 

Every individual in the family $\sum\, e_k\, X_kf$
is determined by the associated values of $e_1, \dots, e_r$; 
for the sake of graphic clarity, we hence imagine 
the $e_k$ as right-angled point coordinates of an 
$r$-times extended manifold.
The points of this manifold then represent all finite transformations
$\sum\, e_k\, X_kf$; hence, they are transformed
by our known linear homogeneous group:
\[
e_k'
=
\sum_{j=1}^r\,\rho_{kj}(a)\,e_j.
\]

At the same time, the origin of coordinates $e_1 = 0$, \dots, $e_r =
0$, i.e. the image-point \deutsch{Bildpunkt} of the identity
transformation $x_i' = x_i$, obviously remains invariant.

Every finite transformation $e_k^0$ belongs to a completely determined
one-term group, whose transformations are defined by the equations:
\[
\frac{e_1}{e_1^0}
=\cdots=
\frac{e_r}{e_r^0}\,;
\]
but according to our interpretation, these equations represent a 
straight line
through $e_k = 0$, that is to say:

\plainstatement{Every one-term subgroup of the group $X_1f, \dots,
X_rf$ is represented, in the space $e_1, \dots, e_r$, by a 
straight line which
passes through the origin of coordinates $e_k = 0$; conversely, every
straight
line through the origin of coordinates represents such a one-term
group.}

Every other subgroup of the $r$-term group $X_1f, \dots, X_rf$
consists of one-term subgroups and these one-term groups are
determined by the infinitesimal transformations that the subgroup in
question contains. But according to Chap.~\ref{kapitel-11},
Proposition~6, p.~\pageref{Satz-6-S-211}, the infinitesimal
transformations of a subgroup of the group $X_1f, \dots, X_rf$ can be
defined by means of \emphasis{linear homogeneous} equations between
$e_1, \dots, e_r$; these equations naturally define also the one-term
groups which belong to the subgroup, hence they actually define the
finite transformations of the subgroup in question. Expressed
differently:

\plainstatement{
Every $m$-term subgroup of the group $X_1f, \dots, X_rf$ is
represented, in the space $e_1, \dots, e_r$, by a 
straight\footnote{\,
The ``straight'' character just 
means that the manifold in question is a linear
subspace of the space $e_1, \dots, e_r$.
} 
\deutsch{eben} \label{straight-eben}
$m$-times
extended manifold which passes through the origin of coordinates: $e_1
= 0$, \dots, $e_r = 0$.}

Of course, the converse does not hold true in general; it occurs only
very exceptionally that \emphasis{every} straight manifold through the
origin of coordinates represents a subgroup.

The adjoint group $e_k' = \sum\, \rho_{ kj} (a)\, e_j$ now transforms
linearly the points $e_k$. If an arbitrary point $\overline{ e}_k$
remains invariant by all transformations of the group, and in
consequence of that, also every point $m \, \overline{ e}_k = {\rm
Const.}\, \overline{ e}_k$, then according to what precedes, this
means that the transformations of the one-term group $\sum\,
\overline{ e}_k\, X_kf$ are interchangeable with all transformations
of the group $X_1f, \dots, X_rf$. \label{S-280}

Every \emphasis{invariant} \label{S-280-ter}
subgroup of the group $X_1f, \dots, X_rf$
is represented by a straight manifold containing the origin of
coordinates which keeps its position through all transformations $e_k'
= \sum\, \rho_{ kj} (a)\, e_j$. On the other hand, every manifold
through the origin of coordinates which remains invariant by the
adjoint group represents an invariant subgroup of $X_1f, \dots, X_rf$
(cf. Chap.~\ref{kapitel-15}, Proposition~5,
p.~\pageref{Satz-5-S-261}).

Now, let $M$ be any straight manifold through the origin of
coordinates which represents a subgroup of the $G_r$: $X_1f, \dots,
X_rf$. If an arbitrary transformation of the adjoint group is executed
on $M$, this gives a new straight manifold which, likewise, represents
a subgroup of the $G_r$. \label{S-280-bis}
Obviously, this new subgroup is equivalent
\deutsch{ähnlich}, through a transformation of the $G_r$, to the
initial one, and if we follow the process of the theory of
substitutions, we can express this as follows: \emphasis{the two
subgroups just discussed are \terminology{conjugate}
\deutsch{gleichberechtigt}\footnote{\,
Literally in German: they are equal, they are considered on the same
basis, or they have the same rights.
} 
inside the group $X_1f, \dots, X_rf$}.

The totality of all subgroups which, inside the group $G_r$, are
conjugate to the subgroup which is represented by $M$ is represented
by a family of straight manifolds, namely by the family that one
receives as soon as one executes on $M$ all transformations of the
adjoint group. Two different manifolds of this family naturally
represent conjugate subgroups.

From this we conclude that every invariant subgroup of the $G_r$ is
conjugate only to itself inside the $G_r$.

Since the family of all subgroups $g$ which are conjugate to a given
one result from the latter by the execution of all transformations
$e_k' = \sum\, \rho_{ kj} (a)\, e_j$, then this family itself must be
reproduced by all transformations $e_k'' = \sum\, \rho_{kj}(b)\, e_j'$.
Because if one executes the transformations $e_k' = \sum\, 
\rho_{kj}(a)\, e_j$ and $e_k'' = \sum\, \rho_{kj}(a)\, e_j'$
one after the other, then one obtains the same result
as if one would have applied all transformations:
\[
e_k''
=
\sum\,\rho_{kj}(c)\,e_j
\ \ \ \ \ \ \ \ \ \ \ \ \ 
{\scriptstyle{(c_k\,=\,\varphi_k\,(a,\,\,b))}}
\]
to the initially given subgroup; in the two cases, one obtains the
said family.

Now, if all conjugate subgroups $g$ have a continuous number of
mutually common transformations, then the totality of all these
transformations are represented by a straight manifold in the space
$e_1, \dots, e_r$. Naturally, this straight manifold remains invariant
by all transformations of the adjoint group and hence, according to
what was said above, it represents an invariant subgroup of the
$G_r$. Consequently, the following holds true. \label{S-281}

\def\thetheorem{50}\begin{theorem}
If the subgroups of a group $G_r$ which are conjugate inside the $G_r$
to a determined subgroup have a family of transformations in common,
then the totality of these transformations forms a subgroup
invariant in the $G_r$.
\end{theorem}

\sectionengellie{\S\,\,\,79.}

We distribute \deutsch{eintheilen} the subgroups of an arbitrary
$r$-term group $G_r$ in different classes which we call
\terminology{types \deutsch{Typen} \label{S-281-bis}
of subgroups of the $G_r$}. We
reckon as belonging to the \emphasis{same} type the groups which are
mutually conjugate inside the $G_r$; the groups which are not
conjugate inside the $G_r$ are reckoned as belonging to
\emphasis{different} types.

If one knows any subgroup of the $G_r$, then at once, one can
determine all subgroups which belong to the same type. Thanks to this
fact, the enumeration \deutsch{Aufzählung} of the subgroups of a given
group is essentially facilitated, since one clearly does not have to
write down all subgroups, but rather, one needs only
to enumerate the different types of subgroups by indicating
a representative for every type, hence a subgroup which belongs
to the type in question. 

We also speak of different types for the finite transformations of a
group. We reckon as belonging to the same type two finite
transformations: $e_1^0 \, X_1f + \cdots + e_r^0 \, X_r f$ and
$\overline{ e}_1 \, X_1f + \cdots + \overline{ e}_r\, X_rf$ of the
$r$-term group $X_1f, \dots, X_rf$ if and only if they are conjugate
inside the $G_r$, that is to say: when the adjoint group contains a
transformation which transfers the point $e_1^0, \dots, e_r^0$ to the
point $\overline{ e}_1, \dots, \overline{ e}_r$. It is clear, should
it be mentioned, that in the concerned transformation $e_k' = \sum\,
\gamma_{ kj}\, e_j$ of the adjoint group the determinant of the
coefficients $\gamma_{ kj}$ should not vanish.

At present, we want to take up again the question about the types of
subgroups of a given group, though not in complete generality;
instead, we want at least to show how one has to proceed in order to
find the extant types of one-term groups and of finite
transformations.

At first, we ask for all types of finite transformations.

Let $e_1^0 \, X_1f + \cdots + e_r^0 \, X_rf$ be any finite
transformation of the group $X_1f, \dots, X_rf$. If, on the the point
$e_1^0, \dots, e_r^0$, we execute all transformations of the adjoint
group, we obtain the image-points \deutsch{Bildpunkte} of all the
finite transformations of our group which are conjugate to $\sum\,
e_k^0\, X_kf$ and thus, belong to the same type as $\sum\, e_k^0\,
X_kf$. According to Chap.~\ref{kapitel-14},
p.~\pageref{Theorem-36-S-225}, the totality of all these points forms
a manifold invariant by the adjoint group and to be precise, a
so-called \emphasis{smallest} invariant manifold, as was said at that
time.

We can consider the finite equations of the adjoint group as known;
consequently, we are in a position to indicate without integration the
equations of the just mentioned smallest invariant manifold
(cf. Chap.~\ref{kapitel-14}, Theorem~37,
p.~\pageref{Theorem-37-S-226}). Now, since every such
\emphasis{smallest} invariant manifold represents the totality of all
finite transformations which belong to a certain type, then with this,
all types of finite transformations of our group are found. As a
result, the following holds true.

\def\theproposition{3}\begin{proposition}
If an $r$-term group $x_i' = f_i ( x_1, \dots, x_n, \, a_1, \dots,
a_r)$ with the $r$ independent infinitesimal transformations $X_1f,
\dots, X_rf$ is presented, then one finds in the following way all
types of finite transformations $e_1\, X_1f + \cdots + e_r\, X_rf$ of
this group: one sets up the finite equations $e_k' = \sum\, \rho_{ kj}
(a)\, e_j$ of the adjoint group of the group $X_1f, \dots, X_rf$ and
one then determines, in the space of the $e$, the smallest manifolds
which remain invariant by the adjoint group; these manifolds represent
the demanded types.
\end{proposition}

With all of that, one does not forget that only the transformations of
the adjoint group for which the determinant of the coefficients does
not vanish are permitted. About this point, one may compare with what
was said in Chap.~\ref{kapitel-14}, p.~\pageref{S-225-bis}.

Next, we seek all types of one-term subgroups, or, what is the same:
all \emphasis{types of infinitesimal transformations} of our group.

Two one-term groups $\sum\, e_k^0\, X_kf$ and $\sum\, \overline{
e}_k\, X_kf$ are conjugate inside the $G_r$ when there is, in the
adjoint group, a transformation which transfers the straight line:
\[
\frac{e_1}{e_1^0}
=\cdots=
\frac{e_r}{e_r^0}
\] 
to the straight line:
\[
\frac{e_1}{\overline{e}_1}
=\cdots=
\frac{e_r}{\overline{e}_r}\,;
\]
these two straight lines passing through the origin of coordinates are
indeed the images of the one-term subgroups in question. Hence if we
imagine that all transformations of the adjoint group are executed on
the first one of these two straight lines, we obtain all straight
lines passing through the origin of coordinates that represent
one-term subgroups which are conjugate to the group $\sum\, e_k^0\,
X_kf$. Naturally, the totality of all these straight lines forms a
manifold invariant by the adjoint group, namely the smallest invariant
manifold to which the straight line: $e_1 / e_1^0 = \cdots = e_r /
e_r^0$ belongs. It is clear that this manifold is represented by a
system of equations which is \emphasis{homogeneous} in the variables
$e_1, \dots, e_r$.

Conversely, it stands to reason that every system of equations
homogeneous in the $e$ which admits the adjoint group represents
an invariant family of one-term groups. Now, since every 
system of equations homogeneous in the $e$ is characterized
to be homogeneous by the fact that it admits all
transformations of the form:
\def\theequation{6}\begin{equation}
e_1'=\lambda\,e_1,
\,\,\,\dots,\,\,\,
e_r'=\lambda\,e_r,
\end{equation}
it follows that we obtain all invariant families of one-term groups by
looking up at all manifolds of the space $e_1, \dots, e_r$ which,
aside from the transformations of the adjoint group, yet also admit
all transformations of the form~\thetag{ 6}. In particular, if we
seek all smallest invariant manifolds of this constitution, we clearly
obtain all existing types of one-term subgroups.

The transformations~\thetag{ 6} form a one-term group whose
infinitesimal transformations read:
\[
Ef
=
\sum_{k=1}^r\,e_k\,
\frac{\partial f}{\partial e_k}.
\]
If we add $Ef$ to the infinitesimal transformations $E_1f, \dots,
E_rf$ of the adjoint group, we again obtain the infinitesimal
transformations of a group; indeed, the expressions $\leftbracket
E_k,\,E \rightbracket$ all vanish identically, as their computation
shows. Visibly, everything amounts to the determination of the
smallest manifolds which remain invariant by the group $E_1f, \dots,
E_rf, Ef$. But according to the instructions in
Chap~\ref{kapitel-14}, p.~\pageref{S-225-ter}, this determination can
be accomplished, since together with the finite equations of the
adjoint group, the finite equations of the group just mentioned are
also known without effort.
Hence we have the

\def\theproposition{4}\begin{proposition}
If an $r$-term group $x_i' = f_i ( x_1, \dots, x_n, \, a_1, \dots,
a_r)$ with the $r$ independent infinitesimal transformations $X_1f,
\dots, X_rf$ is presented, then one finds as follows all types of
one-term subgroups of this group, or what is the same, all types of
infinitesimal transformations $E_1f, \dots, E_rf$: one sets up the
infinitesimal transformations of the adjoint group of the group $X_1f,
\dots, X_rf$, then one computes the finite equations of the group
which is generated by the $r+1$ infinitesimal transformations $E_1f,
\dots, E_rf$ and:
\[
Ef
=
\sum_{k=1}^r\,e_k\,
\frac{\partial f}{\partial e_k}
\]
and lastly, one determines, in the space $e_1, \dots, e_r$, the
smallest manifolds invariant by the just defined group; these
manifolds represent the demanded types.
\end{proposition}

In what precedes, it is shown how one can determine all types of
finite transformations and all types of one-term subgroups of a given
$r$-term group. Now in a couple of words, we will yet consider
somewhat more precisely the connection which exists between these two
problems; we will see that the settlement of one of these two problems
facilitates each time the settlement of the other.

On the first hand, 
we assume that we already know all types of \emphasis{finite
transformations} of the group $X_1f, \dots, X_rf$, so that all
smallest manifolds of the space $e_1, \dots, e_r$ which admit the
adjoint group $E_1f, \dots, E_rf$ are known to us. Then how one must
proceed in order to find the smallest manifolds invariant by the group
$Ef$, $E_1f, \dots, E_rf$?

It is clear that every manifold invariant by the group $Ef$, $E_1f,
\dots, E_rf$ admits also the adjoint group; consequently, every sought
manifold either must be one of the already known manifolds, or it must
contain at least one of the known manifolds. Hence in order to find
all the sought manifolds, we only need to take the known manifolds one
after the other and for each of them, to look up at the smallest
manifold invariant by the group $Ef$, $E_1f, \dots, E_rf$ in which it
is contained.

Let:
\def\theequation{7}\begin{equation}
W_1(e_1,\dots,e_r)=0,
\,\,\,\dots,\,\,\,
W_m(e_1,\dots,e_r)=0,
\end{equation}
or shortly $M$, be one of the known manifolds which admits the adjoint
group $E_1f, \dots, E_rf$. Now, how one finds the smallest manifold
which admits the group $Ef$, $E_1f, \dots, E_rf$ and which in addition
comprises the manifold $M$?

The sought manifold necessarily contains the origin of coordinates
$e_1 = 0$, \dots, $e_r = 0$ and moreover, it consists in nothing but
straight lines passing through it, hence it certainly contains the
manifold $M'$ which is formed of the straight lines between the points
of $M$ and the origin of coordinates. Now, if we can prove that $M'$
admits the infinitesimal transformations $Ef$, $E_1f, \dots, E_rf$,
then at the same time, we have proved that $M'$ is the sought
manifold.

Visibly, the equations of $M'$ are obtained by eliminating the
parameter $\tau$ from the equations:
\def\theequation{8}\begin{equation}
W_1(e_1\,\tau,\,\dots\,,e_r\,\tau)=0,
\,\,\,\dots,\,\,\,
W_m(e_1\,\tau,\,\dots,\,e_r\,\tau)=0.
\end{equation}
Consequently, $M'$ can be interpreted as the totality of all the
$\infty^1$ manifolds that are represented by the equations~\thetag{ 8}
with the arbitrary parameter $\tau$. But the totality of the
manifolds~\thetag{ 8} obviously admits the infinitesimal
transformation $Ef$, since the $\infty^1$ systems of
equations~\thetag{ 8} are permuted with each other by the finite
transformations:
\[
e_1'=\lambda\,e_1,
\,\,\,\dots,\,\,\,
e_r'=\lambda\,e_r
\]
of the one-term group $Ef$. Furthermore, it is easy to see that each
one of the individual systems of equations~\thetag{ 8} allows the
infinitesimal transformations $E_1f, \dots, E_rf$. Indeed, because
the system of equations~\thetag{ 7} allows the infinitesimal
transformations:
\[
E_\mu f
=
\sum_{k,\,\,j}^{1\cdots r}\,
c_{j\mu k}\,e_j\,
\frac{\partial f}{\partial e_k}
\ \ \ \ \ \ \ \ \ \ \ \ \ {\scriptstyle{(\mu\,=\,1\,\cdots\,r)}}
\]
then the system~\thetag{ 8} with the parameter $\tau$ admits the
transformations:
\[
\sum_{k,\,\,j}^{1\cdots r}\,
c_{j\mu k}\,e_j\,\tau\,
\frac{\partial f}{\partial (e_k\tau)}
=
E_\mu f,
\]
that is to say, it admits $E_1f, \dots, E_rf$ themselves, whichever
value $\tau$ can have. 

From this, we conclude that the totality of the $\infty^1$
manifolds~\thetag{ 8} admits the infinitesimal transformations $Ef$,
$E_1f, \dots, E_rf$, hence that the manifold $M'$ which coincides with
this totality really is the sought manifold;
as said, this manifold is analytically represented by the
equations which are obtained by elimination of the parameter $\tau$. 

If we imagine that for \emphasis{every} manifold $M$, the
accompanying manifold $M'$ is formed, then according to
the above, we obtain all types of one-term subgroups of
the group $X_1f, \dots, X_rf$.\,---

On the other hand, we assume that we know all types of
\emphasis{one-term subgroups} of the group $X_1f, \dots, X_rf$, 
and we then seek to determine from them all types of
finite transformations of these groups. 

All smallest manifolds invariant by the group $Ef$, 
$E_1f, \dots, E_rf$ are known to us and we must seek
all smallest manifolds invariant by the group $E_1f, \dots, E_rf$. 
But now, since every sought manifold is contained in
one of the known manifolds, we only have to consider for
itself each individual known manifold and to look up
at manifolds of the demanded constitution that are
located in each such manifold. 

Let the $q$-times extended manifold ${\sf M}$ be one of 
the smallest manifolds invariant by the group $Ef$, 
$E_1f, \dots, E_rf$. Then for the points of ${\sf M}$
(cf. Chap.~\ref{kapitel-14}, p.~\pageref{Theorem-42-S-237}), 
all the $(q+1) \times (q+1)$ determinants, but not all
$q \times q$ ones, of the matrix:
\def\theequation{9}\begin{equation}
\left\vert
\begin{array}{ccccc}
e_1 & \,\cdot\, & \,\cdot\, & \,\cdot\, & e_r
\\
\sum_{k=1}^r\,c_{k11}\,e_k 
& \,\cdot\, & \,\cdot\, & \,\cdot\, &
\sum_{k=1}^r\,c_{k1r}\,e_k
\\
\cdot
& \,\cdot\, & \,\cdot\, & \,\cdot\, &
\cdot
\\
\sum_{k=1}^r\,c_{kr1}\,e_k 
& \,\cdot\, & \,\cdot\, & \,\cdot\, &
\sum_{k=1}^r\,c_{krr}\,e_k
\end{array}
\right\vert
\end{equation}
vanish. 

Now, the question whether ${\sf M}$ decomposes in subsidiary domains
which remain invariant by the group $E_1f, \dots, E_rf$ is settled by
the behaviour of the $q \times q$ subdeterminants of the determinant:
\[
\Delta
=
\left\vert
\begin{array}{ccccc}
\sum_{k=1}^r\,c_{k11}\,e_k
& \,\cdot\, & \,\cdot\, & \,\cdot\, &
\sum_{k=1}^r\,c_{k1r}\,e_k
\\
\cdot
& \,\cdot\, & \,\cdot\, & \,\cdot\, &
\cdot
\\
\sum_{k=1}^r\,c_{kr1}\,e_k 
& \,\cdot\, & \,\cdot\, & \,\cdot\, &
\sum_{k=1}^r\,c_{krr}\,e_k
\end{array}
\right\vert.
\] 

If, for the points of ${\sf M}$, not all the subdeterminants in
question vanish, then no decomposition of ${\sf M}$ in smaller
manifolds invariant by the adjoint group takes place. 

By contrast, if the concerned subdeterminant all vanish for the points
of ${\sf M}$, then ${\sf M}$ decomposes in $\infty^1$ $(q-1)$-times
extended subsidiary domains invariant by the group $E_1f, \dots,
E_rf$; that there are exactly $\infty^1$ such subsidiary domains
follows from the fact that surely for the points of ${\sf M}$, not all
$(q-1) \times (q-1)$ subdeterminants of $\Delta$ vanish, for on the
contrary case, all the $q \times q$ determinants of the
matrix~\thetag{ 9} would vanish simultaneously. Naturally, one can
sets up without integration the equations of the discussed subsidiary
domains, since the finite equations of the group $E_1f, \dots, E_rf$
are known.

\linestop

In the present chapter, we always have up to now considered the
finite transformations of the group $X_1f, \dots, X_rf$
only as individuals and we have interpreted them as
points of an $r$-times extended space.

There is another standpoint which is equally legitimate. We can also
consider as individuals the one-term subgroups, or, what amounts to
the same, the infinitesimal transformations of our group, and
interpret them as points of a now $(r-1)$-times extended space. Then
clearly, we must understand the quantities $e_1, \dots, e_r$ as
\emphasis{homogeneous} coordinates in this space.

Our intention is not to follow further these views; above all, by
considering what was said earlier on, it appears for instance evident
that every $m$-term subgroup of the group $X_1f, \dots, X_rf$ is
represented, in the $(r-1)$-times extended space, by a
smooth\footnote{\,
Namely, the projectivization of a linear subspace (Lie subalgebra) or
``straight'' \deutsch{eben} manifold (cf. \pageref{straight-eben}). }
\deutsch{eben} manifold of $m-1$ dimensions. We only want to derive a
simple proposition which follows by taking as a basis the new
interpretation and which can often be useful for the conceptual
researches about transformation groups.

We imagine that the general finite transformation:
\[
x_i'
=
x_i
+
\frac{t}{1}\,
\sum_{k=1}^r\,\lambda_k^0\,X_k\,x_i
+\cdots
\ \ \ \ \ \ \ \ \ \ \ \ \ {\scriptstyle{(i\,=\,1\,\cdots\,n)}}
\]
of the one-term group $\lambda_1^0\, X_1f + \cdots + \lambda_r^0 \,
X_rf$ is executed on the infinitesimal transformation: $e_1^0 \, X_1f
+ \cdots + e_r^0 \, X_r f$ of our $r$-term group, that is to say, we
imagine that, in place of $x_1, \dots, x_n$, the new variables $x_1',
\dots, x_n'$ are introduced in the expression $e_1^0 \, X_1f + \cdots
+ e_r^0 \, X_rf$. According to p.~\pageref{S-273}, the infinitesimal
transformation: $e_1^0 \, X_1f + \cdots + e_r^0 \, X_rf$ is
transferred at the same time to the $\infty^1$ infinitesimal
transformations: $e_1 \, X_1f + \cdots + e_r\, X_rf$ of our group,
where $e_1, \dots, e_r$ are certain functions of $e_1^0, \dots, e_r^0$
and $t$ that determine themselves through the differential equations:
\def\theequation{4'}\begin{equation}
\frac{\D\,e_s}{\D\,t}
=
-\,\sum_{\nu=1}^r\,
\lambda_\nu^0\,
\sum_{k=1}^r\,c_{\nu ks}\,e_k
\ \ \ \ \ \ \ \ \ \ \ \ \ 
{\scriptstyle{(s\,=\,1\,\cdots\,r)}},
\end{equation}
with the initial conditions: $e_1 = e_1^0$, \dots, $e_r = e_r^0$ for
$t = 0$.

Since every infinitesimal transformation of our group is represented
by a point of the $(r-1)$-times extended space mentioned a short while
ago, we can also obviously say: If all $\infty^1$ transformations of
the one-term group: $\lambda_1^0 \, X_1f + \cdots + \lambda_r^0 \, X_r
f$ are executed on the infinitesimal transformation: 
$e_1^0 \, X_1f + \cdots + e_r^0 \, X_rf$, 
then the image-point of this infinitesimal 
transformation moves on a certain curve of
the space $e_1, \dots, e_r$. 

Now, there is a very simple definition for the tangent
to this curve at the point $e_1^0 \colon \cdots \colon e_r^0$. 
Namely, 
the equations of the tangent in question have the form:
\[
\aligned
e_s\,
\sum_{\nu,\,\,k}^{1\cdots r}\,
c_{\nu k\tau}\,
&
\lambda_\nu^0\,e_k^0
-
e_\tau\,
\sum_{\nu,\,\,k}^{1\cdots r}\,
c_{\nu ks}\,\lambda_\nu^0\,e_k^0
=
0
\\
&
\ \ \ \ \
{\scriptstyle{(s,\,\,\tau\,=\,1\,\cdots\,r)}},
\endaligned
\]
whence the point whose homogeneous coordinates possess the values:
\[
e_s
=
\sum_{\nu,\,\,k}^{1\cdots r}\,
c_{\nu ks}\,\lambda_\nu^0\,e_k^0
\ \ \ \ \ \ \ \ \ \ \ \ \ {\scriptstyle{(s\,=\,1\,\cdots\,r)}}
\]
obviously lie on the tangent. But this point visibly is the
image-point of the infinitesimal transformation:
\[
\bigg\leftbracket
\sum_{\nu=1}^r\,\lambda_\nu^0\,X_\nu f,
\,\,\,\,
\sum_{k=1}^r\,e_k^0\,X_kf
\bigg\rightbracket
=
\sum_{s=1}^r\,
\bigg\{
\sum_{\nu,\,\,k}^{1\cdots r}\,
c_{\nu ks}\,\lambda_\nu^0\,e_k^0
\bigg\}\,
X_sf
\]
which is obtained by \emphasis{combination} of the two infinitesimal
transformations $\sum\, \lambda_\nu^0\, X_\nu f$ and $\sum\, e_k^0\,
X_kf$. Consequently, we have the

\def\theproposition{5}\begin{proposition}
\label{Satz-5-S-288}
If one interprets the $\infty^{ r-1}$ infinitesimal transformations
$e_1\, X_1f + \cdots + e_r\, X_rf$ of an $r$-term group $X_1f, \dots,
X_rf$ as points of an $(r-1)$-times extended space by considering
$e_1, \dots, e_r$ as homogeneous coordinates in this space, then the
following happens: If all transformations of a determined one term
group: $\lambda_1^0 \, X_1f + \cdots + \lambda_r^0 \, X_rf$ are
executed on a determined infinitesimal transformation $e_1^0 \, X_1f +
\cdots + e_r^0 \, X_rf$, then the image-point of the transformation
$e_1^0 \, X_1f + \cdots + e_r^0 \, X_rf$ describes a curve, the
tangent of which at the point $e_1^0 \colon \cdots \colon e_r^0$ may
be obtained by connecting, through a straight line, this point to the
image-point of the infinitesimal transformation $\leftbracket \sum\,
\lambda_\nu^0 \, X_\nu f,\,\, \sum\, e_k^0 \, X_kf \rightbracket$; on
the other hand, if all transformations of the one-term group $\sum\,
e_k^0 \, X_kf$ are executed on the infinitesimal transformation
$\sum\, \lambda_\nu^0 \, X_\nu f$, then the image-point of the
transformation $\sum\, \lambda_\nu^0 \, X_\nu f$ describes a curve,
the tangent of which one obtains by connecting, through a straight
line, this point to the image-point of the infinitesimal
transformation $\leftbracket \sum\, e_k^0 \, X_kf, \, \, \sum\,
\lambda_\nu^0 \, X_\nu f \rightbracket$.
\end{proposition}

\linestop


\chapter{Composition and Isomorphism}
\label{kapitel-17}
\chaptermark{Composition and Isomorphism}

\setcounter{footnote}{0}

\abstract*{??}

Several problems that one can raise about an $r$-term group 
$X_1f, \dots, X_rf$ require, for their solution, 
only the knowledge of the constants $c_{ iks}$ in 
the relations: 
\[
\leftbracket
X_i,\,X_k
\rightbracket
=
\sum_{s=1}^r\,
c_{iks}\,X_sf.
\]

For instance, we have seen that the determination of all subgroups of
the group $X_1f, \dots, X_rf$ depends only on the constants $c_{ iks}$
and that exactly the same also holds true for the determination of all
types of subgroups (cf. Theorem~33, p.~\pageref{Theorem-33-S-210} and
Chap.~\ref{kapitel-16}, p.~\pageref{S-280} and \pageref{S-281}).

It is evident that the constants $c_{ iks}$ actually play the rôle of
certain properties of the group $X_1f, \dots, X_rf$. For the totality
of these properties, we introduce a specific terminology, we call them
the \terminology{composition} \deutsch{Zusammensetzung} of the group,
and we thus say that \emphasis{the constants $c_{ iks}$ in the
relations:
\def\theequation{1}\begin{equation}
\leftbracket
X_i,\,X_k
\rightbracket
=
\sum_{s=1}^r\,c_{iks}\,X_sf
\end{equation}
determine the \terminology{composition}
of the $r$-term group $X_1f, \dots, X_rf$}. 

\sectionengellie{\S\,\,\,80.}

The system of the $c_{ iks}$ which determines the composition of the
$r$-term group $X_1f, \dots, X_rf$ is in turn not completely
determined. Indeed, the individual $c_{ iks}$ receive in general
other numerical values when one chooses, in place of $X_1f, \dots,
X_rf$, any other $r$ independent infinitesimal transformations $e_1\,
X_1f + \cdots + e_r\, X_rf$.

From this, it follows that two different systems of $c_{ iks}$ can
represent, in certain circumstances, the composition of one and the
same group. But how can one recognize that this is the case?

We start from the relations:
\def\theequation{1}\begin{equation}
\leftbracket
X_i,\,X_k
\rightbracket
=
\sum_{s=1}^r\,c_{iks}\,X_sf
\ \ \ \ \ \ \ \ \ \ \ \ \ {\scriptstyle{(i,\,\,k\,=\,1\,\cdots\,r)}}
\end{equation}
which exist between $r$ determined independent infinitesimal
transformations $X_1f, \dots, X_rf$ of our group. We seek the general
form of the relations by which $r$ arbitrary independent infinitesimal
transformations $\mathcal{ X}_1f, \dots, \mathcal{ X}_rf$ of the group
$X_1f, \dots, X_rf$ are linked.

If the concerned relations have the form:
\def\theequation{2}\begin{equation}
\leftbracket
\mathcal{X}_i,\,\mathcal{X}_k
\rightbracket
=
\sum_{s=1}^r\,c_{iks}'\,\mathcal{X}_sf,
\end{equation}
then the system of the constants $c_{ iks}'$ is the most general one
which represents the composition of the group $X_1f, \dots, X_rf$.
Hence the thing is only about the computation of the $c_{ iks}'$.

Since $\mathcal{ X}_1f, \dots, \mathcal{ X}_rf$ are
supposed to be arbitrary
independent infinitesimal transformations 
of the group $X_1f, \dots, X_rf$, 
we have:
\[
\mathcal{X}_kf
=
\sum_{j=1}^r\,h_{kj}\,X_jf
\ \ \ \ \ \ \ \ \ \ \ \ \ {\scriptstyle{(k\,=\,1\,\cdots\,r)}},
\]
where the constants $h_{ kj}$ can take all the possible
values which do not bring to zero 
the determinant:
\[
D
=
{\textstyle{\sum}}\,\pm\,
h_{11}\cdots h_{rr}.
\]

By a calculation, it comes:
\[
\leftbracket
\mathcal{X}_i,\,\mathcal{X}_k
\rightbracket
=
\sum_{j,\,\,\pi}^{1\cdots r}\,
h_{ij}\,h_{k\pi}\,
\leftbracket
X_j,\,X_\pi
\rightbracket
=
\sum_{j,\,\,\pi,\,\,s}^{1\cdots r}\,
h_{ij}\,h_{k\pi}\,c_{j\pi s}\,X_sf\,;
\]
on the other hand, it follows from~\thetag{ 2}:
\[
\leftbracket
\mathcal{X}_i,\,\mathcal{X}_k
\rightbracket
=
\sum_{\pi,\,\,s}^{1\cdots r}\,
h_{\pi s}\,c_{ik\pi}'\,X_sf.
\]

If we compare with each other these two expressions for $\leftbracket
\mathcal{ X}_i, \, \mathcal{ X}_k \rightbracket$ and if we take
account of the fact that $X_1f, \dots, X_rf$ are independent
infinitesimal transformations, we obtain the relations:
\def\theequation{3}\begin{equation}
\sum_{\pi=1}^r\,h_{\pi s}\,c_{ik\pi}'
=
\sum_{j,\,\,\pi}^{1\cdots r}\,
h_{ij}\,h_{k\pi}\,c_{j\pi s}
\ \ \ \ \ \ \ \ \ \ \ \ \ {\scriptstyle{(s\,=\,1\,\cdots\,r)}}.
\end{equation}

Under the assumptions made, these equations can be solved with respect
to the $c_{ ik\pi}'$, hence one has\footnote{\,
Implicitly here, one
sees a standard way to 
represent the matrix which is the inverse of $(h_{ kj})$. 
}: 
\def\theequation{4}\begin{equation}
\label{S-291}
\aligned
c_{ik\rho}'
=
\frac{1}{D}\,
&
\sum_{s=1}^r\,
\bigg\{
\frac{\partial D}{\partial h_{\rho s}}\,
\sum_{j,\,\,\pi}^{1\cdots r}\,
h_{ij}\,h_{k\pi}
\bigg\}\,c_{j\pi s}
\\
& 
\ \ \ \
{\scriptstyle{(i,\,\,k,\,\,\rho\,=\,1\,\cdots\,r)}}
\endaligned
\end{equation}

With these words and according to the above, we have found the general
form of all systems of constants which determine the composition of
the group $X_1f, \dots, X_rf$. At the same time, we have at least
theoretical means to decide whether a given system of constants
$\overline{ c}_{ iks}$ determines the composition of the group $X_1f,
\dots, X_rf$; namely, such a system obviously possesses this property
if and only if one can choose the parameters $h_{kj}$ in such a way
that $c_{ iks}' = \overline{ c}_{ iks}$.\,---

If \emphasis{two} $r$-term groups are given, we can compare their
compositions. Clearly, thanks to the above developments, we are in a
position to decide whether the two groups have one and the same
composition, or have different compositions. Here, we do not need to
pay heed to the number of variables.

\plainstatement{We say that two $r$-term groups which have one and the
same composition are \label{S-291-bis}
\terminology{equally composed}
\deutsch{gleichzusammengesetzt}.}

If:
\[
X_kf
=
\sum_{i=1}^n\,\xi_{ki}(x_1,\dots,x_n)\,
\frac{\partial f}{\partial x_i}
\ \ \ \ \ \ \ \ \ \ \ \ \ {\scriptstyle{(k\,=\,1\,\cdots\,r)}}
\]
are independent infinitesimal transformations of an $r$-term group and
if:
\[
Y_kf
=
\sum_{\mu=1}^m\,\eta_{k\mu}(y_1,\dots,y_m)\,
\frac{\partial f}{\partial y_\mu}
\ \ \ \ \ \ \ \ \ \ \ \ \ {\scriptstyle{(k\,=\,1\,\cdots\,r)}}
\]
are independent infinitesimal transformations of a second $r$-term
group, and in addition, if there are relations:
\[
\leftbracket
X_i,\,X_k
\rightbracket
=
\sum_{s=1}^r\,c_{iks}\,X_sf,
\] 
\emphasis{then obviously, these two groups are equally composed when,
and only when, amongst the infinitesimal transformations $e_1\, Y_1 f
+ \cdots + e_r\, Y_rf$ of the second group, one can indicate $r$
mutually independent transformations $\mathcal{ Y}_1f, \dots, \mathcal{
Y}_rf$ such that the relations:
\[
\leftbracket
\mathcal{Y}_i,\,\mathcal{Y}_k
\rightbracket
=
\sum_{s=1}^r\,c_{iks}\,\mathcal{Y}_sf
\]
hold identically.}

If the relations holding between $Y_1f, \dots, Y_rf$ 
have the form:
\[
\leftbracket
Y_i,\,Y_k
\rightbracket
=
\sum_{s=1}^r\,\overline{c}_{iks}\,Y_sf,
\]
then we can say: the two groups are equally composed when and only
when it is possible to choose the parameters $h_{kj}$ in the
equations~\thetag{ 4} in such a way that every $c_{ iks}'$ is equal to
the corresponding $\overline{ c}_{ iks}$.

\medskip

One can also compare the compositions of two groups which do not have
the same number of parameters. This is made possible by the
introduction of the general concept: 
\label{Begriff-Isomorphismus}
\emphasis{isomorphism}
\deutsch{Isomorphismus}.

\renewcommand{\thefootnote}{\fnsymbol{footnote}}
\emphasis{The $r$-term group $X_1f, \dots, X_rf$:
\label{S-292}
\[
\leftbracket
X_i,\,X_k
\rightbracket
=
\sum_{s=1}^r\,c_{iks}\,X_sf
\]
is said to be \terminology{isomorphic} \deutsch{isomorph} to the $(r -
q)$-term group: $Y_1 f, \dots, Y_{ r-q} f$ when it is possible to
choose $r$ infinitesimal infinitesimal transformations:
\[
\aligned
\mathcal{Y}_kf
=
h_{k1}\,Y_1f
&
+\cdots+
h_{k,\,r-q}\,Y_{r-q}f
\\
& 
{\scriptstyle{(k\,=\,1\,\cdots\,r)}}
\endaligned
\]
in the $(r-q)$ so that not all $(r-q) \times (r-q)$ determinants of
the matrix:
\[
\left\vert
\begin{array}{cccc}
h_{11} & \,\cdot\, & \,\cdot\, & h_{1,r-q}
\\
\cdot & \,\cdot\, & \,\cdot\, & \cdot
\\
\cdot & \,\cdot\, & \,\cdot\, & \cdot
\\
h_{r1} & \,\cdot\, & \,\cdot\, & h_{r,r-q}
\end{array}
\right\vert
\]
and so that at the same time, the relations:
\[
\leftbracket
\mathcal{Y}_i,\,\mathcal{Y}_k
\rightbracket
=
\sum_{s=1}^r\,c_{iks}\,\mathcal{Y}_sf
\]
hold identically.}\footnote[1]{\,
Cf. Volume~III, \cite{enlie1893-17}, p.~701, remarks.
}
\renewcommand{\thefootnote}{\arabic{footnote}}

Let there be an isomorphism in this sense and let $\mathcal{ Y}_1f,
\dots, \mathcal{ Y}_rf$ be already chosen in the indicated way. Then
if we associate to the infinitesimal transformation $e_1\, X_1f +
\cdots + e_r\, X_rf$ of the $r$-term group always the infinitesimal
transformation:
\[
e_1\,\mathcal{Y}_1f
+\cdots+
e_r\,\mathcal{Y}_rf
\]
of the $(r-q)$-term group, whichever values the constants $e_1, \dots,
e_r$ may have, then the following clearly holds true: when ${\sf
Y}_1f$ is the transformation of the $(r-q)$-term group which is
associated to the transformation $\Xi_1 f$ of the other group, and
when, correspondingly, ${\sf Y}_2 f$ is associated to the
transformation $\Xi_2f$, then the transformation $\leftbracket
\Xi_1,\,\Xi_2 \rightbracket$ always corresponds to the transformation
$\leftbracket {\sf Y}_1,\,{\sf Y}_2 \rightbracket$. We express this
more briefly as follows: through the indicated correspondence of the
infinitesimal transformations of the two groups, the groups are
\emphasis{isomorphically related one to another}. Visibly, this
isomorphic condition is completely determined when one knows that, to
$X_1f, \dots, X_rf$, are associated the transformations $\mathcal{
Y}_1f, \dots, \mathcal{ Y}_rf$, respectively.

One makes a distinction between \terminology{holoedric} and
\terminology{meroedric} \emphasis{isomorphisms}. 
\label{S-293} The
\emphasis{holoedric} case occurs when the number $q$, which appears in
the definition of the isomorphism, has the value zero; the
\emphasis{meroedric} case when $q$ is larger than zero.
Correspondingly, one says that the two groups are holoedrically, or
meroedrically, isomorphic one to another.

Visibly, the property of being equally composed \deutsch{die
Eigenschaft des Gleichzusammengesetztseins} of two groups is a special
case of isomorphism; namely, two equally composed groups
are always holoedrically isomorphic, and conversely.

Two meroedrically isomorphic groups are for example
the two: 
\[
\frac{\partial f}{\partial x_1},\ \ \ \ \
x_1\,\frac{\partial f}{\partial x_1},\ \ \ \ \
x_1^2\,\frac{\partial f}{\partial x_1},\ \ \ \ \
\frac{\partial f}{\partial x_2}
\]
and:
\[
\frac{\partial f}{\partial y_1},\ \ \ \ \
y_1\,\frac{\partial f}{\partial y_1},\ \ \ \ \
y_1^2\,\frac{\partial f}{\partial y_1},
\]
with, respectively, four and three parameters.
We obtain that these two groups are meroedrically 
isomorphic when, to the four infinitesimal 
transformations:
\[
X_1f
=
\frac{\partial f}{\partial x_1},\ \ \
X_2f
=
x_1\,\frac{\partial f}{\partial x_1},\ \ \
X_3f
=
x_1^2\,\frac{\partial f}{\partial x_1},\ \ \
X_4f
=
\frac{\partial f}{\partial x_1}
+
\frac{\partial f}{\partial x_2}
\]
of the first, the following four, say, are associated:
\[
\mathcal{Y}_1f
=
\frac{\partial f}{\partial y_1},\ \ \
\mathcal{Y}_2f
=
y_1\,\frac{\partial f}{\partial y_1},\ \ \
\mathcal{Y}_3f
=
y_1^2\,\frac{\partial f}{\partial y_1},\ \ \
\mathcal{Y}_4f
=
\frac{\partial f}{\partial y_1}.
\] 

\renewcommand{\thefootnote}{\fnsymbol{footnote}}
In the theory of substitutions, one also speaks of isomorphic groups,
although the notion of isomorphism \label{S-293-bis} happens to be
apparently different from the one here\footnote[1]{\,
Camille \name{Jordan}, Traité des substitutions, Paris 1870.
}. 
\renewcommand{\thefootnote}{\arabic{footnote}}
Later (cf. Chap.~\ref{kapitel-21}: {\em The Group of Parameters}) we
will convince ourselves that nevertheless, the concept of isomorphism
following from our definition corresponds perfectly to the concept
which one obtains as soon as one translates the definition of the
theory of substitution directly into the theory of finite continuous
groups.

At present, we shall at first derive a few simple consequences from
our definition of isomorphism.

In the preceding chapter, we have seen that to every $r$-term group
$X_1f, \dots, X_rf$ is associated a certain linear homogeneous group,
the adjoint group, as we have named it. From the relations which
exist between the infinitesimal transformations of the adjoint group
(cf. Chap.~\ref{kapitel-16}, p.~\pageref{S-274}), it becomes
immediately evident that the group $X_1f, \dots, X_rf$ is isomorphic
to its adjoint group; however, the two groups are holoedrically
isomorphic only when the group $X_1f, \dots, X_rf$ contains no
excellent infinitesimal transformation, because the adjoint group is
$r$-term only in this case, whereas it always contains less than $r$
parameters in the contrary case. Thus:

\def\thetheorem{55}\begin{theorem}
To every $r$-term group $X_1f, \dots, X_rf$ is associated an
isomorphic linear homogeneous group, namely the adjoint group; this
group is holoedrically isomorphic to the group $X_1f, \dots, X_rf$
only when the latter contains no excellent infinitesimal
transformation.
\end{theorem}

\smallercharacters{
We do not consider here the question of whether to every $r$-term
group which contains excellent infinitesimal transformations one can
also associate an holoedrically isomorphic linear homogeneous group.
Yet by an example, we want to show that this is in any case possible
in many circumstances, also when the given $r$-term group contains
excellent infinitesimal transformations.

Let the $r$-term group $X_1f, \dots, X_rf$ contains precisely $r - m$
independent excellent infinitesimal transformations, and let $X_1f,
\dots, X_rf$ be chosen in such a way that $X_{ m+1}f, \dots, X_rf$ are
excellent infinitesimal transformations; then between $X_1f, \dots,
X_rf$, there exist relations of the form:
\[
\aligned
&
\leftbracket
X_\mu,\,X_\nu
\rightbracket
=
c_{\mu\nu 1}\,X_1f
+\cdots+
c_{\mu\nu r}\,X_rf
\\
&\ \ \ \
\leftbracket
X_\mu,\,X_{m+k}
\rightbracket
=
\leftbracket
X_{m+k},\,X_{m+j}
\rightbracket
=
0
\\
&
\ \ \ \ \ \ \ \ \ \ \
{\scriptstyle{(\mu,\,\,\nu\,=\,1\,\cdots\,m\,;\,\,\,
k,\,\,j\,=\,1\,\cdots\,r\,-\,m)}}.
\endaligned
\]

In the associated adjoint group $E_1f, \dots, E_rf$, there are only
$m$ independent infinitesimal transformations: $E_1f, \dots, E_mf$,
while $E_{ m+1}f, \dots, E_rf$ vanish identically, and hence $E_1f,
\dots, E_mf$ are linked together by the relations:
\[
\leftbracket
E_\mu,\,E_\nu
\rightbracket
=
c_{\mu\nu 1}\,E_1f
+\cdots+
c_{\mu\nu m}\,E_mf.
\]

In particular, if all $c_{ \mu,\nu, m+1}, \dots, c_{\mu\nu r}$ vanish,
one can always indicate an $r$-term linear homogeneous group which is
holoedrically isomorphic to the group $X_1f, \dots, X_rf$. In this
case namely, $X_1f, \dots, X_mf$ actually generate an $m$-term group
to which the group $E_1f, \dots, E_mf$ is holoedrically
isomorphic. Hence if we set:
\[
{\sf E}_{m+1}f
=
e_{r+1}\,\frac{\partial f}{\partial e_{r+1}},
\,\,\,\dots,\,\,\,
{\sf E}_rf
=
e_{2r-m}\,\frac{\partial f}{\partial e_{2r-m}},
\]
then the $r$ independent infinitesimal transformations:
\[
E_1f,\,\dots,\,E_mf,\,\,\,
{\sf E}_{m+1}f,\,\dots,\,{\sf E}_rf
\]
obviously generate a linear homogeneous group which is
holoedrically isomorphic to the group $X_1f, \dots, X_rf$. 

But also in the cases where not all $c_{ \mu, \nu, m+1}, \dots, c_{
\mu \nu r}$ vanish, one can often easily indicate an holoedrically
isomorphic linear homogeneous group.
As an example, we use the three-term group
$X_1f$, $X_2f$, $X_3f$: 
\[
\leftbracket
X_1,\,X_2
\rightbracket
=
X_3f,\ \ \ \ \
\leftbracket
X_1,\,X_3
\rightbracket
=
\leftbracket
X_2,\,X_3
\rightbracket
=
0,
\]
which contains an excellent infinitesimal transformation, 
namely $X_3f$; it is holoedrically isomorphic to the
linear homogeneous group:
\[
{\sf E}_1f
=
\alpha_3\,
\frac{\partial f}{\partial\alpha_1},\ \ \ \ \
{\sf E}_2f
=
\alpha_1\,
\frac{\partial f}{\partial\alpha_2},\ \ \ \ \
{\sf E}_3f
=
\alpha_3\,
\frac{\partial f}{\partial\alpha_2}.
\]

}

As we have seen in Chap.~\ref{kapitel-9}, p.~\pageref{S-170}, 
the constants $c_{ iks}$ in the equations:
\def\theequation{1}\begin{equation}
\leftbracket
X_i,\,X_k
\rightbracket
=
\sum_{s=1}^r\,c_{iks}\,X_sf
\end{equation}
satisfy the relations:
\def\theequation{5}\begin{equation}
\left\{
\aligned
&
\ \ \ \ \ \ \ \ \ \ \ \ \ \ \ \ \ \ 
c_{iks}+c_{kis}=0
\\
&
\sum_{\nu=1}^r\,
\big\{
c_{ik\nu}\,c_{\nu js}+c_{kj\nu}\,c_{\nu is}+c_{ji\nu}\,c_{\nu ks}
\big\}
=
0
\\
&\ \ \ \ \ \ \ \ \ \ \ \ \ \ \ \ \ \ \ \ \ 
{\scriptstyle{(i,\,\,k,\,\,j,\,\,s\,=\,1\,\cdots\,r)}}.
\endaligned\right.
\end{equation}

With the help of these relations, we succeeded, in
Chap.~\ref{kapitel-16}, p.~\pageref{S-274}, to prove that the $r$
infinitesimal transformations:
\def\theequation{6}\begin{equation}
E_\mu f
=
\sum_{k,\,\,j}^{1\cdots r}\,
c_{j\mu k}\,e_j\,
\frac{\partial f}{\partial e_k}
\ \ \ \ \ \ \ \ \ \ \ \ \ {\scriptstyle{(\mu\,=\,1\,\cdots\,r)}}
\end{equation}
of the group adjoint to the group $X_1f, \dots, X_rf$ stand pairwise
in the relationships:
\def\theequation{7}\begin{equation}
\leftbracket
E_\mu,\,E_\nu
\rightbracket
=
\sum_{s=1}^r\,c_{\mu\nu s}\,E_sf.
\end{equation}

But for this proof of the relations~\thetag{ 7}, we have used no more
than the fact that the $c_{ iks}$ satisfied the equations~\thetag{ 5},
namely we have made no use of the fact that we knew $r$ infinitesimal
transformations $X_1f, \dots, X_rf$ which were linked together by the
relations~\thetag{ 1}. Thanks to the cited developments, it is hence
established that the $r$ infinitesimal transformations~\thetag{ 6}
always stand in the relationships~\thetag{ 7} when the $c_{ iks}$
satisfy the equations~\thetag{ 5}.

Consequently, if we know a system of $c_{ iks}$ which satisfies the
relations~\thetag{ 5}, we can immediately 
indicate $r$ linear homogeneous infinitesimal
transformations $E_1f, \dots, E_rf$, namely 
the transformations~\thetag{ 6}, 
which stand pairwise in the relationships:
\[
\leftbracket
E_i,\,E_k
\rightbracket
=
\sum_{s=1}^r\,c_{iks}\,E_sf.
\]

It is evident that the so obtained infinitesimal transformations
$E_1f, \dots, E_rf$ generate a group, and to be precise, a group with
$r$ or less parameters; clearly, they generate a group with exactly
$r$ parameters only when they are mutually independent, hence when it
is impossible that the equations:
\[
g_1\,E_1f
+\cdots+
g_r\,E_rf
=
0,
\]
or the $r^2$ equivalent equations:
\[
g_1\,c_{j1k}
+
g_2\,c_{j2k}
+\cdots+
g_r\,c_{jrk}
\ \ \ \ \ \ \ \ \ \ \ \ \ 
{\scriptstyle{(j,\,\,k\,=\,1\,\cdots\,r)}}
\]
are satisfied by not all vanishing quantities $g_1, \dots, g_r$.

As a result, we have the

\renewcommand{\thefootnote}{\fnsymbol{footnote}}
\def\thetheorem{52}\begin{theorem}\footnote[1]{\,
\name{Lie}, Archiv for Math. og Nat. Vol. 1, 
p. 192, Christiania 1876.
} 
\label{Theorem-52-S-296}
When the constants $c_{ iks}$ ${\scriptstyle{( i,\,\, k,\,\,s\, =
\,1\, \cdots\, r)}}$ possess values such that all relations of the
form:
\def\theequation{5}\begin{equation}
\left\{
\aligned
&
\ \ \ \ \ \ \ \ \ \ \ \ \ \ \ \ \ \ 
c_{iks}+c_{kis}=0
\\
&
\sum_{\nu=1}^r\,
\big\{
c_{ik\nu}\,c_{\nu js}+c_{kj\nu}\,c_{\nu is}+c_{ji\nu}\,c_{\nu ks}
\big\}
=
0
\\
&\ \ \ \ \ \ \ \ \ \ \ \ \ \ \ \ \ \ \ \ \ 
{\scriptstyle{(i,\,\,k,\,\,j,\,\,s\,=\,1\,\cdots\,r)}}
\endaligned\right.
\end{equation}
are satisfied, then the $r$ linear homogeneous
infinitesimal transformations:
\[
E_\mu f
=
\sum_{k,\,\,j}^{1\cdots r}\,
c_{j\mu k}\,e_j\,
\frac{\partial f}{\partial e_k}
\ \ \ \ \ \ \ \ \ \ \ \ \ 
{\scriptstyle{(\mu\,=\,1\,\cdots\,r)}}
\]
stand pairwise in the relationships:
\[
\leftbracket
E_i,\,E_k
\rightbracket
=
\sum_{s=1}^r\,c_{iks}\,E_sf
\ \ \ \ \ \ \ \ \ \ \ \ \ 
{\scriptstyle{(i,\,\,k\,=\,1\,\cdots\,r)}},
\]
and hence, they generate a linear homogeneous group. 
In particular, if the $c_{ iks}$ are constituted so 
that not all $r \times r$ 
determinants, the horizontal series of which have
the form:
\[
\left\vert
c_{j1k} \ c_{j2k} \ \cdots \ c_{jrk}
\right\vert
\ \ \ \ \ \ \ \ \ \ \ \ \ 
{\scriptstyle{(j,\,\,k\,=\,1\,\cdots\,r)}},
\]
vanish, then $E_1f, \dots, E_rf$ are independent infinitesimal
transformations and they generate an $r$-term group whose composition
is determined by the system of the $c_{ iks}$, and which contains no
excellent infinitesimal transformation. In all other cases, the group
generated by $E_1f, \dots, E_rf$ has less than $r$ parameters.
\end{theorem}
\renewcommand{\thefootnote}{\arabic{footnote}}

\sectionengellie{\S\,\,\,81.}

The results of the preceding paragraph suggest the presumption that
actually, every system of $c_{ iks}$ which satisfies the
relations~\thetag{ 5} represents the composition of a certain $r$-term
group. This presumption, corresponds to the truth \deutsch{Diese
Vermuthung entspricht der Wahrheit}, for the following really 
holds\footnote{\,
This is the so-called \terminology{Third Fundamental Theorem}
of Lie's theory, cf. Vol.~III.
}. 

\def\theproposition{1}\begin{proposition}
\label{Proposition-1-S-297}
If the constants $c_{ iks}$
${\scriptstyle{(i,\,\,k,\,\,s\,=\,1\,\cdots\,r)}}$ possess values such
that the relations:
\def\theequation{5}\begin{equation}
\left\{
\aligned
&
\ \ \ \ \ \ \ \ \ \ \ \ \ \ \ \ \ \ 
c_{iks}+c_{kis}=0
\\
&
\sum_{\nu=1}^r\,
\big\{
c_{ik\nu}\,c_{\nu js}+c_{kj\nu}\,c_{\nu is}+c_{ji\nu}\,c_{\nu ks}
\big\}
=
0
\\
&\ \ \ \ \ \ \ \ \ \ \ \ \ \ \ \ \ \ \ \ \ 
{\scriptstyle{(i,\,\,k,\,\,j,\,\,s\,=\,1\,\cdots\,r)}}
\endaligned\right.
\end{equation}
are satisfied, then there are always, in a space of the appropriate
number of dimensions, $r$ independent infinitesimal transformations
$X_1f, \dots, X_rf$ which stand pairwise in the relationships:
\[
\leftbracket
X_i,\,X_k
\rightbracket
=
\sum_{s=1}^r\,
c_{iks}\,X_sf
\]
and hence generate an $r$-term group of the composition $c_{ iks}$.
\end{proposition}

For the time being, we suppress the proof of this important
proposition, in order not to be forced to digress for a long while,
and we will perform this proof only in the next volume. Of course,
until then, we will make the least possible use of this proposition.

From the Proposition~1, it results that the totality of all possible
compositions of $r$-term groups is represented by the totality of all
systems of $c_{ iks}$ which satisfy the equations~\thetag{ 5}. If one
knows all such systems of $c_{ iks}$, then with this, one knows at the
same time all compositions of $r$-term groups.

But now, as we have seen on p.~\pageref{S-291}, there are in general
infinitely many systems of $c_{ iks}$ which represent one and the same
composition; if a system of $c_{ iks}$ is given which represents a
composition, then one finds all systems of $c_{ iks}'$ which represent
the same composition by means of the equations~\thetag{ 4}, in which
it is understood that the $h_{ kj}$ are arbitrary parameters. Hence,
when one knows all systems of $c_{ iks}$ which satisfy the
equations~\thetag{ 5}, a specific research is yet required in order to
find out which of these systems represent different compositions. In
order to be able to execute this research, we must at first consider
somehow more closely the equations~\thetag{ 4}.

For the moment, we disregard the fact that the $c_{ iks}$ are linked
together by some relations; rather, we consider the $c_{ iks}$, and
likewise the $c_{ iks}'$, as variables independent of each other. On
the basis of this conception, it will be shown that the
equations~\thetag{ 4} represent a continuous transformation group in
the variables $c_{ iks}$.

In order to prove the claimed property of the equations~\thetag{ 4},
we will directly execute one after the other two
transformations~\thetag{ 4}, or, what is the same, two
transformations:
\def\theequation{3}\begin{equation}
\aligned
\sum_{\pi=1}^r\,h_{\pi s}\,
&
c_{ik\pi}'
=
\sum_{j,\,\,\pi}^{1\cdots r}\,
h_{ij}\,h_{k\pi}\,c_{j\pi s}
\\
&
{\scriptstyle{(i,\,\,k,\,\,s\,=\,1\,\cdots\,r)}}.
\endaligned
\end{equation}

To begin with, we therefore transfer the $c_{ iks}$ to the $c_{ iks}'$
by means of the transformation~\thetag{ 3} and then the $c_{ iks}'$ to
the $c_{ iks}''$ by means of the transformation:
\def\theequation{3'}\begin{equation}
\sum_{\pi=1}^r\,h_{\pi s}'\,c_{ik\pi}''
=
\sum_{j,\,\,\pi}^{1\cdots r}\,
h_{ij}'\,h_{k\pi}'\,c_{j\pi s}'.
\end{equation}

In this way, we obtain a new transformation, the equations of which
are obtained when the $c_{ iks}'$ are took away from~\thetag{ 3}
and~\thetag{ 3'}. It is to be proved that this new transformation has
the form:
\def\theequation{3''}\begin{equation}
\sum_{\pi=1}^r\,
h_{\pi s}''\,c_{ik\pi}''
=
\sum_{j,\,\,\pi}^{1\cdots r}\,
h_{ij}''\,h_{k\pi}''\,c_{j\pi s}',
\end{equation}
where the $h''$ are functions of only the $h$ and the $h'$.

We multiply~\thetag{ 3'} by $h_{s \sigma}$ and we sum with respect to
$s$; then we receive:
\[
\sum_{\pi,\,\,s}^{1\cdots r}\,
h_{\pi s}'\,h_{s\sigma}\,c_{ik\pi}''
=
\sum_{j,\,\,\pi,\,\,s}^{1\cdots r}\,
h_{ij}'\,h_{k\pi}'\,h_{s\sigma}\,c_{j\pi s}',
\]
or, because of~\thetag{ 3}:
\[
\sum_{\pi,\,\,s}^{1\cdots r}\,
h_{\pi s}'\,h_{s\sigma}\,c_{ik\pi}''
=
\sum_{j,\,\,\pi,\,\,\tau,\,\,\rho}^{1\cdots r}\,
h_{ij}'\,h_{k\pi}'\,h_{j\tau}\,h_{\pi\rho}\,c_{\tau\rho\sigma}.
\]

This is the discussed new transformation; it converts into~\thetag{
3''} when one sets:
\[
h_{\pi\sigma}''
=
\sum_{s=1}^r\,h_{\pi s}'\,h_{s\sigma}.
\]

As a result, it is proved that the transformations~\thetag{ 4}
effectively form a group.

Now, we claim that the transformations of this group leave invariant
the equations~\thetag{ 5}.

Let $c_{ iks}$ be a system of constants which satisfies the
relations~\thetag{ 5}, hence according to Proposition~1, which
represents the composition of a certain $r$-term group. Then as we
know, the system of the $c_{ iks}'$ which is determined by the
relations~\thetag{ 4} represents in the same way a composition, namely
the same composition as the one of the system of the $c_{ iks}$;
consequently, the $c_{ iks}'$ also satisfy relations of the
form~\thetag{ 5}. Thus, the transformations~\thetag{ 4} transfers
every system $c_{ iks}$ which satisfies~\thetag{ 5} to a system $c_{
iks}'$ having the same constitution, that is to say, they leave
invariant the system of equations~\thetag{ 5}, and this was just our
claim.

At present, we interpret the $r^3$ variables $c_{ iks}$ as point
coordinates in a space of $r^3$ dimensions.

In this space, a certain manifold $M$ which is invariant by the
transformations of the group~\thetag{ 4} is sorted by the
equations~\thetag{ 5}. Every point of $M$\,---\,as we can
say\,---\,represents a composition of $r$-term group, and conversely,
every possible composition of $r$-term group is represented by certain
points of $M$. Two different points of $M$ represent one and the same
composition when there is, in the group~\thetag{ 4}, a transformation
which transfers the first point to the other.

Hence, if $P$ is an arbitrary point of $M$, then the totality of all
positions that the point $P$ takes by the transformations of the
group~\thetag{ 4} coincides with the totality of all points which
represent the same composition as $P$. We know from before
(Chap.~\ref{kapitel-14}, p.~\pageref{Theorem-36-S-225}) that that this
totality of points forms a manifold invariant by the group~\thetag{
4}, and to be precise, a so-called smallest invariant manifold.

From this, it follows that one has to proceed as follows in order to
find all different types of $r$-term groups:

One determines all smallest manifolds located in $M$ that remain
invariant by the group~\thetag{ 4}; on \emphasis{each} such manifold,
one chooses an arbitrary point $c_{ iks}$: the systems of values $c_{
iks}$ which belong to the chosen points then represent all types of
different compositions of $r$-term groups.

Since the finite equations of the group~\thetag{ 4} are here, the
determination of the smallest invariant manifolds has to be considered
as an executable operation; it requires only the resolution of
algebraic equations.

We therefore have the

\def\thetheorem{53}\begin{theorem}
\label{Theorem-3-S-300}
The determination of all essentially different compositions
of $r$-term groups requires only algebraic operations.
\end{theorem}

\sectionengellie{\S\,\,\,82.}

Now, let $X_1f, \dots, X_rf$ be an $r$-term group $G_r$
of the composition:
\[
\leftbracket
X_i,\,X_k
\rightbracket
=
\sum_{s=1}^r\,c_{iks}\,X_sf.
\]
Furthermore, let $Y_1f, \dots, Y_{ r-q}f$ be an $(r-q)$-term group
isomorphic to the $G_r$, and to be precise, meroedrically isomorphic,
so that $q$ is therefore larger than zero. We want to denote this
second group shortly by $G_{ r-q}$.

Let the two groups be, in the way indicated on p.~\pageref{S-292}
and~\pageref{S-293}, isomorphically related one to the other; so in the
$G_{ r-q}$, let $r$ infinitesimal transformations $\mathcal{ Y}_1f,
\dots, \mathcal{ Y}_rf$ be chosen which stand in the relationships:
\[
\leftbracket
\mathcal{Y}_i,\,
\mathcal{Y}_k
\rightbracket
=
\sum_{s=1}^r\,c_{iks}\,\mathcal{Y}_sf,
\]
where $\mathcal{ Y}_1f, \dots, \mathcal{ Y}_{ r-q}f$ are mutually
independent, whereas $\mathcal{ Y}_{ r-q+1}f, \dots, \mathcal{ Y}_rf$
are defined by the identities:
\def\theequation{8}\begin{equation}
\aligned
\mathcal{Y}_{r-q+k}f
\equiv
d_{k1}\,
&
\mathcal{Y}_1f
+\cdots+
d_{k,\,r-q}\,\mathcal{Y}_{r-q}f
\\
&
{\scriptstyle{(k\,=\,1\,\cdots\,q)}}.
\endaligned
\end{equation}

Under the assumptions made, one obviously has:
\[
\aligned
\bigg\leftbracket
\mathcal{Y}_if,\,\,\,
&
\mathcal{Y}_{r-q+k}f
-
\sum_{\mu=1}^{r-q}\,d_{k\mu}\,\mathcal{Y}_\mu f
\bigg\rightbracket
\equiv
0
\\
&\ \ \ \
{\scriptstyle{(j\,=\,1\,\cdots\,r\,;\,\,\,
k\,=\,1\,\cdots\,q)}},
\endaligned
\]
or:
\[
\sum_{s=1}^r\,
\bigg\{
c_{j,\,r-q+k,\,s}
-
\sum_{\mu=1}^{r-q}\,
d_{k\mu}\,c_{j\mu s}
\bigg\}\,
\mathcal{Y}_sf
\equiv
0.
\]

If we replace here $\mathcal{ Y}_{ r-q+1}f, \dots, \mathcal{ Y}_rf$ by
their values~\thetag{ 8}, we obtain linear relations between
$\mathcal{ Y}_1f, \dots, \mathcal{ Y}_{ r-q}f$; but obviously, such
linear relations can hold only when the coefficients of every
individual $\mathcal{ Y}_1f, \dots, \mathcal{ Y}_{ r-q}f$ possess
the value zero.

From the vanishing of these coefficients, it follows that
the $rq$ expressions: 
\[
\bigg\leftbracket
X_jf,\,\,\,
X_{r-q+k}f
-
\sum_{\mu=1}^{r-q}\,d_{k\mu}\,X_\mu f
\bigg\rightbracket
=
\sum_{s=1}^r\,
\bigg\{
c_{j,\,r-q+k,\,s}
-
\sum_{\mu=1}^{r-q}\,d_{k\mu}\,c_{j\mu s}
\bigg\}\,
X_sf
\]
can be linearly deduced from the $q$ infinitesimal transformations:
\def\theequation{9}\begin{equation}
X_{r-q+k}f
-
\sum_{\mu=1}^{r-q}\,d_{k\mu}\,X_\mu f
\ \ \ \ \ \ \ \ \ \ \ \ \ {\scriptstyle{(k\,=\,1\,\cdots\,q)}}.
\end{equation}
Expressed differently: 
the $q$ independent infinitesimal transformations~\thetag{ 9}
generate a $q$-term \emphasis{invariant}
subgroup of the group $X_1f, \dots, X_rf$. 

\def\thetheorem{54}\begin{theorem}
\label{Theorem-54-S-301}
If the $(r-q)$-term group $Y_1f, \dots, Y_{ r-q}f$
is isomorphic to the $r$-term group: 
$X_1f, \dots, X_rf$, and if
$\mathcal{ Y}_1f, \dots, \mathcal{ Y}_rf$
are infinitesimal transformations of the
$(r-q)$-term group such that firstly
$\mathcal{ Y}_1f, \dots, \mathcal{ Y}_{ r-q}f$ are mutually
independent while by contrast, $\mathcal{ Y}_{ r-q+1}, \dots, 
\mathcal{ Y}_rf$ can be linearly deduced from 
$\mathcal{ Y}_1f, \dots, \mathcal{ Y}_{ r-q}f$: 
\[
\mathcal{Y}_{r-q+k}f
\equiv
d_{k1}\,\mathcal{Y}_1f
+\cdots+
d_{k,\,r-q}\,\mathcal{Y}_{r-q}f
\ \ \ \ \ \ \ \ \ \ \ \ \ {\scriptstyle{(k\,=\,1\,\cdots\,q)}},
\]
and secondly such that, simultaneously with the relations:
\[
\leftbracket
X_i,\,X_k
\rightbracket
=
\sum_{s=1}^r\,c_{iks}\,X_sf,
\]
the analogous relations:
\[
\leftbracket
\mathcal{Y}_i,\,\mathcal{Y}_k
\rightbracket
=
\sum_{s=1}^r\,c_{iks}\,\mathcal{Y}_sf
\]
hold, then the $q$ infinitesimal transformations:
\[
X_{r-q+k}f
-
d_{k1}\,X_1f
-\cdots-
d_{k,\,r-q}\,X_{r-q}f
\ \ \ \ \ \ \ \ \ \ \ \ \ {\scriptstyle{(k\,=\,1\,\cdots\,q)}}
\]
generate a $q$-term invariant subgroup of the group
$X_1f, \dots, X_rf$. 
\end{theorem}

The $G_r$: $X_1f, \dots, X_rf$ and the $G_{ r-q}$: $Y_1f, \dots, Y_{
r-q}f$ are, as we know, isomorphically related when, to every
infinitesimal transformation of the form: $e_1\, X_1f + \cdots + e_r\,
X_rf$, we associate the infinitesimal transformation:
\[
\sum_{k=1}^r\,e_k\,\mathcal{Y}_kf
=
\sum_{k=1}^{r-q}\,
\bigg\{
e_k
+
\sum_{j=1}^q\,e_{r-q+j}\,d_{jk}
\bigg\}\,\mathcal{Y}_kf.
\]

Through this correspondence, the transformations of the $G_r$ which
belong to the $q$-term group~\thetag{ 9} are the only ones which
correspond to the identically vanishing infinitesimal transformation
of the $G_{ r-q}$. Consequently, in general, to every one-term
subgroup of the $G_r$, there corresponds a completely determined
one-term subgroup of the $G_{ r-q}$, and it is only to the one-term
subgroups of the group~\thetag{ 9} that there correspond no one-term
subgroups of the $G_{ r-q}$, for the associated one-term groups indeed
reduce to the identity transformation. Conversely, to one and the
same one-term subgroup $h_1\, \mathcal{ Y}_1f + \cdots + h_{ r-q} \,
\mathcal{ Y}_{ r-q} f$ of the $G_{ r-q}$, there correspond in total
$\infty^q$ different one-term subgroups of the $G_r$, namely all the
groups of the form:
\[
\sum_{k=1}^{r-q}\,h_k\,X_kf
+
\sum_{j=1}^q\,\lambda_j\,
\bigg\{
X_{r-q+j}f
-
\sum_{\mu=1}^{r-q}\,d_{j\mu}\,X_\mu f
\bigg\}\,,
\]
where $\lambda_1, \dots, \lambda_q$ denote arbitrary constants.

Now, a certain correspondence between the subgroups of the $G_r$ and
the subgroups of the $G_{ r-q}$ actually takes place.

If $m$ arbitrary mutually independent infinitesimal transformations:
\[
l_{\mu 1}\,X_1f
+\cdots+
l_{\mu r}\,X_rf
\ \ \ \ \ \ \ \ \ \ \ \ \ {\scriptstyle{(\mu\,=\,1\,\cdots\,m)}}
\]
generate an $m$-term subgroup of the $G_r$, then
the $m$ infinitesimal transformations:
\[
\aligned
\sum_{k=1}^r\,l_{\mu k}\,\mathcal{Y}_kf
=
\sum_{k=1}^{r-q}\,
&
\bigg\{
l_{\mu k}
+
\sum_{j=1}^q\,l_{\mu,\,r-q+j}\,d_{jk}
\bigg\}\,\mathcal{Y}_kf
\\
&
{\scriptstyle{(\mu\,=\,1\,\cdots\,m)}}
\endaligned
\]
obviously generate a subgroup of the $G_{ r-q}$. This subgroup is at
most $m$-term, and in particular, it is $0$-term, that is to say, it
consists of only the identity transformation when the $m$-term
subgroup of the $G_r$ is contained in the $q$-term 
group~\thetag{ 9}, and in fact clearly, only in this
case. 

Conversely, if $m'$ arbitrary infinitesimal transformations:
\[
l_{\mu 1}\,\mathcal{Y}_1f
+\cdots+
l_{\mu,\,r-q}\,\mathcal{Y}_{r-q}f
\ \ \ \ \ \ \ \ \ \ \ \ \ {\scriptstyle{(\mu\,=\,1\,\cdots\,m')}}
\]
generate an $m'$-term subgroup of the $G_{ r-q}$: $Y_1f, \dots, Y_{
r-q}f$, then the $m'$ infinitesimal transformations:
\[
l_{\mu 1}X_1f
+\cdots+
l_{\mu,\,r-q}\,X_{r-q}f
\ \ \ \ \ \ \ \ \ \ \ \ \ {\scriptstyle{(\mu\,=\,1\,\cdots\,m')}},
\]
together with the $q$ transformations:
\[
X_{r-q+k}f
-
d_{k1}\,X_1f
-\cdots-
d_{k,\,r-q}\,X_{r-q}f
\ \ \ \ \ \ \ \ \ \ \ \ \ {\scriptstyle{(k\,=\,1\,\cdots\,q)}}
\]
always generate an $(m'+q)$-term subgroup of the $G_r$.

From this, we see that to every subgroup of the $G_r$ corresponds a
completely determined subgroup of the $G_{ r-q}$ which, in certain
circumstances, consists of only the identity transformation; and
moreover, we see that to every subgroup of the $G_{ r-q}$ corresponds
at least one subgroup of the $G_r$. If we know all subgroups of the
$G_r$ and if we determine all the subgroups of the $G_{ r-q}$
corresponding to them, we then obtain all subgroups of the $G_{ r-q}$.
Thus, the following holds.

\def\theproposition{2}\begin{proposition}
If one has isomorphically related an $r$-term group of which one
knows all subgroups to an $(r-q)$-term group, then one can also
indicate straightaway all subgroups of the
$(r-q)$-term group. 
\end{proposition}

Let $X_1f, \dots, X_rf$, or shortly $G_r$, be an $r$-term
group of the composition:
\[
\leftbracket
X_i,\,X_k
\rightbracket
=
\sum_{s=1}^r\,c_{iks}\,X_sf.
\]

We want to study what are the different compositions that a group
meroedrically isomorphic to the $G_r$ can have.

Up to now, we know only the following: If the $G_r$ can be
isomorphically related to an $(r-q)$-term group, then there is in the
$G_r$ a completely determined $q$-term invariant subgroup which
corresponds to the identity transformation in the $(r-q)$-term group.
Now, we claim that this proposition can be reversed in the following
way: If the $G_r$ contains a $q$-term invariant subgroup, then there
always is an $(r-q)$-term group $G_{ r-q}$ which is isomorphic to the
$G_r$ and which can be isomorphically related to the $G_r$ in such a
way that the $q$-term invariant group inside the $G_r$ corresponds to
the identity transformation\footnote{\,
The notion of \terminology{quotient group} will in fact not
come up here; instead, Lie uses Proposition~1 above.
} 
in the $G_{ r-q}$.

For reasons of convenience, we imagine that the $r$ independent
infinitesimal transformations $X_1f, \dots, X_rf$ are chosen so that
$X_{ r-q+1}f, \dots, X_rf$ generate the said invariant subgroup. Then
our claim clearly amounts to the fact that, in arbitrary variables
$y_1, y_2, \dots$, there are $r$ infinitesimal transformations
$\mathcal{ Y}_1f, \dots, \mathcal{ Y}_r f$ which satisfy the following
two conditions: firstly, $\mathcal{ Y}_{ r-q+1}f, \dots, \mathcal{
Y}_rf$ vanish identically, while $\mathcal{ Y}_1f, \dots, \mathcal{
Y}_{ r-q}f$ are mutually independent, and secondly, the relations:
\def\theequation{10}\begin{equation}
\leftbracket
\mathcal{Y}_i,\,\mathcal{Y}_k
\rightbracket
=
c_{ik1}\,\mathcal{Y}_1f
+\cdots+
c_{ikr}\,\mathcal{Y}_rf
\ \ \ \ \ \ \ \ \ \ \ \ \ {\scriptstyle{(i,\,\,k\,=\,1\,\cdots\,r)}}
\end{equation}
must hold identically.

Since $X_{ r-q+1}f, \dots, X_rf$ generate a subgroup invariant in the
$G_r$, all $c_{ ik1}$, $c_{ ik2}$, \dots, $c_{ ik, \, r-q}$ in which
at least one of the two indices $i$ and $k$ is larger than $r - q$,
are equal to zero. Thus, if in 
the relations~\thetag{ 10}, we set equal to zero
all the infinitesimal transformations
$\mathcal{ Y}_{ r-q+1}f, \dots, \mathcal{ Y}_rf$, then
all relations for which not both $i$ and $k$ are 
smaller than $r - q+1$ will be identically satisfied, and we
keep only the following relations between
$\mathcal{ Y}_1f, \dots, \mathcal{ Y}_{ r-q}f$: 
\def\theequation{11}\begin{equation}
\leftbracket
\mathcal{Y}_i,\,\mathcal{Y}_k
\rightbracket
=
c_{ik1}\,\mathcal{Y}_1f
+\cdots+
c_{ik,\,r-q}\,\mathcal{Y}_{r-q}f
\ \ \ \ \ \ \ \ \ \ \ \ \ 
{\scriptstyle{(i,\,\,k\,=\,1\,\cdots\,r\,-\,q)}}.
\end{equation}

So, we need only to prove that there are $r - q$ independent
infinitesimal transformations $\mathcal{ Y}_1f, \dots, \mathcal{
Y}_{ r-q}f$ which are linked together by the relations~\thetag{
11}, or, what is the same: that there is an $(r - q)$-term 
group, the composition of which is represented by the system
of the constants $c_{ iks}$
${\scriptstyle{(i,\,\, k,\,\, s\, = \,1\, \cdots\, r\,-\,q)}}$.

In order to prove this, we start from the fact that, for
$i, k, j = 
1, \dots, r-q$, the Jacobi identity:
\[
\big\leftbracket
\leftbracket
X_i,\,X_k
\rightbracket,\,
X_j
\big\rightbracket
+
\big\leftbracket
\leftbracket
X_k,\,X_j
\rightbracket,\,
X_i
\big\rightbracket
+
\big\leftbracket
\leftbracket
X_j,\,X_i
\rightbracket,\,
X_k
\big\rightbracket
=
0,
\]
holds, and it can also be written as:
\[
\aligned
&\ \ \ 
\sum_{\mu=1}^{r-q}\,
\bigg\{
c_{ik\mu}
\leftbracket
X_\mu,\,X_j
\rightbracket
+
c_{kj\mu}
\leftbracket
X_\mu,\,X_i
\rightbracket
+
c_{ji\mu}
\leftbracket
X_\mu,\,X_k
\rightbracket
\bigg\}
\\
&
+
\sum_{\pi=1}^q\,
\Big\leftbracket
X_{r-q+\pi}f,\,\,\,\,
c_{ik,\,r-q+\pi}\,X_jf
+
c_{kj,\,r-q+\pi}\,X_if
+
c_{ji,\,r-q+\pi}\,X_kf
\Big\rightbracket
=
0.
\endaligned
\]

If we develop here the left-hand side and if we take into
consideration that the coefficients of $X_1f, \dots, X_{ r-q}f$
must vanish, we obtain, between the constants $c_{ iks}$
which appear in~\thetag{ 11}, the following relations:
\[
\aligned
\sum_{\mu=1}^{r-q}\,
\Big\{
&
c_{ik\mu}\,c_{\mu j\nu}
+
c_{kj\mu}\,c_{\mu i\nu}
+
c_{ji\mu}\,c_{\mu k\nu}
\Big\}
=
0
\\
&\ \ \ \ \ \ \ \ \ \ \
{\scriptstyle{(i,\,\,k,\,\,j,\,\,\nu\,=\,1\,\cdots\,r\,-\,q)}}.
\endaligned
\] 
These relations show that the $c_{ iks}$ ${\scriptstyle{(i, \,\,k,
\,\,s\, = \,1\, \cdots\, r\,- \,q)}}$ determine a composition in the
sense defined earlier on. As was observed on
p.~\pageref{Proposition-1-S-297}, there surely are $(r-q)$-term groups
the composition of which is determined by the $c_{ iks}$
${\scriptstyle{(i, \,\,k, \,\,s\, = \,1\, \cdots\, r\,- \,q)}}$.

\renewcommand{\thefootnote}{\fnsymbol{footnote}}
As a result, the claim stated above is proved. 
If we yet add what we already knew for a while, 
we obtain the\footnote[1]{\,
We will show later that the cited Proposition~1,
p.~\pageref{Proposition-1-S-297} is not at all indispensable 
\label{S-304}
for the
developments of the text. Cf. \name{Lie}, Archiv for Math. og Nat.
Vol. 10, p.~357, Christiania 1885.
} 
\renewcommand{\thefootnote}{\arabic{footnote}}

\def\theproposition{3}\begin{proposition}
If one knows all invariant subgroups of the $r$-term group 
$X_1f, \dots, X_rf$: 
\[
\leftbracket
X_i,\,X_k
\rightbracket
=
\sum_{s=1}^r\,c_{iks}\,X_sf,
\]
then one can indicate all compositions which a group isomorphic to the
group $X_1f, \dots, X_rf$ can have.
\end{proposition}

In order to really set up the concerned compositions, one
must proceed as follows:

If the $q >0$ independent infinitesimal transformations:
\[
g_{\mu 1}\,X_1f 
+\cdots+
g_{\mu r}\,X_rf
\ \ \ \ \ \ \ \ \ \ \ \ \ {\scriptstyle{(\mu\,=\,1\,\cdots\,q)}}
\]
generate a $q$-term invariant subgroup of the group: 
$X_1f, \dots, X_rf$, or shortly $G_r$, then one sets: 
\[
g_{\mu 1}\,\mathcal{Y}_1f
+\cdots+
g_{\mu r}\,\mathcal{Y}_rf
=
0
\ \ \ \ \ \ \ \ \ \ \ \ \ {\scriptstyle{(\mu\,=\,1\,\cdots\,q)}},
\]
one solves with respect to $q$ of the $r$ expressions
$\mathcal{ Y}_1f, \dots, \mathcal{ Y}_rf$, and
one eliminates them from the relations:
\[
\leftbracket
\mathcal{Y}_i,\,\mathcal{Y}_k
\rightbracket
=
\sum_{s=1}^r\,c_{iks}\,\mathcal{Y}_sf.
\]
Between the $r - q$ expressions remaining amongst the expressions
$\mathcal{ Y}_1f, \dots, \mathcal{ Y}_rf$, one then obtains relations
which define the composition of an $(r - q)$-term group isomorphic to
the $G_r$. If one proceeds in this way for every individual subgroup
of the $G_r$, one obtains all the desired compositions.

Since every group is one's own invariant subgroup, it follows that to
every $r$-term group $G_r$ is associated a meroedrically isomorphic
group, namely the group which is formed by the identity
transformation.  If the $G_r$ is simple (cf. Chap.~\ref{kapitel-15},
p.~\pageref{S-264}), then the identity transformation is evidently the
only group which is meroedrically isomorphic to it.

\sectionengellie{\S\,\,\,83.}

\label{S-305-sq}
In this paragraph, we consider an important case in which groups occur
that are isomorphic to a given group.

Let the $r$-term group:
\[
X_kf
=
\sum_{i=1}^n\,\xi_{ki}(x_1,\dots,x_n)\,
\frac{\partial f}{\partial x_i}
\ \ \ \ \ \ \ \ \ \ \ \ \ {\scriptstyle{(k\,=\,1\,\cdots\,r)}}
\]
of the space $x_1, \dots, x_n$ be imprimitive, and let:
\[
u_1(x_1,\dots,x_n)
=
{\rm const.},
\,\,\,\,\dots,\,\,\,\,
u_{n-q}(x_1,\dots,x_n)
=
{\rm const.}
\]
be a decomposition of the space into $\infty^{ n-q}$
$q$-times extended manifolds invariant by the group.

According to Chap.~\ref{kapitel-7} and
to Chap.~\ref{kapitel-13}, p.~\pageref{S-222}, 
the $n-q$ mutually 
independent functions $u_1, \dots, u_{ n-q}$ satisfy
relations of the form:
\[
X_k\,u_\nu
=
\omega_{k\nu}(u_1,\dots,u_{n-q})
\ \ \ \ \ \ \ \ \ \ \ \ \ 
{\scriptstyle{(k\,=\,1\,\cdots\,r\,;\,\,\,
\nu\,=\,1\,\cdots\,n\,-\,q)}},
\]
hence the $r$ expressions:
\[
\aligned
\sum_{\nu=1}^{n-q}\,
X_k\,u_\nu\,\frac{\partial f}{\partial u_\nu}
&
=
\sum_{\nu=1}^{n-q}\,
\omega_{k\nu}(u_1,\dots,u_{n-q})\,
\frac{\partial f}{\partial u_\nu}
=
\overline{X}_kf
\\
&
\ \ \ \ \ \ \ \ \ \ \ \ \ {\scriptstyle{(k\,=\,1\,\cdots\,r)}}
\endaligned
\]
represent as many infinitesimal transformations in the variables $u_1,
\dots, u_{n-q}$. We claim that $\overline{ X}_1f, \dots, \overline{
X}_rf$ generate a group isomorphic\footnote{\,
This notion of \terminology{reduced} group $\overline{ X}_1f, \dots,
\overline{ X}_rf$ appearing just now is central in the general
algorithm, devised by \name{Lie} and developed later in Vol. III,
towards the classification of all imprimitive transformation groups,
and it also unveils appropriately the true mathematical causality of
the nowadays seemingly unusual notion of isomorphism which was
introduced earlier on, p.~\pageref{Begriff-Isomorphismus}.
} 
to the group $X_1f, \dots, X_rf$.

For the proof, we form:
\[
\aligned
\overline{X}_i
\big(\overline{X}_k(f)\big)
-
\overline{X}_k
\big(\overline{X}_i(f)\big)
&
=
\sum_{\nu=1}^{n-q}\,
\big(
\overline{X}_i\,\omega_{k\nu}
-
\overline{X}_k\,\omega_{i\nu}
\big)\,
\frac{\partial f}{\partial u_\nu}
\\
&
=
\sum_{\nu=1}^{n-q}\,
\big(
X_i\,\omega_{k\nu}-X_k\,\omega_{i\nu}
\big)\,
\frac{\partial f}{\partial u_\nu}\,\,;
\endaligned
\]
now, we have:
\[
X_i\big(X_k(f)\big)
-
X_k\big(X_i(f)\big)
=
\sum_{s=1}^r\,c_{iks}\,X_sf,
\]
or, when we insert $u_\nu$ in place of $f$:
\[
X_i\,\omega_{k\nu}
-
X_k\,\omega_{i\nu}
=
\sum_{s=1}^r\,c_{iks}\,X_s\,u_\nu
=
\sum_{s=1}^r\,c_{iks}\,\omega_{s\nu},
\]
hence, after inserting the found values of
$X_i\, \omega_{ k\nu} - X_k \, \omega_{ i\nu}$, it comes:
\[
\overline{X}_i\big(\overline{X}_k(f)\big)
-
\overline{X}_k\big(\overline{X}_i(f)\big)
=
\sum_{s=1}^r\,\sum_{\nu=1}^{n-q}\,
c_{iks}\,\omega_{s\nu}\,
\frac{\partial f}{\partial u_\nu},
\]
or, what amounts to the same:
\[
\leftbracket
\overline{X}_i,\,\overline{X}_k
\rightbracket
=
\sum_{s=1}^r\,
c_{iks}\,\overline{X}_sf.
\]

But this is what was to be proved.

The group $\overline{ X}_1f, \dots, \overline{ X}_rf$ has a very
simple conceptual meaning.

From the definition of imprimitivity, it follows that the totality of
the $\infty^{ n-q}$ manifolds $u_1 = {\rm const.}$, \dots, $u_{ n-q} =
{\rm const.}$ remains invariant by the group $X_1f, \dots, X_rf$,
hence that the $\infty^{ n-q}$ manifolds are permuted by every
transformation of this group. Consequently, to every transformation
of the group $X_1f, \dots, X_rf$, there corresponds a certain
permutation of our $\infty^{ n-q}$ manifolds, or, what is the same, a
transformation in the $n - q$ variables $u_1, \dots, u_{ n-q}$. It is
clear that the totality of all the so obtained transformations in the
variables $u$ form a group, namely just the group which is generated
by $\overline{ X}_1f, \dots, \overline{ X}_rf$.

Naturally, the group $\overline{ X}_1f, \dots, \overline{ X}_rf$
needs not be holoedrically isomorphic to the initial group, 
but evidently, it is only meroedrically isomorphic to it
when, amongst the infinitesimal transformations $e_1\, \overline{ X}_1f
+\cdots+ e_r\, \overline{ X}_rf$, there is at least one
which vanishes identically without $e_1, \dots, e_r$ being
all zero, hence when at least one amongst the infinitesimal
transformations $e_1\, X_1f + \cdots + e_r\, X_rf$
leaves individually invariant each one
of our $\infty^{ n-q}$ manifolds.

Thus, we have the

\def\theproposition{4}\begin{proposition}
\label{Satz-4-S-307}
If the $r$-term group $X_1f, \dots, X_rf$ of the space $x_1, \dots, 
x_n$ is imprimitive and if the equations:
\[
u_1(x_1,\dots,x_n)
=
{\rm const.},
\,\,\,\,\dots,\,\,\,\,
u_{n-q}(x_1,\dots,x_n)
=
{\rm const.}
\]
represent a decomposition of the space in $\infty^{ n-q}$
$q$-times extended manifolds, then the infinitesimal
transformations: 
\[
\aligned
\sum_{\nu=1}^{n-q}\,
X_k\,u_\nu\,
\frac{\partial f}{\partial u_\nu}
&
=
\sum_{\nu=1}^{n-q}\,
\omega_{k\nu}(u_1,\dots,u_{n-q})\,
\frac{\partial f}{\partial u_\nu}
=
\overline{X}_kf
\\
&
\ \ \ \ \ \ \ \ \
{\scriptstyle{(k\,=\,1\,\cdots\,r)}}
\endaligned
\]
in the variables $u_1, \dots, u_{ n-q}$, generate a group isomorphic
to the group $X_1f, \dots, X_rf$ which indicates in which way the
$\infty^{ n-q}$ manifolds are permuted by the transformations of the
group $X_1f, \dots, X_rf$. If, amongst the infinitesimal
transformations $e_1\, X_1f + \cdots + e_r\, X_rf$, there are exactly
$r - \rho$ independent ones which leave 
\label{S-307}
individually invariant each
one of the $\infty^{ n-q}$ manifolds, then the group $\overline{
X}_1f, \dots, \overline{ X}_rf$ is just $\rho$-term.
\end{proposition}

Using the Theorem~54, p.~\pageref{Theorem-54-S-301}, we yet obtain
the

\def\theproposition{5}\begin{proposition}
If the group $X_1f, \dots, X_rf$ of the space
$x_1, \dots, x_n$ is imprimitive and if:
\[
u_1(x_1,\dots,x_n)
=
{\rm const.},
\,\,\,\,\dots,\,\,\,\,
u_{n-q}(x_1,\dots,x_n)
=
{\rm const.}
\]
is a decomposition of the space in $\infty^{ n-q}$ $q$-times extended
manifolds, then the totality of all infinitesimal transformations
$e_1\, X_1f + \cdots + e_r\, X_rf$ which leave individually invariant
each one of these manifolds generates an invariant subgroup of the
group $X_1f, \dots, X_rf$. 
\end{proposition}

At present, we consider an arbitrary $r$-term group 
$X_1f, \dots, X_rf$ which leaves invariant a manifold 
of the space $x_1, \dots, x_n$. 

According to Chap.~\ref{kapitel-14}, p.~\pageref{Theorem-40-S-233}, 
the points of this manifold are in turn transformed
by a group, the infinitesimal transformations of which
can be immediately indicated when the
equations of the manifold are in resolved form. 
Indeed, if the equations of the manifold read
in the following way:
\[
x_1
=
\varphi_1(x_{n-m+1},\dots,x_n),
\,\,\,\,\dots,\,\,\,\,
x_{n-m}
=
\varphi_{n-m}(x_{n-m+1},\dots,x_n),
\]
and if $x_{ n-m+1}, \dots, x_n$ are chosen as coordinates for
the points of the manifold, then the infinitesimal
transformations of the
group in question possess the
form
(cf. Chap.~\ref{kapitel-14}, p.~\pageref{S-234}):
\[
\aligned
\overline{X}_kf
=
\sum_{\mu=1}^m\,
\xi_{k,\,n-m+\mu}
(\varphi_1,
&
\dots,\varphi_{n-m},\,
x_{n-m+1},\dots,x_n)\,
\frac{\partial f}{\partial x_{n-m+\mu}}
\\
&\ \ \
{\scriptstyle{(k\,=\,1\,\cdots\,r)}}.
\endaligned
\]

What is more, as we have already proved
at that time, the following
relations hold:
\[
\leftbracket
\overline{X}_i,\,\overline{X}_k
\rightbracket
=
\sum_{s=1}^r\,c_{iks}\,\overline{X}_sf.
\]

From this, we see that the group $\overline{ X}_1f, 
\dots, \overline{ X}_rf$ is isomorphic to the group
$X_1f, \dots, X_rf$; in particular, it
is meroedrically isomorphic to it
when, amongst the infinitesimal transformations:
\[
e_1\,\overline{X}_1f
+\cdots+
e_r\,\overline{X}_rf,
\]
there is at least one transformation which vanishes identically, so
that at least one amongst the infinitesimal transformations $e_1\,
X_1f + \cdots + e_r\, X_rf$ leaves fixed every individual point of the
manifold.

We can therefore complement the Theorem~40 in Chap.~\ref{kapitel-14},
p.~\pageref{Theorem-40-S-233} as follows:

\def\theproposition{6}\begin{proposition}
If the $r$-term group $X_1f, \dots, X_rf$ in the variables $x_1, 
\dots, x_n$ leaves invariant the manifold:
\[
x_1
=
\varphi_1(x_{n-m+1},\dots,x_n),
\,\,\,\,\dots,\,\,\,\,
x_{n-m}
=
\varphi_{n-m}(x_{n-m+1},\dots,x_n),
\]
then the reduced infinitesimal transformations:
\[
\aligned
\overline{X}_kf
=
\sum_{\mu=1}^m\,
\xi_{k,\,n-m+\mu}
(\varphi_1,
&
\dots,\varphi_{n-m},\,
x_{n-m+1},\dots,x_n)\,
\frac{\partial f}{\partial x_{n-m+\mu}}
\\
&
\ \ \ \ \ \ 
{\scriptstyle{(k\,=\,1\,\cdots\,r)}}
\endaligned
\]
in the variables $x_{ n-m+1}, \dots, x_n$ generate a group isomorphic
to the $r$-term group which indicates in which way the points of the
manifold are permuted by the transformations of the group $X_1f,
\dots, X_rf$. If, amongst the infinitesimal transformations $e_1\,
X_1f + \cdots + e_r\, X_rf$, there are exactly $r - \rho$ independent
ones which leave invariant every individual point of the manifold,
then the group $\overline{ X}_1 f, \dots, \overline{ X}_rf$ is just
$\rho$-term.
\end{proposition}

In addition, the following also holds true.

\def\theproposition{7}\begin{proposition}
\label{Satz-7-S-309}
If the $r$-term group $X_1f, \dots, X_rf$ of the space $x_1, \dots,
x_n$ leaves invariant a manifold, then the totality of all
infinitesimal transformations $e_1 \, X_1f + \cdots + e_r\, X_rf$
which leave fixed every individual point of this manifold generate an
invariant subgroup of the $r$-term group.
\end{proposition}

\linestop

Let:
\[
X_kf
=
\sum_{i=1}^n\,\xi_{ki}(x_1,\dots,x_n)\,
\frac{\partial f}{\partial x_i}
\ \ \ \ \ \ \ \ \ \ \ \ \ {\scriptstyle{(k\,=\,1\,\cdots\,r)}}
\]
be an $r$-term \emphasis{intransitive} group of the $R_n$, and one
which contains only a \emphasis{discrete} number of invariant
subgroups.

The complete system which is determined by the equations:
\[
X_1f=0,
\,\,\,\dots,\,\,\,
X_rf=0
\]
is $q$-term under these assumptions, where $q < n$, and it possesses
$n - q$ independent solutions. Hence we can imagine that the variables
$x_1, \dots, x_n$ are chosen from the beginning in such a way that
$x_{ q+1}, \dots, x_n$ are solutions of the complete system; if we do
this, then because $\xi_{ k, \, q+j} = X_k\, x_{ q+j}$, all $\xi_{k,
q+1}, \dots, \xi_{kr}$ vanish.

Evidently, each one of the $\infty^{ n-q}$ $q$-times extended
manifolds:
\[
x_{q+1}=a_{q+1},
\,\,\,\dots,\,\,\,
x_n=a_n
\]
remains now invariant by the group $X_1f, \dots, X_rf$; but their
points are permuted and according to Chap.~\ref{kapitel-14},
p.~\pageref{Theorem-40-S-233}, by a group. If we choose $x_1, \dots,
x_q$ as coordinates for the points of the manifold, then the
infinitesimal transformations of the group in question will be:
\[
\aligned
\overline{X}_kf
=
\sum_{i=1}^q\,\xi_{ki}
&
(x_1,\dots,x_q,\,a_{q+1},\dots,a_n)\,
\frac{\partial f}{\partial x_i}
\\
&\ \ \ \ \ \
{\scriptstyle{(k\,=\,1\,\cdots\,r)}}
\endaligned
\]

Whether this group is $r$-term or not remains temporarily uncertain,
and in any case, it is isomorphic to the group $X_1f, \dots, X_rf$.

If the group $\overline{ X}_1f, \dots, \overline{ X}_rf$ if
$\rho$-term ($\rho \leqslant r$), then the group $X_1f, \dots, X_rf$
contains an $(r - \rho)$-term invariant subgroup which leaves fixed
every individual point of the manifold $x_{ q+1} = a_{ q+1}, \dots,
x_n = a_n$ (cf. Propositions~6 and~7). This invariant subgroup cannot
change with the values of $a_{ q+1}, \dots, a_n$, since otherwise,
there would be in the group $X_1f, \dots, X_rf$ a continuous series of
invariant subgroups, and this would contradict our assumption.
Consequently, there exists in the group $X_1f, \dots, X_rf$
an $(r - \rho)$-term invariant subgroup which leaves
fixed all points of an arbitrary manifold amongst the
$\infty^{ n-q}$ manifolds:
\[
x_{q+1}=a_{q+1},
\,\,\,\dots,\,\,\,
x_n=a_n,
\]
hence which actually leaves fixed all points of the space $x_1, \dots,
x_n$. But now, the identity transformation is the only one which
leaves at rest all points of the space $x_1, \dots, x_n$, so one has
$r - \rho = 0$ and $\rho = r$, that is to say: the group 
$\overline{ X}_1f, \dots, \overline{ X}_rf$ is holoedrically
isomorphic to the group $X_1f, \dots, X_rf$. 

As a result, the following holds true. 

\def\theproposition{8}\begin{proposition}
\label{Satz-8-S-310}
If $X_1f, \dots, X_rf$ are independent infinitesimal transformations
of an $r$-term group with the absolute invariants $\Omega_1 ( x_1,
\dots, x_n)$, \dots, $\Omega_{ n-q} ( x_1, \dots, x_n)$, and if there
is only, in this group, a discrete number of invariant subgroups, then
the points of an arbitrary invariant domain: $\Omega_1 = a_1$, \dots,
$\Omega_{ n-q} = a_{ n-q}$ are transformed by a group 
\terminology{holoedrically}
isomorphic to the group $X_1f, \dots, X_rf$.
\end{proposition}

\linestop


\chapter{Finite Groups, the Transformations of Which 
\\
Form Discrete Continuous Families}
\label{kapitel-18}
\chaptermark{Groups with Discrete Families of Transformations}

\setcounter{footnote}{0}

\abstract*{??}

\renewcommand{\thefootnote}{\fnsymbol{footnote}}
So far, we have only occupied ourselves with continuous transformation
groups, hence with groups which are represented by \emphasis{one}
system of equations of the form:
\[
x_i'
=
f_i(x_1,\dots,x_n,\,a_1,\dots,a_r)
\ \ \ \ \ \ \ \ \ \ \ \ \ {\scriptstyle{(i\,=\,1\,\cdots\,n)}}.
\]
In this chapter, we shall also treat briefly the finite groups which
cannot be represented by a single system of equations, but which can
only represented by several such systems; these are the groups about
which we already made a mention in the Introduction,
p.~\pageref{mention-discontinuous}\footnote[1]{\,
\name{Lie}, Verhandlungen der Gesellschaft der Wissenchaften zu 
Christiania, Nr. 12, p.~1, 1883.
} 
\renewcommand{\thefootnote}{\arabic{footnote}}

Thus, we imagine that a series of systems of equations of the form:
\def\theequation{1}\begin{equation}
\aligned
x_i'
=
f_i^{(k)}\big(x_1,\dots,\,
&
x_n,\,a_1^{(k)},\dots,a_{r_k}^{(k)}\big)
\ \ \ \ \ \ \ \ \ \ \ \ \ {\scriptstyle{(i\,=\,1\,\cdots\,n)}}
\\
&
{\scriptstyle{(k\,=\,1,\,2\,\cdots)}}
\endaligned
\end{equation}
is presented, in which each system contains a finite number $r_k$ of
arbitrary parameters $a_1^{ (k)}, \dots, a_{ r_k}^{(k)}$, and we
assume that the totality of all transformations which are represented
by these systems of equations forms a group.

Since each one of the systems of equations~\thetag{ 1} represents a
continuous family of transformations, our group consists of a discrete
number of continuous transformations. It is clear that every
continuous family of transformations of our group either coincides
with one of the families~\thetag{ 1}, or must be contained in one of
these families. Of course, we assume that none of the
families~\thetag{ 1} is contained in one of the remaining families.

It is our intention to develop the foundations of a general theory of
the sorts of groups just defined, but for reasons of simplicity we
will introduce a few restrictions, which, incidentally, are not to be
considered as essential.

Firstly, we make the assumption that the transformations of the
group~\thetag{ 1} are ordered as inverses by pairs.
Hence, although the transformations of the family:
\[
x_i'
=
f_i^{(k)}
\big(
x_1,\dots,x_n,\,a_1^{(k)},\dots,a_{r_k}^{(k)}\big)
\ \ \ \ \ \ \ \ \ \ \ \ \ {\scriptstyle{(i\,=\,1\,\cdots\,n)}}
\]
are not already ordered as inverses by pairs, 
the totality of the associated inverse transformations
is supposed to form a family which belongs to the
group, hence which is contained amongst the families~\thetag{ 1}.

Secondly, we assume that the number of the
families~\thetag{ 1} is finite, say equal to $m$. 
However, when the propositions we derive are also 
true for infinitely many families~\thetag{ 1}, 
we will occasionally point out that this is the case. 

\sectionengellie{\S\,\,\,84.}

To the two assumptions about the group~\thetag{ 1} which we have made
in the introduction of the chapter, we yet want to temporarily add the
third assumption that all families of the group should contain the
same number, say $r$, of essential parameters. In the next paragraph,
we show that this third assumption follows from the first two, and
hence is superfluous.

Let:
\[
x_i'
=
f_i^{(k)}
(x_1,\dots,x_n,\,a_1,\dots,a_r)
\ \ \ \ \ \ \ \ \ \ \ \ \ {\scriptstyle{(i\,=\,1\,\cdots\,n)}}
\]
be one of the $m$ families of $\infty^r$ transformations
of which our group consists, so $k$ is any of the numbers $1, 
\dots, m$. 

By resolution of the equations written above, we obtain a family
of transformations:
\[
x_i
=
F_i^{(k)}
(x_1',\dots,x_n',\,a_1,\dots,a_r)
\ \ \ \ \ \ \ \ \ \ \ \ \ {\scriptstyle{(i\,=\,1\,\cdots\,n)}}
\]
which, under the assumptions made, equally belongs to the group.

Therefore, when the two transformations:
\[
\aligned
x_i
&
=
F_i^{(k)}(x_1',\dots,x_n',\,a_1,\dots,a_r)
\\
x_i''
&
=
f_i^{(k)}
(x_1,\dots,x_n,\,a_1+h_1,\dots,a_r+h_r)
\endaligned
\ \ \ \ \ \ \ \ \ \ \ \ \ {\scriptstyle{(i\,=\,1\,\cdots\,n)}}
\]
are executed one after the other, it again
comes a transformation of our group, namely
the following one:
\def\theequation{2}\begin{equation}
x_i''
=
f_i^{(k)}
\big(
F_1^{(k)}(x',a),\dots,F_n^{(k)}(x',a),\,\,
a_1+h_1,\dots,a_r+h_r
\big)
\ \ \ \ \ \ \ \ \ \ \ \ \ {\scriptstyle{(i\,=\,1\,\cdots\,n)}}.
\end{equation}

Here, we expand the right-hand side with respect
to powers of $h_1, \dots, h_r$ and we find:
\[
x_i''
=
f_i^{(k)}
\big(
F^{(k)}(x',a),\,a
\big)
+
\sum_{j=1}^r\,h_j\,
\bigg[
\frac{\partial f_i^{(k)}(x,a)}{\partial a_j}
\bigg]_{x=F^{(k)}(x',a)}
+\cdots,
\]
where all the left out terms in $h_1, \dots, h_r$ are of
second order and of higher order. But if we take into account
that the two transformations:
\[
x_i
=
F_i^{(k)}(x',a),
\ \ \ \ \ \ \ \ \ \
x_i'
=
f_i^{(k)}(x,a)
\ \ \ \ \ \ \ \ \ \ \ \ \ {\scriptstyle{(i\,=\,1\,\cdots\,n)}}
\]
are inverse to each other, and if in addition, we
yet set for abbreviation:
\def\theequation{3}\begin{equation}
\bigg[
\frac{\partial f_i^{(k)}}{\partial a_j}_{
x=F^{(k)}(x',a)}
\bigg]
=
\eta_{ji}^{(k)}
(x_1',\dots,x_n',\,a_1,\dots,a_r),
\end{equation}
then we see that the transformation just found
has the shape: 
\def\theequation{4}\begin{equation}
x_i''
=
x_i'
+
\sum_{j=1}^r\,h_j\,\eta_{ji}^{(k)}(x',a)
+\cdots
\ \ \ \ \ \ \ \ \ \ \ \ \ {\scriptstyle{(i\,=\,1\,\cdots\,n)}}.
\end{equation}

It is easy to see that there are no
functions $\chi_1 (a), \dots, \chi_r (a)$ independent
of $x_1', \dots, x_n'$ which satisfy the $n$ equations:
\[
\sum_{j=1}^r\,\chi_j(a_1,\dots,a_r)\,
\eta_{ji}^{(k)}(x',a)
=
0
\ \ \ \ \ \ \ \ \ \ \ \ \ {\scriptstyle{(i\,=\,1\,\cdots\,n)}}
\]
identically, without vanishing all. Indeed, if one makes the
substitution $x_\nu' = f_\nu^{ (k)} (x, a)$
in these equations, one obtains the equations:
\[
\sum_{j=1}^r\,\chi_j(a)\,
\frac{\partial f_i^{(k)}(x,a)}{\partial a_j}
=
0
\ \ \ \ \ \ \ \ \ \ \ \ \ {\scriptstyle{(i\,=\,1\,\cdots\,n)}},
\]
which must as well be satisfied identically; 
but according to Chap.~\ref{essential-parameters}, 
Theorem~\ref{Theorem-essential}, p.~\pageref{Theorem-essential},
this is impossible, because the parameters
$a_1, \dots, a_r$ in the transformation
equations $x_i' = f_i^{ (k)} (x, a)$ are
essential.

From this, we conclude that the $r$ infinitesimal transformations:
\def\theequation{5}\begin{equation}
\sum_{i=1}^n\,
\eta_{ji}^{(k)}
(x_1',\dots,x_n',\,a_1,\dots,a_r)\,
\frac{\partial f}{\partial x_i'}
\ \ \ \ \ \ \ \ \ \ \ \ \ {\scriptstyle{(j\,=\,1\,\cdots\,r)}}
\end{equation}
are always mutually independent when $a_1, \dots, a_r$
is a system of values in general position.

By $a_1^0, \dots, a_r^0$, we want to understand a system of
values in general position and we want to set:
\[
\eta_{ji}^{(1)}(x',\,a^0)
=
\xi_{ji}(x_1',\dots,x_n'),
\]
by conferring the special value $1$ to the number $k$.
We will show that all infinitesimal transformations~\thetag{ 5}
can be linearly deduced from the $r$ mutually
independent infinitesimal transformations:
\[
X_j'f
=
\sum_{i=1}^n\,\xi_{ji}(x_1',\dots,x_n')\,
\frac{\partial f}{\partial x_i'}
\ \ \ \ \ \ \ \ \ \ \ \ \ {\scriptstyle{(j\,=\,1\,\cdots\,r)}},
\]
whichever value $k$ can have as one of the integers
$1$, $2$, \dots, $m$ and whichever
value $a_1, \dots, a_r$ can have. 

The proof here has great similarities with the developments
in Chap.~\ref{one-term-groups}, p.~\pageref{S-76} sq. 

We execute two transformations of our group one after
the other, namely at first the transformation:
\def\theequation{6}\begin{equation}
x_i''
=
x_i'
+
\sum_{j=1}^r\,h_j\,\xi_{ji}(x')
+\cdots
\ \ \ \ \ \ \ \ \ \ \ \ \ {\scriptstyle{(i\,=\,1\,\cdots\,n)}}
\end{equation}
which results from the transformation~\thetag{ 4} when one sets $k =
1$ and $a_1 = a_1^0$, \dots, $a_r = a_r^0$, and afterwards, the
transformation:
\[
x_i'''
=
x_i''
+
\sum_{j=1}^r\,\rho_j\,\eta_{ij}^{(k)}(x'',a)
+\cdots
\ \ \ \ \ \ \ \ \ \ \ \ \ {\scriptstyle{(i\,=\,1\,\cdots\,n)}}
\]
which, likewise, possesses the form~\thetag{ 4}. 
In this way, we obtain the following transformation
belonging to our group:
\[
\aligned
x_i'''
=
x_i'
+
\sum_{j=1}^r\,h_j\,
&
\xi_{ji}(x')
+
\sum_{j=1}^r\,\rho_j\,
\eta_{ji}^{(k)}(x',a)
+\cdots
\\
&
\ \ \ \ {\scriptstyle{(i\,=\,1\,\cdots\,n)}},
\endaligned
\]
where the left out terms are of second order and of
higher order in the $2r$ quantities $h_1, \dots, h_r$, 
$\rho_1, \dots, \rho_r$. 

Disregarding the $a$, the $2r$ arbitrary parameters $h_1, \dots, h_r$,
$\rho_1, \dots, \rho_r$ appear in the latter transformation, whereas
our group can contain only transformations with $r$ essential
parameters. From the Proposition~4 of the
Chap.~\ref{one-term-groups}, p.~\pageref{Satz-4-S-65}, it therefore
results that, amongst the $2r$ infinitesimal transformations: \thetag{
5} and $X_1'f, \dots, X_r'f$, only $r$ can be in existence that are
mutually independent. But because $X_1'f, \dots, X_r'f$ are mutually
independent, the infinitesimal transformations~\thetag{ 5} must be
linearly expressible in terms of $X_1'f, \dots, X_r'f$, for all values
$1$, $2$, \dots, $m$ of $k$ and for all values of the $a$.

At present, 
thanks to considerations analogous to those in 
Chap.~?? \Fill, we realize that identities of the form:
\[
\aligned
\sum_{i=1}^n\,\eta_{ji}^{(k)}(x',a)\,
&
\frac{\partial f}{\partial x_i'}
\equiv
\sum_{\pi=1}^r\,\psi_{j\pi}^{(k)}(a_1,\dots,a_r)\,
X_\pi'f
\\
&
{\scriptstyle{(k\,=\,1\,\cdots\,m\,;\,\,\,
j\,=\,1\,\cdots\,r)}}
\endaligned
\]
hold, where the $\psi_{ j\pi}^{ (k)}$ are completely determined
analytic functions of $a_1, \dots, a_r$.

Lastly, if we remember the equations~\thetag{ 3}
which can evidently also be written as:
\[
\eta_{ji}^{(k)}
\big(
f_1^{(k)}(x,a),\dots,f_n^{(k)}(x,a),\,
a_1,\dots,a_r
\big)
\equiv
\frac{\partial f_i^{(k)}(x,a)}{\partial a_j},
\]
and if we compare these equations to the identities:
\[
\eta_{ji}^{(k)}(x',a)
\equiv
\sum_{\pi=1}^r\,
\psi_{j\pi}^{(k)}(a)\,
\xi_{\pi i}(x'),
\]
we then obtain the identities:
\def\theequation{7}\begin{equation}
\aligned
\frac{\partial f_i^{(k)}(x,a)}{\partial a_j}
\equiv
&\,
\sum_{\pi=1}^r\,\psi_{j\pi}^{(k)}(a_1,\dots,a_r)\,
\xi_{\pi i}
\big(f_1^{(k)}(x,a),\,\dots,\,f_n^{(k)}(x,a)\big)
\\
&
\endaligned
\end{equation}

The functions $\xi_{ \pi i}$ here are independent
of the index $k$, but by contrast, the
$\psi_{ j\pi}^{(k)}$ are not.

Thus, we have the

\renewcommand{\thefootnote}{\fnsymbol{footnote}}
\def\thetheorem{55}\begin{theorem}
If the $m$ systems of equations:
\[
\aligned
x_i'
=
f_i^{(k)}(
&
x_1,\dots,x_n,\,a_1,\dots,a_r)
\ \ \ \ \ \ \ \ \ \ \ \ \ {\scriptstyle{(i\,=\,1\,\cdots\,n)}}
\\
&
\ \ \ \ \ \ \ \ \ \
{\scriptstyle{(k\,=\,1\,\cdots\,m)}},
\endaligned
\]
in each of which the $r$ parameters $a_1, \dots, a_r$ are essential,
represent all transformations of a group (and if at the same time, all
these transformations can be ordered as inverses by
pairs)\footnote[1]{\,
As one easily realizes, the Theorem~55 still remains correct when the
words in brackets are crossed out. (Compare to the developments of
the Chaps.~\ref{fundamental-differential}, \ref{one-term-groups}
and~\ref{kapitel-9}.)
}, 
then there are $r$ independent infinitesimal
transformations:
\[
X_jf
=
\sum_{i=1}^n\,\xi_{ji}(x_1,\dots,x_n)\,
\frac{\partial f}{\partial x_i}
\ \ \ \ \ \ \ \ \ \ \ \ \ {\scriptstyle{(j\,=\,1\,\cdots\,r)}}
\]
which stand in such a relationship to the group 
that each family:
\[
x_i'
=
f_i^{(k)}
(x_1,\dots,x_n,\,a_1,\dots,a_r)
\ \ \ \ \ \ \ \ \ \ \ \ \ {\scriptstyle{(i\,=\,1\,\cdots\,n)}}
\]
satisfies differential equations of the form:
\def\theequation{7'}\begin{equation}
\aligned
\frac{\partial f_i^{(k)}}{\partial a_j}
=
\sum_{\pi=1}^r\,
&
\psi_{j\pi}^{(k)}(a_1,\dots,a_r)\,
\xi_{\pi i}(x_1',\dots,x_n')
\\
&
{\scriptstyle{(i\,=\,1\,\cdots\,n\,;\,\,\,
j\,=\,1\,\cdots\,r)}}.
\endaligned
\end{equation}
\end{theorem}
\renewcommand{\thefootnote}{\arabic{footnote}}

From this, it follows at first that, according to the Theorems~21,
p.~\pageref{Theorem-21-S-149} and~24, p.~\pageref{Theorem-24-S-158},
the $r$ infinitesimal transformations $X_1f, \dots, X_rf$ generate an
$r$-term group.

Furthermore, it is clear that the Theorem~25 in Chap.~\ref{kapitel-9},
p.~\pageref{Theorem-25-S-160} finds an application to each one of the
families $x_i' = f_i^{ (k)} ( x, a)$: every transformation $x_i' =
f_i^{ (k)} (x,a)$ whose parameters $a_1, \dots, a_r$ lie in a certain
neighbourhood of $\overline{ a}_1, \dots, \overline{ a}_r$ can be
obtained by executing firstly the transformation $\overline{ x}_i =
f_i^{ (k)} (x, \, \overline{ a})$ and afterwards, a certain
transformation of the $r$-term group $X_1f, \dots, X_rf$.

\label{S-315-sq}
Now, since our group contains all transformations of the form~\thetag{
6}:
\[
x_i''
=
x_i'
+
\sum_{j=1}^r\,h_j\,\xi_{ji}(x')
+\cdots
\ \ \ \ \ \ \ \ \ \ \ \ \ {\scriptstyle{(i\,=\,1\,\cdots\,n)}}
\]
and since one of these transformations is the identity transformation,
namely the one with the parameters $h_1 = 0$, \dots, $h_r = 0$, then
the identity transformation must occur in one of the families $x_i' =
f_i^{ (k)} (x, a)$. Consequently, one amongst these families is just
the $r$-term group $X_1f, \dots, X_rf$. Of course, this is the only
$r$-term group generated by infinitesimal transformations which is
contained in the group $x_i' = f_i^{ (k)} ( x, a)$.

As a result, we have gained the following theorem:

\def\thetheorem{56}\begin{theorem}
\label{Theorem-56-S-315}
Every group:
\[
\aligned
x_i'
=
f_i^{(k)}(
&
x_1,\dots,x_n,\,a_1,\dots,a_r)
\ \ \ \ \ \ \ \ \ \ \ \ \ {\scriptstyle{(i\,=\,1\,\cdots\,n)}}
\\
&
\ \ \ \ \ \ \ 
{\scriptstyle{(k\,=\,1\,\cdots\,m)}}
\endaligned
\]
having the constitution indicated in the preceding theorem contains
one, and only one, $r$-term group with transformations inverse by
pairs. This $r$-term group is generated by the infinitesimal
transformations $X_1f, \dots, X_rf$ defined in the previous theorem,
its $\infty^r$ transformations form one of the $m$ families $x_i' =
f_i^{ (k)} (x, a)$ and in addition, they stand in the following
relationship with respect to each one of the remaining $m - 1$
families: If $\overline{ x}_i = f_i^{ (k)} (x, \overline{ a})$ is an
arbitrary transformation of the family $x_i' = f_i^{ (k)} (x, a)$,
then every other transformation of this family whose parameters $a_1,
\dots, a_r$ lie in a certain neighbourhood of $\overline{ a}_1, \dots,
\overline{ a}_r$ can be obtained by executing firstly the
transformation $\overline{ x}_i = f_i^{ (k)} (x, \, \overline{ a})$,
and afterwards, a certain transformation $x_i' = \omega_i ( \overline{
x}_1, \dots, \overline{ x}_n)$ of the $r$-term group $X_1f, \dots,
X_rf$.
\end{theorem}

For example, we consider the group consisting of the two families of
$\infty^1$ transformations that are represented by the two systems of
equations:
\[
\aligned
x'
&
=
x\,\cos a-y\,\sin a,
\ \ \ \ \ \ \ \ \ \ 
y'
=
x\,\sin a+y\,\cos a
\\
x'
&
=
x\,\cos a+y\,\sin a,
\ \ \ \ \ \ \ \ \ \ 
y'
=
x\,\sin a-y\,\cos a.
\endaligned
\]
For the two families, one obtains by differentiation with respect to
$a$:
\[
\frac{\D\,x'}{\D\,a}
=
-\,y',
\ \ \ \ \ \ \ \ \ \ 
\frac{\D\,y'}{\D\,a}
=
x',
\]
so that in the present case, the functions $\psi_{ j\pi}^{ (k)}$
mentioned in the Theorem~55 are independent of the index $k$.

The first one of the above two families is a one-term group which is
generated by the infinitesimal transformation $y\, \partial f /
\partial x - x\, \partial f / \partial y$. The general transformation
of the second family will be obtained when one executes firstly the
transformation:
\[
\overline{x}
=
x,
\ \ \ \ \ \ \ \ \ \ 
\overline{y}
=
-\,y,
\]
and afterwards the transformation:
\[
x'
=
\overline{x}\,\cos a
-
\overline{y}\,\sin a,
\ \ \ \ \ \ \ \ \ \ 
y'
=
\overline{x}\,\sin a
+
\overline{y}\,\cos a
\]
of the one-term group $y\, \partial f / \partial x - x\, \partial f /
\partial y$, hence the general transformation of the first
family.\,---

It should not remain unmentioned that the two Theorems~55 and~56 also
remain yet valid when the concerned group consists of an
\emphasis{infinite} number of discrete continuous families which all
contain the same number of essential parameters.\,---

\smallercharacters{

Before we go further, yet a few not unimportant remarks.

Let:
\[
x_i'
=
f_i(x_1,\dots,x_n,\,a_1,\dots,a_r)
\ \ \ \ \ \ \ \ \ \ \ \ \ {\scriptstyle{(i\,=\,1\,\cdots\,n)}}
\]
be an $r$-term continuous group which does not contain the identity
transformation. Then according to Theorem~25, Chap.~\ref{kapitel-9},
p.~\pageref{Theorem-25-S-160}, there is an $r$-term group $X_1f,
\dots, X_rf$ with the identity transformation and with pairwise
inverse transformations which stands in the following relationship to
the group $x_i' = f_i ( x, a)$, or shortly $\mathfrak{ G}$: If one
executes at first a transformation $\overline{ x}_i = f_i ( x, \,
\overline{ a})$ of $\mathfrak{ G}$ and afterwards a transformation:
\[
x_i'
=
\overline{x}_i
+
\sum_{k=1}^r\,
e_k\,\overline{X}_k\,\overline{x}_i
+\cdots
\ \ \ \ \ \ \ \ \ \ \ \ \ {\scriptstyle{(i\,=\,1\,\cdots\,n)}}
\]
of the group $X_1f, \dots, X_rf$, then one always obtains a
transformation of $\mathfrak{ G}$.

At present, one can easily prove that one then also always obtains a
transformation of $\mathfrak{ G}$ when one firstly executes a
transformation of the group $X_1f, \dots, X_rf$ and afterwards a
transformation of the group $\mathfrak{ G}$. We do not want to spend
time in order to produce this proof in details, and we only want to
remark that for this proof, one may employ considerations 
completely similar to those of Chap.~\ref{one-term-groups}
(Cf. Theorem~58). 

Now, if we take together these two relationships between $\mathfrak{
G}$ and the group $X_1f, \dots, X_rf$ and in addition, if we take into
account that we have to deal with two groups, then we realize
immediately that the transformations of the group $\mathfrak{ G}$ and
of the group $X_1f, \dots, X_rf$, when combined, again form a group,
but to be precise, a group the transformations of which cannot be
ordered as inverses by pairs.

At present, we also attract in the circle of our considerations yet
the family of the transformations:
\[
x_i'
=
F_i(x_1,\dots,x_n,\,a_1,\dots,a_r)
\ \ \ \ \ \ \ \ \ \ \ \ \ {\scriptstyle{(i\,=\,1\,\cdots\,n)}}
\]
which are inverses of the transformations $x_i' = f_i (x, a)$.
According to Theorem~2\footnote{\,
This theorem simply states that if transformation equations $x_i '= f
( x, a)$ are stable by composition, the same holds for the inverse
transformations $x_i = F_i ( x', a)$ modulo possible shrinkings of
domains.
} 
(\cite{enlie1888-18}, p.~19), this family of transformations also
forms a group, that may be called $\mathfrak{ G}'$. We shall show
that the transformations of the three groups $\mathfrak{ G}$,
$\mathfrak{ G}'$, $X_1f, \dots, X_rf$, when taken together, form again
a group, and naturally, a group with pairwise inverse transformations.

Let $T$ be the general symbol of a transformation of $\mathfrak{ G}$, 
whence $T^{ -1}$ is the general symbol of a transformation 
of $\mathfrak{ G}'$; by $S$, it will always be understood a
transformation of the group $X_1f, \dots, X_rf$. 

We already know that all transformations $T$ and $S$ taken
together form a group and that the transformations
$T^{ -1}$ taken for themselves do the same. From this,
we realize the existence of relations which have
the following form:
\[
\aligned
&
T_\alpha\,T_\beta
=
T_\gamma,
\ \ \ \ \ \ \
S_\lambda\,S_\mu
=
S_\nu,
\ \ \ \ \ \ \
T_\beta^{-1}\,T_\alpha^{-1}
=
T_\gamma^{-1}
\\
&
\ \ \ \ \ \ \ \ \ \ \ \ \ \ \
T_\alpha\,S_\lambda
=
T_\pi,
\ \ \ \ \ \ \
S_\lambda\,T_\alpha
=
T_\rho.
\endaligned
\]

We can also write as follows the second series of these relations:
\[
S_\lambda^{-1}\,T_\alpha^{-1}
=
T_\pi^{-1},
\ \ \ \ \ \ \
T_\alpha^{-1}\,S_\lambda^{-1}
=
T_\rho^{-1}.
\]
Now, since the group of the $S$ consists of pairwise inverse
transformations, it is immediately clear that the $T^{ -1}$ together
with the $S$ form a group. Moreover, we have:
\[
\aligned
T_\alpha^{-1}\,T_\pi
&
=
T_\alpha^{-1}\,T_\alpha\,S_\lambda
=
S_\lambda,
\\
T_\alpha\,T_\rho^{-1}
&
=
T_\alpha\,T_\alpha^{-1}\,S_\lambda^{-1}
=
S_\lambda^{-1},
\endaligned
\]
and therefore, the totality of all $S$, $T$, $T^{ -1}$ also forms a
group.

With these words, the promised proof is brought: the transformations
of the three groups $\mathfrak{ G}$, $\mathfrak{ G}'$, $X_1f, \dots,
X_rf$ form a group together, and naturally, a group with pairwise
inverse transformations.

}

\sectionengellie{\S\,\,\,85.}

At present, we take up a standpoint more general than in the previous
paragraph. We drop the special assumption\footnote{\,
(i.e. the third one)
} 
made there and we only maintain the two settlements which were done in
the Introduction.

Thus, we consider a group $G$ which consists of $m$ discrete families
with, respectively, $r_1$, $r_2$, \dots, $r_m$ essential parameters
and which, for each one of its transformations, also contains the
inverse transformation. We will prove that the numbers $r_1$, $r_2$,
\dots, $r_m$ are all mutually equal. Then from this, it follows that
the assumption: $r_1 = r_2 = \cdots = r_m$ made in the previous
paragraph was not a restriction.

We execute two transformations of the group one after the other, 
firstly a transformation:
\[
x_i'
=
f_i^{(k)}
(x_1,\dots,x_n,\,a_1,\dots,a_{r_k})
\ \ \ \ \ \ \ \ \ \ \ \ \ {\scriptstyle{(i\,=\,1\,\cdots\,n)}}
\]
of a family with $r_k$ parameters, and secondly a transformation:
\[
x_i''
=
f_i^{(j)}
(x_1',\dots,x_n',\,b_1,\dots,b_{r_j})
\ \ \ \ \ \ \ \ \ \ \ \ \ {\scriptstyle{(i\,=\,1\,\cdots\,n)}}
\]
of a family with $r_j$ parameters. 

In this way, we find a transformation:
\[
x_i''
=
f_i^{(j)}
\big(
f_1^{(k)}(x,a),\,\dots,\,f_n^{(k)}(x,a),\,
b_1,\dots,b_{r_j}
\big)
\]
which belongs to our group and which formally contains $r_k + r_j$
arbitrary parameters. Now, if the largest number amongst the numbers
$r_1$, $r_2$, \dots, $r_m$ has the value $r$, then amongst these $r_k
+ r_j$ arbitrary parameters, there are no more than $r$ which are
essential, but also no less essential ones than what indicates the
largest of the two numbers $r_k$ and $r_j$. So in particular, if both
numbers $r_k$ and $r_j$ are equal to $r$, then the last written
transformation contains exactly $r$ essential parameters.
Consequently, all families of our group which contain exactly $r$
essential parameters already form a group $\Gamma$ when taken
together.

Now, we can immediately apply the Theorem~56 of the previous paragraph
to the group $\Gamma$. From this, we see that $\Gamma$ and hence also
$G$ contains an $r$-term group: 
\[
x_i'
=
f_i(x_1,\dots,x_n,\,a_1,\dots,a_r)
\ \ \ \ \ \ \ \ \ \ \ \ \ {\scriptstyle{(i\,=\,1\,\cdots\,n)}}
\]
which is generated by $r$ independent infinitesimal transformations.
Thus, if we at first execute the transformation $x_i' = f_i ( x, a)$
and then the transformation:
\def\theequation{8}\begin{equation}
x_i''
=
f_i^{(j)}
(x_1',\dots,x_n',\,b_1,\dots,b_{r_j}),
\end{equation}
we obtain again a transformation of the group $G$, namely the
following one:
\def\theequation{9}\begin{equation}
x_i''
=
f_i^{(j)}
\big(
f_1(x,a),\,\dots,\,f_n(x,a),\,\,b_1,\dots,b_{r_j}
\big)
\ \ \ \ \ \ \ \ \ \ \ \ \ {\scriptstyle{(i\,=\,1\,\cdots\,n)}}.
\end{equation}

Of the $r + r_j$ parameters of this transformation, exactly $r$ are
essential, so the equations just written represent a continuous family
of $\infty^r$ transformations of the group $G$. But since the group
$x_i' = f_i ( x, a)$ contains the identity transformation, there are
special values of the parameters $a_1, \dots, a_r$ for which the
functions $f_1 ( x, a)$, \dots, $f_n ( x, a)$ reduce to $x_1, \dots,
x_n$, respectively. Consequently, the family of the $\infty^{ r_j}$
transformations~\thetag{ 8} is contained in the family of the
$\infty^r$ transformations~\thetag{ 9}. According to the remarks
made in the introduction of the chapter, this is possible
only when the two families coincide, hence when $r_j$ is
equal to $r$. 

As a result, it is proved that the numbers $r_1$, $r_2$, \dots, $r_m$
are all really equal one to another. Consequently, we have the

\def\thetheorem{57}\begin{theorem}
If a group whose transformations are pairwise inverse one to another
consists of $m$ continuous families of transformations and if each one
of these families contains only a finite number of arbitrary
parameters, then the families all have the same number of essential
parameters.
\end{theorem}

Besides, this theorem still remains valid also when the number of the
families of which the group consists is infinitely large, when each
one of these infinitely many families contains just a finite number
$\rho_k$ of arbitrary parameters and when at the same time, amongst
all the numbers $\rho_k$, a largest one is extant.

\sectionengellie{\S\,\,\,86.}

As up to now, let $G$ consist of $m$ discrete, continuous families of
$\infty^r$ transformations; in addition, we assume that the
transformations of $G$ are mutually inverse by pairs.

According to the Theorem~57 in the preceding paragraph, there is in
the group $G$ one and only one $r$-term group generated by $r$
independent infinitesimal transformations. Now, if $S$ is the symbol
of the general transformation of this $r$-term group and $T$ is the
symbol of an arbitrary transformation of $G$, then in the same way:
\[
T^{-1}\,S\,T
\]
is the symbol of the general infinitesimal transformation of an
$r$-term group generated by infinitesimal transformations. Since this
new group is contained in $G$, it must coincide with the group of all
$S$; according to the terminology \deutsch{Terminologie} introduced in
the Chap.~\ref{kapitel-15}, p.~\pageref{S-261}, we can also express
this as: the discussed $r$-term group remains invariant by every
transformation $T$. Thus the following holds true.

\def\thetheorem{58}\begin{theorem}
\label{Theorem-58-S-320}
If a group $G$ with pairwise inverse transformations consists
of several families of transformations, then the largest
group generated by infinitesimal transformations which is
contained in $G$ remains invariant by every transformation
of $G$. 
\end{theorem}

From this theorem, it comes how one can construct groups which consist
of several continuous families of $\infty^r$ transformations.

Let $X_1f, \dots, X_rf$ be an $r$-term group in the variables $x_1,
\dots, x_n$ and let again $S$ be the symbol of the general
transformation of this group.

Now, when a group with pairwise inverse transformations contains all
transformations of the $r$-term group $X_1f, \dots, X_rf$, and in
addition yet contains a finite number, say $m-1$, of discrete families
of $\infty^r$ transformations, then in consequence of Theorem~56,
p.~\pageref{Theorem-56-S-315}, it possesses the form:
\def\theequation{10}\begin{equation}
T_0\,S,\ \ 
T_1\,S,
\,\,\,\dots,\,\,\,
T_{m-1}S.
\end{equation}

Here, $T_0$ means the identity transformation and $T_1, \dots, T_{
m-1}$ are, according to the last theorem, constituted in such a way
that their totality leaves invariant all transformations $S$. This
property of the $T_\nu$ expresses analytically by the fact that for
every transformation $S_k$ of the group $X_1f, \dots, X_rf$, a
relation of the form:
\[
T_\nu^{-1}\,S_k\,T_\nu
=
S_j
\]
exists, where the transformation $S_j$ again belongs to the group
$X_1f, \dots, X_rf$. Incidentally, such a relation also holds true
when $\nu$ is equal to zero, for indeed one then has $S_k = S_j$.

Since $T_1, \dots, T_{ m-1}$ themselves belong to the
transformations~\thetag{ 10}, the totality of all
transformations~\thetag{ 10} can then be a group only when all
transformations $T_\mu \, T_\nu$ also belong to this totality. Hence,
aside from the above relations, the $T_i$ must also satisfy yet
relations of the form:
\[
T_\mu\,T_\nu
=
T_\pi\,S_\tau,
\]
where $\mu$ and $\nu$ denote arbitrary numbers amongst the numbers
$1$, $2$, \dots, $m-1$, while $\pi$ runs through the values $0$, $1$,
\dots, $m-1$. If one the numbers $\mu$, $\nu$, say $\mu$, is equal to
zero, then actually, there is already relation of the form indicated,
since indeed $T_\pi = T_\nu$ and $S_\tau$ is, just as $T_0$, the
identity transformation.

On the other hand, if the transformations $T_\mu$ possess the 
properties indicated, then for all values $0$, $1$, \dots, 
$m-1$ of the two numbers $\mu$ and $\nu$, there exist
relations of the form:
\[
T_\mu\,S_k\,T_\nu\,S_l
=
T_\mu\,T_\nu\,S_j\,S_l
=
T_\pi\,S_\tau\,S_j\,S_l
=
T_\pi\,S_\rho,
\]
and therefore, the totality of all transformations~\thetag{ 10} forms
a group. It is easy to see that in any case and in general, the
transformations of such a group order as inverses by pairs.

At present, it is yet to be asked how the transformations $T_\mu$ must
be constituted in order that the $m$ families~\thetag{ 10} are all
distinct from each other.

Obviously, all the families~\thetag{ 10} are distinct from each other
when no transformation of the group belongs simultaneously to two of
these families. On the other hand, if any two of these families, say:
$T_\mu \, S$ and $T_\nu \, S$, have a transformation in common, then
they are identical, because from the existence of a relation of the
form:
\[
T_\mu\,S_k
=
T_\nu\,S_l,
\]
it immediately follows:
\[
T_\nu
=
T_\mu\,S_k\,S_l^{-1},
\]
hence the family of the transformations $T_\nu \, S$ has the form:
\[
T_\mu\,S_k\,S_l^{-1}\,S,
\]
that is to say, it is identical to the family: $T_\mu \, S$.

Thus, for the $m$ families~\thetag{ 10} to be distinct from each
other, it is necessary and sufficient that no two of the
transformations $T_0$, $T_1$, \dots, $T_{ m-1}$ be linked by a
relation of the form:
\[
T_\nu
=
T_\mu\,S_j
\ \ \ \ \ \ \ \ \ \ \ \ \ 
{\scriptstyle{(\nu\,\neq\,\mu)}}.
\]

We summarize the gained result in the

\def\thetheorem{59}\begin{theorem}
If $S$ is the symbol of the general transformation of the $r$-term
group $X_1f, \dots, X_rf$, and if moreover, $T_1, \dots, T_{ m-1}$
are transformations which leave invariant the group $X_1f, \dots, X_rf$
and which, in addition, are linked together and jointly with 
the identity transformation $T_0$ by relations of the form:
\[
T_\mu\,T_\nu
=
T_\pi\,S,
\]
but not by relations of the form:
\[
T_\nu
=
T_\mu\,S_j
\ \ \ \ \ \ \ \ \ \ \ \ \ 
{\scriptstyle{(\nu\,\neq\,\mu)}},
\]
then the totality of all transformations:
\[
T_0\,S,\ \ 
T_1\,S,
\,\,\,\dots,\,\,\,
T_{m-1}\,S
\]
form a group with pairwise inverse transformations which consists of
$m$ discrete continuous families of $\infty^r$ transformations and
which, at the same time, contains all transformations of the group
$X_1f, \dots, X_rf$. If one chooses the transformations $T_1, \dots,
T_{ m-1}$ in all possible ways, then one obtains all groups having the
constitution indicated.
\end{theorem}

In the next chapter, we give a general method for the determination of
all transformations which leave invariant a given group $X_1f, \dots,
X_rf$.

\medskip

If one has two different systems of transformations $T_1, \dots, T_{
m-1}$, say $T_1, \dots, T_{ m-1}$ and $T_1', \dots, T_{ m-1}'$, then
obviously, the two groups:
\[
\aligned
&
T_0\,S,\ \
T_1\,S,
\,\,\,\dots,\,\,\,
T_{m-1}\,S
\\
&
T_0\,S,\ \ 
T_1',\ \
\,\,\,\dots,\,\,\,
T_{m-1}'\,S
\endaligned
\]
are always distinct when and only when it is not possible to represent
each one of the transformations $T_1', \dots, T_{ m-1}'$ in the form:
\[
T_\mu'
=
T_{i_\mu}\,S_{k_\mu}.
\]

Needless to say, one can frequently arrange that the $m$
transformations $T_0$, $T_1$, \dots, $T_{ m-1}$
already form a discontinuous group for themselves.

\smallskip{\sf Example.}
The $n$ infinitesimal transformations:
\[
\frac{\partial f}{\partial x_1},
\,\,\,\dots,\,\,\,
\frac{\partial f}{\partial x_n}
\]
generate an $r$-term group. The totality of all transformations
that leave this group invariant forms a finite continuous
group which is generated by the $n + n^2$ infinitesimal
transformations:
\[
\frac{\partial f}{\partial x_i},\ \ \
x_i\,\frac{\partial f}{\partial x_k}
\ \ \ \ \ \ \ \ \ \ \ \ \ {\scriptstyle{(i,\,\,k\,=\,1\,\cdots\,n)}}
\]
Now, if amongst the $\infty^{ nn}$ transformations:
\[
x_i'
=
a_{i1}\,x_1
+\cdots+
a_{in}\,x_n
\ \ \ \ \ \ \ \ \ \ \ \ \ {\scriptstyle{(i\,=\,1\,\cdots\,n)}}
\]
of the group:
\[
x_i\,\frac{\partial f}{\partial x_k}
\ \ \ \ \ \ \ \ \ \ \ \ \ {\scriptstyle{(i,\,\,k\,=\,1\,\cdots\,n)}},
\]
one chooses $m$ arbitrary transformations that form a discontinuous
group as transformations $T_0$, $T_1$, \dots, $T_{ m-1}$, and if one
sets for $S$ the general transformation:
\[
x_1'
=
x_1+a_1,
\,\,\,\dots,\,\,\,
x_n'
=
x_n
+
a_n
\]
of the group $\partial f / \partial x_1$, \dots, $\partial f /
\partial x_n$, then one always obtains a group which consists of $m$
discrete families and which comprises all $\infty^n$ transformations
of the group $\partial f / \partial x_1$, \dots, $\partial f /
\partial x_n$.

\medskip

Naturally, the Theorem~58, p.~\pageref{Theorem-58-S-320} also holds
true when the group $G$ consists of infinitely many families of
$\infty^r$ transformations. Hence if one wants to construct such a
group, one only has to seek infinitely many discrete transformations:
\[
T_1,\,\,
T_2,\,\,
\dots
\] 
which leave invariant the group $X_1f, \dots, X_rf$ and
which, in addition, satisfy pairwise relations of the form:
\[
T_\mu\,T_\nu
=
T_\pi\,S_\tau,
\]
but by contrast, which are neither mutually, nor together 
with the identity transformation, linked by relations of the form:
\[
T_\mu
=
T_\nu\,S_j.
\]
The totality of all transformations:
\[
T_0\,S,\ \
T_1\,S,\ \
T_2\,S,\,\,\,\dots
\]
then forms a group $G$ which comprises the group $X_1f, \dots, X_rf$
and which consists of infinitely many different families of
$\infty^r$ transformations. 

But the transformations of the group $G$ found in this way 
are in general not ordered as inverses by pairs; 
in order that they enjoy this property, 
each one of the transformations $T_1$, $T_2$, \dots, 
must also satisfy, aside from the relations indicated above, 
yet relations of the form:
\[
T_\mu^{-1}
=
T_{k_\mu}\,S_{j_\mu}.
\]

\sectionengellie{\S\,\,\,87.}

At present, we still make a few observations about the 
\emphasis{invariant} subgroups as we have considered them
in the previous three paragraphs.

Let $G$ be a group which consists of the $m$ discrete families:
\[
\aligned
x_i^{(k)}
=
f_i^{(k)}
\big(x_1,\dots,x_n,\,
&
a_1^{(k)},\dots,a_r^{(k)}\big)
\ \ \ \ \ \ \ \ \ \ \ \ \ {\scriptstyle{(i\,=\,1\,\cdots\,n)}}
\\
&
{\scriptstyle{(k\,=\,1\,\cdots\,m)}}
\endaligned
\]
of $\infty^r$ transformations and which, for each one of its
transformations, also contains the inverse transformation.
In particular, let:
\[
x_i'
=
f_i^{(1)}
(x_1,\dots,x_n,\,a_1,\dots,a_r)
\ \ \ \ \ \ \ \ \ \ \ \ \ {\scriptstyle{(i\,=\,1\,\cdots\,n)}}
\]
be the $r$-term group generated by $r$ infinitesimal 
transformations: $X_1f, \dots, X_rf$ which is contained
in $G$. 

In accordance with Chap.~\ref{kapitel-6}, p.~\pageref{S-95}, 
we say that every function which admits all transformations of $G$, 
hence which satisfies the $m$ equations:
\[
\mathfrak{U}\big(x_1^{(k)},\dots,x_n^{(k)}\big)
=
\mathfrak{U}(x_1,\dots,x_n)
\ \ \ \ \ \ \ \ \ \ \ \ \
{\scriptstyle{(k\,=\,1\,\cdots\,m)}},
\]
is an \terminology{invariant} of $G$. We want to show how one
can find the invariants of $G$. 

\medskip

Evidently, every invariant of $G$ is at the same time an invariant of
the $r$-term group $X_1f, \dots, X_rf$ and therefore, it is a solution
of the complete system which is determined by the equations:
\[
X_1f=0,
\,\,\,\dots,\,\,\,
X_rf=0.
\]

If this complete system is $n$-term, then the group $X_1f, \dots,
X_rf$ and therefore also the group $G$ possess in fact no
invariant. So, we assume that the said complete system is $(n-q)$-term
and we denote by $u_1, \dots, u_q$ any $q$ of its solutions that are
independent.

According to Theorem~58, the group $X_1f, \dots, X_rf$ remains
invariant by all transformations of $G$; consequently, the
$(n-q)$-term complete system determined by the equations:
\[
X_1f=0,
\,\,\,\dots,\,\,\,
X_rf=0
\]
also admits all transformations of $G$. Consequently (cf. Chap.~8,
p.~\pageref{S-138}), the solutions $u_1, \dots, u_q$ of this complete
system satisfy relations of the shape:
\def\theequation{11}\begin{equation}
\aligned
u_j
\big(x_1^{(k)},\dots,x_n^{(k)}\big)
&
=
\omega_j^{(k)}
\big(
u_1(x),\dots,u_q(x),\,
a_1^{(k)},\dots,a_r^{(k)}
\big)
\\ 
& \ \ \ \
{\scriptstyle{(j\,=\,1\,\cdots\,q\,;\,\,\,
k\,=\,1\,\cdots\,m)}}.
\endaligned
\end{equation}

It can be shown here that the functions $\omega_j^{ (k)}$ are all free
of the parameters $a_1^{ (k)}, \dots, a_r^{ (k)}$.

By $\overline{ a}_1^{ (k)}, \dots, \overline{ a}_r^{ (k)}$, we want to
denote an arbitrary fixed system of values. If the system of values
$a_1^{ (k)}, \dots, a_r^{ (k)}$ lies in a certain neighbourhood of
$\overline{ a}_1^{ (k)}, \dots, \overline{ a}_r^{ (k)}$, then
according to Theorem~56, p.~\pageref{Theorem-56-S-315}, the
transformation:
\[
x_i^{(k)}
=
f_i^{(k)}
\big(x_1,\dots,x_n,\,
a_1^{(k)},\dots,a_r^{(k)}\big)
\]
can be obtained by executing firstly the transformation:
\[
\overline{x}_i^{(k)}
=
f_i^{(k)}
\big(x_1,\dots,x_n,\,
\overline{a}_1^{(k)},\dots,\overline{a}_r^{(k)}\big)
\]
and afterwards, a certain transformation:
\[
x_i^{(k)}
=
f_i^{(1)}
\big(\overline{x}_1^{(k)},\dots,
\overline{x}_n^{(k)},\,
\alpha_1,\dots,\alpha_r\big)
\]
of the group $X_1f, \dots, X_rf$. So we have:
\[
x_i^{(k)}
=
f_i^{(1)}
\big(
f_1^{(k)}(x,\,\overline{a}^{(k)}),
\,\,\dots,\,\,
f_n^{(k)}(x,\,\overline{a}^{(k)}),\,\,
\alpha_1,\dots,\alpha_r
\big).
\]

Now, $u_1, \dots, u_q$ are invariants of the group $X_1f, 
\dots, X_rf$ and therefore, they satisfy relations of the shape:
\[
u_j\big(x_1^{(k)},\dots,x_n^{(k)}\big)
=
u_j
\big(\overline{x}_1^{(k)},\dots,
\overline{x}_n^{(k)}\big).
\]

On the other hand, we have:
\[
u_j\big(
\overline{x}_1^{(k)},\dots,\overline{x}_n^{(k)}\big)
=
\overline{\omega}_j^{(k)}
\big(u_1(x),\dots,u_q(x)\big),
\]
where the $\overline{ \omega}_j^{ (k)}$ depend only upon $u_1, \dots,
u_q$ and contain no arbitrary parameters, since $\overline{ a}_1^{
(k)}, \dots, \overline{ a}_r^{ (k)}$ are numerical constants, indeed.
Thus we obtain:
\def\theequation{11'}\begin{equation}
\aligned
u_j\big(x_1^{(k)},\dots,x_n^{(k)}\big)
=
\overline{\omega}_j^{(k)}
\big(u_1(x),\dots,u_q(x)\big),
\endaligned
\end{equation}
and with this, it is proved that the functions $\omega_j^{ (k)}$ in
the equations~\thetag{ 11} are effectively free of the parameters
$a_1^{ (k)}, \dots, a_r^{ (k)}$.

According to what has bee said above, every invariant $\mathfrak{ U}
(x_1, \dots, x_n)$ of the group $G$ satisfies $m$ equations of the
shape:
\[
\mathfrak{U}
\big(x_1^{(k)},\dots,x_n^{(k)}\big)
=
\mathfrak{U}(x_1,\dots,x_n)
\ \ \ \ \ \ \ \ \ \ \ \ \ {\scriptstyle{(k\,=\,1\,\cdots\,m)}}\,;
\]
but since it is in addition a function of only $u_1, \dots, u_q$, 
say:
\[
\mathfrak{U}(x_1,\dots,x_n)
=
J(u_1,\dots,u_q),
\]
then at the same time, it satisfies the $m$ relations:
\def\theequation{12}\begin{equation}
\aligned
J
\big(
\omega_1^{(k)}(u_1,\dots,u_q),
&
\,\,\dots,\,\,
\omega_q^{(k)}(u_1,\dots,u_q)
\big)
=
J(u_1,\dots,u_q)
\\
&
\ \ \ \
{\scriptstyle{(k\,=\,1\,\cdots\,m)}}.
\endaligned
\end{equation}

Conversely, every function $J ( u_1, \dots, u_q)$ that satisfies the
$m$ functional equations just written obviously is an invariant of the
group $G$. Consequently, in order to find all invariants of $G$, we
only need to fulfill these functional equations in the most general
way.

The problem of determining all solutions of the functional
equations~\thetag{ 12} is visibly identical to the problem of
determining all functions of $u_1, \dots, u_q$ which admit the $m$
transformations:
\def\theequation{13}\begin{equation}
\aligned
u_j'
&
=
\omega_j^{(k)}(u_1,\dots,u_q)
\ \ \ \ \ \ \ \ \ \ \ \ \ {\scriptstyle{(j\,=\,1\,\cdots\,q)}}
\\
&
\ \ \ \ \ \ \ \ \ \ \ \ \ 
{\scriptstyle{(k\,=\,1\,\cdots\,m)}}.
\endaligned
\end{equation}
But these $m$ transformations form a discontinuous group, as one
realizes without difficulty from the group property
\deutsch{Gruppeneigenschaft} of $G$. Thus, our problem stated at the
outset of the paragraph is lead back to a problem of the theory of the
discontinuous groups.

We summarize the gained result in the

\def\thetheorem{60}\begin{theorem}
If a group $G$ (whose transformations are pairwise inverse)
consists of several discrete families of $\infty^r$ transformations:
\[
\aligned
x_i^{(k)}
&
=
f_i^{(k)}
\big(x_1,\dots,x_n,\,\,a_1^{(k)},\,\dots,\,a_r^{(k)}\big)
\ \ \ \ \ \ \ \ \ \ \ \ \ {\scriptstyle{(i\,=\,1\,\cdots\,n)}}
\\
&
\ \ \ \ \ \ \ \ \ \ \ \ \ \ \ \ \ \ \ \ \ 
{\scriptstyle{(k\,=\,1,\,\,2\,\cdots)}},
\endaligned
\]
then all invariants of $G$ are at the same time also invariants of the
$r$-term continuous group: $X_1f, \dots, X_rf$ determined by $G$. If
one knows the invariants of the latter group, hence if one knows any
$q$ arbitrary independent solutions $u_1, \dots, u_q$ of the $(n -
q)$-term complete system which is determined by the equations:
\[
X_1f=0,
\,\,\,\dots,\,\,\,
X_rf=0,
\]
the one finds the invariants of $G$ in the following way: one forms at
first the relations:
\[
\aligned
u_j\big(x_1^{(k)},\dots,x_n^{(k)}\big)
&
=
\omega_j^{(k)}
\big(u_1(x),\dots,u_q(x)\big)
\ \ \ \ \ \ \ \ \ \ \ \ \ {\scriptstyle{(j\,=\,1\,\cdots\,q)}}
\\
&
\ \ \ \ \ \
{\scriptstyle{(k\,=\,1,\,\,2\,\cdots)}},
\endaligned
\]
which, under the assumptions made, exist, and in which the $\omega_j^{
(k)}$ depend only upon the two indices $j$ and $k$; afterwards, one
determines all functions of $u_1, \dots, u_q$ which admit the
discontinuous group formed by the transformations:
\[
\aligned
u_j'
&
=
\omega_j^{(k)}(u_1,\dots,u_q)
\ \ \ \ \ \ \ \ \ \ \ \ \ {\scriptstyle{(j\,=\,1\,\cdots\,q)}}
\\
&
\ \ \ \ \ \ \ \ \ \ \ \ \ \ \ \ \ \
{\scriptstyle{(k\,=\,1,\,\,2\,\cdots)}}.
\endaligned
\]
The concerned functions are the invariants of the group $G$.
\end{theorem}

A similar theorem clearly holds true when the group $G$ consists of
infinitely many continuous families of transformations.

\smallercharacters{

\renewcommand{\thefootnote}{\fnsymbol{footnote}}
One can propose to oneself the problem of finding
all invariants that a given family of transformations:
\[
x_i'
=
\varphi_i(x_1,\dots,x_n,\,a_1,\dots,a_r)
\ \ \ \ \ \ \ \ \ \ \ \ \ {\scriptstyle{(i\,=\,1\,\cdots\,n)}}
\]
possesses, or all invariants that are common to several such
families.\footnote[1]{\,
Cf. \name{Lie}, Berichte der K. Sächs. Ges. d. W., 1. August 1887.
} 
The family, respectively the families, can here be completely
arbitrary and they need not belong to a finite group.
\renewcommand{\thefootnote}{\arabic{footnote}}

We do not intend to treat exhaustively this problem; let it only be
remarked that the concerned invariants are solutions, though not
arbitrary solutions, of a certain complete system that can be easily
be indicated. Indeed, since the sought invariants, aside from the
given transformations, also obviously admit yet the associated inverse
transformations, then one can very easily set up certain infinitesimal
transformations by which they remain invariant as well. In general,
these infinitesimal transformations contain arbitrary elements, an in
particular, certain parameters; when set equal to zero, these
infinitesimal transformations provide linear partial differential
equations which must be satisfied by the sought invariants. Now, it
is always possible to set up the smallest complete system which
embraces all these differential equations. If one knows a system of
solutions $u_1$, $u_2$, \dots\, of this complete system, then one
forms an arbitrary function $\Omega ( u_1, u_2, \dots)$ of them, one
executes on it the general transformation of the given family and one
determines $\Omega$ in the most general way in order that $\Omega$
behaves invariantly.

}

\linestop


\chapter{Theory of the Similarity 
\deutsch{Aehnlichkeit} of $r$-term Groups}
\label{kapitel-19}
\chaptermark{Theory of the Similarity of $r$-term Groups}

\setcounter{footnote}{0}

\abstract*{??}

It is often of the utmost importance to answer the question whether a
given $r$-term group $x_i' = f_i ( x_1, \dots, x_s, \, a_1, \dots,
a_r)$ of the $s$-times extended space is \emphasis{similar}
\deutsch{ähnlich} to another given $r$-term group $y_i' = F_i ( y_1,
\dots, y_s, \, b_1, \dots, b_r)$ of the same space, hence whether one
can introduce, in place of the $x$ and of the $a$, new
\emphasis{variables}: $y_1, \dots, y_s$ and new \emphasis{parameters}:
$b_1, \dots, b_r$ so that the first group converts into the second
group (Chap.~\ref{fundamental-differential}, p.~\pageref{S-24}). If
one knows, in a given case, that such a transfer of one of the groups
to the other is possible, then a second question raises itself: how
one accomplishes the concerned transfer in the most general way?

In the present chapter, we provide means for answering the two
questions.

To begin with, we show that the first one of the two questions can be
replaced by the following more simple question: under which conditions
does there exist a transformation:
\[
y_i
=
\Phi(x_1,\dots,x_s)
\ \ \ \ \ \ \ \ \ \ \ \ \ {\scriptstyle{(i\,=\,1\,\cdots\,s)}}
\]
of such a nature that $r$ arbitrary independent infinitesimal
transformations of the group $x_i' = f_i ( x, a)$ are transferred to
infinitesimal transformations of the group $y_i' = F_i ( y, b)$ by the
introduction of the variables $y_1, \dots, y_s$? We settle this
simpler question by setting up certain conditions which are necessary
for the existence of a transformation $y_i = \Phi_i ( x)$ of the
demanded constitution, and which prove to also be sufficient. At the
same time, we shall see that all possibly existing transformations
$y_i = \Phi_i ( x)$ of the demanded constitution can be determined by
integrating complete systems. With that, the second one of the two
questions stated above will then also be answered.

\sectionengellie{\S\,\,\,88.}

Let the two $r$-term groups: $x_i' = f_i ( x, a)$ and $y_i' = F_i ( y,
b)$ be similar to each other and to be precise, let the first one be
transferred to the second one when the new variables $y_i = \Phi_i (
x_1, \dots, x_s)$ are introduced in place of $x_1, \dots, x_s$, and
when the new parameters $b_k = \beta_k ( a_1, \dots, a_r)$ are
introduced in place of $a_1, \dots, a_r$. Furthermore, 
let: 
\[
X_kf
=
\sum_{i=1}^s\,\xi_{ki}(x_1,\dots,x_s)\,
\frac{\partial f}{\partial x_i}
\ \ \ \ \ \ \ \ \ \ \ \ \ {\scriptstyle{(k\,=\,1\,\cdots\,r)}}
\]
be $r$ arbitrary independent infinitesimal transformations
of the group $x_i' = f_i ( x, a)$; by the introduction of the new
variables $y_1, \dots, y_s$, let them take the form:
\[
\aligned
X_kf
=
\sum_{i=1}^s\,X_k\,y_i\,
\frac{\partial f}{\partial y_i}
&
=
\sum_{i=1}^s\,\mathfrak{y}_{ki}(y_1,\dots,y_s)\,
\frac{\partial f}{\partial y_i}
=
\mathfrak{Y}_kf
\\
&\ \ \
{\scriptstyle{(k\,=\,1\,\cdots\,r)}}.
\endaligned
\]

We decompose the transition from the group $x_i' = f_i ( x, a)$ to the
group $y_i' = F_i ( y, b)$ in a series of steps.

At first, we bring the group $x_i' = f_i ( x, a)$ to the canonical
form:
\def\theequation{1}\begin{equation}
x_i'
=
x_i
+
\sum_{k=1}^r\,e_k\,\xi_{ki}(x_1,\dots,x_s)
+\cdots
\ \ \ \ \ \ \ \ \ \ \ \ \ {\scriptstyle{(i\,=\,1\,\cdots\,s)}}
\end{equation}
by introducing for $a_1, \dots, a_r$ certain (needless to say)
perfectly determined functions of $e_1, \dots, e_r$. Then evidently,
we must obtain from~\thetag{ 1} the equations $y_i' = F_i ( y, b)$
when we introduce the new variables $y_1, \dots, y_s$ in place of the
$x$, and when in addition, we insert for $e_1, \dots, e_r$ some
completely determined functions of $b_1, \dots, b_r$.

But now, according to Chap.~\ref{one-term-groups}, 
p.~\pageref{S-58}, after the introduction of the variables $y_1,
\dots, y_s$, the equations~\thetag{ 1} take the shape:
\def\theequation{1'}\begin{equation} y_i' = y_i +
\sum_{k=1}^r\,e_k\,\mathfrak{y}_{ki}(y_1,\dots,y_s) +\cdots \ \ \ \ \
\ \ \ \ \ \ \ \ {\scriptstyle{(i\,=\,1\,\cdots\,s)}},
\end{equation}
hence the equations~\thetag{ 1'} must become identical to the
equations $y_i' = F_i ( y,b)$ when one expresses $e_1, \dots, e_r$ in
terms of $b_1, \dots, b_r$ in the indicated way. Consequently, the
equations~\thetag{ 1'} are a form of the group $y_i' = F_i ( y,b)$,
and to be precise, a canonical form, as an examination teaches
\deutsch{wie der Augenschein lehrt}. In other words: $\mathfrak{
Y}_1f, \dots, \mathfrak{ Y}_rf$ are independent infinitesimal
transformations of this group.

Thus, if after the introduction of the new variables $y_1, \dots, y_s$
and of the new parameters $b_1, \dots, b_r$, the $r$-term group $x_i'
= f_i ( x_1, \dots, x_s, \, a_1, \dots, a_r)$ converts into the group
$y_i' = F_i ( y_1, \dots, y_s, \, b_1, \dots, b_r)$, then the
infinitesimal transformations of the first group convert into the
infinitesimal transformations of the second group, also after the
introduction of the variables $y_1, \dots, y_s$.

Obviously, the converse also holds true: if the two groups $x_i' = f_i
( x, a)$ and $y_i' = F_i ( y, b)$ stand in the mutual relationship
that the infinitesimal transformations of the first one are
transferred to the infinitesimal transformations of the second one
after the introduction of the new variables $y_1, \dots, y_s$, then
one can always transfer the group $x_i' = f_i ( x,a)$ to the group
$y_i' = F_i ( y, b)$ by means of appropriate choices of the variables
and of the parameters. Indeed, by the introduction of the $y$, the
canonical form~\thetag{ 1} of the group $x_i' = f_i( x, a)$ is
transferred to the canonical form~\thetag{ 1'} of the group $y_i' =
F_i ( y, b)$.

As a result, we have the

\def\thetheorem{61}\begin{theorem}
Two $r$-term groups in the same number of variables are similar to
each other if and only if it is possible to transfer any $r$ arbitrary
independent infinitesimal transformations of the first one to
infinitesimal transformations of the second one by the introduction of
new variables.
\end{theorem}

When the question is to examine whether the two $r$-term groups $x_i'
= f_i ( x, a)$ and $y_i = F_i ( y, b)$ are similar to each other, then
one needs only to consider \deutsch{ins Auge fassen} the infinitesimal
transformations of the two groups and to ask whether they can be
transferred to each other.

From the above developments, one can yet immediately draw another more
important conclusion.

We know that between the infinitesimal transformations
$X_1f, \dots, X_rf$ of the group $x_i' = f_i ( x,a)$, there exist
relations of the form:
\[
X_i\big(X_k(f)\big)
-
X_k\big(X_i(f)\big)
=
\leftbracket
X_i,\,X_k
\rightbracket
=
\sum_{\sigma=1}^r\,c_{ik\sigma}\,X_\sigma f.
\]
According to Chap.~\ref{kapitel-5}, Proposition~2,
p.~\pageref{Satz-2-S-84}, after the introduction of the new variables
$y_1, \dots, y_s$, these relations receive the form:
\[
\mathfrak{Y}_i
\big(\mathfrak{Y}_k(f)\big)
-
\mathfrak{Y}_k\big(\mathfrak{Y}_i(f)\big)
=
\leftbracket
\mathfrak{Y}_i,\,\mathfrak{Y}_k
\rightbracket
=
\sum_{\sigma=1}^r\,c_{ik\sigma}\,\mathfrak{Y}_\sigma f,
\] 
so the $r$ independent infinitesimal relations $\mathfrak{ Y}_1f,
\dots, \mathfrak{ Y}_rf$ of the group $y_i' = F_i ( y, b)$ are linked
together by exactly the same relations as the infinitesimal
transformations $X_1f, \dots, X_rf$ of the group $x_i' = f_i ( x, a)$.

According to the way of expressing introduced in
Chap.~\ref{kapitel-17}, p.~\pageref{S-291-bis} and~\pageref{S-293}, we
can therefore say:

\def\thetheorem{62}\begin{theorem}
If two $r$-term groups in the same number of variables are similar to
each other, then they are also equally composed, or, what is the same,
holoedrically isomorphic.
\end{theorem} 

At the same time, it is clear that the transformation $y_i = \Phi_i (
x_1, \dots, x_s)$ establishes a holoedrically isomorphic relationship
between the two groups: $x_i' = f_i ( x, a)$ and $y_i' = F_i ( y, b)$
since it associates to the $r$ independent infinitesimal
transformations $X_1f, \dots, X_rf$ of the first group the $r$
independent transformations $\mathfrak{ Y}_1f, \dots, \mathfrak{
Y}_rf$ of the other, and since through this correspondence, the two
groups are obviously related to each other in a holoedrically
\label{S-330}
isomorphic way.

\sectionengellie{\S\,\,\,89.} 

\label{S-330-bis}
At present, we imagine that two arbitrary $r$-term group
in the same number of variables are presented; let their 
infinitesimal transformations be:
\[
X_kf
=
\sum_{i=1}^s\,\xi_{ki}(x_1,\dots,x_s)\,
\frac{\partial f}{\partial x_i}
\ \ \ \ \ \ \ \ \ \ \ \ \ {\scriptstyle{(k\,=\,1\,\cdots\,r)}}
\]
and:
\[
Z_kf
=
\sum_{i=1}^s\,\zeta_{ki}(x_1,\dots,x_s)\,
\frac{\partial f}{\partial y_i}
\ \ \ \ \ \ \ \ \ \ \ \ \ {\scriptstyle{(k\,=\,1\,\cdots\,r)}}.
\]

We ask: are these two groups similar to each other, or not?

Our answer to this question must obviously fall in the negative sense
as soon as the two groups are not equally composed; indeed, according
to Theorem~62, only the equally composed groups can be similar to each
other. Consequently, we need only to occupy ourselves with the case when
the two groups are equally composed; the question whether this case
really happens may always be settled thanks to an algebraic
discussion \deutsch{algebraische Discussion}.

According to that, we assume from now on that the two presented groups
are \emphasis{equally composed}.

According to Theorem~61, the two equally composed groups $X_1f, \dots,
X_rf$ and $Z_1f, \dots, Z_rf$ are similar to each other if and only if
there is a transformation $y_i = \Phi_i ( x_1, \dots, x_s)$ of such a
constitution that $X_1f, \dots, X_rf$ convert into infinitesimal
transformations of the group $Z_1f, \dots, Z_rf$ after the
introduction of the new variables $y_1, \dots, y_s$. If we combine
this with the remark at the end of the previous paragraph, we realize
that the two groups are similar to each other if and only if a
holoedrically isomorphic relation can be established between them so
that it is possible, by the introduction of appropriate new variables:
$y_i = \Phi_i ( x_1, \dots, x_s)$ to transfer the $r$ infinitesimal
transformations $X_1f, \dots, X_rf$ precisely to the infinitesimal
transformations $\mathfrak{ Y}_1f, \dots, \mathfrak{ Y}_rf$ of the
group $Z_1f, \dots, Z_rf$ that are associated through the isomorphic
relation.

The next step that we must climb towards the answer to the question we
stated is therefore to relate the two groups in the most general
holoedrically isomorphic way.\label{S-331}

In the group $Z_1f, \dots, Z_rf$, we choose in the most general 
way $r$ independent infinitesimal transformations:
\[
Y_kf
=
\sum_{j=1}^r\,g_{kj}\,Z_jf
\ \ \ \ \ \ \ \ \ \ \ \ \ {\scriptstyle{(k\,=\,1\,\cdots\,r)}}
\]
such that together with the relations:
\def\theequation{2}\begin{equation}
\leftbracket
X_i,\,X_k
\rightbracket
=
\sum_{\sigma=1}^r\,c_{ik\sigma}\,X_\sigma f,
\end{equation}
there are at the same time the relations:
\def\theequation{2'}\begin{equation}
\leftbracket
Y_i,\,Y_k
\rightbracket
=
\sum_{\sigma=1}^r\,c_{ik\sigma}\,Y_\sigma f.
\end{equation}
If this is occurs\,---\,only algebraic operations are required for
that\,---, then we associate the infinitesimal transformations $Y_1f,
\dots, Y_rf$ to $X_1f, \dots, X_rf$, respectively, and we obtain that
the two groups are holoedrically isomorphically related to each other
in the most general way.

In the present chapter, by the $g_{kj}$, we understand everywhere
\emphasis{the most general system of constants which satisfies the
requirement stated just now}.

Here, it is to be remarked that in general, the $g_{ kj}$ depend upon
arbitrary elements, once upon arbitrary parameters and then
upon certain arbitrarinesses \deutsch{Willkürlichkeiten} which are
caused by the algebraic operations that are necessary for the
determination of the $g_{ kj}$; indeed, it is thinkable that there are
several discrete families of systems of values $g_{ kj}$ which possess
the demanded constitution.

At present, the question is whether, amongst the reciprocal isomorphic
relationships between the two groups found in this way, there is one
which has the property indicated above. In other words: is it possible
to specialize the arbitrary elements which occur in the coefficients
$g_{ kj}$ in such a way that $X_1f, \dots, X_rf$ can be transferred to
$Y_1f, \dots, Y_rf$, respectively, by means of the introduction of
appropriate new variables: $y_i = \Phi_i ( x_1, \dots, x_s)$?

When, but only when, this question must have been answered, one may
conclude that the two groups $X_1f, \dots, X_rf$ and $Z_1f, \dots,
Z_rf$ are similar to each other.

Let $X_1f, \dots, X_nf$ ($n \leqslant r$) be linked together
by no linear relation of the form:
\[
\chi_1(x_1,\dots,x_s)\,X_1f
+\cdots+
\chi_n(x_1,\dots,x_s)\,X_nf
=
0,
\]
while by contrast, $X_{n+1}f, \dots, X_rf$ can be
linearly expressed in terms of $X_1f, \dots, X_nf$:
\def\theequation{3}\begin{equation}
\label{S-332}
\aligned
X_{n+k}f
&
\equiv
\varphi_{k1}(x_1,\dots,x_s)\,X_1f
+\cdots+
\varphi_{kn}(x_1,\dots,x_s)\,X_nf
\\
&
\ \ \ \ \ \ \ \ \ \ \ \ \ \ \ \ \ \ \ \ \ \ \ \
{\scriptstyle{(k\,=\,1\,\cdots\,r\,-\,n)}}.
\endaligned
\end{equation}

Now, if the transformation $y_i = \Phi_i ( x_1, \dots, x_s)$ is
constituted in such a way that, after introduction of the $y$, $X_1f,
\dots, X_rf$ are transferred to certain infinitesimal transformations
$\mathfrak{ Y}_1f, \dots, \mathfrak{ Y}_rf$ which belong to the group
$Z_1f, \dots, Z_rf$, then naturally, $\mathfrak{ Y}_1f, \dots,
\mathfrak{ Y}_n f$ are not linked together by a linear relation of the
form:
\[
\psi_1(y_1,\dots,y_s)\,\mathfrak{Y}_1f
+\cdots+
\psi_n(y_1,\dots,y_s)\,\mathfrak{Y}_nf
=
0\,;
\]
by contrast, we visibly obtain for $\mathfrak{ Y}_{ n+1}f, \dots,
\mathfrak{ Y}_rf$ expressions of the shape:
\[
\mathfrak{Y}_{n+k}f
\equiv
\sum_{\nu=1}^n\,
\overline{\varphi}_{k\nu}(y_1,\dots,y_s)\,
\mathfrak{Y}_\nu f
\ \ \ \ \ \ \ \ \ \ \ \ \
{\scriptstyle{(k\,=\,1\,\cdots\,r\,-\,n)}},
\]
in which the $\overline{ \varphi}_{ k\nu} (y)$ 
come into existence after the introduction
of the variables $y$ in place of the $x$, so that 
the $n( r - n)$ equations:
\[
\overline{\varphi}_{k\nu}(y_1,\dots,y_s)
=
\varphi_{k\nu}(x_1,\dots,x_s)
\]
are hence identities after the substitution: 
$y_i = \Phi_i ( x_1, \dots, x_s)$. 

From this, we conclude that:

\plainstatement{
If there is no linear relation between $X_1f, \dots, X_nf$, while $X_{
n+1}f, \dots, X_rf$ express linearly in terms of $X_1f, \dots, X_nf$
by virtue of the relations~\thetag{ 3}, then the two equally composed
groups $X_1f, \dots, X_rf$ and $Z_1f, \dots, Z_rf$ can be similar only
when the arbitrary elements in the coefficients $g_{ kj}$ defined
above can be chosen in such a way that the infinitesimal
transformations:\label{S-333}
\[
Y_kf
=
\sum_{j=1}^r\,g_{kj}\,Z_jf
\ \ \ \ \ \ \ \ \ \ \ \ \ {\scriptstyle{(k\,=\,1\,\cdots\,r)}}
\]
possess the following properties: firstly, $Y_1f, \dots, Y_nf$ are
linked by no linear relation, while by contrast, $Y_{ n+1}f, \dots,
Y_rf$ express linearly in terms of $Y_1f, \dots, Y_nf$:
\def\theequation{3'}\begin{equation}
Y_{n+k}f
\equiv
\sum_{\nu=1}^n\,\psi_{k\nu}(y_1,\dots,y_s)\,Y_\nu f
\ \ \ \ \ \ \ \ \ \ \ \ \ 
{\scriptstyle{(k\,=\,1\,\cdots\,r\,-\,n)}},
\end{equation}
and secondly, the $n( r - n)$ equations:
\label{S-334-bis}
\def\theequation{4}\begin{equation}
\varphi_{k\nu}(x_1,\dots,x_s)
-
\psi_{k\nu}(y_1,\dots,y_s)
=
0
\ \ \ \ \ \ \ \ \ \ \ \ \ 
{\scriptstyle{(k\,=\,1\,\cdots\,r\,-\,n\,;\,\,\,
\nu\,=\,1\,\cdots\,n)}},
\end{equation}
are compatible with each other and they give
relations neither between the $x$ alone 
nor between the $y$ alone.}

These conditions are \emphasis{necessary} for the similarity
between the two equally composed groups $X_1f, \dots, X_rf$ 
and $Z_1f, \dots, Z_rf$. We claim that the same conditions
are \emphasis{sufficient}. 
\label{S-333-bis} Said more precisely, we claim:
\emphasis{when the said conditions are satisfied, 
then there always is a transformation $y_i = \Phi_i ( x_1, 
\dots, x_s)$ which transfers the infinitesimal
transformations $X_1f, \dots, X_rf$ to, respectively, 
$Y_1f, \dots, Y_rf$, so the two groups
$X_1f, \dots, X_rf$ and $Z_1f, \dots, Z_rf$ are 
similar to each other}. 

The proof of this claim will be produced 
while we develop a method which conducts to the
determination of a transformation having the
indicated constitution.

\medskip

Our more present standpoint is therefore the following:

In the $s$ variables $x_1, \dots, x_s$, let an 
$r$-term group:
\[
X_kf
=
\sum_{i=1}^s\,\xi_{ki}(x_1,\dots,x_s)\,
\frac{\partial f}{\partial x_i}
\ \ \ \ \ \ \ \ \ \ \ \ \ {\scriptstyle{(k\,=\,1\,\cdots\,r)}}
\]
be presented, the composition of which is determined by the relations:
\[
\leftbracket
X_i,\,X_k
\rightbracket
=
\sum_{\sigma=1}^r\,c_{ik\sigma}\,X_\sigma f.
\]
Between $X_1f, \dots, X_nf$, there is at the same time no linear
relation of the form:
\[
\chi_1(x_1,\dots,x_s)\,X_1f
+\cdots+
\chi_n(x_1,\dots,x_s)\,X_nf
=
0,
\]
while by contrast $X_{ n+1}f, \dots, X_rf$ express themselves in terms
of $X_1f, \dots, X_nf$:
\def\theequation{3}\begin{equation}
X_{n+k}f
\equiv
\sum_{\nu=1}^n\,\varphi_{k\nu}(x_1,\dots,x_s)\,X_\nu f
\ \ \ \ \ \ \ \ \ \ \ \ \ 
{\scriptstyle{(k\,=\,1\,\cdots\,n\,-\,r)}}.
\end{equation}

Furthermore, let us be given an $r$-term group:
\[
Z_kf
=
\sum_{i=1}^s\,\zeta_{ki}(y_1,\dots,y_s)\,
\frac{\partial f}{\partial y_i}
\ \ \ \ \ \ \ \ \ \ \ \ \ {\scriptstyle{(k\,=\,1\,\cdots\,r)}}
\]
which is equally composed with the group $X_1f, \dots, X_rf$, 
and let $r$ independent infinitesimal transformations:
\[
\label{S-334-quart}
Y_kf
=
\sum_{j=1}^r\,\overline{g}_{kj}\,Z_jf
=
\sum_{i=1}^s\,\eta_{ki}(y_1,\dots,y_s)\,
\frac{\partial f}{\partial y_i}
\ \ \ \ \ \ \ \ \ \ \ \ \ {\scriptstyle{(k\,=\,1\,\cdots\,r)}}
\]
in this group be chosen in such a way that firstly, the
relations:
\[
\leftbracket
Y_i,\,Y_k
\rightbracket
=
\sum_{\sigma=1}^r\,c_{ik\sigma}\,Y_\sigma f
\]
are identically satisfied, and such that secondly, 
$Y_1f, \dots, Y_nf$ are linked together by no relation 
of the form:
\[
\psi_1(y_1,\dots,y_s)\,Y_1f
+\cdots+
\psi_n(y_1,\dots,y_s)\,Y_nf
=
0,
\]
while by contrast $Y_{ n+1}f, \dots, Y_rf$ 
express themselves in the following way:
\def\theequation{3'}\begin{equation}
Y_{n+k}f
\equiv
\sum_{\nu=1}^n\,\psi_{k\nu}(y_1,\dots,y_s)\,Y_\nu f
\ \ \ \ \ \ \ \ \ \ \ \ \ 
{\scriptstyle{(k\,=\,1\,\cdots\,n\,-\,r)}},
\end{equation}
and such that, lastly, the $n( r - n)$ equations:
\def\theequation{4}\begin{equation}
\varphi_{k\nu}(x_1,\dots,x_s)
-
\psi_{k\nu}(y_1,\dots,y_s)
=
0
\ \ \ \ \ \ \ \ \ \ \ \ \ 
{\scriptstyle{(k\,=\,1\,\cdots\,r\,-\,n\,;\,\,\,
\nu\,=\,1\,\cdots\,n)}}
\end{equation}
are compatible with each other and give relations neither between the
$x$ alone, nor between the $y$ alone.

\plainstatement{To seek a transformation: \label{S-334-ter}
\def\theequation{5}\begin{equation}
y_i
=
\Phi_i(x_1,\dots,x_s)
\ \ \ \ \ \ \ \ \ \ \ \ \ {\scriptstyle{(i\,=\,1\,\cdots\,s)}}
\end{equation}
which transfers $X_1f, \dots, X_rf$ to $Y_1f, \dots, Y_rf$,
respectively.}

We can add:

\plainstatement{The sought transformation is constituted in such a way
that the equations~\thetag{ 4} are identities after the substitution:
$y_1 = \Phi_1 ( x)$, \dots, $y_s = \Phi_s ( x)$.}

With these words, the problem which is to be settled
is enunciated. \label{S-334}

\sectionengellie{\S\,\,\,90.}

Before we attack in its complete generality the problem stated at the
end of the preceding paragraph, we want to consider a special case, 
the settlement of which turns out to be substantially simpler;
we mean the case $n = r$ which, obviously, can occur only
when $s$ is at least equal to $r$. 

Thus, let the entire number $n$ defined above be equal to $r$.

It is clear that in this case, neither between $X_1f, \dots, X_rf$, 
nor between $Z_1f, \dots, Z_rf$ there exists a linear relation.
\label{S-335-sq}
From this, it follows that the $r$ infinitesimal transformations:
\[
Y_kf
=
\sum_{j=1}^r\,g_{kj}\,Z_jf
\ \ \ \ \ \ \ \ \ \ \ \ \ {\scriptstyle{(k\,=\,1\,\cdots\,r)}}
\]
satisfy by themselves the conditions indicated on p.~\pageref{S-333},
without it begin necessary to specialize 
further the arbitrary elements contained
in the $g_{ kj}$. Indeed, firstly there is no linear relation
between $Y_1f, \dots, Y_rf$, and secondly, the equations~\thetag{ 4}
reduce to the identity $0 = 0$, hence they are certainly
compatible to each other and they produce relations
neither between the $x$ alone, nor between the $y$ alone. 

Thus, if the claim stated on p.~\pageref{S-333-bis} is correct, 
our two $r$-term groups must be similar, and to be
precise, there must exist a transformation 
$y_i = \Phi_i ( x_1, \dots, x_s)$ which transfers
$X_1f, \dots, X_rf$ to $Y_1f, \dots, Y_rf$, respectively.
We seek to determine such a transformation.

When we introduce the variables $y_1, \dots, y_s$ in $X_kf$ by
means of the transformation:
\def\theequation{5}\begin{equation}
y_i
=
\Phi_i(x_1,\dots,x_s)
\ \ \ \ \ \ \ \ \ \ \ \ \ {\scriptstyle{(i\,=\,1\,\cdots\,s)}},
\end{equation}
we obtain:
\[
X_kf
=
\sum_{i=1}^s\,X_k\,y_i\,
\frac{\partial f}{\partial y_i}
=
\sum_{i=1}^s\,X_k\,\Phi_i\,
\frac{\partial f}{\partial y_i}\,;
\]
so by comparing to:
\[
Y_kf
=
\sum_{i=1}^s\,\eta_{ki}(y_1,\dots,y_s)\,
\frac{\partial f}{\partial y_i}
=
\sum_{i=1}^s\,Y_k\,y_i\,
\frac{\partial f}{\partial y_i},
\]
we realize that the transformation~\thetag{ 5}
transfers $X_1f, \dots, X_rf$ to $Y_1f, \dots, Y_rf$, 
respectively, if and only if the $r\, s$ equations:
\def\theequation{6}\begin{equation}
Y_k\,y_i
-
X_k\,\Phi_i
=
0
\ \ \ \ \ \ \ \ \ \ \ \ \ 
{\scriptstyle{(k\,=\,1\,\cdots\,r\,;\,\,\,
i\,=\,1\,\cdots\,s)}}
\end{equation}
become identities after the substitution: 
$y_1 = \Phi_1 ( x)$, \dots, $y_s = \Phi_s ( x)$. 

Now, if we set: 
\[
X_kf
+
Y_kf
=
\Omega_kf
\ \ \ \ \ \ \ \ \ \ \ \ \ {\scriptstyle{(k\,=\,1\,\cdots\,r)}},
\]
we have:
\[
\Omega_k(y_i-\Phi_i)
=
Y_k\,y_i
-
X_k\,\Phi_i\,;
\]
thus, when the equations~\thetag{ 5} represent a transformation having
the constitution demanded, the $r\, s$ expressions $\Omega_k ( y_i -
\Phi_i)$ all vanish by means of~\thetag{ 5}, or, what amounts to the
same: the system of equations~\thetag{ 5} admits the $r$ infinitesimal
transformations $\Omega_1f, \dots, \Omega_rf$
(cf. Chap.~\ref{kapitel-7}, p.~\pageref{S-109-sq} sq.).

On the other hand, the following obviously holds true: every system of
equations solvable with respect to $x_1, \dots, x_s$ of the
form~\thetag{ 5} which admits the $r$ infinitesimal transformations
$\Omega_1f, \dots, \Omega_rf$ represents a transformation which
transfers $X_1f, \dots, X_rf$ to $Y_1f, \dots, Y_rf$, respectively.

At present, if we take into account that the $r$ infinitesimal
transformations $\Omega_1f, \dots, \Omega_rf$ stand pairwise in the
relationships:
\[
\leftbracket
\Omega_i,\,\Omega_k
\rightbracket
=
\leftbracket
X_i,\,X_k
\rightbracket
+
\leftbracket
Y_i,\,Y_k
\rightbracket
=
\sum_{\sigma=1}^r\,c_{ik\sigma}\,\Omega_\sigma f,
\]
and therefore generate an $r$-term group in the $2s$ variables $x_1,
\dots, x_s$, $y_1, \dots, y_s$, we can thus say:

\plainstatement{The totality of all transformations~\thetag{ 5} which
transfer $X_1f, \dots, X_rf$ to $Y_1f, \dots, Y_rf$, respectively, is
identical to the totality of all systems of equations solvable with
respect to $x_1, \dots, x_n$ of the form~\thetag{ 5} which admit the
$r$-term group $\Omega_1f, \dots, \Omega_rf$.}

The determination of a transformation of the discussed constitution is
therefore lead back to the determination of a certain system of
equations in the $2s$ variables $x_1, \dots, x_s$, $y_1, \dots, y_s$;
this system of equations must possess the following properties:
\emphasis{it must consist of $s$ independent equations, it must be
solvable both with respect to $x_1, \dots, x_s$ and with respect to
$y_1, \dots, y_s$, and lastly, it must admit the $r$-term group
$\Omega_1f, \dots, \Omega_rf$}.

For the resolution of this new problem, we can base ourselves on the
developments of the Chap.~\ref{kapitel-14}.

Under the assumptions made, the $r \times r$ determinants of the
matrix:
\[
\left\vert
\begin{array}{cccc}
\xi_{11}(x) & \,\cdot\, & \,\cdot\, & \xi_{1s}(x)
\\
\cdot & \,\cdot\, & \,\cdot\, & \cdot
\\
\xi_{r1}(x) & \,\cdot\, & \,\cdot\, & \xi_{rs}(x)
\end{array}
\right\vert
\]
do not all vanish identically and even less all $r \times
r$ determinants of the matrix:
\def\theequation{7}\begin{equation}
\left\vert
\begin{array}{cccccccc}
\xi_{11}(x) & \,\cdot\, & \,\cdot\, & \xi_{1s}(x) &
\eta_{11}(y) & \,\cdot\, & \,\cdot\, & \eta_{1s}(y)
\\
\cdot & \,\cdot\, & \,\cdot\, & \cdot &
\cdot & \,\cdot\, & \,\cdot\, & \cdot
\\
\xi_{r1}(x) & \,\cdot\, & \,\cdot\, & \xi_{rs}(x) &
\eta_{r1}(y) & \,\cdot\, & \,\cdot\, & \eta_{rs}(y)
\end{array}
\right\vert.
\end{equation}
Hence when a system of equations brings to zero all $r \times r$
determinants of the matrix~\thetag{ 7}, it must necessarily contain
relations between the $x$ alone.

Our problem is the determination of a system of equations of the form:
\def\theequation{8}\begin{equation}
y_1-\Phi_1(x_1,\dots,x_s)=0,
\,\,\,\dots,\,\,\,
y_s-\Phi_s(x_1,\dots,x_s)=0
\end{equation}
which admits the group $X_k f + Y_kf$, and which is at the same time
solvable with respect to $x_1, \dots, x_s$.

A system of equations of this nature certainly contains no relation
between $x_1, \dots, x_s$ only, so it does not bring to zero all $r
\times r$ determinants of the matrix~\thetag{ 7} and according to
Theorem~17 in Chap.~\ref{kapitel-7}, p.~\pageref{Theorem-17-S-123}, it
can be brought to a such a form that it contains only relations
between the solutions of the complete system:
\def\theequation{9}\begin{equation}
\Omega_kf
=
X_kf+Y_kf
=
0
\ \ \ \ \ \ \ \ \ \ \ \ \ {\scriptstyle{(k\,=\,1\,\cdots\,r)}}.
\end{equation}
Consequently, our problem will be solved when we will succeed to
determine $s$ independent relations between the solutions of this
complete system that are solvable with respect to $x_1, \dots, x_s$
and with respect to $y_1, \dots, y_s$ as well.

The complete system~\thetag{ 9} possesses $2s - r$ independent
solutions which we can evidently choose in such a way that $s - r$ of
them, say:
\[
u_1(x_1,\dots,x_s),
\,\,\dots,\,\,
u_{s-r}(x_1,\dots,x_s)
\]
depend only upon the $x$, so that $s - r$ other:
\[
v_1(y_1,\dots,y_s),
\,\,\dots,\,\,
v_{s-r}(y_1,\dots,y_s)
\]
depend only upon the $y$, while by contrast, the $r$ remaining ones:
\[
w_1(x_1,\dots,x_s,\,y_1,\dots,y_s),
\,\,\,\dots,\,\,\,
w_r(x_1,\dots,x_s,\,y_1,\dots,y_s)
\]
must contain certain $x$ and certain $y$ as well. Here, the $s$
functions $u_1, \dots, u_{ s-r}$, $w_1, \dots, w_r$ are mutually
independent as far as $x_1, \dots, x_s$ are concerned, and the $s$
functions $v_1, \dots, v_{ s-r}$, $w_1, \dots, w_r$ are so too, as far
as $y_1, \dots, y_s$ are concerned; from this, it follows that the $r$
equations of the complete system are solvable with respect to $r$ of
the differential quotients $\partial f / \partial y_1$, \dots,
$\partial f / \partial y_s$, and with respect to $r$ of the
differential quotients $\partial f / \partial x_1$, \dots, $\partial f
/ \partial x_s$ as well (cf. Chap.~\ref{kapitel-5}, Theorem~12,
p.~\pageref{Theorem-12-S-91}).

Now, when $s$ mutually independent relations between the $u$, $v$, $w$
are solvable both with respect to $x_1, \dots, x_s$ and with respect
to $y_1, \dots, y_s$? Clearly, when and only when they can be solved
both with respect to $u_1, \dots, u_{ s - r}$, $w_1, \dots, w_r$ and
with respect to $v_1, \dots, v_{ s-r}$, $w_1, \dots, w_r$, hence when
they can be brought to the form:
\def\theequation{10}\begin{equation}
\aligned
v_1
&
=
\mathfrak{F}_1(u_1,\dots,u_{s-r}),
\,\,\,\dots,\,\,\,
v_{s-r}
=
\mathfrak{F}_{s-r}(u_1,\dots,u_{s-r}),
\\
w_1
&
=
\mathfrak{G}_1(u_1,\dots,u_{s-r}),
\,\,\,\dots,\,\,\,\,\,\,\,
w_r
=
\mathfrak{G}_r(u_1,\dots,u_{s-r}),
\endaligned
\end{equation}
where $\mathfrak{ F}_1, \dots, \mathfrak{ F}_{ s-r}$ denote arbitrary
\emphasis{mutually independent} functions of their arguments, while
the functions $\mathfrak{ G}_1, \dots, \mathfrak{ G}_r$ are submitted
to absolutely no restriction.

The equations~\thetag{ 10} represent the most general system of
equations which consists of $s$ independent equations, which admits
the group $\Omega_1f, \dots, \Omega_rf$ and which can be solved both
with respect to $x_1, \dots, x_s$ and with respect to $y_1, \dots,
y_s$; at the same time, they represent the most general transformation
between $x_1, \dots, x_s$ and $y_1, \dots, y_s$ which transfers $X_1f,
\dots, X_rf$ to $Y_1f, \dots, Y_rf$, respectively. It is therefore
proved that there exist transformations which achieve this transfer,
hence that our two groups are effectively similar to each other.

But even more: the equations~\thetag{ 10} actually represent the most
general transformation which transfers the group $X_1f, \dots, X_rf$
to the group $Z_1f, \dots, Z_rf$.

In fact, if $y_i = \Psi_i ( x_1, \dots, x_s)$ is an arbitrary
transformation which transfers the one group to the other, then it
converts $X_1f, \dots, X_rf$ into certain infinitesimal
transformations $\mathfrak{ Y}_1f, \dots, \mathfrak{ Y}_rf$ of the
group $Z_1f, \dots, Z_rf$ which stand pairwise in the relationships:
\[
\leftbracket
\mathfrak{Y}_i,\,\mathfrak{Y}_k
\rightbracket
=
\sum_{\sigma=1}^r\,c_{ik\sigma}\,
\mathfrak{Y}_\sigma f,
\]
and which can hence be obtained from the $r$ infinitesimal
transformations:
\[
Y_kf
=
\sum_{j=1}^r\,g_{kj}\,Z_jf
\ \ \ \ \ \ \ \ \ \ \ \ \ {\scriptstyle{(k\,=\,1\,\cdots\,r)}}
\]
when one specializes in an appropriate way the arbitrary elements
appearing in the $g_{ kj}$. But now, all transformations which
convert $X_1f, \dots, X_rf$ into $Y_1f, \dots, Y_rf$, respectively,
are contained in the form~\thetag{ 10}, so in particular, the
transformation $y_i = \Psi_i ( x_1, \dots, x_n)$ is also contained in
this form.

Now, by summarizing the gained result, we can say:

\def\thetheorem{63}\begin{theorem}
If the two $r$-term groups:
\[
X_kf
=
\sum_{i=1}^s\,\xi_{ki}(x_1,\dots,x_s)\,
\frac{\partial f}{\partial x_i}
\ \ \ \ \ \ \ \ \ \ \ \ \ {\scriptstyle{(k\,=\,1\,\cdots\,r)}}
\]
and:
\[
Z_kf
=
\sum_{i=1}^s\,\zeta_{ki}(y_1,\dots,y_s)\,
\frac{\partial f}{\partial y_i}
\ \ \ \ \ \ \ \ \ \ \ \ \ {\scriptstyle{(k\,=\,1\,\cdots\,r)}}
\]
are equally composed and if neither $X_1f, \dots, X_rf$, nor $Z_1f,
\dots, Z_rf$ are linked together by linear relations, then the two
groups are also similar to each other. On obtains the most general
transformation which transfers the one group to the other in the
following way: One chooses the $r^2$ constants $g_{ kj}$ in the most
general way so that the $r$ infinitesimal transformations:
\[
Y_kf
=
\sum_{j=1}^r\,g_{kj}\,Z_jf
\ \ \ \ \ \ \ \ \ \ \ \ \ {\scriptstyle{(k\,=\,1\,\cdots\,r)}}
\]
are mutually independent and so that, together with the relations:
\[
\leftbracket
X_i,\,X_k
\rightbracket
=
\sum_{\sigma=1}^r\,c_{ik\sigma}\,X_\sigma f,
\]
there are at the same time the relations:
\[
\leftbracket
Y_i,\,Y_k
\rightbracket
=
\sum_{\sigma=1}^r\,c_{ik\sigma}\,Y_\sigma f\,;
\]
afterwards, one forms the $r$-term complete system:
\[
X_kf+Y_kf
=
0
\ \ \ \ \ \ \ \ \ \ \ \ \ {\scriptstyle{(k\,=\,1\,\cdots\,r)}}
\]
in the $2s$ variables $x_1, \dots, x_s$, $y_1, \dots, y_s$
and one determines $2s -r$ independent solutions
of it, namely $s - r$ independent solutions:
\[
u_1(x_1,\dots,x_s),
\,\,\dots,\,\,
u_{s-r}(x_1,\dots,x_s)
\]
which contain only the $x$, plus $s - r$ independent solutions:
\[
v_1(y_1,\dots,y_s),
\,\,\dots,\,\,
v_{s-r}(y_1,\dots,y_s)
\]
which contain only the $y$, and $r$ solutions:
\[
w_1(x_1,\dots,x_s,\,y_1,\dots,y_s),
\,\,\,\dots,\,\,\,
w_r(x_1,\dots,x_s,\,y_1,\dots,y_s)
\]
which are mutually independent and which are independent of $u_1,
\dots, u_{ s-r}$, $v_1, \dots, v_{ s-r}$; if this takes place, then
the system of equations:
\[
\aligned
v_1
&
=
\mathfrak{F}_1(u_1,\dots,u_{s-r}),
\,\,\,\dots,\,\,\,
v_{s-r}
=
\mathfrak{F}_{s-r}(u_1,\dots,u_{s-r}),
\\
w_1
&
=
\mathfrak{G}_1(u_1,\dots,u_{s-r}),
\,\,\,\dots,\,\,\,\,\,\,\,
w_r
=
\mathfrak{G}_r(u_1,\dots,u_{s-r}),
\endaligned
\]
represents the demanded transformation; here, $\mathfrak{ G}_1, \dots,
\mathfrak{ G}_r$ are perfectly arbitrary functions of their arguments;
by contrast, $\mathfrak{ F}_1, \dots, \mathfrak{ F}_{ s-r}$ are
subjected to the restriction that they must be mutually 
independent.
\end{theorem}

From this, it results in particular the

\def\theproposition{1}\begin{proposition}
If the $r \leqslant s$ independent infinitesimal transformations:
\[
X_kf
=
\sum_{i=1}^s\,\xi_{ki}(x_1,\dots,x_s)\,
\frac{\partial f}{\partial x_i}
\ \ \ \ \ \ \ \ \ \ \ \ \ {\scriptstyle{(k\,=\,1\,\cdots\,r)}}
\]
stand pairwise in the relationships:
\[
\leftbracket
X_i,\,X_k
\rightbracket
=
0
\ \ \ \ \ \ \ \ \ \ \ \ \ {\scriptstyle{(i,\,\,k\,=\,1\,\cdots\,r)}},
\]
without being, however, linked together by a linear relation
of the form:
\[
\sum_{k=1}^r\,\chi_k(x_1,\dots,x_s)\,X_kf
=
0, 
\]
then they generate an $r$-term group which is
similar to the group of transformations:
\[
Y_1f
=
\frac{\partial f}{\partial y_1},
\,\,\,\dots,\,\,\,
Y_rf
=
\frac{\partial f}{\partial y_r}.
\]
\end{proposition}

A case of the utmost importance is when the two numbers $s$ and $r$
are equal to each other, so that the two groups $X_1f, \dots, X_rf$
and $Z_1f, \dots, Z_rf$ are transitive, or more precisely: simply
transitive (cf. Chap.~\ref{kapitel-13}, p.~\pageref{S-212}).
We want to enunciate the Theorem~63 for this special case:

\def\thetheorem{64}\begin{theorem}
\label{Theorem-64-S-340}
Two simply transitive equally composed groups in the same number of
variables are always also similar to each other.
If:
\[
X_kf
=
\sum_{i=1}^r\,\xi_{ki}(x_1,\dots,x_r)\,
\frac{\partial f}{\partial x_i}
\ \ \ \ \ \ \ \ \ \ \ \ \ {\scriptstyle{(k\,=\,1\,\cdots\,r)}}
\]
and:
\[
Z_kf
=
\sum_{i=1}^r\,\zeta_{ki}(y_1,\dots,y_r)\,
\frac{\partial f}{\partial y_i}
\ \ \ \ \ \ \ \ \ \ \ \ \ {\scriptstyle{(k\,=\,1\,\cdots\,r)}}
\]
are the infinitesimal transformations of the two groups, then one
finds in the following way the most general transformation which
transfers the one group to the other: One chooses the $r^2$ constants
$g_{ kj}$ in the most general way so that the $r$ infinitesimal
transformations:
\[
Y_kf
=
\sum_{j=1}^r\,g_{kj}\,Z_jf
\ \ \ \ \ \ \ \ \ \ \ \ \ {\scriptstyle{(k\,=\,1\,\cdots\,r)}}
\]
are mutually independent and so that, together with the relations:
\[
\leftbracket
X_i,\,X_k
\rightbracket
=
\sum_{\sigma=1}^r\,c_{ik\sigma}\,X_\sigma f,
\]
there are at the same time the relations:
\[
\leftbracket
Y_i,\,Y_k
\rightbracket
=
\sum_{\sigma=1}^r\,
c_{ik\sigma}\,Y_\sigma f\,;
\]
moreover, one forms the $r$-term complete system:
\[
X_kf
+
Y_kf
=
0
\ \ \ \ \ \ \ \ \ \ \ \ \ {\scriptstyle{(k\,=\,1\,\cdots\,r)}}
\]
in the $2r$ variables $x_1, \dots, x_r$, $y_1, \dots, y_r$ 
and one determines $r$ arbitrary independent solutions:
\[
w_1(x_1,\dots,x_r,\,y_1,\dots,y_r),
\,\,\,\dots,\,\,\,
w_r(x_1,\dots,x_r,\,y_1,\dots,y_r)
\]
of it; then the $r$ equations:
\[
w_1
=
a_1,
\,\,\,\dots,\,\,\,
w_r
=
a_r
\]
with the $r$ arbitrary constants $a_1, \dots, a_r$ represent the most
general transformation having the constitution demanded.
\end{theorem}

\sectionengellie{\S\,\,\,91.}

At present, we turn to the treatment of the general problem that we
have stated at the end of \S\,\,89 (p.~\pageref{S-334}).

At first, we can prove exactly as in the previous paragraph that every
transformation: $y_i = \Phi_i ( x_1, \dots, x_s)$ which transfers
$X_1f, \dots, X_rf$ to $Y_1f, \dots, Y_rf$, respectively, represents a
system of equations which admits the $r$-term group $\Omega_k f = X_kf
+ Y_kf$ and that, on the other hand, every system of equations
solvable with respect to $x_1, \dots, x_s$:
\[
y_i
=
\Phi_i(x_1,\dots,x_s)
\]
which admits the group $\Omega_1f, \dots, \Omega_r f$ represents a
transformation which transfers $X_1f, \dots, X_rf$ to $Y_1f, \dots,
Y_rf$, respectively.

Now, according to p.~\pageref{S-334-bis}, every transformation
$y_i = \Phi_i ( x_1, \dots, x_s)$ which transfers
$X_1f, \dots, X_rf$ to $Y_1f, \dots, Y_rf$, respectively,
is constituted in such a way that the equations~\thetag{ 4}
become identities after the substitution: 
$y_1 = \Phi_1 ( x)$, \dots, $y_s = \Phi_s ( x)$. Consequently, 
we can also enunciate as follows the
problem formulated on p.~\pageref{S-334-ter}:

\shortplainstatement{\label{S-341}
To seek, in the $2s$ variables $x_1, \dots, x_s$, $y_1, \dots, y_s$, a
system of equations which admits the $r$-term group $\Omega_1f, \dots,
\Omega_rf$, which consists of exactly $s$ independent equations that
are solvable both with respect to $x_1, \dots, x_s$ and to $y_1,
\dots, y_s$, and lastly, which embraces\footnote{\,
See the footnote on p.~\pageref{explain-embrace}.
} 
the $n ( r - n)$ equations:
\def\theequation{4}\begin{equation}
\varphi_{k\nu}(x_1,\dots,x_s)
-
\psi_{k\nu}(x_1,\dots,x_s)
=
0
\ \ \ \ \ \ \ \ \ \ \ \ \ 
{\scriptstyle{(k\,=\,1\,\cdots\,r\,-\,n\,;\,\,\,
\nu\,=\,1\,\cdots\,n)}}.
\end{equation}} 

For the solution of this problem, it is if great importance
that \emphasis{the system of equations~\thetag{ 4}
admits in turn the $r$-term group $\Omega_1f, \dots, \Omega_r f$}.

In order to prove this, we imagine that the matrix
which is associated to the infinitesimal transformations
$\Omega_1f, \dots, \Omega_rf$:
\def\theequation{11}\begin{equation}
\left\vert
\begin{array}{cccccccc}
\xi_{11}(x) & \,\cdot\, & \,\cdot\, & \xi_{1s}(x) &
\eta_{11}(y) & \,\cdot\, & \,\cdot\, & \eta_{1s}(y)
\\
\cdot & \,\cdot\, & \,\cdot\, & \cdot &
\cdot & \,\cdot\, & \,\cdot\, & \cdot
\\
\xi_{r1}(x) & \,\cdot\, & \,\cdot\, & \xi_{rs}(x) &
\eta_{r1}(y) & \,\cdot\, & \,\cdot\, & \eta_{rs}(y)
\end{array}
\right\vert
\end{equation}
is written down. We will show that the equations~\thetag{ 4} are a
system of equations which is obtained by setting equal to zero all
$(n+1) \times (n+1)$ determinants of the matrix~\thetag{ 11}. With
this, according to Theorem~39, Chap.~\ref{kapitel-14},
p.~\pageref{Theorem-39-S-228}, it will be proved that the system of
equations~\thetag{ 4} admits the group $\Omega_1f, \dots, \Omega_rf$.

Amongst the $(n+1) \times (n+1)$ determinants of the matrix~\thetag{
11}, there are in particular those of the form:
\[
\Delta
=
\left\vert
\begin{array}{ccccc}
\xi_{1k_1}(x) & \,\cdot\, & \,\cdot\, &
\xi_{1k_n}(x) & \eta_{1\sigma}(y)
\\
\cdot & \,\cdot\, & \,\cdot\, & \cdot & \cdot
\\
\xi_{nk_1}(x) & \,\cdot\, & \,\cdot\, &
\xi_{nk_n}(x) & \eta_{n\sigma}(y)
\\
\xi_{n+j,\,k_1}(x) & \,\cdot\, & \,\cdot\, &
\xi_{n+j,\,k_n}(x) & \eta_{n+j,\,\sigma}(y)
\end{array}
\right\vert.
\]
If we replace in $\Delta$ the members of the last horizontal row by
their values from~\thetag{ 3} and~\thetag{ 3'}, namely by the
following values:
\[
\xi_{n+j,\,k_\mu}(x)
\equiv
\sum_{\nu=1}^n\,\varphi_{j\nu}(x)\,
\xi_{\nu,\,k_\mu}(x),
\ \ \ \ \ \ \
\eta_{n+j,\,\sigma}(y)
\equiv
\sum_{\nu=1}^n\,
\psi_{j\nu}(y)\,\eta_{\nu\sigma}(y),
\]
and if we subtract from the last horizontal row the first $n$ rows,
after we have multiplied them before by $\varphi_{ j1} (x), \dots,
\varphi_{ jn} (x)$, respectively, then we receive:
\[
\Delta
=
\sum\,\pm\,\xi_{1k_1}(x)\cdots\,\xi_{nk_n}(x)\,
\sum_{\nu=1}^n\,\eta_{\nu\sigma}(y)
\big\{
\psi_{j\nu}(y)
-
\varphi_{j\nu}(x)
\big\}.
\]

Here, under the assumptions made earlier on, the determinants
of the form:
\[
D
=
\sum\,\pm\,\xi_{1k_1}(x)\cdots\,\xi_{nk_n}(x)
\]
do not all vanish identically, and likewise not all determinants:
\[
\mathfrak{D}
=
\sum\,\pm\,\eta_{1k_1}(y)\cdots\,\eta_{nk_n}(y)
\] 
vanish identically.

Obviously, a system of equations that brings to zero all determinants
$\Delta$ must either contain all equations of the form $D = 0$, or it
must contain the $s ( r - n)$ equations:
\[
\sum_{\nu=1}^n\,
\eta_{\nu\sigma}(y)\,
\big\{
\psi_{j\nu}(y)
-
\varphi_{j\nu}(x)
\big\}
=
0
\ \ \ \ \ \ \ \ \ \ \ \ \ 
{\scriptstyle{(\sigma\,=\,1\,\cdots\,s\,;\,\,\,
j\,=\,1\,\cdots\,r\,-\,n)}}\,;
\]
in the latter case, it embraces either all equations of the form
$\mathfrak{ D} = 0$, or the $n ( r - n)$ equations:
\def\theequation{4}\begin{equation}
\varphi_{j\nu}(x)
-
\psi_{j\nu}(y)
=
0
\ \ \ \ \ \ \ \ \ \ \ \ \ 
{\scriptstyle{(j\,=\,1\,\cdots\,r\,-\,n\,;\,\,\,
\nu\,=\,1\,\cdots\,n)}}.
\end{equation}

In the latter one of these three cases, we now
observe what follows:

The system of equation~\thetag{ 4} brings to zero not only all
determinants $\Delta$, but actually also all $(n+1) \times (n+1)$
determinants of the matrix~\thetag{ 11}. One realizes this immediately
when one writes the matrix under the form:
\[
\left\vert
\begin{array}{cccccccc}
\xi_{11}(x) & \,\cdot\, & \,\cdot\, & \xi_{1s}(x) &
\eta_{11}(y) & \,\cdot\, & \,\cdot\, & \eta_{1s}(y)
\\
\cdot & \,\cdot\, & \,\cdot\, & \cdot &
\cdot & \,\cdot\, & \,\cdot\, & \cdot
\\
\sum_\nu^{1\cdots n}\,\varphi_{1\nu}\,\xi_{\nu 1} & 
\,\cdot\, & \,\cdot\, & 
\sum_\nu^{1\cdots n}\,\varphi_{1\nu}\,\xi_{\nu s}\ &
\ \sum_\nu^{1\cdots n}\,\psi_{1\nu}\,\eta_{\nu 1} &
 \,\cdot\, & \,\cdot\, &
\sum_\nu^{1\cdots n}\,\psi_{1\nu}\,\eta_{\nu s}
\\
\cdot & \,\cdot\, & \,\cdot\, & \cdot &
\cdot & \,\cdot\, & \,\cdot\, & \cdot
\\
\sum_\nu^{1\cdots n}\,\varphi_{r-n,\nu}\,\xi_{\nu 1} & 
\,\cdot\, & \,\cdot\, & 
\sum_\nu^{1\cdots n}\,\varphi_{r-n,\nu}\,\xi_{\nu s}\ &
\ \sum_\nu^{1\cdots n}\,\psi_{r-n,\nu}\,\eta_{\nu 1} &
 \,\cdot\, & \,\cdot\, &
\sum_\nu^{1\cdots n}\,\psi_{r-n,\nu}\,\eta_{\nu s}
\end{array}
\right\vert,
\]
and when one makes afterwards the substitution: $\psi_{ k\nu} ( y) =
\varphi_{ k\nu}(x)$; the $(n+1) \times (n+1)$ determinants of the so
obtained matrix are all identically zero. Consequently, \thetag{ 4}
belongs to the systems of equations that one obtains by setting equal
to zero all $(n+1) \times (n+1)$ determinants of the matrix~\thetag{
11}, and therefore, according to the theorem cited above,
it admits the $r$-term group $\Omega_1f, \dots, \Omega_rf$. 

\smallskip

The important property of the system of equations~\thetag{ 4}
just proved can also yet be realized in another, somehow
more direct way. 

According to p.~\pageref{S-109-sq} and to p.~\pageref{S-223}, the
system of equations~\thetag{ 4} admits in any case the $r$-term group
$\Omega_1f, \dots, \Omega_rf$ only when all equations of the form:
\[
\Omega_j
\big(\varphi_{k\nu}(x)
-
\psi_{k\nu}(y)\big)
=
X_j\,\varphi_{k\nu}(x)
-
Y_j\,\psi_{k\nu}(y)
=
0
\]
are a consequence of~\thetag{ 4}. That this condition is
satisfied in the present case can be easily verified.

For $j = 1, \dots, r$ and $k = 1, \dots, r-n$, we have:
\[
\aligned
\leftbracket
X_j,\,X_{n+k}
\rightbracket
&
=
\bigg\leftbracket
X_jf,\,\,
\sum_{\nu=1}^n\,\varphi_{k\nu}\,X_\nu f
\bigg\rightbracket
\\
&
=
\sum_{\nu=1}^n\,X_j\,\varphi_{k\nu}\,X_\nu f
+
\sum_{\nu=1}^n\,\varphi_{k\nu}
\leftbracket
X_j,\,X_\nu
\rightbracket.
\endaligned
\] 

Moreover, we have in general:
\[
\leftbracket
X_j,\,X_\mu
\rightbracket
=
\sum_{\pi=1}^r\,c_{j\mu\pi}\,X_\pi f
=
\sum_{\nu=1}^n\,
\bigg\{
c_{j\mu\nu}
+
\sum_{\tau=1}^{r-n}\,
c_{j\mu,\,n+\tau}\,\varphi_{\tau\nu}
\bigg\}\,
X_\nu f.
\]

If we insert these values in the preceding equation and if in
addition, we take into account that $X_1f, \dots, X_nf$ are not linked
together by a linear relation, we then find:
\def\theequation{12}\begin{equation}
\left\{
\aligned
X_j\,\varphi_{k\nu}
&
=
c_{j,\,n+k,\,\nu}
+
\sum_{\tau=1}^{r-n}\,
c_{j,\,n+k,\,n+\tau}\,\varphi_{\tau\nu}
\\
&
\ \ \ \ \
-
\sum_{\mu=1}^n\,\varphi_{k\mu}
\bigg\{
c_{j\mu\nu}
+
\sum_{\tau=1}^{r-n}\,
c_{j,\,\mu,\,n+\tau}\,\varphi_{\tau\nu}
\bigg\}.
\endaligned\right.
\end{equation}
Consequently, we have:
\[
\label{S-343-sq}
X_j\,\varphi_{k\nu}
=
\Pi_{jk\nu}
(\varphi_{11},\,\varphi_{12},\,\dots)
\ \ \ \ \ \ \ \ \ \ \ \ \ 
{\scriptstyle{(j\,=\,1\,\cdots\,r\,;\,\,\,
k\,=\,1\,\cdots\,r\,-\,n\,;\,\,\,
\nu\,=\,1\,\cdots\,r)}},
\]
and a completely similar computation gives:
\[
Y_j\,\psi_{k\nu}
=
\Pi_{jk\nu}
(\psi_{11},\,\psi_{12},\,\dots),
\]
where in the two cases $\Pi_{ jk\nu}$ denote the same functions of
their arguments.

At present, we obtain:
\[
\Omega_j(\varphi_{k\nu}-\psi_{k\nu})
=
\Pi_{jk\nu}(\varphi_{11},\,\varphi_{12},\,\dots)
-
\Pi_{jk\nu}(\psi_{11},\,\psi_{12},\,\dots),
\]
from which it is to be seen that the expressions $\Omega_j ( \varphi_{
k\nu} - \psi_{ k\nu})$ effectively vanish by virtue
of~\thetag{ 4}. 

In general, the $n ( r - n)$ equations $\varphi_{ k\nu} (x) = \psi_{
k\nu} (x)$ are not be mutually independent, but rather, they can be
replaced by a smaller number of mutually independent equations, say by
the following $s - \rho \leqslant s$ following ones:
\def\theequation{13}\begin{equation}
\varphi_k(x_1,\dots,x_s)
=
\psi_k(y_1,\dots,y_s)
\ \ \ \ \ \ \ \ \ \ \ \ \ 
{\scriptstyle{(k\,=\,1\,\cdots\,s\,-\,\rho)}}.
\end{equation}
Naturally, each $X_j \, \varphi_k$ will then be
a function of $\varphi_1, \dots, \varphi_{ s - \rho}$ alone:
\def\theequation{14}\begin{equation}
X_j\,\varphi_k
=
\pi_{jk}(\varphi_1,\dots,\varphi_{s-\rho})
\ \ \ \ \ \ \ \ \ \ \ \ \ 
{\scriptstyle{(j\,=\,1\,\cdots\,r\,;\,\,\,
k\,=\,1\,\cdots\,s\,-\,\rho)}},
\end{equation}
and every $Y_j\, \psi_k$ will be the same function of $\psi_1, \dots,
\psi_{ s - \rho}$:
\def\theequation{14'}\begin{equation}
Y_j\,\psi_k
=
\pi_{jk}(\psi_1,\dots,\psi_{s-\rho})
\ \ \ \ \ \ \ \ \ \ \ \ \ 
{\scriptstyle{(j\,=\,1\,\cdots\,r\,;\,\,\,
k\,=\,1\,\cdots\,s\,-\,\rho)}},
\end{equation}
hence all $\Omega_j ( \varphi_k - \psi_k)$ vanish by virtue of the
system of equations: $\varphi_1 = \psi_1$, \dots, $\varphi_{ s - \rho}
= \psi_{ s - \rho}$. But since this system of equations is presented
under a form which satisfies the requirement set on 
p.~\pageref{S-107-sq}
sq., then according to p.~\pageref{Theorem-14-S-112} and to
p.~\pageref{S-223}, we conclude that it admits the $r$-term group
$\Omega_1 f, \dots, \Omega_r f$.

\plainstatement{Now, if $s - \rho = s$, so $\rho = 0$, 
\label{S-344}
then the $s$
equations $\varphi_1 = \psi_1$, \dots, $\varphi_s = \psi_s$ actually
represent already a transformation in the variables $x_1, \dots, x_s$,
$y_1, \dots, y_s$. This transformation obviously transfers $X_1f,
\dots, X_rf$ to $Y_1f, \dots, Y_rf$, respectively, and at the same
time, it is the most general transformation that does this.}

Thus, the case $s - \rho = s$ is settled from the beginning,
but by contrast, the case $s - \rho < s$ requires a closer study.

\medskip

In order to simplify the additional considerations, we want at first
to introduce appropriate new variables in place of the variables $x$
and $y$.

\label{S-344-sq}
The $n$ mutually independent equations: $X_1f = 0$, \dots, $X_n f = 0$
form an $n$-term complete system in the $s$ variables $x_1, \dots,
x_s$, hence they have $s - n$ independent solutions in common, and
likewise, the $n$ equations: $Y_1f = 0$, \dots, $Y_n f = 0$ in the
variables $y_1, \dots, y_s$ have exactly $s - n$ independent solutions
in common.

It stands to reason to simplify the infinitesimal transformations
$X_1f, \dots, X_rf$ and $Y_1f, \dots, Y_rf$ by introducing, in place
of the $x$, new independent variables of which $s-n$ are independent
solutions of the complete system: $X_1f = 0$, \dots, $X_n f = 0$, and
by introducing, in place of the $y$, new independent variables of
which $s - n$ are independent solutions of the complete system: $Y_1f
= 0$, \dots, $Y_n f = 0$.

On the other hand, it stands to reason to simplify the
equations~\thetag{ 13} by introducing the functions $\varphi_1 (x),
\dots, \varphi_{ s - \rho}(x)$ and $\psi_1 ( y), \dots, \psi_{ s-\rho}
(y)$ as new variables.

We want to attempt to combine the two simplifications as far as
possible.

At first, we start by introducing appropriate new variables in place
of the $x$.

Amongst the solutions of the complete system $X_1f = 0$, \dots, $X_n f
= 0$, there can be some which can be expressed in terms of the
functions $\varphi_1 ( x)$, \dots, $\varphi_{ s - \rho} (x)$ alone;
all solutions $F ( \varphi_1, \dots, \varphi_{ s - \rho})$ of this
nature determine themselves from the $n$ differential equations:
\def\theequation{15}\begin{equation}
\sum_{\nu=1}^{s-\rho}\,
X_k\,\varphi_\nu\,
\frac{\partial F}{\partial\varphi_\nu}
=
\sum_{\nu=1}^{s-\rho}\,\pi_{k\nu}
(\varphi_1,\dots,\varphi_{s-\rho})\,
\frac{\partial F}{\partial\varphi_\nu}
=
0
\ \ \ \ \ \ \ \ \ \ \ \ \ {\scriptstyle{(k\,=\,1\,\cdots\,n)}}
\end{equation}
in the $s - \rho$ variables $\varphi_1, \dots, \varphi_{ s - \rho}$.
We want to suppose that these equations possess exactly $s - q
\leqslant s - \rho$ independent solutions, say the following ones:
\[
\label{S-345}
\mathfrak{U}_1(\varphi_1,\dots,\varphi_{s-\rho}),
\,\,\,\dots,\,\,\,
\mathfrak{U}_{s-q}(\varphi_1,\dots,\varphi_{s-\rho}).
\]

Under these assumptions:
\[
\mathfrak{U}_1
\big(\varphi_1(x),\dots,\varphi_{s-\rho}(x)\big)
=
u_1(x),
\,\,\,\dots,\,\,\,
\mathfrak{U}_{s-q}
\big(\varphi_1(x),\dots,\varphi_{s-\rho}(x)\big)
=
u_{s-q}(x)
\]
are obviously independent solutions of the complete system: $X_1f =
0$, \dots, $X_nf = 0$ which can be expressed in terms of $\varphi_1 (
x), \dots, \varphi_{ s - \rho} (x)$ alone, and such that every other
solution of the same constitution is a function of $u_1 ( x), \dots,
u_{ s - q} (x)$ only. Naturally, we also have at the same time the
inequation: $s - q \leqslant s - n$, hence none of the two numbers
$\rho$ and $n$ is larger than $q$.

Now, let:
\[
u_{s-q+1}
=
u_{s-q+1}(x),
\,\,\,\dots,\,\,\,
u_{s-n}
=
u_{s-n}(x)
\]
be $q - n$ arbitrary mutually independent solutions of the complete
system $X_1f = 0$, \dots, $X_n f = 0$ that are also independent of
$u_1 ( x), \dots, u_{ s - q} (x)$. We will show that the $s - \rho + q
- n$ functions: $u_{ s - q + 1} (x)$, \dots, $u_{ s-n} (x)$,
$\varphi_1 ( x), \dots, \varphi_{ s - \rho} (x)$ are mutually
independent.
\label{S-345-bis}

Since $\varphi_1 ( x), \dots, \varphi_{ s - \rho} (x)$ are mutually
independent, there are, amongst the functions $\varphi_1 ( x)$, \dots,
$\varphi_{ s - \rho} (x)$, $u_{ s-q+1}(x)$, \dots, $u_{ s-n}(x)$ at
least $s - \rho$, say exactly $s - \rho + q - n - h$ that are mutually
independent, where $0 \leqslant h \leqslant q - n$. We want to suppose
that precisely $\varphi_1 ( x)$, \dots, $\varphi_{ s - \rho} (x)$,
$u_{ s-q+ h + 1}(x)$, \dots, $u_{ s-n}(x)$ are mutually independent,
while $u_{ s - q + 1} (x)$, \dots, $u_{ s - q + h} (x)$ can be
expressed in terms of $\varphi_1 ( x)$, \dots, $\varphi_{ s - \rho}
(x)$, $u_{ s - q + h + 1} (x)$, \dots, $u_{ s - n}$ only.  Thus, there
must exist, between the quantities $u_{ s - q + 1}, \dots, u_{ s-n}$,
$\varphi_1, \dots, \varphi_{ s - \rho}$ relations of the form:
\def\theequation{16}\begin{equation}
u_{s-q+j}
=
\chi_j
\big(
u_{s-q+h+1},\,\dots,\,u_{s-n},\,\,
\varphi_1,\,\dots,\,\varphi_{s-\rho}
\big)
\ \ \ \ \ \ \ \ \ \ \ \ \ 
{\scriptstyle{(j\,=\,1\,\cdots\,h)}}
\end{equation}
which reduce to identities after the substitution:
\def\theequation{17}\begin{equation}
\left\{
\aligned
u_{s-q+1}
&
=
u_{s-q+1}(x),
\,\,\,\dots,\,\,\,
u_{s-n}
=
u_{s-n}(x)
\\
\varphi_1
&
=
\varphi_1(x),
\,\,\,\dots,\,\,\,
\varphi_{s-\rho}
=
\varphi_{s-\rho}(x).
\endaligned\right.
\end{equation}

If we interpret the substitution~\thetag{ 17} 
by the sign $[ \, \, \, ]$, we then have:
\[
[\,u_{s-q+j}-\chi_j\,]
\equiv
0
\ \ \ \ \ \ \ \ \ \ \ \ \ 
{\scriptstyle{(j\,=\,1\,\cdots\,h)}},
\]
whence it comes:
\[
X_k\,
[\,u_{s-q+j}-\chi_j\,]
=
-\,
\sum_{\nu=1}^{s-\rho}\,
[\,\pi_{k\nu}(\varphi_1,\dots,\varphi_{s-\rho})\,]\,\,
\bigg[\,
\frac{\partial\chi_j}{\partial\varphi_\nu}
\,\bigg]
\equiv
0
\ \ \ \ \ \ \ \ \ \ \ \ \ 
{\scriptstyle{(k\,=\,1\,\cdots\,n\,;\,\,\,
j\,=\,1\,\cdots\,h)}},
\]
that is to say: all the expressions:
\def\theequation{18}\begin{equation}
\sum_{\nu=1}^{s-\rho}\,
\pi_{k\nu}(\varphi_1,\dots,\varphi_{s-\rho})\,
\frac{\partial\chi_j}{\partial\varphi_\nu}
\ \ \ \ \ \ \ \ \ \ \ \ \ 
{\scriptstyle{(k\,=\,1\,\cdots\,n\,;\,\,\,
j\,=\,1\,\cdots\,h)}}
\end{equation}
vanish identically after the substitution~\thetag{ 17}. But now,
these expressions are all free of $u_{ s - q +1 }, \dots, u_{ s - q +
h}$, so if they were not identically zero and they would vanish
identically after the substitution~\thetag{ 17}, then the functions
$u_{ s - q + h + 1} (x)$, \dots, $u_{ s-n}( x)$, $\varphi_1 (x)$,
\dots, $\varphi_{s - \rho} (x)$ would not be mutually independent, but
this is in contradiction to the assumption. Consequently, the
expressions~\thetag{ 18} are in fact identically zero, or, what
amounts to the same, the functions $\chi_1, \dots, \chi_h$ are
solutions of the $n$ differential equations~\thetag{ 15} in the $s -
\rho$ variables $\varphi_1, \dots, \varphi_{ s - \rho}$. From this,
it results that $\varphi_1, \dots, \varphi_{ s - \rho}$ appear in the
$\chi_j$ only in the combination: $\mathfrak{ U}_1 ( \varphi_1, \dots,
\varphi_{ s - \rho})$, \dots, $\mathfrak{ U}_{ s - q} (\varphi_1,
\dots, \varphi_{ s - \rho})$, so that the $h$ equations~\thetag{ 16}
can be replaced by $h$ relations of the form:
\def\theequation{19}\begin{equation}
u_{s-q+j}
=
\overline{\chi}_j
\big(
u_{s-q+h+1},\,\dots,\,u_{s-n},\,\,u_1,\dots,u_{s-q}
\big)
\ \ \ \ \ \ \ \ \ \ \ \ \ {\scriptstyle{(j\,=\,1\,\cdots\,h)}},
\end{equation}
which in turn reduce to identities after the substitution:
\[
u_1=u_1(x),
\,\,\,\dots,\,\,\,
u_{s-n}=u_{s-n}(x).
\]

Obviously, relations of the form~\thetag{ 19} cannot exist, for $u_1,
\dots, u_{ s-n}$ are independent solutions of the complete system
$X_1f = 0$, \dots, $X_nf = 0$; consequently, $h$ is equal to zero. As
a result, it is proved that the $s - \rho + q - n$ functions
$\varphi_1 (x), \dots, \varphi_{ s - \rho}(x)$, $u_{ s- q + 1}(x)$,
\dots, $u_{ s-n} (x)$ really are independent of each other.

We therefore see that, between the $s - n + s - \rho$ quantities
$u_1, \dots, u_{ s-n}$, $\varphi_1, \dots, \varphi_{ s - \rho}$, 
no other relations exist that are independent from the $s-q$ relations:
\def\theequation{20}\begin{equation}
u_1
=
\mathfrak{U}_1(\varphi_1,\dots,\varphi_{s-\rho}),
\,\,\,\dots,\,\,\,
u_{s-q}
=
\mathfrak{U}_{s-q}(\varphi_1,\dots,\varphi_{s-\rho})
\end{equation}
which reduce to identities after the substitution:
\[
u_1
=
u_1(x),
\,\,\,\dots,\,\,\,
u_{s-q}
=
u_{s-q}(x),
\ \ \ \ \
\varphi_1
=
\varphi_1(x),
\,\,\,\dots,\,\,\,
\varphi_{s-\rho}
=
\varphi_{s-\rho}(x).
\]

From what has been just said, it results that, amongst the $s - n + s
- \rho$ functions $u_1(x)$, \dots, $u_{ s-n}(x)$, $\varphi_1(x)$,
\dots, $\varphi_{ s - \rho}(x)$, there exist exactly $s - n + s - \rho
- (s - q) = s - n + q - \rho$ that are mutually independent, namely
for instance the $s - n$ functions: $u_1 (x), \dots, u_{ s-n} (x)$ to
which $q - \rho \geqslant 0$ of the functions $\varphi_1 (x), \dots,
\varphi_{ s - \rho}(x)$ are yet added. If we agree that the
equations~\thetag{ 20} can be resolved precisely with respect to
$\varphi_{ q - \rho +1}, \dots, \varphi_{ s - \rho}$, we can conclude
that precisely the $s - n + q - \rho$ functions: $u_1(x)$, \dots, $u_{
s-n}(x)$, $\varphi_1 (x)$, \dots, $\varphi_{ q - \rho}(x)$ are
mutually independent. Here, the number $q - \rho$, or shortly $m$,
is certainly not larger than $n$, since the sum $s - n + q - \rho$
can naturally not exceed the number $s$ of the variables $x$.

At present, we have gone so far that we can introduce new independent
variables $x_1', \dots, x_s'$ in place of $x_1, \dots, x_s$; we choose
them in the following way:

We set simply:
\[
x_{q+1}'
=
u_1(x_1,\dots,x_s),
\,\,\,\dots,\,\,\,
x_s'
=
u_{s-q}(x_1,\dots,x_s),
\]
furthermore:
\[
x_{n+1}'
=
u_{s-q+1}(x_1,\dots,x_s),
\,\,\,\dots,\,\,\,
x_q'
=
u_{s-n}(x_1,\dots,x_s),
\]
and:
\[
x_1'
=
\varphi_1(x_1,\dots,x_s),
\,\,\,\dots,\,\,\,
x_m'
=
\varphi_m(x_1,\dots,x_s),
\]
where $m = q - \rho$ is not larger than $n$; in addition, 
we yet set: 
\[
x_{m+1}'
=
\lambda_1(x_1,\dots,x_s),
\,\,\,\dots,\,\,\,
x_n'
=
\lambda_{n-m}(x_1,\dots,x_s),
\]
where $\lambda_1 (x), \dots, \lambda_{ n-m} (x)$ denote arbitrary
mutually independent functions that are also independent of $u_1 (x)$,
\dots, $u_{ s-n} (x)$, $\varphi_1 (x)$, \dots, $\varphi_m (x)$.

In a completely similar way, we introduce new independent variables
in place of $y_1, \dots, y_s$. 

We form the $s -q$ functions:
\[
\mathfrak{U}_1
\big(\psi_1(y),\,\dots,\,\psi_{s-\rho}(y)\big)
=
v_1(y),
\,\,\,\dots,\,\,\,
\mathfrak{U}_{s-q}
\big(\psi_1(y),\,\dots,\,\psi_{s-\rho}(y)\big)
=
v_{s-q}(y)
\]
that are evidently independent solutions of the $n$-term complete
system: $Y_1 f = 0$, \dots, $Y_n f = 0$; moreover, we determine any $q
- n$ arbitrary mutually independent solutions: $v_{ s -q +1} (y)$,
\dots, $v_{ s-n} (y)$ of the same complete system that are also
independent of $v_1 (y)$, \dots, $v_{ s-q} (y)$. Then it is clear
that the $s -n + q - \rho$ functions:
\[
v_1(y),\,\dots,\,v_{s-n}(y),\ \
\psi_1(y),\,\dots,\,\psi_{q-\rho}(y)
\]
are mutually independent.

We now choose the new variables $y_1', \dots, y_s'$ in 
exactly the same way as the variables $x'$ a short
while ago.

We set simply:
\[
y_{q+1}'
=
v_1(y_1,\dots,y_s),
\,\,\,\dots,\,\,\,
y_s'
=
v_{s-q}(y_1,\dots,y_s),
\]
furthermore:
\[
y_{n+1}'
=
v_{s-q+1}(y_1,\dots,y_s),
\,\,\,\dots,\,\,\,
y_q'
=
v_{s-n}(y_1,\dots,y_s),
\]
and:
\[
y_1'
=
\psi_1(y_1,\dots,y_s),
\,\,\,\dots,\,\,\,
y_m'
=
\psi_m(y_1,\dots,y_s)\,;
\]
in addition, we yet set:
\[
y_{m+1}'
=
\Lambda_1(y_1,\dots,y_s),
\,\,\,\dots,\,\,\,
y_n'
=
\Lambda_{n-m}(y_1,\dots,y_s),
\]
where $\Lambda_1, \dots, \Lambda_{ n-m}$ are arbitrary mutually
independent functions that are also independent of: $v_1 (y)$, \dots,
$v_{ s-n} (y)$, $\psi_1 (y)$, \dots, $\psi_m (y)$.

At present, we introduce the new variables $x'$ and $y'$ in the
infinitesimal transformations $X_kf$, $Y_kf$ and in the
equations~\thetag{13}.

Since all $X_k\, x_{ n+1}'$, \dots, $X_k\, x_s'$
vanish identically and
since all $X_k\, \varphi_1$, \dots, $X_k\, \varphi_m$
depend only on $\varphi_1, \dots, \varphi_{ s - \rho}$, 
that is to say, only on
$x_1', \dots, x_m'$, $x_{ q+1}', \dots, x_s'$, the $X_k f$
receive the following form in 
the variables $x_1', \dots, x_s'$: 
\[
\aligned
X_kf
&
=
\sum_{\mu=1}^m\,\mathfrak{x}_{k\mu}
\big(
x_1',\dots,x_m',\,\,x_{q+1}',\dots,x_s'
\big)\,
\frac{\partial f}{\partial x_\mu'}
+
\\
&
\ \ \ \ \
+
\sum_{j=1}^{n-m}\,\mathfrak{x}_{k,\,m+j}
\big(
x_1',\dots,x_m',\dots,x_n',\dots,x_q',\dots,x_s'
\big)\,
\frac{\partial f}{\partial x_{m+j}'}
=
\Xi_kf.
\endaligned
\]
In the same way, we have:
\[
\aligned
Y_kf
&
=
\sum_{\mu=1}^m\,\mathfrak{x}_{k\mu}
\big(
y_1',\dots,y_m',\,\,y_{q+1}',\dots,y_s'
\big)\,
\frac{\partial f}{\partial y_\mu'}
+
\\
&
\ \ \ \ \
+
\sum_{j=1}^{n-m}\,
\mathfrak{y}_{k,\,m+j}
\big(
y_1',\dots,y_m',\dots,y_n',\dots,y_q',\dots,y_s'
\big)\,
\frac{\partial f}{\partial y_{m+j}'}
=
H_kf.
\endaligned
\]
Here, the $\mathfrak{ x}_{ k\mu} (y_1', \dots, y_m', \, y_{ q+1}',
\dots, y_s')$ denote the same functions of their arguments as the
$\mathfrak{ x}_{ k\mu} (x_1', \dots, x_m', \, x_{ q+1}', \dots,
x_s')$. Indeed, it results from the equations~\thetag{ 14}
and~\thetag{ 14'} that $Y_k\, y_\mu'$ is the same function of $y_1',
\dots, y_m'$, $y_{ q+1}', \dots, y_s'$ as $X_k\, x_\mu'$ is of $x_1',
\dots, x_m'$, $x_{ q+1}', \dots, x_s'$, where it is understood that
$\mu$ is an arbitrary number amongst $1, 2, \dots, m$.

\medskip

On the other hand, we must determine which form the system of
equations~\thetag{ 13} receives in the new variables.

\label{S-349-sq}
The system of equations can evidently be replaced by the following one:
\[
\aligned
& 
\ \ \ \ \ \ \ \ \ \
\varphi_1-\psi_1=0,
\,\,\,\dots,\,\,\,
\varphi_m-\psi_m=0,
\\
&
\ \ \ \
\mathfrak{U}_1(\varphi_1,\dots,\varphi_{s-\rho})
-
\mathfrak{U}_1(\psi_1,\dots,\psi_{s-\rho})
=
0,\,\,\,\dots,
\\
&
\mathfrak{U}_{s-q}(\varphi_1,\dots,\varphi_{s-\rho})
-
\mathfrak{U}_{s-q}(\psi_1,\dots,\psi_{s-\rho})
=
0.
\endaligned
\]
If we introduce our new variables in this system, we obviously obtain
the simple system:
\def\theequation{21}\begin{equation}
\left\{
\aligned
x_1'-y_1'=0,
\,\,\,\dots,\,\,\,
x_m'-y_m'=0,
\\
x_{q+1}'-y_{q+1}'=0,
\,\,\,\dots,\,\,\,
x_s'-y_s'=0\,;
\endaligned\right.
\end{equation}
so the system of equations~\thetag{ 13}, or, what is the same, the
system of equations~\thetag{ 4}, can be brought to this form after the
introduction of the new variables.

Finally, if we remember that neither $X_1f, \dots, 
X_nf$, nor $Y_1f, \dots, Y_nf$ are linked by linear
relations and that all $(n+1) \times
(n+1)$ determinants of the matrix~\thetag{ 11}
vanish by means of~\thetag{ 4} while
not all $n \times n$ determinants do, we
then realize that in the same way, 
neither $\Xi_1 f, \dots, \Xi_nf$, nor
$H_1f, \dots, H_nf$ are linked by linear relations, 
and that all $(n+1) \times (n+1)$
determinants, but no all $n \times n$ ones, 
of the matrix which can be formed with the
coefficients of the differential quotients
of $f$ in the $r$ infinitesimal transformations
$\Omega_k f = \Xi_k f + H_k f$
vanish by means of~\thetag{ 21}. 

\medskip

Thanks to the preceding developments, the problem
stated in the outset of the paragraph, p.~\pageref{S-341}, 
is lead back to the following simpler problem: 

\plainstatement{
To seek, in the $2s$ variables $x_1', \dots, x_s'$, 
$y_1', \dots, y_s'$, a system of equations
which admits the $r$-term group:
\[
\Omega_kf
=
\Xi_kf+H_kf
\ \ \ \ \ \ \ \ \ \ \ \ \ {\scriptstyle{(k\,=\,1\,\cdots\,r)}}
\]
and in addition which consists of $s$ independent equations that are
solvable with respect to $x_1', \dots, x_s'$ and to $y_1', \dots,
y_s'$ as well, and lastly, which comprises the $s - q + m$
equations~\thetag{ 21}.}

In order to solve this new problem, we remember
Chap.~\ref{kapitel-14}, p.~\pageref{Theorem-40-S-233} sq.; from what
was said at that time, we deduce that every system of equations which
admits the group $\Xi_k f + H_k f = \Omega_kf$ and which comprises at
the same time the equations~\thetag{ 21} can be obtained by adding to
the equations~\thetag{ 21} a system of equations in the $s + q - m$
variables $x_1'$, \dots, $x_m'$, \dots, $x_n'$, \dots, $x_q'$, \dots,
$x_s'$, \dots, $y_{ m+1}', \dots, y_q'$ that admits the $r$-term
group:
\[
\aligned
\overline{\Omega}_kf
&
=
\sum_{\mu=1}^m\,
\mathfrak{x}_{k\mu}
\big(
x_1',\dots,x_m',\,\,x_{q+1}',\dots,x_s'
\big)\,
\frac{\partial f}{\partial x_\mu'}
+
\\
&
\ \ \ \ \
+
\sum_{j=1}^{n-m}\,
\mathfrak{x}_{k,\,m+j}
\big(
x_1',\dots,x_m',\dots,x_n',\dots,x_q',\dots,x_s'
\big)\,
\frac{\partial f}{\partial x_{m+j}'}
+
\\
&
\ \ \ \ \
+
\sum_{j=1}^{n-m}\,
\mathfrak{y}_{k,\,m+j}
\big(
x_1',\dots,x_m',\,\,y_{m+1}',\dots,y_n',\dots,y_q',\,\,
x_{q+1}',\dots,x_s'
\big)\,
\frac{\partial f}{\partial y_{m+j}'}.
\endaligned
\]

Here, $\overline{ \Omega}_kf$ is obtained by leaving out all terms
with $\partial f / \partial y_1'$, \dots, $\partial f / \partial y_m'$
in $\Omega_k f = \Xi_k f + H_kf$ and by making everywhere the
substitution:
\[
y_1'=x_1',
\,\,\,\dots,\,\,\,
y_m'=x_m',\ \ \
y_{q+1}'=x_{q+1}',
\,\,\,\dots,\,\,\,
y_s'=x_s'
\]
by means of~\thetag{ 21} in the remaining terms. This formation of
$\overline{ \Omega}_kf$ shows that in the matrix formed with
$\overline{ \Omega}_1f, \dots, \overline{ \Omega}_rf$:
\def\theequation{22}\begin{equation}
\left\vert
\begin{array}{cccccccc}
\mathfrak{x}_{11}(x') & \,\cdot\, & \,\cdot\, &
\mathfrak{x}_{1n}(x')\, & \,\mathfrak{y}_{1,\,m+1}(x',y') 
& \,\cdot\, & \,\cdot\, &
\mathfrak{y}_{1n}(x',y')
\\
\cdot & \,\cdot\, & \,\cdot\, & \cdot & \cdot 
& \,\cdot\, & \,\cdot\, & \cdot
\\
\mathfrak{x}_{r1}(x') & \,\cdot\, & \,\cdot\, &
\mathfrak{x}_{rn}(x')\, & \,\mathfrak{y}_{r,\,m+1}(x',y') 
& \,\cdot\, & \,\cdot\, &
\mathfrak{y}_{rn}(x',y')
\end{array}
\right\vert,
\end{equation}
all $(n+1) \times (n+1)$ determinants vanish identically, but not all
$n \times n$ ones. If we yet add that neither $\Xi_1f, \dots, \Xi_nf$
nor $H_1f, \dots, H_nf$ are linked by linear relations, we realize
immediately that in particular the two $n \times n$ determinants:
\def\theequation{23}\begin{equation}
\sum\,\pm\,\mathfrak{x}_{11}(x')\cdots\,\mathfrak{x}_{nn}(x')
\end{equation}
and:
\def\theequation{24}\begin{equation}
\left\vert
\begin{array}{cccccccc}
\mathfrak{x}_{11}(x') & \,\cdot\, & \,\cdot\, &
\mathfrak{x}_{1m}(x')\, & \,\mathfrak{y}_{1,\,m+1}(x',y') 
& \,\cdot\, & \,\cdot\, &
\mathfrak{y}_{1n}(x',y')
\\
\cdot & \,\cdot\, & \,\cdot\, & \cdot & \cdot 
& \,\cdot\, & \,\cdot\, & \cdot
\\
\mathfrak{x}_{n1}(x') & \,\cdot\, & \,\cdot\, &
\mathfrak{x}_{nm}(x')\, & \,\mathfrak{y}_{n,\,m+1}(x',y') 
& \,\cdot\, & \,\cdot\, &
\mathfrak{y}_{nn}(x',y')
\end{array}
\right\vert
\end{equation}
are not identically zero.

But now, the matter for us is not to find \emphasis{all} systems of
equations which comprise~\thetag{ 21} and which admit the group
$\Omega_1f, \dots, \Omega_rf$, and the matter is only to determine
systems of equations of this sort which consist of exactly $s$
independent equations that are solvable both with respect to $x_1',
\dots, x_s'$ and to $y_1', \dots, y_s'$. Consequently, we do not have
to set up \emphasis{all} systems of equations in $x_1', \dots, x_s'$,
$y_{ m+1}', \dots, y_q'$ which admit $\overline{ \Omega}_1f, \dots,
\overline{ \Omega}_rf$, but only the systems which consist of exactly
$s - (m+s-q) = q -m$ independent equations and which in addition are
solvable both with respect to $x_{ m+1}', \dots, x_n', \dots, x_q'$
and with respect to $y_{ m+1}', \dots, y_n', \dots, y_q'$.

If a system of equation in $x_1', \dots, x_s'$, $y_{ m+1}', \dots,
y_q'$ satisfies the requirement just stated, then it can be brought to
the form:
\def\theequation{25}\begin{equation}
y_{m+j}'
=
\Phi_{m+j}
(x_1',\dots,x_m',\dots,x_n',\dots,x_q',\dots,x_s')
\ \ \ \ \ \ \ \ \ \ \ \ \ 
{\scriptstyle{(j\,=\,1\,\cdots\,q\,-\,m)}},
\end{equation}
where the functions $\Phi_{ m+1}, \dots, \Phi_q$ are in turn
independent relatively to $x_{ m+1}', \dots, x_q'$. Now, it is clear
that a system of equations of the form~\thetag{ 25} cannot bring to
zero all $n \times n$ determinants of the matrix~\thetag{ 22}, since
in any case, the determinant~\thetag{ 23} cannot be equal to zero by
virtue of~\thetag{ 25}. On the other hand, since all $(n+1) \times
(n+1)$ determinants of~\thetag{ 22} vanish identically, then according
to Chap.~\ref{kapitel-14}, p.~\pageref{S-230}, it comes that every
system of equations of the form~\thetag{ 25} which admits the group
$\overline{ \Omega}_1f, \dots, \overline{ \Omega}_rf$ is represented
by relations between the solutions of the equations: $\overline{
\Omega}_1f = 0$, \dots, $\overline{ \Omega}_r f = 0$, or, what is the
same, by relations between the solutions of the $n$-term complete
system: $\overline{ \Omega}_1f = 0$, \dots, $\overline{ \Omega}_n f =
0$.

The $n$-term complete system $\overline{ \Omega}_1f = 0$, \dots,
$\overline{ \Omega}_n f = 0$ contains $s + q - m$ independent
variables and possesses therefore $s - n + q - m$ independent
solutions; one can immediately indicate $s - n + q - n$ independent
solutions, namely: $x_{ n+1}'$, \dots, $x_q'$, \dots, $x_s'$, $y_{
n+1}'$, \dots, $y_q'$, while the $n - m$ lacking ones must be
determined by integration and can obviously be brought to the form:
\[
\omega_1
(x_1',\dots,x_m',\dots,x_s',\,y_{m+1}',\dots,y_q'),
\,\,\,\dots,\,\,\,
\omega_{n-m}
(x_1',\dots,x_m',\dots,x_s',\,y_{m+1}',\dots,y_q').
\]
Since the two determinants~\thetag{ 23} and~\thetag{ 24} do not vanish
identically, the equations of our complete system are solvable both
with respect to $\partial f / \partial x_1'$, \dots, $\partial f /
\partial x_m'$, \dots, $\partial f / \partial x_n'$ and with respect
to $\partial f / \partial x_1'$, \dots, $\partial f / \partial x_m'$,
$\partial f / \partial y_{ m+1}'$, \dots, $\partial f / \partial y_n'$
and therefore, its $s - n + q - m$ independent solutions $x_{ n+1}',
\dots, x_s'$, $y_{ n+1}', \dots, y_q'$, $\omega_1, \dots, \omega_{
n-m}$ are mutually independent both relatively to $x_{ n+1}'$, \dots,
$x_q'$, \dots, $x_s'$, $y_{ m+1}', \dots, y_q'$ and relatively to $x_{
m+1}'$, $\dots$, $x_n'$, \dots, $x_q'$, \dots, $x_s'$, $y_{ n+1}',
\dots, y_q'$ (cf. Theorem~12, p.~\pageref{Theorem-12-S-91}).
Consequently, the functions $\omega_1, \dots, \omega_{ n-m}$ are
mutually independent both relatively to $y_{ m+1}', \dots, y_n'$ and
relatively to $x_{ m+1}', \dots, x_n'$.

Every system of equations~\thetag{ 25} which satisfies the stated
requirements is represented by relations between the solutions $x_{
n+1}', \dots, x_s'$, $y_{ n+1}', \dots, y_q'$, $\omega_1, \dots,
\omega_{ n-m}$ and to be precise, by $q - m$ relations that are
solvable both with respect to $y_{ m+1}', \dots, y_q'$ and with
respect to $x_{ m+1}', \dots, x_q'$. Visibly, these relations must be
solvable both with respect to $\omega_1, \dots, \omega_{ n-m}$, $x_{
n+1}', \dots, x_q'$ and with respect to $\omega_1, \dots, \omega_{
n-m}$, $y_{ n+1}', \dots, y_q'$, so they must have the form:
\def\theequation{26}\begin{equation}
\left\{
\aligned
&
\omega_\mu
(x_1',\dots,x_m',·\,x_{m+1}',\dots,x_s',\,y_{m+1}',\dots,y_q')
=
\chi_\mu(x_{n+1}',\dots,x_q',\dots,x_s')
\\
&
\ \ \ \ \ \ \ \ \ \ \ \ \ \ \ \ \ \ \ \ \ \ \ \ \ \ 
\ \ \ \ \ \ \ \ \ \ \ \ \ \ \ \ 
{\scriptstyle{(\mu\,=\,1\,\cdots\,n\,-\,m)}}
\\
&
y_{n+1}'
=
\Pi_1(x_{n+1}',\dots,x_q',\dots,x_s'),
\,\,\,\dots,\,\,\,
y_q'
=
\Pi_{q-n}(x_{n+1}',\dots,x_q',\dots,x_s'),
\endaligned\right.
\end{equation}
where $\Pi_1, \dots, \Pi_{ q-n}$ are mutually independent
relatively to $x_{ n+1}', \dots, x_q'$. 

Conversely, every system of equations of the form~\thetag{ 26} in
which $\Pi_1, \dots, \Pi_{ q - n}$ are mutually independent relatively
to $x_{ n+1}', \dots, x_q'$ is solvable both with respect to $y_{
m+1}', \dots, y_q'$ and with respect to $x_{ m+1}', \dots, x_q'$ and
since in addition, it admits the group $\overline{ \Omega}_1f, \dots,
\overline{ \Omega}_rf$, then it possesses all properties which the
sought system of equations in the variables $x_1', \dots, x_s'$, $y_{
m+1}', \dots, y_q'$ should have.

From this, we conclude that the equations~\thetag{ 26} represent the
most general system of equations which admits the group $\overline{
\Omega}_1f, \dots, \overline{ \Omega}_rf$, which consists of exactly
$q -m$ independent equations, and which is solvable both with respect
to $y_{ m+1}', \dots, y_q'$ and with respect to $x_{ m+1}', \dots,
x_q'$; about it, $\chi_1, \dots, \chi_{ n-m}$ are absolutely arbitrary
functions of their arguments and $\Pi_1, \dots, \Pi_{ q - n}$ as well,
though with the restriction for the latter that they must be mutually
independent relatively to $x_{ n+1}', \dots, x_q'$.

\medskip

Thus, if we add the system of equation~\thetag{ 26} to the
equations~\thetag{ 21}, we obtain the most general system of equations
in $x_1', \dots, x_s'$, $y_1', \dots, y_s'$ which admits the group
$\Omega_k = \Xi_kf + H_kf$, which comprises the equations~\thetag{
21}, which consists of $s$ independent equations and which is solvable
both with respect to $x_1', \dots, x_s'$ and with respect to $y_1',
\dots, y_s'$. Finally, if, in this system of equations, we express
the variables $x'$ and $y'$ in terms of the initial variables $x$ and
$y$, we obtain the most general system of equations which admits the
group $\Omega_k f = X_k f + Y_k f$, which comprises the
equations~\thetag{ 4}, which consists of $s$ independent equations and
which is solvable both with respect to $x_1, \dots, x_s$ and with
respect to $y_1, \dots, y_s$. In other words: we obtain the most
general transformation $y_i = \Phi_i ( x_1, \dots, x_s)$ which
transfers $X_1f, \dots, X_rf$ to $Y_1f, \dots, Y_rf$, respectively.

As a result, the problem stated on p.~\pageref{S-334} is settled, and
it has been achieved even more than what was actually required there,
because we not only know a transformation of the demanded
constitution, we know all these transformations. At the same time, the
claim stated on p.~\pageref{S-333-bis} is proved, namely it is proved
that the two groups $X_1f, \dots, X_rf$ and $Z_1f, \dots, Z_rf$ are
really similar to each other, under the assumptions made there.

Lastly, on the basis of the developments preceding, we can yet
determine the most general transformation which transfers the group
$X_1f, \dots, X_rf$ to the group $Z_1f, \dots, Z_rf$; indeed, if we
remember the considerations of p.~\pageref{S-331} sq., and if we
combine them with the result gained just now, we then see immediately
that the transformation in question can be found in the following way:
In the infinitesimal transformations:
\[
Y_kf
=
\sum_{j=1}^r\,
\overline{g}_{kj}\,Z_jf
\ \ \ \ \ \ \ \ \ \ \ \ \ {\scriptstyle{(k\,=\,1\,\cdots\,r)}}
\]
defined on p.~\pageref{S-334-quart}, one chooses the constants
$\overline{ g}_{kj}$ in the most general way and afterwards, following
the method given by us, one determines the most general transformation
which transfers $X_1f, \dots, X_rf$ to $Y_1f, \dots, Y_rf$,
respectively; this then is at the same time the most general
transformation which actually transfers the group $X_1f, \dots, X_rf$
to the group $Z_1f, \dots, Z_rf$.

As the equations~\thetag{ 26} show, the transformation found in this
way contains $n - m + q - n = q - m$ arbitrary functions of $s - n$
arguments, and in addition, certain arbitrary elements which come from
the $\overline{ g}_{ kj}$; these are firstly certain arbitrary
parameters, and secondly, certain arbitrarinesses which come from the
fact that the $\overline{ g}_{ kj}$ are determined by algebraic
operations. From this, it follows that the said transformation cannot
in all cases be represented by a single system of equations.\,---

At present, we summarize our results. 

At first, we have the

\renewcommand{\thefootnote}{\fnsymbol{footnote}}
\def\thetheorem{65}\begin{theorem}
\label{Theorem-65-S-353}
Two $r$-term groups:
\[
X_kf
=
\sum_{i=1}^s\,\xi_{ki}(x_1,\dots,x_s)\,
\frac{\partial f}{\partial x_i}
\ \ \ \ \ \ \ \ \ \ \ \ \ {\scriptstyle{(k\,=\,1\,\cdots\,r)}}
\]
and:
\[
Z_kf
=
\sum_{i=1}^s\,\zeta_{ki}(y_1,\dots,y_s)\,
\frac{\partial f}{\partial y_i}
\ \ \ \ \ \ \ \ \ \ \ \ \ {\scriptstyle{(k\,=\,1\,\cdots\,r)}}
\]
in the same number of variables are similar to each other
if and only if the following conditions are satisfied:

\terminology{Firstly:} the two groups must be equally composed; 
so, if the relations:
\[
\leftbracket
X_i,\,X_k
\rightbracket
=
\sum_{\sigma=1}^r\,
c_{ik\sigma}\,X_\sigma f
\]
hold, it must be possible to determine $r^2$ constants
$g_{ kj}$ such that the $r$ 
infinitesimal transformations:
\[
\mathfrak{Y}_kf
=
\sum_{j=1}^r\,g_{kj}\,Z_jf
\ \ \ \ \ \ \ \ \ \ \ \ \ {\scriptstyle{(k\,=\,1\,\cdots\,r)}}
\]
are mutually independent and such that the relations:
\[
\leftbracket
\mathfrak{Y}_i,\,\mathfrak{Y}_k
\rightbracket
=
\sum_{\sigma=1}^r\,c_{ik\sigma}\,
\mathfrak{Y}_\sigma f
\]
are identically satisfied. 

\terminology{Secondly:} if $X_1f, \dots, X_rf$ are constituted
in such a way that, $X_1f,\dots, X_nf$ (say) are linked together by 
no relation of the form:
\[
\chi_1(x_1,\dots,x_s)\,X_1f
+\cdots+
\chi_n(x_1,\dots,x_s)\,X_nf
=
0,
\]
while by contrast $X_{ n+1}f, \dots, X_rf$ express themselves linearly
in terms of $X_1f, \dots, X_nf$:
\[
X_{n+k}f
\equiv
\sum_{\nu=1}^n\,\varphi_{k\nu}(x_1,\dots,x_s)\,X_\nu f
\ \ \ \ \ \ \ \ \ \ \ \ \ 
{\scriptstyle{(k\,=\,1\,\cdots\,r\,-\,n)}},
\]
them amongst the systems of $g_{ kj}$ which satisfy the 
above requirements, there must exist at least one, 
say the system: $g_{ kj} = \overline{ g}_{ kj}$, which
is constituted in such a way that, of the
$r$ infinitesimal transformations:
\[
Y_kf
=
\sum_{j=1}^r\,\overline{g}_{kj}\,Z_jf
\ \ \ \ \ \ \ \ \ \ \ \ \ {\scriptstyle{(k\,=\,1\,\cdots\,r)}},
\]
the first $n$ ones are linked together by no linear relation
of the form:
\[
\psi_1(y_1,\dots,y_s)\,Y_1f
+\cdots+
\psi_n(y_1,\dots,y_s)\,Y_nf
=
0,
\]
while $Y_{ n+1}f, \dots, Y_rf$ express themselves 
in terms of $Y_1f, \dots, Y_nf$: 
\[
Y_{n+k}f
\equiv
\sum_{\nu=1}^n\,\psi_{k\nu}(y_1,\dots,y_s)\,Y_\nu f
\ \ \ \ \ \ \ \ \ \ \ \ \ 
{\scriptstyle{(k\,=\,1\,\cdots\,r\,-\,n)}},
\]
and such that in addition, the $n ( r - n)$ equations:
\[
\varphi_{k\nu}(x_1,\dots,x_s)
-
\psi_{k\nu}(y_1,\dots,y_s)
=
0
\ \ \ \ \ \ \ \ \ \ \ \ \ 
{\scriptstyle{(k\,=\,1\,\cdots\,r\,-\,n\,;\,\,\,
\nu\,=\,1\,\cdots\,n)}}
\]
neither contradict themselves mutually, 
nor produce relations between the $x$ alone or between
the $y$ alone.\footnote[1]{\,
\name{Lie}, Archiv for Math. og Naturv. Vols. 3 and 4, 
Christiania 1878 and 1879; Math. Ann. Vol. XXV, pp.~96--107.
}
\end{theorem}
\renewcommand{\thefootnote}{\arabic{footnote}}

Moreover, we yet have the

\def\theproposition{2}\begin{proposition}
\label{Satz-2-S-354}
If the two $r$-term groups $X_1f, \dots, X_rf$ and $Z_1f, \dots, 
Z_rf$ are similar to each other and if the $r$
infinitesimal transformations:
\[
Y_kf
=
\sum_{j=1}^r\,\overline{g}_{kj}\,Z_jf
\ \ \ \ \ \ \ \ \ \ \ \ \ {\scriptstyle{(k\,=\,1\,\cdots\,r)}}
\]
are chosen as is indicated by Theorem~65, then 
there always exists at least one transformation: 
\[
y_1
=
\Phi_1(x_1,\dots,x_s),
\,\,\,\dots,\,\,\,
y_s
=
\Phi_s(x_1,\dots,x_s)
\]
which transfers $X_1f, \dots, X_rf$ to $Y_1f, \dots, Y_rf$,
respectively. One can set up each transformation of this nature as
soon as one has integrated certain complete systems. One obtains the
most general transformation which actually transfers the group $X_1f,
\dots, X_rf$ to the group $Z_1f, \dots, Z_rf$ by choosing the
constants $\overline{ g}_{ kj}$ in the $Y_kf$ in the most general way,
and afterwards, by seeking the most general transformation which
transfers $X_1f, \dots, X_rf$ to $Y_1f, \dots, Y_rf$, respectively.
\end{proposition}

The complete systems which are spoken of in this proposition
do not at all appear in one case, namely in the case
where the number $\rho$ defined earlier on is equal to zero, 
hence when amongst the $n ( r - n)$ functions
$\varphi_{ k\nu} (x_1, \dots, x_s)$, there are
exactly $s$ which are mutually independent. 
Indeed, we have already remarked on p.~\pageref{S-344}
that in this case, the equations:
\[
\varphi_{k\nu}(x_1,\dots,x_s)
-
\psi_{k\nu}(y_1,\dots,y_s)
=
0
\ \ \ \ \ \ \ \ \ \ \ \ \ 
{\scriptstyle{(k\,=\,1\,\cdots\,r\,-\,n\,;\,\,\,
\nu\,=\,1\,\cdots\,n)}}
\]
are solvable both with respect to $y_1, \dots, y_s$ and with respect
to $x_1, \dots, x_s$, and that they represent the most general
transformation which transfers $X_1f, \dots, X_rf$ to $Y_1f, \dots,
Y_rf$, respectively.

On the other hand, there are cases in which the integration of the
mentioned complete systems is executable, for instance, this is always
so when the finite equations of both groups $X_1f, \dots, X_rf$ and
$Y_1f, \dots, Y_rf$ are known; however, we cannot
get involved in this sort of questions.
\label{S-355} 

\sectionengellie{\S\,\,\,92.}

Let the two $r$-term groups $X_1f, \dots, X_rf$ and $Z_1f, \dots,
Z_rf$ be similar to each other, so that the transformations $X_1f,
\dots, X_rf$ convert into the infinitesimal transformations of the
other group.

If we now choose $r$ independent infinitesimal transformations: 
\[
Y_kf
=
\sum_{j=1}^r\,g_{kj}\,Z_jf
\ \ \ \ \ \ \ \ \ \ \ \ \ {\scriptstyle{(k\,=\,1\,\cdots\,r)}}
\]
in the group $Z_1f, \dots, Z_rf$, then according to what precedes,
there exists a transformation which transfers $X_1f, \dots, X_rf$ to
$Y_1f, \dots, Y_rf$, respectively, if and only if $Y_1f, \dots, Y_rf$
possess the properties indicated in Theorem~65.

We assume that $Y_1f, \dots, Y_rf$ satisfy this requirement and that
$y_i = \Phi_i (x_1, \dots, x_s)$ is a transformation which transfers
$X_1f, \dots, X_rf$ to $Y_1f, \dots, Y_rf$, respectively. Then by the
transformation $y_i = \Phi_i (x)$, the general infinitesimal
transformation: $e_1\, X_1f + \cdots + e_r\, X_rf$ receives the form:
$e_1\, Y_1f + \cdots + e_r\, Y_rf$, hence our transformation
associates to every infinitesimal transformation of the group $X_1f,
\dots, X_rf$ a completely determined infinitesimal transformation of
the group $Z_1f, \dots, Z_rf$, and conversely. The univalent and
invertible relationships which is established in this way between the
infinitesimal transformations of the two groups is, according to
p.~\pageref{S-330}, a holoedrically isomorphic one.

Now, let $x_1^0, \dots, x_s^0$ be a point in general position, i.e. a
point for which the transformations $X_1f, \dots, X_nf$ produce $n$
independent directions, so that the functions $\varphi_{ k \nu } (x)$
in the identities:
\def\theequation{3}\begin{equation}
X_{n+k}f
\equiv
\sum_{\nu=1}^n\,\varphi_{k\nu}(x_1,\dots,x_n)\,X_\nu f
\ \ \ \ \ \ \ \ \ \ \ \ \ {\scriptstyle{(k\,=\,1\,\cdots\,r\,-\,n)}}
\end{equation}
behave regularly. According to Chap.~\ref{kapitel-11},
p.~\pageref{S-203-ter}, the infinitesimal transformations $e_1\, X_1f
+ \cdots + e_r\, X_rf$ that leave invariant the point $x_1^0, \dots,
x_s^0$ have the form:
\def\theequation{27}\begin{equation}
\sum_{k=1}^{r-n}\,
\varepsilon_k\,
\bigg\{
X_{q+k}f
-
\sum_{\nu=1}^n\,\varphi_{k\nu}(x_1^0,\dots,x_s^0)\,X_\nu f
\bigg\},
\end{equation}
where it is understood that $\varepsilon_1, \dots, \varepsilon_{ r-n}$
are arbitrary parameters, and all these infinitesimal transformations
generate an $(r - n)$-term subgroup of the group $X_1f, \dots, X_rf$.

If, in the identities~\thetag{ 3}, we introduce the new variables
$y_1, \dots, y_s$ by means of the transformation: $y_i = \Phi_i ( x_1,
\dots, x_s)$, then we obtain between $Y_1f, \dots, Y_rf$ the
identities:
\def\theequation{3'}\begin{equation}
Y_{n+k}f
\equiv
\sum_{\nu=1}^n\,\psi_{k\nu}(y_1,\dots,y_s)\,Y_\nu f
\ \ \ \ \ \ \ \ \ \ \ \ \ {\scriptstyle{(k\,=\,1\,\cdots\,r\,-\,n)}}.
\end{equation}
Hence with $y_i^0 = \Phi_i ( x_1^0, \dots, x_s^0)$, after the
introduction of the variables $y$, the infinitesimal
transformation~\thetag{ 27} is transferred to:
\def\theequation{27'}\begin{equation}
\sum_{k=1}^{r-n}\,\varepsilon_k\,
\bigg\{
Y_{q+k}
-
\sum_{\nu=1}^n\,\psi_{k\nu}(y_1^0,\dots y_s^0)\,Y_\nu f
\bigg\},
\end{equation}
that is to say, to the most general infinitesimal transformation
$e_1\, Y_1f + \cdots + e_r\, Y_rf$ which leaves fixed the point:
$y_1^0, \dots, y_s^0$ in general position. At the same time, all
infinitesimal transformations of the form~\thetag{ 27'} naturally
generate an $(r-n)$-term subgroup of the group $Z_1f, \dots, Z_rf$.

In that, we have an important property of the holoedrically isomorphic
relationship which is established between the two groups by the
transformation: $y_i = \Phi_i ( x_1, \dots, x_s)$. Indeed, to the
most general subgroup of $X_1f, \dots, X_rf$ which leaves fixed an
arbitrarily chosen point: $x_1^0, \dots, x_s^0$ in general position,
this holoedrically isomorphic relationship always associates the most
general subgroup of $Y_1f, \dots, Y_rf$ which leaves fixed a point:
$y_1^0, \dots, y_s^0$ in general position. Exactly the same
association is found in the reverse direction; in other words: when
the point: $x_1^0, \dots, x_s^0$ runs through all possible positions,
then the point: $y_1^0, \dots, y_s^0$ also runs through all possible
positions.

\label{S-357-sq}
Now conversely, for two $r$-term groups in the same number of
variables to be similar to each other, then obviously, one must be
able to produce between them a holoedrically isomorphic relationship
of the constitution just described. We claim that this necessary
condition is at the same time also sufficient; we will show that the
two groups are really similar, when such a holoedrically isomorphic
relationship can be produced between them.

In fact, let $Z_1f, \dots, Z_rf$ be any $r$-term group which can be
related in a holoedrically isomorphic way to the group $X_1f, \dots,
X_rf$ in the said manner; let $e_1\, Y_1f + \cdots + e_r\, Y_rf$ be
the infinitesimal transformation of the group $Z_1f, \dots, Z_rf$
which is associated to the general infinitesimal transformation $e_1\,
X_1f + \cdots + e_r\, X_rf$ of the group $X_1f, \dots, X_rf$ through
the concerned holoedrically isomorphic relationship.

If $X_1f, \dots, X_nf$ are linked together by no linear relations,
while $X_{ n+1}f, \dots, X_rf$ can be expressed by means of $X_1f,
\dots, X_nf$ in the known way, then the most general infinitesimal
transformation contained in the group $X_1f, \dots, X_rf$ which leaves
invariant the point: $x_1^0, \dots, x_s^0$ in general position reads
as follows:
\def\theequation{27}\begin{equation}
\sum_{k=1}^{r-n}\,\varepsilon_k\,
\bigg\{
X_{n+k}f
-
\sum_{\nu=1}^n\,
\varphi_{k\nu}(x_1^0,\dots,x_s^0)\,X_\nu f
\bigg\}.
\end{equation}
Under the assumptions made, to it corresponds, in the group 
$Z_1f, \dots, Z_rf$, the infinitesimal transformation:
\def\theequation{28}\begin{equation}
\sum_{k=1}^{r-n}\,\varepsilon_k\,
\bigg\{
Y_{n+k}f
-
\sum_{\nu=1}^n\,
\varphi_{k\nu}(x_1^0,\dots,x_s^0)\,Y_\nu f
\bigg\},
\end{equation}
which now is in turn the most general transformation of the group
$Z_1f, \dots, Z_rf$ which leaves at rest a certain point: $y_1^0,
\dots, y_s^0$ in general position. Here, if the point $x_1^0, \dots,
x_s^0$ runs through all possible positions, then the point $y_1^0,
\dots, y_s^0$ does the same.

Since~\thetag{ 28} is the most general infinitesimal transformation
$e_1\, Y_1f + \cdots + e_r\, Y_rf$ which leaves invariant the point:
$y_1^0, \dots, y_s^0$ in general position, $Y_1f, \dots, Y_nf$ can be
linked together by no linear relation, and by contrast, $Y_{ n+1}f,
\dots, Y_rf$ must be expressible by means of $Y_1f, \dots, Y_nf$ in
the known way. From this, we deduce that the most general
infinitesimal transformation $e_1\, Y_1f + \cdots + e_r\, Y_rf$ which
leaves fixed the point: $y_1^0, \dots, y_s^0$ can also be represented
by the following expression:
\def\theequation{28'}\begin{equation}
\sum_{k=1}^{r-n}\,\varepsilon_k'\,
\bigg\{
Y_{n+k}f
-
\sum_{\nu=1}^n\,\psi_{k\nu}(y_1^0,\dots,y_s^0)\,Y_\nu f
\bigg\}.
\end{equation}

Evidently, every infinitesimal transformation contained in the
expression~\thetag{ 28} is identical to one of the infinitesimal
transformations~\thetag{ 28'}, so for arbitrarily chosen
$\varepsilon_1, \dots, \varepsilon_{ r -n}$, it must
always be possible to determine $\varepsilon_1', 
\dots, \varepsilon_{ r - n}'$ so that the equation:
\[
\sum_{k=1}^{r-n}\,
(\varepsilon_k-\varepsilon_k')\,Y_{n+k}f
-
\sum_{\nu=1}^n\,
\bigg\{
\sum_{k=1}^{r-n}\,
\big(
\varepsilon_k\,\varphi_{k\nu}(x^0)
-
\varepsilon_k'\,\psi_{k\nu}(y^0)
\big)
\bigg\}\,Y_\nu f
=
0
\]
is identically satisfied.

Because $Y_1f, \dots, Y_rf$ are independent infinitesimal
transformations, the equation just written decomposes in the following
ones:
\def\theequation{29}\begin{equation}
\aligned
\varepsilon_k-\varepsilon_k'
=
0,
&
\ \ \ \ \ \ \ \
\sum_{j=1}^{r-n}\,
\big(
\varepsilon_j\,\varphi_{j\nu}(x^0)
-
\varepsilon_j'\,\psi_{j\nu}(y^0)
\big)
=
0
\\
&
\ 
{\scriptstyle{(k\,=\,1\,\cdots\,r\,-\,n\,;\,\,\,
\nu\,=\,1\,\cdots\,n)}}.
\endaligned
\end{equation}

From this, it comes immediately:
\[
\varepsilon_1'=\varepsilon_1,
\,\,\,\dots,\,\,\,
\varepsilon_{r-n}'=\varepsilon_{r-n},
\]
and in addition, we obtain thanks to the arbitrariness of the 
$\varepsilon$, yet the equations:
\[
\varphi_{k\nu}(x_1^0,\dots,x_s^0)
-
\psi_{k\nu}(y_1^0,\dots,y_s^0)
=
0
\ \ \ \ \ \ \ \ \ \ \ \ \
{\scriptstyle{(k\,=\,1\,\cdots\,r\,-\,n\,;\,\,\,
\nu\,=\,1\,\cdots\,n)}},
\]
which must therefore, under the assumptions made, hold for the two
points $x_1^0, \dots, x_s^0$ and $y_1^0, \dots, y_s^0$. At present, if
we still remember that the point $y_1^0, \dots, y_s^0$ runs through
all possible positions, as soon as $x_1^0, \dots, x_s^0$ does this,
then we realize immediately that, under the assumptions made, the $n (
r - n)$ equations:
\def\theequation{4}\begin{equation}
\varphi_{k\nu}(x_1,\dots,x_s)
-
\psi_{k\nu}(y_1,\dots,y_s)
=
0
\ \ \ \ \ \ \ \ \ \ \ \ \
{\scriptstyle{(k\,=\,1\,\cdots\,r\,-\,n\,;\,\,\,
\nu\,=\,1\,\cdots\,n)}}
\end{equation}
are compatible with each other and produce relations neither between
the $x$ alone, not between the $y$ alone.

According to Theorem~65, p.~\pageref{Theorem-65-S-353}, it is thus
proved that the two groups $X_1f, \dots, X_rf$ and $Z_1f, \dots, Z_rf$
are similar to each other, and this is just what we wanted to prove.

We therefore have the

\def\theproposition{3}\begin{proposition}
\label{Satz-3-S-359}
Two $r$-term groups $G$ and $\Gamma$ in the same number of variables
are similar to each other if and only if it is possible to relate them
in a holoedrically isomorphic way so that the most general subgroup of
$G$ which leaves invariant a determined point in general position
always corresponds, in whichever way the point may be chosen, to the
most general subgroup of $\Gamma$ which leaves invariant a certain
point in general position, and so that the same correspondence also
holds in the reverse direction.
\end{proposition}

But at the same time, it is also proved that, under the assumption
made right now, there exists a transformation: $y_i = \Phi_i ( x_1,
\dots, x_s)$ which transfers $X_1f, \dots, X_rf$ to $Y_1f, \dots,
Y_rf$, respectively. We can therefore also state the following
somewhat more specific proposition:

\def\theproposition{4}\begin{proposition}
If the $r$ independent infinitesimal transformations:
\[
X_kf
=
\sum_{i=1}^s\,\xi_{ki}(x_1,\dots,x_s)\,
\frac{\partial f}{\partial x_i}
\ \ \ \ \ \ \ \ \ \ \ \ \ {\scriptstyle{(k\,=\,1\,\cdots\,r)}}
\]
generate an $r$-term group $G$, and if the $r$ independent
infinitesimal transformations:
\[
Y_kf
=
\sum_{i=1}^s\,\eta_{ki}(y_1,\dots,y_s)\,
\frac{\partial f}{\partial y_i}
\ \ \ \ \ \ \ \ \ \ \ \ \ {\scriptstyle{(k\,=\,1\,\cdots\,r)}}
\]
generate an $r$-term group $\Gamma$, then there is a transformation:
$y_i = \Phi_i ( x_1, \dots,x_s)$ which transfers $X_1f, \dots, X_rf$
to $Y_1f, \dots, Y_rf$, respectively, if and only if the following
conditions are satisfied:

\terminology{Firstly:} if $G$ contains exactly $r - n$ independent
infinitesimal transformations which leave invariant an arbitrarily
chosen point in general position, then $\Gamma$ must
also contain exactly $r - n$ independent infinitesimal
transformations having this constitution. 

\terminology{Secondly:} if one associates to every infinitesimal
transformation $e_1\, X_1f + \cdots + e_r\, X_rf$ of $G$ the
infinitesimal transformation $e_1\, Y_1f + \cdots + e_r\, Y_rf$ of
$\Gamma$, then the two groups must be related to each other in a
holoedrically isomorphic way, in the manner indicated by the previous
proposition.
\end{proposition}

Because the most general subgroup of $G$ which leaves invariant a
determined point in general position is completely defined by this
point, and moreover, because the holoedrically isomorphic relationship
between $G$ and $\Gamma$ mentioned several times associates to every
subgroup of $G$ of this kind a subgroup of $\Gamma$ constituted in the
same way, it follows that this relationship between $G$ and $\Gamma$
also establishes a correspondence between the points $x_1, \dots, x_s$
and the points $y_1, \dots, y_s$; however, this correspondence is in
general infinitely multivalent, for to every point $x_1, \dots, x_s$
there obviously correspond all points $y_1, \dots, y_s$ which satisfy
the equations~\thetag{ 4}, and inversely. As a result, this agrees
with the fact that there are in general infinitely many
transformations which transfer $X_1f, \dots, X_rf$ to $Y_1f, \dots,
Y_rf$, respectively.

For \emphasis{transitive} groups, the criterion of similarity
enunciated in Proposition~3 turns out to be particularly simple. We
shall see later (Chap.~\ref{kapitel-21}) that two $r$-term transitive
groups $G$ and $\Gamma$ in the same number, say $s$, of variables are
already similar to each other when it is possible to relate them in a
holoedrically isomorphic way so that, to a single $(r - s)$-term
subgroup of $G$ which leaves fixed one point in general position there
corresponds an $(r - s)$-term subgroup of $\Gamma$ having the same
constitution.

\medskip

If an $r$-term group:
\[
X_kf
=
\sum_{i=1}^s\,\xi_{ki}(x_1,\dots,x_s)\,
\frac{\partial f}{\partial x_i}
\ \ \ \ \ \ \ \ \ \ \ \ \ {\scriptstyle{(k\,=\,1\,\cdots\,r)}}
\]
having the composition:
\[
\leftbracket
X_i,\,X_k
\rightbracket
=
\sum_{\sigma=1}^r\,c_{ik\sigma}\,X_\sigma f
\]
is presented, then one can ask for all transformations:
\[
x_i'
=
\Phi_i(x_1,\dots,x_s)
\ \ \ \ \ \ \ \ \ \ \ \ \ {\scriptstyle{(i\,=\,1\,\cdots\,s)}}
\]
that leave it invariant, namely one can ask for all transformations
through which the group is similar to itself.

Thanks to the developments of the previous paragraphs, we are in a
position to determine all the transformations in question.

To begin with, we relate in the most general holoedrically isomorphic
way the group $X_1f, \dots, X_rf$ to itself, hence we choose in the
most general way $r$ independent infinitesimal transformations:
\[
\Xi_kf
=
\sum_{j=1}^r\,g_{kj}\,X_jf
\ \ \ \ \ \ \ \ \ \ \ \ \ {\scriptstyle{(k\,=\,1\,\cdots\,r)}}
\]
that stand pairwise in the relationships:
\[
\leftbracket
\Xi_i,\,\Xi_k
\rightbracket
=
\sum_{\sigma=1}^r\,c_{ik\sigma}\,\Xi_\sigma f.
\]
Afterwards, we specialize the arbitrary elements contained in the $g_{
kj}$ in such a way that the following conditions are satisfied: when
there are no relations between $X_1f, \dots, X_nf$, while $X_{ n+1}f,
\dots, X_rf$ express themselves by means of $X_1f, \dots, X_nf$:
\[
X_{n+k}f
\equiv
\sum_{\nu=1}^n\,
\varphi_{k\nu}(x_1,\dots,x_s)\,X_\nu f
\ \ \ \ \ \ \ \ \ \ \ \ \ {\scriptstyle{(k\,=\,1\,\cdots\,r\,-\,n)}}, 
\]
then firstly, also $\Xi_1 f, \dots, \Xi_n f$ should be linked together
by no linear relation, while by contrast $\Xi_{ n+1} f, \dots, \Xi_rf$
also express themselves by means of $\Xi_1f, \dots, \Xi_n f$:
\[
\Xi_{n+k}f
\equiv
\sum_{\nu=1}^n\,\psi_{k\nu}(x_1,\dots,x_s)\,\Xi_\nu f
\ \ \ \ \ \ \ \ \ \ \ \ \ {\scriptstyle{(k\,=\,1\,\cdots\,r\,-\,n)}},
\]
and secondly, the $n ( r - n)$ equations: 
\[
\varphi_{k\nu}(x_1',\dots,x_s')
=
\psi_{k\nu}(x_1,\dots,x_s)
\ \ \ \ \ \ \ \ \ \ \ \ \ 
{\scriptstyle{(k\,=\,1\,\cdots\,r\,-\,n\,;\,\,\,
\nu\,=\,1\,\cdots\,n)}}
\]
should be mutually compatible, and they should produce relations
neither between the $x$ alone, nor between the $x'$ alone.

If all of this is realized, then following the introduction of
\S\,\,91, we determine the most general transformation $x_i' = \Phi_i
( x_1, \dots, x_s)$ which transfers $\Xi_1 f, \dots, \Xi_rf$ to
$X_1'f, \dots, X_r'f$, respectively, where:
\[
X_k'f
=
\sum_{i=1}^r\,\xi_{ki}(x_1',\dots,x_s')\,
\frac{\partial f}{\partial x_i'}.
\]
The transformation in question is then the most general one by which
the group $X_1f, \dots, X_rf$ is transferred to itself.

It is clear that the totality of all transformations which leave
invariant
\label{S-361} 
the group $X_1f, \dots, X_rf$ forms a group by itself,
namely \emphasis{the largest group in which $X_1f, \dots, X_rf$ is
contained as an invariant subgroup}. This group can be finite or
infinite, continuous or not continuous, but in all circumstances its
transformations order as inverses by pairs, because when a
transformation $x_i' = \Phi_i ( x_1, \dots, x_s)$ leaves the group
$X_1f, \dots, X_rf$ invariant, then the associated inverse
transformation also does this.

If by chance the just defined group consists of a finite number of
different families of transformations and if at the same time each one
of these families contains only a finite number of arbitrary
parameters, then according to Chap.~\ref{kapitel-18},
p.~\pageref{S-315-sq} sq., there is in the concerned group a family of
transformations which constitutes a finite continuous group; then this
family is the largest continuous group in which $X_1f, \dots, X_rf$ is
contained as an invariant subgroup.

Up to now, we have spoken of similarity only for groups which contain
the same number of variables. But about that, it was not excluded
that certain of the variables were absolutely not transformed by the
concerned groups, so that they actually did not appear in the
infinitesimal transformations.

Now, for the groups which do not contain the same number of variables,
one can also speak of similarity; indeed, one can always complete the
number of variables in one of the groups so that one adds a necessary
number of variables that are absolutely not transformed by the
concerned group. Then one has two groups in the same number of
variables and one can examine whether they are similar to each other,
or not.

In the sequel, unless the contrary is expressly stressed, we shall
actually interpret the concept of similarity in the original, narrower
sense.

\sectionengellie{\S\,\,\,93.}

In order to illustrate the general theory of similarity by an example,
we will examine whether the two three-term groups in two independent
variables:
\[
\aligned
X_1f
&
=
\frac{\partial f}{\partial x_1},
\ \ \ \ \ \ \ \
X_2f
=
x_1\,\frac{\partial f}{\partial x_1}
+
x_2\,\frac{\partial f}{\partial x_2},
\\
X_3f
&
=
x_1^2\,\frac{\partial f}{\partial x_1}
+
(2\,x_1x_2+C\,x_2^2)\,
\frac{\partial f}{\partial x_2}
\endaligned
\]
and:
\[
\aligned
Y_1f
&
=
\frac{\partial f}{\partial y_1}
+
\frac{\partial f}{\partial y_2},
\ \ \ \ \ \ \ \ \ \
Y_2f
=
y_1\,\frac{\partial f}{\partial y_1}
+
y_2\,\frac{\partial f}{\partial y_2},
\\
&
Y_3f
=
y_1^2\,\frac{\partial f}{\partial y_1}
+
y_2^2\,\frac{\partial f}{\partial y_2}
\endaligned
\]
are similar to each other.

As one has remarked, the $Y_kf$ here are already chosen
in such a way that one has at the same time:
\[
\leftbracket
X_1,\,X_2
\rightbracket
=
X_1f,
\ \ \ \ \ \ \
\leftbracket
X_1,\,X_3
\rightbracket
=
2\,X_2f,
\ \ \ \ \ \ \
\leftbracket
X_2,\,X_3
\rightbracket
=
X_3f
\]
and:
\[
\leftbracket
Y_1,\,Y_2
\rightbracket
=
Y_1f,
\ \ \ \ \ \ \
\leftbracket
Y_1,\,Y_3
\rightbracket
=
2\,Y_2f,
\ \ \ \ \ \ \
\leftbracket
Y_2,\,Y_3
\rightbracket
=
Y_3f
\]
Although the $Y_kf$ are not chosen in the most general
way so that the shown relations hold, it is nevertheless
not necessary to do this, because the result would
not be modified by this. 

We find:
\[
X_3f
\equiv
-\,(x_1^2+C\,x_1x_2)\,X_1f
+
(2x_1+C\,x_2)\,X_2f
\]
and:
\[
Y_3f
\equiv
-\,y_1y_2\,Y_1f
+
(y_1+y_2)\,Y_2f,
\]
whence one must have: 
\[
y_1y_2
=
x_1^2+C\,x_1x_2,
\ \ \ \ \ \ \ \
y_1+y_2=2\,x_1+C\,x_2.
\]
As long as the constant $C$ does not vanish, these equations determine
a transformation, and consequently, according to our general theory,
the two groups are similar to each other in the case $C \neq 0$. But
if $C = 0$, it comes a relation between $y_1$ and $y_2$ alone, hence
in any case, there exists no transformation which transfers $X_1f$,
$X_2f$, $X_3f$ to $Y_1f$, $Y_2f$, $Y_3f$, respectively; one can easily
convince oneself that in this case, the two groups are actually not
similar to each other.

\sectionengellie{\S\,\,\,94.}

Subsequently to the theory of the similarity of $r$-term groups, we
yet want to briefly treat a somewhat more general question and to
indicate its solution.

We imagine, in the variables $x_1, \dots, x_s$, that any $p$
infinitesimal transformations: $X_1f, \dots, X_pf$ are presented,
hence not exactly some which generate a finite group, and likewise, we
imagine in $y_1, \dots, y_s$ that any $p$ infinitesimal
transformations $Y_1f, \dots, Y_pf$ are presented. We ask under which
conditions there is a transformation $y_i = \Phi_i ( x_1, \dots, x_s)$
which transfers $X_1f, \dots, X_pf$ to $Y_1f, \dots, Y_pf$,
respectively. Here, we do not demand that $X_1f, \dots, X_pf$ should
be independent infinitesimal transformations, since this assumption
would in fact not interfere with the generality of the considerations
following, but it would complicate the presentation, because it would
always have to be taken into account.

At first, we can lead back the above general problem to the special
case case where the independent equations amongst the equations: $X_1f
= 0$, \dots, $X_p f = 0$ form a complete system.

Indeed, if there is a transformation which transfers $X_1f, \dots,
X_pf$ to $Y_1f, \dots, Y_pf$, respectively, then according to
Chap.~\ref{one-term-groups}, p.~\pageref{Satz-2-S-84}, every
expression:
\[
X_k\big(X_j(f)\big)
-
X_j\big(X_k(f)\big)
=
\leftbracket
X_k,\,X_j
\rightbracket
\]
also converts into the corresponding expression:
\[
Y_k\big(Y_j(f)\big)
-
Y_j\big(Y_k(f)\big)
=
\leftbracket
Y_k,\,Y_j
\rightbracket.
\]
Hence, if the independent equations amongst the equations $X_1f = 0$,
\dots, $X_pf = 0$ do actually not form already a complete system, then
to $X_1f, \dots, X_pf$, we can yet add all the expressions
$\leftbracket X_k,\, X_j \rightbracket$, and also, we must only add to
$Y_1f, \dots, Y_pf$ all the expressions $\leftbracket Y_k,\, Y_j
\rightbracket$. The question whether there exists a transformation of
the demanded constitution then amounts to the question whether there
exists a transformation with transfers $X_1f, \dots, X_pf$,
$\leftbracket X_k, \, X_j \rightbracket$ to $Y_1f, \dots, Y_pf$,
$\leftbracket Y_k,\, Y_j \rightbracket$, respectively, where one has
to set for $k$ and $j$ the numbers $1, 2, \dots, p$ one after the
other.

Now, if the independent equations amongst the equations: $X_1f = 0$,
\dots, $X_p f = 0$, $\leftbracket X_k,\,X_j \rightbracket = 0$ also do
not form a complete system, then we yet add the expressions
$\big\leftbracket X_i,\, \leftbracket X_k,\, X_j \rightbracket
\big\rightbracket$ and the $\big\leftbracket \leftbracket X_i,\,X_k
\rightbracket,\, \leftbracket X_j,\,X_l \rightbracket
\big\rightbracket$, and the corresponding expressions in the $Yf$ as
well. If we continue in this way, we obtain at the end that our
initial problem is lead back to the following one:

\plainstatement{In the variables $x_1, \dots, x_s$, let $r$
infinitesimal transformations: $X_1f, \dots, X_rf$ be presented, of
which $n \leqslant r$, say $X_1f, \dots, X_nf$, are linked together by
no linear relation, while $X_{ n+1}f, \dots, X_rf$ can be linearly
expressed in terms of $X_1f, \dots, X_nf$:
\[
X_{n+k}f
\equiv
\sum_{\nu=1}^n\,\varphi_{k\nu}(x_1,\dots,x_s)\,X_\nu f
\ \ \ \ \ \ \ \ \ \ \ \ \ {\scriptstyle{(k\,=\,1\,\cdots\,n\,-\,r)}}\,;
\]
in addition, let relations of the form:
\[
\leftbracket
X_k,\,X_j
\rightbracket
=
\sum_{\nu=1}^n\,
\varphi_{kj\nu}(x_1,\dots,x_s)\,X_\nu f
\ \ \ \ \ \ \ \ \ \ \ \ \ {\scriptstyle{(k,\,\,j\,=\,1\,\cdots\,r)}}
\]
yet hold, so that the independent equations amongst the equations
$X_1f = 0$, \dots, $X_rf = 0$ form an $n$-term complete system.
Moreover, in the variables $y_1, \dots, y_s$, let $r$ infinitesimal
transformations $Y_1f, \dots, Y_rf$ be presented. To determine whether
there is a transformation $y_i = \Phi_i ( x_1, \dots, x_s)$ which
transfers $X_1f, \dots, X_rf$ to $Y_1 f, \dots, Y_rf$, respectively.}

For a transformation of the constitution demanded here to exist,
$Y_1f, \dots, Y_rf$ should naturally not be linked together by a
linear relation, while $Y_{ n+1}f, \dots, Y_rf$ should be linearly
expressible in terms of $Y_1f, \dots, Y_nf$:
\[
Y_{n+k}f
\equiv
\sum_{\nu=1}^n\,\psi_{k\nu}(y_1,\dots,y_s)\,Y_\nu f
\ \ \ \ \ \ \ \ \ \ \ \ \ {\scriptstyle{(k\,=\,1\,\cdots\,r\,-\,n)}}, 
\]
and there should be relations of the form: 
\[
\leftbracket
Y_k,\,Y_j
\rightbracket
=
\sum_{\nu=1}^n\,\psi_{kj\nu}(y_1,\dots,y_s)\,Y_\nu f
\ \ \ \ \ \ \ \ \ \ \ \ \ {\scriptstyle{(k,\,\,j\,=\,1\,\cdots\,r)}}.
\]
In addition, the equations:
\def\theequation{30}\begin{equation}
\left\{
\aligned
\varphi_{k\nu}(x_1,\dots,x_s)
-
\psi_{k\nu}(y_1,\dots,y_s)
&
=
0
\ \ \ \ \ \ \ \ \ \ \ \ \ 
{\scriptstyle{(k\,=\,1\,\cdots\,r\,-\,n\,;\,\,\,
\nu\,=\,1\,\cdots\,n)}}
\\
\varphi_{kj\nu}(x_1,\dots,x_s)
-
\psi_{kj\nu}(y_1,\dots,y_s)
&
=
0
\ \ \ \ \ \ \ \ \ \ \ \ \ 
{\scriptstyle{(k,\,\,j\,=\,1\,\cdots\,r\,;\,\,\,
\nu\,=\,1\,\cdots\,n)}}
\endaligned\right.
\end{equation}
should yet neither contradict mutually, nor conduct to relations
between the $x$ alone or the $y$ alone, since these equations will
obviously reduce to identities after the substitution: $y_i = \Phi_i (
x_1, \dots, x_s)$, when the transformation: $y_i = \Phi_i ( x)$
converts $X_1f, \dots, X_rf$ into $Y_1f, \dots, Y_rf$, respectively.

We want to assume that all these conditions are satisfied
and that all the equations~\thetag{ 30} reduce to the
$\rho$ mutually independent equations:
\def\theequation{31}\begin{equation}
\varphi_1(x)-\psi_1(x)=0,
\,\,\,\dots,\,\,\,
\varphi_\rho(x)-\psi_\rho(x)=0.
\end{equation}

Thanks to considerations completely similar to those of
p.~\pageref{S-335-sq} sq., we realize that the determination of a
transformation which transfers $X_1f, \dots, X_rf$ to $Y_1f, \dots,
Y_rf$, respectively, amounts to determining a system of equations in
the $2s$ variables $x_1, \dots, x_s$, $y_1, \dots, y_s$, and to be
precise, a system of equations having the following constitution: it
must admit the $r$ infinitesimal transformations: $\Omega_k f = X_k f
+ Y_k f$, it must consist of exactly $s$ independent equations, it
must be solvable both with respect to $x_1, \dots, x_s$ and with
respect to $y_1, \dots, y_s$, and lastly, it must comprise the $\rho$
equations~\thetag{ 31}.

A system of equations which comprises the $\rho$ equations~\thetag{ 31}
and which admits the $r$ infinitesimal transformations
$\Omega_kf$ embraces at the same time
all the $r \, \rho$ equations:
\[
\aligned
\Omega_k
&
\big(\varphi_j(x)-\psi_j(x)\big)
=
X_k\,\varphi_j(x)
-
Y_k\,\psi_j(y)
=
0
\\
&
\ \ \ \ \ \ \ \ \ \ \ \ \ \ \ \ \
{\scriptstyle{(k\,=\,1\,\cdots\,r\,;\,\,\,
j\,=\,1\,\cdots\,\rho)}}.
\endaligned
\]
If these equations are not a consequence of~\thetag{ 31}, we can again
deduce from them new equations which must be contained in the sought
system of equations, and so on. If we proceed 
in this way, then at the end, we
must come either to relations which contradict each other, or to
relations between the $x$ alone, or to to relations between the $y$
alone, or lastly, to a system of $\sigma \leqslant s$ independent
equations:
\def\theequation{32}\begin{equation}
\varphi_1(x)-\psi_1(y)=0,
\,\,\,\dots,\,\,\,
\varphi_\sigma(x)-\psi_\sigma(y)
=
0
\ \ \ \ \ \ \ \ \ \ \ \ \ 
{\scriptstyle{(\sigma\,\geqslant\,\rho)}}
\end{equation}
which possesses the following two properties: it produces no relation
between the $x$ or the $y$ alone, and it admits the $r$ infinitesimal
transformations $\Omega_k f$, so that each one of the $r \, \sigma$
equations:
\[
\Omega_k\big(\varphi_j(x)-\psi_j(y)\big)
=
0
\ \ \ \ \ \ \ \ \ \ \ \ \ \ \ \ \
{\scriptstyle{(k\,=\,1\,\cdots\,r\,;\,\,\,
j\,=\,1\,\cdots\,\sigma)}}.
\]
is a consequence of~\thetag{ 32}. 

Evidently, there can be a transformation which transfers $X_1f, \dots,
X_rf$ to $Y_1f, \dots, Y_rf$, respectively, only when we are conducted
to a system of equations~\thetag{ 32} having the constitution defined
just now by means of the indicated operations. So we need to consider
only this case.

If the entire number $\sigma$ is precisely equal to $s$, then the
system of equations~\thetag{ 32} taken for itself represents a
transformation which achieves the demanded transfer and to be precise,
it is obviously the only transformation which does this. We will show
that also in the case $\sigma < s$, a transformation exists which
transfers $X_1f, \dots, X_rf$ to $Y_1f, \dots, Y_rf$, respectively;
here, we produce the proof of that by indicating a method which
conducts to the determination of a transformation having the demanded
constitution.

It is clear that the equations~\thetag{ 32} neither cancel the
independence of the equations: $X_1f = 0$, \dots, $X_nf = 0$, nor they
cancel the independence of the equations: $Y_1f = 0$, \dots, $Y_n f =
0$, for they produce relations neither between the $x$ alone,
nor between the $y$ alone. 

Besides, it is to be observed that relations of the form:
\[
\aligned
\leftbracket
\Omega_k,\,\Omega_j
\rightbracket
&
=
\sum_{\nu=1}^n\,\varphi_{kj\nu}(x)\,\Omega_\nu f
=
\sum_{\nu=1}^n\,\psi_{kj\nu}(y)\,\Omega_\nu f
\\
&
\ \ \ \ \ \ \ \ \ \ \ \ \ 
{\scriptstyle{(k,\,\,j\,=\,1,\,\,2\,\cdots\,n)}}
\endaligned
\]
hold, in which the coefficients $\varphi_{ kj \nu} (x) = \psi_{ kj\nu}
(y)$ behave regularly in general for the systems of values of the
system of equations: $\varphi_1 - \psi_1 = 0$, \dots, $\varphi_\sigma
- \psi_\sigma = 0$. Thus, the case settled in Theorem~19,
p.~\pageref{Theorem-19-S-132} is present here.

As in p.~\pageref{S-344-sq} sq., we introduce the solutions of the
$n$-term complete system $X_1f = 0$, \dots, $X_nf = 0$ as new $x$ and
those of the complete system $Y_1f = 0$, \dots, $Y_n f = 0$ as new
$y$, and on the occasion, exactly as we did at that time, we have to
make a distinction between the solutions which can be expressed in
terms of $\varphi_1 (x), \dots, \varphi_\sigma (x)$ or, respectively,
in terms of $\psi_1 ( y), \dots, \psi_\sigma (y)$, and the solutions
which are independent of the $\varphi$, or respectively, of the
$\psi$. In this way, we simplify the shape of the equations~\thetag{
32}, and then, exactly as on p.~\pageref{S-349-sq} sq., we can
determine a system of equations which admits $\Omega_1f, \dots, 
\Omega_rf$, which contains exactly $s$ independent equations, 
which is solvable both with respect to $x_1, \dots, x_s$ 
and with respect to $y_1, \dots, y_s$, and lastly, which
embraces the equations~\thetag{ 32}. 
Obviously, the obtained system of equations represents
a transformation which transfers $X_1f, \dots, X_rf$
to $Y_1f, \dots, Y_rf$, respectively. 

Thus, we have the following Theorem.

\renewcommand{\thefootnote}{\fnsymbol{footnote}}
\def\thetheorem{66}\begin{theorem}
If, in the variables $x_1, \dots, x_s$, $p$ infinitesimal
transformations $X_1f, \dots, X_pf$ are presented and if, in the
variables $y_1, \dots, y_s$, $p$ infinitesimal transformations $Y_1f,
\dots, Y_pf$ are also presented, then one can always decide, by means
of differentiations and of eliminations, whether there is a
transformation $y_i = \Phi_i ( x_1, \dots, x_s)$ which transfers
$X_1f, \dots, X_rf$ to $Y_1f , \dots, Y_rf$, respectively;
if there is such a transformation, then one can 
determine the most general transformation
which accomplishes the concerned transfer, as
soon as one has integrated certain complete systems.\footnote[1]{\,
\name{Lie}, Archiv for Math. og Naturv., Vol. 3, 
p.~125, Christiania 1878.
}
\end{theorem}
\renewcommand{\thefootnote}{\arabic{footnote}}

\linestop


\chapter{Groups, the Transformations of Which
\\
Are Interchangeable With All Transformations
\\
of a Given Group
}
\label{kapitel-20}
\chaptermark{Transformations Which Are Interchangeable With All
Transformations of a Group}

\setcounter{footnote}{0}

\abstract*{??}

Thanks to the developments of the \S\S\,\,89, 90, 91,
pp.~\pageref{S-330-bis} up to \pageref{S-355}, we are in the position
to determine all transformations which leave invariant a given
$r$-term group:
\[
X_kf
=
\sum_{i=1}^s\,\xi_{ki}(x_1,\dots,x_s)\,
\frac{\partial f}{\partial x_i}
\ \ \ \ \ \ \ \ \ \ \ \ \ {\scriptstyle{(k\,=\,1\,\cdots\,r)}}.
\]
Amongst all transformations of this nature, we now want yet to pick
out those which in addition possess the property \emphasis{of leaving
invariant every individual transformation of the group $X_1f, \dots,
X_rf$}, and we want to occupy ourselves somehow more closely with
these transformations.

When a transformation $T$ leaves invariant every individual 
transformation $S$ of the group $X_1f, \dots, X_rf$, then
according to p.~\pageref{S-258}, it stands with respect
to $S$ in the relationship:
\[
T^{-1}\,S\,T
=
S,
\]
or, what is the same, in the relationship:
\[
S\,T
=
TS,
\]
hence it is interchangeable with $T$. Thus, we can also characterize
the transformations just defined in the following way: \emphasis{they
are the transformations which are interchangeable with all
transformations of the group $X_1f, \dots, X_rf$}.

\sectionengellie{\S\,\,\,95.}

According to Chap.~\ref{kapitel-15}, p.~\pageref{S-255}, the
expression $e_1\, X_1f + \cdots + e_r\, X_rf$ can be regarded as the
general symbol of a transformation of the group $X_1f, \dots, X_rf$.
Consequently, a transformation $x_i' = \Phi_i ( x_1, \dots, x_s)$ will
leave invariant every individual transformation of the group $X_1f,
\dots, X_rf$ when, for arbitrary choice of the $e$, it leaves
invariant the expression $e_1 \, X_1f + \cdots + e_r\, X_rf$, hence
when the expression $e_1 \, X_1f + \cdots + e_r\, X_rf$ takes the
form: $e_1\, X_1'f + \cdots + e_r \, X_r'f$ after the introduction of
the new variables $x_i' = \Phi_i ( x)$, where:
\[
X_k'f
=
\sum_{i=1}^s\,
\xi_{ki}(x_1',\dots,x_s')\,
\frac{\partial f}{\partial x_i'}.
\]
Here, for this to hold, it is necessary and sufficient that the $r$
infinitesimal transformations: $X_1f, \dots, X_rf$ receive the form
$X_1'f, \dots, X_r'f$, respectively, after the introduction of the new
variables $x_i'$.

It is clear that there are transformations $x_i' = \Phi_i (x)$ having
the constitution demanded; indeed, the identity transformation $x_i' =
x_i$ is such a transformation; besides, this results from the
Proposition~2 of the preceding chapter (p.~\pageref{Satz-2-S-354}),
because this proposition shows that there are transformations which
transfer $X_1f, \dots, X_rf$ to $X_1'f, \dots, X_r'f$, respectively.
In addition, it follows from the developments achieved at that time
that the most general transformation of the demanded constitution is
represented by the most general system of equations:
\[
x_1'
=
\Phi_1(x_1,\dots,x_s),
\,\,\,\dots,\,\,\,
x_s'
=
\Phi_s(x_1,\dots,x_s)
\]
which admits the $r$-term group: $X_1f + X_1'f$, \dots, $X_rf + X_r'f$
in the $2\, s$ variables $x_1, \dots, x_s$, $x_1', \dots, x_s'$.

If one executes, one after the other, two transformations which leave
invariant every individual transformation of the group: $X_1f, \dots,
X_rf$, then obviously, one always obtains a transformation which does
the same; consequently, 
\label{S-368} the totality of all transformations of this
constitution forms a group $G$. This group can be discontinuous, and
it can even reduce to the identity transformation; it can consist of
several discrete families of which each one contains only a finite
number of arbitrary parameters, it can be infinite; but 
its transformations are always ordered as inverses by pairs, 
since if a transformation leaves invariant
all transformations of the group: $X_1f, \dots, X_rf$, 
then the associated inverse transformation
naturally possesses the same property. 

\smallercharacters{
If the just defined group $G$ contains only a finite number of
arbitrary parameters, then it belongs to the category of groups which
was discussed in Chap.~\ref{kapitel-18}, and according to Theorem~56,
p.~\pageref{Theorem-56-S-315}, it certainly comprises a finite
continuous subgroups generated by infinitesimal transformations. On
the other hand, if the group $G$ is infinite, then thanks to
considerations similar to those of Chap.~\ref{kapitel-18}, it can be
proved that it comprises one-term groups, and in fact, infinitely many
such groups that are generated by infinitely many independent
infinitesimal transformations.

}

At present, we take up directly the problem of determining 
all one-term groups which are contained in the group $G$. 

According to Chap.~\ref{kapitel-15}, p.~\pageref{S-258}
and~\ref{S-259}, the $r$ expressions
$X_1f, \dots, X_rf$ remain invariant by all transformations
of the one-term group $Zf$ when the $r$ relations:
\[
\leftbracket
X_k,\,Z
\rightbracket
=
X_k\big(Z(f)\big)
-
Z\big(X_k(f)\big)
=
0
\ \ \ \ \ \ \ \ \ \ \ \ \ {\scriptstyle{(k\,=\,1\,\cdots\,r)}}
\] 
are identically satisfied, hence when the infinitesimal
transformation 
$Zf$ is interchangeable with all transformations
of the group $X_1f, \dots, X_rf$. Thus, the determination
of all one-term groups having the constitution
defined a short while ago amounts to the
determination of the most general 
infinitesimal transformation $Zf$ which is
interchangeable with all $X_kf$. 

If one has two infinitesimal transformations
$Z_1f$ and $Z_2f$ which are interchangeable
with all $X_kf$, then all
expressions
$\leftbracket Z_1, \, X_k \rightbracket$
and $\leftbracket Z_2,\, X_k \rightbracket$
vanish identically, whence the
Jacobi identity:
\[
\big\leftbracket
\leftbracket
Z_1,\,Z_2\rightbracket,\,X_k
\big\rightbracket
+
\big\leftbracket
\leftbracket
Z_2,\,X_k\rightbracket,\,Z_1
\big\rightbracket
+
\big\leftbracket
\leftbracket
X_k,\,Z_1\rightbracket,\,Z_2
\big\rightbracket
\equiv
0
\]
reduces to:
\[
\big\leftbracket
\leftbracket
Z_1,\,Z_2\rightbracket,\,X_k
\big\rightbracket
\equiv
0.
\]

Thus, the following holds:

\def\theproposition{1}\begin{proposition}
\label{Satz-1-S-369}
If the two infinitesimal transformations $Z_1f$ and $Z_2f$ are
interchangeable with all infinitesimal transformations of the $r$-term
group $X_1f, \dots, X_rf$, then the transformation $\leftbracket Z_1,
\, Z_2 \rightbracket$ also is so.
\end{proposition}

On the other hand, every infinitesimal transformation $a\, Z_1f + b\,
Z_2f$ is at the same time interchangeable with $X_1f, \dots, X_rf$,
whichever values the constants $a$ and $b$ may have. Hence if by
chance there is only a finite number, say $q$, of independent
infinitesimal transformations $Z_1f, \dots, Z_qf$ which are
interchangeable with $X_1f, \dots, X_rf$, then the most general
infinitesimal transformation having the same constitution has the
form: $\lambda_1 \, Z_1f + \cdots + \lambda_q \, Z_qf$, where it is
understood that $\lambda_1, \dots, \lambda_q$ are arbitrary
parameters. Then because of Proposition~1, there must exist relations
of the form:
\[
\leftbracket
Z_i,\,Z_k
\rightbracket
=
\sum_{\sigma=1}^q\,c_{ik\sigma}'\,Z_\sigma f, 
\]
so that \emphasis{$Z_1f, \dots, Z_q f$ 
generate a $q$-term group}. 

\medskip

At present, we seek to determine directly the most general
infinitesimal transformation:
\[
Zf
=
\sum_{i=1}^s\,\zeta_i(x_1,\dots,x_s)\,
\frac{\partial f}{\partial x_i}
\]
which is interchangeable with all infinitesimal 
transformations of the group $X_1f, \dots, X_rf$. 

The $r$ condition-equations
\deutsch{Bedingungsgleichungen}:
\[
\leftbracket
X_1,\,Z
\rightbracket
=
0,
\,\,\,\dots,\,\,\,
\leftbracket
X_r,\,Z
\rightbracket
=
0
\]
decompose immediately in the following $r\, s$ equations:
\[
X_k\,\zeta_i
=
Z\,\xi_{ki}
\ \ \ \ \ \ \ \ \ \ \ \ \ 
{\scriptstyle{(k\,=\,1\,\cdots\,r\,;\,\,\,
i\,=\,1\,\cdots\,s)}},
\]
or, if written in more length: 
\def\theequation{1}\begin{equation}
\label{S-370}
\sum_{\nu=1}^s\,\xi_{k\nu}(x)\,
\frac{\partial\zeta_i}{\partial x_\nu}
=
\sum_{\nu=1}^s\,
\frac{\partial\xi_{ki}}{\partial x_\nu}\,\zeta_\nu
\ \ \ \ \ \ \ \ \ \ \ \ \ 
{\scriptstyle{(k\,=\,1\,\cdots\,r\,;\,\,\,
i\,=\,1\,\cdots\,s)}}.
\end{equation}
The question is to determine the most general solutions
$\zeta_1, \dots, \zeta_s$ to 
these differential equations.

Let:
\def\theequation{2}\begin{equation}
\zeta_i
=
\omega_i(x_1,\dots,x_s)
\ \ \ \ \ \ \ \ \ \ \ \ \ {\scriptstyle{(i\,=\,1\,\cdots\,n)}}
\end{equation}
be any system of solutions of~\thetag{ 1}, 
whence all the expressions:
\[
X_k\,\omega_i
-
\sum_{\nu=1}^s\,
\frac{\partial\xi_{ki}}{\partial x_\nu}\,\zeta_\nu
\]
vanish identically after the substitution: $\zeta_1 = \omega_1 (x)$,
\dots, $\zeta_s = \omega_s (x)$; in other words: the system of
equations~\thetag{ 2} in the $2\, s$ variables $x_1, \dots, x_s$,
$\zeta_1, \dots, \zeta_s$ admits the $r$ infinitesimal
transformations:
\[
\label{S-370-sq}
W_k\,F
=
X_kf
+
\sum_{i=1}^s\,
\bigg\{
\sum_{\nu=1}^s\,
\frac{\partial\xi_{ki}}{\partial x_\nu}\,
\zeta_\nu
\bigg\}\,
\frac{\partial f}{\partial\zeta_i}
\ \ \ \ \ \ \ \ \ \ \ \ \ {\scriptstyle{(k\,=\,1\,\cdots\,r)}}.
\]
Conversely, if a system of equations of the form~\thetag{ 2}
admits the infinitesimal transformations $W_1f, \dots, 
W_rf$, then $\omega_1 (x), \dots, \omega_s (x)$ are
obviously solutions of the differential equations~\thetag{ 1}. 
Consequently, the integration of the differential 
equations~\thetag{ 1} is equivalent to the
determination of the most general system of equations~\thetag{ 2}
which admits the infinitesimal transformations $W_1f, \dots, 
W_rf$. 

Every system of equations which admits $W_1f, \dots, W_rf$
also allows the infinitesimal transformation
$W_k (W_j (f)) - W_j ( W_k (f)) = 
\leftbracket W_k,\, W_j \rightbracket$; we
compute it. 

We have:
\[
\aligned
\leftbracket
W_k,\,W_j
\rightbracket
=
\leftbracket
X_k,\,X_j
\rightbracket
&
+
\sum_{i,\,\,\mu,\,\,\nu}^{1\cdots s}\,
\bigg\{
\xi_{k\nu}\,
\frac{\partial^2\xi_{ji}}{\partial x_\mu\partial x_\nu}
-
\xi_{j\nu}\,
\frac{\partial^2\xi_{ki}}{\partial x_\mu\partial x_\nu}
\bigg\}\,
\zeta_\mu\,
\frac{\partial f}{\partial\zeta_i}
\\
&
+
\sum_{i,\,\,\mu,\,\,\nu}^{1\cdots s}\,
\bigg\{
\frac{\partial\xi_{k\nu}}{\partial x_\mu}\,
\frac{\partial\xi_{ji}}{\partial x_\nu}
-
\frac{\partial\xi_{j\nu}}{\partial x_\mu}\,
\frac{\partial\xi_{ki}}{\partial x_\nu}
\bigg\}\,
\zeta_\mu\,
\frac{\partial f}{\partial\zeta_i},
\endaligned
\]
and here, the right-hand side can be written:
\[
\leftbracket
X_k,\,X_j
\rightbracket
+
\sum_{i,\,\,\mu}^{1\cdots s}\,
\frac{\partial}{\partial x_\mu}\,
\sum_{\nu=1}^s\,
\bigg\{
\xi_{k\nu}\,
\frac{\partial\xi_{ji}}{\partial x_\nu}
-
\xi_{j\nu}\,
\frac{\partial\xi_{ki}}{\partial x_\nu}
\bigg\}\,
\zeta_\mu\,
\frac{\partial f}{\partial\zeta_i}.
\]
But now, there are relations of the form:
\[
\leftbracket
X_k,\,X_j
\rightbracket
=
\sum_{\sigma=1}^r\,c_{kj\sigma}\,X_\sigma f,
\]
from which it follows:
\[
\sum_{\nu=1}^s\,
\bigg\{
\xi_{k\nu}\,
\frac{\partial\xi_{ji}}{\partial x_\nu}
-
\xi_{j\nu}\,
\frac{\partial\xi_{ki}}{\partial x_\nu}
\bigg\}
=
\sum_{\sigma=1}^r\,c_{kj\sigma}\,\xi_{\sigma i},
\]
and therefore, we get simply:
\[
\leftbracket
W_k,\,W_j
\rightbracket
=
\sum_{\sigma=1}^r\,
c_{kj\sigma}\,W_\sigma f.
\]

From this, we see that $W_1f, \dots, W_rf$ generate an $r$-term group
in the $2s$ variables $x_1, \dots, x_s$, $\zeta_1, \dots, \zeta_s$;
but this was closely presumed.

\medskip

There can be $n$, amongst the infinitesimal transformations
$X_1f, \dots, X_rf$, say $X_1f, \dots, X_nf$ that are
linked together by no linear relation of the
form:
\[
\chi_1(x_1,\dots,x_s)\,X_1f
+\cdots+
\chi_n(x_1,\dots,x_s)\,X_nf
=
0,
\]
while $X_{ n+1}f, \dots, X_rf$ can be expressed in the
following way: 
\def\theequation{3}\begin{equation}
X_{n+k}f
\equiv
\sum_{\nu=1}^n\,\varphi_{k\nu}(x_1,\dots,x_s)\,
X_\nu f
\ \ \ \ \ \ \ \ \ \ \ \ \ 
{\scriptstyle{(k\,=\,1\,\cdots\,r\,-\,n)}}.
\end{equation}
Under these assumptions, the differential equations~\thetag{ 1}
can also be written as:
\[
\aligned
&
\ \ \ \ \ \ \ \ \ \ \ \ \ 
X_\nu\,\zeta_i
-
Z\,\xi_{\nu i}
=
0
\ \ \ \ \ \ \ \ \ \ \ \ \ 
{\scriptstyle{(\nu\,=\,1\,\cdots\,n\,;\,\,\,
i\,=\,1\,\cdots\,s)}}
\\
&
\sum_{\nu=1}^n\,\varphi_{k\nu}(x)\,X_\nu\,\zeta_i
-
Z\,\xi_{n+k,\,i}
=
0
\ \ \ \ \ \ \ \ \ \ \ \ \ 
{\scriptstyle{(k\,=\,1\,\cdots\,r\,-\,n\,;\,\,\,
i\,=\,1\,\cdots\,s)}}.
\endaligned
\]
Hence if the expressions $X_1\, \zeta_1, \dots, 
X_n\, \zeta_i$ are took away, we obtain
between $\zeta_1, \dots, \zeta_s$ the finite equations:
\def\theequation{4}\begin{equation}
\aligned
\sum_{\pi=1}^s\,
&
\bigg\{
\frac{\partial\xi_{n+k,\,i}}{\partial x_\pi}
-
\sum_{\nu=1}^n\,\varphi_{k\nu}\,
\frac{\partial\xi_{\nu i}}{\partial x_\pi}
\bigg\}\,
\zeta_\pi
=
0
\\
&\ \ \ \ \ \ \ \ 
{\scriptstyle{(k\,=\,1\,\cdots\,r\,-\,n\,;\,\,\,
i\,=\,1\,\cdots\,s)}}.
\endaligned
\end{equation}

It stands to reason that every system of equations of the
form~\thetag{ 2} which admits the group $W_1f, \dots, W_rf$ must
comprise the equations~\thetag{ 4}.

The equations~\thetag{ 4} are linear and homogeneous in $\zeta_1,
\dots, \zeta_s$; so if amongst them, one finds exactly $s$ that are
mutually independent, then $\zeta_1 = 0$, \dots, $\zeta_s = 0$ is the
only system of solutions which satisfies the equations.
In this case, there is only one system of equations of the
form~\thetag{ 2} which admits the group $W_1f, \dots, W_rf$, 
namely the system of equations: $\zeta_1 = 0$, \dots, 
$\zeta_s = 0$, hence there is no infinitesimal transformation
$Zf$ which is interchangeable with all $X_kf$. 

Differently, assume that the equations~\thetag{ 4} reduce to less
than $s$, say to \label{S-372}
$m < s$ independent equations. We will
see that in this case, aside from the
useless system $\zeta_1 = 0$, \dots, $\zeta_s = 0$, there
are yet also other systems of the form~\thetag{ 2}
which admit the group $W_1f, \dots, W_rf$. 

\medskip 

Above all, we observe that the system of equations~\thetag{ 4} admits
the group $W_1f, \dots, W_rf$.

In order to prove this, we write down the matrix which is associated
to the infinitesimal transformations $W_1f, \dots, W_rf$:
\def\theequation{5}\begin{equation}
\left\vert
\begin{array}{cccccccc}
\xi_{11} & \,\cdot\, & \,\cdot\, & \xi_{1s}\, & \,
\sum_{\nu=1}^s\frac{\partial\xi_{11}}{\partial x_\nu}\zeta_\nu
& \,\cdot\, & \,\cdot\, & 
\sum_{\nu=1}^s\frac{\partial\xi_{1s}}{\partial x_\nu}\zeta\nu
\\
\cdot & \,\cdot\, & \,\cdot\, & \cdot & \cdot &
\,\cdot\, & \,\cdot\, & \cdot
\\
\xi_{r1} & \,\cdot\, & \,\cdot\, & \xi_{rs}\, & \,
\sum_{\nu=1}^s\frac{\partial\xi_{r1}}{\partial x_\nu}\zeta_\nu
& \,\cdot\, & \,\cdot\, & 
\sum_{\nu=1}^s\frac{\partial\xi_{rs}}{\partial x_\nu}\zeta\nu
\end{array}
\right\vert.
\end{equation}
We will show that~\thetag{ 4} belongs to the system of equations that
one obtains by setting equal to zero all $(n+1) \times (n+1)$
determinants of this matrix. As a result, according to
Chap.~\ref{kapitel-14}, Theorem~39, p.~\pageref{Theorem-39-S-228}, it
will be proved that~\thetag{ 4} admits the the group $W_1f, \dots,
W_rf$.

Amongst the $(n+1) \times (n+1)$ determinants of the matrix~\thetag{
5}, there are those of the form:
\[
D
=
\left\vert
\begin{array}{ccccc}
\xi_{1k_1} & \,\cdot\, & \,\cdot\, & \xi_{1k_n}\, & \,
\sum_{\nu=1}^s\frac{\partial\xi_{1\sigma}}{\partial x_\nu}\zeta_\nu
\\
\cdot & \,\cdot\, & \,\cdot\, & \cdot & \cdot
\\
\xi_{nk_1} & \,\cdot\, & \,\cdot\, & \xi_{nk_n}\, & \,
\sum_{\nu=1}^s\frac{\partial\xi_{n\sigma}}{\partial x_\nu}\zeta_\nu
\\
\cdot & \,\cdot\, & \,\cdot\, & \cdot & \cdot
\\
\xi_{n+j,\,k_1} & \,\cdot\, & \,\cdot\, & \xi_{n+j,\,k_n}\, & \,
\sum_{\nu=1}^s\frac{\partial\xi_{n+j,\,\sigma}}{\partial x_\nu}\zeta_\nu
\end{array}
\right\vert.
\]
If we use the identities following from~\thetag{ 3}:
\[
\xi_{n+j,\,\pi}
\equiv
\sum_{\nu=1}^n\,\varphi_{j\nu}(x)\,\xi_{\nu\pi}
\ \ \ \ \ \ \ \ \ \ \ \ \ 
{\scriptstyle{(\pi\,=\,1\,\cdots\,s\,;\,\,\,
j\,=\,1\,\cdots\,r\,-\,n)}},
\]
then we can obviously write $D$ as: 
\[
D
=
\sum\,\pm\,
\xi_{1k_1}\cdots\,\xi_{nk_n}\,
\sum_{\nu=1}^s\,
\bigg\{
\frac{\partial\xi_{n+j,\,\sigma}}{\partial x_\nu}
-
\sum_{\tau=1}^n\,\varphi_{j\tau}\,
\frac{\partial\xi_{\tau\sigma}}{\partial x_\nu}
\bigg\}\,\zeta_\nu.
\]
Now, under the assumptions made, not all determinants of the form:
\[
\sum\,\pm\,\xi_{1k_1}\cdots\,\xi_{nk_n}
\]
vanish identically; so if we set equal to zero all determinants of the
form $D$, we obtain either relations between the $x_1, \dots, x_s$
alone, or we obtain the system of equations~\thetag{ 4}. But as it is
easy to see, this system of equations actually brings to zero all
$(n+1) \times (n+1)$ determinants of the matrix~\thetag{ 5}, hence it
admits the group $W_1f, \dots, W_rf$.

The determination of the most general system of equations~\thetag{ 2}
which admits the group $W_1f, \dots, W_rf$ and which comprises the
equations~\thetag{ 4} can now be executed on the basis of
Chap.~\ref{kapitel-14}, p.~\pageref{S-236} up to p.~\pageref{S-238}.

The system of equations~\thetag{ 4} brings to zero all $(n+1) \times
(n+1)$ determinants of the matrix~\thetag{ 5}, but not all $n \times
n$ determinants; likewise, a system of equations of the form~\thetag{
2} cannot make zero all $n \times n$ determinants of the
matrix~\thetag{ 5}. Thus, we proceed in the following way: We solve
the equations~\thetag{ 4} with respect to $m$ of the quantities
$\zeta_1, \dots, \zeta_s$, say with respect to $\zeta_1, \dots,
\zeta_m$, afterwards, following the introduction of the cited
developments, we form the reduced infinitesimal transformations
$\overline{ W}_1f, \dots, \overline{ W}_rf$ in the $2\, s - m$
variables $x_1, \dots, x_s$, $\zeta_{ m+1}, \dots, \zeta_s$, and
lastly, we determine $2\, s - m - n$ arbitrary independent solutions
of the $n$-term complete system which is formed by the $n$ equations:
$\overline{ W}_1f = 0$, \dots, $\overline{ W}_n f = 0$.

The $n$ equations $\overline{ W}_1f = 0$, \dots, $\overline{ W}_n f =
0$ are solvable with respect to $n$ of the differential quotients
$\partial f / \partial x_1$, \dots, $\partial f / \partial x_s$, hence
its $2\, s - m - n$ independent solutions are mutually independent
relatively to $s - n$ of the $x$ and to the variables $\zeta_{ m+1},
\dots, \zeta_s$ (cf. Chap.~\ref{kapitel-15}, Theorem~12,
p.~\pageref{Theorem-12-S-91}). Now, amongst the solutions of the
complete system $\overline{ W}_k f = 0$, there are exactly $s - n$
independent ones which satisfy at the same time the $n$-term complete
system: $X_1f = 0$, \dots, $X_nf = 0$, and hence depend only upon the
$x$, and they can be called:
\[
\mathfrak{u}_1(x_1,\dots,x_s),
\,\,\,\dots,\,\,\,
\mathfrak{u}_{s-n}(x_1,\dots,x_s).
\]
So, if: 
\[
\mathfrak{B}_1
(\zeta_{m+1},\dots,\zeta_s,\,x_1,\dots,x_s),
\,\,\,\dots,\,\,\,
\mathfrak{B}_{s-m}
(\zeta_{m+1},\dots,\zeta_s,\,x_1,\dots,x_s)
\]
are $s - m$ arbitrary mutually independent, and independent of the
$\mathfrak{ u}$, solutions of the complete system $\overline{ W}_k f =
0$, then necessarily, these solutions are mutually independent
relatively to $\zeta_{ m+1}, \dots, \zeta_s$.

At present, we obtain the most general system of equations~\thetag{ 2}
which admits the group $W_1f, \dots, W_rf$ by adding to the
equations~\thetag{ 4}, in the most general way, $s - m$ mutually
independent relations between $\mathfrak{ u}_1, \dots, \mathfrak{ u}_{
s-n}$, $\mathfrak{ B}_1, \dots, \mathfrak{ B}_{ s-m}$ that are
solvable with respect to $\zeta_{ m+1}, \dots, \zeta_s$. It is clear
that these relations must be solvable with respect to $\mathfrak{
B}_1, \dots, \mathfrak{ B}_{ s-m}$, so that they can be brought to the
form:
\def\theequation{6}\begin{equation}
\mathfrak{B}_\mu
(\zeta_{m+1},\dots,\zeta_s,\,x_1,\dots,x_s)
=
\Omega_\mu
\big(\mathfrak{u}_1(x),\,\dots,\,\mathfrak{u}_{s-n}(x)\big)
\ \ \ \ \ \ \ \ \ \ \ \ \ 
{\scriptstyle{(\mu\,=\,1\,\cdots\,s\,-\,m)}}.
\end{equation}
Here, the $\Omega_\mu$ are subject to no restriction at all, and they
are absolutely arbitrary functions of their arguments.

Thus, if we add the equations~\thetag{ 6} to~\thetag{ 4} and if we
solve, what is always possible, all of them with respect to $\zeta_1,
\dots, \zeta_s$, then we obtain the most general system of
equations~\thetag{ 2} which admits the group $W_1f, \dots, W_rf$ and
as a result also, the most general system of solutions to the
differential equations~\thetag{ 1}. Visibly, this most general system
of solutions contains $s - m$ arbitrary functions of $\mathfrak{ u}_1,
\dots, \mathfrak{ u}_{ s - n}$.

But now, the differential equations~\thetag{ 1} are linear
and homogeneous in the unknowns $\zeta_1, \dots, \zeta_s$; 
thus, it can be concluded that its most general 
system of solutions $\zeta_1, \dots, \zeta_s$ can be
deduced from
$s - m$ particular systems of solutions:
\[
\zeta_{\mu 1}(x),
\,\dots,\,
\zeta_{\mu s}(x)
\ \ \ \ \ \ \ \ \ \ \ \ \ 
{\scriptstyle{(\mu\,=\,1\,\cdots\,s\,-\,m)}}
\]
in the following way:
\[
\aligned
\zeta_i
=
\chi_1
(\mathfrak{u}_1,\dots,\,\mathfrak{u}_{s-n})\,
\zeta_{1i}
&
+\cdots+
\chi_{s-m}
(\mathfrak{u}_1,\dots,\,\mathfrak{u}_{s-n})\,
\zeta_{s-m,\,i}
\\
&\ \ \
{\scriptstyle{(i\,=\,1\,\cdots\,s)}},
\endaligned
\]
where the $\chi$ are completely arbitrary functions of the $\mathfrak{
u}$; naturally, the particular system of solutions in question must be
constituted in such a way that there are no $s - m$ functions $\psi_1
( \mathfrak{ u}_1, \dots, \mathfrak{ u}_{ s - n})$, \dots, $\psi_{ s -
m} ( \mathfrak{ u}_1, \dots, \mathfrak{ u}_{ s-n})$ which satisfy
identically the $s$ equations:
\[
\psi_1(\mathfrak{u})\,\zeta_{1i}
+\cdots+
\psi_{s-m}(\mathfrak{u})\,\zeta_{s-m,\,i}
=
0
\ \ \ \ \ \ \ \ \ \ \ \ \ {\scriptstyle{(i\,=\,1\,\cdots\,s)}}.
\]
It can even be proved that actually, there are no $s - m$ functions
$\Psi_1 ( x_1, \dots, x_s)$, \dots, $\Psi_{ s - m} (x_1, \dots, x_s)$
which satisfy identically the $s - m$ equations:
\[
\Psi_1(x)\,\zeta_{1i}
+\cdots+
\Psi_{s-m}\,\zeta_{s-m,\,i}
=
0
\ \ \ \ \ \ \ \ \ \ \ \ \ 
{\scriptstyle{(i\,=\,1\,\cdots\,s\,-\,m)}}.
\]
Indeed, the $s - m$ equations: 
\[
\zeta_{m+\sigma}
=
\sum_{\mu=1}^{s-m}\,\chi_\mu(\mathfrak{u})\,
\zeta_{\mu,\,m+\sigma}
\ \ \ \ \ \ \ \ \ \ \ \ \ 
{\scriptstyle{(\sigma\,=\,1\,\cdots\,s\,-\,m)}}
\]
are obviously equivalent to the equations~\thetag{ 6}, 
hence they must be solvable with respect to
$\chi_1, \dots, \chi_{ s-m}$ and the determinant:
\[
\sum\,\pm\,\zeta_{1,\,m+1}\cdots\,
\zeta_{s-m,\,m+s-m}
\]
should not vanish identically. 

According to these preparations, we can finally determine the form
that the most general infinitesimal transformation $Zf$
interchangeable with $X_1f, \dots, X_rf$ possesses.
The concerned transformation reads:
\[
Zf
=
\sum_{\mu=1}^{s-m}\,
\chi_\mu(\mathfrak{u}_1,\dots,\mathfrak{u}_{s-m})\,
Z_\mu f, 
\]
where the $s - m$ infinitesimal transformations:
\[
Z_\mu f
=
\sum_{i=1}^s\,\zeta_{\mu i}(x_1,\dots,x_s)\,
\frac{\partial f}{\partial x_i}
\ \ \ \ \ \ \ \ \ \ \ \ \ 
{\scriptstyle{(\mu\,=\,1\,\cdots\,s\,-\,m)}}
\]
are linked together by no linear relation of the form:
\[
\Psi_1(x_1,\dots,x_s)\,Z_1f
+\cdots+
\Psi_{s-m}(x_1,\dots,x_s)\,Z_{s-m}f
=
0.
\]
In addition, since according to Proposition~1,
p.~\pageref{Satz-1-S-369}, every infinitesimal transformation
$\leftbracket Z_\mu,\, Z_\nu \rightbracket$ is also interchangeable
with $X_1f, \dots, X_rf$, there are relations of the specific form:
\[
\leftbracket
Z_\mu,\,Z_\nu
\rightbracket
=
\sum_{\pi=1}^{s-m}\,
\omega_{\mu\nu\pi}
(\mathfrak{u}_1,\dots,\mathfrak{u}_{s-m})\,
Z_\pi f
\ \ \ \ \ \ \ \ \ \ \ \ \ 
{\scriptstyle{(\mu,\,\,\nu\,=\,1\,\cdots\,s\,-\,m)}}.
\]

We thus see: There always exists a transformation which is
interchangeable with $X_1f, \dots, X_rf$ when and only when the number
$m$ defined on p.~\pageref{S-372} is smaller than the number $s$ of
the variables $x$. If $m < s$ and if at the same time, $n < s$, then
the group $X_1f, \dots, X_r f$ is intransitive, so the most general
infinitesimal transformation $Zf$ interchangeable with $X_1f, \dots,
X_rf$ depends upon arbitrary functions. If by contrast $n = s$, then
the group $X_1f, \dots, X_rf$ is transitive, so the most general
transformation $Zf$ interchangeable with $X_1f, \dots, X_rf$ can be
linearly deduced from $s - m$ independent infinitesimal
transformations: $Z_1f, \dots, Z_{ s-m} f$; according to a remark
made earlier on (p.~\pageref{Satz-1-S-369}), the concerned infinitesimal
transformations then generate an $(s-m)$-term group.\,---

We know that an infinitesimal transformation interchangeable
with $X_1f, \dots, X_rf$ exists only when the equations~\thetag{ 4}
reduce to less than $s$ independent equations. 
We can give a somewhat more lucid interpretation of this
condition by remembering the identities~\thetag{ 3}, 
or the equivalent identities:
\[
\xi_{n+k,\,i}
-
\sum_{\nu=1}^n\,\varphi_{k\nu}(x)\,
\xi_{\nu i}
\equiv
0
\ \ \ \ \ \ \ \ \ \ \ \ \ 
{\scriptstyle{(k\,=\,1\,\cdots\,r\,-\,n\,;\,\,\,
i\,=\,1\,\cdots\,s)}},
\]
that define the functions $\varphi_{ k\nu}$. 
Indeed, if we differentiate with respect to $x_j$ the
identities just written, we obtain the following
identities:
\[
\aligned
\frac{\partial\xi_{n+k,\,i}}{\partial x_j}
&
-
\sum_{\nu=1}^n\,
\bigg\{
\varphi_{k\nu}\,
\frac{\partial\xi_{\nu i}}{\partial x_j}
+
\xi_{\nu i}\,
\frac{\partial\varphi_{k\nu}}{\partial x_j}
\bigg\}
\equiv 
0
\\
&\ \ \
{\scriptstyle{(k\,=\,1\,\cdots\,r\,-\,n\,;\,\,\,
i\,=\,1,\,\,j\,\cdots\,s)}},
\endaligned
\]
by virtue of which the equations~\thetag{ 4} can 
be replaced by the equivalent equations:
\[
\sum_{\nu=1}^n\,
\xi_{\nu i}\,
\sum_{j=1}^s\,\zeta_j\,
\frac{\partial\varphi_{k\nu}}{\partial x_j}
=
0
\ \ \ \ \ \ \ \ \ \ \ \ \ 
{\scriptstyle{(i\,=\,1\,\cdots\,s\,;\,\,\,
k\,=\,1\,\cdots\,r\,-\,n)}}.
\]
But since not all determinants of the form $\sum\, \pm \, \xi_{ k_1}
\cdots\, \xi_{ k_n}$ vanish identically, then in turn, the latter
equations are equivalent to:
\def\theequation{4'}\begin{equation}
\sum_{j=1}^s\,\zeta_j\,
\frac{\partial\varphi_{k\nu}}{\partial x_j}
=
0
\ \ \ \ \ \ \ \ \ \ \ \ \ 
{\scriptstyle{(k\,=\,1\,\cdots\,r\,-\,n\,;\,\,\,
\nu\,=\,1\,\cdots\,n)}}.
\end{equation}

Thus, there is a transformation $Zf$ interchangeable with $X_1f,
\dots, X_rf$ only when the linear equations~\thetag{ 4'} in the
$\zeta_j$ reduce to less than $s$ independent equations.

The equations~\thetag{ 4'} are more clearly arranged than the
equations~\thetag{ 4}, and they have in addition a simple meaning, for
they express that each one of the $n ( r - n)$ functions $\varphi_{ k
\nu} (x_1, \dots, x_s)$ admits all the infinitesimal transformations
$Zf$. 

If, amongst the equations~\thetag{ 4}, there are exactly $m$ that are
mutually independent, then naturally, amongst the equations~\thetag{
4'}, there are also exactly $m$ that are mutually independent, and
therefore, $m$ is nothing but the number of independent functions
amongst the $n ( r - n)$ functions $\varphi_{ k\nu} (x)$.

We recapitulate the gained result:

\renewcommand{\thefootnote}{\fnsymbol{footnote}}
\def\thetheorem{67}\begin{theorem}
\label{Theorem-67-S-376}
If, amongst the $r$ infinitesimal transformations:
\[
X_kf
=
\sum_{i=1}^s\,\xi_{ki}(x_1,\dots,x_s)\,
\frac{\partial f}{\partial x_i}
\ \ \ \ \ \ \ \ \ \ \ \ \ {\scriptstyle{(k\,=\,1\,\cdots\,r)}},
\]
of an $r$-term group, $X_1f, \dots, X_nf$, say, are linked together by
no linear relation, while $X_{ n+1}f, \dots, X_rf$ can be expressed
linearly in terms of $X_1f, \dots, X_nf$:
\[
X_{n+j}f
\equiv
\sum_{\nu=1}^n\,
\varphi_{j\nu}(x_1,\dots,x_s)\,X_\nu f
\ \ \ \ \ \ \ \ \ \ \ \ \ {\scriptstyle{(j\,=\,1\,\cdots\,r\,-\,n)}},
\]
and if amongst the $n ( r - n)$ functions $\varphi_{ k\nu}$, there
exist exactly $s$ that are mutually independent, then there is no
infinitesimal transformation which is interchangeable with all $X_kf$;
by contrast, if amongst the functions $\varphi_{ k\nu}$, there exist
less than $s$, say only $m$, that are mutually independent, then there
are infinitesimal transformations which are interchangeable with
$X_1f, \dots, X_rf$, and to be precise, the most general
infinitesimal transformation $Zf$ of this nature possesses
the form:
\[
Zf
=
\psi_1(\mathfrak{u}_1,\dots,\mathfrak{u}_{s-n})\,Z_1f
+\cdots+
\psi_{s-m}(\mathfrak{u}_1,\dots,\mathfrak{u}_{s-n})\,
Z_{s-m}f,
\]
where $\mathfrak{u}_1, \dots, \mathfrak{ u}_{ s -n}$ denote
independent solutions of the $n$-term complete system: $X_1 f = 0$,
\dots, $X_nf = 0$, where furthermore $\psi_1, \dots, \psi_{ s - m}$
mean arbitrary functions of their arguments, and lastly, 
where the infinitesimal transformations:
\[
Z_\mu f
=
\sum_{i=1}^s\,
\zeta_{\mu i}(x_1,\dots,x_s)\,
\frac{\partial f}{\partial x_i}
\ \ \ \ \ \ \ \ \ \ \ \ \ 
{\scriptstyle{(\mu\,=\,1\,\cdots\,s\,-\,m)}}
\]
interchangeable with $X_1f, \dots, X_rf$ and linked together by no
linear relation stand pairwise in relationships of the
form\footnote{\,
\name{Lie} suggested this general theorem in 
the Math. Ann. Vol. XXV.
}:

\[
\leftbracket
Z_\mu,\,Z_\nu
\rightbracket
=
\sum_{\pi=1}^{s-m}\,
\omega_{\mu\nu\pi}
(\mathfrak{u}_1,\dots,\mathfrak{u}_{s-m})\,
Z_\pi f.
\]
\end{theorem}
\renewcommand{\thefootnote}{\arabic{footnote}}

\smallercharacters{
Besides, a part of the result stated in this proposition follows
immediately from the developments of the previous chapter. Indeed, if
the number of independent functions amongst the functions $\varphi_{
k\nu } (x_1, \dots, x_s)$ is exactly equal to $s$, then according to
p.~\pageref{S-344}, the equations:
\[
\varphi_{k\nu}(x_1',\dots,x_s')
=
\varphi_{k\nu}(x_1,\dots,x_s)
\ \ \ \ \ \ \ \ \ \ \ \ \ 
{\scriptstyle{(k\,=\,1\,\cdots\,r\,-\,n\,;\,\,\,
\nu\,=\,1\,\cdots\,n)}}
\]
represent the most general transformation
which transfers $X_1f, \dots, X_rf$ to $X_1'f, \dots, X_r'f$,
respectively. By contrast, if amongst the $\varphi_{ k \nu}$, there
are less than $s$ that are mutually independent, then according to
Theorem~65, p.~\pageref{Theorem-65-S-353}, there is a continuous
amount \deutsch{Menge} of transformations which leaves invariant every
$X_k f$; we remarked already on p.~\pageref{S-368} that the totality
of all such transformations forms a group, about which it can be
directly seen that it comprises one-term groups.

}

It seems to be adequate to yet state expressly the following
proposition:

\def\theproposition{2}\begin{proposition}
If there actually exists an infinitesimal transformation which is
interchangeable with all infinitesimal transformations $X_1f, \dots,
X_rf$ of an $r$-term group of the space $x_1, \dots, x_s$, then the
most general infinitesimal transformation interchangeable with $X_1f,
\dots, X_rf$ contains arbitrary functions, as soon as the group $X_1f,
\dots, X_rf$ is intransitive, but by contrast, it
contains only arbitrary
parameters, as soon as the group is transitive.
\end{proposition}

\sectionengellie{\S\,\,\,96.}

Of outstanding significance is the case where the group $X_kf$ is
simply transitive, thus the case $s = r = n$.

If one would desire, in this case, to know the most general
transformation $x_i' = \Phi_i ( x_1, \dots, x_n)$ by virtue of which
every $X_kf$ takes the form $X_k' f$, then one would only have to seek
$n$ independent solutions $\Omega_1, \dots, \Omega_n$ of the $n$-term
complete system:
\[
X_kf+X_k'f
=
0
\ \ \ \ \ \ \ \ \ \ \ \ \ {\scriptstyle{(k\,=\,1\,\cdots\,n)}}
\]
and to set these solutions equal to 
arbitrary constants $a_1, \dots, a_n$.
The equations:
\[
\Omega_k
(x_1,\dots,x_n,\,x_1',\dots,x_n')
=
a_k
\ \ \ \ \ \ \ \ \ \ \ \ \ {\scriptstyle{(k\,=\,1\,\cdots\,n)}}
\]
are then solvable both with respect to the $x$ and
with respect to the $x'$, and they represent the
transformation demanded. 

From the beginning, we know that the totality of all transformations
$\Omega_k = a_k$ forms a group. At present, we see that this group is
$n$-term and simply transitive; this results immediately from the form
in which the group is presented, resolved with respect to its $n$
parameters.

The group:
\[
\Omega_1(x,x')=a_1,
\,\,\,\dots,\,\,\,
\Omega_n(x,x')=a_n
\]
contains the identity transformation and $n$ independent
infinitesimal transformations:
\[
Z_kf
=
\sum_{i=1}^n\,\zeta_{ki}(x_1,\dots,x_n)\,
\frac{\partial f}{\partial x_i}
\ \ \ \ \ \ \ \ \ \ \ \ \ {\scriptstyle{(k\,=\,1\,\cdots\,n)}},
\]
and this results from the developments of the preceding paragraph; but
this is also clear in itself, since the transformations of the group
are ordered as inverses by pairs, and hence the Theorem~56 in
Chap.~\ref{kapitel-18}, p.~\pageref{Theorem-56-S-315} finds an
application. Besides, it follows from the transitivity of the group
$Z_1f, \dots, Z_nf$ that $Z_1f, \dots, Z_nf$ are linked together by no
linear relation of the form $\sum\, \chi_i (x_1, \dots, x_s)\, Z_if =
0$.

\medskip

Between the two \emphasis{simply transitive} groups $X_kf$ and $Z_if$,
there is a full relationship of reciprocity
\deutsch{Reciprocitätsverhältniss}.
If the $X_kf$ are given, then the general infinitesimal
transformation $Zf = \sum\, e_i\, Z_if$
is completely defined by the $n$ equations: 
\[
X_k\big(Z(f)\big)
-
Z\big(X_k(f)\big)
=
0
\ \ \ \ \ \ \ \ \ \ \ \ \ {\scriptstyle{(k\,=\,1\,\cdots\,n)}}\,;
\]
on the other hand, if the $Z_if$ are given then the
$n$ equations:
\[
X\big(Z_i(f)\big)
-
Z_i\big(X(f)\big)
=
0
\ \ \ \ \ \ \ \ \ \ \ \ \ {\scriptstyle{(i\,=\,1\,\cdots\,n)}}
\]
determine in the same way the general infinitesimal 
transformation:
\[
Xf
=
\sum\,e_k\,X_kf.
\]

But with that, the peculiar relationship in which the two groups stand
is not yet written down exhaustively. Indeed, as we will now show, the
two groups are also yet equally composed; since both are simply
transitive, it then follows immediately that they are similar to each
other (cf. Chap.~\ref{kapitel-19}, Theorem~64,
p.~\pageref{Theorem-64-S-340}).

In the infinitesimal transformations $X_kf$ and $Z_if$, we imagine
that the $\xi_{ k\nu}$ and the $\zeta_{ i\nu}$ are expanded in powers
of $x_1, \dots, x_n$, and we assume at the same time that $x_1 = 0$,
\dots, $x_n = 0$ is a point in general position. According to
Chap.~\ref{kapitel-13}, p.~\pageref{Satz-4-S-217}, since the group
$X_kf$ is simply transitive, it contains exactly $n$ infinitesimal
transformations of zeroth order in $x_1, \dots, x_n$ out of which no
infinitesimal transformation of first order, or of higher order, can
be linearly deduced.

We can therefore imagine that $X_1f, \dots, X_nf$ are replaced
\label{S-379}
by $n$ other independent infinitesimal transformations
$\mathfrak{ X}_1f, \dots, \mathfrak{ X}_nf$
which have the form: 
\[
\aligned
\mathfrak{X}_kf
=
\frac{\partial f}{\partial x_k}
&
+
\sum_{\mu,\,\,\nu}^{1\cdots n}\,
h_{k\mu\nu}\,x_\mu\,
\frac{\partial f}{\partial x_\nu}
+\cdots
\\
&
\ \ \ {\scriptstyle{(k\,=\,1\,\cdots\,n)}}, 
\endaligned
\]
after leaving out the terms of second order and
of higher order. In the same way, we can replace
$Z_1f, \dots, Z_nf$ by $n$ independent infinitesimal
transformations of the form:
\[
\aligned
\mathfrak{Z}_if
=
-\,\frac{\partial f}{\partial x_i}
&
+
\sum_{\mu,\,\nu}^{1\cdots n}\,
l_{i\mu\nu}\,x_\mu\,
\frac{\partial f}{\partial x_\nu}
+\cdots
\\
&
\ \ \ {\scriptstyle{(i\,=\,1\,\cdots\,n)}}.
\endaligned
\]

After these preparations, we form the equations
$\leftbracket \mathfrak{X}_k,\, \mathfrak{Z}_i \rightbracket = 0$\,;
the same equation takes the form: 
\[
\sum_{\nu=1}^n\,
(l_{ik\nu}+h_{ki\nu})\,
\frac{\partial f}{\partial x_\nu}
+\cdots
=
0, 
\]
from which it follows:
\[
l_{ik\nu}
=
-\,h_{ki\nu}.
\]
Computations of the same sort yield:
\[
\leftbracket
\mathfrak{X}_k,\,\mathfrak{X}_j
\rightbracket
=
\sum_{\nu=1}^n\,
(h_{jk\nu}-h_{kj\nu})\,
\frac{\partial f}{\partial x_\nu}
+\cdots
\]
and:
\[
\aligned
\leftbracket
\mathfrak{Z}_k,\,\mathfrak{Z}_j
\rightbracket
&
=
\sum_{\nu=1}^n\,
(-\,l_{jk\nu}+l_{kj\nu})\,
\frac{\partial f}{\partial x_\nu}
+\cdots
\\
&
=
\sum_{\nu=1}^n\,
(h_{kj\nu}-h_{jk\nu})\,
\frac{\partial f}{\partial x_\nu}
+\cdots.
\endaligned
\]
On the other hand, we have: 
\[
\leftbracket
\mathfrak{X}_k,\,\mathfrak{X}_j
\rightbracket
=
\sum_{\nu=1}^n\,\mathfrak{c}_{kj\nu}\,
\mathfrak{X}_\nu f,
\ \ \ \ \ \ \ \ \ \ \ \
\leftbracket
\mathfrak{Z}_k,\,\mathfrak{Z}_i
\rightbracket
=
\sum_{\nu=1}^n\,\mathfrak{c}_{kj\nu}'\,
\mathfrak{Z}_\nu f\,;
\]
here, we insert the expressions just found for: $\mathfrak{ X}_\nu f$,
$\mathfrak{ Z}_\nu f$, $\leftbracket \mathfrak{ X}_k,\, \mathfrak{
X}_j \rightbracket$, $\leftbracket \mathfrak{ Z}_k,\, \mathfrak{ Z}_j
\rightbracket$ and afterwards, we make the substitution: $x_1 = 0$,
\dots, $x_n = 0$; it then comes:
\[
\mathfrak{c}_{kj\nu}
=
h_{jk\nu}-h_{kj\nu}
=
\mathfrak{c}_{kj\nu}'.
\]

Thus, our two groups are effectively equally composed, and in
consequence of that, as already remarked above, they are similar to
each other. As a result, the following holds:

\renewcommand{\thefootnote}{\fnsymbol{footnote}}
\def\thetheorem{68}\begin{theorem}
\label{Theorem-68-S-380}
If $X_1f, \dots, X_nf$ are independent infinitesimal transformations
of a simply transitive group in the variables $x_1, \dots, x_n$, then
the $n$ equations $\leftbracket X_k, \, Z \rightbracket = 0$ define
the general infinitesimal transformation $Zf$ of a second simply
transitive group $Z_1f, \dots, Z_nf$ which has the same composition as
the group $X_1f, \dots, X_nf$ and which is at the same time similar to
it. The relationship between these two simply transitive groups is a
reciprocal relationship: each one of the two groups consists of the
totality of all one-term groups whose transformations are
interchangeable with all transformations of the other
group.\footnote[1]{\,
\name{Lie} communicated the Theorem~68 in the Gesellschaft
der Wissenschaften zu Christiania in Nov. 1882 and in May 1883; 
cf. also the Math. Ann. Vol. XXV, p.~107 sq.
}
\end{theorem}
\renewcommand{\thefootnote}{\arabic{footnote}}

It is convenient to call the groups $X_kf$ and $Z_if$ as
\terminology{reciprocal} \deutsch{reciproke} transformation groups, or
always, the one as the reciprocal group of the other.

If we remember from Chap.~\ref{kapitel-16}, p.~\pageref{S-276}, that
the excellent infinitesimal transformations $e_1\, X_1f + \cdots +
e_r\, X_rf$ of the group $X_1f, \dots, X_nf$ are defined by the $n$
equations $\leftbracket X_k, \, Y \rightbracket = 0$, then we can yet
state the following proposition:

\def\theproposition{3}\begin{proposition}
The common infinitesimal transformations of two reciprocal simply
transitive groups are at the same time the excellent infinitesimal
transformations of both groups.
\end{proposition}

Let the two $n$-term groups $X_1f, \dots, X_nf$ and $Z_1f, \dots,
Z_nf$ in the $n$ variables $x_1, \dots, x_n$ be simply transitive and
reciprocal to each other. Then, if after the introduction of new
variables $x_1', \dots, x_n'$, \label{S-381}
the $X_kf$ are transferred to $X_k'f$
and the $Z_kf$ to $Z_k'f$, the two simply transitive groups $X_1'f,
\dots, X_n'f$ and $Z_1' f, \dots, Z_n'f$ are also reciprocal to each
other; this follows immediately from the relations: $\leftbracket
X_k',\, Z_j' \rightbracket = \leftbracket X_k,\, Z_j \rightbracket
\equiv 0$.

From this, it follows in particular the

\def\theproposition{4}\begin{proposition}
If an $n$-term simply transitive group in $n$ variables remains
invariant by a transformation, then at the same time, its reciprocal
simply transitive group remains invariant.
\end{proposition}

From this proposition it lastly follows the next one:

\def\theproposition{5}\begin{proposition}
The largest group of the space $x_1, \dots, x_n$ in which an $n$-term
simply transitive group of this space is contained as an invariant
subgroup coincides with the largest subgroup in which the reciprocal
group of this simply transitive group is contained as an invariant
subgroup.
\end{proposition}

\medskip

We want to give a couple of simple examples to the preceding
general developments about simply transitive groups. 

The group:
\[
\frac{\partial f}{\partial x},
\ \ \ \ \ 
x\,\frac{\partial f}{\partial x}
+
y\,\frac{\partial f}{\partial y}
\]
of the place $x$, $y$ is simply transitive. 
The finite equations of the reciprocal group are obtained
by integration of the complete system: 
\[
\frac{\partial f}{\partial x}
+
\frac{\partial f}{\partial x'}
=
0,
\ \ \ \ \ \ \
x\,\frac{\partial f}{\partial x}
+
y\,\frac{\partial f}{\partial y}
+
x'\,\frac{\partial f}{\partial x'}
+
y'\,\frac{\partial f}{\partial y'}
=
0
\]
under the following form:
\[
\frac{x'-x}{y}
=
a,
\ \ \ \ \ \ \
\frac{y'}{y}
=
b,
\]
hence solved with respect to $x'$ and $y'$:
\[
x'
=
x+a\,y,
\ \ \ \ \ \ \
y'
=
b\,y.
\]
The infinitesimal transformations of the reciprocal group are
therefore: 
\[
y\,\frac{\partial f}{\partial x},
\ \ \ \ \ \ \
y\,\frac{\partial f}{\partial y}.
\]

An interesting example is provided by the six-term projective group of
a nondegenerate surface of second order in ordinary space. Indeed,
this group contains two three-term simply transitive and reciprocal to
each other groups, of which the first group leaves fixed all 
generatrices of the one family, while the other group 
leaves fixed all generatrices of the other family. 

If $z - xy = 0$ is the equation of the surface, 
then: 
\[
\aligned
X_1f
&
=
\frac{\partial f}{\partial x}
+
y\,\frac{\partial f}{\partial z},
\ \ \ \ \ \ \
X_2f
=
x\,\frac{\partial f}{\partial x}
+
z\,\frac{\partial f}{\partial z},
\\
X_3f
&
=
x^2\,\frac{\partial f}{\partial x}
+
(xy-z)\,\frac{\partial f}{\partial y}
+
xz\,\frac{\partial f}{\partial z}
\endaligned
\]
is one of the two simply transitive groups and:
\[
\aligned
Z_1f
&
=
\frac{\partial f}{\partial y}
+
x\,\frac{\partial f}{\partial z},
\ \ \ \ \ \ \
Z_2f
=
y\,\frac{\partial f}{\partial y}
+
z\,\frac{\partial f}{\partial z}
\\
Z_3f
&
=
(xy-z)\,\frac{\partial f}{\partial x}
+
y^2\,\frac{\partial f}{\partial y}
+
yz\,\frac{\partial f}{\partial z}
\endaligned
\]
is the other.

The two reciprocal groups are naturally similar to each other; but it
is to be observed that they are similar through a
\emphasis{projective} transformation, namely through every projective
transformation which interchanges the two families of generatrices of
the surface.

\sectionengellie{\S\,\,\,97.}

We proceed to the general studies about reciprocal simply transitive
groups.

Let $X_1f, \dots, X_nf$ and $Z_1f, \dots, Z_nf$ be two simply
transitive and reciprocal groups in the variables $x_1, \dots, x_n$.

If $n$ would be equal to $1$, then the two groups would be identical
to each other, as one easily convinces oneself; we therefore want to
assume that $n$ is larger than $1$. Then the group $Z_1f, \dots,
Z_nf$ certainly contains subgroups. If $Z_1f, \dots, Z_mf$ ($m < n$)
is such a subgroup, then the $m$ equations:
\[
Z_1f=0,
\,\,\,\dots,\,\,\,
Z_mf=0
\]
form an $m$-term complete system which, as it
results from the identities:
\[
\leftbracket
X_i,\,Z_1
\rightbracket
\equiv
0,
\,\,\,\dots,\,\,\,
\leftbracket
X_i,\,Z_m
\rightbracket
\equiv
0
\ \ \ \ \ \ \ \ \ \ \ \ \ {\scriptstyle{(i\,=\,1\,\cdots\,n)}},
\]
admits the group $X_1f, \dots, X_nf$ (cf. Chap.~\ref{kapitel-8},
Theorem~20, p.~\pageref{Theorem-20-S-140}). Consequently, the group
$X_1f, \dots, X_nf$ is \emphasis{imprimitive}.

If $u_1, \dots, u_{ n-m}$ are independent solutions of the
complete system $Z_1f = 0$, \dots, 
$Z_mf = 0$, then according to
Chap.~\ref{kapitel-8}, Proposition~1, 
\pageref{Satz-1-S-139}, there are relations of the form:
\[
X_i\,u_\nu
=
\omega_{i\nu}(u_1,\dots,u_{n-m})
\ \ \ \ \ \ \ \ \ \ \ \ \ 
{\scriptstyle{(i\,=\,1\,\cdots\,n\,;\,\,\,
\nu\,=\,1\,\cdots\,n\,-\,m)}},
\]
and hence (cf. p.~\pageref{S-143}) the $\infty^{ n-m}$ 
$m$-times extended manifolds: 
\[
u_1
=
{\rm const.},
\,\,\,\dots,\,\,\,
u_{n-m}
=
{\rm const.}
\]
are mutually permuted by the group $X_1f, \dots, X_rf$. 

We therefore have the:

\def\theproposition{6}\begin{proposition}
\label{Satz-6-S-383}
Every simply transitive group: $X_1f, \dots, X_nf$ in $n > 1$
variables: $x_1, \dots, x_n$ is imprimitive; if $Z_1f, \dots, Z_n f$ is
the associated reciprocal group and if $Z_1f, \dots, Z_mf$ is an
arbitrary subgroup of it with the invariants $u_1, \dots, u_{ n-m}$,
then the $m$-term complete system $Z_1f = 0$, \dots, $Z_mf = 0$ admits
the group $X_1f, \dots, X_nf$ and the $\infty^{ n-m}$ $m$-times
extended manifolds: $u_1 = {\rm const.}$, \dots, $u_{ n-m} = {\rm
const.}$ are mutually permuted by this group.
\end{proposition}

The preceding proposition shows that every $m$-term subgroup of the
group $Z_1f, \dots, Z_nf$ provides a completely determined
decomposition of the space $x_1, \dots, x_n$ in $\infty^{ n-m}$
$m$-times extended manifolds invariant by the group $X_1f, \dots,
X_nf$. Hence, if we imagine that all subgroups of the group $Z_1f,
\dots, Z_nf$ are determined, and that all the associated invariants
are computed, then we obtain infinitely many decompositions of the
space invariant by the group $X_1f, \dots, X_nf$. It can be proved
\label{S-383}
that all existing invariant decompositions are found in this way, so
that the Proposition~6 can be reversed. 

Indeed, let $Y_1f = 0$, \dots, $Y_m f = 0$ be any $m$-term complete
system which admits the group $X_1f, \dots, X_nf$, and let $u_1,
\dots, u_{ n-m}$ be independent solutions of this complete system, so
that the family of the $\infty^{ n-m}$ $m$-times extended manifolds:
\[
u_1
=
{\rm const.},
\,\,\,\dots,\,\,\,
u_{n-m}
=
{\rm const.}
\] 
represent a decomposition of the space $x_1, \dots, x_n$ invariant by
the group $X_1f, \dots, X_nf$. We claim that the group $Z_1f, \dots,
Z_nf$ contains a completely determined $m$-term subgroup which leaves
individually invariant each one of these $\infty^{ n-m}$ manifolds;
this is just the reversion of the Proposition~6.

At first, we introduce the functions $u_1, \dots, u_{ n-m}$ and $m$
arbitrary mutually independent functions: $v_1, \dots, v_m$ of the $x$
that are also independent of $u_1, \dots, u_{ n-m}$ as variables in
our reciprocal groups, and we get:
\[
X_kf
=
\sum_{\nu=1}^{n-m}\,
\omega_{k\nu}(u_1,\dots,u_{n-m})\,
\frac{\partial f}{\partial u_\nu}
+
\sum_{\mu=1}^m\,
X_k\,v_\mu\,
\frac{\partial f}{\partial v_\mu}
\ \ \ \ \ \ \ \ \ \ \ \ \ {\scriptstyle{(k\,=\,1\,\cdots\,n)}}
\]
and:
\[
Z_kf
=
\sum_{\nu=1}^{n-m}\,
Z_k\,u_\nu\,
\frac{\partial f}{\partial u_\nu}
+
\sum_{\mu=1}^m\,
Z_k\,v_\mu\,
\frac{\partial f}{\partial v_\mu}
\ \ \ \ \ \ \ \ \ \ \ \ \ {\scriptstyle{(k\,=\,1\,\cdots\,n)}},
\]
where the $X_k\, v_\mu$, $Z_k\, u_\nu$ and $Z_k\, v_\mu$ are certain
functions of the $u$ and the $v$. Our claim now obviously amounts to
the fact that $m$ independent infinitesimal transformations should be
linearly deducible from $Z_1f, \dots, Z_nf$ which should absolutely
not transform $u_1, \dots, u_{ n-m}$, hence in which the coefficients
of $\partial f / \partial u_1$, \dots, $\partial f / \partial u_{
n-m}$ should all be zero.

In order to be able to prove this, we must compute the coefficients of
$\partial f / \partial u_1$, \dots, $\partial f / \partial u_{ n-m}$
in the general infinitesimal transformation $Zf = e_1\, Z_1f + \cdots
+ e_n \, Z_nf$ of the group $Z_1f, \dots, Z_nf$.

The infinitesimal transformation $Zf$ is
completely determined by the relations:
\[
X_1\big(Z(f)\big)
-
Z\big(X_1(f)\big)
=0,
\,\,\,\dots,\,\,\,
X_n\big(Z(f)\big)
-
Z\big(X_n(f)\big)
=
0.
\]
If we replace $f$ in these equations by 
$u_\nu$, we obtain the relations:
\[
X_1(Z\,u_\nu)
-
Z\,\omega_{1\nu}
=
0,
\,\,\,\dots,\,\,\,
X_n(Z\,u_\nu)
-
Z\,\omega_{n\nu}
=
0
\ \ \ \ \ \ \ \ \ \ \ \ \ 
{\scriptstyle{(\nu\,=\,1\,\cdots\,n\,-\,m)}}.
\]
Consequently, the functions \label{S-384}
$Z\, u_1, \dots, Z\, u_{ n- m}$
are solutions of the differential equations:
\def\theequation{7}\begin{equation}
X_k\,\rho_\nu
-
\sum_{\mu=1}^{n-m}\,
\frac{\partial\omega_{k\nu}}{\partial u_\mu}\,
\rho_\mu
=
0
\ \ \ \ \ \ \ \ \ \ \ \ \ 
{\scriptstyle{(k\,=\,1\,\cdots\,n\,;\,\,\,
\nu\,=\,1\,\cdots\,n\,-\,m)}},
\end{equation}
in which $\rho_1, \dots, \rho_{ n-m}$ are to be seen as the unknown
functions.

If:
\def\theequation{8}\begin{equation}
\rho_1
=
\psi_1
(u_1,\dots,u_{n-m},\,v_1,\dots,v_m),
\,\,\,\dots,\,\,\,
\rho_{n-m}
=
\psi_{n-m}
(u_1,\dots,u_{n-m},\,v_1,\dots,v_m)
\end{equation}
is an arbitrary system of solutions to the differential
equations~\thetag{ 7}, then all the $n ( n - m)$ expressions:
\[
X_k\,\psi_\nu
-
\sum_{\mu=1}^{n-m}\,
\frac{\partial\omega_{k\nu}}{\partial u_\mu}\,
\rho_\mu
\]
vanish identically after the substitution: $\rho_1 = \psi_1$, \dots,
$\rho_{ n-m} = \psi_{ n-m}$. Hence if we interpret the
equations~\thetag{ 8} as a system of equations in the $n + n - m$
variables: $u_1, \dots, u_{ n-m}$, $v_1, \dots, v_m$, $\rho_1, \dots,
\rho_{ n-m}$, we realize immediately that this system of equations
admits the $n$ infinitesimal transformations:
\def\theequation{9}\begin{equation}
U_kf
=
X_kf
+
\sum_{\nu=1}^{n-m}\,
\bigg\{
\sum_{\mu=1}^{n-m}\,
\frac{\partial\omega_{k\nu}}{\partial u_\mu}\,\rho_\mu
\bigg\}\,
\frac{\partial f}{\partial\rho_\nu}
\ \ \ \ \ \ \ \ \ \ \ \ \ {\scriptstyle{(k\,=\,1\,\cdots\,n)}}.
\end{equation}
The converse is also clear: if one known an arbitrary system of
equations of the form~\thetag{ 8} which admits the infinitesimal
transformations $U_1f, \dots, U_nf$, then one also knows a system of
solutions of the differential equations~\thetag{ 7}, since the
functions $\psi_1, \dots, \psi_{ n-m}$ are such a system.

From this, it results that the determination of the most general
system of solutions to the differential equations~\thetag{ 7} amounts
to determining, in the $2\, n - m$ variables $u$, $v$, $\rho$, the most
general system of equations~\thetag{ 8} that admits the infinitesimal
transformations $U_1f, \dots, U_nf$. 

Every system of equations which admits $U_1f, \dots, U_nf$ also 
allows all the infinitesimal transformations
$\leftbracket U_i,\, U_k \rightbracket$. 
By a calculation, we find: 
\[
\aligned
\leftbracket
U_i,\,U_k
\rightbracket
=
\leftbracket
X_i,\,X_k
\rightbracket
&
+
\sum_{\mu,\,\,\nu,\,\,\pi}^{1\cdots\,n-m}\,
\bigg(
\omega_{i\pi}\,
\frac{\partial^2\omega_{k\nu}}{\partial u_\mu\partial u_\pi}
-
\omega_{k\pi}\,
\frac{\partial^2\omega_{i\nu}}{\partial u_\mu\partial u_\pi}
\bigg)\,\rho_\mu\,
\frac{\partial f}{\partial\rho_\nu}
\\
&
+
\sum_{\mu,\,\,\nu,\,\,\pi}^{1\cdots\,n-m}\,
\bigg(
\frac{\partial\omega_{i\nu}}{\partial u_\mu}\,
\frac{\partial\omega_{k\pi}}{\partial u_\nu}
-
\frac{\partial\omega_{k\nu}}{\partial u_\mu}\,
\frac{\partial\omega_{i\pi}}{\partial u_\nu}
\bigg)\,\rho_\mu\,
\frac{\partial f}{\partial\rho_\pi},
\endaligned
\]
or else, if written differently:
\[
\leftbracket
U_i,\,U_k
\rightbracket
=
\leftbracket
X_i,\,X_k
\rightbracket
+
\sum_{\mu,\,\,\nu}^{1\cdots\,n-m}\,
\frac{\partial}{\partial u_\mu}\,
\big(
X_i\,\omega_{k\nu}-X_k\,\omega_{i\nu}
\big)\,
\rho_\mu\,
\frac{\partial f}{\partial\rho_\nu}.
\]
Now, since $X_1f, \dots, X_nf$ generate a group, there 
exist relations of the form: 
\[
\leftbracket
X_i,\,X_k
\rightbracket
=
X_i\big(X_k(f)\big)
-
X_k\big(X_i(f)\big)
=
\sum_{\sigma=1}^n\,
c_{ik\sigma}\,X_\sigma f,
\]
hence in particular, we have: 
\[
\aligned
X_i\big(X_k\,u_\nu\big)
-
X_k\big(X_i\,u_\nu\big)
&
=
X_i\,\omega_{k\nu}
-
X_k\,\omega_{i\nu}
\\
&
=
\sum_{\sigma=1}^n\,c_{ik\sigma}\,\omega_{\sigma\nu}, 
\endaligned
\]
and it comes:
\[
\leftbracket
U_i,\,U_k
\rightbracket
=
\sum_{\sigma=1}^n\,c_{ik\sigma}\,U_\sigma f.
\]
Consequently, the infinitesimal transformations $U_1f, \dots, U_nf$
generate an $n$-term group in the $2\, n - m$ variables $u_1, \dots,
u_{ n-m}$, $v_1, \dots, v_m$, $\rho_1, \dots, \rho_{ n-m}$.

Thus, the question at present is to determine the most general system
of equations of the form~\thetag{ 8} that admits the group $U_1f,
\dots, U_nf$. Since $X_1f, \dots, X_nf$ are linked together by no
linear relation, not all $n \times n$ determinants of the matrix which
is associated to $U_1f, \dots, U_nf$ vanish, and likewise also, not
all these determinants vanish by virtue of a system of equations of
the form~\thetag{ 8}. It follows therefore from
Chap.~\ref{kapitel-14}, Theorem~42, p.~\pageref{Theorem-42-S-237} that
every system of equations of the form~\thetag{ 8} which admits the
group $U_1f, \dots, U_nf$ is represented by relations between the
solutions of the $n$-term complete system $U_1f = 0$, \dots, $U_nf =
0$. Now, this complete system possesses exactly $n - m$
independent solutions, say:
\[
\Psi_\mu
(u_1,\dots,u_{n-m},\,v_1,\dots,v_m,\,\rho_1,\dots,\rho_{n-m})
\ \ \ \ \ \ \ \ \ \ \ \ \ {\scriptstyle{(k\,=\,1\,\cdots\,n\,-\,m)}}
\]
and to be precise, $\Psi_1, \dots, \Psi_{ n-m}$ are mutually
independent relatively to $\rho_1, \dots, \rho_{ n-m}$, because the
complete system is solvable with respect to the differential quotients
$\partial f / \partial u_1$, \dots, $\partial f / \partial u_{ n-m}$,
$\partial f / \partial v_1$, \dots, $\partial f / \partial v_m$
(cf. Chap.~\ref{kapitel-5}, Theorem~12, p.~\pageref{Theorem-12-S-91}).
From this, we conclude that the most general system of
equations~\thetag{ 8} that admits the group $U_1f, \dots, U_nf$ can be
given the form:
\[
\Psi_1=C_1,
\,\,\,\dots,\,\,\,
\Psi_{n-m}=C_{n-m},
\]
where $C_1, \dots, C_{ n-m}$ denote arbitrary constants. 

If we solve the system of equations just found
with respect to $\rho_1, \dots, \rho_{ n-m}$, which 
is always possible, then we obtain the most
general system of solutions to the differential equations~\thetag{ 7}, 
hence we see that this most general system of solutions
contains exactly $n - m$ arbitrary, essential constants. 
Now, since the differential equations~\thetag{ 7} are
linear and homogeneous in the unknowns $\rho_1, \dots, \rho_{ n-m}$, 
the said most general system of solutions must 
receive the form:
\def\theequation{10}\begin{equation}
\aligned
\rho_\mu
=
C_1'\,\psi_\mu^{(1)}
&
+
C_2'\,\psi_\mu^{(2)}
+\cdots+
C_{n-m}'\,\psi_\mu^{(n-m)}
\\
&\ \ \
{\scriptstyle{(\mu\,=\,1\,\cdots\,n\,-\,m)}},
\endaligned
\end{equation}
where the $n - m$ systems of functions:
\[
\psi_1^{(\nu)},\ \
\psi_2^{(\nu)},
\,\,\,\dots,\,\,\,
\psi_{n-m}^{(\nu)}
\ \ \ \ \ \ \ \ \ \ \ \ \ 
{\scriptstyle{(\nu\,=\,1\,\cdots\,n\,-\,m)}}
\]
represent the same number $n - m$ of linearly independent
systems of solutions to the differential
equations~\thetag{ 7} that are free of arbitrary 
constants, and where the $C'$ are arbitrary 
constants. Of course, the determinant
of the $\psi$ does not vanish, for the equations~\thetag{ 10}
must be solvable with respect to $C_1', \dots, C_{ n-m}'$. 

The functions: $Z\, u_1, \dots, Z\, u_{ n-m}$ are, 
according to p.~\pageref{S-384}, solutions
of the differential equations~\thetag{ 7}, hence they have the form:
\[
Z\,u_\mu
=
\overline{C}_1\,\psi_\mu^{(1)}
+\cdots+
\overline{C}_{n-m}\,\psi_\mu^{(n-m)}
\ \ \ \ \ \ \ \ \ \ \ \ \ 
{\scriptstyle{(\mu\,=\,1\,\cdots\,n\,-\,m)}}.
\]
Here, the $\overline{ C}_\nu$ are constants about which we temporarily
do not know anything more precise; it could still be thinkable that
they are linked together by linear relations.

From the values of the $Z\, u_\mu$, it follows that the general
infinitesimal transformation $Zf$ of the group: $Z_1f, \dots, Z_nf$
can be given the following representation:
\[
Zf
=
\sum_{\mu,\,\,\nu}^{1\cdots\,n-m}\,
\overline{C}_\mu\,\psi_\nu^{(\mu)}\,
\frac{\partial f}{\partial u_\nu}
+
\sum_{\mu=1}^m\,Z\,v_\mu\,
\frac{\partial f}{\partial v_\mu}.
\]
Now, since the $n - m$ expressions:
\[
\sum_{\nu=1}^{n-m}\,
\psi_\nu^{(\mu)}\,
\frac{\partial f}{\partial u_\nu}
\ \ \ \ \ \ \ \ \ \ \ \ \ 
{\scriptstyle{(\mu\,=\,1\,\cdots\,n\,-\,m)}}
\]
represent independent infinitesimal transformations, then from $Z_1f,
\dots, Z_nf$, one can obviously deduce linearly at least $m$, hence
say exactly $m + \varepsilon$, independent infinitesimal
transformations in which the coefficients of $\partial f / \partial
u_1$, \dots, $\partial f / \partial u_{ n-m}$ are equal to zero.
These $m + \varepsilon$ infinitesimal transformations naturally
generate an $(m+ \varepsilon)$-term subgroup of the group: $Z_1f,
\dots, Z_nf$, and in fact, a subgroup which leaves fixed each one of
the $\infty^{ n-m}$ $m$-times extended manifolds $u_1 = {\rm const.}$,
\dots, $u_{ n-m} = {\rm const.}$ But this is only possible when the
entire number $\varepsilon$ is equal to zero, since if $\varepsilon$
would be $>0$, the group $Z_1f, \dots, Z_nf$ could not be simply
transitive.

As a result, the claim made on p.~\pageref{S-383} is proved
and we can therefore state the following proposition:

\def\theproposition{7}\begin{proposition}
\label{Satz-7-S-387}
If $X_1f, \dots, X_nf$ and $Z_1f, \dots, Z_nf$ are two reciprocal
simply transitive groups in the $n$ variables $x_1, \dots, x_n$, and
if:
\[
u_1(x_1,\dots,x_n)
=
{\rm const.},
\,\,\,\dots,\,\,\,
u_{n-m}(x_1,\dots,x_n)
=
{\rm const.}
\]
is an arbitrary decomposition of the space $x_1, \dots, x_n$ in
$\infty^{ n-m}$ $m$-times extended manifolds that is invariant by the
group $X_1f, \dots, X_nf$, then the group $Z_1f, \dots, Z_nf$ always
contains an $m$-term subgroup which leaves individually fixed each one
of these $\infty^{ n-m}$ manifolds.
\end{proposition}

By combining this proposition with the Proposition~6, 
p.~\pageref{Satz-6-S-383}, we obtain the:

\def\thetheorem{69}\begin{theorem}
\label{Theorem-69-S-387}
If the $n$-term group $X_1f, \dots, X_nf$ in the $n$ variables $x_1,
\dots, x_n$ is simply transitive, then one finds all $m$-term complete
systems that this groups admits, or, what is the same, all invariant
decompositions of the space $x_1, \dots, x_n$ in $\infty^{ n-m}$
$m$-times extended manifolds, in the following way: One determines at
first the simply transitive group: $Z_1f, \dots, Z_nf$ which is
reciprocal to $X_1f, \dots, X_nf$ and one sets up all $m$-term
subgroups of the former; if, say:
\[
{\sf Z}_\mu f
=
g_{\mu 1}\,Z_1f
+\cdots+
g_{\mu n}\,Z_nf
\ \ \ \ \ \ \ \ \ \ \ \ \ {\scriptstyle{(\mu\,=\,1\,\cdots\,m)}}
\]
is one of the found subgroups, then the equations ${\sf Z}_1f = 0$,
\dots, ${\sf Z}_m f = 0$ represent one of the sought complete systems
and they determine a decomposition of the space $x_1, \dots, x_n$ in
$\infty^{ n-m}$ $m$-times extended manifolds that is invariant by the
group $X_1f, \dots, X_nf$; if, for each one of the found subgroups,
one forms the $m$-term complete system which the subgroup provides,
then one obtains all $m$-term complete systems that the group $X_1f,
\dots, X_nf$ admits. If one undertakes the indicated study for each
one of the numbers $m = 1, 2, \dots, n-1$, then one actually obtains
all complete systems that the group $X_1f, \dots, X_nf$ admits, 
and therefore at the same time, all decompositions of
the space $x_1, \dots, x_n$ that are invariant by this group.
\end{theorem}

The above theorem contains a solution to the problem of determining
all possible ways in which a given \emphasis{simply transitive} group
can be imprimitive.

\medskip

\label{S-388-sq}
At present, let the equations:
\[
u_1
=
{\rm const.},
\,\,\,\dots,\,\,\,
u_{n-m}
=
{\rm const.}
\]
again represent an arbitrary decomposition of the space $x_1, \dots,
x_n$ in $\infty^{ n-m}$ $m$-times extended manifolds invariant by the
group $X_1f, \dots, X_nf$, so that hence, when $u_1, \dots, u_{ n-m}$,
together with appropriate functions $v_1, \dots, v_m$, are introduced
as new variables, $X_1f, \dots, X_nf$ receive the form:
\[
X_kf
=
\sum_{\nu=1}^{n-m}\,
\omega_{k\nu}(u_1,\dots,u_{n-m})\,
\frac{\partial f}{\partial u_\nu}
+
\sum_{\mu=1}^m\,
\overline{\xi}_{k\mu}
(u_1,\dots,u_{n-m},\,v_1,\dots,v_m)\,
\frac{\partial f}{\partial v_\mu}.
\]

Here, not all $(n - m) \times (n - m)$ determinants of the matrix:
\[
\label{S-388}
\left\vert
\begin{array}{cccc}
\omega_{11} & \,\cdot\, & \,\cdot\, & \omega_{1,\,n-m}
\\
\cdot & \,\cdot\, & \,\cdot\, & \cdot
\\
\omega_{n1} & \,\cdot\, & \,\cdot\, & \omega_{n,\,n-m}
\end{array}
\right\vert
\]
can vanish identically, since otherwise, $X_1f, \dots, X_nf$ would be
linked together by a linear relation, which is contrary to the
assumption. So, if by $u_1^0, \dots, u_{ n-m}^0$, we understand a
general system of values, then according to Chap.~\ref{kapitel-13},
p.~\pageref{S-222-bis}, the group $X_1f, \dots, X_nf$ contains exactly
$\infty^{ m-1}$ different infinitesimal transformations which leave
invariant the system of equations:
\[
u_1
=
u_1^0,
\,\,\,\dots,\,\,\,
u_{n-m}
=
u_{n-m}^0\,;
\]
naturally, these infinitesimal transformations then generate an
$m$-term \label{S-388-bis}
subgroup, the most general subgroup of the group $X_1f,
\dots, X_nf$ which leaves fixed the $m$-times extended manifold $u_1 =
u_1^0$, \dots, $u_{ n-m} = u_{ n-m}^0$, or shortly $M$.

On the other hand, the said $m$-times extended manifold $M$ now admits
also $m$ independent infinitesimal transformations of the reciprocal
group $Z_1f, \dots, Z_nf$, for according to Proposition~7,
p.~\pageref{Satz-7-S-387}, this group contains an $m$-term subgroup
which leaves individually fixed each one of the $\infty^{ n-m}$
manifolds:
\[
u_1
=
{\rm const.},
\,\,\,\dots,\,\,\,
u_{n-m}
=
{\rm const.}
\]
In addition, it is clear that $M$ cannot admit a larger subgroup of
the group $Z_1f, \dots, Z_nf$, for $u_1^0, \dots, u_{ n-m}^0$ is
supposed to be a general system of values.

From this, we see that $M$ allows exactly $m$ infinitesimal
transformations both from the two reciprocal groups $X_1f, \dots,
X_nf$ and $Z_1f, \dots, Z_nf$, hence from each one, a completely
determined $m$-term subgroup.

We can make these two $m$-term subgroups visible by choosing an
arbitrary system of values: $v_1^0, \dots, v_m^0$ in general position
and by expanding the infinitesimal transformations of our two
reciprocal groups in powers of: $u_1 - u_1^0$, \dots, $u_{ n-m} - u_{
n-m}^0$, $v_1 - v_1^0$, \dots, $v_m - v_m^0$. Indeed, if we disregard
the terms of first order and of higher order, then similarly as on
p.~\pageref{S-379}, we can replace the infinitesimal transformations:
$X_1f, \dots, X_nf$ by $n$ transformations of the form:
\[
\aligned
&
\mathfrak{X}_1f
=
\frac{\partial f}{\partial u_1}
+\cdots,
\,\,\,\dots,\,\,\,
\mathfrak{X}_{n-m}f
=
\frac{\partial f}{\partial u_{n-m}}
+\cdots,
\\
&
\mathfrak{X}_{n-m+1}f
=
\frac{\partial f}{\partial v_1}
+\cdots,
\,\,\,\dots,\,\,\,
\mathfrak{X}_nf
=
\frac{\partial f}{\partial v_m}
+\cdots,
\endaligned
\]
and we can also replace $Z_1f, \dots, Z_nf$ by $n$ transformations
of the form:
\[
\aligned
&
\mathfrak{Z}_1f
=
-\,\frac{\partial f}{\partial u_1}
+\cdots,
\,\,\,\dots,\,\,\,
\mathfrak{Z}_{n-m}f
=
-\,\frac{\partial f}{\partial u_{n-m}}
+\cdots,
\\
&
\mathfrak{Z}_{n-m+1}f
=
-\,\frac{\partial f}{\partial v_1}
+\cdots,
\,\,\,\dots,\,\,\,
\mathfrak{Z}_nf
=
-\,\frac{\partial f}{\partial v_m}
+\cdots.
\endaligned
\]
Here, $\mathfrak{ X}_{ n-m+1}f, \dots, \mathfrak{ X}_nf$ are obviously
independent infinitesimal transformations which leave invariant the
manifold: $u_1 = u_1^0$, \dots, $u_{ n-m} = u_{ n-m}^0$, and
$\mathfrak{ Z}_{ n-m + 1}f, \dots, \mathfrak{ Z}_nf$ are independent
infinitesimal transformations which do the same; thus, $\mathfrak{
X}_{ n-m+1}f, \dots, \mathfrak{ X}_nf$ and $\mathfrak{ Z}_{ n-m+1}f,
\dots, \mathfrak{ Z}_nf$ are the two $m$-term subgroups about which we
have spoken just now.

From this representation of the two subgroup, 
one can yet derive a few noticeable conclusions. 

Between $\mathfrak{ X}_1f, \dots, \mathfrak{ X}_nf$, there are
relations of the form:
\[
\leftbracket
\mathfrak{X}_i,\,\mathfrak{X}_k
\rightbracket
=
\sum_{\nu=1}^n\,\mathfrak{c}_{ik\nu}\,\mathfrak{X}_\nu f,
\]
and according to p.~\pageref{S-379} sq., between 
$\mathfrak{ Z}_1f, \dots, \mathfrak{ Z}_nf$, there are 
the same relations:
\[
\leftbracket
\mathfrak{Z}_i,\,\mathfrak{Z}_k
\rightbracket
=
\sum_{\nu=1}^n\,\mathfrak{c}_{ik\nu}\,
\mathfrak{Z}_\nu f,
\]
with the same constants $\mathfrak{ c}_{ ik\nu}$. From this, it comes
that the two subgroups: $\mathfrak{ X}_{ n - m + 1}f, \dots,
\mathfrak{ X}_nf$ and $\mathfrak{ Z}_{ n-m +1}f, \dots, \mathfrak{
Z}_nf$ are equally composed. But it even comes more. Indeed, since
the two groups: $X_1f, \dots, X_nf$ and: $Z_1f, \dots, Z_nf$ are
related to each other in a holoedrically isomorphic way, when to every
infinitesimal transformation: $e_1\, \mathfrak{ X}_1f + \cdots + e_r\,
\mathfrak{ X}_rf$ is associated the infinitesimal transformation:
$e_1\, \mathfrak{ Z}_1 f + \cdots + e_r\, \mathfrak{ Z}_rf$, and since
through this association, the subgroup: $\mathfrak{ X}_{ n-m+1} f,
\dots, \mathfrak{ X}_rf$ obviously corresponds to the subgroup:
$\mathfrak{ Z}_{ n-m+1}f, \dots, \mathfrak{ Z}_rf$, then it becomes
evident that the two reciprocal groups can be related to each other in
a holoedrically isomorphic way so that the two $m$-term subgroups
which leave $M$ invariant correspond to each other.

If we summarize the results of the pp.~\pageref{S-388-sq} sq., we then
have the:

\def\theproposition{8}\begin{proposition}
\label{Satz-8-S-390}
If an $m$-times extended manifold of the space $x_1, \dots, x_n$
admits exactly $m$ independent infinitesimal transformations, and
hence an $m$-term subgroup, of a simply transitive group $X_1f, \dots,
X_nf$ of this space, then it admits at the same time exactly $m$
independent infinitesimal transformations, and hence an $m$-term
subgroup, of the simply transitive group: $Z_1f, \dots, Z_nf$ which is
reciprocal to $X_1f, \dots, X_nf$. The two $m$-term subgroups defined
in this way are equally composed and it is possible to relate the two
simply transitive \label{S-390}
reciprocal groups in a holoedrically isomorphic way
so that these $m$-term subgroups correspond to each other.
\end{proposition}

\sectionengellie{\S\,\,\,98.}

\label{S-390-bis}
The largest portion of the results of the preceding paragraph can be
derived by means of simple conceptual considerations. We shall now
undertake this, and at the same time, we shall yet gain a few further
results.

As up to now, let $X_1f, \dots, X_nf$ and $Z_1f, \dots, Z_nf$ be two
reciprocal simply transitive groups in the variables $x_1, \dots,
x_n$; moreover, let $Z_1f, \dots, Z_m f$ be again an arbitrary
$m$-term subgroup of the group $Z_1f, \dots, Z_nf$ and let $u_1,
\dots, u_{ n-m}$ be its invariants. Let the letter $S$ be the general
symbol of a transformation of the group $Z_1f, \dots, Z_mf$, and
lastly, let $T$ be an arbitrarily chosen transformation of the group
$X_1f, \dots, X_nf$.

Now, if $P$ is any point of the space $x_1, \dots, x_n$, then:
\[
(P')
=
(P)\,S
\]
is the general symbol of a point on the manifold:
\[
u_1
=
{\rm const.},
\,\,\,\dots,\,\,\,
u_{n-m}
=
{\rm const.}
\]
which passes through the point $P$. Furthermore, since the
transformations $S$ and $T$ are mutually interchangeable, we have:
\[
(P')\,T
=
(P)\,S\,T
=
(P)\,T\,S,
\]
hence, if we denote by $\Pi$ the point $(P) \, T$, we have:
\def\theequation{11}\begin{equation}
(P')\,T
=
(\Pi)\,S.
\end{equation}
Here, $(\Pi)\, S$ is the general symbol of a point on the manifold:
$u_\nu = {\rm const.}$ passing though $\Pi$. Consequently, our
symbolic equation~\thetag{ 11} says that the transformation $T$, hence
actually every transformation of the group $X_1f, \dots, X_nf$,
permutes the $\infty^{ n-m}$ manifolds: $u_1 = {\rm const.}$, \dots,
$u_{ n-m} = {\rm const.}$, by transferring each one of these manifolds
to a manifold of the same family.

As a result, the Proposition~6, p.~\pageref{Satz-6-S-383}, is derived.

But we can also prove the converse of this proposition by means of
such conceptual considerations.

Let us imagine that an arbitrary decomposition of the space $x_1,
\dots, x_n$ in $\infty^{ n-m}$ $m$-times extended manifolds invariant
by the group $X_1f, \dots, X_nf$ is given, and let us assume that $M$
is one of these $\infty^{ n-m}$ manifolds. By $P$ and $P'$, let us
understand two arbitrary points of $M$, and by $T$, an arbitrary
transformation of the group $X_1f, \dots, X_nf$.

By the execution of $T$, let the point 
$P$ be transferred to $\Pi$, so we have:
\[
(\Pi)
=
(P)\,T\,;
\]
on the other hand, there always is in the reciprocal group 
$Z_1f, \dots, Z_nf$ one, and also only one, transformation
which transfers $P$ to $P'$: 
\[
(P')
=
(P)\,S.
\]
Because of:
\[
(P)\,S\,T
=
(P)\,T\,S,
\]
we also have:
\def\theequation{12}\begin{equation}
(P')\,T
=
(\Pi)\,S.
\end{equation}

We must attempt to interpret this equation.

To begin with, we assume that the point $\Pi$ also belongs to the
manifold $M$. In this case, the transformation $T$ has the property
of leaving $M$ invariant. Indeed, $T$ permutes the mentioned
$\infty^{ n-m}$ manifolds, but on the other hand, $T$ transfers a
point of $M$, namely $P$, to a point of $M$, namely the point $\Pi$;
thus, $T$ must transfer all points of $M$ to points of $M$, hence it
must leave $M$ invariant. Now, since $P'$ lies on $M$, then the point
$(P') \, T$ also belongs to the manifold $M$, and because of~\thetag{
12}, the point $(\Pi)\, S$ too; but by appropriate choice of $T$,
$\Pi$ can be an arbitrary point of $M$, hence the transformation $S$
also transfers every point of $M$ to a point of $M$, so that it also
leaves invariant the manifold $M$.

On the other hand, we can suppose that $\Pi$ is any point of an
arbitrary other manifold amongst the $\infty^{ n-m}$ in question; if
we make this assumption, then we immediately realize from~\thetag{ 12}
that $S$ leaves at rest the concerned manifold.

With these words, it is proved that the group $Z_1f, \dots, Z_nf$
contains a transformation, namely the transformation $S$, which leaves
individually fixed each one of our $\infty^{ n-m}$ manifolds. But
there are evidently $\infty^m$ different such transformations, since
after fixing $P$, the point $P'$ can yet be chosen inside $M$ in
$\infty^m$ different ways. However, there are no more than $\infty^m$
transformations of this sort in the group $Z_1f, \dots, Z_nf$, because
this group is simply transitive; consequently, the $\infty^m$ existing
transformations generate an $m$-term subgroup of the group $Z_1f,
\dots, Z_nf$.

As a result, the Proposition~7, p.~\pageref{Satz-7-S-387}, is proved.

Obviously, the manifold $M$ admits, aside from the $\infty^m$
transformations of the group $Z_1f, \dots, Z_nf$, also yet $\infty^m$
transformations of the group $X_1f, \dots, X_nf$, which in turn form
an $m$-term subgroup of this group. This is a result which
is stated in Proposition~8, p.~\pageref{Satz-8-S-390}.

Something essentially new arises when the number 
$m$ in the developments just conducted is chosen
equal to $n$. Up to now, this case did not come into
consideration, because to it, there corresponds
no decomposition
of the space $x_1, \dots, x_n$.

If $m$ is equal to $n$, then the manifold $M$
coincides with the space $x_1, \dots, x_n$ itself; hence
$P$ and $P'$ are arbitrary points of the space, 
and by appropriate choice of $P$ and $P'$, 
$S$ can be any transformation of the reciprocal 
group $Z_1f, \dots, Z_nf$. 

If we choose $P$ and $P'$ fixed, the transformation $S$ is
completely determined; next, if $T$ is an arbitrary 
transformation of the group $X_1f, \dots, X_nf$, we have: 
\[
(P)\,S\,T
=
(P)\,T\,S,
\]
or, because $(P') = (P)\, S$: 
\def\theequation{13}\begin{equation}
(P')\,T
=
(P)\,T\,S.
\end{equation}

Here, by an appropriate choice of the transformation $T$, the
point $(P)\, T$ can be brought to coincidence with any 
arbitrary point of the space $x_1, \dots, x_n$; 
the same holds for the point $(P')\, T$. Consequently, 
thanks to the equation~\thetag{ 13}, 
we are in a position to indicate, for every point
$\mathfrak{ P}$ of the space, the new position
$\mathfrak{ P}'$ that it gets by the transformation $S$; 
we need only to determine the transformation ${\sf T}$ of
the group $X_1f, \dots, X_nf$ which transfers
$P$ to $\mathfrak{ P}$, hence which satisfies
the symbolic equation:
\[
(P)\,{\sf T}
=
(\mathfrak{P}).
\]
Then we have:
\[
(\mathfrak{P}')
=
(\mathfrak{P})\,S
=
(P)\,{\sf T}\,S,
\]
whence:
\[
(\mathfrak{P}')
=
(P')\,{\sf T}.
\]

Now, if we let the point $\mathfrak{ P}$ take all possible positions,
or, what is the same, if we set for ${\sf T}$ one after the other all
transformations of the group $X_1f, \dots, X_nf$, then we obtain that,
to every point of the space $x_1, \dots, x_n$ is associated a
completely determined other point, hence we obtain a transformation of
the space $x_1, \dots, x_n$, namely the transformation $S$. Lastly, if
we yet choose the points $P$ and $P'$ in all possible ways, then we
obviously obtain all transformations of the group $Z_1f, \dots, Z_nf$.

Thus, we can say:

\plainstatement{
If two points $P$ and $P'$ of the space $x_1, \dots, x_n$ are
transformed in cogredient way \deutsch{in cogredienter Weise} by means
of each one of the $\infty^n$ transformations of a simply transitive
group $X_1f, \dots, X_nf$, then the transformation which transfers
each one of the $\infty^n$ positions taken by $P$ to the corresponding
position of the point $P'$ belongs to the group $Z_1f, \dots, Z_nf$
reciprocal to the group $X_1f, \dots, X_nf$. If one chooses the
points $P$ and $P'$ in all possible ways, then one obtains all
transformations of the group $Z_1f, \dots, Z_nf$.}

It goes without saying that here, the points $P$ and $P'$ are always
to be understood as points in general position, or said more
precisely, points which lie on no manifold invariant by the group
$X_1f, \dots, X_nf$.

If one knows the finite equations of the group $X_1f, \dots, X_nf$,
then one can use the construction just found for the transformations
of the group $Z_1f, \dots, Z_nf$ in order to set up the finite
equations of this group.

Let:
\[
x_i' 
= 
f_i(x_1,\dots,x_n,\,a_1,\dots,a_n)
\ \ \ \ \ \ \ \ \ \ \ \ \ {\scriptstyle{(i\,=\,1\,\cdots\,n)}}
\]
be the finite equations of the group $X_1f, \dots, X_nf$. If one
calls the coordinates of the point say $x_1^0, \dots, x_n^0$, and
those of the point $P'$ say $u_1^0, \dots, u_n^0$, then by the
$\infty^n$ transformations of the group $X_1f, \dots, X_nf$, $P$
receives the $\infty^n$ different positions:
\[
y_i
=
f_i(x_1^0,\dots,x_n^0,\,a_1,\dots,a_n)
\ \ \ \ \ \ \ \ \ \ \ \ \ {\scriptstyle{(i\,=\,1\,\cdots\,n)}}
\]
and $P'$ receives the $\infty^n$ positions:
\[
y_i'
=
f_i(u_1^0,\dots,u_n^0,\,a_1,\dots,a_n)
\ \ \ \ \ \ \ \ \ \ \ \ \ {\scriptstyle{(i\,=\,1\,\cdots\,n)}}.
\]
Every system of values of the $a$ provides positions for $P$ and $P'$
which correspond to each other; hence if we eliminate from the
equations $y_i = f_i ( x^0, a)$ and $y_i' = f_i ( u^0, a)$ the
parameters $a$, we obtain the equations:
\def\theequation{14}\begin{equation}
y_i'
=
\mathfrak{F}_i
(y_1,\dots,y_n,\,x_1^0,\dots,x_n^0,\,u_1^0,\dots,u_n^0)
\ \ \ \ \ \ \ \ \ \ \ \ \ {\scriptstyle{(i\,=\,1\,\cdots\,n)}}
\end{equation}
of a transformation of the group $Z_1f, \dots, Z_nf$, namely the
transformation which transfers the point $P$ to $P'$. Lastly, if we
let $x_1^0, \dots, x_n^0$ and $u_1^0, \dots, u_n^0$ take all possible
values, we obtain all transformations of the group $Z_1f, \dots, Z_nf$.

In the equations~\thetag{ 14} of the group $Z_1f, \dots, Z_nf$ just
found, there appear $2\, n$ arbitrary parameters; however, this is
only fictitious, only $n$ of these parameters are essential. Indeed,
one can obtain every individual transformation of the group $Z_1f,
\dots, Z_nf$ in $\infty^n$ different ways, since one can always choose
arbitrarily the point $P$, whereas the point $P'$ is determined by the
concerned transformation after the fixed choice of $P$.

From this, it follows that one can also derive in this way all
transformations of the group $Z_1f, \dots, Z_nf$, by choosing the
point $P$ fixed once for all, and only by letting the point $P'$ take
all possible positions; that is to say, one can insert determined
numbers for the quantities $x_1^0, \dots, x_n^0$ and one only needs to
interpret $u_1^0, \dots, u_n^0$ as arbitrary parameters.

Thus, the following holds:

\def\theproposition{9}\begin{proposition}
If the finite equations:
\[
x_i'
=
f_i(x_1,\dots,x_n,\,a_1,\dots,a_n)
\ \ \ \ \ \ \ \ \ \ \ \ \ {\scriptstyle{(i\,=\,1\,\cdots\,n)}}
\]
of a simply transitive $n$-term group of the space $x_1, \dots, x_n$
are presented, then one finds the equations of the reciprocal simply
transitive group in the following way:

In the equations:
\[
y_i
=
f_i(x_1^0,\dots,x_n^0,\,a_1,\dots,a_n)
\ \ \ \ \ \ \ \ \ \ \ \ \ {\scriptstyle{(i\,=\,1\,\cdots\,n)}},
\]
one confers a fixed value to the $x^0$ and afterwards, one eliminates
the $n$ quantities $a_1, \dots, a_n$ from these equations and from the
equations:
\[
y_i'
=
f_i(u_1^0,\dots,u_n^0,\,a_1,\dots,a_n)
\ \ \ \ \ \ \ \ \ \ \ \ \ {\scriptstyle{(i\,=\,1\,\cdots\,n)}}\,;
\]
the resulting equations:
\[
y_i'
=
\mathfrak{F}_i(y_1,\dots,y_n,\,
x_1^0,\dots,x_n^0,\,u_1^0,\dots,u_n^0)
\ \ \ \ \ \ \ \ \ \ \ \ \ {\scriptstyle{(i\,=\,1\,\cdots\,n)}}
\]
with the $n$ arbitrary parameters are the equations
of the reciprocal group. The assumption is that
the $x_k^0$ are chosen so that the 
point $x_1^0, \dots, x_n^0$ lies on no manifold
which remains invariant by the group: 
$x_i' = f_i ( x, a)$. \label{S-395-bis}
\end{proposition}

\sectionengellie{\S\,\,\,99.}

The Proposition~7, p.~\pageref{Satz-7-S-387} is a special case of a
general proposition that also holds true for certain groups which are
not simply transitive. We now want to derive this general
proposition; on the occasion, we obtain at the same time a new proof
for the Proposition~7.

Let $X_1f, \dots, X_nf$ be an $n$-term group in the variables
$x_1, \dots, x_n$ and let the number $n$ be not larger than
$s$. In addition, we make the assumption that 
$X_1f, \dots, X_nf$ are linked together by no
linear relation of the form:
\[
\chi_1(x_1,\dots,x_s)\,X_1f
+\cdots+
\chi_n(x_1,\dots,x_s)\,X_nf
=
0.
\]

The proposition to be proved amounts to the following: if one knows an
arbitrary
$m$-term complete system: $\mathfrak{ Y}_1 f = 0$, \dots, $\mathfrak{
Y}_m f = 0$ that the group $X_1f, \dots, X_nf$ admits, then one can
always bring the complete system to a form: $Y_1f = 0$, \dots, $Y_mf =
0$ such that \label{S-395}
\emphasis{the infinitesimal transformations $Y_1f, \dots,
Y_mf$ are all interchangeable with $X_1f, \dots, X_nf$}, and that in
the relations:
\[
\leftbracket
Y_i,\,Y_k
\rightbracket
=
\sum_{\nu=1}^m\,
\tau_{ik\nu}(x_1,\dots,x_s)\,Y_\nu f
\]
which hold between $Y_1f, \dots, Y_mf$, the $\tau_{ ik \nu}$ are all
solutions of the $n$-term complete system: $X_1f = 0$, \dots, $X_nf =
0$.

In the special case $s = n$, where the group $X_1f, \dots, X_nf$ is
simply transitive, the infinitesimal transformations $Y_1f, \dots,
Y_mf$ obviously belong to the simply transitive group $Z_1f, \dots,
Z_nf$ reciprocal to $X_1f, \dots, X_nf$; moreover, since the $n$-term
complete system $X_1f = 0$, \dots, $X_nf = 0$ possesses in this case
no other solutions than $f = {\rm const.}$, the functions $\tau_{
ik\nu}$ are then plain constants so that $Y_1f, \dots, Y_mf$ generate
an $m$-term subgroup of the group $Z_1f, \dots, Z_nf$. We therefore
have the Proposition~7, p.~\pageref{Satz-7-S-387}.

However, we deal at present with the general case. 

We therefore imagine that an $m$-term complete system:
\[
\mathfrak{Y}_1f=0,
\,\,\,\dots,\,\,\,
\mathfrak{Y}_mf=0
\]
is presented which admits the group $X_1f, \dots, X_nf$,
so that the $Xf$ and $\mathfrak{ Y}f$ are linked together
by relations of the form:
\def\theequation{15}\begin{equation}
\leftbracket
\mathfrak{Y}_\mu,\,X_k
\rightbracket
=
\sum_{\nu=1}^m\,\alpha_{\mu k\nu}(x_1,\dots,x_s)\,
\mathfrak{Y}_\nu f
\ \ \ \ \ \ \ \ \ \ \ \ \ 
{\scriptstyle{(k\,=\,1\,\cdots\,n\,;\,\,\,
\mu\,=\,1\,\cdots\,m)}}
\end{equation}
(cf. Chap.~\ref{kapitel-13}, p.~\pageref{S-221-bis}).

To begin with, we now attempt to determine $m$ functions
$\rho_1, \dots, \rho_m$ of the $x$ so that the infinitesimal
transformation: 
\[\label{S-396}
Yf
=
\sum_{\mu=1}^m\,\rho_\mu(x_1,\dots,x_s)\,
\mathfrak{Y}_\mu f
\]
is interchangeable with all the $n$ infinitesimal transformations
$X_kf$. We therefore have to satisfy the $n$ equations:
\[
\aligned
\leftbracket
X_k,\,Y
\rightbracket
&
=
\sum_{\mu=1}^m\,
X_k\,\rho_\mu\,\mathfrak{Y}_\mu f
+
\sum_{\mu=1}^m\,
\rho_\mu\,
\leftbracket
X_k,\,\mathfrak{Y}_\mu
\rightbracket
\\
&
=
\sum_{\mu=1}^m\,
\bigg\{
X_k\,\rho_\mu
-
\sum_{\nu=1}^m\,\alpha_{\nu k\mu}\,\rho_\nu
\bigg\}\,
\mathfrak{Y}_\mu f
=
0,
\endaligned
\]
or, because $\mathfrak{ Y}_1f, \dots, \mathfrak{ Y}_m f$
cannot be linked together by linear relations, 
the following $m \, n$ relations:
\def\theequation{16}\begin{equation}
X_k\,\rho_\mu
-
\sum_{\nu=1}^m\,\alpha_{\nu k\mu}\,\rho_\nu
=
0
\ \ \ \ \ \ \ \ \ \ \ \ \ 
{\scriptstyle{(k\,=\,1\,\cdots\,n\,;\,\,\,
\mu\,=\,1\,\cdots\,m)}}.
\end{equation}
These are differential equations by means of which the $\rho$ are 
to be determined.

If:
\def\theequation{17}\begin{equation}
\rho_1
=
P_1(x_1,\dots,x_s),
\,\,\,\dots,\,\,\,
\rho_m
=
P_m(x_1,\dots,x_s)
\end{equation}
is a system of solutions of the differential equations~\thetag{ 16}, 
then the equations:
\[
X_k\,P_\mu
-
\sum_{\nu=1}^m\,\alpha_{\nu k\mu}\,\rho_\nu
=
0
\]
are satisfied identically after the substitution: 
$\rho_1 = P_1$, \dots, $\rho_m = P_m$.
Expressed differently: the system of equations~\thetag{ 17}
in the $s + m$ variables $x_1, \dots, x_s$, $\rho_1, \dots, \rho_m$
admits the infinitesimal transformations:
\[
\mathfrak{U}_kf
=
X_kf
+
\sum_{\mu=1}^m\,
\bigg\{
\sum_{\nu=1}^m\,\alpha_{\nu k\mu}\,\rho_\nu
\bigg\}\,
\frac{\partial f}{\partial\rho_\mu}
\ \ \ \ \ \ \ \ \ \ \ \ \ {\scriptstyle{(k\,=\,1\,\cdots\,n)}}.
\]
Conversely, if a system of equations of the form~\thetag{ 17}
admits the infinitesimal transformations 
$\mathfrak{ U}_1f, \dots, \mathfrak{ U}_nf$, then the
functions $P_1, \dots, P_m$ are obviously solutions
of the differential equations~\thetag{ 16}.

From this, we see that the integration of the differential 
equations~\thetag{ 16} amounts to determining, in the
$s + m$ variables $x$, $\rho$, the most
general system of equations~\thetag{ 17}, 
which admits the infinitesimal transformations
$\mathfrak{ U}_1f, \dots, \mathfrak{ U}_nf$. 

The system of equations to be determined also admits
the infinitesimal transformations: $\leftbracket
\mathfrak{ U}_i,\, \mathfrak{ U}_k \rightbracket$. 

By calculation, we find: 
\[
\aligned
\leftbracket
\mathfrak{U}_i,\,\mathfrak{U}_k
\rightbracket
=
\leftbracket
X_i,\,X_k
\rightbracket
&
+
\sum_{\mu,\,\,\nu}^{1\cdots\,m}
\big\{
X_i\,\alpha_{\nu k\mu}
-
X_k\,\alpha_{\nu i\mu}
\big\}\,
\rho_\nu\,
\frac{\partial f}{\partial\rho_\mu}
\\
&
+
\sum_{\mu,\,\,\nu,\,\,\pi}^{1\cdots\,m}\,
\big\{
\alpha_{\nu i\pi}\,\alpha_{\pi k\mu}
-
\alpha_{\nu k\pi}\,\alpha_{\pi i\mu}
\big\}\,
\rho_\nu\,
\frac{\partial f}{\partial\rho_\mu}.
\endaligned
\]
Here, in order to simplify the right-hand side, 
we form the Jacobi identity (cf. Chap.~\ref{kapitel-5}, 
p.~\pageref{jacobi-identity}):
\[
\big\leftbracket
\mathfrak{Y}_\nu,\,
\leftbracket 
X_i,\,X_k
\rightbracket
\big\rightbracket
+
\big\leftbracket
X_i,\,
\leftbracket 
X_k,\,\mathfrak{Y}_\nu
\big\rightbracket
+
\rightbracket
\big\leftbracket
X_k,\,
\leftbracket 
\mathfrak{Y}_\nu,\,X_i
\rightbracket
\big\rightbracket
=
0,
\]
which, using~\thetag{ 15}, can be written as:
\[
\big\leftbracket
\mathfrak{Y}_\nu,\,
\leftbracket 
X_i,\,X_k
\rightbracket
\big\rightbracket
=
\sum_{\mu=1}^m
\big\{
\leftbracket
X_if,\,\alpha_{\nu k\mu}\,\mathfrak{Y}_\mu f
\rightbracket
-
\leftbracket
X_kf,\,\alpha_{\nu i\mu}\,\mathfrak{Y}_\mu f
\rightbracket
\big\},
\]
or:
\[
\aligned
\big\leftbracket
\mathfrak{Y}_\nu,\,
\leftbracket 
X_i,\,X_k
\rightbracket
\big\rightbracket
&
=
\sum_{\mu=1}^m\,
\big\{
X_i\,\alpha_{\nu k\mu}
-
X_k\,\alpha_{\nu i\mu}
\big\}\,\mathfrak{Y}_\mu f
\\
&
\ \ \ \ \
+
\sum_{\mu,\,\,\pi}^{1\cdots\,m}\,
\big\{
\alpha_{\nu i\mu}\,\alpha_{\mu k\pi}
-
\alpha_{\nu k\mu}\,\alpha_{\mu i\pi}
\big\}\,\mathfrak{Y}_\pi f.
\endaligned
\]
But now, $X_1f, \dots, X_nf$ generate an $n$-term group, whence: 
\[
\leftbracket
X_i,\,X_k
\rightbracket
=
\sum_{\sigma=1}^n\,c_{ik\sigma}\,X_\sigma f,
\]
from which it follows:
\[
\big\leftbracket
\mathfrak{Y}_\nu,\,
\leftbracket 
X_i,\,X_k
\rightbracket
\big\rightbracket
=
\sum_{\mu=1}^m\,
\bigg\{
\sum_{\sigma=1}^n\,c_{ik\sigma}\,\alpha_{\nu\sigma\mu}
\bigg\}\,\mathfrak{Y}_\mu f.
\]
If we still take into account that $\mathfrak{ Y}_1f, 
\dots, \mathfrak{ Y}_mf$ are linked together by no 
linear relations, we obtain: 
\[
\aligned
X_i\,\alpha_{\nu k\mu}
-
X_k\,\alpha_{\nu i\mu}
&
+
\sum_{\pi=1}^m\,
\big\{
\alpha_{\nu i\pi}\,\alpha_{\pi k\mu}
-
\alpha_{\nu k\pi}\,\alpha_{\pi i\mu}
\big\}
\\
&
=
\sum_{\sigma=1}^n\,c_{ik\sigma}\,\alpha_{\nu\sigma\mu}.
\endaligned
\]
By inserting these values in the expression
found above for $\leftbracket \mathfrak{ U}_i, \, 
\mathfrak{ U}_k \rightbracket$, it comes:
\[
\leftbracket
\mathfrak{U}_i,\,\mathfrak{U}_k
\rightbracket
=
\sum_{\sigma=1}^n\,
c_{ik\sigma}\,\mathfrak{U}_\sigma f,
\]
so $\mathfrak{ U}_1 f, \dots, \mathfrak{ U}_n f$ generate an $n$-term
group in the $s + m$ variables $x$, $\rho$.

At present, the question is to determine the most general system of
equations~\thetag{ 17} which admits the $n$-term group $\mathfrak{
U}_1f, \dots, \mathfrak{ U}_nf$.

Since $X_1f, \dots, X_nf$ are linked together by no linear relation,
not all $n \times n$ determinants vanish identically in the matrix
associated to $\mathfrak{ U}_1f, \dots, \mathfrak{ U}_nf$; but
obviously, these $n \times n$ determinants cannot all be equal to
zero, even by virtue of a system of equations of the form~\thetag{
17}. Consequently, every system of equations of the form~\thetag{ 17}
which admits the group $\mathfrak{ U}_1f, \dots, \mathfrak{ U}_nf$ is
represented by relations between the solutions of the $n$-term
complete system: $\mathfrak{ U}_1 f = 0$, \dots, $\mathfrak{ U}_nf =
0$ (cf. Chap.~\ref{kapitel-14}, Theorem~42,
p.~\pageref{Theorem-42-S-237}).

The $n$-term complete system $\mathfrak{ U}_1f = 0$, 
\dots, $\mathfrak{ U}_n f = 0$ possesses $s + m - n$ independent
solutions; $s - n$ of these solutions can be chosen in such 
a way that they are free of the $\rho$ and 
depend only upon the $x$, they are the independent
solutions of the $n$-term complete system: 
$X_1f = 0$, \dots, $X_nf = 0$, and we can call them:
\[
\mathfrak{v}_1(x_1,\dots,x_s),
\,\,\dots,\,\,
\mathfrak{v}_{s-n}(x_1,\dots,x_s).
\]
Furthermore, let: 
\[
\Omega_1(\rho_1,\dots,\rho_m,\,x_1,\dots,x_s),
\,\,\dots,\,\,
\Omega_m(\rho_1,\dots,\rho_m,\,x_1,\dots,x_s)
\]
be $m$ arbitrary mutually independent solutions of the complete
system: $\mathfrak{ U}_1 f = 0$, \dots, $\mathfrak{ U}_nf = 0$ that
are also independent of the $\mathfrak{ v}$.

Since the equations: $\mathfrak{ U}_1f = 0$, \dots, $\mathfrak{ U}_nf
= 0$ are solvable with respect to $n$ of the differential quotients
$\partial f / \partial x_1$, \dots, $\partial f / \partial x_s$, the
$s - n + m$ functions $\mathfrak{ v}_1, \dots, \mathfrak{ v}_{ s-n}$,
$\Omega_1, \dots, \Omega_m$ must be mutually independent relatively to
$s - n$ of the variables $x_1, \dots, x_s$, and relatively to $\rho_1,
\dots, \rho_m$. Consequently, $\Omega_1, \dots, \Omega_m$ are
mutually independent relatively to $\rho_1, \dots, \rho_m$.

\medskip
k
After these preparations, we can identify the
most general system of equations which admits the
group $\mathfrak{ U}_1f, \dots, \mathfrak{ U}_nf$ and
which can be given the form~\thetag{ 17}. 

The concerned system of equations consists of $m$ relations
between $\mathfrak{ v}_1, \dots, \mathfrak{ v}_{ s-n}$, 
$\Omega_1, \dots, \Omega_m$ and is solvable with 
respect to $\rho_1, \dots, \rho_m$; 
consequently, it is solvable with respect to $\Omega_1, \dots, 
\Omega_m$ and has the form:
\def\theequation{18}\begin{equation}
\Omega_\mu
(\rho_1,\dots,\rho_m,\,x_1,\dots,x_s)
=
\Psi_\mu
\big(
\mathfrak{v}_1(x),\,\dots,\,
\mathfrak{v}_{s-n}(x)
\big)
\ \ \ \ \ \ \ \ \ \ \ \ \ 
{\scriptstyle{(\mu\,=\,1\,\cdots\,m)}},
\end{equation}
where the $\Psi$ are absolutely arbitrary functions of their
arguments. If we solve this system of equations with respect to
$\rho_1, \dots, \rho_m$, then for $\rho_1, \dots, \rho_m$, we receive
expressions that represent the most general system of solutions to the
differential equations~\thetag{ 16}.

Thus, the most general system of solutions~\thetag{ 16} contains $m$
arbitrary functions of $\mathfrak{ v}_1, \dots, \mathfrak{ v}_{
s-n}$. Now, since the equations~\thetag{ 16} in the unknowns $\rho_1,
\dots, \rho_m$ are linear and homogeneous, one can conclude that the
most general system of solutions $\rho_1, \dots, \rho_m$ to it can be
deduced from $m$ particular systems of solutions:
\[
P_1^{(\mu)}(x_1,\dots,x_s),
\,\,\dots,\,\,
P_m^{(\mu)}(x_1,\dots,x_s)
\ \ \ \ \ \ \ \ \ \ \ \ \ 
{\scriptstyle{(\mu\,=\,1\,\cdots\,m)}}
\]
in the following way:
\def\theequation{19}\begin{equation}
\aligned
\rho_\nu
=
\chi_1(\mathfrak{v}_1,\dots,\mathfrak{v}_{s-n})\,
&
P_\nu^{(1)}
+\cdots+
\chi_m(\mathfrak{v}_1,\dots,\mathfrak{v}_{s-n})\,
P_\nu^{(m)}
\\
&
{\scriptstyle{(\mu\,=\,1\,\cdots\,m)}},
\endaligned
\end{equation}
where it is understood that the $\chi$ are arbitrary functions of
their arguments. Naturally, the $m$ particular systems of solutions
must be constituted in such a way that it is not possible to determine
$m$ functions $\psi_1, \dots, \psi_m$ of $\mathfrak{ v}_1, \dots,
\mathfrak{ v}_{ s-n}$ such that the $m$ equations:
\[
\sum_{\mu=1}^m\,
\psi_\mu(\mathfrak{v}_1,\dots,\mathfrak{v}_{s-n})\,
P_\nu^{(\mu)}
=
0
\ \ \ \ \ \ \ \ \ \ \ \ \ 
{\scriptstyle{(\nu\,=\,1\,\cdots\,m)}}
\]
are satisfied simultaneously.

Lastly, it is still to be observed that the equations~\thetag{ 19} are
solvable with respect to $\chi_1, \dots, \chi_m$, so that the
determinant of the $P_\nu^{ (\mu)}$ does not vanish
identically. Indeed, according to what precedes, one must be able to
bring the equations~\thetag{ 19} to the form~\thetag{ 18}, and this is
obviously impossible when the determinant of the $P_\nu^{ (\mu)}$
vanishes.

\medskip

At present, we can write down the most general transformation: 
\[
Yf
=
\sum_{\mu=1}^m\,\rho_\mu\,\mathfrak{Y}_\mu f
\]
which is interchangeable with $X_1f, \dots, X_nf$ (cf. above,
p.~\pageref{S-396}). Its form is:
\[
Yf
=
\sum_{\nu=1}^m\,
\chi_\nu
(\mathfrak{v}_1,\dots,\mathfrak{v}_{s-n})\,
\sum_{\mu=1}^m\,
P_\mu^{(\nu)}
(x_1,\dots,x_s)\,
\mathfrak{Y}_\mu f,
\]
or:
\[
Yf
=
\sum_{\nu=1}^m\,\chi_\nu
(\mathfrak{v}_1,\dots,\mathfrak{v}_{s-n})\,Y_\nu f,
\]
if we set: 
\[
\sum_{\mu=1}^m\,P_\mu^{(\nu)}\,
\mathfrak{Y}_\mu f
=
Y_\nu f
\ \ \ \ \ \ \ \ \ \ \ \ \ 
{\scriptstyle{(\nu\,=\,1\,\cdots\,m)}}.
\]

Here obviously, $Y_1f, \dots, Y_mf$ are all interchangeable with
$X_1f, \dots, X_nf$; furthermore, the equations: $Y_1f = 0$, \dots,
$Y_m f = 0$ are equivalent to the equations: $\mathfrak{ Y}_1f = 0$,
\dots, $\mathfrak{ Y}_m f = 0$, and hence they in turn form an
$m$-term complete system that admits the group $X_1f, \dots, X_nf$.

As a consequence of that, between $Y_1f, \dots, Y_mf$, there 
are relations of the form:
\[
\leftbracket
Y_\mu,\,Y_\nu
\rightbracket
=
\sum_{\pi=1}^m\,
\tau_{\mu\nu\pi}(x_1,\dots,x_s)\,
Y_\pi f.
\]
But from the Jacobi identity:
\[
\big\leftbracket
X_k,\,
\leftbracket
Y_\mu,\,Y_\nu
\rightbracket
\big\rightbracket
+
\big\leftbracket
Y_\mu,\,
\leftbracket
Y_\nu,\,X_k
\rightbracket
\big\rightbracket
+
\big\leftbracket
Y_\nu,\,
\leftbracket
X_k,\,Y_\mu
\rightbracket
\big\rightbracket
\equiv
0\,,
\]
in which the last two terms are identically zero, it
results immediately:
\[
\big\leftbracket
X_k,\,
\leftbracket
Y_\mu,\,Y_\nu
\rightbracket
\big\rightbracket
\equiv
0.
\]
We must therefore have:
\[
\sum_{\pi=1}^m\,X_k\,\tau_{\mu\nu\pi}\,Y_\pi f
\equiv 0,
\]
or, because $Y_1f, \dots, Y_mf$ are linked together by 
no linear relation:
\[
X_k\,\tau_{\mu\nu\pi}
\equiv
0
\ \ \ \ \ \ \ \ \ \ \ \ \ {\scriptstyle{(k\,=\,1\,\cdots\,n)}},
\]
that is to say: the $\tau_{ \mu \nu \pi}$ are all solutions of the
complete system: $X_1f = 0$, \dots, $X_nf = 0$, they are functions of
$\mathfrak{ v}_1, \dots, \mathfrak{ v}_{ s-n}$ alone.

The complete system: $Y_1f = 0$, \dots, $Y_mf = 0$ possesses all the
properties indicated on p~\pageref{S-395}; we can therefore enunciate
the statement:

\renewcommand{\thefootnote}{\fnsymbol{footnote}}
\def\thetheorem{70}\begin{theorem}
If an $m$-term complete system: $\mathfrak{ Y}_1f = 0$, \dots,
$\mathfrak{ Y}_mf = 0$ in the $s$ variables $x_1, \dots, x_s$ admits
the $n$-term group $X_1f, \dots, X_nf$ and if this group is
constituted in such a way that between $X_1f, \dots, X_nf$, there is
no linear relation of the form:
\[
\chi_1(x_1,\dots,x_s)\,X_1f
+\cdots+
\chi_n(x_1,\dots,x_s)\,X_nf
=
0,
\]
then it is possible to determine $m^2$ functions $P_\mu^{ (\nu)}
(x_1, \dots, x_s)$ with not identically vanishing determinant
such that the $m$ infinitesimal transformations:
\[
Y_\nu f
=
\sum_{\mu=1}^m\,P_\mu^{(\nu)}(x_1,\dots,x_s)\,
\mathfrak{Y}_\mu f
\ \ \ \ \ \ \ \ \ \ \ \ \ {\scriptstyle{(\nu\,=\,1\,\cdots\,m)}}
\]
are all interchangeable with $X_1f, \dots, X_nf$. In turn, the
equations: $\mathfrak{ Y}_1f = 0$, \dots, $\mathfrak{ Y}_mf = 0$
equivalent to $Y_1f = 0$, \dots, $Y_m f = 0$ form an $m$-term complete
system that admits the group $X_1f, \dots, X_nf$. Lastly, between
$Y_1f, \dots, Y_mf$, there are relations of the specific form:
\[
\leftbracket
Y_\mu,\,Y_\nu
\rightbracket
=
\sum_{\pi=1}^m\,
\vartheta_{\mu\nu\pi}
(\mathfrak{v}_1,\dots,\mathfrak{v}_{s-n})\,
Y_\pi f, 
\]
where $\mathfrak{ v}_1, \dots, \mathfrak{ v}_{ s-n}$ are independent
solutions of the $n$-term complete system: $X_1f = 0$, \dots, $X_nf =
0$.\footnote [1]{\,
During the summer semester 1887, in a lecture about the general theory
of integration of the differential equations that admit a finite
continuous group, \name{Lie} developed the Theorem~70 which comprises
the Theorem~69 as a special case.}
\end{theorem}
\renewcommand{\thefootnote}{\arabic{footnote}}

\linestop


\chapter{The Group of Parameters}
\label{kapitel-21}
\chaptermark{The Group of Parameters}

\setcounter{footnote}{0}

\abstract*{??}

If one executes three transformations of the space
$x_1, \dots, x_n$ subsequently, say the following ones: 
\def\theequation{1}\begin{equation}
\left\{
\aligned
x_i'
&
=
f_i(x_1,\dots,x_n)
\ \ \ \ \ \ \ \ \ \ \ \ \ {\scriptstyle{(i\,=\,1\,\cdots\,n)}}
\\
x_i''
&
=
g_i(x_1',\dots,x_n')
\ \ \ \ \ \ \ \ \ \ \ \ \ {\scriptstyle{(i\,=\,1\,\cdots\,n)}}
\\
x_i'''
&
=
h_i(x_1'',\dots,x_n'')
\ \ \ \ \ \ \ \ \ \ \ \ \ {\scriptstyle{(i\,=\,1\,\cdots\,n)}},
\endaligned\right.
\end{equation}
the one gets a new transformation:
\[
x_i'''
=
\omega_i(x_1,\dots,x_n)
\ \ \ \ \ \ \ \ \ \ \ \ \ {\scriptstyle{(i\,=\,1\,\cdots\,n)}}
\]
of the space $x_1, \dots, x_n$. 

The equations of the new transformation are obtained when the $2\, n$
variables: $x_1', \dots, x_n'$, $x_1'', \dots, x_n''$ are eliminated
from the $3\, n$ equations~\thetag{ 1}.
Clearly, one can execute this elimination in two 
different ways, since one can begin
either by taking away the $x'$, or by taking away the
$x''$. In the first case, one obtains firstly between 
the $x$ and the $x''$ the relations: 
\[
x_i''
=
g_i\big(f_1(x),\,\dots,\,f_n(x)\big)
\ \ \ \ \ \ \ \ \ \ \ \ \ {\scriptstyle{(i\,=\,1\,\cdots\,n)}}, 
\]
and afterwards, one has to insert these values of
$x_1'', \dots, x_n''$ in the equations:
\[
x_i'''
=
h_i(x_1'',\dots,x_n'')
\ \ \ \ \ \ \ \ \ \ \ \ \ {\scriptstyle{(i\,=\,1\,\cdots\,n)}}.
\]
In the second case, one obtains firstly between the
$x'$ and the $x'''$ the relations:
\[
x_i'''
=
h_i
\big(g_1(x'),\,\dots,\,g_n(x')\big)
\ \ \ \ \ \ \ \ \ \ \ \ \ {\scriptstyle{(i\,=\,1\,\cdots\,n)}},
\]
and then one has yet to replace the $x'$ by their values:
\[
x_i'
=
f_i(x_1,\dots,x_n)
\ \ \ \ \ \ \ \ \ \ \ \ \ {\scriptstyle{(i\,=\,1\,\cdots\,n)}}.
\]

The obvious observation that in both ways indicated above, one obtains
in the two cases the same transformation as a final result receives a
content \deutsch{Inhalt} when one interprets the transformation as an
operation and when one applies the symbolism of substitution theory.

For the three substitutions~\thetag{ 1}, we want to introduce the
symbols: $S$, $T$, $U$ one after the other, so that the
transformation:
\[
x_i''
=
g_i
\big(f_1(x),\,\dots,\,f_n(x)\big)
\ \ \ \ \ \ \ \ \ \ \ \ \ {\scriptstyle{(i\,=\,1\,\cdots\,n)}}
\]
receives the symbol: $S\, T$, the transformation:
\[
x_i'''
=
h_i
\big(g_1(x'),\,\dots,\,g_n(x')\big)
\ \ \ \ \ \ \ \ \ \ \ \ \ {\scriptstyle{(i\,=\,1\,\cdots\,n)}},
\]
the symbol: $T\, U$ and lastly, the transformation: $x_i''' = \omega_i
( x_1, \dots, x_n)$, the symbol: $S\, T\, U$. We can then express the
two ways of forming the transformation: $x_i''' = \omega_i ( x_1,
\dots, x_n)$ discussed above by saying that this transformation is
obtained both when one executes at first the transformation: $S\, T$
and next the transformation $U$, and when one executes at first the
transformation $S$, and next the transformation $T\, U$.

Symbolically, these facts can be expressed by the equations:
\def\theequation{2}\begin{equation}
(S\,T)\,U
=
S\,(T\,U)
=
S\,T\,U,
\end{equation}
and its known that they say that the operations $S$, $T$, $U$ satisfy
the so-called \terminology{associative rule} \deutsch{associative
Gesetz}.

We can therefore say:

\plainstatement{The transformations of an $n$-times
extended space are operations for which the law
of associativity holds true.}

At present, let us consider the specific case where
$S$, $T$, $U$ are arbitrary transformations of an $r$-term group:
\[
x_i'
=
f_i(x_1,\dots,x_n,\,a_1,\dots,a_r)
\ \ \ \ \ \ \ \ \ \ \ \ \ {\scriptstyle{(i\,=\,1\,\cdots\,n)}}.
\]
In this case, the transformations $S\, T$, $T\, U$
and $S\, T\, U$ naturally belong also to the
group in question. We now want to see what
can be concluded from the validity of
the associative law. 

\sectionengellie{\S\,\,\,100.}

\label{S-402-sq}
Let the transformations $S$, $T$, $U$ of our group be:
\def\theequation{3}\begin{equation}
\left\{
\aligned
S\colon\ \ \ \ \ \
x_i'
&
=
f_i(x_1,\dots,x_n,\,a_1,\dots,a_r)
\\
T\colon\ \ \ \ \
x_i''
&
=
f_i(x_1',\dots,x_n',\,b_1,\dots,b_r)
\\
U\colon\ \ \ \
x_i'''
&
=
f_i(x_1'',\dots,x_n'',\,c_1,\dots,c_r).
\endaligned\right.
\end{equation}

For the transformation $S\, T$, there result from 
this in the known way equations of the form:
\[
x_i''
=
f_i
\big(x_1,\dots,x_n,\,\,
\varphi_1(a,b),\,\dots,\,\varphi_r(a,b)\big),
\]
where the functions $\varphi_1 ( a, b)$, \dots, $\varphi_r ( a, b)$
are, according to Chap.~\ref{fundamental-differential},
p.~\pageref{A-2}, mutually independent both relatively to 
$a_1, \dots, a_r$ and relatively to 
$b_1, \dots, b_r$. 

Furthermore, the equations of the transformation
$(S\, T)\, U$ are: 
\def\theequation{4}\begin{equation}
x_i'''
=
f_i
\big(x_1,\dots,x_n,\,\,
\varphi_1\big(\varphi(a,b),\,c\big),
\,\,\dots,\,\,
\varphi_r\big(\varphi(a,b),\,c\big)\big).
\end{equation}

On the other hand, we have for the transformation 
$T\, U$ the equations: 
\[
x_i'''
=
f_i
\big(
x_1',\dots,x_n',\,\,
\varphi_1(b,c),\,\dots,\,\varphi_r(b,c)
\big),
\]
hence for $S\, (T\, U)$ the following ones:
\def\theequation{4'}\begin{equation}
x_i'''
=
f_i
\big(
x_1,\dots,x_n,\,\,
\varphi_1\big(a,\,\varphi(b,c)\big),
\,\,\dots,\,\,
\varphi_r\big(a,\,\varphi(b,c)\big)
\big).
\end{equation}

Now, it holds: $(S\, T)\, U = S\, (T\, U)$, whence the
two transformations~\thetag{ 4} and~\thetag{ 4'}
must be identical to each other; by comparison of the
two parameters in the two transformations, we
therefore obtain the following relations:
\[
\aligned
\varphi_k
\big(\varphi_1(a,b),
\,\,\dots,\,\,
\varphi_r(a,b),\,\,
&
c_1,\dots,c_r\big)
=
\varphi_k\big(
a_1,\dots,a_r,\,\,
\varphi_1(b,c),
\,\,\dots,\,\,
\varphi_r(b,c)\big)
\\
&
\ \ \ \ \ \ \ \ \ \ \ 
{\scriptstyle{(k\,=\,1\,\cdots\,r)}},
\endaligned
\]
or more shortly:
\def\theequation{5}\begin{equation}
\varphi_k
\big(\varphi(a,b),\,c\big)
=
\varphi_k\big(a,\,\varphi(b,c)\big)
\ \ \ \ \ \ \ \ \ \ \ \ \ {\scriptstyle{(k\,=\,1\,\cdots\,r)}}.
\end{equation}

Thus, 
the functions $\varphi_1, \dots, \varphi_r$ must identically 
satisfy these relations for all values of
the $a$, $b$, $c$. 

The equations~\thetag{ 5} express the law of associativity for three
arbitrary transformations of the group $x_i' = f_i ( x, a)$. But they
can yet be interpreted in another way; namely, they say that the
equations:
\def\theequation{6}\begin{equation}\label{S-404}
a_k'
=
\varphi_k(a_1,\dots,a_r,\,b_1,\dots,b_r)
\ \ \ \ \ \ \ \ \ \ \ \ \ {\scriptstyle{(k\,=\,1\,\cdots\,r)}}
\end{equation}
in the variables $a_1, \dots, a_r$ represent a group, 
and to be precise, a group with the $r$ parameters
$b_1, \dots, b_r$. 

Indeed, if we execute two transformations~\thetag{ 6}, say:
\[
a_k'
=
\varphi_k(a_1,\dots,a_r,\,b_1,\dots,b_r)
\]
and:
\[
a_k''
=
\varphi_k(a_1',\dots,a_r',\,c_1,\dots,c_r)
\]
one after the other, we obtain the transformation:
\[
a_k''
=
\varphi_k\big(
\varphi_1(a,b),
\,\,\dots,\,\,
\varphi_r(a,b),\,\,
c_1,\dots,c_r\big)
\]
which, by virtue of~\thetag{ 5}, takes the shape: 
\[
a_k''
=
\varphi_k
\big(
a_1,\dots,a_r,\,\,
\varphi_1(b,c),
\,\,\dots,\,\,
\varphi_r(b,c)\big),
\]
and hence belongs likewise to the transformations~\thetag{ 6}.
As a result, our assertion is proved. 

\plainstatement{We want to call the group:
\def\theequation{6}\begin{equation}
a_k'
=
\varphi_k(a_1,\dots,a_r,\,b_1,\dots,b_r)
\ \ \ \ \ \ \ \ \ \ \ \ \ {\scriptstyle{(k\,=\,1\,\cdots\,r)}}
\end{equation}
the \terminology{parameter group} \deutsch{Parametergruppe} of the
group $x_i' = f_i (x_1, \dots, x_n, \, a_1, \dots, a_r)$.}

According to the already cited page~\pageref{A-2}, 
the equations $a_k' = \varphi_k ( a, b)$ are solvable
with respect to the $r$ parameters $b_1, 
\dots, b_r$: 
\[
b_k
=
\psi_k
(a_1,\dots,a_r,\,
a_1',\dots,a_r')
\ \ \ \ \ \ \ \ \ \ \ \ \ {\scriptstyle{(k\,=\,1\,\cdots\,r)}}.
\]
From this, we conclude that the parameter group 
is $r$-term, is transitive, and in fact, 
simply transitive. In addition, the equations~\thetag{ 5}
show that the parameter group is its own
parameter group. 

We can therefore say in summary: 

\def\thetheorem{71}\begin{theorem}
\label{Theorem-71-S-404}
If the functions $f_i ( x_1, \dots, x_n, \, 
a_1, \dots, a_r)$ in the equations $x_i' = f_i ( x, a)$
of an $r$-term group satisfy the functional equations: 
\[
\aligned
f_i\big(
f_1(x,a),
\,\,\dots,\,\,
f_n(x,a),\,\,
&
b_1,\dots,b_r
\big)
=
f_i\big(
x_1,\dots,x_n,\,\,
\varphi_1(a,b),
\,\,\dots,\,\,
\varphi_r(a,b)
\big)
\\
&
\ \ \ \ \ \ \ \ \ \ \ \ \ 
{\scriptstyle{(i\,=\,1\,\cdots\,n)}},
\endaligned
\]
then the $r$ relations:
\[
a_i'
=
\varphi_i(a_1,\dots,a_r,\,b_1,\dots,b_r)
\ \ \ \ \ \ \ \ \ \ \ \ \ 
{\scriptstyle{(i\,=\,1\,\cdots\,r)}}
\]
determine an $r$-term group between the 
$2\, r$ variables $a$ and $a'$: the
parameter group of the original group. 
This parameter group is simply transitive and
is its own parameter group.
\end{theorem}

\sectionengellie{\S\,\,\,101.}

If the $r$-term group: $x_i' = f_i ( x_1, \dots, x_n,\, 
a_1, \dots, a_r)$
is generated by the $r$ 
independent infinitesimal transformations: 
\[
X_kf
=
\sum_{i=1}^n\,\xi_{ki}(x_1,\dots,x_n)\,
\frac{\partial f}{\partial x_i}
\ \ \ \ \ \ \ \ \ \ \ \ \ {\scriptstyle{(k\,=\,1\,\cdots\,r)}},
\]
then according to Chap.~\ref{kapitel-9},
p.~\pageref{Theorem-24-S-158}, its transformations are ordered
together as inverses by pairs. Evidently, the transformations of the
associated parameter group: $a_k' = \varphi_k ( a,b)$ are then also
ordered together as inverses by pairs, whence according to
Chap.~\ref{kapitel-9}, p.~\pageref{S-169} above, the parameter group
contains exactly $r$ independent infinitesimal transformations and is
generated by them.

We can also establish this important property of the parameter group
in the following way, where we find at the same time the infinitesimal
transformations of this group.

In the equations: $x_i' = f_i (x_1, \dots, x_n, \, a_1, \dots, a_r)$,
if the $x_i'$ are considered as functions of the $x$ and of the $a$,
then according to Theorem~3, p.~\pageref{Theorem-3-S-33}, there are
differential equations of the form:
\def\theequation{7}\begin{equation}\label{S-405}
\frac{\partial x_i'}{\partial a_k}
=
\sum_{j=1}^r\,\psi_{kj}(a_1,\dots,a_r)\,
\xi_{ji}(x_1',\dots,x_n')
\ \ \ \ \ \ \ \ \ \ \ \ \ 
{\scriptstyle{(i\,=\,1\,\cdots\,n\,;\,\,
k\,=\,1\,\cdots\,r)}},
\end{equation}
and they can also be written:
\def\theequation{7'}\begin{equation}
\xi_{ji}(x_1',\dots,x_n')
=
\sum_{k=1}^r\,
\alpha_{jk}(a_1,\dots,a_r)\,
\frac{\partial x_i'}{\partial a_k}
\ \ \ \ \ \ \ \ \ \ \ \ \ 
{\scriptstyle{(j\,=\,1\,\cdots\,r\,;\,\,
i\,=\,1\,\cdots\,n)}}.
\end{equation}

Here, according to the Chap.~\ref{fundamental-differential},
Proposition~\ref{economic-differential},
p.~\pageref{economic-differential},
the $\xi_{ ji}$ and the $\alpha_{ jk}$
are determined by the equations:
\[
\xi_{ji}(x_1',\dots,x_n')
=
\bigg[
\frac{\partial x_i'}{\partial b_j}
\bigg]_{b=b^0},
\ \ \ \ \ \ \ \ \ \ \ \ \ \ \
\alpha_{jk}(a_1,\dots,a_r)
=
\bigg[
\frac{\partial a_k}{\partial b_j}
\bigg]_{b=b^0},
\]
whose meaning has been explained in the place indicated.

For the group $a_k' = \varphi_k (a, b)$, 
one naturally obtains analogous differential equations:
\[
\frac{\partial a_i'}{\partial b_k}
=
\sum_{j=1}^r\,\psi_{kj}(b_1,\dots,b_r)\,
\overline{\psi}_{ji}(a_1',\dots,a_r')
\ \ \ \ \ \ \ \ \ \ \ \ \ 
{\scriptstyle{(i,\,\,k\,=\,1\,\cdots\,r)}}
\] 
in which the $\psi_{ kj}$ have the same signification
as in~\thetag{ 7}, while the 
$\overline{ \psi}_{ ji} (a)$ are determined by:
\[
\overline{\psi}_{ji}(a_1,\dots,a_r)
=
\bigg[
\frac{\partial a_i}{\partial b_j}
\bigg]_{b=b^0}.
\]
\emphasis{In consequence of that, the 
$\overline{ \psi}_{ ji}$ are the same functions
of their arguments as the $\alpha_{ ji}$}; 
in correspondence to the
formulas~\thetag{ 7} and~\thetag{ 7'}, 
we therefore obtain the two following ones: 
\def\theequation{8}\begin{equation}
\frac{\partial a_i'}{\partial b_k}
=
\sum_{j=1}^r\,\psi_{kj}(b_1,\dots,b_r)\,
\alpha_{ji}(a_1',\dots,a_r')
\ \ \ \ \ \ \ \ \ \ \ \ \ 
{\scriptstyle{(i,\,\,k\,=\,1\,\cdots\,r)}}
\end{equation}
and:
\def\theequation{8'}\begin{equation}
\alpha_{ji}
(a_1',\dots,a_r')
=
\sum_{k=1}^r\,
\alpha_{jk}(b_1,\dots,b_r)\,
\frac{\partial a_i'}{\partial b_k}
\ \ \ \ \ \ \ \ \ \ \ \ \ 
{\scriptstyle{(j,\,\,i\,=\,1\,\cdots\,r)}}.
\end{equation}

Now, the group: $x_i' = f_i ( x_1, \dots, x_n, \, 
a_1, \dots, a_r)$ contains the identity transformation and in fact, we
can always suppose that the parameters $a_1^0, \dots, a_r^0$ of the
identity transformation lie in the region $(\!( a )\!)$ defined on
p.~\pageref{S-16}. According to p.~\pageref{psi-k-j-e}\footnote{\,
According to the footnote, this determinant is equal to
$(-1)^r$. 
}, 
the determinant: 
\[
\sum\,\pm\,
\psi_{11}(a^0)\cdots\,\psi_{rr}(a^0)
\] 
is then certainly different from zero.

On the other hand, it is clear that, also in the family of the
transformations: $a_k' = \varphi_k ( a, b)$, the identity
transformation: $a_1' = a_1$, \dots, $a_r' = a_r$ appears, and even,
that the transformation:
\[
a_k'
=
\varphi_k(a_1,\dots,a_r,\,a_1^0,\dots,a_r^0)
\ \ \ \ \ \ \ \ \ \ \ \ \ {\scriptstyle{(k\,=\,1\,\cdots\,r)}}
\]
is the identity transformation. Thus, if we take into account the
existence of the differential equations~\thetag{ 8} and if we apply
the Theorem~9, p.~\pageref{Theorem-9-S-72}, we realize immediately
that the family of the $\infty^r$ transformations: $a_k' =
\varphi_k ( a, b)$ coincides with the family of the $\infty^{ r-1}$
one-term groups:
\[
\sum_{k=1}^r\,\lambda_k\,
\sum_{i=1}^r\,\alpha_{ki}
(a_1,\dots,a_r)\,
\frac{\partial f}{\partial a_i}.
\]

Consequently, the group: $a_k' = \varphi_k ( a, b)$ is
generated by the $r$ infinitesimal transformations:
\[
A_kf
=
\sum_{j=1}^r\,\alpha_{kj}(a_1,\dots,a_r)\,
\frac{\partial f}{\partial a_j}
\ \ \ \ \ \ \ \ \ \ \ \ \ {\scriptstyle{(k\,=\,1\,\cdots\,r)}}
\]

Naturally, these infinitesimal transformations are
independent of each other, since there exists
between them absolutely no relation of the form:
\[
\chi_1(a_1,\dots,a_r)\,A_1f
+\cdots+
\chi_r(a_1,\dots,a_r)\,A_rf
=
0,
\]
as we have already realized in Chap.~\ref{fundamental-differential},
Theorem~3, p.~\pageref{Theorem-3-S-33}; 
on the other hand, it can also be concluded from this that
the parameter group is simply 
transitive (Theorem~71). 

\medskip

In what precedes, we have not only proved that the
group: $a_k' = \varphi_k ( a, b)$ is generated
by infinitesimal transformations, but we have also found
these infinitesimal transformations themselves. 
From this, we can immediately deduce a new important
property of the group: $a_k' = \varphi_k ( a, b)$. 

The $r$ infinitesimal transformations: $X_1f, \dots, X_rf$ 
of the group: $x_i' = f_i ( x, a)$ are linked together
by relations of the form:
\[
\leftbracket
X_i,\,X_k
\rightbracket
=
\sum_{s=1}^r\,c_{iks}\,X_sf,
\]
and according to Theorem~21, p.~\pageref{Theorem-21-S-149}, 
there are, between $A_1f, \dots, A_rf$, 
relations of the same form:
\[
\leftbracket
A_i,\,A_k
\rightbracket
=
\sum_{s=1}^r\,c_{iks}\,A_sf, 
\]
with the same constants $c_{ iks}$. By applying the way of expressing
introduced in Chap.~\ref{kapitel-17}, p.~\ref{S-291-bis} and
p.~\pageref{S-293}, we can hence state the following theorem:

\renewcommand{\thefootnote}{\fnsymbol{footnote}}
\def\thetheorem{72}\begin{theorem}
\label{Theorem-72-S-407}
Every $r$-term group: 
\[
x_i'
=
f_i(x_1,\dots,x_n,\,a_1,\dots,a_r)
\ \ \ \ \ \ \ \ \ \ \ \ \ {\scriptstyle{(i\,=\,1\,\cdots\,n)}}
\]
is equally composed with its parameter group:
\[
a_k'
=
\varphi_k(a_1,\dots,a_r,\,b_1,\dots,b_r)
\ \ \ \ \ \ \ \ \ \ \ \ \ 
{\scriptstyle{(k\,=\,1\,\cdots\,r)}},
\]
or, what is the same, it is holoedrically isomorphic to it. 
If the infinitesimal transformations of the group: 
$x_i' = f_i ( x, a)$ satisfy the relations:
\[
\xi_{ji}(x_1',\dots,x_n')
=
\sum_{k=1}^r\,\alpha_{jk}(a_1,\dots,a_r)\,
\frac{\partial x_i'}{\partial a_k}
\ \ \ \ \ \ \ \ \ \ \ \ \ 
{\scriptstyle{(i\,=\,1\,\cdots\,n\,;\,\,\,
j\,=\,1\,\cdots\,r)}}, 
\]
then the $r$ expressions: 
\[
A_jf
=
\sum_{k=1}^r\,
\alpha_{jk}(a_1,\dots,a_r)\,
\frac{\partial f}{\partial a_k}
\ \ \ \ \ \ \ \ \ \ \ \ \ {\scriptstyle{(j\,=\,1\,\cdots\,r)}}
\]
represent $r$ independent infinitesimal transformations
of the parameter group.\footnote[1]{\,
\name{Lie}, Gesellschaft d. W. zu Christiania 1884, no.~15.
} 
\end{theorem}
\renewcommand{\thefootnote}{\arabic{footnote}}

\sectionengellie{\S\,102.}

As usual, let:
\[
X_kf
=
\sum_{i=1}^n\,\xi_{ki}(x_1,\dots,x_n)\,
\frac{\partial f}{\partial x_i}
\ \ \ \ \ \ \ \ \ \ \ \ \ {\scriptstyle{(k\,=\,1\,\cdots\,r)}}
\]
be independent infinitesimal transformations
of an $r$-term group in the variables $x_1, \dots, x_n$. 

Under the guidance of Chap.~\ref{kapitel-9}, 
p.~\pageref{S-156} and p.~\pageref{S-157}, 
one determines $r$ independent infinitesimal transformations:
\[
\overline{A}_kf
=
\sum_{j=1}^r\,
\overline{\alpha}_{kj}(a_1,\dots,a_r)\,
\frac{\partial f}{\partial a_j}
\ \ \ \ \ \ \ \ \ \ \ \ \ {\scriptstyle{(k\,=\,1\,\cdots\,r)}}
\]
which are linked together by no linear relation 
of the form:
\[
\chi_1(a_1,\dots,a_r)\,\overline{A}_1f
+\cdots+
\chi_r(a_1,\dots,a_r)\,\overline{A}_rf
=
0,
\]
but which stand pairwise in the relationships:
\[
\leftbracket\,
\overline{A}_i,\,\overline{A}_k
\rightbracket
=
\sum_{s=1}^r\,c_{iks}\,\overline{A}_sf\,;
\]
in other words, one determines $r$ independent infinitesimal
transformations of any simply transitive $r$-term group which has the
same composition as the group $X_1f, \dots, X_rf$.

Afterwards, one selects any system of values: $\overline{ a}_1, \dots,
\overline{ a}_r$ in the neighbourhood of which the $\overline{
\alpha}_{ ki} (a_1, \dots, a_r)$ behave regularly and for which the
determinant $\sum\, \pm \, \overline{ \alpha}_{ 11} \cdots\,
\overline{ \alpha}_{ rr}$ does not vanish, and then, relatively to
this system of values, one determines the general solutions of the
$r$-term complete system:
\def\theequation{9}\begin{equation}
\aligned
\sum_{i=1}^n\,\xi_{ki}(x')\,
\frac{\partial f}{\partial x_i'}
+
\sum_{j=1}^r\,
&
\overline{\alpha}_{kj}(a)\,
\frac{\partial f}{\partial a_j}
=
X_k'f
+
\overline{A}_kf
=
0
\\
&
\ \ \ \ \
{\scriptstyle{(k\,=\,1\,\cdots\,r)}}
\endaligned
\end{equation}
in the $n + r$ variables: $x_1', \dots, x_n'$, $a_1, \dots, a_r$.

Now, if $\overline{ F}_1 ( x_1', \dots, x_n', \, a_1, \dots, a_r)$,
\dots, $\overline{ F}_n ( x', \, a)$ are the general solutions in
question, one sets:
\[
x_i
=
\overline{F}_i
(x_1',\dots,x_n',\,a_1,\dots,a_r)
\ \ \ \ \ \ \ \ \ \ \ \ \ {\scriptstyle{(i\,=\,1\,\cdots\,n)}}\,;
\]
by resolution with respect to $x_1', \dots, x_n'$, one
obtains equations of the form:
\[
x_i'
=
\overline{f}_i
(x_1,\dots,x_n,\,a_1,\dots,a_r)
\ \ \ \ \ \ \ \ \ \ \ \ \ {\scriptstyle{(i\,=\,1\,\cdots\,n)}}
\]
which, according to Theorem~23, p.~\pageref{Theorem-23-S-154}, 
represent the finite transformations of an $r$-term group, 
namely of the group: $X_1f, \dots, X_rf$ itself. 

In the mentioned theorem, we already observed that the
equations: $x_i' = \overline{ f}_i ( x, a)$
can be brought to the form: 
\[
x_i'
=
x_i
+
\sum_{k=1}^r\,e_k\,X_k\,x_i
+
\sum_{k,\,\,j}^{1\cdots\,r}\,
\frac{e_k\,e_j}{1\cdot 2}\,
X_k\,X_j\,x_i
+\cdots
\ \ \ \ \ \ \ \ \ \ \ \ \ {\scriptstyle{(i\,=\,1\,\cdots\,n)}}
\]
after the introduction of new parameters, so that these equations
represent the finite equations of the $r$-term group which is
generated by the infinitesimal transformations: $X_1f, \dots, X_rf$.

It has no influence whether we choose the group $\overline{ A}_1f,
\dots, \overline{ A}_rf$ or any other group amongst the infinitely
many simply transitive groups which have the same composition as the
group $X_1f, \dots, X_rf$, and it is completely indifferent whether we
choose the system of values $\overline{ a}_1, \dots, \overline{ a}_r$
or any other system of values: in the way indicated, we always obtain
an analytic representation for the finite transformations of the
group: $X_1f, \dots, X_rf$.

We also made this remark already in Chap.~\ref{kapitel-9}, though
naturally not with the same words. But at present, we have reached a
point where we can see directly why one always obtains the equations
of the same group for various choices of the group: $\overline{ A}_1f,
\dots, \overline{ A}_rf$ and of the systems of values: $\overline{
a}_1, \dots, \overline{ a}_r$.

\medskip

Indeed, let:
\[
\mathfrak{A}_kf
=
\sum_{j=1}^r\,\beta_{kj}
(\mathfrak{a}_1,\dots,\mathfrak{a}_r)\,
\frac{\partial f}{\partial\mathfrak{a}_j}
\ \ \ \ \ \ \ \ \ \ \ \ \ {\scriptstyle{(k\,=\,1\,\cdots\,r)}}
\]
be any other simply transitive group
equally composed with the group: $X_1f, \dots, X_rf$ for
which the relations:
\[
\leftbracket
\mathfrak{A}_i,\,\mathfrak{A}_k
\rightbracket
=
\sum_{s=1}^r\,c_{iks}\,\mathfrak{A}_sf
\]
are identically satisfied, and moreover, let: $\mathfrak{ a}_1^0,
\dots, \mathfrak{ a}_r^0$ be any system of values in the neighbourhood
of which all $\beta_{ kj} ( \mathfrak{ a})$ behave regularly and for
which the determinant $\sum\, \pm\, \beta_{ 11} \cdots\, \beta_{ rr}$
does not vanish.

Because the two groups: $\overline{ A}_1 f, \dots, \overline{ A}_rf$
and $\mathfrak{ A}_1f, \dots, \mathfrak{ A}_rf$ are equally composed
and are both simply transitive, then (Theorem~64,
p.~\pageref{Theorem-64-S-340}) there are $\infty^r$ different
transformations:
\[
a_k
=
\lambda_k
\big(
\mathfrak{a}_1,\dots,\mathfrak{a}_r,\,C_1,\dots,C_r
\big)
\ \ \ \ \ \ \ \ \ \ \ \ \ {\scriptstyle{(k\,=\,1\,\cdots\,r)}}
\]
which transfer $\mathfrak{ A}_1f, \dots, \mathfrak{ A}_rf$ to:
$\overline{ A}_1f, \dots, \overline{ A}_rf$, respectively. The
equations of these transformations are solvable with respect to the
arbitrary parameters: $C_1, \dots, C_r$, so in particular, one can
choose $C_1, \dots, C_r$ in such a way that the equations:
\[
\overline{a}_k
=
\lambda_k(\mathfrak{a}_1^0,\dots,\mathfrak{a}_r^0,\,
C_1,\dots,C_r)
\ \ \ \ \ \ \ \ \ \ \ \ \ {\scriptstyle{(k\,=\,1\,\cdots\,r)}}
\]
are satisfied.

If: $C_1^0, \dots, C_r^0$ are the values of $C_1, \dots,
C_r$ obtained in this way, and if one sets:
\[
\lambda_k(\mathfrak{a}_1,\dots,\mathfrak{a}_r,\,
C_1^0,\dots,C_r^0)
=
\pi_k(\mathfrak{a}_1,\dots,\mathfrak{a}_r)
\ \ \ \ \ \ \ \ \ \ \ \ \ {\scriptstyle{(k\,=\,1\,\cdots\,r)}},
\]
then the equations:
\def\theequation{10}\begin{equation}
a_k
=
\pi_k(\mathfrak{a}_1,\dots,\mathfrak{a}_r)
\ \ \ \ \ \ \ \ \ \ \ \ \ {\scriptstyle{(k\,=\,1\,\cdots\,r)}}
\end{equation}
represent a transformation which transfers $\mathfrak{ A}_1f, \dots,
\mathfrak{ A}_rf$ to: $\overline{ A}_1f, \dots, \overline{ A}_rf$,
respectively, and which in addition transfers the system of values:
$\mathfrak{ a}_1^0, \dots, \mathfrak{ a}_r^0$ to the system of values:
$\overline{ a}_1, \dots, \overline{ a}_r$.

From this, it results that we obtain the general solutions of the
complete system:
\def\theequation{9'}\begin{equation}
X_k'f
+
\mathfrak{A}_kf
=
0
\ \ \ \ \ \ \ \ \ \ \ \ \ {\scriptstyle{(k\,=\,1\,\cdots\,r)}}
\end{equation}
relatively to: $\mathfrak{ a}_1 = \mathfrak{ a}_1^0$, \dots,
$\mathfrak{ a}_r = \mathfrak{ a}_r^0$ when we make the substitution:
$a_1 = \pi_1 ( \mathfrak{ a})$, \dots, $a_r = \pi_r ( \mathfrak{ a})$
in the general solutions:
\[
\overline{F}_i
(x_1',\dots,x_n',\,a_1,\dots,a_r)
\ \ \ \ \ \ \ \ \ \ \ \ \ {\scriptstyle{(i\,=\,1\,\cdots\,n)}}
\]
of the complete system~\thetag{ 9}. Therefore, if we had used,
instead of the simply transitive group: $\overline{ A}_1f, \dots,
\overline{ A}_rf$, the group: $\mathfrak{ A}_1f, \dots, \mathfrak{
A}_rf$, and if we had used, instead of the system of values:
$\overline{ a}_1, \dots, \overline{ a}_r$, the system of values:
$\mathfrak{ a}_1^0, \dots, \mathfrak{ a}_r^0$, then instead of the
equations: $x_i' = \overline{ f}_i (x_1, \dots, x_n, \, a_1, \dots,
a_r)$, we would have received the equations:
\[
x_i'
=
\overline{f}_i
\big(
x_1,\dots,x_n,\,
\pi_1(\mathfrak{a}),\,\dots,\,\pi_r(\mathfrak{a})\big)
\ \ \ \ \ \ \ \ \ \ \ \ \ {\scriptstyle{(i\,=\,1\,\cdots\,n)}}.
\]
But evidently, these equations are transferred to the former ones when
the new parameters $a_1, \dots, a_r$ are introduced in place of
$\mathfrak{ a}_1, \dots, \mathfrak{ a}_r$ by virtue of the
equations~\thetag{ 10}.

\medskip

As in the paragraphs 100 and~101, p.~\pageref{S-402-sq} sq., 
we want again to start with a determined form:
\[
x_i'
=
f_i(x_1,\dots,x_n,\,a_1,\dots,a_r)
\ \ \ \ \ \ \ \ \ \ \ \ \ {\scriptstyle{(i\,=\,1\,\cdots\,n)}}
\]
of the group: $X_1f, \dots, X_rf$, and as before, we let $a_1^0,
\dots, a_r^0$ denote the parameters attached to the identity
transformation: $x_i' = x_i$.

Then it is easy to identify the complete system, the integration of
which conducts precisely to the equations: $x_i' = f_i ( x, a)$.

The complete system in question has simply the form:
\[
X_k'
+
A_kf
=
0
\ \ \ \ \ \ \ \ \ \ \ \ \ {\scriptstyle{(k\,=\,1\,\cdots\,r)}},
\]
where, in the infinitesimal transformations:
\[
A_kf
=
\sum_{j=1}^r\,\alpha_{kj}
(a_1,\dots,a_r)\,
\frac{\partial f}{\partial a_j}
\ \ \ \ \ \ \ \ \ \ \ \ \ {\scriptstyle{(k\,=\,1\,\cdots\,r)}},
\]
the functions $\alpha_{ kj} ( a)$ are the
same as in the differential equations~\thetag{ 7'}. 
If: 
\[
F_i(x_1',\dots,x_n',\,a_1,\dots,a_r)
\ \ \ \ \ \ \ \ \ \ \ \ \ {\scriptstyle{(i\,=\,1\,\cdots\,n)}}
\]
denote the general solutions to the complete system: $X_k' f + A_kf =
0$ relatively to: $a_1 = a_1^0$, \dots, $a_r = a_r^0$, then by
resolution with respect to $x_1', \dots, x_n'$, the equations:
\[
x_i
=
F_i
(x_1',\dots,x_n',\,a_1,\dots,a_r)
\ \ \ \ \ \ \ \ \ \ \ \ \ {\scriptstyle{(i\,=\,1\,\cdots\,n)}}
\]
give the equations: $x_i' = f_i ( x, a)$ exactly. All of this follows
from the developments of the Chaps.~\ref{fundamental-differential}
and~\ref{kapitel-9}.

We now apply these considerations to the parameter group associated to
the group: $x_i' = f_i ( x, a)$, whose equations, according to
p.~\pageref{S-404}, have the form:
\[
a_k'
=
\varphi_k(a_1,\dots,a_r,\,b_1,\dots,b_r)
\ \ \ \ \ \ \ \ \ \ \ \ \ {\scriptstyle{(k\,=\,1\,\cdots\,r)}}
\]
and whose identity transformation possesses the parameters: 
$b_1 = a_1^0$, \dots, $b_r = a_r^0$. 

The complete system, through the introduction of which
the equations: $a_k' = \varphi_k ( a, b)$ of the
parameter group can be found, visibly reads:
\[
\sum_{j=1}^r\,\alpha_{kj}(a_1',\dots,a_r')\,
\frac{\partial f}{\partial a_j'}
+
\sum_{j=1}^r\,\alpha_{kj}(b_1,\dots,b_r)\,
\frac{\partial f}{\partial b_j}
=
0
\ \ \ \ \ \ \ \ \ \ \ \ \ {\scriptstyle{(k\,=\,1\,\cdots\,r)}}.
\]
If we determine the general solutions:
\[
H_j(a_1',\dots,a_r',\,b_1,\dots,b_r)
\ \ \ \ \ \ \ \ \ \ \ \ \ {\scriptstyle{(j\,=\,1\,\cdots\,r)}}
\]
of this complete system relatively to: 
$b_k = a_k^0$ and if we solve afterwards the equations:
\[
a_j
=
H_j
(a_1',\dots,a_r',\,b_1,\dots,b_r)
\ \ \ \ \ \ \ \ \ \ \ \ \ {\scriptstyle{(j\,=\,1\,\cdots\,r)}}
\]
with respect to $a_1', \dots, a_r'$, we obtain: 
$a_k' = \varphi_k ( a, b)$. 

Thus, the following holds: 

\def\thetheorem{73}\begin{theorem}
\label{Theorem-73-S-411}
If one knows the infinitesimal transformations:
\[
A_kf
=
\sum_{j=1}^r\,\alpha_{kj}(a_1,\dots,a_r)\,
\frac{\partial f}{\partial a_j}
\ \ \ \ \ \ \ \ \ \ \ \ \ {\scriptstyle{(k\,=\,1\,\cdots\,r)}}
\]
of the parameter group of an $r$-term group and if one knows that the
identity transformation of this group possesses the parameters:
$a_1^0, \dots, a_r^0$, whence the identity transformation of the
parameter group also possesses the parameters: $a_1^0, \dots, a_r^0$,
then one finds the finite equations of the parameter group in the
following way: One determines the general solutions of the complete
system:
\[
\sum_{j=1}^r\,\alpha_{kj}(a')\,
\frac{\partial f}{\partial a_j'}
+
\sum_{j=1}^r\,\alpha_{kj}(b)\,
\frac{\partial f}{\partial b_j}
=
0
\ \ \ \ \ \ \ \ \ \ \ \ \ {\scriptstyle{(k\,=\,1\,\cdots\,r)}}
\]
relatively to: $b_1 = a_1^0$, \dots, $b_r = a_r^0$; 
if: 
\[
H_j
(a_1',\dots,a_r',\,b_1,\dots,b_r)
\ \ \ \ \ \ \ \ \ \ \ \ \ {\scriptstyle{(j\,=\,1\,\cdots\,r)}}
\]
are these general solutions, then one obtains the sought equations of
the parameter group by solving the $r$ equations:
\[
a_j
=
H_j
(a_1',\dots,a_r',\,b_1,\dots,b_r)
\ \ \ \ \ \ \ \ \ \ \ \ \ {\scriptstyle{(j\,=\,1\,\cdots\,r)}}
\]
with respect to $a_1', \dots, a_r'$. 
\end{theorem}

\medskip

If, in an $r$-term group $x_i' = f_i ( x, a)$, one introduces new
parameters $\mathfrak{ a}_k$ in place of the $a$, then group receives
a new form to which a new parameter group is also naturally
associated.

The connection between the new and the old parameter groups is very
simple. Indeed, if the new parameters are determined by the equations:
\[
a_j
=
\lambda_j
(\mathfrak{a}_1,\dots,\mathfrak{a}_r)
\ \ \ \ \ \ \ \ \ \ \ \ \ {\scriptstyle{(j\,=\,1\,\cdots\,r)}},
\]
then the new form of the group: $x_i' = 
f_i ( x,a)$ is the following one:
\[
x_i'
=
f_i
\big(
x_1,\dots,x_n,\,
\lambda_1(\mathfrak{a}),\,\dots,\,\lambda_r(\mathfrak{a})
\big)
\ \ \ \ \ \ \ \ \ \ \ \ \ {\scriptstyle{(i\,=\,1\,\cdots\,n)}},
\]
and the new parameter group reads:
\[
\aligned
\lambda_k(\mathfrak{a}_1',\dots,\mathfrak{a}_r')
=
\varphi_k
&
\big(
\lambda_1(\mathfrak{a}),\,\dots,\,\lambda_r(\mathfrak{a}),\,\,
\lambda_1(\mathfrak{b}),\,\dots,\,\lambda_r(\mathfrak{b})
\big)
\\
&
\ \ \ \ \ \
{\scriptstyle{(k\,=\,1\,\cdots\,r)}}
\endaligned
\]
while the old one was $a_k' = \varphi_k ( a, b)$. 
Consequently, the new one comes from the old one
when by executing the substitution $a_j = \lambda_j ( 
\mathfrak{ a})$ both on the $a$ and on the $b$, 
that is to say, by inserting, for the $a'$, $a$, $b$, 
the following values:
\[
a_j'
=
\lambda_j(\mathfrak{a}'),
\ \ \ \ \
a_j
=
\lambda_j(\mathfrak{a}),
\ \ \ \ \
b_j
=
\lambda_j(\mathfrak{b})
\ \ \ \ \ \ \ \ \ \ \ \ \ {\scriptstyle{(j\,=\,1\,\cdots\,r)}}
\]
in the equations: $a_k' = \varphi_k ( a, b)$. 

\sectionengellie{\S\,\,\,103.}

According to Theorem~72, p.~\pageref{Theorem-72-S-407}, every $r$-term
group has the same composition as its parameter group, and
consequently, the $r$-term group which have the same parameter group
are equally composed with each other.

Now, we claim that conversely, two equally composed $r$-term groups
can always be brought, after the introduction of new parameters, to a
form such that they two have the same parameter group.

In order to prove this claim, we imagine that the
infinitesimal transformations of the two groups are given; 
let the ones of the first group be:
\[
X_kf
=
\sum_{i=1}^n\,\xi_{ki}(x_1,\dots,x_n)\,
\frac{\partial f}{\partial x_i}
\ \ \ \ \ \ \ \ \ \ \ \ \ {\scriptstyle{(k\,=\,1\,\cdots\,r)}},
\]
and the ones of the other group be:
\[
Y_kf
=
\sum_{\mu=1}^m\,
\eta_{k\mu}(y_1,\dots,y_m)\,
\frac{\partial f}{\partial y_\mu}
\ \ \ \ \ \ \ \ \ \ \ \ \ {\scriptstyle{(k\,=\,1\,\cdots\,r)}}.
\]
Since the two groups are equally composed, 
we can assume that the $X_kf$ and the $Y_kf$
are already chosen in such a way that, 
simultaneously with the relations:
\[
\leftbracket
X_i,\,X_k
\rightbracket
=
\sum_{s=1}^r\,c_{iks}\,X_sf,
\]
the relations: 
\[
\leftbracket
Y_i,\,Y_k
\rightbracket
=
\sum_{s=1}^r\,c_{iks}\,Y_sf
\]
hold.

Lastly, we yet imagine that $r$ independent
infinitesimal transformations:
\[
A_kf
=
\sum_{j=1}^r\,
\alpha_{kj}(a_1,\dots,a_r)\,
\frac{\partial f}{\partial a_j}
\ \ \ \ \ \ \ \ \ \ \ \ \ {\scriptstyle{(k\,=\,1\,\cdots\,r)}}
\]
are given which generate a simply transitive $r$-term group
equally composed with these groups and which
stand pairwise in the relationships:
\[
\leftbracket
A_i,\,A_k
\rightbracket
=
\sum_{s=1}^r\,c_{iks}\,A_sf.
\]

At present, we form the complete system:
\[
X_k'f
+
A_kf
=
0
\ \ \ \ \ \ \ \ \ \ \ \ \ {\scriptstyle{(k\,=\,1\,\cdots\,r)}}
\]
and we determine its general solutions:
\[
F_i(x_1',\dots,x_n',\,a_1,\dots,a_r)
\ \ \ \ \ \ \ \ \ \ \ \ \ {\scriptstyle{(i\,=\,1\,\cdots\,n)}}
\]
relatively to an arbitrary system of values: $a_1 = a_1^0$, \dots,
$a_r = a_r^0$; furthermore, we form the complete system:
\[
Y_k'f
+
A_kf
=
0
\ \ \ \ \ \ \ \ \ \ \ \ \ {\scriptstyle{(k\,=\,1\,\cdots\,r)}}
\]
and we determine its general solutions:
\[
\mathfrak{F}_\mu
(y_1',\dots,y_m',\,a_1,\dots,a_r)
\ \ \ \ \ \ \ \ \ \ \ \ \ {\scriptstyle{(i\,=\,1\,\cdots\,n)}}
\]
relatively to the same system of values: $a_k = a_k^0$. 

In addition, we yet solve the $n$ equations: $x_i = 
F_i ( x', a)$ with respect to $x_1', \dots, x_n'$: 
\def\theequation{11}\begin{equation}
x_i'
=
f_i(x_1,\dots,x_n,\,a_1,\dots,a_r)
\ \ \ \ \ \ \ \ \ \ \ \ \ {\scriptstyle{(i\,=\,1\,\cdots\,n)}},
\end{equation}
and likewise, the $m$ equations: $y_\mu = \mathfrak{ F}_\mu ( y', a)$
with respect to $y_1', \dots, y_m'$:
\def\theequation{11'}\begin{equation}
y_\mu'
=
\mathfrak{f}_\mu
(y_1,\dots,y_m,\,a_1,\dots,a_r)
\ \ \ \ \ \ \ \ \ \ \ \ \ {\scriptstyle{(\mu\,=\,1\,\cdots\,n)}}.
\end{equation}

In whichever form the finite equations of the group: $X_1f, \dots,
X_rf$ are present, then obviously, they can be brought to the
form~\thetag{ 11} by introducing new parameters, and likewise, in
whichever form the finite equations of the group: $Y_1f, \dots, Y_rf$
are present, they can always be given the form~\thetag{ 11'} by
introducing new parameters.

But the Theorem~73 p.~\pageref{Theorem-73-S-411} can be applied to the
two groups~\thetag{ 11} and~\thetag{ 11'}. Indeed, in the two groups,
the identity transformation has the parameters: $a_1^0, \dots, a_r^0$
and for the two groups, the associated parameter groups contain the
$r$ independent infinitesimal transformations: $A_1f, \dots, A_rf$,
and consequently, according to the mentioned theorem, one can indicate
for each one of the two groups the associated parameter group, and
these parameter groups are the same for both of them.

As a result, the claim stated above is proved.\,---

Thus, it is at present established that two groups which have the same
parameter group are equally composed, and on the other hand, that two
groups that which are equally composed, can be brought, by introducing
new parameters, to a form in which they two have the same parameter
group. Thus, we can say:

\def\thetheorem{74}\begin{theorem}
\label{Theorem-74-S-414}
Two $r$-term groups:
\[
x_i'
=
f_i(x_1,\dots,x_n,\,a_1,\dots,a_r)
\ \ \ \ \ \ \ \ \ \ \ \ \ {\scriptstyle{(i\,=\,1\,\cdots\,n)}}
\]
and:
\[
y_\mu'
=
g_\mu(y_1,\dots,y_m,\,
\mathfrak{a}_1,\dots,\mathfrak{a}_r)
\ \ \ \ \ \ \ \ \ \ \ \ \ {\scriptstyle{(\mu\,=\,1\,\cdots\,m)}}
\]
are equally composed if and only if it is possible
to represent the parameters: $\mathfrak{ a}_1, 
\dots, \mathfrak{ a}_r$ as independent functions
of the $a$: 
\[
\mathfrak{a}_k
=
\chi_k(a_1,\dots,a_r)
\ \ \ \ \ \ \ \ \ \ \ \ \ {\scriptstyle{(k\,=\,1\,\cdots\,r)}}
\]
in such a way that the parameter group of the group:
\[
y_\mu'
=
g_\mu
\big(
y_1,\dots,y_m,\,
\chi_1(a),\dots,\chi_r(a)
\big)
\ \ \ \ \ \ \ \ \ \ \ \ \ {\scriptstyle{(\mu\,=\,1\,\cdots\,m)}}
\]
coincides with the parameter group of the group: 
$x_i' = f_i ( x, a)$. 
\end{theorem}

It is particularly noticeable that our two equally composed groups:
$X_1f, \dots, X_rf$ and: $Y_1f, \dots, Y_rf$ then have the same
parameter group when one writes their finite equations under the
canonical form:
\def\theequation{12}\begin{equation}
x_i'
=
x_i
+
\sum_{k=1}^r\,e_k\,X_k\,x_i
+
\sum_{k,\,\,j}^{1\cdots\,r}\,
\frac{e_k\,e_j}{1\cdot 2}\,
X_k\,X_j\,x_i
+\cdots
\ \ \ \ \ \ \ \ \ \ \ \ \ {\scriptstyle{(i\,=\,1\,\cdots\,n)}}
\end{equation}
and, respectively:
\def\theequation{12'}\begin{equation}
y_\mu'
=
y_\mu
+
\sum_{k=1}^r\,e_k\,Y_k\,y_\mu
+
\sum_{k,\,\,j}^{1\cdots\,r}\,
\frac{e_k\,e_j}{1\cdot 2}\,
Y_k\,Y_j\,y_\mu
+\cdots
\ \ \ \ \ \ \ \ \ \ \ \ \ {\scriptstyle{(\mu\,=\,1\,\cdots\,m)}}.
\end{equation}

In fact, from the developments in Chap.~\ref{one-term-groups},
Sect.~\ref{generation-by-one-term}, it follows that the
equations~\thetag{ 11} receive the form~\thetag{ 12} after the
substitution:
\def\theequation{13}\begin{equation}\label{S-414}
\aligned
a_\nu
=
a_\nu^0
&
+
\sum_{k=1}^r\,e_k\,
\big[A_k\,a_\nu\big]_{a=a^0}
+
\sum_{k,\,\,j}^{1\cdots\,r}\,
\frac{e_k\,e_j}{1\cdot 2}\,
\big[A_k\,A_j\,a_\nu\big]_{a=a^0}
+\cdots
\\
&
\ \ \ \ \ \ \ \ \ \ \ \ \ 
\ \ \ \ \ \ \ \ \ \ \ \ \ 
{\scriptstyle{(\nu\,=\,1\,\cdots\,r)}}
\endaligned
\end{equation}
and that the equations~\thetag{ 11'} are transferred
to~\thetag{ 12'} by the same substitution. 
Now, since the equations~\thetag{ 13} represent
a transformation between the parameters: 
$a_1, \dots, a_r$ and $e_1, \dots, e_r$ and
since the two groups~\thetag{ 11} and~\thetag{ 11'}
have their parameter group in common, it follows
that one and the same parameter group 
is associated to the two groups~\thetag{ 12} and~\thetag{ 12'}. 

Thus, the following holds:

\def\theproposition{1}\begin{proposition}
\label{Satz-1-S-414}
If the $r$ independent infinitesimal transformations:
\[
X_kf
=
\sum_{i=1}^n\,\xi_{ki}(x_1,\dots,x_n)\,
\frac{\partial f}{\partial x_i}
\ \ \ \ \ \ \ \ \ \ \ \ \ {\scriptstyle{(k\,=\,1\,\cdots\,r)}}
\]
stand pairwise in the relationships:
\[
\leftbracket
X_i,\,X_k
\rightbracket
=
\sum_{s=1}^r\,c_{iks}\,X_sf
\]
and if, on the other hand, the $r$ independent
infinitesimal transformations:
\[
Y_kf
=
\sum_{\mu=1}^m\,\eta_{k\mu}(y_1,\dots,y_m)\,
\frac{\partial f}{\partial y_\mu}
\ \ \ \ \ \ \ \ \ \ \ \ \ {\scriptstyle{(k\,=\,1\,\cdots\,r)}}
\]
stand pairwise in the same relationships:
\[
\leftbracket
Y_i,\,Y_k
\rightbracket
=
\sum_{s=1}^r\,c_{iks}\,Y_sf,
\]
then the two equally composed groups:
\[
x_i'
=
x_i
+
\sum_{k=1}^r\,e_k\,X_k\,x_i
+
\sum_{k,\,\,j}^{1\cdots\,r}\,
\frac{e_k\,e_j}{1\cdot 2}\,
X_k\,X_j\,x_i
+\cdots
\ \ \ \ \ \ \ \ \ \ \ \ \ {\scriptstyle{(i\,=\,1\,\cdots\,n)}}
\]
and:
\[
y_\mu'
=
y_\mu
+
\sum_{k=1}^r\,e_k\,Y_k\,y_\mu
+
\sum_{k,\,\,j}^{1\cdots\,r}\,
\frac{e_k\,e_j}{1\cdot 2}\,
Y_k\,Y_j\,y_\mu
+\cdots
\ \ \ \ \ \ \ \ \ \ \ \ \ {\scriptstyle{(\mu\,=\,1\,\cdots\,m)}}
\]
possess one and the same parameter group.
\end{proposition}

At present, we imagine that two arbitrary equally composed
$r$-term groups: 
\[
x_i'
=
f_i(x_1,\dots,x_n,\,a_1,\dots,a_r)
\ \ \ \ \ \ \ \ \ \ \ \ \ {\scriptstyle{(i\,=\,1\,\cdots\,n)}}
\]
and:
\[
y_\mu'
=
\mathfrak{f}_\mu(y_1,\dots,y_m,\,a_1,\dots,a_r)
\ \ \ \ \ \ \ \ \ \ \ \ \ {\scriptstyle{(\mu\,=\,1\,\cdots\,m)}}
\]
are presented which already have a form such that the parameter group
for both is the same.

Let the infinitesimal transformations of the parameter group in
question be:
\[
A_kf
=
\sum_{j=1}^r\,\alpha_{kj}
(a_1,\dots,a_r)\,
\frac{\partial f}{\partial a_j}
\ \ \ \ \ \ \ \ \ \ \ \ \ {\scriptstyle{(k\,=\,1\,\cdots\,r)}}
\]
and let them be linked together by the relations:
\[
\leftbracket
A_i,\,A_k
\rightbracket
=
\sum_{s=1}^r\,c_{iks}\,A_sf.
\]

According to p.~\pageref{S-405}, $x_1', \dots, x_n'$, when
considered as functions of $x_1, \dots, x_n$, $a_1, \dots, a_r$,
satisfy certain differential equations of the form:
\[
\xi_{ji}(x_1',\dots,x_n')
=
\sum_{k=1}^r\,
\alpha_{jk}(a_1,\dots,a_r)\,
\frac{\partial x_i'}{\partial a_k}
\ \ \ \ \ \ \ \ \ \ \ \ \ 
{\scriptstyle{(j\,=\,1\,\cdots\,r\,;\,\,\,
i\,=\,1\,\cdots\,n)}}.
\]
Here, the $\xi_{ ji} (x')$ are completely determined functions,
because by resolution of the equations: 
$x_i' = f_i ( x,a)$, we get, say: $x_i = F_i ( x', a)$, 
whence it holds identically:
\[
\xi_{ji}(x_1',\dots,x_n')
=
\sum_{k=1}^r\,
\alpha_{jk}(a_1,\dots,a_r)\,
\bigg[
\frac{\partial f_i(x,a)}{\partial a_k}
\bigg]_{x=F(x',a)},
\]
so that we can therefore find the $\xi_{ ji} (x')$ without
difficulty, when we want. 

According to Theorem~21, p.~\pageref{Theorem-21-S-149}, 
the $r$ expressions:
\[
X_kf
=
\sum_{i=1}^n\,
\xi_{ki}(x_1,\dots,x_n)\,
\frac{\partial f}{\partial x_i}
\ \ \ \ \ \ \ \ \ \ \ \ \ {\scriptstyle{(k\,=\,1\,\cdots\,r)}}
\]
are independent infinitesimal transformations
and they are linked together by the relations:
\[
\leftbracket
X_i,\,X_k
\rightbracket
=
\sum_{s=1}^r\,c_{iks}\,X_sf\,;
\]
naturally, they generate the group: $x_i' = 
f_i ( x, a)$. 

For the group: $y_\mu ' = \mathfrak{ f}_\mu
(y, a)$, there are
in a corresponding way differential equations of the form:
\[
\eta_{j\mu}(y_1',\dots,y_m')
=
\sum_{k=1}^r\,
\alpha_{jk}(a_1,\dots,a_r)\,
\frac{\partial y_\mu'}{\partial a_k}
\ \ \ \ \ \ \ \ \ \ \ \ \ 
{\scriptstyle{(j\,=\,1\,\cdots\,r\,;\,\,\,
\mu\,=\,1\,\cdots\,m)}},
\]
where the $\eta_{ j\mu} (y')$ are completely determined
functions. The $r$ expressions:
\[
Y_kf
=
\sum_{\mu=1}^m\,
\eta_{k\mu}(y_1,\dots,y_m)\,
\frac{\partial f}{\partial y_\mu}
\ \ \ \ \ \ \ \ \ \ \ \ \ {\scriptstyle{(k\,=\,1\,\cdots\,r)}}
\]
are independent infinitesimal transformations and they are
linked together by the relations:
\[
\leftbracket
Y_i,\,Y_k
\rightbracket
=
\sum_{s=1}^r\,
c_{iks}\,Y_sf\,;
\]
naturally, they generate the group: $y_\mu ' = \mathfrak{ f}_\mu
(y, a)$. 

At present, exactly as on p.~\pageref{S-414}, we recognize that, 
after the substitution: 
\[
a_\nu
=
a_\nu^0
+
\sum_{k=1}^r\,
e_k\,\big[A_k\,a_\nu\big]_{a=a^0}
+\cdots
\ \ \ \ \ \ \ \ \ \ \ \ \ {\scriptstyle{(\nu\,=\,1\,\cdots\,r)}},
\]
the
groups: $x_i' = f_i ( x, a)$ and: $y_\mu' = \mathfrak{ f}_\mu 
(y, a)$ receive the forms: 
\[
x_i'
=
x_i
+
\sum_{k=1}^r\,e_k\,X_k\,x_i
+\cdots
\ \ \ \ \ \ \ \ \ \ \ \ \ {\scriptstyle{(i\,=\,1\,\cdots\,n)}}
\]
and, respectively:
\[
y_\mu'
=
y_\mu
+
\sum_{k=1}^r\,
e_k\,Y_k\,y_\mu
+\cdots
\ \ \ \ \ \ \ \ \ \ \ \ \ {\scriptstyle{(\mu\,=\,1\,\cdots\,m)}}.
\]
Therefore, we get the

\def\theproposition{2}\begin{proposition}
\label{Satz-2-S-417}
If the two $r$-term groups:
\[
x_i'
=
f_i(x_1,\dots,x_n,\,a_1,\dots,a_r)
\ \ \ \ \ \ \ \ \ \ \ \ \ {\scriptstyle{(i\,=\,1\,\cdots\,n)}}
\]
and:
\[
y_\mu'
=
\mathfrak{f}_\mu
(y_1,\dots,y_m,\,a_1,\dots,a_r)
\ \ \ \ \ \ \ \ \ \ \ \ \ {\scriptstyle{(\mu\,=\,1\,\cdots\,m)}}
\]
have the same parameter group, then it is possible
to introduce, in place of the $a$, new parameters: 
$e_1, \dots, e_r$ so that the two groups
receive the forms: 
\[
x_i'
=
x_i
+
\sum_{k=1}^r\,e_k\,X_k\,x_i
+\cdots
\ \ \ \ \ \ \ \ \ \ \ \ \ {\scriptstyle{(i\,=\,1\,\cdots\,n)}}
\]
and, respectively:
\[
y_\mu'
=
y_\mu
+
\sum_{k=1}^r\,
e_k\,Y_k\,y_\mu
+\cdots
\ \ \ \ \ \ \ \ \ \ \ \ \ {\scriptstyle{(\mu\,=\,1\,\cdots\,m)}}\,;
\]
here, $X_1f, \dots, X_rf$ and $Y_1f, \dots, Y_rf$ are independent
infinitesimal transformations of the two groups such that,
simultaneously with the relations:
\[
\leftbracket
X_i,\,X_k
\rightbracket
=
\sum_{s=1}^r\,c_{iks}\,X_sf
\]
there hold the relations:
\[
\leftbracket
Y_i,\,Y_k
\rightbracket
=
\sum_{s=1}^r\,c_{iks}\,Y_sf,
\]
so that the two groups are related to each other in a holoedrically
isomorphic way when the infinitesimal transformation $e_1\, Y_1 f +
\cdots + e_r\, Y_rf$ is associated to every infinitesimal
transformation: $e_1\, X_1f + \cdots + e_r\, X_rf$.
\end{proposition}

Now, about the two groups: $x_i' = f_i ( x, a)$ and: $y_\mu' =
\mathfrak{ f}_\mu ( y, a)$, we make the same assumptions as in
Proposition~2; in addition, we want to assume yet that: $a_k' =
\varphi_k ( a, b)$ are the finite equations of their common parameter
group.

Under this assumption, the following obviously holds true of each one
of the two groups: If two transformations of the group which have the
parameters: $a_1, \dots, a_r$ and: $b_1, \dots, b_r$, respectively,
are executed one after the other, then the resulting transformation
belongs to the group and it possesses the parameters: $\varphi_1 (
a,b)$, \dots, $\varphi_r ( a, b)$.

We can express this fact somehow differently if we mutually associate
the transformations of the two groups in a way so that every
transformation of the one group corresponds to the transformation of
the other group which has the same parameters. Indeed, we can then
say: if $S$ is a transformation of the one group and if $\mathfrak{
S}$ is the corresponding transformation of the other group, and
moreover, if $T$ is a second transformation of the one group and if
$\mathfrak{ T}$ the corresponding transformation of the other group,
then the transformation $S\, T$ of the one group corresponds to the
transformation $\mathfrak{ S}\, \mathfrak{ T}$ in the other group.

Such a mutual association of the transformations of both groups is
then, according to Theorem~74, p.~\pageref{Theorem-74-S-414} always
possible when and only when the two groups are equally composed. So we
have the:

\renewcommand{\thefootnote}{\fnsymbol{footnote}}
\def\thetheorem{75}\begin{theorem}
\label{Theorem-75-S-418}
Two $r$-term groups are equally composed if and only if it is possible
to relate the transformations of the one group to the transformations
of the other group in a univalent invertible way so that the following
holds true: If, in the one group, one executes two transformations one
after the other and if, in the other group, one executes one after the
other and in the same order the corresponding transformations, then
the transformation that one obtains in the one group corresponds to
the transformation that one obtains in the other group.\footnote[1]{\,
Cf. \name{Lie}, Archiv for Mathematik og Naturvidenskab 1876; 
Math. Ann. Vol. XXV, p.~77; G. d. W. zu Christiania, 1884, no.~15.
} 
\end{theorem}
\renewcommand{\thefootnote}{\arabic{footnote}}

The above considerations provide a new important result when they are
applied to the Proposition~1, p.~\pageref{Satz-1-S-414}. In order to
be able to state this result under the most simple form, we remember
two different things: firstly, that the two groups: $X_1f, \dots,
X_rf$ and: $Y_1f, \dots, Y_rf$ are related to each other in a
holoedrically isomorphic way when to every infinitesimal
transformation: $e_1\, X_1f + \cdots + e_r\, X_rf$ is associated the
infinitesimal transformation: $e_1\, Y_1f + \cdots + e_r\, Y_rf$, and
secondly, that the expression: $e_1\, X_1f + \cdots + e_r\, X_rf$ can
also be interpreted as the symbol of a finite transformation of the
group: $X_1f, \dots, X_rf$ (cf. Chap.~\ref{kapitel-17},
p.~\pageref{S-293} and Chap.~\ref{kapitel-15}, p.~\pageref{S-255}).
By taking this into account, we can say:

\def\theproposition{3}\begin{proposition}
\label{Satz-3-S-418}
Let the two equally composed $r$-term groups:
\[
X_kf
=
\sum_{i=1}^n\,\xi_{ki}(x_1,\dots,x_n)\,
\frac{\partial f}{\partial x_i}
\ \ \ \ \ \ \ \ \ \ \ \ \ {\scriptstyle{(k\,=\,1\,\cdots\,r)}}
\]
and:
\[
Y_kf
=
\sum_{\mu=1}^m\,\eta_{k\mu}(y_1,\dots,y_m)\,
\frac{\partial f}{\partial y_\mu}
\ \ \ \ \ \ \ \ \ \ \ \ \ {\scriptstyle{(k\,=\,1\,\cdots\,r)}}
\]
be related to each other in a holoedrically isomorphic way
when one associates to every infinitesimal transformation: 
$e_1\, X_1f + \cdots + e_r\, X_rf$ the infinitesimal
transformation: $e_1\, Y_1f + \cdots + e_r\, Y_rf$ so that, 
simultaneously with the relations:
\[
\leftbracket
X_i,\,X_k
\rightbracket
=
\sum_{s=1}^r\,c_{iks}\,X_sf
\]
there hold the relations:
\[
\leftbracket
Y_i,\,Y_k
\rightbracket
=
\sum_{s=1}^r\,c_{iks}\,Y_sf,
\]
Then, if one interprets the expressions: $\sum\, e_k\, X_kf$ and
$\sum\, e_k\, Y_kf$ as the general symbols of the finite
transformations of the two groups: $X_1f, \dots, X_rf$ and $Y_1f,
\dots, Y_rf$, the following holds true: If the two transformations:
$\sum\, e_k\, X_kf$ and $\sum\, e_k' \, X_kf$ of the group: $X_1f,
\dots, X_rf$ give, when executed one after the other, the
transformation: $\sum\, e_k''\, X_kf$, then the two transformations:
$\sum\, e_k\, Y_kf$ and $\sum\, e_k'\, Y_kf$ give the transformation:
$\sum\, e_k''\, Y_kf$ when executed one after the other.
\end{proposition}

Thus, if one has related two equally composed $r$-term groups in a
holoedrically isomorphic way in the sense of Chap.~\ref{kapitel-17},
p.~\pageref{S-293}, then at the same time thanks to this, one has
established a univalent invertible relationship between the finite
transformations of the two groups, as is written down in Theorem~75.

But the converse also holds true: If one has produced a univalent
invertible relationship between the transformations of two $r$-term
equally composed groups which has property written down in Theorem~75,
then at the same time thanks to this, the two groups are related to
each other in a holoedrically isomorphic way.

In fact, let: $x_i' = f_i ( x_1, \dots, x_n,\, a_1, \dots, a_r)$ be
the one group and let: $y_\mu ' = \mathfrak{ f}_\mu ( y_1, \dots, y_m,
\, a_1, \dots, a_r)$ be the transformation of the other group that
corresponds to the transformation: $x_i' = f_i ( x, a)$. Then the two
groups: $x_i' = f_i ( x,a)$ and $y_\mu ' = \mathfrak{ f}_\mu ( y, a)$
obviously have one and the same parameter group and hence, according
to Proposition~2, p.~\pageref{Satz-2-S-417}, by introducing
appropriate new parameters: $e_1, \dots, e_r$, they
can be given the following forms:
\[
x_i'
=
x_i
+
\sum_{k=1}^r\,e_k\,X_k\,x_i
+\cdots
\ \ \ \ \ \ \ \ \ \ \ \ \ {\scriptstyle{(i\,=\,1\,\cdots\,n)}}
\]
and:
\[
y_\mu'
=
y_\mu
+
\sum_{k=1}^r\,e_k\,Y_k\,y_\mu
+\cdots
\ \ \ \ \ \ \ \ \ \ \ \ \ {\scriptstyle{(\mu\,=\,1\,\cdots\,m)}}.
\]
From this, it results that the said univalent invertible relationship
between the transformations of the two groups amounts to the fact
that, to every finite transformation: $e_1\, X_1f + \cdots + e_r\,
X_rf$ of the one group is associated the finite transformation: $e_1\,
Y_1f + \cdots + e_r\, Y_rf$ of the other group. Therefore at the same
time, to every infinitesimal transformation: $\sum\, e_k\, X_kf$ is
associated the infinitesimal transformation: $\sum\, e_k\, Y_kf$, and
as a result, according to Proposition~2, a holoedrically isomorphic
relationship between the two groups is effectively established.

\medskip

\label{S-420-sq}
In the theory of substitutions, one defines the holoedric isomorphism
of two groups and the holoedrically isomorphic relationship between
two groups differently than what we have done in
Chap.~\ref{kapitel-17}. There, one says that two groups with the same
number of substitutions are ``equally composed'' or ``holoedrically
isomorphic'' when one can produce, between the transformations of the
two groups, a univalent invertible relationship which has the property
written in the Theorem~75, p.~\pageref{Theorem-75-S-418}; if such a
relationship between two holoedrically isomorphic groups is really
established, then one says that the two groups are ``related to each
other in a holoedrically isomorphic way''.

But already on p.~\pageref{S-293-bis}, we observed that materially
\deutsch{materiell}, our definition of the concept in question
corresponds precisely to the one which is usual in the theory of
substitutions, as far as such a correspondence may in any case be
possible between domains so different as the theory of substitutions
and the theory of transformation groups.

Our last developments show that the claim made on
p.~\pageref{S-293-bis} is correct. From Theorem~75,
p.~\pageref{Theorem-75-S-418}, it results that two $r$-term groups
which are holoedrically isomorphic in the sense of
Chap.~\ref{kapitel-17}, p.~\pageref{S-293} must also be called
holoedrically isomorphic in the sense of the theory of substitutions,
and conversely. From Proposition~3, p.~\pageref{Satz-3-S-418} and
from the remarks following, it is evident that two $r$-term groups
which are related to each other in a holoedrically isomorphic way in
the sense of Chap.~\ref{kapitel-17}, p.~\pageref{S-293}, are also
holoedrically isomorphic in the sense of substitution theory.

\smallercharacters{

It still remains to prove that our definition of the meroedric
isomorphism (cf. Chap.~\ref{kapitel-17}, p.~\pageref{S-293}) also
corresponds to the definition of meroedric isomorphism that is given
by the theory of substitutions.

Let $X_1f, \dots, X_{ r-q}f, \dots, X_rf$ and $Y_1f, \dots,
Y_{ r-q}f$ be two meroedrically isomorphic groups and let
$X_{ r- q + 1}f, \dots, X_rf$ be precisely the invariant
subgroup of the $r$-term group which corresponds to the
identity transformation of the $(r-q)$-term group. 
Lastly, let $A_1f, \dots, A_rf$ be a simply transitive
group in the variables $a_1, \dots, a_r$ which has the
same composition as $X_1f, \dots, X_rf$. Here, we want to 
assume that $a_1, \dots, a_{ r-q}$ are solutions of the
complete system:
\[
A_{r-q+1}f=0,
\,\,\,\dots,\,\,\,
A_rf=0.
\]

After this choice of the variables, $A_1f, \dots, A_{r-q}f$
possess the form:
\[
\aligned
A_kf
=
\sum_{i=1}^{r-q}\,
&
\alpha_{ki}(a_1,\dots,a_{r-q})\,
\frac{\partial f}{\partial a_i}
+
\sum_{j=r-q+1}^r\,
\beta_{kj}(a_1,\dots,a_r)\,
\frac{\partial f}{\partial a_j}
\\
&
\ \ \ \ \ \ \ \ \ \ \ \ \ \ \ \ \ \ \ \ \ \ \ 
{\scriptstyle{(k\,=\,1\,\cdots\,r\,-\,q)}},
\endaligned
\]
while $A_{ r-q+1}f, \dots, A_rf$ have the form:
\[
A_kf
=
\sum_{j=r-q+1}^r\,\beta_{kj}(a_1,\dots,a_r)\,
\frac{\partial f}{\partial a_j}
\ \ \ \ \ \ \ \ \ \ \ \ \ 
{\scriptstyle{(k\,=\,1\,\cdots\,r\,-\,q\,+\,1)}}.
\]
At the same time, it is clear that the reduced infinitesimal
transformations:
\[
\overline{A}_kf
=
\sum_{i=1}^{r-q}\,
\alpha_{ki}(a_1,\dots,a_{r-q})\,
\frac{\partial f}{\partial a_i}
\ \ \ \ \ \ \ \ \ \ \ \ \  
{\scriptstyle{(k\,=\,1\,\cdots\,r\,-\,q)}}
\]
generate a simply transitive group equally composed with the
$(r-q)$-term group $Y_1f, \dots, Y_{ r-q}f$.

If one denotes by $B_kf$, $\overline{ B}_kf$ the infinitesimal
transformations $A_kf$, $\overline{ A}_kf$ written down in the
variables $b$ instead of the variables $a$, then according to earlier
discussions, one finds the parameter group of the $(r-q)$-term group
$Y_1f, \dots, Y_{ r-q}f$ by integrating the complete system:
\[
\overline{A}_kf
+
\overline{B}_kf
=
0
\ \ \ \ \ \ \ \ \ \ \ \ \  
{\scriptstyle{(k\,=\,1\,\cdots\,r\,-\,q)}},
\]
and also the parameter group of the $r$-term group
$X_1f, \dots, X_rf$ by integrating the complete system:
\[
A_kf
+
B_kf
=
0
\ \ \ \ \ \ \ \ \ \ \ \ \ {\scriptstyle{(k\,=\,1\,\cdots\,r)}}.
\]

If the equations of the parameter group associated with the
$(r-q)$-term group obtained in this way are, say:
\[
a_k'
=
\varphi_k(a_1,\dots,a_{r-q},\,b_1,\dots,b_{r-q})
\ \ \ \ \ \ \ \ \ \ \ \ \  
{\scriptstyle{(k\,=\,1\,\cdots\,r\,-\,q)}},
\]
then it is clear that the parameter group associated
with the $r$-term group is representable by equations of the form:
\[
\aligned
a_k'
&
=
\varphi_k(a_1,\dots,a_{r-q},\,b_1,\dots,b_{r-q})
\ \ \ \ \ \ \ \ \ \ \ \ \  
{\scriptstyle{(k\,=\,1\,\cdots\,r\,-\,q)}}
\\
a_i'
&
=
\psi_i(a_1,\dots,a_r,\,b_1,\dots,b_r)
\ \ \ \ \ \ \ \ \ \ \ \ \ 
{\scriptstyle{(i\,=\,r\,-\,q\,+\,1\,\cdots\,r)}}.
\endaligned
\]
From this, it follows that the groups $X_1f, \dots, X_rf$ and $Y_1f,
\dots, Y_{ r-q}f$ are meroedrically isomorphic in the sense of
substitution theory.

\medskip

At present, we want to derive the result obtained just now 
yet in another way. However, we do not consider it to be necessary
to conduct in details this second, in itself noticeable, method, 
because it presents great analogies with the preceding
developments.
 
We imagine that the transformation equations:
\[
x_i'
=
f_i(x_1,\dots,x_n,\,a_1,\dots,a_r)
\ \ \ \ \ \ \ \ \ \ \ \ \ {\scriptstyle{(i\,=\,1\,\cdots\,n)}}
\]
are presented, but we leave undecided whether the $a$ are
essential parameters or not. Now, if the $f_i$ satisfy 
differential equations of the form:
\[
\frac{\partial f_i}{\partial a_k}
=
\sum_{j=1}^r\,\psi_{kj}(a_1,\dots,a_r)\,\xi_{ji}(f_1,\dots,f_n)
\]
which can be resolved with respect to the $\xi_{ ji}$:
\[
\xi_{ji}(f)
=
\sum_{k=1}^r\,\alpha_{jk}(a_1,\dots,a_r)\,
\frac{\partial f_i}{\partial a_k},
\]
then we realize easily (cf. Chap.~\ref{one-term-groups},
Sect.~\ref{application-theorem-9}) that every transformation: $x_i' =
f_i ( x, a)$ whose parameters $a$ lie in a certain neighbourhood of
$\overline{ a}_1, \dots, \overline{ a}_r$ can be obtained by executing
firstly the transformation: $\overline{ x}_i = f_i ( x, \, \overline{
a})$, and then a certain transformation:
\[
x_i'
=
w_i
(\overline{x}_1,\dots,\overline{x}_n,\,
\lambda_1,\dots,\lambda_r)
\]
of a one-term group:
\[
\lambda_1\,X_1f
+\cdots+
\lambda_r\,X_rf
=
\sum_{k=1}^r\,\lambda_k\,
\sum_{i=1}^n\,\xi_{ki}(x)\,
\frac{\partial f}{\partial x_i}.
\]
In addition, we find that the concerned values of the parameters
$\lambda$ only depend upon the form of the functions $\alpha_{ 
jk} (a)$, and on the two systems of values $\overline{ a}_k$
and $a_k$ as well.

Under the assumptions made, when one sets:
\[
\sum_{k=1}^r\,
\alpha_{jk}(a)\,
\frac{\partial f}{\partial a_k}
=
A_kf,
\ \ \ \ \ \ \ \ \ \ \ \
\sum_{i=1}^n\,\xi_{ki}(x')\,
\frac{\partial f}{\partial x_i'}
=
X_k'f, 
\]
one now obtains (compare with the pages~\pageref{S-146-sq} sq.)
that relations of the form:
\[
\aligned
\leftbracket
X_i',\,X_k'
\rightbracket
&
=
\sum_{s=1}^r\,
\vartheta_{iks}(x_1',\dots,x_n',\,a_1,\dots,a_r)\,
X_s'f
\\
\leftbracket
A_i,\,A_k
\rightbracket
&
=
\sum_{s=1}^r\,
\vartheta_{iks}(x_1',\dots,x_n',\,a_1,\dots,a_r)\,
A_sf
\endaligned
\]
hold, in which the $\vartheta_{ iks}$ are even independent of the
$x'$. However, contrary to the earlier analogous developments, it
cannot be proved now that in the \emphasis{two} latter equations, the
$\vartheta$ can be set equal to absolute constants.
But when we consider only the equations:
\[
\leftbracket
X_i',\,X_k'
\rightbracket
=
\sum_{s=1}^r\,
\vartheta_{iks}(a)\,X_s'f, 
\]
it is clear that by particularizing the $a$, they
provide relations of the form:
\[
\leftbracket
X_i',\,X_k'
\rightbracket
=
\sum_{s=1}^r\,c_{iks}\,X_s'f
\]
in which the $c_{ iks}$ denote constants. Hence it is also sure under
our present assumptions that all finite transformations $\lambda_1 \,
X_1f + \cdots + \lambda_r \, X_rf$ form a group which possesses the
same number of essential parameters as the family $x_i' = f_i ( x,
a)$.

\medskip

On the other hand, let $X_1' f, \dots, X_r'f$ denote $r$ infinitesimal
transformation that are not necessarily independent and which stand in
the relationships:
\[
\leftbracket
X_i',\,X_k'
\rightbracket
=
\sum_{s=1}^r\,c_{iks}\,X_s'f\,;
\]
furthermore, let: 
\[
A_kf
=
\sum_{i=1}^r\,\alpha_{ki}(a_1,\dots,a_r)\,
\frac{\partial f}{\partial a_i}
\ \ \ \ \ \ \ \ \ \ \ \ \ {\scriptstyle{(k\,=\,1\,\cdots\,r)}}
\]
be $r$ independent infinitesimal transformations of a simply
transitive group whose composition is given by the equations:
\[
\leftbracket
A_i,\,A_k
\rightbracket
=
\sum_{s=1}^r\,c_{iks}\,A_sf.
\]

At present, we form the $r$-term complete system:
\[
X_k'f+A_kf
=
0
\ \ \ \ \ \ \ \ \ \ \ \ \ {\scriptstyle{(k\,=\,1\,\cdots\,r)}},
\]
we compute its general solutions $F_1, \dots, F_n$ relatively
to a suitable system of values $a_1 = a_1^0$, \dots, 
$a_r = a_r^0$ and we set: $x_1 = F_1$, \dots, 
$x_n = F_n$. Afterwards, the equations resulting
by resolution: 
\[
x_k'
=
f_k(x_1,\dots,x_n,\,a_1,\dots,a_r)
\]
determine a family of at most $\infty^r$ transformations
which obviously comprises the identity transformation 
$x_k' = x_k$. 

Now, by proceeding exactly as on p.~\pageref{S-152}, we obtain
at first the system of equations:
\[
\sum_{\mu=1}^r\,\alpha_{j\mu}(a)\,
\frac{\partial x_\nu'}{\partial a_\mu}
=
\xi_{j\nu}(x'), 
\]
and then from it, by resolution:
\[
\frac{\partial x_\nu'}{\partial a_\mu}
=
\sum_{j=1}^r\,\psi_{\mu j}(a)\,\xi_{j\nu}(x').
\]

From this, we conclude: \emphasis{firstly}, that
the family $x_i' = f_i ( x, a)$ consists
of the transformations of all one-term groups
$\sum\, \lambda_k\, X_kf$, 
\emphasis{secondly}, that all these 
transformations form a \emphasis{group} with 
at most $r$ essential parameters, and lastly
\emphasis{thirdly} that from the two transformations:
\[
x_i'
=
f_i(x,a),
\ \ \ \ \ \ \ \ \ \ \
x_i''
=
f_i(x',b),
\]
a third transformation $x_i'' = f_i ( x, c)$ comes into existence
whose parameters: $c_1 = \varphi_1 ( a, b)$, \dots, $c_r = \varphi_r (
a, b)$ are determined by these two systems of values $a_k$, $b_i$ and
by the form of the functions $\alpha_{ kj} ( a)$.

But with this, the result obtained earlier on is derived in a new
manner, without that we need to enter further details. 

\medskip

If a given $r$-term group $X_1f, \dots, X_{ r-q}f, \dots, X_rf$
contains a known invariant subgroup, say $X_{ r-q+1}f, \dots, X_rf$,
then at present, we can easily indicate a meroedrically isomorphic $(r
- q)$-term group, the identity transformation of which corresponds to
the said invariant subgroup (cf. p.~\pageref{S-304}, footnote). In
fact, one forms a simply transitive group $A_1f, \dots, A_{ r-q}f,
\dots, A_rf$ in the variables $a_1, \dots, a_r$ which is equally
composed with the $r$-term group $X_1f, \dots, X_rf$. At the same
time, we can assume as earlier on that $a_1, \dots, a_{ r-q}$ are
invariants of the $q$-term group $A_{ r - q + 1}f, \dots, A_rf$.
Then if we again set: 
\[
\aligned
A_kf
=
\sum_{i=1}^{r-q}\,
&
\alpha_{ki}(a_1,\dots,a_{r-q})\,
\frac{\partial f}{\partial a_i}
+
\sum_{j=r-q+1}^r\,\beta_{kj}(a)\,
\frac{\partial f}{\partial a_j}
\\
&
\ \ \ \ \ \ \ \ \ \ \ \ \ \ \ \ \
{\scriptstyle{(k\,=\,1\,\cdots\,r\,-\,q)}},
\endaligned
\]
then the reduced infinitesimal transformations:
\[
\overline{A}_kf
=
\sum_{i=1}^{r-q}\,
\alpha_{ki}(a_1,\dots,a_{r-q})\,
\frac{\partial f}{\partial a_i}
\ \ \ \ \ \ \ \ \ \ \ \ \ 
{\scriptstyle{(k\,=\,1\,\cdots\,r\,-\,q)}}
\]
obviously generate an $(r-q)$-term group having the constitution
demanded.

\renewcommand{\thefootnote}{\fnsymbol{footnote}}
In addition, it results from these developments that \emphasis{every
proposition about the composition of $(r-q)$-term groups produces
without effort a proposition about the composition of $r$-term groups
with a $q$-term invariant subgroup}. This general principle which has
its analogue in the theory of substitutions will be exploited in the
third volume.\footnote[1]{\,
Cf. \name{Lie}, Math. Ann. Vol. XXV, p.~137.
}
\renewcommand{\thefootnote}{\arabic{footnote}}

}

\sectionengellie{\S\,\,\,104.}

In Chap.~\ref{kapitel-19}, Proposition~3, p.~\pageref{Satz-3-S-359},
we gave a remarkable form to the criterion for the similarity of
equally composed groups; namely, we showed that two equally composed
groups in the same number of variables are similar to each other if
and only if they can be related to each other in a holoedrically
isomorphic, completely specific way.

Already at that time, we announced that the criterion in question
could yet be simplified in an essential way as soon as two concerned
groups are \emphasis{transitive}; we announced that the following
theorem holds:

\renewcommand{\thefootnote}{\fnsymbol{footnote}}
\def\thetheorem{76}\begin{theorem}
\label{Theorem-76-S-425}
Two equally composed transitive groups:
\[
X_kf
=
\sum_{i=1}^n\,\xi_{ki}(x_1,\dots,x_n)\,
\frac{\partial f}{\partial x_i}
\ \ \ \ \ \ \ \ \ \ \ \ \ {\scriptstyle{(k\,=\,1\,\cdots\,r)}}
\]
and:
\[
Z_kf
=
\sum_{i=1}^n\,\zeta_{ki}(y_1,\dots,y_n)\,
\frac{\partial f}{\partial y_i}
\ \ \ \ \ \ \ \ \ \ \ \ \ {\scriptstyle{(k\,=\,1\,\cdots\,r)}}
\]
in the same number of variables are similar to each other if and only
if the following condition is satisfied: If one chooses a determined
point $x_1^0, \dots, x_n^0$ which lies on no manifold invariant by the
group: $X_1f, \dots, X_rf$, then it must be possible to relate the two
groups one to another in a holoedrically isomorphic way so that the
largest subgroup contained in the group: $X_1f, \dots, X_rf$ which
leaves invariant the point: $x_1^0, \dots, x_n^0$ corresponds to the
largest subgroup of the group: $Z_1f, \dots, Z_rf$ which leaves at
rest a certain point: $y_1^0, \dots, y_n^0$.\footnote[1]{\,
Cf. \name{Lie}, Archiv for Math. og Naturv.
Christiania 1885, p.~388 and p.~389.
} 
\end{theorem}
\renewcommand{\thefootnote}{\arabic{footnote}}

From Proposition~3, p.~\pageref{Satz-3-S-359}, it is clear that this
condition for the similarity of the two groups is necessary; hence we
need only to prove that it is also sufficient.

So, we imagine that the assumptions of the theorem are satisfied,
namely we imagine that in the group: $Z_1f, \dots, Z_rf$, independent
infinitesimal transformations: $Y_1f, \dots, Y_rf$ are chosen so that
our two groups are related to each other in the holoedrically
isomorphic way described in Theorem~76, when to every infinitesimal
transformation: $e_1\, X_1f + \cdots + e_r\, X_rf$ is associated the
infinitesimal transformation: $e_1\, Y_1f + \cdots + e_r\, Y_rf$.

According to the developments of the preceding paragraph, by means of
the concerned holoedrically isomorphic relationship, a univalent
invertible relationship between the finite transformations of the two
groups is also produced. We can describe the latter relationship
simply by interpreting, in the known way, $\sum\, e_k\, X_kf$ and
$\sum\, e_k\, Y_kf$ as symbols of the finite transformations of our
groups and in addition, by denoting for convenience the finite
transformation $\sum\, e_k\, X_kf$ shortly by: $T_{ (e)}$, and as
well, the transformation: $\sum\, e_k\, Y_kf$ shortly by: $\mathfrak{
T}_{ (e)}$.

Indeed, under these assumptions, the holoedrically isomorphic
relationship between our two groups associates the transformation
$\mathfrak{ T}_{ (e)}$ to the transformation $T_{ (e)}$, and
conversely, and this association is constituted in such a way that the
two transformations: $T_{ (e)} \, T_{ (e')}$ and $\mathfrak{ T}_{
(e)}\, \mathfrak{ T}_{ (e')}$ are related to each other, where $e_1,
\dots, e_r$ and $e_1', \dots, e_r'$ denote completely arbitrary
systems of values.

The point: $x_1^0, \dots, x_n^0$ remains invariant by exactly $r - n$
independent infinitesimal transformations of the group $X_1f, \dots,
X_rf$, and consequently, it admits exactly $\infty^{ r - n}$ different
finite transformations of this group, transformations which form, as
is known, an $(r - n)$-term subgroup. For the most general
transformation $T_{ (e)}$ which leaves at rest the point: $x_1^0,
\dots, x_n^0$, we want to introduce the symbol: $S_{ (a_1, \dots, a_{
r-n})}$, or shortly: $S_{ (a)}$; here, by $a_1, \dots, a_{ r-n}$ we
understand arbitrary parameters.

We denote by $\mathfrak{ S}_{ (a_1, \dots, a_{ r-n})}$, or shortly by
$\mathfrak{ S}_{ (a)}$, the transformations $\mathfrak{ T}_{ (e)}$ of
the group: $Z_1f, \dots, Z_rf$ which are associated to the
transformations $S_{ (a)}$ of the group: $X_1f, \dots, X_rf$; under
the assumptions made, the $\infty^{ r - n}$ transformations
$\mathfrak{ S}_{ (a)}$ then form the largest subgroup contained in the
group: $Z_1f, \dots, Z_rf$ which leaves invariant the point: $y_1^0,
\dots, y_n^0$. This is the reason why the point: $y_1^0, \dots, y_n^0$
belongs to no manifold which remains invariant by the group: $Z_1f,
\dots, Z_rf$.

\medskip

After these preparations, we conduct the following reflections.

By execution of the transformation $T_{ (e)}$ the point $x_i^0$ is
transferred to the point: $(x_i^0)\, T_{ (e)}$, whose position
naturally depends on the values of the parameters $e_1, \dots,
e_r$. According to Chap.~\ref{kapitel-14}, Proposition~1,
p.~\pageref{Satz-1-S-227}, this new point in turn admits exactly
$\infty^{ r-n}$ transformations of the group: $X_1f, \dots, X_rf$,
namely all transformations of the form:
\[
T_{(e)}^{-1}\,S_{(a)}\,T_{(e)}
\]
with the $r - n$ arbitrary parameters: $a_1, \dots, a_{ r - n}$.

On the other hand, the point $y_i^0$ is transferred by the
transformation $\mathfrak{ T}_{ (e)}$ to the point: $(y_i^0)\,
\mathfrak{ T}_{ (e)}$ which, naturally, admits exactly $\infty^{ r -
n}$ transformations of the group: $Z_1f, \dots, Z_rf$, namely all
transformations of the form:
\[
\mathfrak{T}_{(e)}^{-1}\,
\mathfrak{S}_{(a)}\,
\mathfrak{T}_{(e)}.
\]

Now, this is visibly the transformations of the group: $Z_1f, \dots,
Z_rf$ which is associated to the transformation: $T_{ (e)}^{ -1}\, S_{
(a)}\, T_{ (e)}$ of the group: $X_1f, \dots, X_rf$; so, we see that,
by means of our holoedrically isomorphic relationship, exactly the
same holds true for the points: $(x_i^0)\, T_{ (e)}$ and $(y_i^0)\,
\mathfrak{ T}_{ (e)}$ as for the points: $x_i^0$ and $y_i^0$; namely,
to the largest subgroup contained in the group $X_1f, \dots, X_rf$
which leaves invariant the point: $(x_i^0)\, T_{ (e)}$ there
corresponds the largest subgroup of the group: $Z_1f, \dots, Z_rf$ by
which the point: $(y_i^0)\, \mathfrak{ T}_{ (e)}$ remains fixed.

Now, we give to the parameters $e_1, \dots, e_r$ in the
transformations $T_{ (e)}$ and $\mathfrak{ T}_{ (e)}$ all possible
values, gradually. Then because of the transitivity of the group:
$X_1f, \dots, X_rf$, the point: $(x_i^0)\, T_{ (e)}$ is transferred
gradually to all points of the space $x_1, \dots, x_n$ that lie on no
invariant manifold, hence to all points in general position. At the
same time, the point: $(y_i^0)\, \mathfrak{ T}_{ (e)}$ is transferred
to all points in general position in the space $y_1, \dots, y_n$.

With these words, it is proved that the holoedrically isomorphic
relationship between our two groups possesses the properties put
together in Proposition~3, p.~\pageref{Satz-3-S-359}; consequently,
according to that proposition, the group: $X_1f, \dots, X_rf$ is
similar to the group: $Z_1f, \dots, Z_rf$ and the correctness of the
Theorem~76, p.~\pageref{Theorem-76-S-425} is established.

\medskip

Besides, for the proof of the Theorem~76, one does absolutely
not need to go further back to the cited proposition. 
Rather, it can be directly concluded from the 
developments above that a transformation exists which
transfers the infinitesimal transformations: $X_1f, \dots, X_rf$
to $Y_1f, \dots, Y_rf$, respectively.

We have seen that through our holoedrically isomorphic relationship,
the point: $(y_i^0)\, \mathfrak{ T}_{ (e)}$ of the space $y_1, \dots,
y_n$ is associated to the point: $(x_i^0)\, T_{ (e)}$ of the space
$x_1, \dots, x_n$, and to be precise, in this way, a univalent
invertible relationship is produced between the points of the two
spaces in general position, inside a certain region. Now, if we
interpret the $x_i$ and the $y_i$ as point coordinates of one and the
same $n$-times extended space: $R_n$, and both with respect to the
same system of coordinates, then there is a completely determined
transformation ${\sf T}$ of the space $R_n$ which always transfers the
points with the coordinates: $(x_i^0)\, T_{ (e)}$ to the point with
the coordinates: $(y_i^0)\, \mathfrak{ T}_{ (e)}$. We will prove that
this transformation ${\sf T}$ transfers the infinitesimal
transformations: $X_1f, \dots, X_rf$ to: $Y_1f, \dots, Y_rf$,
respectively, so that our two $r$-term groups are similar to each
other precisely thanks to the transformation ${\sf T}$.

The transformation ${\sf T}$ satisfies all the infinitely
many symbolic equations:
\def\theequation{14}\begin{equation}
(x_i^0)\,T_{(e')}\,{\sf T}
=
(y_i^0)\,\mathfrak{T}_{(e')},
\end{equation}
with the $r$ arbitrary parameters $e_1', \dots, e_r'$, and it is even
defined by these equations. From this, it results that ${\sf T}$
satisfies at the same time the symbolic equations:
\[
(x_i^0)\,T_{(e')}\,T_{(e)}\,{\sf T}
=
(y_i^0)\,\mathfrak{T}_{(e')}\,\mathfrak{T}_{(e)},
\]
whichever values the parameters $e$ and $e'$ may have. 

If we compare these equations with the equations~\thetag{ 14}
which we can obviously also write as:
\[
(x_i^0)\,T_{(e')}
=
(y_i^0)\,\mathfrak{T}_{(e')}\,{\sf T}^{-1},
\]
then we obtain:
\[
(y_i^0)\,\mathfrak{T}_{(e')}\,{\sf T}^{-1}\,T_{(e)}\,{\sf T}
=
(y_i^0)\,\mathfrak{T}_{(e')}\,\mathfrak{T}_{(e)},
\]
or, what is the same:
\def\theequation{15}\begin{equation}
(y_i^0)\,\mathfrak{T}_{(e')}\,{\sf T}^{-1}\,T_{(e)}\,
{\sf T}\,\mathfrak{T}_{(e)}^{-1}
=
(y_i^0)\,\mathfrak{T}_{(e')}.
\end{equation}

By appropriate choice of $e_1', \dots, e_r'$, the point: $(y_i^0)\,
\mathfrak{ T}_{ (e')}$ can be brought in coincidence with every point
in general position in the $R_m$, hence the equations~\thetag{ 15}
express that every transformation:
\def\theequation{16}\begin{equation} {\sf T}^{-1}\,T_{(e)}\,{\sf
T}\,\mathfrak{T}_{(e)}^{-1}
\end{equation}
leaves fixed all points in general position in the $R_n$. But
obviously, this is possible only when all transformations~\thetag{ 16}
coincide with the identity transformation.

With these words, it is proved that the following $\infty^r$ symbolic
equations hold:
\def\theequation{17}\begin{equation}
{\sf T}^{-1}\,T_{(e)}\,{\sf T}
=
\mathfrak{T}_{(e)},
\end{equation}
hence that the transformation ${\sf T}$ transfers every transformation
$T_{ (e)}$ of the group: $X_1f, \dots, X_rf$ to the corresponding
transformation $\mathfrak{ T}_{ (e)}$ of the group: $Z_1f, \dots,
Z_rf$. In other words, the expression: $e_1\, X_1f + \cdots + e_r\,
X_rf$ converts into the expression: $e_1\, Y_1f + \cdots + e_r\, Y_rf$
after the execution of the transformation ${\sf T}$, or, what is the
same, the $r$ infinitesimal transformations: $X_1f, \dots, X_rf$ are
transferred, by the execution of ${\sf T}$, to $Y_1f, \dots, Y_rf$,
respectively.

Consequently, the two groups: $X_1f, \dots, X_rf$ and $Z_1f, \dots,
Z_rf$ are similar to each other and the Theorem~76,
p.~\pageref{Theorem-76-S-425} is at present proved, independently of
the Proposition~3, p.~\pageref{Satz-3-S-359}.

\sectionengellie{\S\,\,\,105.}

In \S\,\,100, we saw that the same equations defined there: 
$a_k' = \varphi_k ( a, b)$ represented a group, and
in fact, we realized that between the $r$ functions
$\varphi_k ( a, b)$, the $r$ identities:
\def\theequation{5}\begin{equation}
\varphi_k
\big(
\varphi_1(a,b),\,\dots,\,\varphi_r(a,b),\,\,
c_1,\dots,c_r
\big)
\equiv
\varphi_k
\big(
a_1,\dots,a_r,\,\,
\varphi_1(b,c),\,\dots,\varphi_r(b,c)\big)
\end{equation}
held true.

Now, one convinces oneself in the same way that the equations:
\def\theequation{18}\begin{equation}
a_k'
=
\varphi_k(b_1,\dots,b_r,\,a_1,\dots,a_r)
\ \ \ \ \ \ \ \ \ \ \ \ \ {\scriptstyle{(k\,=\,1\,\cdots\,r)}}
\end{equation}
in the variables $a$ also represent a group with
the parameters $b$ and to be precise, a
simply transitive $r$-term group.

At present, a few remarks about this new group.

This group possesses the remarkable property that each
one of its transformations is interchangeable with every 
transformation of the parameter group: $a_k' = 
\varphi_k (a,b)$. 
Indeed, if one executes at first the transformation~\thetag{ 18}
and then any transformation:
\[
a_k''
=
\varphi_k(a_1',\dots,a_r',\,c_1,\dots,c_r)
\ \ \ \ \ \ \ \ \ \ \ \ \ {\scriptstyle{(k\,=\,1\,\cdots\,r)}}
\]
of the parameter group: $a_k' = \varphi_k ( a, b)$, 
then one obtains the transformation:
\[
a_k''
=
\varphi_k
\big(
\varphi_1(b,a),\,\dots,\varphi_r(b,a),\,\,
c_1,\dots,c_r
\big)
\ \ \ \ \ \ \ \ \ \ \ \ \ {\scriptstyle{(k\,=\,1\,\cdots\,r)}}
\]
which, because of the identities~\thetag{ 5}, can be
brought to the form:
\[
a_k''
=
\varphi_k
\big(b_1,\dots,b_r,\,
\varphi_1(a,c),\,\dots,\,\varphi_r(a,c)
\big)
\ \ \ \ \ \ \ \ \ \ \ \ \ {\scriptstyle{(k\,=\,1\,\cdots\,r)}}.
\]
But this transformation can be obtained
by executing at first the transformation:
\[
a_k'
=
\varphi_k(a_1,\dots,a_r,\,c_1,\dots,c_r)
\]
of the parameter group and then the transformation: 
\[
a_k''
=
\varphi_k(b_1,\dots,b_r,\,a_1',\dots,a_r')
\]
of the group~\thetag{ 18}.

\renewcommand{\thefootnote}{\fnsymbol{footnote}}
Thus, the group~\thetag{ 18} is nothing else than
the reciprocal simply transitive group associated to 
the parameter group\footnote[1]{\,
The observation that the two groups~\thetag{ 6} and~\thetag{ 18}
discussed in the text stand in a relationship such that
the transformations of the one are interchangeable
with the transformations of the other, 
is due to \name{Engel}.
}; 
\renewcommand{\thefootnote}{\arabic{footnote}}
according to Theorem~68, p.~\pageref{Theorem-68-S-380}, 
it is equally composed with the parameter group, and
even similar to it; besides, it is naturally also equally
composed with the group: $x_i' = f_i ( x, a)$ itself.

\linestop


\chapter{The Determination of All $r$-term Groups}
\label{kapitel-22}
\chaptermark{The Determination of All $r$-term Groups}

\setcounter{footnote}{0}

\abstract*{??}

\label{S-429-sq}
Already in Chap.~\ref{kapitel-17}, p.~\pageref{Proposition-1-S-297},
we have emphasized that every system of constants $c_{ iks}$ which
satisfies the relations:
\def\theequation{1}\begin{equation}
\left\{
\aligned
&
\ \ \ \ \ \ \ \ \ \ \ \ \ \ \ \ \ \ 
c_{iks}+c_{kis}=0
\\
&
\sum_{\nu=1}^r\,
\big(
c_{ik\nu}\,c_{\nu js}+c_{kj\nu}\,c_{\nu is}+c_{ji\nu}\,c_{\nu ks}
\big)
=
0
\\
&\ \ \ \ \ \ \ \ \ \ \ \ \ \ \ \ \ \ \ \ \ 
{\scriptstyle{(i,\,\,k,\,\,j,\,\,s\,=\,1\,\cdots\,r)}}.
\endaligned\right.
\end{equation}
represents a possible composition of $r$-term group and that there
always are $r$-term groups whose composition is determined just by
this system of $c_{ iks}$. As we remarked at that time, the proof for
this will be first provided in full generality in the second volume;
there, we will imagine that an arbitrary system of $c_{ iks}$ having
the said constitution is presented and we will prove that in order to
find the infinitesimal transformations of an $r$-term group having the
composition $c_{ iks}$, only the integration of simultaneous systems
of ordinary differential equations is required in any case. The
finite equations of this group can likewise be obtained by integrating
ordinary differential equations, according to
Chaps.~\ref{one-term-groups} and~\ref{kapitel-9}.

Now in the present chapter, we show two kinds of things:

\terminology{Firstly}, we imagine that the finite equations: 
\[
x_i'
=
f_i(x_1,\dots,x_n,\,a_1,\dots,a_r)
\ \ \ \ \ \ \ \ \ \ \ \ \ {\scriptstyle{(i\,=\,1\,\cdots\,n)}}
\]
of an $r$-term group are presented and we show that by integrating
simultaneous systems of ordinary differential equations, one can in
every case find all $r$-term \emphasis{transitive} groups which are
equally composed with the group: $x_i' = f_i ( x, a)$.

\terminology{Secondly}, \label{S-430} we show that one can determine all
\emphasis{intransitive} $r$-term groups without integration as soon as
one knows all \emphasis{transitive} groups with $r$ or less
parameters.

If we combine these results with what was said above and if we yet add
that according to Theorem~53, p.~\pageref{Theorem-3-S-300}, the
determination of all possible compositions of groups with given number
of parameters requires only algebraic operations, then we immediately
realize what follows:

\plainstatement{If the number $r$ is given, then, aside from
executable operations, the determination of all $r$-term groups
requires at most the integration of simultaneous systems of ordinary
differential equations.}

Naturally, the question whether the integration of the appearing
differential equations is executable or not is of great importance.
However, we can enter this question neither in this chapter, nor
in the concerned place of the second volume, since the answer
to this question presupposes the theory of integration
which cannot be developed wholly in a work about
transformation groups, and which would rather demand a
separate treatment.

\sectionengellie{\S\,\,\,106.}

In this paragraph, we put together various method for the
determination of \emphasis{simply transitive} groups, partly because
these method will find applications in the sequel, partly because they
are noticeable in themselves. 

At first, let the equations:
\def\theequation{2}\begin{equation}
x_i'
=
f_i(x_1,\dots,x_n,\,a_1,\dots,a_r)
\ \ \ \ \ \ \ \ \ \ \ \ \ {\scriptstyle{(i\,=\,1\,\cdots\,n)}}
\end{equation}
of an $r$-term group be presented. We seek the finite equations of a
simply transitive group which is equally composed with the presented
group.

According to Chap.~\ref{kapitel-21}, p.~\pageref{Theorem-71-S-404},
the parameter group of the group~\thetag{ 2} is simply transitive and
is equally composed with the group~\thetag{ 2}. Now, the finite
equations:
\[
a_k'
=
\varphi_k(a_1,\dots,a_r,\,b_1,\dots,b_r)
\ \ \ \ \ \ \ \ \ \ \ \ \ {\scriptstyle{(k\,=\,1\,\cdots\,r)}}
\]
of this group can be found by means of executable operations, namely
by resolving finite equations. Thus, we can set up the finite
equations of \emphasis{one} simply transitive group having the
constitution demanded, namely just the said parameter group.

We therefore have the:

\def\thetheorem{77}\begin{theorem}
\label{Theorem-77-S-431}
If the finite equations of an $r$-term group are presented, then one
can find the equations of an equally composed simply transitive group
by means of executable operations.
\end{theorem}

Obviously, together with this one simply transitive group, all other
simply transitive groups which are equally composed with the group:
$x_i' = f_i (x,a)$ are also given, for as a consequence of Theorem~64,
p.~\pageref{Theorem-64-S-340}, all these groups are similar to each
other.

On the other hand, we assume that not the finite equations,
but instead only the infinitesimal transformations:
\[
X_kf
=
\sum_{i=1}^n\,\xi_{ki}(x_1,\dots,x_n)\,
\frac{\partial f}{\partial x_i}
\ \ \ \ \ \ \ \ \ \ \ \ \ {\scriptstyle{(k\,=\,1\,\cdots\,r)}}
\]
of an $r$-term group are presented, and we seek to determine
the infinitesimal transformations of an equally composed
simply transitive group. 

Visibly, we have already solved this problem in Chap.~\ref{kapitel-9}, 
p.~\pageref{S-156} and~\pageref{S-157}, though it
was interpreted differently there, because at that
time, we did not have yet the concepts of
simply transitive group and of being equally composed
\deutsch{Gleichzusammengesetztseins}. 
At present, we can state as follows our previous solution 
to the problem: 

\def\theproposition{1}\begin{proposition}
\label{Satz-1-S-431}
If the infinitesimal transformations:
\[
X_kf
=
\sum_{i=1}^n\,\xi_{ki}(x_1,\dots,x_n)\,
\frac{\partial f}{\partial x_i}
\ \ \ \ \ \ \ \ \ \ \ \ \ {\scriptstyle{(k\,=\,1\,\cdots\,r)}}
\]
of an $r$-term group are presented, then one finds
the infinitesimal transformations of an 
equally composed simply transitive group in the following
way: One sets:
\[
\aligned
X_k^{(\mu)}
=
\sum_{i=1}^n\,
&
\xi_{ki}
\big(x_1^{(\mu)},\,\dots,\,x_n^{(\mu)}\big)\,
\frac{\partial f}{\partial x_i^{(\mu)}}
\ \ \ \ \ \ \ \ \ \ \ \ \ {\scriptstyle{(k\,=\,1\,\cdots\,r)}}
\\
& \ \ \ \ \ \ \ \
{\scriptstyle{(\mu\,=\,1,\,2\,\cdots\,r\,-\,1)}},
\endaligned
\]
one forms the $r$ infinitesimal transformations:
\[
W_kf
=
X_kf
+
X_k^{(1)}f
+\cdots+
X_k^{(r-1)}f
\ \ \ \ \ \ \ \ \ \ \ \ \ {\scriptstyle{(k\,=\,1\,\cdots\,r)}}
\]
and one determines $rn - r$ arbitrary independent solutions: 
$u_1, \dots, u_{ nr - r}$ of the $r$-term complete system:
\[
W_1f=0,
\,\,\,\dots,\,\,\,
W_rf=0\,;
\]
then one introduces $u_1, \dots, u_{ nr-r}$ together with
$r$ appropriate functions: $y_1, \dots, y_r$ of the $n\, r$
variables $x_i^{ (\mu)}$ as new independent variables and
in this way, one obtains the infinitesimal transformations:
\[
W_kf
=
\sum_{j=1}^r\,\eta_{kj}
(y_1,\dots,y_r,\,u_1,\dots,u_{nr-r})\,
\frac{\partial f}{\partial y_j}
\ \ \ \ \ \ \ \ \ \ \ \ \ {\scriptstyle{(k\,=\,1\,\cdots\,r)}},
\]
and if in addition one understands numerical constants
by $u_1^0, \dots, u_{ nr - r}^0$, then the $r$ independent
infinitesimal transformations:
\[
\mathfrak{W}_kf
=
\sum_{j=1}^r\,\eta_{kj}
(y_1,\dots,y_r,\,u_1^0,\dots,u_{nr-r}^0)\,
\frac{\partial f}{\partial y_j}
\ \ \ \ \ \ \ \ \ \ \ \ \ {\scriptstyle{(k\,=\,1\,\cdots\,r)}}
\]
generate a simply transitive group equally composed with the
group: $X_1f, \dots, X_rf$. If: $X_1f, \dots, X_rf$ are linked
together by the relations:
\[
\leftbracket
X_i,\,X_k
\rightbracket
=
\sum_{s=1}^r\,c_{iks}\,X_sf,
\]
then the same relations hold between $\mathfrak{ W}_1f, \dots,
\mathfrak{ W}_rf$:
\[
\leftbracket
\mathfrak{W}_i,\,\mathfrak{W}_k
\rightbracket
=
\sum_{s=1}^r\,c_{iks}\,\mathfrak{W}_sf.
\]
\end{proposition}

\label{S-432}
The process which is described in the above proposition requires the
integration of ordinary differential equations, and naturally, this
integration is not always executable. Theoretically, this is
completely indifferent, since in the present chapter, we do not manage
without integration. Nevertheless, in the developments of the next
paragraph where the question is to really set up simply transitive
groups, we will nowhere make use of the described process, and we will
apply it only in one place, where it holds, in order to prove the
existence of one simply transitive group with certain properties.

Lastly, we still want to take up the standpoint where we imagine only
that a composition of an $r$-term group is given, hence a system of
$c_{ iks}$ which satisfies the equations~\thetag{ 1}. We will show
that it is at least in very many cases possible, by means of
executable operations, to set up the finite equations of a 
simply transitive group of the composition $c_{ iks}$. 

For this, we use the Theorem~52, p.~\pageref{Theorem-52-S-296}. 
According to it, the $r$ infinitesimal transformations:
\[
E_\mu f
=
\sum_{k,\,\,j}^{1\cdots\, r}\,
c_{j\mu k}\,e_j\,
\frac{\partial f}{\partial e_k}
\ \ \ \ \ \ \ \ \ \ \ \ \ {\scriptstyle{(k\,=\,1\,\cdots\,r)}}
\]
stand pairwise in the relationships:
\[
\leftbracket
E_i,\,E_k
\rightbracket
=
\sum_{s=1}^r\,c_{iks}\,E_sf,
\]
hence they generate a linear homogeneous group in the variables:
$e_1, \dots, e_r$. This group is in particular $r$-term
when not all $r \times r$ determinants vanish whose
horizontal rows possess the form:
\[
\big\vert
c_{j1k}\,\,\,c_{j2k}
\,\,\,\cdots\,\,\,
c_{jrk}
\big\vert
\ \ \ \ \ \ \ \ \ \ \ \ \ 
{\scriptstyle{(j,\,\,k\,=\,1\,\cdots\,r)}}
\]
and it has then evidently the composition $c_{ iks}$.

Now, the finite equations of the group: $E_1f, \dots, E_rf$ can be set
up by means of executable operations (cf. p.~\pageref{S-273-bis});
hence when these $r \times r$ determinants are not all equal to zero,
we can indicate the finite equations of an $r$-term group having the
composition $c_{ iks}$. But from this, it immediately follows
(Theorem~77, p.~\pageref{Theorem-77-S-431}) that we can also indicate
the finite equations of a simply transitive group having the
composition $c_{ iks}$.

With these words, we have gained the proposition standing next:
 
\def\theproposition{2}\begin{proposition}
If a system of constants $c_{ iks}$ is presented which
satisfies the equations:
\def\theequation{1}\begin{equation}
\left\{
\aligned
&
\ \ \ \ \ \ \ \ \ \ \ \ \ \ \ \ \ \ 
c_{iks}+c_{kis}=0
\\
&
\sum_{\nu=1}^r\,
\big(
c_{ik\nu}\,c_{\nu js}+c_{kj\nu}\,c_{\nu is}+c_{ji\nu}\,c_{\nu ks}
\big)
=
0
\\
&\ \ \ \ \ \ \ \ \ \ \ \ \ \ \ \ \ \ \ \ \ 
{\scriptstyle{(i,\,\,k,\,\,j,\,\,s\,=\,1\,\cdots\,r)}}
\endaligned\right.
\end{equation}
and which is constituted in such a way that not all
$r \times r$ determinants vanish whose horizontal rows have
the form:
\[
\big\vert
c_{j1k}\,\,\,c_{j2k}
\,\,\,\cdots\,\,\,
c_{jrk}
\big\vert
\ \ \ \ \ \ \ \ \ \ \ \ \ 
{\scriptstyle{(j,\,\,k\,=\,1\,\cdots\,r)}},
\]
then by means of \terminology{executable} operations, one can always
find the finite equations of a simply transitive group having the
composition $c_{ iks}$.
\end{proposition}

Besides, the same can also be proved in a completely analogous way for
other systems of $c_{ iks}$, but we do not want to go further in that
direction. Only a few more observations.

If $r$ infinitesimal transformations: $X_1f, \dots, X_rf$ are
presented which generate an $r$-term group, then the system of $c_{
iks}$ which is associated to the group in question is given at the
same time. Now, if this system of $c_{ iks}$ is constituted
in such a way that not all determinants considered in Proposition~2
vanish, then visibly, the group: $X_1f, \dots, X_rf$ contains
no excellent infinitesimal transformation (cf. p.~\pageref{S-276}), 
hence we realize that the following holds:

\plainstatement{If $r$ independent infinitesimal transformations:
$X_1f, \dots, X_rf$ are presented which generate an $r$-term group, 
then when this group contains no excellent infinitesimal
transformation, one can set up in every case, by means
of executable operations, the finite equations of a
simply transitive group which is equally composed with the 
group: $X_1f, \dots, X_rf$.}

\sectionengellie{\S\,\,\,107.}

At present, we tackle the first one of the two problems 
the solution of which we have promised in the introduction
of the chapter.

\plainstatement{So, we imagine that the finite equations:
\def\theequation{2}\begin{equation}
x_i'
=
f_i(x_1,\dots,x_n,\,a_1,\dots,a_r)
\ \ \ \ \ \ \ \ \ \ \ \ \ {\scriptstyle{(i\,=\,1\,\cdots\,n)}}
\end{equation}
of an $r$-term group are presented and we take up the
problem of determining all \terminology{transitive} $r$-term
groups which are equally composed with the presented group.}

When an $r$-term group $\Gamma$ of any arbitrary space is transitive
and has the same composition as the group: $x_i' = f_i ( x, a)$, the
same two properties evidently hold true for all groups of the same
space that are similar to $\Gamma$. Hence, we classify the sought
groups in classes by reckoning as belonging to one and the same class
all groups of the demanded constitution which contain an equal number
of parameter and which in addition are similar to each other. We call
the totality of all groups which belong to such a class a
\label{S-434-sq}
\terminology{type} of transitive group of given composition.

This classification of the sought groups has the advantage that it is
not necessary for us to really write down all the sought
groups. Indeed, if \emphasis{one} group of the constitution demanded
is known, then at the same time, all groups similar to it are known,
hence all groups which belong to the same type. Consequently, it
suffices completely that we enumerate how many different types there
are of the sought groups and that we indicate for every individual
type a representative, hence a group which belongs to the concerned
type.

In order to settle our problem, we take the following route:

\plainstatement{At first,
\label{S-435} we indicate a process which provides the
transitive groups equally composed with the group: $x_i' = f_i ( x,
a)$. Afterwards, we conduct the proof that, by means of this process,
one can obtain every group having the concerned constitution, so that
notably, one finds at least one representative for every type of such
groups. Lastly, we give criteria to decide whether two different
groups obtained by our process are similar to each other, or not,
hence to decide whether they belong to the same type, or not.
Thanks to this, we then become in a position to have a
view of the existing mutually distinct types and at the same
time, for each one of these types, to have a representative.}

This is the program that will be carried out in the sequel.

\medskip

Let:
\[
Y_kf
=
\sum_{j=1}^r\,\eta_{kj}(y_1,\dots,y_r)\,
\frac{\partial f}{\partial y_j}
\ \ \ \ \ \ \ \ \ \ \ \ \ {\scriptstyle{(k\,=\,1\,\cdots\,r)}}
\]
be independent infinitesimal transformations of a simply transitive
group which is equally composed with the group: $x_i' = f_i ( x,
a)$. Moreover, let:
\[
Z_kf
=
\sum_{j=1}^r\,\zeta_{kj}(y_1,\dots,y_r)\,
\frac{\partial f}{\partial y_j}
\ \ \ \ \ \ \ \ \ \ \ \ \ {\scriptstyle{(k\,=\,1\,\cdots\,r)}}
\]
be independent infinitesimal transformations of the associated
reciprocal group which, according to Theorem~68,
p.~\pageref{Theorem-68-S-380}, is simply transitive as well, and is
equally composed with the group: $x_i' = f_i ( x,a)$.

Evidently, all these assumptions can be satisfied. For instance, two
simply transitive groups having the constitution indicated are the
parameter group of the group: $x_i' = f_i ( x,a)$ defined in the
previous chapter:
\[
a_k'
=
\varphi_k
(a_1,\dots,a_r,\,b_1,\dots,b_r)
\ \ \ \ \ \ \ \ \ \ \ \ \ {\scriptstyle{(k\,=\,1\,\cdots\,r)}}
\]
and its reciprocal group:
\[
a_k'
=
\varphi_k(b_1,\dots,b_r,\,a_1,\dots,a_r)
\ \ \ \ \ \ \ \ \ \ \ \ \ {\scriptstyle{(k\,=\,1\,\cdots\,r)}}
\]
and the finite equations of both may even be set up by means of
executable operations.

After these preparations, we turn to the explanation of the process
announced above which produces transitive groups that are equally
composed with the group: $x_i' = f_i ( x, a)$.

In Chap.~\ref{kapitel-17}, p.~\pageref{S-305-sq}, we showed how every
decomposition of the space invariant by a group can be used in order
to set up an isomorphic group. We want to apply this to the group:
$Y_1f, \dots, Y_rf$.

On the basis of Theorem~69, p.~\pageref{Theorem-69-S-387}, we seek an
arbitrary decomposition of the space $y_1, \dots, y_r$ invariant by
the group $Y_1f, \dots, Y_rf$. By means of algebraic operations, we
determine and $m$-term subgroup of the group: $Z_1f, \dots, Z_rf$.
If:
\def\theequation{3}\begin{equation}
{\sf Z}_\mu f
=
\varepsilon_{\mu_1}\,Z_1f
+\cdots+
\varepsilon_{\mu_r}\,Z_rf
\ \ \ \ \ \ \ \ \ \ \ \ \ {\scriptstyle{(\mu\,=\,1\,\cdots\,m)}}
\end{equation}
are independent infinitesimal transformations of this subgroup, we
form the $m$-term complete system:
\[
{\sf Z}_1f=0,
\,\,\,\dots,\,\,\,
{\sf Z}_mf=0
\]
and we compute, by integrating it, $r - m$ arbitrary invariants:
\[
u_1(y_1,\dots,y_r),
\,\,\,\dots,\,\,\,
u_{r-m}(y_1,\dots,y_r)
\]
of the group: ${\sf Z}_1f, \dots, {\sf Z}_mf$. Then the
equations:
\def\theequation{4}\begin{equation}
u_1(y_1,\dots,y_r)
=
{\rm const.},
\,\,\,\dots,\,\,\,
u_{r-m}(y_1,\dots,y_r)
=
{\rm const.}
\end{equation}
determine a decomposition of the space $y_1, \dots, y_r$ in $\infty^{
r - m}$ $m$-times extended manifolds invariant by the group: $Y_1f,
\dots, Y_rf$.

Now, we introduce $u_1 (y), \dots, u_{ r-m} (y)$ together with 
$m$ other appropriate functions: $v_1, \dots, v_m$ of the $y$
as new variables in $Y_1f, \dots, Y_rf$ and we obtain:
\[
\aligned
Y_kf
=
\sum_{\nu=1}^{r-m}\,
\omega_{k\nu}(u_1,\dots,u_{r-m})\,
&
\frac{\partial f}{\partial u_\nu}
+
\sum_{\mu=1}^m\,
w_{k\mu}(u_1,\dots,u_{r-m},\,
v_1,\dots,v_m)\,
\frac{\partial f}{\partial v_\mu}
\\
& 
\ \ \ \ \ {\scriptstyle{(k\,=\,1\,\cdots\,r)}}.
\endaligned
\]
From this, we lastly form the reduced infinitesimal transformations:
\[
U_kf
=
\sum_{\nu=1}^{r-m}\,
\omega_{k\nu}(u_1,\dots,u_{r-m})\,
\frac{\partial f}{\partial u_\nu}
\ \ \ \ \ \ \ \ \ \ \ \ \ {\scriptstyle{(k\,=\,1\,\cdots\,r)}}.
\]
According to Chap.~\ref{kapitel-17}, Proposition~4,
p.~\pageref{Satz-4-S-307}, they generate a group in the variables:
$u_1, \dots, u_{ r-m}$ which is isomorphic with the group: $Y_1f, \dots,
Y_rf$.

If we had chosen, instead of $u_1 (y), \dots, u_{ r-m} (y)$, any other
independent invariants: $u_1'(y), \dots, u_{ r-m}'(y)$ of the group:
${\sf Z}_1f, \dots, {\sf Z}_mf$, then in place of the group: $U_1f,
\dots, U_rf$, we would have obtained another group in the variables
$u_1', \dots, u_{ r-m}'$, but this group would visibly be similar to
the group: $U_1f, \dots, U_rf$, because $u_1'(y), \dots, u_{ r-m}'
(y)$ are independent functions of: $u_1 (y), \dots, u_{ r-m} (y)$. If
we would replace $u_1 (y), \dots, u_{ r-m} (y)$ by the most general
system of $r - m$ independent invariants of the group: ${\sf Z}_1f,
\dots, {\sf Z}_mf$, then we would obtain the most general group in $r
- m$ variables which is similar to the group: $U_1f, \dots, U_rf$.

From the fact that the group: $Y_1f, \dots, Y_rf$ is simply
transitive, it follows, as we have observed already in
Chap.~\ref{kapitel-20}, p.~\pageref{S-388}, that not all $(r - m)
\times (r - m)$ determinants of the matrix:
\[
\left\vert
\begin{array}{cccc}
\omega_{11}(u) & \,\cdot\, & \,\cdot\, & \omega_{1,\,r-m}(u)
\\
\cdot & \,\cdot\, & \,\cdot\, & \cdot
\\
\cdot & \,\cdot\, & \,\cdot\, & \cdot
\\
\omega_{r1}(u) & \,\cdot\, & \,\cdot\, & \omega_{r,\,r-m}(u)
\end{array}
\right\vert
\]
vanish identically; expressed differently: it results that the group:
$U_1f, \dots, U_{ r-m}f$ in the $r - m$ variables $u$ is transitive.

\medskip

We therefore have a method for setting up transitive groups which are
isomorphic with the group: $Y_1f, \dots, Y_rf$. However, this does not
suffices, for we demand transitive groups that are equally composed
with the group: $Y_1f, \dots, Y_rf$, hence are holoedrically
isomorphic to it. Thus, the group: $U_1f, \dots, U_rf$ is useful for
us only when it is $r$-term. Under which conditions is it so?

Since the group: $U_1f, \dots, U_rf$ is transitive, it contains at
least $r - m$ essential parameters, hence it will in general contain
exactly $r - l$ essential parameters, where $0 \leqslant l \leqslant
m$. Then according to p.~\pageref{S-307}, in the group: $Y_1f, \dots,
Y_rf$, there are exactly $l$ independent infinitesimal transformations
which leave individually fixed each one of the $\infty^{ r-m}$
manifolds~\thetag{ 4}, and these infinitesimal transformations
generate an $l$-term invariant subgroup of the group: $Y_1f, \dots,
Y_rf$. For reasons of brevity, we want to denote the concerned
subgroup by $g$.

If ${\sf M}$ is an arbitrary manifold amongst the generally located
manifolds~\thetag{ 4}, then ${\sf M}$ admits exactly $m$ independent
infinitesimal transformations: $e_1\, Y_1f + \cdots + e_r\, Y_rf$
which generate an $m$-term subgroup $\gamma$ of the group: $Y_1f,
\dots, Y_rf$ (cf. page~\pageref{S-388-bis}). Naturally, the invariant
subgroup $g$ is contained in this group $\gamma$. On the other hand,
${\sf M}$ admits exactly $m$ independent infinitesimal transformations
of the reciprocal group: $Z_1f, \dots, Z_rf$, namely: ${\sf Z}_1f,
\dots, {\sf Z}_mf$, which also generate an $m$-term group. Therefore,
according to Chap.~\ref{kapitel-20}, p.~\pageref{S-390}, one can
relate the two simply transitive groups: $Y_1f, \dots, Y_rf$ and
$Z_1f, \dots, Z_rf$ to each other in a holoedrically isomorphic way so
that the subgroup $\gamma$ corresponds to the subgroup: ${\sf Z}_1f,
\dots, {\sf Z}_mf$. On the occasion, to the invariant subgroup $g$,
there visibly corresponds an $l$-term invariant subgroup $g'$ of the
group: $Z_1f, \dots, Z_rf$, and in fact, $g'$ is contained in the
subgroup: ${\sf Z}_1f, \dots, {\sf Z}_mf$.

We therefore see: when the group: $U_1f, \dots, U_rf$
is exactly $(r - l)$-term, then the group: ${\sf Z}_1f, 
\dots, {\sf Z}_mf$ contains an $l$-term subgroup $g'$
which is invariant in the group: $Z_1f, \dots, Z_rf$. 

\medskip

Conversely: when the group: ${\sf Z}_1f, \dots, {\sf Z}_mf$
contains an $l$-term subgroup:
\[
{\sf Z}_\lambda'f
=
h_{\lambda 1}\,{\sf Z}_1f
+\cdots+
h_{\lambda m}\,{\sf Z}_mf
\ \ \ \ \ \ \ \ \ \ \ \ \ 
{\scriptstyle{(\lambda\,=\,1\,\cdots\,l)}}
\]
which is invariant in the group: $Z_1f, \dots, Z_rf$, then the group:
$U_1f, \dots, U_rf$ can at most be $(r - l)$-term.

Indeed, under the assumptions just made, the $l$ mutually
independent equations:
\def\theequation{5}\begin{equation}
{\sf Z}_1'f=0,
\,\,\,\dots,\,\,\,
{\sf Z}_l'f=0
\end{equation}
form an $l$-term complete system with $r-l$ independent
solutions:
\[
\psi_1(y_1,\dots,y_r),
\,\,\,\dots,\,\,\,
\psi_{r-l}(y_1,\dots,y_r).
\]
Furthermore, there are relations of the form:
\[
\aligned
\leftbracket
{\sf Z}_k,\,{\sf Z}_\lambda'
\rightbracket
&
=
h_{k\lambda 1}\,{\sf Z}_1'f
+\cdots+
h_{k\lambda l}\,{\sf Z}_l'f
\\
& \ \ \ \
{\scriptstyle{(k\,=\,1\,\cdots\,r\,;\,\,\,
\lambda\,=\,1\,\cdots\,l)}}
\endaligned
\]
which express that the complete system~\thetag{ 5} admits the
group: $Z_1f, \dots, Z_rf$. Consequently, the equations:
\def\theequation{6}\begin{equation}
\psi_1(y_1,\dots,y_r)
=
{\rm const.},
\,\,\,\dots,\,\,\,
\psi_{r-l}(y_1,\dots,y_r)
=
{\rm const.}
\end{equation}
represent a decomposition of the space $y_1, \dots, y_r$ invariant by
the group $Z_1f, \dots, Z_rf$, and to be precise, a decomposition in
$\infty^{ r - l}$ $l$-times extended manifolds.

These $\infty^{ r-l}$ manifolds stand in a very simple relationship
with respect to the $\infty^{ r-m}$ manifolds:
\def\theequation{4}\begin{equation}
u_1(y_1,\dots,y_r)
=
{\rm const.},
\,\,\,\dots,\,\,\,
u_{r-m}(y_1,\dots,y_r)
=
{\rm const.},
\end{equation}
namely each one of the manifolds~\thetag{ 4} consists of
$\infty^{ m - l}$ different manifolds~\thetag{ 6}. 
This follows without effort from the fact that
$u_1 ( y), \dots, u_{ r - m} (y)$, as solutions of the complete
system:
\[
{\sf Z}_1f=0,
\,\,\,\dots,\,\,\,
{\sf Z}_mf=0,
\]
satisfy simultaneously the complete system~\thetag{ 5}, and can
therefore be represented as functions of: $\psi_1 (y), \dots, \psi_{
r-l} (y)$.

Now, according to Chap.~\ref{kapitel-20}, p.~\pageref{Satz-7-S-387},
the reciprocal group: $Y_1f, \dots, Y_rf$ contains exactly $l$
independent infinitesimal transformations which leave individually
fixed each one of the $\infty^{ r-l}$ manifolds~\thetag{ 6}, hence it
obviously contains at least $l$ independent infinitesimal
transformations which leave individually fixed each one of the
$\infty^{ r-m}$ manifolds~\thetag{ 4}; but from this, it immediately
follows that, under the assumption made, the group: $U_1f, \dots, U_rf$
can be at most $(r-l)$-term, as we claimed above.

\medskip

Thanks to the preceding developments, it is proved that the group:
$U_1f, \dots, U_rf$ is $(r-l)$-term if and only if the group: ${\sf
Z}_1f, \dots, {\sf Z}_mf$ contains an $l$-term subgroup invariant in
the group: $Z_1f, \dots, Z_rf$, but no larger subgroup of the same
nature. In particular, the group: $U_1f, \dots, U_rf$ is $r$-term if
and only if the group: ${\sf Z}_1f, \dots, {\sf Z}_mf$ contains, aside
from the identity transformation, no subgroup invariant in the group
$Z_1f, \dots, Z_rf$. Here, the word `subgroup' is to be understood in
its widest sense, so that one then also has to consider the group:
${\sf Z}_1f, \dots, {\sf Z}_mf$ itself as a subgroup contained in it.

By summarizing the gained result, we can at present say:

\def\thetheorem{78}\begin{theorem}
\label{Theorem-78-S-439}
If the two $r$-term groups:
\[
Y_kf
=
\sum_{j=1}^r\,\eta_{kj}(y_1,\dots,y_r)\,
\frac{\partial f}{\partial y_j}
\ \ \ \ \ \ \ \ \ \ \ \ \ {\scriptstyle{(k\,=\,1\,\cdots\,r)}}
\]
and:
\[
Z_kf
=
\sum_{j=1}^r\,\zeta_{kj}(y_1,\dots,y_r)\,
\frac{\partial f}{\partial y_j}
\ \ \ \ \ \ \ \ \ \ \ \ \ {\scriptstyle{(k\,=\,1\,\cdots\,r)}}
\]
are simply transitive and reciprocal to each other, if, 
moreover: 
\[
{\sf Z}_\mu f
=
\varepsilon_{\mu_1}\,Z_1f
+\cdots+
\varepsilon_{\mu r}\,Z_rf
\ \ \ \ \ \ \ \ \ \ \ \ \ {\scriptstyle{(\mu\,=\,1\,\cdots\,m)}}
\]
is an $m$-term subgroup of the group: $Z_1f, \dots, Z_rf$, and
if:
\[
u_1(y_1,\dots,y_r),
\,\,\dots,\,\,
u_{r-m}(y_1,\dots,y_r)
\]
are independent invariants of this subgroup, then the $r$ 
infinitesimal transformations:
\[
\aligned
\sum_{\nu=1}^{r-m}\,
Y_k\,u_\nu\,
\frac{\partial f}{\partial u_\nu}
&
=
\sum_{\nu=1}^{r-m}\,
\omega_{k\nu}(u_1,\dots,u_{r-m})\,
\frac{\partial f}{\partial u_\nu}
=
U_kf
\\
&
\ \ \ \ \ \ \ \ \ \ \ \ \
{\scriptstyle{(k\,=\,1\,\cdots\,r)}}
\endaligned
\]
in the $r - m$ variables: $u_1, \dots, u_{ r-m}$ generate a transitive
group, isomorphic with the group: $Y_1f, \dots, Y_rf$. This group is
$(r - l)$-term when there is in the group: ${\sf Z}_1f, \dots, {\sf
Z}_mf$ an $l$-term subgroup, but no larger subgroup, which is
invariant in the group: $Z_1f, \dots, Z_rf$. In particular, it is
$r$-term and equally composed with the group: $Y_1f, \dots, Y_rf$ if
and only if the group: ${\sf Z}_1f, \dots, {\sf Z}_mf$ neither is
invariant itself in the group: $Z_1f, \dots, Z_rf$, nor contains,
aside from the identity transformation, another subgroup invariant in
the group: $Z_1f, \dots, Z_rf$.

If one replaces $u_1(y), \dots, u_{ r-m} (y)$ by the most general
system: $u_1'(y), \dots, u_{ r-m}'(y)$ of $r-m$ independent invariants
of the group: ${\sf Z}_1f, \dots, {\sf Z}_mf$, then in place of the
group: $U_1f, \dots,U_rf$, one obtains the most general group in $r -
m$ variables similar to it. In particular, if the group: $U_1f, \dots,
U_rf$ is $r$-term, then one obtains in this way the most general
transitive group equally composed with the group: $Y_1f, \dots, Y_rf$
that belongs to the same type as the group: $U_1f, \dots, U_rf$.
\end{theorem}

In addition, from the developments used for the proof of this theorem,
it yet results the following

\def\theproposition{3}\begin{proposition}
If $Y_1f, \dots, Y_rf$ and $Z_1f, \dots, Z_rf$ are two reciprocal
simply transitive groups and if: $Z_1f, \dots, Z_lf$ is an invariant
$l$-term subgroup of the second group, then the invariants of this
$l$-term group can also be defined as the invariants of a certain
$l$-term group which is contained as an invariant subgroup in the
$r$-term group: $Y_1f, \dots, Y_rf$.
\end{proposition}

\smallskip

Thanks to the Theorem~78, the first part of the program stated on
p.~\pageref{S-435} is settled, and we are in possession of a process
which provides transitive groups equally composed with the group: $x_i
' = f_i ( x, a)$. We now come to the second part, namely to the proof
that every group of this sort can be found thanks to our process.

In $r - m$ variables, let an arbitrary transitive group be presented
which is equally composed with the group: $x_i' = f_i ( x, a)$; let
its infinitesimal transformations be:
\[
\mathfrak{X}_kf
=
\sum_{\nu=1}^{r-m}\,
\mathfrak{x}_{k\nu}(z_1,\dots,z_{r-m})\,
\frac{\partial f}{\partial z_\nu}
\ \ \ \ \ \ \ \ \ \ \ \ \ {\scriptstyle{(k\,=\,1\,\cdots\,r)}}
\]

At first, we prove that amongst the simply transitive groups of the
same composition, there is in any case one which can be obtained from
the group: $\mathfrak{ X}_1f, \dots, \mathfrak{ X}_rf$ in a way
completely analogous to the way in which the group: $U_1f, \dots,
U_rf$ was obtained from the simply transitive group: $Y_1f, \dots,
Y_rf$.

To this end, under the guidance of Proposition~1,
p.~\pageref{Satz-1-S-431}, we form, in the $r ( r -m)$ variables: $z$,
$z^{ (1)}$, \dots, $z^{ (r-1)}$, the $r$ infinitesimal
transformations:
\[
W_kf
=
\mathfrak{X}_kf
+
\mathfrak{X}_k^{(1)}f
+\cdots+
\mathfrak{X}_k^{(r-1)}f
\ \ \ \ \ \ \ \ \ \ \ \ \ {\scriptstyle{(k\,=\,1\,\cdots\,r)}}.
\]
Now, if $\varphi_1, \dots, \varphi_R$ are $(r-m)r - r = R$ independent
solutions of the $r$-term complete system:
\def\theequation{7}\begin{equation}
W_1f=0,
\,\,\,\dots,\,\,\,
W_rf=0,
\end{equation}
then we introduce them, together with $z_1, \dots, z_{ r-m}$ and yet
together with $m$ functions $\overline{ z}_1, \dots, \overline{ z}_m$
of the $z$, $z^{ (1)}$, \dots, $z^{ (r-1)}$, as new variables. This is
possible, since the equations~\thetag{ 7} are solvable with respect to
$\partial f / \partial z_1$, \dots, $\partial f / \partial z_{ r-m}$
because of the transitivity of the group: $\mathfrak{ X}_1f, \dots,
\mathfrak{ X}_rf$, whence $z_1, \dots, z_{ r-m}$ are independent of
the functions: $\varphi_1, \dots, \varphi_R$. In the new variables,
$W_1f, \dots, W_rf$ receive the form:
\[
\aligned
W_kf
=
\mathfrak{X}_kf
+
\sum_{\mu=1}^m\,
&
w_{k\mu}(z_1,\dots,z_{r-m},\,
\overline{z}_1,\dots,\overline{z}_m,\,
\varphi_1,\dots,\varphi_R)\,
\frac{\partial f}{\partial\overline{z}_\mu}
\\
&
\ \ \ \ \ \ \ \ \ \ \ \ \ 
{\scriptstyle{(k\,=\,1\,\cdots\,r)}},
\endaligned
\]
and here, when we confer to the $\varphi$ suitable fixed
values $\varphi_1^0, \dots, \varphi_R^0$ and when we set:
\[
w_k
(z_1,\dots,z_{r-m},\,
\overline{z}_1,\dots,\overline{z}_m,\,
\varphi_1^0,\dots,\varphi_R^0)
=
w_k^0
(z_1,\dots,z_{r-m},\,
\overline{z}_1,\dots,\overline{z}_m),
\]
then:
\[
\mathfrak{W}_kf
=
\mathfrak{X}_kf
+
\sum_{\mu=1}^m\,
w_{k\mu}^0
(z_1,\dots,z_{r-m},\,
\overline{z}_1,\dots,\overline{z}_m)\,
\frac{\partial f}{\partial\overline{z}_\mu}
\ \ \ \ \ \ \ \ \ \ \ \ \ 
{\scriptstyle{(k\,=\,1\,\cdots\,r)}}
\]
are independent infinitesimal transformations of a simply transitive
group which is equally composed with the group: $\mathfrak{ X}_1f,
\dots, \mathfrak{ X}_rf$ and therefore, also equally composed with the
group: $x_i' = f_i ( x, a)$.

With this, we have found a simply transitive group having the
constitution indicated earlier on. 

Indeed, the equations:
\[
z_1
=
{\rm const.},
\,\,\,\dots,\,\,\,
z_{r-m}
=
{\rm const.}
\]
obviously determine a decomposition of the space: $z_1, \dots, z_{
r-m}$, $\overline{ z}_1, \dots, \overline{ z}_m$ in $\infty^{ r-m}$
$m$-times extended manifolds invariant by the group: $\mathfrak{
W}_1f, \dots, \mathfrak{ W}_rf$. The group: $\mathfrak{ X}_1f, \dots,
\mathfrak{ X}_rf$ indicates in which way these $\infty^{ r - m}$
manifolds are permuted by the group: $\mathfrak{ W}_1f, \dots,
\mathfrak{ W}_rf$. So, between the two groups: $\mathfrak{ X}_1f,
\dots, \mathfrak{ X}_rf$ and $\mathfrak{ W}_1f, \dots, \mathfrak{
W}_rf$, there is a relationship completely analogous to the one above
between the two groups: $U_1f, \dots, U_rf$ and $Y_1f, \dots, Y_rf$.

At present, there is no difficulty to prove that the group:
$\mathfrak{ X}_1f, \dots, \mathfrak{ X}_rf$ is similar to one of the
groups that we obtain when we apply the process described in
Theorem~78, p.~\pageref{Theorem-78-S-439} to two determined
simply transitive reciprocal groups: $Y_1f, \dots, Y_rf$ 
and $Z_1f, \dots, Z_rf$ having the concerned composition.

Let: ${\sf B}_1f, \dots, {\sf B}_rf$ be the simply transitive
group reciprocal to: $\mathfrak{ W}_1f, \dots, 
\mathfrak{ W}_rf$. Naturally, this group is equally composed
with the group: $Z_1f, \dots, Z_rf$ and
hence, also similar to it (cf. Chap.~\ref{kapitel-19}, 
Theorem~64, p.~\pageref{Theorem-64-S-340}).

We want to assume that the transformation:
\def\theequation{8}\begin{equation}
\left\{
\aligned
z_1
&
=
z_1(y_1,\dots,y_r),
\,\,\,\dots,\,\,\,
z_{r-m}
=
z_{r-m}(y_1,\dots,y_r)
\\
\overline{z}_1
&
=
\overline{z}_1(y_1,\dots,y_r),
\,\,\,\dots,\,\,\,
\overline{z}_m
=
\overline{z}_m(y_1,\dots,y_r)
\endaligned\right.
\end{equation}
transfers the group: $Z_1f, \dots, Z_rf$ to the group: ${\sf B}_1f,
\dots, {\sf B}_rf$. Then according to Chap.~\ref{kapitel-20},
p.~\pageref{S-381}, through the same transformation, the group: $Y_1f,
\dots, Y_rf$ is transferred at the same time to the group: $\mathfrak{
W}_1f, \dots, \mathfrak{ W}_rf$, hence by virtue of~\thetag{ 8}, there
are relations of the form:
\[
Y_kf
=
\sum_{j=1}^r\,
\mathfrak{h}_{kj}\,\mathfrak{W}_jf
\ \ \ \ \ \ \ \ \ \ \ \ \ {\scriptstyle{(k\,=\,1\,\cdots\,r)}},
\]
where the determinant of the constants $\mathfrak{ h}_{ kj}$ does not
vanish. But from these relations, the following relations immediately
follow:
\def\theequation{9}\begin{equation}
\sum_{\nu=1}^{r-m}\,Y_k\,z_\nu\,
\frac{\partial f}{\partial z_\nu}
=
\sum_{j=1}^r\,\mathfrak{h}_{kj}\,
\mathfrak{X}_jf
\ \ \ \ \ \ \ \ \ \ \ \ \ {\scriptstyle{(k\,=\,1\,\cdots\,r)}},
\end{equation}
and likewise, they hold identically by virtue of~\thetag{ 8}.

This is the reason why the equations:
\[
z_1(y_1,\dots,y_r)
=
{\rm const.},
\,\,\,\dots,\,\,\,
z_{r-m}(y_1,\dots,y_r)
=
{\rm const.}
\]
represent a decomposition of the space $y_1, \dots, y_r$ invariant by
the group: $Y_1f, \dots, Y_rf$, or, what amount to the same, the
reason why $z_1 ( y), \dots, z_{ r-m} (y)$ are independent invariants
of a completely determined $m$-term subgroup $\mathfrak{ g}$ of the
group $Z_1f, \dots, Z_rf$. Consequently, we obtain the group:
$\mathfrak{ X}_1f, \dots, \mathfrak{ X}_rf$ thanks to the process
described in Theorem~78, p.~\pageref{Theorem-78-S-439} when we arrange
ourselves as follows: As group: ${\sf Z}_1f, \dots, {\sf Z}_mf$, we
choose the $m$-term subgroup $\mathfrak{ g}$ of the group: $Z_1f,
\dots, Z_rf$ just defined, and as functions: $u_1 (y), \dots, u_{ r-m}
(y)$, we choose the just said invariants: $z_1 (y), \dots, z_{ r-m}
(y)$ of the group $\mathfrak{ g}$. Indeed, under these assumptions,
the relations~\thetag{ 9} hold true, in which the right-hand side
expressions are independent infinitesimal transformations of the
group: $\mathfrak{ X}_1f, \dots, \mathfrak{ X}_rf$.

As a result, it is proved that, thanks to the process which is
described in Theorem~78, p.~\pageref{Theorem-78-S-439}, one can find
every transitive group isomorphic with the group: $x_i' = f_i (x,
a)$. We can therefore enunciate the following theorem:

\renewcommand{\thefootnote}{\fnsymbol{footnote}}
\def\thetheorem{79}\begin{theorem}
\label{Theorem-79-S-443}
If the finite equations:
\[
x_i'
=
f_i(x_1,\dots,x_n,\,a_1,\dots,a_r)
\ \ \ \ \ \ \ \ \ \ \ \ \ {\scriptstyle{(i\,=\,1\,\cdots\,n)}}
\]
of an $r$-term group are presented, then one finds in the following
way all transitive groups that are equally composed
with this group:

To begin with, one determines, which requires only executable
operations, two $r$-term simply transitive groups:
\[
Y_kf
=
\sum_{j=1}^r\,\eta_{kj}(y_1,\dots,y_r)\,
\frac{\partial f}{\partial y_j}
\ \ \ \ \ \ \ \ \ \ \ \ \ {\scriptstyle{(k\,=\,1\,\cdots\,r)}}
\]
and:
\[
Z_kf
=
\sum_{j=1}^r\,\zeta_{kj}(y_1,\dots,y_r)\,
\frac{\partial f}{\partial y_j}
\ \ \ \ \ \ \ \ \ \ \ \ \ {\scriptstyle{(k\,=\,1\,\cdots\,r)}}
\]
that are reciprocal to each other and are equally composed with the
group presented. Afterwards, by means of algebraic operations, one
sets up all subgroups of the group: $Z_1f, \dots, Z_rf$ which neither
are invariant in this group, nor contain, aside from the identity
transformation, a subgroup invariant in the group: $Z_1f, \dots,
Z_rf$. Under the guidance of Theorem~78, each one of the found
subgroups produces all transitive groups equally composed with the
group: $x_i' = f_i ( x,a)$ that belong to a certain type. If one
determines these groups for each one of the found subgroups, then one
obtains all transitive groups that are equally composed with the
group: $x_i' = f_i ( x, a)$.\footnote[1]{\,
\name{Lie}, Archiv for Math., Vol. 10, Christiania 1885, and
Verh. der Gesellsch. d. W. zu Chr. a. 1884; Berichte der Kgl. Sächs.
G. d. W., 1888.
} 
\end{theorem}
\renewcommand{\thefootnote}{\arabic{footnote}}

The second part of our program is now carried out, and we are in the
position to identify all types of transitive groups equally composed
with the group: $x_i' = f_i ( x, a)$. Every subgroup of the group:
$Z_1f, \dots, Z_rf$ which has the property mentioned in the Theorem~79
provides us with such a type. It yet remains to decide when 
two different subgroups of the group: $Z_1f, \dots, Z_rf$
produce different types, and when they produce the same type.

Let the $m$ independent infinitesimal transformations:
\label{S-443-sq}
\[
{\sf Z}_\mu f
=
\sum_{k=1}^r\,
\varepsilon_{\mu k}\,Z_kf
\ \ \ \ \ \ \ \ \ \ \ \ \ {\scriptstyle{(\mu\,=\,1\,\cdots\,m)}}
\]
generate an $m$-term group which neither is invariant in the
group: $Z_1f, \dots, Z_rf$, nor contains, aside from the
identity transformation, another subgroup invariant in the group:
$Z_1f, \dots, Z_rf$. 

\label{S-444-445}
Let the functions:
\[
u_1(y_1,\dots,y_r),
\,\,\dots,\,\,
u_{r-m}(y_1,\dots,y_r)
\]
be independent invariants of the group: ${\sf Z}_1f, \dots, 
{\sf Z}_mf$. 
Under these assumptions, according to Theorem~79, the group:
\[
U_kf
=
\sum_{\nu=1}^{r-m}\,
Y_k\,u_\nu\,
\frac{\partial f}{\partial u_\nu}
=
\sum_{\nu=1}^{r-m}\,
\omega_{k\nu}(u_1,\dots,u_{r-m})\,
\frac{\partial f}{\partial u_\nu}
\ \ \ \ \ \ \ \ \ \ \ \ \ {\scriptstyle{(k\,=\,1\,\cdots\,r)}}
\]
in the $r - m$ variables: $u_1, \dots, u_{ r-m}$ has the same
composition as the group: $x_i' = f_i ( x,a)$, and is in addition
transitive, hence it is a representative of the type of groups which
corresponds to the subgroup: ${\sf Z}_1f, \dots, {\sf Z}_mf$.

If another subgroup of the group: $Z_1f, \dots, Z_rf$ is supposed to
provide the same type of group, then 
this subgroup must evidently be $m$-term,
for it is only in this case that it can provide transitive groups in
$r - m$ variables that are equally composed with the group: $x_i' =
f_i ( x,a)$.

So, we assume that:
\[
\mathfrak{Z}_\mu f
=
\sum_{k=1}^r\,\mathfrak{e}_{\mu k}\,Z_kf
\ \ \ \ \ \ \ \ \ \ \ \ \ {\scriptstyle{(\mu\,=\,1\,\cdots\,m)}}
\]
is another $m$-term subgroup of the group: $Z_1f, \dots, Z_rf$, 
and that this subgroup too neither is invariant in the
group: $Z_1f, \dots, Z_rf$, nor contains, aside from the identity
transformation, another subgroup invariant in the group:
$Z_1f, \dots, Z_rf$. Let independent invariants of the
group: $\mathfrak{ Z}_1f, \dots, \mathfrak{ Z}_mf$ be:
\[
\mathfrak{u}_1(y_1,\dots,y_r),
\,\,\dots,\,\,
\mathfrak{u}_{r-m}(y_1,\dots,y_r).
\]

Under these assumptions, the group:
\[
\mathfrak{U}_kf
=
\sum_{\nu=1}^{r-m}\,
Y_k\,\mathfrak{u}_\nu\,
\frac{\partial f}{\partial\mathfrak{u}_\nu}
=
\sum_{\nu=1}^{r-m}\,
\mathfrak{o}_{k\nu}
(\mathfrak{u}_1,\dots,\mathfrak{u}_{r-m})\,
\frac{\partial f}{\partial\mathfrak{u}_\nu}
\ \ \ \ \ \ \ \ \ \ \ \ \ {\scriptstyle{(k\,=\,1\,\cdots\,r)}}
\]
in the $r - m$ variables: $\mathfrak{ u}_1, \dots, \mathfrak{ u}_{
r-m}$ is transitive and has the same composition as the group: $x_i' =
f_i ( x, a)$.

The question whether the two subgroups: ${\sf Z}_1f, \dots, {\sf
Z}_mf$ and $\mathfrak{ Z}_1f, \dots, \mathfrak{ Z}_mf$ of the group:
$Z_1f, \dots, Z_rf$ provide the same type of group, or not, can now
obviously be expressed also as follows: Are the two groups: $U_1f,
\dots, U_rf$ and: $\mathfrak{ U}_1f, \dots, \mathfrak{ U}_rf$ similar
to each other, or not?

The two groups about which it is at present question are transitive;
the question whether they are similar, or not similar, can be decided
on the basis of Theorem~76, p.~\pageref{Theorem-76-S-425}.

From this theorem, we see that the groups: $U_1f, \dots, U_rf$ and:
$\mathfrak{ U}_1f, \dots, \mathfrak{ U}_rf$ are similar to each other
if and only if it is possible to relate them to each other in a
holoedrically isomorphic way so that the following condition is
satisfied: The most general subgroup of the group: $U_1f, \dots, U_rf$
which leaves invariant an arbitrarily chosen, but determined, 
system of values: $u_1 = u_1^0$, \dots, $u_{ r-m} = 
u_{ r-m}^0$ in general position must correspond to the most
general subgroup of the group: $\mathfrak{ U}_1f, \dots, 
\mathfrak{ U}_r f$ which leaves invariant a certain system
of values: $\mathfrak{ u}_1 = \mathfrak{ u}_1^0$, \dots, 
$\mathfrak{ u}_{ r-m} = \mathfrak{ u}_{ r-m}^0$ in general
position. 

Now, the groups: $U_1f, \dots, U_rf$ and: $\mathfrak{ U}_1f, \dots, 
\mathfrak{ U}_rf$ are, because of their derivation, 
both related to the group: $Y_1f, \dots, Y_rf$ in a holoedrically
isomorphic way, hence we can express the necessary and
sufficient criterion for their similarity also obviously as follows:
\emphasis{Similarity happens to hold when and only when
it is possible to relate the group: $Y_1f, \dots, Y_rf$ to 
itself in a holoedrically isomorphic way so that
its largest subgroup $G$ which fixes the manifold:
\def\theequation{10}\begin{equation}
u_1(y_1,\dots,y_r)
=
u_1^0,
\,\,\,\dots,\,\,\,
u_{r-m}(y_1,\dots,y_r)
=
u_{r-m}^0
\end{equation}
corresponds to its largest subgroup $\mathfrak{ G}$ which fixes the
manifold:}
\def\theequation{11}\begin{equation}
\mathfrak{u}_1(y_1,\dots,y_r)
=
\mathfrak{u}_1^0,
\,\,\,\dots,\,\,\,
\mathfrak{u}_{r-m}(y_1,\dots,y_r)
=
\mathfrak{u}_{r-m}^0.
\end{equation}

The criterion found with these words can be brought to another
remarkable form. In fact, we consider the group: $Z_1f, \dots, Z_rf$
reciprocal to the group: $Y_1f, \dots, Y_rf$. In it, ${\sf Z}_1f,
\dots, {\sf Z}_mf$ is the largest subgroup which leaves invariant the
manifold~\thetag{ 10}, and $\mathfrak{ Z}_1f, \dots, \mathfrak{ Z}_mf$
is the largest subgroup which leaves invariant the manifold~\thetag{
11}. According to Chap.~\ref{kapitel-20}, Proposition~8,
p.~\pageref{Satz-8-S-390}, one can then relate the two groups: $Y_1f,
\dots, Y_rf$ and: $Z_1f, \dots, Z_rf$ in a holoedrically isomorphic
way so that the subgroup $G$ corresponds to the subgroup: ${\sf Z}_1f,
\dots, {\sf Z}_mf$; but one can also related them together in a
holoedrically isomorphic way so that the subgroup $\mathfrak{ G}$
corresponds to the subgroup: $\mathfrak{ Z}_1f, \dots, \mathfrak{
Z}_mf$.

From this, it results that the group: $Y_1f, \dots, Y_rf$ can be
related to itself in a holoedrically isomorphic way as described just
above if and only if it is possible to relate the group: $Z_1f, \dots,
Z_rf$ to itself in a holoedrically isomorphic way so that the
subgroup: ${\sf Z}_1f, \dots, {\sf Z}_mf$ corresponds to the subgroup:
$\mathfrak{ Z}_1f, \dots, \mathfrak{ Z}_mf$.

We therefore have the:

\def\thetheorem{80}\begin{theorem}
\label{Theorem-80-S-445}
If the two $r$-term groups:
\[
Y_kf
=
\sum_{j=1}^r\,
\eta_{kj}(y_1,\dots,y_r)\,
\frac{\partial f}{\partial y_j}
\ \ \ \ \ \ \ \ \ \ \ \ \ {\scriptstyle{(k\,=\,1\,\cdots\,r)}}
\]
and:
\[
Z_kf
=
\sum_{j=1}^r\,\zeta_{kj}(y_1,\dots,y_r)\,
\frac{\partial f}{\partial y_j}
\ \ \ \ \ \ \ \ \ \ \ \ \ {\scriptstyle{(k\,=\,1\,\cdots\,r)}}
\]
are simply transitive and reciprocal to each other, 
if moreover:
\[
{\sf Z}_\mu f
=
\sum_{k=1}^r\,\varepsilon_{\mu k}\,Z_kf
\ \ \ \ \ \ \ \ \ \ \ \ \ {\scriptstyle{(\mu\,=\,1\,\cdots\,m)}}
\]
and:
\[
\mathfrak{Z}_\mu f
=
\sum_{k=1}^r\,\mathfrak{e}_{\mu k}\,
Z_kf
\ \ \ \ \ \ \ \ \ \ \ \ \ {\scriptstyle{(\mu\,=\,1\,\cdots\,m)}}
\]
are two $m$-term subgroups of the group: $Z_1f, \dots, Z_rf$ which
both neither are invariant in this group, nor contain, aside from the
identity transformation, another subgroup invariant in this group, and
lastly, if: $u_1 ( y), \dots, u_{ r-m} (y)$ and: $\mathfrak{ u}_1 (
y), \dots, \mathfrak{ u}_{ r-m} (y)$ are independent invariants of
these two $m$-term subgroups, respectively, then the two transitive,
both equally composed with the group $Y_1f, \dots, Y_rf$, groups:
\[
U_kf
=
\sum_{\nu=1}^{r-m}\,Y_k\,u_\nu\,
\frac{\partial f}{\partial u_\nu}
=
\sum_{\nu=1}^{r-m}\,
\omega_{k\nu}(u_1,\dots,u_{r-m})\,
\frac{\partial f}{\partial u_\nu}
\ \ \ \ \ \ \ \ \ \ \ \ \ {\scriptstyle{(k\,=\,1\,\cdots\,r)}}
\]
and:
\[
\mathfrak{U}_kf
=
\sum_{\nu=1}^{r-m}\,Y_k\,
\mathfrak{u}_\nu\,
\frac{\partial f}{\partial\mathfrak{u}_\nu}
=
\sum_{\nu=1}^{r-m}\,
\mathfrak{o}_{k\nu}
(\mathfrak{u}_1,\dots,\mathfrak{u}_{r-m})\,
\frac{\partial f}{\partial\mathfrak{u}_\nu}
\ \ \ \ \ \ \ \ \ \ \ \ \ {\scriptstyle{(k\,=\,1\,\cdots\,r)}}
\]
are similar to each other if and only if it is possible to relate the
group: $Z_1f, \dots, Z_rf$ to itself in a holoedrically isomorphic way
so that the subgroup: ${\sf Z}_1f, \dots, {\sf Z}_mf$ corresponds to
the subgroup: $\mathfrak{ Z}_1f, \dots, \mathfrak{ Z}_mf$.
\end{theorem}

Thanks to this theorem, the last part of the program stated on
p.~\pageref{S-435} is now also settled. At present, we can decide
whether two different subgroups of the group: $Z_1f, \dots, Z_rf$
provide, or do not provide, different types of transitive groups
equally composed with the group: $x_i' = f_i ( x, a)$. Clearly, for
that, only a research about the subgroups of the group $Z_1f, \dots,
Z_rf$ is required, or, what is the same, about the subgroups of the
group: $x_i' = f_i ( x, a)$.

We recapitulate the necessary operations in a theorem:

\def\thetheorem{81}\begin{theorem}
If the finite equations, or the infinitesimal transformations of an
$r$-term group $\Gamma$ are presented, one can determine in the
following way how many different types of transitive groups having the
same composition as $\Gamma$ there are: One determines all $m$-term
subgroups of $\Gamma$, but one excludes those which either are
invariant in $\Gamma$ or do contain a subgroup invariant in $\Gamma$
different from the identity transformation. One distributes the found
$m$-term subgroups in classes by reckoning that two subgroups always
belong to the same class when it is possible to relate $\Gamma$ to
itself in a holoedrically isomorphic way so that the two subgroups
correspond to each other. To each class of $m$-term subgroups obtained
in this way there corresponds a completely determined type of
transitive groups in $r - m$ variables equally composed with $\Gamma$;
to different classes there correspond different types. If one
undertakes this study for each one of the numbers: $m = 0, 1, 2,
\dots, r-1$, then one can have a view \deutsch{übersehen} of all
existing types.
\end{theorem}

\label{S-447-sq} 
We yet observe what follows: The operations required in the Theorem~81
are all executable, even when the group $\Gamma$ is not given, and
when only its composition is given. Also in this case, only
executable operations are then necessary. Since the number of the
subgroups of $\Gamma$ only depends upon arbitrary parameters, the
number of existing types depends at most upon arbitrary parameters.
In particular, there is only one single type of simply transitive
groups which are equally composed with the group $\Gamma$. But this
already results from the developments of the previous paragraph.

By combining the two Theorems~81 and~78, we yet obtain the following:

\renewcommand{\thefootnote}{\fnsymbol{footnote}}
\def\thetheorem{82}\begin{theorem}
If the finite equations: 
\[
x_i'
=
f_i(x_1,\dots,x_n,\,a_1,\dots,a_r)
\ \ \ \ \ \ \ \ \ \ \ \ \ {\scriptstyle{(i\,=\,1\,\cdots\,n)}}
\]
of an $r$-term group are presented, then the determination of all
equally composed transitive groups requires in all cases, while
disregarding executable operations, only the integration of
simultaneous systems of ordinary differential
equations.\footnote[1]{\,
Cf. \name{Lie}, Math. Annalen Vol. XVI, p.~528.
} 
\end{theorem}
\renewcommand{\thefootnote}{\arabic{footnote}}

If one wants to list all $r$-term transitive groups in $n$ variables,
then one determines all compositions of $r$-term groups and one then
seeks, for each composition, the associated types of transitive groups
in $n$ variables. \label{S-448}
For given $r$ (and $n$), all these types decompose
in a bounded number of kinds \deutsch{Gattung} so that the types of a
kind have the same analytic representation. However, the analytic
expressions for all types of a kind contain certain parameters the
number of which we can always imagine to be lowered down to a minimum.
That such parameters occur stems from two different facts: firstly,
from the fact that, to a given $r > 2$, an unbounded number of
different compositions is associated; secondly from the fact that in
general, a given $r$-term group contains an unbounded number of
subgroups which produce nothing but different types. We do not
consider it to be appropriate here to pursue these considerations
further.

\sectionengellie{\S\,\,\,108.}

Once again, we move back to the standpoint we took on
p.~\pageref{S-443-sq} sq.

\label{S-448-sq}
We had used two $m$-term subgroups: ${\sf Z}_1f, \dots, {\sf Z}_mf$
and: $\mathfrak{ Z}_1f, \dots, \mathfrak{ Z}_mf$ of the group: $Z_1f,
\dots, Z_rf$ in order to produce transitive groups equally composed
with the group: $Y_1f, \dots, Y_rf$, and we have found the two groups:
$U_1f, \dots, U_rf$ and: $\mathfrak{ U}_1f, \dots, \mathfrak{ U}_rf$.
At present, the question arises: under which conditions are
these two group similar to each other?

We have answered this question at that time, when we based 
ourselves on the Theorem~76, p.~\pageref{Theorem-76-S-425}, 
and in this way, we have obtained the Theorem~80, 
p.~\pageref{Theorem-80-S-445}.
At present, we want to take up again the question and
to attempt to answer it without using the Theorem~76.

\medskip

Evidently, we are close to presume that the groups: $U_1f, \dots, U_rf$
and: $\mathfrak{ U}_1f, \dots, \mathfrak{ U}_rf$ are in any
case similar
to each other when there is a transformation:
\def\theequation{12}\begin{equation}
\overline{y}_k
=
\Omega_k(y_1,\dots,y_r)
\ \ \ \ \ \ \ \ \ \ \ \ \ {\scriptstyle{(k\,=\,1\,\cdots\,r)}}
\end{equation}
which converts the subgroup: ${\sf Z}_1f, \dots, {\sf Z}_mf$
into the subgroup: $\mathfrak{ Z}_1f, \dots, \mathfrak{ Z}_mf$
and which transfers at the same time the group: $Z_1f, \dots, Z_rf$
into itself. We will show that this presumption
corresponds to the truth.

Thus, let~\thetag{ 12} be a transformation which possesses
the indicated properties. Through this transformation,
the invariants of the group: ${\sf Z}_1f, \dots, {\sf Z}_mf$
are obviously transferred to the invariants of the group:
$\mathfrak{ Z}_1f, \dots, \mathfrak{ Z}_mf$, hence
we have by virtue of~\thetag{ 12}:
\[
\mathfrak{u}_\nu
(\overline{y}_1,\dots,\overline{y}_r)
=
\chi_\nu
\big(
u_1(y),\,\dots,\,u_{r-m}(y)
\big)
\ \ \ \ \ \ \ \ \ \ \ \ \ 
{\scriptstyle{(\nu\,=\,1\,\cdots\,r\,-\,m)}},
\]
where the functions: $\chi_1, \dots, \chi_{ r-m}$ are
absolutely determined and are mutually
independent relatively to: $u_1 (y), \dots, 
u_{ r-m} (y)$. 
On the other hand, according to
Chap.~\ref{kapitel-20}, p.~\pageref{S-381}, 
through the transformation~\thetag{ 12}, 
not only the group: $Z_1f, \dots, Z_rf$
is transferred into itself, but also
its reciprocal group: $Y_1f, \dots, Y_rf$, so we have:
\[
\sum_{j=1}^r\,\eta_{kj}
(\overline{y}_1,\dots,\overline{y}_r)\,
\frac{\partial f}{\partial\overline{y}_j}
=
\overline{Y}_kf
=
\sum_{j=1}^r\,h_{kj}\,Y_jf
\ \ \ \ \ \ \ \ \ \ \ \ \ {\scriptstyle{(k\,=\,1\,\cdots\,r)}},
\]
where the $h_{ kj}$ are constants, the determinant of which does not
vanish.

Now, if we set for $f$, in the equations just written, an arbitrary
function $F$ of: $u_1 (y), \dots, u_{ r-m} (y)$, or, what is the same,
a function of: $\mathfrak{ u}_1 ( \overline{ y})$, \dots, $\mathfrak{
u}_{ r-m} ( \overline{ y})$, it then comes:
\[
\overline{Y}_kF
=
\sum_{\nu=1}^{r-m}\,
\mathfrak{o}_{k\nu}
(\mathfrak{u}_1,\dots,\mathfrak{u}_{r-m})\,
\frac{\partial F}{\partial\mathfrak{u}_\nu}
=
\sum_{j=1}^r\,h_{kj}\,
\sum_{\nu=1}^{r-m}\,
\omega_{j\nu}
(u_1,\dots,u_{r-m})\,
\frac{\partial F}{\partial u_\nu}.
\]
In other words, the two groups: $U_1f, \dots, U_rf$ and: $\mathfrak{
U}_1f, \dots, \mathfrak{ U}_rf$ are similar to each other: obviously:
\def\theequation{13}\begin{equation}
\mathfrak{u}_\nu
=
\chi_\nu(u_1,\dots,u_{r-m})
\ \ \ \ \ \ \ \ \ \ \ \ \ 
{\scriptstyle{(\nu\,=\,1\,\cdots\,r\,-\,m)}}
\end{equation}
is a transformation which transfers the one group
to the other. 

As a result, the presumption enunciated above is proved.

\medskip

But it can be proved that the converse also holds true: When the two
groups: $U_1f, \dots, U_rf$ and: $\mathfrak{ U}_1f, \dots, \mathfrak{
U}_rf$ are similar to each other, then there always is a
transformation which transfers the subgroup: ${\sf Z}_1f, \dots, {\sf
Z}_mf$ to the subgroup: $\mathfrak{ Z}_1f, \dots,
\mathfrak{ Z}_mf$ and which transfers at the same
time the group: $Z_1f, \dots, Z_rf$ to itself.

Thus, let the two groups: $U_kf$ and: $\mathfrak{ U}_kf$ be similar to
each other and let:
\def\theequation{14}\begin{equation}
\mathfrak{u}_\nu
=
\psi_\nu(u_1,\dots,u_{r-m})
\ \ \ \ \ \ \ \ \ \ \ \ \ 
{\scriptstyle{(\nu\,=\,1\,\cdots\,r\,-\,m)}}
\end{equation}
be a transformation which transfers the one group to the other, 
so that we have:
\def\theequation{15}\begin{equation}
\mathfrak{U}_kf
=
\sum_{j=1}^r\,\delta_{kj}\,U_jf
=
U_k'f
\ \ \ \ \ \ \ \ \ \ \ \ \ {\scriptstyle{(k\,=\,1\,\cdots\,r)}}.
\end{equation}
Here, by $\delta_{ kj}$, it is to be understood constants the
determinant of which does not vanish. If $\mathfrak{ U}_1f, \dots,
\mathfrak{ U}_rf$ are linked together by relations of the form:
\[
\leftbracket
\mathfrak{U}_i,\,\mathfrak{U}_k
\rightbracket
=
\sum_{s=1}^r\,c_{iks}\,\mathfrak{U}_sf,
\]
then naturally, $U_1'f, \dots, U_r'f$ are linked together by 
the relations:
\[
\leftbracket
U_i',\,U_k'
\rightbracket
=
\sum_{s=1}^r\,c_{iks}\,U_s'f.
\]

Above, we have seen that every transformation~\thetag{ 12}
which leaves invariant the group: $Z_1f, \dots, Z_rf$ 
and which transfers the subgroup of the ${\sf Z}_\mu f$
to the subgroup of the $\mathfrak{ Z}_\mu f$ provides a
completely determined transformation~\thetag{ 13}
which transfers the group of the $U_kf$ to the group
of the $\mathfrak{ U}_kf$. 
Under the present assumptions, we already know a transformation
which accomplishes the latter transfer, namely the 
transformation~\thetag{ 14}. Therefore, we
attempt to determine a transformation:
\def\theequation{16}\begin{equation}
\overline{y}_k
=
O_k(y_1,\dots,y_r)
\ \ \ \ \ \ \ \ \ \ \ \ \ {\scriptstyle{(k\,=\,1\,\cdots\,r)}}
\end{equation}
which leaves invariant the group: $Z_1f, \dots, Z_rf$, 
which converts the subgroup: ${\sf Z}_1f, \dots, {\sf Z}_mf$
into the subgroup: $\mathfrak{ Z}_1f, \dots,
\mathfrak{ Z}_mf$ and lastly, which provides exactly the
transformation~\thetag{ 14}.

It is clear that every transformation~\thetag{ 16}
of the kind demanded must be constituted in such a way 
that its equations embrace the $r-m$ mutually independent
equations:
\def\theequation{14'}\begin{equation}
\mathfrak{u}_\nu
(\overline{y}_1,\dots,\overline{y}_r)
=
\psi_\nu
\big(
u_1(y),\,\dots,\,u_{r-m}(y)
\big)
\ \ \ \ \ \ \ \ \ \ \ \ \ 
{\scriptstyle{(\nu\,=\,1\,\cdots\,r\,-\,m)}}.
\end{equation}
If it satisfies this condition, and in addition, if it yet leaves
invariant the group: $Z_1f, \dots, Z_rf$, then it satisfies all our
demands. Indeed, on the first hand, it transfers the invariants of
the group: ${\sf Z}_1f, \dots, {\sf Z}_mf$ to the invariants of the
group: $\mathfrak{ Z}_1f, \dots, \mathfrak{ Z}_mf$, whence it converts
the first one of these two groups into the second one, and on the
other hand, it visibly produces the transformation~\thetag{ 14} by
virtue of which the two groups: $U_1f, \dots, U_rf$ and $\mathfrak{
U}_1f, \dots, \mathfrak{ U}_rf$ are similar to each other.

But now, whether we require of the transformation~\thetag{ 16} that it
leaves invariant the group: $Z_1f, \dots, Z_rf$, or whether we require
that it transfers the group: $Y_1f, \dots, Y_rf$ into itself, this
obviously is completely indifferent. We can therefore
interpret our problem as follows:

\plainstatement{To seek a transformation~\thetag{ 16}
which leaves invariant the group: $Y_1f, \dots, Y_rf$ and which is
constituted in such a way that its equations embrace the
equations~\thetag{ 14'}.}

\medskip

For the sake of simplification, we introduce new variables.

We choose $m$ arbitrary mutually independent functions: $v_1 ( y),
\dots, v_m ( y)$ that are also independent of: $u_1 ( y), \dots, u_{
r-m} (y)$, and moreover, we choose $m$ arbitrary mutually independent
functions: $\mathfrak{ v}_1 ( \overline{ y}), \dots,
\mathfrak{ v}_m ( \overline{ y})$
that are also independent of $\mathfrak{ u}_1 ( \overline{ y}), \dots,
\mathfrak{ u}_{ r-m} ( \overline{ y})$.
We introduce the functions: $u_1 ( y), \dots, u_{ r-m} (y)$, 
$v_1 ( y), \dots, v_m ( y)$ as new variables in place of: 
$y_1, \dots, y_r$ and the functions: 
$\mathfrak{ u}_1 ( \overline{ y}), \dots, 
\mathfrak{ u}_{ r- m} ( \overline{ y})$, 
$\mathfrak{ v}_1 ( \overline{ y}), \dots, 
\mathfrak{ v}_m ( \overline{ y})$ in place of:
$\overline{ y}_1, \dots, \overline{ y}_r$. 

In the new variables, the sought transformation~\thetag{ 16} 
necessarily receives the form:
\def\theequation{16'}\begin{equation}
\left\{
\aligned
\mathfrak{u}_\nu
&
=
\psi_\nu(u_1,\dots,u_{r-m})
\ \ \ \ \ \ \ \ \ \ \ \ \ 
{\scriptstyle{(\nu\,=\,1\,\cdots\,r\,-\,m)}}
\\
\mathfrak{v}_\mu
&
=
\Psi_\mu(u_1,\dots,u_{r-m},\,v_1,\dots,v_m)
\ \ \ \ \ \ \ \ \ \ \ \ \ 
{\scriptstyle{(\mu\,=\,1\,\cdots\,m)}},
\endaligned\right.
\end{equation}
where it is already taken account of the fact that the
present equations must comprise the equations~\thetag{ 14}.

But what do we have in place of the requirement that
the transformation~\thetag{ 16} should leave invariant
the group: $Y_1f, \dots, Y_rf$?

Clearly, in the new variables: $u_1, \dots, u_{ r-m}$,
$v_1, \dots, v_m$, the group: $Y_1f, \dots, Y_rf$ receives
the form:
\[
U_kf
+
\sum_{\mu=1}^m\,
w_{k\mu}
(u_1,\dots,u_{r-m},\,v_1,\dots,v_m)\,
\frac{\partial f}{\partial v_\mu}
=
U_kf
+
V_kf
\ \ \ \ \ \ \ \ \ \ \ \ \ {\scriptstyle{(k\,=\,1\,\cdots\,r)}}.
\]
On the other hand, after the introduction of
$\mathfrak{ u}_1, \dots, \mathfrak{ u}_{ r-m}$, 
$\mathfrak{ v}_1, \dots, \mathfrak{ v}_m$, the
infinitesimal transformations:
\[
\overline{Y}_kf
=
\sum_{j=1}^r\,\eta_{kj}
(\overline{y}_1,\dots,\overline{y}_r)\,
\frac{\partial f}{\partial\overline{y}_j}
\ \ \ \ \ \ \ \ \ \ \ \ \ {\scriptstyle{(k\,=\,1\,\cdots\,r)}}
\]
are transferred to:
\[
\mathfrak{U}_kf
+
\sum_{\mu=1}^m\,
\mathfrak{w}_{k\mu}
(\mathfrak{u}_1,\dots,\mathfrak{u}_{r-m},\,
\mathfrak{v}_1,\dots,\mathfrak{v}_m)\,
\frac{\partial f}{\partial\mathfrak{v}_\mu}
=
\mathfrak{U}_kf
+
\mathfrak{V}_kf
\ \ \ \ \ \ \ \ \ \ \ \ \ {\scriptstyle{(k\,=\,1\,\cdots\,r)}}.
\]
Consequently, we must require of the transformation~\thetag{ 16'}
that it transfers the group: 
$U_1f + V_1f$, \dots, $U_r f + V_rf$ to the group:
$\mathfrak{ U}_1f + \mathfrak{ V}_1f$, \dots, 
$\mathfrak{ U}_rf + \mathfrak{ V}_rf$; 
we must attempt to determine the functions:
$\Psi_1, \dots, \Psi_m$ accordingly. 

As the equation~\thetag{ 15} show, the $r$ independent infinitesimal
transformations:
\[
\sum_{j=1}^r\,\delta_{1j}\,U_jf,
\,\,\,\dots,\,\,\,
\sum_{j=1}^r\,\delta_{rj}\,U_jf
\]
are transferred by the transformation:
\def\theequation{14}\begin{equation}
\mathfrak{u}_\nu
=
\psi_\nu(u_1,\dots,u_{r-m})
\ \ \ \ \ \ \ \ \ \ \ \ \ 
{\scriptstyle{(\nu\,=\,1\,\cdots\,r\,-\,m)}}
\end{equation}
to the transformations:
\[
\mathfrak{U}_1f,
\,\dots,\,
\mathfrak{U}_rf,
\]
respectively. Hence, if a transformation of the form~\thetag{ 16'}
is supposed to convert the group of the $U_kf + V_kf$
to the group of the $\mathfrak{ U}_kf + \mathfrak{ V}_kf$, 
then through it, the $r$ independent infinitesimal transformations:
\[
\sum_{j=1}^r\,\delta_{1j}\,(U_jf+V_jf),
\,\,\,\dots,\,\,\,
\sum_{j=1}^r\,\delta_{rj}\,(U_jf+V_jf)
\]
are transferred to:
\[
\mathfrak{U}_1f
+
\mathfrak{V}_1f,
\,\,\,\dots,\,\,\,
\mathfrak{U}_rf
+
\mathfrak{V}_rf.
\]
This condition is necessary, and simultaneously also
sufficient.

\medskip

Thanks to the same considerations as in Chap.~\ref{kapitel-19}, 
p.~\pageref{S-335-sq} sq., we now recognize that every 
transformation~\thetag{ 16'} having the constitution demanded
can also be defined as a system of equations
in the $2\, r$ variables $u$, $v$, 
$\mathfrak{ u}$, $\mathfrak{ v}$ which possesses the form:
\def\theequation{16'}\begin{equation}
\left\{
\aligned
\mathfrak{u}_\nu
&
=
\psi_\nu(u_1,\dots,u_{r-m})
\ \ \ \ \ \ \ \ \ \ \ \ \ 
{\scriptstyle{(\nu\,=\,1\,\cdots\,r\,-\,m)}}
\\
\mathfrak{v}_\mu
&
=
\Psi_\mu(u_1,\dots,u_{r-m},\,v_1,\dots,v_m)
\ \ \ \ \ \ \ \ \ \ \ \ \ 
{\scriptstyle{(\mu\,=\,1\,\cdots\,m)}},
\endaligned\right.
\end{equation}
which admits the $r$-term group:
\[
W_kf
=
\mathfrak{U}_kf
+
\mathfrak{V}_kf
+
\sum_{j=1}^r\,\delta_{kj}\,
(U_jf+V_jf)
\ \ \ \ \ \ \ \ \ \ \ \ \ {\scriptstyle{(k\,=\,1\,\cdots\,r)}}
\]
and which is solvable with respect to $u_1, \dots, u_{ r-m}$, 
$v_1, \dots, v_m$.

According to our assumption, the system of equations:
\def\theequation{14}\begin{equation}
\mathfrak{u}_\nu
=
\psi_\nu(u_1,\dots,u_{r-m})
\ \ \ \ \ \ \ \ \ \ \ \ \ 
{\scriptstyle{(\nu\,=\,1\,\cdots\,r\,-\,m)}}
\end{equation}
represents a transformation which transfers the $r$
independent infinitesimal transformations:
\[
\sum_{j=1}^r\,\delta_{kj}\,U_jf
\ \ \ \ \ \ \ \ \ \ \ \ \ {\scriptstyle{(k\,=\,1\,\cdots\,r)}}
\]
to: $\mathfrak{ U}_1f, \dots, \mathfrak{ U}_rf$, 
respectively, whence it admits the $r$-term group:
\[
\mathfrak{U}_kf
+
\sum_{j=1}^r\,\delta_{kj}\,U_jf
\ \ \ \ \ \ \ \ \ \ \ \ \ {\scriptstyle{(k\,=\,1\,\cdots\,r)}},
\]
and consequently, also the group: $W_1f, \dots, W_rf$.

We can therefore use the developments of Chap.~\ref{kapitel-14}, 
pp.~\pageref{Theorem-40-S-233}--\pageref{S-236} in order to 
find a system
of equations~\thetag{ 16'} having the constitution demanded.

To begin with, from $W_1f, \dots, W_rf$, we form certain reduced
infinitesimal transformations: $\mathfrak{ W}_1f, \dots, 
\mathfrak{ W}_rf$ by leaving out all terms 
with the differential quotients: $\partial f / \partial 
\mathfrak{ u}_1$, \dots, $\partial f / \partial 
\mathfrak{ u}_{ r-m}$ and by making the substitution: 
$\mathfrak{ u}_1 = \psi_1 ( u)$, \dots, 
$\mathfrak{ u}_{ r-m} = \psi_{ r-m} (u)$ in all
the terms remaining.
If this substitution is indicated by the sign: 
$[ \,\,\, ]$, then the $\mathfrak{ W}_kf$ read as follows:
\[
\mathfrak{W}_kf
=
\sum_{\mu=1}^m\,
\big[
\mathfrak{w}_{k\mu}
(\mathfrak{u}_1,\dots,\mathfrak{u}_{r-m},\,
\mathfrak{v}_1,\dots,\mathfrak{v}_m
\big]\,
\frac{\partial f}{\partial\mathfrak{v}_\mu}
+
\sum_{j=1}^r\,
\delta_{kj}\,
(U_jf+V_jf),
\]
or, more shortly:
\[
\mathfrak{W}_kf
=
\big[
\mathfrak{V}_kf
\big]
+
\sum_{j=1}^r\,
\delta_{kj}\,
(U_jf+V_jf)
\ \ \ \ \ \ \ \ \ \ \ \ \ {\scriptstyle{(k\,=\,1\,\cdots\,r)}}.
\]
Naturally, $\mathfrak{ W}_1f, \dots, \mathfrak{ W}_rf$ generate a
group in the $r + m$ variables $\mathfrak{ v}_1, \dots, \mathfrak{
v}_m$, $u_1, \dots, u_{ r-m}$, $v_1, \dots, v_m$ and in our case, a
group which is evidently $r$-term.

At present, we determine in the $\mathfrak{ v}$, $u$, $v$ a system of
equations of the form:
\def\theequation{17}\begin{equation}
\mathfrak{v}_\mu
=
\Psi_\mu
(u_1,\dots,u_{r-m},\,
v_1,\dots,v_m)
\ \ \ \ \ \ \ \ \ \ \ \ \ {\scriptstyle{(\mu\,=\,1\,\cdots\,m)}}
\end{equation}
which admits the group: $\mathfrak{ W}_1f, \dots, \mathfrak{ W}_rf$
and which is solvable with respect to $v_1, \dots, v_m$.
Lastly, when we add this system to the equations~\thetag{ 14}, 
we then obtain a system of equations of the form~\thetag{ 16'}
which possesses the properties explained above. 

Since the $r$-term group:
\[
\sum_{j=1}^r\,\delta_{kj}\,
(U_jf+V_jf)
\ \ \ \ \ \ \ \ \ \ \ \ \ {\scriptstyle{(k\,=\,1\,\cdots\,r)}}
\]
is simply transitive, then in the matrix which can be formed
with the coefficients of $\mathfrak{ W}_1f, \dots, \mathfrak{ W}_rf$, 
it is certain that not all $r \times r$ determinants vanish
identically, and they can even less vanish all by virtue of a
system of equations of the form~\thetag{ 17}.
Consequently, every system of equations of the form~\thetag{ 17}
which admits the group: $\mathfrak{ W}_1f, \dots, \mathfrak{ W}_rf$
can be represented by $m$ independent relations between
$m$ arbitrary independent solutions of the $r$-term
complete system:
\def\theequation{18}\begin{equation}
\mathfrak{W}_1f=0,
\,\,\,\dots,\,\,\,
\mathfrak{W}_rf=0
\end{equation}
in the $r + m$ variables $\mathfrak{ v}_1, \dots, 
\mathfrak{ v}_m$, $u_1, \dots, u_{ r-m}$, 
$v_1, \dots, v_m$. 

The equations~\thetag{ 18} are obviously solvable with respect to:
$\partial f / \partial u_1$, \dots, 
$\partial f / \partial u_{ r-m}$, 
$\partial f / \partial v_1$, \dots, 
$\partial f / \partial v_m$; on the other hand, they
are also solvable with respect to: 
$\partial f / \partial u_1$, \dots, 
$\partial f / \partial u_{ r-m}$, 
$\partial f / \partial \mathfrak{ v}_1$, \dots,
$\partial f / \partial \mathfrak{ v}_m$, because if we
introduce the new variables: $u_1, \dots, u_{ r-m}$, 
$\mathfrak{ v}_1, \dots, \mathfrak{ v}_m$ in place of
$\mathfrak{ u}_1, \dots, \mathfrak{ u}_{ r-m}$, 
$\mathfrak{ v}_1, \dots, \mathfrak{ v}_m$
by means of the equations:
\def\theequation{14}\begin{equation}
\mathfrak{u}_\mu
=
\psi_\mu(u_1,\dots,u_{r-m})
\ \ \ \ \ \ \ \ \ \ \ \ \ 
{\scriptstyle{(\mu\,=\,1\,\cdots\,r\,-\,m)}}
\end{equation}
in the $r$ infinitesimal transformations:
\[
\mathfrak{U}_1f
+
\mathfrak{V}_1f,
\,\,\,\dots,\,\,\,
\mathfrak{U}_rf
+
\mathfrak{V}_rf,
\]
then we obtain the infinitesimal transformations:
\[
\sum_{j=1}^r\,\delta_{kj}\,U_jf
+
\big[\mathfrak{V}_kf\big]
\ \ \ \ \ \ \ \ \ \ \ \ \ {\scriptstyle{(k\,=\,1\,\cdots\,r)}}
\]
which in turn generate therefore a simply transitive group in 
the variables: $u_1, \dots, u_{ r-m}$, 
$\mathfrak{ v}_1, \dots, \mathfrak{ v}_m$.

Consequently, if:
\[
P_\mu
\big(
\mathfrak{v}_1,\dots,\mathfrak{v}_m,\,
u_1,\dots,u_{r-m},\,
v_1,\dots,v_m
\big)
\ \ \ \ \ \ \ \ \ \ \ \ \ 
{\scriptstyle{(\mu\,=\,1\,\cdots\,m)}}
\]
are $m$ arbitrary independent solutions of the complete 
system~\thetag{ 18}, then these solutions are mutually
independent both relatively to $\mathfrak{ v}_1, 
\dots, \mathfrak{ v}_m$ and relatively to 
$v_1, \dots, v_m$ (cf. Theorem~12, p.~\pageref{Theorem-12-S-91}).

From this, it results that the most general system of equations
of the form~\thetag{ 17} which can be resolved with respect
to $v_1, \dots, v_m$ and which admits the group: 
$\mathfrak{ W}_1f, \dots, \mathfrak{ W}_rf$ can be obtained
by solving the $m$ equations:
\def\theequation{19}\begin{equation}
P_\mu
(\mathfrak{v}_1,\dots,\mathfrak{v}_m,\,
u_1,\dots,u_{r-m},\,
v_1,\dots,v_m)
=
{\rm const.}
\ \ \ \ \ \ \ \ \ \ \ \ \ {\scriptstyle{(\mu\,=\,1\,\cdots\,m)}}
\end{equation}
with respect to $\mathfrak{ v}_1, \dots, \mathfrak{ v}_m$.

\medskip

At present, we can immediately indicate a transformation,
and in fact the most general transformation~\thetag{ 16'}, which
transfers the infinitesimal transformations: 
\[
\sum_{j=1}^r\,\delta_{1j}\,
(U_j+V_jf),
\,\,\,\dots,\,\,\,
\sum_{j=1}^r\,\delta_{rj}\,
(U_jf+V_jf)
\]
to the infinitesimal transformations:
\[
\mathfrak{U}_1f
+
\mathfrak{V}_1f,
\,\,\dots,\,\,
\mathfrak{U}_rf
+
\mathfrak{V}_rf,
\]
respectively; this transformation is simply represented by the
equations~\thetag{ 14} and~\thetag{ 19} together. Lastly, if we
introduce again the variables: $y_1, \dots, y_r$, $\overline{ y}_1,
\dots, \overline{ y}_r$ in~\thetag{ 14} and in~\thetag{ 19}, we obtain
a transformation which leaves invariant the group: $Y_1f, \dots, Y_rf$
and whose equations embrace the equations~\thetag{ 14'}; in other
words, we obtain a transformation which leaves the group: $Z_1f,
\dots, Z_rf$ invariant and which transfers the subgroup: ${\sf Z}_1f,
\dots, {\sf Z}_mf$ to the subgroup: $\mathfrak{ Z}_1f, \dots,
\mathfrak{ Z}_mf$.

With these words, it is proved that there always exists a
transformations having the constitution just described, as soon as the
two groups: $U_1f, \dots, U_rf$ and $\mathfrak{ U}_1f, \dots,
\mathfrak{ U}_rf$ are similar to each other. But since the similarity
of the two groups follows from the existence of such a transformation,
according to p.~\pageref{S-448-sq} sq., we can say:

\plainstatement{
The two groups: $U_1f, \dots, U_rf$ and $\mathfrak{ U}_1f,
\dots, \mathfrak{ U}_rf$ are similar to each other
if and only if there is a transformation which leaves invariant the
group: $Z_1f, \dots, Z_rf$ and which transfers the subgroup: ${\sf
Z}_1f, \dots, {\sf Z}_mf$ to the subgroup: $\mathfrak{ Z}_1f, \dots,
\mathfrak{ Z}_mf$.}

Now, the group: $Z_1f, \dots, Z_rf$ is simply transitive, hence it is
clear that there always is a transformation of this sort when and only
when the group: $Z_1f, \dots, Z_rf$ can be related to itself in a
holoedrically isomorphic way so that the two subgroups: ${\sf Z}_1f,
\dots, {\sf Z}_mf$ and $\mathfrak{ Z}_1f, \dots, \mathfrak{ Z}_mf$
correspond to each other. Consequently, for the similarity of the
groups: $U_1f, \dots, U_rf$ and $\mathfrak{ U}_1f,
\dots, \mathfrak{ U}_rf$, we find exactly the same criterion
as we had expressed in Theorem~80, p.~\pageref{Theorem-80-S-445}.

Besides, the preceding developments can also be used to derive a new
proof of the Theorem~76 in Chap.~\ref{kapitel-21},
p.~\pageref{Theorem-76-S-425}.

At first, thanks to considerations completely similar to the ones on
p.~\pageref{S-444-445}, it can be proved that the transitive groups:
$U_1f, \dots, U_rf$ and $\mathfrak{ U}_1f, \dots, \mathfrak{ U}_rf$
can be related to each other in a the holoedrically isomorphic way
described in Theorem~76 if and only if it is possible to relate the
group: $Z_1f, \dots, Z_rf$ to itself in a holoedrically isomorphic way
in such a way that the subgroups: ${\sf Z}_1f, \dots, {\sf Z}_mf$ and
$\mathfrak{ Z}_1f, \dots, \mathfrak{ Z}_mf$ correspond to each other.
Afterwards, it follows from the preceding developments
that the conditions of the Theorem~76 for the similarity
of the two groups: $U_1f, \dots, U_rf$ and
$\mathfrak{ U}_1f, \dots, \mathfrak{ U}_rf$
are necessary and sufficient. 

\sectionengellie{\S\,\,\,109.}

Now, we turn to the second one of the two problems, the settlement of
which was announced in the introduction of the chapter (on
p.~\pageref{S-430}): \emphasis{to the determination of all $r$-term
intransitive groups}; as was already said at that time, we want to
undertake the determination in question under the assumption that all
transitive groups with $r$ or less parameters are given. Since all
transitive groups with an equal number of parameters can be ordered in
classes according to their composition and moreover, since all
transitive groups having one and the same composition decompose in a
series of types (cf. p.~\pageref{S-434-sq} sq.), we can precise our
assumption somehow more exactly by supposing \emphasis{firstly} that
all possible compositions of a group with $r$ or less parameters are
known and \emphasis{secondly} by imagining that for each one of these
compositions, all possible types of transitive groups having the
concerned composition are given.

To begin with, we consider an arbitrary $r$-term intransitive group.

If $X_1f, \dots, X_rf$ are independent infinitesimal 
transformations of an $r$-term group of the space 
$x_1, \dots, x_n$, the $r$ equations:
\def\theequation{20}\begin{equation}
X_1f=0,
\,\,\,\dots,\,\,\,
X_rf=0
\end{equation}
have a certain number, say exactly $n - l > 0$, of
independent solutions in common.
Hence we can imagine that the variables $x_1, \dots, x_n$ are
chosen from the beginning in such a way that
$x_{ l+1}, \dots, x_n$ are such independent solutions. 
Then $X_1f, \dots, X_rf$ will receive the form:
\[
X_kf
=
\sum_{\lambda=1}^l\,
\xi_{k\lambda}(x_1,\dots,x_l,\,
x_{l+1},\dots,x_n)\,
\frac{\partial f}{\partial x_\lambda}
\ \ \ \ \ \ \ \ \ \ \ \ \ {\scriptstyle{(k\,=\,1\,\cdots\,r)}}
\]
where now naturally, not all $l \times l$ determinants of the
matrix:
\[
\left\vert
\begin{array}{cccc}
\xi_{11}(x) & \,\cdot\, & \,\cdot\, & \xi_{1l}(x)
\\
\cdot & \,\cdot\, & \,\cdot\, & \cdot
\\
\cdot & \,\cdot\, & \,\cdot\, & \cdot
\\
\xi_{r1}(x) & \,\cdot\, & \,\cdot\, & \xi_{rl}(x)
\end{array}
\right\vert
\]
vanish identically, since otherwise, the equations~\thetag{ 20}
would have more than $n - l$ independent solutions in common.

If the number $r$, which is at least equal to $l$, would
be exactly equal to $l$, then $X_1f, \dots, X_rf$ would
be linked together by no relation of the form:
\[
\chi_1(x_{l+1},\dots,x_n)\,X_1f
+\cdots+
\chi_r(x_{l+1},\dots,x_n)\,X_rf
=
0\,;
\]
but now, $r$ needs not be equal to $l$, hence relations of the form
just described can also very well exist without that all the functions
$\chi_1, \dots, \chi_r$ vanish. So we want to assume that $X_1f,
\dots, X_mf$, say, are linked together by no such relation, while by
contrast, $X_{ m+1}f, \dots, X_rf$ may be expressed as follows in
terms of $X_1f, \dots, X_mf$:
\def\theequation{21}\begin{equation}
X_{m+\nu}f
\equiv
\sum_{\mu=1}^m\,
\vartheta_{\nu\mu}(x_{l+1},\dots,x_n)\,
X_\mu f
\ \ \ \ \ \ \ \ \ \ \ \ \ 
{\scriptstyle{(\nu\,=\,1\,\cdots\,r\,-\,m)}}.
\end{equation}
Here of course, $m$ satisfies the inequations: $l \leqslant m 
\leqslant r$. 

By combination of the equation~\thetag{ 21} with the relations:
\[
\leftbracket
X_i,\,X_k
\rightbracket
=
\sum_{s=1}^r\,c_{iks}\,X_sf
\ \ \ \ \ \ \ \ \ \ \ \ \ 
{\scriptstyle{(i,\,\,k\,=\,1\,\cdots\,r)}}
\]
which hold true in all circumstances, we yet recognize that
$X_1f, \dots, X_mf$ stand pairwise in the relationships:
\[
\aligned
\leftbracket
X_\lambda,\,X_\mu
\rightbracket
=
\sum_{\pi=1}^m\,
\bigg\{
c_{\lambda\mu\pi}
&
+
\sum_{\nu=1}^{r-m}\,
c_{\lambda,\,\mu,\,m+\nu}\,
\vartheta_{\nu\pi}
(x_{l+1},\dots,x_n)
\bigg\}\,X_\pi f
\\
& \ \ \ \ \ \ \ 
{\scriptstyle{(\lambda,\,\,\mu\,=\,1\,\cdots\,m)}}.
\endaligned
\]

\smallskip

Now, if the variables $x_{ l+1}, \dots, x_n$ are replaced by arbitrary
constants: $a_{ l+1}, \dots, a_n$ and if 
$x_1, \dots, x_l$ only are still considered as variables, then 
it is clear that the $r$ infinitesimal transformations:
\[
\overline{X}_kf
=
\sum_{\lambda=1}^l\,
\xi_{k\lambda}
(x_1,\dots,x_l,\,a_{l+1},\dots,a_n)\,
\frac{\partial f}{\partial x_\lambda}
\ \ \ \ \ \ \ \ \ \ \ \ \ {\scriptstyle{(k\,=\,1\,\cdots\,r)}}
\]
in the $l$ variables $x_1, \dots, x_l$ are not anymore independent
of each other, but can be linearly deduced from the $m$
independent infinitesimal transformations:
\[
\overline{X}_\mu f
=
\sum_{\lambda=1}^l\,
\xi_{\mu\lambda}
(x_1,\dots,x_l,\,
a_{l+1},\dots,a_n)\,
\frac{\partial f}{\partial x_\lambda}
\ \ \ \ \ \ \ \ \ \ \ \ \ {\scriptstyle{(\mu\,=\,1\,\cdots\,m)}}.
\]
The $m$ infinitesimal transformations $\overline{ X}_1f, \dots, 
\overline{ X}_mf$ are in turn obviously linked together by the
relations:
\def\theequation{22}\begin{equation}
\leftbracket
\overline{X}_\lambda,\,
\overline{X}_\mu
\rightbracket
=
\sum_{\pi=1}^m\,
\bigg\{
c_{\lambda\mu\pi}
+
\sum_{\nu=1}^{r-m}\,
c_{\lambda,\,\mu,\,m+\nu}\,
\vartheta_{\nu\pi}(a_{l+1},\dots,a_n)
\bigg\}\,
\overline{X}_\pi f,
\end{equation}
and consequently, whichever values one may
confer to the parameters $a_{ l+1}, \dots, a_n$, they always generate
an $m$-term group in the variables: $x_1, \dots, x_l$, and of course,
a transitive group.

Now, by conferring to the parameters $a_{ l+1}, \dots, a_n$ all
possible values gradually, one obtains $\infty^{ n-l}$ $m$-term groups
in $l$ variables. On the occasion, it is thinkable, though not
necessary, that these $\infty^{ n-l}$ groups are similar to each
other. If this is not the case, then these groups order themselves in
$\infty^{ n - l - \sigma}$ families, each one consisting in
$\infty^\sigma$ groups, and to be precise, in such a way that two
groups in the same family are similar, while by contrast, two groups
belonging to two different families are not similar.

\plainstatement{In all circumstances, our $\infty^{ n-l}$
groups belong to the same kind of type
\deutsch{Typengattung}; in the latter case, this kind 
depends upon essential parameters
(cf. p.~\pageref{S-447-sq} sq.).}

At present, one may overview how one can find all intransitive
$r$-term groups in $n$ variables. One chooses two numbers $l$ and $m$
so that $l \leqslant m \leqslant r$ and so that in addition $l < n$,
and one then forms all kinds \deutsch{Gattungen} of transitive
$m$-term groups in $l$ variables. If $Y_1f, \dots, Y_mf$:
\[
Y_kf
=
\sum_{i=1}^l\,
\eta_{ki}(x_1,\dots,x_l,\,\alpha_1,\alpha_2,\dots)\,
\frac{\partial f}{\partial x_i}
\ \ \ \ \ \ \ \ \ \ \ \ \ {\scriptstyle{(k\,=\,1\,\cdots\,m)}}
\]
is such a kind with the essential parameters $\alpha_1$, 
$\alpha_2$, \dots, then one interprets these parameters
as unknown functions of $x_{ l+1}, \dots, x_n$, one
sets:
\[
X_kf
=
\sum_{i=1}^m\,\beta_{ki}(x_{l+1},\dots,x_n)\,Y_if
\ \ \ \ \ \ \ \ \ \ \ \ \ {\scriptstyle{(k\,=\,1\,\cdots\,r)}},
\]
and lastly, one attempts to choose the yet undetermined functions
$\alpha_j$, $\beta_{ ki}$ of $x_{ l+1}, \dots, x_n$ in the most
general way so that $X_1f, \dots, X_rf$ become independent
infinitesimal transformations of an $r$-term group. This requirement
conducts to certain \emphasis{finite relations} between the $\alpha$,
the $\beta$ and the $c_{ iks}$ of the sought $r$-term group which must
be satisfied in the most general way. The expressions of the
infinitesimal transformations determined in this way contain, apart
from certain arbitrary constants $c_{ iks}$, yet certain arbitrary
functions of the invariants of the group.

We therefore have at first the

\renewcommand{\thefootnote}{\fnsymbol{footnote}}
\def\thetheorem{83}\begin{theorem}
The determination of all intransitive $r$-term groups: $X_1f, \dots, 
X_rf$ in $n$ variables requires, as soon as all transitive
groups with $r$ or less parameters are found, no integration, 
but only executable operations.\footnote[1]{\,
\name{Lie}, Archiv for Math. og Naturv. Vol.~10, Christiania 1885; 
Math. Ann. Vol. XVI, p.~528, 1880.
}
\end{theorem}
\renewcommand{\thefootnote}{\arabic{footnote}}

It is to be observed that the computations necessary for the
determination of all intransitive groups depend only upon
the numbers $r$, $l$ and $m$. By contrast, the number $n$
plays no direct rôle. Consequently:

\def\thetheorem{84}\begin{theorem}
\label{Theorem-84-S-458}
The determination of all $r$-term groups in an arbitrary number
of variables can be led back to the determination of all 
$r$-term groups in $r$ or less variables. 
\end{theorem}

\smallercharacters{
Yet a brief remark about the determination of all intransitive
$r$-term groups of given composition.

If a group of the concerned composition contains only a finite number
of invariant subgroups, then $m$ must be equal to $r$ (Proposition~8,
p.~\pageref{Satz-8-S-310}), and consequently, the settlement of the
problem formulated just now amounts without effort to the
determination of all transitive groups of the concerned composition.
By contrast, if there occur infinitely many invariant subgroups,
$m$ can be smaller than $r$; then according to the developments just
mentioned, the $m$-term group $\overline{ X}_1f, \dots, \overline{
X}_mf$ discussed above must be meroedrically isomorphic
to the sought $r$-term group. We do not want to undertake to show
more closely how the problem can be settled in this case. 

}

\linestop


\chapter{Invariant Families of Manifolds}
\label{kapitel-23}
\chaptermark{Invariant Families of Manifolds}

\setcounter{footnote}{0}

\abstract*{??}

If $x_1, \dots, x_n$ are point coordinates of an $n$-times
extended space, then a family of manifolds of this space
is represented by equations of the form:
\def\theequation{1}\begin{equation}
\Omega_1(x_1,\dots,x_n,\,l_1,\dots,l_m)
=
0,
\,\,\,\dots,\,\,\,
\Omega_{n-q}(x_1,\dots,x_n,\,l_1,\dots,l_m)
=
0,
\end{equation}
in which, aside from the variables $x_1, \dots, x_n$, 
yet certain parameters: $l_1, \dots, l_m$ are present.

If one executes an arbitrary transformation:
\[
x_i'
=
f_i(x_1,\dots,x_n)
\ \ \ \ \ \ \ \ \ \ \ \ \ {\scriptstyle{(i\,=\,1\,\cdots\,n)}}
\]
of the space $x_1, \dots, x_n$, then each one of the
manifolds~\thetag{ 1} is transferred to a new manifold, hence the
whole family~\thetag{ 1} converts into a new family a new family of
manifolds. One obtains the equations of this new family when one takes
away $x_1, \dots, x_n$ from~\thetag{ 1} with the help of the $n$
equations: $x_i' = f_i ( x)$. Now in particular, if the new family of
manifolds coincides with the original family~\thetag{ 1}, whence every
manifold of the family~\thetag{ 1} is transferred by the
transformation: $x_i' = f_i (x)$ to a manifold of the same family,
then we say \emphasis{that the family~\thetag{ 1}
\terminology{admits} the transformation in question, 
or that it \terminology{remains invariant} by it}.

If a family of manifolds in the space $x_1, \dots, x_n$ admits all
transformations of an $r$-term group, then we say
\emphasis{that it \terminology{admits} the group in question}.

Examples of invariant families of manifolds in the space $x_1, \dots,
x_n$ have already appeared to us several times; every intransitive
group decomposes the space into an invariant family of individually
invariant manifolds (Chap.~\ref{kapitel-13},
pp.~\pageref{S-215}--\pageref{S-216}); every imprimitive group
decomposes the space in an invariant family of manifolds that it
permutes (p.~\pageref{S-220-sq} sq.); also, every manifold
individually invariant by a group may be interpreted as an invariant
family of manifolds, namely as a family parametrized by a point.

In what follows, we now consider a completely arbitrary family of
manifolds. At first, we study under which conditions this family
admits a single given transformation, or a given $r$-term
group. Afterwards, we imagine that a group is given by which the
family remains invariant and we determine the law according to which
the manifolds of the family are permuted with each other by the
transformations of this group. In this way, we obtain a new process to
set up the groups which are isomorphic with a given group. Finally,
we give a method for finding all families of manifolds invariant by a
given group.

\sectionengellie{\S\,\,\,110.}

Let the equations:
\def\theequation{1}\begin{equation}
\Omega_k(x_1,\dots,x_n,\,l_1,\dots,l_m)
=
0
\ \ \ \ \ \ \ \ \ \ \ \ \ 
{\scriptstyle{(k\,=\,1\,\cdots\,n\,-\,q)}}
\end{equation}
with the $m$ arbitrary parameters $l_1, \dots, l_m$ represent
an arbitrary family of manifolds.

If $l_1, \dots, l_m$ are absolutely arbitrary parameters, 
then naturally, one should not be able to eliminate all the $x$
from~\thetag{ 1}; 
therefore, the equations~\thetag{ 1} must be solvable
with respect to $n - q$ of the variables
$x_1, \dots, x_n$. 

By contrast, it is not excluded that the $l$ can be eliminated
from~\thetag{ 1}, and this says nothing but,
that some relations between the $x$ alone can be 
derived from~\thetag{ 1}. Only the number of independent
equations free of the $l$ which follow from~\thetag{ 1}
must be smaller than $n - q$, since otherwise, 
the parameters $l_1, \dots, l_r$ would be only apparent 
and the equations~\thetag{ 1} would therefore represent
not a family of manifolds, but a single manifold.

Before we study how the family of manifolds~\thetag{ 1}
behaves relatively to transformations of the $x$, we
must first make a few remarks about the nature of the 
equations~\thetag{ 1}. 

The equations~\thetag{ 1} contain $m$ parameters $l_1, \dots, l_m$; 
if we let these parameters take all possible values, then
we obtain $\infty^m$ different systems of values: 
$l_1, \dots, l_m$, but not necessarily $\infty^m$
different manifolds. It must therefore be 
determined under which conditions the given system of equations
represents exactly $\infty^m$ different manifolds, or in 
other words: one must give a criterion to determine
whether the parameters: $l_1, \dots, l_m$
in the equations~\thetag{ 1} are
\emphasis{essential}, or not.

In order to find such a criterion, we imagine that the equations: 
$\Omega_k = 0$ are solved with respect to $n - q$ of the $x$, 
say with respect to $x_{ q+1}, \dots, x_n$:
\def\theequation{2}\begin{equation}
\aligned
x_{q+k}
&
=
\psi_{q+k}
(x_1,\dots,x_q,\,l_1,\dots,l_m)
\\
&
\ \ \ \ \ \ \ \ \ \ \ \ \ 
{\scriptstyle{(k\,=\,1\,\cdots\,n\,-\,q)}}
\endaligned
\end{equation}
and we imagine moreover that the functions $\psi_{ q+k}$ are
expanded with respect to the powers of: $x_1 - x_1^0$, \dots, 
$x_q - x_q^0$ in the neighbourhood of an arbitrary system
of values: $x_1^0, \dots, x_q^0$. The coefficients of the
expansion, whose number is naturally infinitely large
in general, will be analytic functions of $l_1, \dots, l_m$
and we may call them: 
\[
\Lambda_j(l_1,\dots,l_m)
\ \ \ \ \ \ \ \ \ \ \ \ \ 
{\scriptstyle{(j\,=\,1,\,\,2\,\cdots\,)}}
\]
The question amounts just to how many independent functions are extant
amongst all the functions $\Lambda_1$, $\Lambda_2$, \dots.

Indeed, if amongst the $\Lambda$, there are exactly $l$ that are
mutually independent functions\,---\,there are anyway surely no more
than $l$\,---, then to the $\infty^m$ different systems of values:
$l_1, \dots, l_m$, there obviously correspond also $\infty^m$
different systems of values: $\Lambda_1$, $\Lambda_2$, \dots, and
therefore $\infty^m$ different manifolds~\thetag{ 2}, that is to say,
the parameters are essential in the equations~\thetag{ 2}, and hence
also in the equations~\thetag{ 1}.

Otherwise, assume that amongst the functions $\Lambda_1$, $\Lambda_2$,
\dots, there are not $m$
but less functions, say only $m - h$, that are mutually
independent. In this case, all the $\Lambda_j$ can be expressed in
terms of $m - h$ of them, say in terms of: $\Lambda_1', \dots,
\Lambda_{ m-h}'$, which naturally must be mutually independent. To the
$\infty^m$ different systems of values $l_1, \dots, l_m$, 
there correspond therefore $\infty^{ m-h}$ different systems
of values $\Lambda_1', \dots, \Lambda_{ m-h}'$, and
$\infty^{ m-h}$ different systems of values $\Lambda_1$, 
$\Lambda_2$, \dots, so that the equations~\thetag{ 2}, 
and thus also the equations~\thetag{ 1}, represent
only $\infty^{ m-h}$ different manifolds.
This comes to expression in the clearest way
when one observes that the functions $\psi_{ q + k} (x,l)$
contain the parameters $l_1, \dots, l_m$ only in the
combinations: $\Lambda_1', \dots, \Lambda_{ m-h}'$, so 
that the equations~\thetag{ 2} possess the form:
\def\theequation{2'}\begin{equation}
x_{q+k}
=
\overline{\psi}_{q+k}
(x_1,\dots,x_q,\,
\Lambda_1',\dots,\Lambda_{m-h}')
\ \ \ \ \ \ \ \ \ \ \ \ \ 
{\scriptstyle{(k\,=\,1\,\cdots\,n\,-\,q)}}.
\end{equation}
From this, it indeed results that we can 
introduce precisely $\Lambda_1', \dots, \Lambda_{ m-h}'$
as new parameters in place of $l_1, \dots, l_n$, a process by which
the number of arbitrary parameters
appearing in~\thetag{ 2} is lowered to $m - h$.

Thus, we can say:

\plainstatement{The equations:
\[
\Omega_k(x_1,\dots,x_n,\,l_1,\dots,l_m)
=
0
\ \ \ \ \ \ \ \ \ \ \ \ \ 
{\scriptstyle{(k\,=\,1\,\cdots\,n\,-\,q)}}
\]
represent $\infty^m$ different manifolds only when it is not possible
to indicate $m - h < m$ functions $\pi_1, \dots, \pi_{ m-h}$ of $l_1,
\dots, l_m$ such that, in the resolved equations:
\[
x_{q+k}
=
\psi_{q+k}(x_1,\dots,x_q,\,l_1,\dots,l_m)
\ \ \ \ \ \ \ \ \ \ \ \ \ 
{\scriptstyle{(k\,=\,1\,\cdots\,n\,-\,q)}},
\]
the functions $\psi_{ q+1}, \dots, \psi_n$ can be expressed in terms
of $x_1, \dots, x_n$ and of $\pi_1, \dots, \pi_{ m-h}$ only. By
contrast, if it is possible to indicate $m - h < m$ functions
$\pi_\mu$ having this constitution, then the equations $\Omega_k = 0$
represent at most $\infty^{ m-h}$ manifolds and the parameters $l_1,
\dots, l_m$ are hence not essential.}

\smallercharacters{
Here, yet a brief remark.

If the functions: $\Lambda_1$, $\Lambda_2$, \dots\, discussed above
take the values: $\Lambda_1^0$, $\Lambda_2^0$ after the substitution:
$l_1 = l_1^0$, \dots, $l_m = l_m^0$, then the equations:
\[
\Lambda_j(l_1,\dots,l_m)
=
\Lambda_j^0
\ \ \ \ \ \ \ \ \ \ \ \ \ 
{\scriptstyle{(j\,=\,1,\,\,2\,\cdots\,)}}
\]
define the totality of all systems of values $l_1, \dots, l_m$ which,
when inserted in~\thetag{ 2} or in~\thetag{ 1}, provide the same
manifold as the system of values: $l_1^0, \dots, l_m^0$. Now, if
amongst the functions $\Lambda_1$, $\Lambda_2$, \dots, there are $m$
functions that are mutually independent, then the parameters $l_1,
\dots, l_m$ are all essential, whence the following obviously holds:
Around every system of values $l_1^0, \dots, l_m^0$ in general
position, one can delimit a region such that two distinct systems of
values $l_1, \dots, l_m$ of the concerned region always provide two
manifolds that are also distinct.

}

The definition for essentiality \deutsch{Wesentlichkeit} and for
inessentiality \deutsch{Nichtwesentlichkeit} of the parameters $l_1,
\dots, l_m$ (respectively) given above is satisfied only when the
equations: $\Omega_1 = 0$, \dots, $\Omega_q = 0$ are already solved
with respect to $n - q$ of the $x$. However, it is desirable, in
principle and also for what follows, to reshape this definition so
that it fits also to a not resolved system of equations $\Omega_k =
0$.

There is no difficulty for doing that.

Let the parameters $l_1, \dots, l_m$ in the equations:
\def\theequation{2}\begin{equation}
x_{q+k}
=
\psi_{q+k}
(x_1,\dots,x_q,\,l_1,\dots,l_m)
\ \ \ \ \ \ \ \ \ \ \ \ \ 
{\scriptstyle{(k\,=\,1\,\cdots\,n\,-\,q)}}
\end{equation}
be not essential, whence one can indicate $m - h < m$ functions:
$\pi_1 (l), \dots, \pi_{ m-h} (l)$ so that $\psi_{ q+1} (x, l)$,
\dots, $\psi_n ( x,l)$ can be expressed in terms of $x_1, \dots, x_q$
and of $\pi_1 (l), \dots, \pi_{ m-h} (l)$ alone. 
Evidently, there is then at least one linear partial differential 
equation: 
\[
Lf
=
\sum_{\mu=1}^m\,
\lambda_\mu(l_1,\dots,l_m)\,
\frac{\partial f}{\partial l_\mu}
=
0
\]
with the coefficients: $\lambda_1 ( l), \dots, \lambda_\mu ( l)$
free of the $x$ which is satisfied identically by all functions:
$\pi_1 ( l), \dots, \pi_{ m - h} (l)$, and hence also
by all the functions: $x_{ q+i} - \psi_{ q + i} ( x, l)$. 
Thus, we can also say (cf. Chap.~\ref{kapitel-7},
Theorem~15, p.~\pageref{Theorem-15-S-117}): 

When the parameters $l_1, \dots, l_m$ are not essential in~\thetag{
2}, then the system of equations~\thetag{ 2}, interpreted as a system
of equations in the variables: $x_1, \dots, x_n$, $l_1, \dots, l_m$,
admits an infinitesimal transformation $Lf$ in the variables $l_1,
\dots, l_m$ alone.

But the converse also holds true: When the system of equation~\thetag{
2}, interpreted as a system of equations in the variables $x_1, \dots,
x_n$, $l_1, \dots, l_m$, admits an infinitesimal transformation $Lf$
in the $l$ alone, then the parameters $l_1, \dots, l_m$ are not
essential in the system of equations; it is indeed immediately clear
that, under the assumption made, $\psi_{ q+1} ( x, l)$, \dots, $\psi_n
( x,l)$ are solutions of the partial differential equation: $Lf = 0$,
whence it is possible to lower the number of arbitrary parameters
appearing in~\thetag{ 2}.

If we bear in mind that the system of equations~\thetag{ 2}
is only another form of the system of equations~\thetag{ 1}, 
then we realize without effort
(cf. p.~\pageref{Satz-2-S-111}) that the following proposition holds:

\def\theproposition{1}\begin{proposition}
A system of equations:
\def\theequation{1}\begin{equation}
\Omega_k
(x_1,\dots,x_n,\,l_1,\dots,l_m)
=
0
\ \ \ \ \ \ \ \ \ \ \ \ \ 
{\scriptstyle{(k\,=\,1\,\cdots\,n\,-\,q)}}
\end{equation}
with the $l_1, \dots, l_m$ parameters which is solvable with respect
to $n - q$ of the variables $x_1, \dots, x_n$ represents $\infty^m$
different manifolds of the space $x_1, \dots, x_n$ if and only if,
when it is regarded as a system of equations in the 
$n + m$ variables: $x_1, \dots, x_n$, $l_1, \dots, l_m$, 
it admits no infinitesimal transformation:
\[
Lf
=
\sum_{\mu=1}^m\,
\lambda_\mu(l_1,\dots,l_m)\,
\frac{\partial f}{\partial l_\mu}
\]
in the variables $l$ alone.
\end{proposition}

\sectionengellie{\S\,\,\,111.}

\label{S-463-sq}
In the space $x_1, \dots, x_n$, let at present a family of $\infty^m$
different manifolds be determined by the $n - q$ equations: $\Omega_k
( x, l) = 0$, or by the equally good \deutsch{gleichwertig} equations:
\def\theequation{2}\begin{equation}
x_{q+k}
=
\psi_{q+k}
(x_1,\dots,x_n,\,l_1,\dots,l_m)
\ \ \ \ \ \ \ \ \ \ \ \ \ 
{\scriptstyle{(k\,=\,1\,\cdots\,n\,-\,q)}}.
\end{equation}

If this family is supposed to remain invariant by the transformation:
$x_i' = f_i ( x_1, \dots, x_n)$, then every manifold of the family
must be transferred by the concerned transformation again into a
manifold of the family. Hence, if by $l_1', \dots, l_m'$ we
understand the parameters of the manifold of the family into which the
manifold with the parameters $l_1, \dots, l_m$ 
\label{S-463} is transferred by
the transformation: $x_i' = f_i ( x)$, then after the introduction of
the new variables: $x_1' = f_1 ( x)$, \dots, $x_n' = f_n (x)$, the
equations~\thetag{ 2} must receive the form:
\[
x_{q+k}'
=
\psi_{q+k}
(x_1',\dots,x_q',\,l_1',\dots,l_m')
\ \ \ \ \ \ \ \ \ \ \ \ \ 
{\scriptstyle{(k\,=\,1\,\cdots\,n\,-\,q)}},
\]
where the parameters $l_1', \dots, l_m'$ 
depend naturally only upon the $l$.

But now, after the introduction of the $x'$, the
equation~\thetag{ 2} evidently take up the form:
\[
x_{q+k}'
=
\Psi_{q+k}
(x_1',\dots,x_q',\,l_1,\dots,l_m)
\ \ \ \ \ \ \ \ \ \ \ \ \ 
{\scriptstyle{(k\,=\,1\,\cdots\,n\,-\,q)}}\,;
\]
thus, if the family~\thetag{ 2} is supposed to 
remain invariant by the transformation: $x_i' = f_i ( x)$, 
then it must be possible to satisfy the $n - q$
equations:
\def\theequation{3}\begin{equation}
\psi_{q+k}
(x_1',\dots,x_q',\,l_1',\dots,l_m')
=
\Psi_{q+k}(x_1',\dots,x_q',\,l_1,\dots,l_m)
\ \ \ \ \ \ \ \ \ \ \ \ \ 
{\scriptstyle{(k\,=\,1\,\cdots\,n\,-\,q)}}
\end{equation}
independently of the values of the variables $x_1', \dots, x_q'$.

If one expands the two sides of~\thetag{ 3}, in the neighbourhood of an
arbitrary system of values ${x_1'}^0, \dots, {x_n'}^0$, with respect to
the powers of: $x_1' - {x_1'}^0$, \dots, $x_n' - {x_n'}^0$, if one
compares the coefficients, and if one takes into account that $l_1',
\dots, l_m'$ are essential parameters, then one realizes that $l_1',
\dots, l_m'$ must be entirely determined functions of $l_1, \dots,
l_m$: 
\[
l_\mu'
=
\chi_\mu(l_1,\dots,l_m)
\ \ \ \ \ \ \ \ \ \ \ \ \ 
{\scriptstyle{(\mu\,=\,1\,\cdots\,m)}}
\]

Conversely, the $l$ must also naturally be representable as functions
of the $l'$, because through the transition from the $x'$ to the $x$,
the family of our manifolds must also remain unchanged.

We therefore see that the equations:
\[
x_i'
=
f_i(x_1,\dots,x_n)
\ \ \ \ \ \ \ 
{\scriptstyle{(i\,=\,1\,\cdots\,n)}},
\ \ \ \ \ \ \ \ \
l_\mu'
=
\chi_\mu(l_1,\dots,l_m)
\ \ \ \ \ \ \ 
{\scriptstyle{(\mu\,=\,1\,\cdots\,m)}}
\]
taken together represent a transformation in the $n + m$ variables:
$x_1, \dots, x_n$, $l_1, \dots, l_m$ that leaves invariant the system
of equations: $x_{ q + k} = \psi_{ q + k} (x, l$ in these $n + m$
variables. Consequently, we can say:

\shortplainstatement{
A family of $\infty^m$ manifolds:
\[
\Omega_k
(x_1,\dots,x_n,\,l_1,\dots,l_m)
=
0
\ \ \ \ \ \ \ \ \ \ \ \ \ 
{\scriptstyle{(k\,=\,1\,\cdots\,n\,-\,q)}}
\]
in the space $x_1, \dots, x_n$ admits the transformation:
\[
x_i'
=
f_i(x_1,\dots,x_n)
\ \ \ \ \ \ \ \ \ \ \ \ \ {\scriptstyle{(i\,=\,1\,\cdots\,n)}}
\]
if and only if it is possible to add to this transformation
in the $x$ a corresponding transformation:
\[
\label{S-464}
l_\mu'
=
\chi_\mu(l_1,\dots,l_m)
\ \ \ \ \ \ \ \ \ \ \ \ \
{\scriptstyle{(\mu\,=\,1\,\cdots\,m)}} 
\]
in the $l$ in such a way that the system of equations:
$\Omega_k ( x, l) = 0$ in the $n + m$ variables $x_1, \dots,
x_n$, $l_1, \dots, l_m$ allows the transformation:
\[
x_i'
=
f_i(x_1,\dots,x_n),
\ \ \ \ \ \ \
l_\mu'
=
\chi_\mu(l_1,\dots,l_m).
\]}

From the above considerations, it becomes clear that the
transformation: $l_\mu ' = \chi_\mu (l)$, when it actually exists,
is the only one of its kind; indeed, it is completely 
determined when the transformation $x_i' = f_i (x)$ is
known. The transformation: $l_\mu ' = \chi_\mu (l)$ 
therefore contains no arbitrary parameters. 

If the family of the $\infty^m$ manifolds: $\Omega_k (x, l) = 0$
allows two different transformations:
\[
x_i'
=
f_i(x_1,\dots,x_n)\,;
\ \ \ \ \ \
x_i''
=
F_i(x_1',\dots,x_n'),
\]
then obviously, it admits also the transformation:
\[
x_i''
=
F_i
\big(
f_1(x),\,\dots,\,f_n(x)
\big),
\] 
which comes into existence by executing each one 
of the two transformations one after the other. 
From this, we conclude:

\plainstatement{The totality of all transformations $x_i' = f_i ( x_1,
\dots, x_n)$ which leave invariant a family of $\infty^m$ manifolds of
the space $x_1, \dots, x_n$ forms a group.}

Naturally, this group needs neither be finite, nor continuous; in
complete generality, we can only say: its transformations are
ordered as inverses by pairs. Hence in particular, when it contains
only a finite number of arbitrary parameters, then it belongs to the
category of groups which was discussed in Chap.~\ref{kapitel-18}, and
according to Theorem~56, p.~\pageref{Theorem-56-S-315}, it contains a
completely determined finite continuous subgroup generated by
infinitesimal transformations. Evidently, this subgroup is the
largest continuous subgroup by which the family $\Omega_k ( x, l) = 0$
remains invariant.

We now turn ourselves to the consideration of finite continuous
groups which leave invariant the family of the $\infty^m$ manifolds:
\def\theequation{1}\begin{equation}
\Omega_k(x_1,\dots,x_n,\,l_1,\dots,l_m)
\ \ \ \ \ \ \ \ \ \ \ \ \ 
{\scriptstyle{(k\,=\,1\,\cdots\,n\,-\,q)}}\,;
\end{equation}
however, we restrict ourselves at first for reasons of simplicity
to the case of a one-term group having 
the concerned constitution.

Let the family of the $\infty^m$ manifolds~\thetag{ 1}
admit all transformations:
\def\thequation{4}\begin{equation}
x_i'
=
f_i(x_1,\dots,x_n,\,\varepsilon)
\ \ \ \ \ \ \ \ \ \ \ \ \ {\scriptstyle{(i\,=\,1\,\cdots\,n)}}
\end{equation}
of the one-term group:
\[
Xf
=
\sum_{i=1}^n\,\xi_{ki}(x_1,\dots,x_n)\,
\frac{\partial f}{\partial x_i}.
\]
The transformation in the $l$ which, according to p.~\pageref{S-464},
corresponds to the general transformation: $x_i' = f_i ( x, a)$, can
be read:
\def\theequation{4'}\begin{equation}
l_\mu'
=
\chi_\mu(l_1,\dots,l_m,\,\varepsilon)
\ \ \ \ \ \ \ \ \ \ \ \ \ {\scriptstyle{(\mu\,=\,1\,\cdots\,m)}}.
\end{equation}

It is easy to see that the totality of all transformations
of the form: 
\def\theequation{4''}\begin{equation}
\left\{
\aligned
x_i'
&
=
f_i(x_1,\dots,x_n,\,\varepsilon)
\ \ \ \ \ \ \ \ \ \ \ \ \ {\scriptstyle{(i\,=\,1\,\cdots\,n)}}
\\
l_\mu'
&
=
\chi_\mu(l_1,\dots,l_m,\,\varepsilon)
\ \ \ \ \ \ \ \ \ \ \ \ \ {\scriptstyle{(\mu\,=\,1\,\cdots\,m)}}
\endaligned\right.
\end{equation}
forms in turn a one-term group in the $n + m$ variables:
$x_1, \dots, x_n$, $l_1, \dots, l_m$.

In fact, according to p.~\pageref{S-464}, 
the transformations~\thetag{ 4''} leave invariant
the system of equations~\thetag{ 1};
hence if one executes at first the transformation~\thetag{ 4''}
and afterwards a transformation of the same form with the
parameter $\varepsilon_1$, then one gets a transformation:
\[
\aligned
x_i''
&
=
f_i\big(
f_1(x,\varepsilon),\,\dots,\,f_n(x,\varepsilon),\,
\varepsilon_1\big)
\ \ \ \ \ \ \ \ \ \ \ \ \ {\scriptstyle{(i\,=\,1\,\cdots\,n)}}
\\
l_\mu''
&
=
\chi_\mu\big(
\chi_1(l,\varepsilon),\,\dots,\,\chi_m(l,\varepsilon),\,
\varepsilon_1\big)
\ \ \ \ \ \ \ \ \ \ \ \ \ {\scriptstyle{(\mu\,=\,1\,\cdots\,m)}}
\endaligned
\]
which, likewise, leaves invariant the system of equations~\thetag{ 1}.
Now, the transformation:
\def\theequation{5}\begin{equation}
x_i''
=
f_i\big(
f_1(x,\varepsilon),\,\dots,\,f_n(x,\varepsilon),\,
\varepsilon_1\big)
\end{equation}
belongs to the one-term group $Xf$ and can hence
be brought to the form:
\[
x_i''
=
f_i(x_1,\dots,x_n,\,\varepsilon_2)
\ \ \ \ \ \ \ \ \ \ \ \ \ {\scriptstyle{(i\,=\,1\,\cdots\,n)}},
\]
where $\varepsilon_2$ depends only on $\varepsilon$ and
on $\varepsilon_1$. Consequently, 
the transformation in the $l$ which corresponds to the
transformation~\thetag{ 5} has the shape:
\[
l_\mu''
=
\chi_\mu(l_1,\dots,l_m,\,\varepsilon_2)
\ \ \ \ \ \ \ \ \ \ \ \ \ {\scriptstyle{(\mu\,=\,1\,\cdots\,m)}}
\]
and we therefore deduce:
\[
\aligned
\chi_\mu\big(
\chi_1(l,\varepsilon),\,\dots,\,
&
\chi_m(l,\varepsilon),\,\varepsilon_1)
=
\chi_\mu(l_1,\dots,l_m,\,\varepsilon_2)
\\
&
\ \ \ \ \
{\scriptstyle{(\mu\,=\,1\,\cdots\,m)}}.
\endaligned
\]

As a result, it is proved that the equations~\thetag{ 4''} effectively
represent a group in the $n + m$ variables: $x_1, \dots, x_n$, $l_1,
\dots, l_m$, and to be precise, a group which possesses the same
parameter group as the given group: $x_i' = f_i ( x, \varepsilon)$.
At the same time, it is clear that the equations~\thetag{ 4'} taken
for themselves also represent a group in the variables: $l_1, \dots,
l_m$, however not necessarily a one-term group, because it is
thinkable that the parameter $\varepsilon$ in the equations~\thetag{
4'} is completely missing.

\smallercharacters{

The transformations of the group~\thetag{ 4} order as inverses by
pairs, hence the same visibly holds true for the transformations of
the group~\thetag{ 4'}. From this, it comes
(cf. Chap.~\ref{kapitel-9}, p.~\pageref{S-169} above) that the
group~\thetag{ 4''} also contains the identity transformation, and in
addition, an infinitesimal transformation:
\[
\sum_{i=1}^n\,
\xi_i(x_1,\dots,x_n)\,
\frac{\partial f}{\partial x_i}
+
\sum_{\mu=1}^m\,\lambda_\mu
(l_1,\dots,l_m)\,
\frac{\partial f}{\partial\lambda_\mu}
=
Xf
+
Lf,
\]
by which it is generated.

}

We recapitulate the gained result in the:

\def\theproposition{2}\begin{proposition}
If the family of the $\infty^m$ manifolds:
\[
\Omega_k(x_1,\dots,x_n,\,l_1,\dots,l_m)
=
0
\ \ \ \ \ \ \ \ \ \ \ \ \ 
{\scriptstyle{(k\,=\,1\,\cdots\,n\,-\,q)}}
\]
of the space $x_1, \dots, x_n$ admits all transformations:
\[
x_i'
=
f_i(x_1,\dots,x_n,\,\varepsilon)
\ \ \ \ \ \ \ \ \ \ \ \ \ {\scriptstyle{(i\,=\,1\,\cdots\,n)}}
\]
of the one-term group:
\[
Xf
=
\sum_{i=1}^n\,\xi_i(x_1,\dots,x_n)\,
\frac{\partial f}{\partial x_i},
\]
then the corresponding transformations:
\[
\aligned
x_i'
&
=
f_i(x_1,\dots,x_n,\,\varepsilon)
\ \ \ \ \ \ \ \ \ \ \ \ \ {\scriptstyle{(i\,=\,1\,\cdots\,n)}}
\\
l_\mu'
&
=
\chi_\mu(l_1,\dots,l_m,\,\varepsilon)
\ \ \ \ \ \ \ \ \ \ \ \ \ {\scriptstyle{(\mu\,=\,1\,\cdots\,m)}}
\endaligned
\]
which leave invariant the system of equations: $\Omega_1 ( x, l) = 0$,
\dots, $\Omega_{ n-q} (x, l) = 0$ in the $n + m$ variables: $x_1,
\dots, x_n$, $l_1, \dots, l_m$ form a one-term group with an
infinitesimal transformation of the shape:
\[
\sum_{i=1}^n\,
\xi_i(x_1,\dots,x_n)\,
\frac{\partial f}{\partial x_i}
+
\sum_{\mu=1}^m\,\lambda_\mu
(l_1,\dots,l_m)\,
\frac{\partial f}{\partial\lambda_\mu}
=
Xf
+
Lf.
\]
\end{proposition}

At present, we set up the following definition:

\smallskip{\sf Definition.}
\emphasis{A family of $\infty^m$ manifolds:
\[
\Omega_k(x_1,\dots,x_n,\,l_1,\dots,l_m)
=
0
\ \ \ \ \ \ \ \ \ \ \ \ \ 
{\scriptstyle{(k\,=\,1\,\cdots\,n\,-\,q)}}
\]
of the space $x_1, \dots, x_n$ admits the infinitesimal
transformation:
\[
Xf
=
\sum_{i=1}^n\,\xi_i(x_1,\dots,x_n)\,
\frac{\partial f}{\partial x_i}
\]
when there is, in $l_1, \dots, l_m$, an infinitesimal
transformation:
\[
Lf
=
\sum_{\mu=1}^m\,
\lambda_\mu(l_1,\dots,l_m)\,
\frac{\partial f}{\partial l_\mu}
\]
such that the system of equations: $\Omega_1 ( x, l) = 0$, \dots,
$\Omega_{ n-q} ( x, l) = 0$ in the $n + m$ variables: $x_1, \dots,
x_n$, $l_1, \dots, l_m$ admits the infinitesimal transformation $Xf +
Lf$.}

It yet remains here to consider the question whether the infinitesimal
transformation $Lf$ is completely determined by the transformation
$Xf$. One easily realizes that this is the case; indeed, if the system
of equations: $\Omega_k ( x, l) = 0$ admits the two infinitesimal
transformations: $Xl + Lf$ and $Xf + \mathfrak{ L}f$, then it admits
at the same time the transformation:
\[
Xf+Lf
-
(Xf+\mathfrak{L}f)
=
Lf
-
\mathfrak{L}f\,;
\]
but since the parameters of the family: $\Omega_k ( x, l) = 0$ are
essential, the expression: $Lf - \mathfrak{ L}f$ must vanish
identically, hence the transformation $\mathfrak{ L}f$ cannot be
distinct from the transformation $Lf$.

By taking as a basis the above definition, we can obviously express
the content of Proposition~2 as follows:

If the family of the $\infty^m$ manifolds:
\[
\Omega_k(x_1,\dots,x_n,\,l_1,\dots,l_m)
=
0
\ \ \ \ \ \ \ \ \ \ \ \ \ {\scriptstyle{(k\,=\,1\,\cdots\,n\,-\,q)}}
\]
admits the one-term group $Xf$, then at the same time, 
it admits the infinitesimal transformation $Xf$. 

Conversely: when the family of the $\infty^m$ manifolds:
$\Omega_k ( x, l) = 0$ admits the infinitesimal
transformation $Xf$, then it also admits the 
one-term group $Xf$. 

Indeed, under the assumption made, if the system of 
equations: $\Omega_k ( x, l) = 0$ in the $n + m$ variables
$x$, $l$ admits an infinitesimal transformation of the form:
\[
Xf+Lf
=
\sum_{i=1}^n\,\xi_i(x_1,\dots,x_n)\,
\frac{\partial f}{\partial x_i}
+
\sum_{\mu=1}^m\,\lambda_\mu
(l_1,\dots,l_m)\,
\frac{\partial f}{\partial l_\mu},
\]
then\footnote{\,
See Chap.~\ref{kapitel-6}, Theorem~14, 
p.~\pageref{Theorem-14-S-112} for this general fact.
} 
it admits at the same time all transformations of the one-term group:
$Xf + Lf$, and therefore, the family of the $\infty^m$ manifolds:
$\Omega_k ( x, l) = 0$ in the space $x_1, \dots, x_n$ allows all
transformations of the one-term group $Xf$.

With these words, we have prove the:

\def\theproposition{3}\begin{proposition}
\label{Satz-3-S-468}
The family of the $\infty^m$ manifolds:
\[
\Omega_k(x_1,\dots,x_n,\,l_1,\dots,l_m)
=
0
\ \ \ \ \ \ \ \ \ \ \ \ \ {\scriptstyle{(k\,=\,1\,\cdots\,n\,-\,q)}}
\]
in the space $x_1, \dots, x_n$ admits the one-term group $Xf$ if and
only if it admits the infinitesimal transformation $Xf$.
\end{proposition}

If one wants to know whether the family of the $\infty^m$ manifolds:
$\Omega_k ( x, l) = 0$ admits a given infinitesimal transformation
$Xf$, then one will at first solve the equations $\Omega_k ( x, l) =
0$ with respect to $n - q$ of the $x$:
\def\theequation{2}\begin{equation}
x_{q+k}
=
\psi_{q+k}
(x_1,\dots,x_q,\,l_1,\dots,l_m)
\ \ \ \ \ \ \ \ \ \ \ \ \ {\scriptstyle{(k\,=\,1\,\cdots\,n\,-\,q)}},
\end{equation}
and afterwards, one will attempt to determine $m$ functions:
$\lambda_1 ( l), \dots, \lambda_m ( l)$ of the $l$ so that the system
of equations~\thetag{ 2} in the $n + m$ variables $x$, $l$ admits the
infinitesimal transformation:
\[
Xf+Lf
=
\sum_{i=1}^n\,\xi_i(x_1,\dots,x_n)\,
\frac{\partial f}{\partial x_i}
+
\sum_{\mu=1}^m\,\lambda_\mu(l_1,\dots,l_m)\,
\frac{\partial f}{\partial l_\mu}.
\]

If one denotes the substitution: $x_{ q+1} = \psi_{ q+1}$, \dots,
$x_n = \psi_n$ by the sign 
\label{S-468} $[ \,\,\, ]$, then one obviously 
obtains for $\lambda_1 (l)$, \dots, $\lambda_m (l)$ the following
equations:
\[
\aligned
\sum_{\mu=1}^m\,\lambda_\mu(l)\,
&
\frac{\partial\psi_{q+k}}{\partial l_\mu}
=
\big[\xi_{q+k}\big]
-
\sum_{j=1}^q\,
\big[\xi_j\big]\,
\frac{\partial\psi_{q+k}}{\partial x_j}
\\
&
\ \ \ \ \ \ \ \ \ \ \ {\scriptstyle{(k\,=\,1\,\cdots\,n\,-\,q)}},
\endaligned
\]
which must be satisfied independently of the values of the
variables $x_1, \dots, x_q$. Theoretically, 
there is absolutely no difficulty to decide whether this is
possible. One finds either that there is no system of functions
$\lambda_\mu ( l)$ which has the constitution demanded, 
or one finds a system of this sort, but then also only one
such system, 
for the infinitesimal transformation $Lf$ is indeed,
when it actually exists, completely determined by $Xf$. 

At present, we will prove that the following proposition holds:

\def\theproposition{4}\begin{proposition}
\label{Satz-4-S-469}
If the family of the $\infty^m$ manifolds:
\[
\Omega_k(x_1,\dots,x_n,\,l_1,\dots,l_m)
=
0
\ \ \ \ \ \ \ \ \ \ \ {\scriptstyle{(k\,=\,1\,\cdots\,n\,-\,q)}}
\]
in the space $x_1, \dots, x_n$ admits the two infinitesimal 
transformations:
\[
X_1f
=
\sum_{i=1}^n\,\xi_{1i}(x)\,
\frac{\partial f}{\partial x_i},
\ \ \ \ \ \ \ \ \ \ \ \
X_2f
=
\sum_{i=1}^n\,\xi_{2i}(x)\,
\frac{\partial f}{\partial x_i},
\]
then it admits not only every infinitesimal transformation:
\[
e_1\,X_1f+e_2\,X_2f
\]
which can be linearly deduced from $X_1f$ and $X_2f$, but
also the transformation: 
$\leftbracket X_1, \, X_2 \rightbracket$.
\end{proposition}

Under the assumptions that are made in the proposition, there are two
infinitesimal transformations $L_1f$ and $L_2f$ in the $l$ alone that
are constituted in such a way that the system of equations: $\Omega_k
( x, l) = 0$ in the $n + m$ variables $x$, $l$ allows the two
infinitesimal transformations $X_1f + L_1f$, $X_2f + L_2f$. According
to Chap.~\ref{kapitel-7}, Proposition~5, p.~\pageref{Satz-5-S-118},
the system of equations $\Omega_k ( x, l) = 0$ then also admits the
infinitesimal transformation: $\leftbracket X_1, \, X_2 \rightbracket
+ \leftbracket L_1, \, L_2 \rightbracket$ coming into existence by
combination; but this precisely says that the family of the $\infty^m$
manifolds $\Omega_k ( x, l) = 0$ also admits $\leftbracket X_1, X_2
\rightbracket$. As a result, our proposition is proved.

The considerations that we have applied in the proof
just conducted also give yet what follows:

\plainstatement{
If the family of the $\infty^m$ manifolds $\Omega_k ( x, l) = 0$
admits the two infinitesimal transformations $X_1f$ and $X_2f$, and if
$L_1f$ and $L_2f$, respectively, are the corresponding infinitesimal
transformations in $l_1, \dots, l_m$ alone, then the infinitesimal
transformation $\leftbracket L_1, \, L_2 \rightbracket$
corresponds to the transformation 
$\leftbracket X_1, \, X_2 \rightbracket$.}

\smallskip

At present, we assume that the family of the $\infty^m$
manifolds: $\Omega_k ( x, l) = 0$ in the space: $x_1, \dots, x_n$
allows all transformations of an arbitrary $r$-term group:
\def\theequation{6}\begin{equation}
x_i'
=
f_i(x_1,\dots,x_n,\,a_1,\dots,a_r)
\ \ \ \ \ \ \ \ \ \ \ \ \ {\scriptstyle{(i\,=\,1\,\cdots\,n)}}.
\end{equation}
Let the group in question be generated by the $r$ independent
infinitesimal transformations:
\[
X_kf
=
\sum_{i=1}^n\,\xi_{ki}(x_1,\dots,x_n)\,
\frac{\partial f}{\partial x_i},
\]
so that between: $X_1f, \dots, X_nf$, there are relations of the
known form:
\[
\leftbracket
X_i,\,X_k
\rightbracket
=
\sum_{s=1}^r\,c_{iks}\,X_sf
\ \ \ \ \ \ \ \ \ \ \ \ \ {\scriptstyle{(i,\,\,k\,=\,1\,\cdots\,r)}}.
\]

If:
\def\theequation{6'}\begin{equation}
l_\mu'
=
\chi_\mu(l_1,\dots,l_m,\,a_1,\dots,a_r)
\ \ \ \ \ \ \ \ \ \ \ \ \ {\scriptstyle{(\mu\,=\,1\,\cdots\,m)}}
\end{equation}
is the transformation in the $l$ which corresponds to the general
transformation: $x_i' = f_i ( x,a)$ of the group: $X_1f, \dots, X_rf$,
then the equations~\thetag{ 6} and~\thetag{ 6'} taken together
represent an $r$-term group in the $n + m$ variables: $x_1, \dots,
x_n$, $l_1, \dots, l_m$.

In fact, if $T_{ (a_1, \dots, a_r)}$, or shortly $T_{ (a)}$, and $T_{
(b_1, \dots, b_r)}$, or shortly $T_{ (b)}$, are two transformations of
the group: $x_i' = f_i ( x,a)$, then when executed one after the
other, they produce in the known way a transformation: $T_{ (a)} \,
T_{ (b)} = T_{ (c)}$ of the same group, where the parameters $c_1,
\dots, c_r$ depend only on the $a$ and on the $b$.

On the other hand, if $S_{ (a)}$ and $S_{ (b)}$ are the
transformations~\thetag{ 6'} which correspond, respectively, to the
transformations $T_{ (a)}$ and $T_{ (b)}$, then one obviously obtains
the transformation in the $l$ corresponding to the transformation $T_{
(a)} \, T_{ (b)}$ when one executes the two transformations $S_{ (a)}$
and $S_{ (b)}$ one after the other, that is to say: the transformation
$S_{ (a)}\, S_{ (b)}$ corresponds to the transformation $T_{ (a)}\,
T_{ (b)}$. But now, we have: $T_{ (a)}\, T_{ (b)} = T_{ (c)}$ and the
transformation $S_{ (c)}$ in the $l$ corresponds to the transformation
$T_{ (c)}$, whence we must have: $S_{ (a)}\, S_{ (b)} = S_{ (c)}$.

With these words, it is proved that the equations~\thetag{ 6}
and~\thetag{ 6'} really represent an $r$-term group in the variables:
$x_1, \dots, x_n$, $l_1, \dots, l_m$ and to be precise, a group which
is holoedrically isomorphic with the group: $x_i' = f_i ( x, a)$;
indeed, both groups visibly have one and the same parameter group
(cf. Chap.~\ref{kapitel-21}, p.~\pageref{S-402-sq} sq.).

At the same time, it is proved that the equations~\thetag{ 6'} in turn
also represent a group in the variables $l_1, \dots, l_m$, and in
fact, a group \label{S-470} isomorphic with the group: $x_i' = f_i (
x, a)$ (cf. Chap.~\ref{kapitel-21}, p.~\pageref{S-420-sq} sq.), as it
results from the symbolic relations holding simultaneously:
\[
T_{(a)}\,T_{(b)}
=
T_{(c)}
\]
and:
\[
S_{(a)}\,S_{(b)}
=
S_{(c)}.
\]
If one associates to every transformation of the group~\thetag{ 6} the
transformation of the group~\thetag{ 6'} which is determined by it,
then one obtains the two groups~\thetag{ 6} and~\thetag{ 6'} related
to each other in an isomorphic way.

The isomorphism of the two groups~\thetag{ 6} and~\thetag{ 6'} needs
not at all be holoedric, and in certain circumstances, the
group~\thetag{ 6'} can even reduce to the identity transformation,
namely when the group $x' = f( x,a)$ leaves individually untouched
each one of the $\infty^m$ manifolds: $\Omega_k ( x, l) = 0$.

\medskip

We will show that the group which is represented by the joint
equations~\thetag{ 6} and~\thetag{ 6'} is generated by $r$ independent
infinitesimal transformations.

Since the family of the $\infty^m$ manifolds: $\Omega_k ( x, l)$
allows all transformations of the group: $x_i' = f_i ( x, a)$, it
admits in particular each one of the $r$ one-term groups: $X_1f,
\dots, X_rf$, hence according to Proposition~3,
p.~\pageref{Satz-3-S-468}, it also admits each one of the $r$
infinitesimal transformations: $X_1f, \dots, X_rf$. From this, it
results that to every infinitesimal transformation $X_kf$ is
associated a completely determined infinitesimal transformation:
\[
L_kf
=
\sum_{\mu=1}^m\,
\lambda_{k\mu}(l_1,\dots,l_m)\,
\frac{\partial f}{\partial l_\mu}
\]
of such a constitution that the system of equations: $\Omega_k ( x, l)
= 0$ in the $n + m$ variables: $x_1, \dots, x_n$, $l_1, \dots, l_m$
admits the infinitesimal transformation: $X_kf + L_kf$.

Naturally, the system of equations: $\Omega_k ( x, l) = 0$ admits at
the same time every infinitesimal transformation: $e_1\, (X_1f + L_1f)
+ \cdots + e_r\, (X_rf + L_rf)$, and in consequence of that, also
every one-term group: $e_1\, (X_1f + L_1f) + \cdots + e_r\, (X_rf +
L_rf)$. But since the group: $x_i' = f_i ( x, a)$ consists of the
totality of all one-term groups: $e_1\, X_1f + \cdots + e_r\, X_rf$,
then the group represented by the equations~\thetag{ 6} and~\thetag{
6'} must obviously be identical to the totality of all one-term
groups: $\sum\, e_k\, X_kf + \sum\, e_k\, L_kf$, hence it must be
generated by the $r$ independent infinitesimal transformations: $X_kf
+ L_kf$, was what to be shown.

Now, from this, it follows that two arbitrary infinitesimal
transformations amongst the $X_k + L_kf$ must satisfy relations of the
form:
\[
\leftbracket
X_i+L_if,\,\,
X_k+L_kf
\rightbracket
=
\sum_{s=1}^r\,
c_{iks}'\,
(X_sf+L_sf).
\]
By verifying this directly, we obtain a new proof of the fact that all
finite transformations: $x_i' = f_i ( x, a)$, $l_\mu' = \chi_\mu ( l,
a)$ form a group. But at the same time, we realize that $c_{ iks}' =
c_{ iks}$, which is therefore coherent with the fact that our new
group in the $x$ and $l$ possesses the same parameter group as the
given group: $x_i' = f_i ( x, a)$.

The system of equations: $\Omega_k ( x, l) = 0$ admits, simultaneously
with the infinitesimal transformations: $X_1f + L_1f$, \dots, $X_rf +
L_rf$, the transformations:
\[
\aligned
\leftbracket
X_i,\,X_k
\rightbracket
&
+
\leftbracket
L_i,\,L_k
\rightbracket
=
\sum_{s=1}^r\,c_{iks}\,X_sf
+
\leftbracket
L_i,\,L_k
\rightbracket
\\
&
\ \ \ \ \ \ \ \ \ \ \ \ \ 
{\scriptstyle{(i,\,\,k\,=\,1\,\cdots\,r)}},
\endaligned
\]
and therefore also the following ones:
\[
\leftbracket
X_i,\,X_k
\rightbracket
+
\leftbracket
L_i,\,L_k
\rightbracket
-
\sum_{s=1}^r\,c_{iks}\,
(X_sf+L_sf)
=
\leftbracket
L_i,\,L_k
\rightbracket
-
\sum_{s=1}^r\,c_{iks}\,L_sf
\]
in the variables $l_1, \dots, l_m$ alone. But because of the
constitution of the system of equations: $\Omega_k ( x, l) = 0$, this
is possible only when the infinitesimal transformations just written
vanish identically, hence when the relations:
\[
\leftbracket
L_i,\,L_k
\rightbracket
=
\sum_{s=1}^r\,c_{iks}\,L_sf
\]
hold true. From this, it results immediately:
\[
\leftbracket
L_i,\,L_k
\rightbracket
+
\leftbracket
X_i,\,X_k
\rightbracket
=
\sum_{s=1}^r\,c_{iks}\,
(X_sf+L_sf),
\]
whence the said property of the joint equations~\thetag{ 6}
and~\thetag{ 6'} is proved.

According to the above, it goes without saying that the group~\thetag{
6'} is generated by the $r$ infinitesimal transformations: $L_1f,
\dots, L_rf$. At the same time, in what precedes, we have a new proof
of the fact that the equations~\thetag{ 6'} represent a group in the
variables $l_1, \dots, l_m$ which is isomorphic with the group:
$X_1f, \dots, X_rf$. 

At present, the isomorphic relationship between the two
groups~\thetag{ 6} and~\thetag{ 6'} mentioned on p.~\pageref{S-470}
can also be defined by saying that to every infinitesimal
transformation: $e_1\, X_1f + \cdots + e_r\, X_rf$ is associated
the infinitesimal transformation: $e_1\, L_1f + \cdots + 
e_r \, L_rf$ determined by it in the $l$. 

We therefore have the:

\def\theproposition{5}\begin{proposition}
\label{Satz-5-S-472}
If the family of the $\infty^m$ manifolds:
\[
\Omega_k(x_1,\dots,x_n,\,l_1,\dots,l_m)
=
0
\ \ \ \ \ \ \ \ \ \ \ \ \ {\scriptstyle{(k\,=\,1\,\cdots\,n\,-\,q)}}
\]
in the space: $x_1, \dots, x_n$ admits the $r$ independent 
infinitesimal transformations:
\[
X_kf
=
\sum_{i=1}^n\,\xi_{ki}(x_1,\dots,x_n)\,
\frac{\partial f}{\partial x_i}
\ \ \ \ \ \ \ \ \ \ \ \ \ {\scriptstyle{(k\,=\,1\,\cdots\,r)}}
\]
of an $r$-term group having the composition:
\[
\leftbracket
X_i,\,X_k
\rightbracket
=
\sum_{s=1}^r\,c_{iks}\,X_sf
\ \ \ \ \ \ \ \ \ \ \ \ \ 
{\scriptstyle{(i,\,\,k\,=\,1\,\cdots\,r)}},
\]
then to every $X_kf$ there corresponds a completely determined
infinitesimal transformation:
\[
L_kf
=
\sum_{\mu=1}^m\,
\lambda_{k\mu}(l_1,\dots,l_m)\,
\frac{\partial f}{\partial l_\mu}
\]
of such a constitution that the system of equations: $\Omega_k ( x, l)
= 0$ in the $n + m$ variables: $x_1, \dots, x_n$, $l_1, \dots, l_m$
admits the $r$ infinitesimal transformations: $X_kf + L_kf$, so the
$r$ infinitesimal transformations $L_kf$ stand pairwise in the
relationships:
\[
\leftbracket
L_i,\,L_k
\rightbracket
=
\sum_{s=1}^r\,c_{iks}\,L_sf
\ \ \ \ \ \ \ \ \ \ \ \ \ 
{\scriptstyle{(i,\,\,k\,=\,1\,\cdots\,r)}},
\]
that is to say, they generate a group isomorphic with the group:
$X_1f, \dots, X_rf$.
\end{proposition}

Here, the following proposition may yet be expressly stated and proved
independently:

\def\theproposition{6}\begin{proposition}
If the $r$-term group: $X_1f, \dots, X_rf$ of the space $x_1, \dots,
x_n$ contains exactly $p \leqslant r$ independent infinitesimal
transformations which leave invariant the $\infty^m$ manifolds:
\[
\Omega_k(x_1,\dots,x_n,\,l_1,\dots,l_m)
=
0
\ \ \ \ \ \ \ \ \ \ \ \ \ {\scriptstyle{(k\,=\,1\,\cdots\,n\,-\,q)}},
\]
then these $p$ independent infinitesimal transformations generate a
$p$-term subgroup of the group: $X_1f, \dots, X_rf$.
\end{proposition}

The proof is very simple. Let:
\[
\Xi_\pi f
=
\sum_{j=1}^r\,g_{\pi j}\,X_jf
\ \ \ \ \ \ \ \ \ \ \ \ \ {\scriptstyle{(\pi\,=\,1\,\cdots\,p)}}
\]
be such independent infinitesimal transformations of the group: $X_1f,
\dots, X_rf$ which leave invariant the family: $\Omega_k ( x, l) = 0$,
so that every other infinitesimal transformation: $e_1\, X_1f + \cdots
+ e_r\, X_rf$ which does the same can be linearly deduced from
$\Xi_1f, \dots, \Xi_pf$. Then according to Proposition~4,
p.~\pageref{Satz-4-S-469}, the family: $\Omega_k ( x,l) = 0$ also
admits every infinitesimal transformation: $\leftbracket \Xi_\mu, \,
\Xi_\nu \rightbracket$ of the group: $X_1f, \dots, X_rf$, and
consequently, there are relations of the form:
\[
\leftbracket
\Xi_\mu,\,\Xi_\nu
\rightbracket
=
\sum_{\pi=1}^p\,g_{\mu\nu\pi}\,\Xi_\pi f
\ \ \ \ \ \ \ \ \ \ \ \ \ 
{\scriptstyle{(\mu,\,\,\nu\,=\,1\,\cdots\,p)}}
\]
in which the $g_{ \mu \nu \pi}$ denote constants. From this,
it results that $\Xi_1f, \dots, \Xi_pf$ effectively generate
a $p$-term group.\,---

If the family of the $\infty^m$ manifolds:
\[
\Omega_k
(x_1,\dots,x_n,\,l_1,\dots,l_m)
=
0
\ \ \ \ \ \ \ \ \ \ \ \ \ {\scriptstyle{(k\,=\,1\,\cdots\,n\,-\,q)}}
\]
is presented, and if in addition an arbitrary $r$-term group:
$X_1f, \dots, X_rf$ of the space $x_1, \dots, x_n$ is also
presented, one can ask how many independent infinitesimal
transformations of the group: $X_1f, \dots, X_rf$ 
does the family presented admit. We now overview how this
question can be answered.

Indeed, if the family: $\Omega_k ( x, l) = 0$ is supposed 
to admit an infinitesimal transformation of the form:
$e_1\, X_1f + \cdots + e_r\, X_rf$, then it must be possible
to indicate an infinitesimal transformation: $Lf$ in the
$l$ alone such that the system of equations: $\Omega_k ( x, l) = 0$
in the $n + m$ variables: $x_1, \dots, x_n$, $l_1, \dots, l_m$
admits the infinitesimal transformation: $\sum\, e_k\, X_kf + 
Lf$. Hence, if we imagine that the equations: $\Omega_k ( x, l) = 0$
are resolved with respect to $n - q$ of the $x$:
\[
x_{q+k}
=
\psi_{q+k}
(x_1,\dots,x_q,\,l_1,\dots,l_m)
\ \ \ \ \ \ \ \ \ \ \ \ \ {\scriptstyle{(k\,=\,1\,\cdots\,n\,-\,q)}},
\]
and if, as on p.~\pageref{S-468}, we denote the substitution: 
$x_{ q+1} = \psi_{ q+1}$, \dots, $x_n = \psi_n$ by the sign
$[ \,\,\, ]$, then we only have to determine the constants
$e_1, \dots, e_r$ and the functions: $\lambda_1 ( l), \dots, 
\lambda_m ( l)$ in the most general way so that the $n - q$ equations:
\[
\aligned
\sum_{\mu=1}^m\,\lambda_\mu(l)\,
\frac{\partial\psi_{q+k}}{\partial l_\mu}
&
=
\sum_{\sigma=1}^r\,e_\sigma\,
\bigg\{
\big[\xi_{\sigma,\,q+k}\big]
-
\sum_{j=1}^q\,
\big[\xi_{\sigma j}\big]\,
\frac{\partial\psi_{q+k}}{\partial x_j}
\bigg\}
\\
&
\ \ \ \ \ \ {\scriptstyle{(k\,=\,1\,\cdots\,n\,-\,q)}}
\endaligned
\]
are satisfied identically, independently of the values of the
variables $x_1, \dots, x_q$. In this way, we find the most general
infinitesimal transformation: $e_1\, X_1f + \cdots + 
e_r\, X_rf$ which leaves invariant the family: $\Omega_k ( x, l) = 0$.

\sectionengellie{\S\,\,\,112.}

When the family of the $\infty^m$ manifolds:
\def\theequation{1}\begin{equation}
\Omega_k(x_1,\dots,x_n,\,l_1,\dots,l_m)
=
0
\ \ \ \ \ \ \ \ \ \ \ \ \ {\scriptstyle{(k\,=\,1\,\cdots\,n\,-\,q)}}
\end{equation}
admits the transformation: $x_i' = f_i ( x_1, \dots, x_n)$, then there
is, as we have seen on p.~\pageref{S-463-sq} sq., a completely
determined transformation: $l_\mu' = \chi_\mu ( l_1, \dots, l_m)$ of
such a constitution that the system of equations: $\Omega_k ( x, l) =
0$ in the $n + m$ variables $x$, $l$ admits the transformation:
\[
\left\{
\aligned
x_i'
&
=
f_i(x_1,\dots,x_n)
\ \ \ \ \ \ \ \ \ \ \ \ \ {\scriptstyle{(i\,=\,1\,\cdots\,n)}}
\\
l_\mu'
&
=
\chi_\mu(l_1,\dots,l_m)
\ \ \ \ \ \ \ \ \ \ \ \ \ {\scriptstyle{(\mu\,=\,1\,\cdots\,m)}}.
\endaligned\right.
\]

Already on p.~\pageref{S-463}, we observed that the equations: $l_\mu'
= \chi_\mu ( l)$ determine the parameters of the manifold of the
invariant family~\thetag{ 1} into which the manifold with the
parameters $l_1, \dots, l_m$ is transferred by the transformation:
$x_i' = f_i ( x)$. Hence, when we interpret $l_1, \dots, l_m$
virtually as coordinates of individual manifolds of the
family~\thetag{ 1}, the equations: $l_\mu' = \chi_\mu ( l)$ indicate
the law according to which the manifolds of our invariant family are
permuted with each other by the transformation: $x_i' = f_i ( x)$.

For instance, if a special system of values: $l_1^0, \dots, l_m^0$
admits the transformation:
\[
l_\mu'
=
\chi_\mu(l_1,\dots,l_m)
\ \ \ \ \ \ \ \ \ \ \ \ \ {\scriptstyle{(\mu\,=\,1\,\cdots\,m)}},
\]
then at the same time, the manifold:
\[
\Omega_k
(x_1,\dots,x_n,\,l_1^0,\dots,l_m^0)
=
0
\ \ \ \ \ \ \ \ \ \ \ \ \ {\scriptstyle{(k\,=\,1\,\cdots\,n\,-\,q)}}
\]
admits the transformation:
\[
x_i'
=
f_i(x_1,\dots,x_n)
\ \ \ \ \ \ \ \ \ \ \ \ \ {\scriptstyle{(i\,=\,1\,\cdots\,n)}}.
\]

In fact, the manifold: $\Omega_k ( x, l^0) = 0$ is transferred, by the
execution of the transformation: $x_i' = f_i ( x)$, to the new
manifold:
\[
\Omega_k\big(
x_1',\dots,x_n',\,
\chi_1(l^0),\,\dots,\,\chi_m(l^0)\big)
=
0
\ \ \ \ \ \ \ \ \ \ \ \ \ {\scriptstyle{(k\,=\,1\,\cdots\,n\,-\,q)}}\,;
\]
but under the assumption made, we have:
\[
\chi_\mu(l_1^0,\dots,l_m^0)
=
l_\mu^0
\ \ \ \ \ \ \ \ \ \ \ \ \ {\scriptstyle{(\mu\,=\,1\,\cdots\,m)}},
\]
hence the new manifold coincides with the old one and the manifold:
$\Omega_k ( x, l^0) = 0$ really remains invariant.

However by contrast, the converse does not always hold true. If an
arbitrary special manifold:
\[
\Omega_k(x_1,\dots,x_n,\,l_1^0,\dots,l_m^0)
=
0
\ \ \ \ \ \ \ \ \ \ \ \ \ {\scriptstyle{(k\,=\,1\,\cdots\,n\,-\,q)}}
\]
of the family~\thetag{ 1} admits the transformation: $x_i' = f_i(x)$,
then it does not follow from that with necessity that the system of
values: $l_1^0, \dots, l_m^0$ admits the transformation: $l_\mu' =
\chi_\mu ( l)$. Indeed, it is thinkable that, to the system of values:
$l_1^0, \dots, l_m^0$, there is associated a continuous series of
systems of values: $l_1, \dots, l_m$ which provide the same manifold
as the system of values: $l_1^0, \dots, l_m^0$, when inserted
in~\thetag{ 1}; in this case, from the invariance of the manifold:
$\Omega_k ( x, l^0) = 0$, it only follows that the individual systems
of values: $l_1, \dots, l_m$ of the series just defined are permuted
with each other by the transformation: $l_\mu' = \chi_\mu (l)$, but
not that the system of values: $l_1^0, \dots, l_m^0$ remains invariant
by the transformation in question. \emphasis{Nevertheless, if the
$l_k$ do not have special values but general values, then the
manifold: $\Omega_k ( x, l) = 0$ admits the transformation: $x_i' =
f_i (x, a)$ only when the system of values $l_k$ allows the
corresponding transformation: $l_\mu' = \chi_\mu ( l)$, and also,
always in this case too.}

Now, we consider the general case where the family of the $\infty^m$
manifolds:
\def\theequation{1}\begin{equation}
\Omega_k(x_1,\dots,x_n,\,l_1,\dots,l_m)
=
0
\ \ \ \ \ \ \ \ \ \ \ \ \ {\scriptstyle{(k\,=\,1\,\cdots\,n\,-\,q)}}
\end{equation}
admits the $r$-term group:
\[
x_i'
=
f_i(x_1,\dots,x_n,\,a_1,\dots,a_r)
\ \ \ \ \ \ \ \ \ \ \ \ \ {\scriptstyle{(i\,=\,1\,\cdots\,n)}}
\]
with the $r$ independent infinitesimal transformations: $X_1f, \dots, 
X_rf$. Let the group in the $l$ which corresponds to the group: $x_i'
= f_i ( x,a)$ have the form:
\[
l_\mu'
=
\chi_\mu(l_1,\dots,l_m,\,a_1,\dots,a_r)
\ \ \ \ \ \ \ \ \ \ \ \ \ {\scriptstyle{(\mu\,=\,1\,\cdots\,m)}},
\]
and let it be generated by the $r$ infinitesimal transformations:
$L_1f, \dots, L_rf$ which in turn correspond naturally to the
infinitesimal transformations: $X_1f, \dots, X_rf$, respectively.

At first, the question is how one decides whether the infinitesimal
transformation: $e_1^0\, X_1f + \cdots + e_r^0 \, X_rf$ leaves
invariant a determined manifold contained in the family~\thetag{ 1}:
\def\theequation{7}\begin{equation}
\Omega_k(x_1,\dots,x_n,\,
l_1^0,\dots,l_m^0)
=
0
\ \ \ \ \ \ \ \ \ \ \ \ \ {\scriptstyle{(k\,=\,1\,\cdots\,n\,-\,q)}}.
\end{equation}

It is easy to see that the manifold~\thetag{ 7} admits in any case the
infinitesimal transformation: $e_1^0 \, X_1f + \cdots + e_r^0 \, X_rf$
when the system of values: $l_1^0, \dots, l_m^0$ admits the
infinitesimal transformation: $e_1^0 \, L_1f + \cdots + e_r^0\, L_rf$.

In fact, the system of equations~\thetag{ 1} in the $m + n$ variables:
$x_1, \dots, x_n$, $l_1, \dots, l_m$ admits, under the assumptions
made, the infinitesimal transformation: $\sum\, e_j^0 \, (X_j f +
L_jf)$; the $n - q$ expressions:
\[
\sum_{j=1}^r\,e_j^0\,
\big(
X_j\,\Omega_k+L_j\,\Omega_k
\big)
=
\sum_{j=1}^r\,e_j^0\,
\bigg\{
\sum_{i=1}^n\,\xi_{ji}(x)\,
\frac{\partial\Omega_k}{\partial x_i}
+
\sum_{\mu=1}^m\,\lambda_{j\mu}(l)\,
\frac{\partial\Omega_k}{\partial l_\mu}
\bigg\}
\]
therefore vanish all by virtue of: $\Omega_1 ( x, l) = 0$, \dots,
$\Omega_{ n -q} ( x, l) = 0$. This also holds true when we set: $l_1 =
l_1^0$, \dots, $l_m = l_m^0$; but now, the system of values $l_1^0,
\dots, l_m^0$ admits the infinitesimal transformation: $\sum\, e_j^0\,
L_jf$ and hence, the $m$ expressions:
\[
e_1^0\,\lambda_{1\mu}(l^0)
+\cdots+
e_r^0\,\lambda_{r\mu}(l^0)
\ \ \ \ \ \ \ \ \ \ \ \ \ {\scriptstyle{(\mu\,=\,1\,\cdots\,m)}}
\]
vanish all. Consequently, the $n - q$ expressions:
\[
\sum_{i=1}^n\,
\bigg\{
\sum_{j=1}^r\,e_j^0\,\xi_{ji}(x)
\bigg\}\,
\frac{\partial\Omega_k(x,l^0)}{\partial x_i}
\ \ \ \ \ \ \ \ \ \ \ \ \ {\scriptstyle{(k\,=\,1\,\cdots\,n\,-\,q)}}
\]
vanish all by virtue of~\thetag{ 7}, that is to say, the
manifold~\thetag{ 7} really admits the infinitesimal transformation:
$\sum\, e_j^0\, X_jf$.

Nevertheless, the sufficient criterion found here with these words is
not necessary.

Indeed, if our family of manifolds $\Omega_k ( x, l) = 0$ remains
invariant by the group: $X_1f, \dots, X_rf$, and if at the same time,
the special manifold: $\Omega_k ( x, l^0) = 0$ is supposed to admit
the infinitesimal transformation: $e_1^0\, X_1f + \cdots + e_r^0 \,
X_rf$, then for that, it is only necessary that the $n - q$
expressions:
\[
\sum_{\mu=1}^m\,\sum_{j=1}^r\,
e_j^0\,\lambda_{j\mu}(l^0)\,
\frac{\partial\Omega_k(x,l^0)}{\partial l_\mu}
\]
be equal to zero by virtue of the system of equations: $\Omega_k ( x,
l^0) = 0$; but it is not at all necessary that the $r$ expressions:
$\sum\, e_j^0\, \lambda_{ j\mu} ( l^0)$ themselves vanish.

\smallskip

Thus, if one wants to know whether every manifold: $\Omega_k ( x, l) 
= 0$ of the invariant family admits one or several infinitesimal
transformations: $\sum\, e_k\, X_kf$, and if in addition, one wants
for every manifold $l_k$ to find the concerned infinitesimal
transformation, then one has to determined the $e_k$ in the most
general way as functions of the $l$ so that the equations:
\def\theequation{8}\begin{equation}
\sum_{\mu=1}^m\,\sum_{j=1}^r\,
e_j\,\lambda_{j\mu}(l)\,
\frac{\partial\Omega(x,l)}{\partial l_\mu}
=
0
\end{equation}
are a consequence of the system of equations: $\Omega_k ( x, l) = 0$. 
But since, according to an assumption made earlier on, the $l$ are
essential parameters, the system of equations: $\Omega_k ( x, l) = 0$
in the $n + m$ variables $x$, $l$ can admit no not identically
vanishing transformation:
\[
\sum_{\mu=1}^m\,\sum_{j=1}^r\,
e_j(l)\,\lambda_{j\mu}(l)\,
\frac{\partial f}{\partial l_\mu},
\]
whence it follows with necessity that:
\[
e_1\,\lambda_{1\mu}(l)
+\cdots+
e_r\,\lambda_{r\mu}(l)
=
0
\ \ \ \ \ \ \ \ \ \ \ \ \ 
{\scriptstyle{(\mu\,=\,1,\,\,2\,\cdots\,m)}}.
\]

We therefore get the following:

\plainstatement{If: $\sum\, e_k^0 \, L_kf$ is the most general
infinitesimal transformation contained in the group: $L_1f, \dots,
L_rf$ which leaves invariant the system of values: $l_1^0, \dots,
l_m^0$ \terminology{located in general position}, then: $\sum\,
e_k^0\, X_kf$ is the most general infinitesimal transformation
contained in the group: $X_1f, \dots, X_rf$ which leaves invariant the
manifold~\thetag{ 7} located in general position.}

We can assume that amongst the $r$ infinitesimal transformations:
$L_1f, \dots, L_rf$, exactly $m - p$, say: $L_1f, \dots, L_{ m-p} f$,
are linked together by no linear relation of the form:
\[
\alpha_1(l_1,\dots,l_m)\,L_1f
+\cdots+
\alpha_{m-p}(l_1,\dots,l_m)\,L_{m-p}f
=
0,
\]
while by contrast: $L_{m-p+1}f, \dots, L_rf$ can be expressed linearly
in terms of $L_1f, \dots, L_{ m-p}f$:
\[
L_{m-p+j}f
\equiv
\sum_{\mu=1}^{m-p}\,
\alpha_{j\mu}(l_1,\dots,l_m)\,L_\mu f
\ \ \ \ \ \ \ \ \ \ \ \ \ 
{\scriptstyle{(j\,=\,1\,\cdots\,r\,-\,m\,+\,p)}}.
\]
We can always insure that this assumption holds, even when the
infinitesimal transformations: $L_1f, \dots, L_rf$ are
not mutually independent, which can very well occur.

It is now clear that the $r - m + p$ infinitesimal transformations:
\def\theequation{9}\begin{equation}
L_{m-p+j}f
-
\sum_{\mu=1}^{m-p}\,
\alpha_{j\mu}(l_1^0,\dots,l_m^0)\,
L_\mu f
\ \ \ \ \ \ \ \ \ \ \ \ \ 
{\scriptstyle{(j\,=\,1\,\cdots\,r\,-\,m\,+\,p)}}
\end{equation}
leave invariant the system of values: $l_1^0, \dots, l_m^0$; moreover,
because $l_1^0, \dots, l_m^0$ is a system of values in general
position, one realizes immediately that every infinitesimal
transformation: $\sum\, e_k\, L_kf$ by which the system of values
remains invariant can be linearly deduced from the the $r - m + p$
transformations~\thetag{ 9}. Consequently, the most general
infinitesimal transformation: $\sum\, e_k\, X_kf$ which leaves
invariant the manifold~\thetag{ 7} located in general position can be
linearly deduced from the $r - m + p$ transformations:
\def\theequation{10}\begin{equation}
X_{m-p+j}f
-
\sum_{\mu=1}^{m-p}\,
\alpha_{j\mu}(l_1^0,\dots,l_m^0)\,
X_\mu f
\ \ \ \ \ \ \ \ \ \ \ \ \ 
{\scriptstyle{(j\,=\,1\,\cdots\,r\,-\,m\,+\,p)}}.
\end{equation}

Of course, the infinitesimal transformations~\thetag{ 10}
are mutually independent and they generate an
$(r - m + p)$-term group, namely the most general subgroup
contained in the group: $X_1f, \dots, X_rf$ by which the
manifold~\thetag{ 7} remains invariant.

In particular, if the group: $L_1f, \dots, L_rf$ is transitive, 
then the entire number $p$ has the value zero; so in this
case, every manifold of the family~\thetag{ 1} located
in general position admits exactly $r - m$ independent
infinitesimal transformations of the group: $X_1f, \dots, X_rf$. 

We therefore have the

\def\theproposition{7}\begin{proposition}
\label{Satz-7-S-478}
If the family of the $\infty^m$ manifolds: 
\[
\Omega_k(x_1,\dots,x_n,\,l_1,\dots,l_m)
=
0
\ \ \ \ \ \ \ \ \ \ \ \ \ {\scriptstyle{(k\,=\,1\,\cdots\,n\,-\,q)}}
\]
admits the $r$-term group: $X_1f, \dots, X_rf$, if moreover
$L_1f, \dots, L_rf$ are the infinitesimal transformations
in the $l$ alone which correspond to the $X_jf$ and if lastly
$L_1f, \dots, L_rf$ are linked together by exactly $r - m + p$
independent relations of the form: $\sum\, \beta_j (l)\, L_jf = 0$, 
then the group: $X_1f, \dots, X_rf$ contains exactly $r - m + p$
independent infinitesimal transformations which leave invariant
a manifold $l_1^0, \dots, l_m^0$ in general position.
These transformations generate an $(r - m + p)$-term group.
If the group: $L_1f, \dots, L_rf$ in the $m$ variables $l_1, \dots, 
l_m$ is transitive, then every manifold of the family $\Omega_k
(x, l) = 0$ located in general position admits exactly $r - m$
independent infinitesimal transformations of the group: $X_1f, 
\dots, X_rf$ and these infinitesimal transformations generate
an $(r - m)$-term subgroup.
\end{proposition}

\label{S-479-sq}
We yet add here the obvious remark that the group: $L_1f, \dots, 
L_rf$ in $l_1, \dots, l_m$ is transitive if and only if
every manifold of the family $\Omega_k ( x, l) = 0$ located
in general position can be transferred, by means of at least
one transformation of the group: $X_1f, \dots, X_rf$, to 
every other manifold.

\sectionengellie{\S\,\,\,113.}

\label{S-479}
In the $n$-times extended space $x_1, \dots, x_n$, let an $r$-term
group: $X_1f, \dots, X_rf$ be presented, and in addition, 
let an arbitrary manifold, which we want to denote by $M$, 
be given. We assume that in the equations of $M$,
no arbitrary parameters of any kind appear.

If all $\infty^r$ transformations of the group: $X_1f, \dots, 
X_rf$ are executed on the manifold $M$, then this manifold is
transferred to a series of new manifolds. We will prove that the
totality of all these manifolds remains invariant by the
group: $X_1f, \dots, X_rf$, and that it forms a family invariant
by the group.

Let $M'$ be an arbitrary manifold which belongs to the totality
just said, and let $T_1$ be a transformation of the group:
$X_1f, \dots, X_rf$ which transfers $M$ to $M'$, whence
there is the symbolic equation:
\[
(M)\,T_1
=
(M').
\]
Now, if $T$ is an arbitrary transformation of the group, 
we have:
\[
(M')\,T
=
(M')\,T_1\,T
=
(M)\,T_2,
\]
where the transformation $T_2$ again belongs to the group;
consequently, the manifold $M'$ is transferred by the transformation
$T$ to another manifold of the totality in question. Since this holds
for every manifold $M'$ of the totality, we see that the manifolds
belonging to the totality are permuted with each other by $T$, and
hence actually, by all transformations of the group: $X_1f, \dots,
X_rf$, hence we see that the totality of manifolds defined above
effectively forms a family invariant by the group: $X_1f, \dots,
X_rf$.

It is easy to see that every manifold of this invariant family can be
transferred, by means of at least one transformation of the group, to
every other manifold of the family. Indeed, if:
\[
(M')
=
(M)\,T_1,
\ \ \ \ \ \ \ \ \ \ \
(M'')
=
(M)\,T_2,
\]
then $(M) = (M')\, T_1^{ -1}$, whence:
\[
(M'')
=
(M')\,T_1^{-1}\,T_2,
\]
whence the claim is proved. At the same time, it results from this
that one also obtains the discussed family of manifolds when one
executes all $\infty^r$ transformations of the group on an arbitrary
manifold of the family.

We therefore have the

\def\theproposition{8}\begin{proposition}
If one executes all $\infty^r$ transformations of an $r$-term group:
$X_1f, \dots, X_rf$ of the $R_n$ on a given manifold of this space,
then the totality of all positions that the manifold takes on the
occasion forms a family of manifolds invariant by the group. Every
manifold of this family can be transferred, by means of at least one
transformation of the group, to any other manifold of the family.
The family can be derived from each one of its manifolds in 
the same way as it is deduced from the initially given manifold.
\end{proposition}

\renewcommand{\thefootnote}{\fnsymbol{footnote}}
The manifold $M$ will admit a certain number, that we assume to 
be exactly equal to $r -m$, of infinitesimal transformations
of the group: $X_1f, \dots, X_rf$.\footnote[1]{\,
The totality of all \emphasis{finite} transformations of the group:
$X_1f, \dots, X_rf$ which leave invariant a manifold $M$ always forms
a subgroup (Theorem~32, p.~\pageref{Theorem-32-S-208}), but of course,
this group needs not be a finite continuous group. Nevertheless, 
in the following developments of the text, when we make the
implicit assumption that this subgroup is generated by
infinitesimal transformations, that is not to be interpreted as
an essential restriction, because we can indeed suitably narrow down 
the region $(\!( a)\!)$ introduced on p.~\pageref{S-16}.
}
\renewcommand{\thefootnote}{\arabic{footnote}}
These transformations then generate an $(r - m)$-term subgroup
(cf. Theorem~31, p.~\pageref{Theorem-31-S-207}).

Now, let $S$ be the general symbol of a transformation of this
subgroup, so: $(M)\, S = (M)$; moreover, let $T_1$ be an arbitrary
transformation of the group: $X_1f, \dots, X_rf$ and let $M$ be
transferred to the new position $M'$ by $T_1$: $(M') = (M)\, T_1$.
Then it is easy to indicate \emphasis{all} transformations 
${\sf T}$ of the group $X_kf$ which transfer $M$ to $M'$. 

Indeed, one has: 
\[
(M)\,{\sf T}
=
(M')
=
(M)\,T_1,
\]
whence:
\[
(M)\,{\sf T}\,T_1^{-1}
=
(M),
\]
hence ${\sf T}\, T_1^{ -1}$ is a transformation $S$ and $S\, T_1$ is
the general symbol of all transformations of the group $X_kf$ which
transfer $M$ to $M'$. But there are as many transformations $S\, T_1$
as there are different transformations $S$, that is to say, $\infty^{
r-m}$.

One finds in a similar way all transformations $\mathfrak{ S}$ of our
group which leave $M'$ invariant. From $(M')\, \mathfrak{ S} = (M')$
and $(M)\, T_1 = (M')$, one obtains:
\[
(M)\,T_1\,\mathfrak{S}\,T_1^{-1}
=
(M),
\]
whence:
\[
T_1\,\mathfrak{S}\,T_1^{-1}
=
S,
\ \ \ \ \ \ \ \ \ \
\mathfrak{S}
=
T_1^{-1}\,S\,T_1.
\]
Likewise, there are $\infty^{ r-m}$ such different transformations and
their totality forms an $(r - m)$-term subgroup which is conjugate to
the group $S$ inside the group $X_kf$. 

We summarize this result as follows:

\def\theproposition{9}\begin{proposition}
If a manifold $M$ of the $R_n$ admits exactly $r - m$ independent
infinitesimal transformations of the $r$-term group: $X_1f, \dots,
X_rf$, or briefly $G_r$, if moreover $S$ is the general symbol of the
$\infty^{ r-m}$ finite transformations of the $(r - m)$-term subgroup
which is generated by these $r - m$ infinitesimal transformations, and
lastly, if $T$ is an arbitrary transformation of the $G_r$: $X_1f,
\dots, X_rf$ and if $M$ takes the new position $M'$ after the
execution of $T$, then the $G_r$ contains exactly $\infty^{ r-m}$
different transformations which likewise transfer $M$ to $M'$ and the
general symbol of these transformations is: $S\, T$; in addition, the
$G_r$ contains exactly $\infty^{ r - m}$ transformations which leave
$M'$ invariant, these transformations have $T^{ -1} \, S \, T$ for a
general symbol and they form an $(r - m)$-term subgroup which is
conjugate to the group of the $S$ inside the $G_r$.
\end{proposition}

If we imagine that the $\infty^r$ transformations of the group: $x_i'
= f_i ( x,a)$ are executed on the equations of the manifold $M$, then
we obtain the analytic expression of the discussed invariant family of
manifolds. Formally, this expression contains the $r$ parameters $a_1,
\dots, a_r$, but these parameters need not be all essential. We now
want to determine the number $m'$ of the essential parameters amongst
the parameters $a_1, \dots, a_r$.

Our invariant family consists of $\infty^{ m'}$ different manifolds
and each one of these manifolds can be transferred to every other
manifold by means of a transformation of the group: $X_1f, \dots,
X_rf$. According to Proposition~7, p.~\pageref{Satz-7-S-478}, each
individual manifold amongst the $\infty^{ m'}$ manifolds then admits
exactly $r - m'$ independent infinitesimal transformations of the
$G_r$; but from what precedes, we know that under the assumptions
made, each one of these manifolds admits exactly $\infty^{ r-m}$
finite transformations of the $G_r$, and consequently, we have $m' =
m$ and amongst the $r$ parameters: $a_1, \dots, a_r$, there are
exactly $m$ that are essential.

We therefore have the

\def\theproposition{10}\begin{proposition}
If a manifold of the $n$-times extended space $R_n$ admits exactly $r
- m$ independent infinitesimal transformations of an $r$-term group:
$X_1f, \dots, X_rf$ of this space, then this manifold takes exactly
$\infty^m$ different positions by the $\infty^r$ transformations of
this group.
\end{proposition}
\label{S-482}

\sectionengellie{\S\,\,\,114.}

In the equations of our invariant family of manifolds, the parameters
need not, as said, be all essential, but we can always imagine that
$m \leqslant r$ functions $l_1, \dots, l_m$ of the $a$ are introduced
as new parameters so that the equations of our family are given
the form:
\[
\Omega_k(x_1,\dots,x_n,\,l_1,\dots,l_m)
=
0
\ \ \ \ \ \ \ \ \ \ \ \ \ 
{\scriptstyle{(k\,=\,1\,\cdots\,n\,-\,q)}},
\]
where now $l_1, \dots, l_m$ are essential parameters.

The family $\Omega_k = 0$ remains invariant by the group: $X_1f, \dots,
X_rf$ while its individual manifolds are permuted with each other.
The way how the manifolds are permuted is indicated by 
the group: $L_1f, \dots, L_rf$ in the $l$ alone which, as shown
earlier on, is completely determined by the group: $X_1f, \dots, 
X_rf$. 

The group $L_kf$ in the variables $l_1, \dots, l_m$ is isomorphic
with the group $X_kf$, hence it has at most $r$ essential parameters;
on the other hand it is certainly, under the assumption made, 
transitive (cf. p.~\pageref{S-479-sq} sq.), hence 
it has at least $m$ essential parameters. For us, it only remains
to indicate a simple criterion for determining
how many essential parameters the group $L_kf$ really has. 

Let the group: $L_1f, \dots, L_rf$ be $\rho$-term, where $0 \leqslant
\rho \leqslant r$. Because it is meroedrically isomorphic with the
$G_r$: $X_1f, \dots, X_rf$, then there must exist in the $G_r$
an $(r - \rho)$-term invariant subgroup which corresponds to the
identity transformation in the group $L_kf$ (cf. Theorem~54, 
p.~\pageref{Theorem-54-S-301}). This $(r - \rho)$-term invariant
subgroup of the $G_r$ then leaves individually fixed each one
of the $\infty^m$ manifolds: $\Omega_k ( x, l) = 0$, hence it is
contained in the $(r - m)$-term subgroup $g_{ r-m}$ of the $G_r$
which leaves invariant the manifold $M$ discussed above and at
the same time, it is contained in all $(r - m)$-term subgroups
contained in the $G_r$ which are conjugate, inside the $G_r$, 
to the $g_{ r - m}$ just mentioned.

Conversely, if this $g_{ r - m}$ contains an $(r - \rho)$-term
subgroup which is invariant in the $G_r$, then this subgroup is at the
same time contained in all subgroups of the $G_r$ that are conjugate
to the $g_{ r-m}$, hence it leaves untouched every individual
manifold: $\Omega_k ( x, l) = 0$, and to it, there corresponds the
identity transformation in the group: $L_1f, \dots, L_rf$.

In order to decide how many parameters the group: $L_1f, \dots, L_rf$
contains, we therefore have only to look up at the largest group
contained in the $g_{ r-m}$ which is invariant in the $G_r$. When the
group in question is exactly $(r - \rho)$-term ($\rho \geqslant m$),
then the group: $L_1f, \dots, L_rf$ is exactly $\rho$-term.

We do not want to state this result as a specific proposition, but
instead, we want to recapitulate all the results of the previous two
paragraphs in a theorem.

\renewcommand{\thefootnote}{\fnsymbol{footnote}}
\def\thetheorem{85}\begin{theorem}
\label{Theorem-85-S-483}
If one has an $r$-term group: $X_1f, \dots, X_rf$, or briefly $G_r$,
of the space $x_1, \dots, x_n$ and if one has an arbitrary manifold
$M$ which allows exactly $r - m$ independent infinitesimal
transformations of the $G_r$ and hence which also admits the $(r -
m)$-term subgroup $g_{ r-m}$ generated by these infinitesimal
transformations, then through the $\infty^r$ transformations of the
$G_r$, $M$ takes in total $\infty^m$ different positions the totality
of which remains invariant relatively to the group $G_r$. If one marks
the individual positions of $M$ by means of $m$ parameters: $l_1,
\dots, l_m$, then one obtains a certain group in $l_1, \dots, l_m$:
\[
l_\mu'
=
\chi_\mu(l_1,\dots,l_m;\,a_1,\dots,a_r)
\ \ \ \ \ \ \ \ \ \ \ \ \ {\scriptstyle{(\mu\,=\,1\,\cdots\,m)}}
\]
which indicates in which way the individual positions of $M$ are
permuted with each other by the group: $X_1f, \dots, X_rf$. This
group in the $l$ is transitive and isomorphic with the group: $X_1f,
\dots, X_rf$. If the largest subgroup contained in the $g_{ r-m}$
which is invariant in the $G_r$ is exactly $(r - \rho)$-term, then the
group: $l_\mu' = \chi_\mu ( l, a)$ has exactly $\rho$ essential
parameters. In particular, if the $G_r$ is simple, then the group in
the $l$ is always $r$-term and holoedrically isomorphic to the $G_r$,
with the only exception of the case $m = 0$, in which the group:
$l_\mu ' = \chi_\mu ( l, a)$ consists only of the identity
transformation.\footnote[1]{\,
\name{Lie}, Archiv for Mathematik og Naturv., Vol. 10, Christiania
1885.
}
\end{theorem}
\renewcommand{\thefootnote}{\arabic{footnote}}

The number $r - \rho$ mentioned in the theorem may have each one of
the values: $0$, $1$, \dots, $r - m$; if $r - \rho = 0$, then the $(r
- \rho)$-term subgroup of the $g_{ r-m}$ consists of the identity
transformation, hence the group: $l_\mu ' = \chi_\mu ( l, a)$ is
holoedrically isomorphic to the $G_r$; if $r - \rho = r - m$, then the
$g_{ r-m}$ itself is invariant in the $G_r$ and the group: $l_\mu' =
\chi_\mu ( l, a)$ is only $m$-term.

\sectionengellie{\S\,\,\,115.}

In the last but one paragraph we gave a method to find the families of
manifolds which remain invariant by a given $r$-term group: $X_1f,
\dots, X_rf$. The invariant families that we obtained in this way were
distinguished by the fact that every manifold of a family of this sort
could be transferred, by means of at least one transformation of the
group: $X_1f, \dots, X_rf$, to every other manifold of the family.

At present, we want to generalize the discussed method so that
it produces \emphasis{all} families of manifolds invariant by the
group $X_1f, \dots, X_rf$. 

For this, we are conducted to the obvious observation that
the considerations of the pp.~\pageref{S-479}--\pageref{S-482}
also remain valid when the equations of the manifold $M$ contain
arbitrary parameters, that is to say: when in place
of an individual manifold $M$ we use directly a complete family
of manifolds. Thus, we can also proceed in the following way in order
to find invariant families of manifolds:

We take an arbitrary family:
\def\theequation{11}\begin{equation}
V_k(x_1,\dots,x_n,\,u_1,\dots,u_h)
=
0
\ \ \ \ \ \ \ \ \ \ \ \ \ 
{\scriptstyle{(k\,=\,1\,\cdots\,n\,-\,q)}}
\end{equation}
of $\infty^h$ manifolds and we execute on it all $\infty^r$
transformations of the group: $X_1f, \dots, X_rf$; the totality
of all manifolds that we obtain in this way always forms 
a family invariant by the group: $X_1f, \dots, X_rf$. 

\emphasis{It is clear that we obtain all families invariant
by the group: $X_1f, \dots, X_rf$ when we choose the family~\thetag{
11} in all possible ways}. Indeed, if an arbitrary family invariant
by the group is presented, say the family:
\def\theequation{1}\begin{equation}
\Omega_k(x_1,\dots,x_n,\,l_1,\dots,l_m)
=
0
\ \ \ \ \ \ \ \ \ \ \ \ \ 
{\scriptstyle{(k\,=\,1\,\cdots\,n\,-\,q)}},
\end{equation}
then this family can in any case be obtained when
we just choose, as the family~\thetag{ 11}, the family~\thetag{ 1}.
Besides, one can easily indicate infinitely many other 
families~\thetag{ 11} out of which precisely the family~\thetag{ 1} is 
obtained, but we do not want to spend time on this. 

But it yet remains to answer a question. 

If the family~\thetag{ 11} is presented and if all the transformations
of the group: $X_1f, \dots, X_rf$ are executed on it, then the
equations of the invariant family which comes into existence in 
this way obviously have the form:
\def\theequation{12}\begin{equation}
W_k(x_1,\dots,x_n,\,u_1,\dots,u_h,\,a_1,\dots,a_r)
=
0
\ \ \ \ \ \ \ \ \ \ \ \ \ 
{\scriptstyle{(k\,=\,1\,\cdots\,n\,-\,q)}},
\end{equation}
and therefore, they formally contain $h + r$ arbitrary parameters, 
namely: $u_1, \dots, u_h$, $a_1, \dots, a_r$. How many parameters 
amongst these parameters are essential?

Every generally positioned manifold of the family~\thetag{ 11} takes
by the group: $X_1f, \dots, X_rf$ a certain number, say $\infty^p$, of
different positions; of these $\infty^p$ positions, there is a certain
number, say $\infty^o$, of different positions, which again belong to
the family~\thetag{ 11}. In this manner, the complete family~\thetag{
11} is decomposed in $\infty^{ h - o} $ different subfamilies of
$\infty^o$ manifolds, in such a way that every manifold of the
family~\thetag{ 11} can always be transferred to every manifold which
belongs to one and the same subfamily by means of at least one
transformation of the group: $X_1f, \dots, X_rf$, and such that every
manifold of the family~\thetag{ 11}, as soon as it remains inside the
family~\thetag{ 11} by a transformation of the group: $X_1f, \dots,
X_rf$, remains at the same time in the subfamily to which it belongs.

Now, if we imagine that all transformations of the group: $X_1f,
\dots, X_rf$ are executed on two arbitrary manifolds which belong to
the \emphasis{same} subfamily, we obviously obtain in the two cases
the same family of $\infty^p$ manifolds; on the other hand, if we
imagine that all transformations of the group are executed on two
manifolds of the family~\thetag{ 11} which belong to
\emphasis{distinct} subfamilies, we obtain two distinct families of
$\infty^p$ manifolds which have absolutely no manifold in common.

Hence, when we choose on each one of the $\infty^{ h - o}$ subfamilies
of the family~\thetag{ 11} a manifold and when we execute
all transformations of the group: $X_1f, \dots, X_rf$ on 
the so obtained $\infty^{ h - o}$ manifolds, we receive
$\infty^{ h - o}$ distinct families of $\infty^p$ manifolds, 
in total $\infty^{ h - o + p}$ different manifolds. At
the same time, it is clear that in this way, we obtain exactly
the same family as when we execute all transformations of
our group on the manifolds~\thetag{ 11} themselves. 

As a result, it is proved that the family~\thetag{ 12} consists
of $\infty^{h - o + p}$ different manifolds, hence that amongst
the $h + r$ parameters of the equations~\thetag{12}, exactly
$h - o + p$ are essential. 

\smallercharacters{
We yet want to indicate how one has to proceed in order to determine
the numbers $p$ and $o$ mentioned above.

The number $r - p$ is evidently the number of the infinitesimal
transformations: $e_1\, X_1f + \cdots + e_r\, X_rf$ which leave
invariant a generally positioned manifold of the family~\thetag{
11}. Now, the most general infinitesimal transformation: $\sum\, e_j\,
X_jf$ which leaves invariant the manifold:
\def\theequation{13}\begin{equation}
V_k(x_1,\dots,x_n,\,u_1,\dots,u_h)
=
0
\ \ \ \ \ \ \ \ \ \ \ \ \ 
{\scriptstyle{(k\,=\,1\,\cdots\,n\,-\,q)}}
\end{equation}
necessarily has the form: 
\def\theequation{14}\begin{equation}
\sum_{j=1}^r\,e_j(u_1,\dots,u_h)\,X_jf,
\end{equation}
where the $e_j ( u_1, \dots, u_h)$ are functions of the $u$. So we
need only to determine the functions $e_j ( u)$ in the most general
way so that the system of equations~\thetag{ 13} in the variables:
$x_1, \dots, x_n$ admits the infinitesimal transformation~\thetag{
14}.

If we imagine that the equations~\thetag{ 13} are resolved with
respect to $n - q$ of the $x$:
\[
x_{q+k}
=
\omega_{q+k}(x_1,\dots,x_q,\,u_1,\dots,u_h)
\ \ \ \ \ \ \ \ \ \ \ \ \ 
{\scriptstyle{(k\,=\,1\,\cdots\,n\,-\,q)}},
\]
and if we denote the substitution: $x_{ q+1} = \omega_{ q+1}$, \dots,
$x_n = \omega_n$ by the sign $[ \,\,\, ]$, we visibly obtain for the
functions $e_j ( u)$ the equations:
\def\theequation{15}\begin{equation}
\aligned
\sum_{j=1}^r\,e_j(u_1,\dots,u_h)\,
&
\bigg\{
\big[\xi_{j,\,q+k}\big]
-
\sum_{\pi=1}^q\,
\big[\xi_{j\pi}\big]\,
\frac{\partial\omega_{q+k}}{\partial x_\pi}
\bigg\}
=
0
\\
&
\ \ \ \ \ \ \
{\scriptstyle{(k\,=\,1\,\cdots\,n\,-\,q)}}
\endaligned
\end{equation}
which must be satisfied independently of the values of $x_1, \dots,
x_n$.

If we have determined from these equations the $e_j ( u)$ in the most
general way, we know the most general infinitesimal transformation
$\sum\, e_j\, X_jf$ which leaves invariant the manifold~\thetag{ 13}
and from this, we can immediately deduce the number of the independent
infinitesimal transformations $\sum\, e_j\, X_jf$ of this sort.

For the determination of the number $o$, we proceed as follows:

We seek at first the most general infinitesimal transformation:
$\sum\, \varepsilon_j ( u_1, \dots, u_h)\, X_jf$ that transfers the
manifold~\thetag{ 13} to an infinitely neighbouring manifold, or,
expressed differently: we seek the most general infinitesimal
transformation:
\[
\sum_{j=1}^r\,
\varepsilon_j(u_1,\dots,u_h)\,X_jf
+
\sum_{\sigma=1}^h\,
\Phi_\sigma(u_1,\dots,u_h)\,
\frac{\partial f}{\partial u_\sigma}
\] 
which leaves invariant the system of equations~\thetag{ 13} in the
$n + h$ variables: $x_1, \dots, x_n$, $u_1, \dots, u_h$.

If we keep the notation chosen above, the functions: 
$\varepsilon_j ( u)$ and $\Phi_\sigma ( u)$ are obviously
defined by the equations:
\def\theequation{15'}\begin{equation}
\aligned
\sum_{j=1}^r\,\varepsilon_j(u)\,
\bigg\{
\big[\xi_{j,\,q+k}\big]
&
-
\sum_{\pi=1}^q\,
\big[\xi_{j\pi}\big]\,
\frac{\partial\omega_{q+k}}{\partial x_\pi}
\bigg\}
=
\sum_{\sigma=1}^h\,
\Phi_\sigma(u)\,
\frac{\partial\omega_{q+k}}{\partial u_\sigma}
\\
&
\ \ \ \ \ \ \
{\scriptstyle{(k\,=\,1\,\cdots\,n\,-\,q)}}
\endaligned
\end{equation}
which they must satisfy independently of the values of the
variables $x_1, \dots, x_q$. 

We imagine that the $\varepsilon_j ( u)$ and the $\Phi_\sigma ( u)$
are determined in the most general way from these equations 
and we form the expression:
\[
\sum_{\sigma=1}^h\,\Phi_\sigma(u_1,\dots,u_h)\,
\frac{\partial f}{\partial u_\sigma}.
\]
Evidently, this expression can be deduced from a completely determined
number, say $h' \leqslant h$, of expressions:
\[
\sum_{\sigma=1}^h\,
\Phi_{\tau\sigma}(u_1,\dots,u_h)\,
\frac{\partial f}{\partial u_\sigma}
\ \ \ \ \ \ \ \ \ \ \ \ \ {\scriptstyle{(\tau\,=\,1\,\cdots\,h')}}
\]
by means of additions and of multiplications by functions of the $u$.
According to the nature of things, the $h'$ equations:
\[
\sum_{\sigma=1}^h\,\Phi_{\tau\sigma}
(u_1,\dots,u_h)\,
\frac{\partial f}{\partial u_\sigma}
=
0
\ \ \ \ \ \ \ \ \ \ \ \ \ {\scriptstyle{(\tau\,=\,1\,\cdots\,h')}}
\]
form an $h'$-term complete system with $h - h'$ independent solutions:
\[
w_1(u_1,\dots,u_h),
\,\,\,\dots,\,\,\,
w_{h-h'}(u_1,\dots,u_h),
\]
and it is clear that the equations:
\[
w_1={\rm const.},
\,\,\,\dots,\,\,\,
w_{h-h'}={\rm const.}
\]
determine the $\infty^{ h - o}$ subfamilies in which the
family~\thetag{ 11} can be decomposed, as we saw above.
Consequently, we have: $h - o = h - h'$, hence
$o = h'$.

}

It should not remain unmentioned that the developments of the present
chapter can yet be generalized.

One can for instance, instead of starting from an individual manifold,
start from a discrete number of manifolds or even more generally:
instead of starting from an individual family of manifolds, one can
start from several such families. We call a number of manifolds
briefly a \terminology{figure} \deutsch{Figur}.

\sectionengellie{\S\,\,\,116.}

The Theorem~85, p.~\pageref{Theorem-85-S-483} contains a method for
the determination of transitive groups which are isomorphic with a
given $r$-term group; we have already explained a method for the
determination of \emphasis{all} groups of this sort in
Chap.~\ref{kapitel-22}, p.~\pageref{S-434-sq} sq. At present, we will
show that our new method fundamentally amounts finally to the old one,
and we will in this way come to state an important result found
earlier on, in a new, much more general form.

Let:
\def\theequation{16}\begin{equation}
V_k(x_1,\dots,x_n)
=
0
\ \ \ \ \ \ \ \ \ \ \ \ \ 
{\scriptstyle{(k\,=\,1\,\cdots\,n\,-\,q)}}
\end{equation}
be an arbitrary manifold and let:
\[
x_i'
=
f_i(x_1,\dots,x_n,\,a_1,\dots,a_r)
\ \ \ \ \ \ \ \ \ \ \ \ \ {\scriptstyle{(i\,=\,1\,\cdots\,n)}}
\]
be an arbitrary $r$-term group generated by $r$ independent
infinitesimal transformations. By resolution with respect to the $x$,
the equations: $x_i' = f_i ( x, a)$ may give:
\[
x_i
=
F_i(x_1',\dots,x_n',\,a_1,\dots,a_r)
\ \ \ \ \ \ \ \ \ \ \ \ \ {\scriptstyle{(i\,=\,1\,\cdots\,n)}}.
\]
Lastly, we want yet to assume that the manifold~\thetag{ 16} admits
exactly $r - m$ independent infinitesimal transformations of the
group: $x_i' = f_i ( x, a)$.

If we execute on the manifold~\thetag{ 16} the general transformation:
$x_i' = f_i ( x, a)$ of our group, then according to Theorem~85,
p.~\pageref{Theorem-85-S-483}, we obtain a family of $\infty^m$
manifolds invariant by the group. The equations of this family are:
\def\theequation{17}\begin{equation}
\aligned
V_k\big(
F_1(x',a),\,\dots,\,
&
F_n(x',a)\big)
=
W_k(x_1',\dots,x_n',\,a_1,\dots,a_r)
=
0
\\
&
\ \ \ \ \ \ \ \ \
{\scriptstyle{(k\,=\,1\,\cdots\,n\,-\,q)}},
\endaligned
\end{equation}
hence they formally contain $r$ arbitrary parameters. But amongst
these $r$ parameters, only $m$ are essential, hence it is possible to
indicate $m$ independent functions: $\omega_1 ( a), \dots, \omega_m (
a)$ of the $a$ so that the equations~\thetag{ 17} take the form:
\def\theequation{17'}\begin{equation}
\Omega_k\big(
x_1',\dots,x_n',\,\omega_1(a),\dots,\omega_m(a)\big)
=
0
\ \ \ \ \ \ \ \ \ \ \ \ \ 
{\scriptstyle{(k\,=\,1\,\cdots\,n\,-\,q)}}.
\end{equation}
Here, we can lastly introduce: $l_1 = \omega_1 ( a)$, \dots, $l_m =
\omega_m ( a)$ as new parameters in place of the $a$; then in the
system of equations:
\def\theequation{18}\begin{equation}
\Omega_k(x_1',\dots,x_n',\,l_1,\dots,l_m)
=
0
\ \ \ \ \ \ \ \ \ \ \ \ \ 
{\scriptstyle{(k\,=\,1\,\cdots\,n\,-\,q)}}
\end{equation}
which comes into existence in this way, the parameters $l_1, \dots,
l_m$ are essential.

We find a new representation of our invariant family when we execute
an arbitrary transformation: $x_i'' = f_i ( x', b)$ of our group on
the system of equations~\thetag{ 18}. According to Theorem~85, 
p.~\pageref{Theorem-85-S-483}, \thetag{ 18} receives on the
occasion the form:
\def\theequation{18'}\begin{equation}
\Omega_k(x_1'',\dots,x_n'',\,l_1',\dots,l_m')
=
0
\ \ \ \ \ \ \ \ \ \ \ \ \ 
{\scriptstyle{(k\,=\,1\,\cdots\,n\,-\,q)}},
\end{equation}
where the $l'$ are completely determined functions of the $l$
and of the $b$:
\def\theequation{19}\begin{equation}
l_\mu'
=
\chi_\mu(l_1,\dots,l_m,\,b_1,\dots,b_r)
\ \ \ \ \ \ \ \ \ \ \ \ \ {\scriptstyle{(\mu\,=\,1\,\cdots\,m)}}.
\end{equation}
We must obtain this representation of our family when we execute
the transformation:
\[
x_i''
=
f_i\big(
f_1(x,a),\,\dots,\,f_n(x,a),\,\,b_1,\dots,b_r\big)
=
f_i(x_1,\dots,x_n,\,a_1',\dots,a_r')
\]
directly on the manifold~\thetag{ 16}, where the $a'$ are completely
determined functions of the $a$ and of the $b$:
\def\theequation{20}\begin{equation}
a_k'
=
\varphi_k(a_1,\dots,a_r,\,b_1,\dots,b_r)
\ \ \ \ \ \ \ \ \ \ \ \ \ {\scriptstyle{(k\,=\,1\,\cdots\,r)}}.
\end{equation}
If we do that, we receive the equations of our family at first
under the form:
\[
\aligned
V_k\big(F_1(x'',a'),\,\dots,\,
&
F_n(x'',a')\big)
=
W_k(x_1'',\dots,x_n'',\,a_1',\dots,a_r')
=
0
\\
&
\ \ \ \ \ \ \ \ \ \ \
{\scriptstyle{(k\,=\,1\,\cdots\,n\,-\,q)}},
\endaligned
\]
which we can also obviously write under the form:
\[
\Omega_k\big(
x_1'',\dots,x_n'',\,
\omega_1(a'),\dots,\omega_m(a')\big)
=
0
\ \ \ \ \ \ \ \ \ \ \
{\scriptstyle{(k\,=\,1\,\cdots\,n\,-\,q)}}.
\]
But since these equations must agree with the equations~\thetag{ 18'},
it results that the parameters $l'$ are linked to the $a'$ by the
relations:
\[
l_\mu'
=
\omega_\mu(a_1',\dots,a_r')
\ \ \ \ \ \ \ \ \ \ \ \ \ {\scriptstyle{(\mu\,=\,1\,\cdots\,m)}}.
\]

Because of the equations~\thetag{ 19}, we therefore have:
\[
\omega_\mu(a')
=
\chi_\mu(l_1,\dots,l_m,\,b_1,\dots,b_r)
\ \ \ \ \ \ \ \ \ \ \ \ \ {\scriptstyle{(\mu\,=\,1\,\cdots\,m)}},
\]
or:
\def\theequation{21}\begin{equation}
\omega_\mu(a')
=
\chi_\mu\big(
\omega_1(a),\dots,\omega_m(a),\,b_1,\dots,b_r)
\ \ \ \ \ \ \ \ \ \ \ \ \ {\scriptstyle{(\mu\,=\,1\,\cdots\,m)}}.
\end{equation}
If we make here the substitution: $a_1' = \varphi_1 ( a, b)$, \dots,
$a_r' = \varphi_r ( a, b)$, then we must obtain only identities,
because the parameters $a_1, \dots, a_r$, $b_1, \dots, b_r$ are
absolutely arbitrary, and hence, are linked together by no relation.

At present, if we remember that the equations~\thetag{ 20} represent a
simply transitive group in the variables $a_1, \dots, a_r$ equally
composed with the group: $x_i' = f_i ( x,a)$, namely the associated
parameter group (Chap.~\ref{kapitel-21}, p.~\pageref{S-404}), and that
the equations~\thetag{ 19} represent a transitive group isomorphic
with the group: $x_i' = f_i ( x, a)$, then we realize immediately what
follows:

The equations:
\def\theequation{22}\begin{equation}
\omega_1(a_1,\dots,a_r)
=
{\rm const.},
\,\,\,\dots,\,\,\,
\omega_m(a_1,\dots,a_r)
=
{\rm const.}
\end{equation}
represent a decomposition of the $r$-times extended space $a_1, \dots,
a_r$ invariant by the group~\thetag{ 20}, and to be precise, 
a decomposition in $\infty^m$ $(r - m)$-times extended manifolds.
But the group:
\def\theequation{19}\begin{equation}
l_\mu'
=
\chi_\mu(l_1,\dots,l_m,\,b_1,\dots,b_r)
\ \ \ \ \ \ \ \ \ \ \ \ \ {\scriptstyle{(\mu\,=\,1\,\cdots\,m)}}
\end{equation}
indicates in which way the $\infty^m$ manifolds~\thetag{ 22}
are permuted with each other by the transformations of the simply
transitive group~\thetag{ 20}.

We therefore see that the group~\thetag{ 19} can be derived from the
simply transitive group: $a_k' = \varphi_k ( a, b)$ according to the
rules of the preceding chapter (p.~\pageref{S-435} sq.).

We can use this fact in order to decide under which conditions do we
obtain two groups~\thetag{ 19} similar to each other when we start
from two different manifolds~\thetag{ 16}. Thanks to similar
considerations, we realize that the following statement holds, of
which Theorem~80 in Chap.~\ref{kapitel-22},
p.~\pageref{Theorem-80-S-445} is only a special case, fundamentally:

\def\thetheorem{86}\begin{theorem}
If, in the space $x_1, \dots, x_n$, an $r$-term group: $X_1f, \dots,
X_rf$ is presented, if in addition, two manifolds $M$ and $M'$ are
given, and if one executes all $\infty^r$ transformations of the
group: $X_1f, \dots, X_rf$ on each one of these two manifolds, then
the individual manifolds of the two invariant families that one
obtains in this way are transformed by two transitive groups
isomorphic with the group: $X_1f, \dots, X_rf$ which are similar to
each other when and only when the following two conditions are
satisfied:

\terminology{Firstly}, the two manifolds $M$ and $M'$ must admit the
same number of independent infinitesimal transformations of the form:
$e_1\, X_1f + \cdots + e_r\, X_rf$, and:

\terminology{Secondly}, it must be possible to relate the group:
$X_1f, \dots, X_rf$ to itself in a holoedrically isomorphic way so
that, to every infinitesimal transformation which fixes the one
manifold, there corresponds an infinitesimal transformation which
leaves invariant the other manifold.
\end{theorem}

It goes without saying that this theorem remains also valid
yet when one replaces the two manifolds by two figures.

\sectionengellie{\S\,\,\,117.}

\renewcommand{\thefootnote}{\fnsymbol{footnote}}
Specially worthy of note is the case where one has a manifold, or a
figure, which admits absolutely no infinitesimal transformation of the
$G_r$: $X_1f, \dots, X_rf$. The group~\thetag{ 19} isomorphic to $G_r$
then has the form:
\[
l_k'
=
\chi_k(l_1,\dots,l_r,\,a_1,\dots,a_r)
\ \ \ \ \ \ \ \ \ \ \ \ \ {\scriptstyle{(k\,=\,1\,\cdots\,r)}},
\]
it is simply transitive and hence holoedrically isomorphic with 
the $G_r$.\footnote[1]{\,
\name{Lie}, Gesellsch. d. W. zu Christiania, 1884.
}
\renewcommand{\thefootnote}{\arabic{footnote}}

We therefore have here a general method in order to set up the simply
transitive groups that are equally composed with a given $r$-term
group.

The method applied in the preceding chapter, Proposition~1,
page~\pageref{Satz-1-S-431}, is only a special case of the present
more general method. Indeed at that time\,---\,as we can now express
this\,---\,we used as a figure the totality of $r$ different points of
the $R_n$.

When no special assumption is made about the position of these points,
then the figure consisting of them can allow none of the infinitesimal
transformations of the group $X_kf$. Since the group $X_kf$ certainly
leaves invariant no point in general position, there can be in it at
most $r - 1$ independent infinitesimal transformations which fix such
a point; from these possible $r - 1$ infinitesimal transformations,
one can linearly deduce again at most $r - 2$ independent
infinitesimal transformations by which a second point in general
position yet remains untouched, and so on; one realizes at the end
that there is no infinitesimal transformation in the group by which
$r$ points in general position remain simultaneously invariant.

Hence if we take a figure which consists of $r$ points of this sort,
and if we execute on them the $\infty^r$ transformations of the group:
$X_1f, \dots, X_rf$, then the figure takes $\infty^r$ different
positions, the totality of which remains invariant by the group. These
$\infty^r$ positions are transformed by a simply transitive group.

If we really set up this simply transitive group, we obtain exactly
the same group as the one obtained thanks to the method of the
preceding chapter.

\sectionengellie{\S\,\,\,118.}

It appears to be desirable to illustrate the general developments
of the \S\S\,\,114 and 116 by means of a specific example. However, 
we restrict ourselves here to giving indications, and we
leave to the reader the effective execution of the concerned
simple computations. 

The most general projective group of the plane which leaves
invariant the nondegenerate conic section: $x^2 - 2y = 0$ is 
three-term and contains the following three independent
infinitesimal transformations:
\[
\aligned
X_1f
=
\frac{\partial f}{\partial x}
&
+
x\,\frac{\partial f}{\partial y},
\ \ \ \ \ \ \ \ \ \
X_2f
=
x\,\frac{\partial f}{\partial x}
+
2y\,\frac{\partial f}{\partial y}
\\
&
X_3f
=
(x^2-y)\,\frac{\partial f}{\partial x}
+
xy\,\frac{\partial f}{\partial y}.
\endaligned
\]

One sees immediately that this group\,---\,that we want to
call $G_3$ for short\,---\,is transitive and that its composition
is determined through the relations:
\[
\leftbracket
X_1,\,X_2
\rightbracket
=
X_1f,\ \ \ \ \ \ \
\leftbracket
X_1,\,X_3
\rightbracket
=
X_2f,\ \ \ \ \ \ \
\leftbracket
X_2,\,X_3
\rightbracket
=
X_3f.
\]
From this, it follows by taking account of Chap.~\ref{kapitel-15},
Proposition~8, p.~\pageref{Satz-8-S-263} that the $G_3$ contains no
two-term invariant subgroup; one convinces oneself that it does not
contain either any one-term invariant subgroup. Consequently, the
$G_3$ is \emphasis{simple} (Chap.~\ref{kapitel-15},
p.~\pageref{S-264}).

Every tangent to the fixed conic section $x^2 - 2y = 0$ admits exactly
two independent infinitesimal transformations of the group; besides,
it can be proved that the tangents are the only curves of the plane
which possess this property. In the same way, the conic sections
which enter in contact \deutsch{berühren} in two points with the fixed
conic section are the only curves which admit one and only one
infinitesimal transformation of the $G_3$; as a conic section which
enters twice in contact, one must also certainly count every line of
the plane which cuts \deutsch{schneidet} the conic section in two
separate points. Lastly, it is clear that every point not lying on the
conic section admits one and only one infinitesimal transformation of
the group.

Now, if one chooses as manifold $M$ an arbitrary other curve, hence a
curve which admits absolutely no infinitesimal transformation of the
$G_3$, then by the group, this curve takes $\infty^3$ different
positions and the totality of these positions is obviously transformed
by means of a three-term group. One therefore finds a simply
transitive group of the $R_3$ equally composed with the original
$G_3$. All groups that one obtains in this way are similar to each
other; one amongst them is for instance the three-term group of all
projective transformations of the $R_3$ which leaves invariant a
winding curve of third order.

If one introduces as manifold $M$ a conic section (irreducible or
decomposable) entering in contact in two separate points [with the
fixed conic section], then one obtains a group holoedrically
isomorphic to the $G_3$ in a twice-extended manifold; all groups
obtained in this way are similar to each other. By contrast, if one
uses as a manifold $M$ a conic section having four points entering in
contact [with the fixed conic section], one obtains a completely
different type of three-term group of a twice-extended manifold.

Lastly, if one introduces as manifold $M$ a tangent to the fixed conic
section, then one obtains a three-term group in a once-extended
manifold which is similar to the general projective group of the
straight line.

\sectionengellie{\S\,\,\,119.}

Let a linear homogeneous group:
\[
X_kf
=
\sum_{\mu,\,\,\nu}^{1\cdots\,n}\,
a_{k\mu\nu}\,x_\mu\,
\frac{\partial f}{\partial x_\nu}
\ \ \ \ \ \ \ \ \ \ \ \ \ {\scriptstyle{(k\,=\,1\,\cdots\,r)}}
\]
of the $R_n$ be presented. This group leaves invariant the family
of the $\infty^n$ planes $M_{ n-1}$ of the $R_n$:
\[
u_1\,x_1
+\cdots+
u_n\,x_n
+
1
=
0.
\]
We want to set up the group corresponding to it in the parameters
$u_1, \dots, u_n$. 

According to \S\,\,111, the infinitesimal transformations:
\[
U_kf
=
\sum_{\nu=1}^n\,v_{k\nu}(u_1,\dots,u_n)\,
\frac{\partial f}{\partial u_\nu}
\]
of the sought group are to be determined so that the equation
$\sum\, u_\nu \, x_\nu + 1 = 0$ admits the infinitesimal
transformation $X_kf + U_kf$. So, the $r$ expressions:
\[
\sum_{\mu,\,\,\nu}^{1\cdots\,n}\,
a_{k\mu\nu}\,x_\mu\,u_\nu
+
\sum_{\nu}^{1\cdots\,n}\,
x_\nu\,v_{k\nu}
\]
must vanish by means of $\sum\, u_\nu\, x_\nu + 1 = 0$. 
This is possible only when they vanish identically, that is to
say when $v_{ k\nu} = -\, \sum_\mu \, a_{ k\nu \mu}\, u_\mu$.

\medskip

We therefore find:
\[
U_kf
=
-\,\sum_{\mu,\,\,\nu}^{1\cdots\,n}\,
a_{k\nu\mu}\,u_\mu\,
\frac{\partial f}{\partial u_\nu}.
\]
Consequently, the group $U_kf$ is linear homogeneous too; we know from
the beginning that it is isomorphic with the group $X_kf$.

We call the group: $U_1f, \dots, U_rf$ the group
\terminology{dualistic} \deutsch{dualistisch} to the group: $X_1f,
\dots, X_rf$.

In order to give an example, we want to look up at the group 
dualistic to the adjoint group of an $r$-term group having the
composition:
\[
\leftbracket
Y_i,\,Y_k
\rightbracket
=
\sum_{s=1}^r\,c_{iks}\,Y_sf.
\]

\renewcommand{\thefootnote}{\fnsymbol{footnote}}
The adjoint group reads (cf. p.~\pageref{S-275}):
\[
\sum_{i,\,\,s}^{1\cdots\,r}\,
c_{iks}\,e_i\,
\frac{\partial f}{\partial e_s}
\ \ \ \ \ \ \ \ \ \ \ \ \ {\scriptstyle{(k\,=\,1\,\cdots\,r)}},
\]
whence the group dualistic to it reads\footnote[1]{\,
\name{Lie}, Math. Ann. Vol. XVI, p.~496, cf. also Archiv for
Mathematik, Vol. 1, Christiania 1876.
}: 
\[
\sum_{i,\,\,s}^{1\cdots\,r}\,
c_{kis}\,\varepsilon_s\,
\frac{\partial f}{\partial\varepsilon_i}
\ \ \ \ \ \ \ \ \ \ \ \ \ {\scriptstyle{(k\,=\,1\,\cdots\,r)}},
\]
where the relation: $c_{ iks} = -\, c_{ kis}$ is used.
\renewcommand{\thefootnote}{\arabic{footnote}}

Similar considerations can actually be made for all projective
groups of the $R_n$, since they all leave invariant the family
of the $\infty^n$ straight, $(n-1)$-times extended manifolds of
the $R_n$. However, we do not want to enter these considerations,
and rather, we refer to the next volume in which the concept
of duality \deutsch{Dualität} is considered under a more
general point of view, namely as a special case of the general
concept of contact transformation
\deutsch{Berührungstransformation}.

\sectionengellie{\S\,\,\,120.}

\label{S-494-sq}
Finally, we yet want to consider an important example of a more
general nature.

We assume that we know all $q$-term subgroups of the $G_r$: $X_1f,
\dots, X_rf$. Then the question is still to decide what are the
different \emphasis{types} of such $q$-term subgroups.

In Chap.~\ref{kapitel-16}, p.~\pageref{S-281-bis}, we already have
explained the concept of ``types of subgroups''; according to that, we
reckon two $q$-term subgroups as belonging to the same type when they
are conjugate to each other inside the $G_r$; of all subgroups which
belong to the same type, we therefore need only to indicate a single
one, and in addition, all these subgroups are perfectly determined by
this only one.

We know that every subgroup of the group: $X_1f, \dots, X_rf$ is
represented by a series of linear homogeneous relations between the
parameters: $e_1, \dots, e_r$ of the general infinitesimal
transformation: $e_1\, X_1f + \cdots + e_r\, X_rf$ (cf.
Chap.~\ref{kapitel-12}, p.~\pageref{Satz-6-S-211}). Moreover, we know
that two subgroups are conjugate to each other inside the group:
$X_1f, \dots, X_rf$ if and only if the system of equations between
the $e$ which represents the one subgroup can be transferred, by means
of a transformation of the adjoint group, to the system of equations
which represents the other subgroup (Chap.~\ref{kapitel-16}, 
p.~\pageref{S-280-bis}).

Now according to our assumption, we know all $q$-term subgroups of
the $G_r$, hence we know all systems of $r - q$ independent
linear homogeneous equations between the $e$: 
\[
\sum_{j=1}^r\,h_{kj}\,e_j
=
0
\ \ \ \ \ \ \ \ \ \ \ \ \ {\scriptstyle{(k\,=\,1\,\cdots\,r\,-\,q)}}
\]
which represent $q$-term subgroups. 

For reasons of simplicity, amongst all these systems of equations,
we want to take those which can be resolved with respect to 
$e_{ q+1}, \dots, e_r$, hence which can be brought to the form:
\def\theequation{23}\begin{equation}
\aligned
e_{q+k}
&
=
g_{q+k,\,1}\,e_1
+\cdots+
g_{q+k,\,q}\,e_q
\\
&
\ \ \ \ \ \ \ \ \ \ \ 
{\scriptstyle{(k\,=\,1\,\cdots\,r\,-\,q)}}\,;
\endaligned
\end{equation}
we leave the remaining ones which cannot be brought to this form,
because they could naturally be treated in exactly the same way as
those of the form~\thetag{ 23}.

All systems of values $g_{ q+k,\, j}$ which, when inserted in~\thetag{
23}, provide $q$-term subgroups are defined by means of certain
equations between the $g_{ q+k, \, j}$; however in general, it is not
possible to represent all these systems of values by means of a single
system of equations between the $g$, and rather, a discrete number of
such systems of equations will be necessary if one wants to have all
$q$-term subgroups which are contained in the form~\thetag{ 23}.
Naturally, two different systems of equations of this sort then
provide nothing but different types of $q$-term subgroups.

We restrict ourselves to an arbitrary system of equations 
amongst the concerned systems of equations, say the following one:
\def\theequation{24}\begin{equation}
\Omega_\mu
(g_{q+1,\,1},\,g_{q+1,\,2},\,\dots,\,g_{r,\,q})
=
0
\ \ \ \ \ \ \ \ \ \ \ \ \ 
{\scriptstyle{(\mu\,=\,1,\,\,2\,\cdots\,)}},
\end{equation}
and we now want to see what types of $q$-term subgroups does this
system determine. The equations~\thetag{ 23} determine, when the $g_{
q + k, \, j}$ are completely arbitrary, the family of all straight
$q$-times extended manifolds of the space $e_1, \dots, e_r$ which pass
through the point: $e_1 = 0$, \dots, $e_r = 0$. Of course, this family
of manifolds remains invariant by the adjoint group:
\def\theequation{25}\begin{equation}
e_k'
=
\sum_{j=1}^r\,
\rho_{kj}(a_1,\dots,a_r)\,e_j
\ \ \ \ \ \ \ \ \ \ \ \ \ {\scriptstyle{(k\,=\,1\,\cdots\,r)}}
\end{equation}
of the group: $X_1f, \dots, X_rf$. Hence, if we execute the 
transformation~\thetag{ 25} on the system of equations~\thetag{ 23}, 
we obtain a system of equations in the $e'$ of the corresponding
form:
\[
e_{q+k}'
=
\sum_{j=1}^q\,g_{q+k,\,j}'\,e_j'
\ \ \ \ \ \ \ \ \ \ \ \ \ {\scriptstyle{(k\,=\,1\,\cdots\,r\,-\,q)}},
\]
where the $g'$ are linear homogeneous functions of the $g$
with coefficients which depend upon the $a$:
\def\theequation{26}\begin{equation}
\aligned
g_{q+k,\,j}'
&
=
\sum_{\mu=1}^{r-q}\,\sum_{\nu=1}^q\,
\alpha_{kj\mu\nu}(a_1,\dots,a_r)\,
g_{q+\mu,\,\nu}
\\
&\ \ \ \ \
{\scriptstyle{(k\,=\,1\,\cdots\,r\,-\,q\,;\,\,\,
j\,=\,1\,\cdots\,q)}}.
\endaligned
\end{equation}

According to p.~\pageref{Satz-4-S-469} sq., the equations~\thetag{ 26}
determine a group in the variables $g$. The system of
equations~\thetag{ 24} remains invariant by this group, because every
system of values $g_{ q+k, \, j}$ which provides a subgroup is
naturally transferred to a system of values $g_{ q + k, \, j}'$ which
determines a subgroup; but since the group~\thetag{ 26} is continuous,
it leaves individually invariant all discrete regions of systems of
values $g_{ q + k, j}$ of this sort, hence in particular also the
system of equations:
\def\theequation{24}\begin{equation}
\Omega_\mu(g_{q+1,\,1},\dots,g_{r,\,q})
=
0
\ \ \ \ \ \ \ \ \ \ \ \ \ 
{\scriptstyle{(\mu\,=\,1,\,\,2\,\cdots\,)}}.
\end{equation}

Now, the question is how the systems of values~\thetag{ 24}
are transformed by the group~\thetag{ 26}, and whether every 
system of values can be transferred to every other, or not.

This question receives a graphic sense when we imagine that
the $g_{ q + k, \, j}$ are point coordinates in a space of
$q ( r - q)$ dimensions. Indeed, the equations~\thetag{ 24}
then represent a certain manifold in this space which remains
invariant by the group~\thetag{ 26}. Each point of the
manifold belongs to a certain smallest invariant subsidiary 
domain of the manifold and the points of such a subsidiary 
domain represent nothing but conjugate $q$-term subgroups of the
$G_r$, and to be precise, all the $q$-term subgroups of the
$G_r$ which belong to one and the same type.

Besides, one must draw attention on the fact that only the points $g_{
q + k,\, j}$ of such a smallest invariant subsidiary domain that
belong in turn to no smaller invariant subsidiary domain are to be
counted, because it is only when one delimits the subsidiary domain in
this way that each one of its points can be transferred to any other
point by means of a nondegenerate transformation of the group~\thetag{
26}.

Hence, if one wants to determine all types of $q$-term subgroups which
are contained in the system of equations~\thetag{ 24}, then one only
has to look up at all smallest invariant subsidiary domains of the
manifold $\Omega_\mu = 0$. Every such subsidiary domain determines all
subgroups which belong to the same type, and an arbitrary point of the
subsidiary domain provides a group which can be chosen as a
representative of the concerned type.

Besides, in order to be able to determine the discussed invariant
subsidiary domains, one does not at all need to assume that the finite
equations~\thetag{ 25} of the adjoint group are known; it suffices
that one has the infinitesimal transformations of this group, because
one can then immediately indicate the infinitesimal transformations of
the group~\thetag{ 26} and afterwards, following the rules of
Chap.~\ref{kapitel-14}, one can determine the desired invariant
subsidiary domains; in the case present here, this determination
requires only executable operations.

One observes that the preceding developments remain also applicable
when one knows not an $r$-term \emphasis{group}, but only a possible
\emphasis{composition} of such a group, hence a system of $c_{ iks}$
which satisfies the known relations:
\[
\left\{\aligned
&
\ \ \ \ \ \ \ \ \ \ \ \ \ \ \ \ \ \ 
c_{iks}+c_{kis}
=
0
\\
&
\sum_{\nu=1}^r\,
\big\{
c_{ik\nu}\,c_{\nu js}
+
c_{kj\nu}\,c_{\nu is}
+
c_{ji\nu}\,c_{\nu ks}
\big\}
=
0
\\
&
\ \ \ \ \ \ \ \ \ \ \ \ \ \ \ \ \ \ \ \ 
{\scriptstyle{(i,\,\,k,\,\,j,\,\,s\,=\,1\,\cdots\,r)}}.
\endaligned\right.
\]

If the $G_r$: $X_1f, \dots, X_rf$ is invariant in a larger group
$\mathfrak{ G}$, then the question is often whether two subgroups of
the $G_r$ are conjugate in this larger group $\mathfrak{ G}_r$, or
not. In this case, one can also define differently the concept of
``type of subgroup of the $G_r$'', by reckoning two subgroups of the
$G_r$ as being distinct only when they are not conjugate to each other
also in the $\mathfrak{ G}$.

If one wants to determine all types, in this sense, of subgroups of
the $G_r$, then this task presents no special difficulty. Indeed, the
problem in question is obviously a part of the more general problem of
determining all types of subgroups of the $\mathfrak{ G}$, when the
word ``type'' is understood in the sense of Chap.~\ref{kapitel-16},
p.~\pageref{S-281-bis}.

This study is particularly important when the group $\mathfrak{ G}$
is actually the largest subgroup of the $R_n$ in which the $G_r$: 
$X_1f, \dots, X_rf$ is invariant.

\linestop


\chapter{Systatic and Asystatic Transformation Groups}
\label{kapitel-24}
\chaptermark{Systatic and Asystatic Transformation Groups}

\setcounter{footnote}{0}

\abstract*{??}

In the $s$-times extended space $x_1, \dots, x_s$, let an $r$-term
group:
\[
X_kf
=
\sum_{i=1}^s\,\xi_{ki}(x_1,\dots,x_n)\,
\frac{\partial f}{\partial x_i}
\ \ \ \ \ \ \ \ \ \ \ \ \ {\scriptstyle{(k\,=\,1\,\cdots\,r)}},
\]
or shortly $G_r$, be presented. Of the $r$ infinitesimal 
transformations $X_1f, \dots, X_rf$, let there be precisely $n$, say:
$X_1f, \dots, X_nf$, which are linked together by no linear relation, 
while by contrast $X_{ n+1}f, \dots, X_rf$ can be expressed
as follows:
\def\theequation{1}\begin{equation}
\label{S-497}
X_{n+j}f
\equiv
\sum_{\nu=1}^n\,\varphi_{j\nu}(x_1,\dots,x_n)\,
X_\nu f
\ \ \ \ \ \ \ \ \ \ \ \ \ {\scriptstyle{(j\,=\,1\,\cdots\,r\,-\,n)}}.
\end{equation}

Now, if $x_1^0, \dots, x_s^0$ is a point for which not all
$n \times n$ determinants of the matrix:
\def\theequation{2}\begin{equation}
\left\vert
\begin{array}{cccc}
\xi_{11}(x) & \,\cdot\, & \,\cdot\, & \xi_{1s}(x)
\\
\cdot & \,\cdot\, & \,\cdot\, & \cdot
\\
\xi_{n1}(x) & \,\cdot\, & \,\cdot\, & \xi_{ns}(x)
\end{array}
\right\vert
\end{equation}
vanish, then according to Chap.~\ref{kapitel-11},
p.~\pageref{Satz-7-S-203}, there are in the $G_r$ exactly $r - n$
independent infinitesimal transformations whose power series
expansions with respect to the $x_i - x_i^0$ contain no term of zeroth
order, but only terms of first order or of higher order; the point
$x_1^0, \dots, x_s^0$ therefore admits exactly $r - n$ independent
infinitesimal transformations of the $G_r$ which generate an $(r -
n)$-term subgroup $G_{ r -n}$ of the $G_r$
(cf. Chap.~\ref{kapitel-12}, p.~\pageref{Satz-2-S-205}). The
infinitesimal transformations of this $G_{ r -n}$ can, according to
p.~\pageref{S-203-ter}, be linearly deduced from the $r - n$
independent transformations:
\[
X_{n+j}f
-
\sum_{\nu=1}^n\,
\varphi_{j\nu}(x_1^0,\dots,x_s^0)\,
X_\nu f
\ \ \ \ \ \ \ \ \ \ \ \ \ {\scriptstyle{(j\,=\,1\,\cdots\,r\,-\,n)}}.
\]

In the sequel, we want to briefly call a point $x_1^0, \dots, 
x_s^0$ for which not all $n \times n$ determinants of the 
matrix~\thetag{ 2} vanish, a point \emphasis{in general position}.

\sectionengellie{\S\,\,\,121.}

According to the above, to every point $x_1^0, \dots, x_s^0$ in
general position is associated a completely determined $(r - n)$-term
subgroup of the $G_r$, namely the most general subgroup of the $G_r$
by which it remains invariant.

If we let the point $x^0$ change its position, we receive a new $(r -
n)$-term subgroup of the $G_r$, and since there are $\infty^s$
different points, we receive in total $\infty^s$ subgroups of this
sort; however, it is not said that we obtain $\infty^s$ different
subgroups.

If for example $r = n$, then the discussed $\infty^s$ subgroups
coincide all, namely they reduce all to the identity transformation,
since the $G_r$ contains absolutely no infinitesimal transformation
which leaves at rest a point $x^0$ in general position.

But disregarding also this special case, it can happen that to the
$\infty^s$ points of the space $x_1, \dots, x_s$, only $\infty^{
s-1}$, or less, different groups of the said constitution are
associated; evidently, this will always occur in any case when
\emphasis{there is a continuous family of individual points which
simultaneously keep their positions by the $(r - n)$-term subgroup}.

At present, we want to look up at the analytic conditions under 
which such a phenomenon occurs. At first, we take up the question: 
when do two points in general position remain invariant by the
same $(r - n)$-term subgroup of the $G_r$?

The answer to this question has a great similarity with the
considerations in Chap.~\ref{kapitel-19}, p.~\pageref{S-357-sq} sq.

Let the one point be $x_1^0, \dots, x_s^0$ and let us call $G_{ r-n}$
the associated $(r - n)$-term subgroup of the $G_r$; then the general
infinitesimal transformation of the $G_{ r - n}$ reads:
\[
\sum_{j=1}^{r-n}\,\varepsilon_j\,
\bigg(
X_{n+j}f
-
\sum_{\nu=1}^n\,\varphi_{j\nu}^0\,X_\nu f
\bigg),
\]
where it is understood that the $\varepsilon$ are arbitrary 
parameters.

Let the other point be: $\overline{ x}_1, \dots, \overline{ x}_s$, and
let the general infinitesimal transformation of the subgroup
$\overline{ G}_{ r - n}$ associated to it then be:
\[
\sum_{j=1}^{r-n}\,
\overline{\varepsilon}_j\,
\bigg(
X_{n+j}f
-
\sum_{\nu=1}^n\,
\overline{\varphi}_{j\nu}\,X_\nu f
\bigg).
\]
Now, if the two points are supposed to remain invariant by the same
$(r - n)$-term subgroup, then $G_{ r - n}$ and $\overline{ G}_{ r-n}$
do coincide; for this, it is necessary and sufficient that all
infinitesimal transformations of the one belong to those of the other
group, and conversely; when expressed analytically, one must be able
to satisfy identically the equation:
\[
\sum_{j=1}^{r-n}\,\varepsilon_j\,
\bigg\{
X_{n+j}f
-
\sum_{\nu=1}^n\,\varphi_{j\nu}^0\,X_\nu f
\bigg\}
=
\sum_{j=1}^{r-n}\,
\overline{\varepsilon}_j\,
\bigg\{
X_{n+j}f
-
\sum_{\nu=1}^n\,
\overline{\varphi}_{j\nu}\,X_\nu f
\bigg\}
\]
for arbitrarily chosen $\varepsilon$ thanks to suitable values of the
$\overline{ \varepsilon}$, and also for arbitrarily chosen $\overline{
\varepsilon}$ thanks to suitable values of the $\varepsilon$.

The latter equation can also be written:
\[
\sum_{j=1}^{r-n}\,
(\varepsilon_j-\overline{\varepsilon}_j)\,
X_{n+j}f
-
\sum_{\nu=1}^n\,\sum_{j=1}^{r-n}\,
(\varepsilon_j\,\varphi_{j\nu}^0
-
\overline{\varepsilon}_j\,\overline{\varphi}_{j\nu})\,
X_\nu f
=
0,
\]
hence, because of the independence of the infinitesimal
transformations $X_1f, \dots, X_rf$, if can hold only if $\overline{
\varepsilon}_j = \varepsilon_j$; in addition, since the
$\varepsilon_j$ are absolutely arbitrary, we obtain:
\[
\varphi_{j\nu}(x_1^0,\dots,x_s^0)
=
\varphi_{j\nu}(\overline{x}_1,\dots,\overline{x}_s)
\ \ \ \ \ \ \ \ \ \ \ \ \ 
{\scriptstyle{(j\,=\,1\,\cdots\,r\,-\,n\,;\,\,\,
\nu\,=\,1\,\cdots\,n)}}.
\]

Consequently, the two $(r - n)$-term subgroups of the $G_r$ which are
associated to two distinct points in general position are identical
with each other if and only if each one of the $n ( r - n)$ functions
$\varphi_{ j\nu}$ takes the same numerical values for the one point 
as for the other point.

Hence, if we want to know all points in general position which, 
under the group $G_{ r-n}$, 
keep their positions simultaneously with the point $x_1^0, \dots, 
x_s^0$, then we only have to determine all systems of values
$x$ which satisfy the equations:
\def\theequation{3}\begin{equation}
\varphi_{j\nu}(x_1,\dots,x_s)
=
\varphi_{j\nu}^0
\ \ \ \ \ \ \ \ \ \ \ \ \ 
{\scriptstyle{(j\,=\,1\,\cdots\,r\,-\,n\,;\,\,\,
\nu\,=\,1\,\cdots\,n)}}\,;
\end{equation}
every such system of values provides a point having the
constitution demanded. 

Here, two cases have to be distinguished.

\terminology{Firstly} the number of mutually independent
functions amongst the $n ( r - n)$ functions $\varphi_{ k\nu} (x)$
can be equal to $s$ exactly. In this case, the $(r - n)$-term
subgroup $G_{ r-n}$ which fixes the point: $x_1^0, \dots, x_s^0$
leaves untouched yet at most a discrete number of points 
in general position.

\terminology{Secondly} the number of independent functions amongst the
$\varphi_{ k\nu} (x)$ can be smaller than $s$. In this case, there is
a continuous manifold of points in general position which remain all
invariant by the group $G_{ r - n}$; to every point of the manifold in
question is then associated the same $(r - n)$-term subgroup of the
$G_r$ as to the point $x_1^0, \dots, x_s^0$. At the same time, the
point $x_1^0, \dots, x_s^0$ evidently lies inside the manifold, that
is to say: there are, in the manifold, also points that are infinitely
close to the point: $x_1^0, \dots, x_s^0$.

We assume that the second case happens, so that amongst the $n ( r -
n)$ functions $\varphi_{ k\nu} (x)$, there are only $s - \rho < s$
that are mutually independent, and we may call them $\varphi_1 ( x),
\dots, \varphi_{ s - \rho} (x)$.

Under this assumption, all $\varphi_{ k\nu} ( x)$ can be expressed by
means of $\varphi_1 ( x), \dots, \varphi_{ s - \rho} ( x)$ alone, and
the equations~\thetag{ 3} can be replaced by the $s - \rho$ mutually
independent equations:
\def\theequation{3'}\begin{equation}
\varphi_1(x_1,\dots,x_s)
=
\varphi_1(x_1^0,\dots,x_s^0),
\,\,\,\dots,\,\,\,
\varphi_{s-\rho}(x_1,\dots,x_s)
=
\varphi_{s-\rho}(x_1^0,\dots,x_s^0).
\end{equation}
We therefore see that to every generally located point: $x_1^0, \dots,
x_s^0$ of the space is associated a completely determined $\rho$-times
extended manifold~\thetag{ 3'} which is formed of the totality of all
points to which is associated the same $(r - n)$-term subgroup of the
group: $X_1f, \dots, X_rf$ as to the point: $x_1^0, \dots, x_s^0$.

\smallercharacters{

It is possible that the equations:
\[
\varphi_{j\nu}(x_1,\dots,x_s)
=
\varphi_{j\nu}(x_1^0,\dots,x_s^0)
\ \ \ \ \ \ \ \ \ \ \ \ \ 
{\scriptstyle{(j\,=\,1\,\cdots\,r\,-\,n\,;\,\,\,
\nu\,=\,1\,\cdots\,n)}}
\]
represent, for every system of values $x_1^0, \dots, x_s^0$, a 
manifold which decomposes in a discrete number of different manifolds.

An example is provided by the three-term group:
\[
X_1f
=
\frac{\partial f}{\partial x_1},
\ \ \ \ \ \ \ \ \
X_2f
=
\frac{1}{\cos x_2}\,
\frac{\partial f}{\partial x_2},
\ \ \ \ \ \ \ \ \
X_3f
=
\tan x_2\,
\frac{\partial f}{\partial x_2}.
\]
Here, we have: 
\[
X_3f
\equiv
\sin x_2\,X_2f,
\]
while $X_1f$ and $X_2f$ are linked by no linear relation.
So, if one fixes the point $x_1^0$, $x_2^0$, all points whose
coordinates $x_1$, $x_2$ satisfy the equation:
\[
\sin x_2
=
\sin x_2^0
\]
also remain fixed, that is to say: simultaneously with the point
$x_1^0$, $x_2^0$, every point which lies in one of the infinitely
many lines:
\[
x_2
=
x_2^0+2k\,\pi
\]
parallel to the $x_1$-axis keeps its position, where it is 
understood that $k$ is an arbitrary, positive or negative, 
entire number.

}

Since, under the assumption made above, the group: $X_1f, \dots, X_rf$
associates to every point: $x_1^0, \dots, x_s^0$ a $\rho$-times
extended manifold passing through it, then the whole space $x_1,
\dots, x_s$ obviously decomposes in a family of $\infty^{ s - \rho}$
$\rho$-times extended manifolds:
\def\theequation{4}\begin{equation}
\varphi_1(x_1,\dots,x_s)
=
{\rm const.},
\,\,\,\dots,\,\,\,
\varphi_{s-\rho}(x_1,\dots,x_s)
=
{\rm const.},
\end{equation}
and to be precise, in such a way that all transformations of the group
which fix an arbitrarily chosen point in general position do leave at
rest all points of the manifold~\thetag{ 4} that passes through this
point.

Certainly, it is to be remarked here that one can speak of a real
decomposition of the space only when the number $\rho$ is smaller than
$s$; if $\rho = s$, no real decomposition of the space occurs, since
every transformation of our group which fixes a point in general
position actually leaves invariant all points of the space; in other
words: the identity transformation is the only transformation of the
group which leaves at rest a point in general position.

The preceding considerations give an occasion for an important
division \deutsch{Eintheilung} of all $r$-term groups $X_1f, \dots,
X_rf$ of the space $x_1, \dots, x_s$, and to be precise, for a
division in two different classes.

\renewcommand{\thefootnote}{\fnsymbol{footnote}}
\plainstatement{If a group of the space $x_1, \dots, x_s$ is
constituted in such a way that all its transformations which leave
invariant a point in general position do simultaneously fix all points
of a continuous manifold passing through this point, then we reckon
this group as belonging to the one class and we call them
\terminology{systatic} \deutsch{systatisch}. But we reckon all the
remaining groups, hence those which are not systatic, as belonging to
the other class, and we call them \terminology{asystatic}
\deutsch{asystatisch}.\footnote[1]{\,
The concepts of systatic and asystatic groups, and the theory of these
groups as well, stem from \name{Lie} (Ges. d. W. zu Christiania 1884,
Archiv for Math. Vol. 10, Christiania 1885). The expressive
terminologies: systatic ``with leaving fixed''
\deutsch{mitstehendlassend} and asystatic ``without leaving fixed''
\deutsch{nichtmitstehendlassend} are from \name{Engel}. As for the
rest, the terminology characteristic of the present work was
introduced by \name{Lie}.
}}
\renewcommand{\thefootnote}{\arabic{footnote}}

Using this terminology, we can express the gained result in the
following way: 

\def\thetheorem{87}\begin{theorem}
\label{Theorem-87-S-502}
If the independent infinitesimal transformations $X_1f, \dots, X_rf$
of an $r$-term group in $s$ variables $x_1, \dots, x_s$ are linked
together by $r - n$ linear relations of the form:
\[
X_{n+k}f
\equiv
\sum_{\nu=1}^n\,
\varphi_{k\nu}(x_1,\dots,x_s)\,X_\nu f
\ \ \ \ \ \ \ \ \ \ \ \ \ {\scriptstyle{(k\,=\,1\,\cdots\,r\,-\,n)}},
\]
while between $X_1f, \dots, X_nf$ alone no relation of this sort
holds, then the group is systatic when amongst the $n ( r - n)$
functions $\varphi_{ k\nu} ( x_1, \dots, x_s)$, there are less than
$s$ that are mutually independent; by contrast, if amongst the
functions $\varphi_{ k \nu} ( x_1, \dots, x_s)$, there are $s$
functions that are mutually independent, then the group is asystatic.
\end{theorem}

\sectionengellie{\S\,\,\,122.}

We want to maintain all the assumptions that we have made in the
introduction of the chapter about the $r$-term group: $X_1f, \dots,
X_rf$ of the space $x_1, \dots, x_s$, and at present, we only want to
add yet the assumption that the group is supposed to be systatic.
Thus, we assume that, amongst the $n ( r - n)$ functions $\varphi_{
k\nu} ( x_1, \dots, x_s)$, only $0 \leqslant s - \rho < s$ are
mutually independent, and as above, we may call them $\varphi_1 ( x),
\dots, \varphi_{ s - \rho} ( x)$.

We consider at first the case where the number $s - \rho$ has the
value zero.

If the number $s - \rho$ vanishes, then $n$ is obviously equal to 
$r$, that is to say the $r$ independent infinitesimal transformations:
$X_1f, \dots, X_rf$ are linked together by no linear relation 
of the form:
\[
\chi_1(x_1,\dots,x_s)\,X_1f
+\cdots+
\chi_r(x_1,\dots,x_s)\,X_rf
=
0.
\]
From this, it results that the number $r$ is in any case not larger
than the number $s$ of the variables $x$, hence that the group:
$X_1f, \dots, X_rf$ is either intransitive, or at most simply 
transitive. Hence in the two cases, the group: $X_1f, \dots, X_rf$ is
imprimitive (cf. Chap.~\ref{kapitel-13}, p.~\pageref{S-221} and
Chap.~\ref{kapitel-20}, Proposition~6, p.~\pageref{Satz-6-S-383}.

So we see that the systatic group: $X_1f, \dots, X_rf$ is always
imprimitive when the entire number $s - \rho$ has the value zero.

We now turn ourselves to the case $s - \rho > 0$.

In this case, the equations:
\def\theequation{4}\begin{equation}
\varphi_1(x_1,\dots,x_s)
=
{\rm const.},
\,\,\,\dots,\,\,\,
\varphi_{s-\rho}(x_1,\dots,x_s)
=
{\rm const.}
\end{equation}
provide a real decomposition of the space $x_1, \dots, x_s$ in
$\infty^{ s - \rho}$ $\rho$-times extended manifolds, and in fact, as
we have seen above, a decomposition which is completely determined by
the group: $X_1f, \dots, X_rf$ and which stands in a completely
characteristic relationship to this group. We are very close to
presume that this decomposition remains invariant by the group: $X_1f,
\dots, X_rf$; from this, it would then follow that the systatic group:
$X_1f, \dots, X_rf$ is imprimitive also in the case $s - \rho > 0$
(Chap.~\ref{kapitel-13}, p.~\pageref{S-220-sq}).

On can very easily see that the presumption expressed just now 
corresponds to the truth. In fact, according to Chap~\ref{kapitel-19}, 
p.~\pageref{S-343-sq} sq., the $r ( s - \rho)$ expressions: 
$X_k\, \varphi_1$, \dots, $X_k\, \varphi_{ s - \rho}$ can
be expressed as functions of $\varphi_1, \dots, \varphi_{ s - \rho}$
alone:
\[
X_k\,\varphi_j
=
\pi_{jk}(\varphi_1,\dots,\varphi_{s-\rho})
\ \ \ \ \ \ \ \ \ \ \ \ \ 
{\scriptstyle{(k\,=\,1\,\cdots\,r\,;\,\,\,
j\,=\,1\,\cdots\,s\,-\,\rho)}}.
\]
In this (cf. Chap.~\ref{kapitel-8}, p.~\pageref{Satz-1-S-139} and
\pageref{S-143}) lies the reason why the decomposition~\thetag{ 4}
admits the $r$ infinitesimal transformations $X_kf$, and therefore
actually, the complete group: $X_1f, \dots, X_rf$. 

The developments just carried out prove that a systatic
group of the space $x_1, \dots, x_s$ is always imprimitive. 
We therefore have the

\def\thetheorem{88}\begin{theorem}
Every systatic group is imprimitive.
\end{theorem}

It is not superfluous to establish, also by means of conceptual
considerations, that in the case $s - \rho > 0$, the
decomposition~\thetag{ 4} remains invariant by the systatic
group: $X_1f, \dots, X_rf$. 

Let us denote by $M$ an arbitrary manifold amongst the $\infty^{
s - \rho}$ $\rho$-times extended manifolds~\thetag{ 4}, let
$P$ be the general symbol of a point of the manifold $M$, 
let $S$ be the general symbol of the $\infty^{ r - n}$
transformations of our group which leave untouched
all the points of $M$, and lastly, let us understand by $T$
an arbitrary transformation of our group.

If we execute the transformation $T$ on $M$, we obtain a certain
$\rho$-times extended manifold $M'$, the $\infty^\rho$ points
$P'$ of which are defined by the equation:
\[
(P')
=
(P)\,T.
\]
Now, since every point $P$ remains invariant by all $\infty^{ r -n}$
transformations $S$ of our group, it is clear that every point $P'$
keeps its position by the $\infty^{ r - n}$ transformations $T^{ -1}
\, S \, T$, which belong as well to our group. From this, it results
that $M'$ also belongs to the $\infty^{ s - \rho}$ manifolds~\thetag{
4}; consequently, it is proved that the $\infty^{ s - \rho}$ 
manifolds~\thetag{ 4} are permuted with each other by every 
transformation of the group: $X_1f, \dots, X_rf$, so that
the decomposition~\thetag{ 4} effectively remains invariant
by our group.

Since according to Theorem~88, every systatic group is imprimitive,
every primitive group must be asystatic, inversely.
But there are also imprimitive groups which are asystatic, for
instance the four-term group:
\[
\frac{\partial f}{\partial x},
\ \ \ \ \ \ \ \
\frac{\partial f}{\partial y},
\ \ \ \ \ \ \ \
x\,\frac{\partial f}{\partial y},
\ \ \ \ \ \ \ \
y\,\frac{\partial f}{\partial y}
\]
of the plane $x$, $y$. This group is imprimitive, since
it leaves invariant the family of straight lines: 
$x = {\rm const.}$, but it is at the same time asystatic,
because from the two identities:
\[
x\,\frac{\partial f}{\partial y}
\equiv
x\,\frac{\partial f}{\partial y},
\ \ \ \ \ \ \ \
y\,\frac{\partial f}{\partial y}
\equiv
y\,\frac{\partial f}{\partial y},
\]
it becomes clear that the functions $\varphi_{ k\nu}$ associated
to the group are nothing but $x$ and $y$ themselves; but
they are obviously independent of each other. The three-term
intransitive group:
\[
\frac{\partial f}{\partial y},
\ \ \ \ \ \ \ \
x\,\frac{\partial f}{\partial y},
\ \ \ \ \ \ \ \
y\,\frac{\partial f}{\partial y}
\]
shows that there are even intransitive asystatic groups.

\sectionengellie{\S\,\,\,123.}

If one knows $r$ independent infinitesimal transformations: $X_1f,
\dots, X_rf$ of an $r$-term group of the space $x_1, \dots, x_s$, then
according to Theorem~87, p.~\pageref{Theorem-87-S-502}, one can easily
decide whether the concerned group is systatic or not.

For this purpose, one forms the matrix:
\[
\left\vert
\begin{array}{cccc}
\xi_{11}(x) & \,\cdot\, & \,\cdot\, & \xi_{1s}(x)
\\
\cdot & \,\cdot\, & \,\cdot\, & \cdot
\\
\xi_{r1}(x) & \,\cdot\, & \,\cdot\, & \xi_{rs}(x)
\end{array}
\right\vert
\]
and one studies its determinants. If all the $( n+1) \times (n+1)$
determinants vanish identically, but not all $n \times n$
determinants, and in particular, not all $n \times n$ determinants
that can be formed with the topmost $n$ rows of the matrix, then there
are identities of the form:
\def\theequation{1}\begin{equation}
X_{n+k}f
\equiv
\sum_{\nu=1}^n\,
\varphi_{k\nu}(x_1,\dots,x_s)\,
X_\nu f
\ \ \ \ \ \ \ \ \ \ \ \ \ 
{\scriptstyle{(k\,=\,1\,\cdots\,r\,-\,n)}},
\end{equation}
while $X_1f, \dots, X_nf$ are linked together by no linear relation.
If one has set up the identities~\thetag{ 1}, then one determines the
number of independent functions amongst the $n ( r - n)$ functions
$\varphi_{ k\nu} ( x_1, \dots, x_s)$.

However, in order to be able to decide whether a determined $r$-term
group is systatic or not, one needs absolutely not to know the finite
expressions for the infinitesimal transformations of the group, and
rather, it suffices that one knows the \emphasis{defining equations}
of the group. This sufficiency is based on the fact that, as soon as
the defining equations of the group are presented, one can always
indicate an unrestricted integrable system of total differential
equations, the only integral functions of which are just the 
$\varphi_{ k \nu} (x)$ and the functions of them. Indeed, from
this it clearly turns out that one can determine the number
of independent functions amongst the $\varphi_{ k\nu} (x)$ 
without knowing the $\varphi_{ k\nu} (x)$ themselves.

At present, we will derive this important result.

The independent equations amongst the equations:
\def\theequation{5}\begin{equation}
d\varphi_{k\nu}
=
\sum_{i=1}^s\,
\frac{\partial\varphi_{k\nu}(x)}{\partial x_i}\,
\D\,x_i
=
0
\ \ \ \ \ \ \ \ \ \ \ \ \ 
{\scriptstyle{(k\,=\,1\,\cdots\,r\,-\,n\,;\,\,\,
\nu\,=\,1\,\cdots\,n)}}
\end{equation}
form an unrestricted integrable system of total differential
equations, the only integral functions of which are the 
$\varphi_{ k\nu} (x)$, and the functions of them
(cf. Chap.~\ref{kapitel-5}, p.~\pageref{S-91-sq} sq.).
We start from this system of total differential equations.

Because of~\thetag{ 1}, we have identically:
\[
\xi_{n+j,\,i}
-
\sum_{\nu=1}^n\,\varphi_{j\nu}\,\xi_{\nu i}
\equiv
0
\ \ \ \ \ \ \ \ \ \ \ \ \ 
{\scriptstyle{(j\,=\,1\,\cdots\,r\,-\,n\,;\,\,\,
i\,=\,1\,\cdots\,s)}},
\]
whence it comes by differentiation:
\[
d\xi_{n+j,\,i}
-
\sum_{\nu=1}^n\,\varphi_{j\nu}\,d\xi_{\nu i}
\equiv
\sum_{\nu=1}^n\,\xi_{\nu i}\,d\varphi_{j\nu},
\]
or, since according to the assumption, not all $n \times n$
determinants of the matrix:
\[
\left\vert
\begin{array}{cccc}
\xi_{11}(x) & \,\cdot\, & \,\cdot\, & \xi_{1s}(x)
\\
\cdot & \,\cdot\, & \,\cdot\, & \cdot
\\
\xi_{r1}(x) & \,\cdot\, & \,\cdot\, & \xi_{rs}(x)
\end{array}
\right\vert
\]
vanish, it yet comes:
\[
\aligned
d\varphi_{j\pi}
&
\equiv
\sum_{i=1}^s\,\chi_{\pi i}(x_1,\dots,x_s)\,
\bigg\{
d\xi_{n+j,\,i}
-
\sum_{\nu=1}^n\,\varphi_{j\nu}\,d\xi_{\nu i}
\bigg\}
\\
&
\ \ \ \ \ \ \ \ \ \ \ \ \ \
{\scriptstyle{(j\,=\,1\,\cdots\,r\,-\,n\,;\,\,\,
\pi\,=\,1\,\cdots\,n)}}.
\endaligned
\]
From this, it follows that the system of the total differential 
equations~\thetag{ 5} can be replaced by the following one:
\def\theequation{6}\begin{equation}
\sum_{\pi=1}^s\,
\bigg\{
\frac{\partial\xi_{n+j,\,i}}{\partial x_\pi}
-
\sum_{\nu=1}^n\,\varphi_{j\nu}\,
\frac{\partial\xi_{\nu i}}{\partial x_\pi}
\bigg\}\,\D\,x_\pi
=
0
\ \ \ \ \ \ \ \ \ \ \ \ \ \
{\scriptstyle{(j\,=\,1\,\cdots\,r\,-\,n\,;\,\,\,
i\,=\,1\,\cdots\,s)}}.
\end{equation}
Of course, the independent equations amongst the equations~\thetag{ 6}
form an unrestricted integrable system of total differential
equations, the integral functions of which are just the $\varphi_{
k\nu} (x)$, and the functions of them.

However, it is not possible now to set up the individual
equations~\thetag{ 6} when one only knows the defining equations of
the group: $X_1f, \dots, X_rf$; by contrast, it is possible to replace
the system of equations~\thetag{ 6} by another system which, aside
from the $\D\, x_\pi$, contains only the coefficients of the defining
equations. We arrive at this in the following way.

By $x_1^0, \dots, x_s^0$, we understand a point in general position.
According to p.~\pageref{S-203}, the most general infinitesimal
transformation: $e_1\, X_1f + \cdots + e_r\, X_rf$
whose power series expansion with respect to the
$x_i - x_i^0$ contains only terms of first order or
of higher order has the shape:
\[
\sum_{j=1}^{r-n}\,
\mathfrak{e}_j\,
\bigg(
X_{n+j}f
-
\sum_{\nu=1}^n\,
\varphi_{j\nu}(x_1^0,\dots,x_s^0)\,X_\nu f
\bigg),
\]
where it is understood that the $\mathfrak{ e}_j$ are
arbitrary parameters. But when we really execute the
power series expansion with respect to the $x_i - x_i^0$
and when at the same time we take into account only the
terms of first order, this expression receives the form:
\def\theequation{7}\begin{equation}
\sum_{j=1}^{r-n}\,\mathfrak{e}_j\,
\sum_{i,\,\,\pi}^{1\cdots\,s}\,
\bigg\{
\bigg[
\frac{\partial\xi_{n+j,\,i}}{\partial x_\pi}
\bigg]_{x=x^0}
-
\sum_{\nu=1}^n\,\varphi_{j\nu}(x^0)\,
\bigg[
\frac{\partial\xi_{\nu i}}{\partial x_\pi}
\bigg]_{x=x^0}
\bigg\}\,
(x_\pi-x_\pi^0)\,
\frac{\partial f}{\partial x_i}
+\cdots.
\end{equation}

But we can compute the terms of first order in the expression~\thetag{
7} from the defining equations of the group, and to be precise,
without integration. Indeed, according to Chap.~\ref{kapitel-11},
p.~\pageref{S-188-sq} sq., the most general infinitesimal
transformation whose power series expansion with respect to the $x_i -
x_i^0$ contains only terms of first order or of higher order has the
form:
\[
\sum_{i,\,\,\pi}^{1\cdots\,s}\,
g_{i\pi}'(x_\pi-x_\pi^0)\,
\frac{\partial f}{\partial x_\pi}
+\cdots,
\]
where a certain number of the $s^2$ quantities $g_{ \pi i}'$ which we
denoted by $\varepsilon_1 - \nu_1$ at that time were arbitrary, while
the $s^2 - \varepsilon_1 + \nu_1$ remaining ones were linear
homogeneous functions of these $\varepsilon_1 - \nu_1$ quantities with
coefficients which could be computed immediately from the defining
equations. We therefore obtain, when we start from the defining
equations, the following representation for the expression~\thetag{
7}:
\def\theequation{7'}\begin{equation}
\sum_{j=1}^{\varepsilon_1-\nu_1}\,
\mathfrak{e}_j'\,
\sum_{i,\,\,\pi}^{1\cdots\,s}\,
\alpha_{j\pi i}(x_1^0,\dots,x_s^0)\,
(x_\pi-x_\pi^0)\,
\frac{\partial f}{\partial x_i}
+\cdots,
\end{equation}
where the $\mathfrak{ e}_j'$ denote arbitrary parameters, while the
$\alpha_{ j\pi i} (x^0)$ are completely determined analytic functions
of the $x^0$ which can, as said above, be computed from the defining
equations.

Because~\thetag{ 7} and~\thetag{ 7'} are only different
representations of the same infinitesimal transformation, the factors
of $(x_\pi - x_\pi^0) \, \partial f / \partial x_i$ in the two
expressions must be equal to each other, that is to say, there are the
following $s^2$ relations:
\[
\sum_{j=1}^{r-n}\,\mathfrak{e}_j\,
\bigg[
\frac{\partial\xi_{n+j,\,i}}{\partial x_\pi}
-
\sum_{\nu=1}^n\,\varphi_{j\nu}(x)\,
\frac{\partial\xi_{\nu i}}{\partial x_\pi}
\bigg]_{x=x^0}
=
\sum_{j=1}^{\varepsilon_1-\nu_1}\,
\mathfrak{e}_j'\,
\alpha_{j\pi i}(x^0)
\ \ \ \ \ \ \ \ \ \ \ \ \ 
{\scriptstyle{(i,\,\,\pi\,=\,1\,\cdots\,s)}}.
\]
For arbitrarily chosen $\mathfrak{ e}_j$, one must always be able to
satisfy these relations thanks to suitable choices of the values of
the $\mathfrak{ e}_j'$, and for arbitrarily chosen $\mathfrak{ e}_j'$,
always thanks to suitable choices of the values of the $\mathfrak{
e}_j$.

All of this holds for every point $x_1^0, \dots, x_s^0$ in general
position; hence this also holds true when we consider the $x^0$
as variables and when we substitute them by $x_1, \dots, x_n$.
Consequently, when we have chosen the $\mathfrak{ e}$ in 
completely arbitrary way as functions of the $x$, we
can always determine the $\mathfrak{ e}'$ as functions
of the $x$ in such a way that the equations:
\def\theequation{7''}\begin{equation}
\aligned
\sum_{j=1}^{r-n}\,
\mathfrak{e}_k\,
\bigg\{
\frac{\partial\xi_{n+j,\,i}}{\partial x_\pi}
&
-
\sum_{\nu=1}^n\,\varphi_{j\nu}(x)\,
\frac{\partial\xi_{\nu i}}{\partial x_\pi}
\bigg\}
=
\sum_{j=1}^{\varepsilon_1-\nu_1}\,
\mathfrak{e}_j'\,\alpha_{j\pi i}(x)
\\
& \ \ \ \ \ 
{\scriptstyle{(i,\,\,\pi\,=\,1\,\cdots\,s)}}
\endaligned
\end{equation}
are identically satisfied, and when we have chosen the $\mathfrak{
e}'$ in completely arbitrary way as functions of the $x$, we can
always satisfy~\thetag{ 7''} identically thanks to suitable
functions $\mathfrak{ e}_1, \dots, \mathfrak{ e}_{ r -n}$ of the $x$.

Now, one can obviously replace the equations~\thetag{ 6} by the
following $s$ equations:
\def\theequation{6'}\begin{equation}
\sum_{j=1}^{r-n}\,\mathfrak{e}_j\,
\sum_{\pi=1}^s\,
\bigg\{
\frac{\partial\xi_{n+j,\,i}}{\partial x_\pi}
-
\sum_{\nu=1}^n\,\varphi_{j\nu}\,
\frac{\partial\xi_{\nu i}}{\partial x_\pi}
\bigg\}\,\D\,x_\pi
=
0
\ \ \ \ \ \ \ \ \ \ \ \ \ {\scriptstyle{(i\,=\,1\,\cdots\,s)}},
\end{equation}
provided only that one regards the $\mathfrak{ e}$ as arbitrary
functions of the $x$ in them. Hence from what has been said above, it
follows that the totality of all equations~\thetag{ 6} is equivalent
to the totality of all equations of the form:
\[
\sum_{j=1}^{\varepsilon_1-\nu_1}\,
\mathfrak{e}_j'
\sum_{\pi=1}^s\,
\alpha_{j\pi i}(x_1,\dots,x_s)\,
\D\,x_\pi
=
0
\ \ \ \ \ \ \ \ \ \ \ \ \ {\scriptstyle{(i\,=\,1\,\cdots\,s)}},
\]
in which the $\mathfrak{ e}'$ are to be interpreted as arbitrary
functions of the $x$. Lastly, the latter equations can evidently be
replaced by the $(\varepsilon_1 - \nu_1)s$ equations:
\def\theequation{8}\begin{equation}
\sum_{\pi=1}^s\,
\alpha_{j\pi i}(x_1,\dots,x_s)\,
\D\,x_\pi
=
0
\ \ \ \ \ \ \ \ \ \ \ \ \ 
{\scriptstyle{(j\,=\,1\,\cdots\,\varepsilon_1\,-\,\nu_1\,;\,\,\,
i\,=\,1\,\cdots\,s)}}.
\end{equation}

With these words, it is proved that the two systems of total
differential equations: \thetag{ 6} and~\thetag{ 8} are equivalent to
each other; thus, it results that the independent equations amongst
the equations~\thetag{ 8} form an unrestricted integrable system of
total differential equations, and to be precise, a system, the only
integral functions of which are the $\varphi_{ k\nu} ( x)$ and the
functions of them.

We can therefore say:

\renewcommand{\thefootnote}{\fnsymbol{footnote}}
\def\thetheorem{89}\begin{theorem}
If the defining equations of an $r$-term group of the space $x_1,
\dots, x_s$ are presented, then one decides in the following way
whether the concerned group is systatic or not:

One understands by $x_1^0, \dots, x_s^0$ an arbitrary point in the
neighbourhood of which the coefficients of the resolved defining
equations behave regularly and one determines the terms 
of zeroth order and of first order in the power series expansion
of the general infinitesimal transformation of the group 
with respect to the powers of $x_1 - x_1^0$, \dots, 
$x_s - x_s^0$. Afterwards, one searches for the terms of first
order in the most general infinitesimal transformation
of the group which contains no term of zeroth order.
These terms will have the form:
\[
\sum_{j=1}^{\varepsilon_1-\nu_1}\,
\mathfrak{e}_j'\,
\sum_{i,\,\,\pi}^{1\cdots\, s}\,
\alpha_{j\pi i}(x_1^0,\dots,x_s^0)\,
(x_\pi-x_\pi^0)\,
\frac{\partial f}{\partial x_i},
\]
where the $\mathfrak{ e}_j'$ denote arbitrary parameters, 
while the $\alpha_{ j\pi i} (x^0)$ are completely determined
analytic functions of the $x^0$ and can be computed without
integration from the coefficients of the defining equations. Now, 
one forms the system of the total differential equations:
\def\theequation{8}\begin{equation}
\sum_{\pi=1}^s\,
\alpha_{j\pi i}(x_1,\dots,x_s)\,\D\,x_\pi
=
0
\ \ \ \ \ \ \ \ \ \ \ \ \ 
{\scriptstyle{(j\,=\,1\,\cdots\,\varepsilon_1\,-\,\nu_1\,;\,\,\,
i\,=\,1\,\cdots\,s)}}
\end{equation}
and one determines the number $s - \rho$ of the independent
equations amongst these equations. If $s - \rho < s$, then
the group is systatic, but if $s - \rho = s$, the group is
asystatic.\footnote[1]{\,
\name{Lie}, Archiv for Math., Vol. 10, Christiania 1885.
}
\end{theorem}
\renewcommand{\thefootnote}{\arabic{footnote}}

We can add:

\def\theproposition{1}\begin{proposition}
The $s - \rho$ mutually independent equations amongst the differential
equations~\thetag{ 8} form an unrestricted integrable system with $s -
\rho$ independent integral functions: $\varphi_1 ( x), \dots,
\varphi_{ s - \rho} (x)$. These integral functions stand in the
following relationship to the group: $X_1f, \dots, X_rf$:

If, amongst the $r$ independent infinitesimal transformations: $X_1f,
\dots, X_rf$ there are exactly $n$, say $X_1f, \dots, X_nf$, that are
linked together by no linear relation, while $X_{ n+1}f, \dots, X_rf$
can be expressed linearly in terms of $X_1f, \dots, X_nf$:
\[
X_{n+k}f
\equiv
\sum_{\nu=1}^n\,\varphi_{k\nu}(x_1,\dots,x_s)\,
X_\nu f
\ \ \ \ \ \ \ \ \ \ \ \ \ {\scriptstyle{(k\,=\,1\,\cdots\,r\,-\,n)}},
\]
then all $n ( r - n)$ functions $\varphi_{ k\nu} (x)$ can be 
expressed in terms of $\varphi_1 (x), \dots, \varphi_{ s - \rho}
(x)$ alone.
\end{proposition}

One even does not need to know the defining equations themselves in
order to be able to decide whether a determined group is systatic or
not. \emphasis{For this, one only needs to know the initial terms in
the power series expansions of the infinitesimal transformations of
the group in the neighbourhood of an individual point: $x_1^0, \dots,
x_s^0$ in general position.} Indeed, if one knows these initial terms,
one can obviously compute the numerical values $\alpha_{ j\pi i}^0$
that the functions $\alpha_{ j\pi i} ( x)$ take for $x_1 = x_1^0$,
\dots, $x_s = x_s^0$. The number $s - \rho$ defined above then is
nothing but the number of the mutually independent equations amongst
the linear equations in the $\D\,x_1, \dots, \D\,x_s$:
\[
\sum_{\pi=1}^s\,\alpha_{j\pi i}^0\,\D\,x_\pi
=
0
\ \ \ \ \ \ \ \ \ \ \ \ \ 
{\scriptstyle{(j\,=\,1\,\cdots\,\varepsilon_1\,-\,\nu_1\,;\,\,\,
i\,=\,1\,\cdots\,s)}}.
\]

The system of the total differential equations has a very simple
conceptual meaning. Indeed, as one easily realizes directly, it
defines all the points: $x_1 + \D\,x_1$, \dots, $x_s + \D\,x_s$ 
infinitely
close to the point $x_1, \dots, x_s$ that remain invariant by all
transformations of our group which leave at rest the point $x_1,
\dots, x_s$. Here lies the inner reason \deutsch{innere Grund} why
the $\varphi_{ k\nu} (x)$ and the functions of them are the only
integral functions of the system~\thetag{ 6} or~\thetag{ 8}, because
indeed, the equations:
\[
\varphi_{k\nu}(y_1,\dots,y_s)
=
\varphi_{k\nu}(x_1,\dots,x_s)
\ \ \ \ \ \ \ \ \ \ \ \ \ 
{\scriptstyle{(k\,=\,1\,\cdots\,r\,-\,n\,;\,\,\,
\nu\,=\,1\,\cdots\,n)}}
\]
define all points $y_1, \dots, y_s$ which remain invariant
simultaneously with the point $x_1, \dots, x_s$.

\sectionengellie{\S\,\,\,124.}

We have found that the $r$-term group: $X_1f, \dots, X_rf$ of the
space $x_1, \dots, x_s$ is asystatic or is systatic according to
whether there are exactly $s$, or less than $s$, mutually independent
functions amongst the functions $\varphi_{ k\nu} (x_1, \dots, x_s)$
defined on p.~\pageref{S-332} and~\pageref{S-497}. Now, according to
Chap.~\ref{kapitel-20}, Theorem~67, p.~\pageref{Theorem-67-S-376},
there always is an infinitesimal transformation $Zf$ which is
interchangeable with all $X_kf$ when, and only when, the number of
independent functions amongst the $\varphi_{ k\nu} (x)$ is smaller
than $s$. Consequently, we can also say:

\renewcommand{\thefootnote}{\fnsymbol{footnote}}
\def\theproposition{2}\begin{proposition}
The $r$-term group: $X_1f, \dots, X_rf$ of the space $x_1, \dots, 
x_s$ is systatic if and only if there is an infinitesimal
transformation $Zf$ which is interchangeable with all $X_kf$; 
if there is no such infinitesimal transformation, 
the group: $X_1f, \dots, X_rf$ is asystatic.\footnote[1]{\,
\name{Lie}, Archiv for Math., Vol. 10, p.~377, Christiania 1885.
}
\end{proposition}
\renewcommand{\thefootnote}{\arabic{footnote}}

Since the excellent infinitesimal transformations of a group are 
interchangeable with all the other infinitesimal transformations
of the group, we have in addition:

\def\theproposition{3}\begin{proposition}
Every group which contains one or several excellent infinitesimal
transformations is systatic.
\end{proposition}

Thus for such groups, one realizes already from the composition that
they are systatic.

Finally, if we remember that the adjoint group of a group without
excellent infinitesimal transformation contains $r$ essential
parameters (cf. Theorem~49, p.~\pageref{Theorem-49-S-277}), we see
that the following proposition holds:

\def\theproposition{4}\begin{proposition}
The adjoint group of an asystatic group is always $r$-term.
\end{proposition}

Proposition~2 and Proposition~4 are the generalizations of the
Propositions~1 and~2 of the Chap.~\ref{kapitel-16}
(p.~\pageref{Satz-1-S-277}) announced at that time.

Already from the developments of the Chap.~\ref{kapitel-20} we could
have taken as an opportunity a division of all groups in two classes;
in the first class, we had reckoned every $r$-term group: $X_1f,
\dots, X_rf$ for which there is at least one infinitesimal
transformation interchangeable with all $X_kf$, and in the other
class, all the remaining groups. According to what was said above, it
is clear that this division would coincide with our present division
of the groups in systatic and asystatic groups; the first class would
consist of all systatic group, and the second class, of all asystatic
groups. In what follows, we want to explain this fact by means of
conceptual considerations and at the same time, we want to derive new
important results.

Let $\Theta$ be a transformation which is interchangeable with all
transformations of a given $r$-term group: $X_1f, \dots, X_rf$ of the
space $x_1, \dots, x_s$, and moreover, let $S$ be the general symbol
of the transformations of the group: $X_1f, \dots, X_rf$ which leave
invariant an arbitrary, generally positioned point $P$ of the space
$x_1, \dots, x_s$.

If, $P$ is transferred to the point $P_1$ by the execution of the
transformation $\Theta$, then obviously, $\Theta^{ -1} S \, \Theta$ is
the general symbol of all transformations of the group: $X_1f, \dots,
X_rf$ which leave invariant the point $P_1$. But now, since $\Theta$
is interchangeable with all transformations of the group: $X_1f,
\dots, X_rf$, then the totality of all transformations $\Theta^{ -1} S
\, \Theta$ is identical to the totality of all transformations $S$,
hence we see that all transformations of the group: $X_1f, \dots,
X_rf$ which fix the point $P$ do also leave at rest the point $P_1$.

We now apply this to the case where there is a continuous family of
transformations $\Theta$ which are interchangeable with all
transformations of the group: $X_1f, \dots, X_rf$.

Since $P$ is a point in general position, then by the execution of the
transformations $\Theta$, it takes a continuous series of different
positions. But as we have seen just now, each one of these positions
remains invariant by all transformations $S$, and consequently, the
group: $X_1f, \dots, X_rf$ is systatic.

\medskip

As a result, it is proved that the $r$-term group: $X_1f, \dots, X_rf$
of the space $x_1, \dots, x_s$ is in any case systatic when there is
an infinitesimal transformation $Zf$ interchangeable with all
$X_kf$. It yet remains to show that the converse also holds true,
namely that for every systatic group: $X_1f, \dots, X_rf$, one can
indicate a continuous family of transformations which are
interchangeable with all transformations $X_1f, \dots, X_rf$.

Thus, we imagine that a systatic $r$-term group: $X_1f, \dots, X_rf$
of the space $x_1, \dots, x_s$ is given. We will indicate a
construction that provides infinitely many transformations which are
interchangeable with all transformations of this group.

Every transformation $\Theta$ which is interchangeable with all
transformations of the group: $X_1f, \dots, X_rf$ transfers every
point $P$ of the space to a point $P_1$ which admits exactly the same
transformations of the group as the point $P$; this is what we showed
above. \emphasis{Hence, we choose two arbitrary points $P$ and 
\label{S-511} $P_1$
which admit the same transformations of our group, and we attempt to
determine a transformation $\Theta$ which is interchangeable with all
transformations of our group and which in addition transfers $P$ to
$P_1$.}

Let $P'$ be a point to which $P$ can be transferred by means of a
transformation $T$ of the group: $X_1f, \dots, X_rf$; furthermore,
let, as earlier on, $S$ be the general symbol of all transformations
of this group which leave invariant $P$ and hence $P_1$ too. Then
(Chap.~\ref{kapitel-14}, Proposition~1, p.~\pageref{Satz-1-S-227}),
$S\, T$ is the general symbol of all transformations of our group
which transfer $P$ to $P'$. Obviously, $P_1$ takes the same position
by all these transformations $S\, T$, for one indeed has:
\[
(P_1)\,S\,T
=
(P_1)\,T.
\]

At present, by assuming the existence of \emphasis{one} transformation
$\Theta$ having the constitution just demanded, we can easily see that
every point $P' = (P)\, T$ receives a completely determined new
position by all possible $\Theta$; indeed, this follows immediately
from the equations:
\[
(P')\,\Theta
=
(P)\,T\,\Theta
=
(P)\,\Theta\,T
=
(P_1)\,T.
\]
\emphasis{Hence if $T$ is an arbitrary transformation of our group,
then by every transformation $\Theta$ which actually exists, the point
$(P)\, T$ takes the new position $(P_1)\, T$.}

We add that we obtain in this way no overdetermination
\deutsch{Überbestimmung} of the new position of the point $P'$.
Indeed, if we replace in the latter equations the transformation $T$
by an arbitrary other transformation of the group: $X_1f, \dots, X_rf$
which transfers in the same way $P$ to $P'$, hence if write $S\, T$ in
place of $T$, then it comes again:
\[
(P')\,\Theta
=
(P)\,S\,T\,\Theta
=
(P)\,\Theta\,S\,T
=
(P_1)\,S\,T
=
(P_1)\,T.
\]

We consider at first the special case where the systatic group: $X_1f,
\dots, X_rf$ is \emphasis{transitive}.

When the group: $X_1f, \dots, X_rf$ is transitive, by a suitable
choice of $T$, the point $(P)\, T$ can be brought to coincidence with
every other point $(P_1)\, T$ of the space; hence, if we associate to
every point $(P)\, T$ of the space the point $(P_1)\, T$, a completely
determined transformation $\Theta'$ is defined in this way. If we yet
succeed to prove that $\Theta'$ is interchangeable with all
transformations of our group, then it is clear that $\Theta'$
possesses all properties which we have required of the transformation
$\Theta$, and that $\Theta'$ is the only transformation $\Theta$ which
actually exists.

The fact that $\Theta'$ is really interchangeable with all
transformations of our group can be easily proved. Indeed, we have:
\[
(P)\,T\,\Theta'
=
(P_1)\,T
=
(P)\,\Theta'\,T,
\]
where $T$ means a completely arbitrary transformation of our group.
Hence if we understand in the same way by ${\sf T}$ a completely
arbitrary transformation of our group, we obtain:
\[
(P)\,T\,{\sf T}\,\Theta'
=
(P)\,\Theta'\,T\,{\sf T}
=
(P)\,T\,\Theta'\,{\sf T},
\]
and therefore, the transformation: ${\sf T}\, \Theta' \, {\sf T}^{-1}
{\Theta'}^{ -1}$ leaves invariant the point $(P)\, T$, that is to say,
every point of the space. From this, it follows that ${\sf T}\,
\Theta' \, {\sf T}^{ -1} {\Theta '}^{ -1}$ is the identity
transformation, that is to say: $\Theta'$ is really interchangeable
with all transformations of our group.

Thus, when the systatic group: $X_1f, \dots, X_rf$ is transitive,
to every pair of points $P$, $P_1$ having the constitution
defined above there corresponds one and only one transformation
interchangeable with all transformations of the group.
If one chooses the pair of points in all possible ways,
one obtains infinitely many such transformations, and it is
easy to determine how many: In any case, when $\infty^\rho$
different points remain untouched by all transformations
of the group which leave invariant an arbitrarily chosen
point, then there are exactly $\infty^\rho$ different transformations
which are interchangeable with all transformations of the group:
$X_1f, \dots, X_rf$. This is coherent with Theorem~67,
p.~\pageref{Theorem-67-S-376}, for because of the assumption made
above, amongst the $n ( r - n)$ functions $\varphi_{ k\nu} ( x)$,
there are exactly $s - \rho$ that are mutually independent, hence
there are exactly $\rho$ independent infinitesimal transformations
$Z_1f, \dots, Z_\rho f$ which are interchangeable with all $X_kf$.

We can therefore state the following proposition:

\def\theproposition{5}\begin{proposition}
If the $r$-term group: $X_1f, \dots, X_rf$ of the space $x_1, \dots,
x_s$ is \terminology{transitive}, and if $P$ and $P_1$ are two points
which admit exactly the same infinitesimal transformations of the
group: $X_1f, \dots, X_rf$, then there is one and only one
transformation $\Theta$ which is interchangeable with all
transformations of the group: $X_1f, \dots, X_rf$ and which transfers
$P$ to $P_1$. If one understands by $T$ the general symbol of a
transformation of the group: $X_1f, \dots, X_rf$, then $\Theta$ can be
defined as the transformation which transfers every point $(P)\, T$ to
the point $(P_1)\, T$. If each time exactly $\infty^\rho$ different
points admit precisely the same infinitesimal transformations of the
group: $X_1f, \dots, X_rf$, then there are exactly $\infty^\rho$
different transformations which are interchangeable with all
transformations of the group $X_1f, \dots, X_rf$.
\end{proposition}

The developments about simply transitive groups that we have given
in Chap.~\ref{kapitel-20}, p.~\pageref{S-390-bis}--\pageref{S-395-bis}
are obviously contained as a special case of the developments
carried out just now.

\smallercharacters{At present, we turn to the case where
the $r$-term systatic group $X_1f, \dots, X_rf$
is \emphasis{intransitive}.
However, we want to be brief here.

If the systatic group: $X_1f, \dots, X_rf$ is intransitive, then
there is not only a single transformation which transfers
the point $P$ to the point defined on p.~\pageref{S-511}
and which is interchangeable with all transformations of our group, 
and rather, there are infinitely many different transformations
of this kind. We will indicate how one finds such transformations.

The intransitive systatic group: $X_1f, \dots, X_rf$ determines
several invariant decompositions of the space $x_1, \dots, x_s$. 

A first decomposition is represented by the $s - \rho < s$
equations:
\def\theequation{a}\begin{equation}
\varphi_1(x_1,\dots,x_s)
=
{\rm const.},
\,\,\,\dots,\,\,\,
\varphi_{s-\rho}(x_1,\dots,x_s)
=
{\rm const.}
\end{equation}
A second decomposition is determined by $s - n > 0$ arbitrary 
independent solutions: $u_1 ( x), \dots, u_{ s - n} (x)$
of the $n$-term complete system: $X_1f = 0$, \dots, $X_nf = 0$; 
the analytic expression of this decomposition reads:
\def\theequation{b}\begin{equation}
u_1(x_1,\dots,x_s)
=
{\rm const.},
\,\,\,\dots,\,\,\,
u_{s-n}(x_1,\dots,x_s)
=
{\rm const.}
\end{equation}
In what follows, by $M_\rho$, we always understand one of the
$\infty^{ s - \rho}$ $\rho$-times extended manifolds~\thetag{ a}, 
and by ${\sf M}_n$, we understand one of the $\infty^{ s - n}$
$n$-times extended manifolds~\thetag{ b}.

Amongst the solutions of the complete system: $X_1f = 0$, \dots, 
$X_nf = 0$, there is a certain number, say $s - q \leqslant s 
- n$, which can be expressed in terms of $\varphi_1 ( x), \dots, 
\varphi_{ s - \rho} ( x)$ alone; we want to assume that
$u_1 ( x), \dots, u_{ s - q} (x)$ are such solutions, so that 
$s - q \leqslant s - \rho$ relations of the form:
\[
u_1(x)
=
\mathfrak{U}_1
\big(\varphi_1(x),\dots,\varphi_{s-\rho}(x)\big),
\,\,\,\dots,\,\,\,
u_{s-q}(x)
=
\mathfrak{U}_{s-q}
\big(\varphi_1(x),\dots,\varphi_{s-\rho}(x)\big)
\]
hold, and therefore, every solution of the complete system: $X_1f =
0$, \dots, $X_nf = 0$ which can be expressed in terms of $\varphi_1 (
x), \dots, \varphi_{ s - \rho} (x)$ alone is a function of $u_1 ( x),
\dots, u_{ s - q} (x)$ (cf. Chap.~\ref{kapitel-19},
p.~\pageref{S-345}). Then the $s - q$ equations:
\def\theequation{c}\begin{equation}
u_1(x_1,\dots,x_s)
=
{\rm const.},
\,\,\,\dots,\,\,\,
u_{s-q}(x_1,\dots,x_s)
=
{\rm const.}
\end{equation}
represent a third decomposition invariant by the group: $X_1f, \dots,
X_rf$. The individual manifolds of this decomposition visibly are the
smallest manifolds which consist both of the $M_\rho$ and of the ${\sf
M}_n$ (cf. Chap.~\ref{kapitel-8}, p.~\pageref{S-146}). By $\mathfrak{
M}_q$, we always understand in what follows one of the $\infty^{ s -
q}$ $q$-times extended manifolds~\thetag{ c}.

Lastly, a fourth invariant decomposition is determined by the manifold
sections \deutsch{Schnittmannigfaltigkeiten} of the $M_\rho$ and the
${\sf M}_n$ (Chap.~\ref{kapitel-8}, p.~\pageref{S-145}), and this
decomposition is obviously determined by the $s - \rho + q - n$
equations:
\def\theequation{d}\begin{equation}
\left\{
\aligned
&
\varphi_1(x)
=
{\rm const.},
\,\,\,\dots,\,\,\,
\varphi_{s-\rho}(x)
=
{\rm const.}
\\
&
u_{s-q+1}(x)
=
{\rm const.},
\,\,\,\dots,\,\,\,
u_{s-n}(x)
=
{\rm const.},
\endaligned\right.
\end{equation}
which, according to Chap.~\ref{kapitel-19}, p.~\pageref{S-345-bis}
sq., are independent of each other. In what follows, by $N_{ \rho + n
- q}$, we want to always understand one of the $\infty^{ s - \rho + q
- n}$ $(\rho + n - q)$-times extended manifolds~\thetag{ d}.

Now, in order to find a transformation $\Theta$ which is
interchangeable with all transformations of the group: $X_1f, \dots,
X_rf$, we proceed in the following way:

Inside every $\mathfrak{ M}_q$, we associate to every ${\sf M}_n$
another ${\sf M}_n$ which we may call ${\sf M}_n'$, and to be precise,
we make this association according to an arbitrary analytic
law. Afterwards, on each one of the $\infty^{ s - n}$ ${\sf M}_n$, we
choose an arbitrary point $P$, and to each one of the $\infty^{ s-n}$
chosen points, we associate an arbitrary point $P_1$ on the $N_{ \rho +
n - q}$ in which the $M_\rho$ passing through the point cuts the ${\sf
M}_n'$ which corresponds to the ${\sf M}_n$ passing through the point.

There is one and only one transformation $\Theta'$ which transfers the
$\infty^{ s - n}$ chosen points $P$ to the point $P_1$ corresponding
to them. This transformation is defined by the symbolic equation:
\[
(P)\,T\,\Theta
=
(P_1)\,T, 
\]
in which $P$ is the general symbol of the $\infty^{ s - n}$ chosen
points, while $T$ is the general symbol of the $\infty^r$
transformations of our group.

One convinces oneself easily that the transformation $\Theta$ just
defined is interchangeable with all transformations of the group:
$X_1f, \dots, X_rf$ and that one obtains all transformations $\Theta$
of this constitution when one chooses in the most general way the
arbitrary elements which are contained in the definition of $\Theta$.

We need not to spend time for proving that; let it only be remarked
that when one sets up the analytic expression of the transformation
$\Theta$, one sees immediately that the number of the arbitrary
functions appearing in this expression and the number of the arguments
appearing in these functions agree with the Theorem~67,
p.~\pageref{Theorem-67-S-376}.

Finally, yet a remark which applies both to the intransitive and to
the the transitive systatic groups: When the manifold:
\[
\varphi_{k\nu}(x_1,\dots,x_s)
=
\varphi_{k\nu}(x_1^0,\dots,x_s^0)
\ \ \ \ \ \ \ \ \ \ \ \ \ 
{\scriptstyle{(k\,=\,1\,\cdots\,r\,-\,n\,;\,\,\,
\nu\,=\,1\,\cdots\,n)}}
\]
decomposes for every system of values $x_1^0, \dots, x_s^0$
in several discrete manifolds, then the totality of all 
transformations which are interchangeable with the transformations
of the group: $X_1f, \dots, X_rf$ also decompose in several
discrete families.

}

\sectionengellie{\S\,\,\,125.}

The functions $\varphi_{ k\nu} (x)$ which decide whether a group is
systatic or not have also played a great rôle already in the chapter
about the similarity of $r$-term groups. At present, we want to go
back to the developments of that time, and we want to complete them in
a certain direction.

If, in the same number of variables, two $r$-term groups are
presented:
\[
X_kf
=
\sum_{i=1}^s\,\xi_{ki}(x_1,\dots,x_s)\,
\frac{\partial f}{\partial x_i}
\ \ \ \ \ \ \ \ \ \ \ \ \ {\scriptstyle{(k\,=\,1\,\cdots\,r)}}
\]
and:
\[
Y_kf
=
\sum_{i=1}^s\,\eta_{ki}(y_1,\dots,y_s)\,
\frac{\partial f}{\partial y_i}
\ \ \ \ \ \ \ \ \ \ \ \ \ {\scriptstyle{(k\,=\,1\,\cdots\,r)}},
\]
and if at the same time the relations: 
\[
\leftbracket
X_i,\,X_k
\rightbracket
=
\sum_{\sigma=1}^r\,c_{iks}\,X_\sigma f
\ \ \ \ \ \ \ \ \ \
\text{\rm and}
\ \ \ \ \ \ \ \ \ \
\leftbracket
Y_i,\,Y_k
\rightbracket
=
\sum_{\sigma=1}^r\,c_{iks}\,Y_\sigma f,
\]
hold with the same constants $c_{ iks}$ in the two cases, then
according to Chap.~\ref{kapitel-19}, p.~\pageref{Theorem-65-S-353} sq.,
there is a transformation:
\[
y_\nu
=
\Phi_\nu(x_1,\dots,x_s)
\ \ \ \ \ \ \ \ \ \ \ \ \ {\scriptstyle{(\nu\,=\,1\,\cdots\,s)}}
\]
which transfers $X_1f, \dots, X_rf$ to $Y_1f, \dots, Y_rf$,
respectively, if an only if the following conditions are satisfied:
if, between $X_1f, \dots, X_rf$, there are relations of the form:
\[
X_{n+k}f
=
\sum_{\nu=1}^n\,\varphi_{k\nu}(x_1,\dots,x_s)\,
X_\nu f
\ \ \ \ \ \ \ \ \ \ \ \ \ {\scriptstyle{(k\,=\,1\,\cdots\,r\,-\,n)}},
\]
while $X_1f, \dots, X_rf$ are linked together by no linear relation, 
then between $Y_1f, \dots, Y_rf$, there must exist analogous relations:
\[
Y_{n+k}f
=
\sum_{\nu=1}^n\,\psi_{k\nu}(y_1,\dots,y_s)\,
Y_\nu f
\ \ \ \ \ \ \ \ \ \ \ \ \ {\scriptstyle{(k\,=\,1\,\cdots\,r\,-\,n)}}
\]
but $Y_1f, \dots, Y_nf$ should not be linked together by linear 
relations; in addition, the $n ( r - n)$ equations:
\def\theequation{e}\begin{equation}
\varphi_{k\nu}(x_1,\dots,x_s)
=
\psi_{k\nu}(y_1,\dots,y_s)
\ \ \ \ \ \ \ \ \ \ \ \ \ 
{\scriptstyle{(k\,=\,1\,\cdots\,r\,-\,n\,;\,\,\,
\nu\,=\,1\,\cdots\,n)}}
\end{equation}
should neither contradict with each other, nor provide relations
between the $x$ alone or the $y$ alone.

If, amongst the functions $\varphi_{ k\nu} (x)$, there would be
present less than $s$ that are mutually independent, then the
determination of a transformation which transfers $X_1f, \dots, X_rf$
to $Y_1f, \dots, Y_rf$ would require certain integrations; by
contrast, if the number of independent functions $\varphi_{ k\nu} (x)$
would be equal to $s$, the equations~\thetag{ e} would represent by
themselves a transformation transferring $X_1f, \dots, X_rf$ to $Y_1f,
\dots, Y_rf$, respectively, and in fact, the most general
transformation of this nature. Hence, if we remember that in the
latter case the group: $X_1f, \dots, X_rf$ is asystatic, and naturally
also the group: $Y_1f, \dots, Y_rf$, then we obtain the:

\def\theproposition{6}\begin{proposition}
If one knows that two $r$-term asystatic groups in $s$ variables
are similar and if one has already chosen, in each one of the
two groups, $r$ infinitesimal transformations:
\[
X_kf
=
\sum_{i=1}^s\,\xi_{ki}(x_1,\dots,x_s)\,
\frac{\partial f}{\partial x_i}
\ \ \ \ \ \ \ \ \ \ \ \ \ {\scriptstyle{(k\,=\,1\,\cdots\,r)}}
\]
and:
\[
Y_kf
=
\sum_{i=1}^r\,\eta_{ki}(y_1,\dots,y_s)\,
\frac{\partial f}{\partial y_i}
\ \ \ \ \ \ \ \ \ \ \ \ \ {\scriptstyle{(k\,=\,1\,\cdots\,r)}}
\]
such that there exists a transformation: $y_i = \Phi_i ( x_1, \dots,
x_s)$ which transfers $X_1f, \dots, X_rf$ to $Y_1f, \dots, Y_rf$,
respectively, then one can set up without integration the most
general transformation that achieves the concerned transfer; this most
general transformation contains neither arbitrary functions, nor
arbitrary parameters.
\end{proposition}

From this, it follows that one can find without integration the
most general transformation which actually transfers the asystatic
group $X_1f, \dots, X_rf$ to the group: $Y_1f, \dots, Y_rf$ 
similar to it. To this end, one has one has to proceed as follows:

One determines in the group: $Y_1f, \dots, Y_rf$ in the most
general way $r$ independent infinitesimal transformations:
\[
{\sf Y}_jf
=
\sum_{k=1}^r\,\overline{g}_{jk}\,Y_kf
\ \ \ \ \ \ \ \ \ \ \ \ \ {\scriptstyle{(j\,=\,1\,\cdots\,r)}}
\]
such that firstly, the relations:
\[
\leftbracket
{\sf Y}_i,\,{\sf Y}_k
\rightbracket
=
\sum_{\sigma=1}^r\,c_{iks}\,{\sf Y}_\sigma f
\]
hold, and secondly such that there exists a transformation which
transfers $X_1f, \dots, X_rf$ to $Y_1f, \dots, Y_rf$,
respectively. Then according to what has been said a short while ago,
one can find without integration the most general transformation which
achieves the transfer in question, and as a result, one obtains at the
same time the most general transformation which converts the group:
$X_1f, \dots, X_rf$ into the group: $Y_1f, \dots, Y_rf$. Obviously,
this transformation contains only arbitrary parameters.

In particular, if the group: $Y_1f, \dots, Y_rf$ coincides with the
group: $X_1f, \dots, X_rf$, then in the way indicated, one obtains all
transformations which leave invariant the group: $X_1f, \dots,
X_rf$. According to Chap.~\ref{kapitel-19}, p.~\pageref{S-361}, the
totality of all these transformations forms a group, and in fact in
our case, visibly a finite group. Thus:

\renewcommand{\thefootnote}{\fnsymbol{footnote}}
\def\thetheorem{90}\begin{theorem}
The largest subgroup in which an $r$-term asystatic group: $X_1f,
\dots, X_rf$ of the space $x_1, \dots, x_s$ is contained as an
invariant subgroup contains only a finite number of parameters. One
can find the finite equations of this group without integration as
soon as the infinitesimal transformations of the group: $X_1f, \dots,
X_rf$ are given.\footnote[1]{\,
\name{Lie}, Archiv for Math. Vol. 10, p.~378, Christiania 1885.
}
\end{theorem}
\renewcommand{\thefootnote}{\arabic{footnote}}

It is of importance that one can also find without integration
the finite equations of the asystatic group: $X_1f, \dots, X_rf$ 
itself, as soon as its infinitesimal transformations are given.

One simply sets up the finite equations:
\def\theequation{f}\begin{equation}
e_k'
=
\sum_{j=1}^r\,\psi_{kj}
(\varepsilon_1,\dots,\varepsilon_r)\,e_j
\ \ \ \ \ \ \ \ \ \ \ \ \ {\scriptstyle{(k\,=\,1\,\cdots\,r)}}
\end{equation}
of the adjoint group associated to the group: $X_1f, \dots, X_rf$;
this demands only executable operations (cf. Chap.~\ref{kapitel-16},
p.~\pageref{S-273-bis}). Next, if the sought finite equations of the
group: $X_1f, \dots, X_rf$ have the form: $x_i' = f_i ( x_1, \dots,
x_n,\, \varepsilon_1, \dots, \varepsilon_r)$, then according to
Theorem~48, p.~\pageref{Theorem-48-S-275}, after the introduction of
the new variables: $x_i' = f_i ( x, \varepsilon)$, the infinitesimal
transformations $X_kf$ take the form:
\[
X_kf
=
\sum_{j=1}^r\,\psi_{jk}
(\varepsilon_1,\dots,\varepsilon_r)\,
X_j'f
\ \ \ \ \ \ \ \ \ \ \ \ \ {\scriptstyle{(k\,=\,1\,\cdots\,r)}},
\]
where, as usual, we have set:
\[
\sum_{i=1}^n\,\xi_{ki}(x_1',\dots,x_s')\,
\frac{\partial f}{\partial x_i'}
=
X_k'f.
\]
Now, since the group: $X_1f, \dots, X_rf$ is asystatic, there is
a completely determined transformation between the $x$ and the $x'$
that transfers $X_1f, \dots, X_rf$ to:
\[
\sum_{j=1}^r\,\psi_{jk}
(\varepsilon_1,\dots,\varepsilon_r)\,
X_j'f
\ \ \ \ \ \ \ \ \ \ \ \ \ {\scriptstyle{(k\,=\,1\,\cdots\,r)}},
\]
respectively. If one computes this transformation according to
the former rules, one finds the sought equations $x_i' = f_i ( x,
\varepsilon )$. 

In particular, if the equations~\thetag{ f} are a canonical form of
the adjoint group (cf. Chap.~\ref{kapitel-9}, p.~\pageref{S-171}),
then evidently, one obtains the finite equations of the group: $X_1f,
\dots, X_rf$ also in canonical form. We therefore have the

\def\theproposition{7}\begin{proposition}
If one knows the infinitesimal transformations of an asystatic group
of the space $x_1, \dots, x_s$, then one can always find the finite
equations of this group by means of executable operations and to be
precise, in canonical form.
\end{proposition}

There exist yet more general cases for which the finite equations
of an $r$-term group, the infinitesimal transformations of which
one knows, can be determined without integration. However, we 
do not want to be involved further in such questions, and we
only want to remark that the determination of the finite equations
succeeds, amongst other circumstances, when there is no infinitesimal
transformation interchangeable with all $X_kf$ which does not 
belongs to the group $X_1f, \dots, X_rf$. 

\sectionengellie{\S\,\,\,126.}

Let $X_1f, \dots, X_rf$, or shortly $G_r$, be an $r$-term
\emphasis{systatic} group of the space $x_1, \dots, x_s$, and let $G_{
r - n}$ be the $(r-n)$-term subgroup of the $G_r$ which is associated
to a determined point $x_1^0, \dots, x_s^0$ in general position. The
manifold:
\[
\varphi_{k\nu}(x_1,\dots,x_s)
=
\varphi_{k\nu}(x_1^0,\dots,x_s^0)
\ \ \ \ \ \ \ \ \ \ \ \ \ 
{\scriptstyle{(k\,=\,1\,\cdots\,r\,-\,n\,;\,\,\,
\nu\,=\,1\,\cdots\,n)}}
\]
which consists of all points invariant by the $G_{ r - n}$ may be
denoted by $M$.

Since the $G_{ r - n}$ fixes all points of $M$, it naturally leaves
invariant $M$ itself; but is it thinkable that the $G_r$ contains
transformations which also leave invariant the manifold $M$ without
fixing all of its points. We want to assume that the largest subgroup
of the $G_r$ which leaves $M$ invariant contains exactly $r - l$
parameters and we want to call this subgroup $G_{ r - l}$.

Of course, the $G_{ r - n}$ is either identical to the $G_{ r - l}$ or
contained in it as a subgroup. The latter case occurs always when the
$G_r$ is transitive; indeed, in this case, the point $x_1^0, \dots,
x_s^0$ can be transferred to all points of $M$ by means of suitable
transformations of the $G_r$, and since every transformation of the
$G_r$ which transfers $x_1^0, \dots, x_s^0$ to another point of $M$
visibly leaves invariant the manifold $M$, the $G_r$ contains a
continuous family of transformations which leave invariant
$M$ without fixing all of its points; consequently, in the case
of a transitive group $G_r$, the number $r - l$ is surely larger
than $r - n$. But if the $G_r$ is intransitive, then it is very 
well possible that all transformations of the $G_r$ which leave
$M$ invariant also fix all points of $M$, so that $r - l = r - n$.
This is shown for example by the three-term intransitive systatic
group:
\[
\frac{\partial f}{\partial x_2},
\ \ \ \ \ \ \ \
x_2\,\frac{\partial f}{\partial x_2},
\ \ \ \ \ \ \ \
x_2^2\,\frac{\partial f}{\partial x_2}
\]
of the plane $x_1$, $x_2$.

It is easy to see that the $G_{ r - n}$ is invariant in the $G_{ r -
l}$. In fact, the $G_{ r - n}$ is generated by all infinitesimal
transformations of the $G_{ r - l}$ which leave untouched all the
points of $M$; but according to Chap.~\ref{kapitel-17}, Proposition~7,
p.~\pageref{Satz-7-S-309}, these infinitesimal transformations
generate an invariant subgroup of the $G_{ r - n}$. We therefore
have the:

\def\theproposition{8}\begin{proposition}
\label{Satz-8-S-519}
If: $X_1f, \dots, X_rf$, or $G_r$, is an $r$-term \terminology{systatic}
group of the space $x_1, \dots, x_s$, if moreover $P$ is a point in 
general position, and lastly, if $M$ is the manifold of all points
that remain invariant by all transformations of the $G_r$
which fix $P$, then the largest subgroup of the $G_r$ which leaves
$P$ invariant is either identical to the largest subgroup
which leaves $M$ at rest, or is contained in this subgroup as 
an invariant subgroup. The first case can occur only when the $G_r$ is
intransitive; when the $G_r$ is transitive, the second case
always occurs.
\end{proposition}

On the other hand, if one knows an arbitrary $r$-term group $G_r$ of
the space $x_1, \dots, x_s$ and if one knows that the largest subgroup
$G_{ r - n}$ of the $G_r$ which leaves invariant an arbitrary point
$P$ in general position is invariant in a yet larger subgroup $G_{ r -
h}$ with $r - h > r -n$ parameters, then one can conclude that the
$G_r$ belongs to the class of the systatic groups.

In fact, the point $P$ admits exactly $r - n$ independent
infinitesimal transformations of the $G_{ r -h}$ and hence
(Chap.~\ref{kapitel-23}, Theorem~85, p.~\pageref{Theorem-85-S-483}),
it takes, by all $\infty^{ r - h}$ transformations of the $G_{ r -
h}$, exactly $\infty^{ n - h}$ different positions, where the number
$n -h$, under the assumption made, is at least equal to $1$. Now,
since the $G_{ r - n}$ is invariant in the $G_{ r -h}$, then to each
of these $\infty^{ n -h}$ positions of $P$ is associated exactly the
same $(r - n)$-term subgroup of the $G_r$ as to the point $P$, that is
to say: the $G_r$ is effectively systatic. Thus:

\def\theproposition{9}\begin{proposition}
If the $r$-term group $G_r$ of the space $x_1, \dots, x_s$ is
constituted in such a way that its largest subgroup $G_{ r -n}$ which
leaves invariant an arbitrary point in general position is invariant
in a larger subgroup of the $G_r$, or even in the $G_r$ itself, then
the $G_r$ belongs to the class of the systatic groups.
\end{proposition}

From this proposition, it follows immediately that for an $r$-term
asystatic group of the space $x_1, \dots, x_s$, the subgroup
associated to a point in general position is invariant neither in the
$G_r$ itself, nor in a larger subgroup of the $G_r$. But conversely,
if the $G_r$ is constituted in such a way that the subgroup which is
associated to a point in general position is invariant in no larger
subgroup, then it needs not be asystatic for this reason, and it
necessarily so only when it is at the same time also transitive; this
follows immediately from the Proposition~8,
p.~\pageref{Satz-8-S-519}. We can therefore say:

\def\theproposition{10}\begin{proposition}
An $r$-term \terminology{transitive} group $G_r$ of the space $x_1,
\dots, x_s$ is asystatic if and only if the subgroup associated
to a point in general position is invariant in no larger subgroup.
\end{proposition}

In Chap.~\ref{kapitel-22}, Theorem~79, p.~\pageref{Theorem-79-S-443},
we provided a method for the determination of all
\emphasis{transitive} groups which are equally composed with a given
$r$-term group $\Gamma$. At present, we can specialize this method so
that it provides in particular all \emphasis{transitive asystatic}
groups which are equally composed with the group $\Gamma$.

One determines all subgroups of $\Gamma$ which neither contain a
subgroup invariant in $\Gamma$, nor are invariant in a larger subgroup
of $\Gamma$. Each one of these subgroups provides, when one
proceeds according to the rules of Theorem~78, a transitive asystatic
group which is equally composed with $\Gamma$, and in fact,
one finds in this way all transitive asystatic groups having the
concerned composition.

\smallercharacters{
Taking the Proposition~10 as a basis, one can answer the question
whether a given transitive group is systatic, or asystatic.
One can set up a similar statement by means of which one can
answer the question \emphasis{whether a given transitive group is
primitive, or not}. 

Let the $r$-term \emphasis{transitive} group $G_r$ of the space
$x_1, \dots, x_s$ be imprimitive, and let:
\[
u_1(x_1,\dots,x_s)
=
{\rm const.},
\,\,\,\dots,\,\,\,
u_{s-m}(x_1,\dots,x_s)
=
{\rm const.}
\ \ \ \ \ \ \ \ \ \ \ \ \ 
{\scriptstyle{(0\,<\,m\,<\,s)}}
\]
be a decomposition of the space in $\infty^{ s -m}$ $m$-times
extended manifolds which is invariant by the group.

Since the $G_r$ is transitive, all its transformations which leave
invariant a point $x_1^0, \dots, x_s^0$ in general position form 
an $( r - s)$-term subgroup $G_{ r -s}$, and on the other hand, 
all its transformations which leave invariant the $m$-times extended
manifold:
\[
u_1(x_1,\dots,x_s)
=
u_1(x^0),
\,\,\,\dots,\,\,\,
u_{s-m}(x_1,\dots,x_s)
=
u_{s-m}(x^0)
\]
form an $(r - s + m)$-term subgroup $G_{ r - s + m}$ (cf. 
Chap.~\ref{kapitel-23}, p.~\pageref{Satz-7-S-478}). Here naturally,
the $G_{ r -s}$ is contained as subgroup in the $G_{ r - s + m}$, 
which in turn obviously is an actual subgroup of the $G_r$. 

\medskip

At present, we imagine conversely that an arbitrary $r$-term transitive
group $\mathfrak{ G}_r$ of the space $x_1, \dots, x_s$ is presented,
and we assume that the $(r - s)$-term subgroup $\mathfrak{ G}_{ r -s}$
of the $\mathfrak{ G}_r$ which is associated to a point $P$ in general
position is contained in a larger subgroup:
\[
\mathfrak{G}_{r-s+h}
\ \ \ \ \ \ \ \ \ \ \ \ \ 
{\scriptstyle{(r\,-\,s\,<\,r\,-\,s\,+\,h\,<\,r)}}.
\]

If all transformations of the $\mathfrak{ G}_{ r - s + h}$ are
executed on $P$, then the point takes $\infty^h$ different positions.
These $\infty^h$ positions form an $h$-times extended manifold $M$
which remains invariant by the $\mathfrak{ G}_{ r - s + h}$
(cf. Chap.~\ref{kapitel-23}, p.~\pageref{Theorem-85-S-483}). Lastly,
if we execute on $M$ all transformations of the $\mathfrak{ G}_r$, we
obtain $\infty^{ s - h}$ different $h$-times extended manifolds which
determine a decomposition of the space $x_1, \dots, x_s$.

In fact, since the manifold $M$ remains invariant by the $\mathfrak{
G}_{ r - s + h}$, then by the $\infty^r$ transformations of the
$\mathfrak{ G}_r$, it takes at most $\infty^{ s -h}$ different
positions (cf. Chap.~\ref{kapitel-23}, p.~\pageref{Theorem-85-S-483}),
and on the other hand, thanks to the transitivity of the $\mathfrak{
G}_r$, it takes at least $\infty^{ s - h}$ positions, hence it
receives by the $\mathfrak{ G}_r$ exactly $\infty^{ s - h}$ different
positions which fill exactly the space and hence determine a
decomposition. It follows from Theorem~85,
p.~\pageref{Theorem-85-S-483}, that this decomposition remains
invariant by the $\mathfrak{ G}_r$.

From this, it results that the $\mathfrak{ G}_r$, under the
assumptions made, is imprimitive. If we combine this result with the
one gained above, we then obtain at first the:

\def\thetheorem{91}\begin{theorem}
\label{Theorem-91-S-521}
An $r$-term transitive group $G_r$ of the space $x_1, \dots, x_s$ is
primitive if and only if the $(r - s)$-term subgroup $G_{ r - s}$
which is associated to a point in general position is contained in no
larger subgroup of the $G_r$.
\end{theorem}

But even more, we obviously obtain at the same time a method for the
determination of all possible decompositions of the space $x_1, \dots,
x_s$ which remain invariant by a given transitive group of this space:

\def\thetheorem{92}\begin{theorem}
If an $r$-term transitive group $G_r$ of the space $x_1, \dots, x_s$
is presented, then one finds all possible decompositions of the space
invariant by the group in the following way:

One determines at first the $(r - s)$-term subgroup $G_{ r - s}$ of
the $G_r$ by which an arbitrary point $P$ in general position remains
invariant, a point which does not lie on any manifold invariant by the
$G_r$. Afterwards, one looks up at all subgroups of the $G_r$ which
contain the $G_{ r - s}$. If $G_{ r - s + h}$ is one of these
subgroups, then one executes all transformations of the $G_{ r - s +
h}$ on $P$; in the process, $P$ takes $\infty^h$ different positions
which form an $h$-times extended manifold $M$; now, if one executes
on $M$ all transformations of the $G_r$, then $M$ takes $\infty^{ s
- h}$ different positions which determine a decomposition of the
space invariant by the $G_r$. If one treats in this way all subgroups
of the $G_r$ which comprise the $G_{ r - s}$, one obtains all
decompositions invariant by the $G_r$. 
\end{theorem}

At present, we can also indicate how one has to specialize the method
explained in Chap.~\ref{kapitel-22}, Theorem~79,
p.~\pageref{Theorem-79-S-443} in order to find all
\emphasis{primitive} groups which are equally composed with a given
$r$-term group $\Gamma$.

Since all primitive groups are transitive, one has to proceed in the
following way:

One determines all subgroups of $\Gamma$ which contain no subgroup
invariant in $\Gamma$ and which are contained in no larger subgroup of
$\Gamma$. Each of these subgroups provides, when one proceeds
according to the rules of Theorem~78, a primitive group equally
composed with $\Gamma$ and in fact, one obtains in this way all
primitive groups of this sort.\,---

The Theorem~92 shows that one can find without integration all
decompositions invariant by a \emphasis{transitive} group when the
finite equations of the group are known. Now, the same also holds true
for intransitive groups, although in order to see this, one needs
considerably longer considerations that might
not be advisable to develop here.

}

\linestop


\chapter{Differential Invariants}
\label{kapitel-25}
\chaptermark{Differential Invariants}

\setcounter{footnote}{0}

\abstract*{??}

\label{S-522-sq}
In $n$ variables, we imagine that a transformation is presented:
\def\theequation{1}\begin{equation}
y_i
=
f_i(x_1,\dots,x_n)
\ \ \ \ \ \ \ \ \ \ \ \ \ {\scriptstyle{(i\,=\,1\,\cdots\,n)}},
\end{equation}
and we imagine that this transformation is executed on a system
of $n - q$ independent equations:
\def\theequation{2}\begin{equation}
\Omega_k(x_1,\dots,x_n)
=
0
\ \ \ \ \ \ \ \ \ \ \ \ \ 
{\scriptstyle{(k\,=\,1\,\cdots\,n\,-\,q)}}\,;
\end{equation}
we obtain in this way a new system of equations:
\def\theequation{2'}\begin{equation}
\overline{\Omega}(y_1,\dots,y_n)
=
0
\ \ \ \ \ \ \ \ \ \ \ \ \ 
{\scriptstyle{(k\,=\,1\,\cdots\,n\,-\,q)}}
\end{equation}
in the $y$.

Now, by virtue of the equations~\thetag{ 2}, $n - q$ of the variables
$x_1, \dots, x_n$ can be represented as functions of the $q$ remaining
ones and these $n - q$ naturally possess certain differential
quotients: $p_1$, $p_2$, \dots, with respect to the $q$ remaining.
On the other hand, by virtue of the equations~\thetag{ 2'}, $n - q$ of
the variables $y_1, \dots, y_n$ can be represented as functions of the
$q$ remaining ones and these $n - q$ possess in their turn certain
differential quotients: $p_1'$, $p_2'$, \dots, with respect to the $q$
remaining.

Between the two series of differential quotients so defined, there
exists a certain connection which is essentially independent of the
form of the $n - q$ relations~\thetag{ 2}. In the course of the
chapter, we will explain thoroughly this known connection, and we
recall here that the $p'$ from the first order up to the $m$-th order
can be represented as functions of $x_1, \dots, x_n$ and of the $p$
from the first order up to the $m$-th order, when the $p$ are
interpreted as functions of $x_1', \dots, x_n'$ and of the $p'$. From
this, it results that one can derive from the transformation~\thetag{
1} a new transformation which, aside from the $x$, also transforms the
differential quotients from the first order up to the $m$-th order.
We want to say that \emphasis{this new transformation is obtained by
\terminology{prolongation} \deutsch{Erweiterung} of the
transformation~\thetag{ 1}}.

Visibly, the prolongation of the transformation~\thetag{ 1} can take
place in several very varied ways, because it left just as one likes
how many and which ones of the variables $x_1, \dots, x_n$ will be
seen as functions of the others. One can even increase yet the number
of possibilities by taking in addition auxiliary variables: $t_1$,
$t_2$, \dots, that are absolutely not transformed by the
transformation~\thetag{ 1}, so that one adds the identity
transformation in the auxiliary variables: $t_1$, $t_2$, \dots, to the
equations of the transformation~\thetag{ 1}.

If an $r$-term group is presented, we can prolong \deutsch{erweitern}
all its $\infty^r$ transformations and obtain in this way $\infty^r$
prolonged transformations; the latter transformations constitute in
their turn, as we will see, an $r$-term group which is equally
composed with the original group.

\sectionengellie{\S\,\,\,127.}

We consider at first a special simple sort of the prolongation.

In the transformation:
\def\theequation{1}\begin{equation}
\label{S-524-sq}
y_i
=
f_i(x_1,\dots,x_n)
\ \ \ \ \ \ \ \ \ \ \ \ \ {\scriptstyle{(i\,=\,1\,\cdots\,n)}},
\end{equation}
we consider the variables $x_1, \dots, x_n$ as functions of an 
auxiliary variable $t$ which is absolutely not transformed
by the transformation~\thetag{ 1}. Obviously, $y_1, \dots, y_n$
are then to be interpreted also as functions of $t$; hence if we set:
\[
\frac{\D\,x_i}{\D\,t}
=
x_i^{(1)},
\ \ \ \ \ \ \ \ \ \
\frac{\D\,y_i}{\D\,t}
=
y_i^{(1)},
\]
then it follows from~\thetag{ 1} by differentiation with respect to
$t$:
\def\theequation{3}\begin{equation}
y_i^{(1)}
=
\sum_{\nu=1}^n\,
\frac{\partial f_i}{\partial x_\nu}\,
x_\nu^{(1)}
\ \ \ \ \ \ \ \ \ \ \ \ \ {\scriptstyle{(i\,=\,1\,\cdots\,n)}}.
\end{equation}
If we take together~\thetag{ 1} and~\thetag{ 3}, we obtain
a transformation in the $2\, n$ variables $x_i$ and
$x_i^{ (1)}$. 

The transformation~\thetag{ 1}, \thetag{ 3} is in a certain 
sense the simplest transformation that one can derive
from~\thetag{ 1} by prolongation.

If we apply the special prolongation written just now to all
$\infty^r$ transformations of the $r$-term group: $X_1f, \dots, 
X_rf$, or to:
\def\theequation{4}\begin{equation}
y_i
=
f_i(x_1,\dots,x_n,\,a_1,\dots,a_r)
\ \ \ \ \ \ \ \ \ \ \ \ \ {\scriptstyle{(i\,=\,1\,\cdots\,n)}},
\end{equation}
then we obtain $\infty^r$ prolonged transformations which have
evidently the form:
\def\theequation{5}\begin{equation}
y_i
=
f_i(x_1,\dots,x_n,\,a_1,\dots,a_r),
\ \ \ \ \ \ \ \
y_i^{(1)}
=
\sum_{\nu=1}^n\,
\frac{\partial f_i(x,a)}{\partial x_\nu}\,
x_\nu^{(1)}
\ \ \ \ \ \ \ \ \ {\scriptstyle{(i\,=\,1\,\cdots\,n)}}.
\end{equation}

It can be proved that the equations~\thetag{ 5} represent an
$r$-term group in the $2\, n$ variables $x_i$, $x_i^{ (1)}$.

In fact, the equations:
\[
y_i
=
f_i(x_1,\dots,x_n,\,a_1,\dots,a_r)
\ \ \ \ \ 
\text{\rm and}
\ \ \ \ \
z_i
=
f_i(y_1,\dots,y_n,\,b_1,\dots,b_r)
\]
have by assumption as a consequence:
\[
z_i
=
f_i(x_1,\dots,x_n,\,c_1,\dots,c_r),
\]
where the $c$ depend only on the $a$ and on the $b$. 
Therefore, from:
\[
y_i^{(1)}
=
\sum_{\nu=1}^n\,
\frac{\partial f_i(x,a)}{\partial x_\nu}\,
x_\nu^{(1)}
\ \ \ \ \ 
\text{\rm and}
\ \ \ \ \
z_j^{(1)}
=
\sum_{i=1}^n\,
\frac{\partial f_j(y,b)}{\partial y_i}\,
y_i^{(1)},
\]
it follows by elimination of the $y_i^{ (1)}$:
\[
\aligned
z_j^{(1)}
&
=
\sum_{i,\,\,\nu}^{1\cdots\,n}\,
\frac{\partial f_j(y,b)}{\partial y_i}\,
\frac{\partial f_i(x,a)}{\partial x_\nu}\,
x_\nu^{(1)}
\\
&
=
\sum_{\nu=1}^n\,
\frac{\partial f_j(y,b)}{\partial x_\nu}\,
x_\nu^{(1)}
=
\sum_{\nu=1}^n\,
\frac{\partial f_j(x,c)}{\partial x_\nu}\,
x_\nu^{(1)},
\endaligned
\]
whence the announced proof is achieved.

In addition, from what has been said, it results that the new
group possesses the same
\label{S-525} parameter group as the original group:
\[
y_i
=
f_i(x_1,\dots,x_n,\,a_1,\dots,a_r).
\]

The new group:
\def\theequation{5}\begin{equation}
y_k
=
f_k(x,a),
\ \ \ \ \ \ \ \ \ \
y_k^{(1)}
=
\sum_{i=1}^n\,
\frac{\partial f_k}{\partial x_i}\,
x_i^{(1)}
\ \ \ \ \ \ \ \ \ \ \ \ \ {\scriptstyle{(k\,=\,1\,\cdots\,n)}}
\end{equation}
which shows a first example of a \emphasis{prolonged
group} has a very simple conceptual meaning.

Indeed, if one interprets $x_1, \dots, x_n$ as ordinary point 
coordinates of an $n$-times extended space, then one can
interpret the $2n$ quantities: $x_1, \dots, x_n$, 
$x_1^{ (1)}, \dots, x_n^{ (1)}$ as the determination pieces
\deutsch{Bestimmungsstücke} of a line element \deutsch{Linienelement};
about it, $x_1^{ (1)}, \dots, x_n^{ (1)}$ are homogeneous
coordinates in the domain of the $\infty^{ n-1}$ directions
which pass through the point $x_1, \dots, x_n$. 

\plainstatement{The new group~\thetag{ 5} 
\label{S-525-bis}
in the $x$, $x^{ (1)}$
indicates in which way the line elements of the space $x_1, \dots, x_n$
of the original group: $y_i = f_i ( x, a)$ are permuted
with each other.}

Since the group: $y_i = f_i ( x, a)$ is generated by the $r$ 
infinitesimal transformations: $X_1f, \dots, X_rf$, its finite
equations in canonical form read:
\[
y_k
=
x_k
+
\sum_{j=1}^r\,e_j\,\xi_{jk}
+\cdots
\ \ \ \ \ \ \ \ \ \ \ \ \ {\scriptstyle{(i\,=\,1\,\cdots\,n)}}.
\]
If we set this form of the group: $y_i = f_i ( x, a)$ as fundamental
in order to make up the prolonged group~\thetag{ 5}, then the 
equations of this group receive the shape:
\def\theequation{6}\begin{equation}
\aligned
y_k
=
x_k
+
\sum_{j=1}^r\,e_j\,\xi_{jk}
&
+\cdots,
\ \ \ \ \ \ \ \
y_k^{(1)}
=
x_k^{(1)}
+
\sum_{j=1}^r\,e_j\,\xi_{jk}^{(1)}
+\cdots
\\
&
\ \ \ \ \ \ \ \
{\scriptstyle{(k\,=\,1\,\cdots\,n)}},
\endaligned
\end{equation}
where we have set for abbreviation:
\[
\sum_{i=1}^n\,
\frac{\partial\xi_{jk}}{\partial x_i}\,
x_i^{(1)}
=
\xi_{jk}^{(1)}.
\]

From this, we realize immediately that the prolonged group contains
the identity transformation, and at the same time, we come to the
presumption that it is generated by the $r$ infinitesimal
transformations:
\[
X_k^{(1)}f
=
\sum_{i=1}^n\,\xi_{ki}\,
\frac{\partial f}{\partial x_i}
+
\sum_{i=1}^n\,\xi_{ki}^{(1)}\,
\frac{\partial f}{\partial x_i^{(1)}}
\ \ \ \ \ \ \ \ \ \ \ \ \ {\scriptstyle{(k\,=\,1\,\cdots\,r)}}.
\]
The correctness of this presumption can be established in the 
following way:

According to Theorem~3, p.~\pageref{Theorem-3-S-33},
the $n$ functions: $y_i = f_i ( x, a)$ satisfy differential 
equations of the form:
\def\theequation{7}\begin{equation}
\frac{\partial y_i}{\partial a_k}
=
\sum_{j=1}^r\,\psi_{kj}(a_1,\dots,a_r)\,
\xi_{ji}(y_1,\dots,y_n)
\ \ \ \ \ \ \ \ \ \ \ \ \ 
{\scriptstyle{(i\,=\,1\,\cdots\,n\,;\,\,\,
k\,=\,1\,\cdots\,r)}}
\end{equation}
If we differentiate these equations with respect to $t$ and
if take into account that the $a$ are independent of the $t$, 
it comes:
\def\theequation{8}\begin{equation}
\frac{\partial y_i^{(1)}}{\partial a_k}
=
\sum_{j=1}^r\,\psi_{kj}(a)\,
\sum_{\nu=1}^n\,
\frac{\partial \xi_{ji}(y)}{\partial y_\nu}\,y_\nu^{(1)}
\ \ \ \ \ \ \ \ \ \ \ \ \ 
{\scriptstyle{(i\,=\,1\,\cdots\,n\,;\,\,\,
k\,=\,1\,\cdots\,r)}}.
\end{equation}
Now, if $a_1^0, \dots, a_r^0$ are the parameters of the identity 
transformation in the group: $y_i = f_i ( x, a)$, then the determinant:
$\sum\, \pm\, \psi_{ 11} ( a) \cdots\, \psi_{ rr} (a)$ does not vanish
for: $a_1 = a_1^0$, \dots, $a_r = a_r^0$; but $a_1^0, \dots, a_r^0$
are also the parameters of the identity transformation in the
prolonged group~\thetag{ 5}; consequently, from the equations~\thetag{
7} and~\thetag{ 8}, we can conclude in the known way (cf. 
Chap.~\ref{one-term-groups}, p.~\pageref{application-theorem-9}
sq.), that the prolonged group is really generated by the
$r$ independent infinitesimal transformations $X_k^{ (1)} f$. 
We call the $X_k^{ (1)} f$ the
\emphasis{prolonged infinitesimal transformations}.

As the independent infinitesimal transformations of an $r$-term group, 
the $X_k^{ (1)}f$ satisfy pairwise relations of the form: 
\[
X_k^{(1)}\,X_j^{(1)}f
-
X_j^{(1)}\,X_k^{(1)}f
=
\sum_{s=1}^r\,c_{kjs}'\,X_s^{(1)}f.
\]
Here, the coefficients of the $\partial f /\partial
x_i$ must coincide on the right and on the left, hence 
one must have: $X_k X_j f - X_j X_k f = 
\sum_s\, c_{ kj s}' \, X_sf$, and consequently, one has:
\[
\sum_{s=1}^r\,
(c_{kjs}'-c_{kjs})\,X_sf
=
0,
\]
or, because of the independence of $X_1f, \dots, X_rf$:
\[
c_{kjs}'
=
c_{kjs}
\ \ \ \ \ \ \ \ \ \ \ \ \ 
{\scriptstyle{(k,\,\,j,\,\,s\,=\,1\,\cdots\,r)}}.
\]
Thus, the prolonged group~\thetag{ 5} is holoedrically isomorphic
to the original group, which is coherent with the fact that, 
according to p.~\pageref{S-525}, both groups have the same
parameter group.\,---

We give yet a second direct proof of the result just found.

If $X_kf$ and $X_jf$ are \emphasis{arbitrary} infinitesimal
transformations and if the associated two prolonged infinitesimal
transformations are denoted $X_k^{ (1)} f$ and
$X_j^{ (1)} f$ as earlier on, it follows that:
\[
\aligned
X_k^{(1)}X_j^{(1)}f
-
X_j^{(1)}X_k^{(1)}f
&
=
X_kX_jf
-
X_jX_kf
\\
&
\ \ \ \ \
+
\sum_{\pi,\,\,i,\,\,\nu}^{1\cdots\,n}\,
\bigg\{
\xi_{ki}\,
\frac{\partial^2\xi_{j\pi}}{\partial x_i\partial x_\nu}
-
\xi_{ji}\,
\frac{\partial^2\xi_{k\pi}}{\partial x_i\partial x_\nu}
\bigg\}\,x_\nu^{(1)}\,
\frac{\partial f}{\partial x_\pi^{(1)}}
\\
&
\ \ \ \ \
+
\sum_{\pi,\,\,i,\,\,\nu}^{1\cdots\,n}\,
\bigg\{
\frac{\partial\xi_{ki}}{\partial x_\nu}\,
\frac{\partial\xi_{j\pi}}{\partial x_i}
-
\frac{\partial\xi_{ji}}{\partial x_\nu}\,
\frac{\partial \xi_{k\pi}}{\partial x_i}
\bigg\}\,
x_\nu^{(1)}\,
\frac{\partial f}{\partial x_\pi^{(1)}},
\endaligned
\]
or that:
\def\theequation{9}\begin{equation}
\left\{
\aligned
X_k^{(1)}X_j^{(1)}f
-
X_j^{(1)}X_k^{(1)}f
&
=
\sum_{\pi,\,\,i}^{1\cdots\,n}\,
\bigg\{
\xi_{ki}\,\frac{\partial\xi_{j\pi}}{\partial x_i}
-
\xi_{ji}\,\frac{\partial\xi_{k\pi}}{\partial x_i}
\bigg\}\,
\frac{\partial f}{\partial x_\pi}
\\
&
\ \ \ \ \
+
\sum_{\pi,\,\,i}^{1\cdots\,n}\,
\frac{\D}{\D\,t}\,
\bigg\{
\xi_{ki}\,\frac{\partial\xi_{j\pi}}{\partial x_i}
-
\xi_{ji}\,\frac{\partial\xi_{k\pi}}{\partial x_i}
\bigg\}\,
\frac{\partial f}{\partial x_\pi^{(1)}}.
\endaligned\right.
\end{equation}

\renewcommand{\thefootnote}{\fnsymbol{footnote}}
In particular, if we assume that the $X_kf$ are infinitesimal
transformations of an $r$-term group, hence that:
\[
X_kX_jf
-
X_jX_kf
=
\sum_{s=1}^r\,c_{kjs}\,X_sf,
\]
then it follows that:
\[
\aligned
X_k^{(1)}X_j^{(1)}f
-
X_j^{(1)}X_k^{(1)}f
&
=
\sum_{s=1}^r\,c_{kjs}\,
\sum_{\pi=1}^n\,\xi_{s\pi}\,
\frac{\partial f}{\partial x_\pi}
\\
&
\ \ \ \ \
+
\sum_{s=1}^r\,c_{kjs}\,
\sum_{\pi=1}^n\,
\frac{\D}{\D\,t}\,
\xi_{s\pi}\,
\frac{\partial f}{\partial x_\pi}
=
\sum_{s=1}^r\,c_{kjs}\,X_s^{(1)}f,
\endaligned
\]
was what to be shown.\footnote[1]{\,
The preceding analytic developments present great similarities
with certain developments of Chap.~\ref{kapitel-20} (cf. 
p.~\pageref{S-370-sq} sq.). The inner reason of this
connection lies in the fact that the quantities
$\zeta$ on p.~\pageref{S-370} are nothing but certain differential
quotients: $\delta y / \delta t$ of the $y$ with respect to $t$.
}
\renewcommand{\thefootnote}{\arabic{footnote}}

The general formula~\thetag{ 9} in the computations executed just
now is of special interest and it can be expressed in words
as follows:

\def\theproposition{1}\begin{proposition}
\label{Satz-1-S-527}
If, for two arbitrary infinitesimal transformations:
\[
X_1f
=
\sum_{i=1}^n\,\xi_{1i}\,\frac{\partial f}{\partial x_i},
\ \ \ \ \ \ \ \ \ \ \ \
X_2f
=
\sum_{i=1}^n\,\xi_{2i}\,\frac{\partial f}{\partial x_i},
\]
one forms the prolonged infinitesimal transformations:
\[
X_k^{(1)}f
=
\sum_{i=1}^n\,\xi_{ki}\,
\frac{\partial f}{\partial x_i}
+
\sum_{i=1}^n\,\xi_{ki}^{(1)}\,
\frac{\partial f}{\partial x_i^{(1)}}
\ \ \ \ \ \ \ \ \ \ \ \ \ {\scriptstyle{(k\,=\,1,\,\,2)}},
\]
then $X_1^{ (1)} X_2^{ (1)} f - X_2^{(1)} X_1^{ (1)} f$
is the prolonged infinitesimal transformation associated to
$X_1X_2 f - X_2X_1f$. 
\end{proposition}

\sectionengellie{\S\,\,\,128.}

At present, we ask how one can decide whether the \emphasis{equation}:
\[
\sum_{\nu=1}^n\,U_\nu(x_1,\dots,x_n)\,x_\nu^{(1)}
=
0
\]
admits every finite transformation of the prolonged one-term group:
\[
X^{(1)}f
=
\sum_{i=1}^n\,\xi_i\,
\frac{\partial f}{\partial x_i}
+
\sum_{i,\,\,\nu}^{1\cdots\,n}\,
\frac{\partial\xi_i}{\partial x_\nu}\,
x_\nu^{(1)}\,
\frac{\partial f}{\partial x_i^{(1)}}.
\]

According to Chap.~\ref{kapitel-7}, p.~\pageref{Theorem-14-S-112}, this
holds true if and if and only if $X^{ (1)} ( \sum\, U_\nu \, x_\nu^{
(1)} )$ vanishes by virtue of $\sum\, U_\nu \, x_\nu^{ (1)} = 0$.
By computation, one finds: 
\[
X^{(1)}
\bigg(
\sum_{\nu=1}^n\,U_\nu\,x_\nu^{(1)}
\bigg)
=
\sum_{\nu=1}^n\,X\,U_\nu\,x_\nu^{(1)}
+
\sum_{\nu=1}^n\,U_\nu\,\sum_{i=1}^n\,
\frac{\partial\xi_\nu}{\partial x_i}\,x_i^{(1)},
\]
hence an expression which is linear in the $x_i^{ (1)}$. 
\emphasis{Consequently, the equation: $\sum\, U_\nu\, x_\nu^{ (1)} = 0$
admits the one-term group $X^{ (1)}f$ when and only
when a relation of the form:
\[
X^{(1)}
\bigg(
\sum_{\nu=1}^n\,U_\nu\,x_\nu^{(1)}
\bigg)
=
\rho\,
\sum_{\nu=1}^n\,U_\nu\,x_\nu^{(1)}
\]
holds, where it is understood that $\rho$ is a function of
$x_1, \dots, x_n$ alone.}

Incidentally, it may be observed that the \emphasis{expression}
$\sum\, U_\nu \, x_\nu^{ (1)}$ remains always invariant by every
finite transformation of the prolonged group $X^{ (1)} f$ if and only
if the expression: $X^{ (1)} \big( \sum\, U_\nu \, x_\nu^{ (1)} \big)$
vanishes identically.

Lastly, a system of $m$ equations:
\def\theequation{10}\begin{equation}
\sum_{i=1}^n\,U_{ki}(x_1,\dots,x_n)\,x_i^{(1)}
=
0
\ \ \ \ \ \ \ \ \ \ \ \ \ {\scriptstyle{(k\,=\,1\,\cdots\,m)}}
\end{equation}
will admit all transformations of the prolonged group $X^{ (1)}f$
if and only if every expression: $X^{ (1)}\, 
\big( \sum\, U_{ ki}\, x_i^{ (1)} \big)$
vanishes by virtue of the system, hence when relations of the form:
\[
X^{(1)}
\bigg(
\sum_{i=1}^n\,U_{ki}\,x_i^{(1)}
\bigg)
=
\sum_{j=1}^m\,\rho_{kj}\,
\sum_{i=1}^n\,U_{ji}\,x_i^{(1)}
\ \ \ \ \ \ \ \ \ \ \ \ \ {\scriptstyle{(k\,=\,1\,\cdots\,m)}}
\]
hold, where the $\rho_{ kj}$ denote functions of the $x$ alone.

If we combine the Proposition~5 of Chap.~\ref{kapitel-7}, 
p.~\pageref{Satz-5-S-118} with the proposition proved
above (p.~\pageref{Satz-1-S-527}), it follows immediately the:

\def\theproposition{2}\begin{proposition}
If a system of equations of the form:
\def\theequation{10}\begin{equation}
\sum_{i=1}^n\,U_{ki}(x_1,\dots,x_n)\,x_i^{(1)}
=
0
\ \ \ \ \ \ \ \ \ \ \ \ \ {\scriptstyle{(k\,=\,1\,\cdots\,m)}}
\end{equation}
admits the two one-term groups: $X_1^{ (1)} f$ and $X_2^{ (1)} f$
which are derived from $X_1f$ and $X_2f$, respectively, by means of
the special prolongation defined above, then it admits at the same
time the one-term group: $X_1^{ (1)} X_2^{ (1)} f - X_2^{ (1)} X_1^{
(1)} f$, which is derived by means of prolongation from the group:
$X_1 X_2 f - X_2 X_1 f$.
\end{proposition}

\smallercharacters{

We have obtained the preceding proposition by applying the
Proposition~5 of Chap.~\ref{kapitel-7} to the special case of a system
of equations of the form~\thetag{ 10}. But it must be mentioned that
our present proposition can be proved in a way substantially easier
than the general Proposition~5 of the Chap.~\ref{kapitel-7}. Indeed,
one can convince oneself without difficulty by a calculation that a
system of equations~\thetag{ 10} which admits $X_1^{ (1)} f$ and
$X_2^{ (1)}f$ also allows $X_1^{ (1)} X_2^{ (1)} f - X_2^{ (1)} X_1^{
(1)} f$.

}

We call an equation of the form:
\[
\sum_{i=1}^n\,U_i(x_1,\dots,x_n)\,\D\,x_i
=
0
\]
a \terminology{Pfaffian equation} \deutsch{Pfaffsche Gleichung},
and we denote its left-hand side a
\terminology{Pfaffian expression} \deutsch{Pfaffsche Ausdruck}.

From the developments of the present chapter done up to now,
it results important propositions about systems of Pfaffian 
equations:
\def\theequation{10'}\begin{equation}
\sum_{i=1}^n\,U_{ki}(x_1,\dots,x_n)\,\D\,x_i
=
0
\ \ \ \ \ \ \ \ \ \ \ \ \ {\scriptstyle{(k\,=\,1\,\cdots\,m)}}.
\end{equation}

The execution of the transformation: $y_i = f_i ( x_1, \dots, x_n)$ on
the differential expression: $\sum_i\, U_{ ki} (x_1, \dots, x_n)\,
\D\, x_i$ obviously happens in the way that one transforms
the $2\, n$ quantities $x_1, \dots, x_n$, $\D\, x_1, \dots, \D\, x_n$
by means of the transformation:
\[
y_i
=
f_i(x_1,\dots,x_n),
\ \ \ \ \ \ \ \
\d\,y_i
=
\sum_{\nu=1}^n\,
\frac{\partial f_i}{\partial x_\nu}\,
\D\,x_\nu
\ \ \ \ \ \ \ \ \ \ \ \ \ {\scriptstyle{(i\,=\,1\,\cdots\,n)}}.
\]
From this, it follows that the system of the Pfaffian equations:
\def\theequation{10'}\begin{equation}
\sum_{i=1}^n\,U_{ki}(x_1,\dots,x_n)\,\D\,x_i
=
0
\ \ \ \ \ \ \ \ \ \ \ \ \ {\scriptstyle{(k\,=\,1\,\cdots\,m)}}
\end{equation}
remains invariant by the one-term group $Xf$ if and only if 
the system of equations:
\[
\sum_{i=1}^n\,U_{ki}(x_1,\dots,x_n)\,x_i^{(1)}
=
0
\ \ \ \ \ \ \ \ \ \ \ \ \ {\scriptstyle{(k\,=\,1\,\cdots\,m)}}
\]
admits the one-term prolonged group:
\[
X^{(1)}f
=
\sum_{i=1}^n\,\xi_i\,
\frac{\partial f}{\partial x_i}
+
\sum_{i=1}^n\,\xi_i^{(1)}\,
\frac{\partial f}{\partial x_i^{(1)}},
\]
hence when $m$ relations of the form:
\[
\aligned
\sum_{i=1}^n\,X\,U_{ki}\,x_i^{(1)}
&
+
\sum_{i=1}^n\,U_{ki}\,\sum_{\nu=1}^n\,
\frac{\partial\xi_i}{\partial x_\nu}\,
x_\nu^{(1)}
=
\sum_{j=1}^m\,\rho_{kj}(x)\,
\sum_{i=1}^n\,U_{ji}\,x_i^{(1)}
\\
&
\ \ \ \ \ \ \ \ \ \ \ \ \ {\scriptstyle{(k\,=\,1\,\cdots\,m)}},
\endaligned
\]
or, what amounts to the same, $m$ relations of the form:
\[
\aligned
\sum_{i=1}^n\,X\,U_{ki}\,\D\,x_i
&
+
\sum_{i=1}^n\,U_{ki}\,\D\,(X\,x_i)
=
\sum_{j=1}^m\,\rho_{kj}(x)\,
\sum_{i=1}^n\,U_{ji}\,\D\,x_i
\\
&
\ \ \ \ \ \ \ \ \ \ \ \ \ {\scriptstyle{(k\,=\,1\,\cdots\,m)}}.
\endaligned
\]
Hence if we introduce the abbreviated way of writing:
\[
\sum_{i=1}^n\,X\,U_{ki}\,\D\,x_i
+
\sum_{i=1}^n\,U_{ki}\,\D\,(X\,x_i)
=
X\,
\bigg(
\sum_{i=1}^n\,U_{ki}\,\D\,x_i
\bigg),
\]
we obtain the

\def\theproposition{3}\begin{proposition}
If the system of the $m$ Pfaffian equations:
\def\theequation{10'}\begin{equation}
\sum_{i=1}^n\,U_{ki}(x_1,\dots,x_n)\,\D\,x_i
=
0
\ \ \ \ \ \ \ \ \ \ \ \ \ {\scriptstyle{(k\,=\,1\,\cdots\,m)}}.
\end{equation}
is supposed to admit all transformations: $y_i = x_i + \varepsilon\,
X\, x_i + \cdots$ of the one-term group $Xf$, hence if for every value
of $\varepsilon$, relations of the form:
\[
\sum_{i=1}^n\,U_{ki}(y_1,\dots,y_n)\,\D\,y_i
=
\sum_{j=1}^m\,
\omega_{kj}(x_1,\dots,x_n,\,\varepsilon)\,
\sum_{i=1}^n\,U_{ji}(x)\,\D\,x_i
\ \ \ \ \ \ \ \ \ \ \ \ \ {\scriptstyle{(k\,=\,1\,\cdots\,m)}}
\]
are supposed to hold, then for this, it is necessary and sufficient
that the $m$ expressions $X ( \sum\, U_{ ki}\, \D\, x_i)$ can be
represented under the form:
\[
X\bigg(
\sum_{i=1}^n\,U_{ki}\,\D\,x_i
\bigg)
=
\sum_{j=1}^m\,\rho_{kj}(x)\,
\sum_{i=1}^n\,U_{ji}\,\D\,x_i
\ \ \ \ \ \ \ \ \ \ \ \ \ {\scriptstyle{(k\,=\,1\,\cdots\,m)}}.
\]
\end{proposition}

Lastly, if we introduce the language: \emphasis{the system of the
Pfaffian equations~\thetag{ 10'} admits the infinitesimal transformation
$Xf$ when $m$ relations of the form:
\[
X\bigg(
\sum_{i=1}^n\,U_{ki}\,\D\,x_i
\bigg)
=
\sum_{j=1}^m\,\rho_{kj}(x)\,
\sum_{i=1}^n\,U_{ji}\,\D\,x_i
\]
hold}, then we can also state the latter proposition in the 
following way:

\plainstatement{The system of the $m$ Pfaffian
\label{S-531} equations~\thetag{
10'} admits all transformations of the one-term group $Xf$ if and 
only if it admits the infinitesimal transformation $Xf$.}

In addition, from Proposition~2, it follows easily:

\def\thetheorem{93}\begin{theorem}
\label{Theorem-93-S-531}
If a system of $m$ Pfaffian equations:
\[
\sum_{i=1}^n\,U_{ki}(x_1,\dots,x_n)\,
\D\,x_i
=
0
\ \ \ \ \ \ \ \ \ \ \ \ \ {\scriptstyle{(k\,=\,1\,\cdots\,m)}}
\]
admits the two one-term groups: $X_1f$ and $X_2f$, then at the
same time, it also admits the one-term group: $X_1 X_2 f - X_2 X_2 f$.
\end{theorem}

It is yet to be remarked that the Pfaffian \emphasis{expression}:
$\sum_i\, U_i ( x_1, \dots, x_n)\, \D\, x_i$ remains invariant
by every transformation of the one-term group $Xf$
when the expression: $X ( \sum\, U_i \, \D\, x_i)$ vanishes
identically.

\smallercharacters{

Earlier on (p.~\pageref{S-53} sq.), we have agreed that in place
of $\xi_i$, we also want to write occasionally $\delta x_i / \delta t$
and in place of $Xf$, also $\delta f / \delta t$. In a similar
way, in place of $X ( \sum\, U_i \, \D\, x_i)$, we also want to 
write: $\frac{ \delta}{ \delta t}\, 
\sum\, U_i\, \D\, x_i$ as well; then we have the equation:
\[
\delta
\bigg(
\sum_{i=1}^n\,U_i\,\D\,x_i
\bigg)
=
\sum_{i=1}^n\,\delta\,U_i\,\D\,x_i
+
\sum_{i=1}^n\,U_i\,\D\,\delta x_i,
\]
an expression which occurs with the same meaning as in the 
Calculus of Variations.

}

\sectionengellie{\S\,\,\,129.}

The propositions of the preceding paragraph perfectly suffice in order
derive in complete generality the theory of the prolongation of a
group by means of the addition of differential quotients. However, we
want first to consider a special case, before we tackle the treatment
of the general case. At the same time, we go back to known geometrical
concepts of the thrice-extended space.

By $x, y, z$, we may understand ordinary point coordinates of the
thrice-extended space and moreover, let:
\def\theequation{11}\begin{equation}
x'
=
\Xi(x,y,z),
\ \ \ \ \ \ \
y'
=
{\sf H},
\ \ \ \ \ \ \
z'
=
{\sf Z}
\end{equation}
be a transformation of this space.

By the transformation~\thetag{ 11}, all points $x, y, z$ are
transformed, that is to say, they are transferred to the new
positions: $x', y', z'$. At the same time, all surfaces take new
positions: every surface: $\chi ( x, y, z) = 0$ is transferred to a
new surface: $\psi ( x', y', z') = 0$, the equation of which is
obtained by means of elimination of $x, y, z$ from the
equations~\thetag{ 11} of the transformation in combination with $\chi
= 0$.

It is in the nature of things that the transformation~\thetag{ 11}
transfer surfaces which enter in contact to surfaces which stand in
the same relationship, at least in general.
Hence, if by the expression \terminology{surface element}
\deutsch{Flächenelement} we call \emphasis{the totality 
of a point $x, y, z$ located on the surface and of the tangential
plane:
\[
z_1-z
=
p\,(x_1-x)+q\,(y_1-y)
\ \ \ \ \ \ \ \ \ \ \
\bigg(
p
=
\frac{\partial z}{\partial x},
\ \ \ \ \
q
=
\frac{\partial z}{\partial y}
\bigg)
\]
passing through it}, then we can say that our transformation~\thetag{
11} converts every surface element $x, y, z, p, q$ to a new surface
element $x', y', z', p', q'$.
In other words, there must exist certain equations:
\def\theequation{12}\begin{equation}
p'
=
\Pi(x,y,z,p,q),
\ \ \ \ \ \ \ \ \ 
q'
=
{\sf K}(x,y,z,p,q)
\end{equation}
which we will moreover really set up below. The equations~\thetag{ 11}
and~\thetag{ 12} taken together represent a transformation which comes
into existence by prolongation of the transformation~\thetag{ 11}.

At present, we assume that the $\infty^r$ transformations:
\def\theequation{13}\begin{equation}
x'
=
\Xi(x,y,z,\,a_1,\dots,a_r),
\ \ \ \ \ \ \
y'
=
{\sf H},
\ \ \ \ \ \ \
z'
=
{\sf Z}
\end{equation}
are given which represent an $r$-term group of point transformations
of the space, and we imagine that for each of these $\infty^r$
transformations, the prolonged transformation just defined is set up.
We claim that the $\infty^r$ prolonged transformations:
\def\theequation{14}\begin{equation}
\left\{
\aligned
x'
&
=
\Xi(x,y,z,\,a_1,\dots,a_r),
\ \ \ \ \ \ \
y'
=
{\sf H},
\ \ \ \ \ \ \
z'
=
{\sf Z}
\\
p'
&
=
\Pi(x,y,z,p,q,\,a_1,\dots,a_r),
\ \ \ \ \ \ \
q'
=
{\sf H}
\endaligned\right.
\end{equation}
that we obtain in this way do form an $r$-term \label{S-532}
group in the variables
$x, y, z, p, q$.

For the proof, we interpret the $\infty^r$ transformations~\thetag{
13} as operations and we observe that these operations permute with
each other both the points $x, y, z$ and the surface elements $x, y,
z, p, q$. When regarded as permutations of the points, these
operations form a group, and consequently, when interpreted as
permutations of the surface elements, they must also form a group,
whence the fact that the equations~\thetag{ 14} represent a group
finds its analytic expression.

We will make more precise the considerations made up to now by
developing them in a purely analytic way.

Let:
\def\theequation{11}\begin{equation}
x'
=
\Xi(x,y,z),
\ \ \ \ \ \ \ \ \
y'
=
{\sf H}(x,y,z),
\ \ \ \ \ \ \ \ \
z'
=
{\sf Z}(x,y,z)
\end{equation}
be a transformation between the variables $x, y, z$ and
$x', y', z'$. If we imagine $z$ as an arbitrarily chosen function:
$z = \varphi ( x, y)$ of $x$ and $y$, there exist partial differential
quotients of first order of $z$ with respect to $x$ and $y$; these
are defined by means of the equation:
\[
\D\,z-p\,\D\,x-q\,\D\,y
=
0.
\]
But on the other hand, our transformation converts every relationship
of dependence \deutsch{Abhängigskeitsverhältniss}: $z = \varphi ( x, y)$
between $x, y, z$ in just such a relationship between $z', x', y'$
which can in general be given the form: $z' = \overline{ \varphi}
(x', y')$; hence $z'$ also has two partial differential quotients
$p'$ and $q'$ which, in their turn, satisfy the condition:
\[
\D\,z'
-
p'\,\D\,x'
-
q'\,\D\,y'
=
0.
\]

If, in the equation just written, we substitute $z', x', y'$ with 
their values ${\sf Z}$, $\Xi$, ${\sf H}$ and if we organize the
result with respect to $\D\, z$, $\D\, x$, $\D\, y$, it comes:
\[
\aligned
\bigg(
\frac{\partial{\sf Z}}{\partial z}
-
p'\,\frac{\partial\Xi}{\partial z}
-
q'\,\frac{\partial{\sf H}}{\partial z}
\bigg)\,\D\,z
&
+
\bigg(
\frac{\partial{\sf Z}}{\partial x}
-
p'\,\frac{\partial\Xi}{\partial x}
-
q'\,\frac{\partial{\sf H}}{\partial x}
\bigg)\,\D\,x
\\
&
+
\bigg(
\frac{\partial{\sf Z}}{\partial y}
-
p'\,\frac{\partial\Xi}{\partial y}
-
q'\,\frac{\partial{\sf H}}{\partial y}
\bigg)\,\D\,y
=
0,
\endaligned
\]
or, because of $\D\, z = p\, \D\, x + q\, \D\, y$:
\def\theequation{15}\begin{equation}
\left\{
\aligned
&\ \ \ \ \,
\bigg\{
\frac{\partial{\sf Z}}{\partial x}
-
p'\,\frac{\partial\Xi}{\partial x}
-
q'\,\frac{\partial{\sf H}}{\partial x}
+
p\,\bigg(
\frac{\partial{\sf Z}}{\partial z}
-
p'\,\frac{\partial\Xi}{\partial z}
-
q'\,\frac{\partial{\sf H}}{\partial z}
\bigg)
\bigg\}\,\D\,x
\\
&
+
\bigg\{
\frac{\partial{\sf Z}}{\partial y}
-
p'\,\frac{\partial\Xi}{\partial y}
-
q'\,\frac{\partial{\sf H}}{\partial y}
+
q\,\bigg(
\frac{\partial{\sf Z}}{\partial z}
-
p'\,\frac{\partial\Xi}{\partial z}
-
q'\,\frac{\partial{\sf H}}{\partial z}
\bigg)
\bigg\}\,\D\,y
=
0.
\endaligned\right.
\end{equation}

In the equation~\thetag{ 15}, the two factors of $\D\,x$ and of $\D\,
y$ must vanish, because $\D\, x$ and $\D\, y$ are not linked together
by a relation. We therefore obtain the two equations:
\def\theequation{16}\begin{equation}
\left\{
\aligned
p'\,\bigg(
\frac{\partial\Xi}{\partial x}
+
p\,\frac{\partial\Xi}{\partial z}
\bigg)
+
q'\,\bigg(
\frac{\partial{\sf H}}{\partial x}
+
p\,\frac{\partial{\sf H}}{\partial z}
\bigg)
&
=
\frac{\partial{\sf Z}}{\partial x}
+
p\,\frac{\partial{\sf Z}}{\partial z}
\\
p'\,\bigg(
\frac{\partial\Xi}{\partial y}
+
q\,\frac{\partial\Xi}{\partial z}
\bigg)
+
q'\,\bigg(
\frac{\partial{\sf H}}{\partial y}
+
q\,\frac{\partial{\sf H}}{\partial z}
\bigg)
&
=
\frac{\partial{\sf Z}}{\partial y}
+
q\,\frac{\partial{\sf Z}}{\partial z}.
\endaligned\right.
\end{equation}

These equations can be solved with respect to $p'$ and $q'$. Indeed, 
if the determinant:
\def\theequation{17}\begin{equation}
\left\vert
\begin{array}{cc}
\frac{\partial\Xi}{\partial x}+p\,\frac{\partial\Xi}{\partial z}
\,\,\,&\,\,\,
\frac{\partial\Xi}{\partial y}+q\,\frac{\partial\Xi}{\partial z}
\\
\frac{\partial{\sf H}}{\partial x}+p\,\frac{\partial{\sf H}}{\partial z}
\,\,\,&\,\,\,
\frac{\partial{\sf H}}{\partial y}+q\,\frac{\partial{\sf H}}{\partial z}
\end{array}
\right\vert
\end{equation}
would vanish for arbitrary functions $z$ of $x$ and $y$, then
it would actually vanish identically, that is to say, for every 
value of the variables $x, y, z, p, q$. Evidently, this case
can occur only when all $2 \times 2$ determinants of the matrix:
\[
\left\vert
\begin{array}{ccc}
\frac{\partial\Xi}{\partial x} 
\,\,&\,\,
\frac{\partial\Xi}{\partial y}
\,\,&\,\,
\frac{\partial\Xi}{\partial z}
\\
\frac{\partial{\sf H}}{\partial x} 
\,\,&\,\,
\frac{\partial{\sf H}}{\partial y}
\,\,&\,\,
\frac{\partial{\sf H}}{\partial z}
\end{array}
\right\vert
\]
vanish identically. But this is excluded from the beginning.

We therefore see that the equations~\thetag{ 16} are in general 
solvable with respect to $p'$ and $q'$, whichever function of
$x$ and $y$ can $z$ be, and that the resolution is not possible
only when the function $z = \varphi ( x, y)$ satisfies the partial 
differential equation which arises by setting to zero the 
determinant~\thetag{ 17}. By really executing the resolution of
the equations~\thetag{ 16}, we obtain for $p'$ and $q'$ completely
determined functions of $x, y, z, p, q$:
\[
p'
=
\Pi(x,y,z,p,q),
\ \ \ \ \ \ \ \ \
q'
=
{\sf K}(x,y,z,p,q).
\]
\emphasis{This determination is generally valid, because we have made
no special assumption on the function $z = \varphi ( x, y)$.}

Besides, all of that is known long since.

The equation~\thetag{ 15} shows that it is possible to determine a
quantity $\mathfrak{ a}$ in such a way that the relation:
\def\theequation{18}\begin{equation}
\D\,{\sf Z}
-
p'\,\D\,\Xi
-
q'\,\D\,{\sf H}
=
\mathfrak{a}\,
(\D\,z-p\,\D\,x-q\,\D\,y)
\end{equation}
holds identically in $\D\,x$, $\D\, y$, $\D\, z$. Indeed, if one
expands~\thetag{ 18} with respect to $\D\, x$, $\D\, y$, $\D\, z$, one
obtains at first by comparing the factors of $\D\, z$ in both
sides:
\[
\mathfrak{a}
=
\frac{\partial{\sf Z}}{\partial z}
-
p'\,\frac{\partial\Xi}{\partial z}
-
q'\,\frac{\partial{\sf H}}{\partial z}\,;
\]
but if one sets this value in the equation~\thetag{ 18}, then this
equation converts into~\thetag{ 15}.

\renewcommand{\thefootnote}{\fnsymbol{footnote}}
In its turn, the equation~\thetag{ 18} has a very simple meaning:
\emphasis{it expresses that the prolonged transformation:
\[
x'=\Xi,
\ \ \ \ \ \ \
y'={\sf H},
\ \ \ \ \ \ \
z'={\sf Z},
\ \ \ \ \ \ \
p'=\Pi,
\ \ \ \ \ \ \
q'={\sf K}
\]
leaves invariant the Pfaffian equation: $\D\, z - p\, \D\, x - q\, 
\D\, y = 0$}.\footnote[1]{\,
\name{Lie}, Göttinger Nachr. 1872, p.~480, Verhandl. d. G. d. W. zu
Christiania 1873; Math. Ann. Vol. VIII.
}
\renewcommand{\thefootnote}{\arabic{footnote}}

For what follows, it is very important to observe that the
prolonged transformation: $x' = \Xi$, \dots, $q' = {\sf K}$
is perfectly determined by this property. In other words: if one
knows $\Xi$, ${\sf H}$, ${\sf Z}$, then $\Pi$ and ${\sf K}$ are
uniquely determined by the requirement that the transformation:
$x' = \Xi$, \dots, $q' = {\sf K}$ should leave invariant the
Pfaffian equation: $\D\, z - p\, \D\, x - q\, \D\, y = 0$.

\medskip

At present, we again pass to the consideration of an $r$-term group:
\def\theequation{19}\begin{equation}
x'
=
\Xi(x,y,z,\,a_1,\dots,a_r),
\ \ \ \ \ \ \ \ \
y'
=
{\sf H},
\ \ \ \ \ \ \ \ \
z'
=
{\sf Z},
\end{equation}
and we imagine that each one of the $\infty^r$ transformations of
this group is prolonged in the way indicated formerly by adding
the equations:
\def\theequation{20}\begin{equation}
p'
=
\Pi(x,y,z,p,q,\,a_1,\dots,a_r),
\ \ \ \ \ \ \ \ \
q'
=
{\sf K}.
\end{equation}
Then we already know (cf. p.~\pageref{S-532}) that the 
equations~\thetag{ 19} and~\thetag{ 20} taken together represent
an $r$-term group, but at present, we want to prove this also 
in an analytic way.

If the equations:
\[
x''
=
\Xi(x',y',z',\,b_1,\dots,b_r),
\ \ \ \ \ \ \ \ \
y''
=
{\sf H},
\ \ \ \ \ \ \ \ \
z''
=
{\sf Z}
\]
are combined with~\thetag{ 19}, then one obtains in the known way:
\[
x''
=
\Xi(x,y,z,\,c_1,\dots,c_r),
\ \ \ \ \ \ \ \ \
y''
=
{\sf H},
\ \ \ \ \ \ \ \ \
z''
=
{\sf Z},
\]
where $c$ depends only on the $a$ and on the $b$. But on the other hand,
the two equations:
\[
\aligned
\D\,z'-p'\,\D\,x'-q'\,\D\,y'
&
=
\mathfrak{a}\,
(d\,z-p\,\D\,x-q\,\D\,y)
\\
\D\,z''-p''\,\D\,x''-q''\,\D\,y''
&
=
\mathfrak{a}'\,
(\D\,z'-p'\,\D\,x'-q'\,\D\,y')
\endaligned
\]
give the analogous equation:
\[
\D\,z''-p''\,\D\,x''-q''\,\D\,y''
=
\mathfrak{a}\,\mathfrak{a'}\,
(\D\,z-p\,\D\,x-q\,\D\,y),
\]
which, according to the observation made above, shows that $p''$
and $q''$ depend on $x, y, z, p, q$ and on the $c$ in exactly the
same way as $p'$ and $q'$ depend on $x, y, z, p, q$ and on the $a$. 
We therefore have the

\def\theproposition{4}\begin{proposition}
\label{Satz-4-S-535}
If the equations:
\def\theequation{19}\begin{equation}
x'
=
\Xi(x,y,z,\,a_1,\dots,a_r),
\ \ \ \ \ \ \ \ \
y'
=
{\sf H},
\ \ \ \ \ \ \ \ \
z'
=
{\sf Z}
\end{equation}
represent an $r$-term group and if one determines the functions
$\Pi ( x, y, z, \, a_1, \dots, a_r)$, ${\sf K} ( x, y, z, \, 
a_1, \dots, a_r)$ in such a way that the prolonged transformation
equations:
\def\theequation{21}\begin{equation}
x'
=
\Xi,
\ \ \ \ \ \ \ \ \
y'
=
{\sf H},
\ \ \ \ \ \ \ \ \
z'
=
{\sf Z},
\ \ \ \ \ \ \ \ \
p'
=
\Pi,
\ \ \ \ \ \ \ \ \
q'
=
{\sf K}
\end{equation}
leave invariant the Pfaffian equation: $\D\, z - p\, \D\, x - q\, \D\,
y = 0$, so that a relation of the form:
\[
\D\,z'-p'\,\D\,x'-q'\,\D\,y'
=
\mathfrak{a}
(x,y,z,p,q,\,a_1,\dots,a_r)\,
(\D\,z-p\,\D\,x-q\,\D\,y)
\]
holds, then these transformation equations represent in the same
way an $r$-term group and to be precise, a group which has the
same parameter group as the original group.
\end{proposition}

We now claim: If the $r$-term group~\thetag{ 19} is generated
by $r$ independent infinitesimal transformations\,---\,and as
always, we naturally assume this also here\,---, then at the same
time, the same holds true of the prolonged group~\thetag{ 21}.

We will prove this claim at first for the simple case: $r = 1$. 
Let therefore $r = 1$ and let the group~\thetag{ 19} be generated
by the infinitesimal transformation:
\[
Xf
=
\xi\,\frac{\partial f}{\partial x}
+
\eta\,\frac{\partial f}{\partial y}
+
\zeta\,\frac{\partial f}{\partial z}.
\]
In order to prove that the prolonged group is generated
by an infinitesimal transformation, we only need to show
that from $Xf$, a prolonged infinitesimal transformation:
\[
X^{(1)}f
=
\xi\,\frac{\partial f}{\partial x}
+
\eta\,\frac{\partial f}{\partial y}
+
\zeta\,\frac{\partial f}{\partial z}
+
\pi\,\frac{\partial f}{\partial p}
+
\chi\,\frac{\partial f}{\partial q}
\]
can be derived which leaves invariant the Pfaffian equation:
\def\theequation{22}\begin{equation}
\D\,z-p\,\D\,x-q\,\D\,y
=
0.
\end{equation}
If this is proved, then it is clear that the one-term group in the
variables $x, y, z, p, q$ generated by $X^{ (1)}f$ is identical to the
one-term group~\thetag{ 21}, for this last group then leaves invariant
the Pfaffian equation~\thetag{ 22} and comes into existence by
prolongation of the group $Xf$ which, according to the assumption, is
just the group~\thetag{ 19}.

According to p.~\pageref{S-531}, the infinitesimal transformation $X^{
(1)} f$ leaves invariant the Pfaffian equation~\thetag{ 22} if and
only if $\pi$ and $\chi$ satisfy an equation of the form:
\[
\D\,\zeta-p\,\D\,\xi-q\,\D\,\eta
-\pi\,\D\,x-\chi\,\D\,y
=
\mathfrak{b}\,
(\D\,z-p\,\D\,x-q\,\D\,y),
\]
where it is understood that $\mathfrak{ b}$ is a function of $x, y, z,
p, q$. This equation decomposes in three equations:
\[
\aligned
\mathfrak{b}
&
=
\frac{\partial\zeta}{\partial z}
-
p\,\frac{\partial\xi}{\partial z}
-
q\,\frac{\partial\eta}{\partial z}
\\
\pi
&
=
\frac{\partial\zeta}{\partial x}
-
p\,\frac{\partial\xi}{\partial x}
-
q\,\frac{\partial\eta}{\partial x}
+
\mathfrak{b}\,p
\\
\chi
&
=
\frac{\partial\zeta}{\partial y}
-
p\,\frac{\partial\xi}{\partial y}
-
q\,\frac{\partial\eta}{\partial y}
+
\mathfrak{b}\,q,
\endaligned
\]
whence it comes for $\pi$ and $\chi$ the completely determined values:
\[
\pi
=
\frac{\D\,\zeta}{\D\,x}
-
p\,\frac{\D\,\xi}{\D\,x}
-
q\,\frac{\D\,\eta}{\D\,x},
\ \ \ \ \ \ \ \ \ \ \ \ \ \
\chi
=
\frac{\D\,\zeta}{\D\,y}
-
p\,\frac{\D\,\xi}{\D\,y}
-
q\,\frac{\D\,\eta}{\D\,y},
\]
where for brevity, we have set:
\[
\frac{\partial\Phi(x,y,z)}{\partial x}
+
p\,\frac{\partial\Phi(x,y,z)}{\partial z}
=
\frac{\D\,\Phi}{\D\,x},
\ \ \ \ \ \ \ \ \ \ \ \ \ \
\frac{\partial\Phi(x,y,z)}{\partial y}
+
q\,\frac{\partial\Phi(x,y,z)}{\partial z}
=
\frac{\D\,\Phi}{\D\,y}.
\]

With these words, the existence of the infinitesimal transformation
$X^{ (1)}f$ is proved, and also at the the same time, the correctness
for $r = 1$ of the claim stated above.

Now, let $r$ be arbitrary. Next, if $e_1\, X_1f + \cdots + e_r\, X_rf$
is an arbitrary infinitesimal transformation of the group~\thetag{
19}, then evidently, the one-term group which is generated by the
prolonged infinitesimal transformation: $e_1\, X_1^{ (1)} f + \cdots +
e_r\, X_r^{ (1)} f$ belongs to the group~\thetag{ 21}, because the
concerned one-term group comes indeed into existence by prolongation
of the one-term group: $e_1\, X_1f + \cdots + e_r\, X_rf$. From this,
it results that the group~\thetag{ 21} contains the $\infty^{ r-1}$
one-term groups: $e_1\, X_1^{ (1)}f + \cdots + e_r\, X_r^{ (1)} f$, so
that it is generated by the $r$ independent infinitesimal
transformations: $X_1^{(1)} f, \dots, X_r^{ (1)}f$.

Of course, the $X_k^{ (1)}f$ satisfy pairwise relations of the form:
\[
X_i^{(1)}X_k^{(1)}f-X_k^{(1)}X_i^{(1)}f
=
\sum_{s=1}^r\,c_{iks}'\,X_s^{(1)}f\,;
\]
we will verify this by means of a computation, and hence in addition,
we will recognize that the $c_{ iks}'$ here have the same values as
the $c_{ iks}$ in the relations:
\[
X_iX_kf-X_kX_if
=
\sum_{s=1}^r\,c_{iks}\,X_sf.
\]

According to Theorem~93, p.~\pageref{Theorem-93-S-531}, the Pfaffian
equation: $\D\, z - p\, \D\, x - q\, \D\, y = 0$ admits,
simultaneously with the two infinitesimal transformations: $X_i^{
(1)}f$ and $X_k^{ (1)}f$ also the following: $X_i^{ (1)} X_k^{ (1)} f
- X_k^{ (1)} X_i^{ (1)} f$, which has visibly the form:
\[
X_i^{(1)}X_k^{(1)}f
-
X_k^{(1)}X_i^{(1)}f
=
X_iX_kf-X_kX_if
+
\alpha\,\frac{\partial f}{\partial p}
+
\beta\,\frac{\partial f}{\partial q}.
\]
From this, it results that the infinitesimal transformation:
$X_i^{ (1)} X_k^{ (1)} f - X_k^{ (1)} X_i^{ (1)} f$ is obtained
by prolongation of:
\[
X_iX_kf-X_kX_if
=
\sum_{s=1}^r\,c_{iks}\,X_sf,
\]
so that the relations:
\[
X_i^{(1)}X_k^{(1)}f-X_k^{(1)}X_i^{(1)}f
=
\sum_{s=1}^r\,c_{iks}\,X_s^{(1)}f
\]
hold.

From this, we yet see that the prolonged group~\thetag{ 21} is
holoedrically isomorphic to the group~\thetag{ 19}, which
is coherent with the fact that the two groups have the same
parameter group (Proposition~4, p.~\pageref{Satz-4-S-535}).

\medskip

The group:
\def\theequation{19}\begin{equation}
x'
=
\Xi(x,y,z,\,a_1,\dots,a_r),
\ \ \ \ \ \ \
y'={\sf H},
\ \ \ \ \ \ \
z'={\sf Z}
\end{equation}
can still be prolonged further, namely by taking also differential
quotients of order higher than the first order. Here, we only want
to consider also the prolongation by means of the addition
of differential quotients of
second order. 

Since $z$ is considered as a function of $x$ and $y$, apart
from $p$ and $q$, we have to take also account of the three
differential quotients:
\[
\frac{\partial^2z}{\partial x^2}
=
r,
\ \ \ \ \ \ \
\frac{\partial^2z}{\partial x\partial y}
=
s,
\ \ \ \ \ \ \
\frac{\partial^2z}{\partial y^2}
=
t.
\]

Through the transformations of our group, $x', y', z'$ depend
in the known way upon $x, y, z$, and according to what precedes,
$p', q'$ depend likewise upon $x, y, z, p, q$. As should at
present be shown, $r', s', t'$ can also be represented as functions
of $x, y, z, p, q, r, s, t$:
\[
r'
=
{\sf P}(x,y,z,p,q,r,s,t;\,a_1,\dots,a_r),
\ \ \ \ \ \ \
s'
=
\Sigma,
\ \ \ \ \ \ \
t'
=
{\sf T}.
\]

The quantities $r', s', t'$ are defined by:
\[
\D\,p'-r'\,\D\,x'-s'\,\D\,y'
=
0,
\ \ \ \ \ \ \ \ \ \
\D\,q'-s'\,\D\,x'-t'\,\D\,y'
=
0.
\]

Here, if the $\D\, x'$, $\D\, y'$, $\D\, p'$, $\D\, q'$ are expanded
with respect to $\D\, x$, $\D\, y$, $\D\, z$, $\D\, p$, $\D\, q$, and
next, if the values for $\D\, z$, $\D\, p$, $\D\, q$ are inserted
from:
\[
\D\,z-p\,\D\,x-q\,\D\,y=0,
\ \ \ \ \
\D\,p-r\,\D\,x-s\,\D\,y=0,
\ \ \ \ \
\D\,q-s\,\D\,x-t\,\D\,y=0,
\]
then one obtains two equations of the form: $A\, \D\, x + B\, \D\, y =
0$ whose coefficients, apart from $x', y', z'$ and their differential
quotients, yet contain also $p'$, $q'$, $r'$, $s'$, $t'$, $p$, $q$,
$r$, $s$, $t$. Since $\D\, x$ and $\D\, y$ are independent of
each other, $A$ and $B$ must in the two cases vanish individually,
and hence we find the four equations:
\[
\aligned
r'\bigg(
\frac{\partial\Xi}{\partial x}
+
p\,\frac{\partial\Xi}{\partial z}
\bigg)
+
s'\bigg(
\frac{\partial{\sf H}}{\partial x}
+
p\,\frac{\partial{\sf H}}{\partial z}
\bigg)
&
=
\frac{\partial\Pi}{\partial x}
+
p\,\frac{\partial\Pi}{\partial z}
+
r\,\frac{\partial\Pi}{\partial p}
+
s\,\frac{\partial\Pi}{\partial q}
\\
r'\bigg(
\frac{\partial\Xi}{\partial y}
+
q\,\frac{\partial\Xi}{\partial z}
\bigg)
+
s'\bigg(
\frac{\partial{\sf H}}{\partial y}
+
q\,\frac{\partial{\sf H}}{\partial z}
\bigg)
&
=
\frac{\partial\Pi}{\partial y}
+
q\,\frac{\partial\Pi}{\partial z}
+
s\,\frac{\partial\Pi}{\partial p}
+
t\,\frac{\partial\Pi}{\partial q}
\\
& 
\text{and so on}.
\endaligned
\]

Of these equations, the first two are solvable with respect to $r'$
and $s'$, the last two with respect to $s'$ and $t'$, exactly as the
equations~\thetag{ 16} were solvable with respect to $p'$ and $q'$.
Still, the question is only whether the two values which are found
in this way for $s'$ coincide with each other. But as is known, this
is indeed the case, since otherwise, we would receive a relation between
$x$, $y$, $z$, $p$, $q$, $r$, $s$, $t$ not holding identically
which should be satisfied identically after the substitution: 
$z = \varphi ( x, y)$, and this is impossible, because $\varphi$
is submitted to absolutely no restriction.

Similarly as before, we realize also here that the quantities
$r', s', t'$ are defined uniquely as functions of $x$, $y$, $z$, 
$p$, $q$, $r$, $s$, $t$ by the condition that equations of the form:
\[
\aligned
\D\,p'-r'\,\D\,x'-s'\,\D\,y'
&
=
\alpha_1(\D\,z-p\,\D\,x-q\,\D\,y)
+
\beta_1(\D\,p-r\,\D\,x-s\,\D\,y)
\\
&
\ \ \ \ \ 
+
\gamma_1(\D\,q-s\,\D\,x-t\,\D\,y)
\\
\D\,q'-s'\,\D\,x'-t'\,\D\,y'
&
=
\alpha_2(\D\,z-p\,\D\,x-q\,\D\,y)
+
\beta_2(\D\,p-r\,\D\,x-s\,\D\,y)
\\
&
\ \ \ \ \ 
+
\gamma_2(\D\,q-s\,\D\,x-t\,\D\,y)
\endaligned
\]
hold identically in $\D\, p$, $\D\, q$, $\D\, z$, $\D\, y$, $\D\, x$.
For this, if we remember the former identity:
\[
\D\,z'-p'\,\D\,x'-q'\,\D\,y'
=
\mathfrak{a}(\D\,z-p\,\D\,x-q\,\D\,y)
\]
which determine $p'$ and $q'$, we can say that for given
$\Xi$, ${\sf H}$, ${\sf Z}$, the transformation
equations:
\def\theequation{23}\begin{equation}
x'
=
\Xi(x,y,z,a_1,\dots,a_r),
\,\,\,\dots,\,\,\,
t'
=
{\sf T}(x,y,z,p,q,r,s,t,\,a_1,\dots,a_r)
\end{equation}
are completely and uniquely determined
by the condition that they should leave invariant the system of
Pfaffian equations:
\def\theequation{24}\begin{equation}
\D\,z-p\,\D\,x-q\,\D\,y
=
0,
\ \ \ \ \ \ \
\D\,p-r\,\D\,x-s\,\D\,y
=
0,
\ \ \ \ \ \ \
\D\,q-s\,\D\,x-t\,\D\,y
=
0.
\end{equation}

From this, it becomes clear that after the succession of two 
transformations~\thetag{ 23}, a transformation comes into
existence with again leaves invariant the system~\thetag{ 24}.
But by assumption, the equations:
\[
\aligned
x'
&
=
\Xi(x,y,z,\,a_1,\dots,a_r),
\ \ \ \ \ \ \ \ \ \ \ \ \ \
y'
=
{\sf H},
\ \ \ \ \ \ \
z'
=
{\sf Z}
\\
x''
&
=
\Xi(x',y',z',\,b_1,\dots,b_r),
\ \ \ \ \ \ \ \ \ \
y''
=
{\sf H},
\ \ \ \ \ \ \
z''
=
{\sf Z}
\endaligned
\]
have as a consequence:
\def\theequation{19'}\begin{equation}
x''
=
\Xi(x,y,z,\,c_1,\dots,c_r)
\ \ \ \ \ \ \ \ \ \ \ \ \ \ \
y''
=
{\sf H},
\ \ \ \ \ \ \
z''
=
{\sf Z},
\end{equation}
where the $c$ depend only on the $a$ and on the $b$. Therefore, the two 
transformations:
\[
\aligned
x'
&
=
\Xi(x,y,z,\,a_1,\dots,a_r),
\,\,\,\dots,\,\,\,
\ \ \ \
t'
=
{\sf T}(x,y,z,p,q,r,s,t,\,a_1,\dots,a_r)
\\
x''
&
=
\Xi(x',y',z',\,b_1,\dots,b_r),
\,\,\,\dots,\,\,\,
t''
=
{\sf T}(x',y',z',p',q',r',s',t',\,b_1,\dots,b_r)
\endaligned
\]
when executed one after the other, give a transformation
which results from~\thetag{ 19'} by the same prolongation
as the one by which~\thetag{ 23} comes into existence from:
\def\theequation{19}\begin{equation}
x'
=
\Xi(x,y,z,\,a_1,\dots,a_r),
\ \ \ \ \ \ \ \ \ \
y'
=
{\sf H},
\ \ \ \ \ \ \ \ \ \
z'
=
{\sf Z},
\end{equation}
a transformation therefore which belongs in the same way to the
family~\thetag{ 23}. Thus, the transformations~\thetag{ 23} form
an $r$-term group. 

We will convince ourselves directly that the group~\thetag{ 23}
is generated by $r$ independent infinitesimal transformations.

To begin with, we again assume that $r = 1$ and that the group~\thetag{
19} is generated by the infinitesimal transformation $Xf$.
Then obviously, we need only to prove that there is a prolonged
infinitesimal transformation:
\[
X^{(2)}f
=
Xf
+
\pi\,\frac{\partial f}{\partial p}
+
\chi\,\frac{\partial f}{\partial q}
+
\rho\,\frac{\partial f}{\partial r}
+
\sigma\,\frac{\partial f}{\partial s}
+
\tau\,\frac{\partial f}{\partial t}
\]
which leaves invariant the system of the Pfaffian equations~\thetag{
24}. But there is no difficulty to do that; for $\pi$ and $\chi$, we
find the same values as before, and for $\rho$, $\sigma$, $\tau$, we
obtain in a similar way expressions which depend linearly and
homogeneously in the $\xi$, $\eta$, $\zeta$ and their differential
quotients of first order and of second order.

Also, we find here at first four equations for $\rho$, $\sigma$,
$\tau$; but it can easily be proved that these equations are
compatible with each other. We postpone the realization of this proof
to the consideration of the general case, which appears to be clearer
than the case at hand here.

From the existence of $X^{ (2)} f$, it naturally follows that the
group~\thetag{ 23} is just generated by $X^{ (2)} f$ in the case
$r = 1$. 

One realizes in the same way that for an arbitrary $r$, the
group~\thetag{ 23} is generated by the $r$ infinitesimal 
transformations: $X_1^{ (2)}f, \dots, X_r^{ (2)}f$ which 
are obtained by prolongation of the infinitesimal transformations:
$X_1f, \dots, X_rf$. As a result, it is proved at the same
time that relations of the form:
\[
X_i^{(2)}X_k^{(2)}f
-
X_k^{(2)}X_i^{(2)}f
=
\sum_{s=1}^r\,c_{iks}''\,X_s^{(2)}f
\]
hold.

It can be seen easily that $c_{ iks}'' = c_{ iks}$. In fact, 
together with $X_i^{ (2)} f$ and $X_k^{ (2)}f$, the system of the
Pfaffian equations~\thetag{ 24} also admits the infinitesimal 
transformation: $X_i^{ (2)} X_k^{ (2)} f - X_k^{ (2)} X_i^{ (2)} f
= \big \leftbracket X_i^{ (2)}, \, X_k^{ (2)} 
\big \rightbracket$, but this
transformation obviously has the form:
\[
\big\leftbracket
X_i^{(2)},\,X_k^{(2)}
\big\rightbracket
=
\leftbracket
X_i,\,X_k
\rightbracket
+
\alpha\,
\frac{\partial f}{\partial p}
+
\beta\,
\frac{\partial f}{\partial q}
+
\lambda\,
\frac{\partial f}{\partial r}
+
\mu\,
\frac{\partial f}{\partial s}
+
\nu\,
\frac{\partial f}{\partial t},
\]
hence it comes into existence by prolongation of the transformation:
\[
\leftbracket
X_i,\,X_k
\rightbracket
=
\sum_{s=1}^r\,c_{iks}\,X_sf,
\]
and it can be represented in the following way:
\[
\big\leftbracket
X_i^{(2)},\,X_k^{(2)}
\big\rightbracket
=
\sum_{s=1}^r\,c_{iks}\,
X_s^{(2)}f.
\]

Here lies the reason why the prolonged group~\thetag{ 23} is
equally composed with the original group~\thetag{ 19}.

\sectionengellie{\S\,\,\,130.}

After the realization of the special studies of the preceding 
paragraph, the general theory of the prolongation of a finite
continuous group by addition of differential quotients will cause
us no more difficulty.

To begin with, we consider an individual transformation in the $n + m$
variables: $x_1, \dots, x_n$, $z_1, \dots, z_m$, say the following:
\def\theequation{25}\begin{equation}
\aligned
x_i'
&
=
f_i(x_1,\dots,x_n,\,z_1,\dots,z_m)
\ \ \ \ \ \ \ \ \ \ \ \ \ {\scriptstyle{(i\,=\,1\,\cdots\,n)}}
\\
z_k'
&
=
F_k(x_1,\dots,x_n,\,z_1,\dots,z_m)
\ \ \ \ \ \ \ \ \ \ \ \ \ {\scriptstyle{(k\,=\,1\,\cdots\,m)}}.
\endaligned
\end{equation}
If we want to prolong this transformation by adding differential 
quotients, then we must at first agree on how many and which of
the variables $x_1, \dots, x_n$, $z_1, \dots, z_m$ should
be considered as independent, and which ones should be considered
as dependent. This can occur in very diverse ways and to each 
such possible way there corresponds a completely determined
prolongation of the transformation~\thetag{ 25}.

\plainstatement{In the sequel, we will always consider $x_1, \dots, x_n$
as independent variables and $z_1, \dots, z_m$ as functions of 
$x_1, \dots, x_n$, but which can be chosen arbitrarily. Under this
assumption, $x_1', \dots, x_n'$ are in general mutually independent,
while $z_1', \dots, z_m'$ are functions of $x_1', \dots, x_n'$.}

For the differential quotients of the $z$ with respect to the $x$
and of the $z'$ with respect to the $x'$, we introduce the following
notation:
\[
\frac{\partial z_\nu}{\partial x_k}
=
z_{\nu,\,k},
\ \ \ \ \ \ \
\frac{\partial^{\alpha_1+\cdots+\alpha_n}z_\mu}{
\partial x_1^{\alpha_1}\cdots\,\partial x_n^{\alpha_n}}
=
z_{\mu,\,\alpha_1,\dots,\alpha_n},
\ \ \ \ \ \ \
\frac{\partial^{\alpha_1+\cdots+\alpha_n}z_\mu'}{
\partial{x_1'}^{\alpha_1}\cdots\,\partial{x_n'}^{\alpha_n}}
=
z_{\mu,\,\alpha_1,\dots,\alpha_n}',
\]
and we claim that $z_{ \mu, \alpha_1, \dots, \alpha_n}'$
can be expressed by means of $x_1, \dots, x_n$, $z_1, \dots, z_m$, 
and by means of the differential quotients $z_{ \nu, \, \beta_1, 
\dots, \beta_n}$ of first order up to the 
$(\alpha_1 + \cdots + \alpha_n)$-th order:
\def\theequation{26}\begin{equation}
\aligned
z_{\mu,\,\alpha_1,\dots,\alpha_n}'
&
=
F_{\mu,\,\alpha_1,\dots,\alpha_n}
\big(
x_1,\dots,x_n,\,z_1,\dots,z_m,\,z_{\nu,\,\beta_1,\dots,\beta_n}
\big)
\\
&\
{\scriptstyle{(\nu\,=\,1\,\cdots\,m\,;\,\,\,
\beta_1\,+\,\cdots\,+\,\beta_n\,\leqslant\,
\alpha_1\,+\,\cdots\,+\,\alpha_n)}}.
\endaligned
\end{equation}

If we set $\alpha_1 + \cdots + \alpha_n = N$, then the existence of
equations of the form~\thetag{ 26} is clear for $N = 0$; in order to
establish this existence for an arbitrary $N$, we therefore need only
to show that equations of the form~\thetag{ 26} hold also for
$\alpha_1 + \cdots + \alpha_n = N + 1$ as soon as such equations hold
for $\alpha_1 + \cdots + \alpha_n \leqslant N$.

Thus, let the functions $F_{ \mu, \, \alpha_1, \dots, \alpha_n}$
($\alpha_1 + \cdots + \alpha_n \leqslant N$) be known; then the values
of the $z_{\mu,\, \alpha_1, \dots, \alpha_n}'$ ($\alpha_1 + \cdots +
\alpha_n = N + 1$) are to be determined from the equations:
\def\theequation{27}\begin{equation}
\D\,z_{\mu,\,\alpha_1,\dots,\alpha_n}'
-
\sum_{i=1}^n\,
z_{\mu,\,\alpha_1,\dots,\alpha_i+1,\dots,\alpha_n}'\,\D\,x_i'
=
0
\ \ \ \ \ \ \ \ \ \ \ \ \ 
{\scriptstyle{(\alpha_1\,+\,\cdots\,+\,\alpha_n\,=\,N)}},
\end{equation}
and to be precise, \thetag{ 27} must hold identically by virtue of the
system of equations:
\def\theequation{28}\begin{equation}
\aligned
&
\D\,z_{\nu,\,\beta_1,\dots,\beta_n}
-
\sum_{i=1}^n\,
z_{\nu,\,\beta_1,\dots,\beta_i+1,\dots,\beta_n}\,\D\,x_i
=
0
\\
& \ \ \ \ \ \ \ \ \ \ \ \ 
{\scriptstyle{(\nu\,=\,1\,\cdots\,m\,;\,\,\,
0\,\leqslant\,\beta_1\,+\,\cdots\,+\,\beta_n\,\leqslant\,N)}},
\endaligned
\end{equation}
while $\D\, x_1, \dots, \D\, x_n$ are fully independent of each other.
For $z_{ \mu, \, \alpha_1, \dots, \alpha_i + 1, \dots, \alpha_n}'$,
we therefore obtain the equations:
\def\theequation{29}\begin{equation}
\sum_{i=1}^n\,
z_{\mu,\,\alpha_1,\dots,\alpha_i+1,\dots,\alpha_n}'\,
\bigg\{
\frac{\partial f_i}{\partial x_k}
+
\sum_{\nu=1}^m\,z_{\nu,\,k}\,
\frac{\partial f_i}{\partial z_\nu}
\bigg\}
=
\frac{\D\,F_{\mu,\,\alpha_1,\dots,\alpha_n}}{\D\,x_k}
\ \ \ \ \ \ \ \ \ \ \ \ \ {\scriptstyle{(k\,=\,1\,\cdots\,n)}},
\end{equation}
where the right-hand side means the complete differential quotient:
\[
\aligned
\frac{\partial F_{\mu,\,\alpha_1,\dots,\alpha_n}}{
\partial x_k}
&
+
\sum_{\nu=1}^m\,z_{\nu,\,k}\,
\frac{\partial F_{\mu,\,\alpha_1,\dots,\alpha_n}}{
\partial z_\nu}
+
\sum_{\nu=1}^m\,\sum_\beta\,
z_{\nu,\,\beta_1,\dots,\beta_k+1,\dots,\beta_n}\,
\frac{\partial F_{\mu,\,\alpha_1,\dots,\alpha_n}}{
\partial z_{\nu,\,\beta_1,\dots,\beta_n}}
\\
&
\ \ \ \ \ \ \ \ \ \ \ \ \ \ \ \ \ 
{\scriptstyle{(1\,<\,\beta_1\,+\,\cdots\,+\,\beta_n\,
\leqslant\,\alpha_1\,+\,\cdots\,+\,\alpha_n)}}
\endaligned
\]
of $F_{ \mu, \, \alpha_1, \dots, \alpha_n}$ with respect to $x_k$.

If the equations~\thetag{ 29} would not be solvable with respect
to the $z_{ \mu, \, \alpha_1, \dots, \alpha_i + 1, \dots, \alpha_n }'$
($i = 1, \dots, n$), then the determinant:
\[
\left\vert
\begin{array}{cccc}
\frac{\partial f_1}{\partial x_1}
+
\sum_{\mu=1}^m\,z_{\mu,\,1}\,\frac{\partial f_1}{\partial z_\mu}
\,\,\,&\,\,\,\cdot\,\,\,&\,\,\,\cdot\,\,\,&\,\,\,
\frac{\partial f_1}{\partial x_n}
+
\sum_{\mu=1}^m\,z_{\mu,\,n}\,\frac{\partial f_1}{\partial z_\mu}
\\
\cdot
\,\,\,&\,\,\,\cdot\,\,\,&\,\,\,\cdot\,\,\,&\,\,\,
\cdot
\\
\cdot
\,\,\,&\,\,\,\cdot\,\,\,&\,\,\,\cdot\,\,\,&\,\,\,
\cdot
\\
\frac{\partial f_n}{\partial x_1}
+
\sum_{\mu=1}^m\,z_{\mu,\,1}\,\frac{\partial f_n}{\partial z_\mu}
\,\,\,&\,\,\,\cdot\,\,\,&\,\,\,\cdot\,\,\,&\,\,\,
\frac{\partial f_n}{\partial x_n}
+
\sum_{\mu=1}^m\,z_{\mu,\,n}\,\frac{\partial f_n}{\partial z_\mu}
\end{array}
\right\vert
\]
would necessarily be zero, whichever functions of $x_1, \dots, x_n$
one could substitute for $z_1, \dots, z_m$, that is to say: this
determinant should vanish identically for every value of the 
variables $x_i$, $z_\mu$, $z_{ \mu, i}$. One convinces oneself 
easily that this can happen only when all $m \times m$ determinants
of the matrix:
\[
\left\vert
\begin{array}{cccccccc}
\frac{\partial f_1}{\partial x_1} 
\,\,&\,\,\cdot\,\,&\,\,\cdot\,\,&\,\,
\frac{\partial f_1}{\partial x_n}
\,\,&\,\,
\frac{\partial f_1}{\partial z_1}
\,\,&\,\,\cdot\,\,&\,\,\cdot\,\,&\,\,
\frac{\partial f_1}{\partial z_m}
\\
\cdot
\,\,&\,\,\cdot\,\,&\,\,\cdot\,\,&\,\,
\cdot
\,\,&\,\,
\cdot
\,\,&\,\,\cdot\,\,&\,\,\cdot\,\,&\,\,
\cdot
\\
\frac{\partial f_n}{\partial x_1} 
\,\,&\,\,\cdot\,\,&\,\,\cdot\,\,&\,\,
\frac{\partial f_n}{\partial x_n}
\,\,&\,\,
\frac{\partial f_n}{\partial z_1}
\,\,&\,\,\cdot\,\,&\,\,\cdot\,\,&\,\,
\frac{\partial f_n}{\partial z_m}
\end{array}
\right\vert
\]
vanishes identically, hence when the functions $f_1, \dots, f_n$ are
not independent of each other. But since this is excluded, it follows
that the equations~\thetag{ 29} are solvable with respect to the
$z_{\mu, \, \alpha_1, \dots, \alpha_i + 1, \dots, \alpha_n }'$ ($i =
1, \dots, n$).

It still remains to eliminate an objection. Apparently, the
equations~\thetag{ 29} provide in general different values for the
differential quotients $z_{ \mu, \, \alpha_1, \dots, \alpha_i + 1,
\dots, \alpha_n }'$. But in reality, this is only fictitious, because
otherwise this would give certain not identical relations between the
$x_i$, the $z_\mu$ and their differential quotients which should hold
for completely arbitrary functions $z_1, \dots, z_m$ of $x_1, \dots,
x_n$, which is impossible.

We therefore see: the equations~\thetag{ 29} determine all $z_{ \mu,
\, \alpha_1, \dots, \alpha_i + 1, \dots, \alpha_n}'$ completely and
uniquely, and consequently, there are, under the assumptions made,
equations of the form~\thetag{ 26} also for $\alpha_1 + \cdots +
\alpha_n = N +1$, whence their general existence is established. But
in addition, we yet realize that the transformation:
\def\theequation{30}\begin{equation}
\aligned
&
x_i'
=
f_i(x_1,\dots,x_n,\,z_1,\dots,z_m)
\ \ \ \ \ \ \ \ \
z_\mu'
=
F_\mu(x_1,\dots,x_n,\,z_1,\dots,z_m)
\\
&
z_{\mu,\,\alpha_1,\dots,\alpha_n}'
=
F_{\mu\,\alpha_1,\dots,\alpha_n}
\big(
x_1,\dots,x_n,\,z_1,\dots,z_m,\,z_{\nu,\,\beta_1,\dots,\beta_n}
\big)
\ \ \ \ \ \ \ \ \ \ \ \ \ 
{\scriptstyle{(\sum_\chi\,\beta_\chi\,
\leqslant\,\sum_\chi\,\alpha_\chi)}}
\\
&
\ \ \ \ \ \ \ \ \ \ \ \ \ \ \ \ \ \ \ \ \ \ \ \ \ \ \ \ \ \ \ 
{\scriptstyle{(i\,=\,1\,\cdots\,n\,;\,\,\,
\mu\,=\,1\,\cdots\,m\,;\,\,\,
0\,<\,\alpha_1\,+\,\cdots\,+\,\alpha_n\,\leqslant\,N)}}
\endaligned
\end{equation}
leaves invariant the system of equations:
\def\theequation{31}\begin{equation}
\D\,z_{\mu,\,\beta_1,\dots,\beta_n}
-
\sum_{i=1}^n\,z_{\mu,\,\beta_1,\dots,\beta_i+1,\dots,\beta_n}\,
\D\,x_i
=
0
\ \ \ \ \ \ \ \ \ \ \ \ \ 
{\scriptstyle{(\mu\,=\,1\,\cdots\,m\,;\,\,\,
0\,\leqslant\,\sum_\chi\,\beta_\chi\,<\,N)}}
\end{equation}
and that it is completely defined by this property. Lastly, it is yet
also clear that the succession of two transformations~\thetag{ 30}
which leave invariant the Pfaffian system of equations~\thetag{ 31}
again provides a transformation having this constitution.

\medskip

At present, we assume that the original transformation
equations~\thetag{ 25} contain a certain number, say $r$, of
parameters:
\def\theequation{32}\begin{equation}
\left\{
\aligned
x_i'
&
=
f_i(x_1,\dots,x_n,\,z_1,\dots,z_m,\,a_1,\dots,a_r)
\ \ \ \ \ \ \ \ \ \ \ \ \ {\scriptstyle{(i\,=\,1\,\cdots\,n)}}
\\
z_\mu'
&
=
F_\mu(x_1,\dots,x_n,\,z_1,\dots,z_m,\,a_1,\dots,a_r)
\ \ \ \ \ \ \ \ \ \ \ \ \ {\scriptstyle{(\mu\,=\,1\,\cdots\,m)}}
\endaligned\right.
\end{equation}
and that they represent an $r$-term group generated by $r$ independent
infinitesimal transformations. From~\thetag{ 32}, we imagine that all
equations of the form~\thetag{ 26} are derived, in which $\alpha_1 +
\cdots + \alpha_n \leqslant N$; we claim that these equations taken
together with~\thetag{ 32}:
\def\theequation{33}\begin{equation}
\left\{
\aligned
&
x_i'
=
f_i(x,z,a),
\ \ \ \ \ \ \ \ \ \ \ \ \
z_\mu'
=
F_\mu(x,z,a)
\\
&
z_{\mu,\,\alpha_1,\dots,\alpha_n}'
=
F_{\mu,\,\alpha_1,\dots,\alpha_n}
(x,z,z_{\nu,\,\beta_1,\dots,\beta_n},\,a)
\\
&
\ \ \ \ \ \ 
{\scriptstyle{(i\,=\,1\,\cdots\,n\,;\,\,\,
\mu\,=\,1\,\cdots\,m\,;\,\,\,
0\,<\,\alpha_1\,+\,\cdots\,+\,\alpha_n\,\leqslant\,N)}}
\endaligned\right.
\end{equation}
represent again an $r$-term group.

The proof is very simple.

The transformations:
\[
\aligned
x_i'
&
=
f_i(x,z,a),
\ \ \ \ \ \ \ \ \ \ \
z_k'
=
F_k(x,z,a)
\\
x_i''
&
=
f_i(x',z',b),
\ \ \ \ \ \ \ \ \
z_k''
=
F_k(x',z',b)
\endaligned
\]
executed one after the other give a transformation:
\[
x_i''
=
f_i(x,z,c),
\ \ \ \ \ \ \ \ \
z_k''
=
F_k(x,z,c)
\]
in which the $c$ are certain functions of the $a$ and of the $b$.
According to what has been said above, the $z_{ \mu, \, \alpha_1, 
\dots, \alpha_n}''$ express in terms of the $x_i$, $z_k$, 
$z_{ \nu, \, \beta_1, \dots, \beta_n}$ and of the $c_j$
in exactly the same way as the $z_{ \mu, \, \alpha_1, \dots, \alpha_n}'$
express in terms of the $x_i$, $z_k$, $z_{ \nu,\, \beta_1, \dots, 
\beta_n}$ and of the $a$. As a result, our claim is proved.

\medskip

It still remains to prove that the prolonged group~\thetag{ 33} is
generated by $r$ infinitesimal transformations, just as the original
group~\thetag{ 32}. In order to be able to perform this proof, 
we make at first the following considerations:

We start from an arbitrary infinitesimal transformation:
\[
Xf
=
\sum_{i=1}^n\,\xi_i(x,z)\,
\frac{\partial f}{\partial x_i}
+
\sum_{\mu=1}^m\,\zeta_\mu(x,z)\,
\frac{\partial f}{\partial z_\mu}
\]
and we attempt to form a prolonged infinitesimal transformation
from it:
\[
\aligned
X^{(N)}f
&
=
Xf
+
\sum_{\nu=1}^m\,\sum_\alpha\,
\zeta_{\nu,\,\alpha_1,\dots,\alpha_n}\,
\frac{\partial f}{\partial
z_{\nu,\,\alpha_1,\dots,\alpha_n}}
\\
&
\ \ \ \ \ \ \ \ \ \ \ 
{\scriptstyle{(0\,<\,\alpha_1\,+\,\cdots\,+\,\alpha_n\,\leqslant\,N)}}
\endaligned
\]
which leaves invariant the Pfaffian system of equations:
\def\theequation{31}\begin{equation}
\aligned
&
\D\,z_{\nu,\,\beta_1,\dots,\beta_n}
-
\sum_{i=1}^n\,
z_{\nu,\,\beta_1,\dots,\beta_i+1,\dots,\beta_n}\,
\D\,x_i
=
0
\\
&
\ \ \ \ \ \ \ \ \ \ \ \ \ 
{\scriptstyle{(\nu\,=\,1\,\cdots\,m\,;\,\,\,
0\,\leqslant\,\beta_1\,+\,\cdots\,+\,\beta_n\,<\,N)}}
\endaligned
\end{equation}
For $N = 0$ and $N = 1$, there certainly exists an infinitesimal
transformation $X^{ (N)}f$ of the constitution just described, 
and this does not require any justification.
Hence we can conduct the general proof for the existence of
$X^{ (N)}f$ in the way that we show that, as soon as 
$X^{ (N-1)}f$ and $X^{ (N)}f$ exist, then $X^{ (N+1)}f$ also 
exists.

So the assumption is that $X^{ (N-1)}f$ and $X^{ (N)}f$ exist.
Now, if there is an infinitesimal transformation $X^{ (N+1)}f$
which leaves invariant the system of equations:
\def\theequation{34}\begin{equation}
\aligned
&
\D\,z_{\nu,\beta_1,\dots,\beta_n}
-
\sum_{i=1}^n\,
z_{\nu,\,\beta_1,\dots,\beta_i+1,\dots,\beta_n}\,
\D\,x_i
=
0
\\
&
\ \ \ \ \ \ \ \ \ 
{\scriptstyle{(\nu\,=\,1\,\cdots\,m\,;\,\,\,
0\,\leqslant\,\beta_1\,+\,\cdots\,+\,\beta_n\,<\,N\,+\,1)}}
\endaligned
\end{equation}
then the coefficients still unknown:
\[
\frac{\delta z_{\mu,\,\gamma_1,\dots,\gamma_n}}{
\delta t}
=
\zeta_{\mu,\,\gamma_1,\dots,\gamma_n}
\ \ \ \ \ \ \ \ \ \ \ \ \ 
{\scriptstyle{(\gamma_1\,+\,\cdots\,+\,\gamma_n\,=\,N\,+\,1)}}
\]
must satisfy certain equations of the form:
\[
\aligned
&
\D\,\zeta_{\mu,\,\alpha_1,\dots,\alpha_n}
-
\sum_{i=1}^n\,
z_{\mu,\,\alpha_1,\dots,\alpha_i+1,\dots,\alpha_n}\,
\D\,\xi_i
-
\sum_{i=1}^n\,
\zeta_{\mu,\,\alpha_1,\dots,\alpha_i+1,\dots,\alpha_n}\,
\D\,x_i
\\
&
=
\sum_{\nu=1}^m\,\sum_\beta\,
{\sf P}_{\nu,\,\beta_1,\dots,\beta_n}\,
\bigg\{
\D\,z_{\nu,\,\beta_1,\dots,\beta_n}
-
\sum_{i=1}^n\,
z_{\nu,\,\beta_1,\dots,\beta_i+1,\dots,\beta_n}\,
\D\,x_i
\bigg\}
\\
&
\ \ \ \ \ \ \ \ \ \ \ \ \ \ \ \ \ \ 
{\scriptstyle{(\alpha_1\,+\,\cdots\,+\,\alpha_n\,=\,N\,;\,\,\,
0\,\leqslant\,\beta_1\,+\,\cdots\,+\,\beta_n\,<\,N\,+\,1)}},
\endaligned
\]
and to be precise, independently of the differentials: 
$\D\, x_i$, $\D\, z_{ \nu, \, \beta_1, \dots, \beta_n}$.

\medskip

\renewcommand{\thefootnote}{\fnsymbol{footnote}}
From this at first, the ${\sf P}_{\nu,\, \beta_1, \dots, \beta_n}$ 
determine themselves uniquely and it remains only equations between
the mutually independent differentials: $\D\, x_1, \dots, \D\, x_n$.
Hence if one inserts the values of the ${\sf P}_{ \nu, \beta_1, \dots,
\beta_n}$ and if one compares the coefficients of the $\D\, x_i$ in 
both sides, one obtains for $\zeta_{ \mu, \, \alpha_1, \dots, \alpha_i
+ 1, \dots, \alpha_n}$ the following expression:
\def\theequation{35}\begin{equation}
\aligned
&
\frac{\delta}{\delta t}\,
z_{\mu,\,\alpha_1,\dots,\alpha_i+1,\dots,\alpha_n}
=
\zeta_{\mu,\,\alpha_1,\dots,\alpha_i+1,\dots,\alpha_n}
\\
&
\frac{\D\,\zeta_{\mu,\,\alpha_1,\dots,\alpha_n}}{\D\,x_i}
-
\sum_{j=1}^n\,
z_{\mu,\,\alpha_1,\dots,\alpha_j+1,\dots,\alpha_n}\,
\frac{\D\,\xi_j}{\D\,x_i}
\\
&
\ \ \ \ \ \ \ \ \ \ \ \ \ \ \ \ \ \ 
{\scriptstyle{(\alpha_1\,+\,\cdots\,+\,\alpha_n\,=\,N)}},
\endaligned
\end{equation}
where $\D\, / \D\, x_i$ denotes a complete differential quotient 
with respect to $x_i$.\footnote[1]{\,
The formula~\thetag{ 35} is fundamentally identical with a formula
due to \name{Poisson} in the Calculus of Variations.
}
\renewcommand{\thefootnote}{\arabic{footnote}}

But now, it yet remains a difficulty; indeed, in general, we obtain
for each $\zeta_{ \mu, \, \alpha_1, \dots, \alpha_i + 1, \dots, 
\alpha_n}$ a series of apparently different expressions.

The expression in the right-hand side of~\thetag{ 35} is the
derivation of $\zeta_{ \mu, \, \alpha_1, \dots, \alpha_i + 1, \dots, 
\alpha_n}$ from the value of:
\[
\frac{\delta}{\delta t}\,
\frac{\partial z_{\mu,\,\alpha_1,\dots,\alpha_n}}{\partial x_i}
=
\frac{\delta}{\delta t}\,
z_{\mu,\,\alpha_1,\dots,\alpha_i+1,\dots,\alpha_n}.
\]
But on the other hand, we also have:
\def\theequation{36}\begin{equation}
\aligned
&
\frac{\delta}{\delta t}\,
z_{\mu,\,\alpha_1,\dots,\alpha_i+1,\dots,\alpha_n}
=
\frac{\delta}{\delta t}\,
\frac{\partial z_{\mu,\,\alpha_1,\dots,\alpha_h-1,\dots,
\alpha_i+1,\dots,\alpha_n}}{\partial x_h}
\\
&
=
\frac{\D\,\zeta_{\mu,\,\alpha_1,\dots,\alpha_h-1,\dots,
\alpha_i+1,\dots,\alpha_n}}{\D\,x_h}
-
\sum_{j=1}^n\,
z_{\mu,\,\alpha_1,\dots,\alpha_j+1,\dots,
\alpha_h-1,\dots,\alpha_i+1,\dots,\alpha_n}\,
\frac{\D\,\xi_j}{\D\,x_h},
\endaligned
\end{equation}
where $h$ denotes an arbitrary number amongst $1, 2, \dots, n$
different from $i$. Thus, all possibilities are exhausted. Therefore,
it only remains to prove yet that the latter value of $\frac{ \delta}{
\delta t}\, z_{\mu, \, \alpha_1, \dots, \alpha_i + 1, \dots,
\alpha_n}$ coincides with the value~\thetag{ 35}.

In order to prove this, we remember that the following equations hold:
\[
\aligned
&
\zeta_{\mu,\,\alpha_1,\dots,\alpha_n}
=
\frac{\delta}{\delta t}\,
\frac{\partial z_{\mu,\,\alpha_1,\dots,\alpha_h-1,\dots,\alpha_n}}{
\partial x_h}
\\
&
=
\frac{\D\,\zeta_{\mu,\,\alpha_1,\dots,\alpha_h-1,\dots,\alpha_n}}{
\D\,x_h}
-
\sum_{j=1}^n\,z_{\mu,\,\alpha_1,\dots,\alpha_j+1,\dots,
\alpha_h-1,\dots,\alpha_n}\,
\frac{\D\,\xi_j}{\D\,x_h},
\endaligned
\]
and:
\[
\aligned
&
\zeta_{\mu,\,\alpha_1,\dots,\alpha_h-1,\dots,\alpha_i+1,\dots,\alpha_n}
=
\frac{\delta}{\delta t}\,
\frac{\partial z_{\mu,\,\alpha_1,\dots,\alpha_h-1,\dots,\alpha_n}}{
\partial x_i}
\\
&
=
\frac{\D\,\zeta_{\mu,\,\alpha_1,\dots,\alpha_h-1,\dots,\alpha_n}}{
\D\,x_i}
-
\sum_{j=1}^n\,
z_{\mu,\,\alpha_1,\dots,\alpha_j+1,\dots,\alpha_h-1,\dots,\alpha_n}\,
\frac{\D\,\xi_j}{\D\,x_i}.
\endaligned
\]

If we insert the value of $\zeta_{ \mu, \, \alpha_1, \dots, \alpha_n}$
in~\thetag{ 35}, we obtain:
\[
\aligned
&
\frac{\delta}{\delta t}\,
\frac{\partial z_{\mu,\,\alpha_1,\dots,\alpha_n}}{\partial x_i}
=
\frac{\D}{\D\,x_i}\,
\frac{\D}{\D\,x_h}\,
\zeta_{\mu,\,\alpha_1,\dots,\alpha_h-1,\dots,\alpha_n}
-
\sum_{j=1}^n\,
z_{\mu,\,\alpha_1,\dots,\alpha_j+1,\dots,\alpha_h-1,\dots,
\alpha_i+1,\dots,\alpha_n}\,
\frac{\D\,\xi_j}{\D\,x_h}
\\
&
-
\sum_{j=1}^n\,
z_{\mu,\,\alpha_1,\dots,\alpha_j+1,\dots,\alpha_n}\,
\frac{\D\,\xi_j}{\D\,x_i}
-
\sum_{j=1}^n\,
z_{\mu,\,\alpha_1,\dots,\alpha_j+1,\dots,\alpha_h-1,\dots,\alpha_n}\,
\frac{\D}{\D\,x_i}\,
\frac{\D\,\xi_j}{\D\,x_h}.
\endaligned
\]
On the other hand, if we insert the value of $\zeta_{ \mu, \, \alpha_1,
\dots, \alpha_h - 1, \dots, \alpha_i + 1, \dots, \alpha_n}$ 
in~\thetag{ 36} and if we take into account that $\frac{ \D}{ \D\, 
x_i}\, \frac{ \D}{ \D\, x_h} = \frac{ \D}{\D\, x_h}\, \frac{ \D}{
\D\, x_i}$, we find: 
\[
\frac{\delta}{\delta t}\,
\frac{\partial z_{\mu,\,\alpha_1,\dots,\alpha_h-1,\dots,\alpha_i+1,
\dots,\alpha_n}}{\partial x_h}
=
\frac{\delta}{\delta t}\,
\frac{z_{\mu,\,\alpha_1,\dots,\alpha_n}}{\partial x_i}.
\]
But this is what was to be shown.

At present, it is shown that under the assumptions made, the $\zeta_{
\mu, \, \gamma_1, \dots, \gamma_n}$ ($\gamma_1 + \cdots + \gamma_n = N
+ 1$) really exist and are uniquely determined; consequently, it is
certain that, according to what was said above, to every $Xf$ and for
every value of $N$, there is associated a completely determined
prolonged infinitesimal transformation $X^{ (N)}f$.

The coefficients of $X^{ (N)}f$ are obviously linear and homogeneous
in the $\xi_i$, $\zeta_k$ and their partial differential quotients
with respect to the $x$ and the $z$. Hence if $X_if$ and $X_jf$ are
two infinitesimal transformations of the form $Xf$, and moreover, if
$X_i^{ (N)}f$ and $X_j^{ (N)}f$ are the infinitesimal transformations
prolonged in the way indicated, then:
\[
c_i\,X_i^{(N)}f
+
c_j\,X_j^{(N)}f
\]
results from $c_i\, X_if + c_j\, X_jf$ by means of the prolongation
in question. In addition, since $X_i^{ (N)} X_j^{ (N)} f - 
X_j^{ (N)} X_i^{ (N)}f$ leaves invariant the Pfaffian system of
equations~\thetag{ 31} simultaneously with $X_i^{ (N)}f$ and
$X_j^{ (N)}f$, it follows that:
\[
X_i^{(N)}X_j^{(N)}f
-
X_j^{(N)}X_i^{(N)}f
\]
must come into existence from $X_iX_jf - X_jX_if$ by means of this
prolongation.

Lastly, let $X_1f, \dots, X_rf$ be independent infinitesimal
transformations of the $r$-term group~\thetag{ 32} so that
they satisfy the known relations:
\[
X_iX_jf-X_jX_if
=
\sum_{s=1}^r\,c_{ijs}\,X_sf.
\]

If we form the prolonged infinitesimal transformations: $X_1^{ (N)}f,
\dots, X_r^{ (N)}f$, then $X_i^{ (N)} X_j^{ (N)}f - X_j^{ (N)} X_i^{
(N)}f$ also comes from $X_i X_jf - X_j X_if$ by means of this
prolongation, and consequently, it comes from $\sum\, c_{ iks}\, X_s$,
which again means that we have:
\[
X_i^{(N)}X_j^{(N)}f
-
X_j^{(N)}X_i^{(N)}f
=
\sum_{s=1}^r\,c_{ijs}\,X_s^{(N)}f.
\]
Therefore, the $r$ infinitesimal transformations $X_i^{ (N)}f$
generate, for every value of $N$, a group equally composed with 
the group $X_if$; the former group is obviously identical to
the group~\thetag{ 33} discussed earlier on which was obtained
by prolongating the finite equations~\thetag{ 32} of the group:
$X_1f, \dots, X_rf$.

\renewcommand{\thefootnote}{\fnsymbol{footnote}}
\def\thetheorem{94}\begin{theorem}
If the $\infty^r$ transformations:
\[
\aligned
x_i'
&
=
f_i(x_1,\dots,x_n,\,z_1,\dots,z_m,\,a_1,\dots,a_r)
\ \ \ \ \ \ \ \ \ \ \ \ \ {\scriptstyle{(i\,=\,1\,\cdots\,n)}}
\\
z_k'
&
=
F_k(x_1,\dots,x_n,\,z_1,\dots,z_m,\,a_1,\dots,a_r)
\ \ \ \ \ \ \ \ \ \ \ \ \ {\scriptstyle{(k\,=\,1\,\cdots\,m)}}
\endaligned
\]
in the variables $x_1, \dots, x_n$, $z_1, \dots, z_m$ form an $r$-term
group and if one considers the $z_k$ as functions of the $x_i$ which
can be chosen arbitrarily, then the differential quotients of the $z_k$
with respect to the $x_i$ are also subjected to transformations.
If one takes together all differential quotients from the first
order up to, say, the $N$-th order, then one obtains certain equations:
\[
\aligned
z_{\mu,\,\alpha_1,\dots,\alpha_n}'
&
=
F_{\mu,\,\alpha_1,\dots,\alpha_n}
(x_1,\dots,x_n,\,z_1,\dots,z_m,\,
z_{\nu,\,\beta_1,\dots,\beta_n}\,;\,\,a_1,\dots,a_r)
\\
&
\ \ \ \ \ \ \ \ \ \ \ \ \ \ 
{\scriptstyle{(\beta_1\,+\,\cdots\,+\,\beta_n\,\leqslant\,
\alpha_1\,+\,\cdots\,+\,\alpha_n\,\leqslant\,N)}}
\endaligned
\]
which, when joined to the equations of the original group, represent
an $r$-term group equally composed with the original 
group.\footnote[1]{\,
\name{Lie}, Math. Annalen Vol. XXIV, 1884; Archiv for Math.,
Christiania 1883.
}
\end{theorem}
\renewcommand{\thefootnote}{\arabic{footnote}}

Above, we already mentioned that every given group can be prolonged
in very many different ways; indeed, it is left just as one likes
which variables one wants to choose as the independent ones. 

In addition, one can, from the beginning, substitute the given group
for a group equally composed with it by adding a certain number
of variables: $t_1, \dots, t_\sigma$ which are absolutely not 
transformed by the group, or said differently, that are transformed
only by the identity transformation:
\[
t_1'=t_1,
\,\,\,\dots,\,\,\,
t_\sigma'=t_\sigma.
\]
Now, if one regards as the independent variables an arbitrary number
amongst the original variables and amongst the $t_i$, one can add
differential quotients and prolong; one always comes to an equally
composed group. 

As one sees, the number of possibilities is very large here.

\sectionengellie{\S\,\,\,131.}

The theory of the invariants of an arbitrary group developed earlier
on in Chap.~\ref{kapitel-13} can be immediately applied to our prolonged
groups. 

Since one can always choose the number $N$ so large that the 
infinitesimal transformations $X_i^{ (N)}f$ contain more than $r$ 
independent variables, one can always arrange that the $r$ equations
$X_k^{ (N)}f = 0$ form a complete system with one or more solutions.
These solutions are functions of the $x$, of the $z$ and of the
differential quotients of the latter, they admit every finite
transformation of the prolonged group $X_k^{ (N)}f$ and hence
are absolute invariants of this group; they shall be called
the \terminology{differential invariants} 
\deutsch{Differentialinvarianten} of the original group.

\plainstatement{A function $\Omega$ of $x_1, \dots, x_n$, $z_1, \dots,
z_m$ and of the differential quotients of the $z$ with respect to
the $x$ is called a \terminology{differential invariant}
of the $r$-term group:
\[
\aligned
x_i'
&
=
f_i(x_1,\dots,x_n,\,z_1,\dots,z_m\,;\,\,a_1,\dots,a_r)
\\
z_k'
&
=
F_k(x_1,\dots,x_n,\,z_1,\dots,z_m\,;\,\,a_1,\dots,a_r)
\endaligned
\]
when a relation of the form:
\[
\Omega\big(x_1',\dots,x_n',\,z_1',\dots,z_m',\,
z_{\mu,\,\alpha_1,\dots,\alpha_n}'\big)
=
\Omega(x_1,\dots,x_n,\,z_1,\dots,z_m,\,
z_{\mu,\,\alpha_1,\dots,\alpha_n})
\]
holds identically.}

Since we can choose $N$ arbitrarily large, we have instantly:

\renewcommand{\thefootnote}{\fnsymbol{footnote}}
\def\thetheorem{95}\begin{theorem}
Every continuous transformation group: $X_1f, \dots, X_rf$ determines
an infinite series of differential invariants which define themselves
as solutions of complete systems.\footnote{\,
\name{Lie}, Gesellsch. d. W. zu Christiania 1882; Math. Ann. Vol.
XXIV, 1884.
}
\end{theorem}
\renewcommand{\thefootnote}{\arabic{footnote}}

If one knows the finite equations of the group: $X_1f, \dots, X_rf$,
then in the way explained above, one finds the finite equations of the
prolonged group and afterwards, under the guidance of
Chap.~\ref{kapitel-13}, the differential invariants of any order
without integration. But in general, this method for the determination
of the differential invariants is not practically applicable. In most
cases, the direct integration of the complete system: $X_k^{ (N)}f =
0$ is preferable.

We shall not enter these considerations here. Still, we should only
observe that from sufficiently well known differential invariants, one
can derive new differential invariants by differentiation and by
formation of determinants.

\smallercharacters{

In the variables $x_1, \dots, x_n$, $z_1, \dots, z_m$, if a group is
represented by \emphasis{several} systems of equations, each one with
$r$ parameters, then naturally, there are in the same way prolonged
groups whose differential invariants are those of the
original group. The former general developments not only show
that each such group possesses differential invariants, but also 
at the same time, they show how these differential invariants can be
found.

\sectionengellie{\S\,\,\,132.}

One can also ask for possible \emphasis{systems of differential
equations} which remain invariant by a given group. The determination
of a system of this sort can obviously be carried out by prolonging
the concerned group in a suitable way and by applying the developments
of Chap.~\ref{kapitel-14}, on the basis of which all systems of
equations invariant by the group can be determined. Each system of
equations found in this way then represents a system of differential
equations invariant by the original group.

However, in each individual case, it must yet be specially studied
whether the concerned system of differential equations satisfies the
condition of integrability.\,---

Conversely, one can imagine that a system of differential equations is
given\,---\,integrable or not integrable\,---\,and one can ask the
question whether this system admits a given group. The answer to this
question also presents no difficulty at present. One only has to
prolong the given group in the right way and afterwards, to study
whether the given system of equations admits the prolonged group;
according to Chap.~\ref{kapitel-7}, this study can be conducted
without integration.

\sectionengellie{\S\,\,\,133.}

In order to give a simple application of the preceding theory, we want
to seek the conditions under which a system of differential equations
of the form:
\[
A_k\varphi
=
\sum_{i=1}^n\,\alpha_{ki}(x_1,\dots,x_n)\,
\frac{\partial\varphi}{\partial x_i}
=
0
\ \ \ \ \ \ \ \ \ \ \ \ \ {\scriptstyle{(k\,=\,1\,\cdots\,q)}}
\]
admits the $r$-term group:
\[
X_jf
=
\sum_{i=1}^n\,\xi_{ji}(x_1,\dots,x_n)\,
\frac{\partial f}{\partial x_i}
\ \ \ \ \ \ \ \ \ \ \ \ \ {\scriptstyle{(j\,=\,1\,\cdots\,r)}}
\]
in the variables $x_1, \dots, x_n$; here naturally, the $q$ equations 
$A_k \varphi = 0$ are assumed to be mutually independent. 

In order to be able to answer the question raised, to the variables
$x_1, \dots, x_n$ of the group $X_kf$, we must add yet the variable
$\varphi$ which is not at all transformed by the group. The $x$ are to
be considered as independent variables, $\varphi$ as dependent and the
group: $X_1f, \dots, X_rf$ is thus to be prolonged by adding the $n$
differential quotients:
\[
\frac{\partial\varphi}{\partial x_i}
=
\varphi_i
\ \ \ \ \ \ \ \ \ \ \ \ \ {\scriptstyle{(i\,=\,1\,\cdots\,n)}}.
\]
Afterwards, one has to examine whether the system of equations:
$\sum_k \, \alpha_{ ki} \, \varphi_i = 0$ allows the prolonged group.

To begin with, we compute the infinitely small increment $\delta \,
\varphi_i$ that $\varphi_i$ is given by the infinitesimal
transformation $X_jf$. This increment is to be determined so that
the expression:
\[
\delta\,
\bigg(
\D\,\varphi
-
\sum_{\nu=1}^n\,\varphi_\nu\,\D\,x_\nu
\bigg)
=
\D\,\delta\,\varphi
-
\sum_{\nu=1}^n\,
\big\{
\varphi_\nu\,\D\,\delta\,x_\nu
+
\delta\,\varphi_\nu\,\D\,x_\nu
\big\}
\]
vanishes by virtue of $\D\, \varphi = \sum\, \varphi_\nu\, \D\, x_\nu$.
But since, as remarked above, $\delta \varphi$ is zero, we obtain for
the $\delta \, \varphi_\nu$ the equation:
\[
\sum_{\nu=1}^n\,
\big\{
\varphi_\nu\,\D\,\xi_{j\nu}\,\delta\,t
+
\delta\,\varphi_\nu\,\D\,x_\nu
\big\}
=
0
\]
which must hold identically. Consequently, we have:
\[
\delta\,\varphi_\nu
=
-\,\sum_{\mu=1}^n\,
\frac{\partial\xi_{j\mu}}{\partial x_\nu}\,
\varphi_\mu\,\delta\,t,
\]
and the prolonged infinitesimal transformation $X_jf$ has the shape:
\[
X_j^{(1)}f
=
\sum_{i=1}^n\,\xi_{ji}\,
\frac{\partial f}{\partial x_i}
-
\sum_{i=1}^n\,
\bigg\{
\sum_{\mu=1}^n\,
\frac{\partial\xi_{j\mu}}{\partial x_i}\,
\varphi_\mu
\bigg\}\,
\frac{\partial f}{\partial\varphi_i}.
\]

Now, if the system of equations: $\sum_i\, \alpha_{ ki}(x)\, \varphi_i
= 0$ in the $2\, n$ variables $x_1, \dots, x_n$, $\varphi_1, \dots, 
\varphi_n$ is supposed to admit the infinitesimal transformation
$X_j^{ (1)}f$, then all the $q$ expressions:
\[
X_j^{(1)}
\bigg(
\sum_{i=1}^n\,\alpha_{ki}\,\varphi_i
\bigg)
=
\sum_{i=1}^n\,X_j\,\alpha_{ki}\,\varphi_i
-
\sum_{\mu=1}^n\,
\bigg\{
\sum_{i=1}^n\,\alpha_{ki}\,
\frac{\partial\xi_{j\mu}}{\partial x_i}
\bigg\}\,
\varphi_\mu
\ \ \ \ \ \ \ \ \ \ \ \ \ {\scriptstyle{(k\,=\,1\,\cdots\,q)}}
\]
must vanish by virtue of the system of equations. This necessary, and 
at the same time sufficient, condition is then satisfied only when 
relations of the form:
\[
X_j^{(1)}
\bigg(
\sum_{i=1}^n\,\alpha_{ki}\,\varphi_i
\bigg)
=
\sum_{\sigma=1}^q\,\rho_{jk\sigma}(x)\,
\sum_{i=1}^n\,\alpha_{\sigma i}\,\varphi_i
\]
hold, where the $\rho_{ jk\sigma}$ denote functions of $x_1, \dots, x_n$
only which do not depend on the $\varphi_i$.

With these words, we have found the desired conditions; they can be 
written:
\[
\sum_{i=1}^r\,
\big(
X_j\,\alpha_{ki}
-
A_k\,\xi_{ji}
\big)\,\varphi_i
=
\sum_{\sigma=1}^q\,
\rho_{jk\sigma}\,
\sum_{i=1}^n\,\alpha_{\sigma i}\,\varphi_i,
\]
or, if we insert again $\partial \varphi / \partial x_i$ in
place of $\varphi_i$:
\[
X_j\big(A_k(\varphi)\big)
-
A_k\big(X_j(\varphi)\big)
=
\sum_{\sigma=1}^q\,\rho_{jk\sigma}(x)\,A_\sigma\varphi.
\]

Relations of this sort must be satisfied for every arbitrary function
$\varphi ( x_1, \dots, x_n)$. The system of the $q$ linear partial
differential equations: $A_1 \varphi = 0$, \dots, $A_q \varphi = 0$
always remains invariant by every transformation of the group: $X_1f,
\dots, X_rf$ when these relations hold, and only when they hold.

In the case where the $q$ equations: $A_k \varphi = 0$ form a $q$-term
complete system, this result is not new for us. Indeed, in this
special case, we have already indicated in Chap.~\ref{kapitel-8},
Theorem~20, p.~\pageref{Theorem-20-S-140} the necessary and sufficient
condition just found. However, our present developments accomplish
more than the developments done at that time, for we have shown at
present that the criterion in question holds generally, also when the
equations: $A_k \varphi = 0$ do not form a $q$-term complete system.

\linestop

The origins of the theory of the differential invariants goes back 
long ago; indeed, mathematicians of the previous century already 
have considered and integrated the differential invariants associated
to several specially simple groups. 

For instance, it is known long since that every differential equation
of first order between $x$ and $y$ in which one variable, say $y$, does
not explicitly appear, can be integrated by quadrature; but obviously:
\[
f\bigg(
x,\,\,\frac{\D\,y}{\D\,x}
\bigg)
=
0
=
f(x,\,y')
\]
is nothing but the most general differential equation of first order
that admits the one-term group:
\def\theequation{37}\begin{equation}
\mathfrak{x}
=
x,
\ \ \ \ \ \ \ \ \
\mathfrak{y}
=
y+a
\end{equation}
with the parameter $a$. The most general integral equation
\deutsch{Integralgleichung} can be deduced from a particular integral
equation, and in fact, as we can now say, in such a way that one
executes the general transformation of the one-term group~\thetag{ 37}
on the particular solution in question; thanks to this, one indeed
obtains the equation: $\mathfrak{ y} = \varphi ( \mathfrak{ x}) + a$
with the arbitrary constant $a$.

Furthermore, the homogeneous differential equation:
\[
f\bigg(
\frac{y}{x},\,\,
\frac{\D\,y}{\D\,x}
\bigg)
=
0
=
f
\bigg(
\frac{y}{x},\,y'
\bigg)
\]
is the general form of a differential equation of first order
which admits the one-term group: $\mathfrak{ x} = a\, x$, $\mathfrak{ y}
= a \, y$; here, $f( y / x, \, y')$ is the most general first order
differential invariant associated to this group.

It has been observed since a long time that every particular integral
equation: $F ( x, y) = 0$ of a homogeneous differential equation: $f = 
0$ can, by executing the general transformation: $\mathfrak{ x} = a\,
x$, $\mathfrak{ y} = a\, y$, be transferred to the general integral
equation:
\[
F\bigg(
\frac{\mathfrak{x}}{a},\,\,
\frac{\mathfrak{y}}{a}
\bigg)
=
0.
\]
However, the equation $xy' - y = 0$ must be disregarded here.

A third example is provided by differential equations of the form:
\[
f\bigg(
\frac{\D^m\,y}{\D\,x^m},\,\,
\frac{\D^{m+1}y}{\D\,x^{m+1}}
\bigg)
=
0
=
f\big(
y^{(m)},\,
y^{(m+1)}
\big)\,;
\]
nevertheless, it is not necessary here to write down the group
of all equations of this form.

In the invariant theory of linear transformations, there often appear
true differential invariants relatively to all linear transformations.
They have been considered by \name{Cayley} for the first time.
Nevertheless, it is to be observed on the occasion \emphasis{firstly}
that the differential invariant of \name{Cayley} are not the simplest
ones which are associated to the general linear homogeneous group, and
\emphasis{secondly}, that \name{Cayley} has not considered invariants
of \emphasis{differential} equations, and even less has integrated
such equations.

In a prized essay \deutsch{Preisschrift} achieved in 1867 and
published in 1871 (Determination of a special minimal surface,
Akad. d. W. zu Berlin), \name{H. Schwarz} considered differential
equations of the form:
\[
J
=
\frac{y'\,y'''-\frac{3}{2}\,{y''}^2}{{y'}^2}
-
f(f)
=
0
\]
which, as he himself indicated, had already appeared occasionally
apud \name{Lagrange}. \name{Schwarz} observed that the most general
solution can be deduced from every particular solution: $y = \varphi
(x)$, namely the former has the form:
\[
y
=
\frac{a+b\,\varphi}{c+d\,\varphi},
\]
with the arbitrary constants: $a \colon b \colon c \colon d$. Thus, 
as we can say, the expression $J$ is a differential invariant, 
and to be precise, the most general third order differential invariant 
of the group:
\[
\mathfrak{x}
=
x,
\ \ \ \ \ \ \ \ \ \ \
\mathfrak{y}
=
\frac{a+b\,y}{c+d\,y}.
\]

But now, although all these special theories are undoubtedly valuable,
it is however still to be remarked that the inner connection 
between them, the general principle from which they flow
\deutsch{fliessen}, was missed by the mathematicians. They have not
observed that differential invariants are associated to \emphasis{every}
finite continuous group.

\renewcommand{\thefootnote}{\fnsymbol{footnote}}
In the years 1869--1871, \name{Lie} occupied himself with differential
equations which allow interchangeable infinitesimal transformations and
in 1874, he published a work already announced in 1872 about a
\emphasis{general integration theory of ordinary differential
equations which admit an arbitrary continuous group of 
transformations}.\footnote[1]{\,
Verhandlungen der Gesellsch. d. W. zu Christiania 1870--1874; Math.
Ann. Vol. V, XI; Gött. Nachr. 1874.
}

Afterwards, \name{Halphen} computed the simplest differential invariants
relatively to all projective transformations and he gave in addition 
nice applications to the theory of the \emphasis{linear} differential
equations.\footnote[2]{\,
Thèse sur les invariants différentiels 1878; Journal de l'\'Ecole
Pol. 1880. Cf. also Comptes Rendus Vol. 81, 1875, p.~1053; Journal
de Liouville Novemb. 1876. Mémoire sur la réduction des équat. diff.
lin. aux formes intégrables 1880--1883.
}

After that, \name{Lie} developed in the years 1882--1884 a
\emphasis{general theory of the differential invariants} of the
finite and infinite continuous groups, where he specially 
occupied himself with finite groups in two variables.\footnote[3]{\,
Verh. d. Gesellsch. d. W. zu Christiania, 1882, 1883 and February
1875; Archiv for Math. 1882, 1883; Math. Ann. Vol. XXIV, 1884.
} 
As far as it relates to finite groups, this general theory is explained
in what precedes.
\renewcommand{\thefootnote}{\arabic{footnote}}

Finally, after 1884, \name{Sylvester} and many other English and
American mathematicians have published detailed but \emphasis{special}
studies about differential invariants.

}

\linestop


\chapter{The General Projective Group}
\label{kapitel-26}
\chaptermark{The General Projective Group}

\setcounter{footnote}{0}

\abstract*{??}

The equations:
\def\theequation{1}\begin{equation}
x_\nu'
=
\frac{a_{1\nu}\,x_1+\cdots+a_{n\nu}\,x_n+a_{n+1,\nu}}{
a_{1,n+1}\,x_1+\cdots+a_{n,n+1}\,x_n+a_{n+1,n+1}}
\ \ \ \ \ \ \ \ \ \ \ \ \
{\scriptstyle{(\nu\,=\,1\,\cdots\,n)}}
\end{equation}
determine a group, as one easily convinces oneself, the so-called
\terminology{general projective group}\, of the manifold $x_1, \dots,
x_n$. In the present chapter, we want to study somehow more closely
this important group, which is also called the group of all
\terminology{collineations}\, 
\deutsch{Collineationen}
of the space $x_1, \dots, x_n$, by
focusing our attention especially on its subgroups.

\sectionengellie{\S\,\,\,134.}

The $( n + 1)^2$ parameters $a$ are not all essential: there indeed
appears just their ratios; one of the parameters, best $a_{ n + 1, n+
1}$, can hence be set equal to 1. The values of the parameters are
subjected to the restriction that the substitution determinant
\deutsch{Substitutionsdeterminant} $\sum \, \pm \, a_{ 11} \cdots a_{
n+1, n+1}$ should not be equal to zero, because simultaneously with
it, the functional determinant \deutsch{Functionaldeterminant}:
$\sum\, \pm \, \frac{ \partial x_1'}{\partial x_1} \dots\, \frac{
\partial x_n'}{ \partial x_n}$ would also vanish.

The identity transformation is contained in our group, it corresponds
to the parameter values:
\[
a_{\nu\nu}=1,\ \ \ \ \ \ \
a_{\mu\nu}=0
\ \ \ \ \ \ \ \ \ \ \ \ \
{\scriptstyle{(\mu,\,\nu\,=\,1\,\cdots\,n\,+\,1,\,\,\mu\,\neq\,\nu)}},
\]
for which indeed it comes $x_i ' = x_i$. As a consequence of that,
one obtains the infinitesimal transformations of the group by giving
to the $a_{ \mu \nu}$ the values:
\[
a_{\nu\nu}
=
1+\omega_{\nu\nu},\ \ \ \ \ \ \
a_{n+1,n+1}
=
1,\ \ \ \ \ \ \
a_{\mu\nu}
=
\omega_{\mu\nu},
\]
where the $\omega_{ \mu \nu}$ mean infinitesimal quantities. Thus one
finds:
\[
x_\nu'
=
\bigg(
x_\nu
+
\sum_{\mu=1}^n\,
\omega_{\mu\nu}\,x_\mu
+
\omega_{n+1,\nu}
\bigg)
\bigg(
1
-
\sum_{\mu=1}^n\,
\omega_{\mu,n+1}\,x_\mu
+\cdots
\bigg),
\]
or by leaving out the quantities of second or higher order:
\[
x_\nu'-x_\nu
=
\sum_{\mu=1}^n\,
\omega_{\mu\nu}\,x_\mu
+
\omega_{n+1,\nu}
-
x_\nu\,\sum_{\mu=1}^n\,
\omega_{\mu,n+1}\,x_\mu.
\]
If one sets here all the $\omega_{ \mu \nu}$ with the exception of a
single one equal to zero, then one recognizes gradually that our
group contains the $n\, ( n + 2)$ independent infinitesimal
transformations:
\def\theequation{2}\begin{equation}
\frac{\partial f}{\partial x_i},\ \ \ \ \ \ \
x_i\,\frac{\partial f}{\partial x_k},\ \ \ \ \ \ \
x_i\,\sum_{j=1}^n\,x_j\,\frac{\partial f}{\partial x_j}
\ \ \ \ \ \ \ \ \ \ \ \ \
{\scriptstyle{(i,\,\,k\,=\,1\,\cdots\,n)}}.
\end{equation}

\plainstatement{The general projective group of the $n$-times extended
space $x_1, \dots, x_n$ therefore contains $n \, ( n + 2)$ essential
parameters and is generated by infinitesimal transformations. The
analytic expressions of the latter behave regularly for every point of
the space. }

\emphasis{From now on, we will as a rule write $p_i$ for 
$\partial f / \partial x_i$}. In addition, for reasons of convenience,
we want in this chapter to introduce the abbreviations:
\[
x_i\,p_k
=
T_{ik},
\ \ \ \ \ \ \ \ \ \ \
x_i\,\sum_{k=1}^n\,x_k\,p_k
=
P_i.
\]
Lastly, we still want to agree that $\varepsilon_{ ik}$ should mean
zero whenever $i$ and $k$ are distinct from each other, whereas by
contrast $\varepsilon_{ ii}$ is supposed to have the value 1; this is
a fixing of notation that we have already adopted from time to
time. On such a basis, we can write as follows the relations which
come out through combination of the infinitesimal transformations
$p_i$, $T_{ ik}$, $P_i$:
\[
\aligned
\big\leftbracket p_i,\,p_k\big\rightbracket
&
=
0,\ \ \ \ \ \ \ 
\big\leftbracket P_i,\,P_k\big\rightbracket
=
0,\ \ \ \ \ \ \
\big\leftbracket p_i,\,P_k\big\rightbracket
=
T_{ki}
+
\varepsilon_{ik}\,\sum_{\nu=1}^n\,T_{\nu\nu},
\\
&
\big\leftbracket p_i,\,T_{k\nu}\big\rightbracket
=
\varepsilon_{ik}\,p_\nu,\ \ \ \ \ \ \ \
\big\leftbracket P_i,\,T_{k\nu}\big\rightbracket
=
-\,\varepsilon_{i\nu}\,P_k,
\\
&
\ \ \ \ \ \ \ \ \ \ \ \
\rule[-3pt]{0pt}{21pt}
\big\leftbracket T_{ik},\,T_{\mu\nu}\big\rightbracket
=
\varepsilon_{k\mu}\,T_{i\nu}
-
\varepsilon_{\nu i}\,T_{\mu k}.
\endaligned
\]

One easily convinces oneself that these relations remain unchanged
when one substitutes in them the $p_i$, $T_{ ik}$ and $P_i$ by the
respective expressions standing under them in the pattern:
\def\theequation{3}\begin{equation}
\aligned
&
p_i,\ \ \ \ \ \
T_{ik},\ \ \ \
P_i
\\
&
P_i,\ \ \
-\,T_{ki},\ \ \ \
p_i.
\endaligned
\end{equation}
\emphasis{Thus in this way, the general projective group is
related to itself in a holoedrically isomorphic way.}

\smallskip

\smallercharacters{

One could presume that there is a transformation:
$x_i' = \Phi_i ( x_1, \dots, x_n)$ which transfers the infinitesimal
transformations:
\[
p_i,\ \ \ \ \ \ \
x_ip_k,\ \ \ \ \ \ \ \
x_i\,\sum_{k=1}^n\,x_k\,p_k
\]
to, respectively:
\[
x_i'\,\sum_{k=1}^n\,x_k'\,p_k',\ \ \ \ \
-\,x_k'\,p_i',\ \ \ \ 
p_i'.
\]
But there is no such transformation, and this, for simple reasons,
because the $n$ infinitesimal transformations $p_1, \dots, p_n$
generate an $n$-term \emphasis{transitive} group, while: $x_1' \, \sum\,
x_k' p_k'$, \dots, $x_n' \, \sum \, x_k' p_k'$ generate an $n$-term
\emphasis{intransitive} group.

}

First in the next Volume we will learn to see the full signification
of this important property of the projective group, when the concept
of contact transformation \deutsch{Berührungstransformation} and
especially the duality will be introduced.

\medskip

The general infinitesimal transformation:
\[
\sum_{i=1}^n\,a_i\,p_i
+
\sum_{i,\,\,k=1}^n\,b_{ik}\,T_{ik}
+
\sum_{i=1}^n\,c_i\,P_i
\]
of our group is expanded in powers of $x_1, \dots, x_n$ and visibly
contains only terms of zeroth, first and second order in the $x$. One
realizes easily that the group contains $n$ independent infinitesimal
transformations of zeroth order in $x$, out of which no infinitesimal
transformation of first or second order in the $x$ can be deduced
linearly. For instance, $p_1, \dots, p_n$ are $n$ such infinitesimal
transformations. From this it follows that {\it the general projective
group is transitive}.

Besides, there are $n^2$ infinitesimal transformations of first order
in the $x_i$, for instance all $x_ip_k = T_{ ik}$, out of which no 
transformation
of second order can be deduced linearly. Finally it yet arises $n$
transformations of second order in the $x$:
\[
x_i\,\sum_{k=1}^n\,x_kp_k
=
P_i.
\]
In agreement with the Proposition~9 of the Chap.~\ref{kapitel-15} on
p.~\pageref{Satz-9-S-264}, the $P_i$ are pairwise interchangeable and
in addition, the $T_{ ik}$ together with the $P_i$ generate a subgroup
in which the group of the $P_i$ is contained as an invariant
subgroup.

As one sees, and also as it follows from our remark above about the
relationship between the $p_i$ and the $P_i$, the $p_i$ are also
interchangeable in pairs and they generate together with the $T_{ ik}$ a
subgroup in which the group of the $p_i$ is invariant.

\sectionengellie{\S\,\,\,135.}

For the most important subgroups of the general projective group, it
is advisable to employ special names. If, in the general
expression~\thetag{ 1} of a projective transformation, one lets the
denominator reduce to 1, then one gets a \emphasis{linear}
transformation:
\[
x_\nu'
=
a_{1\nu}\,x_1
+\cdots+
a_{n\nu}\,x_n
+
a_{n+1,\nu}
\ \ \ \ \ \ \ \ \ \ \ \ \
{\scriptstyle{(\nu\,=\,1\,\cdots\,n)}};
\]
all transformations of this form constitute the so-called
\terminology{general linear} group. We have already indicated at the
end of the previous paragraph the infinitesimal transformations of
this group; they are deduced by linear combination from the the
following $n ( n+1)$ ones:
\[
p_i,\ \ \ \ \
x_i\,p_k
\ \ \ \ \ \ \ \ \ \ \ \ \
{\scriptstyle{(i,\,k\,=\,1\,\cdots\,n)}}.
\]

If one interprets $x_1, \dots, x_n$ as coordinates of an $n$-times
extended space $R_n$ and if one translates the language of 
the ordinary space, then one can say that the general linear group
consists of all projective transformations which leave invariant the
infinitely far $(n-1)$-times extended straight manifold, or briefly, the
\terminology{infinitely far plane} \deutsch{unendlich ferne ebene}
$M_{ n - 1}$.

Next, if one remembers that by execution of two finite linear
transformations one after the other, the substitution determinants:
$\sum\, \pm a_{ 11} \cdots a_{ nn}$ multiply with each other, 
then one realizes
without difficulty that the totality of all linear transformations
whose determinant equals 1 constitutes a subgroup, and in fact, an
invariant subgroup, which we want to call the \terminology{special
linear group}. One finds easily that, as the $n \,( n + 1) - 1$
independent infinitesimal transformations of this group, the following
ones can be chosen:
\[
p_i,\ \ \ \ \
x_ip_k,\ \ \ \ \ \
x_ip_i-x_kp_k
\ \ \ \ \ \ \ \ \ \ \ \ \
{\scriptstyle{(i\,\gtrless\,k)}}.
\]

If, amongst all linear transformations, one restricts oneself to those
homogeneous in $x$, then one obtains the \terminology{general linear
homogeneous group}:
\[
x_\nu'
=
a_{1\nu}\,x_1
+\cdots+
a_{n\nu}\,x_n
\ \ \ \ \ \ \ \ \ \ \ \ \
{\scriptstyle{(\nu\,=\,1\,\cdots\,n)}},
\]
all infinitesimal transformations of which possess the form: $\sum \,
b_{ ik}\, x_i p_k$ and hence can be linearly deduced from the $n^2$
transformations: $x_i p_k$. Also this group visibly contains an
invariant subgroup, the \terminology{special linear homogeneous
group}, for which: $\sum\, \pm a_{ 11} \cdots a_{ nn}$ has the value
1. The $n^2 - 1$ infinitesimal transformations of the latter group
are:
\[
x_ip_k,\ \ \ \ \ \ \
x_ip_i-x_kp_k
\ \ \ \ \ \ \ \ \ \ \ \ \ \ \
{\scriptstyle{(i\,\gtrless\,k)}}\,;
\]
hence, the general infinitesimal transformation of the group in
question has the form: $\sum_{ i, k}\, \alpha_{ ik}\, x_i p_k$, where
the $n^2$ arbitrary constants $\alpha_{ ik}$ are only subjected to the
condition $\sum \, \alpha_{ ii} = 0$.

\medskip

Since the expression: $\big\leftbracket x_i p_k, \, \sum_j\, x_j p_j
\big\rightbracket$ always vanishes, it is obvious that the last two
named groups are \emphasis{systatic} and consequently
\emphasis{imprimitive}. Indeed, if one sets:
\[
\frac{x_i}{x_n}
=
y_i,
\ \ \ \ \ \ \
\frac{x_i'}{x_n'}
=
y_i'
\ \ \ \ \ \ \ \ \ \ \ \ \
{\scriptstyle{(i\,=\,1\,\cdots\,n-1)}},
\]
then one receives:
\[
y_\nu'
=
\frac{a_{1\nu}\,y_1+\cdots+a_{n-1,\nu}\,y_{n-1}+a_{n\nu}}{
a_{1,n}\,y_1+\cdots+a_{n-1,n}\,y_{n-1}+a_{n,n}}
\ \ \ \ \ \ \ \ \ \ \ \ \
{\scriptstyle{(\nu\,=\,1\,\cdots\,n-1)}}.
\]
From this, it results that in both cases the $y$ are transformed by the
$(n^2 - 1)$-term general projective group of the $( n - 1)$-times
extended manifold $y_1, \dots, y_{ n - 1}$. Consequently, this group
is isomorphic with the general linear homogeneous group of an $n$-times
extended manifold and with the special linear homogeneous group as
well, \label{S-558-bis}
though the isomorphism is holoedric only for the special linear
homogeneous group, since this one contains $n^2 - 1$ parameters.

\def\thetheorem{96}\begin{theorem}
The special linear homogeneous group:
\[
x_ip_k,\ \ \ \ \ \ 
x_ip_i-x_kp_k
\ \ \ \ \ \ \ \ \ \ \ \ \
{\scriptstyle{(i\,\gtrless\,k\,=\,1\,\cdots\,n)}}
\]
in the variables $x_1, \dots, x_n$ is imprimitive and holoedrically
isomorphic with the general projective group of an $(n - 1)$-times
extended manifold.
\end{theorem}

The formally simplest infinitesimal transformations of the general
projective group are $p_1, \dots, p_n$; these generate, as already
observed, a group actually: the group of all
\terminology{translations}:
\[
x_i'
=
x_i
+
a_i
\ \ \ \ \ \ \ \ \ \ \ \ \ \ \ \
{\scriptstyle{(i\,=\,1\,\cdots\,n)}},
\] 
which obviously is simply transitive.

In fact, $m$ arbitrary infinitesimal translations, 
for instance $p_1, \dots, p_m$, always generate an $m$-term
group. For all of these groups, the following
holds:

\def\theproposition{1}\begin{proposition}
\label{S-558}
All $m$-term groups of translations are conjugate to each other inside
the general projective group, and even already inside the general linear
group.
\end{proposition}

Indeed, the $m$ independent infinitesimal transformations of such a
group always have the form:
\[
\sum_{\nu=1}^n\,b_{\mu\nu}\,p_\nu
\ \ \ \ \ \ \ \ \ \ \ \ \
{\scriptstyle{(\mu\,=\,1\,\cdots\,m)}},
\]
where not all $m\times m$ determinants of the $b_{\mu \nu}$ vanish.

But we can very easily show that by means of some linear transformation,
new variables $x_1', \dots, x_n'$ can be introduced for which one has:
\[
p_\mu'
=
\sum_{\nu=1}^n\,b_{\mu\nu}\,p_\nu
\ \ \ \ \ \ \ \ \ \ \ \ \
{\scriptstyle{(\mu\,=\,1\,\cdots\,m)}}.
\]
In fact, let $p_\mu' = p_1 \, \frac{ \partial x_1}{\partial x_\mu'} +
\cdots + p_n\, \frac{ \partial x_n}{ \partial x_\mu'}$; then we only
need to set:
\[
\frac{\partial x_\nu}{\partial x_\mu'}
=
b_{\mu\nu}
\ \ \ \ \ \ \ \ \ \ \ \ \
{\scriptstyle{(\nu\,=\,1\,\cdots\,n;\,\,\,\mu\,=\,1\,\cdots\,m)}},
\]
while the $\frac{ \partial x_\nu}{\partial x_{ m+1}'}, \dots, \frac{
\partial x_\nu}{ \partial x_n'}$ remain arbitrary. We can give to
these last ones some values such that the equations:
\[
x_\nu
=
\sum_{\mu=1}^n\,b_{\mu\nu}\,x_\mu'
+
\sum_{\pi=m+1}^n\,c_{\pi\nu}\,x_\pi'
\ \ \ \ \ \ \ \ \ \ \ \ \
{\scriptstyle{(\nu\,=\,1\,\cdots\,n)}}
\]
determine a transformation, and then, this transformation transfers the
given group of translations to the group $p_1', \dots, p_m'$. From
this, our proposition follows immediately.

\smallskip

\smallercharacters{We 
want to at least indicate a second proof of the same proposition.
As already observed, the general linear group leaves invariant the
infinitely far plane $M_{ n-1}$, and in fact, it is even the most
general projective group of this nature. Now, every infinitesimal
translation is directed by an infinitely far point and is completely
determined by this point; \emphasis{every $m$-term group of translations
can therefore be represented by an $m$-times extended, infinitely far,
straight manifold $M_m$}. But two infinitely far straight $M_m$
always can be transferred one to the other by a linear transformation
which leaves invariant the infinitely far plane. Consequently, all
$m$-term groups of translations are conjugate to each other inside the
general linear group, and in the same way, inside the general
projective group. 

}

The correspondence indicated earlier on which takes place between
the $p_i$ and the $P_i$ yields, as we prove instantly, the

\def\theproposition{2}\begin{proposition}
\label{Satz-2-S-559}
All $m$-term groups, whose infinitesimal transformations possess the
form $\sum \, e_i \, P_i$, are conjugate to each other inside the
general projective group.
\end{proposition}

For the proof, we start from the fact that two subgroups are
conjugate inside a group $G_r$ when the one can be, by means of a
transformation of the adjoint group of $G_r$, transferred to the
other; here, we have to imagine the subgroup as a straight manifold 
in the space $e_1, \dots, e_r$ which is transformed by the adjoint group
(cf. Chap.~\ref{kapitel-16}, 
p.~\pageref{S-280-bis}). If we now
write the transformations of the projective group firstly in the
sequence $p_i$, $T_{ ik}$, $P_i$ and secondly 
in the sequence $P_i$, $-\, T_{
ki}$, $p_i$, then in the two cases we get the same adjoint group. But
since two $m$-term groups of translations can always be transferred
one to the other by the adjoint group, this must also always be the
case with two $m$-term groups whose infinitesimal transformations can
be deduced linearly from the $P_i$. Furthermore, it even immediately
comes out that two $m$-term groups of this sort are already conjugate
to each other inside the group $P_i$, $T_{ ik}$. With that, our
proposition is proved.

\sectionengellie{\S\,\,\,136.}

We consider now one after the other the general projective group, the
general linear group and the general
linear homogeneous group, and to be
precise, we want to examine whether there are invariant subgroups and
which one are contained in these three groups.

To begin with, the general projective group. Let:
\[
\label{S-560}
S
=
\sum_{i=1}^n\,\alpha_i\,p_i
+
\sum_{i=1}^n\,\sum_{k=1}^n\,\beta_{ik}\,x_ip_k
+
\sum_{i=1}^n\,\gamma_i\,x_i\,\sum_{k=1}^n\,x_kp_k
\]
be an infinitesimal transformation of an invariant subgroup; then
necessarily $\leftbracket p_\nu, \, S \rightbracket$ and $\big
\leftbracket p_\mu, \, \leftbracket p_\nu, \, S \rightbracket \big
\rightbracket$ are also transformations of the same
subgroup. Consequently, in our invariant subgroup, there would
certainly appear an infinitesimal translation $\sum\, \rho_i \,
p_i$.

But because all infinitesimal translations are conjugate to each other
inside the general projective group, they would all appear.
Furthermore, since it is invariant, the subgroup would necessarily
contain all transformations: $\big\leftbracket p_i, \, x_i \, \sum_j
\, x_j p_j \big\rightbracket$, or after computation:
\[
x_i\,p_k
\ \ \ \ \ \ \
{\scriptstyle{(i\,\gtrless\,k)}},
\ \ \ \ \ \ \ \ \ \ \ \ \
x_i\,p_i
+
\sum_{j=1}^n\,x_jp_j.
\]
Adding the $n$ transformations: $x_ip_i + \sum_j \, x_j p_j$, one
obtains: $(n+1) \, \sum\, x_j p_j$, hence $x_i p_i$ and therefore 
actually all $x_i p_k$. Finally, the invariant subgroup would yet
contain all transformations: $\big\leftbracket x_i p_i, \, x_i \,
\sum_k\, x_k p_k \big\rightbracket$, hence all $x_i \, \sum_k \, x_k
p_k$ and thus it would be identical to the general projective group
itself. Thus, our first result is:

\renewcommand{\thefootnote}{\fnsymbol{footnote}}
\def\thetheorem{97}\begin{theorem}
The general projective group in $n$ variables is simple.\footnote[1]{\,
Lie, Math. Ann., Vol.~XXV, p.~130.
}
\end{theorem}
\renewcommand{\thefootnote}{\arabic{footnote}}

Correspondingly, the special linear homogeneous group:
\def\theequation{4}\begin{equation}
x_ip_k,
\ \ \ \ \ \ \ \
x_ip_i-x_kp_k
\ \ \ \ \ \ \ \ \ \ \ \ \
{\scriptstyle{(i\,\gtrless\,k)}}
\end{equation}
is also simple.

\smallskip

The general linear homogeneous group with the $n^2$ infinitesimal
transformations $x_i p_k$ contains, as we saw above, an
invariant subgroup with $n^2 - 1$ parameters, namely the
group~\thetag{ 4} just named.

If there is yet a second invariant subgroup, then this subgroup
obviously cannot comprise the group~\thetag{ 4}, and in the same way,
it even cannot have infinitesimal transformations in common with the
same group, since such transformations would constitute an invariant
subgroup in the simple group~\thetag{ 4} (cf. Proposition~10 of the
Chap.~\ref{kapitel-15} on p.~\pageref{Satz-10-S-264}). Taking the
Proposition~7 of Chap.~\ref{kapitel-12} on p.~\pageref{Satz-7-S-211}
into account, it hence follows that a possible second invariant
subgroup can contain only one infinitesimal transformation, and to be
precise, one of the form:
\[
\sum_{i=1}^n\,x_ip_i
+
\sum_{i,\,\,k}^{1\cdots\,n}\,
\alpha_{ik}\,x_ip_k
\ \ \ \ \ \ \ \ \ \ \ \ \
{\scriptstyle{\big(\sum_{i=1}^n\,\,\alpha_{ii}\,=\,0\big)}}.
\]
In addition, according to Proposition~11 of Chap.~\ref{kapitel-15} on
p.~\pageref{Satz-11-S-264}, the same transformation must be
interchangeable with every transformation of the group~\thetag{ 4}, from
which it follows that the transformation:
\[
\sum_{i,\,\,k}^{1\cdots\,r}\,
\alpha_{ik}\,x_ip_k
\ \ \ \ \ \ \ \ \ \ \ \
{\scriptstyle{\big(\sum_{i=1}^n\,\,\alpha_{ii}\,=\,0\big)}}
\]
must be excellent inside the group~\thetag{ 4}.
But there is no such transformation, whence all the $\alpha_{ ik}$
vanish and it shows up that $x_1 p_1 + \cdots + x_n p_n$ and~\thetag{
4} are the only two invariant subgroups of the group $x_i p_k$.

\def\thetheorem{98}\begin{theorem}
The general linear homogeneous group $x_i p_k$ in $n$ variables
contains only two invariant subgroups, namely the special linear
homogeneous group and the one-term group: $x_1 p_1 + \cdots + x_n
p_n$.
\end{theorem}

At present, one easily manages to set up all invariant groups of the
general linear group. Let:
\[
S
=
\sum_{i=1}^n\,\alpha_i\,p_i
+
\sum_{i=1}^n\,\sum_{k=1}^n\,\beta_{ik}\,x_ip_k
\]
be a transformation of such a subgroup. Then together with $S$, also
$\leftbracket p_j, S \rightbracket$ belongs to the invariant subgroup;
hence this subgroup certainly contains a translation, and because of
Proposition~1, p.~\pageref{S-558}, it contains all of them. The
smallest invariant subgroup therefore consists of the translations
themselves; every other subgroup must contain, aside from the
translations, yet a series of infinitesimal transformations of the
form: $\sum_i \, \sum_k\, \alpha_{ ik }\, x_i p_k$. But these last
transformations visibly generate an invariant subgroup, the linear
homogeneous group $x_i p_k$. So we find:

\renewcommand{\thefootnote}{\fnsymbol{footnote}}
\def\thetheorem{99}\begin{theorem}
The general linear group: $p_i$, $x_i p_k$ contains only three
invariant subgroups\footnote[1]{\,
Lie, Math. Ann., Vol.~XXV, p.~130.
}, 
namely the three ones:
\[
p_i
\ \ \ \ \ \ \ \
p_i,\ \
x_1p_1+\cdots+x_np_n
\ \ \ \ \ \ \ \
p_i,\ \
x_ip_k,\ \
x_ip_i-x_kp_k
\ \ \ \ \ \ \ \ \ \ \ \ \
{\scriptstyle{(i\,\gtrless\,k)}},
\]
with respectively $n$, $n+1$ and $n^2 + n - 1$ parameters.
\end{theorem}
\renewcommand{\thefootnote}{\arabic{footnote}}

If, as already done several times, we employ the terminology which is
common for the ordinary space, we can say: the three invariant
subgroups of the general linear group are firstly the group of all
translations, secondly the group of all similitudes
\deutsch{Aehnlichkeitstransformationen}: $(x_1 - x_1^0) p_1 + \cdots
+ (x_n - x_n^0) p_n$, and lastly the most general linear group which
leaves all volumes unchanged.

\sectionengellie{\S\,\,\,137.}

Before we pass to determining the largest subgroups of the general 
projective group, we put on ahead an observation which will find
several applications in the sequel.

Let $X_1f, \dots, X_rf$ be independent infinitesimal transformations
of an $r$-term group $G_r$ and let a family of $\infty^{ r - m - 1}$
infinitesimal transformations of this group be determined by $m$
independent equations of the form:
\def\theequation{5}\begin{equation}
\sum_{j=1}^r\,\alpha_{kj}\,e_j
=
0
\ \ \ \ \ \ \ \ \ \ \ \ \ {\scriptstyle{(k\,=\,1\,\cdots\,m)}}.
\end{equation}
Assume furthermore 
that one knows from some reason that amongst the infinitesimal
transformations: $e_1\, X_1f + \cdots + e_r\, X_rf$, no infinitesimal
transformation of the form: $e_1\, X_1f + \cdots + e_m X_mf$ is
contained in this family. Then at first, it can be deduced that the 
equations~\thetag{ 5} are solvable with respect to $e_1, \dots, e_m$, 
for if one sets: $e_{ m+1} = \cdots = e_r = 0$ in these equations,
it must follow: $e_1 = \cdots = e_m = 0$, which happens only when
the determinant: $\sum\, \pm \, \alpha_{ 11} \cdots\, \alpha_{ mm}$
does not vanish. Hence if one chooses $e_{ m+1}, \dots, e_r$ arbitrary,
though not all zero, it follows that $e_1, \dots, e_m$ receive
determined values and therefore, the family contains $r - m$ mutually
independent infinitesimal transformations of the form:
\[
X_{m+j}f
+
e_{1j}\,X_1f
+\cdots+
e_{mj}\,X_mf
\ \ \ \ \ \ \ \ \ \ \ \ \ {\scriptstyle{(j\,=\,1\,\cdots\,r\,-\,m)}}.
\]
Thus, the following holds:

\def\theproposition{3}\begin{proposition}
\label{Satz-3-S-562}
If, amongst the infinitesimal transformations: $\sum\, e_k\, X_kf$ of
the $r$-term group: $X_1f, \dots, X_rf$, a family is sorted by means
of $m$ independent linear equations:
\[
\sum_{j=1}^r\,\alpha_{kj}\,e_j
=
0
\ \ \ \ \ \ \ \ \ \ \ \ \ {\scriptstyle{(k\,=\,1\,\cdots\,m)}}
\]
which embraces no infinitesimal transformation of the form: $e_1\, 
X_1f + \cdots + e_m\, X_mf$, then this family contains $r - m$ 
infinitesimal transformations of the shape:
\[
X_{m+j}
+
\sum_{\nu=1}^m\,e_{j\nu}\,X_\nu f
\ \ \ \ \ \ \ \ \ \ \ \ \ {\scriptstyle{(j\,=\,1\,\cdots\,r\,-\,m)}}.
\]
\end{proposition}

\sectionengellie{\S\,\,\,138.}

After these preparations, we turn ourselves specially to the general
projective group. We denote for its number $n \, ( n+2)$ of parameters
shortly by $N$ and we seek at first all subgroups with more than $N -
n$ parameters, hence some $G_{ N - m}$ for which $m < n$. Naturally,
this way of putting the question is meaningful only when the number
$n$ is larger than $1$.

According to an observation made earlier on
(cf. Chap.~\ref{kapitel-12}, Proposition~7,
p.~\pageref{Satz-7-S-211}), the sought $G_{ N-m}$ has at least $n - m$
independent infinitesimal transformations in common with the $n$-term
group $p_1, \dots, p_n$. So $G_{ N-m}$ contains in any case $n - m$
independent infinitesimal translations. If it contains no more than $n
- m$ translations, then thanks to Proposition~1, we can assume that
$p_{ m+1}, \dots, p_n$ are these translations, while no translation of
the form: $e_1\, p_1 + \cdots + e_m\, p_m$ is extant. From the
Proposition~3, it then follows that there appears a transformation:
\[
x_{m+1}\,p_1
+e_1\,p_1+\cdots+e_m\,p_m,
\]
but by combination with $p_{ m+1}$, it would give $p_1$, which would
be a contradiction. Hence in our $G_{ N - m}$ there are surely more
than $n - m$, say $n - q$ ($q < m$) infinitesimal translations, and we
want to assume that these are: $p_{ q+1}, \dots, p_m, \dots, p_n$; by
contrast, when $q >0$, there appear no translations of the form:
$e_1\, p_1 + \cdots + e_q \, p_q$. Now in any case
(cf. Chap.~\ref{kapitel-12}, Proposition~7,
p.~\pageref{Satz-7-S-211}), there is in the $G_{ N - m}$
\emphasis{one} transformation:
\[
\lambda_n\,x_n\,p_1
+\cdots+
\lambda_{q+1}\,x_{q+1}\,p_1
+e_q\,p_q
+\cdots+
e_1\,p_1
\]
in which, according to what precedes, the $\lambda$ cannot vanish all.
By combination with one of the translations $p_{ q+1}, \dots, p_n$, 
we therefore obtain at least once $p_1$, but this was excluded.
Hence the number $q$ cannot be bigger than zero, so $q = 0$ and the
sought $G_{ N-m}$ therefore comprises, when $m < n$, all translations.

Thanks to completely analogous considerations, one realizes that the
$G_{ N - m}$ ($m < n$) must contain all transformations $P_i$. These
considerations coincide even literally with those undertaken just now,
when one replaces the $p_i$, $x_ip_k$, $P_i$ by $P_i$, $-\, x_kp_i$,
$p_i$, respectively, and when one relates to Proposition~2,
p.~\pageref{Satz-2-S-559}.

Thus, our $G_{ N - m}$ contains all $p_i$ and all $P_i$
simultaneously, but then as was already shown earlier on some occasion
(on page~\pageref{S-560}), it contains yet also all $x_i p_k$ and is
hence identical to the general projective group itself. Consequently:

\def\thetheorem{100}\begin{theorem}
The general projective group of the manifold $x_1, \dots, x_n$
contains no subgroup with more than $n \, ( n+2) - n = 
n \, ( n+1)$ parameters.
\end{theorem}

\sectionengellie{\S\,\,\,139.}

At present, the question is to determine all subgroups contained in
the general projective group having $N - n = n \, (n+1)$ parameters.
In order to be able to settle completely this problem, we must treat
individually a series of various possibilities.

At first, we seek all $n \, (n+1)$-term subgroups which contain no
infinitesimal translation $\sum\, e_k\, p_k$. According to
Proposition~3, p.~\pageref{Satz-3-S-562}, there surely exists a
transformation:
\[
U
=
\sum_{i=1}^n\,x_i\,p_i
-
\sum_{i=1}^n\alpha_i\,p_i
=
\sum_{i=1}^n\,
(x_i-\alpha_i)\,p_i
\] 
and for every value of $i$ and $k$, there is in the same way a 
transformation of the form:
\[
T
=
(x_i-\alpha_i)\,p_k
+
\sum_{j=1}^n\,\beta_{ikj}\,p_j.
\]
By combination of the two infinitesimal transformations $U$ and $T$,
we obtain the expression: $\leftbracket U, \, T \rightbracket = 
-\, \sum_j\, \beta_{ ikj}\, p_j$, and because our group contains
no infinitesimal transformation of this form, all $\beta_{ ikj}$
must vanish.

Lastly, there are yet $n$ infinitesimal transformations:
\[
P_i
+
\sum_{k=1}^n\,\gamma_{ik}\,p_k,
\]
or, what amounts to the same:
\[
P_i'
=
(x_i-\alpha_i)\,
\sum_{j=1}^n\,
(x_j-\alpha_j)\,p_j
+
\sum_{k=1}^n\,\delta_{ik}\,p_k.
\]
If we make the combination of $P_i'$ with $\sum\, (x_k - \alpha_k)\,
p_k = U$, we obtain:
\[
\leftbracket
U,\,P_i'
\rightbracket
=
(x_i-\alpha_i)\,
\sum_{j=1}^n\,(x_j-\alpha_j)\,p_j
-
\sum_{k=1}^n\,\delta_{ik}\,p_k,
\]
so that all $\delta_{ ik}$ vanish. 

Thus, in the sought $n \, ( n+1)$-term subgroups, the following
$n \, (n+1)$ independent infinitesimal transformations must appear:
\def\theequation{6}\begin{equation}
(x_i-\alpha_i)\,p_k,
\ \ \ \ \ \ \ 
(x_i-\alpha_i)\,
\sum_{j=1}^n\,(x_j-\alpha_j)\,p_j
\ \ \ \ \ \ \ \ \ \ \ \ \ 
{\scriptstyle{(i,\,\,k\,=\,1\,\cdots\,n)}}.
\end{equation}
Thanks to pairwise combinations, one easily convinces oneself that
these infinitesimal transformations effectively generate an $n\, (
n+1)$-term group. Besides, this also follows from the fact that all
the infinitesimal transformations~\thetag{ 6} leave invariant the
point $x_i = \alpha_i$ lying in the domain of the finite. Indeed,
they are mutually independent and their number equals $n \, ( n+1)$,
that is to say, exactly equal to the number of independent
infinitesimal transformations there are in the $n\, (n+2)$-term
projective group that leave invariant the point $x_i =
\alpha_i$. According to Proposition~2, p.~\pageref{Satz-2-S-205}, the
infinitesimal transformations~\thetag{ 6} therefore generate an $n \,
(n+1)$-term group.

As a result, every $n \, (n+1)$-term projective group of the $R_n$ in
which no infinitesimal translation $\sum\, e_k\, p_k$ appears consists
of all projective transformations that fix a point located in the
domain of the finite.

If, in the above computations, we would actually have written $P_i$, 
$-\, x_ip_k$, $p_i$ in place of $p_i$, $x_k\, p_i$, $P_i$, 
respectively, then we would have found all $n\, (n+1)$-term
subgroups which contain no infinitesimal transformation 
$\sum\, e_k\, P_k$. We can therefore make exactly the same 
exchange in the expressions~\thetag{ 6} and we obtain in this way:
\[
x_k\,p_i
+
\alpha_i\,P_k,
\ \ \ \ \ \ \ 
p_i
+
\alpha_i\,\sum_{j=1}^n\,
x_j\,p_j
+
\sum_{j=1}^n\,
\alpha_j\,(x_j\,p_i+\alpha_i\,P_j).
\]
Here, we may clearly take away the the term $\sum\, \alpha_j\, 
(x_j\, p_i + \alpha_i \, P_j)$ and we therefore obtain the general
form of the $n \, ( n+1)$-term subgroups which contain
no transformation $\sum\, e_k\, P_k$ as follows:
\def\theequation{7}\begin{equation}
p_i
+
\alpha_i\,
\sum_{j=1}^n\,x_j\,p_j,
\ \ \ \ \ \ \ \ \ \
x_i\,p_k
+
\alpha_k\,P_i.
\end{equation}
The fact that these infinitesimal transformation generate a group
follows from their derivation, but naturally, one could also 
corroborate this directly.

If all $\alpha_i$ vanish, then we have the general linear group
already discussed earlier on which leaves invariant the plane
$M_{ n-1}$ located at infinity. Hence we are very close to presume
that in the general case where not all $\alpha_i$ are equal to zero,
there exists in the same way a plane $M_{ n-1}\colon$
$\lambda_1 \, x_1 + \cdots + \lambda_n \, x_n + \lambda = 0$
which admits all infinitesimal transformations~\thetag{ 7}.

By execution of the infinitesimal transformations: $p_i + \alpha_i\, 
\sum\, x_j\, p_j$, we obtain the following conditions for the 
$\lambda_i$:
\[
\lambda_i
+
\alpha_i\,\sum_{j=1}^n\,\lambda_j\,x_j
=
0
=
\lambda_i
-
\alpha_i\,\lambda. 
\]
The quantity $\lambda$ can hence in any case not vanish and may be
set equal to $1$; so if there actually is an invariant plane 
$M_{ n-1}$, this plane can only have the form: $\alpha_1\, x_1 + 
\cdots + \alpha_n\, x_n + 1 = 0$. In fact, this latter plane yet
admits also the infinitesimal transformations: $x_i \, p_k + \alpha_k
\, P_i$. 

Every subgroup~\thetag{ 7} therefore leaves invariant a plane $M_{ n
- 1}$ of the $R_n$ and is in addition the most general projective
group of the $R_n$ which leaves at rest the plane $M_{ n-1}$ in 
question. Every further infinitesimal projective transformation which
would leave invariant the plane: $M_{ n-1}\colon$
$\alpha_1\, x_1 + \cdots + \alpha_n\, x_n + 1 = 0$ could indeed
be given the form:
\[
\sum_{k=1}^n\,e_k\,P_k
=
\sum_{k=1}^n\,e_k\,x_k\,
\sum_{i=1}^n\,x_i\,p_i.
\] 
But if one executes this infinitesimal transformation on the 
$M_{ n-1}$, it comes:
\[
\sum_{k=1}^n\,e_k\,x_k\,
\sum_{i=1}^n\,\alpha_i\,x_i
=
0
=
-\,\sum_{k=1}^n\,e_k\,x_k,
\]
whence $e_1 = \cdots = e_n = 0$. 

Thus, if an $n\, (n+1)$-term projective group of the $R_n$ contains no
infinitesimal transformation $\sum\, e_k\, P_k$, it consists of all
projective transformations which leave invariant a plane $M_{ n-1}$
not passing through the origin of coordinates.

From the previous results, we can yet derive, thanks to a simple
transformation, a few other results which will be useful to us in the
future. Indeed, if we transfer to infinity the former origin of
coordinates by means of the collineation:
\def\theequation{7'}\begin{equation}
x_1
=
\frac{1}{x_1'},
\ \
x_2
=
\frac{x_2'}{x_1'},
\,\,\,\dots,\,\,\,
x_n
=
\frac{x_n'}{x_1'},
\end{equation}
we then obtain:
\[
p_1'
=
\sum_{i=1}^n\,p_i\,
\frac{\partial x_i}{\partial x_1'}
=
-\,
\frac{1}{{x_1'}^2}\,
\bigg(
p_1
+
\sum_{i=2}^n\,x_i'\,p_i
\bigg),
\]
and hence:
\[
\aligned
p_1'
&
=
-\,x_1\,\sum_{i=1}^n\,x_i\,p_i
=
-\,P_1,
\\
x_k'\,p_1'
&
=
-\,x_k\,\sum_{i=1}^n\,x_i\,p_i
=
-\,P_k
\ \ \ \ \ \ \ \ \ \ \ \ \ {\scriptstyle{(k\,=\,2\,\cdots\,n)}}.
\endaligned
\]
In the same way, we obtain:
\[
x_1'\,\sum_{i=1}^n\,x_i'\,p_i'
=
-\,p_1,
\ \ \ \ \ \ \
x_1'\,p_k'
=
-\,p_k
\ \ \ \ \ \ \ \ \ \ \ \ \ {\scriptstyle{(k\,=\,2\,\cdots\,n)}},
\]
as every infinitesimal projective transformation is actually
transferred to a transformation of the same kind after the
introduction of the $x'$ (cf. Chap.~\ref{fundamental-differential},
Proposition~4, p.~\pageref{Satz-4-S-81}).

From this, we see: our collineation~\thetag{ 7'} converts every $n\,
(n+1)$-term projective group in which no infinitesimal transformation
$\sum\, e_k\, P_k$ appears into a projective group which contains no
transformation $e_1\, p_1 + e_2\, x_2 \, p_1 + \cdots + e_n\, x_n \,
p_1$. In the same way, every projective group free of all $\sum\,
e_k\, p_k$ is transferred to a projective group which contains no
transformation: $e_1\, P_1 + e_2\, x_1\, p_2 + \cdots + e_n\, x_1 \,
p_n$.

Thus, if in an $n\, (n+1)$-term projective group, there are no
infinitesimal transformations: $e_1\, p_1 + e_2\, x_2\, p_1 + \cdots +
e_n\, x_n\, p_1$, then this group consists of all projective
transformations which leave invariant a certain plane $M_{ n-1}$. By
contrast, if in the group there are no infinitesimal transformations:
$e_1\, P_1 + e_2\, x_1\, p_2 + \cdots + e_n\, x_1\, p_n$, then this
group consists of all projective transformations which leave invariant
a certain point.

The general projective group possesses the property of being able to
transfer \emphasis{every} point of the $R_n$ to \emphasis{every} other
point, and \emphasis{every} plane $M_{ n-1}$ to \emphasis{every} other
plane. From this, it results that every $n \, (n+1)$-term projective
group of the $R_n$ which leaves invariant a point is conjugate, inside
the general projective group of the $R_n$, to every other group of the
same sort, and likewise, it follows that every $n \, (n+1)$-term
projective group of the $R_n$ which leaves invariant a plane $M_{
n-1}$ is conjugate to every other projective group of this sort.

Finally, we can at present look up at all $n\, (n+1)$-term groups of
the $R_n$.

Since we know all groups that contain no translation, it only remains
to find the groups in which some infinitesimal translations appear. We
assume that there are exactly $q$ independent such infinitesimal
translations, say $p_n, \dots, p_{ n - q + 1}$, so that, when $n - q$
is $>0$, no translation of the form: $e_1\, p_1 + \cdots + e_{ n-q}\,
p_{ n - q}$ is extant.

Then in our group, there surely exists an infinitesimal transformation
of the form:
\def\theequation{8}\begin{equation}
\sum_{i=1}^q\,x_{n-q+i}\,
\sum_{k=1}^{n-q}\,
\lambda_{ik}\,p_k
+
e_1\,p_1
+\cdots+
e_{n-q}\,p_{n-q},
\end{equation}
only when the number $(n - q) (q + 1)$ of terms contained in this
infinitesimal transformation is larger than $n$.

This is the case only when $(n-q) (q+1) - n = q\, (n - q - 1)$ is
larger than zero, from which it follows that we must temporarily
disregard the cases $q = n - 1$ and $q = n$. But if we assume that $q$
is smaller than $n - 1$ and if we make a combination of the
transformation~\thetag{ 8}\,---\,in which obviously not all $\lambda_{
ik}$ vanish\,---\,with each one of the extant translations $p_{
n-q+1}, \dots, p_n$, we then obtain in all circumstances a not
identically vanishing transformation of the form $\mu_1 \, p_1 +
\cdots + \mu_{ n-q}\, p_{n-q}$, and this is a contradiction. Thus, the
number $q$ cannot be smaller than $n - 1$.

Therefore, if an $n\, (n+1)$-term projective group contains
\emphasis{one} infinitesimal translation $\sum\, e_k\, p_k$, it
contains at least $n - 1$ independent translations.

In this result we can again, as so often before, replace the $p_i$ by
the $P_i$ and find that there always exist $n - 1$ independent
transformations $\sum\, e_k\, P_k$ in every group of the said
sort as soon as there exists only a single transformation of this
form. 

At present, we seek all $n\, (n+1)$-term projective groups with 
exactly $n - 1$ independent infinitesimal translations $\sum\, e_k\,
p_k$, say with $p_2, \dots, p_n$. No such group can contain a 
transformation of the form:
\[
e_1\,p_1
+
e_2\,x_2\,p_1
+\cdots+
e_n\,x_n\,p_1,
\]
because by combination with $p_2, \dots, p_n$, we would obtain $p_1$,
what is excluded. According to what has been seen above, all these
groups belong to the category of the $n\, (n+1)$-term projective
groups which leave invariant a plane $M_{ n-1}$. Correspondingly, 
it results that every $n\, (n+1)$-term group which contains exactly
$n - 1$ independent transformations $\sum\, e_k\, P_k$ leaves
invariant a point. 

It still remains to determine the $n \, (n+1)$-term projective groups
which contain all the $n$ translations $p_1, \dots, p_n$. Besides,
the $P_i$ cannot yet appear all, because otherwise, the group would
coincide with the general projective group itself. 

Therefore, according to what has been said above, there are only
two possibilities: either there is absolutely no transformation 
$\sum\, e_k\, P_k$, or there are $n - 1$ independent such 
transformations. Both cases are already settled above. 

As a result, our study is brought to a conclusion. The result is the
following:

\renewcommand{\thefootnote}{\fnsymbol{footnote}}
\def\thetheorem{101}\begin{theorem}
The largest subgroups of the general projective group of an $n$-times
extended manifold contain $n\, (n+1)$ parameters. Each such subgroup
consists of all projective transformations which leave invariant
either a plane $M_{ n-1}$, or a point. In the first case, the
subgroup is conjugate inside the general projective group to the
general linear group $p_i$, $x_ip_k$, and in the second
case, to the group $x_kp_i$, $P_k$.\footnote[1]{\,
\name{Lie}, Math. Ann. Vol. XXV, p.~130.
}
\end{theorem}
\renewcommand{\thefootnote}{\arabic{footnote}}

Because the groups of the one category come from the groups of the
other category through the exchange of $p_i$, $x_ip_k$, $P_i$ with
$P_i$, $-\, x_k p_i$, $p_i$, all $n\, (n+1)$-term subgroups
of the general projective group are holoedrically isomorphic 
to each other. From what precedes, it follows in addition the

\def\theproposition{4}\begin{proposition}
The general projective group of the space $x_1, \dots, x_n$
can be related to itself in a holoedrically isomorphic way
so that the largest subgroups which leave invariant one point
correspond each time to the largest subgroups which leave
invariant a plane $M_{ n-1}$.
\end{proposition}

If $n = 1$, the difference between point and plane $M_{ n-1}$
disappears; the general three-term
projective group of the once-extended
manifold therefore contains only {\em one} category of two-term
subgroups, and all of these are
conjugate in the three-term projective
group.

According to the previous developments, the general projective
group of an $n$-times extended space contains $n\, (n+1)$
independent infinitesimal transformations which leave
at rest a given point, and to be precise, these infinitesimal
transformations
generate a subgroup which is contained in no larger subgroup.
From this, it follows (Chap.~\ref{kapitel-24}, 
Theorem~91, p.~\pageref{Theorem-91-S-521})
that \emphasis{the general projective group is primitive and
all the more asystatic.}

\sectionengellie{\S\,\,\,140.}

We derive here a few general considerations about the
determination of all subgroups of the general linear group
$p_i$, $x_ip_k$ 
${\scriptstyle{(i,\,\,k\,=\,1\,\cdots\,n)}}$
of the $R_n$.

The general infinitesimal transformation of a linear group of the
$R_n$ can be written:
\def\theequation{9}\begin{equation}
\left\{
\aligned
&
\sum_{i,\,\,k}^{i\,\gtrless\,k}\,
a_{ik}\,x_i\,p_k
+
\sum_{i=1}^{n-1}\,
b_i\,(x_ip_i-x_np_n)
\\
&
+
c\,\sum_{k=1}^n\,x_k\,p_k
+
\sum_{k=1}^n\,d_k\,p_k.
\endaligned\right.
\end{equation}
If one makes combination between two infinitesimal transformations
of this shape, one obtains a transformation:
\[
\sum_{i,\,\,k}^{i\,\gtrless\,k}\,
A_{ik}\,x_i\,p_k
+
\sum_{i=1}^{n-1}\,
B_i\,(x_ip_i-x_np_n)
+
C\,
\sum_{k=1}^n\,x_k\,p_k
+
c\,\sum_{k=1}^n\,D_k\,p_k,
\]
in which the $A_{ ik}$, $B_i$ and $C$ ($C = 0$) depend only 
of $a_{ ik}$, $b_i$ and $c$. Consequently, the reduced infinitesimal
transformation:
\def\theequation{10}\begin{equation}
\sum_{i,\,\,k}^{i\,\gtrless\,k}\,
a_{ik}\,x_i\,p_k
+
\sum_{i=1}^{n-1}\,b_i\,
(x_ip_i-x_np_n)
+
c\,\sum_{k=1}^n\,x_k\,p_k
\end{equation}
is in turn the general infinitesimal transformation of a linear
homogeneous group in $x_1, \dots, x_n$.

From this, we realize that the problem of determining all subgroups of
the general linear group decomposes in two problems which must be
settled one after the other. At first, all subgroups of the general
linear homogeneous group $x_i\, p_k$
${\scriptstyle{(i,\,\,k\,=\,1\,\cdots\,n)}}$ have to be sought;
afterwards, to the infinitesimal transformations of each of the found
groups, one must add terms $\sum\, \beta_k\, p_k$ in the most general
way so that one again obtains a group. Thus, if $X_1f, \dots, X_rf$
is one of the found linear homogeneous groups, one has
to determine all groups of the form:
\[
\aligned
X_kf
&
+
\sum_{i=1}^n\,\alpha_{ki}\,p_i,
\ \ \ \ \ \ \ \ \ \ \ \ \ \ \
\sum_{i=1}^n\,\beta_{\mu i}\,p_i
\\
&
\ \ \ \ \
{\scriptstyle{(k\,=\,1\,\cdots\,r\,;\,\,\,
\mu\,=\,1\,\cdots\,m,\,\,\,\,m\,\leqslant\,n)}}.
\endaligned
\]
We do not want here to say anything more about the further treatment
of these two reduced problems; rather, we refer to the third volume
where the detailed studies about the projective groups of the plane
and of the thrice-extended space will be brought.

By contrast, we do not want to neglect to draw attention to the
geometrical signification that the decomposition just said of the
problem in question has.

To this end, we imagine that the group~\thetag{ 9} is prolonged, by
regarding, as in Chap.~\ref{kapitel-25}, p.~\pageref{S-524-sq} sq.,
the $x$ as functions of an auxiliary variable $t$ and by taking with
them the differential quotients $\D\, x_i / \D\, t = x_i'$. In the
process, we obtain the group:
\[
\aligned
&
\sum\,a_{ik}\,x_i\,p_k
+
\sum\,b_i\,(x_ip_i-x_np_n)
+
c\,\sum\,x_k\,p_k
+
\sum\,d_k\,p_k
\\
&
+
\sum\,a_{ik}\,x_i'\,p_k'
+
\sum\,b_i\,(x_i'p_i'-x_n'p_n')
+
c\,\sum\,x_k'\,p_k',
\endaligned
\]
in which the terms in the $x_i'$ determine for themselves a group,
namely precisely the group~\thetag{ 10} just found. But now, as we
have seen loc. cit., the $x_i'$ can be interpreted as homogeneous
coordinates of the directions through the point $x_1, \dots, x_n$ of
the $R_n$. The fact that the $x_i'$ are transformed for themselves by
the above group therefore means nothing else than the fact that
parallel lines are transferred to parallel lines by every linear
transformation of the $x$; the directions which a determined system of
values $x_i'$ associates to all points of the $R_n$ are indeed
parallel to each other. But every bundle \deutsch{Bündel} of parallel
lines provides a completely determined point on the plane $M_{ n-1}$
at the infinity of the $R_n$, hence \emphasis{$x_1', \dots, x_n'$ can
be virtually interpreted as homogeneous coordinates of the points on
the infinitely far plane $M_{ n-1}$} and the group:
\def\theequation{11}\begin{equation}
\sum\,a_{ik}\,x_i'\,p_k'
+
\sum\,b_i\,(x_i'p_i'-x_n'p_n')
+
c\,\sum\,x_k'\,p_k'
\end{equation}
therefore indicates how the infinitely far points of the $R_n$ are
transformed by the group~\thetag{ 9}. At the same time, it is yet to
be observed that the infinitesimal transformation $\sum\, x_k' \,
p_k'$ leaves fixed all infinitely far points, so that these points are
transformed by the last group exactly as if $c$ would be zero.

At present, we have the inner reason for the decomposition indicated
above of the problem of determining all linear groups of the
$R_n$. The groups in question are simply thought to be distributed in
classes, and in each class are reckoned all the groups for which the
group~\thetag{ 11} is the same, so that the infinitely far points of
the $R_n$ are transformed in the same way (cf. for 
this purpose Theorem~40,
p.~\pageref{Theorem-40-S-233}).

\sectionengellie{\S\,\,\,141.}

In order to give at least \emphasis{one} application of the
general considerations developed just above, we want to 
determine all linear groups of the $R_n$ which transform
the infinitely far points of the $R_n$ in the most general 
way. For all of these groups, the associated group~\thetag{ 11}
has the form:
\[
x_i'\,p_k',
\ \ \ \
x_i'\,p_i'-x_k'\,p_k'
\ \ \ \ \ \ \ \ \ \ \ \ \ 
{\scriptstyle{(i,\,\,k\,=\,1\,\cdots\,n,\,\,\,\,
i\,\gtrless\,k)}},
\]
where in certain circumstances, the transformation $x_1'\, p_1' +
\cdots + x_n'\,p_n'$ can yet occur, which leaves individually fixed
all infinitely far points. The infinitely far points are thus
transformed by an $(n^2 - 1)$-term group and to be precise, by the
general projective group of an $(n-1)$-times extended space.

Each one of the sought groups must contain $n^2 - 1$ infinitesimal
transformations out of which none can be linearly deduced which
leaves invariant all infinitely far points, hence none which
possesses the form: $\gamma\, \sum_j\, x_j\, p_j + \sum\, 
\gamma_k\, p_k$. Therefore, the group surely contains
$n^2 - 1$ infinitesimal transformations of the form:
\def\theequation{12}\begin{equation}
\left\{
\aligned
&
x_i\,p_k
+
\alpha_{ik}\,\sum_{j=1}^n\,x_j\,p_j
+
\sum_{\nu=1}^n\,\beta_{ik\nu}\,p_\nu
\ \ \ \ \ \ \ \ \ \ \ \ \ 
{\scriptstyle{(i\,\gtrless\,k)}}
\\
&
x_i\,p_i-x_n\,p_n
+
\alpha_i\,\sum_{j=1}^n\,x_j\,p_j
+
\sum_{\nu=1}^n\,\beta_{i\nu}\,p_\nu.
\endaligned\right.
\end{equation}
In addition, one or several infinitesimal transformations
of the form: $\gamma\, \sum_j\, x_j\, p_j + 
\sum_\nu \, \gamma_\nu\, p_\nu$ can yet occur.

If a group of the demanded sort contains a translation, then it
contains all translations. Indeed, if $p_1 + e_2\, p_2 + 
\cdots + e_n\, p_n$ is the translation in question, we
make a combination of it with:
\[
x_1\,p_k
+
\alpha_{1k}\,\sum_{j=1}^n\,x_j\,p_j
+
\sum_{\nu=1}^n\,\beta_{1k\nu}\,p_\nu
\ \ \ \ \ \ \ \ \ \ \ \ \ {\scriptstyle{(k\,=\,2\,\cdots\,n)}}
\]
and we obtain in this way $p_2, \dots, p_n$ and therefore all 
$p_i$. 

Consequently, we want at first to assume that all translations
appear. Then if there is yet the transformation
$\sum\, x_j\, p_j$, we have the general linear group itself. 
By contrast, if the transformation $\sum\, x_j\, p_j$ is not
extant, we obtain by combination of the transformations~\thetag{
12} in which we can set beforehand
the $\beta_{ ik\nu}$ and $\beta_{ i\nu}$ equal to zero:
\[
\bigg\leftbracket
x_i\,p_k
+
\alpha_{ik}\,\sum_{j=1}^n\,x_j\,p_j,\,\,\,
x_k\,p_i
+
\alpha_{ki}\,\sum_{j=1}^n\,x_j\,p_j
\bigg\rightbracket
=
x_i\,p_i
-
x_k\,p_k,
\]
so that all the $\alpha_i$ are zero. But in addition, it comes:
\[
\bigg\leftbracket
x_i\,p_i-x_k\,p_k,\,\,\,
x_i\,p_k+\alpha_{ik}\,\sum_{j=1}^n\,x_j\,p_j
\bigg\rightbracket
=
2\,x_i\,p_k.
\]
The concerned group is therefore the special linear group.

If in the sought group there is absolutely no translation, 
but by contrast a transformation:
\def\theequation{13}\begin{equation}
\sum_{j=1}^n\,x_j\,p_j
+
\sum_{\nu=1}^n\,\gamma_\nu\,p_\nu
=
\sum_{j=1}^n\,(x_j+\gamma_j)\,p_j,
\end{equation}
then all $\alpha_{ ik}$ and all $\alpha_i$ can be set equal to zero.
If we yet write the infinitesimal transformations~\thetag{ 12} 
in the form:
\[
\aligned
&
\ \ \ \ \ \ \ \ \ \ \ \ \ 
(x_i+\gamma_i)\,p_k
+
\sum_{\nu=1}^n\,\beta_{ik\nu}'\,p_\nu
\ \ \ \ \ \ \ \ \ \ \ \ \ 
{\scriptstyle{(i\,\gtrless\,k)}}
\\
&
(x_i+\gamma_i)\,p_i
-
(x_n+\gamma_n)\,p_n
+
\sum_{\nu=1}^n\,\beta_{i\nu}'\,p_\nu,
\endaligned
\]
we then realize immediately by combination with $\sum\, (x_i + \gamma_i
)\, p_i$ that all $\beta_{ ik\nu}'$ and $\beta_{ i\nu}'$ vanish and
we thus find the group:
\[
(x_i+\gamma_i)\,p_k
\ \ \ \ \ \ \ \ \ \ \ \ \ 
{\scriptstyle{(i,\,\,k\,=\,1\,\cdots\,n)}}.
\]

Lastly, if there also occurs no transformation of the form~\thetag{ 
13}, then one obtains at first by combination from~\thetag{ 12}
that all $\alpha_{ ik}$ and $\alpha_i$ are equal to zero. 
Furthermore, it comes:
\[
\bigg\leftbracket
x_i\,p_k
+
\sum_{\nu=1}^n\,\beta_{ik\nu}\,p_\nu,\,\,\,
x_k\,p_i
+
\sum_{\nu=1}^n\,\beta_{ki\nu}\,p_\nu
\bigg\rightbracket
=
x_i\,p_i
-
x_k\,p_k
+
\beta_{ikk}\,p_i
-
\beta_{kii}\,p_k,
\]
and in addition:
\[
\aligned
\bigg\leftbracket
(x_i+\beta_{ikk})\,p_i
&
-
(x_k+\beta_{kii})\,p_k,\,\,\,
x_i\,p_k
+
\sum_{\nu=1}^n\,\beta_{ik\nu}\,p_\nu
\bigg\rightbracket
\\
&
=
2\,x_i\,p_k
+
2\,\beta_{ikk}\,p_k
-
\beta_{iki}\,p_i,
\endaligned
\]
so that all $\beta_{ ik\nu}$ vanish with the exception of the $\beta_{
ikk}$ and $\beta_{ kii}$, by means of which also the $\beta_{ i\nu}$
can be expressed.

Now, if $n > 2$, the $\beta_{ ikk}$ could vary with $k$; however, this
is not the case. Indeed, if we replace $i$ and $k$ firstly by $k$, $j$
and secondly by $j$, $i$ in:
\[
(x_i+\beta_{ikk})\,p_i
-
(x_k+\beta_{kii})\,p_k
\]
and if we add together the three obtained infinitesimal
transformations, then the sum must vanish, since no translation should
occur. So we obtain: $\beta_{ ikk} = \beta_{ ijj}$ and so on. Thus, if
we write shortly $\beta_i$ for $\beta_{ ikk}$, we have the group:
\[
(x_i+\beta_i)\,p_k,
\ \ \ \ \ \ \
(x_i+\beta_i)\,p_i
-
(x_k+\beta_k)\,p_k
\ \ \ \ \ \ \ \ \ \ \ \ \ 
{\scriptstyle{(i\,\gtrless\,k)}}.
\]

With that, all cases are settled. If, in the two latter forms of
groups \deutsch{Gruppenformen}, we yet introduce $x_i + \gamma_i$ and
$x_i + \beta_i$ as new $x_i$, respectively, we may recapitulate
our result as follows:

\renewcommand{\thefootnote}{\fnsymbol{footnote}}
\def\thetheorem{102}\begin{theorem}
The general linear group in $n$ variables contains only three
different sorts of subgroups which transform the points of the
infinitely far plane in an $(n^2 - 1)$-term way, as does the general
linear group itself: firstly, the special linear group and secondly
all groups that are conjugate to the two homogeneous
groups:\footnote[1]{\,
\name{Lie}, Archiv for Math. og Nat. Vol. IX, p.~103 and 104,
Christiania 1884.
}
\[
x_i\,p_k\,;
\ \ \ \ \ \ \ \ \ \
x_i\,p_k\ \ \
x_i\,p_k-x_k\,p_k
\ \ \ \ \ \ \ \ \ \ \ \ \ 
{\scriptstyle{(i\,\gtrless\,k)}}.
\]
\end{theorem}
\renewcommand{\thefootnote}{\arabic{footnote}}

Here, we thus have a characteristic property which is common to
all these groups already known to us. By contrast, the
distinctive marks of the four groups in question are briefly the
following: The general linear group leaves invariant the 
ratios of all volumes; the special linear group leaves invariant
all volumes. The general and the special linear homogeneous
\label{S-574}
groups differ from the general and the special linear groups, 
respectively, in that they yet leave invariant the point $x_i = 0$.

\sectionengellie{\S\,\,\,142.}

In the beginning of the previous paragraph, we have seen that the
determination of all linear groups of the $R_n$ is essentially
produced as soon as all linear homogeneous projective groups
of the $R_n$ are determined. There is no special difficulty
to settle this last problem when one knows all projective groups
of the $R_{ n-1}$. At present, we yet want to show that.

As we know, the general linear homogeneous group of the $R_n$: 
$x_i\, p_k$ ${\scriptstyle{(i,\,\,k\,=\,1\,\cdots\,n)}}$
contains an invariant subgroup with $n^2 - 1$ parameters, namely 
the special linear homogeneous group:
\[
x_i\,p_k,
\ \ \ \ \ \
x_i\,p_i-x_k\,p_k
\ \ \ \ \ \ \ \ \ \ \ \ \ 
{\scriptstyle{(i,\,\,k\,=\,1\,\cdots\,n,\,\,\,
i\,\gtrless\,k)}}.
\]

This last group is equally composed with the general projective group
of the $R_{ n-1}$, so its subgroups can immediately be written down
when all projective groups of the $R_{ n-1}$ are known (cf. about 
that the next chapter). After that, one finds the subgroups of the 
group $x_i\, p_k$ thanks to the following considerations:

An $r$-term subgroup $G_r$ of the group $x_i\, p_k$ is either
contained at the same time in the group $x_i\, p_k$, $x_i\, p_i - 
x_k\, p_k$ ${\scriptstyle{(i\,\gtrless\,k)}}$ or it is not. In the 
first case it would be already known, and in the second case, 
according to Proposition~7, p.~\pageref{Satz-7-S-211},
it would have a $G_{ r-1}$ in common with the special linear
homogeneous group. Hence, in order to find all linear homogeneous
groups $G_r$ of this sort, we need only to add, to every $G_{ r-1}
\colon$ $X_1f, \dots, X_{ r-1}f$ of the group $x_i\, p_k$, 
$x_i\, p_i - x_k\, p_k$
${\scriptstyle{(i\,\gtrless\,k)}}$, an infinitesimal transformation
of the form:
\[
{\sf Y}f
=
\sum_{i=1}^n\,x_i\,p_i
+
\sum_{k,\,\,j}^{1\cdots\,n}\,
\alpha_{kj}\,x_k\,p_j
\ \ \ \ \ \ \ \ \ \ \ \ \ 
{\scriptstyle{\left(\sum_k\,\,\alpha_{kk}\,=\,0\right)}},
\]
and to determine, by combination with $X_1f, \dots, X_{ r-1}f$, 
all values of the $\alpha_{ kj}$ which produce a group. Here, 
it is to be observed that the $(r-1)$-term group $X_1f, \dots, 
X_{ r-1}f$ must obviously be invariant in the sought $r$-term group;
we therefore find, so to say, the most general values of the 
$\alpha_{ kj}$ when we seek the most general linear homogeneous
infinitesimal transformation ${\sf Y}f$ which leaves invariant the
given $(r-1)$-term group. Besides, we always find \emphasis{one}
$r$-term group, namely when we choose all $\alpha_{ kj}$ equal to
zero.

By taking into consideration the infinitesimal transformations $X_1f,
\dots, X_{ r-1}f$, one realizes that $r - 1$ of the $\alpha_{ kj}$ can
be made equal to zero. Hence the smaller the number $r$ is, the more
some constants must be determined. But for small values of $r$, the
following method is otherwise often more convenient:

Indeed, the $r$-term subgroups of the group $x_i\, p_k$ can yet in
another way be distributed in two categories; firstly, in the
category of subgroups which contain the transformation $\sum\, x_i\,
p_i$\,---\,they can, under the assumptions made, be written
down instantly\,---\, and secondly, in the category of subgroups
which do not contain the transformation $x_1\, p_1 + \cdots + x_n\, 
p_n$. The $r$ infinitesimal transformations of a group from the
latter category must have the form:
\def\theequation{14}\begin{equation}
X_kf
+
\alpha_k\,\sum_{i=1}^n\,x_i\,p_i,
\end{equation}
where the $X_kf$ represent infinitesimal transformations of the
special linear homogeneous group. Because the transformation $\sum\,
x_i\, p_i$ is interchangeable with all $X_kf$, for the execution of
the bracket operation \deutsch{Klammeroperation}, it is completely
indifferent whether the $\alpha_i$ vanish or not; hence $X_1f, \dots,
X_rf$ must themselves generate a group and to be precise, an $r$-term
subgroup of the group $x_i\, p_k$, $x_i\, p_i - x_k\, p_k$
${\scriptstyle{(i\,\gtrless\,k)}}$, one of those which we assume as
known. Therefore, it only remains to determine 
the $\alpha_k$ in the most general way so that the infinitesimal
transformations~\thetag{ 14} generate a group. In certain circumstances,
one can realize easily that certain of the $\alpha_k$ must vanish; 
indeed, if there is an equation of the form:
\[
\leftbracket
X_i,\,X_k
\rightbracket
=
X_jf,
\]
then $\alpha_j$ must necessarily be zero.

If all $\leftbracket X_i, \, X_k \rightbracket$ generate a $\rho$-term
group (Chap.~\ref{kapitel-15}, Proposition~6,
p.~\pageref{Satz-6-S-261}), then we can assume that the infinitesimal
transformations are chosen so that all $\leftbracket X_i, \, X_k
\rightbracket$ can be linearly deduced from $X_1f, \dots, X_\rho f$.
Then we have $\alpha_1 = \cdots = \alpha_\rho = 0$, while all other
$\alpha_i$ can be different from zero. It must be specially studied
in every individual case whether the different values of these
parameters produce different \emphasis{types} of linear homogeneous
groups, or in other words, whether the concerned parameters are
essential. The settlement of this problem for $n = 2$ and for
$n = 3$ appears in the third volume. 

\sectionengellie{\S\,\,\,143.}

If an arbitrary linear homogeneous group $G_r$ is presented which
is not contained in the special linear homogeneous group, then
as we have already observed above, this $G_r$ comprises an 
\emphasis{invariant} $(r-1)$-term subgroup whose infinitesimal 
transformations are characterized by the fact that they possess the
form:
\[
\sum_{i,\,\,k}^{i\,\gtrless\,k}\,
a_{ik}\,x_i\,p_k
+
\sum_{i=1}^{n-1}\,a_i\,
(x_i\,p_i-x_n\,p_n).
\]

If we now apply this observation to the adjoint associated to an 
\emphasis{arbitrary} $r$-term group: $X_1f, \dots, X_rf$:
\[
E_\nu f
=
\sum_{k,\,\,s}^{1\cdots\,r}\,
c_{k\nu s}\,e_k\,
\frac{\partial f}{\partial e_s}
\ \ \ \ \ \ \ \ \ \ \ \ \ {\scriptstyle{(\nu\,=\,1\,\cdots\,r)}},
\]
then we recognize immediately that this adjoint
group, when the $r$ sums $\sum_k\, 
c_{ k\nu k}$ do not all vanish, contains an invariant subgroup whose
infinitesimal transformations: $\lambda_1\, E_1f + \cdots + 
\lambda_r \, E_rf$ are defined by means of the condition-equations:
\[
\sum_{\nu=1}^r\,\lambda_\nu\,
\sum_{k=1}^r\,c_{k\nu k}
=
0.
\]

Lastly, if we remember that every group is isomorphic with its adjoint
group, we obtain the

\renewcommand{\thefootnote}{\fnsymbol{footnote}}
\def\theproposition{5}\begin{proposition}
If $r$ independent infinitesimal transformations $X_1f, \dots, X_rf$
of an $r$-term group stand pairwise in the relationships: 
$\leftbracket X_i,\, X_k \rightbracket = \sum_s\, c_{ iks}\, X_sf$
and if at least one of the $r$ sums $\sum_k\, c_{\nu kk}$ is different
from zero, then all infinitesimal transformations $\lambda_1\, X_1f
+ \cdots + \lambda_r\, X_rf$ which satisfy the condition:
\[
\sum_{\nu=1}^r\,\lambda_\nu\,
\sum_{k=1}^r\,c_{\nu kk}
=
0
\]
generate an invariant $(r - 1)$-term subgroup.\footnote[1]{\,
\name{Lie}, Archiv for Math., Vol. IX, p.~89, Christiania 1884; 
Fortschritte der Mathematik, Vol. XVI, p.~325.
}
\end{proposition}
\renewcommand{\thefootnote}{\arabic{footnote}}

\smallercharacters{

For a thorough study of the general projective group, one must
naturally devote a special attention to its adjoint group $\sum\, 
e_k\, E_kf$ and to the associated invariant systems of equations
in the $e_k$ as well. Here only two brief, but important observations.

If, as usual, one interprets the $e_k$ as homogeneous point
coordinates of a space with $n^2 + 2\, n - 1$ dimensions, then amongst
all invariant manifolds of this space, there is a determined one whose
dimension number possesses the smallest value. This important manifold
consists of all points $e_k$ which represent either translations or
transformations conjugate to translations. This manifold is contained
in no even \deutsch{eben} manifold of the space in question. Besides,
the known classification of all projective transformations naturally
provides without effort \emphasis{all} invariant manifolds of the
space $e_1, \dots, e_r$.

For every infinitesimal projective transformation which is conjugate
to a translation, all points of a plane $M_{ n-1}$ in the space $x_k$
remain invariant and simultaneously, all planes $M_{ n-1}$ which pass
through a certain point of this $M_{ n-1}$. Every transformation of
this sort is completely determined by means of the firstly said plane
$M_{ n-1}$ and by the distinguished invariant point of this plane.

If $n = 2$, then as we can say, every projective transformation of the
plane $x_1$, $x_2$ which is conjugate to an infinitesimal translation
is completely represented by means of a \emphasis{line
element}. Accordingly, in thrice-extended space $x_1$, $x_2$, $x_3$,
every projective transformation conjugate to an infinitesimal 
translation is represented by a \emphasis{surface element}, 
and so on. These observations will be exploited in the third
volume. 

}

\linestop


\chapter{Linear Homogeneous Groups}
\label{kapitel-27}
\chaptermark{Linear Homogeneous Groups}

\setcounter{footnote}{0}

\abstract*{??}

In the previous chapter, p.~\pageref{S-574}, we have characterized the
general linear homogeneous group in $n$ variables $x_1, \dots, x_n$ as
the most general projective group of the $n$-times extended space
$x_1, \dots, x_n$, or shortly $R_n$, that leaves invariant the
infinitely far plane $M_{ n-1}$ and simultaneously the point $x_1 =
\cdots = x_n = 0$. This group receives another meaning when one
interprets $x_1, \dots, x_n$ as homogeneous coordinates in an
$(n-1)$-times extended space $R_{ n-1}$. In the present chapter, this
interpretation shall be led at the foundation.

\sectionengellie{\S\,\,\,144.}

We imagine that the transition from the ordinary Cartesian coordinates
$y_1, \dots, y_{ n-1}$ of the $(n-1)$-times extended space
$R_{ n-1}$ to the homogeneous coordinates $x_1, \dots, x_n$ of the
same space is procured by means of the equations:
\[
y_k
=
\frac{x_k}{x_n}
\ \ \ \ \ \ \ \ \ \ \ \ \ 
{\scriptstyle{(k\,=\,1\,\cdots\,n\,-\,1)}}.
\]

To the $n^2$-term general linear homogeneous group:
\def\theequation{1}\begin{equation}
x_i'
=
\sum_{k=1}^n\,\alpha_{ik}\,x_k
\ \ \ \ \ \ \ \ \ \ \ \ \ {\scriptstyle{(i\,=\,1\,\cdots\,n)}}
\end{equation}
then corresponds, in the variables $y_1, \dots, y_{ n-1}$, the
meroedrically isomorphic group:
\def\theequation{2}\begin{equation}
y_i'
=
\frac{\sum_{k=1}^{n-1}\,\alpha_{ik}\,y_k+\alpha_{in}}{
\sum_{k=1}^{n-1}\,\alpha_{nk}\,y_k+\alpha_{nn}}
\ \ \ \ \ \ \ \ \ \ \ \ \ 
{\scriptstyle{(k\,=\,1\,\cdots\,n\,-\,1)}},
\end{equation}
the $(n^2 - 1)$-term general projective group of the $R_{ n-1}$.
Thus, to every linear homogeneous transformation~\thetag{ 1}
corresponds a single projective transformation~\thetag{ 2}, hence a
completely determined collineation of the $R_{ n-1}$, whereas
conversely, to every projective transformation~\thetag{ 2}, there
correspond in total $\infty^1$ different linear homogeneous
transformations~\thetag{ 1}.

In addition, one can also establish a univalent invertible association
between linear homogeneous transformations in $x_1, \dots, x_n$ and
projective transformations in $y_1, \dots, y_{ n-1}$ when one submits
the constants $\alpha_{ ik}$ to the condition $\sum\, \pm \, \alpha_{
11} \cdots\, \alpha_{ nn} = 1$, hence when one considers instead of
the general the special linear homogeneous group which is
holoedrically isomorphic to the general projective group~\thetag{ 2}
(cf. Chap.~\ref{kapitel-26}, p.~\pageref{S-558-bis}). We want to
develop in details this association for the infinitesimal
transformations of the two groups.

The special linear homogeneous group in $x_1, \dots, x_n$ contains
the following $n^2 - 1$ independent infinitesimal transformations:
\def\theequation{3}\begin{equation}
x_i\,p_k,
\ \ \ \ \ \ \ \
x_i\,p_i-x_k\,p_k
\ \ \ \ \ \ \ \ \ \ \ \ \ 
{\scriptstyle{(i\,\gtrless\,k)}}.
\end{equation}
In order to find the corresponding infinitesimal transformations
in $y_1, \dots, y_{n-1}$, we only have to compute, for each
of the individual infinitesimal transformations just written, 
the increment:
\[
\delta\,y_i
=
\frac{x_n\,\delta\,x_i-x_i\,\delta\,x_n}{x_n^2}
\ \ \ \ \ \ \ \ \ \ \ \ \ {\scriptstyle{(i\,=\,1\,\cdots\,n\,-\,1)}}.
\]
In this way, we find the following table:
\def\theequation{4}\begin{equation}
\left\{
\aligned
&
x_n\,p_k
\equiv
q_k,
\ \ \ \ \ \ \
x_k\,p_n
\equiv
-\,y_k\,(y_1\,q_1+\cdots+y_{n-1}\,q_{n-1})
\\
&
x_k\,p_k-x_n\,p_n
\equiv
y_k\,q_k
+
y_1\,q_1
+\cdots+
y_{n-1}\,q_{n-1},
\ \ \ \ \ \ \
x_i\,p_k
\equiv
y_i\,q_k
\\
&
\ \ \ \ \ \ \ \ \ \ \ \ \ \ \ \ \ \ \ \ \ \ \ \ \ \ \ \ \ \ 
{\scriptstyle{(i,\,\,k\,=\,1\,\cdots\,n\,-\,1\,;\,\,\,
i\,\gtrless\,k)}},
\endaligned\right.
\end{equation}
where $q_i$ is written for $\partial f / \partial y_i$. This table
also provides inversely the infinitesimal transformations of the
special linear homogeneous group~\thetag{ 3} corresponding to 
every infinitesimal transformation of the projective group~\thetag{ 2};
indeed, from the equations~\thetag{ 4}, we obtain without effort:
\[
\aligned
n\,(y_1\,q_1+\cdots+y_{n-1}\,q_{n-1})
&
\equiv
x_1\,p_1+\cdots+x_n\,p_n
-
n\,x_n\,p_n
\\
n\,y_k\,q_k
&
\equiv
n\,x_k\,p_k
-
(x_1\,p_1+\cdots+x_n\,p_n).
\endaligned
\]

At present, it is also easy to indicate which $\infty^1$ infinitesimal
transformations of the general linear homogeneous group~\thetag{ 1}
correspond to a given infinitesimal transformation of the projective
group~\thetag{ 2}. Indeed, the infinitesimal transformation $x_1\, p_1
+ \cdots + x_n\, p_n$ reduces in the variables $y_1, \dots, y_{ n- 1}$
to the identity, for all increases of the $y_k$ vanish. Consequently,
if the infinitesimal transformation $Yf$ of the projective
group~\thetag{ 2} corresponds to the infinitesimal transformation $Xf$
of the special homogeneous group~\thetag{ 3}, then all $\infty^1$
infinitesimal transformations of the general homogeneous
group~\thetag{ 1} which correspond to $Yf$ are contained in the
expression $Xf + c\, (x_1\, p_1 + \cdots + x_n\, p_n)$, where $c$
denotes an arbitrary constant.

\medskip

If a system of equations:
\[
\Omega_k(y_1,\dots,y_{n-1})
=
0
\ \ \ \ \ \ \ \ \ \ \ \ \ {\scriptstyle{(k\,=\,1\,\cdots\,m)}}
\]
in $y_1, \dots, y_{ n-1}$ admits a finite or an infinitesimal
projective transformation, then naturally, the corresponding system of
equations in the $x$:
\[
\Omega_k
\bigg(
\frac{x_1}{x_n},
\,\,\dots,\,\,
\frac{x_{n-1}}{x_n}
\bigg)
=
0
\ \ \ \ \ \ \ \ \ \ \ \ \ {\scriptstyle{(k\,=\,1\,\cdots\,m)}}
\]
admits the corresponding finite or infinitesimal transformation
of the group~\thetag{ 3}; but in addition, because of its homogeneity,
it also admits yet the infinitesimal transformation: $x_1\, p_1 + 
\cdots + x_n\, p_n$. 

Conversely, every system of equations in $x_1, \dots, x_n$ which
admits $x_1\, p_1 + \cdots + x_n\, p_n$ is homogeneous. But now, when
we write a projective group of the $R_{ n-1}$ in the homogeneous
variables $x_1, \dots, x_n$, we are concerned only with the ratios of
the $x$, hence also, we are concerned only with systems of equations
that are homogeneous in the $x$. So, when we have written in a
homogeneous way the infinitesimal transformations of a projective
group of the $R_{ n-1}$ with the help of the table~\thetag{ 4} and
when we want to look up at the associated invariant systems of
equations, that is why we will always add yet the infinitesimal
transformation $x_1\, p_1 + \cdots + x_n\, p_n$. The group in $x_1,
\dots, x_n$ obtained in this way is the true analytic representation
in homogeneous variables of the concerned projective group of the $R_{
n-1}$.

For the study of a projective group, if one also wants to take the 
infinite into consideration, then 
one must write the group in homogeneous
variables. 

\sectionengellie{\S\,\,\,145.}

In the preceding paragraph, we have shown that the infinitesimal
projective transformations of the $R_{ n-1}$ can be replaced
by infinitesimal linear homogeneous transformations in $n$ variables.
At present, we want to imagine that an arbitrary transformation
of this sort in $x_1, \dots, x_n$ is presented, say:
\[
Xf
=
\sum_{i,\,\,k}^{1\cdots\,n}\,
a_{ki}\,x_i\,p_k,
\]
and we want to submit it to a closer examination. Namely, we want
to look for plane manifolds of the $R_{ n-1}$ which admit the
infinitesimal transformation in question. In this way, we will
succeed to show that $Xf$ can always be given a certain canonical
form after the introduction of appropriate new variables: 
$x_i' = \sum\, c_{ ik}\, x_k$. 

If a plane: $M_{ n-2}\colon \sum\, c_i\, x_i = 0$ of the $R_{ n-1}$
admits the infinitesimal transformation $Xf$, then according to 
Theorem~14, p.~\pageref{Theorem-14-S-112}, it admits at the same
time all finite transformations of the associated one-term group. 
Now, the $M_{ n-2}$ admits the infinitesimal transformation if and
only if the expression $X ( \sum\, c_i\, x_i)$ vanishes simultaneously
with $\sum\, c_i\, x_i$. Since $X ( \sum\, c_i\, x_i)$ is linear
in the $x_i$, this condition amounts to the fact that 
a relation of the form:
\[
X\bigg(
\sum_{k=1}^n\,c_k\,x_k
\bigg)
=
\rho\,\sum_{i=1}^n\,c_i\,x_i
\]
holds identically, where $\rho$ denotes a constant. The 
condition-equation following from this: 
\[
\sum_{i=1}^n\,
\bigg(
\sum_{k=1}^n\,a_{ki}\,c_k
-
\rho\,c_i
\bigg)\,x_i
=
0
\]
decomposes in the $n$ equations:
\def\theequation{5}\begin{equation}
a_{1i}\,c_1
+\cdots+
(a_{ii}-\rho)\,c_i
+\cdots+
a_{ni}\,c_n
=
0
\ \ \ \ \ \ \ \ \ \ \ \ \ {\scriptstyle{(i\,=\,1\,\cdots\,n)}},
\end{equation}
which can be satisfied only when the determinant:
\def\theequation{6}\begin{equation}
\left\vert
\begin{array}{cccccc}
a_{11}-\rho \,&\, a_{21} \,&\, \cdot \,&\, \cdot \,&\, \cdot \,&\, a_{n1}
\\
a_{12} \,&\, a_{22}-\rho \,&\, \cdot \,&\, \cdot \,&\, \cdot \,&\, a_{n2}
\\ 
\cdot \,&\, \cdot \,&\, \cdot \,&\, \cdot \,&\, \cdot \,&\, \cdot
\\
a_{1n} \,&\, a_{2n} \,&\, \cdot \,&\, \cdot \,&\, \cdot \,&\, 
a_{nn}-\rho
\end{array}
\right\vert
=
\Delta(\rho)
\end{equation}
vanishes. This produces for $\rho$ an equation of degree $n$
with surely $n$ roots, amongst which some can be multiple, though.
Therefore in all circumstances, there is at least one
plane $M_{ n-2} \colon$ $c_1\, x_1 + \cdots + c_n \, x_n = 0$
which remains invariant by the one-term group $Xf$. 

\smallercharacters{

It is known that one sees in exactly the same way that every
\emphasis{finite} projective transformation, or, when written
homogeneously, every transformation $x_i' = \sum\, b_{ ik}\, x_k$
likewise leaves fixed at least one plane $M_{ n-2}$. In addition, this
follows from the fact that every transformation $x_i' = \sum\, b_{
ik}\, x_k$ is associated to a one-term group $Xf$.

}

If the equation $\Delta ( \rho) = 0$ found above has exactly $n$
different roots $\rho$, then there are in total $n$ separate planes
$M_{ n-2}$ which remain invariant by the group $Xf$; indeed, two
distinct roots $\rho_1$ and $\rho_2$ of the equation for $\rho$ always
provide, because of the form of the equations~\thetag{ 5}, two
distinct systems of values $c_1 \colon c_2 \colon \cdots\, \colon
c_n$. By contrast, if there are multiple roots $\rho$, then various
cases can occur. If, for an $m$-fold root $\rho$, the
determinant~\thetag{ 6} vanishes itself, but not all of its $(n - 1)
\times (n - 1)$ determinants, then for the concerned value of $\rho$,
exactly $n - 1$ of the equations remain independent of each other and
the ratios of the $c_i$ are then determined all. To the $m$-fold root
is hence associated only a single invariant plane $M_{ n-2}$, but this
plane counts $m$ times. By contrast, if for an $m$-fold root not only
the determinant~\thetag{ 6} itself is equal to zero, but if all its
$(n-1) \times (n-1)$, \dots, $(n - q + 1) \times (n - q + 1)$
subdeterminants also vanish, whereas not all $( n - q) \times ( n -
q)$ subdeterminants do ($q \leqslant m$), then the equations~\thetag{ 5}
reduce to exactly $n - q$ independent equations and amongst the ratios
of the $c_i$, there remain $q - 1$ that can be chosen
arbitrarily. Thus in this case, the $m$-fold root $\rho$ gives
a family of $\infty^{ q - 1}$ planes $M_{ n-2}$ which remain 
individually invariant. 

It is easy to realize that a root of the equation $\Delta ( \rho) = 0$
is at least $q$-fold when all $(n - q + 1) \times ( n - q + 1)$
subdeterminants of $\Delta ( \rho)$ vanish for this root. Indeed, the
differential quotients of order $(q-1)$ of $\Delta ( \rho)$
with respect to $\rho$ express themselves as sums of the
$(n - q + 1) \times ( n - q + 1)$ subdeterminants of $\Delta ( \rho)$.

\medskip

If we would assume that the duality theory is already known at this
place, then we could immediately conclude that the infinitesimal
transformation $Xf$ leaves invariant at least one \emphasis{point} in
the $R_{ n-1}$. But we prefer to prove this also directly,
particularly because we gain on the occasion a deeper insight in the
circumstances.

In the homogeneous variables $x_1, \dots, x_n$, a point is represented
by $n - 1$ equations of the form $x_i\, x_k^0 - x_k x_i^0$; so the
point will admit all transformations of the one-term group $Xf$ when
$x_k^0\, X x_i - x_i^0\, X x_k$ vanishes by virtue of the equations
$x_i \, x_k^0 - x_k\, x_i^0 = 0$, hence when $n$ relations of the
form:
\[
Xx_i
=
\sigma\,x_i
\ \ \ \ \ \ \ \ \ \ \ \ \ {\scriptstyle{(i\,=\,1\,\cdots\,n)}}
\]
hold. For the $x_i$ and for $\sigma$, we therefore obtain the $n$
condition-equations:
\[
a_{i1}\,x_1
+\cdots+
(a_{ii}-\sigma)\,x_i
+\cdots+
a_{in}\,x_n
=
0
\ \ \ \ \ \ \ \ \ \ \ \ \ {\scriptstyle{(i\,=\,1\,\cdots\,n)}}.
\]
If we disregard the nonsensical solution $x_1 = \cdots = x_n = 0$, 
$\sigma$ must be a root of the equation:
\[
\Delta(\sigma)
=
\left\vert
\begin{array}{cccccc}
a_{11}-\sigma \,&\, a_{12} \,&\, \cdot \,&\, \cdot \,&\, \cdot \,&\, 
a_{1n}
\\
a_{21} \,&\, a_{22}-\sigma \,&\, \cdot \,&\, \cdot \,&\, \cdot \,&\, 
a_{2n}
\\ 
\cdot \,&\, \cdot \,&\, \cdot \,&\, \cdot \,&\, \cdot \,&\, \cdot
\\
a_{n1} \,&\, a_{n2} \,&\, \cdot \,&\, \cdot \,&\, \cdot \,&\, 
a_{nn}-\sigma
\end{array}
\right\vert
\]
and each such root produces an invariant point. The determination of
the points which stay fixed by $Xf$ therefore conducts to the same
algebraic equation as does the determination of the invariant plane:
$M_{ n-2}\colon$ $c_1\, x_1 + \cdots + c_n\, x_n = 0$.

Thus, if this equation of degree $n$ has $n$ different roots, then in
the $R_{ n-1} \colon$ $x_1 \colon x_2 \colon \cdots\, \colon x_n$, not
only $n$ different planes $M_{ n-2} \colon$ $\sum\, c_k\, x_k = 0$
remain invariant, but also simultaneously, $n$ separate points. In the
process, these $n$ points do not lie all in one and the same plane
$M_{ n-2}$, because if $n$ points stay fixed in a plane $M_{ n-2}$,
then necessarily, infinitely many points of the $M_{ n-2}$ keep their
positions, which is excluded under the assumption made. We can
express briefly this property of the $n$ invariant points by saying
that a nondegenerate $n$-frame \deutsch{$n$-Flach} remains invariant
by $Xf$.

If there appear multiple roots, then two cases must again be
distinguished. If, for an $m$-fold root, not all $(n - 1) \times
(n-1)$ subdeterminants of the determinant~\thetag{ 6} vanish, then
this root gives an invariant plane $M_{ n-2}$ counting $m$ times and
an invariant point counting $m$ times. By contrast, if, for the root
in question, all $(n - 1) \times ( n - 1)$, \dots, $(n - q + 1) \times
( n - q + 1)$ subdeterminants are equal to zero ($q \leqslant m$)
without all $( n - q) \times (n - q)$ subdeterminants vanishing, then
to this root is associated a family of $\infty^{ q - 1}$ individually
invariant planes $M_{ n-2}$ and a plane manifold of $\infty^{ q-1}$
individually invariant points. Therefore, the following holds

\def\theproposition{1}\begin{proposition}
Every infinitesimal transformation:
\[
\sum_{i,\,\,k}^{1\cdots\,n}\,
a_{ki}\,x_i\,p_k
\]
in the homogeneous variables $x_1, \dots, x_n$, or, what is the same,
every infinitesimal projective transformation in $n - 1$ variables:
\[
\frac{x_1}{x_n},
\,\,\dots,\,\,
\frac{x_{n-1}}{x_n},
\]
leaves invariant a series of points $x_1 \colon x_2 \colon \cdots\,
\colon x_n$ and a series of planes $M_{ n-2} \colon$ $c_1\, x_1 +
\cdots + c_n\, x_n = 0$. The points which stay fixed fill a finite
number, and to be precise at most $n$, separate manifolds. Likewise,
the invariant planes $M_{ n-2}$ form a finite number, and to be
precise at most $n$, separate linear pencils.
\end{proposition}

If the infinitesimal transformation $Xf$ has the form: $\sum\, x_k\, 
p_k$, then it actually leaves invariant all points $x_1 \colon x_2
\colon \cdots\, \colon x_n$ and naturally also, all planes
$M_{ n-2} \colon$ $\sum\, c_i\, x_i = 0$. 

\medskip

From what has been found up to now, one can now draw further 
conclusions. At first, we consider once more the special case
where the equation of degree $n$ discussed above has $n$ different
roots. Then there are $n$ separate invariant planes $M_{ n-2} 
\colon$ $\sum_i\, c_{ ki}\, x_i = 0$ which form a true
$n$-frame according to what precedes. We can therefore introduce:
\[
x_k'
=
\sum_{i=1}^n\,c_{ki}\,x_i
\ \ \ \ \ \ \ \ \ \ \ \ \ {\scriptstyle{(k\,=\,1\,\cdots\,n)}}
\]
as new homogeneous variables and in the process, we must 
obtain an infinitesimal transformation in $x_1', \dots, x_n'$
which leaves invariant the $n$ equation $x_k' = 0$, hence
which possesses the form:
\[
Xf
=
a_1'\,x_1'\,p_1'
+\cdots+
a_n'\,x_n'\,p_n'.
\]
Under the assumptions made, $Xf$ can be brought to this canonical
form. Naturally, no two of the quantities $a_1', \dots, a_n'$
are equal to each other here, for the equation:
\[
(a_1'-\rho)\cdots\,(a_n'-\rho)
=
0
\]
must obviously have $n$ distinct roots.

Now, similar canonical forms of $Xf$ also exist when the equation for
$\rho$ possesses multiple roots. However, we do not want to get
involved in the consideration of them, and we only want to show
that there is a canonical form to which \emphasis{every}
infinitesimal transformation:
\[
Xf
=
\sum_{i,\,\,k}^{1\cdots\,n}\,
a_{ki}\,x_i\,p_k
\]
can be brought thanks to an appropriate change of variables: $x_k' =
\sum\, h_{ ki}\, x_i$, that is to say hence, thanks to an appropriate
collineation of the $R_{ n-1}$, completely without paying heed to the
constitution of the equation $\Delta ( \rho) = 0$.

Since $Xf$ always leaves a point invariant, we can imagine that our
coordinates are chosen so that the point: $x_1 = \cdots = x_{ n-1} =
0$ stay fixed. On the occasion, we find:
\def\theequation{7}\begin{equation}
Xf
=
\sum_{k=1}^{n-1}\,\sum_{i=1}^{n-1}\,
a_{ik}'\,x_k\,p_i
+
\sum_{k=1}^n\,a_{nk}'\,x_k\,p_n
=
X'f
+
\sum_{k=1}^n\,a_{nk}'\,x_k\,p_n.
\end{equation}
To the linear homogeneous infinitesimal transformation $X'f$ in the $n
- 1$ variables $x_1, \dots, x_{ n-1}$ we can apply the same process as
the one which provided us with the reduction of $Xf$ to the
form~\thetag{ 7}, and we obtain in this way:
\[
Xf
=
\sum_{k=1}^{n-2}\,\sum_{i=1}^{n-2}\,
a_{ik}''\,x_k\,p_i
+
\sum_{k=1}^{n-1}\,a_{n-1,\,k}''\,x_k\,p_{n-1}
+
\sum_{k=1}^n\,a_{nk}'\,x_k\,p_n.
\]
Here, we can again treat in an analogous way the first term of the
right-hand side. Finally, we obtain the

\def\thetheorem{103}\begin{theorem}
In every linear homogeneous infinitesimal transformation with $n$
variables, one can introduce as new independent variables
$n$ linear homogeneous functions of these variables so that
the concerned infinitesimal transformation takes the form:
\[
a_{11}\,x_1\,p_1
+
(a_{21}\,x_1+a_{22}\,x_2)\,p_2
+\cdots+
(a_{n1}\,x_1+\cdots+a_{nn}\,x_n)\,p_n.
\]
\end{theorem}

It is of a certain interest to interpret the way in which this
result has been gained.

Since $Xf$ leaves a point invariant in any case, if we chose above
such a point as corner of coordinates \deutsch{Coordinateneckpunkt}:
$x_1 = \cdots = x_{ n-1} = 0$. The ratios $x_1 \colon x_2 \colon
\cdots\, \colon x_{ n -1}$ then represent the straight lines through
the chosen point; simultaneously with the points $x_1 \colon x_2
\colon \cdots \colon x_n$, these straight lines are permuted with each
other by the transformation $Xf$, and to be precise, by means of the
linear transformation $X'f$ in the $n - 1$ variables $x_1, \dots, x_{
n-1}$. But according to what precedes, a system of ratios $x_1 \colon
\cdots\, \colon x_{ n-1}$ must remain unchanged by $X'f$, that is to
say, a straight line through the point just said. If we choose this
line as edge $x_1 = \cdots = x_{ n-2} = 0$ of our coordinate system,
then the ratios $x_1 \colon \cdots\, \colon x_{ n-2}$ represent the
planes $M_2$ passing through this edge. In the same way, these $M_2$
are permuted with each other by $Xf$; one amongst them surely remains
invariant and gives again a closer determination of the coordinate
system, and so on.

In this way, one realizes that the coordinate system can be chosen so
that $Xf$ receives the normal form indicated above. However, we want
to recapitulate yet in a specific proposition the considerations just
made, since they express a general property of the infinitesimal
projective transformations of the $R_{ n-1}$:

\def\theproposition{2}\begin{proposition}
By every infinitesimal projective transformation of the $R_{ n-1}$
there remain invariant: at least one point; through every invariant
point: at least one straight line; \dots; and lastly, through
every invariant plane $M_{ n-2}$: at least one plane $M_{ n-2}$.
\end{proposition}

\sectionengellie{\S\,\,\,146.}

In the essence of things, the results of the previous paragraph are
known long since and they coincide in the main thing with the
reduction to a normal form, due to \name{Cauchy}, of a system of
linear ordinary differential equations with constant coefficients.

In what follows, we will generalize in an essential way the
developments conducted up to now, but beforehand, we must insert a few
observations which, strictly speaking, are certainly subordinate to
the general developments in Chap.~\ref{kapitel-29}, p.~\pageref{S-479}
sq.

Let an arbitrary $r$-term group $G_r$ with an $(r - m)$-term
\emphasis{invariant} subgroup $G_{ r - m}$ be presented, and 
let $T$ denote an arbitrary transformation of the $G_r$, while
$S$ denotes an arbitrary transformation of the $G_{ r - m}$. 
Then according to the assumption, there are certain equations
of the form:
\[
T^{-1}\,S\,T
=
S_1,
\ \ \ \ \ \ \
T\,S_1\,T^{-1}
=
S,
\]
where $S_1$ is again a transformation of the $G_{ r - m}$, and
in fact, a completely arbitrary one, provided only that one
chooses the $S$ appropriately.

Now, if every transformation $S$ leaves invariant a certain 
point figure
\deutsch{Punktfigur} $M$, so that:
\[
(M)\,S
=
(M),
\]
then the following equation also holds:
\[
(M)\,T\,T^{-1}\,S\,T
=
(M)\,T,
\]
and it shows that the figure $(M)\, T$ admits every transformation
$T^{ -1}\, S \, T$, hence actually, every transformation $S$.

\def\theproposition{3}\begin{proposition}
If $G_r$ is an $r$-term group, if $G_{ r-m}$ is an invariant
subgroup of it, and lastly, if $M$ is a point figure which admits
all transformations of the $G_{ r - m}$, then every position
that $M$ takes by means of a transformation $G_r$ also
remains invariant by all transformations of the $G_{ r - m}$.
\end{proposition}

Now in particular, let the two groups $G_r$ and $G_{ r - m}$ be
projective groups of the $R_{ n-1}$ and let them be written down in
$n$ homogeneous variables. Let the figure $M$ be a point. Then every
individual point invariant by all $S$ is transferred, by the execution
of all transformations $T$, only to points which again admit all
transformations $S$. Consequently, the \emphasis{totality of all}
points invariant by the group $G_{ r - m}$ also remains invariant by
the group $G_r$, though in general, the individual points of this
totality are permuted with each other by the $G_r$.

Above, we saw that all the points which remain invariant by a one-term
projective group can be ordered in at most $n$ mutually distinct plane
manifolds. Naturally, this also holds true for the points which keep
their position by our $G_{ r - m}$. Hence, let $M_1$, $M_2$, \dots, 
$M_\rho$ ($\rho \leqslant n$) be distinct plane manifolds all points
of which remain fixed by the $G_{ r - m}$, whereas, out of
these manifolds, there are no points invariant by the $G_{ r - m}$. 
Now, if the group $G_r$ were discontinuous\,---\,the considerations
made remain valid also for this case\,---, then the manifolds
$M_1, \dots, M_\rho$ could be permuted with each other by the
$G_r$. Not so when the $G_r$ is continuous, hence when it is
generated by infinitesimal transformations, because if for
instance $M_1$ would take new positions by means of the transformations
contained in the $G_r$, then these positions would form a
continuous family. But now, we have just seen that $M_1$ can at most
be given the finitely different positions $M_1, \dots,
M_\rho$. Consequently, no transformation of the $G_r$ can change the
position of $M_1$. In the same way, every other of the $\rho$ $M_k$
naturally stays at its place by all transformations of the $G_r$.

At present, we want to add yet the assumption that $m$ has the value
$1$. Let $X_1f, \dots, X_{ r - 1}f$ be independent infinitesimal
transformations of the $G_{ r-1}$ invariant in $G_r$, and let $X_1f,
\dots, X_{ r - 1}f$, $Yf$ be the transformations of the $G_r$;
moreover, let $x_1 = 0$, \dots, $x_q = 0$ be one of the plane
manifolds $M_1, \dots, M_\rho$ all points of which remain invariant by
$X_1f, \dots, X_{ r - 1}f$.

Since, as we know, $Yf$ surely leaves invariant the manifold $x_1 = 
0$, \dots, $x_q = 0$, it must have the form:
\[
Yf
=
\sum_{i,\,\,k}^{1\cdots\,q}\,
a_{ki}\,x_i\,p_k
+
\sum_{j=q+1}^n\,\sum_{\nu=1}^n\,
a_{j\nu}\,x_\nu\,p_j.
\]
The points of the manifold in question are transformed by $Yf$
(cf. Theorem~40, p.~\pageref{Theorem-40-S-233}) and
to be precise, by means of the transformation which is comes
from $Yf$ after the substitution $x_1 = \cdots = x_q = 0$, namely:
\def\theequation{8}\begin{equation}
\sum_{j=q+1}^n\,
\sum_{\nu=q+1}^n\,
a_{j\nu}\,x_\nu\,p_j\,;
\end{equation}
here, the variables $x_{ q+1}, \dots, x_n$ are to be considered as
coordinates for the points of the manifold. Now, the
transformation~\thetag{ 8} leaves invariant at least one point $x_{
q+1}^0 \colon \cdots\, \colon x_n^0$ \emphasis{inside} the manifold
$x_1 = \cdots = x_q = 0$, further, it leaves invariant a straight line
passing through this point, and so on.

Saying this, we have the

\def\thetheorem{104}\begin{theorem}
If an $(r + 1)$-term projective group $X_1f, \dots, X_rf$, $Yf$ of the
$R_{ n-1}$ contains an $r$-term invariant subgroup $X_1f, \dots, X_rf$
and if there are points of the $R_{ n-1}$ which remain invariant by
all $X_kf$, then every plane manifold consisting of such invariant
points which is not contained in a larger manifold of this sort also
admits all transformations of the $(r + 1)$-term group. Besides, the
points of each such manifold are are permuted with each other, but so 
that at least one of these points keeps its position by all
transformations of the $(r + 1)$-term group.
\end{theorem}

Now, we will apply this theorem to a special category of linear
homogeneous groups.

Let $X_1f, \dots, X_rf$ be the infinitesimal transformations of an
$r$-term linear homogeneous group, and let, as is always possible,
$X_1f, \dots, X_\rho f$ generate a $\rho$-term subgroup which is
invariant in the $( \rho + 1 )$-term subgroup $X_1f, \dots, X_{ \rho
+1}f$, when $\rho$ is smaller than $r$. Analytically, these
assumptions find their expressions in certain relations of the form:
\[
\leftbracket
X_i,\,X_{i+k}
\rightbracket
=
\sum_{s=1}^{i+k-1}\,
c_{iks}\,X_sf
\ \ \ \ \ \ \ \ \ \ \ \ \ 
{\scriptstyle{(i\,=\,1\,\cdots\,r\,-\,1\,;\,\,\,
k\,=\,1\,\cdots\,r\,-\,i)}}.
\]

Now, if we interpret as before $x_1, \dots, x_n$ as homogeneous
coordinates of an $R_{ n-1}$, we then see that there always is
in this $R_{ n-1}$ a point invariant by $X_1f$, 
moreover that there always is a point invariant by
$X_1f$ and by $X_2f$ as well, and finally, a point
which actually admits all transformations $X_1f, \dots, X_rf$, 
hence for which relations of the form:
\[
X_k\,x_i
=
\alpha_k\,x_i
\ \ \ \ \ \ \ \ \ \ \ \ \ 
{\scriptstyle{(k\,=\,1\,\cdots\,r\,;\,\,\,
i\,=\,1\,\cdots\,n)}}
\]
hold. If we choose the variables $x_i$ so that $x_1 = \cdots = 
x_{ n-1} = 0$ is this invariant point, then in each of the 
$r$ expressions:
\[
X_kf
=
\xi_{k1}\,p_1
+\cdots+
\xi_{kn}\,p_n,
\]
the $n - 1$ first coefficients $\xi_{ k1}, \dots, \xi_{k, \, n-1}$
depend only upon $x_1, \dots, x_{ n-1}$. 
Consequently, the reduced expressions:
\[
X_k'f
=
\xi_{k1}\,p_1
+\cdots+
\xi_{k,\,n-1}\,p_{n-1}
\]
stand again in the relationships:
\[
\leftbracket
X_i',\,X_{i+k}'
\rightbracket
=
\sum_{s=1}^{i+k-1}\,
c_{iks}\,X_s'f.
\]
However, $X_1'f, \dots, X_r'f$ need not be mutually independent
anymore, but nevertheless, we can apply to them the same
considerations as above for $X_1f, \dots, X_rf$, because
on the occasion, it was made absolutely no use of the independence
of the $X_kf$. Similarly as above, we can hence imagine that
the variables $x_1, \dots, x_{ n-1}$ are chosen in such a way
that all $\xi_{ k1}, \dots, \xi_{ k, \, n-2}$ depend only
on $x_1, \dots, x_{ n-2}$. At present, the expressions:
\[
X_k''f
=
\xi_{k1}\,p_1
+\cdots+
\xi_{k,\,n-2}\,p_{n-2}
\]
can be treated in exactly the same way, and so on.

We therefore obtain the

\renewcommand{\thefootnote}{\fnsymbol{footnote}}
\def\thetheorem{105}\begin{theorem}
If $X_1f, \dots, X_rf$ are independent infinitesimal transformations
of an $r$-term linear homogeneous group in the variables $x_1, \dots,
x_n$ and if relations of the specific form:
\[
\leftbracket
X_i,\,X_{i+k}
\rightbracket
=
\sum_{s=1}^{i+k-1}\,
c_{iks}\,X_sf
\ \ \ \ \ \ \ \ \ \ \ \ \ 
{\scriptstyle{(i\,=\,1\,\cdots\,r\,-\,1\,;\,\,\,
k\,=\,1\,\cdots\,r\,-\,i)}}
\]
hold, then one can always introduce linear homogeneous functions
of $x_1, \dots, x_n$ as new independent variables so that all $X_kf$ 
simultaneously receive the canonical form:\footnote{\,
\name{Lie}, Archiv for Math., Vol. 3, p.~110 and p.~111, Theorem~3; 
Christiania 1878.
}
\[
a_{k11}\,x_1\,p_1
+
(a_{k21}\,x_1+a_{k22}\,x_2)\,p_2
+\cdots+
(a_{kn1}\,x_1+\cdots+a_{knn}\,x_n)\,p_n.
\]
\end{theorem}
\renewcommand{\thefootnote}{\arabic{footnote}}

On the other hand, if we remember that $X_1f, \dots, X_rf$ 
is also a projective group of the $R_{ n-1}$, we can say:

\def\theproposition{4}\begin{proposition}
\label{Satz-4-S-589}
If $X_1f, \dots, X_rf$ is an $r$-term projective group of the $R_{
n-1}$ having the specific composition:
\[
\leftbracket
X_i,\,X_{i+k}
\rightbracket
=
\sum_{s=1}^{i+k-1}\,
c_{iks}\,X_sf
\ \ \ \ \ \ \ \ \ \ \ \ \ 
{\scriptstyle{(i\,=\,1\,\cdots\,r\,-\,1\,;\,\,\,
k\,=\,1\,\cdots\,r\,-\,i)}},
\]
then in the $R_{ n-1}$, there is at least one point $M_0$ invariant by
the group; through every invariant point, there passes at least one
invariant straight line; \dots; through every invariant plane $M_{
n-3}$, there passes at least one invariant plane $M_{ n-2}$. In
certain circumstances, several (infinitely many) such series of
invariant manifolds $M_0$, $M_1$, \dots, $M_{ n-2}$ are associated to
the group.
\end{proposition}

If we maintain the assumptions of this last proposition and if we
assume in addition that in the space $x_1 \colon x_2 \colon \cdots\,
x_n$, already several invariant plane manifolds ${\sf M}_{ \rho_1}, 
\dots, {\sf M}_{ \rho_q}$ are known, about which each one
is contained in the one following next, then amongst the series
of invariant manifolds $M_0$, $M_1$, \dots, $M_{ n-2}$
mentioned in Proposition~4, there obviously is 
at least one
for which $M_{ \rho_1}$ coincides with ${\sf M}_{ \rho_1}$,
at least one
for which $M_{ \rho_2}$ coincides with ${\sf M}_{ \rho_2}$, \dots,
at least
for which $M_{ \rho_q}$ coincides with ${\sf M}_{ \rho_q}$.

\sectionengellie{\S\,\,\,147.}

Apparently, the studies developed above certainly have a considerably
special character, but this special character is the very reason why 
they nevertheless possess a general signification for any finite 
continuous group, because for every such group, an isomorphic
linear homogeneous group can be shown, namely the associated adjoint
group (Chap.~\ref{kapitel-16}).

Amongst other things, we can employ our theory mentioned above
in order to prove that every group with more than two parameters
contains two-term subgroups, and that every group with more than
three parameters contains three-term subgroups.

Let $X_1f, \dots, X_rf$ be an arbitrary $r$-term group and let:
\[
E_kf
=
\sum_{i,\,\,j}^{1\cdots\,r}\,
c_{jki}\,e_j\,
\frac{\partial f}{\partial e_i}
\ \ \ \ \ \ \ \ \ \ \ \ \ {\scriptstyle{(k\,=\,1\,\cdots\,r)}}
\]
be the associated adjoint group, so that the relations:
\[
\leftbracket
X_i,\,X_k
\rightbracket
=
\sum_{s=1}^r\,c_{iks}\,X_sf,
\ \ \ \ \ \ \ \ \ \
\leftbracket
E_i,\,E_k
\rightbracket
=
\sum_{s=1}^r\,c_{iks}\,E_sf
\]
hold simultaneously.

We claim that in any case, $X_1f$ belongs to \emphasis{one} two-term
group, hence that an infinitesimal transformation $\lambda_2\, X_2f
+ \cdots + \lambda_r \, X_rf$ exists for which one has:
\[
\leftbracket
X_1,\,
\lambda_2\,X_2
+\cdots+
\lambda_r\,X_r
\rightbracket
=
\rho\,X_1f
+
\mu\,\sum_{k=2}^r\,\lambda_k\,X_kf.
\]
This condition represents itself in the form:
\[
\sum_{k=2}^r\,\lambda_k\,
\sum_{s=1}^r\,c_{1ks}\,X_sf
=
\rho\,X_1f
+
\mu\,\sum_{k=2}^r\,\lambda_k\,X_kf,
\]
and it decomposes itself in the $r$ equations:
\def\theequation{9}\begin{equation}
\left\{
\aligned
\sum_{k=2}^r\,\lambda_k\,c_{1k1}
&
=
\rho
\\
\sum_{k=2}^r\,\lambda_k\,c_{1ks}
&
=
\mu\,\lambda_s
\ \ \ \ \ \ \ \ \ \ \ \ \ {\scriptstyle{(s\,=\,2\,\cdots\,r)}}.
\endaligned\right.
\end{equation}
Here, this obviously amounts to satisfying only the last $r - 1$
equations; but one sees directly in the simplest way that this 
is possible. Indeed, one obtains for $\mu$ the equation:
\[
\left\vert
\begin{array}{cccccc}
c_{122}-\mu \,&\, c_{132} \,&\, \cdot \,&\, \cdot \,&\, \cdot \,&\,
c_{1r2}
\\
c_{123} \,&\, c_{133}-\mu \,&\, \cdot \,&\, \cdot \,&\, \cdot \,&\,
c_{1r3}
\\
\cdot \,&\, \cdot \,&\, \cdot \,&\, \cdot \,&\, \cdot \,&\, \cdot
\\
c_{12r} \,&\, c_{13r} \,&\, \cdot \,&\, \cdot \,&\, \cdot \,&\,
c_{1rr}-\mu
\end{array}
\right\vert
\]
which always possesses roots; for this reason, there always exists a
system of $\lambda_k$ not all vanishing which satisfies the above
conditions and which also determines $\rho$.

But the fact that the equations~\thetag{ 9} may be satisfied is also
an immediate consequence of our theory mentioned above; although this
observation seems hardly necessary here, where the relationships are
so simple, we nevertheless do not want to miss this, because for the
more general cases that are to be treated next, we will not make it
without this theory.

In the infinitesimal transformation:
\[
E_1f
=
\sum_{i=1}^r\,\sum_{k=1}^r\,c_{k1i}\,e_k\,
\frac{\partial f}{\partial e_i}
=
\sum_{i=1}^r\,\varepsilon_i\,
\frac{\partial f}{\partial e_i}
\]
of the adjoint group, all $\varepsilon_i$ are free of $e_1$, since
$c_{ 11i}$ is always zero. Now, the cut linear homogeneous
infinitesimal transformation $\varepsilon_2 \, p_2 + \cdots +
\varepsilon_r \, p_r$ in the $r - 1$ variables $e_2, \dots, e_r$
surely leaves invariant a system $e_2 \colon \cdots\, \colon e_r$; so
it is possible to satisfy the equations:
\[
\sum_{k=2}^r\,c_{k1i}\,e_k
=
\sigma\,e_i
\ \ \ \ \ \ \ \ \ \ \ \ \ {\scriptstyle{(i\,=\,2\,\cdots\,r)}},
\]
but these equations differ absolutely not from the last $r - 1$
equations~\thetag{ 9}, since one has $c_{ k 1 i} = - \, c_{
1ki}$ indeed.

At present, we can state more precisely our claim above that every
group with more than two parameters contains two-term subgroups, in the
following way:

\renewcommand{\thefootnote}{\fnsymbol{footnote}}
\def\theproposition{5}\begin{proposition}
Every infinitesimal transformation of a group with more than two
parameters is contained in at least one two-term 
subgroup.\footnote[1]{\,
\name{Lie}, Archiv for Math. Vol. 1, p.~192. Christiania 1876.
}
\end{proposition}
\renewcommand{\thefootnote}{\arabic{footnote}}

At present, we want to assume that $X_1f$ and $X_2f$ generate a
two-term subgroup, so that a relation of the form:
\[
\leftbracket
X_1,\,X_2
\rightbracket
=
c_{121}\,X_1f
+
c_{122}\,X_2f
\]
holds; we claim that, as soon as $r$ is larger than $3$, there also
always exists a three-term subgroup in which $X_1f$ and $X_2f$ are
contained.

At first, of the constants $c_{ 121}$ and $c_{ 122}$, when the two are
not already zero, we can always make one equal to zero. In fact, if,
say, $c_{ 121}$ is different from zero, then we introduce $X_1f +
\frac{ c_{ 122}}{ c_{ 121}} \, X_2f$ as new $X_1f$ and we obtain:
\[
\leftbracket
X_1,\,X_2
\rightbracket
=
c_{121}\,X_1f.
\]
Thus, every two-term group can be brought to this form.

If $\lambda_3 \, X_3 f + \cdots + \lambda_r \, X_r$ is supposed to
generate a three-term group together with $X_1f$ and $X_2f$, then one
must have:
\def\theequation{10}\begin{equation}
\left\{
\aligned
\big\leftbracket
X_1,\,\,
\lambda_3\,X_3
+\cdots+
\lambda_r\,X_r
\big\rightbracket
&
=
\alpha_1\,X_1f
+
\alpha_2\,X_2f
+
\mu\,(\lambda_3\,X_3f+\cdots+\lambda_r\,X_rf)
\\
\big\leftbracket
X_2,\,\,
\lambda_3\,X_3
+\cdots+
\lambda_r\,X_r
\big\rightbracket
&
=
\beta_1\,X_1f
+
\beta_2\,X_2f
+
\nu\,(\lambda_3\,X_3f+\cdots+\lambda_r\,X_rf).
\endaligned\right.
\end{equation}
From this, we obtain the following condition-equations:
\def\theequation{11}\begin{equation}
\left\{
\aligned
\sum_{k=3}^r\,\lambda_k\,c_{1k1}
&
=
\alpha_1,
\ \ \ \ \ \ \ \ \ \ \ \ \ \ \ \ \ \
\sum_{k=3}^r\,\lambda_k\,c_{1k2}
=
\alpha_2,
\\
\sum_{k=3}^r\,\lambda_k\,c_{2k1}
&
=
\beta_1,
\ \ \ \ \ \ \ \ \ \ \ \ \ \ \ \ \ \
\sum_{k=3}^r\,\lambda_k\,c_{2k2}
=
\beta_2,
\\
\sum_{k=3}^r\,c_{1ks}\,\lambda_k
&
=
\mu\,\lambda_s,
\ \ \ \ \ \ \ \ \ \ \ \ \ \ 
\sum_{k=3}^r\,c_{2ks}\,\lambda_k
=
\nu\,\lambda_s
\ \ \ \ \ \ \ \ \ \ \ \ \
{\scriptstyle{(s\,=\,3\,\cdots\,r)}}.
\endaligned\right.
\end{equation}

As one has observed, it is only necessary to prove that the $2\, ( r -
2)$ equations in the last row can be satisfied, since the remaining
equations can always be satisfied afterwards.

In order to settle this question, we form the infinitesimal
transformations: 
\[
\aligned
E_1f
&
=
\sum_{i=1}^r\,\sum_{k=1}^r\,
c_{k1i}\,e_k\,
\frac{\partial f}{\partial e_i}
=
\sum_{i=1}^r\,\varepsilon_{1i}\,
\frac{\partial f}{\partial e_i}
\\
E_2f
&
=
\sum_{i=1}^r\,\sum_{k=1}^r\,
c_{k2i}\,e_k\,
\frac{\partial f}{\partial e_i}
=
\sum_{i=1}^r\,\varepsilon_{2i}\,
\frac{\partial f}{\partial e_i},
\endaligned
\]
which, as we know, stand in the relationship:
\[
\leftbracket
E_1,\,E_2
\rightbracket
=
c_{121}\,E_1f.
\]
Since all $c_{ 11i}$, $c_{ 22i}$ and the $c_{ 122}, \dots, 
c_{ 12r}$ are equal to zero, both $e_1$ and $e_2$ do not
appear at all 
in $\varepsilon_{ 13}, \dots, \varepsilon_{ 1r}$, $\varepsilon_{ 23},
\dots, \varepsilon_{ 2r}$. The reduced expressions:
\[
\overline{E}_1f
=
\sum_{i=3}^r\,\varepsilon_{1i}\,
\frac{\partial f}{\partial e_i},
\ \ \ \ \ \ \ \ \ \ \ \ \ \
\overline{E}_2f
=
\sum_{i=3}^r\,\varepsilon_{2i}\,
\frac{\partial f}{\partial e_i}
\]
therefore are linear homogeneous infinitesimal transformations
in $e_3, \dots, e_r$ and they satisfy the relation:
\[
\big\leftbracket
\overline{E}_1,\,\overline{E}_2
\big\rightbracket
=
c_{121}\,\overline{E}_1f.
\]
According to what precedes, it follows from this that there exists
at least one system $e_3 \colon \cdots\, \colon e_r$ which admits
the two infinitesimal transformations $\overline{ E}_1f$ and
$\overline{ E}_2f$, hence which satisfies conditions of the form:
\[
\sum_{k=3}^r\,
c_{k1i}\,e_k
=
\sigma\,e_i,
\ \ \ \ \ \ \ \ \ \ \ \ \ \
\sum_{k=3}^r\,c_{k2i}\,e_k
=
\tau\,e_i
\ \ \ \ \ \ \ \ \ \ \ \ \ {\scriptstyle{(i\,=\,3\,\cdots\,r)}}.
\]
But this is exactly what was to be proved, for these last
equations are nothing else than the last equations~\thetag{ 11}, 
about which the question was whether they could be satisfied.
Naturally, the quantities $\alpha_1$, $\alpha_2$, $\beta_1$, 
$\beta_2$ are also determined together with the $\lambda_i$.

Since we can choose in all circumstances
$\lambda_3, \dots, \lambda_r$ so that equations of the
form~\thetag{ 3} hold, we can say:

\renewcommand{\thefootnote}{\fnsymbol{footnote}}
\def\thetheorem{106}\begin{theorem}
Every infinitesimal transformation and likewise every 
two-term subgroup of a group with more than three
parameters is contained in at least one three-term 
subgroup.\footnote[1]{\,
\name{Lie}, Archiv for Math. Vol. 1, p.~193, Vol. 3, pp.~114--116, 
Christiania 1876 and 1878.
}
\end{theorem}
\renewcommand{\thefootnote}{\arabic{footnote}}

One could be conducted to presume that every three-term subgroup
is also contained in at least one four-term subgroup, and so on, 
but this presumption is not confirmed. Our process of proof could
be employed further only when \emphasis{every} three-term group
$X_1f$, $X_2f$, $X_3f$ could be brought to the form:
\[
\aligned
&
\leftbracket
X_1,\,X_2
\rightbracket
=
c_{121}\,X_1f,
\ \ \ \ \ \ \ \ \ \ \ \ \
\leftbracket
X_1,\,X_3
\rightbracket
=
c_{131}\,X_1f
+
c_{132}\,X_2f
\\
&
\ \ \ \ \ \ \ \ \ \ \ \ \
\leftbracket
X_2,\,X_3
\rightbracket
=
c_{231}\,X_1f
+
c_{232}\,X_2f,
\endaligned
\]
hence not only when $X_1f$ would be invariant in the group $X_1f$, 
$X_2f$, but when this last group would also be invariant in the
entire three-term group. 

But for all three-term groups of the composition:
\[
\leftbracket
X_1,\,X_2
\rightbracket
=
X_1f,
\ \ \ \ \ \ \ \ \ \ 
\leftbracket
X_1,\,X_3
\rightbracket
=
2\,X_2f,
\ \ \ \ \ \ \ \ \ \ 
\leftbracket
X_2,\,X_3
\rightbracket
=
X_3f,
\]
this is \emphasis{not} the case (Chap.~\ref{kapitel-15}, Proposition~12,
p.~\pageref{Satz-12-S-270}).

Still, we want to briefly dwell on a special case in which an 
$m$-term subgroup is really contained in an $(m+1)$-term subgroup.

Let an arbitrary $r$-term group $X_1f, \dots, X_rf$ be presented
which contains an $m$-term subgroup $X_1f, \dots, X_mf$ having
the characteristic composition:
\[
\leftbracket
X_i,\,X_{i+k}
\rightbracket
=
\sum_{s=1}^{i+k-1}\,c_{i,\,i+k,\,s}\,X_sf
\ \ \ \ \ \ \ \ \ \ \ \ \ 
{\scriptstyle{(i\,<\,m,\,\,\,i\,+\,k\,\leqslant\,m)}}
\]
already mentioned. We claim that there always exists an 
$(m+1)$-term subgroup of the form:
\[
X_1f,\dots,X_mf,
\ \ \ \ \ \ \
\lambda_{m+1}\,X_{m+1}f
+\cdots+
\lambda_r\,X_rf.
\]

Our assertion amounts to the fact that certain relations of the form:
\[
\aligned
\big\leftbracket
X_j,\,\,\lambda_{m+1}\,X_{m+1}
+\cdots+
&
\lambda_r\,X_r
\big\rightbracket
=
\sum_{k=1}^m\,\alpha_{jk}\,X_kf
+
\mu_j\sum_{s=m+1}^r\,\lambda_s\,X_sf
\\
&\ \ \ 
{\scriptstyle{(j\,=\,1\,\cdots\,m)}}
\endaligned
\]
hold. By decomposition, it comes:
\def\theequation{12}\begin{equation}
\left\{
\aligned
\sum_{i=m+1}^r\,\lambda_i\,c_{jik}
&
=
\alpha_{jk}
\ \ \ \ \ \ \ \ \ \ \ \ \ 
{\scriptstyle{(j,\,\,k\,=\,1\,\cdots\,m)}}
\\
\sum_{i=m+1}^r\,\lambda_i\,c_{jis}
&
=
\mu_j\,\lambda_s
\ \ \ \ \ \ \ \ \ \ \ \ \ 
{\scriptstyle{(j\,=\,1\,\cdots\,m\,;\,\,\,
s\,=\,m\,+\,1\,\cdots\,r)}}.
\endaligned\right.
\end{equation}

Thus, the question is whether the last $m \, ( r - m)$ equations
can be satisfied; then always, the first $m^2$ equations can be
satisfied. 

In order to settle this question, we form the $m$ infinitesimal 
transformations:
\[
\aligned
E_kf
&
=
\sum_{i=1}^r\,\sum_{j=1}^r\,
c_{jki}\,e_j\,
\frac{\partial f}{\partial e_i}
=
\sum_{i=1}^r\,
\varepsilon_{ki}\,
\frac{\partial f}{\partial e_i}
\\
&
\ \ \ \ \ \ \ \ \ \ \ \ \ \
{\scriptstyle{(k\,=\,1\,\cdots\,m)}}
\endaligned
\]
which satisfy in pairs the relations:
\[
\leftbracket
E_i,\,E_{i+k}
\rightbracket
=
\sum_{s=1}^{i+k-1}\,
c_{i,\,i+k,\,s}\,E_sf
\ \ \ \ \ \ \ \ \ \ \ \ \ \
{\scriptstyle{(i\,+\,k\,\leqslant\,m)}}.
\]
We observe that all $c_{ j, \, k, \, m+1}$, \dots, $c_{ jkr}$
vanish for which $j$ and $k$ are smaller than $m+1$, and
from this we deduce that in $E_1f, \dots, E_mf$, all
coefficients $\varepsilon_{ k, \, m+1}, \dots, \varepsilon_{ kr}$
are free of $e_1, \dots, e_m$. The reduced infinitesimal 
transformations:
\[
\overline{E}_kf
=
\sum_{i=m+1}^r\,\varepsilon_{ki}\,
\frac{\partial f}{\partial e_i}
\ \ \ \ \ \ \ \ \ \ \ \ \ \
{\scriptstyle{(k\,=\,1\,\cdots\,m)}}
\]
in the variables $e_{ m+1}, \dots, e_r$ therefore stand pairwise in the
relationships:
\[
\big\leftbracket
\overline{E}_i,\,\overline{E}_k
\big\rightbracket
=
\sum_{s=1}^{i+k-1}\,
c_{i,\,i+k,\,s}\,\overline{E}_sf
\ \ \ \ \ \ \ \ \ \ \ \ \ \
{\scriptstyle{(i\,+\,k\,\leqslant\,m)}}
\]
and consequently, according to Proposition~4,
p.~\pageref{Satz-4-S-589}, there is a system $e_{ m+1} \colon \cdots\,
\colon e_r$ which admits all infinitesimal transformations $\overline{
E}_1f, \dots, \overline{ E}_mf$. But this says nothing but that it is
possible to satisfy the equations:
\[
\sum_{j=m+1}^r\,
c_{jks}\,e_j
=
\sigma_k\,e_s
\ \ \ \ \ \ \ \ \ \ \ \ \ \
{\scriptstyle{(k\,=\,1\,\cdots\,m\,;\,\,\,
s\,=\,m\,+\,1\,\cdots\,r)}}.
\]
But these equations are exactly the same as the equations~\thetag{ 12}
found above and thus, everything we wanted to show is effectively
proved.

\renewcommand{\thefootnote}{\fnsymbol{footnote}}
\def\thetheorem{107}\begin{theorem}
If, in an $r$-term group $X_1f, \dots, X_rf$, there is an $m$-term
subgroup $X_1f, \dots, X_mf$ having the specific composition:
\[
\leftbracket
X_i,\,X_{i+k}
\rightbracket
=
\sum_{s=1}^{i+k-1}\,
c_{i,\,i+k,\,s}\,X_sf
\ \ \ \ \ \ \ \ \ \ \ \ \ \
{\scriptstyle{(i\,<\,m,\,\,\,i\,+\,k\,\leqslant\,m)}},
\]
then this $m$-term subgroup is always contained in at least
one $(m + 1)$-term subgroup.\footnote[1]{\,
\name{Lie}, Archiv for Math., Vol. 3, pp.~114--116, 
Christiania 1878.
}
\end{theorem}
\renewcommand{\thefootnote}{\arabic{footnote}}

\sectionengellie{\S\,\,\,148.}

We yet want to derive the theorem stated just now also thanks to
conceptual considerations, by interpreting, as in
Chap.~\ref{kapitel-16}, the infinitesimal transformations $e_1\, X_1f
+ \cdots + e_r\, X_rf$ of our group as points of an $(r - 1)$-times
extended space with the homogeneous coordinates $e_1, \dots, e_r$.

Since the group $X_1f, \dots, X_mf, \dots, X_rf$ is isomorphic with
its adjoint group: $E_1f, \dots, E_mf, \dots, E_rf$, then under the
assumptions made, $E_1f, \dots, E_mf$ satisfy relations of the form:
\[
\leftbracket
E_i,\,E_{i+k}
\rightbracket
=
\sum_{s=1}^{i+k-1}\,
c_{i,\,i+k,\,s}\,E_sf
\ \ \ \ \ \ \ \ \ \ \ \ \ 
{\scriptstyle{(i\,+\,k\,\leqslant\,m)}}, 
\]
hence they generate a subgroup which, in the space $e_k$, is represented
by the $(m - 1)$-times extended plane manifold:
\[
e_{m+1}=0,
\,\,\,\dots,\,\,\,
e_r=0.
\]
This plane manifold naturally admits the subgroup $E_1f, \dots, 
E_mf$, and the same obviously holds true of the \emphasis{family}
of all $m$-times extended plane manifolds:
\[
\frac{e_{m+1}}{\mathfrak{e}_{m+1}}
=
\frac{e_{m+2}}{\mathfrak{e}_{m+2}}
=
\cdots
=
\frac{e_r}{\mathfrak{e}_r}
\]
which pass through the manifold $e_{ m+1} = 0$, \dots, $e_r = 0$.
But since the parameters $\mathfrak{ e}_{ m+1} \colon \cdots\,
\colon \mathfrak{ e}_r$ are transformed by a linear homogeneous
group which is isomorphic with the subgroup $E_1f, \dots, E_mf$
(Chap.~\ref{kapitel-23}, Proposition~5, p.~\pageref{Satz-5-S-472}),
then amongst the $m$-times extended manifolds of our invariant
family, there is at least one which remains invariant by the
subgroup $E_1f, \dots, E_mf$. This $m$-times extended manifold
is the image of an $(m+1)$-term subgroup of the group: 
$X_1f, \dots, X_rf$ 
(cf. Proposition~5, p.~\pageref{Satz-5-S-288}).

The conceptual considerations made just now which have again conducted
us to Theorem~107 are in essence identical to the analytical
considerations developed earlier on. However, the synthetical
explanation is more transparent than the analytical one \deutsch{Doch
ist die synthetische Begründung durchsichtiger als die analytische}.

\medskip

For the studies about the composition of transformation groups, 
it is actually advisable to set as fundamental the interpretation
of all infinitesimal transformations $e_1\, X_1f + \cdots + 
e_r\, X_rf$ as the points of a space which is transformed by the
linear homogeneous adjoint group. Thanks to a example, 
we will yet put in light the fruitfulness and the simplicity of
this method which shall also find multiple applications in 
the third volume, and at the same time, we will derive
a new remarkable statement.

We consider an $r$-term group $X_1f, \dots, X_rf$, the infinitesimal
transformations of which are linked together by relations of
the form:
\[
\leftbracket
X_i,\,X_{i+k}
\rightbracket
=
\sum_{s=1}^{i+k-1}\,
c_{i,\,i+k,\,s}\,
X_sf.
\]
The infinitesimal transformations $E_1f, \dots, E_rf$ of the 
adjoint group then satisfy the analogous equations:
\[
\leftbracket
E_i,\,E_{i+k}
\rightbracket
=
\sum_{s=1}^{i+k-1}\,
c_{i,\,i+k,\,s}\,E_sf.
\]

Consequently (Proposition~4, p.~\pageref{Satz-4-S-589}), the space
$e_k$ contains at least one point invariant by the adjoint group.
Through every point of this sort, there passes at least
one invariant straight line $M_1$, through every such straight
line, there passes at least one invariant plane $M_2$, and so on.

Now, if we interpret the points $e_k$ as infinitesimal transformations,
we see that our $r$-term group: $X_1f, \dots, X_rf$ contains at
least one invariant one-term subgroup; next, that every 
invariant one-term subgroup is contained in at least
one invariant two-term subgroup; next, that every invariant
two-term subgroup is contained in at least one invariant three-term
subgroup; and so on (cf. p.~\pageref{S-280-ter}).

Thus, the following holds true.

\renewcommand{\thefootnote}{\fnsymbol{footnote}}
\def\thetheorem{108}\begin{theorem}
If an $r$-term group contains $r$ independent infinitesimal 
transformations $Y_1f, \dots, Y_rf$ which satisfy relations
of the form:
\[
\leftbracket
Y_i,\,Y_{i+k}
\rightbracket
=
\sum_{s=1}^{i+k-1}\,
c_{i,\,i+k,\,s}\,Y_sf,
\]
then at the same time, it contains $r$ independent infinitesimal
transformations $Z_1f, \dots, Z_rf$ between which relations
of the form: 
\[
\leftbracket
Z_i,\,Z_{i+k}
\rightbracket
=
\sum_{s=1}^{i+k-1}\,
c_{i,\,i+k,\,s}\,Z_sf,
\]
hold; then $Z_1f, \dots, Z_if$ generate for every $i < r$ an 
$i$-term subgroup $\mathfrak{ G}_i$ and to be precise, every
$\mathfrak{ G}_i$ is invariant in every $\mathfrak{ G}_{ i+k}$
and in the group itself $Y_1f, \dots, Y_rf$ as well.\footnote[1]{\,
\name{Lie}, Archiv for Math., Vol. 3, p.~112 and p.~113, Christiania
1878; cf. also Vol. IX, pp.~79--82.
}
\end{theorem}
\renewcommand{\thefootnote}{\arabic{footnote}}

If we maintain the assumptions of this theorem and if we assume in
addition that by chance, several invariant subgroups, say $G_{
\rho_1}$, $G_{ \rho_2}$, \dots, $G_{ \rho_q}$ are known, of which each
one comprises the one following next, then it becomes evident (cf. the
concluding remarks of the \S\,\,146) that the subgroups $\mathfrak{
G}_i$ discussed in Theorem~108 can be chosen in such a way that
$\mathfrak{ G}_{ \rho_1}$ coincides with $G_{ \rho_1}$, that
$\mathfrak{ G}_{ \rho_2}$ coincides with $G_{ \rho_2}$, \dots, that
$\mathfrak{ G}_{ \rho_q}$ coincides with $G_{ \rho_q}$.

On the other hand, when an $r$-term group: $X_1f, \dots, X_rf$ of the
specific composition considered here is presented, it is always
possible to derive certain invariant subgroups by differentiation.
Indeed, all $\leftbracket X_i, \, X_k \rightbracket$ generate a
\emphasis{firstly derived} invariant subgroup with say $r_1 < r$
independent infinitesimal transformations. If, as in
Chap.~\ref{kapitel-15}, p.~\pageref{S-266}, we denote them by $X_1'f,
\dots, X_{ r_1}'f$, then all $\leftbracket X_i', \, X_k'
\rightbracket$ generate a \emphasis{secondly derived} invariant
subgroup of the $r$-term group with, say, $r_2 < r_1$ independent
infinitesimal transformations $X_1'' f, \dots, X_{ r_2}'' f$, and so
on.

It is therefore possible to bring our group to a form $U_1f, \dots, 
U_i f, \dots, U_rf$ such that \emphasis{firstly} for every $i$, 
the infinitesimal transformations $U_1f, \dots, U_if$ generate
an invariant subgroup, such that \emphasis{secondly} all 
$\leftbracket U_i, \, U_k \rightbracket$ can be linearly
deduced from $U_1f, \dots, U_{ r_1}f$, such that \emphasis{thirdly} all
$\big \leftbracket \leftbracket U_i, \, U_k \rightbracket,\,\,
\leftbracket U_\alpha, \, U_\beta \rightbracket \big \rightbracket$
can be linearly deduced from $U_1f, \dots, U_{ r_2}f$, and so on.

\linestop


\chapter{Approach \deutsch{Ansatz} towards the Determination
\\
of All Finite Continuous Groups
\\
of the $n$-times Extended Space
}
\label{kapitel-28}
\chaptermark{Approach towards the Determination of All Finite
Groups of the $R_n$}

\setcounter{footnote}{0}

\abstract*{??}

It is unlikely that one becomes close to be in the position of
determining all finite continuous transformation groups; indeed, it is
even uncertain whether this will ever succeed. Therefore, instead of the
general problem to determine \emphasis{all} finite continuous groups,
one would do well to tackle at first more special problems which
concern the determination of certain categories of finite continuous
groups. Namely, more special problems of this kind are the following 
three:

\begin{itemize}

\smallskip\item[{\sf Firstly}]
\ \ \ \ \ \ \ \ \ \ \ \ \ \ \
the determination of all $r$-term groups in $n$ variables.

\smallskip\item[{\sf Secondly}]
\ \ \ \ \ \ \ \ \ \ \ \ \ \ \
the determination of all $r$-term groups in general.

\smallskip\item[{\sf Thirdly}]
\ \ \ \ \ \ \ \ \ \ \ \ \ \ \
the determination of all finite continuous groups in $n$ variables.

\end{itemize}\smallskip

In Chap.~\ref{kapitel-22}, p.~\pageref{S-429-sq} sq., we have shown
that the settlement of the first of these problems, aside from
executable operations, requires in any case only the integration of
simultaneous systems of ordinary differential equations. Moreover, we
found that the second of our three problems can be led back to the
first one (Theorem~84, p.~\pageref{Theorem-84-S-458}).

By contrast, the third problem is not at all settled by means of the
developments of the Chap.~\ref{kapitel-22}. Namely if $n > 1$, then
for every value of $r$, how large can it be though, there always are
$r$-term groups in $n$ variables. For instance, if one chooses $r$
functions $F_1, \dots, F_r$ of $x_1, \dots, x_{ n-1}$ so that between
them, no linear relation of the form:
\[
c_1\,F_1
+\cdots+
c_r\,F_r
=
0
\]
with constant coefficients holds, then the $r$ infinitesimal
transformations:
\[
F_1(x_1,\dots,x_{n-1})\,
\frac{\partial f}{\partial x_n},
\,\,\,\dots,\,\,\,
F_r(x_1,\dots,x_{n-1})\,
\frac{\partial f}{\partial x_n}
\]
are mutually independent and pairwise exchangeable, hence they
generate an $r$-term group in the $n$ variables $x_1, \dots, x_n$.

In spite of the important results of the Chap.~\ref{kapitel-22}, the
third one of the problems indicated therefore still awaits for a
solution. This is why we will tackle this problem in the present
chapter and at least provide an approach towards its settlement.

\emphasis{In the sequel, we decompose the problem of determining all
groups of an $n$-times extended space $R_n$ in a series of individual
problems which are independent of each other.} We make it to thanks to
a natural division of the groups of the $R_n$ in classes which are
selected in such a way that the groups of some given class can be
determined without it to be necessary to know any of the groups of the
remaining classes. Admittedly, we cannot indicate a general method
which accomplishes in every case the determination of all groups
of a class; nonetheless, we provide important statements about
the groups belonging to a given class and on the other hand, we develop
a series of considerations which facilitate the determination of all
groups in a class; in the next chapter, this shall be illustrated by
special applications. For these general discussions, we essentially
restrict ourselves to transitive groups, because our classification
precisely is of specific practical meaning for the transitive groups.

This finds good reasons in the fact that later (in the third Volume),
for the determination of all finite continuous groups in one or two or
three variables, we will exploit only in part the developments of the
present chapter.  For the determination of the \emphasis{primitive}
groups of a space, the process explained here is firmly effective; not
so for the determination of the imprimitive groups, and rather, it is
advisable to take a different route. Every imprimitive group of the
$R_n$ decomposes this space in an invariant family of manifolds and
transforms the manifolds of this family by means of an isomorphic
group in less than $n$ variables. From this, it results that one would
do well to tackle at first the determination of all imprimitive groups
of the $R_n$, when one already know all groups in less than $n$
variables. Applying this fundamental principle, we shall undertake in
Volume~III the determination of all imprimitive groups of the $R_2$
and of certain imprimitive groups of the $R_3$.

\sectionengellie{\S\,\,\,149.}

Let:
\[
X_kf
=
\sum_{i=1}^n\,\xi_{ki}(x_1,\dots,x_n)\,
\frac{\partial f}{\partial x_i}
\ \ \ \ \ \ \ \ \ \ \ \ \ {\scriptstyle{(k\,=\,1\,\cdots\,r)}}
\]
be a completely arbitrary $r$-term group of the $n$-times extended
space $x_1, \dots, x_n$, or briefly of the $R_n$.

Under the guidance of Chap.~\ref{kapitel-25}, p.~\pageref{S-522-sq},
we prolong this group by viewing $x_1, \dots, x_n$ as functions of an
auxiliary variable $t$ which is absolutely not transformed by our
group and by considering that the differential quotients: $\D\, x_i /
\D\, t = x_i'$ are transformed by the group.  In the $2\, n$
variables: $x_1, \dots, x_n$, $x_1', \dots, x_n'$, we then obtain the
prolonged group:
\[
\overline{X}_kf
=
\sum_{i=1}^n\,
\xi_{ki}(x)\,
\frac{\partial f}{\partial x_i}
+
\sum_{i=1}^n\,
\bigg(
\sum_{\nu=1}^n\,
\frac{\partial\xi_{ki}}{\partial x_\nu}\,x_\nu'
\bigg)\,
\frac{\partial f}{\partial x_i'}
\ \ \ \ \ \ \ \ \ \ \ \ \ {\scriptstyle{(k\,=\,1\,\cdots\,r)}}.
\]
As we know, this prolonged group indicates in which way the $\infty^{
2n-1}$ line elements: $x_1, \dots, x_n$, $x_1' \colon x_2' \colon
\cdots\, \colon x_n'$ of the space $x_1, \dots, x_n$ are
permuted with each other by the group: $X_1f, \dots, X_rf$
(cf. p.~\pageref{S-525-bis}).

Every point: $x_1^0, \dots, x_n^0$ in general position remains
invariant by a completely determined number of independent
infinitesimal transformations: $e_1\, X_1f + \cdots + e_r\, X_rf$ 
and to be precise, at least by $r - n$ such transformations, 
and at most by $r - 1$. We want to denote by $r - q$ the number
of these independent infinitesimal transformations and
to suppose that: $X_1^0f, \dots, X_{ r - q}^0 f$ are such independent
infinitesimal transformations; then they certainly generate
an $(r - q)$-term subgroup of the group: $X_1f, \dots, X_rf$
(cf. Chap.~\ref{kapitel-12}, p.~\pageref{Satz-1-S-205} sq.). 

In the series expansions of $X_1^0f, \dots, X_{ r - q}^0f$ with 
respect to the powers of the $x_i - x_i^0$, all terms of order
zero are naturally lacking and there appear only terms of 
first or higher order:
\[
X_k^0f
=
\sum_{\nu=1}^n\,
\bigg\{
\sum_{i=1}^n\,
\alpha_{ki\nu}(x_1^0,\dots,x_n^0)\,
(x_i-x_i^0)
+\cdots
\bigg\}\,
\frac{\partial f}{\partial x_\nu}
\ \ \ \ \ \ \ \ \ \ \ \ \ 
{\scriptstyle{(k\,=\,1\,\cdots\,r\,-\,q)}}.
\]

The group: $X_1^0f, \dots, X_{ r - q}^0f$ leaves the point $x_1^0, 
\dots, x_n^0$ invariant, but it permutes the line elements which
pass through this point; how? this is what the associated
prolonged group shows:
\[
\aligned
\overline{X}_k^0f
&
=
\sum_{\nu=1}^n\,
\bigg\{
\sum_{i=1}^n\,
\alpha_{ki\nu}(x_1^0,\dots,x_n^0)\,
(x_i-x_i^0)
+\cdots
\bigg\}\,
\frac{\partial f}{\partial x_\nu}
\\
&
\ \ \ \ \
+
\sum_{\nu=1}^n\,
\bigg\{
\sum_{i=1}^n\,\alpha_{ki\nu}(x_1^0,\dots,x_n^0)\,x_n'
+\cdots
\bigg\}\,
\frac{\partial f}{\partial x_\nu'}
\ \ \ \ \ \ \ \ \ \ \ \ \ 
{\scriptstyle{(k\,=\,1\,\cdots\,r\,-\,q)}},
\endaligned
\]
and it transforms the $\infty^{ 2n-1}$ line elements of the space
$x_1, \dots, x_n$ in exactly the same way as the group: $X_1^0 f,
\dots, X_{ r - q}^0f$. Hence, if one wants to restrict oneself to the
line elements which pass through the point $x_1^0, \dots, x_n^0$, and
to disregard the remaining ones, then under the guidance of
Chap.~\ref{kapitel-14}, p.~\pageref{S-233-sq} sq., one has to leave
out, in the $\overline{ X}_k^0f$, all terms with the differential
quotients of $f$ with respect to $x_1, \dots, x_n$ and to make the
substitution: $x_1 = x_1^0$, \dots, $x_n = x_n^0$ in the terms
remaining. The so obtained reduced infinitesimal transformations:
\def\theequation{1}\begin{equation}
L_kf
=
\sum_{i,\,\,\nu}^{1\cdots\,n}\,
\alpha_{ki\nu}(x_1^0,\dots,x_n^0)\,
x_i'\,
\frac{\partial f}{\partial x_\nu'}
\ \ \ \ \ \ \ \ \ \ \ \ \ 
{\scriptstyle{(k\,=\,1\,\cdots\,r\,-\,q)}}
\end{equation}
generate \emphasis{a linear homogeneous group in the $n$ variables
$x_1', \dots, x_n'$}; this group, which is isomorphic with
the group: $\overline{ X}_1^0f, \dots, \overline{ X}_{ r - q}^0f$, 
and naturally also, with the group: $X_1^0f, \dots, X_{ r - q}^0f$, 
indicates in which way the two groups said just now transform the
$\infty^{ n-1}$ line elements: $x_1' \colon x_2' \colon \cdots\, 
x_n'$ through the point $x_1^0, \dots, x_n^0$. 

Visibly, 
the linear homogeneous group: $L_1f, \dots, L_{ r - q}f$ is perfectly
determined by the terms of first order in the power series expansions
of $X_1^0f, \dots, X_{ r - q}^0f$, and therefore, it contains
as many essential parameters as the group: $X_1^0f, \dots, X_{ r - 
q}^0f$ contains independent infinitesimal transformations of first
order out of which no transformation of second order, or of
higher order, can be linearly deduced. From this, it follows
that one can set up the group: $L_1f, \dots, L_{ r - q}f$
as soon as one knows how many independent infinitesimal
transformations of first order the group: $X_1, \dots, X_r$
contains in the neighbourhood of $x_1^0, \dots, x_n^0$
out of which no transformation of second or higher order can
be linearly deduced, and in addition, as soon as one knows
the terms of first order in the power series expansions of these
infinitesimal transformations. Conversely, as soon as one knows
the group: $L_1f, \dots, L_rf$, one can indicate the number of and
the terms of those independent infinitesimal transformations of first
order in the group: $X_1f, \dots, X_rf$ out of which no transformation
of second or higher order in the $x_i - x_i^0$ can be linearly 
deduced.

\smallercharacters{

One can derive the linear homogeneous group: $L_1f, \dots, 
L_{ r - q}f$ also in the following way:

Since the question is about the way in which the directions through
the fixed point: $x_1^0, \dots, x_n^0$ are transformed, one may 
substitute the group: $X_1^0f, \dots, X_{ r - q}^0f$ for the
following:
\[
\sum_{i,\,\,\nu}^{1\cdots\,n}\,
\alpha_{ki\nu}(x_1^0,\dots,x_n^0)\,
(x_i-x_i^0)\,
\frac{\partial f}{\partial x_\nu}
\ \ \ \ \ \ \ \ \ \ \ \ \ 
{\scriptstyle{(k\,=\,1\,\cdots\,r\,-\,q)}}
\]
which is obtained by leaving out all terms of second and higher
order. Now, $x_1 - x_1^0$, \dots, $x_n - x_n^0$ can here directly
conceived as homogeneous coordinates of the line elements through
the point $x_1^0, \dots, x_n^0$, whence the group:
\[
\sum_{i,\,\,\nu}^{1\cdots\,n}\,
\alpha_{ki\nu}(x_1^0,\dots,x_n^0)\,
(x_i-x_i^0)\,
\frac{\partial f}{\partial(x_\nu-x_\nu^0)}
\ \ \ \ \ \ \ \ \ \ \ \ \ 
{\scriptstyle{(k\,=\,1\,\cdots\,r\,-\,q)}}
\]
indicates how these line elements are transformed. This is coherent
with the above.

}

We have seen that the $r$-term group: $X_1f, \dots, X_rf$ 
associates to every point $x_1^0, \dots, x_n^0$ in general position
a completely determined linear homogeneous group~\thetag{ 1}
which, though, turns out to be different for different points.
At present, we study how this linear homogeneous group behaves after
the introduction of some new variables.

In place of $x_1, \dots, x_n$, we introduce the new variables:
\def\theequation{2}\begin{equation}
y_i
=
y_i^0
+
\sum_{\nu=1}^n\,
a_{i\nu}(x_\nu-x_\nu^0)
+\cdots
\ \ \ \ \ \ \ \ \ \ \ \ \ {\scriptstyle{(i\,=\,1\,\cdots\,n)}}
\end{equation}
which are ordinary power series of $x_1 - x_1^0$, \dots, $x_n - x_n^0$;
as always, we assume on the occasion that the determinant:
\[
\sum\,\pm\,a_{11}\cdots\,a_{nn}
\]
is different from zero, so that inversely also, the $x_i$ are ordinary
power series of $y_1 - y_1^0$, \dots, $y_n - y_n^0$.

In $y_1, \dots, y_n$, let the $X_k^0$ receive the form:
\[
Y_k^0f
=
\sum_{\nu=1}^n\,
\bigg\{
\sum_{i=1}^n\,\beta_{ki\nu}\,(y_i-y_i^0)
+\cdots
\bigg\}\,
\frac{\partial f}{\partial y_\nu}
\ \ \ \ \ \ \ \ \ \ \ \ \ 
{\scriptstyle{(k\,=\,1\,\cdots\,r\,-\,q)}}.
\]
Then it is clear at first that the group which comes into existence
from $X_1f, \dots, X_rf$ after the introduction of the new variables
$y_1, \dots, y_n$ associates to the point $y_1^0, \dots, y_n^0$ the
linear homogeneous group:
\def\theequation{1'}\begin{equation}
\mathfrak{L}_kf
=
\sum_{i,\,\,\nu}^{1\cdots\,n}\,
\beta_{ki\nu}\,y_i'\,
\frac{\partial f}{\partial y_\nu'}
\ \ \ \ \ \ \ \ \ \ \ \ \ 
{\scriptstyle{(k\,=\,1\,\cdots\,r\,-\,q)}}.
\end{equation}
But on the other hand, it is clear (cf. Chap.~\ref{kapitel-11}, 
p.~\pageref{S-196-sq} sq.) that the terms of first order:
\[
\sum_{i,\,\,\nu}^{1\cdots\,r}\,
\beta_{ki\nu}\,(y_i-y_i^0)\,
\frac{\partial f}{\partial y_\nu}
\]
of $Y_k^0$ can be obtained from the terms of first order:
\[
\sum_{i,\,\,\nu}^{1\cdots\,n}\,
\alpha_{ki\nu}\,(x_i-x_i^0)\,
\frac{\partial f}{\partial x_\nu}
\]
of $X_k^0f$ after the introduction of the new variables:
\[
y_i
=
y_i^0
+
\sum_{\nu=1}^n\,
a_{i\nu}\,(x_\nu-x_\nu^0)
\ \ \ \ \ \ \ \ \ \ \ \ \ {\scriptstyle{(i\,=\,1\,\cdots\,n)}}.
\]
Consequently, it follows that the linear homogeneous group~\thetag{ 1'}
comes into existence from the linear homogeneous group~\thetag{ 1}
when one introduces in~\thetag{ 1}, by means of the linear
homogeneous transformation:
\[
y_i'
=
\sum_{\nu=1}^n\,a_{i\nu}\,x_\nu'
\ \ \ \ \ \ \ \ \ \ \ \ \ {\scriptstyle{(i\,=\,1\,\cdots\,n)}},
\]
the new variables: $y_1', \dots, y_n'$ in place of the $x'$.

In this lies the reason why the linear homogeneous group~\thetag{ 1}
is essentially independent from the analytic representation of the
group: $X_1f, \dots, X_rf$, that is to say, from the choice of the
variables; indeed, if, by means of a transformation~\thetag{ 2},
one introduces new variables in the group: $X_1f, \dots, X_rf$,
then the linear homogeneous group~\thetag{ 1} converts into another
linear homogeneous group~\thetag{ 1'} which is conjugate to~\thetag{
1} inside the general linear homogeneous group of the $R_n$ (cf. 
Chap.~\ref{kapitel-16}, p.~\pageref{S-280-bis}).
Thanks to an appropriate choice of the constants $a_{ i\nu}$ in the
transformation~\thetag{ 2}, one can obviously insure that the 
group~\thetag{ 1'} associated to the point $y_k^0$ becomes an
arbitrary group conjugate to~\thetag{ 1}.

Now, we specially assume that the transformation~\thetag{ 2} belongs
to the group: $X_1f, \dots, X_rf$ itself. In this case, \thetag{ 1'}
is visibly the linear homogeneous group that the group: $X_1f, \dots, 
X_rf$ associates to the point: $y_1^0, \dots, y_n^0$, and consequently,
\thetag{ 1'} comes from~\thetag{ 1} when one replaces the $x^0$ in
the $\alpha_{ ki\nu} ( x_1^0, \dots, x_n^0)$ by the $y^0$. But since
the two groups~\thetag{ 1} and~\thetag{ 1'} are conjugate inside
the general linear homogeneous group of the $R_n$, we have the
following

\def\thetheorem{109}\begin{theorem}
Every $r$-term group: $X_1f, \dots, X_rf$ of the space $x_1, \dots,
x_n$ associates to every point $x_1^0, \dots, x_n^0$ in general
position a completely determined linear homogeneous group of the $R_n$
which indicates in which way the line elements through this point are
transformed, as soon as the point is fixed.  To those points which can
be transferred one to another by means of transformations of the
group: $X_1f, \dots, X_rf$ are associated linear homogeneous groups
which are conjugate inside the general linear homogeneous group of the
$R_n$. In particular, if the group: $X_1f, \dots, X_rf$ is transitive,
then to all points which lie in no invariant manifold, it associates
linear homogeneous groups that are conjugate to each other inside the
general linear homogeneous group.
\end{theorem}

The above theorem which recapitulates the most important result up to
now, provides now the classification of all groups of the $R_n$
announced in the introduction.

We consider at first the \emphasis{transitive} groups.

We reckon as belonging to the same class two transitive groups $G$ and
$\Gamma$ when the linear homogeneous group that $G$ associates to an
arbitrary point in general position is conjugate to the linear
homogeneous group that $\Gamma$ associates to such an arbitrary point.
In the opposite case, we reckon $G$ and $\Gamma$ as belonging to
different classes.

Thus, we differentiate as many classes of transitive groups of the
$R_n$ as there are \emphasis{types} of subgroups of the general linear
homogeneous group of the $R_n$ (cf. p.~\pageref{S-281-bis}).  Later,
we will see that to each one of such classes, there belongs in any
case a transitive group of the $R_n$.

If two transitive groups $G$ and $\Gamma$ of the space $x_1, \dots,
x_n$ belong to the same class, then in the neighbourhood of any point
$x_1^0, \dots, x_n^0$ in general position, they obviously contain the
same number of independent infinitesimal transformations of first
order in the $x_i - x_i^0$ out of which no transformation of second or
of higher order can be linearly deduced. In addition, since one can
always, by introducing new variables, reshape $\Gamma$ so that it
associates to the point $x_1^0, \dots, x_n^0$ exactly the same linear
homogeneous group as does $G$, then in all cases, one can insure that
the terms of first order in the infinitesimal transformations of first
order in question are the same for the two groups. In addition,
because they are transitive, $G$ and $\Gamma$ contain, in the
neighbourhood of $x_1^0, \dots, x_n^0$, $n$ independent infinitesimal
transformations of zeroth order out of which no transformation of
first or higher order can be linearly deduced; by contrast, the
numbers of terms and the initial terms of second, third, \dots orders
can very well be different for $\Gamma$ and for $G$.  \emphasis{Here
lies the reason why two transitive groups of the space $x_1, \dots,
x_n$ which belong to the same class need absolutely not have the same
number of parameters.}

\smallercharacters{

We now turn ourselves to the \emphasis{intransitive} groups.

To every transitive group of the $R_n$ was associated a completely
determined type of linear homogeneous group of the $R_n$; for the
intransitive groups, this is in general not the case.  Every
intransitive group $G$ of the $R_n$ decomposes this space in a family
of $\infty^{ n-q}$ ($0 < q < n$) individually invariant $q$-times
extended manifolds $M_q$, but so that the points of every individual
$M_q$ are transformed transitively (cf. Chap.~\ref{kapitel-13},
p.~\pageref{S-216}). From this, it follows that $G$ associates to all
points of one and the same $M_q$ always conjugate linear homogeneous
groups, but not necessarily to the points of different $M_q$.

In general, our $\infty^{ n - q}$ $M_q$ are gathered in continuous
families so that conjugate linear homogeneous groups are associated
only to the points which belong to the $M_q$ in the same family.  If
each such family consists of exactly $\infty^m$ $M_q$, then the whole
$R_n$ decomposes in $\infty^{ n - q + m}$ $(q + m)$-times extended
manifolds ${\sf M}_{ q + m}$ and to every ${\sf M}_{ q + m}$ is
associated a completely determined type of linear homogeneous groups,
while to different ${\sf M}_{ q+m}$ are associated different types
too. So to our group is associated a family of $\infty^{ n - q - m}$
different types; the totality of all these types can naturally be
represented by certain analytic expressions with $n - q - m$ essential
arbitrary parameters. We can also express this as follows: all the
concerned types belong to the same kind of types
\deutsch{Typengattung} (cf. Chap.~\ref{kapitel-22},
p.~\pageref{S-448}).

Now, we reckon two intransitive groups of the $R_n$ as belonging to
the same class when the same kind of type of linear homogeneous
groups of the $R_n$ is associated to both of them.

}

\sectionengellie{\S\,\,\,105.}

We can use the classification of all groups of the $R_n$ just
described in order to provide an approach \deutsch{Ansatz} towards the
determination of these groups. But we will only undertake this for the
\emphasis{transitive} groups.

If we imagine the variables chosen so that the origin of coordinates:
$x_1 = 0$, \dots, $x_n = 0$ is a point in general position, then
every transitive group of the $R_n$ contains, in the the neighbourhood
of the origin of coordinates, $n$ independent infinitesimal 
transformations of zeroth order in the $x_i$:
\[
T_i^{(0)}
=
\frac{\partial f}{\partial x_i}
+\cdots
\ \ \ \ \ \ \ \ \ \ \ \ \ {\scriptstyle{(i\,=\,1\,\cdots\,n)}}
\]
out of which no transformation of first or higher order can be
linearly deduced.

Furthermore, every transitive group of the $R_n$ contains in 
general also certain infinitesimal transformations of first
order in the $x_i$ that depend, according to what precede,
on the class to which the group belongs. Now, since a completely
determined class of transitive groups of the $R_n$ is 
associated to every type of linear homogeneous groups of this
space, we want to choose any such type and to restrict ourselves
to the consideration of those transitive groups which
belong to the corresponding class.

Let the $m_1$-term group:
\def\theequation{3}\begin{equation}
\sum_{i,\,\,\nu}^{1\cdots\,n}\,
\alpha_{ji\nu}\,x_i'\,
\frac{\partial f}{\partial x_\nu'}
\ \ \ \ \ \ \ \ \ \ \ \ \ 
{\scriptstyle{(j\,=\,1\,\cdots\,m_1\,;\,\,\,
0\,\leqslant\,m_1\,\leqslant\,n^2)}}
\end{equation}
be a representative of the chosen type of linear homogeneous groups.
Then according to the above, every transitive group of the $R_n$
which belongs to the corresponding class can, by means of an 
appropriate choice of variables, be brought to a form such that
in the neighbourhood of: $x_1 = 0$, \dots, $x_n = 0$, 
it contains the following $m_1$ independent infinitesimal
transformations of first order:
\[
T_j^{(1)}
=
\sum_{\nu=1}^n\,
\bigg\{
\sum_{i=1}^n\,\alpha_{ji\nu}\,x_i
+\cdots
\bigg\}\,
\frac{\partial f}{\partial x_\nu}
\ \ \ \ \ \ \ \ \ \ \ \ \ 
{\scriptstyle{(j\,=\,1\,\cdots\,m_1)}}.
\]
These $m_1$ infinitesimal transformations $T_j^{ (1)}$ are
constituted in such a
way that out of them, no transformation of
second or higher order can be linearly deduced and on the other
hand, such that in every first order infinitesimal transformation
of the group, the terms of first order can be linearly deduced
from the terms of first order in $T_1^{ (1)}, \dots, T_{ m_1}^{
(1)}$.

Besides, we already see now that in all circumstances, there is
at least one transitive group which belongs to the class
chosen by us; indeed, \emphasis{one} such group can be
immediately indicated, namely the $( n + m_1)$-term group:
\[
\frac{\partial f}{\partial x_1},
\,\,\dots,\,\,
\frac{\partial f}{\partial x_n},
\ \ \ \ \ \ \ \ \ \ \
\sum_{i,\,\,\nu}^{1\cdots\,n}\,
\alpha_{ji\nu}\,x_i\,
\frac{\partial f}{\partial x_\nu}
\ \ \ \ \ \ \ \ \ \ \ \ \ 
{\scriptstyle{(j\,=\,1\,\cdots\,m_1)}},
\]
and it is obtained by leaving out all terms of first, second, 
higher orders in the $P_i$ and in the $T_j^{ (1)}$.

\renewcommand{\thefootnote}{\fnsymbol{footnote}}
Aside from the infinitesimal transformations of zeroth and first order
already indicated, every group which belongs to our class can contain,
in the neighbourhood of: $x_1 = 0$, \dots, $x_n = 0$, yet a certain
number of independent infinitesimal transformations of second order,
out of which no transformation of third or higher order can be
linearly deduced, and moreover, a certain number of infinitesimal
transformations of third, fourth, \dots orders; but according to
Theorem~29, p.~\pageref{Theorem-29-S-192}, there always is an entire
number $s \geqslant 1$ characteristic to the group of such a nature
that the group contains infinitesimal transformations of second,
third, \dots, $s$-th orders, while by contrast, it contains no
transformations of orders $(s+1)$ or higher.\footnote[1]{\,
What we denoted by $s$ in Theorem~29 is called here $s + 1$.
} 
From this, it follows that the totality of all groups of our class
decomposes in a series of subclasses: to each value of $s$ there
corresponds a subclass.
\renewcommand{\thefootnote}{\arabic{footnote}}

Since there are infinitely many entire numbers $s$ which are
$\geqslant 1$, the number of the subclasses just defined is infinitely
large, though every subclass needs absolutely not be represented
effectively by a group. We will show how it can be decided for every
individual value of $s$ whether some groups belong to the concerned
subclass; for all that, the question fundamentally is about values of
$s$ larger than $1$ only, for we already know that the subclass: $s =
1$ contains some groups.

Let $s_0 \geqslant 1$ be an arbitrarily chosen but completely
determined entire number; we ask whether there are groups in our class
which belong to the subclass: $s = s_0$.

If there are groups of this sort, let for instance $G$ be one of them.
If, for $k = 0, 1, 2, \dots, s_0$, we leave out from the infinitesimal
transformations of order $k$ of $G$ all terms of orders $(k+1)$ and
higher, then we obviously obtain independent infinitesimal
transformations which generate a \label{S-607}
certain group $\Gamma$ and to be
precise, a group which belongs to the subclass $s = s_0$ just as $G$.

From this, it follows that the subclass: $s = s_0$, as soon as it
actually comprises groups, contains at least one group $\Gamma$
having the following specific constitution: $\Gamma$ is generated
by the $n + m_1$ infinitesimal transformations of zeroth and 
of first orders:
\[
\aligned
&
{\sf T}_i^{(0)}
=
\frac{\partial f}{\partial x_i},
\ \ \ \ \ \ \
{\sf T}_j^{(1)}
=
\sum_{i,\,\,\nu}^{1\cdots\,n}\,
\alpha_{ji\nu}\,x_i\,
\frac{\partial f}{\partial x_\nu}
\\
&
\ \ \ \ \ \ \ \ \ \ \ \ \ 
{\scriptstyle{(i\,=\,1\,\cdots\,n\,;\,\,\,
j\,=\,1\,\cdots\,m_1)}},
\endaligned
\]
and in addition, by $m_2$ independent infinitesimal transformations
of second order, by $m_3$ independent ones of third order, 
\dots, by $m_{ s_0}$ independent ones of $s_0$-th order; the general 
form of these transformations is:
\[
\aligned
&
{\sf T}_{i_k}^{(k)}
=
\sum_{\nu=1}^n\,
\xi_{i_k,\,\nu}^{(k)}
(x_1,\dots,x_n)\,
\frac{\partial f}{\partial x_\nu}
\\
&
\ \ \ \ \ \
{\scriptstyle{(k\,=\,2,\,3\,\cdots\,s_0\,;\,\,\,
i_k\,=\,1,\,2\,\cdots\,m_k)}},
\endaligned
\]
where the $\xi^{ (k)}$ are entire $k$-th order homogeneous functions
of their arguments. At the same time, it results that to every group
$G$ which belongs to the subclass $s = s_0$ is associated a completely
determined group $\Gamma$ having the constitution just described.

By executable operations, it can be decided whether there is a group
$\Gamma$ which possesses the properties just described; by means of
executable operations, one can even determine all possibly existing
groups $\Gamma$.

In fact, the number of all possible systems: $m_2$, $m_3$, \dots, $m_{
s_0}$ is at first finite. Furthermore, if one has chosen such a
system, one can always determine, in the most general way and by means
of algebraic operations, $m_2 + \cdots + m_{ s_0}$ independent
infinitesimal transformations:
\[
{\sf T}_{i_k}^{(k)}
\ \ \ \ \ \ \ \ \ \ \ \ \
{\scriptstyle{(k\,=\,2,\,3\,\cdots\,s_0\,;\,\,\,
i_k\,=\,1,\,2\,\cdots\,m_k)}}
\] 
which, together with the ${\sf T}_i^{ (0)}$, ${\sf T}_j^{ (1)}$
generate an $(n + m_1 + \cdots + m_{ s_0})$-term group.  Indeed, for
that, one only needs to determine in the most general way the
coefficients in the functions $\xi^{ (k)}$ so that every
transformation:
\[
\big\leftbracket
{\sf T}_{i_k}^{(k)},\,\,
{\sf T}_{j_\mu}^{(\mu)}
\big\rightbracket
\ \ \ \ \ \ \ \ \ \ \ \ \
{\scriptstyle{(k,\,\,\mu\,=\,2,\,3\,\cdots\,s_0)}}
\]
may be linearly deduced from:
\[
{\sf T}_\pi^{(k+\mu-1)}
\ \ \ \ \ \ \ \ \ \ \ \ \ 
{\scriptstyle{(\pi\,=\,1\,\cdots\,m_{k+\mu-1})}},
\]
as soon as $k + \mu - 1 \leqslant s_0$, and vanishes identically as
soon as $k + \mu - 1 > s_0$. It is clear on the occasion that one
obtains only algebraic equations for the unknown coefficients.

\smallskip

In what precedes, it is shown that for every \emphasis{individual}
entire number $s_0 \geqslant 1$, it can be realized by means of 
executable operations whether there are groups of our class which
belong to the subclass: $s = s_0$. But we possess no general method
which accomplishes this for \emphasis{all} entire numbers $s > 1$
\emphasis{in one stroke}. Only in special cases, for special
constitutions of the linear homogeneous group~\thetag{ 3} did
we succeed to recognize how many and which ones of the
infinitely many subclasses are represented by some groups. 
On the occasion, it happens that a maximum exists for the number
$s$, so that only the classes whose number $s$ does not 
exceed a certain maximum really contain some groups (cf. the
Chap.~\ref{kapitel-29}); however, for the existence of 
such a maximum, we do not have a general criterion. Nevertheless,
we believe that it is possible to set up such a criterion.

In consequence of that, we will restrict ourselves to explain
how one can find all groups of our class which belong to a 
determined subclass, say in the subclass: $s = s_0 \geqslant 1$. 

According to p.~\pageref{S-607}, to every group of the subclass:
$s = s_0$ is associated a completely determined group $\Gamma$ of
the same subclass with infinitesimal transformations which have a
specific form described above. Now, since all groups
$\Gamma$ of this sort which belong to the subclass: $s = s_0$
can be determined, as we know, 
by means of executable operations, we yet need
only to show in which way one can find the groups of the subclass: 
$s = s_0$ to which are associated an arbitrarily chosen concerned
group $\Gamma$.

Let:
\def\theequation{4}\begin{equation}
\left\{
\aligned
{\sf T}_i^{(0)}
=
\frac{\partial f}{\partial x_i},
\ \ \ \ \ \ \ \
{\sf T}_j^{(1)}
&
=
\sum_{\nu,\,\,\pi}^{1\cdots\,n}\,
\alpha_{j\nu\pi}\,x_\nu\,
\frac{\partial f}{\partial x_\pi}
=
\sum_{\pi=1}^n\,\xi_{j\pi}^{(1)}(x_1,\dots,x_n)\,
\frac{\partial f}{\partial x_\pi}
\\
{\sf T}_{i_k}^{(k)}
&
=
\sum_{\nu=1}^n\,\xi_{i_k,\,\nu}^{(k)}(x_1,\dots,x_n)\,
\frac{\partial f}{\partial x_\nu}
\\
&
\!\!\!\!\!\!\!\!\!\!\!\!\!\!\!\!\!\!
{\scriptstyle{(i\,=\,1\,\cdots\,n\,;\,\,\,
j=\,1\,\cdots\,m_1\,;\,\,\,
i_k\,=\,1\,\cdots\,m_k\,;\,\,\,
k\,=\,2\,\cdots\,s_0)}}
\endaligned\right.
\end{equation}
be the $n + m_1 + \cdots + m_{ s_O}$ independent infinitesimal
transformations of an arbitrary group amongst the discussed
groups $\Gamma$, and let the composition of this 
group be determined by the relations:
\def\theequation{5}\begin{equation}
\aligned
&
\big\leftbracket
{\sf T}_{i_k}^{(k)},\,\,
{\sf T}_{j_\mu}^{(\mu)}
\big\rightbracket
=
\sum_{\pi=1}^{m_{k+\mu-1}}\,
c_{i_k\,j_\mu\,\pi}\,
{\sf T}_\pi^{(k+\mu-1)}
\\
&
{\scriptstyle{(k,\,\,\mu\,=\,0,\,\,1\,\cdots\,s_0\,;\,\,\,
i_k\,=\,1\,\cdots\,m_k\,;\,\,\,
j_\mu\,=\,1\,\cdots\,m_\mu\,;\,\,\,
m_0\,=\,n)}},
\endaligned
\end{equation}
in which the $c$ of the right-hand side are to be considered as
known and in particular, do vanish as soon as $k + \mu - 1$
exceeds the number $s_0$.

Every group $\mathfrak{ G}$ belonging to the subclass: $s = s_0$
to which is associated the group~\thetag{ 4} contains 
$n + m_1 + \cdots + m_{ s_0}$ parameters and is generated
by the same number of independent infinitesimal transformations; 
these transformations have the form:
\def\theequation{6}\begin{equation}
\aligned
T_i^{(0)}
=
\frac{\partial f}{\partial x_i}
&
+\cdots,
\ \ \ \ \ \ \
T_{i_k}^{(k)}
=
\sum_{\nu=1}^n\,\xi_{i_k,\,\nu}^{(k)}
(x_1,\dots,x_n)\,
\frac{\partial f}{\partial x_\nu}
+\cdots
\\
&
\ \ \
{\scriptstyle{(i\,=\,1\,\cdots\,n\,;\,\,\,
i_k\,=\,1\,\cdots\,m_k\,;\,\,\,
k\,=\,1\,\cdots\,s_0)}},
\endaligned
\end{equation}
where, generally, the left out terms are of higher order than those
written. \emphasis{The question is nothing but to determine all
groups, the infinitesimal transformations of which have the 
form~\thetag{ 6}.}

Evidently, 
the composition of a group with the infinitesimal 
transformations~\thetag{ 6} is represented by relations of
the form:
\def\theequation{7}\begin{equation}
\aligned
&
\left\{
\aligned
\big\leftbracket
T_{i_k}^{(k)},\,\,
T_{j_\mu}^{(\mu)}
\big\rightbracket
&
=
\sum_{\pi=1}^{m_{k+\mu-1}}\,
c_{i_k\,j_\mu\,\pi}\,\,
T_\pi^{(k+\mu-1)}
\\
&\ \ \ \ \
+
\sum_{\tau=k+\mu}^{s_0}\,
\sum_{\pi_\tau}^{1\cdots\,m_\tau}\,
C_{i_k\,j_\mu\,\pi_\tau}\,
T_{\pi_\tau}^{(\tau)}
\endaligned\right.
\\
&
\ \ \ \ \ \
{\scriptstyle{(k,\,\,\mu\,=\,0,\,\,1\,\cdots\,s_0\,;\,\,\,
i_k\,=\,1\,\cdots\,m_k\,;\,\,\,
j_\mu\,=\,1\,\cdots\,m_\mu\,;\,\,\,
m_0\,=\,n)}}.
\endaligned
\end{equation}
Here, the $C$ are certain constants which, as is known, satisfy
relations derived from the Jacobi identity
(cf. Chap.~\ref{kapitel-9}, p.~\pageref{S-169-sq} sq.).

We imagine that the concerned relations between the $C$ are set
up and that the most general system of $C$ which satisfies them
is computed. Since all the relations are algebraic, this computation
requires only executable operations. In addition, the form of the 
relations~\thetag{ 7} can be simplified by replacing every 
infinitesimal transformation of $k$-th order: $T_{ i_k}^{ (k)}$
by another transformation of order $k$:
\def\theequation{8}\begin{equation}
\label{S-610}
\aligned
&
\mathfrak{T}_{i_k}^{(k)}
=
T_{i_k}^{(k)}
+
\sum_{\tau=k+1}^{s_0}\,
\sum_{\pi_\tau}^{1\cdots\,m_\tau}\,
{\sf P}_{i_k\,\pi_\tau}\,
T_{\pi_\tau}^{(\tau)}
\\
&
\ \ \ \ \ \ \ \ \ \ \ 
{\scriptstyle{(i_k\,=\,1\,\cdots\,m_k\,;\,\,\,
k\,=\,0,\,\,1\,\cdots\,s_0)}},
\endaligned
\end{equation}
where it is understood that the ${\sf P}$ are arbitrary numerical
quantities. Indeed, if one introduces the $\mathfrak{ T}$ in 
place of the $T$, one obtains in place of~\thetag{ 7} relations
of the form:
\def\theequation{7'}\begin{equation}
\left\{
\aligned
\big\leftbracket
\mathfrak{T}_{i_k}^{(k)},\,\,
\mathfrak{T}_{j_\mu}^{(\mu)}
\big\rightbracket
&
=
\sum_{\pi=1}^{m_{k+\mu-1}}\,
c_{i_k\,j_\mu\,\pi}\,
\mathfrak{T}_\pi^{(k+\mu-1)}
\\
&
\ \ \ \ \
+
\sum_{\tau=k+\mu}^{s_0}\,
\sum_{\pi_\tau}^{1\cdots\,m_\tau}\,
\mathfrak{C}_{i_k\,j_\mu\,\pi_\tau}\,
\mathfrak{T}_{\pi_\tau}^{(\tau)},
\endaligned\right.
\end{equation}
where, between the $\mathfrak{ C}$ and the $C$, a connection holds
which can be easily indicated. Now, one will provide the ${\sf P}$, 
which are perfectly arbitrary in 
order that the coefficients $\mathfrak{ C}$
receive the simplest possible numerical values; thanks to this, 
one achieves a certain simplification.

If, in the equations~\thetag{ 8}, all ${\sf P}$ are chosen fixed, 
then we say that the infinitesimal transformation of $k$-th 
order: $\mathfrak{ T}_{ i_k}^{ (k)}$ is \emphasis{normalized}.

If one knows all systems of $C$ which satisfy the relations
mentioned a short while ago, then one therefore knows at the same time
all compositions that a group of the form~\thetag{ 6} can possibly 
have. Still, the question is only whether, for each of the so
defined compositions, there are groups of the form~\thetag{ 6}
which have precisely the composition in question, and how
these groups can be found, in case they exist. 

Thus, we imagine that an arbitrary system of values $C$
is given in the equations~\thetag{ 7} which satisfies
the relations discussed, so that the equations~\thetag{ 7}
represent a possible composition of an $(n + m_1 + \cdots + 
m_{ s_0})$-term group (cf. 
p.~\pageref{Proposition-1-S-297}).

At first, we realize easily that if there actually are groups
of the form~\thetag{ 6} which have the composition~\thetag{ 7}, 
then they all are similar to each other, whence they all belong
to the same type of transitive groups of the space $x_1, \dots, x_n$
(Chap.~\ref{kapitel-22}, p.~\pageref{S-434-sq}).

In fact, if we have two groups $\mathfrak{ G}$ and $\mathfrak{ G}'$
which both possess the form~\thetag{ 6} and the composition~\thetag{
7}, then these two groups can obviously be related to each other
in a holoedrically isomorphic way thanks to a choice of their
infinitesimal transformations so that the largest subgroup
of $\mathfrak{ G}$ which leaves invariant the point: $x_1 = 0$, \dots,
$x_n = 0$ corresponds to the largest subgroup of $\mathfrak{ G}'$
which fixes this point. Now, since under the assumptions made, the 
point: $x_1 = 0$, \dots, $x_n = 0$ is a point in general position, 
then according to Theorem~76, p.~\pageref{Theorem-76-S-425}, these
two \emphasis{transitive} group $\mathfrak{ G}$ and $\mathfrak{ G}'$
are similar to each other. But this is what was to be proved.

Moreover, we will show that there always are groups of the 
form~\thetag{ 6} which possess the composition~\thetag{ 7}.

We determine right at the front in a space $R_N$ of $N = n + m_1 
+ \cdots + m_{ s_0}$ dimensions a simply transitive group:
\def\theequation{9}\begin{equation}
W_if,
\ \ \ \ \
W_{i_k}^{(k)}f
\ \ \ \ \ \ \ \ \ \ \ \ \ 
{\scriptstyle{(i\,=\,1\,\cdots\,n\,;\,\,\,
k\,=\,1\,\cdots\,s_0\,;\,\,\,
i_k\,=\,1\,\cdots\,m_k)}}
\end{equation}
of the composition~\thetag{ 7}; according to Chap.~\ref{kapitel-22}, 
pp.~\pageref{S-429-sq}--\pageref{S-432}, this is always possible
and this requires at most the integration of ordinary differential 
equations. Afterwards, we choose any manifold $\mathfrak{ M}$
of the $R_N$ which admits the $(N - n)$-term subgroup:
\def\theequation{10}\begin{equation}
W_{i_k}^{(k)}f
\ \ \ \ \ \ \ \ \ \ \ \ \ 
{\scriptstyle{(k\,=\,1\,\cdots\,s_0\,;\,\,\,
i_k\,=\,1\,\cdots\,m_k)}},
\end{equation}
of the group~\thetag{ 9},
but no larger subgroup; this property is possessed by every 
characteristic manifold (p.~\pageref{S-101}) of the $(N-n)$-term
complete system:
\[
W_{i_k}^{(k)}f
=
0
\ \ \ \ \ \ \ \ \ \ \ \ \ 
{\scriptstyle{(k\,=\,1\,\cdots\,s_0\,;\,\,\,
i_k\,=\,1\,\cdots\,m_k)}}.
\]
Through the $\infty^N$ transformations of the group~\thetag{ 9},
$\mathfrak{ M}$ takes precisely $\infty^n$ different positions 
whose totality forms an invariant family. If we characterize
the individual manifolds of this family by means of $n$ coordinates
$x_1, \dots, x_n$, we obtain a transitive group in $x_1, \dots, 
x_n$:
\def\theequation{11}\begin{equation}
X_i^{(0)}f,
\ \ \ \ \
X_{i_k}^{(k)}f
\ \ \ \ \ \ \ \ \ \ \ \ \ 
{\scriptstyle{(i\,=\,1\,\cdots\,n\,;\,\,\,
k\,=\,1\,\cdots\,s_0\,;\,\,\,
i_k\,=\,1\,\cdots\,m_k)}}
\end{equation}
that indicates in which way the manifolds of our invariant family
are permuted with each other by the group~\thetag{ 9}. This new
group is isomorphic with the group~\thetag{ 9} so that relations
of the form:
\def\theequation{7''}\begin{equation}
\aligned
\big\leftbracket
&
X_{i_k}^{(k)},\,\,
X_{j_\mu}^{(\mu)}
\big\rightbracket
=
\sum_{\pi=1}^{m_{k+\mu-1}}\,
c_{i_k\,j_\mu\,\pi}\,
X_\pi^{(k+\mu-1)}f
+
\sum_{\tau=k+\mu}^{s_0}\,
\sum_{\pi_\tau}^{1\cdots\,m_\tau}\,
C_{i_k\,j_\mu\,\pi_\tau}\,
X_{\pi_\tau}^{(\tau)}f
\\
&
\ \ \ \ \ \ \ \ \ \ \ \ \ \ \ \ \ \ 
{\scriptstyle{(k,\,\,\mu\,=\,0,\,\,1\,\cdots\,s_0\,;\,\,\,
i_k\,=\,1\,\cdots\,m_k\,;\,\,\,
j_\mu\,=\,1\,\cdots\,m_\mu\,;\,\,\,
m_0\,=\,n)}}.
\endaligned
\end{equation}
hold (cf. Theorem~85, p.~\pageref{Theorem-85-S-483})

It can be easily proved that the infinitesimal transformations of the 
group~\thetag{ 11} receive the form~\thetag{ 6} after an appropriate
choice of the variables $x_1, \dots, x_n$. 

In order to conduct the concerned proof, we imagine above all that
the variables are chosen so that $\mathfrak{ M}$ receives the 
coordinates: $x_1 = 0$, \dots, $x_n = 0$. Then obviously, all
infinitesimal transformations:
\[
\ \ \ \ \
X_{i_k}^{(k)}f
\ \ \ \ \ \ \ \ \ \ \ \ \ 
{\scriptstyle{(k\,=\,1\,\cdots\,s_0\,;\,\,\,
i_k\,=\,1\,\cdots\,m_k)}}
\]
in the $x_\nu$ are of first or higher order, while by contrast:
$X_1^0f, \dots, X_n^0f$ are independent infinitesimal transformations
of zeroth order, out of which no transformation of first or
higher order can be linearly deduced; this last fact follows from
the transitivity of the group~\thetag{ 11}. It is therefore clear
that, notwithstanding our assumption just made, we can choose
the variables $x_1, \dots, x_n$ so that $X_1^0f, \dots, X_n^0f$ receive
the form:
\[
X_i^{(0)}f
=
\frac{\partial f}{\partial x_i}
+\cdots
\ \ \ \ \ \ \ \ \ \ \ \ \ {\scriptstyle{(i\,=\,1\,\cdots\,n)}}.
\]

After these preparations, we want at first to determine the
initial terms in the power series expansions of the $m_1$ 
infinitesimal transformations $X_j^{ (1)}f$. 

We have:
\[
X_j^{(1)}f
=
\sum_{\nu=1}^n\,\zeta_{j\nu}^{(1)}
(x_1,\dots,x_n)\,
\frac{\partial f}{\partial x_\nu}
+\cdots
\ \ \ \ \ \ \ \ \ \ \ \ \ {\scriptstyle{(j\,=\,1\,\cdots\,m_1)}},
\]
where it is understood that the $\zeta^{ (1)}$ are linear homogeneous
functions of their arguments. If we insert this expression in the 
$n$ relations:
\[
\big\leftbracket
X_i^{(0)},\,\,X_j^{(1)}
\big\rightbracket
=
\sum_{\pi=1}^n\,c_{ij\pi}\,X_\pi^{(0)}f
+
\sum_{\tau=1}^{s_0}\,
\sum_{\pi_\tau}^{1\cdots\,m_\tau}\,
C_{ij\pi_\tau}\,
X_{\pi_\tau}^{(\tau)}f
\ \ \ \ \ \ \ \ \ \ \ \ \ {\scriptstyle{(i\,=\,1\,\cdots\,n)}}
\]
and if we compare the terms of zeroth order in the two sides, it
then comes:
\[
\sum_{\nu=1}^n\,
\frac{\partial\zeta_{j\nu}^{(1)}}{\partial x_i}\,
\frac{\partial f}{\partial x_\nu}
=
\sum_{\pi=1}^n\,c_{ij\pi}\,
\frac{\partial f}{\partial x_\pi}
\ \ \ \ \ \ \ \ \ \ \ \ \ {\scriptstyle{(i\,=\,1\,\cdots\,n)}},
\]
hence all first order differential quotients of $\zeta_{
j\nu}^{ (1)}$ are completely determined, and because of the known
equation:
\[
\sum_{i=1}^n\,x_i\,
\frac{\partial\zeta_{j\nu}^{(1)}}{\partial x_i}
=
\zeta_{j\nu}^{(1)}, 
\]
also $\zeta_{ \nu}^{ (1)}$ itself is completely determined.
Now, since $\xi_{ j\nu}^{ (1)}$ obviously satisfies all differential
equations which we have found just now, it results that: $\zeta_{
j\nu}^{ (1)} = \xi_{ j\nu}^{ (1)}$, and consequently:
\[
X_j^{(1)}f
=
\sum_{\nu=1}^n\,\xi_{j\nu}^{(1)}
(x_1,\dots,x_n)\,
\frac{\partial f}{\partial x_\nu}
+\cdots
\ \ \ \ \ \ \ \ \ \ \ \ \ {\scriptstyle{(j\,=\,1\,\cdots\,m_1)}}.
\]

In the same way, we find generally:
\[
X_{i_k}^{(k)}f
=
\sum_{\nu=1}^n\,
\xi_{i_k\,\nu}^{(k)}
(x_1,\dots,x_n)\,
\frac{\partial f}{\partial x_\nu}
+\cdots
\ \ \ \ \ \ \ \ \ \ \ \ \ 
{\scriptstyle{(k\,=\,1\,\cdots\,s_0\,;\,\,\,
i_k\,=\,1\,\cdots\,m_k)}},
\]
or in other words: we find that the infinitesimal transformations
of the group~\thetag{ 11} have the form~\thetag{ 6}. Evidently, 
here lies the reason why the $N$ infinitesimal transformations~\thetag{
11} are mutually independent, so that the group~\thetag{ 11} is
$N$-term and has the same composition as the group~\thetag{ 9}.

\smallercharacters{

Since, as we have seen just now, the group~\thetag{ 11} is $N$-term, 
the group:
\def\theequation{9}\begin{equation}
W_i^{(0)}f,
\ \ \ \ \
W_{i_k}^{(k)}f
\ \ \ \ \ \ \ \ \ \ \ \ \ 
{\scriptstyle{(i\,=\,1\,\cdots\,n\,;\,\,\,
k\,=\,1\,\cdots\,s_0\,;\,\,\,
i_k\,=\,1\,\cdots\,m_k)}}
\end{equation}
can contain no invariant subgroup which belongs to the group:
\def\theequation{10}\begin{equation}
W_{i_k}^{(k)}f
\ \ \ \ \ \ \ \ \ \ \ \ \ 
{\scriptstyle{(k\,=\,1\,\cdots\,s_0\,;\,\,\,
i_k\,=\,1\,\cdots\,m_k)}}
\end{equation}
(Theorem~85, p.~\pageref{Theorem-85-S-483}). But one can also
realize directly that there is no such invariant subgroup of 
the group~\thetag{ 9}; in this way, by taking the Theorem~85 into
consideration, one obtains a new proof of the fact that the 
group~\thetag{ 11} has $N$ essential parameters.

Every invariant subgroup of~\thetag{ 9} which belongs to the 
group~\thetag{ 10} obviously contains an infinitesimal transformation:
\[
\mathfrak{W}f
=
\sum_{\tau=1}^{m_p}\,
\rho_\tau\,W_\tau^{(p)}f
+
\sum_{k=p+1}^{s_0}\,
\sum_{\tau_k}^{1\cdots\,m_k}\,
\sigma_{\tau_k}\,
W_{\tau_k}^{(k)}f
\ \ \ \ \ \ \ \ \ \ \ \ \ 
{\scriptstyle{(p\,\geqslant\,1)}},
\]
in which not all the $m_p$ quantities $\rho_\tau$ vanish. 
Simultaneously with $\mathfrak{ W}f$, there appear in $\mathfrak{ g}$
also the $n$ infinitesimal transformations:
\[
\leftbracket
W_1^0,\,\mathfrak{W}
\rightbracket,
\,\,\cdots,\,\,
\leftbracket
W_n^0,\,\mathfrak{W}
\rightbracket\,;
\]
consequently, as one sees from the composition~\thetag{ 7} of the
group~\thetag{ 9}, there surely appears in $\mathfrak{ g}$
a transformation of the form:
\[
\mathfrak{W}'f
=
\sum_{\tau=1}^{m_p-1}\,
\rho_\tau'\,
W_\tau^{(p-1)}f
+
\sum_{k=p}^{s_0}\,
\sum_{\tau_k}^{1\cdots\,m_k}\,
\sigma_{\tau_k}'\,
W_{\tau_k}^{(k)}f,
\]
in which not all the $m_{ p-1}$ quantities $\rho_\tau'$ are zero.  In
this way, one realizes finally that $\mathfrak{ g}$ must contain an
infinitesimal transformation which does not belong to the
group~\thetag{ 10}; but this contradicts the assumption that
$\mathfrak{ g}$ should be contained in the group~\thetag{ 10}, hence
there is no group $\mathfrak{ g}$ having the supposed constitution.

}

The developments achieved up to now provide the

\renewcommand{\thefootnote}{\fnsymbol{footnote}}
\def\thetheorem{110}\begin{theorem}
\label{Theorem-110-S-614}
In the variables $x_1, \dots, x_n$, let a transitive group be
presented which, in the neighbourhood of the point: $x_1 = 0$, \dots,
$x_n = 0$ in general position contains the following infinitesimal
transformations: firstly, $n$ zeroth order independent transformations
of the form:
\[
{\sf T}_1^{(0)}
=
\frac{\partial f}{\partial x_1},
\,\,\,\dots,\,\,\,
{\sf T}_n^{(0)}
=
\frac{\partial f}{\partial x_n},
\]
and secondly, for every $k = 1, 2, \dots, s$, $m_k > 0$ independent
$k$-th order transformations of the form:
\[
{\sf T}_{i_k}^{(k)}
=
\sum_{\nu=1}^n\,
\xi_{i_k\,\nu}^{(k)}
(x_1,\dots,x_n)\,
\frac{\partial f}{\partial x_\nu}
\ \ \ \ \ \ \ \ \ \ \ \ \ 
{\scriptstyle{(k\,=\,1\,\cdots\,s\,;\,\,\,
i_k\,=\,1\,\cdots\,m_k)}},
\]
where the $\xi^{ (k)}$ denote completely homogeneous functions of
order $k$. Also, let the composition of the group be determined
by the relations:
\[
\aligned
&
\ 
\big\leftbracket
{\sf T}_{i_k}^{(k)},\,\,
{\sf T}_{j_\mu}^{(\mu)}
\big\rightbracket
=
\sum_{\pi=1}^{m_{k+\mu-1}}\,
c_{i_k\,j_\mu\,\pi}\,
{\sf T}_\pi^{(k+\mu-1)}
\\
&
{\scriptstyle{(k,\,\,\mu\,=\,0,\,\,1\,\cdots\,s\,;\,\,\,
i_k\,=\,1\,\cdots\,m_k\,;\,\,\,
j_\mu\,=\,1\,\cdots\,m_\mu\,;\,\,\,
m_0\,=\,n)}},
\endaligned
\]
where the $c$ in the right-hand side vanish all as soon as $k + \mu
- 1$ is larger than $s$. Then one finds in th following
way all $(n + m_1 + \cdots + m_s)$-term groups of the space
$x_1, \dots, x_n$, the infinitesimal transformations of
which possess the form:

\def\theequation{6}\begin{equation}
\aligned
T_i^{(0)}
=
\frac{\partial f}{\partial x_i}
&
+\cdots,
\ \ \ \ \ \ \
T_{i_k}^{(k)}
=
\sum_{\nu=1}^n\,\xi_{i_k,\,\nu}^{(k)}
(x_1,\dots,x_n)\,
\frac{\partial f}{\partial x_\nu}
+\cdots
\\
&
\ \ \
{\scriptstyle{(i\,=\,1\,\cdots\,n\,;\,\,\,
i_k\,=\,1\,\cdots\,m_k\,;\,\,\,
k\,=\,1\,\cdots\,s)}},
\endaligned
\end{equation}
in the neighbourhood of: $x_1 = 0$, \dots, $x_n = 0$.

One determines the constants $C$ in the equations:
\def\theequation{7}\begin{equation}
\aligned
&
\left\{
\aligned
\big\leftbracket
T_{i_k}^{(k)},\,\,
T_{j_\mu}^{(\mu)}
\big\rightbracket
&
=
\sum_{\pi=1}^{m_{k+\mu-1}}\,
c_{i_k\,j_\mu\,\pi}\,\,
T_\pi^{(k+\mu-1)}
\\
&\ \ \ \ \
+
\sum_{\tau=k+\mu}^{s_0}\,
\sum_{\pi_\tau}^{1\cdots\,m_\tau}\,
C_{i_k\,j_\mu\,\pi_\tau}\,
T_{\pi_\tau}^{(\tau)}
\endaligned\right.
\\
&
\ \ \ \ \ \
{\scriptstyle{(k,\,\,\mu\,=\,0,\,\,1\,\cdots\,s\,;\,\,\,
i_k\,=\,1\,\cdots\,m_k\,;\,\,\,
j_\mu\,=\,1\,\cdots\,m_\mu\,;\,\,\,
m_0\,=\,n)}}
\endaligned
\end{equation}
in the most general way so that they satisfy the relations 
following from the Jacobi identity. If this takes place, then
the equations~\thetag{ 7} represent all compositions that the
sought groups may have. To every individual composition amongst
these compositions correspond groups of the form~\thetag{ 6}
which are all similar
to each other and which can be found in any case by integrating
ordinary differential equations.\footnote[1]{\,
\name{Lie}, Archiv for Math. Vol. X, pp.~381--389. 1885.
}
\end{theorem}
\renewcommand{\thefootnote}{\arabic{footnote}}

\sectionengellie{\S\,\,\,151.}

In the preceding paragraphs, we have reduced
the problem of determining all transitive groups of the 
$R_n$ to the following four problems:

\smallskip
$\mathfrak{A}$. \emphasis{To find all types of linear homogeneous 
groups in $n$ variables.}

\smallskip
If one has solved this problem, then in the neighbourhood
of any point in general position, one knows the possible forms of the
initial terms in all infinitesimal transformations of first order
which can appear in a transitive group of the $R_n$. 

\smallskip
$\mathfrak{B}$. \emphasis{Assuming that, in the neighbourhood of a
point in general position, the initial terms of the
first order infinitesimal
transformations of a transitive group of the $R_n$ are
given, to determine all possible forms of the initial terms in the
infinitesimal transformations of second and higher order.}

\smallskip
$\mathfrak{C}$. \emphasis{To determine all compositions that a 
transitive group of the $R_n$ can have, a group which, in the 
neighbourhood of a point in general position contains certain 
infinitesimal transformations of first, second, \dots, $s$-th order
with given initial terms, and which by contrast, contains no 
transformation of $(s + 1)$-th or higher order.}

\smallskip
$\mathfrak{D}$. \emphasis{Assuming that one of the compositions
found in the preceding problem is given, to set up a transitive
group which possesses this composition and whose infinitesimal
transformations of first, \dots, $s$-th order have the form given
in the preceding problem.}

\smallskip
If one knows \emphasis{one} group having the constitution demanded
in the last problem, then one knows \emphasis{all} such groups,
since these are, according to Theorem~110, similar to each other.

The settlement of the first one amongst these four problems requires
only executable, or said more precisely: only algebraic operations
(cf. Chap.~\ref{kapitel-12}, p.~\pageref{Theorem-33-S-210} and
Chap.~\ref{kapitel-23}, p.~\pageref{S-494-sq} sq.).
By contrast, we did not succeed to reduce the problem $\mathfrak{ B}$
to a finite number of executable operations. The problem
$\mathfrak{ C}$ again requires only algebraic operations.
Lastly, as we have shown in the previous paragraph, 
the problem $\mathfrak{ D}$ can in any case be settled
by integrating ordinary differential equations.

At present, we want to consider a special case which presents
a characteristic simplification and we want to carry it out
in details.

Assume that the $(n + m_1 + \cdots + m_s)$-term group:
\def\theequation{12}\begin{equation}
\aligned
{\sf T}_i^{(0)}
&
=
\frac{\partial f}{\partial x_i},
\ \ \ \ \ \ 
{\sf T}_{i_k}^{(k)}
=
\sum_{\nu=1}^n\,
\xi_{i_k\,\nu}^{(k)}
(x_1,\dots,x_n)\,
\frac{\partial f}{\partial x_\nu}
\\
&
\ \ \ \ \ \ \ \
{\scriptstyle{(i\,=\,1\,\cdots\,n\,;\,\,\,
k\,=\,1\,\cdots\,s\,;\,\,\,
i_k\,=\,1\,\cdots\,m_k)}}
\endaligned
\end{equation}
which we imagined as given in Theorem~110, 
p.~\pageref{Theorem-110-S-614}, contains amongst its first-order
infinitesimal transformations:
\[
e_1\,{\sf T}_1^{(1)}
+\cdots+
e_{m_1}\,{\sf T}_{m_1}^{(1)}
\]
specially one of the form:
\[
x_1\,\frac{\partial f}{\partial x_1}
+\cdots+
x_n\,\frac{\partial f}{\partial x_n}
\]
and let, say, ${\sf T}_{ m_1}^{ (1)}$ have this form.  We will show
that under this assumption, it is always possible to determine all
transitive $(n + m_1 + \cdots + m_s)$-term groups of the space whose
transformations, in the neighbourhood of the point: $x_1 = 0$, \dots,
$x_n = 0$ in general position, have the form:
\def\theequation{12'}\begin{equation}
\aligned
T_i^{(0)}
&
=
\frac{\partial f}{\partial x_i}
+\cdots,
\ \ \ \ \ \ \
T_{i_k}^{(k)}
=
\sum_{\nu=1}^n\,
\xi_{i_k\,\nu}^{(k)}(x_1,\dots,x_n)\,
\frac{\partial f}{\partial x_\nu}
+\cdots
\\
&
\ \ \ \ \ \ \ \ \ \ \ \ \ \ \ \
{\scriptstyle{(i\,=\,1\,\cdots\,n\,;\,\,\,
k\,=\,1\,\cdots\,s\,;\,\,\,
i_k\,=\,1\,\cdots\,m_k)}}.
\endaligned
\end{equation}

For:
\[
T_{m_1}^{(1)}
=
\sum_{\nu=1}^n\,
x_\nu\,
\frac{\partial f}{\partial x_\nu}
+\cdots,
\]
we write $U$, a notation which we shall also employ on later occasions.

At first, it is clear that between $U$ and the infinitesimal
transformations of order $s$ of a group of the form~\thetag{ 12'}, 
the following relations hold:
\[
\big\leftbracket
T_\pi^{(s)},\,U
\big\rightbracket
=
(1-s)\,T_\pi^{(s)}
\ \ \ \ \ \ \ \ \ \ \ \ \ 
{\scriptstyle{(\pi\,=\,1\,\cdots\,m_s)}}.
\]

In the same way, between $U$ and the infinitesimal transformations
of order $(s - 1)$, there are relations of the form:
\[
\big\leftbracket
T_j^{(s-1)},\,U
\big\rightbracket
=
(2-s)\,T_j^{(s-1)}
+
\sum_{\pi=1}^{m_s}\,
K_{j\pi}\,T_\pi^{(s)}
\ \ \ \ \ \ \ \ \ \ \ \ \ 
{\scriptstyle{(j\,=\,1\,\cdots\,m_{s-1})}},
\]
where the $K$ are unknown constants. In order to simplify these
relations, we set (cf. p.~\pageref{S-610}):
\[
\mathfrak{T}_j^{(s-1)}
=
T_j^{(s-1)}
+
\sum_{\pi=1}^{m_s}\,
{\sf P}_{j\pi}\,
T_\pi^{(s)}
\ \ \ \ \ \ \ \ \ \ \ \ \ 
{\scriptstyle{(j\,=\,1\,\cdots\,m_{s-1})}},
\]
and we find:
\[
\big\leftbracket
\mathfrak{T}_j^{(s-1)},\,\,U
\big\rightbracket
=
(2-s)\,\mathfrak{T}_j^{(s-1)}
+
\sum_{\pi=1}^{m_s}\,
\big(
K_{j\pi}
+
(1-s-2+s)\,{\sf P}_{j\pi}
\big)\,T_{j\pi}^{(s)},
\]
hence when we choose the ${\sf P}$ in an appropriate way:
\[
\big\leftbracket
\mathfrak{T}_j^{(s-1)},\,\,U
\big\rightbracket
=
(2-s)\,\mathfrak{T}_j^{(s-1)}
\ \ \ \ \ \ \ \ \ \ \ \ \ 
{\scriptstyle{(j\,=\,1\,\cdots\,m_{s-1})}}.
\]
As a result, the infinitesimal transformations of order $(s - 1)$
of the group~\thetag{ 12'} are \emphasis{normalized}.

In exactly the same way, we can normalize the infinitesimal
transformations of order $(s - 2)$ by setting:
\[
\aligned
\mathfrak{T}_j^{(s-2)}
=
T_j^{(s-2)}
&
+
\sum_{\pi=1}^{m_{s-1}}\,
{\sf P}_{j\pi}'\,\mathfrak{T}_\pi^{(s-1)}
+
\sum_{\pi=1}^{m_s}\,
{\sf P}_{j\pi}''\,T_\pi^{(s)}
\\
& \ \ \ 
{\scriptstyle{(j\,=\,1\,\cdots\,m_{s-2})}},
\endaligned
\]
and dispose appropriately of the ${\sf P}'$ and of the ${\sf P}''$;
in this way, it comes:
\[
\big\leftbracket
\mathfrak{T}_j^{(s-2)},\,\,U
\big\rightbracket
=
(3-s)\,\mathfrak{T}_j^{(s-2)}
\ \ \ \ \ \ \ \ \ \ \ \ \ 
{\scriptstyle{(j\,=\,1\,\cdots\,m_{s-2})}}.
\] 
By proceeding alike, we obtain finally, when we generally write $T$ 
for $\mathfrak{ T}$:
\def\theequation{13}\begin{equation}
\big\leftbracket
T_{i_k}^{(k)},\,\,U
\big\rightbracket
=
(1-k)\,T_{i_k}^{(k)}
\ \ \ \ \ \ \ \ \ \ \ \ \ 
{\scriptstyle{(k\,=\,0,\,\,1\,\cdots\,s\,;\,\,\,
i_k\,=\,1\,\cdots\,m_k\,;\,\,\,
m_0\,=\,n)}}.
\end{equation}

As a result, the infinitesimal transformations of zeroth, first, 
up to $s$-th order are all normalized, \emphasis{except for $U$
itself}.

At present, we remember that because of the composition of the
group~\thetag{ 12'}, between the $T$, there are relations of the
form:
\def\theequation{14}\begin{equation}
\aligned
&
\left\{
\aligned
\big\leftbracket
T_{i_k}^{(k)},\,
T_{j_\mu}^{(\mu)}
\big\rightbracket
&
=
\sum_{\pi=1}^{m_{k+\mu-1}}\,
c_{i_k\,j_\mu\,\pi}\,
T_\pi^{(k+\mu-1)}
\\
&
\ \ \ \ \
+
\sum_{\tau=k+\mu}^s\,\sum_{\pi_\tau}^{1\cdots\,m_\tau}\,
C_{i_k\,j_\mu\,\pi_\tau}\,
T_{\pi_\tau}^{(\tau)}
\endaligned\right.
\\
& 
\ \ \ \ \ \ \ 
{\scriptstyle{(k,\,\,\mu\,=\,0,\,\,1\,\cdots\,s\,;\,\,\,
i_k\,=\,1\,\cdots\,m_k\,;\,\,\,
j_\mu\,=\,1\cdots\,m_\mu\,;\,\,\,
m_0\,=\,n)}},
\endaligned
\end{equation} 
in which the $c$ actually vanish as soon as $k + \mu - 1 > s$. 
In order to determine the unknown constants $C$, we form the
Jacobi identity:
\[
\big\leftbracket
\leftbracket
T_{i_k}^{(k)},\,
T_{j_\mu}^{(\mu)}
\rightbracket,\,\,
U
\big\rightbracket
+
\big\leftbracket
\leftbracket
T_{j_\mu}^{(\mu)},\,
U
\rightbracket,\,\,
T_{i_k}^{(k)}
\big\rightbracket
+
\big\leftbracket
\leftbracket
U,\,
T_{i_k}^{(k)}
\rightbracket,\,\,
T_{j_\mu}^{(\mu)}
\big\rightbracket
=
0,
\]
which, because of~\thetag{ 13}, obviously takes the form:
\[
\big\leftbracket
\leftbracket
T_{i_k}^{(k)},\,
T_{j_\mu}^{(\mu)}
\rightbracket,\,\,
U
\big\rightbracket
=
(2-k-\mu)\,
\big\leftbracket
T_{i_k}^{(k)},\,
T_{j_\mu}^{(\mu)}
\big\rightbracket.
\]
Here, if we insert the expression~\thetag{ 14} and if we use once
more the equations~\thetag{ 13}, it comes:
\[
\aligned
\sum_{\pi=1}^{m_{k+\mu-1}}\,
c_{i_k\,j_\mu\,\pi}\,(2-k-\mu)\,
T_\pi^{(k+\mu-1)}
&
+
\sum_{\tau=k+\mu}^s\,
\sum_{\pi_\tau}^{1\cdots\,m_\tau}\,
C_{i_k\,j_\mu\,\pi_\tau}\,(1-\tau)\,T_{\pi_\tau}^{(\tau)}
\\
&
=
(2-k-\mu)\,
\big\leftbracket
T_{i_k}^{(k)},\,
T_{j_\mu}^{(\mu)}
\big\rightbracket,
\endaligned
\]
or, because the $T$ are independent infinitesimal transformations:
\[
C_{i_k\,j_\mu\,\pi_\tau}\,(\tau+1-k-\mu)
=
0
\ \ \ \ \ \ \ \ \ \ \ \ \ 
{\scriptstyle{(k\,+\,\mu\,\leqslant\,\tau\,\leqslant s)}}.
\]
From this, it results that all the $C$ vanish.

Thus, under the assumption made, all groups of 
the form~\thetag{ 12'} have the same composition as the
group~\thetag{ 12}, hence according to Theorem~110, 
p.~\pageref{Theorem-110-S-614}, they all are similar to each other
and similar to the group~\thetag{ 12}. 

In consequence of that, we can say:

\def\thetheorem{111}\begin{theorem}
If a transitive $(n + m_1 + \cdots + m_s)$-term group in the variables
$x_1, \dots, x_n$ contains, in the neighbourhood of the point 
in general position: $x_1 = 0$, \dots, $x_n = 0$, aside from the 
$n$ independent infinitesimal transformations of zeroth order in the
$x$:
\[
T_i
=
\frac{\partial f}{\partial x_i}
+\cdots
\ \ \ \ \ \ \ \ \ \ \ \ \ {\scriptstyle{(i\,=\,1\,\cdots\,n)}},
\]
yet $m_k$, for $k = 1, 2, \dots, s$, independent infinitesimal
transformations of order $k$ out of which no transformation
of order $k+1$ or higher can be linearly deduced, and if it
specially contains a first order
infinitesimal transformation of the form:
\[
\sum_{\nu=1}^n\,
x_\nu\,\frac{\partial f}{\partial x_\nu}
+\cdots,
\]
then thanks to the introduction of new variables $x_1', \dots, x_n'$, 
it can be brought to the form:
\[
\aligned
T_i
=
{\sf T}_i
&
=
\frac{\partial f}{\partial x_i'},
\ \ \ \ \ \ \ 
T_{i_k}^{(k)}
=
{\sf T}_{i_k}^{(k)}
=
\sum_{\nu=1}^n\,
\xi_{i_k\,\nu}^{(k)}
(x_1',\dots,x_n')\,
\frac{\partial f}{\partial x_\nu'}
\\
&
\ \ \ \ \ \ \ \ \ \ \ \ \ \ 
{\scriptstyle{(i\,=\,1\,\cdots\,n\,;\,\,\,
k\,=\,1\,\cdots\,s\,;\,\,\,
i_k\,=\,1\,\cdots\,m_k)}}\,;
\endaligned
\]
here, the $\xi^{(k)}$ are the entire homogeneous functions of order
$k$ which determine the terms of order $k$ in the infinitesimal
transformations of order $k$:
\[
T_{i_k}^{(k)}
=
\sum_{\nu=1}^n\,\xi_{i_k\,\nu}^{(k)}(x_1,\dots,x_n)\,
\frac{\partial f}{\partial x_\nu}
+\cdots
\ \ \ \ \ \ \ \ \ \ \ \ \ 
{\scriptstyle{(k\,=\,1\,\cdots\,s\,;\,\,\,
i_k\,=\,1\,\cdots\,m_k)}}
\]
of the group.
\end{theorem}

\linestop


\chapter{Characteristic Properties of the Groups 
\\
Which are Equivalent to Certain Projective Groups}
\label{kapitel-29}
\chaptermark{Groups that are Equivalent to Certain Projective Groups}

\setcounter{footnote}{0}

\abstract*{??}

In the preceding chapter, we gave a classification of all
\emphasis{transtitive}\, groups: $X_1 f, \dots, X_r f$ of the $n$-fold
extended space $x_1, \dots, x_n$. We chose a point $x_1^0, \dots,
x_n^0$ in general position and considered all infinitesimal
transformations of the group which leave this point at rest, hence all
transformations whose power series expansion with respect to the $x_i
- x_i^0$ have the form:
\[
\sum_{i,\,\nu\,=\,1}^{n}\,
\alpha_{ji\nu}
(x_1^0,\dots,x_n^0)\cdot
(x_i-x_i^0)\,
\frac{\partial f}{\partial x_\nu}
+
\cdots
\ \ \ \ \ \ \ \ \ \
{\scriptstyle{(j\,=\,1\,,2,\,\dots)}},
\]
where the left out terms are of second or of higher order\footnote{\,
Such terms are systematically written ``$+ \cdots$'' by Engel and Lie.
}. 
Then the linear homogeneous group:
\[
L_jf
=
\sum_{i,\,\nu\,=\,1}^n\,
\alpha_{ji\nu}(x_1^0,\dots,x_n^0)\,x_i'\,
\frac{\partial f}{\partial x_\nu'}
\ \ \ \ \ \ \ \ \ \
{\scriptstyle{(j\,=\,1\,,2,\,\dots)}},
\]
showed in which way the $\infty^{ n - 1}$ line-elements $x_1 ' \!:\!
x_2 ' \!:\! \cdots \!:\! x_n'$ through the point $x_1^0, \dots, x_n^0$
are transformed by those transformations of the group: $X_1 f, \dots,
X_r f$, which leave this point invariant.

\renewcommand{\thefootnote}{\fnsymbol{footnote}}
In the present chapter, to begin with, we solve the problem of
determining all \emphasis{transitive}\, groups: $X_1 f, \dots, X_r f$ of
the space $x_1, \dots, x_n$ for which the linear homogeneous group:
$L_1 f, L_2 f, \dots$ assigned to a point in general position
coincides either with the general linear homogeneous group or with the
special linear homogeneous 
group \footnote[1]{\,
It is easy to see that a group of $R_n$ which assigns the general
or the special linear homogeneous group to a point $x_1^0, \dots,
x_n^0$ in general position, is always transitive. Indeed, in the
neighborhood of $x_1^0, \dots, x_n^0$, the group certainly comprises
an infinitesimal transformation of zeroth order in the $x_i - x_i^0$,
hence a transformation: $\sum\, \alpha_i \, p_i + \cdots$, in which
not all $\alpha_i$ are equally null. In addition, the group surely
comprises $n ( n-1)$ first order transformations of the form:
\[
(x_i-x_i^0)\,p_k
+\cdots
\ \ \ \ \ \ \ \ \ \
{\scriptstyle{(i,\,k\,=\,1\,\cdots\,n;\,\,i\,\neq\,k)}}.
\]
If one makes Combination [bracketting] of the latter with $\sum\,
\alpha_i\, p_i + \cdots$, then one recognizes that $n$ transformations
of the form: $p_1 + \cdots$, \dots, $p_n + \cdots$ do appear, whence
the group is actually transitive.
}. 
Then it comes out the curious
result that every such group: $X_1 f, \dots, X_r f$ is equivalent
either to the general projective group, or to the general linear
group, or to the special linear group of the space $x_1, \dots, x_n$.
\renewcommand{\thefootnote}{\arabic{footnote}}

In addition, we yet show in the last paragraph of the chapter that in
the space $x_1, \dots, x_n$ there is no finite continuous group which
can transfer $m > n + 2$ arbitrarily chosen points in general position
to just the same kind of $m$ points; at the same time, we show that
the general projective group and the groups that are equivalent to it
are the only groups of the space $x_1, \dots, x_n$ which can transfer
$n + 2$ arbitrarily chosen points in general position to just the same
kind of $n + 2$ points.

\sectionengellie{\S\,\,\,152.}

Every transitive group of $R_n$ comprises, in the neighborhood of the
point: $x_1^0, \dots, x_n^0$ in general position, $n$ independent
infinitesimal transformations of zeroth order in the $x_i - x_i^0$:
\[
p_1+\cdots,\ \
p_2+\cdots,\ \
\dots,\ \
p_n+\cdots,
\]
where, according to the fixation of notation on p.~555, 
$p_i$ is written in place of $\frac{ \partial f}{\partial x_i}$.

\smallskip

Now, if $G$ assigns to the point $x_1^0, \dots, x_n^0$ the general
linear homogeneous group as group: $L_1f, L_2 f, \dots$, then it
comprises in the neighborhood of $x_1^0, \dots, x_n^0$ the largest
possible number, namely $n^2$, of such infinitesimal transformations of
first order in $x_1 - x_1^0, \dots, x_n - x_n^0$, out of which no
transformation of second or of higher order can be deduced by linear
combination. These $n^2$ first order transformations 
have the form:
\[
(x_i-x_i^0)\,p_k+\cdots
\ \ \ \ \ \ \ \ \ \
{\scriptstyle{(i,\,k\,=\,1\,\cdots\,n)}}.
\]
If, on the other hand, $G$ assigns to the point $x_1^0, \dots, x_n^0$
the special linear homogeneous group, then it comprises, in the
neighborhood of the point, only $n^2 - 1$ independent first order
infinitesimal transformations out of which no transformation of second
or of higher order can be deduced by linear combination; the same
transformations have the form:
\[
\aligned
(x_i-x_i^0)\,p_k
+\cdots,&
\ \ \
(x_i-x_i^0)\,p_i-(x_k-x_k^0)\,p_k
+\cdots
\\
&
{\scriptstyle{(i,\,k\,=\,1\,\cdots\,n\,;\,\,i\,\gtrless\,k)}}.
\endaligned
\]

Therefore, when we yet choose the point $x_1^0, \dots, x_n^0$ as the
origin of coordinates, we can enunciate as follows the problem
indicated in the introduction of the chapter:

\smallskip\label{S-620}
{\it 
In the variables $x_1, \dots, x_n$, to seek all groups $X_1 f, \dots,
X_r f$, or shortly $G$, which, in the neighborhood of the point in
general position: $x_1 = 0, \dots, x_n = 0$, comprise the following
infinitesimal transformations of zeroth and of first order:
either the $n + n^2$:
\[
p_i+\cdots,\ \ \
x_i\,p_k+\cdots
\ \ \ \ \ \ \ \ \ \
{\scriptstyle{(i,\,k\,=\,1\,\cdots\,n)}};
\]
or the $n + n^2 - 1$:
\[
\aligned
p_i+\cdots,\ \ \
x_i\,p_k
&
+\cdots,\ \ \
x_i\,p_i-x_k\,p_k+\cdots
\\
&
{\scriptstyle{(i,\,k\,=\,1\,\cdots\,n\,;\,\,i\,\gtrless\,k)}}.
\endaligned
\]}

To begin with, the two cases can be treated simultaneously; one must
only, as long as possible, disregard the fact that in the first case,
aside from the transformations which occur in the second case, yet one
transformation appears: $\sum^i\, x_i\, p_i$.

The group $X_1 f, \dots, X_r f$ comprises infinitesimal
transformations whose expansion in power series begins with terms of
second, or of relatively higher order in the $x$, shortly,
infinitesimal transformations of second or higher order, respectively.
We search for the highest order number\footnote{\,
Actually, $s$ will be shown to be equal to 2, not more, and the
corresponding transformations of second order to be necessarily of the
form $x_i ( x_1\, p_1 + \cdots + x_n\, p_n) + \cdots$, $i = 1, \dots,
n$.
} 
$s$ of the existing
transformations and we even look for the determination of such
transformations.

We can suppose this number $s$ to be bigger than 1, since we already know
all the first order infinitesimal transformations that appear. Let
\[
K
=
\vartheta_1\,p_1
+\cdots+
\vartheta_n\,p_n
+\cdots
\]
be a $s$-th order infinitesimal transformation of the group; at the
same time, the $\vartheta$'s are supposed to denote completely
homogeneous functions of order $s$ in $x_1, \dots, x_n$, while the
left out terms are of higher order. Obviously, $\vartheta_1, \dots,
\vartheta_n$ do not vanish all, because there should
not be given any infinitesimal transformation $\sum\, c_k\, X_kf$
whose power series expansion starts with terms of $(s + 1)$-th or of
higher order; so we can assume in any case that $\vartheta_1$ is not
identically null.

Let the lowest power of $x_1$ which appears in $\vartheta_1$ be the
$\alpha_1$-th and let
\[
x_1^{\alpha_1}\,x_2^{\alpha_2}\cdots 
x_n^{\alpha_n}
\ \ \ \ \ \ \ \ \ \ \ \ \ \
{\scriptstyle{(\,\alpha_1\,+\,\cdots\,+\,\alpha_n\,=\,s\,)}}
\]
be a term of $\vartheta_1$ with nonvanishing coefficient. Now, we
make Combination of the transformation $x_1 \, p_2 + \cdots$ with $K$,
and we make Combination of the result once more with $x_1 \, p_2 +
\cdots$, and so on\,\,---\,\,in total $\alpha_2$ times. We make
Combination of the $s$-th order infinitesimal transformation obtained
this way with $x_1 \, p_3 + \cdots$, and we then proceed $\alpha_3$
times one after the other, etc., and finally, we still apply
$\alpha_n$ times $x_1 \, p_n + \cdots$ one after the other. In this
way, we recognize in the end that a transformation of the form:
\[
K'
=
x_1^s\,p_1+\vartheta_2'\,p_2
+\cdots+
\vartheta_n'\,p_n
+\cdots
\]
belongs to our group.

All remaining terms of $\vartheta_1$ are cancelled; indeed, the same
terms contain either also the power $x_1^{ \alpha_1}$ or a higher one,
and in both cases, the power of a variable $x_i$ ($i > 1$) is
certainly not the $\alpha_i$-th, but a lowest one; consequently, the
corresponding term vanishes by making $\alpha_i$ times Combination
with $x_1 \, p_i + \cdots$. As always, the terms of $(s + 1)$-th and
of higher order are, also here, left out of 
consideration\footnote{\,
Here is an 
equivalent reformulation.
Adapting slightly notation, let $\vartheta_1 = \sum_{ \vert \alpha
\vert = s}\, A_\alpha\, x^\alpha$ and choose $\beta \in \N^n$ with
$\vert \beta \vert = s$ so that $\beta_1 = \inf \{ \alpha_1 : \,
A_\alpha \neq 0 \}$. For any other monomial $A_\alpha \, x^\alpha$
with $A_\alpha \neq 0$ we have either $\alpha_1 > \beta_1$ or
$\alpha_1 = \beta_1$. In the first case, since $\alpha_1 + \alpha_2 +
\cdots + \alpha_i + \cdots + \alpha_n = \beta_1 + \beta_2 + \cdots +
\beta_i + \cdots + \beta_n = s$, there must exist a $i$ with $2
\leqslant i \leqslant n$ such that $\alpha_i < \beta_i$. In the second
case, namely if $\alpha_1 = \beta_1$, then $\alpha_2 + \cdots +
\alpha_i + \cdots + \alpha_n = \beta_2 + \cdots + \beta_i + \cdots +
\beta_n$, and again, there must exist a $i$ with $2 \leqslant i
\leqslant n$ such that $\alpha_i < \beta_i$, because otherwise,
$\alpha_2 \geqslant \beta_2$, \dots, $\alpha_n \geqslant \beta_n$
together with $\vert \alpha \vert = \vert \beta \vert$
implies $\alpha_2 = \beta_2$, \dots, $\alpha_n = \beta_n$, hence
$\alpha = \beta$, a contradiction.
} 
(cf. Chap.~\ref{kapitel-11}, 
Theorem~30, p.~\pageref{Theorem-30}).

In order to determine more precisely the form of the transformation
$K'$, we make use of an auxiliary proposition which can be
exploited several times with benefit.

\def\theproposition{1}\begin{proposition}
If the infinitesimal transformations $Cf$ and $B_1 f + \cdots + B_m f$
belong to a group and if furthermore, there are $m$ relations of the
form $\leftbracket C, \, B_k \rightbracket = \varepsilon_k \, B_kf$, where the constants
$\varepsilon_k$ are all distinct one another, then the group comprises
all the $m$ infinitesimal transformations 
$B_1 f, \dots, B_m f$.
\end{proposition}

The proof of this auxiliary proposition is very simple. 
Aside from $B_1 f + \cdots + B_m f$, the group yet obviously comprises
the following infinitesimal transformations:
\[
\aligned
\big\leftbracket
C,\,B_1+\cdots+B_m
\big\rightbracket
&
=
\varepsilon_1\,B_1f
+\cdots+
\varepsilon_m\,B_mf
\\
\big\leftbracket
C,\,
\varepsilon_1\,B_1
+\cdots+
\varepsilon_m\,B_m
\big\rightbracket
&
=
\varepsilon_1^2\,B_1f
+\cdots+
\varepsilon_m^2\,B_mf
\\
\cdots\cdots\cdots\cdots\cdots\cdots\cdots\cdots\cdots
&
\cdots\cdots\cdots\cdots\cdots\cdots\cdots\cdots\cdots
\\
\big\leftbracket
C,\,\,
\varepsilon_1^{m-2}\,B_1
+\cdots+
\varepsilon_m^{m-2}\,B_m
\big\rightbracket
&
=
\varepsilon_1^{m-1}\,B_1f
+\cdots+
\varepsilon_m^{m-1}\,B_mf.
\endaligned
\]
But since the determinant:
\[
\left\vert
\begin{array}{ccccc}
1 & \cdots & 1
\\
\varepsilon_1 & \cdots & \varepsilon_m
\\
\vdots & \ddots & \vdots
\\
\varepsilon_1^{m-1} & \cdots & \varepsilon_m^{m-1}
\end{array}
\right\vert
=
\prod_{i>k}\,
(\varepsilon_i-\varepsilon_k)
\]
does not vanish by assumption, each individual transformation $B_k f$
can be deduced, by multiplication with appropriate constants and by
subsequent addition, from the just found transformations. Thus the
proposition is proved.

\smallskip

In order to apply the same proposition, we now make Combination of the
transformation $K'$, written in detail:
\[
\aligned
x_1^s\,p_1
+
\big\{
{\textstyle{\sum}}\,
A_\beta\,x_1^{\beta_1}
&
\cdots x_n^{\beta_n}
\big\}\,p_2
+\cdots+
\big\{
{\textstyle{\sum}}\,
N_\nu\,x_1^{\nu_1}\cdots x_n^{\nu_n}
\big\}\,p_n
+\cdots
\\
&
{\scriptstyle{(\beta_1\,+\,\cdots\,+\,\beta_n\,=\,\cdots\,=\,
\nu_1\,+\,\cdots\,+\,\nu_n\,=\,s)}}
\endaligned
\]
with the transformation:
\[
x_1\,p_1-x_i\,p_i+\cdots,
\]
on the understanding that $i$ is any of the numbers $2, \dots,
n$. Then one obtains:
\[
\aligned
(s-1)\,
&
x_1^s\,p_1
+
\big\{
{\textstyle{\sum}}\,
A_\beta\,
(\beta_1-\beta_i+\varepsilon_{2i})\,
x_1^{\beta_1}\cdots x_n^{\beta_n}
\big\}\,p_2
+\cdots+
\\
&
+
\big\{
{\textstyle{\sum}}\,
N_\nu\,
(\nu_1-\nu_i+\varepsilon_{ni})\,
x_1^{\nu_1}\cdots x_n^{\nu_n}
\big\}\,p_n
+\cdots,
\endaligned
\]
where $\varepsilon_{ ii }$ has the value 1, while all $\varepsilon_{
ki}$ ($i \neq k$) vanish. Therefore the $s$-th order terms of the
infinitesimal transformation $K'$ have the form $B_1 f + \cdots + B_m
f$ discussed above, where the $B_k f$ are reproduced through
Combination with $x_1 \, p_1 - x_i \, p_i$, 
but with distinct factors\footnote{\,
Here are considerations about derived homogeneous groups.
The next phrase of the text applies a fundamental observation of
Chap.~28 that we must reconstitute just as a preparation. Let $G$ be
a finite continuous (locally) transitive homogeneous group acting on
the $(x_1, \dots, x_n)$-space and suppose that the origin $0$ is the
central point, of course in general position whenever needed. Then
$G$ comprises $n$ transformations of zeroth order:
\[
{\sf T}_i^{(0)}
=
\frac{\partial f}{\partial x_i}
+
\cdots
\ \ \ \ \ \ \ \ \ \
{\scriptstyle{(i\,=\,1\,\cdots\,n)}}
\]
out of which no transformation of first or of higher order can be
deduced by linear combination. Further, $G$ contains a certain
number, say $m_1$ (possibly null), of infinitesimal transformations of
first order:
\[
{\sf T}_j^{(1)}
=
\sum_{\nu=1}^n\,
\bigg\{
\sum_{i=1}^n\,\alpha_{ji\nu}\,x_i+\cdots
\bigg\}\,
\frac{\partial f}{\partial x_\nu}
\ \ \ \ \ \ \ \ \ \
{\scriptstyle{(j\,=\,1\,\cdots\,m_1)}}
\]
that are linearly independent modulo transformations of order
$\geqslant 2$. Generally, for $k = 1, 2, \dots$ up to a finite maximal
order $s_0$, the group $G$ comprises a certain number $m_k$ (possibly
null, and indeed null for $k \geqslant s_0 + 1$ by definition of
$s_0$) of infinitesimal transformations of $k$-th order:
\[
{\sf T}_j^{(k)}
=
\sum_{\nu=1}^n\,
\Big\{
\xi_{j\nu}^{(k)}(x)+\cdots
\Big\}\,
\frac{\partial f}{\partial x_\nu}
\ \ \ \ \ \ \ \ \ \
{\scriptstyle{(j\,=\,1\,\cdots\,m_k)}}
\]
that are linearly independent modulo transformations of order
$\geqslant m_k + 1$, where each $\xi_{ j\nu}^{ ( k)} (x)$ is a
homogeneous polynomial of degree $m_k$ in $x_1, \dots, x_n$, and where
the left out terms are of order $\geqslant m_k + 1$.

By bracketing (cf. Chap.~\ref{kapitel-11}, Theorem~30, 
p.~\pageref{Theorem-30}), we get:
\def\theequation{{\bf a}}\begin{equation}
\big\leftbracket
{\sf T}_{j_1}^{(k_1)},\,
{\sf T}_{j_2}^{(k_2)}
\big\rightbracket
=
\sum_{\nu=1}^n\,
\bigg\{
\xi_{j_1\mu}^{(k_1)}\,
\frac{\partial\xi_{j_2\nu}^{(k_2)}}{\partial x_\mu}
-
\xi_{j_2\mu}^{(k_2)}\,
\frac{\partial\xi_{j_1\nu}^{(k_1)}}{\partial x_\mu}
+
\cdots
\bigg\}\,
\frac{\partial f}{\partial x_\nu},
\end{equation}
where the written terms are of order $k_1 + k_2 - 1$ (they might be
null) and where the left out terms are of order $\geqslant k_1 + k_2$.
The linear combination of infinitesimal transformations to which the
bracket $\big\leftbracket {\sf T}_{ j_1}^{ (k_1)}, \, {\sf T}_{ j_2}^{ ( k_2)}
\big\rightbracket$ should be equal can only contain transformations of order
$\geqslant k_1 + k_2 - 1$, hence may be written:
\def\theequation{{\bf b}}\begin{equation}
\big\leftbracket
{\sf T}_{j_1}^{(k_1)},\,
{\sf T}_{j_2}^{(k_2)}
\big\rightbracket
=
\sum_{l=1}^{m_{k_1+k_2-1}}\,
C_{j_1j_2l}^{k_1k_2}\cdot
{\sf T}_l^{(k_1+k_2-1)}
+
\cdots,
\end{equation}
where the left out terms comprise linear combinations of
transformations of order $\geqslant k_1 + k_2$. In particular, if
$k_1 + k_2 - 1 \geqslant s_0 +1$, the brackets $\big\leftbracket {\sf T}_{ j_1}^{
(k_1)},\, {\sf T}_{ j_2}^{ (k_2)} \big\rightbracket$ must vanish.

For every $k = 0, 1, \dots, s_0$ and every $j = 1, \dots, m_k$, we
define the projected infinitesimal transformation:
\[
\widehat{\sf T}_j^{(k)}
:=
\sum_{\nu=1}^n\,
\xi_{j\nu}^{(k)}(x)\,
\frac{\partial f}{\partial x_\nu},
\]
defined by erasing all remainders,
which is homogeneous of order $m_k$. Equivalently, $\widehat{ \sf
T}_j^{ (k)} = \widehat{ \pi}_k \big( {\sf T}_j^{ (k)} \big)$ if, by
$\widehat{ \pi}_k$, we denote the projection of infinitesimal
transformations (and of analytic functions) onto the vector space of
all homogeneous monomials of degree $k$. We notice that a
reformulation of~\thetag{ {\bf a}} simply says:
\[
\widehat{\pi}_{k_1+k_2-1}
\Big(
\big\leftbracket
{\sf T}_{j_1}^{(k_1)},\,
{\sf T}_{j_2}^{(k_2)}
\big\rightbracket
\Big)
=
\Big\leftbracket
\widehat{\pi}_{k_1}\big({\sf T}_{j_1}^{(k_1)}\big),\,
\widehat{\pi}_{k_2}\big({\sf T}_{j_2}^{(k_2)}\big)
\Big\rightbracket
=
\big\leftbracket
\widehat{\sf T}_{j_1}^{(k_1)},\,
\widehat{\sf T}_{j_2}^{(k_2)}
\big\rightbracket.
\] 
Finally, by taking the $\widehat{ \pi}_{ k_1 + k_2 - 1}$-projection
of both sides of~\thetag{ {\bf b}}, we obtain Lie algebra-type identities:
\[
\aligned
\big\leftbracket
\widehat{\sf T}_{j_1}^{(k_1)},\,
\widehat{\sf T}_{j_2}^{(k_2)}
\big\rightbracket
&
=
\widehat{\pi}_{k_1+k_2-1}
\Big(
\big\leftbracket
{\sf T}_{j_1}^{(k_1)},\,
{\sf T}_{j_2}^{(k_2)}
\big\rightbracket
\Big)
\\
&
=
\widehat{\pi}_{k_1+k_2-1}
\bigg(
\sum_{l=1}^{m_{k_1+k_2-1}}\,
C_{j_1j_2l}^{k_1k_2}\,
{\sf T}_l^{k_1+k_2-1}
+\cdots
\bigg)
\\
&
=
\sum_{l=1}^{m_{k_1+k_2-1}}\,
C_{j_1j_2l}^{k_1k_2}\,
\widehat{\sf T}_l^{k_1+k_2-1}
\endaligned
\]
which show that to the arbitrary initial group $G$ is always
associated a group $\widehat{ G}$ constituted by infinitesimal
transformations $\widehat{ \sf T}_j^{ ( k)}$, $k = 0, 1, \dots, s_0$,
$j = 1, \dots, m_k$, having polynomial, homogeneous coefficients.
}. 

Now, since from the expansion in power series of our group a new
group $\Gamma$ can be derived in such a way that in each power series
expansion, only the terms of lowest order are kept (Chap.~28,
p.~607), then our Proposition~1 shows that each individual $B_k f$
belongs to the group $\Gamma$. Next, there evidently is a $B_k f$ which
embraces the term $x_1^s \, p_1$ and hence is reproduced with the
factor $s - 1$. The remaining terms of this $B_k f$ are defined
through the equations\footnote{\,
Here is an enlightement. 
The derived homogeneous group $\widehat{ G} \equiv \Gamma$ contains
the infinitesimal transformations $\widehat{ D}_i = x_1 \, p_1 - x_i
\, p_i$, for $i = 1, \dots, n$, and also:
\[
\widehat{K}'
=
x_1^s\,p_1
+
\big\{
{\textstyle{\sum}}\,
A_\beta\,x^\beta
\big\}\,p_2
+\cdots+
\big\{
{\textstyle{\sum}}\,
N_\gamma\,x^\gamma
\big\}\,p_n,
\]
all the remainders being suppressed (one must have a group for
Proposition~1 to apply). At first, for {\it fixed}\, $i$, one looks at
the Lie bracket, joint with $\widehat{ D}_i$, of each {\it monomial}\,
homogeneous infinitesimal transformation $\big\{ B_\gamma \, x^\gamma
\big\} \, p_k$, $k \geqslant 2$, $\vert \gamma \vert = s$, which appears
in $\widehat{ K}'$ and which is distinct from $x_1^s \, p_1$:
\[
\big\leftbracket
\widehat{D}_i,\,
\big\{B_\gamma\,x^\gamma\big\}\,p_k
\big\rightbracket
=
(\gamma_1-\gamma_i+\varepsilon_{ki})\,
\big\{B_\gamma\,x^\gamma\big\}\,p_k.
\]
One then collects all the monomials $\big\{ B_\gamma\,
x^\gamma\big\}\, p_k$ of $\widehat{ K}'$ having the {\it same}\,
reproducing factor $\gamma_1 - \gamma_i + \varepsilon_{ ki} = s-1$
{\it as}\, the infinitesimal transformation $x_1^s \, p_1$, and one
calls $\widehat{ B}_1^i$ the corresponding sum (it depends on $i$),
which is a part of $\widehat{ K}'$. On the other hand, there are
finitely many {\it other}\, values of the integers $\gamma_1 -
\gamma_i + \varepsilon_{ ki}$, say $m ( i)$, and one decomposes
accordingly $\widehat{ K}' = \widehat{ B}_1^i + \widehat{ B}_2^i +
\cdots + \widehat{ B}_{ m( i)}^i$. Then Proposition~1 yields that
$\widehat{ B}_1^i$ belongs to $\widehat{ G}$. Other $\widehat{ B}_1^k$
are left out.

Let $i' \neq i$ be another integer and consider bracketing with
$\widehat{ D}_{ i'}$. The same reasoning applied to $\widehat{ B}_1^i$
(instead of $\widehat{ K}'$) yields a decomposition $\widehat{ B}_1^i
= \widehat{ B}_{ 1, 1}^{ i, i'} + \widehat{ B}_{ 1, 2}^{ i, i'} +
\cdots + \widehat{ B}_{ 1, m( i')}^{ i, i'}$ with the first term
$\widehat{ B}_{ 1, 1}^{ i, i'}$ collecting monomials of $\widehat{
B}_1^i$ that are reproduced with the factor $s - 1$; clearly, $x_1^s\,
p_1$ still belongs to $\widehat{ B}_{ 1, 1}^{ i, i'}$. Proposition~1,
yields again that $\widehat{ B}_{ 1, 1}^{ i, i'}$ belongs to
$\widehat{ G}$. By induction, letting $i = 1$, $i' = 2$, \dots\,, $i^{
( n)} = n$, one gets an infinitesimal transformation $\widehat{ B}_{
1, 1, \dots, 1}^{ 1, 2, \dots, n}$ of $\widehat{ G}$\,\,---\,\,denoted
by $B_k f$ in the translated text\,\,---\,\, such that $\big\leftbracket
\widehat{ D}_i, \, \widehat{ B}_{ 1, 1, \dots, 1}^{ 1, 2, \dots, n}
\big\rightbracket = (s-1) \, \widehat{ B}_{ 1, 1, \dots, 1}^{ 1, 2, \dots, n}$,
{\it for all}\, $i = 1, 2, \dots, n$.
}: 
\[
\beta_1-\beta_i+\varepsilon_{2i}
=\cdots=
\nu_1-\nu_i+\varepsilon_{ni}
=
s-1.
\]

Since $\varepsilon_{ 23}, \dots, \varepsilon_{ 2n}$ vanish, from the
same equations, we obtain immediately $\beta_3 = \cdots = \beta_n$ and
consequently:
\[
\beta_1+\beta_2+\cdots+\beta_n
=
\beta_1+\beta_2+(n-2)\,\beta_3
=
s,
\]
whence it yet comes:
\[
\beta_2=\beta_1-s+2,
\ \ \ \ \ \ \ \ \ \ \ \ \ \ \
\beta_3=\beta_1-s+1.
\]
By elimination of $\beta_2$ and $\beta_3$, it follows:
\[
\aligned
&
(\beta_1-s+1)\,n=0
\\
\text{\rm hence}\ \ \ \ \
\beta_1=s-1,
&
\ \ \ \ \ \ \
\beta_2=1,
\ \ \ \ \ \ \
\beta_3=\cdots=\beta_n=0.
\endaligned
\]
In the same way, the $\nu_k$ determine themselves, and so on. In
brief, we realize that our group: $X_1 f, \dots, X_r f$ comprises an
infinitesimal transformation of the form:
\[
K''
=
x_1^s\,p_1
+
A_2\,x_1^{s-1}\,x_2\,p_2
+\cdots+
A_n\,x_1^{s-1}\,x_n\,p_n
+\cdots.
\]

By making Combination of $K''$ with $p_1 + \cdots$, we get:
\[
\aligned
s\,x_1^{s-1}\,p_1
&
+
(s-1)\,A_2\,x_1^{s-2}\,x_2\,p_2
+
\\
&
+\cdots+
(s-1)\,A_n\,x_1^{s-2}\,x_n\,p_n
+\cdots
=
L,
\endaligned
\]
and hence our group comprises an infinitesimal transformation, namely
$\big\leftbracket L, \, K''\big\rightbracket$, which possesses the form:
\[
s\,x_1^{2s-2}\,p_1
+
\eta_2\,p_2
+\cdots+
\eta_n\,p_n
+\cdots.
\]
But now, since $2 s - 2$ should not be larger than $s$ and since, on
the other hand, $s$ is larger than 1, it follows that:
\[
s
=
2,
\]
so that we have:
\[
K''
=
x_1^2\,p_1
+
A_2\,x_1\,x_2\,p_2
+\cdots+
A_n\,x_1\,x_n\,p_n
+\cdots.
\]

Furthermore, it comes:
\[
\big\leftbracket
x_1\,p_i+\cdots,\,
K''
\big\rightbracket
=
(A_i-1)\,x_1^2\,p_i
+\cdots.
\]
If now $A_i$ were different from 1, then we would obtain, one after
the other, both transformations:
\[
\aligned
&
\big\leftbracket
x_1^2\,p_i+\cdots,\,
x_i\,p_1+\cdots
\big\rightbracket
=
x_1^2\,p_1
-
2x_1x_i\,p_i
+\cdots,
\\
&
\big\leftbracket
x_1^2\,p_1-2x_1x_i\,p_i+\cdots,\,
x_1^2\,p_i+\cdots
\big\rightbracket
=
4x_1^3\,p_i
+\cdots.
\endaligned
\]
But since no infinitesimal transformation of third order should be
found, all the $A_i$ are equal to 1. Therefore $K''$ has the form:
\[
x_1\,\big(
x_1\,p_1+x_2\,p_2
+\cdots+
x_n\,p_n
\big)
+\cdots
\]
and thereupon lastly, by making Combination with $x_i \, p_1 +
\cdots$, we find generally:
\[
x_i\,
\big(
x_1\,p_1
+\cdots+
x_n\,p_n
\big)
+\cdots.
\]

\plainstatement{Consequently, 
if a group $X_1 f, \dots, X_r f$ of the nature required
on p.~\pageref{S-620} comprises infinitesimal transformations of order
higher than the first, then it comprises only such transformations
that are of order two, and in fact in any case, $n$ of the form: }
\[
x_i\,
\big(
x_1\,p_1
+\cdots+
x_n\,p_n
\big)
+\cdots
\ \ \ \ \ \ \ \ \ \
{\scriptstyle{(i\,=\,1\,\cdots\,n)}}.
\]

If we add up all the $n$ transformations:
\[
\Big\leftbracket
p_i+\cdots,\,
x_i\,{\textstyle{\sum_{k=1}^n}}\,x_k\,p_k
+\cdots
\Big\rightbracket
=
x_i\,p_i
+
{\textstyle{\sum_{k=1}^n}}\,x_k\,p_k
+\cdots,
\]
together, then we obtain the transformation $x_1\, p_1 + \cdots +
x_n\, p_n + \cdots$. Hence if the group $X_1 f, \dots, X_r f$ contains
infinitesimal transformations of second order, then the associated
linear homogeneous group $L_1 f, L_2 f, \dots$ is the general linear
homogeneous group. Or conversely:

\plainstatement{The 
group $X_1 f, \dots, X_r f$ never contains infinitesimal
transformations of order higher than the first,
when the associated group $L_1 f, L_2 f, \dots$ is 
the special linear homogeneous group}.

For abbreviation, we write the infinitesimal transformation $x_i \,
\big( x_1 \, p_1 + \cdots + x_n\, p_n \big) + \cdots$ in the form $H_i
+ \cdots$, on the understanding that $H_i$ denotes the terms of second
order. Now, it could be thinkable that in one group $X_1 f, \dots,
X_r f$, except the $n$ transformations $H_k + \cdots$, there still
appeared others of second order. Let such a one be for instance:
\[
\tau_1\,p_1
+\cdots+
\tau_n\,p_n
+\cdots
=
{\sf T}
+\cdots,
\]
where the $\tau_i$ mean homogeneous functions of second order in $x_1,
\dots, x_n$ and ${\sf T}$ the sum $\sum \, \tau_k \, p_k$. Then the
expression $\big\leftbracket H_i, \, {\sf T} \big\rightbracket$ must obviously vanish, since
the lowest term of this bracket represents an infinitesimal
transformation of third order. Consequently, $\tau_1, \dots, \tau_n$
satisfy the equation:
\[
\bigg\leftbracket
x_i\,\sum_{k=1}^n\,x_k\,p_k,\,
\sum_{j=1}^n\,\tau_j\,p_j
\bigg\rightbracket
=
0;
\]
this equation decomposes into the following $n$ equations:
\[
x_i\,
\bigg(
\sum_{k=1}^n\,x_k\,
\frac{\partial\tau_j}{\partial x_k}
-
\tau_j
\bigg)
-
x_j\,\tau_i
=
0
\]
and from this it follows, aside from the natural equation:
\[
\sum_{k=1}^n\,x_k\,\frac{\partial\tau_j}{\partial x_k}
=
2\tau_j
\]
yet also $x_i \, \tau_j - x_j\, \tau_i = 0$.
Therefore the $\tau_i$ and ${\sf T}$ have the form:
\[
\tau_i
=
\sum_{k=1}^n\,\alpha_k\,x_kx_i,
\ \ \ \ \ \ \ \ \ \ \ \ \ \ \ \ \ \ \ \ \ \ \ \
{\sf T}
=
\sum_{k=1}^n\,\alpha_k\,H_k.
\]

Thus we have proved that a group: $X_1 f, \dots, X_r f$ of the
indicated constitution can comprise no second order infinitesimal
transformation apart from the $n$ transformations $H_k + \cdots$.

In total, we therefore have the following cases:

\label{S-625}
\plainstatement{If 
the linear homogeneous group $L_1f, L_2f, \dots$ is the
general linear homogeneous group, then the concerned group $X_1 f,
\dots, X_r f$ comprises exactly $n$ infinitesimal transformations of
zeroth order:
\[
P_i
=
p_i
+\cdots
\ \ \ \ \ \ \ \ \ \
{\scriptstyle{(i\,=\,1\,\cdots\,n)}}
\]
and $n^2$ of first order:
\[
T_{ik}
=
x_i\,p_k
+\cdots
\ \ \ \ \ \ \ \ \ \
{\scriptstyle{(i,\,k\,=\,1\,\cdots\,n)}}.
\]
Either transformations of higher order do not occur at all, or there are
extant the following $n$:
\[
S_i
=
x_i\,\big(
x_1\,p_1+\cdots+x_n\,p_n
\big)
+\cdots
\ \ \ \ \ \ \ \ \ \
{\scriptstyle{(i\,=\,1\,\cdots\,n)}}.
\]

If the group: $L_1 f, L_2 f, \dots$ is the special linear homogeneous
group, then the group: $X_1 f, \dots, X_r f$ comprises exactly $n$
infinitesimal transformations of zeroth order:
\[
P_i
=
p_i
+\cdots
\ \ \ \ \ \ \ \ \ \
{\scriptstyle{(i\,=\,1\,\cdots\,n)}}
\]
and in addition, $n^2 - 1$ of first order:
\[
T_{ik}
=
x_i\,p_k
+\cdots,
\ \ \ \ \
T_{ii}-T_{kk}
=
x_i\,p_i-x_k\,p_k
+\cdots
\ \ \ \ \ \ \ 
{\scriptstyle{(i\,\gtrless\,k\,=\,1\,\cdots\,n)}}.
\]
but no one of higher order. }

We will treat these three cases one after the other. At first, we
bring the Relations [brackets] between the infinitesimal
transformations, and afterwards, even these transformations
themselves, to a form which is as simple as possible. On the occasion,
we remark that among the three cases, the first two are already
finished off by the developments of the paragraph~151 on p.~616 sq.
Nonetheless, we maintain that it is advisable to also treat in detail
these two cases.

\sectionengellie{\sf\S\,\,\,153.}

The first case, where $n$ transformations of zeroth order
and $n^2$ of first order do appear, is the simplest one.

We can indicate without effort the Relations between the $n^2$
infinitesimal transformations. They are:
\[
\big\leftbracket
T_{ik},\,T_{\nu\pi}
\big\rightbracket
=
\varepsilon_{k\nu}\,T_{i\pi}
-
\varepsilon_{\pi i}\,T_{\nu k},
\]
where $\varepsilon_{ ik}$ vanishes as soon as $i$ and $k$ are
distinct, whereas $\varepsilon_{ ii}$ has the value 1. In particular,
it is of importance that each expression:
\[
\bigg\leftbracket
T_{ik},\,\,
\sum_{\nu=1}^n\,T_{\nu\nu}
\bigg\rightbracket
\]
vanishes. Further, if, for reasons of abbreviation, we introduce the
symbol:
\[
U
=
\sum_{\nu=1}^n\,T_{\nu\nu},
\]
there exist Relations of the form:
\[
\big\leftbracket
P_i,\,U
\big\rightbracket
=
P_i
+
\sum_{\nu=1}^n\,\sum_{\pi=1}^n\,
\alpha_{\nu\pi}\,T_{\nu\pi},
\]
or, when we introduce as new $P_i$ the right-hand side:
$\big\leftbracket P_i, \, U \big\rightbracket = P_i$.

Lastly, there are Relations of the form:
\[
\aligned
\big\leftbracket
P_i,\,T_{kj}
\big\rightbracket
&
=
\varepsilon_{ik}\,P_j
+
\sum\,\sum\,
\beta_{\nu\pi}\,T_{\nu\pi},
\\
\big\leftbracket
P_i,\,P_k
\big\rightbracket
&
=
\sum\,
\gamma_\nu\,P_\nu
+
\sum\,\sum\,
\delta_{\nu\pi}\,T_{\nu\pi}.
\endaligned
\]
But the Jacobian identities:
\[
\aligned
\big\leftbracket\big\leftbracket P_i,\,T_{kj}\big\rightbracket,\,U\big\rightbracket
-
\big\leftbracket P_i,\,T_{kj}\big\rightbracket
&
=
0,
\\
\big\leftbracket\big\leftbracket P_i,\,P_k\big\rightbracket,\,U\big\rightbracket
-
2\,\big\leftbracket P_i,\,P_k\big\rightbracket
&
=
0
\endaligned
\]
show immediately that all constants $\beta$, $\gamma$, $\delta$ do
vanish, whence it is valid that:
\[
\big\leftbracket P_i,\,T_{kj}\big\rightbracket
=
\varepsilon_{ik}\,P_j,
\ \ \ \ \ \ \ \ \ \ \ \ \ \ \ \
\big\leftbracket P_i,\,P_k\big\rightbracket
=
0.
\]

As a result, all the Relations between the infinitesimal transformations
of our group are known.

The $r$ infinitesimal transformations $P_i = p_i + \cdots$ generate a
simply transitive group which has 
\terminology{the same composition} \footnote{\,
Namely here: they both have the same, in fact
vanishing, Lie brackets.
} 
\deutsch{gleichzusammengesetz ist} 
as the group $p_1', \dots, p_n'$
and is hence also equivalent to it
(Chap.~19, p.~339, 
Prop.~1\footnote{\,
The so-called \terminology{Frobenius theorem}\, is needed in the
phrase just below. Classically (cf. \cite{ stk2000}), one performs a
preliminary reduction to a commuting system of vector fields, reducing
the proof to the following concrete

\smallskip\noindent{\bf Proposition.}
{\em If the $r \leqslant s$ independent
infinitesimal transformations:
\[
X_kf
=
\sum_{i=1}^s\,\xi_{ki}(x_1,\dots,x_s)\,
\frac{\partial f}{\partial x_i}
\ \ \ \ \ \ \ \ \ \
{\scriptstyle{(k\,=\,1\,\cdots\,r)}}
\]
stand in the Relationships:
\[
\big\leftbracket X_i,\,X_k\big\rightbracket
=
0
\ \ \ \ \ \ \ \ \ \
{\scriptstyle{(i,\,k\,=\,1\,\cdots\,r)}}
\]
without though being tied up together through a linear relation of the
form:
\[
\sum_{k=1}^n\,\chi_k(x_1,\dots,x_s)\,X_sf
=0,
\]
then they generate an $r$-term group which is equivalent to the group
of translations: }
\[
Y_1f
=
\frac{\partial f}{\partial y_1},
\dots,\,
Y_rf
=
\frac{\partial f}{\partial y_r}.
\]

\noindent
(Namely, in a neighborhood of a generic point, there is a local
diffeomorphism $x \mapsto y = y (x)$ straightening each $X_k$ to $Y_k
= \frac{ \partial }{ \partial y_k}$.)
}). 
So we can introduce such new variables
$x_1', \dots, x_n'$ achieving that:
\[
P_i
=
p_i'
\ \ \ \ \ \ \ \ \ \
{\scriptstyle{(i\,=\,1\,\cdots\,n)}}.
\]
The form $\xi_1 \, p_1' + \cdots + \xi_n \, p_n'$ which $U$ takes in
the new variables determines itself from the Relations $\big\leftbracket P_i, \,
U \big\rightbracket = P_i$; the same Relations yield:
\[
U
=
\sum_{k=1}^n\,
(x_k'+\alpha_k)\,p_k'
\]
and, when $x_k' + \alpha_k$ is introduced as new $x_k$, which does not
change the form of the $P_i$, one achieves $U = \sum\, x_k\,
p_k$. From the Relations:
\[
\big\leftbracket
P_i,\,T_{kj}
\big\rightbracket
=
\varepsilon_{ik}\,P_j
\]
one finds in the same way:
\[
T_{kj}
=
x_k\,p_j
+
\sum_{\nu=1}^n\,
\alpha_{kj\nu}\,p_\nu;
\]
but since $\big\leftbracket T_{ kj}, \, U \big\rightbracket$ must vanish, all the $\alpha_{
kj\nu}$ are equally null. Therefore we have the group:
\[
P_i
=
p_i,
\ \ \ \ \ \ \ \ \ \
T_{ik}
=
x_i\,p_k
\ \ \ \ \ \ \ \ \ \
{\scriptstyle{(i,\,k\,=\,1\,\cdots\,n)}},
\]
that is to say, {\it all groups which belong to the first case are
equivalent to the general linear group of the manifold $x_1, \dots,
x_n$}.

\sectionengellie{\S\,\,\,154.}

We now come to the second case, where, apart from the $n$
transformations of zeroth order $P_i = p_i + \cdots$ and the $n^2$ of
first order $T_{ ik} = x_i\, p_k + \cdots$, there still appear the $n$
transformations of second order:
\[
S_i
=
x_i\,\big(x_1\,p_1+\cdots+x_n\,p_n)
+\cdots.
\]

The following Relations are obtained without any effort:
\[
\big\leftbracket S_i,\,S_k\big\rightbracket
=
0,
\ \ \ \ \ \ \ \ 
\big\leftbracket T_{ik},\,S_j\big\rightbracket
=
\varepsilon_{kj}\,S_i,
\ \ \ \ \ \ \ \
\big\leftbracket U,\,S_j\big\rightbracket
=
S_j.
\]
Moreover, there are equations of the form:
\[
\big\leftbracket
T_{ik},\,U
\big\rightbracket
=
\sum_{j=1}^n\,\alpha_{ikj}\,S_j.
\]
If at first we suppose that $i$ and $k$ are not all both equal to $n$,
then we can introduce $T_{ ik} + \sum\, \alpha_{ ikj}\, S_j$ as new
$T_{ ik}$ and obtain correspondingly the Relation: $\big\leftbracket T_{ ik},\, U
\big\rightbracket = 0$ for the concerned values of $i$ and $k$. Since in addition
$\big\leftbracket U, \, U \big\rightbracket = 0$ it comes generally:
\[
\big\leftbracket T_{ik},\,U\big\rightbracket
=
0
\ \ \ \ \ \ \ \ \ \
{\scriptstyle{(i,\,k\,=\,1\,\cdots\,n)}}.
\]

Further, one has:
\[
\big\leftbracket
T_{ik},\,T_{\nu\pi}
\big\rightbracket
=
\varepsilon_{k\nu}\,T_{i\pi}
-
\varepsilon_{\pi i}\,T_{\nu k}
+
\sum_{j=1}^n\,\beta_j\,S_j;
\]
but the identity:
\[
\big\leftbracket\big\leftbracket T_{ik},\,T_{\nu\pi}\big\rightbracket,\,U\big\rightbracket
=
0
\]
enables to recognize that all $\beta_j$ vanish, so it holds:\[
\big\leftbracket
T_{ik},\,T_{\nu\pi}
\big\rightbracket
=
\varepsilon_{k\nu}\,T_{i\pi}
-
\varepsilon_{\pi i}\,T_{\nu k}.
\]
From the relation:
\[
\big\leftbracket
P_i,\,U
\big\rightbracket
=
P_i
+
\sum_{\nu=1}^n\,\sum_{\pi=1}^n\,
\gamma_{\nu\pi}\,T_{\nu\pi}
+
\sum_{\nu=1}^n\,\delta_\nu\,S_\nu
\]
when:
\[
P_i
+
\sum_{\nu=1}^n\,\sum_{\pi=1}^n\,
\gamma_{\nu\pi}\,T_{\nu\pi}
+
\frac{1}{2}\,\sum_{\nu=1}^n\,\delta_\nu\,S_\nu
\]
is introduced as new $P_i$, it comes:
\[
\big\leftbracket P_i,\,U\big\rightbracket
=
P_i.
\]

Furthermore, one has:
\[
\big\leftbracket P_i,\,S_k\big\rightbracket
=
\varepsilon_{ik}\,U
+
T_{ki}
+
\sum_{\nu=1}^n\,\lambda_\nu\,S_\nu;
\]
but if we form the identity:
\[
\big\leftbracket\big\leftbracket P_i,\,S_k\big\rightbracket,\,U\big\rightbracket
+
\big\leftbracket P_i,\,S_k\big\rightbracket
-
\big\leftbracket P_i,\,S_k\big\rightbracket
=
0,
\]
we then find that the $\lambda_\nu$ are null; consequently we have:
\[
\big\leftbracket
P_i,\,S_k
\big\rightbracket
=
\varepsilon_{ik}\,U
+
T_{ki}.
\]

Finally, there are relations of the form:
\[
\aligned
\big\leftbracket
P_i,\,T_{kj}
\big\rightbracket
&
=
\varepsilon_{ik}\,P_j
+
\sum\,\sum\,\alpha_{\nu\pi}\,T_{\nu\pi}
+
\sum\,\lambda_\nu\,S_\nu,
\\
\big\leftbracket
P_i,\,P_k
\big\rightbracket
&
=
\sum\,g_\nu\,P_\nu
+
\sum\,\sum\,h_{\nu\pi}\,T_{\nu\pi}
+
\sum\,l_\nu\,S_\nu;
\endaligned
\]
but also here all constants do vanish because of the identities:
\[
\aligned
\big\leftbracket\big\leftbracket P_i,\,T_{kj}\big\rightbracket,\,U\big\rightbracket
-
\big\leftbracket P_i,\,T_{kj}\big\rightbracket
&
=
0,
\\
\big\leftbracket\big\leftbracket P_i,\,P_k\big\rightbracket,\,U\big\rightbracket
-
2\,\big\leftbracket P_i,\,P_k\big\rightbracket
&
=
0.
\endaligned
\]
Therefore it comes:
\[
\big\leftbracket P_i,\,T_{kj}\big\rightbracket
=
\varepsilon_{ik}\,P_j,
\ \ \ \ \ \ \ \ \ \ \ \ \ \
\big\leftbracket P_i,\,P_k\big\rightbracket
=
0.
\]

All the relations between the infinitesimal transformations of our
group are now determined. One will notice that the infinitesimal
transformations $P_i$, $T_{ ik}$ generate a subgroup which possesses
the form considered in the first case and hence takes the form:
\[
P_i
=
p_i,
\ \ \ \ \ \ \ \ \ \ \ \ \ \ \
T_{ik}
=
x_i\,p_k
\ \ \ \ \ \ \ \ \ \
{\scriptstyle{(i,\,k\,=\,1\,\cdots\,n)}}
\]
by an appropriate choice of the variables. In the new
variables, the infinitesimal transformations $S_i$ are, say, equal
to $\sum\, \xi_{ ik}\, p_k$, where the $\xi_{ ik}$ satisfy the
relations:
\[
\big\leftbracket
P_\nu,\,S_i
\big\rightbracket
=
\varepsilon_{\nu i}\,U
+
T_{i\nu},
\ \ \ \ \ \ \ \ \ \ \ \ \ 
\big\leftbracket U,\,S_i\big\rightbracket
=
S_i.
\]
We find from this:
\[
\frac{\partial\xi_{ik}}{\partial x_\nu}
=
\varepsilon_{\nu i}\,x_k
+
\varepsilon_{\nu k}\,x_i,
\ \ \ \ \ \ \ \ \ \ \ \ \ \ \ \ \ \ \
\sum_{\nu=1}^n\,x_\nu\,
\frac{\partial\xi_{ik}}{\partial x_\nu}
=
2\,\xi_{ik},
\]
therefore $\xi_{ ik} = x_i x_k$ and:
\[
S_i
=
x_i\big(x_1\,p_1
+\cdots+
x_n\,p_n
\big).
\]

Consequently we have the group:
\[
P_i
=
p_i,
\ \ \ \ 
T_{ik}
=
x_i\,p_k,
\ \ \ \
S_i
=
x_i\big(
x_1\,p_1
+\cdots+
x_n\,p_n
\big)
\ \ \ \ \ \ \ 
{\scriptstyle{(i,\,k\,=\,1\,\cdots\,n)}};
\]
that is to say, \emphasis{all groups which belong to the
second case are equivalent to the general projective group of the
manifold $x_1, \dots, x_n$}.

\bigskip

\centerline{\sf\S\,\,\,155.}

\medskip

It remains the third and last case with $n$ infinitesimal
transformations of zeroth order: $P_i = p_i + \cdots$ and $n^2 - 1$
of first order:
\[
T_{ik}
=
x_i\,p_k
+\cdots,
\ \ \ \ 
T_{ii}-T_{kk}
=
x_i\,p_i-x_k\,p_k
+\cdots
\ \ \ \ \ \ \
{\scriptstyle{(i,\,\gtrless\,k\,=\,1\,\cdots\,n)}},
\]
while transformations of higher order do not occur.

It is convenient to replace the infinitesimal transformations $T_{ ii}
- T_{ kk}$ by one of the form: $\alpha_1 \, T_{ 11} + \cdots +
\alpha_n\, T_{ nn}$, thinking that the $\alpha_i$ are arbitrary
constants, though subjected to the condition $\sum_i\, \alpha_i = 0$.
Then to begin with, the following Relations hold:
\[
\big\leftbracket
T_{ik},\,T_{\nu\pi}
\big\rightbracket
=
\varepsilon_{k\nu}\,T_{i\pi}
-
\varepsilon_{\pi i}\,T_{\nu k}
\ \ \ \ \ \ \ \ \ \
{\scriptstyle{(i\,\gtrless\,k,\,\,\,\nu\,\gtrless\,\pi)}},
\]
\[
\bigg\leftbracket
T_{ik},\,\,
\sum_{\nu=1}^n\,\alpha_\nu\,T_{\nu\nu}
\bigg\rightbracket
=
(\alpha_k-\alpha_i)\,T_{ik},
\ \ \ \ \ \ 
\bigg\leftbracket
T_{ii}-T_{kk},\,\,
\sum_{\nu=1}^n\,
\alpha_\nu\,T_{\nu\nu}
\bigg\rightbracket
=
0.
\]
Further, there is an equation of the form:
\[
\bigg\leftbracket
P_i,\,\,
\sum_{\nu=1}^n\,\alpha_\nu\,T_{\nu\nu}
\bigg\rightbracket
=
\alpha_i\,P_i
+
\sum_{k=1}^n\,\sum_{j=1}^n\,
\lambda_{ikj}\,T_{kj}
\ \ \ \ \ \ \ \ \ \ \ \ \ \ \ \
{\scriptstyle{(\sum_{k=1}^n\,\lambda_{ikk}\,=\,0)}}.
\]
Therefore if we set:
\[
P_i'
=
P_i
+
\sum_{k=1}^n\,\sum_{j=1}^n\,
l_{ikj}\,T_{kj}
\ \ \ \ \ \ \ \ \ \ \ \ \ \ \ \
{\scriptstyle{(\sum_{k=1}^n\,l_{ikk}\,=\,0)}},
\]
it comes:
\[
\bigg\leftbracket
P_i',\,\,
\sum_{\nu=1}^n\,\alpha_\nu\,T_{\nu\nu}
\bigg\rightbracket
=
\alpha_i\,P_i'
+
\sum_{k=1}^n\,\sum_{j=1}^n\,
\Big\{
\lambda_{ikj}-(\alpha_i+\alpha_k-\alpha_j)
\Big\}\,
T_{kj}.
\]
Now, we imagine that a completely determined system of values is
chosen for $\alpha_1, \dots, \alpha_n$ which satisfies the condition
$\sum\, \alpha_i = 0$ and has in addition the property that no
expression $\alpha_i + \alpha_k - \alpha_j$ vanishes, a demand which
can always be satisfied. Under these assumptions, we can choose the
$l_{ ikj}$ so that the equations:
\[
\lambda_{ikj}
-
(\alpha_i+\alpha_k-\alpha_j)\,l_{ikj}
=
0
\]
are satisfied, where the condition $\sum_k \, l_{ ikk} = 0$ is
automatically fulfilled. Consequently, we obtain after utilizing again
the initial designation:
\[
\bigg\leftbracket
P_i,\,\,
\sum_{\nu=1}^n\,\alpha_\nu\,T_{\nu\nu}
\bigg\rightbracket
=
\alpha_i\,P_i
\ \ \ \ \ \ \ \ \ \ \ \ \ \ \ \
{\scriptstyle{(\sum_{\nu=1}^n\,\alpha_\nu\,=\,0)}}.
\]

If the expression:
\[
\aligned
\bigg\leftbracket
P_i,\,\,
\sum_{k=1}^n\,\beta_k\,T_{kk}
&
\bigg\rightbracket
=
\beta_i\,P_i
+
\sum_{\nu=1}^n\,\sum_{\pi=1}^n\,
g_{\nu\pi}\,T_{\nu\pi}
\\
&
{\scriptstyle{(\sum\,\beta_k\,=\,\sum\,g_{\nu\nu}\,=\,0)}}
\endaligned
\]
is inserted into the identity:
\[
\bigg\leftbracket\Big\leftbracket
P_i,\,\,
\sum_{k=1}^n\,\beta_k\,T_{kk}
\Big\rightbracket,\,\,
\sum_{k=1}^n\,\alpha_k\,T_{kk}
\bigg\rightbracket
-
\alpha_i\,
\bigg\leftbracket
P_i,\,\,
\sum_{k=1}^n\,\beta_k\,T_{kk}
\bigg\rightbracket
=
0,
\]
then it comes:
\[
\sum_{\nu=1}^n\,\sum_{\pi=1}^n\,
(\alpha_\pi-\alpha_\nu-\alpha_i)\,
g_{\nu\pi}\,T_{\nu\pi}
=
0.
\]
Because of the nature of the $\alpha$, it follows from this that all
$g_{ \nu \pi}$ are equally null, and therefore the equation:
\[
\bigg\leftbracket
P_i,\,\,
\sum_{k=1}^n\,\beta_k\,T_{kk}
\bigg\rightbracket
=
\beta_i\,T_i
\]
holds for all systems of values $\beta_1, \dots, \beta_n$ which
satisfy the condition $\sum\, \beta_k = 0$. Furthermore, there is a
relation of the form:
\[
\aligned
\big\leftbracket P_i,\,T_{kj}\big\rightbracket
&
=
\varepsilon_{ik}\,P_j
+
\sum_{\nu=1}^n\,\sum_{\pi=1}^n\,
h_{\nu\pi}\,T_{\nu\pi}
\\
&
\ \ \
{\scriptstyle{(k\,\gtrless\,j,\,\,\sum_\nu\,\,h_{\nu\nu}\,=\,0)}}.
\endaligned
\]
The identity:
\[
\bigg\leftbracket
\big\leftbracket P_i,\,T_{kj}\big\rightbracket,\,\,
\sum_{\tau=1}^n\,\beta_\tau\,T_{\tau\tau}
\bigg\rightbracket
-
(\beta_i+\beta_j-\beta_k)\,
\big\leftbracket P_i,\,T_{kj}\big\rightbracket
=
0
\]
therefore takes the form:
\[
\varepsilon_{ik}\,(\beta_k-\beta_i)\,P_j
+
\sum_{\nu=1}^n\,\sum_{\pi=1}^n\,
\big(
\beta_\pi-\beta_\nu-\beta_i+\beta_k-\beta_j
\big)\,
h_{\nu\pi}\,T_{\nu\pi}
=
0;
\]
consequently the $h_{ \nu \pi}$ must vanish, since the $\beta_\nu$,
while disregarding the condition $\sum_\nu \, \beta_\nu = 0$, are
completely arbitrary:
\[
\big\leftbracket
P_i,\,
T_{kj}
\big\rightbracket
=
\varepsilon_{ik}\,P_j
\ \ \ \ \ \ \ \ \ \ \ \ \ \
{\scriptstyle{(k\,\gtrless\,j)}}.
\]

Finally, the relations:
\[
\big\leftbracket
P_i,\,P_k
\big\rightbracket
=
\sum_{\nu=1}^n\,m_\nu\,P_\nu
+
\sum_{\nu=1}^n\,\sum_{\pi=1}^n\,
m_{\nu\pi}\,T_{\nu\pi}
\ \ \ \ \ \ \ \ \ \ \ \ \ \
{\scriptstyle{(\sum_{\nu=1}^n\,\,m_{\nu\nu}\,=\,0)}}
\]
are still to be examined. By calculating the identity:
\[
\bigg\leftbracket
\big\leftbracket P_i,\,P_k\big\rightbracket,\,\,
\sum_{\tau=1}^n\,\beta_\tau\,T_{\tau\tau}
\bigg\rightbracket
-
(\beta_i+\beta_k)\,
\big\leftbracket P_i,\,P_k\big\rightbracket
=
0,
\]
we find:
\begin{small}
\[
\sum_{\nu=1}^n\,
(\beta_\nu-\beta_i-\beta_k)\,m_\nu\,P_\nu
+
\sum_{\nu=1}^n\,\sum_{\pi=1}^n\,
(\beta_\pi-\beta_\nu-\beta_i-\beta_k)\,m_{\nu\pi}\,T_{\nu\pi}
=
0;
\]
\end{small}

\noindent
because of the arbitrariness of the $\beta_\nu$, all the $m_\nu$ and
$m_{ \nu\pi}$ must be equally null. Consequently, we have:
\[
\big\leftbracket P_i,\,P_k\big\rightbracket
=
0
\]
and we therefore know all Relations between the infinitesimal
transformations of our group.

In the same way as it has been achieved in the first case, we bring
$P_1, \dots, P_n$ by means of an appropriate choice of variables to
the form:
\[
P_i
=
p_i
\ \ \ \ \ \ \ \ \ \
{\scriptstyle{(i\,=\,1\,\cdots\,n)}}.
\]
From this, by proceeding as in the end of the paragraph 153, we
conclude that \emphasis{all groups which belong to the third case are
equivalent to the special linear group of the manifold $x_1, \dots,
x_n$}.

By unifying the found results we therefore obtain the

\renewcommand{\thefootnote}{\fnsymbol{footnote}}
\def\thetheorem{112}\begin{theorem}
\label{Theorem-I-112}
If a transitive group $G$ in $n$ variables is constituted so that all
of its transformations which leave invariant a point in general
position do transform the line-elements passing through the point by
means of the general or of the special linear homogeneous group $L_1
f, L_2 f, \dots$, then $G$ is equivalent either to the general
projective group, or to the general linear group, or to the special
linear group, in $n$ 
variables\footnote[1]{\,
Lie, Archiv for Math., Vol.~3, Christiania 1878; cf. also
Math. Ann. Vol.~XVI, Vol.~XXV and Götting. Nachr. 1874, p.~539.
}. 
\end{theorem}
\renewcommand{\thefootnote}{\arabic{footnote}}

If we call \terminology{$m$-fold transitive}\, 
an $r$-term group of the space
$x_1, \dots, x_n$ when it contains at least one transformation which
transfers any $m$ given points in mutually general position to $m$
other arbitrary given points in general position, then we can now
easily prove {\it firstly}, that $m$ is always $\leqslant n + 2$ and
{\it secondly}, that every group for which $m = n + 2$ is equivalent
to the general projective group of the $n$-fold extended space.

Indeed, if the infinitesimal transformations: $X_1 f, \dots, X_r f$ in
the variables $x_1, \dots, x_n$ generate an $m$-fold transitive group,
then it stands immediately to reason that the linear homogeneous group
$L_1f, L_2 f, \dots$ assigned to a point $x_k^0$ in general position
does transform the $\infty^{ n-1}$ line-elements passing through this
point by means of an $(m-1)$-fold transitive projective group. But
now, since the general projective group of an $(n-1)$-fold extended
space is known to be $(n+1)$-transitive, we then realize that:
\[
m-1
\leqslant
n+1
\ \ \ \ \ \ \ \
\text{\rm and hence:}
\ \ \ \ \
m\leqslant n+2.
\]
Therefore, the following holds.

\def\thetheorem{113}\begin{theorem}
A finite continuous group in $n$ variables is
at most $(n + 2)$-fold transitive.
\end{theorem}

If an $r$-term group: $X_1f, \dots, X_r f$ in $n$ variables is exactly
$(n + 2)$-transitive, then as observed earlier on, the group $L_1f,
L_2f, \dots$ assigned to a point $x_k^0$ in general position is the
general or the special linear homogeneous group in $n$ variables;
consequently, the group: $X_1f, \dots, X_r f$ is equivalent either to
the general projective group, or to the general linear group, or to
the special linear group in $n$ variables. But now, because only the
first-mentioned among these three groups is $( n + 2)$-transitive, we
obtain the

\def\thetheorem{114}\begin{theorem}
If an $r$-term group in $n$ variables is $( n + 2)$-transitive,
then it is equivalent to the general projective group of the $n$-fold
extended space.
\end{theorem}

In the second and in the third volume, a series of studies which are
analogous to those conducted in this chapter shall be, among other
things, undertaken.

\linestop

\backmatter

\end{document}